\documentclass[11pt, prb, onecolumn, floatfix, letterpaper, nofootinbib ]{revtex4-2}

\usepackage[a-1b]{pdfx}

\usepackage{graphicx}
\usepackage{dcolumn}
\usepackage{bm}
\usepackage{multirow}

\usepackage[normalem]{ulem}

\usepackage{amsmath}
\usepackage{amsthm}
\usepackage{amstext}
\usepackage{amssymb}
\usepackage{mathrsfs}
\usepackage{amsfonts}
\usepackage{amsbsy} 

\usepackage[all,matrix,cmtip]{xy}

\usepackage{csquotes}
\MakeOuterQuote{"}

\usepackage{color}

\definecolor{red}{rgb}{1,0,0}
\definecolor{blue}{rgb}{0,0,1}
\definecolor{dblue}{rgb}{0,0,0.4}
\definecolor{green}{rgb}{0,1,0}
\definecolor{dgreen}{rgb}{0,0.4,0}
\definecolor{black}{rgb}{0,0,0}
\definecolor{white}{rgb}{1,1,1}

\definecolor{brn}{rgb}{.8,.4,.0}
\definecolor{redo}{rgb}{1,.5,.0}
\definecolor{ddgrn}{rgb}{0,0.4,0}
\definecolor{dgrn}{rgb}{0,0.55,0}
\definecolor{dbl}{rgb}{0,0,0.5}
\definecolor{grey}{rgb}{0.5,0.5,0.5}

\newcommand{\white}{\color{white}}

\newcommand{\black}{\color{black}}
\newcommand{\grey}{\color{grey}}

\usepackage[bbgreekl]{mathbbol}
\newcommand{\one}{\mathbf{1}}

\newcommand{\Z}{\mathbb{Z}}
\newcommand{\C}{\mathbb{C}}
\newcommand{\R}{\mathbb{R}}
\newcommand{\Q}{\mathbb{Q}}

\newcommand{\N}{\mathbb{N}}
\newcommand{\BO}{\mathbb{O}}

\renewcommand{\v}[1]{\boldsymbol{#1}} 
\renewcommand{\t}[1]{\widetilde{#1}} 
\newcommand{\ii}{\hspace{1pt}\mathrm{i}\hspace{1pt}}
\newcommand{\ee}{\hspace{1pt}\mathrm{e}}

\newcommand{\Rf}[1]{Ref.~\onlinecite{#1}}

\newcommand{\Tr}{{\rm Tr}}

\newcommand{\ie}{{\it i.e.~}} 
 
\newcommand{\etc}{{\it etc.}}

\newcommand{\bpm}{\begin{pmatrix}}
\newcommand{\epm}{\end{pmatrix}}
\newcommand{\bmm}{\begin{matrix}}
\newcommand{\emm}{\end{matrix}}

\newcommand{\cA}{ {\cal A} } 
\newcommand{\cB}{ {\cal B} }
\newcommand{\cC}{ {\cal C} } 
\newcommand{\cD}{ {\cal D} } 
 
\newcommand{\cF}{ {\cal F} }

\newcommand{\cM}{ {\cal M} }

\newcommand{\cR}{ {\cal R} }

\usepackage{euscript}

\newcommand\eC           {\EuScript{C}}
\newcommand\eD           {\EuScript{D}}
\newcommand\eE          {\EuScript{E}}
\newcommand\eF          {\EuScript{F}}

\newcommand\eM          {\EuScript{M}}

\newcommand\eZ         {\EuScript{Z}}

\newcommand{\al}{\alpha} 
\newcommand{\bt}{\beta} 
\newcommand{\del}{\delta} 
 
\newcommand{\eps}{\epsilon} 
\newcommand{\veps}{\varepsilon} 
\newcommand{\ga}{\gamma}

\newcommand{\la}{\lambda} 
 
\newcommand{\om}{\omega} 
 
\renewcommand{\th}{\theta} 
 
\newcommand{\si}{\sigma}

\usepackage{longtable}
\usepackage{tabularx}

\usepackage{hyperref}
\hypersetup{
    colorlinks = true,
    allcolors = blue,
}

\newcommand\cNG{\mathcal{NG}}
\newcommand\onebb{\mathbb{1}}

\newtheorem{thm}{Theorem}[section]
\newtheorem{cor}[thm]{Corollary}
\newtheorem{prop}[thm]{Proposition}
\newtheorem{lem}[thm]{Lemma}

\newtheorem{defn}[thm]{Definition}

\renewcommand\Vec{\mathcal{V}\mathrm{ec}}
\newcommand\Rep{\mathcal{R}\mathrm{ep}}

\newcommand\ord{\operatorname{ord}}
\newcommand\pord{\operatorname{pord}}

\newcommand{\FPdim}{\operatorname{FPdim}}

\newcommand{\frt}{\mathfrak{t}}
\newcommand{\frs}{\mathfrak{s}}

\newcommand{\SL}{\operatorname{SL}_2(\Z)}
\newcommand{\hSL}{{\widehat{\operatorname{SL}_2(\Z)}}}

\newcommand{\qsl}[1]{{\operatorname{SL}_2(\Z / #1 \Z)}}

\newcommand{\id}{\mathrm{id}}
\newcommand{\Gal}{\operatorname{Gal}}
\newcommand{\GL}{\operatorname{GL}}

\newcommand\GQ{\Gal(\bar\Q)}

\newcommand\hs{\hat\sigma}

\newcommand{\irep}[1]{}

\begin{document}

\title{Classification of modular data up to rank 12}

\author{Siu-Hung Ng}
\address{Department of Mathematics,
Louisiana State University,
Baton Rouge, LA 70803, USA}

\author{Eric C. Rowell}
\affiliation{ Department of Mathematics, Texas A\&M University, College
Station, TX 77843, USA }
\affiliation{ School of Mathematics, Leeds University, Leeds, LS2 9JT UK }

\author{Xiao-Gang Wen}
\affiliation{Department of Physics, Massachusetts Institute of Technology,
Cambridge, Massachusetts 02139, USA}

\begin{abstract} 

We use the computer algebra system GAP to classify  modular data up to rank 12.
This extends the previously obtained classification of  modular data up to rank
6.  Our classification includes all the modular data from modular tensor
categories up to rank 12, with a few possible exceptions at rank 12 and levels $5,7$ and $14$.  
Those exceptions are eliminated up to a certain bound by an extensive finite search in place of required infinite search.
Our list contains a few potential unitary modular
data which are not known to correspond to any unitary
modular tensor categories (such as those from Kac-Moody algebra, twisted
quantum doubles of finite group, as well as their Abelian anyon condensations).
It remains to be shown if those potential modular data can be realized by
modular tensor categories or not. We provide some evidence that all may be constructed from centers of near-group
categories or gauging group symmetries of known modular tensor categories,
with the exception of a total of five cases at rank 11 (with $D^2 =1964.590$) and 12 (with $D^2 =3926.660$).  The
classification of modular data corresponds to a classification of modular
tensor categories (up to modular isotopes which are not expected to be present
at low ranks).  The classification of modular tensor categories leads to a
classification of gapped quantum phases of matter in 2-dimensional space for
bosonic lattice systems with no symmetry, as well as a classification of
generalized symmetries in 1-dimensional space.

\end{abstract}

\maketitle

\tableofcontents
\setcounter{tocdepth}{2} 

\section{Introduction}

\subsection{Gapped liquid phases of quantum matter 
and braided fusion higher categories}

Quantum states of matter can be divided into four classes:
\begin{itemize}
\item
Gapped liquid: All excitations have a gap and there are no low energy
excitations.  So the gapped states appear to be trivial at low energies.  Band
insulator and quantum Hall states are examples of gapped liquid states.

\item
Gapped non-liquid: All excitations also have an energy gap.  But in contrast to
gapped liquid, by definition, a gapped non-liquid cannot ``dissolve'' product
states \cite{ZW1490,SM1403}.  Gapped fracton states are examples of gapped
non-liquid states \cite{C0502,H11011962}.

\item
Gapless liquid: There are finitely  many types of gapless low energy
excitations.  Dirac/Weyl semimetal, superfluid, critical point at continuous
phase transition are examples of gapless liquid states.

\item
Gapless non-liquid: There are infinity many types of gapless low energy
excitations.  Fermi metal, Bose metal, {\it etc.} are examples of gapless
non-liquid states.

\end{itemize}

People used to believe that Landau symmetry breaking theory provides a
systematic description of gapped phases of quantum matter.  In this case, group
theory that describes  symmetry breaking patterns provides a mathematical
foundation and classification of spontaneous symmetry breaking states.

The experimental discovery of fractional quantum Hall states suggested that
Landau symmetry breaking theory fails to describe all gapped phases. This led to
the theoretical discovery of a new order in gapped liquid states: topological
order \cite{W8987,W9039,KW9327}, which corresponds to pattern of long range
many-body entanglement \cite{CGW1038}.  But what mathematical theory
systematically describes various topological orders ({\it i.e.} patterns of
long range entanglement)?

\begin{table}[tb] 
\caption{List of topological orders (TO) (up to $E(8)$ invertible topological
order and up to modular isotopes) for bosonic systems with no symmetry in
2-dimensional space, which are classified  by the unitary modular data (UMD)
with increasing rank (the number of anyon types).  The Abelian TOs have only
Abelian anyons (\ie pointed simple objects).  The non-Abelian TOs have at least
one non-Abelian anyon.  The prime TOs cannot be viewed as stacking of two
non-trivial TOs with fewer anyon types.  Such a classification also leads to a
classification of symmetry-TOs in 2-dimensional space, which classifies the
generalized global symmetries in 1-dimensional space, up to holo-equivalence
\cite{KZ200514178}.  Those generalized symmetries include, but can go beyond,
finite-group symmetries (with potential anomalies).  Note that the list
includes all modular data, but also contain a few potential modular data at
rank 11, and 12 whose realizations are unknown (see Table
\ref{tab:exotic}).  } \label{toptable} \centering 
\begin{tabular}{|c | c|c|c|c|c|c|c|c|c|c|c|c|}
\hline
\# of anyon types  (rank)
& 1 & 2 &  3 &  4 &  5 &  6 &  7 &  8 &  9 & 10 & 11 & 12\\
\hline
\hline
\# of TOs (UMD)
& 1 & 4 & 12 & 18 & 10 & 50 & 28 & 64 & 81 & 76 & 44 & 221\\
\hline
\# of prime TOs (prime UMD)
& 1 & 4 & 12 & 8 & 10 & 10 & 28 & 20 & 20 & 40 & 44 & 33\\
\hline
\# of Abelian TOs (pointed UMD) 
& 1 & 2 &  2 &  9 &  2 &  4 &  2 & 20 &  4 &  4 & 2 & 18\\
\hline
\# of non-Abelian TOs (non-pointed UMD) 
& 0 & 2 & 10 &  9 &  8 & 46 & 26 & 44 & 77 & 72 & 42 & 203\\
\hline
\hline
\# of symTOs (UMTC in trivial Witt class) 
& 1 & 0 & 0 & 3 & 0 & 0 & 0 & 6 & 6 & 3 & 0 & 3\\
\hline
\# of finite-group symmetries (with anomaly $\om$) 
& $1_{\Z_1}$ & 0 & 0 & $2_{\Z_2^\om}$ & 0 & 0 & 0 & $6_{S_3^\om}$ & $3_{\Z_3^\om}$ & 0 & 0 & 0 \\
\hline
\end{tabular}
\end{table}

There are two approaches:
\begin{itemize}
\item
\textbf{Ground state based}: the robust degenerate ground states led to the
discovery and physical definition of topological order \cite{W8987,WN9077}.  It
is a striking property of topological order that the degeneracy of ground
states depends on the topology of the closed space where the system lives.  We
note that the degenerate ground states give rise to a vector bundle over the
moduli space of gapped quantum systems \cite{W9039,KW9327,KW200411904}.  The
vector bundle, plus many additional conditions, can form a foundation for a
general theory of topological order.  We refer to this ground-state-based
approach as moduli bundle theory.  We note that the holonomy of the vector
bundle give rise to a projective representation of the mapping class group of
the space on which the  ground states live.  Thus the  degenerate ground states
form a projective representation of the mapping class group.  When the space is
a 2-dimensional torus, the mapping class group is $\SL$ and the representation
(with a particular choice of basis) is called modular data.  We see that the
mapping class group representations and the modular data are the key
ingredients of the moduli bundle theory.

\item
\textbf{Excitation based}: We may also use topological excitations (which is
defined as excitations that cannot be created individually) to describe
topological orders.  Those excitations can fuse and braid. Thus one can use a
non-degenerate braided (higher) fusion category to describe topological
excitations and its associated topological order \cite{K10042307,KW1458}.  In
2-dimensional space, the topological excitations are point-like (called
anyons). They carry fractional Abelian \cite{LM7701,W8257,ASW8422,H8483} or
non-Abelian statistics \cite{W9102,MR9162} described by braid group
representations \cite{W8413,GMS8503}.  Those anyons (and the associated
topological orders in 2-dimensional space) are systematically describe by
modular tenor category theory.  It is interesting to note that modular tensor
categories were first used to systematically describe rational conformal field
theories \cite{MS8977}.  Then topological quantum field theories in
2-dimensional space (which contain anyons) were shown to be closely connected
to rational conformal field theory \cite{W8951}.  In particular, the structure
of modular tenor categories in 1+1D conformal field theory have a natural
interpretation in terms of anyons in 2+1D topological quantum field theory.
This led to the modular tensor category description of anyons (for a review,
see \Rf{K062}) and 2+1D bosonic topological orders (for a review, see
\Rf{W150605768}).

\end{itemize}

To summarize, gapped quantum liquid \cite{ZW1490,SM1403} phases of matter (\ie
topological orders) in $n$-dimensional space are described by moduli bundle
theory \cite{W9039,KW9327,KW200411904} or braided fusion $n-1$-categories with
trivial center \cite{KW1458,KZ170200673}.  For example, for bosonic systems
with no symmetry, there is no non-trivial gapped quantum liquid phases (\ie no
non-trivial topological order) in 1-dimensional space
\cite{FNW9243,CGW1107,SPC1139}, since the mapping class group of a circle,
SL$_1(\Z)$,  is trivial.  The gapped quantum liquid phases (\ie topological
orders) in 2-dimensional space (up to stacking of $E(8)$ invertible topological
orders), are classified by modular data for the torus and generalized modular
data for high genus surfaces \cite{W9039,KW9327,KW200411904,WW190810381}, or
alternatively, by unitary modular tensor categories (UMTC, which are braided
fusion $1$-categories with trivial center).  The gapped quantum liquid phases
in 3+1-dimensions are also classified.  For example, those without emergent
fermions are classified by a finite group $G$ and its group  cohomology classes
$\om \in H^4(G;\mathbb{R}/\Z)$ \cite{LW170404221}.

Because moduli bundle theory and non-degenerate braided fusion (higher)
category theory describe the same physical object -- topological order, in this
paper, we are going apply this connection in 2-dimensional space, and use the
moduli bundle approach to classify modular tensor categories through modular
data.  In other words, we will use modular data  to classify modular tensor
categories, up to modular isotopes.  Here modular isotopes correspond to
different modular tensor categories with the same modular data.  The first
example of modular isotopes is given in \Rf{MS170802796}, the twisted quantum
double $D^\om(\Z_{11}\rtimes \Z_5)$ of rank 49.  There are no known modular
isotopes at rank 12 or less.  Thus,  at low ranks, a classification of modular
data likely corresponds to a classification of modular tensor categories.  A
classification of modular tensor categories in turn gives rise to a
classification of all gapped quantum phases of matter in 2-dimensional space
(see Table \ref{toptable}).

\subsection{Generalized symmetry and braided fusion higher categories in
trivial Witt class}

\begin{table}[tb] 
\caption{Fusion rule for rank-11 $D^2 \approx 1964.590$ potential modular data }
\label{fusion} \centering \scriptsize
\begin{tabularx}{6.5in}{ |c||c|X|X|X|X|X|X|X|X|X|X|}
 \hline 
$\otimes$ & $\onebb$  & $a$  & $b$  & $c$  & $d$  & $e$  & $f$  & $g$  & $h$  & $i$  & $j$ \\ 
\hline 
 \hline 
$\onebb$  & $ \onebb$  & $ a$  & $ b$  & $ c$  & $ d$  & $ e$  & $ f$  & $ g$  & $ h$  & $ i$  & $ j$  \\ 
 \hline 
$a$  & $ a$  & $ \onebb \oplus a \oplus c \oplus d \oplus e \oplus f \oplus i$  & $ b \oplus g \oplus h \oplus i \oplus j$  & $ a \oplus f \oplus g \oplus i \oplus j$  & $ a \oplus d \oplus e \oplus f \oplus g \oplus i \oplus j$  & $ a \oplus d \oplus e \oplus f \oplus g \oplus h \oplus i \oplus j$  & $ a \oplus c \oplus d \oplus e \oplus f \oplus g \oplus h \oplus i \oplus j$  & $ b \oplus c \oplus d \oplus e \oplus f \oplus g \oplus h \oplus i \oplus j$  & $ b \oplus e \oplus f \oplus g \oplus h \oplus i \oplus 2j$  & $ a \oplus b \oplus c \oplus d \oplus e \oplus f \oplus g \oplus h \oplus 2i \oplus 2j$  & $ b \oplus c \oplus d \oplus e \oplus f \oplus g \oplus 2h \oplus 2i \oplus 2j$  \\ 
 \hline 
$b$  & $ b$  & $ b \oplus g \oplus h \oplus i \oplus j$  & $ \onebb \oplus a \oplus b \oplus d \oplus e \oplus g \oplus i$  & $ c \oplus f \oplus h \oplus i \oplus j$  & $ b \oplus d \oplus e \oplus g \oplus h \oplus i \oplus j$  & $ b \oplus d \oplus e \oplus f \oplus g \oplus h \oplus i \oplus j$  & $ c \oplus e \oplus f \oplus g \oplus h \oplus i \oplus 2j$  & $ a \oplus b \oplus d \oplus e \oplus f \oplus g \oplus h \oplus i \oplus j$  & $ a \oplus c \oplus d \oplus e \oplus f \oplus g \oplus h \oplus i \oplus j$  & $ a \oplus b \oplus c \oplus d \oplus e \oplus f \oplus g \oplus h \oplus 2i \oplus 2j$  & $ a \oplus c \oplus d \oplus e \oplus 2f \oplus g \oplus h \oplus 2i \oplus 2j$  \\ 
 \hline 
$c$  & $ c$  & $ a \oplus f \oplus g \oplus i \oplus j$  & $ c \oplus f \oplus h \oplus i \oplus j$  & $ \onebb \oplus b \oplus c \oplus d \oplus e \oplus h \oplus i$  & $ c \oplus d \oplus e \oplus f \oplus h \oplus i \oplus j$  & $ c \oplus d \oplus e \oplus f \oplus g \oplus h \oplus i \oplus j$  & $ a \oplus b \oplus d \oplus e \oplus f \oplus g \oplus h \oplus i \oplus j$  & $ a \oplus e \oplus f \oplus g \oplus h \oplus i \oplus 2j$  & $ b \oplus c \oplus d \oplus e \oplus f \oplus g \oplus h \oplus i \oplus j$  & $ a \oplus b \oplus c \oplus d \oplus e \oplus f \oplus g \oplus h \oplus 2i \oplus 2j$  & $ a \oplus b \oplus d \oplus e \oplus f \oplus 2g \oplus h \oplus 2i \oplus 2j$  \\ 
 \hline 
$d$  & $ d$  & $ a \oplus d \oplus e \oplus f \oplus g \oplus i \oplus j$  & $ b \oplus d \oplus e \oplus g \oplus h \oplus i \oplus j$  & $ c \oplus d \oplus e \oplus f \oplus h \oplus i \oplus j$  & $ \onebb \oplus a \oplus b \oplus c \oplus d \oplus e \oplus f \oplus g \oplus h \oplus i \oplus j$  & $ a \oplus b \oplus c \oplus d \oplus e \oplus f \oplus g \oplus h \oplus 2i \oplus j$  & $ a \oplus c \oplus d \oplus e \oplus f \oplus g \oplus h \oplus 2i \oplus 2j$  & $ a \oplus b \oplus d \oplus e \oplus f \oplus g \oplus h \oplus 2i \oplus 2j$  & $ b \oplus c \oplus d \oplus e \oplus f \oplus g \oplus h \oplus 2i \oplus 2j$  & $ a \oplus b \oplus c \oplus d \oplus 2e \oplus 2f \oplus 2g \oplus 2h \oplus 2i \oplus 2j$  & $ a \oplus b \oplus c \oplus d \oplus e \oplus 2f \oplus 2g \oplus 2h \oplus 2i \oplus 3j$  \\ 
 \hline 
$e$  & $ e$  & $ a \oplus d \oplus e \oplus f \oplus g \oplus h \oplus i \oplus j$  & $ b \oplus d \oplus e \oplus f \oplus g \oplus h \oplus i \oplus j$  & $ c \oplus d \oplus e \oplus f \oplus g \oplus h \oplus i \oplus j$  & $ a \oplus b \oplus c \oplus d \oplus e \oplus f \oplus g \oplus h \oplus 2i \oplus j$  & $ \onebb \oplus a \oplus b \oplus c \oplus d \oplus e \oplus f \oplus g \oplus h \oplus 2i \oplus 2j$  & $ a \oplus b \oplus c \oplus d \oplus e \oplus 2f \oplus g \oplus h \oplus 2i \oplus 2j$  & $ a \oplus b \oplus c \oplus d \oplus e \oplus f \oplus 2g \oplus h \oplus 2i \oplus 2j$  & $ a \oplus b \oplus c \oplus d \oplus e \oplus f \oplus g \oplus 2h \oplus 2i \oplus 2j$  & $ a \oplus b \oplus c \oplus 2d \oplus 2e \oplus 2f \oplus 2g \oplus 2h \oplus 2i \oplus 3j$  & $ a \oplus b \oplus c \oplus d \oplus 2e \oplus 2f \oplus 2g \oplus 2h \oplus 3i \oplus 3j$  \\ 
 \hline 
$f$  & $ f$  & $ a \oplus c \oplus d \oplus e \oplus f \oplus g \oplus h \oplus i \oplus j$  & $ c \oplus e \oplus f \oplus g \oplus h \oplus i \oplus 2j$  & $ a \oplus b \oplus d \oplus e \oplus f \oplus g \oplus h \oplus i \oplus j$  & $ a \oplus c \oplus d \oplus e \oplus f \oplus g \oplus h \oplus 2i \oplus 2j$  & $ a \oplus b \oplus c \oplus d \oplus e \oplus 2f \oplus g \oplus h \oplus 2i \oplus 2j$  & $ \onebb \oplus a \oplus b \oplus c \oplus d \oplus 2e \oplus f \oplus 2g \oplus h \oplus 2i \oplus 2j$  & $ a \oplus b \oplus c \oplus d \oplus e \oplus 2f \oplus g \oplus 2h \oplus 2i \oplus 2j$  & $ a \oplus b \oplus c \oplus d \oplus e \oplus f \oplus 2g \oplus 2h \oplus 2i \oplus 2j$  & $ a \oplus b \oplus c \oplus 2d \oplus 2e \oplus 2f \oplus 2g \oplus 2h \oplus 3i \oplus 3j$  & $ a \oplus 2b \oplus c \oplus 2d \oplus 2e \oplus 2f \oplus 2g \oplus 2h \oplus 3i \oplus 3j$  \\ 
 \hline 
$g$  & $ g$  & $ b \oplus c \oplus d \oplus e \oplus f \oplus g \oplus h \oplus i \oplus j$  & $ a \oplus b \oplus d \oplus e \oplus f \oplus g \oplus h \oplus i \oplus j$  & $ a \oplus e \oplus f \oplus g \oplus h \oplus i \oplus 2j$  & $ a \oplus b \oplus d \oplus e \oplus f \oplus g \oplus h \oplus 2i \oplus 2j$  & $ a \oplus b \oplus c \oplus d \oplus e \oplus f \oplus 2g \oplus h \oplus 2i \oplus 2j$  & $ a \oplus b \oplus c \oplus d \oplus e \oplus 2f \oplus g \oplus 2h \oplus 2i \oplus 2j$  & $ \onebb \oplus a \oplus b \oplus c \oplus d \oplus 2e \oplus f \oplus g \oplus 2h \oplus 2i \oplus 2j$  & $ a \oplus b \oplus c \oplus d \oplus e \oplus 2f \oplus 2g \oplus h \oplus 2i \oplus 2j$  & $ a \oplus b \oplus c \oplus 2d \oplus 2e \oplus 2f \oplus 2g \oplus 2h \oplus 3i \oplus 3j$  & $ a \oplus b \oplus 2c \oplus 2d \oplus 2e \oplus 2f \oplus 2g \oplus 2h \oplus 3i \oplus 3j$  \\ 
 \hline 
$h$  & $ h$  & $ b \oplus e \oplus f \oplus g \oplus h \oplus i \oplus 2j$  & $ a \oplus c \oplus d \oplus e \oplus f \oplus g \oplus h \oplus i \oplus j$  & $ b \oplus c \oplus d \oplus e \oplus f \oplus g \oplus h \oplus i \oplus j$  & $ b \oplus c \oplus d \oplus e \oplus f \oplus g \oplus h \oplus 2i \oplus 2j$  & $ a \oplus b \oplus c \oplus d \oplus e \oplus f \oplus g \oplus 2h \oplus 2i \oplus 2j$  & $ a \oplus b \oplus c \oplus d \oplus e \oplus f \oplus 2g \oplus 2h \oplus 2i \oplus 2j$  & $ a \oplus b \oplus c \oplus d \oplus e \oplus 2f \oplus 2g \oplus h \oplus 2i \oplus 2j$  & $ \onebb \oplus a \oplus b \oplus c \oplus d \oplus 2e \oplus 2f \oplus g \oplus h \oplus 2i \oplus 2j$  & $ a \oplus b \oplus c \oplus 2d \oplus 2e \oplus 2f \oplus 2g \oplus 2h \oplus 3i \oplus 3j$  & $ 2a \oplus b \oplus c \oplus 2d \oplus 2e \oplus 2f \oplus 2g \oplus 2h \oplus 3i \oplus 3j$  \\ 
 \hline 
$i$  & $ i$  & $ a \oplus b \oplus c \oplus d \oplus e \oplus f \oplus g \oplus h \oplus 2i \oplus 2j$  & $ a \oplus b \oplus c \oplus d \oplus e \oplus f \oplus g \oplus h \oplus 2i \oplus 2j$  & $ a \oplus b \oplus c \oplus d \oplus e \oplus f \oplus g \oplus h \oplus 2i \oplus 2j$  & $ a \oplus b \oplus c \oplus d \oplus 2e \oplus 2f \oplus 2g \oplus 2h \oplus 2i \oplus 2j$  & $ a \oplus b \oplus c \oplus 2d \oplus 2e \oplus 2f \oplus 2g \oplus 2h \oplus 2i \oplus 3j$  & $ a \oplus b \oplus c \oplus 2d \oplus 2e \oplus 2f \oplus 2g \oplus 2h \oplus 3i \oplus 3j$  & $ a \oplus b \oplus c \oplus 2d \oplus 2e \oplus 2f \oplus 2g \oplus 2h \oplus 3i \oplus 3j$  & $ a \oplus b \oplus c \oplus 2d \oplus 2e \oplus 2f \oplus 2g \oplus 2h \oplus 3i \oplus 3j$  & $ \onebb \oplus 2a \oplus 2b \oplus 2c \oplus 2d \oplus 2e \oplus 3f \oplus 3g \oplus 3h \oplus 4i \oplus 4j$  & $ 2a \oplus 2b \oplus 2c \oplus 2d \oplus 3e \oplus 3f \oplus 3g \oplus 3h \oplus 4i \oplus 4j$  \\ 
 \hline 
$j$  & $ j$  & $ b \oplus c \oplus d \oplus e \oplus f \oplus g \oplus 2h \oplus 2i \oplus 2j$  & $ a \oplus c \oplus d \oplus e \oplus 2f \oplus g \oplus h \oplus 2i \oplus 2j$  & $ a \oplus b \oplus d \oplus e \oplus f \oplus 2g \oplus h \oplus 2i \oplus 2j$  & $ a \oplus b \oplus c \oplus d \oplus e \oplus 2f \oplus 2g \oplus 2h \oplus 2i \oplus 3j$  & $ a \oplus b \oplus c \oplus d \oplus 2e \oplus 2f \oplus 2g \oplus 2h \oplus 3i \oplus 3j$  & $ a \oplus 2b \oplus c \oplus 2d \oplus 2e \oplus 2f \oplus 2g \oplus 2h \oplus 3i \oplus 3j$  & $ a \oplus b \oplus 2c \oplus 2d \oplus 2e \oplus 2f \oplus 2g \oplus 2h \oplus 3i \oplus 3j$  & $ 2a \oplus b \oplus c \oplus 2d \oplus 2e \oplus 2f \oplus 2g \oplus 2h \oplus 3i \oplus 3j$  & $ 2a \oplus 2b \oplus 2c \oplus 2d \oplus 3e \oplus 3f \oplus 3g \oplus 3h \oplus 4i \oplus 4j$  & $ \onebb \oplus 2a \oplus 2b \oplus 2c \oplus 3d \oplus 3e \oplus 3f \oplus 3g \oplus 3h \oplus 4i \oplus 4j$  \\ 
 \hline 
\end{tabularx}
\end{table}

After a systematic understanding of gapped phases of quantum matter, the next
goal is to have a systematic understanding of gapless liquid phases of quantum
matter.  This is a wide open and much harder problem.  One idea is to use
emergent symmetries at low energies to character those gapless phases, hoping
to obtain a systematic understanding via the maximal emergent symmetries
\cite{CW221214432}.

It became more and more clear that emergent symmetries in gapless quantum
states are generalized symmetries, which can be a combination of ordinary symmetry (described by group),
higher-form symmetry \cite{NOc0605316,NOc0702377,KT13094721,GW14125148}
(described by higher-group with only one non-trivial layer),
higher-group symmetry \cite{KT13094721}, anomalous ordinary symmetry
\cite{H8035,CGL1314,W1313,KT14030617}, anomalous higher symmetry
\cite{KT13094721,GW14125148,TK151102929,T171209542,P180201139,DT180210104,BH180309336,ZW180809394,WW181211968,WW181211967,GW181211959,WW181211955,W181202517},
beyond-anomalous symmetry \cite{CW220303596}, non-invertible 0-symmetry (in
1+1D)
\cite{PZh0011021,CSh0107001,FSh0204148,FSh0607247,FS09095013,BT170402330,CY180204445,TW191202817,KZ191213168,I210315588,Q200509072},
non-invertible higher symmetry (which includes algebraic higher symmetry)
\cite{DR11070495,KZ200308898,KZ200514178,FT220907471}, and/or non-invertible gravitational
anomaly \cite{KW1458,FV14095723,M14107442,KZ150201690,KZ170200673,JW190513279}.
Thus emergent symmetries can go beyond the group and higher group theory
description.  It was proposed
\cite{KZ150201690,JW191213492,LB200304328,KZ200308898,KZ200514178,GK200805960,KZ210703858,KZ201102859,KZ220105726,CW220303596,FT220907471}
that the (holo-equivalent classes of)  generalized symmetries in
$n$-dimensional space can be systematically described by topological order (TO)
with gappable boundary in one higher dimension, or equivalently by
non-degenerate braided fusion $n-1$-categories in trivial Witt class.  Such
topological order with gappable boundary (\ie non-degenerate braided fusion
higher category in trivial Witt class) is referred to as symmetry-TO (symTO).
Thus symTO, replacing group and higher group, describes generalized symmetry,
which can be anomalous or beyond anomalous.

Classification of groups is a major achievement of modern mathematics.  Such a
classification is important since groups describe possible symmetries in our
world. However, from the above discussion, we see that holo-equivalent
symmetries in our quantum world are actually described by non-degenerate
braided fusion higher categories in trivial Witt class.  Thus our
classification of modular tensor categories also leads to a classification of
emergent generalized symmetry for quantum systems in 1-dimensional space (see
Table \ref{toptable}).

\subsection{Summary of results}
\label{sec:sum}

In this paper, we use the GAP computer algebraic system \cite{GAP4} to classify
modular data up to rank 12.  We focus on \emph{prime} modular data, {\it i.e.}
those that are not simply products of two smaller modular data.  Those prime
modular data are given in Section \ref{MDbyGal}, where, for each Galois orbit,
we list one modular data (an unitary one if exists).  We also indicate the MTCs
that realize the modular data, if realizations are known.  We list all unitary
modular data in Appendix \ref{UMD}, and all unitary and non-unitary modular
data in Appendix \ref{allMD}.

From our explicit classifications, we find some exotic potential modular data
which are not realized by Kac-Moody algebras or by twisted quantum doubles, nor
by their Deligne product, their Galois conjugations, their changes of spherical
structure, and their Abelian anyon condensations \cite{LW170107820}.  We show
or provide strong evidence that all but five of the  exotic potential modular
data are indeed modular data that can be realized by center of near group
fusion category or gauging the automorphism of known modular data (followed by
some condensation reductions, see Section \ref{realizations} for details).  

One of the five potential modular data with unknown realization is labeled by
$11_{\frac{32}{5},1964.}^{35,581}$ which has rank 11, central charge $ c =
\frac{32}{5} $ mod 8, and total quantum dimension $D^2 \approx 1964.590$. We do
not know for sure if $11_{\frac{32}{5},1964.}^{35,581}$ is a modular data or
not.  The topological spins and quantum dimensions for such a potential modular
data are given by\\
$s_i =  0,
\frac{2}{35},
\frac{22}{35},
\frac{32}{35},
\frac{1}{5},
0,
\frac{3}{7},
\frac{5}{7},
\frac{6}{7},
\frac{3}{5},
\frac{1}{5} 
$,
\\[2mm]
 $d_i =  1$, {\footnotesize
 $\frac12\Big(5+2\sqrt5+\sqrt{7(5+2\sqrt5)}\Big)$,
 $\frac12\Big(5+2\sqrt5+\sqrt{7(5+2\sqrt5)}\Big)$,
 $\frac12\Big(5+2\sqrt5+\sqrt{7(5+2\sqrt5)}\Big)$,\\ {\white.}\ \ \ \ \ \
 $\frac14\Big(9+5\sqrt5+\sqrt{14(25+11\sqrt5)}\Big)$,
 $\frac14\Big(11+7\sqrt5+\sqrt{14(25+11\sqrt5)}\Big)$,
 $\frac14\Big(15+7\sqrt5+\sqrt{14(25+11\sqrt5)}\Big)$,\\ {\white.}\ \ \ \ \ \ 
 $\frac14\Big(15+7\sqrt5+\sqrt{14(25+11\sqrt5)}\Big)$,
 $\frac14\Big(15+7\sqrt5+\sqrt{14(25+11\sqrt5)}\Big)$,
 $\frac14\Big(21+7\sqrt5+\sqrt{14(65+29\sqrt5)}\Big)$,\\ {\white.}\ \ \ \ \ \
 $\frac14\Big(19+9\sqrt5+\sqrt{14(65+29\sqrt5)}\Big)$}.
\\
The $S$ matrices is given in Section \ref{ss:rank11}, and fusion by Table
\ref{fusion}.  The higher central charges \cite{NW181211234} are
$
1$, $ 
\zeta_5^4$, $ 
-\zeta_5^4$, $ 
-\zeta_5$, $ 
-\zeta_5$, $ 
\frac00$, $ 
-\zeta_5^4$, $ 
\zeta_5$, $ 
-\zeta_5$, $ 
-\zeta_5$, $ 
\frac00$, $ 
\zeta_5^4$, $ 
-\zeta_5^4$, $ 
-\zeta_5$, $ 
-\zeta_5^2$, $ 
\frac00$, $ 
\zeta_5^4$, $ 
-\zeta_5^4$, $ 
-\zeta_5$, $ 
\zeta_5$, $ 
\frac00$, $ 
-\zeta_5^3$, $ 
-\zeta_5^4$, $ 
-\zeta_5$, $ 
\zeta_5$, $ 
\frac00$, $ 
-\zeta_5^4$, $ 
-\zeta_5^4$, $ 
\zeta_5^4$, $ 
-\zeta_5$, $ 
\frac00$, $ 
-\zeta_5^4$, $ 
-\zeta_5^4$, $ 
-\zeta_5$, $ 
\zeta_5
$.
Such a potential modular data has no non-trivial condensable algebra.  There are
three other potential modular data at rank 12 whose realizations are unknown.
Those potential modular data have central charge $ c = 0 $ mod 8 and total
quantum dimension $D^2 \approx 3926.660$ (see Section \ref{ss:rank12}).

\section{The necessary conditions on  modular data}

Modular data is an important invariant of MTC.  In the following, we list
necessary conditions on  modular data.  We will use those conditions, trying to
classify modular data. 
Many conditions are well-known and can be found in, e.g.
\cite{BK}.  
\begin{prop}
\label{p:MD}
The modular data $(S,T)$ of an MTC satisfies:
\begin{enumerate}
\item
$S,T$ are symmetric complex matrices, indexed by $i,j=1,\ldots,
r$.\footnote{The index also labels the simple objects in the MTC, with $i=1$
corresponding to the unit object, and $r$ is the \textbf{rank} of the modular
data and the MTC.}

\item
$T$ is unitary, diagonal, and $T_{11}=1$.

\item
$S_{11}=1$. Let  $d_i = S_{1i} $ and
$D=\sqrt{\sum_{i=1}^{r} d_i d_i^*}$ (the positive root). Then 
\begin{align}
S S^\dag = D^2\id,
\end{align}
and the $d_i\in\mathbb{R}$.
\item
$S_{ij}$ are cyclotomic integers in $\Q_{\ord(T)}$\footnote{Here $\Q_n$ denotes the field $\Q(\zeta_n)$ for a primitive $n$th root of unity $\zeta_n$} \cite{NS10}. 
The ratios
$S_{ij}/S_{1j}$ are cyclotomic integers for all $i,j$ \cite{CG}. Also there is a $j$ such
that $S_{ij}/S_{1j} \in [1,+\infty) $ for all $i$ \cite{ENO}.

\item

Let $\theta_i = T_{ii}$ and $p_\pm = \sum_{i=1}^{r} d_i^2 (\theta_i)^{\pm 1}$.
Then $p_+/p_-$ is a root of unity, and $p_+=D\ee^{\ii 2\pi c/8}$ for some
rational number $c$.\footnote{The \textbf{central charge} $c$ of the modular
data and of the MTC is only defined modulo $8$.} Moreover, the modular data
$(S,T)$ is associated with a projective $\SL$ representation, since:
\begin{align}
\label{SL2rep}
 (ST)^3 = p_+ S^2,\ \ \ \frac{S^2}{D^2} =C,\ \ \ C^2 =\id,
\end{align}
where $C$ is a permutation matrix
satisfying
\begin{align}
 \Tr(C) > 0.
\end{align}

\item
$D$ is a cyclotomic integer.  $D^5/\ord(T)$ is an algebraic  integer, which is
also a cyclotomic integer \cite{E02}.  $D/d_i$ are cyclotomic integers (see
Lemma \ref{Ddi}).  

\item Cauchy Theorem \cite{BNRW}: The set of prime divisors of $\ord(T)$
coincides with the prime divisors of $\mathrm{norm}(D^2)$.\footnote{Here
$\mathrm{norm}(x)$ is the product of the distinct Galois conjugates of the
algebraic number $x$.} The prime divisors of norm$(D)$ and $\ord(T)$ coincide.  
The prime divisors of norm$(D/d_i)$ are part of those of $\ord(T)$.

\item  Verlinde formula $($cf. \cite{V8860}$):$
\begin{equation}\label{Ver0} 
 N^{ij}_k = \frac{1}{D^2}\sum_{l=1}^{r} \frac{S_{li} S_{lj} S_{lk}^*}{ d_l }
\in \N ,  
\end{equation} 
where $i,j, k=1,2,\ldots,r$ and $\N$ is the set of non-negative integers.\footnote{The $N^{ij}_k$ are called the fusion coefficients.} The $N^{ij}_1$ satisfy
\begin{align}
 N^{ij}_1 = C_{ij},
\end{align}
which defines a
charge conjugation $i \to \bar i$ via
\begin{align}
 N^{\bar i j}_1 = \delta_{i j}.
\end{align}

\item 
\label{FScnd}
Let $n \in \N_+$.  The $n^\text{th}$ Frobenius-Schur indicator of the $i$-th
simple object
\begin{equation}
 \label{nunFS}
 \nu_n(i)= D^{-2} \sum_{j, k} N_i^{jk} (d_j\theta_j^n) (d_k\theta_k^n)^*
\end{equation}
 is a cyclotomic integer whose conductor divides $n$ and
$\ord(T)$ \cite{NS07a, NS07b}.  The 1st Frobenius-Schur indicator satisfies $\nu_1(i)=\delta_{i,1}$ while the 2nd
Frobenius-Schur indicator $\nu_2(i)$ satisfies $\nu_2(i)=0$ if $i\neq \bar i$,
and $\nu_2(i)=\pm 1$ if $i = \bar i$ (see \cite{Bantay, NS07a, RSW0777}).

\item The twists, fusion coefficients and $S$-matrix entries satisfy the balancing equation\footnote{This holds more generally in ribbon fusion categories, i.e. premodular categories.}:
\begin{equation}
\label{eq:balancing}
S_{ij}=\sum_k N^{ij}_k \frac{\theta_i\theta_j}{\theta_k} d_k .
\end{equation}


%
%
%
%
%
%
%

\end{enumerate}
\end{prop}

Based on a physics argument, \Rf{KW200411904} conjectured that modular data
also satisfy 
\begin{align}
c\, D_g/2 &\in \Z \   \text{ for } g\geq 3,
&
 D_g &= \sum_i \big(D/d_i\big)^{2(g-1)}, \ \ \ \ g = 0,1,2,3,4,\cdots.
\end{align}
We note that $\frac{D}{d_i}$ are real cyclotomic integers.  The
Galois conjugation of $\frac{D}{d_i}$ just permutes the $i$-indices:
$\si(\frac{D}{d_i}) = \frac{D}{d_{\hat \si(i)}}$.  Thus $D_g$ is a cyclotomic
integer that is invariant under all Galois conjugations, which implies that
$D_g$ is a positive integer.  In fact, $D_g$ is the ground state
degeneracy of the corresponding topological order on genus $g$ Riemann
surface.  We find that both the unitary and non-unitary modular data that we
obtained satisfy the above condition, although the condition is
argued for unitary topological orders.  For $g=2$, the condition $c\, D_2/2 \in
\Z$ is not satisfied. But a weaker condition 
\begin{align}
c\, D_2 \in \Z
\end{align}
is satisfied by all the unitary and non-unitary modular data that we
obtained, as well as by unitary modular data constructed from Kac-Moody algebras up to rank 200.  Also, \Rf{Em0207007} showed that 
\begin{align}
c\, D_1\, D^5/2 =  c\, r\, D^5/2 = \text{cyclotomic integers.}
\end{align}
This generalizes the above result to $g=1$ case.  The unitary and non-unitary modular data that we
obtained actually satisfy a stronger condition
\begin{align}
c\, D^5/2 = \text{cyclotomic integers.}  
\end{align}

From \eqref{SL2rep}, we see that the modular data $(S,T)$ is closely related to
the $\SL$ representations. We are going to use this relation to classify
modular data.  Let us first summarize some important facts about $\SL$
representations.  Let $\frs=\begin{bmatrix} 0&-1\\ 1&
0 \end{bmatrix} $, $\frt =\begin{bmatrix} 1&1\\ 0& 1 \end{bmatrix} $ be the standard generators of $\SL$. This admits the presentation:
\begin{align}
\SL = \langle \frs, \frt\mid \frs^4 = \id, (\frs \frt)^3 = \frs^2\rangle\,.
\end{align}

We note that for any positive integer $n$, the reduction $\Z \to \Z_n$ defines
a surjective group homomorphism  $\pi_n:\SL \to \qsl{n}$. Thus, a
representation of $\qsl{n}$ is also a representation of $\SL$, which will be
called a \emph{congruence} representation of $\SL$ in this paper.  It is
immediate to see that a representation of $\qsl{n}$ is also a $\qsl{mn}$
representation for any positive integer $m$. The smallest positive integer $n$
such that a congruence representation $\rho$ of $\SL$  factors through $\pi_n :
\SL \to \qsl{n}$ is called the \emph{level} of $\rho$. It is known that the
level $n=\ord(\rho(\frt))$ (cf.  \cite[Lem. A.1]{DLN}).  Here $\ord(t)$ is
defined as:
\begin{defn}{\rm
Let $t$ be any matrix over $\C$. The smallest positive integer $n$ such that
$t^n=\id$ is called the \emph{order} of $t$,
and denoted by $\ord(t):=n$. If such integer does not exist, we define
$\ord(t):=\infty$. 
}\end{defn}
Similarly
\begin{defn}{\rm
Let $t$ be any matrix over $\C$. The smallest positive integer $n$ such that
$t^n=\al \id$ for some $\al \in \C^\times$ is called the \emph{projective order} of $t$,
and denoted by $\pord(t):=n$. If such integer does not exist, we define
$\pord(t):=\infty$. 
}\end{defn}

We can organize the finite level irreducible representations of $\SL$ by the
level and the dimension of the representations.  Due to the Chinese remainder
theorem, if the level of a irreducible representation $\rho$ factors as
$n=\prod_i p_i^{k_i}$ where $p_i$ are distinct primes, then $\rho \cong
\bigotimes_i \rho_i$ where $\rho_i$ are level $p_i^{k_i}$ representations.
Thus we can construct all irreducible $\SL$ representations as tensor products
of irreducible $\SL$ representations of prime-power levels, which in turn,
yields a construction of all semisimple $\SL$ representations $\rho$ via direct
sums of the irreducible representations.

Define $\Q_n=\Q(\zeta_{n})$ to be the cyclotomic field of order $n$. For any
positive integer $n$, we can construct a faithful representation $D_n:
\Gal(\Q_n) \to \qsl{n}$, which identifies the Galois group $\Gal(\Q_n) \cong
\Z_n^\times$ with the diagonal subgroup of $\qsl{n}$ \cite[Remark 4.5]{DLN}.
More generally,  for any $\si \in \GQ$, $\si(\Q_n) = \Q_n$ and so there exists
an integer $a$ (unique modulo $n$) such that $\si( \zeta_n )=\zeta_n^a$ and  
\begin{align}
\label{DSLn}
 D_n(\si) :=
\frt^a \frs \frt^b \frs \frt^a \frs^{-1}  
=
\begin{pmatrix}
 a & 0\\
 0 & b\\
\end{pmatrix} \in \qsl{n}\,,
\end{align}
 where $b$ satisfies $ab \equiv 1 \mod n$. If $\rho$ is a  level $n$
representation of $\SL$, the composition 
\begin{equation} 
\label{defDrho1}
    D_\rho(\si):=\rho \circ D_n(\si)
\end{equation}
defines a representation of $\GQ$. We may also write $D_{n}(\si)$ as $D_n(a)$.
Such a representation of Galois group captures the Galois conjugation action on
$\SL$ representations of  modular data, and plays a very important role in our
classification. In other words, in our classification, we look for some $\SL$
representations such that $D_{\rho}(a)$ is a signed permutation matrix.

We also note that the 1-dimensional representations of $\SL$, denoted as
$\hSL$, form a cyclic group of order 12 under tensor product.  We will take
$\chi \in \hSL$ defined by $\chi(\frt) = \zeta_{12}$ to be the generator, where
$\zeta_n:=\ee^{2 \pi \ii /n}$ and $\zeta_n^k:=\ee^{2 \pi \ii k/n}$. Under this
convention, every 1-dimensional representation of $\SL$ is equivalent to
$\chi^\al$ for some integer $\al$, unique modulo 12:
\begin{align}
\chi^\al(\frs) = \overline \zeta_4^\al, \quad \chi^\al(\frt)=\zeta_{12}^\al .
\end{align}

From a modular data, we can obtain a particular type of $\SL$ representations,
called MD representations.  From a MD representation, we can also obtain its
corresponding modular data.  Thus we can classify modular data by classifying
MD representations.  In the following, we describe the detailed relation
between modular data and MD representations, and the necessary conditions for a
$\SL$ representation to be a MD representation.  Many of the following
collection of results on MD representations were proved in \cite{NS10, DLN}:
\begin{prop}
\label{p:MDcond}
Given a modular data $S,T$ of rank $r$, let $\rho_\al$ be any one of its 
\textbf{MD representations}, which is defined as
\begin{align}
\label{Ss2tSs}
 \rho_\al(\frs) = P \overline\zeta_4^\al \ee^{2\pi \ii \frac{c}{8}} \frac{S}{p_+} P^\top, \ \
 \rho_\al(\frt) = P \zeta_{12}^\al \ee^{-2\pi \ii \frac{c}{24}} T P^\top
\ \ \ (\al \in \Z_{12}), 
\end{align}
where $\al = 0,1,\cdots,11$ and $P$ is a permutation matrix.
Then, there exists a rational number $c$ (called central charge), such that
$\rho_\al$ for all $\al$ has the following properties:
\begin{enumerate}
\item $\rho_\al$ is an unitary and symmetric matrix representation of $\SL$ with
level $\ord(\rho_\al(\frt))$, and $\ord(T) \mid \ord(\rho_\al(\frt)) \mid 12
\ord(T)$\,.

\item The conductor of the elements of $\rho_\al(\frs)$ 
divides $\ord(\rho_\al(\frt))$.

\item 
If $\rho_\al$ is equivalent to a direct sum of two $\SL$ representations
\begin{align}
 \rho_\al \cong \rho \oplus \rho',
\end{align}
then the eigenvalues of $\rho(\frt)$ and $\rho'(\frt)$ must overlap.  This
implies that if $\rho_\al \cong \rho \oplus \chi_1 \oplus \dots \oplus
\chi_\ell$ for some 1-dimensional representations $\chi_1, \dots, \chi_\ell$,
then $\chi_1, \cdots \chi_\ell$ are the same 1-dimensional representation.

\item Suppose that $\rho_\al \cong \rho \oplus \ell \chi$ for an irreducible
representation $\rho$ with non-degenerate $\rho(\frt)$, and an 1-dimensional
representation $\chi$.  If $\ell \ne 2\dim(\rho) -1$ or $\ell > 1$, then
$(\rho(\frs) \chi(\frs)^{-1})^2 = \id$. 

\item $\rho_\al$ satisfies
\begin{align}
\rho_\al \not\cong n \rho
\end{align}
for any integer $n >1$ and any  representation $\rho$ such that
$\rho(\frt)$ is non-degenerate.

\item
If $\rho_\al(s)^2=\pm \id$ ({\it i.e.} if the modular data or MTC is self dual),
$\pord(\rho_\al(\frt)) $ is a prime and satisfies $\pord(\rho_\al(\frt)) = 1$ mod
4, then the representation $\rho_\al$ cannot be a direct sum of a
$d$-dimensional irreducible $\SL$ representation and two or more 1-dimensional
$\SL$ representations with $d=(p+1)/2$.  

\item Let $3< p < q$ be prime such that $pq \equiv 3 \mod 4$ and
$\pord(\rho_\al(\frt))=pq$, then the rank $r \ne \frac{p+q}{2}+1$.  Moreover, if
$p > 5$, rank $r > \frac{p+q}{2}+1$.

\item  The number of self dual objects is greater than 0. Thus
\begin{align}
 \Tr(\rho_\al(\frs)^2) \neq 0 .
\end{align}
Since $\Tr(\rho_\al(\frs)^2) \neq 0$, let us introduce
\begin{align}
 C = \frac{\Tr(\rho_\al(\frs)^2)}{|\Tr(\rho_\al(\frs)^2)|}   \rho_\al(\frs)^2.
\end{align}
The above $C$ is the charge conjugation operator of MTC,
{\it i.e.} $C$ is a permutation matrix of order 2.  In particular,
$\Tr(C)$ is the number of  self dual objects. Also, for each eigenvalue $\tilde
\theta$ of $\rho_\al(\frt)$, 
\begin{align}
\Tr_{\tilde \theta}(C) \geq 0,
\end{align}
where $\Tr_{\tilde \theta}$ is the trace in the degenerate
subspace of $\rho_\al(\frt)$ with eigenvalue $\tilde \theta$.

\item 
If the modular data is integral  and $\ord(\rho_\al(\frt))$ = odd, then
\begin{align}
\Tr(C) =\Tr^2(\rho^2_\al(\frs)) = 1,
\end{align}
\ie the unit object is the only self-dual object.

\item 
For any Galois conjugation $\si$ in $\Gal(\Q_{\ord(\rho_\al(\frt))})$, there is a
permutation of the indices, $i \to \hs(i)$, and $\eps_\si(i)\in \{1,-1\}$, such
that
\begin{align}
\label{Galact}
\si \big(\rho_\al(\frs)_{i,j}\big) &
= \eps_\si(i)\rho_\al(\frs)_{\hat \si (i),j} 
= \rho_\al(\frs)_{i,\hat \si (j)}\eps_\si(j) 
\\
\si^2 \big(\rho_\al(\frt)_{i,i}\big) &= \rho_\al(\frt)_{\hat \si (i),\hat \si (i)},
\end{align}
for all $i,j$. 

\item 
For any integer $a$ coprime to $n  = \ord(\rho_\al(\frt))$, we define
\begin{align}
\label{defDrho}
D_{\rho_\al}(a) &:= \rho_\al(\frt^a \frs \frt^b \frs \frt^a \frs^{-1}) 
= D_{\rho_\al}(a+\ord(\rho_\al(\frt))) , 
\nonumber\\
&
\text{ where } ab \equiv 1 \mod \ord(\rho_\al(\frt))\,. 
\end{align}
For any $\si \in \mathrm{Gal}(\Q_n)$, $\si(\zeta_n) = \zeta_n^a$ for some unique integer $a$ modulo $n$. We define
\begin{align}
\label{eq:Drho}
D_{\rho_\al}(\si) := D_{\rho_\al}(a)\,. 
\end{align}
By \cite[Theorem II]{DLN}, $D_{\rho_\al}: \mathrm{Gal}(\Q_n)=(\Z_n)^\times
\to \GL_r(\C)$ is a representation equivalent to the restriction of $\rho_\al$
on the diagonal subgroup of $\qsl{n}$.  $D_{\rho_\al}(\si)$ in \eqref{eq:Drho}
must be a \textbf{signed permutation}
\begin{align}
\label{Drhosihat}
 (D_{\rho_\al}(\si))_{i,j} = \eps_\si(i) \delta_{\hs(i),j}.
\end{align}
and satisfies
\begin{align}
\label{siDrho}
\si(\rho_\al(\frs)) &= D_{\rho_\al}(\si) \rho_\al(\frs) =\rho_\al(\frs)D_{\rho_\al}^\top(\si),
\nonumber\\
\si^2(\rho_\al(\frt)) &= D_{\rho_\al}(\si) \rho_\al(\frt) D_{\rho_\al}^\top(\si)
\end{align}

\item 
There exists a $u$ such that $\rho_\al(\frs)_{uu} \neq 0$ and
\begin{align}
\label{SSN}
&  \rho_\al(\frs)_{ui} \neq 0 , \ \ \ \
 \frac{ \rho_\al(\frs)_{ij} }{ \rho_\al(\frs)_{uu} },\ 
 \frac{ \rho_\al(\frs)_{ij} }{ \rho_\al(\frs)_{uj} }
\in \BO_{\ord(T)}, \ \ \ \
 \frac{ \rho_\al(\frs)_{ij} }{ \rho_\al(\frs)_{i'j'} }
\in \Q_{\ord(T)},
\nonumber\\
& N^{ij}_k = \sum_{l=0}^{r-1} \frac{
\rho_\al(\frs)_{li} \rho_\al(\frs)_{lj} \rho_\al(\frs^{-1})_{lk}}{ \rho_\al(\frs)_{lu} } \in\N. 
\nonumber\\
& \forall \ i,j,k = 1,2,\ldots,r.
\end{align}
(The index $u$ corresponds the unit object of MTC).  Here,  $\Q_{\ord(T)}$ is the field
of cyclotomic number and $\BO_{\ord(T)}$ is the ring of cyclotomic integer.
Also, $ \rho_\al(\frs)_{ui}$ for $i \in \{1,\cdots, r\}$ are either all real or
all imaginary.

\item
Let $n \in \N_+$.  The $n^\text{th}$ Frobenius-Schur indicator of the $i$-th
simple object
\begin{align}
 \label{nunFSapp}
 \nu_n(i)&= 
\sum_{j, k=0}^{r-1} N_i^{jk} \rho_\al(\frs)_{ju}\theta_j^n [\rho_\al(\frs)_{ku}\theta_k^n]^*
=\sum_{j, k=0}^{r-1} N_i^{jk} 
\rho_\al(\frt^n\frs)_{ju} \rho_\al(\frt^{-n}\frs^{-1})_{ku}
\nonumber\\
&=\sum_{j,k,l=0}^{r-1} \frac{
\rho_\al(\frs)_{lj} \rho_\al(\frs)_{lk} \rho^*_\al(\frs)_{li}}{ \rho_\al(\frs)_{lu} }
\rho_\al(\frt^n\frs)_{ju} \rho_\al(\frt^{-n}\frs^{-1})_{ku}
\nonumber\\
&=\sum_{l=0}^{r-1} \frac{
\rho_\al(\frs\frt^n\frs)_{lu} \rho_\al(\frs\frt^{-n}\frs^{-1})_{lu} \rho_\al(\frs^{-1})_{li}}{ \rho_\al(\frs)_{lu} }
\end{align}
is a cyclotomic integer whose conductor divides $n$ and $\ord(T)$.  The 1st
Frobenius-Schur indicator satisfies $\nu_1(i)=\delta_{iu}$ while the 2nd
Frobenius-Schur indicator $\nu_2(i)$ satisfies 
$\nu_2(i)=\pm \rho_\al(\frs^2)_{ii}$
(see \cite{Bantay, NS07b, RSW0777}).

The above condition can also be rewritten as
\begin{align}
\rho_\al(\frs)_{lu} 
\sum_i
\rho_\al(\frs)_{li} \nu_n(i)
=
 \rho_\al(\frs\frt^n\frs)_{lu} \rho_\al(\frs\frt^{-n}\frs^{-1})_{lu} 
\end{align}
Summing over $l$, we obtain
\begin{align}
\sum_i \t C_{ui} \nu_n(i) = \t C_{uu} \ \ \to \ \
 \t C_{uu} \nu_n(u) = \t C_{uu} \ \ \to \ \ \nu_n(u) = 1.
\end{align}

After a signed-diagonal conjugation $V_{ij} = v_i \del_{ij}$ that changes
$\rho_\al(\frs)_{ij}$
to $\rho_\mathrm{pMD} (\frs)_{ij}$,
we find
\begin{align}
\rho_\mathrm{pMD}(\frs)_{lu} 
&
\sum_i
\rho_\mathrm{pMD}(\frs)_{li} \nu_n(i) v_iv_u
=
 \rho_\mathrm{pMD}(\frs\frt^n\frs)_{lu} \rho_\mathrm{pMD}(\frs\frt^{-n}\frs^{-1})_{lu} 
\nonumber\\ &
\rho_\mathrm{pMD}(\frs)_{lu} 
\sum_i
\rho_\mathrm{pMD}(\frs)_{li} \nu^\text{pMD}_n(i) 
=
 \rho_\mathrm{pMD}(\frs\frt^n\frs)_{lu} \rho_\mathrm{pMD}(\frs\frt^{-n}\frs^{-1})_{lu} 
\nonumber\\
\nu^\text{pMD}_n(i) v_iv_u &= \nu_n(i) v_iv_u, \ \ \ 
\nu^\text{pMD}_n(u)=1, \ \ \  
\nu^\text{pMD}_n(u) = \pm  \rho_\mathrm{pMD}(\frs^2)_{ii}
\end{align}

\end{enumerate}
\end{prop}

Here we like to remark that the condition involving central charge
\begin{align}
p_\pm = \sum_{i=1}^{r} d_i^2 \theta_i^{\pm 1},\ \ \
p_\pm=D\ee^{\pm \ii 2\pi c/8},
\end{align}
is not a new condition. It comes from the $\SL$ condition.  First,
\eqref{Ss2tSs} can be rewritten as
\begin{align}
\label{tSs2tSs}
 \rho_\al(\frs) =  \ee^{2\pi \ii \frac{\t c}{8}} \frac{\rho_\al(\frs)/\rho_\al(\frs)_{uu}}{p_+} , \ \
 \rho_\al(\frt) =  \ee^{-2\pi \ii \frac{\t c}{24}} \rho_\al(\frt)/\rho_\al(\frt)_{uu},
\end{align}
for a $\t c \in \Q$.
We find $\rho_\al(\frt)_{uu} =  \ee^{-2\pi \ii \frac{\t c}{24}}$,
which allows us to rewrite \eqref{tSs2tSs} as
\begin{align}
\rho_\al(\frs)_{uu} \rho_\al^3(\frt)_{uu} p_+ 
= 
\rho_\al(\frs)_{uu} \rho_\al^3(\frt)_{uu} 
\sum_{i=1}^{r} 
\Big(\frac{\rho_\al(\frs)_{iu}}{\rho_\al(\frs)_{uu}}\Big)^2  
\frac{\rho_\al(\frt)_{ii}}{\rho_\al(\frt)_{uu}}  
= 1,
\end{align}
and becomes an condition on $\rho_\al(\frs),\rho_\al(\frt)$.
But this is not a new condition, since
the above can be rewritten as
\begin{align}
\rho_\al^2(\frt)_{uu}
\sum_{i=1}^{r} 
\rho_\al(\frs)_{ui}
\rho_\al(\frs)_{iu}
\rho_\al(\frt)_{ii} 
= \rho_\al(\frs)_{uu} 
\ \ \text{ or }\ \
\big(
\rho_\al(\frt)
\rho_\al(\frs)
\rho_\al(\frt)
\rho_\al(\frs)
\rho_\al(\frt)
\big)_{uu}
= \rho_\al(\frs)_{uu},
\end{align}
which is a consequence of $\SL$ representation.

Using the irreducible $\SL$ representations obtained by GAP package SL2Reps, we
can explicitly constructed all unitary representations of $\SL$ (up to unitary
equivalence).  However, this only gives the $\SL$ representations in some
arbitrary basis, not in the basis yielding MD representations (\emph{i.e.}
satisfying Proposition \ref{p:MDcond}), since   MD representations are $\SL$
representations in a particular basis.

We can improve the situation by choosing a basis to make $\rho(\frt)$ diagonal
and $\rho(\frs)$ symmetric.  We can choose more special bases to make the $\SL$
representations closer to the basis of MD representations.  Since we are going
to use several types of bases, let us define these choices:
\begin{defn}
An unitary $\SL$ representations $\tilde\rho$ is called a \textbf{general}
$\SL$ matrix representations if $\tilde\rho(\frt)$ is diagonal \footnote{We
will consider only $\SL$ matrix representations with diagonal $\tilde\rho(\frt)$
in this paper.}.  A general $\SL$ matrix representation $\tilde \rho$ is called
\textbf{symmetric} if $\tilde\rho(\frs)$ is symmetric.  
An general $\SL$ matrix representation $\tilde\rho$ is called
\textbf{irrep-sum} if $\tilde\rho(\frs),\tilde\rho(\frt)$ are matrix-direct sum
of irreducible $\SL$ representations.  An $\SL$ matrix representations
$\tilde\rho$ is called an $\SL$ representation \textbf{of modular data} $S,T$,
if $\tilde\rho$ is unitarily equivalent to an MD
representation of the modular data, {\it i.e.},
\begin{align}
\label{STrhoU1}
\tilde\rho(\frs) = \ee^{-2\pi \ii \frac{\al}{4}}\frac{1}{D}\, US U^\dag, \ \ \ \
\tilde\rho(\frt) = UTU^\dag \ee^{2\pi \ii (\frac{-c}{24} + \frac{\al}{12})},
\end{align}
for some unitary matrix $U$ and $\al \in \Z_{12}$, where $c$ is the central
charge.\footnote{Note that $D^2$ is always positive and $D$ in \eqref{STrhoU1}
is the positive root of $D^2$, even for non-unitary cases.} An $\SL$ matrix
representations $\rho_\mathrm{pMD}$ is called \textbf{pseudo MD (pMD)}
representation, if $\rho_\mathrm{pMD}$ is related to an MD representation of the
modular data via a conjugation of a \textbf{signed diagonal} matrix $V$:\footnote{A
signed diagonal matrix is a diagonal matrix with diagonal elements $\pm 1$.}
\begin{align}
\label{STrhoU}
\rho_\textrm{pMD}(\frs) = \ee^{-2\pi \ii \frac{\al}{4}}\frac{1}{D}\, V S V, \ \ \ \
\rho_\textrm{pMD}(\frt) = VTV \ee^{2\pi \ii (\frac{-c}{24} + \frac{\al}{12})}.
\end{align}
\end{defn}

We find that all irreducible unitary representations of $\SL$ are unitarily
equivalent to symmetric matrix representations of $\SL$.  We will start with
those symmetric matrix representations to obtain a classification of modular
data.

\section{Our strategy }

We will use the following strategy to classify modular data of a given rank.

\begin{enumerate}
\item 
Obtain all the irreducible representations of dimensions up to the rank $r$,
using GAP package SL2Reps created by Siu-Hung Ng, Yilong Wang, and Samuel
Wilson.  Then construct all the dimension-$r$ representations,
$\rho_\textrm{isum}$, from those irreducible representations.

\item Using some conditions (see Section \ref{MDreps}), to reject the
representations that are not unitarily equivalent to MD representations.

\item For the remain representations $\rho_\textrm{isum}$, find
all the unitary matrices $U$ that transform them to pMD representations (see
Section \ref{sec:pMD}):
\begin{align}
\rho_\textrm{pMD} = U\rho_\textrm{isum} U^\dag.
\end{align}
Reject those representations for which the unitary matrices $U$ do not exist.
This is the most difficult step, since we need to find finite solutions from
infinite possibilities of unitary $U$.  

The key is to generate many conditions on the unitary $U$, so that the number of the solutions 
of those condition is finite. To achieve this,  we first
consider all the possible $D_\rho(\si)$'s (see Section \ref{calDrho}) and unit-row index $u$.  Since
$D_\rho(\si)$'s are signed permutations, those possibilities are finite.  Once
$D_\rho(\si)$ and $u$ are known, we can use them to obtain many conditions on
$U$, in addition to the unitary conditions.  

Also, from $D_\rho(\si)$, we can
determine if the corresponding MTC is integral or not.  This allows us to use
two different approaches to handle integral cases (see Section \ref{intcat})
and non-integral cases (see Section \ref{sec:pMD}) separately.  The integral
and non-integral cases are quite different, and require different approaches to
handle them.

Once we know the unit-row index $u$, we can obtain conditions on $U$ from the
second Frobenius-Schur indicator (see \eqref{nunFSapp}). Also
$\textrm{norm}(\rho(\frs)_{ui}) = q_i$ are inverse of integers.
If we can isolate the polynomial conditions that depend only on $q_i$'s (see Section \ref{reducevars}), then 
we can use the generalized Egyptian-fraction method to solve $q_i$'s (see Section \ref{gEFrac}).
Those are important tricks that make our calculation possible.

\item Find all the signed permutations $P_\textrm{sgn}$
that transform the pMD representations to MD representations:
\begin{align}
S &= P_\textrm{sgn} \frac{\rho_\textrm{pMD}(\frs)}{(\rho_\textrm{pMD})_{uu}(\frs)} P_\textrm{sgn}^\top ,
&
T &= P_\textrm{sgn} \frac{\rho_\textrm{pMD}(\frt)}{(\rho_\textrm{pMD})_{uu}(\frt)} P_\textrm{sgn}^\top
.
\end{align}
This step is easy, since the number of possible signed permutations is finite.
We just search all the  signed permutations so that the resulting $S,T$ satisfy
the conditions for modular data.  The signed permutations have the form
$P_\textrm{sgn} = P V_\textrm{sd}$, where $P$ is a permutation matrix and
$V_\textrm{sd}$ is a signed diagonal matrix.  We may fix $P$ and only search for
$V_\textrm{sd}$. From the Verlinde formula \eqref{SSN}, we know how $N^{ij}_k$
transform under the conjugation of $V_\textrm{sd}$.  This allows us to find
$V_\textrm{sd}$ that make $N^{ij}_k$ non-negative.

\end{enumerate}


\section{Candidate representations of modular data from
$\SL$ representations }
\label{MDreps}

We note that different choices of orthogonal basis give rise to different
matrix representations of $\SL$. The modular data $S,T$ are obtained from some
particular choices of the basis.  Some properties of the MD representations of
a modular data do not depend on the choices of basis in the eigenspaces of
$\tilde\rho(\frt)$ (induced by the block-diagonal unitary transformation $U$ in
\eqref{STrhoU} that leaves $\tilde\rho(\frt)$ invariant). Those properties
remain valid for any general $\SL$ representations $\tilde\rho$ of the modular
data.  In the following, we collect the basis-independent conditions on the
$\SL$ matrix representations of modular data.  This will help us to narrow the
list of $\SL$ representations that are related modular data.
\begin{prop}
\label{p:gencond1}
Let $\tilde\rho$ be a general $\SL$ matrix representations of a modular data or
a MTC.  Then $\tilde\rho$ must satisfy the following conditions:
\begin{enumerate}
\item 
If $\tilde\rho$ is a direct sum of two $\SL$ representations
\begin{align}
 \tilde\rho \cong \rho \oplus \rho',
\end{align}
then the diagonals entries
of $\rho(\frt)$ and $\rho'(\frt)$ must overlap.  

\item
Suppose that $\tilde\rho \cong \rho \oplus \ell \chi$ for an irreducible
representation $\rho$ with $\rho(\frt)$ non-degenerate, and a character $\chi$.  If $\ell \neq 1$ and $\ell \ne 2\dim(\rho) -1$, then
$(\rho(\frs) \chi(\frs)^{-1})^2 = \id$. 

\item
If $\tilde\rho(\frs)^2=\pm \id$, and $\pord(\tilde\rho(\frt)) = 1$ mod 4 and is a
prime, then the representation $\t\rho$ cannot be a direct sum of a
$d$-dimensional irreducible $\SL$ representation and two or more 1-dimensional
$\SL$ representations with $d=(\pord(\tilde\rho(\frt)) +1)/2$.  

\item $\tilde\rho$ satisfies
\begin{align}
\tilde\rho \not\cong n \rho
\end{align}
for any integer $n >1$ and any  representation $\rho$ such that $\rho(\frt)$ is non-degenerate.

\item
Let $3< p < q$ be primes such that $pq \equiv 3 \mod 4$ and
$\pord(\rho(\frt))=pq$, then the rank $r \ne \frac{p+q}{2}+1$.  Moreover, if
$p > 5$, rank $r > \frac{p+q}{2}+1$.


\end{enumerate}
\end{prop}

Some other properties of an MD representation do depend on the choice of basis.
To make use of those properties, we can construct some combinations of
$\tilde\rho(\frs)$'s that are invariant under the block-diagonal unitary
transformation $U$.  

The eigenvalues of $\tilde\rho(\frt)$ partition the indices of the basis vectors. To construct the invariant combinations of $\tilde\rho(\frs)$, for any eigenvalue $\tilde\theta$ of $\tilde\rho(\frt)$, let
\begin{align}
\label{Ith}
 I_{\tilde\theta} = \{ i \, \big|\, \tilde\rho(\frt)_{ii}=\tilde\theta\}. 
 \end{align} 
Let $I = I_{\tilde\theta}$, $J = J_{\tilde\theta'}$, $K = K_{\tilde\theta''}$ for some eigenvalues $\tilde\theta$, $\tilde\theta'$, $\tilde\theta''$ of $\tilde\rho(\frt)$.
We see that the following uniform polynomials of $\tilde\rho(\frs)$ are
invariant 
\begin{align}
\label{invcomb}
P_I(\rho(\frs)) =
\Tr \tilde\rho(\frs)_{II} &\equiv \sum_{i\in I} \tilde\rho(\frs)_{ii},
\nonumber\\
P_{IJ}(\rho(\frs)) =
\Tr \tilde\rho(\frs)_{IJ} \tilde\rho(\frs)_{JI} &\equiv 
\sum_{i\in I,j\in J} \tilde\rho(\frs)_{i,j} \tilde\rho(\frs)_{ji},
\\
P_{IJK}(\rho(\frs)) =
\Tr \tilde\rho(\frs)_{IJ} \tilde\rho(\frs)_{JK} \tilde\rho(\frs)_{KI} 
&\equiv \sum_{i\in I,j\in J,k\in K} 
\tilde\rho(\frs)_{i,j} 
\tilde\rho(\frs)_{j,k}
\tilde\rho(\frs)_{k,i}
.
\nonumber 
\end{align}
Certainly we can construction many other invariant uniform polynomials in the
similar way.  Using those invariant uniform polynomials, we have the following
results
\begin{prop}
\label{p:gencond2}
Let $\tilde\rho$ be a general $\SL$ representations of a modular data or a
MTC.  Then following statements hold:
\begin{enumerate}

\item $\tilde\rho(\frs)$ satisfies
\begin{align}
 \Tr(\tilde\rho(\frs)^2)  \in \Z \setminus\{0\} .
\end{align}
Let 
\begin{align}
\label{Crho}
 C = \frac{\Tr(\tilde\rho(\frs)^2)}{|\Tr(\tilde\rho(\frs)^2)|}   \tilde\rho(\frs)^2.
\end{align}
For all $I$, 
\begin{align}
\label{PIC}
 P_I(C) \geq 0.
\end{align}

\item 
\label{cndI}
The conductor of $P_\mathrm{odd}(\tilde\rho(\frs))$ divides
$\ord(\tilde\rho(\frt))$ for all the invariant uniform polynomials
$P_\mathrm{odd}$ with odd powers of $\tilde\rho(\frs)$ (such as $P_I$ and
$P_{IJK}$ in \eqref{invcomb}).  The conductor of
$P_\mathrm{even}(\tilde\rho(\frs))$ divides $\pord(\tilde\rho(\frt))$ for all the
invariant uniform polynomials $P_\mathrm{even}$ with even powers of
$\tilde\rho(\frs)$ (such as $P_{IJ}$ in \eqref{invcomb}).

\item 
 For any Galois conjugation $\si\in\Gal(\Q_{\ord(\rho(\frt))})$, there is a
permutation on the set $\{I\}$, $I \to \hs(I)$, such that
\begin{align}
\si P_{IJ}(\tilde\rho(\frs)) &
= P_{I\hs(J)}(\tilde\rho(\frs))
=P_{\hs(I)J}(\tilde\rho(\frs))
\nonumber\\
\si^2 \big(\tilde\theta_I\big) &= \tilde\theta_{\hat \si (I)},
\end{align}
for all $I,J$.

\item 
For any invariant uniform polynomials $P$ (such as those in \eqref{invcomb})
\begin{align}
\si P\big(\tilde\rho(\frs)\big) = 
 P\big(\si \tilde\rho(\frs)\big) = 
P \big(\tilde\rho(\frt)^a \tilde\rho(\frs) \tilde\rho(\frt)^b \tilde\rho(\frs) \tilde\rho(\frt)^a\big)
\end{align}
where $\si\in\Gal(\Q_{\ord(\t\rho(\frt))})$,
and $a,b$ are given by $\si( \ee^{\ii 2
\pi/\ord(\t\rho(\frt))} ) = \ee^{a \ii 2 \pi/\ord(\t\rho(\frt))}$ and $ab \equiv 1$ mod
$\ord(\t\rho(\frt))$.  
\end{enumerate}
\end{prop}


Instead of constructing invariants, there is another way to make use of the
properties of an MD representation that depend on the choices of basis.  We can
choose a more special basis, so that the basis is closer to the basis that
leads to the MD representation.  For example, we can choose a basis to make
$\tilde\rho(\frs)$ symmetric ({\it i.e.} to make $\tilde\rho$ a symmetric
representation).

Now consider a symmetric $\SL$ matrix representation $\tilde \rho$ of a modular
data or a MTC.  We find that the restriction of the unitary $U$ in
\eqref{STrhoU} on the non-degenerate subspace (see \Rf{NW220314829} Theorem
3.4) must be a signed diagonal matrix.  In this case some properties of MD
representation apply to the blocks of the symmetric representation within the
non-degenerate subspace.  This allows us to obtain
\begin{prop}
\label{p:gal_sym}
Let $\tilde\rho$ be a symmetric $\SL$ representations equivalent to an MD representation.
 Let 
\begin{align}
I_\mathrm{ndeg} &:=\{i \mid \tilde\rho(\frt)_{i,i} \text{ is a
non-degenerate eigenvalue}\},
\end{align}
Then there exists an orthogonal $U$ such that $U \tilde \rho U^\top$ is a pMD representation, and the following statements hold:

\begin{enumerate}
\item 
\label{cndndegI}
The conductor of $(U \tilde\rho(\frs) U^\top)_{i,j} $ divides
$\ord(\tilde\rho(\frt))$ for all $i,j$.  This implies that the conductor of
$(\tilde\rho(\frs))_{i,j} $ divides $\ord(\tilde\rho(\frt))$ for all $i,j \in
I_\mathrm{ndeg}$.

\item 
 For any Galois conjugation $\si$ in $\Gal(\Q_{\ord(\t\rho(\frt))})$, there is a
permutation $i \to \hs(i)$, such that
\begin{align}
\si \big((U \tilde\rho(\frs)U^\top)_{i,j}\big) &
= \eps_\si(i) (U \tilde\rho(\frs) U^\top)_{\hat \si (i),j} 
= (U \tilde\rho(\frs) U^\top)_{i,\hat \si (j)} \eps_\si(j)
\nonumber\\
\si^2 \big( \tilde\rho(\frt)_{i,i}\big) &= 
 \tilde\rho(\frt)_{\hat \si (i),\hat \si (i)},
\end{align}
for all $i,j$, where $\eps_\si(i)\in \{1,-1\}$. 
This implies that
\begin{align}
\si \big( \tilde\rho(\frs)_{i,j}\big) & = \tilde\rho(\frs)_{\hat \si (i),j} 
\ \ \mathrm{or}  \ \
\si \big( \tilde\rho(\frs)_{i,j}\big)  = -  \tilde\rho(\frs)_{\hat \si (i),j} 
\nonumber\\
\si \big( \tilde\rho(\frs)_{i,j}\big) & =  \tilde\rho(\frs) _{i,\hat \si (j)} 
\ \ \mathrm{or}  \ \
\si \big( \tilde\rho(\frs)_{i,j}\big)  = - \tilde\rho(\frs) _{i,\hat \si (j)} 
\end{align}
for all $i,j \in I_\mathrm{ndeg}$.  This also implies that $D_{\t\rho}(\si)$
defined in \eqref{eq:Drho} is a signed permutation matrix in the $I_\mathrm{ndeg}$
block, {\it i.e.} $(D_{\t\rho}(\si))_{i,j}$ for $i,j \in I_\mathrm{ndeg}$ are
matrix elements of a signed permutation matrix.

\item 
\label{sirhondegI}
For all $i,j$, 
\begin{align}
\si\big((U\tilde\rho(\frs) U^\top)_{i,j}\big) = 
\big(U \tilde\rho(\frt)^a \tilde\rho(\frs) \tilde\rho(\frt)^b \tilde\rho(\frs) \tilde\rho(\frt)^a U^\top\big)_{i,j}
\end{align}
where $\si\in\Gal(\Q_{\ord(\t\rho(\frt))})$, and
$a,b$ are given by $\si( \ee^{\ii 2 \pi/\ord(\t\rho(\frt))} ) = \ee^{a \ii 2
\pi/\ord(\t\rho(\frt))}$ and $ab \equiv 1$ mod $\ord(\t\rho(\frt))$.  
This implies that
\begin{align}
\label{sirhondeg}
\si\big((\tilde\rho(\frs) )_{i,j}\big) = 
\big( \tilde\rho(\frt)^a \tilde\rho(\frs) \tilde\rho(\frt)^b \tilde\rho(\frs) \tilde\rho(\frt)^a \big)_{i,j}.
\end{align}
for all $i,j \in I_\mathrm{ndeg}$.

\item
\label{nonzeroI} Both $T$ and $\t\rho(\frt)$ are diagonal, and without loss of
generality, we may assume $\t\rho(\frt)$ is a scalar multiple of $T$. In this
case $U$ in \eqref{STrhoU} is a block diagonal matrix preserving the
eigenspaces of $\t\rho(\frt)$.  Let $I_\mathrm{nonzero} = \{i\}$ be a set of
indices such that the $i^\mathrm{th}$ row of  $U \t\rho(\frs) U^\top$ contains
no zeros for some othorgonal $U$ satisfying $U\t\rho(\frt)U^\top =
\t\rho(\frt)$.  The index for the unit object of MTC must be in
$I_\mathrm{nonzero}$.  Thus $I_\mathrm{nonzero}$ must be nonempty:
\begin{align}
\label{Ineq0}
 I_\mathrm{nonzero} \neq \emptyset.
\end{align}

\item 
Let $I_{\t\theta}$ be a set of indices for an eigenspace $E_{\t\theta}$
of $\tilde\rho(\frt)$
\begin{align}
I_{\t\theta} &:=\{i \mid \tilde\rho(\frt)_{i,i} = \t\theta \}.
\end{align}
Then there exists a $I_{\t\theta}$ such that
\begin{align}
\label{TrCunit}
I_{\t\theta} \cap I_\mathrm{nonzero} \neq \emptyset\ \text{ and }\
\Tr_{E_{\t\theta}} C > 0,
\end{align}
where $C$ is given in \eqref{Crho}.



\end{enumerate}
\end{prop}

It is very helpful to determine if a $\SL$ representation gives rise to an
integral MTC or not.  This is because  integral MTCs satisfy more conditions.
\begin{thm}
\label{intMTCu}
Let $\t\rho$ be a representations of $\SL$.  If $u^\text{th}$ row of $\t\rho$ is
the unit row, and $D_{\t\rho}(\si)_{uu} = \pm 1$ for all $\si \in
\Gal(\Q_{\ord(\t\rho(\frt))}/\Q)$, then $\t\rho$ is either equivalent to an MD
representation $\rho_\mathrm{MD}$ of integral MTC or is not equivalent to any
MD representation.
\end{thm}
The above result comes from Eq. \eqref{Drhosihat}.  
$D_{\t\rho}(\si)_{uu} = \pm 1$ implies that
$D_{\rho_\mathrm{MD}}(\si)_{uu} = \pm 1$, which in turn implies 
that $\hat\si(u)=u$.
Then Eq. \eqref{Galact} implies that, for the corresponding MD representation
$\rho_\mathrm{MD}$ of $\t\rho$,
\begin{align}
\label{Galdi}
\si \Big(\frac{\rho_\mathrm{MD}(\frs)_{ij}}{\rho_\mathrm{MD}(\frs)_{iu}} \Big) &
= \frac{\rho_\mathrm{MD}(\frs)_{\hat \si (i)j}}{\rho_\mathrm{MD}(\frs)_{\si (i)u}} 
\end{align}
for all $i,j$, and hence $\si(d_i) = d_i$ for all $\si \in
\Gal(\Q_{\ord(\t\rho(\frt))}/\Q)$, where
$d_i=\frac{\rho_\mathrm{MD}(\frs)_{ui}}{\rho_\mathrm{MD}(\frs)_{uu}}$ is the
quantum dimension.
Thus, $d_i$ are integral, and the corresponding MTC is integral, if it exists.
If we do not know the unit row, then we have
\begin{thm}
\label{intMTC}
Let $\{\t\theta\}_\mathrm{nonzero}$ is a set the eigenvalues of $\t\rho(\frt)$,
$\t\theta$, such that $\t\theta$ is a 24$^\mathrm{th}$ root of unity and
$I_{\t\theta}$ has a non-empty overlap with $I_\mathrm{nonzero}$.  If
$D_{\t\rho}(\si)_{I_{\t\theta}} = \pm \id$ for all $\si \in
\Gal(\Q_{\ord(\t\rho(\frt))}/\Q)$ and for all $\t\theta$ in
$\{\t\theta\}_\mathrm{nonzero}$, then $\t\rho$ is either equivalent to an MD
representation of integral MTC or is not equivalent to any MD representation.
\end{thm}
We also have 
\begin{thm}
If a $\SL$ representation $\t\rho$ satisfies $\pord(\t\rho(\frt)) = \ord(T) \in
\{2,3,4,6\}$, then $\t\rho$ is either equivalent to an MD representation of
integral MTC or is not equivalent to any MD representation.
\end{thm}

\begin{table}[tb] 
\caption{The numbers of
the candidate irrep-sum $\SL$ representations for each rank.
} 
\label{Nirrep} 
\centering
\begin{tabular}{ |c|c|c|c|c|c|c|c|c|c|c|c| } 
\hline 
rank           & 2 & 3 & 4 &  5 &  6 &  7 &   8 &   9 &   10 &   11 &   12 \\
\hline 
 number of reps & 2 & 4 & 9 & 20 & 57 &106 & 258 & 533 & 1210 & 2374 & 5288 \\ 
\hline 
\end{tabular}
\end{table}

Using GAP System for Computational Discrete Algebra, we obtain a list of
symmetric irrep-sum $\SL$ matrix representations that satisfy the conditions in
Propositions \ref{p:gencond1}, \ref{p:gencond2}, and \ref{p:gal_sym}.  Also,
our list only includes one representive for each orbit generated by Galois
conjugations and tensoring 1-dim $\SL$ representations.
The numbers of
those candidate irrep-sum $\SL$ representations for each rank are given in Table \ref{Nirrep}.

Some of those symmetric irrep-sum $\SL$ matrix representations are
representations of modular data, while others are not.  However, the list
includes all the symmetric irrep-sum $\SL$ matrix representations of modular
data or MTC's.  In the next section, we will use GAP group to determine which
of those irrep-sum representations can give rise to modular data, and which
should be rejected.  However, there are a few cases at rank-12 are hard to
handle.  We have to use extensive search to reject those cases with high
likelyhood.  Those calculations are presented in the first few sections of the
Appendix.

\section{Candidate \lowercase{p}MD representations from 
$\SL$ representations }
\label{sec:pMD}

Our $\SL$ representation $\t\rho$ has a form of direct sum of irreducible
representations: $\t\rho = \rho_\textrm{isum}$.  We have chosen a
special basis in the eigenspaces of a $\SL$ matrix representation
$\rho_\textrm{isum}$ to make $\rho_\textrm{isum}(\frs)$ symmetric.  But such a
special basis is still not special enough to make $\rho_\textrm{isum}$ to be a
MD representation $\rho_\al$.

We can choose a more special basis to make $\rho_\textrm{isum}(\frs^2)$ a
signed permutation matrix, in addition to making $\rho_\textrm{isum}(\frs)$
symmetric.  We know that, for a MD representation $\rho_\al$, $\rho(\frs^2)$ is
a signed permutation matrix. So the new special basis makes
$\rho_\textrm{isum}$ closer to the MD representation $\rho_\al$.

We can choose an even more special basis in the eigenspaces of
$\rho_\textrm{isum}(\frt)$ to make $\rho_\textrm{isum}$ into a pseudo MD
representation that differs from a MD representation $\rho_\al$ only by the
conjugation of a signed diagonal matrix $V_\mathrm{sd}$: $ \rho_\mathrm{pMD} =
V_\mathrm{sd} \rho_\al V_\mathrm{sd}$.  Pseudo MD representation has a property
that the matrix $D_{\rho_\mathrm{pMD}}(\si)$ defined in \eqref{defDrho} are
signed permutations.

We would like to point out that, since both $\rho_\textrm{isum}$ and
$\rho_\mathrm{pMD}$ are symmetric $\SL$ matrix representations that are related
by an unitary transformation, according to Theorem 3.4 in \Rf{NW220314829},
they can be related by an orthogonal transformation.  An generic orthogonal
transformation contains continuous real parameters.  This leads to infinite many
potential pseudo MD representations $\rho_\mathrm{pMD}$ and infinite many
potential MD representations $\rho_\al$.  This makes it impossible to check
one-by-one, if those potential  $\rho_\al$'s are indeed  MD representations.

However, when eigenspaces of $\rho_\textrm{isum}(\frt)$ are all 1-dimensional
(\ie non-degenerate), the orthogonal matrices $U$ that transform
$\rho_\textrm{isum}$ to $\rho_\mathrm{pMD}$ must be an identity matrix, up to
signed diagonal matrices.  This leads to only a finite many potential MD
representations $\rho_\al$.  We can then check each of the  possible
$\rho_\al$'s, to see if it is a MD representation.

Even if some eigenspaces of $\rho_\textrm{isum}(\frt)$ are degenerate, under
certain conditions, the number of orthogonal transformations $U$ that transform
$D_{\rho_\textrm{isum}}(\si)$ into signed permutations
$D_{\rho_\textrm{pMD}}(\si)$ can still be finite.  Let $I_{\th}$ is the set of
indices for the degenerate eigenspace of $\rho_\mathrm{pMD}(\frt)$ with
eigenvalue $\th$.  Let us consider $D_{\rho_\textrm{isum}}(\si)$ in the
$I_{\th}$-block, which are denoted as
$D_{\rho_\textrm{isum}}(\si)_{I_\th,I_\th}$.  If the common eigenspaces for
$D_{\rho_\textrm{isum}}(\si)_{I_\th,I_\th}$ with different $\si$'s are
non-degenerate, then there is only a finite number of orthogonal
transformations in the $I_\th$ block, $U_{I_\th,I_\th}$, that transform
$D_{\rho_\textrm{isum}}(\si)$ into signed permutations
$D_{\rho_\textrm{pMD}}(\si)$.  

\subsection{Cases where
$D_{\rho_\textrm{isum}}(\si)_{I_\th,I_\th}$ are signed diagonal matrices
with non-degenerate common eigenspaces}

For example, when the non-zero $D_{\rho_\textrm{isum}}(\si)_{I_\th,I_\th}$ are
signed diagonal matrices with non-degenerate common eigenspaces, the most
general orthogonal matrices $U_{I_\th,I_\th}$ have only a finite number of
choices.  They must be either an identity matrix, or a combinations of $\pm
45^\circ$ rotations among various pairs of indices, up to signed diagonal
matrices,

Let us consider a concrete case. In a 3-dimensional eigenspace of
$\rho_\textrm{isum}(\frt)$, the non-zero $D_{\rho_\textrm{isum}}(\si)$'s may generate a $3\times 3$
matrix groups $MG$, given by
\begin{align}
\label{MG3x3}
 MG &=\Big\{  
\begin{pmatrix} 1 &0 & 0\\ 0 & 1 & 0\\ 0&0&1\\ \end{pmatrix}, 
\begin{pmatrix} -1 &0 & 0\\ 0 & -1 & 0\\ 0&0&1\\ \end{pmatrix}, 
\begin{pmatrix} 1 &0 & 0\\ 0 & -1 & 0\\ 0&0&1\\ \end{pmatrix}, 
\begin{pmatrix} -1 &0 & 0\\ 0 & 1 & 0\\ 0&0&1\\ \end{pmatrix} 
\Big\}. 
\end{align}
To find the most general orthogonal matrices that transform the above $3\times
3$ matrices in $MG$ into signed permutation matrices, we first show
\begin{thm}
\label{PPsgn}
If $P$ is a permutation matrix with $P^2 = \id$, then $P$ is a direct sum of
$2\times 2$ and $1\times 1$ matrices.  If $P_\mathrm{sgn}$ is a signed
permutation matrix with $P_\mathrm{sgn}^2 = \id$, then $P_\mathrm{sgn}$ is a
direct sum of $2\times 2$ and $1\times 1$ matrices. The $2\times 2$ matrices
are the $\pm 45^\circ$ rotations mentioned above.
\end{thm}
\begin{proof}[Proof of Theorem \ref{PPsgn}] 
If $P$ is a permutation matrix with $P^2 = \id$, $P$ must be a pair-wise
permutation, and thus $P$ is a direct sum of $2\times 2$ and $1\times 1$
matrices.  The reduction from signed permutation matrix to permutation matrix
by ignoring the signs is homomorphism of the matrix product.  If
$P_\mathrm{sgn}$ is a signed permutation matrix with $P_\mathrm{sgn}^2 = \id$,
its reduction given rise to a permutation matrix $P$ with $P^2 = \id$.  Since
$P$ is a direct sum of $2\times 2$ and $1\times 1$ matrices, $P_\mathrm{sgn}$
is also a direct sum of $2\times 2$ and $1\times 1$ matrices.
\end{proof}

Using the above result, we can show that
the most general orthogonal matrices that transform all  $D_{
\rho_\mathrm{isum}}(\si_\mathrm{inv})|_{E_{\t\theta}}$'s into signed
permutations must have one of the following forms 
\begin{align}
\label{U3x3a}
 U &=  
\frac{PV_\mathrm{sd}}{\sqrt 2}  
\begin{pmatrix}
 1& 1 & 0\\
 1&-1 & 0\\
 0& 0 & 1\\
\end{pmatrix} ,
 \text{ or } \
 U =   \frac{PV_\mathrm{sd}}{\sqrt 2} 
 \begin{pmatrix}
-1& 1 & 0\\
 1& 1 & 0\\
 0& 0 & 1\\
\end{pmatrix} ,
\nonumber\\
\text{or }\  U &=  
\frac{PV_\mathrm{sd}}{\sqrt 2}  
\begin{pmatrix}
 1& 0& 1 \\
 0& 1& 0\\
 1& 0&-1 \\
\end{pmatrix} ,
 \text{ or } \
 U =   \frac{PV_\mathrm{sd}}{\sqrt 2} 
\begin{pmatrix}
-1& 0& 1 \\
 0& 1& 0\\
 1& 0& 1 \\
\end{pmatrix} ,
\nonumber\\
\text{or }\  U &=  
\frac{PV_\mathrm{sd}}{\sqrt 2}  
\begin{pmatrix}
 1& 0& 0 \\
 0& 1& 1\\
 0& 1&-1 \\
\end{pmatrix} ,
 \text{ or } \
 U =   \frac{PV_\mathrm{sd}}{\sqrt 2} 
\begin{pmatrix}
 1& 0& 0 \\
 0&-1& 1\\
 0& 1& 1 \\
\end{pmatrix} ,
\nonumber\\
\text{or }\  U &=  
PV_\mathrm{sd}
\begin{pmatrix}
 1& 0& 0 \\
 0& 1& 0\\
 0& 0& 1 \\
\end{pmatrix} .
\end{align}
where $V_\mathrm{sd}$ are signed diagonal matrices, and $P$ are permutation
matrices.  We note that the non-trivial part of $U$ is a $2\times 2$ block for
index $(1,2)$, $(1,3)$, and $(2,3)$.  The $2\times 2$ block has three
possibilities 
\begin{align}
\begin{pmatrix}
 1& 0 \\
 0& 1\\
\end{pmatrix},\ \ \ \
\frac{1}{\sqrt 2} 
\begin{pmatrix}
 1& 1 \\
 1&-1\\
\end{pmatrix}, \ \ \ \
\frac{1}{\sqrt 2} 
\begin{pmatrix}
-1& 1 \\
 1& 1\\
\end{pmatrix}.
\end{align}
This is a general pattern that apply for cases when the non-zero
$D_{\rho_\textrm{isum}}(\si)_{I_\th,I_\th}$ are signed diagonal matrices with
non-degenerate common eigenspaces.

Some times, the non-zero $D_{\rho_\textrm{isum}}(\si)_{I_\th,I_\th}$ are not
signed diagonal matrices.  We need to examine those cases individually.

\subsection{Within a 2-dimensional eigenspace of $\rho_\mathrm{isum}(\frt)$}

In this case, the matrix groups $MG$ generated by non-zero 2-by-2 matrices, $D_{
\rho_\mathrm{isum}}(\si)_{I_\th,I_\th}$, can have several
different forms, for those passing representations. By examine the
computer results, we find that matrix groups $MG$ can be
\begin{align}
 MG &=\Big\{  
\begin{pmatrix} 1 &0\\ 0 & 1\\ \end{pmatrix} 
\Big\}, 
&&\text{ for } \dim(\rho_\mathrm{isum}) \geq 5;
\nonumber\\
 MG &=\Big\{  
\begin{pmatrix} 1 &0\\ 0 & 1\\ \end{pmatrix}, 
-\begin{pmatrix} 1 &0\\ 0 & 1\\ \end{pmatrix} 
\Big\}, 
&&\text{ for } \dim(\rho_\mathrm{isum}) \geq 6.
\end{align}
Those $D_{\rho_\textrm{isum}}(\si)_{I_\th,I_\th}$'s have degenerate common
eigenspaces.  Thus the resulting $U_{I_\th,I_\th}$'s have continuous parameters
and are not finite many.

We also have 
\begin{align}
 MG &=\Big\{  
\begin{pmatrix} 1 &0\\ 0 & 1\\ \end{pmatrix}, 
\begin{pmatrix} 1 &0\\ 0 & -1\\ \end{pmatrix} 
\Big\}, 
&&\text{ for } \dim(\rho_\mathrm{isum}) \geq 4;
\nonumber\\
 MG &=\Big\{  
\begin{pmatrix} 1 &0\\ 0 & 1\\ \end{pmatrix}, 
-\begin{pmatrix} 1 &0\\ 0 & 1\\ \end{pmatrix}, 
\begin{pmatrix} 1 &0\\ 0 & -1\\ \end{pmatrix}, 
\begin{pmatrix} -1 &0\\ 0 & 1\\ \end{pmatrix} 
\Big\}, 
&&\text{ for } \dim(\rho_\mathrm{isum}) \geq 6.
\end{align}
In those two cases
\begin{align}
\label{U2x2a}
 U = \frac1{\sqrt 2} \begin{pmatrix}
 1& 1\\
 1&-1\\
\end{pmatrix}
\ \ \text{ or } \ \
 U = \frac1{\sqrt 2} \begin{pmatrix}
-1& 1\\
 1& 1\\
\end{pmatrix}
\ \ \text{ or } \ \
 U = \begin{pmatrix}
 1& 0\\
 0& 1\\
\end{pmatrix}
\end{align}
will transform all  $D_{
\rho_\mathrm{isum}}(\si_\mathrm{inv})_{I_\th,I_\th}$'s into signed
permutations.  
In general we have
\begin{thm}
\label{MG2x2}
Let
\begin{align}
 MG_2 &=\Big\{  
\begin{pmatrix} 1 &0\\ 0 & 1\\ \end{pmatrix}, 
\begin{pmatrix} 1 &0\\ 0 & -1\\ \end{pmatrix} 
\Big\}, 
\nonumber\\
 MG_4 &=\Big\{  
\begin{pmatrix} 1 &0\\ 0 & 1\\ \end{pmatrix}, 
-\begin{pmatrix} 1 &0\\ 0 & 1\\ \end{pmatrix}, 
\begin{pmatrix} 1 &0\\ 0 & -1\\ \end{pmatrix}, 
\begin{pmatrix} -1 &0\\ 0 & 1\\ \end{pmatrix} 
\Big\}. 
\end{align}
The most general orthogonal matrices that transform all  matrices in $MG_2$ or
$MG_4$ into signed permutations must have one of the following forms 
\begin{align}
\label{U2x2}
 U =  
\frac{PV_\mathrm{sd}}{\sqrt 2}  
\begin{pmatrix}
 1& 1\\
 1&-1\\
\end{pmatrix} ,
 \text{ or } 
 U =   \frac{PV_\mathrm{sd}}{\sqrt 2} 
 \begin{pmatrix}
-1& 1\\
 1& 1\\
\end{pmatrix} ,
 \text{ or } 
 U =  P V_\mathrm{sd} \begin{pmatrix}
 1& 0\\
 0& 1\\
\end{pmatrix}
\end{align}
where $V_\mathrm{sd}$ are signed diagonal matrices, and $P$ are permutation
matrices.  The number of the orthogonal transformations $U$ is finite.
\end{thm}
\begin{proof}[Proof of Theorem \ref{MG2x2}] 
We only need to consider the first matrix group $MG_2$, where the matrix group
is isomorphic to the $\mathbb{Z}_2$ group.  There are only four matrix groups
formed by 2-dimensional signed permutations matrices, that are  isomorphic
$\mathbb{Z}_2$.  The four matrix groups are generated by the following four
generators respectively:
\begin{align}
 \begin{pmatrix}
 1& 0\\
 0&-1\\
\end{pmatrix},\ \
 \begin{pmatrix}
-1& 0\\
 0&1\\
\end{pmatrix},\ \
 \begin{pmatrix}
0&1\\
1&0\\
\end{pmatrix},\ \
 \begin{pmatrix}
0&-1\\
-1&0\\
\end{pmatrix}.
\end{align}
An orthogonal transformation $U$ that transforms $MG$ to one of the above
matrix groups must have a from $U=VU_0$, where $V$ transforms $MG_2$ 
into itself,
and $U_0$ is a fixed orthogonal transformation that transforms $MG_2$ to one of
the above matrix groups.
We can choose $U_0$ to have the following form
\begin{align}
\label{U2x2b}
 U_0 =  
\frac{P}{\sqrt 2}  
\begin{pmatrix}
 1& 1\\
 1&-1\\
\end{pmatrix} ,
 \text{ or } 
 U_0 =   \frac{P}{\sqrt 2} 
 \begin{pmatrix}
-1& 1\\
 1& 1\\
\end{pmatrix} ,
 \text{ or } 
 U_0 =  P  \begin{pmatrix}
 1& 0\\
 0& 1\\
\end{pmatrix}.
\end{align}
To keep $MG$ unchanged
$V$ must satisfy
\begin{align}
 V  
\begin{pmatrix}
 1& 0\\
 0& -1\\
\end{pmatrix}
=
\begin{pmatrix}
 1& 0\\
 0& -1\\
\end{pmatrix}
V.
\end{align}
We find that $V$ must be diagonal. Thus $V$, as an orthogonal matrix, must be
signed diagonal. This gives us the result \eqref{U2x2}.
\end{proof}

If $\dim(\rho_\mathrm{isum}) \geq 8$, it is possible that the matrix group of
$D_{ \rho_\mathrm{isum}}(\si)_{I_\th,I_\th}$'s is generated by the following
non-diagonal matrix 
\begin{align}
\pm \begin{pmatrix}
0 & -1\\
1 &  0\\
\end{pmatrix}
\end{align}
This is because the direct sum decomposition of $\rho_\mathrm{isum}$ contains a
dimension-6 irreducible representation, whose $\rho(\frt)$ has a 2-dimensional
eigenspace.  The representation can give rise to this form of $D_{
\rho_\mathrm{isum}}(\si)_{I_\th,I_\th}$'s.

The eigenvalues of the matrices are $(\ii,-\ii)$.  
The most general  orthogonal matrices that transform all non-zero $D_{
\rho_\mathrm{isum}}(\si)_{I_\th,I_\th}$'s into signed
permutations must have the form 
\begin{align}
 U =  P V_\mathrm{sd} \begin{pmatrix}
 1& 0\\
 0& 1\\
\end{pmatrix}
.
\end{align}

If $\dim(\rho_\mathrm{isum}) \geq 8$, it is also possible that non-zero $D_{
\rho_\mathrm{isum}}(\si)_{I_\th,I_\th}$'s form the following matrix group:
\begin{align}
 \begin{pmatrix}
1& 0 \\ 
0 & 1 \\ 
\end{pmatrix}
,\ \ 
 \begin{pmatrix}
-\frac{1}{2}& -{\frac{\sqrt3}{2}} \\ 
{\frac{\sqrt3}{2}}& -\frac{1}{2} \\ 
\end{pmatrix}
,\ \ 
 \begin{pmatrix}
-\frac{1}{2}& {\frac{\sqrt3}{2}} \\ 
-{\frac{\sqrt3}{2}}& -\frac{1}{2} \\
\end{pmatrix}
\end{align}
This is because the direct sum decomposition of $\rho_\mathrm{isum}$ contains a
dimension-8 irreducible representation whose $\rho(\frt)$ has a 2-dimensional
eigenspace, which gives rise to the this form of $D_{
\rho_\mathrm{isum}}(\si)_{I_\th,I_\th}$'s.  In the above, the eigenvalues of
the later two matrices are $\pm (\ee^{\ii 2\pi/3},\ee^{-\ii 2\pi/3})$.  Since a
permutation of two elements can only have orders 1 or 2, the corresponding $2
\times 2$ signed permutation matrix  can only have eigenvalues $1$, $-1$ or
$\pm \ii$. Any other eigenvalue is not possible.  Thus, there is no orthogonal
matrix that can transform the above two matrices into signed permutations.
Such $\rho_\mathrm{isum}$ is not a representation of any modular data.

We also find cases where a 2-by-2 $D_{\rho_\mathrm{isum}}(\si)_{I_\th,I_\th}$
takes one of the following forms (up to conjugation of signed diagonal
matrices): 
$\begin{pmatrix}
-\frac{1+\sqrt{5}}{4},
& \frac{1}{2}c_{20}^{3} \\ 
-\frac{1}{2}c_{20}^{3},
& -\frac{1+\sqrt{5}}{4} \\ 
\end{pmatrix}
$,  
$\begin{pmatrix}
-\frac{1-\sqrt{5}}{4},
& -\frac{1}{2}c_{20}^{1} \\ 
\frac{1}{2}c_{20}^{1},
& -\frac{1-\sqrt{5}}{4} \\ 
\end{pmatrix}
$,  
$\begin{pmatrix}
\frac{1}{2},
& \frac{\sqrt{3}}{2} \\ 
-\frac{\sqrt{3}}{2},
& \frac{1}{2} \\ 
\end{pmatrix}
$.  
We find $D_{\rho_\mathrm{isum}}^4(\si)_{I_\th,I_\th} \neq \id$ for those
matrices.  Thus, those $D_{\rho_\mathrm{isum}}(\si)_{I_\th,I_\th}$ are not
similar to any 2-by-2 signed permutation.

\subsection{Within a 3-dimensional eigenspace of $\rho_\mathrm{isum}(\frt)$ }

We find cases where a 3-by-3 $D_{\rho_\mathrm{isum}}(\si)_{I_\th,I_\th}$
is given by
$\begin{pmatrix}
\frac{1}{2},
& -\frac{\sqrt{3}}{2},
& 0 \\ 
\frac{\sqrt{3}}{2},
& \frac{1}{2},
& 0 \\ 
0,
& 0,
& 1 \\ 
\end{pmatrix}
$, whose order is 3 (\ie cube to identity).  
Since the trace of the matrix is non-zero,
it cannot be similar to any 3-by-3 order-3 signed permutation matrix 
(since 3-by-3 order-3 signed permutation matrices all have zero trace).

There are also cases where 3-by-3 generators of 
$D_{\rho_\mathrm{isum}}(\si)_{I_\th,I_\th}$'s are given by
$\begin{pmatrix}
1,
& 0,
& 0 \\ 
0,
& -1,
& 0 \\ 
0,
& 0,
& -1 \\ 
\end{pmatrix}
$,  
$\begin{pmatrix}
1,
& 0,
& 0 \\ 
0,
& -\frac{1}{2},
& \frac{\sqrt{3}}{2} \\ 
0,
& -\frac{\sqrt{3}}{2},
& -\frac{1}{2} \\ 
\end{pmatrix}
$.  
The second matrix is of order-3 and can only be similar to signed $(1,2,3)$
permutations.  The most general orthogonal transformations that transform the
second matrix into signed permutations  have a form
\begin{align}
\label{W3by3}
U = P V_\textrm{sd} \begin{pmatrix}
\frac{\sqrt{3}}{3},
& -\frac{\sqrt{6}}{6},
& -\frac{\sqrt{2}}{2} \\ 
\frac{\sqrt{3}}{3},
& \frac{\sqrt{6}}{3},
& 0 \\ 
\frac{\sqrt{3}}{3},
& -\frac{\sqrt{6}}{6},
& \frac{\sqrt{2}}{2} \\ 
\end{pmatrix}
\end{align}
where $V_\textrm{sd}$ is a signed diagonal matrix and $P$ a permutation matrix.
But those orthogonal transformations all fail to
transform the first matrix into signed permutations.
So those cases are rejected.

There are cases where a 3-by-3  
$D_{\rho_\mathrm{isum}}(\si)_{I_\th,I_\th}$'s is 
$\begin{pmatrix}
1,
& 0,
& 0 \\ 
0,
& 0,
& -1 \\ 
0,
& 1,
& 0 \\ 
\end{pmatrix}
$.  
This matrix is of order-4 and can only be similar to signed $(2,3)$
permutations.
The most general orthogonal matrices that conjugate this matrix 
into signed permutation matrices have a form $P V_\mathrm{sd} $.

There are also cases where a 3-by-3  
$D_{\rho_\mathrm{isum}}(\si)_{I_\th,I_\th}$'s is (up to sign and signed permutations)
$\begin{pmatrix}
1,
& 0,
& 0 \\ 
0,
& -\frac{1}{2},
& \frac{\sqrt{3}}{2} \\ 
0,
& -\frac{\sqrt{3}}{2},
& -\frac{1}{2} \\ 
\end{pmatrix}
$.  
The most general orthogonal transformations
that transform this matrix into signed permutations are given by
\eqref{W3by3}.

\subsection{Within a 4-dimensional or 5-dimensional eigenspace of
$\rho_\mathrm{isum}(\frt)$ }

Some cases have 4-by-4  $D_{\rho_\mathrm{isum}}(\si)_{I_\th,I_\th}$'s generated
by the following generators (up to sign and signed permutations)
\begin{align}
\begin{pmatrix}
1,
& 0,
& 0,
& 0 \\ 
0,
& -1,
& 0,
& 0 \\ 
0,
& 0,
& -1,
& 0 \\ 
0,
& 0,
& 0,
& -1 \\ 
\end{pmatrix}
, \ \ \ \ \ \   
\begin{pmatrix}
1,
& 0,
& 0,
& 0 \\ 
0,
& 0,
& 0,
& -1 \\ 
0,
& 0,
& -1,
& 0 \\ 
0,
& 1,
& 0,
& 0 \\ 
\end{pmatrix}
\end{align}
The second matrix is a signed $(1,3)(2,4)$ permutation.
It can also be transformed into signed $(1,2,3,4)$ permutation by
the following matrix
\begin{align}
\label{W4by4a}
U = P V_\textrm{sd} \begin{pmatrix}
\frac{1}{2},
& \frac{\sqrt{2}}{2},
& \frac{1}{2},
& 0 \\ 
\frac{1}{2},
& 0,
& -\frac{1}{2},
& -\frac{\sqrt{2}}{2} \\ 
\frac{1}{2},
& -\frac{\sqrt{2}}{2},
& \frac{1}{2},
& 0 \\ 
\frac{1}{2},
& 0,
& -\frac{1}{2},
& \frac{\sqrt{2}}{2} \\ 
\end{pmatrix}
\end{align}
By this $U$ fails to transform the first matrix into a signed permutation.
Thus, the most general orthogonal transformations that transform the two
matrices into signed permutations have a form $U = P V_\textrm{sd}$.

Some cases have 4-by-4  $D_{\rho_\mathrm{isum}}(\si)_{I_\th,I_\th}$'s generated
by the following generators (up to sign and signed permutations)
\begin{align}
\begin{pmatrix}
1,
& 0,
& 0,
& 0 \\ 
0,
& 0,
& 0,
& -1 \\ 
0,
& 0,
& -1,
& 0 \\ 
0,
& 1,
& 0,
& 0 \\ 
\end{pmatrix}
\end{align}
The most general orthogonal transformations that transform this matrix
into signed permutations have a form $U = P V_\textrm{sd}$ or \eqref{W4by4a}.

Some cases have 5-by-5  $D_{\rho_\mathrm{isum}}(\si)_{I_\th,I_\th}$'s generated
by the following generators (up to sign and signed permutations)
\begin{align}
\begin{pmatrix}
-1, & 0, & 0, & 0, & 0 \\ 
0, & 1, & 0, & 0, & 0 \\ 
0, & 0, & 1, & 0, & 0 \\ 
0, & 0, & 0, & 1, & 0 \\ 
0, & 0, & 0, & 0, & 1 \\ 
\end{pmatrix}
,\ \ \ \ \ \ \
\begin{pmatrix}
1, & 0, & 0, & 0, & 0 \\ 
0, & 1, & 0, & 0, & 0 \\ 
0, & 0, & 0, & 0, & -1 \\ 
0, & 0, & 0, & -1, & 0 \\ 
0, & 0, & 1, & 0, & 0 \\ 
\end{pmatrix}
\end{align}
The most general orthogonal transformations that transform this matrix
into signed permutations have a form $U = P V_\textrm{sd}$ or 
\begin{align}
\label{W5by5}
U = P V_\textrm{sd} \begin{pmatrix}
1 & 0 & 0 & 0 & 0 \\
0 & \frac{1}{2}, & \frac{\sqrt{2}}{2}, & \frac{1}{2}, & 0 \\ 
0 & \frac{1}{2}, & 0, & -\frac{1}{2}, & -\frac{\sqrt{2}}{2} \\ 
0 & \frac{1}{2}, & -\frac{\sqrt{2}}{2}, & \frac{1}{2}, & 0 \\ 
0 & \frac{1}{2}, & 0, & -\frac{1}{2}, & \frac{\sqrt{2}}{2} \\ 
\end{pmatrix}
\end{align}

Some cases have 4-by-4  $D_{\rho_\mathrm{isum}}(\si)_{I_\th,I_\th}$'s generated
by the following generators (up to sign and signed permutations)
\begin{align}
\begin{pmatrix}
-1, & 0, & 0, & 0 \\ 
0, & 1, & 0, & 0 \\ 
0, & 0, & 1, & 0 \\ 
0, & 0, & 0, & 1 \\ 
\end{pmatrix}
,\ \ \ \ \ \
\begin{pmatrix}
1, & 0, & 0, & 0 \\ 
0, & 1, & 0, & 0 \\ 
0, & 0, & -\frac{1}{2}, & \frac{\sqrt{3}}{2} \\ 
0, & 0, & -\frac{\sqrt{3}}{2}, & -\frac{1}{2} \\ 
\end{pmatrix}
\end{align}
The second matrix is of order-6 and can only be similar to signed $(2,3,4)$
permutations.  The most general orthogonal transformations that transform the
second matrix (and the first matrix) into signed permutations  have a form
\begin{align}
\label{W4by4b}
U = P V_\textrm{sd} \begin{pmatrix}
1 & 0 & 0 & 0\\
0 &
\frac{\sqrt{3}}{3},
& -\frac{\sqrt{6}}{6},
& -\frac{\sqrt{2}}{2} \\ 
0 &
\frac{\sqrt{3}}{3},
& \frac{\sqrt{6}}{3},
& 0 \\ 
0 &
\frac{\sqrt{3}}{3},
& -\frac{\sqrt{6}}{6},
& \frac{\sqrt{2}}{2} \\ 
\end{pmatrix}
\end{align}

\subsection{General degenerate cases}
 
In general, the orthogonal matrix $U$ that transform $\rho_\textrm{isum}$ to a
pMD representation $\rho_\textrm{pMD}$ 
\begin{align}
 \rho_\textrm{pMD}(\frt) = U \rho_\textrm{isum}(\frt) U^\top,\ \ \ \
 \rho_\textrm{pMD}(\frs) = U \rho_\textrm{isum}(\frs) U^\top,
\end{align}
contains continuous parameters.  For example, $U$ may have the following form
\begin{align}
\label{Ugen}
U=
\begin{pmatrix}
\sqrt2 & -\sqrt2 & 0 & 0 & 0& 0& 0 &\cdots \\
\sqrt2 & \sqrt2 & 0 & 0 & 0& 0& 0 &\cdots \\
0 & 0 & u_1 & u_2 & 0& 0& 0 &\cdots \\
0 & 0 & u_2 &-u_1 & 0 & 0& 0 &\cdots \\
0 & 0 & 0 &  0 &  u_3 & u_4 & u_5  &\cdots \\
0 & 0 & 0 &  0 & u_6 & u_7 & u_8 &\cdots \\
0 & 0 & 0 &  0 & u_9& u_{10}& u_{11} &\cdots \\
 \vdots &  \vdots & \vdots &  \vdots & \vdots& \vdots& \vdots &\ddots \\
\end{pmatrix}
\ \
\text{ for }
\tilde\rho(\frt)=
\begin{pmatrix}
\t \th_0 & 0 & 0 & 0 &  0& 0& 0 &\cdots \\
0 & \t \th_0 & 0 &  0& 0& 0 &\cdots \\
0 & 0 & \t \th_1 & 0 & 0& 0& 0 &\cdots \\
0 & 0 & 0 & \t \th_1 & 0& 0& 0 &\cdots \\
0 & 0 & 0 &  0 & \t\th_2& 0& 0 &\cdots \\
0 & 0 & 0 &  0 & 0& \t\th_2& 0 &\cdots \\
0 & 0 & 0 &  0 & 0& 0& \t\th_2 &\cdots \\
 \vdots &  \vdots & \vdots &  \vdots & \vdots& \vdots& \vdots &\ddots \\
\end{pmatrix},
\end{align}
where $u_i$'s satisfy the following
 orthogonality conditions:
\begin{align}
 u_1^2+u_2^2 -1 =0,\ \ 
 u_3^2+u_4^2+u_5^2 -1 =0,\ \  
 u_3u_6+u_4u_7+u_5u_8  =0,\ \ \cdots 
\end{align}
where will be called zero conditions.  When $D_{\t\rho}(\si)$ are not signed
diagonal matrices with non-degenerate common eigenspaces, we will assume $U$ to
have the general form \eqref{Ugen}.

Because $u_i$ in $U$ are real numbers, different choices of $u_i$'s will give us
infinitely many potential pMD representations $\rho_\textrm{pMD}$, and hence
infinitely many potential MD representations $\rho_\textrm{MD}$, after a finite
number conjugations by signed diagonal matrices.  Thus we cannot check those MD
representations one by one to see which of them satisfy Proposition \ref{p:MD}.
Thus we need to use additional conditions on
$u_i$.  If we have enough conditions which only allow a finite numbers of
solutions, then we get a finite number of $U$'s which lead to a finite number
potential pMD representations.

We will use the following conditions on pMD representations to obtain equations
on $u_i$'s.  Those conditions on pMD representations are derived from condition
on MD representations (see Proposition \ref{p:MDcond}), by noticing that a pMD
representation is related to a MD representation via
\begin{align}
 (\rho_\mathrm{pMD})_{ij} = v_iv_j(\rho_\mathrm{MD})_{ij},\ \ \ v_i \in \{+1,-1\}
\end{align}
\begin{prop}
\label{pMDcond}
A pMD representation $\rho_\mathrm{pMD}$ has the following properties:
\begin{enumerate}

\item The conductor of the elements of $\rho_\mathrm{pMD}(\frs)$ divides
$\ord(\rho_\mathrm{pMD}(\frt))$.  


%



\item  The number of self dual objects is greater than 0. Thus
\begin{align}
 \Tr(\rho_\mathrm{pMD}^2(\frs)) \neq 0 .
\end{align}
Since $\Tr(\rho_\mathrm{pMD}^2(\frs)) \neq 0$, let us introduce
\begin{align}
\label{Cdef}
 C = \frac{\Tr(\rho_\mathrm{pMD}^2(\frs))}{|\Tr(\rho_\mathrm{pMD}^2(\frs))|}   \rho_\mathrm{pMD}^2(\frs).
\end{align}
The above $C$ is the charge conjugation operator of MTC, {\it i.e.} $C$ is a
permutation matrix of order 2.  In particular, $\Tr(C)$ is the number of  self
dual objects. Also, for each eigenvalue $\tilde \theta$ of
$\rho_\mathrm{pMD}(\frt)$, 
\begin{align}
\Tr_{\tilde \theta}(C) \geq 0,
\end{align}
where $\Tr_{\tilde \theta}$ is the trace in the degenerate subspace of
$\rho_\mathrm{pMD}(\frt)$ with eigenvalue $\tilde \theta$.

\item 
If the modular data is integral  and $\ord(\rho_\mathrm{pMD}(\frt))$ = odd, then
\begin{align}
\Tr(C) =\Big( \Tr\big(\rho^2_\mathrm{pMD}(\frs)\big) \Big)^2 = 1,
\end{align}
\ie the unit object is the only self-dual object.

\item 
For any Galois conjugation $\si$ in $\Gal(\Q_{\ord(\rho_\mathrm{pMD}(\frt))})$, there is a
permutation of the indices, $i \to \hs(i)$, and $\eps_\si(i)\in \{1,-1\}$, such
that
\begin{align}
\label{pMDGalact}
\si \big(\rho_\mathrm{pMD}(\frs)_{i,j}\big) &
= \eps_\si(i)\rho_\mathrm{pMD}(\frs)_{\hat \si (i),j} 
= \rho_\mathrm{pMD}(\frs)_{i,\hat \si (j)}\eps_\si(j) 
\\
\si^2 \big(\rho_\mathrm{pMD}(\frt)_{i,i}\big) &= \rho_\mathrm{pMD}(\frt)_{\hat \si (i),\hat \si (i)},
\end{align}
for all $i,j$. 

\item 
For any integer $a$ coprime to $n  = \ord(\rho_\mathrm{pMD}(\frt))$, we define
\begin{align}
\label{pMDdefDrho}
D_{\rho_\mathrm{pMD}}(a) &:= \rho_\mathrm{pMD}(\frt^a \frs \frt^b \frs \frt^a \frs^{-1}) 
= D_{\rho_\mathrm{pMD}}(a+\ord(\rho_\mathrm{pMD}(\frt))) , 
\nonumber\\
&
\text{ where } ab \equiv 1 \mod \ord(\rho_\mathrm{pMD}(\frt))\,. 
\end{align}
For any $\si \in \mathrm{Gal}(\Q_n)$, $\si(\zeta_n) = \zeta_n^a$ for some unique integer $a$ modulo $n$. We define
\begin{align}
\label{pMDeq:Drho}
D_{\rho_\mathrm{pMD}}(\si) := D_{\rho_\mathrm{pMD}}(a)\,. 
\end{align}
By \cite[Theorem II]{DLN}, $D_{\rho_\mathrm{pMD}}:
\mathrm{Gal}(\Q_n)=(\Z_n)^\times \to \GL_r(\C)$ is a representation equivalent
to the restriction of $\rho_\mathrm{pMD}$ on the diagonal subgroup of
$\qsl{n}$.  $D_{\rho_\mathrm{pMD}}(\si)$ in \eqref{pMDdefDrho} must be a
\textbf{signed permutation}
\begin{align}
\label{pMDDrhosihat}
 (D_{\rho_\mathrm{pMD}}(\si))_{i,j} = \eps_\si(i) \delta_{\hs(i),j}.
\end{align}
and satisfies
\begin{align}
\label{pMDsiDrho}
\si(\rho_\mathrm{pMD}(\frs)) &= D_{\rho_\mathrm{pMD}}(\si) \rho_\mathrm{pMD}(\frs) =\rho_\mathrm{pMD}(\frs)D_{\rho_\mathrm{pMD}}^\top(\si),
\nonumber\\
\si^2(\rho_\mathrm{pMD}(\frt)) &= D_{\rho_\mathrm{pMD}}(\si) \rho_\mathrm{pMD}(\frt) D_{\rho_\mathrm{pMD}}^\top(\si)
\end{align}

\item 
There exists a $u$ such that $\rho_\mathrm{pMD}(\frs)_{uu} \neq 0$ and
\begin{align}
\label{pMDSSN}
&  
 \frac{ \rho_\mathrm{pMD}(\frs)_{ij} }{ \rho_\mathrm{pMD}(\frs)_{uu} },\ 
 \frac{ \rho_\mathrm{pMD}(\frs)_{ij} }{ \rho_\mathrm{pMD}(\frs)_{uj} }
\in \BO_{\ord(T)}, \ \ \ \
 \frac{ \rho_\mathrm{pMD}(\frs)_{ij} }{ \rho_\mathrm{pMD}(\frs)_{i'j'} }
\in \Q_{\ord(T)},
\nonumber\\
&
\rho_\mathrm{pMD}(\frs)_{ui} \neq 0 , \ \ \ \
 \frac{ 1}{ \rho_\mathrm{pMD}(\frs)_{ui} }
\in \BO_{\ord( \rho(\frt))}, 
, 
\nonumber\\
& \t N^{ij}_k = \sum_{l=0}^{r-1} \frac{
\rho_\mathrm{pMD}(\frs)_{li} \rho_\mathrm{pMD}(\frs)_{lj} \rho_\mathrm{pMD}(\frs^{-1})_{lk}}{ \rho_\mathrm{pMD}(\frs)_{lu} } \in\Z
,\ \  \forall \ i,j,k = 1,2,\ldots,r.
\end{align}
($u$ corresponds the unit object of MTC. Also see Lemma \ref{Ddi}.)

\item
Let $n \in \N_+$.  The $n^\text{th}$ pseudo Frobenius-Schur indicator of the
$i$-th simple object
\begin{align}
 \label{pMDnunFSapp}
 \t\nu_n(i)&= 
\sum_{l=1}^{r} \frac{
\rho_\mathrm{pMD}(\frs\frt^n\frs)_{lu} \rho_\mathrm{pMD}(\frs\frt^{-n}\frs^{-1})_{lu} \rho_\mathrm{pMD}(\frs^{-1})_{li}}{ \rho_\mathrm{pMD}(\frs)_{lu} }
\end{align}
is a cyclotomic integer whose conductor divides $n$ and $\ord(T)$.  The 1st
pseudo Frobenius-Schur indicator satisfies $\t\nu_1(i)=\delta_{iu}$ while the
2nd pseudo Frobenius-Schur indicator $\t\nu_2(i)$ satisfies $\t\nu_2(i) \in \{
\rho_\mathrm{pMD}(\frs^2)_{ii},
-\rho_\mathrm{pMD}(\frs^2)_{ii} 
\}$ (see \cite{Bantay, NS07b, RSW0777}).  We also
have the identity $\t \nu_n(u) =1$.  

\item
$\t D = 1/\rho_\mathrm{pMD}(\frs)_{uu}$ is a cyclotomic integer.  $\t
D^5/\ord(T)$ is an algebraic  integer, which is also a cyclotomic integer.  The
prime divisors of norm$(\t D)$ and $\ord(T)$ coincide.  $\t D/\t d_i =
1/\rho_\mathrm{pMD}(\frs)_{ui}$ are cyclotomic integers (see Lemma \ref{Ddi}),
where $\t d_i = \rho_\mathrm{pMD}(\frs)_{ui}/\rho_\mathrm{pMD}(\frs)_{uu}$.
The prime divisors of norm$(\t D/\t d_i)  =
\textrm{norm}(1/\rho_\mathrm{pMD}(\frs)_{ui})$ are part of the prime divisors
of $\ord(T)$. This also implies that the prime divisors of norm$(\t d_i)  =
\textrm{norm}(\rho_\mathrm{pMD}(\frs)_{ui}/\rho_\mathrm{pMD}(\frs)_{uu})$ are
part of the prime divisors of $\ord(T)$.

\item 
If the above conditions are satisfied for choosing an index $u$ as the unit
index, then the above conditions are also satisfied for choosing another index
$u'$ as the unit index, provided that $u$ and $u'$ are related by an Galois
conjugation: $u'=\hat \sigma (u)$

\end{enumerate}
\end{prop}

\begin{lem}
\label{Ddi}
For a pMD representation $\rho_\mathrm{pMD}$,  let $s=\rho_\mathrm{pMD}(\frs)$,
the $u$-th row of $s$ the unit row, $n =\ord(\rho_\mathrm{pMD}(\frt))$, and $N
=\pord(\rho_\mathrm{pMD}(\frt))$.  Then $1/s_{uj}$ is a cyclotomic integer in
$\Q_n$ which divides $1/s_{uu}$ for  all $j$. In particular, the prime factors
of $\mathrm{norm}(1/s_{uj})$ can only be part of those of $N$.
\end{lem}
\begin{proof}
Let $d_i$ and $D^2$ be the quantum dimensions and the total quantum dimension
of the modular tensor category that corresponds to a pMD representation
$\rho_\mathrm{pMD}$.  Note that $\t D = 1/s_{uu} \in
\Q_n$  and $\t d_j =\frac{s_{uj}}{s_{uu}} \in
\Q_N$, which are cyclotomic integers for all
$j$.  $D$ and $\t D$ differ only by some 4-th roots of unity.  $d_i$ and $\t
d_i$ differ only by $\pm 1$.  Since $D/d_j$ is an algebraic integer and is
equal to $1/s_{uj}$ up to a 4-th root of unity,  $1/s_{uj}$ is an algebraic
integer in $\Q_n$ which divides  $1/s_{uu}$ as
algebraic integers. By the Cauchy theorem, $\mathrm{norm}(1/s_{uu})$ has
exactly the same prime factors as those of $N$.
Since $\mathrm{norm}(1/s_{uj}) \mid \mathrm{norm}(1/s_{uu})$, the prime factors
of $\mathrm{norm}(1/s_{uj})$ must be a subcollection of the prime factors of
$N$.
\end{proof}

In our computer algebra calculation, the conditions on $u_i$ are written as
\emph{zero conditions} $f_1(u_i)=0 \text{ and }  f_2(u_i)=0\text{ and }
\cdots$, where $f_k(u_i)$ are multi-variable polynomials. 
Since those zero conditions
must be satisfied simultaneously, we refer to the set of zero conditions as
and-connected zero conditions.

We may have another set zero conditions that must be satisfied simultaneously.
But we only require one of the two sets of zero conditions to be satisfied.
Thus the two sets of zero conditions are connected by ``or'', which give rise
to or-connected sets of and-connected zero conditions.  Some time, we create
and-connected sets of or-connected zero conditions, and we need to convert them
to or-connected sets of and-connected zero conditions.  Our computer code is
designed to manage those different logically connected zero conditions.

When we have two zero or-connected conditions, $g_1(u_i)=0$ or $h_1(u_i)=0$, we
could combine them into one zero condition $f_1(u_i) = g_1(u_i)h_1(u_i) = 0$.
But we will not do it. We will store the zero conditions
as or-connected sets of and-connected zero conditions:
\begin{align}
& [g_1(u_i)=0 \text{ and } f_2(u_i)=0\text{ and } \cdots] 
\nonumber\\
\text{or  } \ \ 
& [h_1(u_i)=0
\text{ and }  f_2(u_i)=0 \text{ and } \cdots] .
\end{align}
Those two sets of zero conditions can be viewed as a factorization of a single
set of zero conditions $f_1(u_i)=0 \text{ and }  f_2(u_i)=0 \text{ and }
\cdots$.  So storing the $u_i$'s conditions as two sets of and-connected zero
conditions allows us to avoid factorizing some algebraic equations (such as
factor $f_1(u_i)=0$ into $g_1(u_i)=0$ or $h_1(u_i)=0$), which is difficult and
unreliable for multi-variable polynomials.  This is a strategy that we use in
computer algebra calculation to construct the conditions on $u_i$'s: \emph{We
try to construct many or-connected sets of and-connected zero conditions.}

We call this stage of calculation as $u$-stage, where we try to factorize a
zero condition $f_1(u_i) = g_1(u_i)h_1(u_i) = 0$ into or-connected zero
conditions: $g_1(u_i)=0$ or $h_1(u_i)=0$.  To help factorizing, we write GAP
code to choose Gr\"obner basis for the and-connected zero conditions so that
some of the  zero conditions to have as few variables as possible.  When a zero
condition contains a single $u_i$ variable and if this variable is known to be
a cyclotonic number of a know conductor, then we can factorize the zero
condition.

If after factorization, we obtain say
\begin{align}
& [u_1=0 \text{ and } u_2- 1=0 \text{ and } u_3^2+u_4^2 - 1 = 0
\cdots] 
\nonumber\\
\text{or } \ \ & [u_1-\frac{\sqrt2}2 =0 \text{ and } u_2-\frac{\sqrt2}2 = 0 \text{ and }
u_3^2+u_4^2 - 1 =0 \cdots], 
\end{align}
then $u_1$ and $u_2$ are solved, but $u_3$ and $u_4$ remain unsolved.  Since
$u_i$ are real, the condition $u_1^2+u_2^2=0$ becomes $u_1=0 \text{ and }
u_2=0$.  The condition $u_1^2+u_2^2+1=0$ leads to a rejection.

If we fail to solve all the variables $u_i$, we replace the matrix elements of
$\rho_\mathrm{pMD}(\frs)$ with $u_i$ variable's by cyclotomic numbers which
contain only finite rational variables $r_i$.  This is because the conductor of
the matrix elements of $\rho_\mathrm{pMD}(\frs)$ is
$\ord(\rho_\mathrm{pMD}(\frt))$ which is finite.  At this stage of rational
variables (called $r$-stage), we can utilize addition conditions that involve
Galois conjugations, which we cannot use at the $u$-stage.  Also, the zero
condition like $r_1^2 -\frac23 =$ will lead to a rejection at $r$-stage.

If we fail to solve all the variable $r_i$, we can go to stage of integer
variables $n_i$, the $n$-stage, (or we can go to the $n$-stage directly from the
$u$-stage).  This is achieved by representing $\rho_\mathrm{pMD}(\frs)$ in
terms of cyclotomic integers:
\begin{align}
\label{KKJ}
&  
K_{ij} \equiv \frac{ \rho_\mathrm{pMD}(\frs)_{ij} }{ \rho_\mathrm{pMD}(\frs)_{uu} }
\in \BO_{\ord(T)},\ \ \ 
\t K_{ij} \equiv \frac{ \rho_\mathrm{pMD}(\frs)_{ij} }{ \rho_\mathrm{pMD}(\frs)_{uj} }
\in \BO_{\ord(T)},\ \ \ 
J_i \equiv \frac{ 1}{ \rho_\mathrm{pMD}(\frs)_{ui} }
\in \BO_{\ord( \rho(\frt))}.
\end{align}
From $K_{ij},\ \t K_{ij},\ J_i$ we can recover $\rho_\mathrm{pMD}(\frs)$.  
Since $K_{ij},\ \t K_{ij},\ J_i$ have finite conductors
$\ord(\rho_\mathrm{pMD}(\frt))$ or $\ord(T)$, we can replace the elements with
$r_i$ or $u_i$ variables by a finite number of integer variables $n_i$.

\begin{prop}
\label{KKQcond}
The two matrices $K,\ \t K$ and the vector $J$ satisfy the following
conditions, which we can use to solve for  integer variables $n_i$:
\begin{enumerate}
\item
From the their definition, we find
\begin{align}
 K_{ij} = K_{ji}, \ \ \ \ K_{iu} = \text{real}, \ \ \ \
 K_{ij} J_{j}  = \t K_{ij} J_u .
\end{align}

\item
For any integer $\si$ coprime to $\ord(\rho_\mathrm{pMD}(\frt))$, we have
\begin{align}
D_{\rho_\mathrm{pMD}}(\si) K J_u &= 
\rho_\mathrm{pMD}(\frt)^\si K \rho_\mathrm{pMD}(\frt)^\ga K \rho_\mathrm{pMD}(\frt)^\si 
\end{align}
where  $\ga$ satisfies $\si \ga  = 1 \mod \ord(\rho_\mathrm{pMD}(\frt))$. 
\begin{align}
 D_{\rho_\mathrm{pMD}}(-1) J_u^2 = K^2 ,
\end{align}
where $D_{\rho_\mathrm{pMD}}(\si)=
D_{\rho_\mathrm{pMD}}(\si+\ord(\rho_\mathrm{pMD}(\frt))) $ and are signed
permutation matrics, which for a representation of the multiplication group
$\Z_{\ord(\rho_\mathrm{pMD}(\frt))}^\times$.

\item The Galois conguations:
\begin{align}
 \si(J_i) &= \sum_{j=0}^{r-1} J_j (D_{\rho_\mathrm{pMD}}(\si)^\top)_{ji}, 
&
 \si(\t K_{ij}) &= \sum_{k=0}^{r-1} \t K_{ik} |(D_{\rho_\mathrm{pMD}}(\si)^\top)_{kj}|
,
\nonumber\\
\si(K_{ij})J_u &= \sum_{k=0}^{r-1}  D_{\rho_\mathrm{pMD}}(\si)_{ik} K_{kj} \si(J_u),
&
\si(\t K_{ij})J_j &= \sum_{k=0}^{r-1}  D_{\rho_\mathrm{pMD}}(\si)_{ik} \t K_{kj} \si(J_j)
.
\end{align}

\item
If we know $|K_{iu}|=1$ for an index $i$, then
\begin{align}
 |K_{ij}|=1 \ \ \text{ for all indices } j .
\end{align}

\item Pseudo characters of $\SL$ representation:
\begin{align}
 J_u \Tr(\rho_\mathrm{isum}(\frs)_{I_\th, I_\th}) &= \Tr(K_{I_\th,I_\th}), \ \ \ \
\nonumber\\
 J_u^2 \Tr\big(
\rho_\mathrm{isum}(\frs)_{I_\th,I_{\th'}}
\rho_\mathrm{isum}(\frs)_{I_{\th'},I_\th} 
\big) 
&=
 \Tr\big(
K_{I_\th,I_{\th'}}
K_{I_{\th'},I_\th} 
\big) 
.
\end{align}
Here, $I_{\th}$ is the set of indices for the degenerate eigenspace of
$\rho_\mathrm{pMD}(\frt)$ with eigenvalue $\th$.  For example, the matrix
$K_{I_\th,I_{\th'}}$ is a block of the $K$-matrix that connects $I_{\th}$ and
$I_{\th'}$.

\item Orthogonality conditions:
\begin{align}
\sum_{i=0}^{r-1}  \frac{1}{|J_i|^2} = 1,\ \ 
\sum_{k=0}^{r-1} 
K_{ik} K_{kj}^* = |J_u|^2 \del_{ij}, \ \ 
\sum_{k=0}^{r-1} \t K_{ki} \t K_{kj}^* = |J_i|^2 \del_{ij}, \ \ 
\sum_{k=0}^{r-1} K_{ik} \t K_{kj}^* = J_u J_i^* \del_{ij}
\end{align}

\item
The pseudo fusion coefficients $\t N^{ij}_k$:
\begin{align}
\t N^{ij}_k &= 
\rho_\mathrm{pMD}(\frs)_{uu}
\rho_\mathrm{pMD}(\frs)_{uu}^*
\sum_{l=0}^{r-1} \frac{
\rho_\mathrm{pMD}(\frs)_{li} }
{ \rho_\mathrm{pMD}(\frs)_{uu} } 
\frac{
\rho_\mathrm{pMD}(\frs)_{lj} 
}{\rho_\mathrm{pMD}(\frs)_{uu}}
\frac{
\rho_\mathrm{pMD}(\frs)_{lk}^*}
{\rho_\mathrm{pMD}(\frs)_{lu}^*}
\nonumber\\
& =
\frac{
\sum_{l=0}^{r-1}
K_{il} K_{jl} \t K_{kl}^*
}
{|J_u|^2}
\in\Z,
\end{align}
where we have used the fact that
$\frac{\rho_\mathrm{pMD}(\frs)_{uu}}{\rho_\mathrm{pMD}(\frs)_{ul}}$ are real,
which leads to $
\frac{\rho_\mathrm{pMD}(\frs)_{uu}}{\rho_\mathrm{pMD}(\frs)_{ul}} =
\frac{\rho_\mathrm{pMD}(\frs)_{uu}^*}{\rho_\mathrm{pMD}(\frs)_{ul}^*} $.

\item
The pseudo Frobenius-Schur indicators $\t \nu_n(i)$:
\begin{align}
 \t\nu_n(i)&= 
|\rho_\mathrm{pMD}(\frs)_{uu}|^4
\sum_{l=1}^{r} 
\frac{ \rho_\mathrm{pMD}(\frs\frt^n\frs)_{lu} }
{(\rho_\mathrm{pMD}(\frs)_{uu})^2}
\frac{ \rho_\mathrm{pMD}(\frs\frt^{-n}\frs^{-1})_{lu} }
{ \rho_\mathrm{pMD}(\frs)_{uu} \rho_\mathrm{pMD}(\frs)_{uu}^* }
\frac{ \rho_\mathrm{pMD}(\frs)^*_{il}}
{ \rho_\mathrm{pMD}(\frs)^*_{ul} }
\nonumber\\
&= 
|J_u|^{-4}
\sum_{l=1}^{r} 
\big(K \rho_\mathrm{pMD}^n(\frt) K\big)_{lu} 
\big(K \rho_\mathrm{pMD}^{-n}(\frt) K^*\big)_{lu} 
\t K^*_{il},
\end{align}
where $\t\nu_1(i)=\delta_{iu}$ and $\t\nu_2(i) \in \{
\rho_\mathrm{pMD}(\frs^2)_{ii}, -\rho_\mathrm{pMD}(\frs^2)_{ii} \}$

\end{enumerate}
\end{prop}

At $n$-stage, there are more tricks to solve zero conditions. For example, a
zero condition $n_1 n_2 = 10$ will lead to or-connected zero conditions: 
\begin{align}
& [n_1 = 1 \text{ and } n_2=10] \ \ \text{ or }\ \ [n_1 = -1 \text{ and } n_2=-10]
\nonumber\\
\text{or } \ \ & [n_1 = 2 \text{ and } n_2=5] \ \ \ \text{ or }\ \ [n_1 = -2 \text{ and } n_2=-5]
\nonumber\\
\text{or } \ \ & [n_1 = 5 \text{ and } n_2=2] \ \ \ \text{ or }\ \ [n_1 = -5 \text{ and } n_2=-2] 
\nonumber\\
\text{or } \ \ & [n_1 = 10 \text{ and } n_2=1] \ \ \text{ or }\ \ [n_1 = -10 \text{ and } n_2=-1].  
\end{align}
Also, a zero condition $n_1^2 -n_1n_2+ n_2^2 = 1$ of elliptic form will lead to
or-connected zero conditions:
\begin{align}
[n_1 = 1 \text{ and } n_2=1] \ \ \text{ or } \ \  
[n_1 = -1 \text{ and } n_2=-1].
\end{align}

\subsection{Some details of computer calculations}

\subsubsection{Determining $D_{\rho_\mathrm{pMD}}(\si)$}
\label{calDrho}

$D_{\rho_\mathrm{pMD}}(\si)$ given by Eq. \eqref{pMDdefDrho} contain the
variables $u_i$.  The non-zero elements and the elements containing variables
form blocks: Let $I_{\th}$ be the indices for the degenerate eigenspace
of $\rho_\mathrm{pMD}(\frt)$ with eigenvalue $\th$.  Only elements in the block
$D_{\rho_\mathrm{pMD}}(\si)_{I_{\th},I_{\th'}}$ is non-zero, and for each fixed
$I_{\th}$, there exists only one $I_{\th'}$ such that
$D_{\rho_\mathrm{pMD}}(\si)_{I_{\th},I_{\th'}}$ is non-zero.  Also a non-zero
block is always a square matrix.  This is because the non-zero block
$D_{\rho_\mathrm{pMD}}(\si)_{I_{\th},I_{\th'}}$ connects eigenvalues $\th$ and
$\th' = \si^2(\th)$, and the number degenerate eigenvalues for $\th$ and
$\th'$, mapped into each other by a Galois conjugation, are the same.

We know that the matrix elements of $D_{\rho_\mathrm{pMD}}(\si)$ can only take
three possible values $0,\pm 1$, since $D_{\rho_\mathrm{pMD}}(\si)$ are signed
permutations.  Thus, there are only a finite number of possible
$D_{\rho_\mathrm{pMD}}(\si)$.  So in the first step, we list all those possible
$D_{\rho_\mathrm{pMD}}(\si)$'s.  In fact, the number of possible sets of
$D_{\rho_\mathrm{pMD}}(\si)$'s are not many, since many matrix elements of
$D_{\rho_\mathrm{pMD}}(\si)$'s are either equal or differ by sign.  When we
assign $0,\pm1$ to  matrix elements of $D_{\rho_\mathrm{pMD}}(\si)$, we include
such correlations.  Also, the resulting $D_{\rho_\mathrm{pMD}}(\si)$'s are
commuting signed permutations, which also reduces the number of possible sets
of $D_{\rho_\mathrm{pMD}}(\si)$'s.  The computation load can be further reduced
if we use the  commuting properties of $D_{\rho_\mathrm{pMD}}(\si)$ as early as
possible during our determination of the matrix elements of
$D_{\rho_\mathrm{pMD}}(\si)$.

The number of the possible $D_{\rho_\mathrm{pMD}}(\si)$'s can be further
reduced.  We note that the transformation $U$ is block diagonal.  The block
that contain $u_i$ variable's has a form $U_{I_{\th},I_{\th}}$, \ie $U$ only
maps the eigenspace $I_{\th}$ into itself (see \eqref{Ugen}).  The
transformation $U$ generated by our code has such a property that the block
$U_{I_{\th},I_{\th}}$ with variable's corresponds to the most general orthogonal
transformations in the eigenspace $I_{\th}$.  Those orthogonal transformations,
acting on various eigenspaces, generate signed permutations of the index
$i=0,1,\cdots,r-1$.  We can define two sets of $D_{\rho_\mathrm{pMD}}(\si)$'s as
equivalent if they are connected by those signed permutations by conjugation.
We only need to pick one representative from each equivalence class.

We pick the representative in the following way.  We first pick a $\si$.  Then
we only require $D_{\rho_\mathrm{pMD}}(\si)$ for that one $\si$ to satisfy some
additional conditions.  For a diagonal block
$D_{\rho_\mathrm{pMD}}(\si)_{I_{\th},I_{\th}}$ where $U_{I_{\th},I_{\th}}$
contain $u_i$ variable's, we require the square matrix $D \equiv
D_{\rho_\mathrm{pMD}}(\si)_{I_{\th},I_{\th}}$ to satisfy the following
conditions (\ie we can use signed permutations $P_\mathrm{sgn}$ and
transformation $D\to P_\mathrm{sgn} D P_\mathrm{sgn}^{-1}$ to make $D$ to
satisfy the following conditions):
\begin{enumerate}

\item 
$D$ is a signed permutation matrix.

\item
$D_{i,j} = 0$ if $j \geq i+2$. 

\item
$D_{i,i+1} \in \{0,1\}$, \ie the possibility of $D_{i,i+1}=-1$ is excluded. 

\item
$D_{ii} \geq D_{i+1,i+1}$.

\end{enumerate}
For example, $D$ may take the following form
\begin{align}
\label{Dmat}
 D =
\begin{pmatrix}
 1 & 0 & 0 & 0 & 0\\
 0 & 0 & 1 & 0 & 0\\
 0 & 0 & 0 & 1 & 0\\
 0 & \pm1 & 0 & 0 & 0\\
 0 & 0 & 0 & 0 & -1\\
\end{pmatrix}
\end{align}
In other words, we can use $U_{I_{\th},I_{\th}}$ to make the upper triangle
part of $D$ non-negative.

Let us define an ordering of the blocks $I_\th$ in such a way that
block-level cyclic permutation $( I_{\th_1}, I_{\th_2},\cdots, I_{\th_n},) $
generated by $D_{\rho_\mathrm{pMD}}(\si)_{I_{\th},I_{\th'}}$ satisfies
\begin{align}
I_{\th_1} < I_{\th_2}, \  
I_{\th_2} < I_{\th_3}, \ \cdots, \ 
I_{\th_n} > I_{\th_1}.  
\end{align}
For an variable-containing off-diagonal block
$D_{\rho_\mathrm{pMD}}(\si)_{I_{\th},I_{\th'}}$ with $I_{\th} > I_{\th'}$ where
$U_{I_{\th},I_{\th}}$ and $U_{I_{\th'},I_{\th'}}$ contain $u_i$ variable's, we
require the square matrix $D \equiv
D_{\rho_\mathrm{pMD}}(\si)_{I_{\th},I_{\th'}}$ to satisfy the following
conditions: (\ie we can use two signed permutations $P_\mathrm{sgn},\t
P_\mathrm{sgn}$ and transformation $D\to P_\mathrm{sgn} D \t P_\mathrm{sgn}$ to
make $D$ to satisfy the following conditions):
\begin{enumerate}

\item
$D$ is an identity matrix.

\end{enumerate}

For an variable-containing off-diagonal block
$D_{\rho_\mathrm{pMD}}(\si)_{I_{\th},I_{\th'}}$ with $I_{\th} < I_{\th'}$ where
$U_{I_{\th},I_{\th}}$ and $U_{I_{\th'},I_{\th'}}$ contain $u_i$ variable's, we
require the square matrix $D \equiv
D_{\rho_\mathrm{pMD}}(\si)_{I_{\th},I_{\th'}}$ to satisfy the following
conditions (\ie we can use signed permutations $P_\mathrm{sgn}$ and
transformation $D\to P_\mathrm{sgn} D P_\mathrm{sgn}^{-1}$ to make $D$ to
satisfy the following conditions):
\begin{enumerate}

\item 
$D$ is a signed permutation matrix.

\item
$D_{i,j} = 0$ if $j \geq i+2$. 

\item
$D_{i,i+1} \in \{0,1\}$, \ie the possibility of $D_{i,i+1}=-1$ is excluded. 

\item
$D_{ii} \geq D_{i+1,i+1}$.

\end{enumerate}
We like to remark that since we have fixed
$D_{\rho_\mathrm{pMD}}(\si)_{I_{\th},I_{\th'}}$ with $I_{\th} > I_{\th'}$ to be
identity, the transformation, $D_{\rho_\mathrm{pMD}}(\si)_{I_{\th},I_{\th'}}
\to P_\mathrm{sgn} D_{\rho_\mathrm{pMD}}(\si)_{I_{\th},I_{\th'}} \t
P_\mathrm{sgn}$, that keep $D_{\rho_\mathrm{pMD}}(\si)_{I_{\th},I_{\th'}}$
unchanged must satisfy $\t P_\mathrm{sgn} = P_\mathrm{sgn}^{-1}$.  Therefore,
the transformation on $D_{\rho_\mathrm{pMD}}(\si)_{I_{\th},I_{\th'}}$ with
$I_{\th} < I_{\th'}$ has a form $D_{\rho_\mathrm{pMD}}(\si)_{I_{\th},I_{\th'}}
\to P_\mathrm{sgn} D_{\rho_\mathrm{pMD}}(\si)_{I_{\th},I_{\th'}}
P_\mathrm{sgn}^{-1}$.
For example, $D$ may take the following form
\begin{align}
 D =
\begin{pmatrix}
 0 & \id & 0 \\
 0 & 0 & \id \\
 \tilde D & 0 & 0 \\
\end{pmatrix}
\end{align}
where $\tilde D$ is given by
\eqref{Dmat}.

After knowing $D_{\rho_\mathrm{pMD}}(\si)$, we can obtain a set of zero
conditions on $u_i$'s from Eq. \eqref{pMDsiDrho}, as well as the zero
conditions from the orthogonality conditions of $U$.

\subsubsection{Determining if the MTC is integral}

We expand the number possible cases (\ie the number of or-connected sets of
and-connected zero conditions) further by choosing possible rows of
$\rho_\mathrm{pMD}(\frs)$ to be the unit row.  Certainly, only the row that
contain no zero's can be a unit row.  

Then for each case, we can check if $(D_{\rho_\mathrm{pMD}}(\si))_{uu} = \pm 1$
for all $\si$. If yes, it means $\rho_\mathrm{pMD}$ must correspond to a
integral MTC, if any.  Also, if $\pord(\rho_\mathrm{pMD}(\frt)) \in
\{2,3,4,6\}$, then $\rho_\mathrm{pMD}$ must correspond to an integral MTC as
well, if any.  Using such a method, we find that $\SL$ representations may give
rise to integral MD's only when the prime divisors of
$\pord(\rho_\mathrm{pMD}(\frt))$ belong to the sets listed in Table
\ref{pfac}.  This result will use to construct integral MTCs (see Section
\ref{intMD} for more details).

\begin{table}[tb] 
\caption{
We have computed $\rho_\mathrm{pMD}(\frt)$ of all $\SL$ representations that may
give rise to non-pointed integral MD's for each rank.  This table lists all the
prime divisors of $\pord(\rho_\mathrm{pMD}(\frt))$ of those
$\rho_\mathrm{pMD}$'s.
} 
\label{pfac} 
\centering
\begin{tabular}{ |c|c| } 
\hline 
rank & prime divisors  \\
\hline 
2 & [\ ] \\
\hline 
3 & [\ ] \\
\hline 
4 & [2] \\
\hline 
5 & [2,3] \\
\hline 
6 & [\ ] \\
\hline 
7 & [\ ] \\
\hline 
8 & [2,3] \\
\hline 
9 & [2,3] \\
\hline 
10 & [2], [2,3], [2,3,5] \\
\hline 
11 & [2], [2,3], [2,3,5], [2,5] \\
\hline 
12 & [2,3], [2,3,5], [2,3,7] \\
\hline \end{tabular}
\end{table}

From now on, for each case, we not only know $\rho_\mathrm{pMD}(\frs)$ and
$\rho_\mathrm{pMD}(\frt)$, we also know a possible $D_{\rho_\mathrm{pMD}}(\si)$
of them, as well as a possible index $u$ for the unit row, together with a set
of and-connected zero conditions of $u_i$.  
This allows us to obtain many additional zero conditions on $u_i$'s for such a
case.

\subsubsection{Conditions from Frobenius-Schur indicators}

Since we know the index $u$ of the unit row, we can compute the second
pseudo Frobenius-Schur indicator from $\rho_\mathrm{pMD}(\frs)$:
\begin{align}
 \t\nu_2(i)&= 
\sum_{l=1}^{r} \frac{
\rho_\mathrm{pMD}(\frs\frt^2\frs)_{lu} \rho_\mathrm{pMD}(\frs\frt^{-2}\frs^{-1})_{lu} \rho_\mathrm{pMD}(\frs^{-1})_{li}}{ \rho_\mathrm{pMD}(\frs)_{lu} }.
\end{align}
The 2nd pseudo Frobenius-Schur indicator $\t\nu_2(i)$ satisfies 
\begin{align}
\t\nu_2(i)=\pm
\rho_\mathrm{pMD}(\frs^2)_{ii},\ \ \ \  \t \nu_2(u) =1.  
\end{align}
The  2nd pseudo Frobenius-Schur indicator give us important zero conditions on
$u_i$'s, which can effectively help us to determine $u_i$'s and $r_i$'s.  However, for more
complicated cases (such as when the number of variables is more then
$~12$), the zero conditions from the Frobenius-Schur indicator can
be too complicated.

We would like to remark that the zero conditions from the Frobenius-Schur
indicator have the form $f_1(u_i)/g_1(u_i)=0$, where $f_1(u_i)$ and $g_1(u_i)$
are multi-variable polynomials.  We convert the condition $f_1(u_i)/g_1(u_i)=0$
to a zero condition $f_1(u_i)=0$.  Such a conversion may create some fake
solutions of $u_i$.  We can rule out those fake solutions later, when we check
if the resulting $(S,T)$ matrices satisfy the conditions for modular data or
not.


\subsubsection{Reduce the number of variables}
\label{reducevars}

In order to factorize zero conditions, it is important to find a Gr\"obner
basis for and-connected zero conditions, such that some zero conditions have
fewest variables, in particular single variable.  We use the following strategy
to find such a Gr\"obner basis.  Consider a set of and-connected zero
conditions $[ f_1(u_i) = 0 \text{ and } f_2(u_i) = 0 \text{ and } f_2(u_i) = 0
\cdots ]$.  We want to use the $i^\text{th}$ zero condition $f_i$ to transform
the $i^\text{th}$ zero condition $f_j$.  Assume $f_i$ is a sum of a few
monomials: $f_i = m_1(u_i)+m_2(u_i)+m_3(u_i)$.  We can use the substitution
$m_1 \to -m_2-m_3$ to transform the zero condition $f_j$.  We can also use the
substitution $m_2 \to -m_1-m_3$ to transform the zero condition $f_j$, \etc.
Among all those transformed $f_j$, we select those with fewer variables than
that of $f_j$, and replace $f_j$ by those transformed zero conditions.  If the
one with fewer variables does not exist, we just keep original zero condition
$f_j$.

We perform the calculation for all the pairs $f_i$ and $f_j$.  We then perform
the above  calculation for a few iterations.  We have tested this approach and
found that this is an effective way to reduce the number of variables in zero
conditions.  This method is very important for our calculations.  In
particular, if we can obtain zero conditions with only single variable and if
the conductor of the variable is small, then there is an effective GAP/Singular
function to factorize the single-variable polynomials.

We can also use the above substitution approach, try to eliminate some
variables by simply reducing the number of variables that we want to eliminate
in zero conditions.

\subsubsection{From $u$-stage to $r$-stage}

We replace the $u_i$-dependent matrix elements of $\rho_\mathrm{pMD}(\frs)$
with cyclotomic numbers of a conductor $\ord(\rho_\mathrm{pMD}(\frt))$. Those
cyclotomic numbers are expressed in terms of $r_i$ variables, where $r_i$ are
the expansion coefficients over a cyclotomic basis.  This way, we express
$\rho_\mathrm{pMD}(\frs)$ in terms of $r_i$ variables.  

Compare $\rho_\mathrm{pMD}(\frs)$ in terms of $u_i$ and
$\rho_\mathrm{pMD}(\frs)$ in terms of $r_i$, we obtain many $u$-$r$ relations.
This allows us to convert some $u_i$ zero conditions to $r_i$ zero conditions,
by trying to eliminate the $u$-variables.  The $u_i$ zero conditions that
cannot be converted to $r_i$ zero conditions will be dropped.  As a result,
$r_i$-dependent $\rho_\mathrm{pMD}(\frs)$ may not be equivalent to the starting
$\SL$ representation $\rho_\textrm{isum}(\frs)$.  To partially fix this
problem, we re-implement the $\SL$ conditions and the $D_\rho$ conditions
\eqref{pMDDrhosihat}, \eqref{pMDsiDrho}, on the $r_i$-dependent
$\rho_\mathrm{pMD}$, to obtain additional $r_i$ zero conditions.  We also
compute some simple $\SL$ characters 
\begin{align}
\label{matchKar}
 \Tr_{\t\th} \big(\rho(\frs) \big), \ \ \ 
 \Tr_{\t\th} \big(\rho(\frs) \rho(\frt)^n \rho(\frs)\big), 
\end{align}
for the $\rho_\mathrm{pMD}$ and $\rho_\textrm{isum}$ representations.  Here
$\Tr_\th$ is the trace in the eigenspace of $\rho(\frt)$ with eigenvalue
$\t\th$.  Matching those simple $\SL$ characters also give us additional $r_i$
zero conditions.  This complete our conversion from $u_i$'s to $r_i$'s.  

\subsubsection{Integer conditions and inverse-pair of integer conditions}
\label{invpair}

In $r$-stage, 
$\frac{1}{ \rho_\textrm{pMD}(\frs)_{ui}}$, 
$\frac{\rho_\textrm{pMD}(\frs)_{ij} }{
\rho_\textrm{pMD}(\frs)_{uj}}$ and $\frac{\rho_\textrm{pMD}(\frs)_{ij} }{
\rho_\textrm{pMD}(\frs)_{uu}}$ are cyclotomic integers in term of the rational
$r_i$ variables. 
We can expand the cyclotomic integers in an integral basis of
cyclotomic numbers.  The expansion coefficients are ratios of polynomials of
$r_i$ with rational coefficients, which must be equal to integers.  This gives
us many $r_i$ integer conditions of form
\begin{align}
\label{rintcond}
 \frac{f_1(r_i)} {g_1(r_i)} \in \Z,\ \
 \frac{f_2(r_i)} {g_2(r_i)} \in \Z, \ \cdots
\end{align}

If two $r_i$ integer conditions, $ \frac{f_1(r_i)} {g_1(r_i)} \in \Z,\
\frac{f_2(r_i)} {g_2(r_i)} \in \Z $, satisfy
\begin{align}
 \frac{f_1(r_i)} {g_1(r_i)} 
 \frac{f_2(r_i)} {g_2(r_i)} = n \in \Z,
\end{align}
we will call them an inverse-pair of integer conditions.  For each inverse-pair
of integer conditions, we can obtain a set of $r_i$ zero conditions: 
\begin{align}
 \frac{f_1(r_i)} {g_1(r_i)} = m_1 \ \text{ or }\ 
 \frac{f_1(r_i)} {g_1(r_i)} = -m_1 \ \text{ or }\ 
 \frac{f_1(r_i)} {g_1(r_i)} = m_2 \ \text{ or }\ 
 \frac{f_1(r_i)} {g_1(r_i)} = -m_2 \ \text{ or }\ 
\cdots,
\end{align}
where $m_1,\ m_2, \cdots$ are
factors of $n$.

We note that when $\rho_\textrm{pMD}(\frs)_{uu}$ contains no variables,
$\frac{1}{ \rho_\textrm{pMD}(\frs)_{ui}}$ and
$\frac{\rho_\textrm{pMD}(\frs)_{ui}}{ \rho_\textrm{pMD}(\frs)_{uu}}$ are both
cyclotomic integers.  They produce inverse-pairs of integer conditions.  In
particular norm$(\frac{1}{ \rho_\textrm{pMD}(\frs)_{ui}})$ and
norm$(\frac{\rho_\textrm{pMD}(\frs)_{ui}}{ \rho_\textrm{pMD}(\frs)_{uu}})$ is an
inverse-pair of integer conditions.

\subsubsection{Solve $\rho_\textrm{pMD}(\frs)_{ui}$ using
generalized Egyptian-fraction method}
\label{gEFrac}

The generalized Egyptian-fraction method described here
is an important workhorse in our calculation.
Let $\{p_i\}$ be the set of prime divisors of pord$(\rho_\textrm{pMD}(\frt))$.
We note that norm$(\frac{1}{\rho_\textrm{pMD}(\frs)_{ui}})$ are integers, whose
prime divisors are contained in $\{p_i\}$.  Similarly,
norm$(\frac{1}{\rho_\textrm{pMD}(\frs)_{uu}})/\prod p_i$ is also an integer,
whose prime divisors are contained in $\{p_i\}$.  We will try to use such
conditions to find possible values of $\rho_\textrm{pMD}(\frs)_{ui}$.

The trick is to introduce variables $v_i$ that are reciprocals of positive
integers, whose prime divisors are contained in $\{p_i\}$.
We use those variables to represent
norm$(\frac{1}{\rho_\textrm{pMD}(\frs)_{ui}})$
and
norm$(\frac{1}{\rho_\textrm{pMD}(\frs)_{uu}})/\prod p_i$:
\begin{align}
 \mathrm{norm}(\rho_\textrm{pMD}(\frs)_{ui}) &= v_i \text{ or } -v_i,
\ \ \text{ for } i\neq u
\nonumber\\
 \mathrm{norm}(\rho_\textrm{pMD}(\frs)_{uu}) \prod p_i &= v_u \text{ or } -v_u .
\end{align}
The choices of $\pm$ signs produce a lot of cases, and we need
handle those cases one by one.

For each case, we combine the zero conditions,
$\mathrm{norm}(\rho_\textrm{pMD}(\frs)_{ui}) \pm v_i=0$ and
$\mathrm{norm}(\rho_\textrm{pMD}(\frs)_{uu})\prod p_i \pm v_u = 0$, of $v_i$'s
and $r_i$'s, with other zero conditions of $r_i$'s.  Then, we use the method
outlined in Section \ref{reducevars} to eliminate the $r_i$ variables, trying
to obtain, as many as possible, the zero conditions that contain only $v_i$'s.

A zero condition with only $v_i$'s has the following form
\begin{align}
\label{ccPcP}
 c_0 
+ c_1^+P_1^+(v_i) + c_2^+P_2^+(v_i) +\cdots
- c_1^-P_1^-(v_i) - c_2^-P_2^-(v_i) +\cdots
\end{align}
where $\ c_i^\pm$ are positive integers, $c_0$ is a non-negative inter, and
$P_i^\pm(v_i)$ are products of $v_i$'s.  We apply the following generalized
Egyptian-fraction method to solve the above zero conditions:
\begin{enumerate}

\item
If $c_0 = 0$, the iteration stops, and we only get a single zero
condition \eqref{ccPcP}.

\item
If $c_0 > 0$, then $c_1^-P_1^-(v_i) + c_2^-P_2^-(v_i) +\cdots \geq c_0$.  Let
$c_{j_\mathrm{large}}^-P_{j_\mathrm{large}}^-(v_i)$ be a term among
$c_j^-P_j^-(v_i)$'s, such that $c_{j_\mathrm{large}}^-P_{j_\mathrm{large}}^-(v_i)
\geq c_0/N^-$, where $N^-$ is the number of terms in $c_1^-P_1^-(v_i) +
c_2^-P_2^-(v_i) +\cdots$.  
Such a term must exit.
Because $P_{j_\mathrm{large}}^-(v_i)$ is the inverse
of a positive integer (whose prime divisors are part of those of
$\pord(\rho_\mathrm{pMD}(\frt))$), $P_{j_\mathrm{large}}^-(v_i)$ has only a finite
number of possible values 
$ P_{j_\mathrm{large}}^-(v_i) =1/n^-$, where
\begin{align}
&
n^-  \in \N,
\ \ \ \
n^- \leq \frac{N^- c_{j_\mathrm{large}}^-}{c_0} , 
\nonumber\\
&
\text{ the prime divisors of }
n^- \text{ are contained in } \{p_i\}
.
\end{align}
For each possibility, the zero condition \eqref{ccPcP} is reduced to a zero
condition with one less terms and $P_{j_\mathrm{large}}^-(v_i) =1/n^-$.

We also need to run through different choices of large term among
$c_j^-P_j^-(v_i)$'s.  At the end, we obtain many or-connected sets of
and-connected zero conditions.  A set of and-connected zero conditions contains
$ P_{j_\mathrm{large}}^-(v_i) =1/n^-$ and the reduced zero condition from
\eqref{ccPcP}, as discussed above.

\item
If $c_0 < 0$, then $c_1^+P_1^+(v_i) + c_2^+P_2^+(v_i) +\cdots \geq -c_0$.  We
repeat the above calculation for $c_1^+P_1^+(v_i) + c_2^+P_2^+(v_i) +\cdots$,
and obtain  many or-connected sets of and-connected zero conditions.  A set of
and-connected zero conditions contains $ P_{j_\mathrm{large}}^+(v_i) =1/n^+$ and
the reduced zero condition from \eqref{ccPcP}.

\end{enumerate}

We perform the above calculation, starting from the simplest zero condition of
$v_i$'s that has fewest terms and a non-zero $c_0$.  Then we handle the next simplest zero condition
of $v_i$'s.  At last, we process the resulting or-connected sets of
and-connected zero conditions of $v_i$'s and $r_i$'s, trying to replace $v_i$'s
by $r_i$'s, and convert the zero conditions of $v_i$'s and $r_i$'s to zero
conditions of $r_i$'s. The remaining zero conditions containing $v_i$'s will be
ignored.  This will give us or-connected sets of and-connected zero conditions
of $r_i$'s.

\subsubsection{Integer conditions from norm($\rho_\textrm{pMD}(\frs)_{uu}$)}
\label{tSuucnd5}

We note that $\frac{1}{\rho_\textrm{pMD}(\frs)_{uu}} = \pm D $ or $\pm \ii D$,
where $u$ is the index for the unit row.  Therefore
norm$(\frac{1}{\rho_\textrm{pMD}(\frs)_{uu}} )$ is an integer whose prime
divisors coincide with the  prime divisors of pord$(\rho_\textrm{pMD}(\frt))$.
At $r$-stage, norm$(\frac{1}{\rho_\textrm{pMD}(\frs)_{uu}} )$ is a ratio of two
polynomials of $r_i$'s.  Sometimes,
norm$(\frac{1}{\rho_\textrm{pMD}(\frs)_{uu}} ) \in \Z$ has only finite numbers
of solutions for rational variables $r_i$'s. We can use this property to solve
$\frac{1}{\rho_\textrm{pMD}(\frs)_{uu}}$.

In particular,
there is a class of cases which is hard to solve.  For those cases,
$\rho_\textrm{pMD}(\frs)_{uu}$ only depends on a single rational variable $r$
and 
\begin{align}
\rho_\textrm{pMD}(\frs)_{uu} = r - \frac{\sqrt{5}}{10},\ \ \ \ \ \
 \text{norm}(\frac{1}{\rho_\textrm{pMD}(\frs)_{uu}} )
=
\frac{1}{r^2 -\frac{1}{20}} \in \Z
\end{align}
Also the prime divisors of pord$(\rho_\textrm{pMD}(\frt))$ are contained in
$[2,3,5]$.  This class of cases can be solved in the following way.

Recall that the ring of algebraic integers of $\Q(\sqrt{5})$ is $\Z[\phi]$
where $\phi=  \frac{1+\sqrt{5}}{2}$. It is well-known that the group $G$ of
invertible elements in $\Z[\phi]$ is isomorphic to $\Z_2 \times \Z$, and is
given by $ \{\pm \phi^n \mid n \in \Z\} $.  $\Z[\sqrt{5}]$ is a subring of
$\Z[\phi]$, and its group  $G'$ of invertible elements, which is a subgroup of
$G$. Since $n=3$ is the smallest positive integer $n$ such that $\phi^n \in
\Z[\sqrt{5}]$, 
\begin{align}
G' = \{\pm \phi^{3n} \mid  n \in \Z\}\,.
\end{align}
Note that $\phi^3= 2+\sqrt{5}$ and norm$(\phi^3)=-1$. Thus, $\phi^6 = 9+4\sqrt{5}$ and norm$(\phi^6)=1$. If $u, v \in \N$ such that $u^2 -5 v^2 = \pm 1$, then 
\begin{equation}\label{eq1}
u+ v\sqrt{5}= \left\{
\begin{array}{ll}
(2+\sqrt{5})^{2n}=(9+4\sqrt{5})^n & \text{ if } u^2 -5 v^2 = 1, \\
(2+\sqrt{5})^{2n+1}=(2+\sqrt{5})(9+4\sqrt{5})^n & \text{ if } u^2 -5 v^2 = -1\,. \\
\end{array} \right.
\end{equation}
for some positive integer $n$.
\begin{lem}\label{l:1}
Let $u, v$ be nonzero integers such that $v>0$, $(u^2-5v^2) \mid 5$ and $u=2^a 3^b 5^c$ for some nonnegative integers $a,b,c$. Then $(u, v)$ can only be one of the following pairs:
\begin{enumerate}
\item[\rm (i)] $ (2, \pm 1), (9, \pm 4)$ which are the solutions of $u^2-5 v^2 = \pm 1$.
\item[\rm (ii)] $(5,\pm 2), (20, \pm 9), (360,\pm 161)$ which are the solutions of $u^2-5 v^2 = \pm 5$.
\end{enumerate}
\end{lem}
\begin{proof}
Since $u^2-5 v^2 \mid 5$, we have two cases: $u^2-5 v^2= \eps$ or $5\eps$ for some $\eps = \pm 1$. 

\noindent (i) $u^2-5 v^2= \eps$. It is immediate to see that $5 \nmid u$ and so $c =0$. Next,
we show that $a \le 2$. Assume to the contrary that $a > 2$. Then $u^2 \equiv 0 \mod 8$, and $v$ must be odd. Thus, $u^2 -5 v^2 \equiv 3  \not \equiv \eps \mod 8$. Therefore, $a \le 2$. 

Now, we show that $b \le 2$. Suppose not. Then $b \ge 3$ which implies $3 \nmid v$ and $u \equiv 0 \mod 27$. Since $u^2 -5 v^2 \equiv -5 \equiv 1 \mod 3$, $\eps =1$. By \eqref{eq1}, we find 
\begin{align}
\label{eq2} 2^a 3^b = u &= \frac{1}{2}\left((9+4 \sqrt{5})^n + (9-4 \sqrt{5})^n\right) =\sum_{0 \le 2 i \le n} \binom{n}{2i} 9^{n-2i} \cdot 4^{2i}\cdot 5^{i} \\
\label{eq3}  & \equiv   \left\{\begin{array}{ll}
4^{n}\cdot 5^{\frac{n}{2}} 
\mod 27  & \text{ if } n \text{ is even,}\\
n \cdot 9\cdot  4^{n-1}\cdot 5^{\frac{n-1}{2}} \mod 27 & \text{ if } n \text{ is odd.}
\end{array}\right.
\end{align}
Since the leftmost expression of \eqref{eq2} is divisible by 27, it follows from \eqref{eq3} that $n$ must be odd and $3 \mid n$. Thus, $n = 3k$ for some positive odd integer $k$. It follows from \eqref{eq2} again
\begin{align}
2^a 3^b = \frac{1}{2}\left((9+4 \sqrt{5})^{3k} + (9-4 \sqrt{5})^{3k}\right) =
\frac{1}{2}\left((9+4 \sqrt{5})^{3} + (9-4 \sqrt{5})^{3}\right) z = 2889 z = 3^3 \cdot 107 \cdot z
\end{align}
for some integer $z$ since $k$ is odd. However, this is a contradiction since $107 \nmid 2^a 3^b$. Therefore, $b \le 2$. 

Now, we can solve for the integral solutions $v$ of the equations $u^2-5 v^2 = \pm 1$ for $u=2^a 3^b$ with $ 0 \le a, b \le 2$. The integral solutions to $u^2-5v^2 =\pm 1$ are $(u,v) = (2, \pm 1)$ and $(9,\pm 4)$.

\noindent (ii) $u^2-5 v^2= 5\eps$. Then $5$ divides $u$ (or $c \ge 1$), and we have $v^2 - 5 \overline u^2 =\pm 1$ where $\overline u = u/5 = 2^a 3^b 5^{c-1}$. 

We first show that $c=1$. Assume to the contrary that $c >1$. Then $\overline u \equiv 0 \mod 5$. By \eqref{eq1},  there is a positive integer $n$ such that 
\begin{equation}
\label{eq4}  \overline u=  \frac{1}{2\sqrt{5}}\left((2+\sqrt{5})^n - (2- \sqrt{5})^n\right) = \sum_{0 \le 2 i+1 \le n} \binom{n}{2i+1} 2^{n-2i-1} \cdot 5^{i} \equiv 
\binom{n}{1} 2^{n-1} \mod 5\,.
\end{equation}
This implies $\binom{n}{1} 2^{n-1} \equiv 0 \mod 5$ , and hence $n=5k$ for some integer $k$. Therefore,
\begin{align}
2^a 3^b 5^{c-1} =\frac{1}{2\sqrt{5}}\left((2+\sqrt{5})^{5k} - (2- \sqrt{5})^{5k}\right) = \frac{1}{2\sqrt{5}}\left((2+\sqrt{5})^{5} - (2- \sqrt{5})^{5}\right) z = 305 z = 5\cdot 31 \cdot z
\end{align}
for some integer $z$, which is absurd. Therefore,  $c=1$, and so $\overline u =2^a 3^b$.

Next, we show that  $a \le 3$ and $b \le 2$.   If $a > 3$, then $v$ is odd. Since $v^2 \not\equiv -1 \mod 16$, $\eps=1$ and so
\begin{equation} \label{eq5}
v^2 -5 \overline u^2 =  1\,.
\end{equation}
Similarly, if $b > 2$, then $3 \nmid v$. Since $v^2 \not\equiv -1 \mod 27$, $\eps=1$.  Therefore, $(\overline u,v)$ can only satisfy \eqref{eq5}. It follows from \eqref{eq1} that
\begin{eqnarray}
\label{eq6} \overline u &=& \frac{1}{2\sqrt{5}}\left((9+4\sqrt{5})^{n} - (9- 4\sqrt{5})^{n}\right) = \sum_{0 \le 2i+1 \le n} \binom{n}{2i+1} 9^{n-2i-1}\cdot 4^{2i+1} \cdot 5^{i} \\ 
\label{eq7}  & \equiv &  \left\{\begin{array}{ll}
\binom{n}{1}9^{n-1} 4  \mod 16 & \\
4^n \cdot 5^{(n-1)/2} \mod 27 & \text{if $n$ is odd;}\\
\binom{n}{1}9 \cdot 4^{n-1}\cdot 5^{(n-1)/2} \mod 27 & \text{if $n$ is even.}\
\end{array}\right. 
\end{eqnarray}

If $a > 3$, $\overline u \equiv 0 \mod 16$. Therefore, $n=4k$ for some integer $k$ by \eqref{eq7}. It follows from \eqref{eq6} that
\begin{align}
2^a 3^b &= \overline u= \frac{1}{2\sqrt{5}}\left((9+4\sqrt{5})^{4k} - (9- 4\sqrt{5})^{4k}\right) = \frac{1}{2\sqrt{5}}\left((9+4\sqrt{5})^{4} - (9- 4\sqrt{5})^{4}\right) \cdot z 
\nonumber\\
&= 2^4 3^27\cdot 23\cdot z
\end{align} 
for some integer $z$ but this is absurd. Therefore, $a \le 3$.

If $b > 2$, $\overline u \equiv 0 \mod 27$. Therefore, $n=6k$ for some integer $k$ by \eqref{eq7}. It follows from \eqref{eq6} that
\begin{align}
2^a 3^b &= \overline u = \frac{1}{2\sqrt{5}}\left((9+4\sqrt{5})^{6k} - (9- 4\sqrt{5})^{6k}\right) = \frac{1}{2\sqrt{5}}\left((9+4\sqrt{5})^{6} - (9- 4\sqrt{5})^{6}\right) \cdot z 
\nonumber\\
&= 2^3 3^3 17\cdot 19\cdot 107\cdot  z
\end{align} 
for some integer $z$ but this is absurd. Therefore, $b \le 2$.

Now, we can solve for  $v$ of the equations $u^2 - 5 v^2 = \pm 5$ for $u=2^a\cdot 3^b\cdot 5$ with $ 0 \le a \le 3,\, 0\le b \le 2$. The integral solutions to $u^2-5v^2 =\pm 5$ are $(u,v) = (5,\pm 2)$, $(20, \pm 9)$  and  $(360,\pm 161)$.
\end{proof}

\begin{prop}\label{p:2}
Let $r \in \Q$ such that $(r^2 - \frac{1}{20})^{-1}=\eps \cdot 2^{\al+1}\cdot 3^\bt \cdot 5^\ga$ for some 
nonnegative integers  $\al,\bt,\ga$ and $\eps = \pm 1$. Then, 
\begin{align}
r \in  \left\{0, \pm \frac{1}{4},\pm \frac{1}{5},  \pm \frac {2}{9},  \pm \frac{9}{40}, \pm \frac{161}{720} \right\}\,.
\end{align}
\end{prop}

\begin{proof}
It is clear that $r=0$ is a solution. Now, we assume $r = \frac{x}{y}$ where $x,y$ are nonzero coprime integers with $y > 0$. Then
\begin{equation}\label{eq:2}
\frac{20 y^2}{y^2-20 x^2}=\eps \cdot 2^{\al+1}\cdot 3^\bt \cdot 5^\ga, \quad \text{or} \quad \frac{10 y^2}{y^2-20 x^2}=\eps \cdot 2^{\al}\cdot 3^\bt \cdot 5^\ga\,.
\end{equation}
Thus, $(10-\eps \cdot 2^{\al}\cdot 3^\bt \cdot 5^\ga)\cdot y^2 = -\eps \cdot 2^{\al+2}\cdot 3^\bt \cdot 5^{\ga+1}\cdot x^2$. Since $x, y$ are coprime, $y^2 \mid 2^{\al+2}\cdot 3^b \cdot 5^{\ga+1}$ or
\begin{equation}
y = 2^{a}\cdot 3^{b}\cdot 5^{c}
\end{equation}
for some nonnegative integers $a, b, c$ such that $2a \le \al+2$, $2b \le \bt$, $2c \le \ga+1$.

Equations \eqref{eq:2} also imply that $y^2-20 x^2$ divides $10 y^2$ and so $y^2-20 x^2$ divides $200 x^2$. Since $x, y$ are coprime, $y^2-20 x^2$ divides $\gcd(200 x^2, 200y^2)=200$. Note that  if $5 \mid   (y^2-20 x^2)$, $5 \mid y$. Since $x,y$ are coprime, $5 \nmid x$ and hence $y^2 - 20 x^2 \not\equiv 0 \mod 25$. Therefore, $(y^2 - 20 x^2) \mid 40$.

(i) Suppose $y$ is odd. Then $y^2-20 x^2$ is an odd divisor of 200. Thus, $(y^2-20 x^2) \mid 5$, and so $(y^2 - 5 \overline x^2)\mid 5$ where $\overline x = 2x$. It follows from  Lemma \ref{l:1} that $(y, \overline x) =(9, \pm 4)$ or $(5, \pm 2)$. Therefore,
$r = \frac{x}{y} = \pm \frac{2}{9}$ or $\pm \frac{1}{5}$.

(ii) Now we assume $y$ is even. Then $y^2-20 x^2$ is also and $x$ must be odd.  Thus, $4 \mid (y^2-20 x^2)$. If $8 \mid y^2 -20 x^2$, then $2 \mid \left( (\frac{y}{2})^2 -5 x^2\right)$ which implies $\frac{y}{2}$ is odd. Therefore, $(\frac{y}{2})^2, 5 x^2 \equiv 1 \mod 4$ and so $4 \mid \left((\frac{y}{2})^2 - 5 x^2\right)$ or $16 \mid (y^2 - 20 x^2)$, which cannot be a divisor of $40$. As a consequence, we find
\begin{align}
y^2 -20 x^2 = \pm 4   \text{ or } \pm 20 \,.
\end{align}
Let $\overline y = y/2$. These diophantine equations are equivalent to $(\overline y^2 -5 x^2)\mid 5$.  It follows from  Lemma \ref{l:1} that $(\overline y, x) =(2, \pm 1)$, $(20, \pm 9)$ or $(360, \pm 161)$. Therefore,
$r = \frac{x}{y} = \pm \frac{1}{4}$,  $\pm \frac{9}{40}$ or $\pm \frac{161}{720}$.
\end{proof}

Because $r$ has only a finite number of solutions, thus
$\rho_\textrm{pMD}(\frs)_{uu}$ only takes a finite number of possible values,
and becomes known. This makes the cases soluble using the methods of
inverse-pair of integer conditions discussed in Section \ref{invpair}.

\subsubsection{Go to $n$-stage}

We can go to $n$-stage from $r$-stage or $u$-stage, by considering $K_{ij},
\tilde K_{ij}, J_i$ introduced in \eqref{KKJ}.  $K_{ij}, \tilde K_{ij}, J_i$ are
cyclotomic integers, and those cyclotomic integers are expressed in terms of
$n_i$ variables,  where the $n_i$ variables are the expansion coefficients of the
cyclotomic integers over an integral basis of cyclotomic numbers.  Those $n_i$
variables satisfy many conditions, which are listed in the Proposition
\ref{KKQcond}.

\subsubsection{Solve bounded $n_i$ zero conditions}

A $n_i$ zero condition is bounded if it has the following form
\begin{align}
\text{const.} - \sum_{\{i,j,\cdots} c_{i,j,\cdots} 
n_i^{a_i}
n_j^{a_j}\cdots = 0,\ \ \ a_i,a_j,\cdots \in 2\N,\ \ c_{i,j,\cdots} > 0.
\end{align}
In this case, $n_i$ that can satisfy the zero condition have a simple finite
range.  We can test each set of $n_i$'s in this finite range to see which set
of $n_i$'s are solutions.

\subsubsection{Solve bounded $n_i$ integer conditions}

A simple example of bounded $n_i$ integer condition is given by
$\frac{n_1^2}{1+n_1+n_1^2} \in \Z$.  The integer condition can be reduced to
$\frac{-1-n_1}{1+n_1+n_1^2} \in \Z$.  The reduced integer condition can be
written as
$\frac{-\frac{1}{n_1^2}-\frac{1}{n_1}}{\frac{1}{n_1^2}+\frac{1}{n_1}+1} \in
\Z$, which has the form that except a constant term in the denominator, all other
terms has negative powers.  In this case, $|n_1|$ cannot be too big to satisfy
the  integer condition. Thus we check $n_1$ in a finite range to see which
values satisfy the integer condition.

Let us consider a more general example of integer condition 
\begin{align}
h(n_1,n_2) = \frac{n_1-n_2}{1+n_1+n_1^2n_2} \in \Z, 
\end{align}
which can be rewritten as $
h(n_1,n_2) = \frac{\frac{1}{n_1n_2} - \frac{1}{n_1^2}}{\frac{1}{n_1^2n_2} +
\frac{1}{n_1n_2}+1} \in \Z$.  Except a ``1'' in the denominator, all other
terms have negative exponents, which indicates that  the integer condition is a
bounded integer condition.  We then introduce
\begin{align}
 h_\mathrm{max} (n_1,n_2) =  
\frac{\frac{1}{|n_1n_2|} 
+ \frac{1}{|n_1^2|}}{-\frac{1}{|n_1^2n_2|} -
\frac{1}{|n_1n_2|}+1} 
\end{align}
which is obtained from $h (n_1,n_2)$ by replace $n_i$ with $|n_i|$, changing
the sign of the terms in the numerators to ``$+$'', and the sign of the terms in
the denominator (except the ``1'' term) to ``$-$'', $h_\mathrm{max}$ satisfies $
h_\mathrm{max} (n_1,n_2) \geq |h (n_1,n_2)| $ if $h_\mathrm{max} (n_1,n_2) > 0$.

We note that when $|n_i|$ are large, $ 0 < h_\mathrm{max} (n_1,n_2) <1$.  In this
case $h (n_1,n_2)$ cannot be integer except $h (n_1,n_2)=0$.  Thus we only need
to check $n_1,n_2$'s in a finite range to see if $h (n_1,n_2)$ is integer or
not.  More specifically, we test $n_1=\pm 1, n_2=\pm 1$.  Then $n_1=\pm 2,
n_2=\pm 1$ and $n_1=\pm 1, n_2=\pm 2$, {\it etc}, until $0 < h_\mathrm{max}
(n_1,n_2) < 1$.  Those values of $n_1,n_2$ that make $h (n_1,n_2) \in \Z$ give
us a set of zero conditions $ [ n_1=\al_1 \text{ and } n_2 =\al_2 ] \text{ or }
[ n_1=\al_1' \text{ and } n_2 =\al_2' ] \text{ or } \cdots$, where $\al_i$ and $\al_i'$ are integral constant.

Since we did not check $n_1=0$ or $n_2=0$, as well as the
possibility that $h (n_1,n_2)=0$, we need to add those possible cases back to the
sets of zero conditions. Thus the final sets of zero conditions are given by $
[h (n_1,n_2)=0] \text{ or } [n_1=0] \text{ or } [n_2=0] \text{ or } [ n_1=\al_1
\text{ and } n_2 =\al_2 ] \text{ or } [ n_1=\al_1' \text{ and } n_2 =\al_2' ]
\text{ or } \cdots$.  Those are the sets of zero conditions obtained from the
bounded integer condition $h(n_1,n_2) \in \Z$.

\section{Classification of integral modular data}
\label{intcat}

\subsection{Mathematical results on integral modular categories}
\label{mathintcat}

Various approaches to classifying \emph{integral} modular categories have been developed, see, eg. \cite{BGHKNNPR,BGNPRW,BR,CP}.
Here we summarize some of these approaches with an eye towards automation.

 Here integral means that the FP-dimension $\FPdim(X_i)$ of each simple object is a (necessarily positive) integer. Let $\eC$ be an integral MTC. Firstly, by \cite[Props. 8.23 and 8.24 ]{ENO} every integral category is pseudo-unitary so that one may assume that the dimensions $d_i$ are positive integers (equal to $\FPdim(X_i)$), by adjusting the spherical structure if necessary.  Assume that $1=d_1\leq d_2\leq \cdots\leq d_r$ and set $x_{r-i+1}:=\frac{\dim(\eC)}{d_i^2}\in \mathbb{N}$.  We have the following:

\begin{enumerate}
\item The $x_i$ satisfy the (Egyptian fraction) Diophantine equation  \begin{equation}\label{eqn:Egypt}
1=\sum_{i=1}^r\frac{1}{x_i}.\end{equation}
\item  By a classical result of Landau \cite{Landau} eqn (\ref{eqn:Egypt}) has finitely many solutions $(x_1,\ldots,x_r)\in\mathbb{N}^r$ for fixed $r$.   Indeed $k\leq x_k\leq u_k(r-k+1)$ where $u_k$ is Sylvester's sequence defined by $u_1:=1$ and $u_k:=u_{k-1}(u_{k-1}+1)$.  
\item For any $x_i,x_j$ we have that $\frac{x_i}{x_j}=\frac{d_j^2}{d_i^2}$ is a square rational number.  In particular if for some prime $p$, $p^k\mid\mid x_i$ and $p^m\mid\mid x_j$ then the parity of $k$ and $m$ are the same (here $p^s\mid\mid b$ means that $p^s\mid b$ and $p^{s+1}\not\mid b$).
\item Let $T$ be the $T$-matrix of $\eC$, and let $\wp:=\{p: p\mid\ord(T)\}$ be the set of primes dividing $\ord(T)$.  By the Cauchy theorem \cite{BNRW} 
$\dim(\eC)=\prod_{p\in\wp}p^{\alpha_p}$  where $1\leq \alpha_p\in\mathbb{N}$.  In particular, if $p\mid x_i$ then $p\in\wp$.
\item Suppose that $|\wp|\leq 2$.  Then there are at least 2 invertible objects, i.e. $1=d_1=d_2$.  If $p$ is the minimal prime in $\wp$ then we have $1=d_1=\cdots =d_p$.    This follows from the fact \cite{ENOsolvable} that any fusion category of dimension $p^aq^b$ for primes $p,q$ is \emph{solvable}, and any solvable fusion category has a non-trivial invertible object.
\end{enumerate}

While the bounds provided by Sylvester's sequence are doubly exponential, in our setting the primes dividing $x_i$ are restricted to the finite set $\wp$.  For example, if $\ord(T)=2^a3^b$ then  $x_i=2^{a_i}3^{b_i}$ for $a_i,b_i$ bounded by $\log_2( u_k(r-k+1))$, a much more practical bound.  For modest sized $r$, the known bounds on the primes in $\wp$ make finding solutions to equation (\ref{eqn:Egypt}), restricted to $\wp$, reasonably efficient.  

In practice we find that it is frequently the case that $|\wp |\leq 2$, so that we are assured of having invertible objects.  For example by \cite{CP,CGP} every \emph{odd dimensional} integral modular category of rank at most  23 must have non-trivial invertible objects--they are all pointed!  The smallest rank for which we are aware of an integral modular category with no non-trivial invertible objects is 22: this comes from the Drinfeld center of $\Rep_{A_5}$.  A category is called \emph{perfect} if it has no non-trivial invertible objects.

Generally, a fusion category $\mathcal{C}$ is \textbf{$G$-graded} if there is a decomposition as abelian categories $\mathcal{C}\cong\bigoplus_g \mathcal{C}_g$ such that if $\mathcal{C}_g\otimes \mathcal{C}_h\subset \mathcal{C}_{gh}$, for some group $G$. The grading is faithful if each $\mathcal{C}_g$ is non-trivial. A useful consequence of a faithful $G$-grading is that $\dim(\mathcal{C}_g)$ is constant for all $g$: it is equal to $\dim(\eC)/|G|$. This induces a partition of the simple objects, and hence the list of dimensions. It is clear that the monoidal unit $\mathbf{1}$ must lie in the trivial component $\mathcal{C}_e$, and therefore if $X\in\mathcal{C}_g$ then $X^*\in\mathcal{C}_{g^{-1}}$.  Thus, if every object is self-dual, any faithful grading group must be an elementary abelian $2$-group (but not conversely, there could be non-self-dual objects in the trivially-graded component).

For a modular category $\mathcal{C}$ the largest faithful grading group (called the universal grading) is isomorphic to the group of isomorphism classes of invertible simple objects (see \cite{tcat} for details) so is abelian in this case.  Moreover the trivial component with respect to the universal grading is the adjoint subcategory. A natural way to understand this grading is by defining, for $\varphi\in \hat{A}$, $\mathcal{C}_\varphi$ to be the abelian subcategory generated by objects $X$ such that $c_{X,z}c_{z,X}=\varphi(z) Id_{z\otimes X}$ for all $z\in A$. 

In general for $\eD\subset\eC$ ribbon categories, \textbf{the centralizer of $\eD$ in $\eC$}, denoted $C_{\eC}(\eD)$, is the subcategory generated by those $Y\in\eC$ so that $c_{Y,X}c_{X,Y}=Id_{X\otimes Y}$ for all $X\in\eD$.  For simple objects $X,Y$ we have that $c_{Y,X}c_{X,Y}=Id_{X\otimes Y}$ if and only if $S_{X,Y}=d_Xd_Y$.  The centralizer of $\eD$ in itself $C_\eD(\eD)$ is sometimes denoted $Sym(\eD)$ (we avoid the somewhat vague notation $\eD^\prime$ that is sometimes found in the literature).
For a MTC $\eC$ it is known that the trivial component $\eC_0$ with respect to the universal grading is the \emph{adjoint} subcategory, which is precisely the centralizer $C_{\eC_{pt}}(\eC)$ of the pointed subcategory.  In particular, the symmetric center of $\eC_0$ is  $Sym(\eC_0):=\eC_{pt}\cap  \eC_0$  and is both pointed and \emph{symmetric}.  Thus each simple object  $Z\in\eC_{pt}\cap  \eC_0$ is either bosonic or fermionic, i.e. $\theta_Z=\pm 1$ and $c_{Z,Z}^2=Id$.  If $Z\in Sym(\eC_0)$ is fermionic so that $\theta_Z=-1$ then for any other object $X\in\eC_0$ we have $Z\otimes X\not\cong X$. Indeed, by \eqref{eq:balancing} for $Z$ invertible in $\eC_0$ with $\theta_Z=-1$ we have 
$\theta_{Z^*}\theta_X S_{Z^*,X}=d_{Z\otimes X}\theta_{Z\otimes X}$ 
 so that $-\theta_X=\theta_{Z\otimes X}$.  We can often use this fact to show that $\eC_{pt}\cap  \eC_0$ is Tannakian, i.e., every simple object is bosonic.  For example, if there is only one simple object in $\eC_0$ of some dimension $d$, then there can be no fermionic invertible objects in $\eC_0$.  We can then condense the bosons in $\eC_0$ to obtain a new MTC $\eD$ (the modularization of $\eC_0$: recall that $C_{\eC_0}(\eC_0)=\eC_{pt}\cap \eC_0$ and see \cite{Brug}) with dimension $\dim(\eC_0)/\dim(\eC_{pt}\cap\eC_0)$.  The modularization is a braided tensor functor $F:\eC_0\rightarrow \eD$. This yields further constraints.

Another valuable fact is the following:
\begin{lem}\cite[Theorem 3.2]{MugerPLMS}\label{lem:dimsofsubcat}
 If   $\eD\subset\eC$ is a subcategory of a modular category $\eC$ then $\dim(\eC)=\dim(\eD)\dim(C_{\eC}(\eD))$.
\end{lem}

\begin{prop}\label{integralprop}
    Let $\eC$ be an integral MTC and $G(\eC)=A$ the group of invertible objects, so that $\eC$ has universal grading group $A$. Then:
 \begin{enumerate}
     \item[(a)] If $\eC_0$ has no non-trivial invertible objects then $\eC_0$ is itself modular, so that $\eC\cong \eC_0\boxtimes \eC_{pt}$.  In particular $\eC$ is not prime if $A$ is non-trivial.
     \item[(b)] If $\eC_0$ has a fermionic invertible object then each $d$ appearing as a dimension of a simple object in $\eC_0$ occurs with even multiplicity.
     \item[(c)] If $\eC_0$ has a fermionic invertible then after condensing the maximal Tannakian subcategory of $\eC_0$ the result is super-modular. Methods found in \cite{BGNPRW} can be brought to bear.
     \item[(d)] If $Sym(\eC_0)$ is Tannakian then every dimension $d$ appearing in $\eC_i$ for $i\neq 0$ must occur with multiplicity $\geq 2$.
     \item[(e)] If $Sym(\eC_0)\cong \hat{G}$ is Tannakian ($G$ is abelian) then the condensation $\eC_0$ by $G$ is modular, i.e. $(\eC_0)_G=\eD$ is modular (the \emph{modularization} of $\eC_0$ in \cite{Brug}).  Denote by $F:\eC_0\rightarrow \eD$ this modularization (or $G$-de-equivariantization) functor. In this case $G$ acts on the simple objects of $\eC_0$ by $X\mapsto g\otimes X$.  If $Stab_G(X)$ is cyclic (for example if $Stab_G(X)$ is trivial) then $F(X)$ is a direct sum of $|Stab_G(X)|$ pairwise non-isomorphic simple objects \cite[Lemme 4.3]{Brug}.
     
 \end{enumerate}
\end{prop}
\begin{proof} 
Proof of (a): Since $Sym(\eC_0)$ is pointed, the hypothesis imply $\eC_0$ is modular.  By \cite[Corollary 3.5]{MugerPLMS} $\eC\cong \eC_0\boxtimes \eC_{pt}$.

Proof of (b): If a fermionic invertible object $z\in\eC_0$ exists and $X\in\eC_0$ has dimension $d$ then $z\otimes X\not\cong X$: otherwise eqn \eqref{eq:balancing} implies $-\theta_Xd_X=\theta_{z^*}\theta_XS_{z^*,X}=d_{z^*\otimes X}\theta_{z^*\otimes X}=d_X\theta_X$, a contradiction.

Proof of (c): This is well-known, see eg. \cite{superrank6}.

Proof of (d): Suppose that some $\eC_i$ with $i\neq 0$ has rank $1$, i.e., has a single simple object $Y$.  Then for any simple $Z\in\eC_{pt}$ we have $Z\otimes Y=Y$.  Now the balancing equation \eqref{eq:balancing}
we have $$\theta_{Y}S_{Z,Y}=\theta_Z\theta_{Y^*}S_{Z,Y}=\theta_{Y}\dim(Y),$$ so that $S_{Z,Y}=\dim(Y)\dim(Z)$ for all $Z\in\eC_{pt}$.  This implies that $Y\in\eC_0$ contrary to assumption.

Proof of (e): This is directly taken from \cite{Brug}.
\end{proof}

These can be implemented in the computer calculations.  To give some idea of how this works we present a few examples by hand.  One useful fact is the following:
\begin{lem}\label{lem:bosoninorder4}
    If $\eD$ is a symmetric pointed category of dimension $2^k$ with $k\geq 2$ then $\eD$ contains an order 2 boson $b\not\cong\one$, i.e. an object such that $b^{\otimes 2}\cong\one$ and $\theta_b=1$.
\end{lem}
\begin{proof}
The twist of each simple (invertible) object $\theta_Z=\pm 1$ by
\eqref{eq:balancing}.  If $Z^2\cong \one$ and $\theta_Z=-1$ for each
$Z\not\cong\one$ we obtain the (modular) 3 fermion theory \cite{RSW0777} as a
subcategory, contradicting $\eD$ symmetric.  So either there is an invertible
object $W$ of order $2^s$ with $s\geq 2$ or a non-trivial invertible object $U$
with $\theta_U=1$. First consider such a $W$.  the balancing equation
\eqref{eq:balancing} gives:
$$1=(\theta_W)(\theta_{W^*})S_{W,W^*}=\theta_{W^2}$$ so that $U=W^2$ is a
nontrivial invertible object with $\theta_U=1$.  Thus we reduce to the second
case.  Let $U$ be chosen non-trivial with minimal order $2^t$.  We must show
$t=1$.  Again, \eqref{eq:balancing} implies:

$$1=\theta_U\theta_{U^*}S_{U,U^*}=\theta_{U^2}.$$  Thus by minimality of $t$ we see that $U^2\cong\one$, so $t=1$.
\end{proof}

The following lemma is straight-forward, but useful.
\begin{lem}\label{balancingtrick}
    Suppose that the trivial component $\eC_0$ of an integral MTC $\eC$ has an invertible object $b$ with twist $\theta_b=1$, and $X\in\eC_i$ for $i\neq 0$ is the unique object of dimension $d$.  Then
$b$ and $X$ centralize each other, i.e. $S_{b,X}=d$. In particular, if $Sym(\eC_0)$ is Tannakian then every component $\eC_i$ must have multiple objects of any given dimension $d$.
\end{lem}
\begin{proof} Since $X$ is the only object in $\eC_i$ of dimension $d$, we have $z\cdot X\cong X$ for any invertible $z\in\eC_0$. 
    The balancing equation \eqref{eq:balancing} yields $$\theta_b\theta_XS_{b,X}=d_{b^*X}\theta_{b^*X}=d\theta_X,$$ so that $S_{b,X}=d$.   The second statement follows from the fact that $\eC_{pt}=C_\eC(\eC_0).$
\end{proof}

\subsection{Applications to some examples}

Later in section \ref{intMD}, we use some of the above results and $\SL$
representations to obtain sets of potential quantum dimensions that include all
the integral modular data (see  Table \ref{dLLs}) via an automated computer
calculation. In the following, we will apply the above general results trying
to rule out some of those sets of potential quantum dimensions.  In particular we can deal with all the unrealizable cases in Table \ref{dLLs} directly.

\begin{enumerate}
    
\item{Rank 7}

Consider an MTC $\eC$ with
$\dim(\eC) = 16$, with dimensions partitioned via the grading as
$[1, 1, 1, 1], [2], [2], [2]$ (cf. Table \ref{dLLs}).   This category does not exist by the results of \cite{BGNPRW}, which classifies all rank $7$ integral categories.  To illustrate the power of our methods we will eliminate this case directly, so that our computational techniques plus the following short argument reproduces the main result of \cite{BGNPRW}.
Suppose $\eC$ is such a MTC. Since $\eC_0$ is a symmetric pointed category it must contain an order 2 boson $b$ by Lemma \ref{lem:bosoninorder4}.  Now by Lemma \ref{balancingtrick} $b$ centralizes the $2$ dimensional simple objects, and is thus central in $\eC$.  This contradicts modularity.

\item{Rank 8}

Suppose we have $\eC$ with
$\dim(\eC)= 36$, with dimensions partitioned by gradings being 
$ [1, 2, 2, 3], [1, 2, 2, 3]$.  By Prop. \ref{integralprop}(c) the trivial component must be modular, which does not exist due to \cite{RSW0777}.

\item{Rank 9}

 \begin{itemize}
     
\item Suppose $\eC$ has
$\dim(\eC) = 144$, and with dimensions partitioned via the grading as
$ [1, 1, 1, 1, 4, 4] ,[6] ,[6] ,[6]$ (cf. Table \ref{dLLs}).  This can be eliminated similarly as above: let $b$ the boson obtained from Lemma \ref{lem:bosoninorder4}.  Now by Lemma \ref{balancingtrick} $b$ is central in $\eC$, contradicting modularity.

\item Suppose $\eC$ is a $144$-dimensional MTC with simple objects of dimensions
$$[1,1,3,3,4,6,6,6,6].$$  Since $\eC$ is faithfully $\Z_2$ graded the
dimensions of the simple objects in $\eC_0$ must be $[1,1,3,3,4,6]$.  Now the
non-trivial invertible $Z\in\eC_0$ must be a boson by Prop.
\ref{integralprop}(b).   Now let $X_1$ be a simple object of dimension $4$ and
$X_2$ a simple object of dimension $6$.  We have that $Z\otimes X_i\cong X_i$
for $i=1,2$.  By Prop. \ref{integralprop}(e) we see that $F(X_i)$ sums of 2
non-isomorphic simple objects, and thus we obtain 4 simple objects of dimension
$2,2,3,3$ in the modularization $\eD$. Similarly, the objects $Y_1,Y_2$ of
dimension $3$ in $\eC_0$ must obey $Z\otimes Y_1\cong Y_2$ and so $F(Y_1)\cong
F(Y_2)$ is a simple object of dimension $3$ in $\eD$.  Finally, $F(Z)\cong
F(\mathbf{1})$, so that $\eD$ has simple objects of dimension $1,2,2,3,3,3$.
But no such modular category exists: it has dimension $36$ and so is solvable
but only 1 invertible object--a contradiction.

\item Suppose there were a MTC  with $\dim(\eC) = 288$, and dimensions
partitioned via the grading as: $[1, 1, 3, 3, 4, 6, 6, 6] ,[12]$.  Since there
is a unique object of dimension $4$ in the trivial component, the non-trivial
invertible object $b$ must be a boson.  As the non-trivial component has only
one simple, this contradicts Prop. \ref{integralprop}(d).
 \end{itemize}
\item{rank 10}
 
\begin{itemize}
    
\item Suppose $\eC$ is a MTC with
$\dim(\eC) = 108$ and dimensions partitioned as  $$[1, 1, 2, 2, 2, 2, 6], [3, 3, 6].$$  By Prop. \ref{integralprop}(b) and (c) this category  does not exist.

\item Suppose that $\eC$ is a MTC with dimensions partitions as: $[1,1,1,1,2,2,2],[4],[4],[4]$ as in Table \ref{dLLs}.  By Lemma \ref{lem:bosoninorder4} there must be a boson $b$, which by Lemma \ref{balancingtrick} must centralize all of the simple objects of dimension $4$ as well as $\eC_0$ contradicting modularity.

 \end{itemize} 

 \item{rank 11}
\begin{itemize}
\item 
Suppose $\eC$ is a MTC of dimension $144$ with dimensions partitioned via the grading as:
$[1, 1, 1, 1, 4, 4], [2, 4, 4], [6], [6]$ (cf. Table \ref{dLLs}).  By Lemma \ref{lem:bosoninorder4} we have a boson $b$ in the trivial component.  Since $b\otimes Y\cong Y$ for each of the objects of dimension $6$ they both centralize $b$.  Thus we find that the centralizer $C_{\eC}(\langle b\rangle)$ has dimension at least $3\cdot 6^2$ since $b$ is centralized by $\eC_0$ as well.  But this contradicts Lemma \ref{lem:dimsofsubcat} since $\dim(\langle b\rangle)\cdot C_{\eC}(\langle b\rangle)\geq 6^3>144.$
\item 
Now suppose that $\eC$ has dimension $144$ with dimensions partitioned as:
$$[ 1, 1, 1, 3, 6 ], [ 4, 4, 4 ], [ 4, 4, 4 ]$$ (cf. Table \ref{dLLs}).
Then clearly the trivial component has a Tannakian subcategory equivalent to $\Rep_{\Z_3}$.  By Prop. \ref{integralprop}(e) the modularization of $\eC_0$ is a rank $7$ category of dimension $16$ with dimensions $[1,1,1,1,2,2,2]$ which we have eliminated already.
\end{itemize}
\item{rank 12}

By now the methods are familiar so we quickly eliminate the following (cf. Table \ref{dLLs}):
 \begin{itemize}
     \item  $ [ 1, 1, 1, 1, 2, 2, 2, 2, 4 ], [ 6 ], [ 6 ], [ 6 ] $ cannot occur as the boson afforded by Lemma \ref{lem:bosoninorder4} would be central, by Lemma \ref{balancingtrick}, contradicting modularity.
     \item $[[ 1, 1, 1, 1, 4, 4 ], [ 3, 3, 3, 3 ], [ 6 ], [ 6 ]$ cannot occur as the boson afforded by Lemma \ref{lem:bosoninorder4} has a centralizer that has dimension more than half of the dimension of the category, contradicting Lemma \ref{lem:dimsofsubcat}.
     \item $[  1, 1, 1, 2, 2, 2, 2, 2, 2, 3 ], [ 6 ], [ 6 ]$ cannot occur as Lemma \ref{balancingtrick} implies that the non-trivial invertibles are central, contradicting modularity.
 \end{itemize}

\end{enumerate}

\subsection{Computer calculation of Integral modular data}
\label{intMD}

The methods outlined in Section \ref{sec:pMD} can be used to find both integral
and non-integral modular data from $\SL$ representations. The results are
summarized in Section \ref{MDbyGal}.  However, those methods are less effective
to integral modular data.  For rank 8 and above, we need to use a different and
more effective approach to find integral modular data.  The new approach
contains a few steps.

In the first step, from a list all possible $\SL$ representations that have
potential to produce modular data that we have obtained before, we can obtain a
list of prime divisors $\wp$ of $\pord(\rho(\frt))$ of those $\SL$
representations (see Table \ref{primeLL}).  This information will be used in
the next step.

In the second step, we note that $x_i \equiv D^2/d_i^2,\ i=1,\ldots,r$ are
integers for integral modular data, which satisfy  the Egyptian fraction
condition \eqref{eqn:Egypt}.

\begin{table}[tb] 
\caption{
Sets of prime divisors of $\pord(\rho(\frt))$ of the $\SL$ representations that
have potential to produce modular data.  For each rank, we drop the sets of
prime divisors which are contained in some other sets of prime divisors at that
rank.
} 
\label{primeLL} 
\centering
\begin{tabular}{ |c|c| } 
\hline 
rank & prime divisors  \\
\hline 
2 & [2], [5] \\
\hline 
3 & [2], [3], [7] \\
\hline 
4 & [2,5], [3] \\
\hline 
5 & [2,3], [2,5], [7], [11] \\
\hline 
6 & [2,3], [2,5], [2,7], [3,5], [3,7], [5,7], [13] \\
\hline 
7 & [2,3], [2,5], [2,7], [3,5], [3,7], [11] \\
\hline 
8 & [2,3,5], [2,7], [3,7], [3,11], [13], [17]  \\
\hline 
9 & [2,3,5], [2,3,7], [2,11], [2,13], [3,11], [19] \\
\hline 
10 & [2,3,5], [2,3,7], [2,11], [2,13], [3,11], [3,13], [5,7], [5,11], [5,13], [7,11], [17] \\
\hline 
11 & [2,3,5], [2,3,7], [2,3,11], [2,13], [2,17], [3,5,7], [3,13], [5,13], [7,11], [19], [23] \\
\hline 
12 & [2,3,5], [2,3,7], [2,3,11], [2,3,13], [2,5,7], [2,5,13], [2,17], [3,5,7], [3,7,11], [3,17], [3,19], [5,17]  \\
\hline \end{tabular}
\end{table}

For each rank $r$, the number of the lists of Egyptian fractions is finite.  We
use the following method to obtain all of them.  First we assume $x_i$ are
ordered $x_1 \leq x_2,\leq \cdots \leq x_r$.  Suppose we have known a partial
list $x_1,\ldots,x_k$.  Then the next integer $x_{k+1}$ in the list satisfies $
\frac{1}{x_{k+1}} \geq \frac{1}{r-k}(1-\sum_{i=1}^k \frac{1}{x_i})$.  Thus
\begin{align}
 x_k \leq x_{k+1} \leq \frac{r-k}{1-\sum_{i=1}^k \frac{1}{x_i}} .
\end{align}
We see that the integer $x_{k+1}$ only has a finite range of choices.  Since
here, $x_i$'s are given by the quantum dimensions of modular tensor category:
$x_i=D^2/d_i^2$, the prime divisors of $x_{k+1}$ must be a subset of prime
divisors of $\pord(\rho(\frt))$.  The range of choices can be greatly  reduced
by using this property.  Since the known prime divisors of $\pord(\rho(\frt))$
form small sets (see Table \ref{primeLL}), this allows us to obtain all
possible lists of Egyptian fractions, $1/x_i$'s, for each set of prime divisors
(of $\pord(\rho(\frt))$) in Table \ref{primeLL}.  (Otherwise, the load of
computer calculation is too much for rank 11, 12.)

From $x_i$'s we can get the quantum dimensions $d_i$ via
\begin{align}
 d_i = \sqrt{\frac{x_r}{x_i}},\ \ \ 
D^2 =\sum_i d_i^2 = x_r.
\end{align}
This way we obtain sets of potential quantum dimensions of
all integral modular data.

When $d_i$'s are viewed as quantum dimensions of an integral modular data, they
must satisfy many additional conditions discussed in  Section
\ref{mathintcat}.  We further reduces the sets of potential quantum dimensions
of integral modular data, by using some of those conditions that can be
automated by computer:
\begin{lem}
\begin{itemize}
\item
$d_i =\sqrt{\frac{x_r}{x_i}}$ must an integer.

\item
The number $N_\mathrm{inv}$ of invertible objects (the number of $d_i=1$'s)
divides $x_r=D^2$.  $\eC$ can be decomposed as an Abelian category into $N_\mathrm{inv}$ components, such
that each component $\eC_i$  has the same total quantum dimension
$D^2/N_\mathrm{inv}$. Either $\eC_0$ contains invertible objects, which must be
bosons or fermions, or $\eC_0$ is modular and hence $\eC$ is not prime.

\item
If $\eC_i$ contains invertible objects, then the set of quantum dimensions in
$\eC_i$ is the same as that in $\eC_0$.

\item
Let $N_\mathrm{inv,0}$ be the number of invertibles in $\eC_0$. Let
$P_\mathrm{inv,0}$ be the product of distinct prime factors in $N_\mathrm{inv,0}$.
Then $P_\mathrm{inv,0}$ divides $d_iN_{d_i}$, where $d_i$ is a quantum
dimension in $\eC_0$ and $N_{d_i}$ is the number of simple objects with
dimension $d_i$.  If one of $N_{d_i}$ is odd, then invertible objects in $\eC_0$
must be all bosons, i.e. $Sym(\eC_0)$ is Tannakian. 
\item
If $|\wp|<3$ and $Sym(\eC_0)$ is Tannakian, then $d_\mathrm{min}/N_\mathrm{inv,0} \leq 1$,
where $d_\mathrm{min}$ is the smallest non-unit dimension in $\eC_0$.

\end{itemize}
\end{lem}
\begin{proof}
The first two statements are immediate.  Let us justify the last three.

Let $Z\in\eC_i$ be invertible.  Then the map $X\mapsto Z\cdot X$ is a bijection between $\eC_0$ and $\eC_i$.

Let $p$ be a prime dividing $N_\mathrm{inv,0}$ so that $\Z_p$ acts on the objects of dimension $d_i$. If $p=2$ and the corresponding invertible object is a fermion, $f$ then $f\cdot X\not\cong X$ for each simple $X$ of dimension $d_i$ so that $N_i$ is even.  If $p>2$ the corresponding invertible object $b$ is bosonic.  If $p\nmid d_i$ then $b\cdot X\not\cong X$ since if we condense $b$ the object $X$ would split into $p$ objects of dimension $d_i/p$.  Thus the orbit of $X$ under the $\Z_p$ action must have size $p$, and thus $N_\mathrm{inv,0}$ is divisible by $p$.

The given hypotheses imply that $\eC$ is solvable, and the de-equivariantization of $\eC_0$ by $Sym(\eC_0)\cong\Rep(A)$ is modular and again solvable. Denote this functor by $F:\eC_0\rightarrow (\eC_0)_A$. By solvability $(\eC_0)_A$ has a non-trivial invertible object $Z$. Such an object must appear as a subobject of $F(X)$ for some non-invertible object $X$ of dimension $d$, since $F(W)=\one$ for all invertible objects in $\eC_0$. Using \cite[Prop. 4.4]{Brug} we find that $1=\dim(Z)=\frac{\nu_{Z,X}d}{|Stab_A(X)|}$ where $Stab_A(X)$ is the $X$-stabilizer subgroup of $A$, and $\nu_{Z,X}$ is the multiplicity of $Z$ in $F(X)$.  Now $|A|=N_\mathrm{inv,0}$, so $1=\frac{\nu_{Z,X}d}{|Stab_A(X)|}\geq \frac{d}{N_\mathrm{inv,0}}\geq \frac{d_\mathrm{min}}{N_\mathrm{inv,0}}$.

\end{proof}

In Table \ref{pfacEF}, we list the prime divisors of $D^2$ of the reduced sets
of the potential quantum dimensions.  Note that the  prime divisors of $D^2$
coincide with the prime divisors of $\pord(\rho(\frt))$, for a modular data
(integral or non-integral).  

\begin{table}[tb] 
\caption{
Sets of prime divisors of $D^2$ (and of $\pord(\rho(\frt))$) of the reduced sets
of the potential quantum dimensions for integral modular data.  }
\label{pfacEF} \centering
\begin{tabular}{ |c|c| } 
\hline 
rank & prime divisors  \\
\hline 
2 & [2] \\
\hline 
3 & [3] \\
\hline 
4 & [2] \\
\hline 
5 & [5] \\
\hline 
6 & [2,3] \\
\hline 
7 & [2], [7] \\
\hline 
8 & [2], [2,3]   \\
\hline 
9 & [2], [2,3] \\
\hline 
10 & [2], [2,3], [2,5], [2,3,5] \\
\hline 
11 & [2], [11], [2,3], [2,5], [2,3,5], [2,3,7]  \\
\hline 
12 & [2,3], [2,3,5], [2,3,7], [2,3,11], [2,3,13]  \\
\hline \end{tabular}
\end{table}

In the third step, we go back to $\SL$ representations.  We find all possible
$D_\rho(\si)$'s for a $\SL$ representation $\rho$ and examine all possible
choices of unit index $u$.  If $D_\rho(\si)_{uu} = \pm 1$ for all $\si$'s, then
the $\SL$ representation $\rho$ has a potential to produce integral modular
data, provided that the prime divisors of
$\pord(\rho(\frt))$ is in Table \ref{pfacEF} for the given rank.  This leads
to Table \ref{pfac}, where the possible sets of prime divisors of
$\pord(\rho(\frt))$ are further reduced.  In this calculation, we also obtain a
list of possible $\rho(\frt)$'s for each possible set of prime divisors of
$\pord(\rho(\frt))$ in Table \ref{pfac}.  This information is useful for the
next step of calculation.

In the fourth step, we compute all possible fusion rings (described by fusion
coefficients $N^{ij}_k$) from a set of potential quantum dimensions $d_i$'s,
via the following equation
\begin{align}
\label{ddNd}
 d_i d_j =\sum_{k=1}^r N^{ij}_k d_k .
\end{align}
Since $N^{ij}_k \geq 0$, the number of solutions for the above equation is
finite.  However, for higher ranks, such as for rank 10,11,12, the calculation
load is too much. We need to find ways to make the calculation doable.

One trick is to use the symmetry of $N^{ij}_k$:
\begin{align}
  N^{ij}_{k} = N^{j\bar{k}}_{i}  = N^{\bar{k}i}_{\bar{j}} 
= N^{\bar{i}k}_{j} = N^{\bar{j}\bar{i}}_{\bar{k}} = N^{k\bar{j}}_{i}
= N^{ji}_{k} = N^{\bar{k}j}_{\bar{i}}  = N^{i\bar{k}}_{\bar{j}} 
= N^{k\bar{i}}_{j} = N^{\bar{i}\bar{j}}_{\bar{k} }= N^{\bar{j}k}_{i}
,
\end{align}
to reduce the number of variable of $N^{ij}_k$.  But to use this trick, we need
to consider all possible charge conjugations $i \to \bar{i}$, and solve
possible fusion rings for each choice of charge conjugation.

Eq. \eqref{ddNd} gives a set of linear equations of the variables.  We also
have a set quadratic equations of variables from
\begin{align}
\sum_m N^{ij}_m N^{mk}_l =\sum_n N^{in}_l  N^{jk}_n.
\end{align}

We first solve the linear  equation with minimal number of variables by search
since all the variables are bounded integers.  During the search, we apply the
quadratic equations containing the searching variables to reduce the search.
After solved some variables, we then solve the linear  equation with minimal
number of remaining variables.  Repeating this process, we can solve all the
variables.

The above trick is not enough. We need more tricks.  When the potential quantum
dimensions contain multiple $1$'s, the corresponding modular data will have
multiple invertible objects.  In this case, the corresponding modular tensor
category is graded by an Abelian group whose order is given by the number of
invertible objects. 
For example, $[d_i] = [1, 1, 2, 2, 2, 2, 3, 3]$ is graded
as $\Z_2$: $\eC_0 = [1, 1, 2, 2, 2, 2]$, $\eC_1 =[3, 3]$.  In this case the fusion
between the invertible objects in the trivial component $\eC_0$ is described by an
Abelian group.  The fusion between the components are also described by an
Abelian group.  Those properties can be used to reduce the calculation.

When the trivial component $\eC_0$ contain more than one invertible object,
those invertible objects form a symmetric fusion category.  Let us first assume
all the invertible objects in $\eC_0$ are bosons.  Then we can condense the
bosons $b_p$ whose order is a prime divisor of $N_\mathrm{inv}$, where
$N_\mathrm{inv}$ is the number of invertible objects in $\eC_0$. In other words,
$(b_p)^p$ is the trivial object if $p$ is a  prime divisor of $N_\mathrm{inv}$.  
The group $\Z_p$ acts on the simple objects $X$ of a given fixed dimension $d$ via $X\mapsto b_p\otimes X$.

The condensation of $b_p$ reduces $\eC_0$ to a ribbon category $\eC_0'$ with
several possible dimension arrays $ [d_1',d_2',\cdots ]$.  The dimension in
$\eC_0'$ are obtained in the following way: first, the number of invertible
objects in $\eC_0'$ is at least $ N'_\mathrm{inv}= N_\mathrm{inv}/p$.  If there are
$p$ $d_i$'s in $\eC_0$ with the same value $d'$, then those $p$ $d_i$'s might
be ``combined'' into a single $d'$ in $\eC_0'$.  Also, if a  quantum dimension
$d_i$ in $\eC_0$ is divisible by $p$, then the single $d_i$ in $\eC_0$ can
``split'' into $p$ degenerate $d_{i_j}'$ in $\eC_0'$, where $d_{i_j}' = d_i/p$
for $j=1,\ldots ,p$.  All the quantum dimensions $d_i$ in $\eC_0$ must either
combine or split, as describe above.  Otherwise, the original quantum
dimensions $[d_i]$ do not correspond to any modular tensor category.  We can
condense all $b_p$ for all the prime divisor of $N_\mathrm{inv}$, using the above
method.

At the end, we obtain a list of possible condensation products $\eC_0'$'s.
Then we try to compute the possible fusion rings for each set of
$[d_1',d_2',\cdots]$.  If there is no valid fusion ring for all the possible
condensation products $\eC_0'$'s, then the original quantum dimensions $[d_i]$
does not corresponds to any modular tensor category where the invertible
objects in the trivial components are all bosons.

If the invertible objects in the trivial components contain fermions, then
$N_\mathrm{inv}$ must be even.  If $N_\mathrm{inv}$ is divisible by 4, then $b_2$
is a boson and we can use the approach described above to obtain a list of
possible condensation products $\eC_0'$'s.  But here we further require $d_i'$
in $\eC_0'$ all have an even degeneracy, since the invertible object in
$\eC_0'$ must contain the uncondensed fermions.  If the even $N_\mathrm{inv}$ is
not divisible by 4, then we condensed all the $b_p$'s, except $b_2$, to obtain a
list of possible condensation products.  If there is no valid fusion ring for
all the possible condensation products $\eC_0'$'s, then the original quantum
dimensions $[d_i]$ does not corresponds to any modular tensor category where
the invertible objects in the trivial components contains fermions.  

\begin{table}[tb] 
\caption{ Lists of all potential quantum dimensions for integral modular tensor
categories.  The quantum dimensions are presented as
$[[d_1,d_2,\cdots],[d_1',d_2'\cdots],\cdots]$, where $[d_1,d_2,\cdots]$ form
the trivial component $\eC_0$, $[d_1',d_2',\cdots]$ form the first non-trvial
component $\eC_1$, {\it etc}.  Each set of quantum dimensions has one or more
valid fusion rings, but may not correspond to modular data.  } \label{dLLs}
\centering
\begin{tabular}{ |c|c| } 
\hline 
rank & quantum dimensions \\
\hline 
2 & $ [ 1 ], [ 1 ]$ \\
\hline 
3 & $ [ 1 ], [ 1 ], [ 1 ]$ \\
\hline 
4 & $  [ 1 ], [ 1 ], [ 1 ], [ 1 ]$ \\
\hline 
5 & $ [ 1 ], [ 1 ], [ 1 ], [ 1 ], [ 1 ]$ \\
\hline 
6 & $ [ 1 ], [ 1 ], [ 1 ], [ 1 ], [ 1 ], [ 1 ]$ \\
\hline 
7 & $ [ 1 ], [ 1 ], [ 1 ], [ 1 ], [ 1 ], [ 1 ], [ 1 ]$ \\
&
 $ [ 1, 1, 1, 1 ], [ 2 ], [ 2 ], [ 2 ]$ \\
\hline 
8 & $ [ 1 ], [ 1 ], [ 1 ], [ 1 ], [ 1 ], [ 1 ], [ 1 ], [ 1 ]
$ \\
&
$ [ 1, 1, 2, 2, 2, 2 ], [ 3, 3 ] $   \\
\hline 
9 & $ [ 1 ], [ 1 ], [ 1 ], [ 1 ], [ 1 ], [ 1 ], [ 1 ], [ 1 ], [ 1 ]$ \\
& 
$ [ 1, 1, 1, 1, 4, 4 ], [ 6 ], [ 6 ], [ 6 ]$\\
\hline 
10 & $ [ 1 ], [ 1 ], [ 1 ], [ 1 ], [ 1 ], [ 1 ], [ 1 ], [ 1 ], [ 1 ], [ 1 ]$ \\
&
$ [ 1, 1, 1, 1, 2, 2, 2 ], [ 4 ], [ 4 ], [ 4 ]$ \\
&
$ [ 1, 1, 1, 3 ], [ 2, 2, 2 ], [ 2, 2, 2 ]$ \\
\hline 
11 & $ [ 1 ], [ 1 ], [ 1 ], [ 1 ], [ 1 ], [ 1 ], [ 1 ], [ 1 ], [ 1 ], [ 1 ], [ 1 ]$ \\
&
$ [ 1, 1, 1, 1, 2 ], [ 2, 2 ], [ 2, 2 ], [ 2, 2 ]$ \\
&
$ [ 1, 1, 1, 1, 4, 4 ], [ 2, 4, 4 ], [ 6 ], [ 6 ] $ \\
&
$ [ 1, 1, 1, 3, 6 ], [ 4, 4, 4 ], [ 4, 4, 4 ]$  \\
\hline 
12 & $ [ 1 ], [ 1 ], [ 1 ], [ 1 ], [ 1 ], [ 1 ], [ 1 ], [ 1 ], [ 1 ], [ 1 ], [ 1 ], [ 1 ]$ \\
&
$ [ 1, 1, 1, 1, 2, 2, 2, 2, 4 ], [ 6 ], [ 6 ], [ 6 ]$ \\
&
$ [ 1, 1, 1, 1, 4, 4 ], [ 3, 3, 3, 3 ], [ 6 ], [ 6 ]$ \\
& 
$ [ 1, 1, 1, 2, 2, 2, 2, 2, 2, 3 ], [ 6 ], [ 6 ]$ \\
\hline \end{tabular}
\end{table}

The above condensation consideration is very effective in ruling out many
invalid potential sets of quantum dimensions.  However, if the quantum
dimensions $[d_1,d_2,\cdots]$ contain only a single invertible object (\ie
the unit), the corresponding modular tensor category will be called perfect. In
this case, our above tricks will not apply. For rank 11 and 12, the computation
of the fusion rings from the perfect quantum dimensions is too much for a
current desktop computer to handle.

To solve this problem, we note that for each set of quantum dimensions,
$d_i$'s, we know $D^2$ and its prime divisors.  For such a set of prime
divisors, we know a list of all possible $\rho(\frt)$'s, and the
corresponding twists $\t\th_i$'s, that may give rise to integral modular data.
However, we do not know the matching between $d_i$'s and $\t\th_i$'s.  Thus, we
need to find all permutations in the indices, $p$'s, such that $(d_i,
\t\th_{p(i)})$ are the quantum dimension and the twist of an object in a modular
tensor category.  In particular, the permutation $p$ must satisfy
\begin{align}
 D\ee^{\ii \phi} = \sum_i d_i^2 \t\th_{p(i)}. 
\end{align}
For all possible $\rho(\frt)$'s, we search all possible permutations
$p$ to satisfy the above condition.  We can reduce the search by noticing that
the twists related by Galois conjugation, $\t\th_i \to \si^2 \t\th_i =
\t\th_{\hat \si(i)} $, correspond to the same value of  quantum
dimensions.

If no valid permutations can be found for all possible $\rho(\frt)$'s, then the
set of quantum dimensions, $d_i$'s, does not correspond to any valid integral
modular data.  This condition is very effective in ruling out many perfect
quantum dimensions.  The remaining sets of quantum dimensions that have not
been ruled out are listed in Table \ref{dLLs}.  We manage to obtain all the
fusion rings for each set of quantum dimensions in Table \ref{dLLs}.

In the fifth step, we compute the possible topological spins $s_i$ ($\ee^{\ii
2\pi s_i}$ are eigenvalues of the $T$-matrix) from the obtained fusion ring
$N^{ij}_k$ via\cite{V8821,AM8841,Em0207007,E2009} 
\begin{align}
\label{Vs}
\sum_r V_{ijkl}^r s_r =0 \text{ mod }1
\end{align} 
where
\begin{align}
\ \ \ \ \ \ 
V_{ijkl}^r &=  
N^{ij}_r N^{kl}_{\bar r}+
N^{il}_r N^{jk}_{\bar r}+
N^{ik}_r N^{jl}_{\bar r}
- ( \del_{ir}+ \del_{jr}+ \del_{kr}+ \del_{lr}) \sum_m N^{ij}_m N^{kl}_{\bar m}
.
\end{align}
There are many sets of solutions of $s_i$'s, which can be calculated via the
Smith normal form of above $V$-matrix.  For each set of solution, $s_i$'s, we
can calculate the $S$-matrix via \eqref{eq:balancing}.  We then check if the
resulting $S,T$ form a valid modular data.  This way, we obtain a complete list
of integral modular data, for rank 12 and below.  The results are summarized in
Section \ref{MDbyGal}.  We remark that our results coincide with those of
\cite{palcoux} for rank 13 and below, although our methods are somewhat
different.

\section{Lists of modular data by Galois orbits}
\label{MDbyGal}

By a theorem of Mueger, if $\eC\subset\eD$ are both modular categories then we
have a factorization $\eC\cong\eD\boxtimes \eD^\prime$ where $\eD^\prime$ is
another modular subcategory.  If $\eD$ has no modular subcategories it is said
to be \emph{prime}.  In this section we list all the prime modular data for
ranks 7,8 9, 10, 11 and 12.  To save space, we group the modular data by Galois
orbits generated by Galois conjugations.  We only list one representative for
each Galois orbit.  If a Galois orbit contains unitary modular data (defined by
quantum dimensions $d_i \geq 1$), we will choose the representative to be an
unitary one.   

A grey entry means that the Galois orbit and the previous Galois orbit are
connected by a change of spherical structure (\ie the two Galois orbits each
contain a modular data, such that the two modular data are connected by a
change of spherical structure).  If a Galois orbit contains no unitary modular
data, and one of the non-unitary modular data is pseudo unitary, we will drop
this Galois orbit. Such a Galois orbit is connected to a Galois orbit with
unitary modular data.  If a Galois orbit contains a modular data that is not
prime, we will also drop this Galois orbit.  The resulting lists are given
below.  The list is ordered by $D^2$.

In the list, the $T$-matrix of the modular data is presented as
$(s_0,s_1,\cdots)$ where $T_{ii}=\ee^{\ii 2\pi s_i}$, $0\leq s_i <1$.  The
$S$-matrix is presented as $(S_{11},S_{12},S_{13},\cdots;\ \
S_{22},S_{23},\cdots)$.  The matrix elements of $S$ are given in terms of the
following cyclotomic numbers:
\begin{align}
\zeta^m_n &=\mathrm{e}^{2\pi \mathrm{i} m/n}, \ \ \ \
 c^m_n = \zeta^m_n+\zeta^{-m}_n, \ \ \ \
 s^m_n = \zeta^m_n-\zeta^{-m}_n, 
\nonumber\\
 \xi^{m}_n &= \xi^{m,1}_n, \ \ \
 \eta^{m}_n = \eta^{m,1}_n,\ \ \
 \chi^{m}_n = \chi^{m,1}_n,\ \ \
 \la^{m}_n = \la^{m,1}_n,
\nonumber\\
 \xi^{m,l}_n &=
 (\zeta^m_{2n}-\zeta^{-m}_{2n})/(\zeta_{2n}^l-\zeta_{2n}^{-l}), \ \ \ 
 \eta^{m,l}_n =
 (\zeta^m_{2n}+\zeta^{-m}_{2n})/(\zeta_{2n}^l+\zeta_{2n}^{-l}),
\nonumber\\
 \chi^{m,l}_n &=
 (\zeta^m_{2n}+\zeta^{-m}_{2n})/(\zeta_{2n}^l-\zeta_{2n}^{-l}), \ \ \ 
 \la^{m,l}_n =
 (\zeta^m_{2n}-\zeta^{-m}_{2n})/(\zeta_{2n}^l+\zeta_{2n}^{-l}).
\end{align}

Each modular data in the list is labeled by $r_{c,D^2}^{\ord(T),\text{fp}}$.
For example, $7_{2,7.}^{7,892}$ labels a  modular tensor category with rank $r=
7$, chiral central charge $c=2$, total quantum dimension $D^2=7.0$, order-$T$
$\ord(T) = 7$, and finger print fp $=892$.  Here the ``finger print'' is given
by the first three digits of $|\sum_i (s_i^2-\frac14) d_i|$, so that distinct
modular tensor categories are more likely to have distinct labels.

In the table, we also list the realizations of each modular data.  Usually, a
realization is given by the modular tensor category of Kac-Moody algebra.  For
example, $SU(5)_{5}$ is the  modular tensor category of  $SU(5)$ level 5
Kac-Moody algebra.  $PSU(3)_{5}$ is the  modular tensor category that is the
non-pointed Deligne factor of $SU(3)_{5}$, {\it i.e.} $SU(3)_{5} = PSU(3)_{5}
\boxtimes \cC(\Z_5, q)$.  We also use $O_n$ to represents the modular tensor
category of the $U(1)_{2n}$ orbifold \cite{DVV8985}.  The  modular tensor
category from twisted quantum double is labeled by $\eD^\om(G)$, where $G$ is a
finite group and $\om$ in the cocycle twist.  The  modular tensor category from
twisted Haagerup-Izumi modular data is labeled by Haag$(n)_m$, $m=0,\pm
1,\cdots, \pm n$\cite{EG10061326}.  Also many modular data are realized as
Abelian anyon condensations \cite{LW170107820} of the modular tensor category
from Kac-Moody algebra and/or twisted quantum double.  Two constructions
closely related to Abelian anyon condensations called \emph{zesting}
\cite{16foldway,DGPRZ} and the \emph{condensed fiber product} \cite{CFP} are
also useful.  Most of the potential modular data in the lists are realized by
modular tensor categories, and are indeed modular data.  There are a few
potential modular data whose realizations are not known or not sure, which will
be discussed on Section \ref{unknownrealization}.


\subsection{Rank 7}\label{ss:rank7}

{\small

\noindent1. $7_{2,7.}^{7,892}$ \irep{48}:\ \ 
$d_i$ = ($1.0$,
$1.0$,
$1.0$,
$1.0$,
$1.0$,
$1.0$,
$1.0$) 

\vskip 0.7ex
\hangindent=3em \hangafter=1
$D^2= 7.0 = 
7$

\vskip 0.7ex
\hangindent=3em \hangafter=1
$T = ( 0,
\frac{1}{7},
\frac{1}{7},
\frac{2}{7},
\frac{2}{7},
\frac{4}{7},
\frac{4}{7} )
$,

\vskip 0.7ex
\hangindent=3em \hangafter=1
$S$ = ($ 1$,
$ 1$,
$ 1$,
$ 1$,
$ 1$,
$ 1$,
$ 1$;\ \ 
$ -\zeta_{14}^{3}$,
$ \zeta_{7}^{2}$,
$ -\zeta_{14}^{5}$,
$ \zeta_{7}^{1}$,
$ -\zeta_{14}^{1}$,
$ \zeta_{7}^{3}$;\ \ 
$ -\zeta_{14}^{3}$,
$ \zeta_{7}^{1}$,
$ -\zeta_{14}^{5}$,
$ \zeta_{7}^{3}$,
$ -\zeta_{14}^{1}$;\ \ 
$ \zeta_{7}^{3}$,
$ -\zeta_{14}^{1}$,
$ \zeta_{7}^{2}$,
$ -\zeta_{14}^{3}$;\ \ 
$ \zeta_{7}^{3}$,
$ -\zeta_{14}^{3}$,
$ \zeta_{7}^{2}$;\ \ 
$ -\zeta_{14}^{5}$,
$ \zeta_{7}^{1}$;\ \ 
$ -\zeta_{14}^{5}$)

Realization: $U(7)_1$.

\vskip 1ex

\noindent2. $7_{\frac{27}{4},27.31}^{32,396}$ \irep{96}:\ \ 
$d_i$ = ($1.0$,
$1.0$,
$1.847$,
$1.847$,
$2.414$,
$2.414$,
$2.613$) 

\vskip 0.7ex
\hangindent=3em \hangafter=1
$D^2= 27.313 = 
16+8\sqrt{2}$

\vskip 0.7ex
\hangindent=3em \hangafter=1
$T = ( 0,
\frac{1}{2},
\frac{1}{32},
\frac{1}{32},
\frac{1}{4},
\frac{3}{4},
\frac{21}{32} )
$,

\vskip 0.7ex
\hangindent=3em \hangafter=1
$S$ = ($ 1$,
$ 1$,
$ c_{16}^{1}$,
$ c_{16}^{1}$,
$ 1+\sqrt{2}$,
$ 1+\sqrt{2}$,
$ c^{1}_{16}
+c^{3}_{16}
$;\ \ 
$ 1$,
$ -c_{16}^{1}$,
$ -c_{16}^{1}$,
$ 1+\sqrt{2}$,
$ 1+\sqrt{2}$,
$ -c^{1}_{16}
-c^{3}_{16}
$;\ \ 
$(-c^{1}_{16}
-c^{3}_{16}
)\mathrm{i}$,
$(c^{1}_{16}
+c^{3}_{16}
)\mathrm{i}$,
$ -c_{16}^{1}$,
$ c_{16}^{1}$,
$0$;\ \ 
$(-c^{1}_{16}
-c^{3}_{16}
)\mathrm{i}$,
$ -c_{16}^{1}$,
$ c_{16}^{1}$,
$0$;\ \ 
$ -1$,
$ -1$,
$ c^{1}_{16}
+c^{3}_{16}
$;\ \ 
$ -1$,
$ -c^{1}_{16}
-c^{3}_{16}
$;\ \ 
$0$)

Realization:
Abelian anyon condensation of
$SU(6)_2$ or
$Sp(12)_1$ or
$\overline{SU(2)}_6$.   One of 4 $\Z_2$-zestings of the spin modular category $\overline{SU(2)}_6$, see \cite{16fold}.

\vskip 1ex

\noindent3. $7_{\frac{9}{4},27.31}^{32,918}$ \irep{97}:\ \ 
$d_i$ = ($1.0$,
$1.0$,
$1.847$,
$1.847$,
$2.414$,
$2.414$,
$2.613$) 

\vskip 0.7ex
\hangindent=3em \hangafter=1
$D^2= 27.313 = 
16+8\sqrt{2}$

\vskip 0.7ex
\hangindent=3em \hangafter=1
$T = ( 0,
\frac{1}{2},
\frac{3}{32},
\frac{3}{32},
\frac{1}{4},
\frac{3}{4},
\frac{15}{32} )
$,

\vskip 0.7ex
\hangindent=3em \hangafter=1
$S$ = ($ 1$,
$ 1$,
$ c_{16}^{1}$,
$ c_{16}^{1}$,
$ 1+\sqrt{2}$,
$ 1+\sqrt{2}$,
$ c^{1}_{16}
+c^{3}_{16}
$;\ \ 
$ 1$,
$ -c_{16}^{1}$,
$ -c_{16}^{1}$,
$ 1+\sqrt{2}$,
$ 1+\sqrt{2}$,
$ -c^{1}_{16}
-c^{3}_{16}
$;\ \ 
$ c^{1}_{16}
+c^{3}_{16}
$,
$ -c^{1}_{16}
-c^{3}_{16}
$,
$ c_{16}^{1}$,
$ -c_{16}^{1}$,
$0$;\ \ 
$ c^{1}_{16}
+c^{3}_{16}
$,
$ c_{16}^{1}$,
$ -c_{16}^{1}$,
$0$;\ \ 
$ -1$,
$ -1$,
$ -c^{1}_{16}
-c^{3}_{16}
$;\ \ 
$ -1$,
$ c^{1}_{16}
+c^{3}_{16}
$;\ \ 
$0$)

Realization:
Abelian anyon condensation of
$SU(2)_6$ or
$\overline{SU(6)}_2$ or
$\overline{Sp(12)}_1$.  One of 4 $\Z_2$-zestings of the spin modular category $\overline{SU(2)}_6$, see \cite{16fold}.

\vskip 1ex

\noindent4. $7_{\frac{31}{4},27.31}^{32,159}$ \irep{97}:\ \ 
$d_i$ = ($1.0$,
$1.0$,
$1.847$,
$1.847$,
$2.414$,
$2.414$,
$2.613$) 

\vskip 0.7ex
\hangindent=3em \hangafter=1
$D^2= 27.313 = 
16+8\sqrt{2}$

\vskip 0.7ex
\hangindent=3em \hangafter=1
$T = ( 0,
\frac{1}{2},
\frac{5}{32},
\frac{5}{32},
\frac{1}{4},
\frac{3}{4},
\frac{25}{32} )
$,

\vskip 0.7ex
\hangindent=3em \hangafter=1
$S$ = ($ 1$,
$ 1$,
$ c_{16}^{1}$,
$ c_{16}^{1}$,
$ 1+\sqrt{2}$,
$ 1+\sqrt{2}$,
$ c^{1}_{16}
+c^{3}_{16}
$;\ \ 
$ 1$,
$ -c_{16}^{1}$,
$ -c_{16}^{1}$,
$ 1+\sqrt{2}$,
$ 1+\sqrt{2}$,
$ -c^{1}_{16}
-c^{3}_{16}
$;\ \ 
$ -c^{1}_{16}
-c^{3}_{16}
$,
$ c^{1}_{16}
+c^{3}_{16}
$,
$ -c_{16}^{1}$,
$ c_{16}^{1}$,
$0$;\ \ 
$ -c^{1}_{16}
-c^{3}_{16}
$,
$ -c_{16}^{1}$,
$ c_{16}^{1}$,
$0$;\ \ 
$ -1$,
$ -1$,
$ c^{1}_{16}
+c^{3}_{16}
$;\ \ 
$ -1$,
$ -c^{1}_{16}
-c^{3}_{16}
$;\ \ 
$0$)

Realization:
Abelian anyon condensation of
$SU(6)_2$ or
$Sp(12)_1$ or
$\overline{SU(2)}_6$.  One of 4 $\Z_2$-zestings of the spin modular category $\overline{SU(2)}_6$, see \cite{16fold}.

\vskip 1ex

\noindent5. $7_{\frac{13}{4},27.31}^{32,427}$ \irep{96}:\ \ 
$d_i$ = ($1.0$,
$1.0$,
$1.847$,
$1.847$,
$2.414$,
$2.414$,
$2.613$) 

\vskip 0.7ex
\hangindent=3em \hangafter=1
$D^2= 27.313 = 
16+8\sqrt{2}$

\vskip 0.7ex
\hangindent=3em \hangafter=1
$T = ( 0,
\frac{1}{2},
\frac{7}{32},
\frac{7}{32},
\frac{1}{4},
\frac{3}{4},
\frac{19}{32} )
$,

\vskip 0.7ex
\hangindent=3em \hangafter=1
$S$ = ($ 1$,
$ 1$,
$ c_{16}^{1}$,
$ c_{16}^{1}$,
$ 1+\sqrt{2}$,
$ 1+\sqrt{2}$,
$ c^{1}_{16}
+c^{3}_{16}
$;\ \ 
$ 1$,
$ -c_{16}^{1}$,
$ -c_{16}^{1}$,
$ 1+\sqrt{2}$,
$ 1+\sqrt{2}$,
$ -c^{1}_{16}
-c^{3}_{16}
$;\ \ 
$(-c^{1}_{16}
-c^{3}_{16}
)\mathrm{i}$,
$(c^{1}_{16}
+c^{3}_{16}
)\mathrm{i}$,
$ c_{16}^{1}$,
$ -c_{16}^{1}$,
$0$;\ \ 
$(-c^{1}_{16}
-c^{3}_{16}
)\mathrm{i}$,
$ c_{16}^{1}$,
$ -c_{16}^{1}$,
$0$;\ \ 
$ -1$,
$ -1$,
$ -c^{1}_{16}
-c^{3}_{16}
$;\ \ 
$ -1$,
$ c^{1}_{16}
+c^{3}_{16}
$;\ \ 
$0$)

Realization:
Abelian anyon condensation of
$SU(2)_6$ or
$\overline{SU(6)}_2$ or
$\overline{Sp(12)}_1$.  One of 4 $\Z_2$-zestings of the spin modular category $\overline{SU(2)}_6$, see \cite{16fold}.

\vskip 1ex

\noindent6. $7_{2,28.}^{56,139}$ \irep{99}:\ \ 
$d_i$ = ($1.0$,
$1.0$,
$2.0$,
$2.0$,
$2.0$,
$2.645$,
$2.645$) 

\vskip 0.7ex
\hangindent=3em \hangafter=1
$D^2= 28.0 = 
28$

\vskip 0.7ex
\hangindent=3em \hangafter=1
$T = ( 0,
0,
\frac{1}{7},
\frac{2}{7},
\frac{4}{7},
\frac{1}{8},
\frac{5}{8} )
$,

\vskip 0.7ex
\hangindent=3em \hangafter=1
$S$ = ($ 1$,
$ 1$,
$ 2$,
$ 2$,
$ 2$,
$ \sqrt{7}$,
$ \sqrt{7}$;\ \ 
$ 1$,
$ 2$,
$ 2$,
$ 2$,
$ -\sqrt{7}$,
$ -\sqrt{7}$;\ \ 
$ 2c_{7}^{2}$,
$ 2c_{7}^{1}$,
$ 2c_{7}^{3}$,
$0$,
$0$;\ \ 
$ 2c_{7}^{3}$,
$ 2c_{7}^{2}$,
$0$,
$0$;\ \ 
$ 2c_{7}^{1}$,
$0$,
$0$;\ \ 
$ \sqrt{7}$,
$ -\sqrt{7}$;\ \ 
$ \sqrt{7}$)

Realization:
$\overline{SO(7)}_2$ or
Abelian anyon condensation of
$O_7$.

\vskip 1ex

\noindent7. $7_{2,28.}^{56,680}$ \irep{99}:\ \ 
$d_i$ = ($1.0$,
$1.0$,
$2.0$,
$2.0$,
$2.0$,
$2.645$,
$2.645$) 

\vskip 0.7ex
\hangindent=3em \hangafter=1
$D^2= 28.0 = 
28$

\vskip 0.7ex
\hangindent=3em \hangafter=1
$T = ( 0,
0,
\frac{1}{7},
\frac{2}{7},
\frac{4}{7},
\frac{3}{8},
\frac{7}{8} )
$,

\vskip 0.7ex
\hangindent=3em \hangafter=1
$S$ = ($ 1$,
$ 1$,
$ 2$,
$ 2$,
$ 2$,
$ \sqrt{7}$,
$ \sqrt{7}$;\ \ 
$ 1$,
$ 2$,
$ 2$,
$ 2$,
$ -\sqrt{7}$,
$ -\sqrt{7}$;\ \ 
$ 2c_{7}^{2}$,
$ 2c_{7}^{1}$,
$ 2c_{7}^{3}$,
$0$,
$0$;\ \ 
$ 2c_{7}^{3}$,
$ 2c_{7}^{2}$,
$0$,
$0$;\ \ 
$ 2c_{7}^{1}$,
$0$,
$0$;\ \ 
$ -\sqrt{7}$,
$ \sqrt{7}$;\ \ 
$ -\sqrt{7}$)

Realization:
Abelian anyon condensation or $\Z_2$-zesting of
$SO(7)_2$ or
$\overline{O}_7$.

\vskip 1ex

\noindent8. $7_{\frac{32}{5},86.75}^{15,205}$ \irep{80}:\ \ 
$d_i$ = ($1.0$,
$1.956$,
$2.827$,
$3.574$,
$4.165$,
$4.574$,
$4.783$) 

\vskip 0.7ex
\hangindent=3em \hangafter=1
$D^2= 86.750 = 
30+15c^{1}_{15}
+15c^{2}_{15}
+15c^{3}_{15}
$

\vskip 0.7ex
\hangindent=3em \hangafter=1
$T = ( 0,
\frac{1}{5},
\frac{13}{15},
0,
\frac{3}{5},
\frac{2}{3},
\frac{1}{5} )
$,

\vskip 0.7ex
\hangindent=3em \hangafter=1
$S$ = ($ 1$,
$ -c_{15}^{7}$,
$ \xi_{15}^{3}$,
$ \xi_{15}^{11}$,
$ \xi_{15}^{5}$,
$ \xi_{15}^{9}$,
$ \xi_{15}^{7}$;\ \ 
$ -\xi_{15}^{11}$,
$ \xi_{15}^{9}$,
$ -\xi_{15}^{7}$,
$ \xi_{15}^{5}$,
$ -\xi_{15}^{3}$,
$ 1$;\ \ 
$ \xi_{15}^{9}$,
$ \xi_{15}^{3}$,
$0$,
$ -\xi_{15}^{3}$,
$ -\xi_{15}^{9}$;\ \ 
$ 1$,
$ -\xi_{15}^{5}$,
$ \xi_{15}^{9}$,
$ c_{15}^{7}$;\ \ 
$ -\xi_{15}^{5}$,
$0$,
$ \xi_{15}^{5}$;\ \ 
$ -\xi_{15}^{9}$,
$ \xi_{15}^{3}$;\ \ 
$ -\xi_{15}^{11}$)

Realization: $PSU(2)_{13}$.

\vskip 1ex

\noindent9. $7_{1,93.25}^{8,230}$ \irep{54}:\ \ 
$d_i$ = ($1.0$,
$2.414$,
$2.414$,
$3.414$,
$3.414$,
$4.828$,
$5.828$) 

\vskip 0.7ex
\hangindent=3em \hangafter=1
$D^2= 93.254 = 
48+32\sqrt{2}$

\vskip 0.7ex
\hangindent=3em \hangafter=1
$T = ( 0,
\frac{1}{2},
\frac{1}{2},
\frac{1}{4},
\frac{1}{4},
\frac{5}{8},
0 )
$,

\vskip 0.7ex
\hangindent=3em \hangafter=1
$S$ = ($ 1$,
$ 1+\sqrt{2}$,
$ 1+\sqrt{2}$,
$ 2+\sqrt{2}$,
$ 2+\sqrt{2}$,
$ 2+2\sqrt{2}$,
$ 3+2\sqrt{2}$;\ \ 
$ -1-2  \zeta^{1}_{8}
-2  \zeta^{2}_{8}
$,
$ -1-2  \zeta^{-1}_{8}
+2\zeta^{2}_{8}
$,
$(-2-\sqrt{2})\mathrm{i}$,
$(2+\sqrt{2})\mathrm{i}$,
$ 2+2\sqrt{2}$,
$ -1-\sqrt{2}$;\ \ 
$ -1-2  \zeta^{1}_{8}
-2  \zeta^{2}_{8}
$,
$(2+\sqrt{2})\mathrm{i}$,
$(-2-\sqrt{2})\mathrm{i}$,
$ 2+2\sqrt{2}$,
$ -1-\sqrt{2}$;\ \ 
$ (2+2\sqrt{2})\zeta_{8}^{3}$,
$ (-2-2\sqrt{2})\zeta_{8}^{1}$,
$0$,
$ 2+\sqrt{2}$;\ \ 
$ (2+2\sqrt{2})\zeta_{8}^{3}$,
$0$,
$ 2+\sqrt{2}$;\ \ 
$0$,
$ -2-2\sqrt{2}$;\ \ 
$ 1$)

Realization:
$PSU(3)_5$.

\vskip 1ex

\noindent10. $7_{\frac{30}{11},135.7}^{11,157}$ \irep{66}:\ \ 
$d_i$ = ($1.0$,
$2.918$,
$3.513$,
$3.513$,
$4.601$,
$5.911$,
$6.742$) 

\vskip 0.7ex
\hangindent=3em \hangafter=1
$D^2= 135.778 = 
55+44c^{1}_{11}
+33c^{2}_{11}
+22c^{3}_{11}
+11c^{4}_{11}
$

\vskip 0.7ex
\hangindent=3em \hangafter=1
$T = ( 0,
\frac{1}{11},
\frac{4}{11},
\frac{4}{11},
\frac{3}{11},
\frac{6}{11},
\frac{10}{11} )
$,

\vskip 0.7ex
\hangindent=3em \hangafter=1
$S$ = ($ 1$,
$ 2+c^{1}_{11}
+c^{2}_{11}
+c^{3}_{11}
+c^{4}_{11}
$,
$ \xi_{11}^{5}$,
$ \xi_{11}^{5}$,
$ 2+2c^{1}_{11}
+c^{2}_{11}
+c^{3}_{11}
+c^{4}_{11}
$,
$ 2+2c^{1}_{11}
+c^{2}_{11}
+c^{3}_{11}
$,
$ 2+2c^{1}_{11}
+2c^{2}_{11}
+c^{3}_{11}
$;\ \ 
$ 2+2c^{1}_{11}
+2c^{2}_{11}
+c^{3}_{11}
$,
$ -\xi_{11}^{5}$,
$ -\xi_{11}^{5}$,
$ 2+2c^{1}_{11}
+c^{2}_{11}
+c^{3}_{11}
$,
$ 1$,
$ -2-2  c^{1}_{11}
-c^{2}_{11}
-c^{3}_{11}
-c^{4}_{11}
$;\ \ 
$ s^{2}_{11}
+2\zeta^{3}_{11}
-\zeta^{-3}_{11}
+\zeta^{4}_{11}
+\zeta^{5}_{11}
$,
$ -1-c^{1}_{11}
-2  \zeta^{2}_{11}
-2  \zeta^{3}_{11}
+\zeta^{-3}_{11}
-\zeta^{4}_{11}
-\zeta^{5}_{11}
$,
$ \xi_{11}^{5}$,
$ -\xi_{11}^{5}$,
$ \xi_{11}^{5}$;\ \ 
$ s^{2}_{11}
+2\zeta^{3}_{11}
-\zeta^{-3}_{11}
+\zeta^{4}_{11}
+\zeta^{5}_{11}
$,
$ \xi_{11}^{5}$,
$ -\xi_{11}^{5}$,
$ \xi_{11}^{5}$;\ \ 
$ -2-c^{1}_{11}
-c^{2}_{11}
-c^{3}_{11}
-c^{4}_{11}
$,
$ -2-2  c^{1}_{11}
-2  c^{2}_{11}
-c^{3}_{11}
$,
$ 1$;\ \ 
$ 2+2c^{1}_{11}
+c^{2}_{11}
+c^{3}_{11}
+c^{4}_{11}
$,
$ 2+c^{1}_{11}
+c^{2}_{11}
+c^{3}_{11}
+c^{4}_{11}
$;\ \ 
$ -2-2  c^{1}_{11}
-c^{2}_{11}
-c^{3}_{11}
$)

Realization:
$PSO(10)_3$

\vskip 1ex 

}

\subsection{Rank 8}\label{ss:rank8}

{\small

\noindent1. $8_{1,8.}^{16,123}$ \irep{1000000}:\ \ 
$d_i$ = ($1.0$,
$1.0$,
$1.0$,
$1.0$,
$1.0$,
$1.0$,
$1.0$,
$1.0$) 

\vskip 0.7ex
\hangindent=3em \hangafter=1
$D^2= 8.0 = 
8$

\vskip 0.7ex
\hangindent=3em \hangafter=1
$T = ( 0,
0,
\frac{1}{4},
\frac{1}{4},
\frac{1}{16},
\frac{1}{16},
\frac{9}{16},
\frac{9}{16} )
$,

\vskip 0.7ex
\hangindent=3em \hangafter=1
$S$ = ($ 1$,
$ 1$,
$ 1$,
$ 1$,
$ 1$,
$ 1$,
$ 1$,
$ 1$;\ \ 
$ 1$,
$ 1$,
$ 1$,
$ -1$,
$ -1$,
$ -1$,
$ -1$;\ \ 
$ -1$,
$ -1$,
$-\mathrm{i}$,
$\mathrm{i}$,
$-\mathrm{i}$,
$\mathrm{i}$;\ \ 
$ -1$,
$\mathrm{i}$,
$-\mathrm{i}$,
$\mathrm{i}$,
$-\mathrm{i}$;\ \ 
$ -\zeta_{8}^{3}$,
$ \zeta_{8}^{1}$,
$ \zeta_{8}^{3}$,
$ -\zeta_{8}^{1}$;\ \ 
$ -\zeta_{8}^{3}$,
$ -\zeta_{8}^{1}$,
$ \zeta_{8}^{3}$;\ \ 
$ -\zeta_{8}^{3}$,
$ \zeta_{8}^{1}$;\ \ 
$ -\zeta_{8}^{3}$)

Realization: $U(8)_1$.

\vskip 1ex

\noindent2. $8_{0,36.}^{6,213}$ \irep{1000000}:\ \ 
$d_i$ = ($1.0$,
$1.0$,
$2.0$,
$2.0$,
$2.0$,
$2.0$,
$3.0$,
$3.0$) 

\vskip 0.7ex
\hangindent=3em \hangafter=1
$D^2= 36.0 = 
36$

\vskip 0.7ex
\hangindent=3em \hangafter=1
$T = ( 0,
0,
0,
0,
\frac{1}{3},
\frac{2}{3},
0,
\frac{1}{2} )
$,

\vskip 0.7ex
\hangindent=3em \hangafter=1
$S$ = ($ 1$,
$ 1$,
$ 2$,
$ 2$,
$ 2$,
$ 2$,
$ 3$,
$ 3$;\ \ 
$ 1$,
$ 2$,
$ 2$,
$ 2$,
$ 2$,
$ -3$,
$ -3$;\ \ 
$ 4$,
$ -2$,
$ -2$,
$ -2$,
$0$,
$0$;\ \ 
$ 4$,
$ -2$,
$ -2$,
$0$,
$0$;\ \ 
$ -2$,
$ 4$,
$0$,
$0$;\ \ 
$ -2$,
$0$,
$0$;\ \ 
$ 3$,
$ -3$;\ \ 
$ 3$)

Realization: $\eD(S_3)$

\vskip 1ex

\noindent3. $8_{4,36.}^{6,102}$ \irep{1000000}:\ \ 
$d_i$ = ($1.0$,
$1.0$,
$2.0$,
$2.0$,
$2.0$,
$2.0$,
$3.0$,
$3.0$) 

\vskip 0.7ex
\hangindent=3em \hangafter=1
$D^2= 36.0 = 
36$

\vskip 0.7ex
\hangindent=3em \hangafter=1
$T = ( 0,
0,
\frac{1}{3},
\frac{1}{3},
\frac{2}{3},
\frac{2}{3},
0,
\frac{1}{2} )
$,

\vskip 0.7ex
\hangindent=3em \hangafter=1
$S$ = ($ 1$,
$ 1$,
$ 2$,
$ 2$,
$ 2$,
$ 2$,
$ 3$,
$ 3$;\ \ 
$ 1$,
$ 2$,
$ 2$,
$ 2$,
$ 2$,
$ -3$,
$ -3$;\ \ 
$ -2$,
$ 4$,
$ -2$,
$ -2$,
$0$,
$0$;\ \ 
$ -2$,
$ -2$,
$ -2$,
$0$,
$0$;\ \ 
$ -2$,
$ 4$,
$0$,
$0$;\ \ 
$ -2$,
$0$,
$0$;\ \ 
$ -3$,
$ 3$;\ \ 
$ -3$)

Realization: condensation reductions of 
$\eZ(\cNG(\Z_3\times\Z_3,0))$ (non-group-theoretical, \cite{GNN}).

\vskip 1ex

\noindent4. $8_{0,36.}^{12,101}$ \irep{1000000}:\ \ 
$d_i$ = ($1.0$,
$1.0$,
$2.0$,
$2.0$,
$2.0$,
$2.0$,
$3.0$,
$3.0$) 

\vskip 0.7ex
\hangindent=3em \hangafter=1
$D^2= 36.0 = 
36$

\vskip 0.7ex
\hangindent=3em \hangafter=1
$T = ( 0,
0,
0,
0,
\frac{1}{3},
\frac{2}{3},
\frac{1}{4},
\frac{3}{4} )
$,

\vskip 0.7ex
\hangindent=3em \hangafter=1
$S$ = ($ 1$,
$ 1$,
$ 2$,
$ 2$,
$ 2$,
$ 2$,
$ 3$,
$ 3$;\ \ 
$ 1$,
$ 2$,
$ 2$,
$ 2$,
$ 2$,
$ -3$,
$ -3$;\ \ 
$ 4$,
$ -2$,
$ -2$,
$ -2$,
$0$,
$0$;\ \ 
$ 4$,
$ -2$,
$ -2$,
$0$,
$0$;\ \ 
$ -2$,
$ 4$,
$0$,
$0$;\ \ 
$ -2$,
$0$,
$0$;\ \ 
$ -3$,
$ 3$;\ \ 
$ -3$)

Realization: $\eD^3(S_3)$.

\vskip 1ex

\noindent5. $8_{4,36.}^{12,972}$ \irep{1000000}:\ \ 
$d_i$ = ($1.0$,
$1.0$,
$2.0$,
$2.0$,
$2.0$,
$2.0$,
$3.0$,
$3.0$) 

\vskip 0.7ex
\hangindent=3em \hangafter=1
$D^2= 36.0 = 
36$

\vskip 0.7ex
\hangindent=3em \hangafter=1
$T = ( 0,
0,
\frac{1}{3},
\frac{1}{3},
\frac{2}{3},
\frac{2}{3},
\frac{1}{4},
\frac{3}{4} )
$,

\vskip 0.7ex
\hangindent=3em \hangafter=1
$S$ = ($ 1$,
$ 1$,
$ 2$,
$ 2$,
$ 2$,
$ 2$,
$ 3$,
$ 3$;\ \ 
$ 1$,
$ 2$,
$ 2$,
$ 2$,
$ 2$,
$ -3$,
$ -3$;\ \ 
$ -2$,
$ 4$,
$ -2$,
$ -2$,
$0$,
$0$;\ \ 
$ -2$,
$ -2$,
$ -2$,
$0$,
$0$;\ \ 
$ -2$,
$ 4$,
$0$,
$0$;\ \ 
$ -2$,
$0$,
$0$;\ \ 
$ 3$,
$ -3$;\ \ 
$ 3$)

Realization: condensation reductions of 
$\eZ(\cNG(\Z_3\times\Z_3,0))$ (non-group-theoretical, \cite{GNN}).

\vskip 1ex

\noindent6. $8_{0,36.}^{18,162}$ \irep{1000000}:\ \ 
$d_i$ = ($1.0$,
$1.0$,
$2.0$,
$2.0$,
$2.0$,
$2.0$,
$3.0$,
$3.0$) 

\vskip 0.7ex
\hangindent=3em \hangafter=1
$D^2= 36.0 = 
36$

\vskip 0.7ex
\hangindent=3em \hangafter=1
$T = ( 0,
0,
0,
\frac{1}{9},
\frac{4}{9},
\frac{7}{9},
0,
\frac{1}{2} )
$,

\vskip 0.7ex
\hangindent=3em \hangafter=1
$S$ = ($ 1$,
$ 1$,
$ 2$,
$ 2$,
$ 2$,
$ 2$,
$ 3$,
$ 3$;\ \ 
$ 1$,
$ 2$,
$ 2$,
$ 2$,
$ 2$,
$ -3$,
$ -3$;\ \ 
$ 4$,
$ -2$,
$ -2$,
$ -2$,
$0$,
$0$;\ \ 
$ 2c_{9}^{2}$,
$ 2c_{9}^{4}$,
$ 2c_{9}^{1}$,
$0$,
$0$;\ \ 
$ 2c_{9}^{1}$,
$ 2c_{9}^{2}$,
$0$,
$0$;\ \ 
$ 2c_{9}^{4}$,
$0$,
$0$;\ \ 
$ 3$,
$ -3$;\ \ 
$ 3$)

Realization: $\eD^4(S_3)$ or $SO(9)_2$.

\vskip 1ex

\noindent7. $8_{0,36.}^{36,495}$ \irep{1000000}:\ \ 
$d_i$ = ($1.0$,
$1.0$,
$2.0$,
$2.0$,
$2.0$,
$2.0$,
$3.0$,
$3.0$) 

\vskip 0.7ex
\hangindent=3em \hangafter=1
$D^2= 36.0 = 
36$

\vskip 0.7ex
\hangindent=3em \hangafter=1
$T = ( 0,
0,
0,
\frac{1}{9},
\frac{4}{9},
\frac{7}{9},
\frac{1}{4},
\frac{3}{4} )
$,

\vskip 0.7ex
\hangindent=3em \hangafter=1
$S$ = ($ 1$,
$ 1$,
$ 2$,
$ 2$,
$ 2$,
$ 2$,
$ 3$,
$ 3$;\ \ 
$ 1$,
$ 2$,
$ 2$,
$ 2$,
$ 2$,
$ -3$,
$ -3$;\ \ 
$ 4$,
$ -2$,
$ -2$,
$ -2$,
$0$,
$0$;\ \ 
$ 2c_{9}^{2}$,
$ 2c_{9}^{4}$,
$ 2c_{9}^{1}$,
$0$,
$0$;\ \ 
$ 2c_{9}^{1}$,
$ 2c_{9}^{2}$,
$0$,
$0$;\ \ 
$ 2c_{9}^{4}$,
$0$,
$0$;\ \ 
$ -3$,
$ 3$;\ \ 
$ -3$)

Realization: $\eD^1(S_3)$.

\vskip 1ex

\noindent8. $8_{\frac{62}{17},125.8}^{17,152}$ \irep{199}:\ \ 
$d_i$ = ($1.0$,
$1.965$,
$2.864$,
$3.666$,
$4.342$,
$4.871$,
$5.234$,
$5.418$) 

\vskip 0.7ex
\hangindent=3em \hangafter=1
$D^2= 125.874 = 
36+28c^{1}_{17}
+21c^{2}_{17}
+15c^{3}_{17}
+10c^{4}_{17}
+6c^{5}_{17}
+3c^{6}_{17}
+c^{7}_{17}
$

\vskip 0.7ex
\hangindent=3em \hangafter=1
$T = ( 0,
\frac{5}{17},
\frac{2}{17},
\frac{8}{17},
\frac{6}{17},
\frac{13}{17},
\frac{12}{17},
\frac{3}{17} )
$,

\vskip 0.7ex
\hangindent=3em \hangafter=1
$S$ = ($ 1$,
$ -c_{17}^{8}$,
$ \xi_{17}^{3}$,
$ \xi_{17}^{13}$,
$ \xi_{17}^{5}$,
$ \xi_{17}^{11}$,
$ \xi_{17}^{7}$,
$ \xi_{17}^{9}$;\ \ 
$ -\xi_{17}^{13}$,
$ \xi_{17}^{11}$,
$ -\xi_{17}^{9}$,
$ \xi_{17}^{7}$,
$ -\xi_{17}^{5}$,
$ \xi_{17}^{3}$,
$ -1$;\ \ 
$ \xi_{17}^{9}$,
$ \xi_{17}^{5}$,
$ -c_{17}^{8}$,
$ -1$,
$ -\xi_{17}^{13}$,
$ -\xi_{17}^{7}$;\ \ 
$ -1$,
$ -\xi_{17}^{3}$,
$ \xi_{17}^{7}$,
$ -\xi_{17}^{11}$,
$ -c_{17}^{8}$;\ \ 
$ -\xi_{17}^{9}$,
$ -\xi_{17}^{13}$,
$ 1$,
$ \xi_{17}^{11}$;\ \ 
$ c_{17}^{8}$,
$ \xi_{17}^{9}$,
$ -\xi_{17}^{3}$;\ \ 
$ -c_{17}^{8}$,
$ -\xi_{17}^{5}$;\ \ 
$ \xi_{17}^{13}$)

Realization: $PSU(2)_{15}$.

\vskip 1ex

\noindent9. $8_{\frac{36}{13},223.6}^{13,370}$ \irep{179}:\ \ 
$d_i$ = ($1.0$,
$2.941$,
$4.148$,
$4.148$,
$4.712$,
$6.209$,
$7.345$,
$8.55$) 

\vskip 0.7ex
\hangindent=3em \hangafter=1
$D^2= 223.689 = 
78+65c^{1}_{13}
+52c^{2}_{13}
+39c^{3}_{13}
+26c^{4}_{13}
+13c^{5}_{13}
$

\vskip 0.7ex
\hangindent=3em \hangafter=1
$T = ( 0,
\frac{1}{13},
\frac{8}{13},
\frac{8}{13},
\frac{3}{13},
\frac{6}{13},
\frac{10}{13},
\frac{2}{13} )
$,

\vskip 0.7ex
\hangindent=3em \hangafter=1
$S$ = ($ 1$,
$ 2+c^{1}_{13}
+c^{2}_{13}
+c^{3}_{13}
+c^{4}_{13}
+c^{5}_{13}
$,
$ \xi_{13}^{7}$,
$ \xi_{13}^{7}$,
$ 2+2c^{1}_{13}
+c^{2}_{13}
+c^{3}_{13}
+c^{4}_{13}
+c^{5}_{13}
$,
$ 2+2c^{1}_{13}
+c^{2}_{13}
+c^{3}_{13}
+c^{4}_{13}
$,
$ 2+2c^{1}_{13}
+2c^{2}_{13}
+c^{3}_{13}
+c^{4}_{13}
$,
$ 2+2c^{1}_{13}
+2c^{2}_{13}
+c^{3}_{13}
$;\ \ 
$ 2+2c^{1}_{13}
+2c^{2}_{13}
+c^{3}_{13}
+c^{4}_{13}
$,
$ -\xi_{13}^{7}$,
$ -\xi_{13}^{7}$,
$ 2+2c^{1}_{13}
+2c^{2}_{13}
+c^{3}_{13}
$,
$ 2+2c^{1}_{13}
+c^{2}_{13}
+c^{3}_{13}
+c^{4}_{13}
+c^{5}_{13}
$,
$ -1$,
$ -2-2  c^{1}_{13}
-c^{2}_{13}
-c^{3}_{13}
-c^{4}_{13}
$;\ \ 
$ -1-c^{1}_{13}
-c^{2}_{13}
+c^{5}_{13}
$,
$ 2+2c^{1}_{13}
+2c^{2}_{13}
+c^{3}_{13}
-c^{5}_{13}
$,
$ \xi_{13}^{7}$,
$ -\xi_{13}^{7}$,
$ \xi_{13}^{7}$,
$ -\xi_{13}^{7}$;\ \ 
$ -1-c^{1}_{13}
-c^{2}_{13}
+c^{5}_{13}
$,
$ \xi_{13}^{7}$,
$ -\xi_{13}^{7}$,
$ \xi_{13}^{7}$,
$ -\xi_{13}^{7}$;\ \ 
$ 1$,
$ -2-2  c^{1}_{13}
-2  c^{2}_{13}
-c^{3}_{13}
-c^{4}_{13}
$,
$ -2-2  c^{1}_{13}
-c^{2}_{13}
-c^{3}_{13}
-c^{4}_{13}
$,
$ 2+c^{1}_{13}
+c^{2}_{13}
+c^{3}_{13}
+c^{4}_{13}
+c^{5}_{13}
$;\ \ 
$ -2-c^{1}_{13}
-c^{2}_{13}
-c^{3}_{13}
-c^{4}_{13}
-c^{5}_{13}
$,
$ 2+2c^{1}_{13}
+2c^{2}_{13}
+c^{3}_{13}
$,
$ 1$;\ \ 
$ -2-c^{1}_{13}
-c^{2}_{13}
-c^{3}_{13}
-c^{4}_{13}
-c^{5}_{13}
$,
$ -2-2  c^{1}_{13}
-c^{2}_{13}
-c^{3}_{13}
-c^{4}_{13}
-c^{5}_{13}
$;\ \ 
$ 2+2c^{1}_{13}
+2c^{2}_{13}
+c^{3}_{13}
+c^{4}_{13}
$)

Realization: $PSO(12)_3$.

\vskip 1ex

\noindent10. $8_{4,308.4}^{15,440}$ \irep{183}:\ \ 
$d_i$ = ($1.0$,
$5.854$,
$5.854$,
$5.854$,
$5.854$,
$6.854$,
$7.854$,
$7.854$) 

\vskip 0.7ex
\hangindent=3em \hangafter=1
$D^2= 308.434 = 
\frac{315+135\sqrt{5}}{2}$

\vskip 0.7ex
\hangindent=3em \hangafter=1
$T = ( 0,
0,
0,
\frac{1}{3},
\frac{2}{3},
0,
\frac{2}{5},
\frac{3}{5} )
$,

\vskip 0.7ex
\hangindent=3em \hangafter=1
$S$ = ($ 1$,
$ \frac{5+3\sqrt{5}}{2}$,
$ \frac{5+3\sqrt{5}}{2}$,
$ \frac{5+3\sqrt{5}}{2}$,
$ \frac{5+3\sqrt{5}}{2}$,
$ \frac{7+3\sqrt{5}}{2}$,
$ \frac{9+3\sqrt{5}}{2}$,
$ \frac{9+3\sqrt{5}}{2}$;\ \ 
$ -5-3\sqrt{5}$,
$ \frac{5+3\sqrt{5}}{2}$,
$ \frac{5+3\sqrt{5}}{2}$,
$ \frac{5+3\sqrt{5}}{2}$,
$ -\frac{5+3\sqrt{5}}{2}$,
$0$,
$0$;\ \ 
$ -5-3\sqrt{5}$,
$ \frac{5+3\sqrt{5}}{2}$,
$ \frac{5+3\sqrt{5}}{2}$,
$ -\frac{5+3\sqrt{5}}{2}$,
$0$,
$0$;\ \ 
$ \frac{5+3\sqrt{5}}{2}$,
$ -5-3\sqrt{5}$,
$ -\frac{5+3\sqrt{5}}{2}$,
$0$,
$0$;\ \ 
$ \frac{5+3\sqrt{5}}{2}$,
$ -\frac{5+3\sqrt{5}}{2}$,
$0$,
$0$;\ \ 
$ 1$,
$ \frac{9+3\sqrt{5}}{2}$,
$ \frac{9+3\sqrt{5}}{2}$;\ \ 
$ \frac{3+3\sqrt{5}}{2}$,
$ -6-3\sqrt{5}$;\ \ 
$ \frac{3+3\sqrt{5}}{2}$)

Realization: condensation reductions of $\eZ(\cNG(\Z_5,5))$ \cite{EvansGannon}.

\vskip 1ex

\noindent11. $8_{0,308.4}^{15,100}$ \irep{188}:\ \ 
$d_i$ = ($1.0$,
$5.854$,
$5.854$,
$5.854$,
$5.854$,
$6.854$,
$7.854$,
$7.854$) 

\vskip 0.7ex
\hangindent=3em \hangafter=1
$D^2= 308.434 = 
\frac{315+135\sqrt{5}}{2}$

\vskip 0.7ex
\hangindent=3em \hangafter=1
$T = ( 0,
\frac{1}{3},
\frac{1}{3},
\frac{2}{3},
\frac{2}{3},
0,
\frac{1}{5},
\frac{4}{5} )
$,

\vskip 0.7ex
\hangindent=3em \hangafter=1
$S$ = ($ 1$,
$ \frac{5+3\sqrt{5}}{2}$,
$ \frac{5+3\sqrt{5}}{2}$,
$ \frac{5+3\sqrt{5}}{2}$,
$ \frac{5+3\sqrt{5}}{2}$,
$ \frac{7+3\sqrt{5}}{2}$,
$ \frac{9+3\sqrt{5}}{2}$,
$ \frac{9+3\sqrt{5}}{2}$;\ \ 
$ \frac{5+3\sqrt{5}}{2}$,
$ -5-3\sqrt{5}$,
$ \frac{5+3\sqrt{5}}{2}$,
$ \frac{5+3\sqrt{5}}{2}$,
$ -\frac{5+3\sqrt{5}}{2}$,
$0$,
$0$;\ \ 
$ \frac{5+3\sqrt{5}}{2}$,
$ \frac{5+3\sqrt{5}}{2}$,
$ \frac{5+3\sqrt{5}}{2}$,
$ -\frac{5+3\sqrt{5}}{2}$,
$0$,
$0$;\ \ 
$ \frac{5+3\sqrt{5}}{2}$,
$ -5-3\sqrt{5}$,
$ -\frac{5+3\sqrt{5}}{2}$,
$0$,
$0$;\ \ 
$ \frac{5+3\sqrt{5}}{2}$,
$ -\frac{5+3\sqrt{5}}{2}$,
$0$,
$0$;\ \ 
$ 1$,
$ \frac{9+3\sqrt{5}}{2}$,
$ \frac{9+3\sqrt{5}}{2}$;\ \ 
$ -6-3\sqrt{5}$,
$ \frac{3+3\sqrt{5}}{2}$;\ \ 
$ -6-3\sqrt{5}$)

Realization: condensation reductions of $\eZ(\cNG(\Z_5,5))$ \cite{EvansGannon}.

\vskip 1ex

\noindent12. $8_{4,308.4}^{45,289}$ \irep{235}:\ \ 
$d_i$ = ($1.0$,
$5.854$,
$5.854$,
$5.854$,
$5.854$,
$6.854$,
$7.854$,
$7.854$) 

\vskip 0.7ex
\hangindent=3em \hangafter=1
$D^2= 308.434 = 
\frac{315+135\sqrt{5}}{2}$

\vskip 0.7ex
\hangindent=3em \hangafter=1
$T = ( 0,
0,
\frac{1}{9},
\frac{4}{9},
\frac{7}{9},
0,
\frac{2}{5},
\frac{3}{5} )
$,

\vskip 0.7ex
\hangindent=3em \hangafter=1
$S$ = ($ 1$,
$ \frac{5+3\sqrt{5}}{2}$,
$ \frac{5+3\sqrt{5}}{2}$,
$ \frac{5+3\sqrt{5}}{2}$,
$ \frac{5+3\sqrt{5}}{2}$,
$ \frac{7+3\sqrt{5}}{2}$,
$ \frac{9+3\sqrt{5}}{2}$,
$ \frac{9+3\sqrt{5}}{2}$;\ \ 
$ -5-3\sqrt{5}$,
$ \frac{5+3\sqrt{5}}{2}$,
$ \frac{5+3\sqrt{5}}{2}$,
$ \frac{5+3\sqrt{5}}{2}$,
$ -\frac{5+3\sqrt{5}}{2}$,
$0$,
$0$;\ \ 
$ -3  c^{1}_{45}
+3c^{4}_{45}
-4  c^{10}_{45}
+3c^{11}_{45}
$,
$ 3c^{2}_{45}
+c^{5}_{45}
+3c^{7}_{45}
+c^{10}_{45}
$,
$ 3c^{1}_{45}
-3  c^{2}_{45}
-3  c^{4}_{45}
-c^{5}_{45}
-3  c^{7}_{45}
+3c^{10}_{45}
-3  c^{11}_{45}
$,
$ -\frac{5+3\sqrt{5}}{2}$,
$0$,
$0$;\ \ 
$ 3c^{1}_{45}
-3  c^{2}_{45}
-3  c^{4}_{45}
-c^{5}_{45}
-3  c^{7}_{45}
+3c^{10}_{45}
-3  c^{11}_{45}
$,
$ -3  c^{1}_{45}
+3c^{4}_{45}
-4  c^{10}_{45}
+3c^{11}_{45}
$,
$ -\frac{5+3\sqrt{5}}{2}$,
$0$,
$0$;\ \ 
$ 3c^{2}_{45}
+c^{5}_{45}
+3c^{7}_{45}
+c^{10}_{45}
$,
$ -\frac{5+3\sqrt{5}}{2}$,
$0$,
$0$;\ \ 
$ 1$,
$ \frac{9+3\sqrt{5}}{2}$,
$ \frac{9+3\sqrt{5}}{2}$;\ \ 
$ \frac{3+3\sqrt{5}}{2}$,
$ -6-3\sqrt{5}$;\ \ 
$ \frac{3+3\sqrt{5}}{2}$)

Realization: condensation reductions of $\eZ(\cNG(\Z_5,5))$ \cite{EvansGannon}.

\vskip 1ex 

}

\subsection{Rank 9}\label{ss:rank9}

{\small

\noindent1. $9_{0,9.}^{9,620}$ \irep{0}:\ \ 
$d_i$ = ($1.0$,
$1.0$,
$1.0$,
$1.0$,
$1.0$,
$1.0$,
$1.0$,
$1.0$,
$1.0$) 

\vskip 0.7ex
\hangindent=3em \hangafter=1
$D^2= 9.0 = 
9$

\vskip 0.7ex
\hangindent=3em \hangafter=1
$T = ( 0,
0,
0,
\frac{1}{9},
\frac{1}{9},
\frac{4}{9},
\frac{4}{9},
\frac{7}{9},
\frac{7}{9} )
$,

\vskip 0.7ex
\hangindent=3em \hangafter=1
$S$ = ($ 1$,
$ 1$,
$ 1$,
$ 1$,
$ 1$,
$ 1$,
$ 1$,
$ 1$,
$ 1$;\ \ 
$ 1$,
$ 1$,
$ -\zeta_{6}^{1}$,
$ \zeta_{3}^{1}$,
$ -\zeta_{6}^{1}$,
$ \zeta_{3}^{1}$,
$ -\zeta_{6}^{1}$,
$ \zeta_{3}^{1}$;\ \ 
$ 1$,
$ \zeta_{3}^{1}$,
$ -\zeta_{6}^{1}$,
$ \zeta_{3}^{1}$,
$ -\zeta_{6}^{1}$,
$ \zeta_{3}^{1}$,
$ -\zeta_{6}^{1}$;\ \ 
$ -\zeta_{18}^{5}$,
$ \zeta_{9}^{2}$,
$ \zeta_{9}^{4}$,
$ -\zeta_{18}^{1}$,
$ \zeta_{9}^{1}$,
$ -\zeta_{18}^{7}$;\ \ 
$ -\zeta_{18}^{5}$,
$ -\zeta_{18}^{1}$,
$ \zeta_{9}^{4}$,
$ -\zeta_{18}^{7}$,
$ \zeta_{9}^{1}$;\ \ 
$ \zeta_{9}^{1}$,
$ -\zeta_{18}^{7}$,
$ -\zeta_{18}^{5}$,
$ \zeta_{9}^{2}$;\ \ 
$ \zeta_{9}^{1}$,
$ \zeta_{9}^{2}$,
$ -\zeta_{18}^{5}$;\ \ 
$ \zeta_{9}^{4}$,
$ -\zeta_{18}^{1}$;\ \ 
$ \zeta_{9}^{4}$)

Realization: $U(9)_1$

\vskip 1ex

\noindent2. $9_{2,44.}^{88,112}$ \irep{496}:\ \ 
$d_i$ = ($1.0$,
$1.0$,
$2.0$,
$2.0$,
$2.0$,
$2.0$,
$2.0$,
$3.316$,
$3.316$) 

\vskip 0.7ex
\hangindent=3em \hangafter=1
$D^2= 44.0 = 
44$

\vskip 0.7ex
\hangindent=3em \hangafter=1
$T = ( 0,
0,
\frac{1}{11},
\frac{3}{11},
\frac{4}{11},
\frac{5}{11},
\frac{9}{11},
\frac{1}{8},
\frac{5}{8} )
$,

\vskip 0.7ex
\hangindent=3em \hangafter=1
$S$ = ($ 1$,
$ 1$,
$ 2$,
$ 2$,
$ 2$,
$ 2$,
$ 2$,
$ \sqrt{11}$,
$ \sqrt{11}$;\ \ 
$ 1$,
$ 2$,
$ 2$,
$ 2$,
$ 2$,
$ 2$,
$ -\sqrt{11}$,
$ -\sqrt{11}$;\ \ 
$ 2c_{11}^{2}$,
$ 2c_{11}^{1}$,
$ 2c_{11}^{4}$,
$ 2c_{11}^{3}$,
$ 2c_{11}^{5}$,
$0$,
$0$;\ \ 
$ 2c_{11}^{5}$,
$ 2c_{11}^{2}$,
$ 2c_{11}^{4}$,
$ 2c_{11}^{3}$,
$0$,
$0$;\ \ 
$ 2c_{11}^{3}$,
$ 2c_{11}^{5}$,
$ 2c_{11}^{1}$,
$0$,
$0$;\ \ 
$ 2c_{11}^{1}$,
$ 2c_{11}^{2}$,
$0$,
$0$;\ \ 
$ 2c_{11}^{4}$,
$0$,
$0$;\ \ 
$ \sqrt{11}$,
$ -\sqrt{11}$;\ \ 
$ \sqrt{11}$)

Realization:
$SO(11)_2$ or
Abelian anyon condensation of $O_{11}$ or $\overline{SO(22)}_2$.

\vskip 1ex

\noindent3. $9_{2,44.}^{88,529}$ \irep{496}:\ \ 
$d_i$ = ($1.0$,
$1.0$,
$2.0$,
$2.0$,
$2.0$,
$2.0$,
$2.0$,
$3.316$,
$3.316$) 

\vskip 0.7ex
\hangindent=3em \hangafter=1
$D^2= 44.0 = 
44$

\vskip 0.7ex
\hangindent=3em \hangafter=1
$T = ( 0,
0,
\frac{1}{11},
\frac{3}{11},
\frac{4}{11},
\frac{5}{11},
\frac{9}{11},
\frac{3}{8},
\frac{7}{8} )
$,

\vskip 0.7ex
\hangindent=3em \hangafter=1
$S$ = ($ 1$,
$ 1$,
$ 2$,
$ 2$,
$ 2$,
$ 2$,
$ 2$,
$ \sqrt{11}$,
$ \sqrt{11}$;\ \ 
$ 1$,
$ 2$,
$ 2$,
$ 2$,
$ 2$,
$ 2$,
$ -\sqrt{11}$,
$ -\sqrt{11}$;\ \ 
$ 2c_{11}^{2}$,
$ 2c_{11}^{1}$,
$ 2c_{11}^{4}$,
$ 2c_{11}^{3}$,
$ 2c_{11}^{5}$,
$0$,
$0$;\ \ 
$ 2c_{11}^{5}$,
$ 2c_{11}^{2}$,
$ 2c_{11}^{4}$,
$ 2c_{11}^{3}$,
$0$,
$0$;\ \ 
$ 2c_{11}^{3}$,
$ 2c_{11}^{5}$,
$ 2c_{11}^{1}$,
$0$,
$0$;\ \ 
$ 2c_{11}^{1}$,
$ 2c_{11}^{2}$,
$0$,
$0$;\ \ 
$ 2c_{11}^{4}$,
$0$,
$0$;\ \ 
$ -\sqrt{11}$,
$ \sqrt{11}$;\ \ 
$ -\sqrt{11}$)

Realization:
Abelian anyon condensation of 
$SO(11)_2$ or
$O_{11}$ or $\overline{SO(22)}_2$.

\vskip 1ex

\noindent4. $9_{\frac{12}{5},52.36}^{40,304}$ \irep{483}:\ \ 
$d_i$ = ($1.0$,
$1.0$,
$1.902$,
$1.902$,
$2.618$,
$2.618$,
$3.77$,
$3.77$,
$3.236$) 

\vskip 0.7ex
\hangindent=3em \hangafter=1
$D^2= 52.360 = 
30+10\sqrt{5}$

\vskip 0.7ex
\hangindent=3em \hangafter=1
$T = ( 0,
0,
\frac{3}{40},
\frac{23}{40},
\frac{1}{5},
\frac{1}{5},
\frac{3}{8},
\frac{7}{8},
\frac{3}{5} )
$,

\vskip 0.7ex
\hangindent=3em \hangafter=1
$S$ = ($ 1$,
$ 1$,
$ c_{20}^{1}$,
$ c_{20}^{1}$,
$ \frac{3+\sqrt{5}}{2}$,
$ \frac{3+\sqrt{5}}{2}$,
$ c^{1}_{20}
+c^{3}_{20}
$,
$ c^{1}_{20}
+c^{3}_{20}
$,
$ 1+\sqrt{5}$;\ \ 
$ 1$,
$ -c_{20}^{1}$,
$ -c_{20}^{1}$,
$ \frac{3+\sqrt{5}}{2}$,
$ \frac{3+\sqrt{5}}{2}$,
$ -c^{1}_{20}
-c^{3}_{20}
$,
$ -c^{1}_{20}
-c^{3}_{20}
$,
$ 1+\sqrt{5}$;\ \ 
$ c^{1}_{20}
+c^{3}_{20}
$,
$ -c^{1}_{20}
-c^{3}_{20}
$,
$ -c^{1}_{20}
-c^{3}_{20}
$,
$ c^{1}_{20}
+c^{3}_{20}
$,
$ c_{20}^{1}$,
$ -c_{20}^{1}$,
$0$;\ \ 
$ c^{1}_{20}
+c^{3}_{20}
$,
$ -c^{1}_{20}
-c^{3}_{20}
$,
$ c^{1}_{20}
+c^{3}_{20}
$,
$ -c_{20}^{1}$,
$ c_{20}^{1}$,
$0$;\ \ 
$ 1$,
$ 1$,
$ c_{20}^{1}$,
$ c_{20}^{1}$,
$ -1-\sqrt{5}$;\ \ 
$ 1$,
$ -c_{20}^{1}$,
$ -c_{20}^{1}$,
$ -1-\sqrt{5}$;\ \ 
$ -c^{1}_{20}
-c^{3}_{20}
$,
$ c^{1}_{20}
+c^{3}_{20}
$,
$0$;\ \ 
$ -c^{1}_{20}
-c^{3}_{20}
$,
$0$;\ \ 
$ 1+\sqrt{5}$)

Realization:
$SU(2)_8$.
Abelian anyon condensation of 
$\overline{SU(8)}_2$ or $\overline{Sp(18)}_1$.

\vskip 1ex

\noindent5. $9_{\frac{28}{5},52.36}^{40,247}$ \irep{483}:\ \ 
$d_i$ = ($1.0$,
$1.0$,
$1.902$,
$1.902$,
$2.618$,
$2.618$,
$3.77$,
$3.77$,
$3.236$) 

\vskip 0.7ex
\hangindent=3em \hangafter=1
$D^2= 52.360 = 
30+10\sqrt{5}$

\vskip 0.7ex
\hangindent=3em \hangafter=1
$T = ( 0,
0,
\frac{7}{40},
\frac{27}{40},
\frac{4}{5},
\frac{4}{5},
\frac{3}{8},
\frac{7}{8},
\frac{2}{5} )
$,

\vskip 0.7ex
\hangindent=3em \hangafter=1
$S$ = ($ 1$,
$ 1$,
$ c_{20}^{1}$,
$ c_{20}^{1}$,
$ \frac{3+\sqrt{5}}{2}$,
$ \frac{3+\sqrt{5}}{2}$,
$ c^{1}_{20}
+c^{3}_{20}
$,
$ c^{1}_{20}
+c^{3}_{20}
$,
$ 1+\sqrt{5}$;\ \ 
$ 1$,
$ -c_{20}^{1}$,
$ -c_{20}^{1}$,
$ \frac{3+\sqrt{5}}{2}$,
$ \frac{3+\sqrt{5}}{2}$,
$ -c^{1}_{20}
-c^{3}_{20}
$,
$ -c^{1}_{20}
-c^{3}_{20}
$,
$ 1+\sqrt{5}$;\ \ 
$ -c^{1}_{20}
-c^{3}_{20}
$,
$ c^{1}_{20}
+c^{3}_{20}
$,
$ -c^{1}_{20}
-c^{3}_{20}
$,
$ c^{1}_{20}
+c^{3}_{20}
$,
$ c_{20}^{1}$,
$ -c_{20}^{1}$,
$0$;\ \ 
$ -c^{1}_{20}
-c^{3}_{20}
$,
$ -c^{1}_{20}
-c^{3}_{20}
$,
$ c^{1}_{20}
+c^{3}_{20}
$,
$ -c_{20}^{1}$,
$ c_{20}^{1}$,
$0$;\ \ 
$ 1$,
$ 1$,
$ c_{20}^{1}$,
$ c_{20}^{1}$,
$ -1-\sqrt{5}$;\ \ 
$ 1$,
$ -c_{20}^{1}$,
$ -c_{20}^{1}$,
$ -1-\sqrt{5}$;\ \ 
$ c^{1}_{20}
+c^{3}_{20}
$,
$ -c^{1}_{20}
-c^{3}_{20}
$,
$0$;\ \ 
$ c^{1}_{20}
+c^{3}_{20}
$,
$0$;\ \ 
$ 1+\sqrt{5}$)

Realization:
Abelian anyon condensation of 
$\overline{SU(2)}_8$ or
$SU(8)_2$ or $Sp(18)_1$.

\vskip 1ex

\noindent6. $9_{\frac{120}{19},175.3}^{19,574}$ \irep{414}:\ \ 
$d_i$ = ($1.0$,
$1.972$,
$2.891$,
$3.731$,
$4.469$,
$5.86$,
$5.563$,
$5.889$,
$6.54$) 

\vskip 0.7ex
\hangindent=3em \hangafter=1
$D^2= 175.332 = 
45+36c^{1}_{19}
+28c^{2}_{19}
+21c^{3}_{19}
+15c^{4}_{19}
+10c^{5}_{19}
+6c^{6}_{19}
+3c^{7}_{19}
+c^{8}_{19}
$

\vskip 0.7ex
\hangindent=3em \hangafter=1
$T = ( 0,
\frac{4}{19},
\frac{17}{19},
\frac{1}{19},
\frac{13}{19},
\frac{15}{19},
\frac{7}{19},
\frac{8}{19},
\frac{18}{19} )
$,

\vskip 0.7ex
\hangindent=3em \hangafter=1
$S$ = ($ 1$,
$ -c_{19}^{9}$,
$ \xi_{19}^{3}$,
$ \xi_{19}^{15}$,
$ \xi_{19}^{5}$,
$ \xi_{19}^{13}$,
$ \xi_{19}^{7}$,
$ \xi_{19}^{11}$,
$ \xi_{19}^{9}$;\ \ 
$ -\xi_{19}^{15}$,
$ \xi_{19}^{13}$,
$ -\xi_{19}^{11}$,
$ \xi_{19}^{9}$,
$ -\xi_{19}^{7}$,
$ \xi_{19}^{5}$,
$ -\xi_{19}^{3}$,
$ 1$;\ \ 
$ \xi_{19}^{9}$,
$ \xi_{19}^{7}$,
$ \xi_{19}^{15}$,
$ 1$,
$ c_{19}^{9}$,
$ -\xi_{19}^{5}$,
$ -\xi_{19}^{11}$;\ \ 
$ -\xi_{19}^{3}$,
$ -1$,
$ \xi_{19}^{5}$,
$ -\xi_{19}^{9}$,
$ \xi_{19}^{13}$,
$ c_{19}^{9}$;\ \ 
$ -\xi_{19}^{13}$,
$ -\xi_{19}^{11}$,
$ -\xi_{19}^{3}$,
$ -c_{19}^{9}$,
$ \xi_{19}^{7}$;\ \ 
$ -c_{19}^{9}$,
$ \xi_{19}^{15}$,
$ -\xi_{19}^{9}$,
$ \xi_{19}^{3}$;\ \ 
$ \xi_{19}^{11}$,
$ 1$,
$ -\xi_{19}^{13}$;\ \ 
$ \xi_{19}^{7}$,
$ -\xi_{19}^{15}$;\ \ 
$ \xi_{19}^{5}$)

Realization: $\overline{PSU(2)}_{17}$. Abelian anyon condensation
of $\overline{SU(2)}_{17}$ or $SU(17)_2$ or $Sp(34)_1$.

\vskip 1ex

\noindent7. $9_{\frac{14}{5},343.2}^{15,715}$ \irep{361}:\ \ 
$d_i$ = ($1.0$,
$2.956$,
$4.783$,
$4.783$,
$4.783$,
$6.401$,
$7.739$,
$8.739$,
$9.357$) 

\vskip 0.7ex
\hangindent=3em \hangafter=1
$D^2= 343.211 = 
105+45c^{1}_{15}
+75c^{2}_{15}
+90c^{3}_{15}
$

\vskip 0.7ex
\hangindent=3em \hangafter=1
$T = ( 0,
\frac{1}{15},
\frac{1}{5},
\frac{13}{15},
\frac{13}{15},
\frac{2}{5},
\frac{2}{3},
0,
\frac{2}{5} )
$,

\vskip 0.7ex
\hangindent=3em \hangafter=1
$S$ = ($ 1$,
$ 1+c^{2}_{15}
+c^{3}_{15}
$,
$ \xi_{15}^{7}$,
$ \xi_{15}^{7}$,
$ \xi_{15}^{7}$,
$ 2+c^{1}_{15}
+c^{2}_{15}
+2c^{3}_{15}
$,
$ 2+c^{1}_{15}
+2c^{2}_{15}
+2c^{3}_{15}
$,
$ 3+c^{1}_{15}
+2c^{2}_{15}
+2c^{3}_{15}
$,
$ 3+c^{1}_{15}
+2c^{2}_{15}
+3c^{3}_{15}
$;\ \ 
$ 2+c^{1}_{15}
+2c^{2}_{15}
+2c^{3}_{15}
$,
$ 2\xi_{15}^{7}$,
$ -\xi_{15}^{7}$,
$ -\xi_{15}^{7}$,
$ 2+c^{1}_{15}
+2c^{2}_{15}
+2c^{3}_{15}
$,
$ 1+c^{2}_{15}
+c^{3}_{15}
$,
$ -1-c^{2}_{15}
-c^{3}_{15}
$,
$ -2-c^{1}_{15}
-2  c^{2}_{15}
-2  c^{3}_{15}
$;\ \ 
$ \xi_{15}^{7}$,
$ \xi_{15}^{7}$,
$ \xi_{15}^{7}$,
$ -\xi_{15}^{7}$,
$ -2\xi_{15}^{7}$,
$ -\xi_{15}^{7}$,
$ \xi_{15}^{7}$;\ \ 
$ 2-4  \zeta^{1}_{15}
-2  \zeta^{-1}_{15}
+\zeta^{2}_{15}
-2  \zeta^{-2}_{15}
-\zeta^{3}_{15}
+5\zeta^{-3}_{15}
-5  \zeta^{4}_{15}
$,
$ -3+3\zeta^{1}_{15}
+\zeta^{-1}_{15}
-2  \zeta^{2}_{15}
+\zeta^{-2}_{15}
-6  \zeta^{-3}_{15}
+5\zeta^{4}_{15}
$,
$ -\xi_{15}^{7}$,
$ \xi_{15}^{7}$,
$ -\xi_{15}^{7}$,
$ \xi_{15}^{7}$;\ \ 
$ 2-4  \zeta^{1}_{15}
-2  \zeta^{-1}_{15}
+\zeta^{2}_{15}
-2  \zeta^{-2}_{15}
-\zeta^{3}_{15}
+5\zeta^{-3}_{15}
-5  \zeta^{4}_{15}
$,
$ -\xi_{15}^{7}$,
$ \xi_{15}^{7}$,
$ -\xi_{15}^{7}$,
$ \xi_{15}^{7}$;\ \ 
$ -3-c^{1}_{15}
-2  c^{2}_{15}
-2  c^{3}_{15}
$,
$ 1+c^{2}_{15}
+c^{3}_{15}
$,
$ 3+c^{1}_{15}
+2c^{2}_{15}
+3c^{3}_{15}
$,
$ -1$;\ \ 
$ 2+c^{1}_{15}
+2c^{2}_{15}
+2c^{3}_{15}
$,
$ -2-c^{1}_{15}
-2  c^{2}_{15}
-2  c^{3}_{15}
$,
$ -1-c^{2}_{15}
-c^{3}_{15}
$;\ \ 
$ 1$,
$ 2+c^{1}_{15}
+c^{2}_{15}
+2c^{3}_{15}
$;\ \ 
$ -3-c^{1}_{15}
-2  c^{2}_{15}
-2  c^{3}_{15}
$)

Realization: $PSO(14)_3$.

\vskip 1ex

\noindent8. $9_{7,475.1}^{24,793}$ \irep{448}:\ \ 
$d_i$ = ($1.0$,
$4.449$,
$4.449$,
$5.449$,
$5.449$,
$8.898$,
$8.898$,
$9.898$,
$10.898$) 

\vskip 0.7ex
\hangindent=3em \hangafter=1
$D^2= 475.151 = 
240+96\sqrt{6}$

\vskip 0.7ex
\hangindent=3em \hangafter=1
$T = ( 0,
\frac{1}{4},
\frac{1}{4},
\frac{1}{2},
\frac{1}{2},
\frac{1}{3},
\frac{7}{12},
0,
\frac{7}{8} )
$,

\vskip 0.7ex
\hangindent=3em \hangafter=1
$S$ = ($ 1$,
$ 2+\sqrt{6}$,
$ 2+\sqrt{6}$,
$ 3+\sqrt{6}$,
$ 3+\sqrt{6}$,
$ 4+2\sqrt{6}$,
$ 4+2\sqrt{6}$,
$ 5+2\sqrt{6}$,
$ 6+2\sqrt{6}$;\ \ 
$ 2\xi_{24}^{7}$,
$ 2+2c^{1}_{24}
-2  c^{2}_{24}
-4  c^{3}_{24}
$,
$ -2  c^{2}_{24}
-3  c^{3}_{24}
$,
$ 2c^{2}_{24}
+3c^{3}_{24}
$,
$ -4-2\sqrt{6}$,
$ 4+2\sqrt{6}$,
$ -2-\sqrt{6}$,
$0$;\ \ 
$ 2\xi_{24}^{7}$,
$ 2c^{2}_{24}
+3c^{3}_{24}
$,
$ -2  c^{2}_{24}
-3  c^{3}_{24}
$,
$ -4-2\sqrt{6}$,
$ 4+2\sqrt{6}$,
$ -2-\sqrt{6}$,
$0$;\ \ 
$ 3+2c^{1}_{24}
-2  c^{2}_{24}
-4  c^{3}_{24}
$,
$ 3+2c^{1}_{24}
+2c^{2}_{24}
+2c^{3}_{24}
$,
$0$,
$0$,
$ 3+\sqrt{6}$,
$ -6-2\sqrt{6}$;\ \ 
$ 3+2c^{1}_{24}
-2  c^{2}_{24}
-4  c^{3}_{24}
$,
$0$,
$0$,
$ 3+\sqrt{6}$,
$ -6-2\sqrt{6}$;\ \ 
$ 4+2\sqrt{6}$,
$ 4+2\sqrt{6}$,
$ -4-2\sqrt{6}$,
$0$;\ \ 
$ -4-2\sqrt{6}$,
$ -4-2\sqrt{6}$,
$0$;\ \ 
$ 1$,
$ 6+2\sqrt{6}$;\ \ 
$0$)

Realization:  $G(2)_4$.

\vskip 1ex

\noindent9. $9_{6,668.5}^{12,567}$ \irep{330}:\ \ 
$d_i$ = ($1.0$,
$6.464$,
$6.464$,
$6.464$,
$6.464$,
$6.464$,
$6.464$,
$13.928$,
$14.928$) 

\vskip 0.7ex
\hangindent=3em \hangafter=1
$D^2= 668.553 = 
336+192\sqrt{3}$

\vskip 0.7ex
\hangindent=3em \hangafter=1
$T = ( 0,
\frac{1}{2},
\frac{1}{2},
\frac{1}{4},
\frac{1}{4},
\frac{1}{4},
\frac{1}{4},
0,
\frac{2}{3} )
$,

\vskip 0.7ex
\hangindent=3em \hangafter=1
$S$ = ($ 1$,
$ 3+2\sqrt{3}$,
$ 3+2\sqrt{3}$,
$ 3+2\sqrt{3}$,
$ 3+2\sqrt{3}$,
$ 3+2\sqrt{3}$,
$ 3+2\sqrt{3}$,
$ 7+4\sqrt{3}$,
$ 8+4\sqrt{3}$;\ \ 
$ -7+4\zeta^{1}_{12}
-8  \zeta^{-1}_{12}
+8\zeta^{2}_{12}
$,
$ 1-8  \zeta^{1}_{12}
+4\zeta^{-1}_{12}
-8  \zeta^{2}_{12}
$,
$ 3+2\sqrt{3}$,
$ 3+2\sqrt{3}$,
$ 3+2\sqrt{3}$,
$ 3+2\sqrt{3}$,
$ -3-2\sqrt{3}$,
$0$;\ \ 
$ -7+4\zeta^{1}_{12}
-8  \zeta^{-1}_{12}
+8\zeta^{2}_{12}
$,
$ 3+2\sqrt{3}$,
$ 3+2\sqrt{3}$,
$ 3+2\sqrt{3}$,
$ 3+2\sqrt{3}$,
$ -3-2\sqrt{3}$,
$0$;\ \ 
$ 9+6\sqrt{3}$,
$ -3-2\sqrt{3}$,
$ -3-2\sqrt{3}$,
$ -3-2\sqrt{3}$,
$ -3-2\sqrt{3}$,
$0$;\ \ 
$ 9+6\sqrt{3}$,
$ -3-2\sqrt{3}$,
$ -3-2\sqrt{3}$,
$ -3-2\sqrt{3}$,
$0$;\ \ 
$ 9+6\sqrt{3}$,
$ -3-2\sqrt{3}$,
$ -3-2\sqrt{3}$,
$0$;\ \ 
$ 9+6\sqrt{3}$,
$ -3-2\sqrt{3}$,
$0$;\ \ 
$ 1$,
$ 8+4\sqrt{3}$;\ \ 
$ -8-4\sqrt{3}$)

Realization: $\Z_3$-algebra condensation of $SU(3)_9$ (see also
\cite{Edie-MichellCAMS}.)

\vskip 1ex 

}

\subsection{Rank 10 }\label{ss:rank10}

{\small

\noindent1. $10_{3,24.}^{48,945}$ \irep{1119}:\ \ 
$d_i$ = ($1.0$,
$1.0$,
$1.0$,
$1.0$,
$1.732$,
$1.732$,
$1.732$,
$1.732$,
$2.0$,
$2.0$) 

\vskip 0.7ex
\hangindent=3em \hangafter=1
$D^2= 24.0 = 
24$

\vskip 0.7ex
\hangindent=3em \hangafter=1
$T = ( 0,
0,
\frac{1}{4},
\frac{1}{4},
\frac{3}{16},
\frac{3}{16},
\frac{11}{16},
\frac{11}{16},
\frac{1}{3},
\frac{7}{12} )
$,

\vskip 0.7ex
\hangindent=3em \hangafter=1
$S$ = ($ 1$,
$ 1$,
$ 1$,
$ 1$,
$ \sqrt{3}$,
$ \sqrt{3}$,
$ \sqrt{3}$,
$ \sqrt{3}$,
$ 2$,
$ 2$;\ \ 
$ 1$,
$ 1$,
$ 1$,
$ -\sqrt{3}$,
$ -\sqrt{3}$,
$ -\sqrt{3}$,
$ -\sqrt{3}$,
$ 2$,
$ 2$;\ \ 
$ -1$,
$ -1$,
$(-\sqrt{3})\mathrm{i}$,
$(\sqrt{3})\mathrm{i}$,
$(-\sqrt{3})\mathrm{i}$,
$(\sqrt{3})\mathrm{i}$,
$ 2$,
$ -2$;\ \ 
$ -1$,
$(\sqrt{3})\mathrm{i}$,
$(-\sqrt{3})\mathrm{i}$,
$(\sqrt{3})\mathrm{i}$,
$(-\sqrt{3})\mathrm{i}$,
$ 2$,
$ -2$;\ \ 
$ -\sqrt{3}\zeta_{8}^{3}$,
$ \sqrt{3}\zeta_{8}^{1}$,
$ \sqrt{3}\zeta_{8}^{3}$,
$ -\sqrt{3}\zeta_{8}^{1}$,
$0$,
$0$;\ \ 
$ -\sqrt{3}\zeta_{8}^{3}$,
$ -\sqrt{3}\zeta_{8}^{1}$,
$ \sqrt{3}\zeta_{8}^{3}$,
$0$,
$0$;\ \ 
$ -\sqrt{3}\zeta_{8}^{3}$,
$ \sqrt{3}\zeta_{8}^{1}$,
$0$,
$0$;\ \ 
$ -\sqrt{3}\zeta_{8}^{3}$,
$0$,
$0$;\ \ 
$ -2$,
$ -2$;\ \ 
$ 2$)

Realization:
Abelian anyon condensation of $SU(2)_4$ or $O_3$
or $\overline{SU(4)}_2$ or $\overline{Sp(8)}_1$ or

\vskip 1ex

\noindent2. $10_{7,24.}^{48,721}$ \irep{1119}:\ \ 
$d_i$ = ($1.0$,
$1.0$,
$1.0$,
$1.0$,
$1.732$,
$1.732$,
$1.732$,
$1.732$,
$2.0$,
$2.0$) 

\vskip 0.7ex
\hangindent=3em \hangafter=1
$D^2= 24.0 = 
24$

\vskip 0.7ex
\hangindent=3em \hangafter=1
$T = ( 0,
0,
\frac{1}{4},
\frac{1}{4},
\frac{3}{16},
\frac{3}{16},
\frac{11}{16},
\frac{11}{16},
\frac{2}{3},
\frac{11}{12} )
$,

\vskip 0.7ex
\hangindent=3em \hangafter=1
$S$ = ($ 1$,
$ 1$,
$ 1$,
$ 1$,
$ \sqrt{3}$,
$ \sqrt{3}$,
$ \sqrt{3}$,
$ \sqrt{3}$,
$ 2$,
$ 2$;\ \ 
$ 1$,
$ 1$,
$ 1$,
$ -\sqrt{3}$,
$ -\sqrt{3}$,
$ -\sqrt{3}$,
$ -\sqrt{3}$,
$ 2$,
$ 2$;\ \ 
$ -1$,
$ -1$,
$(-\sqrt{3})\mathrm{i}$,
$(\sqrt{3})\mathrm{i}$,
$(-\sqrt{3})\mathrm{i}$,
$(\sqrt{3})\mathrm{i}$,
$ 2$,
$ -2$;\ \ 
$ -1$,
$(\sqrt{3})\mathrm{i}$,
$(-\sqrt{3})\mathrm{i}$,
$(\sqrt{3})\mathrm{i}$,
$(-\sqrt{3})\mathrm{i}$,
$ 2$,
$ -2$;\ \ 
$ \sqrt{3}\zeta_{8}^{3}$,
$ -\sqrt{3}\zeta_{8}^{1}$,
$ -\sqrt{3}\zeta_{8}^{3}$,
$ \sqrt{3}\zeta_{8}^{1}$,
$0$,
$0$;\ \ 
$ \sqrt{3}\zeta_{8}^{3}$,
$ \sqrt{3}\zeta_{8}^{1}$,
$ -\sqrt{3}\zeta_{8}^{3}$,
$0$,
$0$;\ \ 
$ \sqrt{3}\zeta_{8}^{3}$,
$ -\sqrt{3}\zeta_{8}^{1}$,
$0$,
$0$;\ \ 
$ \sqrt{3}\zeta_{8}^{3}$,
$0$,
$0$;\ \ 
$ -2$,
$ -2$;\ \ 
$ 2$)

Realization:
Abelian anyon condensation of $SU(4)_2$ or $Sp(8)_1$ or
$\overline{SU(2)}_4$ or $\overline{O}_3$.

\vskip 1ex

\noindent3. $10_{4,36.}^{6,152}$ \irep{0}:\ \ 
$d_i$ = ($1.0$,
$1.0$,
$1.0$,
$2.0$,
$2.0$,
$2.0$,
$2.0$,
$2.0$,
$2.0$,
$3.0$) 

\vskip 0.7ex
\hangindent=3em \hangafter=1
$D^2= 36.0 = 
36$

\vskip 0.7ex
\hangindent=3em \hangafter=1
$T = ( 0,
0,
0,
0,
0,
\frac{1}{3},
\frac{1}{3},
\frac{2}{3},
\frac{2}{3},
\frac{1}{2} )
$,

\vskip 0.7ex
\hangindent=3em \hangafter=1
$S$ = ($ 1$,
$ 1$,
$ 1$,
$ 2$,
$ 2$,
$ 2$,
$ 2$,
$ 2$,
$ 2$,
$ 3$;\ \ 
$ 1$,
$ 1$,
$ -2\zeta_{6}^{1}$,
$ 2\zeta_{3}^{1}$,
$ -2\zeta_{6}^{1}$,
$ 2\zeta_{3}^{1}$,
$ -2\zeta_{6}^{1}$,
$ 2\zeta_{3}^{1}$,
$ 3$;\ \ 
$ 1$,
$ 2\zeta_{3}^{1}$,
$ -2\zeta_{6}^{1}$,
$ 2\zeta_{3}^{1}$,
$ -2\zeta_{6}^{1}$,
$ 2\zeta_{3}^{1}$,
$ -2\zeta_{6}^{1}$,
$ 3$;\ \ 
$ -2$,
$ -2$,
$ 2\zeta_{6}^{1}$,
$ -2\zeta_{3}^{1}$,
$ -2\zeta_{3}^{1}$,
$ 2\zeta_{6}^{1}$,
$0$;\ \ 
$ -2$,
$ -2\zeta_{3}^{1}$,
$ 2\zeta_{6}^{1}$,
$ 2\zeta_{6}^{1}$,
$ -2\zeta_{3}^{1}$,
$0$;\ \ 
$ -2\zeta_{3}^{1}$,
$ 2\zeta_{6}^{1}$,
$ -2$,
$ -2$,
$0$;\ \ 
$ -2\zeta_{3}^{1}$,
$ -2$,
$ -2$,
$0$;\ \ 
$ 2\zeta_{6}^{1}$,
$ -2\zeta_{3}^{1}$,
$0$;\ \ 
$ 2\zeta_{6}^{1}$,
$0$;\ \ 
$ -3$)

Realization:  $SU(3)_3$.

\vskip 1ex

\noindent4. $10_{4,36.}^{18,490}$ \irep{0}:\ \ 
$d_i$ = ($1.0$,
$1.0$,
$1.0$,
$2.0$,
$2.0$,
$2.0$,
$2.0$,
$2.0$,
$2.0$,
$3.0$) 

\vskip 0.7ex
\hangindent=3em \hangafter=1
$D^2= 36.0 = 
36$

\vskip 0.7ex
\hangindent=3em \hangafter=1
$T = ( 0,
0,
0,
\frac{1}{9},
\frac{1}{9},
\frac{4}{9},
\frac{4}{9},
\frac{7}{9},
\frac{7}{9},
\frac{1}{2} )
$,

\vskip 0.7ex
\hangindent=3em \hangafter=1
$S$ = ($ 1$,
$ 1$,
$ 1$,
$ 2$,
$ 2$,
$ 2$,
$ 2$,
$ 2$,
$ 2$,
$ 3$;\ \ 
$ 1$,
$ 1$,
$ -2\zeta_{6}^{1}$,
$ 2\zeta_{3}^{1}$,
$ -2\zeta_{6}^{1}$,
$ 2\zeta_{3}^{1}$,
$ -2\zeta_{6}^{1}$,
$ 2\zeta_{3}^{1}$,
$ 3$;\ \ 
$ 1$,
$ 2\zeta_{3}^{1}$,
$ -2\zeta_{6}^{1}$,
$ 2\zeta_{3}^{1}$,
$ -2\zeta_{6}^{1}$,
$ 2\zeta_{3}^{1}$,
$ -2\zeta_{6}^{1}$,
$ 3$;\ \ 
$ 2\zeta_{18}^{5}$,
$ -2\zeta_{9}^{2}$,
$ -2\zeta_{9}^{4}$,
$ 2\zeta_{18}^{1}$,
$ -2\zeta_{9}^{1}$,
$ 2\zeta_{18}^{7}$,
$0$;\ \ 
$ 2\zeta_{18}^{5}$,
$ 2\zeta_{18}^{1}$,
$ -2\zeta_{9}^{4}$,
$ 2\zeta_{18}^{7}$,
$ -2\zeta_{9}^{1}$,
$0$;\ \ 
$ -2\zeta_{9}^{1}$,
$ 2\zeta_{18}^{7}$,
$ 2\zeta_{18}^{5}$,
$ -2\zeta_{9}^{2}$,
$0$;\ \ 
$ -2\zeta_{9}^{1}$,
$ -2\zeta_{9}^{2}$,
$ 2\zeta_{18}^{5}$,
$0$;\ \ 
$ -2\zeta_{9}^{4}$,
$ 2\zeta_{18}^{1}$,
$0$;\ \ 
$ -2\zeta_{9}^{4}$,
$0$;\ \ 
$ -3$)

Realization: Zesting or Abelian anyon condensation of $SU(3)_3$ \cite{DGPRZ}.

\vskip 1ex

\noindent5. $10_{0,52.}^{26,247}$ \irep{1035}:\ \ 
$d_i$ = ($1.0$,
$1.0$,
$2.0$,
$2.0$,
$2.0$,
$2.0$,
$2.0$,
$2.0$,
$3.605$,
$3.605$) 

\vskip 0.7ex
\hangindent=3em \hangafter=1
$D^2= 52.0 = 
52$

\vskip 0.7ex
\hangindent=3em \hangafter=1
$T = ( 0,
0,
\frac{1}{13},
\frac{3}{13},
\frac{4}{13},
\frac{9}{13},
\frac{10}{13},
\frac{12}{13},
0,
\frac{1}{2} )
$,

\vskip 0.7ex
\hangindent=3em \hangafter=1
$S$ = ($ 1$,
$ 1$,
$ 2$,
$ 2$,
$ 2$,
$ 2$,
$ 2$,
$ 2$,
$ \sqrt{13}$,
$ \sqrt{13}$;\ \ 
$ 1$,
$ 2$,
$ 2$,
$ 2$,
$ 2$,
$ 2$,
$ 2$,
$ -\sqrt{13}$,
$ -\sqrt{13}$;\ \ 
$ 2c_{13}^{2}$,
$ 2c_{13}^{5}$,
$ 2c_{13}^{4}$,
$ 2c_{13}^{6}$,
$ 2c_{13}^{1}$,
$ 2c_{13}^{3}$,
$0$,
$0$;\ \ 
$ 2c_{13}^{6}$,
$ 2c_{13}^{3}$,
$ 2c_{13}^{2}$,
$ 2c_{13}^{4}$,
$ 2c_{13}^{1}$,
$0$,
$0$;\ \ 
$ 2c_{13}^{5}$,
$ 2c_{13}^{1}$,
$ 2c_{13}^{2}$,
$ 2c_{13}^{6}$,
$0$,
$0$;\ \ 
$ 2c_{13}^{5}$,
$ 2c_{13}^{3}$,
$ 2c_{13}^{4}$,
$0$,
$0$;\ \ 
$ 2c_{13}^{6}$,
$ 2c_{13}^{5}$,
$0$,
$0$;\ \ 
$ 2c_{13}^{2}$,
$0$,
$0$;\ \ 
$ \sqrt{13}$,
$ -\sqrt{13}$;\ \ 
$ \sqrt{13}$)

Realization:
Abelian anyon condensation of $SO(26)_2$ or $O_{13}$.

\vskip 1ex

\noindent6. $10_{4,52.}^{26,862}$ \irep{1035}:\ \ 
$d_i$ = ($1.0$,
$1.0$,
$2.0$,
$2.0$,
$2.0$,
$2.0$,
$2.0$,
$2.0$,
$3.605$,
$3.605$) 

\vskip 0.7ex
\hangindent=3em \hangafter=1
$D^2= 52.0 = 
52$

\vskip 0.7ex
\hangindent=3em \hangafter=1
$T = ( 0,
0,
\frac{2}{13},
\frac{5}{13},
\frac{6}{13},
\frac{7}{13},
\frac{8}{13},
\frac{11}{13},
0,
\frac{1}{2} )
$,

\vskip 0.7ex
\hangindent=3em \hangafter=1
$S$ = ($ 1$,
$ 1$,
$ 2$,
$ 2$,
$ 2$,
$ 2$,
$ 2$,
$ 2$,
$ \sqrt{13}$,
$ \sqrt{13}$;\ \ 
$ 1$,
$ 2$,
$ 2$,
$ 2$,
$ 2$,
$ 2$,
$ 2$,
$ -\sqrt{13}$,
$ -\sqrt{13}$;\ \ 
$ 2c_{13}^{4}$,
$ 2c_{13}^{1}$,
$ 2c_{13}^{3}$,
$ 2c_{13}^{2}$,
$ 2c_{13}^{5}$,
$ 2c_{13}^{6}$,
$0$,
$0$;\ \ 
$ 2c_{13}^{3}$,
$ 2c_{13}^{4}$,
$ 2c_{13}^{6}$,
$ 2c_{13}^{2}$,
$ 2c_{13}^{5}$,
$0$,
$0$;\ \ 
$ 2c_{13}^{1}$,
$ 2c_{13}^{5}$,
$ 2c_{13}^{6}$,
$ 2c_{13}^{2}$,
$0$,
$0$;\ \ 
$ 2c_{13}^{1}$,
$ 2c_{13}^{4}$,
$ 2c_{13}^{3}$,
$0$,
$0$;\ \ 
$ 2c_{13}^{3}$,
$ 2c_{13}^{1}$,
$0$,
$0$;\ \ 
$ 2c_{13}^{4}$,
$0$,
$0$;\ \ 
$ -\sqrt{13}$,
$ \sqrt{13}$;\ \ 
$ -\sqrt{13}$)

Realization:
Abelian anyon condensation of $SO(13)_2$.

\vskip 1ex

\noindent7. $10_{0,52.}^{52,110}$ \irep{1120}:\ \ 
$d_i$ = ($1.0$,
$1.0$,
$2.0$,
$2.0$,
$2.0$,
$2.0$,
$2.0$,
$2.0$,
$3.605$,
$3.605$) 

\vskip 0.7ex
\hangindent=3em \hangafter=1
$D^2= 52.0 = 
52$

\vskip 0.7ex
\hangindent=3em \hangafter=1
$T = ( 0,
0,
\frac{1}{13},
\frac{3}{13},
\frac{4}{13},
\frac{9}{13},
\frac{10}{13},
\frac{12}{13},
\frac{1}{4},
\frac{3}{4} )
$,

\vskip 0.7ex
\hangindent=3em \hangafter=1
$S$ = ($ 1$,
$ 1$,
$ 2$,
$ 2$,
$ 2$,
$ 2$,
$ 2$,
$ 2$,
$ \sqrt{13}$,
$ \sqrt{13}$;\ \ 
$ 1$,
$ 2$,
$ 2$,
$ 2$,
$ 2$,
$ 2$,
$ 2$,
$ -\sqrt{13}$,
$ -\sqrt{13}$;\ \ 
$ 2c_{13}^{2}$,
$ 2c_{13}^{5}$,
$ 2c_{13}^{4}$,
$ 2c_{13}^{6}$,
$ 2c_{13}^{1}$,
$ 2c_{13}^{3}$,
$0$,
$0$;\ \ 
$ 2c_{13}^{6}$,
$ 2c_{13}^{3}$,
$ 2c_{13}^{2}$,
$ 2c_{13}^{4}$,
$ 2c_{13}^{1}$,
$0$,
$0$;\ \ 
$ 2c_{13}^{5}$,
$ 2c_{13}^{1}$,
$ 2c_{13}^{2}$,
$ 2c_{13}^{6}$,
$0$,
$0$;\ \ 
$ 2c_{13}^{5}$,
$ 2c_{13}^{3}$,
$ 2c_{13}^{4}$,
$0$,
$0$;\ \ 
$ 2c_{13}^{6}$,
$ 2c_{13}^{5}$,
$0$,
$0$;\ \ 
$ 2c_{13}^{2}$,
$0$,
$0$;\ \ 
$ -\sqrt{13}$,
$ \sqrt{13}$;\ \ 
$ -\sqrt{13}$)

Realization:
Abelian anyon condensation of $SO(26)_2$ or $O_{13}$.

\vskip 1ex

\noindent8. $10_{4,52.}^{52,489}$ \irep{1120}:\ \ 
$d_i$ = ($1.0$,
$1.0$,
$2.0$,
$2.0$,
$2.0$,
$2.0$,
$2.0$,
$2.0$,
$3.605$,
$3.605$) 

\vskip 0.7ex
\hangindent=3em \hangafter=1
$D^2= 52.0 = 
52$

\vskip 0.7ex
\hangindent=3em \hangafter=1
$T = ( 0,
0,
\frac{2}{13},
\frac{5}{13},
\frac{6}{13},
\frac{7}{13},
\frac{8}{13},
\frac{11}{13},
\frac{1}{4},
\frac{3}{4} )
$,

\vskip 0.7ex
\hangindent=3em \hangafter=1
$S$ = ($ 1$,
$ 1$,
$ 2$,
$ 2$,
$ 2$,
$ 2$,
$ 2$,
$ 2$,
$ \sqrt{13}$,
$ \sqrt{13}$;\ \ 
$ 1$,
$ 2$,
$ 2$,
$ 2$,
$ 2$,
$ 2$,
$ 2$,
$ -\sqrt{13}$,
$ -\sqrt{13}$;\ \ 
$ 2c_{13}^{4}$,
$ 2c_{13}^{1}$,
$ 2c_{13}^{3}$,
$ 2c_{13}^{2}$,
$ 2c_{13}^{5}$,
$ 2c_{13}^{6}$,
$0$,
$0$;\ \ 
$ 2c_{13}^{3}$,
$ 2c_{13}^{4}$,
$ 2c_{13}^{6}$,
$ 2c_{13}^{2}$,
$ 2c_{13}^{5}$,
$0$,
$0$;\ \ 
$ 2c_{13}^{1}$,
$ 2c_{13}^{5}$,
$ 2c_{13}^{6}$,
$ 2c_{13}^{2}$,
$0$,
$0$;\ \ 
$ 2c_{13}^{1}$,
$ 2c_{13}^{4}$,
$ 2c_{13}^{3}$,
$0$,
$0$;\ \ 
$ 2c_{13}^{3}$,
$ 2c_{13}^{1}$,
$0$,
$0$;\ \ 
$ 2c_{13}^{4}$,
$0$,
$0$;\ \ 
$ \sqrt{13}$,
$ -\sqrt{13}$;\ \ 
$ \sqrt{13}$)

Realization: $SO(13)_2$.

\vskip 1ex

\noindent9. $10_{6,89.56}^{12,311}$ \irep{680}:\ \ 
$d_i$ = ($1.0$,
$1.0$,
$2.732$,
$2.732$,
$2.732$,
$2.732$,
$2.732$,
$3.732$,
$3.732$,
$4.732$) 

\vskip 0.7ex
\hangindent=3em \hangafter=1
$D^2= 89.569 = 
48+24\sqrt{3}$

\vskip 0.7ex
\hangindent=3em \hangafter=1
$T = ( 0,
\frac{1}{2},
0,
\frac{1}{3},
\frac{1}{3},
\frac{1}{3},
\frac{5}{6},
0,
\frac{1}{2},
\frac{3}{4} )
$,

\vskip 0.7ex
\hangindent=3em \hangafter=1
$S$ = ($ 1$,
$ 1$,
$ 1+\sqrt{3}$,
$ 1+\sqrt{3}$,
$ 1+\sqrt{3}$,
$ 1+\sqrt{3}$,
$ 1+\sqrt{3}$,
$ 2+\sqrt{3}$,
$ 2+\sqrt{3}$,
$ 3+\sqrt{3}$;\ \ 
$ 1$,
$ -1-\sqrt{3}$,
$ 1+\sqrt{3}$,
$ -1-\sqrt{3}$,
$ -1-\sqrt{3}$,
$ 1+\sqrt{3}$,
$ 2+\sqrt{3}$,
$ 2+\sqrt{3}$,
$ -3-\sqrt{3}$;\ \ 
$0$,
$ -2-2\sqrt{3}$,
$0$,
$0$,
$ 2+2\sqrt{3}$,
$ -1-\sqrt{3}$,
$ 1+\sqrt{3}$,
$0$;\ \ 
$ 1+\sqrt{3}$,
$ 1+\sqrt{3}$,
$ 1+\sqrt{3}$,
$ 1+\sqrt{3}$,
$ -1-\sqrt{3}$,
$ -1-\sqrt{3}$,
$0$;\ \ 
$(-3-\sqrt{3})\mathrm{i}$,
$(3+\sqrt{3})\mathrm{i}$,
$ -1-\sqrt{3}$,
$ -1-\sqrt{3}$,
$ 1+\sqrt{3}$,
$0$;\ \ 
$(-3-\sqrt{3})\mathrm{i}$,
$ -1-\sqrt{3}$,
$ -1-\sqrt{3}$,
$ 1+\sqrt{3}$,
$0$;\ \ 
$ 1+\sqrt{3}$,
$ -1-\sqrt{3}$,
$ -1-\sqrt{3}$,
$0$;\ \ 
$ 1$,
$ 1$,
$ 3+\sqrt{3}$;\ \ 
$ 1$,
$ -3-\sqrt{3}$;\ \ 
$0$)

Realization: 
Abelian anyon condensation \cite{LW170107820} of $SO(5)_3$ or $Sp(4)_3$
or $\overline{Sp(6)}_2$.

\vskip 1ex

\noindent10. $10_{0,89.56}^{12,155}$ \irep{587}:\ \ 
$d_i$ = ($1.0$,
$1.0$,
$2.732$,
$2.732$,
$2.732$,
$2.732$,
$2.732$,
$3.732$,
$3.732$,
$4.732$) 

\vskip 0.7ex
\hangindent=3em \hangafter=1
$D^2= 89.569 = 
48+24\sqrt{3}$

\vskip 0.7ex
\hangindent=3em \hangafter=1
$T = ( 0,
\frac{1}{2},
\frac{1}{3},
\frac{1}{4},
\frac{5}{6},
\frac{7}{12},
\frac{7}{12},
0,
\frac{1}{2},
0 )
$,

\vskip 0.7ex
\hangindent=3em \hangafter=1
$S$ = ($ 1$,
$ 1$,
$ 1+\sqrt{3}$,
$ 1+\sqrt{3}$,
$ 1+\sqrt{3}$,
$ 1+\sqrt{3}$,
$ 1+\sqrt{3}$,
$ 2+\sqrt{3}$,
$ 2+\sqrt{3}$,
$ 3+\sqrt{3}$;\ \ 
$ 1$,
$ 1+\sqrt{3}$,
$ -1-\sqrt{3}$,
$ 1+\sqrt{3}$,
$ -1-\sqrt{3}$,
$ -1-\sqrt{3}$,
$ 2+\sqrt{3}$,
$ 2+\sqrt{3}$,
$ -3-\sqrt{3}$;\ \ 
$ 1+\sqrt{3}$,
$ -2-2\sqrt{3}$,
$ 1+\sqrt{3}$,
$ 1+\sqrt{3}$,
$ 1+\sqrt{3}$,
$ -1-\sqrt{3}$,
$ -1-\sqrt{3}$,
$0$;\ \ 
$0$,
$ 2+2\sqrt{3}$,
$0$,
$0$,
$ -1-\sqrt{3}$,
$ 1+\sqrt{3}$,
$0$;\ \ 
$ 1+\sqrt{3}$,
$ -1-\sqrt{3}$,
$ -1-\sqrt{3}$,
$ -1-\sqrt{3}$,
$ -1-\sqrt{3}$,
$0$;\ \ 
$(3+\sqrt{3})\mathrm{i}$,
$(-3-\sqrt{3})\mathrm{i}$,
$ -1-\sqrt{3}$,
$ 1+\sqrt{3}$,
$0$;\ \ 
$(3+\sqrt{3})\mathrm{i}$,
$ -1-\sqrt{3}$,
$ 1+\sqrt{3}$,
$0$;\ \ 
$ 1$,
$ 1$,
$ 3+\sqrt{3}$;\ \ 
$ 1$,
$ -3-\sqrt{3}$;\ \ 
$0$)

Realization: 
Abelian anyon condensation of $SO(5)_3$ or $Sp(4)_3$ or $\overline{Sp(6)}_2$.

\vskip 1ex

\noindent11. $10_{4,89.56}^{12,822}$ \irep{587}:\ \ 
$d_i$ = ($1.0$,
$1.0$,
$2.732$,
$2.732$,
$2.732$,
$2.732$,
$2.732$,
$3.732$,
$3.732$,
$4.732$) 

\vskip 0.7ex
\hangindent=3em \hangafter=1
$D^2= 89.569 = 
48+24\sqrt{3}$

\vskip 0.7ex
\hangindent=3em \hangafter=1
$T = ( 0,
\frac{1}{2},
\frac{1}{3},
\frac{3}{4},
\frac{5}{6},
\frac{1}{12},
\frac{1}{12},
0,
\frac{1}{2},
\frac{1}{2} )
$,

\vskip 0.7ex
\hangindent=3em \hangafter=1
$S$ = ($ 1$,
$ 1$,
$ 1+\sqrt{3}$,
$ 1+\sqrt{3}$,
$ 1+\sqrt{3}$,
$ 1+\sqrt{3}$,
$ 1+\sqrt{3}$,
$ 2+\sqrt{3}$,
$ 2+\sqrt{3}$,
$ 3+\sqrt{3}$;\ \ 
$ 1$,
$ 1+\sqrt{3}$,
$ -1-\sqrt{3}$,
$ 1+\sqrt{3}$,
$ -1-\sqrt{3}$,
$ -1-\sqrt{3}$,
$ 2+\sqrt{3}$,
$ 2+\sqrt{3}$,
$ -3-\sqrt{3}$;\ \ 
$ 1+\sqrt{3}$,
$ -2-2\sqrt{3}$,
$ 1+\sqrt{3}$,
$ 1+\sqrt{3}$,
$ 1+\sqrt{3}$,
$ -1-\sqrt{3}$,
$ -1-\sqrt{3}$,
$0$;\ \ 
$0$,
$ 2+2\sqrt{3}$,
$0$,
$0$,
$ -1-\sqrt{3}$,
$ 1+\sqrt{3}$,
$0$;\ \ 
$ 1+\sqrt{3}$,
$ -1-\sqrt{3}$,
$ -1-\sqrt{3}$,
$ -1-\sqrt{3}$,
$ -1-\sqrt{3}$,
$0$;\ \ 
$(3+\sqrt{3})\mathrm{i}$,
$(-3-\sqrt{3})\mathrm{i}$,
$ -1-\sqrt{3}$,
$ 1+\sqrt{3}$,
$0$;\ \ 
$(3+\sqrt{3})\mathrm{i}$,
$ -1-\sqrt{3}$,
$ 1+\sqrt{3}$,
$0$;\ \ 
$ 1$,
$ 1$,
$ 3+\sqrt{3}$;\ \ 
$ 1$,
$ -3-\sqrt{3}$;\ \ 
$0$)

Realization: 
Abelian anyon condensation of $SO(5)_3$ or $Sp(4)_3$  or $\overline{Sp(6)}_2$.

\vskip 1ex

\noindent12. $10_{2,89.56}^{12,119}$ \irep{680}:\ \ 
$d_i$ = ($1.0$,
$1.0$,
$2.732$,
$2.732$,
$2.732$,
$2.732$,
$2.732$,
$3.732$,
$3.732$,
$4.732$) 

\vskip 0.7ex
\hangindent=3em \hangafter=1
$D^2= 89.569 = 
48+24\sqrt{3}$

\vskip 0.7ex
\hangindent=3em \hangafter=1
$T = ( 0,
\frac{1}{2},
\frac{1}{2},
\frac{1}{3},
\frac{5}{6},
\frac{5}{6},
\frac{5}{6},
0,
\frac{1}{2},
\frac{1}{4} )
$,

\vskip 0.7ex
\hangindent=3em \hangafter=1
$S$ = ($ 1$,
$ 1$,
$ 1+\sqrt{3}$,
$ 1+\sqrt{3}$,
$ 1+\sqrt{3}$,
$ 1+\sqrt{3}$,
$ 1+\sqrt{3}$,
$ 2+\sqrt{3}$,
$ 2+\sqrt{3}$,
$ 3+\sqrt{3}$;\ \ 
$ 1$,
$ -1-\sqrt{3}$,
$ 1+\sqrt{3}$,
$ 1+\sqrt{3}$,
$ -1-\sqrt{3}$,
$ -1-\sqrt{3}$,
$ 2+\sqrt{3}$,
$ 2+\sqrt{3}$,
$ -3-\sqrt{3}$;\ \ 
$0$,
$ -2-2\sqrt{3}$,
$ 2+2\sqrt{3}$,
$0$,
$0$,
$ -1-\sqrt{3}$,
$ 1+\sqrt{3}$,
$0$;\ \ 
$ 1+\sqrt{3}$,
$ 1+\sqrt{3}$,
$ 1+\sqrt{3}$,
$ 1+\sqrt{3}$,
$ -1-\sqrt{3}$,
$ -1-\sqrt{3}$,
$0$;\ \ 
$ 1+\sqrt{3}$,
$ -1-\sqrt{3}$,
$ -1-\sqrt{3}$,
$ -1-\sqrt{3}$,
$ -1-\sqrt{3}$,
$0$;\ \ 
$(-3-\sqrt{3})\mathrm{i}$,
$(3+\sqrt{3})\mathrm{i}$,
$ -1-\sqrt{3}$,
$ 1+\sqrt{3}$,
$0$;\ \ 
$(-3-\sqrt{3})\mathrm{i}$,
$ -1-\sqrt{3}$,
$ 1+\sqrt{3}$,
$0$;\ \ 
$ 1$,
$ 1$,
$ 3+\sqrt{3}$;\ \ 
$ 1$,
$ -3-\sqrt{3}$;\ \ 
$0$)

Realization: 
Abelian anyon condensation of $SO(5)_3$ or $Sp(4)_3$  or $\overline{Sp(6)}_2$.

\vskip 1ex

\noindent13. $10_{7,89.56}^{24,123}$ \irep{978}:\ \ 
$d_i$ = ($1.0$,
$1.0$,
$2.732$,
$2.732$,
$2.732$,
$2.732$,
$2.732$,
$3.732$,
$3.732$,
$4.732$) 

\vskip 0.7ex
\hangindent=3em \hangafter=1
$D^2= 89.569 = 
48+24\sqrt{3}$

\vskip 0.7ex
\hangindent=3em \hangafter=1
$T = ( 0,
\frac{1}{2},
\frac{1}{3},
\frac{5}{6},
\frac{1}{8},
\frac{11}{24},
\frac{11}{24},
0,
\frac{1}{2},
\frac{7}{8} )
$,

\vskip 0.7ex
\hangindent=3em \hangafter=1
$S$ = ($ 1$,
$ 1$,
$ 1+\sqrt{3}$,
$ 1+\sqrt{3}$,
$ 1+\sqrt{3}$,
$ 1+\sqrt{3}$,
$ 1+\sqrt{3}$,
$ 2+\sqrt{3}$,
$ 2+\sqrt{3}$,
$ 3+\sqrt{3}$;\ \ 
$ 1$,
$ 1+\sqrt{3}$,
$ 1+\sqrt{3}$,
$ -1-\sqrt{3}$,
$ -1-\sqrt{3}$,
$ -1-\sqrt{3}$,
$ 2+\sqrt{3}$,
$ 2+\sqrt{3}$,
$ -3-\sqrt{3}$;\ \ 
$ 1+\sqrt{3}$,
$ 1+\sqrt{3}$,
$ -2-2\sqrt{3}$,
$ 1+\sqrt{3}$,
$ 1+\sqrt{3}$,
$ -1-\sqrt{3}$,
$ -1-\sqrt{3}$,
$0$;\ \ 
$ 1+\sqrt{3}$,
$ 2+2\sqrt{3}$,
$ -1-\sqrt{3}$,
$ -1-\sqrt{3}$,
$ -1-\sqrt{3}$,
$ -1-\sqrt{3}$,
$0$;\ \ 
$0$,
$0$,
$0$,
$ -1-\sqrt{3}$,
$ 1+\sqrt{3}$,
$0$;\ \ 
$ -3-\sqrt{3}$,
$ 3+\sqrt{3}$,
$ -1-\sqrt{3}$,
$ 1+\sqrt{3}$,
$0$;\ \ 
$ -3-\sqrt{3}$,
$ -1-\sqrt{3}$,
$ 1+\sqrt{3}$,
$0$;\ \ 
$ 1$,
$ 1$,
$ 3+\sqrt{3}$;\ \ 
$ 1$,
$ -3-\sqrt{3}$;\ \ 
$0$)

Realization: 
Abelian anyon condensation of $SO(5)_3$ or $Sp(4)_3$  or $\overline{Sp(6)}_2$.

\vskip 1ex

\noindent14. $10_{1,89.56}^{24,380}$ \irep{977}:\ \ 
$d_i$ = ($1.0$,
$1.0$,
$2.732$,
$2.732$,
$2.732$,
$2.732$,
$2.732$,
$3.732$,
$3.732$,
$4.732$) 

\vskip 0.7ex
\hangindent=3em \hangafter=1
$D^2= 89.569 = 
48+24\sqrt{3}$

\vskip 0.7ex
\hangindent=3em \hangafter=1
$T = ( 0,
\frac{1}{2},
\frac{1}{3},
\frac{5}{6},
\frac{3}{8},
\frac{17}{24},
\frac{17}{24},
0,
\frac{1}{2},
\frac{1}{8} )
$,

\vskip 0.7ex
\hangindent=3em \hangafter=1
$S$ = ($ 1$,
$ 1$,
$ 1+\sqrt{3}$,
$ 1+\sqrt{3}$,
$ 1+\sqrt{3}$,
$ 1+\sqrt{3}$,
$ 1+\sqrt{3}$,
$ 2+\sqrt{3}$,
$ 2+\sqrt{3}$,
$ 3+\sqrt{3}$;\ \ 
$ 1$,
$ 1+\sqrt{3}$,
$ 1+\sqrt{3}$,
$ -1-\sqrt{3}$,
$ -1-\sqrt{3}$,
$ -1-\sqrt{3}$,
$ 2+\sqrt{3}$,
$ 2+\sqrt{3}$,
$ -3-\sqrt{3}$;\ \ 
$ 1+\sqrt{3}$,
$ 1+\sqrt{3}$,
$ -2-2\sqrt{3}$,
$ 1+\sqrt{3}$,
$ 1+\sqrt{3}$,
$ -1-\sqrt{3}$,
$ -1-\sqrt{3}$,
$0$;\ \ 
$ 1+\sqrt{3}$,
$ 2+2\sqrt{3}$,
$ -1-\sqrt{3}$,
$ -1-\sqrt{3}$,
$ -1-\sqrt{3}$,
$ -1-\sqrt{3}$,
$0$;\ \ 
$0$,
$0$,
$0$,
$ -1-\sqrt{3}$,
$ 1+\sqrt{3}$,
$0$;\ \ 
$ 3+\sqrt{3}$,
$ -3-\sqrt{3}$,
$ -1-\sqrt{3}$,
$ 1+\sqrt{3}$,
$0$;\ \ 
$ 3+\sqrt{3}$,
$ -1-\sqrt{3}$,
$ 1+\sqrt{3}$,
$0$;\ \ 
$ 1$,
$ 1$,
$ 3+\sqrt{3}$;\ \ 
$ 1$,
$ -3-\sqrt{3}$;\ \ 
$0$)

Realization: 
$\overline{Sp(6)}_2$.

\vskip 1ex

\noindent15. $10_{\frac{26}{7},236.3}^{21,145}$ \irep{960}:\ \ 
$d_i$ = ($1.0$,
$1.977$,
$2.911$,
$3.779$,
$4.563$,
$5.245$,
$5.810$,
$6.245$,
$6.541$,
$6.690$) 

\vskip 0.7ex
\hangindent=3em \hangafter=1
$D^2= 236.341 = 
42+42c^{1}_{21}
+42c^{2}_{21}
+21c^{3}_{21}
+21c^{4}_{21}
+21c^{5}_{21}
$

\vskip 0.7ex
\hangindent=3em \hangafter=1
$T = ( 0,
\frac{2}{7},
\frac{2}{21},
\frac{3}{7},
\frac{2}{7},
\frac{2}{3},
\frac{4}{7},
0,
\frac{20}{21},
\frac{3}{7} )
$,

\vskip 0.7ex
\hangindent=3em \hangafter=1
$S$ = ($ 1$,
$ -c_{21}^{10}$,
$ \xi_{21}^{3}$,
$ \xi_{21}^{17}$,
$ \xi_{21}^{5}$,
$ \xi_{21}^{15}$,
$ \xi_{21}^{7}$,
$ \xi_{21}^{13}$,
$ \xi_{21}^{9}$,
$ \xi_{21}^{11}$;\ \ 
$ -\xi_{21}^{17}$,
$ \xi_{21}^{15}$,
$ -\xi_{21}^{13}$,
$ \xi_{21}^{11}$,
$ -\xi_{21}^{9}$,
$ \xi_{21}^{7}$,
$ -\xi_{21}^{5}$,
$ \xi_{21}^{3}$,
$ -1$;\ \ 
$ \xi_{21}^{9}$,
$ \xi_{21}^{9}$,
$ \xi_{21}^{15}$,
$ \xi_{21}^{3}$,
$0$,
$ -\xi_{21}^{3}$,
$ -\xi_{21}^{15}$,
$ -\xi_{21}^{9}$;\ \ 
$ -\xi_{21}^{5}$,
$ 1$,
$ \xi_{21}^{3}$,
$ -\xi_{21}^{7}$,
$ \xi_{21}^{11}$,
$ -\xi_{21}^{15}$,
$ -c_{21}^{10}$;\ \ 
$ -\xi_{21}^{17}$,
$ -\xi_{21}^{9}$,
$ -\xi_{21}^{7}$,
$ c_{21}^{10}$,
$ \xi_{21}^{3}$,
$ \xi_{21}^{13}$;\ \ 
$ \xi_{21}^{15}$,
$0$,
$ -\xi_{21}^{15}$,
$ \xi_{21}^{9}$,
$ -\xi_{21}^{3}$;\ \ 
$ \xi_{21}^{7}$,
$ \xi_{21}^{7}$,
$0$,
$ -\xi_{21}^{7}$;\ \ 
$ 1$,
$ -\xi_{21}^{9}$,
$ \xi_{21}^{17}$;\ \ 
$ -\xi_{21}^{3}$,
$ \xi_{21}^{15}$;\ \ 
$ -\xi_{21}^{5}$)

Realization: $PSU(2)_{19}$.

\vskip 1ex

\noindent16. $10_{\frac{48}{17},499.2}^{17,522}$ \irep{888}:\ \ 
$d_i$ = ($1.0$,
$2.965$,
$4.830$,
$5.418$,
$5.418$,
$6.531$,
$8.9$,
$9.214$,
$10.106$,
$10.653$) 

\vskip 0.7ex
\hangindent=3em \hangafter=1
$D^2= 499.210 = 
136+119c^{1}_{17}
+102c^{2}_{17}
+85c^{3}_{17}
+68c^{4}_{17}
+51c^{5}_{17}
+34c^{6}_{17}
+17c^{7}_{17}
$

\vskip 0.7ex
\hangindent=3em \hangafter=1
$T = ( 0,
\frac{1}{17},
\frac{3}{17},
\frac{2}{17},
\frac{2}{17},
\frac{6}{17},
\frac{10}{17},
\frac{15}{17},
\frac{4}{17},
\frac{11}{17} )
$,

\vskip 0.7ex
\hangindent=3em \hangafter=1
$S$ = ($ 1$,
$ 2+c^{1}_{17}
+c^{2}_{17}
+c^{3}_{17}
+c^{4}_{17}
+c^{5}_{17}
+c^{6}_{17}
+c^{7}_{17}
$,
$ 2+2c^{1}_{17}
+c^{2}_{17}
+c^{3}_{17}
+c^{4}_{17}
+c^{5}_{17}
+c^{6}_{17}
+c^{7}_{17}
$,
$ \xi_{17}^{9}$,
$ \xi_{17}^{9}$,
$ 2+2c^{1}_{17}
+c^{2}_{17}
+c^{3}_{17}
+c^{4}_{17}
+c^{5}_{17}
+c^{6}_{17}
$,
$ 2+2c^{1}_{17}
+2c^{2}_{17}
+c^{3}_{17}
+c^{4}_{17}
+c^{5}_{17}
+c^{6}_{17}
$,
$ 2+2c^{1}_{17}
+2c^{2}_{17}
+c^{3}_{17}
+c^{4}_{17}
+c^{5}_{17}
$,
$ 2+2c^{1}_{17}
+2c^{2}_{17}
+2c^{3}_{17}
+c^{4}_{17}
+c^{5}_{17}
$,
$ 2+2c^{1}_{17}
+2c^{2}_{17}
+2c^{3}_{17}
+c^{4}_{17}
$;\ \ 
$ 2+2c^{1}_{17}
+2c^{2}_{17}
+c^{3}_{17}
+c^{4}_{17}
+c^{5}_{17}
+c^{6}_{17}
$,
$ 2+2c^{1}_{17}
+2c^{2}_{17}
+2c^{3}_{17}
+c^{4}_{17}
$,
$ -\xi_{17}^{9}$,
$ -\xi_{17}^{9}$,
$ 2+2c^{1}_{17}
+2c^{2}_{17}
+2c^{3}_{17}
+c^{4}_{17}
+c^{5}_{17}
$,
$ 2+2c^{1}_{17}
+c^{2}_{17}
+c^{3}_{17}
+c^{4}_{17}
+c^{5}_{17}
+c^{6}_{17}
$,
$ 1$,
$ -2-2  c^{1}_{17}
-c^{2}_{17}
-c^{3}_{17}
-c^{4}_{17}
-c^{5}_{17}
-c^{6}_{17}
-c^{7}_{17}
$,
$ -2-2  c^{1}_{17}
-2  c^{2}_{17}
-c^{3}_{17}
-c^{4}_{17}
-c^{5}_{17}
$;\ \ 
$ 2+2c^{1}_{17}
+2c^{2}_{17}
+c^{3}_{17}
+c^{4}_{17}
+c^{5}_{17}
+c^{6}_{17}
$,
$ \xi_{17}^{9}$,
$ \xi_{17}^{9}$,
$ -1$,
$ -2-2  c^{1}_{17}
-2  c^{2}_{17}
-c^{3}_{17}
-c^{4}_{17}
-c^{5}_{17}
$,
$ -2-2  c^{1}_{17}
-2  c^{2}_{17}
-2  c^{3}_{17}
-c^{4}_{17}
-c^{5}_{17}
$,
$ -2-c^{1}_{17}
-c^{2}_{17}
-c^{3}_{17}
-c^{4}_{17}
-c^{5}_{17}
-c^{6}_{17}
-c^{7}_{17}
$,
$ 2+2c^{1}_{17}
+c^{2}_{17}
+c^{3}_{17}
+c^{4}_{17}
+c^{5}_{17}
+c^{6}_{17}
$;\ \ 
$ 4+3c^{1}_{17}
+3c^{2}_{17}
+2c^{3}_{17}
+2c^{4}_{17}
+2c^{5}_{17}
+c^{6}_{17}
$,
$ -3-2  c^{1}_{17}
-2  c^{2}_{17}
-c^{3}_{17}
-c^{4}_{17}
-2  c^{5}_{17}
-c^{6}_{17}
$,
$ -\xi_{17}^{9}$,
$ \xi_{17}^{9}$,
$ -\xi_{17}^{9}$,
$ \xi_{17}^{9}$,
$ -\xi_{17}^{9}$;\ \ 
$ 4+3c^{1}_{17}
+3c^{2}_{17}
+2c^{3}_{17}
+2c^{4}_{17}
+2c^{5}_{17}
+c^{6}_{17}
$,
$ -\xi_{17}^{9}$,
$ \xi_{17}^{9}$,
$ -\xi_{17}^{9}$,
$ \xi_{17}^{9}$,
$ -\xi_{17}^{9}$;\ \ 
$ -2-2  c^{1}_{17}
-2  c^{2}_{17}
-2  c^{3}_{17}
-c^{4}_{17}
$,
$ -2-2  c^{1}_{17}
-c^{2}_{17}
-c^{3}_{17}
-c^{4}_{17}
-c^{5}_{17}
-c^{6}_{17}
-c^{7}_{17}
$,
$ 2+2c^{1}_{17}
+2c^{2}_{17}
+c^{3}_{17}
+c^{4}_{17}
+c^{5}_{17}
+c^{6}_{17}
$,
$ 2+2c^{1}_{17}
+2c^{2}_{17}
+c^{3}_{17}
+c^{4}_{17}
+c^{5}_{17}
$,
$ -2-c^{1}_{17}
-c^{2}_{17}
-c^{3}_{17}
-c^{4}_{17}
-c^{5}_{17}
-c^{6}_{17}
-c^{7}_{17}
$;\ \ 
$ 2+2c^{1}_{17}
+2c^{2}_{17}
+2c^{3}_{17}
+c^{4}_{17}
+c^{5}_{17}
$,
$ 2+c^{1}_{17}
+c^{2}_{17}
+c^{3}_{17}
+c^{4}_{17}
+c^{5}_{17}
+c^{6}_{17}
+c^{7}_{17}
$,
$ -2-2  c^{1}_{17}
-2  c^{2}_{17}
-2  c^{3}_{17}
-c^{4}_{17}
$,
$ -1$;\ \ 
$ -2-2  c^{1}_{17}
-2  c^{2}_{17}
-2  c^{3}_{17}
-c^{4}_{17}
$,
$ 2+2c^{1}_{17}
+c^{2}_{17}
+c^{3}_{17}
+c^{4}_{17}
+c^{5}_{17}
+c^{6}_{17}
$,
$ 2+2c^{1}_{17}
+c^{2}_{17}
+c^{3}_{17}
+c^{4}_{17}
+c^{5}_{17}
+c^{6}_{17}
+c^{7}_{17}
$;\ \ 
$ 1$,
$ -2-2  c^{1}_{17}
-2  c^{2}_{17}
-c^{3}_{17}
-c^{4}_{17}
-c^{5}_{17}
-c^{6}_{17}
$;\ \ 
$ 2+2c^{1}_{17}
+2c^{2}_{17}
+2c^{3}_{17}
+c^{4}_{17}
+c^{5}_{17}
$)

Realization: $PSO(16)_3$.
Abelian anyon condensation of $E(8)_4$.

\vskip 1ex

\noindent17. $10_{0,537.4}^{14,352}$ \irep{783}:\ \ 
$d_i$ = ($1.0$,
$3.493$,
$4.493$,
$4.493$,
$5.603$,
$5.603$,
$9.97$,
$10.97$,
$10.97$,
$11.591$) 

\vskip 0.7ex
\hangindent=3em \hangafter=1
$D^2= 537.478 = 
308+224c^{1}_{7}
+112c^{2}_{7}
$

\vskip 0.7ex
\hangindent=3em \hangafter=1
$T = ( 0,
0,
\frac{2}{7},
\frac{5}{7},
\frac{3}{7},
\frac{4}{7},
0,
\frac{1}{7},
\frac{6}{7},
\frac{1}{2} )
$,

\vskip 0.7ex
\hangindent=3em \hangafter=1
$S$ = ($ 1$,
$ 1+2c^{1}_{7}
$,
$ 2\xi_{7}^{3}$,
$ 2\xi_{7}^{3}$,
$ 4+2c^{1}_{7}
+2c^{2}_{7}
$,
$ 4+2c^{1}_{7}
+2c^{2}_{7}
$,
$ 5+4c^{1}_{7}
+2c^{2}_{7}
$,
$ 6+4c^{1}_{7}
+2c^{2}_{7}
$,
$ 6+4c^{1}_{7}
+2c^{2}_{7}
$,
$ 5+6c^{1}_{7}
+2c^{2}_{7}
$;\ \ 
$ -5-4  c^{1}_{7}
-2  c^{2}_{7}
$,
$ -4-2  c^{1}_{7}
-2  c^{2}_{7}
$,
$ -4-2  c^{1}_{7}
-2  c^{2}_{7}
$,
$ 6+4c^{1}_{7}
+2c^{2}_{7}
$,
$ 6+4c^{1}_{7}
+2c^{2}_{7}
$,
$ 1$,
$ 2\xi_{7}^{3}$,
$ 2\xi_{7}^{3}$,
$ -5-6  c^{1}_{7}
-2  c^{2}_{7}
$;\ \ 
$ -2\xi_{7}^{3}$,
$ 6+6c^{1}_{7}
+2c^{2}_{7}
$,
$ 4+4c^{1}_{7}
+2c^{2}_{7}
$,
$ -4-2  c^{1}_{7}
-2  c^{2}_{7}
$,
$ 6+4c^{1}_{7}
+2c^{2}_{7}
$,
$ -6-4  c^{1}_{7}
-2  c^{2}_{7}
$,
$ -2c_{7}^{1}$,
$0$;\ \ 
$ -2\xi_{7}^{3}$,
$ -4-2  c^{1}_{7}
-2  c^{2}_{7}
$,
$ 4+4c^{1}_{7}
+2c^{2}_{7}
$,
$ 6+4c^{1}_{7}
+2c^{2}_{7}
$,
$ -2c_{7}^{1}$,
$ -6-4  c^{1}_{7}
-2  c^{2}_{7}
$,
$0$;\ \ 
$ 6+4c^{1}_{7}
+2c^{2}_{7}
$,
$ 2c_{7}^{1}$,
$ -2\xi_{7}^{3}$,
$ -6-6  c^{1}_{7}
-2  c^{2}_{7}
$,
$ 2\xi_{7}^{3}$,
$0$;\ \ 
$ 6+4c^{1}_{7}
+2c^{2}_{7}
$,
$ -2\xi_{7}^{3}$,
$ 2\xi_{7}^{3}$,
$ -6-6  c^{1}_{7}
-2  c^{2}_{7}
$,
$0$;\ \ 
$ -1-2  c^{1}_{7}
$,
$ 4+2c^{1}_{7}
+2c^{2}_{7}
$,
$ 4+2c^{1}_{7}
+2c^{2}_{7}
$,
$ -5-6  c^{1}_{7}
-2  c^{2}_{7}
$;\ \ 
$ -4-2  c^{1}_{7}
-2  c^{2}_{7}
$,
$ 4+4c^{1}_{7}
+2c^{2}_{7}
$,
$0$;\ \ 
$ -4-2  c^{1}_{7}
-2  c^{2}_{7}
$,
$0$;\ \ 
$ 5+6c^{1}_{7}
+2c^{2}_{7}
$)

Realization:  $Sp(6)_3$.

\vskip 1ex

\noindent18. $10_{6,684.3}^{77,298}$ \irep{1138}:\ \ 
$d_i$ = ($1.0$,
$7.887$,
$7.887$,
$7.887$,
$7.887$,
$7.887$,
$8.887$,
$9.887$,
$9.887$,
$9.887$) 

\vskip 0.7ex
\hangindent=3em \hangafter=1
$D^2= 684.336 = 
\frac{693+77\sqrt{77}}{2}$

\vskip 0.7ex
\hangindent=3em \hangafter=1
$T = ( 0,
\frac{1}{11},
\frac{3}{11},
\frac{4}{11},
\frac{5}{11},
\frac{9}{11},
0,
\frac{3}{7},
\frac{5}{7},
\frac{6}{7} )
$,

\vskip 0.7ex
\hangindent=3em \hangafter=1
$S$ = ($ 1$,
$ \frac{7+\sqrt{77}}{2}$,
$ \frac{7+\sqrt{77}}{2}$,
$ \frac{7+\sqrt{77}}{2}$,
$ \frac{7+\sqrt{77}}{2}$,
$ \frac{7+\sqrt{77}}{2}$,
$ \frac{9+\sqrt{77}}{2}$,
$ \frac{11+\sqrt{77}}{2}$,
$ \frac{11+\sqrt{77}}{2}$,
$ \frac{11+\sqrt{77}}{2}$;\ \ 
$ -1-2  c^{1}_{77}
+c^{2}_{77}
-c^{3}_{77}
-c^{4}_{77}
+c^{5}_{77}
+2c^{6}_{77}
-c^{7}_{77}
-c^{8}_{77}
+c^{9}_{77}
-2  c^{10}_{77}
-c^{11}_{77}
+c^{12}_{77}
-4  c^{14}_{77}
-c^{15}_{77}
+3c^{16}_{77}
+2c^{17}_{77}
-c^{18}_{77}
+c^{19}_{77}
-c^{21}_{77}
-c^{22}_{77}
-c^{23}_{77}
+c^{28}_{77}
-c^{29}_{77}
$,
$ -2+2c^{2}_{77}
-2  c^{3}_{77}
-2  c^{4}_{77}
+2c^{5}_{77}
-2  c^{6}_{77}
-6  c^{7}_{77}
-2  c^{8}_{77}
+2c^{9}_{77}
-2  c^{11}_{77}
+2c^{12}_{77}
-c^{14}_{77}
-2  c^{15}_{77}
-2  c^{17}_{77}
-3  c^{18}_{77}
+2c^{19}_{77}
-2  c^{22}_{77}
+2c^{23}_{77}
-c^{26}_{77}
-c^{28}_{77}
-3  c^{29}_{77}
$,
$ 1+2c^{1}_{77}
+c^{6}_{77}
+c^{7}_{77}
-2  c^{9}_{77}
+2c^{10}_{77}
-2  c^{13}_{77}
+c^{14}_{77}
+c^{16}_{77}
+c^{17}_{77}
+2c^{21}_{77}
+2c^{23}_{77}
-2  c^{24}_{77}
-2  c^{28}_{77}
$,
$ -1+c^{1}_{77}
+2c^{9}_{77}
+c^{10}_{77}
+2c^{13}_{77}
-c^{14}_{77}
+2c^{18}_{77}
-4  c^{21}_{77}
+c^{23}_{77}
+2c^{24}_{77}
+2c^{26}_{77}
-c^{28}_{77}
+2c^{29}_{77}
$,
$ 5-2  c^{2}_{77}
+2c^{3}_{77}
+2c^{4}_{77}
-2  c^{5}_{77}
+4c^{7}_{77}
+2c^{8}_{77}
-c^{9}_{77}
+2c^{11}_{77}
-2  c^{12}_{77}
+c^{13}_{77}
+4c^{14}_{77}
+2c^{15}_{77}
-2  c^{16}_{77}
-2  c^{19}_{77}
+3c^{21}_{77}
+2c^{22}_{77}
-2  c^{23}_{77}
+c^{24}_{77}
-2  c^{26}_{77}
+3c^{28}_{77}
$,
$ -\frac{7+\sqrt{77}}{2}$,
$0$,
$0$,
$0$;\ \ 
$ 5-2  c^{2}_{77}
+2c^{3}_{77}
+2c^{4}_{77}
-2  c^{5}_{77}
+4c^{7}_{77}
+2c^{8}_{77}
-c^{9}_{77}
+2c^{11}_{77}
-2  c^{12}_{77}
+c^{13}_{77}
+4c^{14}_{77}
+2c^{15}_{77}
-2  c^{16}_{77}
-2  c^{19}_{77}
+3c^{21}_{77}
+2c^{22}_{77}
-2  c^{23}_{77}
+c^{24}_{77}
-2  c^{26}_{77}
+3c^{28}_{77}
$,
$ -1-2  c^{1}_{77}
+c^{2}_{77}
-c^{3}_{77}
-c^{4}_{77}
+c^{5}_{77}
+2c^{6}_{77}
-c^{7}_{77}
-c^{8}_{77}
+c^{9}_{77}
-2  c^{10}_{77}
-c^{11}_{77}
+c^{12}_{77}
-4  c^{14}_{77}
-c^{15}_{77}
+3c^{16}_{77}
+2c^{17}_{77}
-c^{18}_{77}
+c^{19}_{77}
-c^{21}_{77}
-c^{22}_{77}
-c^{23}_{77}
+c^{28}_{77}
-c^{29}_{77}
$,
$ 1+2c^{1}_{77}
+c^{6}_{77}
+c^{7}_{77}
-2  c^{9}_{77}
+2c^{10}_{77}
-2  c^{13}_{77}
+c^{14}_{77}
+c^{16}_{77}
+c^{17}_{77}
+2c^{21}_{77}
+2c^{23}_{77}
-2  c^{24}_{77}
-2  c^{28}_{77}
$,
$ -1+c^{1}_{77}
+2c^{9}_{77}
+c^{10}_{77}
+2c^{13}_{77}
-c^{14}_{77}
+2c^{18}_{77}
-4  c^{21}_{77}
+c^{23}_{77}
+2c^{24}_{77}
+2c^{26}_{77}
-c^{28}_{77}
+2c^{29}_{77}
$,
$ -\frac{7+\sqrt{77}}{2}$,
$0$,
$0$,
$0$;\ \ 
$ -1+c^{1}_{77}
+2c^{9}_{77}
+c^{10}_{77}
+2c^{13}_{77}
-c^{14}_{77}
+2c^{18}_{77}
-4  c^{21}_{77}
+c^{23}_{77}
+2c^{24}_{77}
+2c^{26}_{77}
-c^{28}_{77}
+2c^{29}_{77}
$,
$ 5-2  c^{2}_{77}
+2c^{3}_{77}
+2c^{4}_{77}
-2  c^{5}_{77}
+4c^{7}_{77}
+2c^{8}_{77}
-c^{9}_{77}
+2c^{11}_{77}
-2  c^{12}_{77}
+c^{13}_{77}
+4c^{14}_{77}
+2c^{15}_{77}
-2  c^{16}_{77}
-2  c^{19}_{77}
+3c^{21}_{77}
+2c^{22}_{77}
-2  c^{23}_{77}
+c^{24}_{77}
-2  c^{26}_{77}
+3c^{28}_{77}
$,
$ -2+2c^{2}_{77}
-2  c^{3}_{77}
-2  c^{4}_{77}
+2c^{5}_{77}
-2  c^{6}_{77}
-6  c^{7}_{77}
-2  c^{8}_{77}
+2c^{9}_{77}
-2  c^{11}_{77}
+2c^{12}_{77}
-c^{14}_{77}
-2  c^{15}_{77}
-2  c^{17}_{77}
-3  c^{18}_{77}
+2c^{19}_{77}
-2  c^{22}_{77}
+2c^{23}_{77}
-c^{26}_{77}
-c^{28}_{77}
-3  c^{29}_{77}
$,
$ -\frac{7+\sqrt{77}}{2}$,
$0$,
$0$,
$0$;\ \ 
$ -2+2c^{2}_{77}
-2  c^{3}_{77}
-2  c^{4}_{77}
+2c^{5}_{77}
-2  c^{6}_{77}
-6  c^{7}_{77}
-2  c^{8}_{77}
+2c^{9}_{77}
-2  c^{11}_{77}
+2c^{12}_{77}
-c^{14}_{77}
-2  c^{15}_{77}
-2  c^{17}_{77}
-3  c^{18}_{77}
+2c^{19}_{77}
-2  c^{22}_{77}
+2c^{23}_{77}
-c^{26}_{77}
-c^{28}_{77}
-3  c^{29}_{77}
$,
$ -1-2  c^{1}_{77}
+c^{2}_{77}
-c^{3}_{77}
-c^{4}_{77}
+c^{5}_{77}
+2c^{6}_{77}
-c^{7}_{77}
-c^{8}_{77}
+c^{9}_{77}
-2  c^{10}_{77}
-c^{11}_{77}
+c^{12}_{77}
-4  c^{14}_{77}
-c^{15}_{77}
+3c^{16}_{77}
+2c^{17}_{77}
-c^{18}_{77}
+c^{19}_{77}
-c^{21}_{77}
-c^{22}_{77}
-c^{23}_{77}
+c^{28}_{77}
-c^{29}_{77}
$,
$ -\frac{7+\sqrt{77}}{2}$,
$0$,
$0$,
$0$;\ \ 
$ 1+2c^{1}_{77}
+c^{6}_{77}
+c^{7}_{77}
-2  c^{9}_{77}
+2c^{10}_{77}
-2  c^{13}_{77}
+c^{14}_{77}
+c^{16}_{77}
+c^{17}_{77}
+2c^{21}_{77}
+2c^{23}_{77}
-2  c^{24}_{77}
-2  c^{28}_{77}
$,
$ -\frac{7+\sqrt{77}}{2}$,
$0$,
$0$,
$0$;\ \ 
$ 1$,
$ \frac{11+\sqrt{77}}{2}$,
$ \frac{11+\sqrt{77}}{2}$,
$ \frac{11+\sqrt{77}}{2}$;\ \ 
$ 1+3c^{4}_{77}
+2c^{7}_{77}
-2  c^{9}_{77}
+c^{10}_{77}
+7c^{11}_{77}
+2c^{15}_{77}
-2  c^{16}_{77}
+c^{17}_{77}
+2c^{18}_{77}
-2  c^{19}_{77}
+c^{22}_{77}
-2  c^{23}_{77}
+c^{24}_{77}
+c^{25}_{77}
+2c^{26}_{77}
+2c^{29}_{77}
$,
$ -1+c^{1}_{77}
+c^{2}_{77}
-c^{3}_{77}
-3  c^{4}_{77}
+c^{5}_{77}
+c^{6}_{77}
-c^{7}_{77}
-c^{8}_{77}
+c^{9}_{77}
-2  c^{10}_{77}
-2  c^{11}_{77}
+c^{12}_{77}
+c^{13}_{77}
-c^{14}_{77}
+c^{16}_{77}
-2  c^{17}_{77}
-c^{18}_{77}
+5c^{22}_{77}
+c^{23}_{77}
-2  c^{24}_{77}
-3  c^{25}_{77}
-c^{29}_{77}
$,
$ -4-2  c^{1}_{77}
-2  c^{2}_{77}
+2c^{3}_{77}
+c^{4}_{77}
-2  c^{5}_{77}
-2  c^{6}_{77}
+c^{7}_{77}
+2c^{8}_{77}
-c^{9}_{77}
-4  c^{11}_{77}
-2  c^{12}_{77}
-2  c^{13}_{77}
+2c^{14}_{77}
-c^{15}_{77}
-c^{16}_{77}
+c^{18}_{77}
+c^{19}_{77}
-5  c^{22}_{77}
-c^{23}_{77}
+2c^{25}_{77}
-c^{26}_{77}
+c^{29}_{77}
$;\ \ 
$ -4-2  c^{1}_{77}
-2  c^{2}_{77}
+2c^{3}_{77}
+c^{4}_{77}
-2  c^{5}_{77}
-2  c^{6}_{77}
+c^{7}_{77}
+2c^{8}_{77}
-c^{9}_{77}
-4  c^{11}_{77}
-2  c^{12}_{77}
-2  c^{13}_{77}
+2c^{14}_{77}
-c^{15}_{77}
-c^{16}_{77}
+c^{18}_{77}
+c^{19}_{77}
-5  c^{22}_{77}
-c^{23}_{77}
+2c^{25}_{77}
-c^{26}_{77}
+c^{29}_{77}
$,
$ 1+3c^{4}_{77}
+2c^{7}_{77}
-2  c^{9}_{77}
+c^{10}_{77}
+7c^{11}_{77}
+2c^{15}_{77}
-2  c^{16}_{77}
+c^{17}_{77}
+2c^{18}_{77}
-2  c^{19}_{77}
+c^{22}_{77}
-2  c^{23}_{77}
+c^{24}_{77}
+c^{25}_{77}
+2c^{26}_{77}
+2c^{29}_{77}
$;\ \ 
$ -1+c^{1}_{77}
+c^{2}_{77}
-c^{3}_{77}
-3  c^{4}_{77}
+c^{5}_{77}
+c^{6}_{77}
-c^{7}_{77}
-c^{8}_{77}
+c^{9}_{77}
-2  c^{10}_{77}
-2  c^{11}_{77}
+c^{12}_{77}
+c^{13}_{77}
-c^{14}_{77}
+c^{16}_{77}
-2  c^{17}_{77}
-c^{18}_{77}
+5c^{22}_{77}
+c^{23}_{77}
-2  c^{24}_{77}
-3  c^{25}_{77}
-c^{29}_{77}
$)

Realization: condensation reductions of $\eZ(\cNG(\Z_7,7))$.

\vskip 1ex

\noindent19. $10_{4,1435.}^{10,168}$ \irep{538}:\ \ 
$d_i$ = ($1.0$,
$9.472$,
$9.472$,
$9.472$,
$9.472$,
$9.472$,
$9.472$,
$16.944$,
$16.944$,
$17.944$) 

\vskip 0.7ex
\hangindent=3em \hangafter=1
$D^2= 1435.541 = 
720+320\sqrt{5}$

\vskip 0.7ex
\hangindent=3em \hangafter=1
$T = ( 0,
\frac{1}{2},
\frac{1}{2},
\frac{1}{2},
\frac{1}{2},
\frac{1}{2},
\frac{1}{2},
\frac{1}{5},
\frac{4}{5},
0 )
$,

\vskip 0.7ex
\hangindent=3em \hangafter=1
$S$ = ($ 1$,
$ 5+2\sqrt{5}$,
$ 5+2\sqrt{5}$,
$ 5+2\sqrt{5}$,
$ 5+2\sqrt{5}$,
$ 5+2\sqrt{5}$,
$ 5+2\sqrt{5}$,
$ 8+4\sqrt{5}$,
$ 8+4\sqrt{5}$,
$ 9+4\sqrt{5}$;\ \ 
$ 15+6\sqrt{5}$,
$ -5-2\sqrt{5}$,
$ -5-2\sqrt{5}$,
$ -5-2\sqrt{5}$,
$ -5-2\sqrt{5}$,
$ -5-2\sqrt{5}$,
$0$,
$0$,
$ 5+2\sqrt{5}$;\ \ 
$ 15+6\sqrt{5}$,
$ -5-2\sqrt{5}$,
$ -5-2\sqrt{5}$,
$ -5-2\sqrt{5}$,
$ -5-2\sqrt{5}$,
$0$,
$0$,
$ 5+2\sqrt{5}$;\ \ 
$ 15+6\sqrt{5}$,
$ -5-2\sqrt{5}$,
$ -5-2\sqrt{5}$,
$ -5-2\sqrt{5}$,
$0$,
$0$,
$ 5+2\sqrt{5}$;\ \ 
$ 15+6\sqrt{5}$,
$ -5-2\sqrt{5}$,
$ -5-2\sqrt{5}$,
$0$,
$0$,
$ 5+2\sqrt{5}$;\ \ 
$ 15+6\sqrt{5}$,
$ -5-2\sqrt{5}$,
$0$,
$0$,
$ 5+2\sqrt{5}$;\ \ 
$ 15+6\sqrt{5}$,
$0$,
$0$,
$ 5+2\sqrt{5}$;\ \ 
$ 14+6\sqrt{5}$,
$ -6-2\sqrt{5}$,
$ -8-4\sqrt{5}$;\ \ 
$ 14+6\sqrt{5}$,
$ -8-4\sqrt{5}$;\ \ 
$ 1$)

Realization: Condensation of $\Z_5$ bosons in $SU(5)_5$, see \cite{Edie-MichellCAMS}.

\vskip 1ex

\noindent20. $10_{0,1435.}^{20,676}$ \irep{948}:\ \ 
$d_i$ = ($1.0$,
$9.472$,
$9.472$,
$9.472$,
$9.472$,
$9.472$,
$9.472$,
$16.944$,
$16.944$,
$17.944$) 

\vskip 0.7ex
\hangindent=3em \hangafter=1
$D^2= 1435.541 = 
720+320\sqrt{5}$

\vskip 0.7ex
\hangindent=3em \hangafter=1
$T = ( 0,
0,
0,
\frac{1}{4},
\frac{1}{4},
\frac{3}{4},
\frac{3}{4},
\frac{2}{5},
\frac{3}{5},
0 )
$,

\vskip 0.7ex
\hangindent=3em \hangafter=1
$S$ = ($ 1$,
$ 5+2\sqrt{5}$,
$ 5+2\sqrt{5}$,
$ 5+2\sqrt{5}$,
$ 5+2\sqrt{5}$,
$ 5+2\sqrt{5}$,
$ 5+2\sqrt{5}$,
$ 8+4\sqrt{5}$,
$ 8+4\sqrt{5}$,
$ 9+4\sqrt{5}$;\ \ 
$ 15+6\sqrt{5}$,
$ -5-2\sqrt{5}$,
$ -5-2\sqrt{5}$,
$ -5-2\sqrt{5}$,
$ -5-2\sqrt{5}$,
$ -5-2\sqrt{5}$,
$0$,
$0$,
$ 5+2\sqrt{5}$;\ \ 
$ 15+6\sqrt{5}$,
$ -5-2\sqrt{5}$,
$ -5-2\sqrt{5}$,
$ -5-2\sqrt{5}$,
$ -5-2\sqrt{5}$,
$0$,
$0$,
$ 5+2\sqrt{5}$;\ \ 
$ -3-6  s^{1}_{20}
-4  c^{2}_{20}
+14s^{3}_{20}
$,
$ -3+6s^{1}_{20}
-4  c^{2}_{20}
-14  s^{3}_{20}
$,
$ 5+2\sqrt{5}$,
$ 5+2\sqrt{5}$,
$0$,
$0$,
$ 5+2\sqrt{5}$;\ \ 
$ -3-6  s^{1}_{20}
-4  c^{2}_{20}
+14s^{3}_{20}
$,
$ 5+2\sqrt{5}$,
$ 5+2\sqrt{5}$,
$0$,
$0$,
$ 5+2\sqrt{5}$;\ \ 
$ -3+6s^{1}_{20}
-4  c^{2}_{20}
-14  s^{3}_{20}
$,
$ -3-6  s^{1}_{20}
-4  c^{2}_{20}
+14s^{3}_{20}
$,
$0$,
$0$,
$ 5+2\sqrt{5}$;\ \ 
$ -3+6s^{1}_{20}
-4  c^{2}_{20}
-14  s^{3}_{20}
$,
$0$,
$0$,
$ 5+2\sqrt{5}$;\ \ 
$ -6-2\sqrt{5}$,
$ 14+6\sqrt{5}$,
$ -8-4\sqrt{5}$;\ \ 
$ -6-2\sqrt{5}$,
$ -8-4\sqrt{5}$;\ \ 
$ 1$)

Realization: Condensation reductions of $\eZ(\cNG(\Z_4\times \Z_4,16))$, see \cite{YuZhang}. 

\vskip 1ex

\noindent21. $10_{\frac{60}{11},43.10}^{11,372}$ \irep{561}:\ \ 
$d_i$ = ($1.0$,
$0.309$,
$1.682$,
$1.830$,
$2.397$,
$2.918$,
$-1.88$,
$-1.309$,
$-2.513$,
$-3.513$) 

\vskip 0.7ex
\hangindent=3em \hangafter=1
$D^2= 43.108 = 
33+11c^{1}_{11}
+11c^{2}_{11}
+11c^{3}_{11}
+11c^{4}_{11}
$

\vskip 0.7ex
\hangindent=3em \hangafter=1
$T = ( 0,
\frac{4}{11},
\frac{10}{11},
\frac{1}{11},
\frac{7}{11},
\frac{5}{11},
\frac{2}{11},
\frac{6}{11},
\frac{3}{11},
\frac{9}{11} )
$,

\vskip 0.7ex
\hangindent=3em \hangafter=1
$S$ = ($ 1$,
$ -1-c^{4}_{11}
$,
$ c_{11}^{1}$,
$ 1+c^{2}_{11}
$,
$ 1+c^{1}_{11}
+c^{3}_{11}
$,
$ 2+c^{1}_{11}
+c^{2}_{11}
+c^{3}_{11}
+c^{4}_{11}
$,
$ -1-c^{1}_{11}
-c^{3}_{11}
-c^{4}_{11}
$,
$ c_{11}^{4}$,
$ -c^{1}_{11}
-c^{2}_{11}
$,
$ -\xi_{11}^{5}$;\ \ 
$ -c_{11}^{4}$,
$ \xi_{11}^{5}$,
$ c_{11}^{1}$,
$ 1+c^{1}_{11}
+c^{3}_{11}
+c^{4}_{11}
$,
$ 1+c^{2}_{11}
$,
$ 2+c^{1}_{11}
+c^{2}_{11}
+c^{3}_{11}
+c^{4}_{11}
$,
$ 1+c^{1}_{11}
+c^{3}_{11}
$,
$ 1$,
$ c^{1}_{11}
+c^{2}_{11}
$;\ \ 
$ -1-c^{1}_{11}
-c^{3}_{11}
-c^{4}_{11}
$,
$ -1-c^{1}_{11}
-c^{3}_{11}
$,
$ -1$,
$ -c_{11}^{4}$,
$ -1-c^{4}_{11}
$,
$ c^{1}_{11}
+c^{2}_{11}
$,
$ 1+c^{2}_{11}
$,
$ -2-c^{1}_{11}
-c^{2}_{11}
-c^{3}_{11}
-c^{4}_{11}
$;\ \ 
$ -c_{11}^{4}$,
$ c^{1}_{11}
+c^{2}_{11}
$,
$ -1-c^{4}_{11}
$,
$ 1$,
$ -\xi_{11}^{5}$,
$ 2+c^{1}_{11}
+c^{2}_{11}
+c^{3}_{11}
+c^{4}_{11}
$,
$ 1+c^{1}_{11}
+c^{3}_{11}
+c^{4}_{11}
$;\ \ 
$ 1+c^{2}_{11}
$,
$ -\xi_{11}^{5}$,
$ -c_{11}^{1}$,
$ 2+c^{1}_{11}
+c^{2}_{11}
+c^{3}_{11}
+c^{4}_{11}
$,
$ c_{11}^{4}$,
$ -1-c^{4}_{11}
$;\ \ 
$ 1$,
$ -c^{1}_{11}
-c^{2}_{11}
$,
$ -c_{11}^{1}$,
$ -1-c^{1}_{11}
-c^{3}_{11}
-c^{4}_{11}
$,
$ 1+c^{1}_{11}
+c^{3}_{11}
$;\ \ 
$ \xi_{11}^{5}$,
$ -1-c^{2}_{11}
$,
$ -1-c^{1}_{11}
-c^{3}_{11}
$,
$ c_{11}^{4}$;\ \ 
$ -1-c^{1}_{11}
-c^{3}_{11}
-c^{4}_{11}
$,
$ 1+c^{4}_{11}
$,
$ 1$;\ \ 
$ \xi_{11}^{5}$,
$ -c_{11}^{1}$;\ \ 
$ 1+c^{2}_{11}
$)

Realization: 
$PSO(5)_{\frac{5}{2}}$, {\it i.e.} the adjoint subcategory of the
non-unitary braided fusion category $SO(5)_{\frac{5}{2}}$ corresponding to
$U_q\mathfrak{so}_5$ with $q=e^{\pi i/11}$, see \cite{RowellMathZ}.

\vskip 1ex 

}

\subsection{Rank 11 }\label{ss:rank11}

{\small

\noindent1. $11_{2,11.}^{11,568}$ \irep{983}:\ \ 
$d_i$ = ($1.0$,
$1.0$,
$1.0$,
$1.0$,
$1.0$,
$1.0$,
$1.0$,
$1.0$,
$1.0$,
$1.0$,
$1.0$) 

\vskip 0.7ex
\hangindent=3em \hangafter=1
$D^2= 11.0 = 
11$

\vskip 0.7ex
\hangindent=3em \hangafter=1
$T = ( 0,
\frac{1}{11},
\frac{1}{11},
\frac{3}{11},
\frac{3}{11},
\frac{4}{11},
\frac{4}{11},
\frac{5}{11},
\frac{5}{11},
\frac{9}{11},
\frac{9}{11} )
$,

\vskip 0.7ex
\hangindent=3em \hangafter=1
$S$ = ($ 1$,
$ 1$,
$ 1$,
$ 1$,
$ 1$,
$ 1$,
$ 1$,
$ 1$,
$ 1$,
$ 1$,
$ 1$;\ \ 
$ -\zeta_{22}^{7}$,
$ \zeta_{11}^{2}$,
$ -\zeta_{22}^{9}$,
$ \zeta_{11}^{1}$,
$ -\zeta_{22}^{3}$,
$ \zeta_{11}^{4}$,
$ -\zeta_{22}^{5}$,
$ \zeta_{11}^{3}$,
$ -\zeta_{22}^{1}$,
$ \zeta_{11}^{5}$;\ \ 
$ -\zeta_{22}^{7}$,
$ \zeta_{11}^{1}$,
$ -\zeta_{22}^{9}$,
$ \zeta_{11}^{4}$,
$ -\zeta_{22}^{3}$,
$ \zeta_{11}^{3}$,
$ -\zeta_{22}^{5}$,
$ \zeta_{11}^{5}$,
$ -\zeta_{22}^{1}$;\ \ 
$ \zeta_{11}^{5}$,
$ -\zeta_{22}^{1}$,
$ -\zeta_{22}^{7}$,
$ \zeta_{11}^{2}$,
$ \zeta_{11}^{4}$,
$ -\zeta_{22}^{3}$,
$ \zeta_{11}^{3}$,
$ -\zeta_{22}^{5}$;\ \ 
$ \zeta_{11}^{5}$,
$ \zeta_{11}^{2}$,
$ -\zeta_{22}^{7}$,
$ -\zeta_{22}^{3}$,
$ \zeta_{11}^{4}$,
$ -\zeta_{22}^{5}$,
$ \zeta_{11}^{3}$;\ \ 
$ \zeta_{11}^{3}$,
$ -\zeta_{22}^{5}$,
$ \zeta_{11}^{5}$,
$ -\zeta_{22}^{1}$,
$ \zeta_{11}^{1}$,
$ -\zeta_{22}^{9}$;\ \ 
$ \zeta_{11}^{3}$,
$ -\zeta_{22}^{1}$,
$ \zeta_{11}^{5}$,
$ -\zeta_{22}^{9}$,
$ \zeta_{11}^{1}$;\ \ 
$ \zeta_{11}^{1}$,
$ -\zeta_{22}^{9}$,
$ -\zeta_{22}^{7}$,
$ \zeta_{11}^{2}$;\ \ 
$ \zeta_{11}^{1}$,
$ \zeta_{11}^{2}$,
$ -\zeta_{22}^{7}$;\ \ 
$ \zeta_{11}^{4}$,
$ -\zeta_{22}^{3}$;\ \ 
$ \zeta_{11}^{4}$)

Realization:
$U(11)_1$.

\vskip 1ex

\noindent2. $11_{1,32.}^{16,245}$ \irep{0}:\ \ 
$d_i$ = ($1.0$,
$1.0$,
$1.0$,
$1.0$,
$2.0$,
$2.0$,
$2.0$,
$2.0$,
$2.0$,
$2.0$,
$2.0$) 

\vskip 0.7ex
\hangindent=3em \hangafter=1
$D^2= 32.0 = 
32$

\vskip 0.7ex
\hangindent=3em \hangafter=1
$T = ( 0,
0,
0,
0,
\frac{1}{4},
\frac{1}{16},
\frac{1}{16},
\frac{1}{16},
\frac{9}{16},
\frac{9}{16},
\frac{9}{16} )
$,

\vskip 0.7ex
\hangindent=3em \hangafter=1
$S$ = ($ 1$,
$ 1$,
$ 1$,
$ 1$,
$ 2$,
$ 2$,
$ 2$,
$ 2$,
$ 2$,
$ 2$,
$ 2$;\ \ 
$ 1$,
$ 1$,
$ 1$,
$ 2$,
$ -2$,
$ -2$,
$ 2$,
$ -2$,
$ -2$,
$ 2$;\ \ 
$ 1$,
$ 1$,
$ 2$,
$ -2$,
$ 2$,
$ -2$,
$ -2$,
$ 2$,
$ -2$;\ \ 
$ 1$,
$ 2$,
$ 2$,
$ -2$,
$ -2$,
$ 2$,
$ -2$,
$ -2$;\ \ 
$ -4$,
$0$,
$0$,
$0$,
$0$,
$0$,
$0$;\ \ 
$ 2\sqrt{2}$,
$0$,
$0$,
$ -2\sqrt{2}$,
$0$,
$0$;\ \ 
$ 2\sqrt{2}$,
$0$,
$0$,
$ -2\sqrt{2}$,
$0$;\ \ 
$ 2\sqrt{2}$,
$0$,
$0$,
$ -2\sqrt{2}$;\ \ 
$ 2\sqrt{2}$,
$0$,
$0$;\ \ 
$ 2\sqrt{2}$,
$0$;\ \ 
$ 2\sqrt{2}$)

Realization: $O_4$ or $\overline{SO(16)}_3$ or 
Abelian anyon condensations of $\eD^3(Q_8)$.

\vskip 1ex

\noindent3. $11_{1,32.}^{16,157}$ \irep{0}:\ \ 
$d_i$ = ($1.0$,
$1.0$,
$1.0$,
$1.0$,
$2.0$,
$2.0$,
$2.0$,
$2.0$,
$2.0$,
$2.0$,
$2.0$) 

\vskip 0.7ex
\hangindent=3em \hangafter=1
$D^2= 32.0 = 
32$

\vskip 0.7ex
\hangindent=3em \hangafter=1
$T = ( 0,
0,
0,
0,
\frac{1}{4},
\frac{1}{16},
\frac{1}{16},
\frac{5}{16},
\frac{9}{16},
\frac{9}{16},
\frac{13}{16} )
$,

\vskip 0.7ex
\hangindent=3em \hangafter=1
$S$ = ($ 1$,
$ 1$,
$ 1$,
$ 1$,
$ 2$,
$ 2$,
$ 2$,
$ 2$,
$ 2$,
$ 2$,
$ 2$;\ \ 
$ 1$,
$ 1$,
$ 1$,
$ 2$,
$ -2$,
$ -2$,
$ 2$,
$ -2$,
$ -2$,
$ 2$;\ \ 
$ 1$,
$ 1$,
$ 2$,
$ -2$,
$ 2$,
$ -2$,
$ -2$,
$ 2$,
$ -2$;\ \ 
$ 1$,
$ 2$,
$ 2$,
$ -2$,
$ -2$,
$ 2$,
$ -2$,
$ -2$;\ \ 
$ -4$,
$0$,
$0$,
$0$,
$0$,
$0$,
$0$;\ \ 
$ 2\sqrt{2}$,
$0$,
$0$,
$ -2\sqrt{2}$,
$0$,
$0$;\ \ 
$ 2\sqrt{2}$,
$0$,
$0$,
$ -2\sqrt{2}$,
$0$;\ \ 
$ -2\sqrt{2}$,
$0$,
$0$,
$ 2\sqrt{2}$;\ \ 
$ 2\sqrt{2}$,
$0$,
$0$;\ \ 
$ 2\sqrt{2}$,
$0$;\ \ 
$ -2\sqrt{2}$)

Realization: 
Abelian anyon condensations of $O_4$ or $\eD^1(Q_8)$ or $\overline{\eD^3(Q_8)}$.

\vskip 1ex

\noindent4. $11_{2,60.}^{120,157}$ \irep{2274}:\ \ 
$d_i$ = ($1.0$,
$1.0$,
$2.0$,
$2.0$,
$2.0$,
$2.0$,
$2.0$,
$2.0$,
$2.0$,
$3.872$,
$3.872$) 

\vskip 0.7ex
\hangindent=3em \hangafter=1
$D^2= 60.0 = 
60$

\vskip 0.7ex
\hangindent=3em \hangafter=1
$T = ( 0,
0,
\frac{1}{3},
\frac{1}{5},
\frac{4}{5},
\frac{2}{15},
\frac{2}{15},
\frac{8}{15},
\frac{8}{15},
\frac{1}{8},
\frac{5}{8} )
$,

\vskip 0.7ex
\hangindent=3em \hangafter=1
$S$ = ($ 1$,
$ 1$,
$ 2$,
$ 2$,
$ 2$,
$ 2$,
$ 2$,
$ 2$,
$ 2$,
$ \sqrt{15}$,
$ \sqrt{15}$;\ \ 
$ 1$,
$ 2$,
$ 2$,
$ 2$,
$ 2$,
$ 2$,
$ 2$,
$ 2$,
$ -\sqrt{15}$,
$ -\sqrt{15}$;\ \ 
$ -2$,
$ 4$,
$ 4$,
$ -2$,
$ -2$,
$ -2$,
$ -2$,
$0$,
$0$;\ \ 
$ -1-\sqrt{5}$,
$ -1+\sqrt{5}$,
$ -1+\sqrt{5}$,
$ -1+\sqrt{5}$,
$ -1-\sqrt{5}$,
$ -1-\sqrt{5}$,
$0$,
$0$;\ \ 
$ -1-\sqrt{5}$,
$ -1-\sqrt{5}$,
$ -1-\sqrt{5}$,
$ -1+\sqrt{5}$,
$ -1+\sqrt{5}$,
$0$,
$0$;\ \ 
$ 2c_{15}^{4}$,
$ 2c_{15}^{1}$,
$ 2c_{15}^{7}$,
$ 2c_{15}^{2}$,
$0$,
$0$;\ \ 
$ 2c_{15}^{4}$,
$ 2c_{15}^{2}$,
$ 2c_{15}^{7}$,
$0$,
$0$;\ \ 
$ 2c_{15}^{1}$,
$ 2c_{15}^{4}$,
$0$,
$0$;\ \ 
$ 2c_{15}^{1}$,
$0$,
$0$;\ \ 
$ \sqrt{15}$,
$ -\sqrt{15}$;\ \ 
$ \sqrt{15}$)

Connected to the orbit of $11_{2,60.}^{120,364}$ via a change of spherical
structure.

Realization:
Galois conjugation of 
$SO(15)_2$.

\vskip 1ex 

{\grey 
\noindent5. $11_{2,60.}^{120,364}$ \irep{2274}:\ \ 
$d_i$ = ($1.0$,
$1.0$,
$2.0$,
$2.0$,
$2.0$,
$2.0$,
$2.0$,
$2.0$,
$2.0$,
$3.872$,
$3.872$) 

\vskip 0.7ex
\hangindent=3em \hangafter=1
$D^2= 60.0 = 
60$

\vskip 0.7ex
\hangindent=3em \hangafter=1
$T = ( 0,
0,
\frac{1}{3},
\frac{1}{5},
\frac{4}{5},
\frac{2}{15},
\frac{2}{15},
\frac{8}{15},
\frac{8}{15},
\frac{3}{8},
\frac{7}{8} )
$,

\vskip 0.7ex
\hangindent=3em \hangafter=1
$S$ = ($ 1$,
$ 1$,
$ 2$,
$ 2$,
$ 2$,
$ 2$,
$ 2$,
$ 2$,
$ 2$,
$ \sqrt{15}$,
$ \sqrt{15}$;\ \ 
$ 1$,
$ 2$,
$ 2$,
$ 2$,
$ 2$,
$ 2$,
$ 2$,
$ 2$,
$ -\sqrt{15}$,
$ -\sqrt{15}$;\ \ 
$ -2$,
$ 4$,
$ 4$,
$ -2$,
$ -2$,
$ -2$,
$ -2$,
$0$,
$0$;\ \ 
$ -1-\sqrt{5}$,
$ -1+\sqrt{5}$,
$ -1+\sqrt{5}$,
$ -1+\sqrt{5}$,
$ -1-\sqrt{5}$,
$ -1-\sqrt{5}$,
$0$,
$0$;\ \ 
$ -1-\sqrt{5}$,
$ -1-\sqrt{5}$,
$ -1-\sqrt{5}$,
$ -1+\sqrt{5}$,
$ -1+\sqrt{5}$,
$0$,
$0$;\ \ 
$ 2c_{15}^{4}$,
$ 2c_{15}^{1}$,
$ 2c_{15}^{7}$,
$ 2c_{15}^{2}$,
$0$,
$0$;\ \ 
$ 2c_{15}^{4}$,
$ 2c_{15}^{2}$,
$ 2c_{15}^{7}$,
$0$,
$0$;\ \ 
$ 2c_{15}^{1}$,
$ 2c_{15}^{4}$,
$0$,
$0$;\ \ 
$ 2c_{15}^{1}$,
$0$,
$0$;\ \ 
$ -\sqrt{15}$,
$ \sqrt{15}$;\ \ 
$ -\sqrt{15}$)

Connected to the orbit of $11_{2,60.}^{120,157}$ via a change of spherical
structure.

Realization:
Galois conjugation of 
Abelian anyon condensation of $SO(30)_2$.
}

\vskip 1ex

\noindent6. $11_{\frac{13}{2},89.56}^{48,108}$ \irep{2191}:\ \ 
$d_i$ = ($1.0$,
$1.0$,
$1.931$,
$1.931$,
$2.732$,
$2.732$,
$3.346$,
$3.346$,
$3.732$,
$3.732$,
$3.863$) 

\vskip 0.7ex
\hangindent=3em \hangafter=1
$D^2= 89.569 = 
48+24\sqrt{3}$

\vskip 0.7ex
\hangindent=3em \hangafter=1
$T = ( 0,
\frac{1}{2},
\frac{1}{16},
\frac{1}{16},
\frac{1}{3},
\frac{5}{6},
\frac{13}{16},
\frac{13}{16},
0,
\frac{1}{2},
\frac{19}{48} )
$,

\vskip 0.7ex
\hangindent=3em \hangafter=1
$S$ = ($ 1$,
$ 1$,
$ c_{24}^{1}$,
$ c_{24}^{1}$,
$ 1+\sqrt{3}$,
$ 1+\sqrt{3}$,
$ \frac{3+3\sqrt{3}}{\sqrt{6}}$,
$ \frac{3+3\sqrt{3}}{\sqrt{6}}$,
$ 2+\sqrt{3}$,
$ 2+\sqrt{3}$,
$ 2c_{24}^{1}$;\ \ 
$ 1$,
$ -c_{24}^{1}$,
$ -c_{24}^{1}$,
$ 1+\sqrt{3}$,
$ 1+\sqrt{3}$,
$ \frac{-3-3\sqrt{3}}{\sqrt{6}}$,
$ \frac{-3-3\sqrt{3}}{\sqrt{6}}$,
$ 2+\sqrt{3}$,
$ 2+\sqrt{3}$,
$ -2c_{24}^{1}$;\ \ 
$(\frac{-3-3\sqrt{3}}{\sqrt{6}})\mathrm{i}$,
$(\frac{3+3\sqrt{3}}{\sqrt{6}})\mathrm{i}$,
$ -2c_{24}^{1}$,
$ 2c_{24}^{1}$,
$(\frac{3+3\sqrt{3}}{\sqrt{6}})\mathrm{i}$,
$(\frac{-3-3\sqrt{3}}{\sqrt{6}})\mathrm{i}$,
$ -c_{24}^{1}$,
$ c_{24}^{1}$,
$0$;\ \ 
$(\frac{-3-3\sqrt{3}}{\sqrt{6}})\mathrm{i}$,
$ -2c_{24}^{1}$,
$ 2c_{24}^{1}$,
$(\frac{-3-3\sqrt{3}}{\sqrt{6}})\mathrm{i}$,
$(\frac{3+3\sqrt{3}}{\sqrt{6}})\mathrm{i}$,
$ -c_{24}^{1}$,
$ c_{24}^{1}$,
$0$;\ \ 
$ 1+\sqrt{3}$,
$ 1+\sqrt{3}$,
$0$,
$0$,
$ -1-\sqrt{3}$,
$ -1-\sqrt{3}$,
$ 2c_{24}^{1}$;\ \ 
$ 1+\sqrt{3}$,
$0$,
$0$,
$ -1-\sqrt{3}$,
$ -1-\sqrt{3}$,
$ -2c_{24}^{1}$;\ \ 
$(\frac{3+3\sqrt{3}}{\sqrt{6}})\mathrm{i}$,
$(\frac{-3-3\sqrt{3}}{\sqrt{6}})\mathrm{i}$,
$ \frac{3+3\sqrt{3}}{\sqrt{6}}$,
$ \frac{-3-3\sqrt{3}}{\sqrt{6}}$,
$0$;\ \ 
$(\frac{3+3\sqrt{3}}{\sqrt{6}})\mathrm{i}$,
$ \frac{3+3\sqrt{3}}{\sqrt{6}}$,
$ \frac{-3-3\sqrt{3}}{\sqrt{6}}$,
$0$;\ \ 
$ 1$,
$ 1$,
$ -2c_{24}^{1}$;\ \ 
$ 1$,
$ 2c_{24}^{1}$;\ \ 
$0$) 

Connected to the orbit of $11_{\frac{7}{2},89.56}^{48,628}$ via a change of spherical
structure.

Realization:
Abelian anyon condensation of
$\overline{SU}(2)_{10}$ or
$SU(10)_2$ or
$Sp(20)_1$.

\vskip 1ex 

{\grey
\noindent7. $11_{\frac{7}{2},89.56}^{48,628}$ \irep{2191}:\ \ 
$d_i$ = ($1.0$,
$1.0$,
$1.931$,
$1.931$,
$2.732$,
$2.732$,
$3.346$,
$3.346$,
$3.732$,
$3.732$,
$3.863$) 

\vskip 0.7ex
\hangindent=3em \hangafter=1
$D^2= 89.569 = 
48+24\sqrt{3}$

\vskip 0.7ex
\hangindent=3em \hangafter=1
$T = ( 0,
\frac{1}{2},
\frac{3}{16},
\frac{3}{16},
\frac{2}{3},
\frac{1}{6},
\frac{7}{16},
\frac{7}{16},
0,
\frac{1}{2},
\frac{41}{48} )
$,

\vskip 0.7ex
\hangindent=3em \hangafter=1
$S$ = ($ 1$,
$ 1$,
$ c_{24}^{1}$,
$ c_{24}^{1}$,
$ 1+\sqrt{3}$,
$ 1+\sqrt{3}$,
$ \frac{3+3\sqrt{3}}{\sqrt{6}}$,
$ \frac{3+3\sqrt{3}}{\sqrt{6}}$,
$ 2+\sqrt{3}$,
$ 2+\sqrt{3}$,
$ 2c_{24}^{1}$;\ \ 
$ 1$,
$ -c_{24}^{1}$,
$ -c_{24}^{1}$,
$ 1+\sqrt{3}$,
$ 1+\sqrt{3}$,
$ \frac{-3-3\sqrt{3}}{\sqrt{6}}$,
$ \frac{-3-3\sqrt{3}}{\sqrt{6}}$,
$ 2+\sqrt{3}$,
$ 2+\sqrt{3}$,
$ -2c_{24}^{1}$;\ \ 
$(\frac{-3-3\sqrt{3}}{\sqrt{6}})\mathrm{i}$,
$(\frac{3+3\sqrt{3}}{\sqrt{6}})\mathrm{i}$,
$ -2c_{24}^{1}$,
$ 2c_{24}^{1}$,
$(\frac{3+3\sqrt{3}}{\sqrt{6}})\mathrm{i}$,
$(\frac{-3-3\sqrt{3}}{\sqrt{6}})\mathrm{i}$,
$ -c_{24}^{1}$,
$ c_{24}^{1}$,
$0$;\ \ 
$(\frac{-3-3\sqrt{3}}{\sqrt{6}})\mathrm{i}$,
$ -2c_{24}^{1}$,
$ 2c_{24}^{1}$,
$(\frac{-3-3\sqrt{3}}{\sqrt{6}})\mathrm{i}$,
$(\frac{3+3\sqrt{3}}{\sqrt{6}})\mathrm{i}$,
$ -c_{24}^{1}$,
$ c_{24}^{1}$,
$0$;\ \ 
$ 1+\sqrt{3}$,
$ 1+\sqrt{3}$,
$0$,
$0$,
$ -1-\sqrt{3}$,
$ -1-\sqrt{3}$,
$ 2c_{24}^{1}$;\ \ 
$ 1+\sqrt{3}$,
$0$,
$0$,
$ -1-\sqrt{3}$,
$ -1-\sqrt{3}$,
$ -2c_{24}^{1}$;\ \ 
$(\frac{3+3\sqrt{3}}{\sqrt{6}})\mathrm{i}$,
$(\frac{-3-3\sqrt{3}}{\sqrt{6}})\mathrm{i}$,
$ \frac{3+3\sqrt{3}}{\sqrt{6}}$,
$ \frac{-3-3\sqrt{3}}{\sqrt{6}}$,
$0$;\ \ 
$(\frac{3+3\sqrt{3}}{\sqrt{6}})\mathrm{i}$,
$ \frac{3+3\sqrt{3}}{\sqrt{6}}$,
$ \frac{-3-3\sqrt{3}}{\sqrt{6}}$,
$0$;\ \ 
$ 1$,
$ 1$,
$ -2c_{24}^{1}$;\ \ 
$ 1$,
$ 2c_{24}^{1}$;\ \ 
$0$)

Connected to the orbit of $11_{\frac{13}{2},89.56}^{48,108}$ via a change of spherical
structure.

Realization:
Abelian anyon condensation of
$SU(2)_{10}$ or
$\overline{SU}(10)_2$ or
$\overline{Sp}(20)_1$.
}

\vskip 1ex

\noindent8. $11_{\frac{5}{2},89.56}^{48,214}$ \irep{2192}:\ \ 
$d_i$ = ($1.0$,
$1.0$,
$1.931$,
$1.931$,
$2.732$,
$2.732$,
$3.346$,
$3.346$,
$3.732$,
$3.732$,
$3.863$) 

\vskip 0.7ex
\hangindent=3em \hangafter=1
$D^2= 89.569 = 
48+24\sqrt{3}$

\vskip 0.7ex
\hangindent=3em \hangafter=1
$T = ( 0,
\frac{1}{2},
\frac{1}{16},
\frac{1}{16},
\frac{2}{3},
\frac{1}{6},
\frac{5}{16},
\frac{5}{16},
0,
\frac{1}{2},
\frac{35}{48} )
$,

\vskip 0.7ex
\hangindent=3em \hangafter=1
$S$ = ($ 1$,
$ 1$,
$ c_{24}^{1}$,
$ c_{24}^{1}$,
$ 1+\sqrt{3}$,
$ 1+\sqrt{3}$,
$ \frac{3+3\sqrt{3}}{\sqrt{6}}$,
$ \frac{3+3\sqrt{3}}{\sqrt{6}}$,
$ 2+\sqrt{3}$,
$ 2+\sqrt{3}$,
$ 2c_{24}^{1}$;\ \ 
$ 1$,
$ -c_{24}^{1}$,
$ -c_{24}^{1}$,
$ 1+\sqrt{3}$,
$ 1+\sqrt{3}$,
$ \frac{-3-3\sqrt{3}}{\sqrt{6}}$,
$ \frac{-3-3\sqrt{3}}{\sqrt{6}}$,
$ 2+\sqrt{3}$,
$ 2+\sqrt{3}$,
$ -2c_{24}^{1}$;\ \ 
$ \frac{3+3\sqrt{3}}{\sqrt{6}}$,
$ \frac{-3-3\sqrt{3}}{\sqrt{6}}$,
$ -2c_{24}^{1}$,
$ 2c_{24}^{1}$,
$ \frac{-3-3\sqrt{3}}{\sqrt{6}}$,
$ \frac{3+3\sqrt{3}}{\sqrt{6}}$,
$ -c_{24}^{1}$,
$ c_{24}^{1}$,
$0$;\ \ 
$ \frac{3+3\sqrt{3}}{\sqrt{6}}$,
$ -2c_{24}^{1}$,
$ 2c_{24}^{1}$,
$ \frac{3+3\sqrt{3}}{\sqrt{6}}$,
$ \frac{-3-3\sqrt{3}}{\sqrt{6}}$,
$ -c_{24}^{1}$,
$ c_{24}^{1}$,
$0$;\ \ 
$ 1+\sqrt{3}$,
$ 1+\sqrt{3}$,
$0$,
$0$,
$ -1-\sqrt{3}$,
$ -1-\sqrt{3}$,
$ 2c_{24}^{1}$;\ \ 
$ 1+\sqrt{3}$,
$0$,
$0$,
$ -1-\sqrt{3}$,
$ -1-\sqrt{3}$,
$ -2c_{24}^{1}$;\ \ 
$ \frac{-3-3\sqrt{3}}{\sqrt{6}}$,
$ \frac{3+3\sqrt{3}}{\sqrt{6}}$,
$ \frac{3+3\sqrt{3}}{\sqrt{6}}$,
$ \frac{-3-3\sqrt{3}}{\sqrt{6}}$,
$0$;\ \ 
$ \frac{-3-3\sqrt{3}}{\sqrt{6}}$,
$ \frac{3+3\sqrt{3}}{\sqrt{6}}$,
$ \frac{-3-3\sqrt{3}}{\sqrt{6}}$,
$0$;\ \ 
$ 1$,
$ 1$,
$ -2c_{24}^{1}$;\ \ 
$ 1$,
$ 2c_{24}^{1}$;\ \ 
$0$)

Connected to the orbit of $11_{\frac{15}{2},89.56}^{48,311}$ via a change of spherical
structure.

Realization:
$SU(2)_{10}$ or
Abelian anyon condensation of
$\overline{SU}(10)_2$ or
$\overline{Sp}(20)_1$.

\vskip 1ex 

{\grey
\noindent9. $11_{\frac{15}{2},89.56}^{48,311}$ \irep{2192}:\ \ 
$d_i$ = ($1.0$,
$1.0$,
$1.931$,
$1.931$,
$2.732$,
$2.732$,
$3.346$,
$3.346$,
$3.732$,
$3.732$,
$3.863$) 

\vskip 0.7ex
\hangindent=3em \hangafter=1
$D^2= 89.569 = 
48+24\sqrt{3}$

\vskip 0.7ex
\hangindent=3em \hangafter=1
$T = ( 0,
\frac{1}{2},
\frac{3}{16},
\frac{3}{16},
\frac{1}{3},
\frac{5}{6},
\frac{15}{16},
\frac{15}{16},
0,
\frac{1}{2},
\frac{25}{48} )
$,

\vskip 0.7ex
\hangindent=3em \hangafter=1
$S$ = ($ 1$,
$ 1$,
$ c_{24}^{1}$,
$ c_{24}^{1}$,
$ 1+\sqrt{3}$,
$ 1+\sqrt{3}$,
$ \frac{3+3\sqrt{3}}{\sqrt{6}}$,
$ \frac{3+3\sqrt{3}}{\sqrt{6}}$,
$ 2+\sqrt{3}$,
$ 2+\sqrt{3}$,
$ 2c_{24}^{1}$;\ \ 
$ 1$,
$ -c_{24}^{1}$,
$ -c_{24}^{1}$,
$ 1+\sqrt{3}$,
$ 1+\sqrt{3}$,
$ \frac{-3-3\sqrt{3}}{\sqrt{6}}$,
$ \frac{-3-3\sqrt{3}}{\sqrt{6}}$,
$ 2+\sqrt{3}$,
$ 2+\sqrt{3}$,
$ -2c_{24}^{1}$;\ \ 
$ \frac{-3-3\sqrt{3}}{\sqrt{6}}$,
$ \frac{3+3\sqrt{3}}{\sqrt{6}}$,
$ -2c_{24}^{1}$,
$ 2c_{24}^{1}$,
$ \frac{-3-3\sqrt{3}}{\sqrt{6}}$,
$ \frac{3+3\sqrt{3}}{\sqrt{6}}$,
$ -c_{24}^{1}$,
$ c_{24}^{1}$,
$0$;\ \ 
$ \frac{-3-3\sqrt{3}}{\sqrt{6}}$,
$ -2c_{24}^{1}$,
$ 2c_{24}^{1}$,
$ \frac{3+3\sqrt{3}}{\sqrt{6}}$,
$ \frac{-3-3\sqrt{3}}{\sqrt{6}}$,
$ -c_{24}^{1}$,
$ c_{24}^{1}$,
$0$;\ \ 
$ 1+\sqrt{3}$,
$ 1+\sqrt{3}$,
$0$,
$0$,
$ -1-\sqrt{3}$,
$ -1-\sqrt{3}$,
$ 2c_{24}^{1}$;\ \ 
$ 1+\sqrt{3}$,
$0$,
$0$,
$ -1-\sqrt{3}$,
$ -1-\sqrt{3}$,
$ -2c_{24}^{1}$;\ \ 
$ \frac{3+3\sqrt{3}}{\sqrt{6}}$,
$ \frac{-3-3\sqrt{3}}{\sqrt{6}}$,
$ \frac{3+3\sqrt{3}}{\sqrt{6}}$,
$ \frac{-3-3\sqrt{3}}{\sqrt{6}}$,
$0$;\ \ 
$ \frac{3+3\sqrt{3}}{\sqrt{6}}$,
$ \frac{3+3\sqrt{3}}{\sqrt{6}}$,
$ \frac{-3-3\sqrt{3}}{\sqrt{6}}$,
$0$;\ \ 
$ 1$,
$ 1$,
$ -2c_{24}^{1}$;\ \ 
$ 1$,
$ 2c_{24}^{1}$;\ \ 
$0$)

Connected to the orbit of $11_{\frac{5}{2},89.56}^{48,214}$ via a change of spherical
structure.

Realization:
Abelian anyon condensation of
$\overline{SU}(2)_{10}$ or
$SU(10)_2$ or
$Sp(20)_1$.
}

\vskip 1ex

\noindent10. $11_{\frac{144}{23},310.1}^{23,306}$ \irep{1836}:\ \ 
$d_i$ = ($1.0$,
$1.981$,
$2.925$,
$3.815$,
$4.634$,
$5.367$,
$5.999$,
$6.520$,
$6.919$,
$7.190$,
$7.326$) 

\vskip 0.7ex
\hangindent=3em \hangafter=1
$D^2= 310.117 = 
66+55c^{1}_{23}
+45c^{2}_{23}
+36c^{3}_{23}
+28c^{4}_{23}
+21c^{5}_{23}
+15c^{6}_{23}
+10c^{7}_{23}
+6c^{8}_{23}
+3c^{9}_{23}
+c^{10}_{23}
$

\vskip 0.7ex
\hangindent=3em \hangafter=1
$T = ( 0,
\frac{5}{23},
\frac{21}{23},
\frac{2}{23},
\frac{17}{23},
\frac{20}{23},
\frac{11}{23},
\frac{13}{23},
\frac{3}{23},
\frac{4}{23},
\frac{16}{23} )
$,

\vskip 0.7ex
\hangindent=3em \hangafter=1
$S$ = ($ 1$,
$ -c_{23}^{11}$,
$ \xi_{23}^{3}$,
$ \xi_{23}^{19}$,
$ \xi_{23}^{5}$,
$ \xi_{23}^{17}$,
$ \xi_{23}^{7}$,
$ \xi_{23}^{15}$,
$ \xi_{23}^{9}$,
$ \xi_{23}^{13}$,
$ \xi_{23}^{11}$;\ \ 
$ -\xi_{23}^{19}$,
$ \xi_{23}^{17}$,
$ -\xi_{23}^{15}$,
$ \xi_{23}^{13}$,
$ -\xi_{23}^{11}$,
$ \xi_{23}^{9}$,
$ -\xi_{23}^{7}$,
$ \xi_{23}^{5}$,
$ -\xi_{23}^{3}$,
$ 1$;\ \ 
$ \xi_{23}^{9}$,
$ \xi_{23}^{11}$,
$ \xi_{23}^{15}$,
$ \xi_{23}^{5}$,
$ -c_{23}^{11}$,
$ -1$,
$ -\xi_{23}^{19}$,
$ -\xi_{23}^{7}$,
$ -\xi_{23}^{13}$;\ \ 
$ -\xi_{23}^{7}$,
$ \xi_{23}^{3}$,
$ 1$,
$ -\xi_{23}^{5}$,
$ \xi_{23}^{9}$,
$ -\xi_{23}^{13}$,
$ \xi_{23}^{17}$,
$ c_{23}^{11}$;\ \ 
$ c_{23}^{11}$,
$ -\xi_{23}^{7}$,
$ -\xi_{23}^{11}$,
$ -\xi_{23}^{17}$,
$ -1$,
$ \xi_{23}^{19}$,
$ \xi_{23}^{9}$;\ \ 
$ \xi_{23}^{13}$,
$ -\xi_{23}^{19}$,
$ c_{23}^{11}$,
$ \xi_{23}^{15}$,
$ -\xi_{23}^{9}$,
$ \xi_{23}^{3}$;\ \ 
$ \xi_{23}^{3}$,
$ \xi_{23}^{13}$,
$ \xi_{23}^{17}$,
$ -1$,
$ -\xi_{23}^{15}$;\ \ 
$ -\xi_{23}^{5}$,
$ -\xi_{23}^{3}$,
$ \xi_{23}^{11}$,
$ -\xi_{23}^{19}$;\ \ 
$ -\xi_{23}^{11}$,
$ c_{23}^{11}$,
$ \xi_{23}^{7}$;\ \ 
$ -\xi_{23}^{15}$,
$ \xi_{23}^{5}$;\ \ 
$ -\xi_{23}^{17}$)

Realization:
Abelian anyon condensation of
$\overline{SU}(2)_{21}$ or
$SU(21)_2$ or
$Sp(42)_1$.

\vskip 1ex

\noindent11. $11_{\frac{54}{19},696.5}^{19,306}$ \irep{1746}:\ \ 
$d_i$ = ($1.0$,
$2.972$,
$4.864$,
$6.54$,
$6.54$,
$6.623$,
$8.201$,
$9.556$,
$10.650$,
$11.453$,
$11.944$) 

\vskip 0.7ex
\hangindent=3em \hangafter=1
$D^2= 696.547 = 
171+152c^{1}_{19}
+133c^{2}_{19}
+114c^{3}_{19}
+95c^{4}_{19}
+76c^{5}_{19}
+57c^{6}_{19}
+38c^{7}_{19}
+19c^{8}_{19}
$

\vskip 0.7ex
\hangindent=3em \hangafter=1
$T = ( 0,
\frac{1}{19},
\frac{3}{19},
\frac{7}{19},
\frac{7}{19},
\frac{6}{19},
\frac{10}{19},
\frac{15}{19},
\frac{2}{19},
\frac{9}{19},
\frac{17}{19} )
$,

\vskip 0.7ex
\hangindent=3em \hangafter=1
$S$ = ($ 1$,
$ 2+c^{1}_{19}
+c^{2}_{19}
+c^{3}_{19}
+c^{4}_{19}
+c^{5}_{19}
+c^{6}_{19}
+c^{7}_{19}
+c^{8}_{19}
$,
$ 2+2c^{1}_{19}
+c^{2}_{19}
+c^{3}_{19}
+c^{4}_{19}
+c^{5}_{19}
+c^{6}_{19}
+c^{7}_{19}
+c^{8}_{19}
$,
$ \xi_{19}^{9}$,
$ \xi_{19}^{9}$,
$ 2+2c^{1}_{19}
+c^{2}_{19}
+c^{3}_{19}
+c^{4}_{19}
+c^{5}_{19}
+c^{6}_{19}
+c^{7}_{19}
$,
$ 2+2c^{1}_{19}
+2c^{2}_{19}
+c^{3}_{19}
+c^{4}_{19}
+c^{5}_{19}
+c^{6}_{19}
+c^{7}_{19}
$,
$ 2+2c^{1}_{19}
+2c^{2}_{19}
+c^{3}_{19}
+c^{4}_{19}
+c^{5}_{19}
+c^{6}_{19}
$,
$ 2+2c^{1}_{19}
+2c^{2}_{19}
+2c^{3}_{19}
+c^{4}_{19}
+c^{5}_{19}
+c^{6}_{19}
$,
$ 2+2c^{1}_{19}
+2c^{2}_{19}
+2c^{3}_{19}
+c^{4}_{19}
+c^{5}_{19}
$,
$ 2+2c^{1}_{19}
+2c^{2}_{19}
+2c^{3}_{19}
+2c^{4}_{19}
+c^{5}_{19}
$;\ \ 
$ 2+2c^{1}_{19}
+2c^{2}_{19}
+c^{3}_{19}
+c^{4}_{19}
+c^{5}_{19}
+c^{6}_{19}
+c^{7}_{19}
$,
$ 2+2c^{1}_{19}
+2c^{2}_{19}
+2c^{3}_{19}
+c^{4}_{19}
+c^{5}_{19}
$,
$ -\xi_{19}^{9}$,
$ -\xi_{19}^{9}$,
$ 2+2c^{1}_{19}
+2c^{2}_{19}
+2c^{3}_{19}
+2c^{4}_{19}
+c^{5}_{19}
$,
$ 2+2c^{1}_{19}
+2c^{2}_{19}
+c^{3}_{19}
+c^{4}_{19}
+c^{5}_{19}
+c^{6}_{19}
$,
$ 2+2c^{1}_{19}
+c^{2}_{19}
+c^{3}_{19}
+c^{4}_{19}
+c^{5}_{19}
+c^{6}_{19}
+c^{7}_{19}
+c^{8}_{19}
$,
$ -1$,
$ -2-2  c^{1}_{19}
-c^{2}_{19}
-c^{3}_{19}
-c^{4}_{19}
-c^{5}_{19}
-c^{6}_{19}
-c^{7}_{19}
$,
$ -2-2  c^{1}_{19}
-2  c^{2}_{19}
-2  c^{3}_{19}
-c^{4}_{19}
-c^{5}_{19}
-c^{6}_{19}
$;\ \ 
$ 2+2c^{1}_{19}
+2c^{2}_{19}
+2c^{3}_{19}
+c^{4}_{19}
+c^{5}_{19}
+c^{6}_{19}
$,
$ \xi_{19}^{9}$,
$ \xi_{19}^{9}$,
$ 2+c^{1}_{19}
+c^{2}_{19}
+c^{3}_{19}
+c^{4}_{19}
+c^{5}_{19}
+c^{6}_{19}
+c^{7}_{19}
+c^{8}_{19}
$,
$ -2-2  c^{1}_{19}
-c^{2}_{19}
-c^{3}_{19}
-c^{4}_{19}
-c^{5}_{19}
-c^{6}_{19}
-c^{7}_{19}
$,
$ -2-2  c^{1}_{19}
-2  c^{2}_{19}
-2  c^{3}_{19}
-2  c^{4}_{19}
-c^{5}_{19}
$,
$ -2-2  c^{1}_{19}
-2  c^{2}_{19}
-c^{3}_{19}
-c^{4}_{19}
-c^{5}_{19}
-c^{6}_{19}
$,
$ -1$,
$ 2+2c^{1}_{19}
+2c^{2}_{19}
+c^{3}_{19}
+c^{4}_{19}
+c^{5}_{19}
+c^{6}_{19}
+c^{7}_{19}
$;\ \ 
$ s^{1}_{19}
+s^{2}_{19}
+s^{3}_{19}
+s^{4}_{19}
+2\zeta^{5}_{19}
-\zeta^{-5}_{19}
+\zeta^{6}_{19}
+2\zeta^{7}_{19}
-\zeta^{-7}_{19}
+2\zeta^{8}_{19}
-\zeta^{-8}_{19}
+\zeta^{9}_{19}
$,
$ -1-2  \zeta^{1}_{19}
-2  \zeta^{2}_{19}
-2  \zeta^{3}_{19}
-2  \zeta^{4}_{19}
-2  \zeta^{5}_{19}
+\zeta^{-5}_{19}
-\zeta^{6}_{19}
-2  \zeta^{7}_{19}
+\zeta^{-7}_{19}
-2  \zeta^{8}_{19}
+\zeta^{-8}_{19}
-\zeta^{9}_{19}
$,
$ -\xi_{19}^{9}$,
$ \xi_{19}^{9}$,
$ -\xi_{19}^{9}$,
$ \xi_{19}^{9}$,
$ -\xi_{19}^{9}$,
$ \xi_{19}^{9}$;\ \ 
$ s^{1}_{19}
+s^{2}_{19}
+s^{3}_{19}
+s^{4}_{19}
+2\zeta^{5}_{19}
-\zeta^{-5}_{19}
+\zeta^{6}_{19}
+2\zeta^{7}_{19}
-\zeta^{-7}_{19}
+2\zeta^{8}_{19}
-\zeta^{-8}_{19}
+\zeta^{9}_{19}
$,
$ -\xi_{19}^{9}$,
$ \xi_{19}^{9}$,
$ -\xi_{19}^{9}$,
$ \xi_{19}^{9}$,
$ -\xi_{19}^{9}$,
$ \xi_{19}^{9}$;\ \ 
$ -2-2  c^{1}_{19}
-2  c^{2}_{19}
-c^{3}_{19}
-c^{4}_{19}
-c^{5}_{19}
-c^{6}_{19}
$,
$ -2-2  c^{1}_{19}
-2  c^{2}_{19}
-2  c^{3}_{19}
-c^{4}_{19}
-c^{5}_{19}
-c^{6}_{19}
$,
$ 1$,
$ 2+2c^{1}_{19}
+2c^{2}_{19}
+2c^{3}_{19}
+c^{4}_{19}
+c^{5}_{19}
$,
$ 2+2c^{1}_{19}
+2c^{2}_{19}
+c^{3}_{19}
+c^{4}_{19}
+c^{5}_{19}
+c^{6}_{19}
+c^{7}_{19}
$,
$ -2-2  c^{1}_{19}
-c^{2}_{19}
-c^{3}_{19}
-c^{4}_{19}
-c^{5}_{19}
-c^{6}_{19}
-c^{7}_{19}
-c^{8}_{19}
$;\ \ 
$ 2+2c^{1}_{19}
+c^{2}_{19}
+c^{3}_{19}
+c^{4}_{19}
+c^{5}_{19}
+c^{6}_{19}
+c^{7}_{19}
+c^{8}_{19}
$,
$ 2+2c^{1}_{19}
+2c^{2}_{19}
+2c^{3}_{19}
+c^{4}_{19}
+c^{5}_{19}
$,
$ -2-c^{1}_{19}
-c^{2}_{19}
-c^{3}_{19}
-c^{4}_{19}
-c^{5}_{19}
-c^{6}_{19}
-c^{7}_{19}
-c^{8}_{19}
$,
$ -2-2  c^{1}_{19}
-2  c^{2}_{19}
-2  c^{3}_{19}
-2  c^{4}_{19}
-c^{5}_{19}
$,
$ 1$;\ \ 
$ -2-2  c^{1}_{19}
-c^{2}_{19}
-c^{3}_{19}
-c^{4}_{19}
-c^{5}_{19}
-c^{6}_{19}
-c^{7}_{19}
$,
$ -2-2  c^{1}_{19}
-2  c^{2}_{19}
-c^{3}_{19}
-c^{4}_{19}
-c^{5}_{19}
-c^{6}_{19}
-c^{7}_{19}
$,
$ 2+2c^{1}_{19}
+2c^{2}_{19}
+2c^{3}_{19}
+c^{4}_{19}
+c^{5}_{19}
+c^{6}_{19}
$,
$ 2+c^{1}_{19}
+c^{2}_{19}
+c^{3}_{19}
+c^{4}_{19}
+c^{5}_{19}
+c^{6}_{19}
+c^{7}_{19}
+c^{8}_{19}
$;\ \ 
$ 2+2c^{1}_{19}
+2c^{2}_{19}
+2c^{3}_{19}
+2c^{4}_{19}
+c^{5}_{19}
$,
$ -2-2  c^{1}_{19}
-c^{2}_{19}
-c^{3}_{19}
-c^{4}_{19}
-c^{5}_{19}
-c^{6}_{19}
-c^{7}_{19}
-c^{8}_{19}
$,
$ -2-2  c^{1}_{19}
-c^{2}_{19}
-c^{3}_{19}
-c^{4}_{19}
-c^{5}_{19}
-c^{6}_{19}
-c^{7}_{19}
$;\ \ 
$ -2-c^{1}_{19}
-c^{2}_{19}
-c^{3}_{19}
-c^{4}_{19}
-c^{5}_{19}
-c^{6}_{19}
-c^{7}_{19}
-c^{8}_{19}
$,
$ 2+2c^{1}_{19}
+2c^{2}_{19}
+c^{3}_{19}
+c^{4}_{19}
+c^{5}_{19}
+c^{6}_{19}
$;\ \ 
$ -2-2  c^{1}_{19}
-2  c^{2}_{19}
-2  c^{3}_{19}
-c^{4}_{19}
-c^{5}_{19}
$)

Realization:
Abelian anyon condensation of
$\overline{SO}(18)_3$. 

\vskip 1ex

\noindent12. $11_{3,1337.}^{48,634}$ \irep{2187}:\ \ 
$d_i$ = ($1.0$,
$6.464$,
$6.464$,
$7.464$,
$7.464$,
$12.928$,
$12.928$,
$12.928$,
$13.928$,
$14.928$,
$14.928$) 

\vskip 0.7ex
\hangindent=3em \hangafter=1
$D^2= 1337.107 = 
672+384\sqrt{3}$

\vskip 0.7ex
\hangindent=3em \hangafter=1
$T = ( 0,
0,
0,
\frac{1}{4},
\frac{1}{4},
\frac{3}{4},
\frac{3}{16},
\frac{11}{16},
0,
\frac{1}{3},
\frac{7}{12} )
$,

\vskip 0.7ex
\hangindent=3em \hangafter=1
$S$ = ($ 1$,
$ 3+2\sqrt{3}$,
$ 3+2\sqrt{3}$,
$ 4+2\sqrt{3}$,
$ 4+2\sqrt{3}$,
$ 6+4\sqrt{3}$,
$ 6+4\sqrt{3}$,
$ 6+4\sqrt{3}$,
$ 7+4\sqrt{3}$,
$ 8+4\sqrt{3}$,
$ 8+4\sqrt{3}$;\ \ 
$ -7+4\zeta^{1}_{12}
-8  \zeta^{-1}_{12}
+8\zeta^{2}_{12}
$,
$ 1-8  \zeta^{1}_{12}
+4\zeta^{-1}_{12}
-8  \zeta^{2}_{12}
$,
$(-6-4\sqrt{3})\mathrm{i}$,
$(6+4\sqrt{3})\mathrm{i}$,
$ -6-4\sqrt{3}$,
$ 6+4\sqrt{3}$,
$ 6+4\sqrt{3}$,
$ -3-2\sqrt{3}$,
$0$,
$0$;\ \ 
$ -7+4\zeta^{1}_{12}
-8  \zeta^{-1}_{12}
+8\zeta^{2}_{12}
$,
$(6+4\sqrt{3})\mathrm{i}$,
$(-6-4\sqrt{3})\mathrm{i}$,
$ -6-4\sqrt{3}$,
$ 6+4\sqrt{3}$,
$ 6+4\sqrt{3}$,
$ -3-2\sqrt{3}$,
$0$,
$0$;\ \ 
$ (-8-4\sqrt{3})\zeta_{6}^{1}$,
$ (8+4\sqrt{3})\zeta_{3}^{1}$,
$0$,
$0$,
$0$,
$ 4+2\sqrt{3}$,
$ 8+4\sqrt{3}$,
$ -8-4\sqrt{3}$;\ \ 
$ (-8-4\sqrt{3})\zeta_{6}^{1}$,
$0$,
$0$,
$0$,
$ 4+2\sqrt{3}$,
$ 8+4\sqrt{3}$,
$ -8-4\sqrt{3}$;\ \ 
$ 12+8\sqrt{3}$,
$0$,
$0$,
$ -6-4\sqrt{3}$,
$0$,
$0$;\ \ 
$ \frac{24+12\sqrt{3}}{\sqrt{6}}$,
$ \frac{-24-12\sqrt{3}}{\sqrt{6}}$,
$ -6-4\sqrt{3}$,
$0$,
$0$;\ \ 
$ \frac{24+12\sqrt{3}}{\sqrt{6}}$,
$ -6-4\sqrt{3}$,
$0$,
$0$;\ \ 
$ 1$,
$ 8+4\sqrt{3}$,
$ 8+4\sqrt{3}$;\ \ 
$ -8-4\sqrt{3}$,
$ -8-4\sqrt{3}$;\ \ 
$ 8+4\sqrt{3}$)

Realization: may be condensation reductions of $\eZ(\cNG(\Z_{12},12))$.

\vskip 1ex

\noindent13. $11_{\frac{32}{5},1964.}^{35,581}$ \irep{2077}:\ \ 
$d_i$ = ($1.0$,
$8.807$,
$8.807$,
$8.807$,
$11.632$,
$13.250$,
$14.250$,
$14.250$,
$14.250$,
$19.822$,
$20.440$) 

\vskip 0.7ex
\hangindent=3em \hangafter=1
$D^2= 1964.590 = 
910-280  c^{1}_{35}
+280c^{2}_{35}
+280c^{3}_{35}
+175c^{4}_{35}
+280c^{5}_{35}
-105  c^{6}_{35}
+490c^{7}_{35}
-280  c^{8}_{35}
+175c^{9}_{35}
+280c^{10}_{35}
$

\vskip 0.7ex
\hangindent=3em \hangafter=1
$T = ( 0,
\frac{2}{35},
\frac{22}{35},
\frac{32}{35},
\frac{1}{5},
0,
\frac{3}{7},
\frac{5}{7},
\frac{6}{7},
\frac{3}{5},
\frac{1}{5} )
$,

\vskip 0.7ex
\hangindent=3em \hangafter=1
$S$ = ($ 1$,
$ 4-c^{1}_{35}
+c^{2}_{35}
+c^{3}_{35}
+c^{4}_{35}
+c^{5}_{35}
+2c^{7}_{35}
-c^{8}_{35}
+c^{9}_{35}
+c^{10}_{35}
$,
$ 4-c^{1}_{35}
+c^{2}_{35}
+c^{3}_{35}
+c^{4}_{35}
+c^{5}_{35}
+2c^{7}_{35}
-c^{8}_{35}
+c^{9}_{35}
+c^{10}_{35}
$,
$ 4-c^{1}_{35}
+c^{2}_{35}
+c^{3}_{35}
+c^{4}_{35}
+c^{5}_{35}
+2c^{7}_{35}
-c^{8}_{35}
+c^{9}_{35}
+c^{10}_{35}
$,
$ 5-2  c^{1}_{35}
+2c^{2}_{35}
+2c^{3}_{35}
+c^{4}_{35}
+2c^{5}_{35}
-c^{6}_{35}
+3c^{7}_{35}
-2  c^{8}_{35}
+c^{9}_{35}
+2c^{10}_{35}
$,
$ 6-2  c^{1}_{35}
+2c^{2}_{35}
+2c^{3}_{35}
+c^{4}_{35}
+2c^{5}_{35}
-c^{6}_{35}
+4c^{7}_{35}
-2  c^{8}_{35}
+c^{9}_{35}
+2c^{10}_{35}
$,
$ 7-2  c^{1}_{35}
+2c^{2}_{35}
+2c^{3}_{35}
+c^{4}_{35}
+2c^{5}_{35}
-c^{6}_{35}
+4c^{7}_{35}
-2  c^{8}_{35}
+c^{9}_{35}
+2c^{10}_{35}
$,
$ 7-2  c^{1}_{35}
+2c^{2}_{35}
+2c^{3}_{35}
+c^{4}_{35}
+2c^{5}_{35}
-c^{6}_{35}
+4c^{7}_{35}
-2  c^{8}_{35}
+c^{9}_{35}
+2c^{10}_{35}
$,
$ 7-2  c^{1}_{35}
+2c^{2}_{35}
+2c^{3}_{35}
+c^{4}_{35}
+2c^{5}_{35}
-c^{6}_{35}
+4c^{7}_{35}
-2  c^{8}_{35}
+c^{9}_{35}
+2c^{10}_{35}
$,
$ 9-3  c^{1}_{35}
+3c^{2}_{35}
+3c^{3}_{35}
+2c^{4}_{35}
+3c^{5}_{35}
-c^{6}_{35}
+4c^{7}_{35}
-3  c^{8}_{35}
+2c^{9}_{35}
+3c^{10}_{35}
$,
$ 9-3  c^{1}_{35}
+3c^{2}_{35}
+3c^{3}_{35}
+2c^{4}_{35}
+3c^{5}_{35}
-c^{6}_{35}
+5c^{7}_{35}
-3  c^{8}_{35}
+2c^{9}_{35}
+3c^{10}_{35}
$;\ \ 
$ 1-c^{1}_{35}
+4c^{2}_{35}
-2  c^{3}_{35}
+c^{4}_{35}
+2c^{5}_{35}
-c^{6}_{35}
+c^{7}_{35}
+2c^{9}_{35}
-3  c^{10}_{35}
+2c^{11}_{35}
$,
$ 1+5c^{1}_{35}
+4c^{3}_{35}
+3c^{4}_{35}
+c^{5}_{35}
+4c^{6}_{35}
+2c^{8}_{35}
+3c^{10}_{35}
+c^{11}_{35}
$,
$ 5-6  c^{1}_{35}
-2  c^{2}_{35}
-3  c^{4}_{35}
-c^{5}_{35}
-4  c^{6}_{35}
+3c^{7}_{35}
-4  c^{8}_{35}
-c^{9}_{35}
+2c^{10}_{35}
-3  c^{11}_{35}
$,
$ 7-2  c^{1}_{35}
+2c^{2}_{35}
+2c^{3}_{35}
+c^{4}_{35}
+2c^{5}_{35}
-c^{6}_{35}
+4c^{7}_{35}
-2  c^{8}_{35}
+c^{9}_{35}
+2c^{10}_{35}
$,
$ 4-c^{1}_{35}
+c^{2}_{35}
+c^{3}_{35}
+c^{4}_{35}
+c^{5}_{35}
+2c^{7}_{35}
-c^{8}_{35}
+c^{9}_{35}
+c^{10}_{35}
$,
$ -1-3  c^{1}_{35}
-2  c^{3}_{35}
-2  c^{4}_{35}
-c^{5}_{35}
-2  c^{6}_{35}
-c^{8}_{35}
-2  c^{10}_{35}
$,
$ -3+4c^{1}_{35}
+c^{2}_{35}
+2c^{4}_{35}
+c^{5}_{35}
+2c^{6}_{35}
-2  c^{7}_{35}
+2c^{8}_{35}
+c^{9}_{35}
-c^{10}_{35}
+2c^{11}_{35}
$,
$ -2  c^{2}_{35}
+c^{3}_{35}
-c^{4}_{35}
-c^{5}_{35}
-2  c^{9}_{35}
+2c^{10}_{35}
-2  c^{11}_{35}
$,
$0$,
$ -7+2c^{1}_{35}
-2  c^{2}_{35}
-2  c^{3}_{35}
-c^{4}_{35}
-2  c^{5}_{35}
+c^{6}_{35}
-4  c^{7}_{35}
+2c^{8}_{35}
-c^{9}_{35}
-2  c^{10}_{35}
$;\ \ 
$ 5-6  c^{1}_{35}
-2  c^{2}_{35}
-3  c^{4}_{35}
-c^{5}_{35}
-4  c^{6}_{35}
+3c^{7}_{35}
-4  c^{8}_{35}
-c^{9}_{35}
+2c^{10}_{35}
-3  c^{11}_{35}
$,
$ 1-c^{1}_{35}
+4c^{2}_{35}
-2  c^{3}_{35}
+c^{4}_{35}
+2c^{5}_{35}
-c^{6}_{35}
+c^{7}_{35}
+2c^{9}_{35}
-3  c^{10}_{35}
+2c^{11}_{35}
$,
$ 7-2  c^{1}_{35}
+2c^{2}_{35}
+2c^{3}_{35}
+c^{4}_{35}
+2c^{5}_{35}
-c^{6}_{35}
+4c^{7}_{35}
-2  c^{8}_{35}
+c^{9}_{35}
+2c^{10}_{35}
$,
$ 4-c^{1}_{35}
+c^{2}_{35}
+c^{3}_{35}
+c^{4}_{35}
+c^{5}_{35}
+2c^{7}_{35}
-c^{8}_{35}
+c^{9}_{35}
+c^{10}_{35}
$,
$ -3+4c^{1}_{35}
+c^{2}_{35}
+2c^{4}_{35}
+c^{5}_{35}
+2c^{6}_{35}
-2  c^{7}_{35}
+2c^{8}_{35}
+c^{9}_{35}
-c^{10}_{35}
+2c^{11}_{35}
$,
$ -2  c^{2}_{35}
+c^{3}_{35}
-c^{4}_{35}
-c^{5}_{35}
-2  c^{9}_{35}
+2c^{10}_{35}
-2  c^{11}_{35}
$,
$ -1-3  c^{1}_{35}
-2  c^{3}_{35}
-2  c^{4}_{35}
-c^{5}_{35}
-2  c^{6}_{35}
-c^{8}_{35}
-2  c^{10}_{35}
$,
$0$,
$ -7+2c^{1}_{35}
-2  c^{2}_{35}
-2  c^{3}_{35}
-c^{4}_{35}
-2  c^{5}_{35}
+c^{6}_{35}
-4  c^{7}_{35}
+2c^{8}_{35}
-c^{9}_{35}
-2  c^{10}_{35}
$;\ \ 
$ 1+5c^{1}_{35}
+4c^{3}_{35}
+3c^{4}_{35}
+c^{5}_{35}
+4c^{6}_{35}
+2c^{8}_{35}
+3c^{10}_{35}
+c^{11}_{35}
$,
$ 7-2  c^{1}_{35}
+2c^{2}_{35}
+2c^{3}_{35}
+c^{4}_{35}
+2c^{5}_{35}
-c^{6}_{35}
+4c^{7}_{35}
-2  c^{8}_{35}
+c^{9}_{35}
+2c^{10}_{35}
$,
$ 4-c^{1}_{35}
+c^{2}_{35}
+c^{3}_{35}
+c^{4}_{35}
+c^{5}_{35}
+2c^{7}_{35}
-c^{8}_{35}
+c^{9}_{35}
+c^{10}_{35}
$,
$ -2  c^{2}_{35}
+c^{3}_{35}
-c^{4}_{35}
-c^{5}_{35}
-2  c^{9}_{35}
+2c^{10}_{35}
-2  c^{11}_{35}
$,
$ -1-3  c^{1}_{35}
-2  c^{3}_{35}
-2  c^{4}_{35}
-c^{5}_{35}
-2  c^{6}_{35}
-c^{8}_{35}
-2  c^{10}_{35}
$,
$ -3+4c^{1}_{35}
+c^{2}_{35}
+2c^{4}_{35}
+c^{5}_{35}
+2c^{6}_{35}
-2  c^{7}_{35}
+2c^{8}_{35}
+c^{9}_{35}
-c^{10}_{35}
+2c^{11}_{35}
$,
$0$,
$ -7+2c^{1}_{35}
-2  c^{2}_{35}
-2  c^{3}_{35}
-c^{4}_{35}
-2  c^{5}_{35}
+c^{6}_{35}
-4  c^{7}_{35}
+2c^{8}_{35}
-c^{9}_{35}
-2  c^{10}_{35}
$;\ \ 
$ -6+2c^{1}_{35}
-2  c^{2}_{35}
-2  c^{3}_{35}
-c^{4}_{35}
-2  c^{5}_{35}
+c^{6}_{35}
-4  c^{7}_{35}
+2c^{8}_{35}
-c^{9}_{35}
-2  c^{10}_{35}
$,
$ -9+3c^{1}_{35}
-3  c^{2}_{35}
-3  c^{3}_{35}
-2  c^{4}_{35}
-3  c^{5}_{35}
+c^{6}_{35}
-5  c^{7}_{35}
+3c^{8}_{35}
-2  c^{9}_{35}
-3  c^{10}_{35}
$,
$ -4+c^{1}_{35}
-c^{2}_{35}
-c^{3}_{35}
-c^{4}_{35}
-c^{5}_{35}
-2  c^{7}_{35}
+c^{8}_{35}
-c^{9}_{35}
-c^{10}_{35}
$,
$ -4+c^{1}_{35}
-c^{2}_{35}
-c^{3}_{35}
-c^{4}_{35}
-c^{5}_{35}
-2  c^{7}_{35}
+c^{8}_{35}
-c^{9}_{35}
-c^{10}_{35}
$,
$ -4+c^{1}_{35}
-c^{2}_{35}
-c^{3}_{35}
-c^{4}_{35}
-c^{5}_{35}
-2  c^{7}_{35}
+c^{8}_{35}
-c^{9}_{35}
-c^{10}_{35}
$,
$ 9-3  c^{1}_{35}
+3c^{2}_{35}
+3c^{3}_{35}
+2c^{4}_{35}
+3c^{5}_{35}
-c^{6}_{35}
+4c^{7}_{35}
-3  c^{8}_{35}
+2c^{9}_{35}
+3c^{10}_{35}
$,
$ 1$;\ \ 
$ 1$,
$ 7-2  c^{1}_{35}
+2c^{2}_{35}
+2c^{3}_{35}
+c^{4}_{35}
+2c^{5}_{35}
-c^{6}_{35}
+4c^{7}_{35}
-2  c^{8}_{35}
+c^{9}_{35}
+2c^{10}_{35}
$,
$ 7-2  c^{1}_{35}
+2c^{2}_{35}
+2c^{3}_{35}
+c^{4}_{35}
+2c^{5}_{35}
-c^{6}_{35}
+4c^{7}_{35}
-2  c^{8}_{35}
+c^{9}_{35}
+2c^{10}_{35}
$,
$ 7-2  c^{1}_{35}
+2c^{2}_{35}
+2c^{3}_{35}
+c^{4}_{35}
+2c^{5}_{35}
-c^{6}_{35}
+4c^{7}_{35}
-2  c^{8}_{35}
+c^{9}_{35}
+2c^{10}_{35}
$,
$ -9+3c^{1}_{35}
-3  c^{2}_{35}
-3  c^{3}_{35}
-2  c^{4}_{35}
-3  c^{5}_{35}
+c^{6}_{35}
-4  c^{7}_{35}
+3c^{8}_{35}
-2  c^{9}_{35}
-3  c^{10}_{35}
$,
$ -5+2c^{1}_{35}
-2  c^{2}_{35}
-2  c^{3}_{35}
-c^{4}_{35}
-2  c^{5}_{35}
+c^{6}_{35}
-3  c^{7}_{35}
+2c^{8}_{35}
-c^{9}_{35}
-2  c^{10}_{35}
$;\ \ 
$ -5+6c^{1}_{35}
+2c^{2}_{35}
+3c^{4}_{35}
+c^{5}_{35}
+4c^{6}_{35}
-3  c^{7}_{35}
+4c^{8}_{35}
+c^{9}_{35}
-2  c^{10}_{35}
+3c^{11}_{35}
$,
$ -1+c^{1}_{35}
-4  c^{2}_{35}
+2c^{3}_{35}
-c^{4}_{35}
-2  c^{5}_{35}
+c^{6}_{35}
-c^{7}_{35}
-2  c^{9}_{35}
+3c^{10}_{35}
-2  c^{11}_{35}
$,
$ -1-5  c^{1}_{35}
-4  c^{3}_{35}
-3  c^{4}_{35}
-c^{5}_{35}
-4  c^{6}_{35}
-2  c^{8}_{35}
-3  c^{10}_{35}
-c^{11}_{35}
$,
$0$,
$ 4-c^{1}_{35}
+c^{2}_{35}
+c^{3}_{35}
+c^{4}_{35}
+c^{5}_{35}
+2c^{7}_{35}
-c^{8}_{35}
+c^{9}_{35}
+c^{10}_{35}
$;\ \ 
$ -1-5  c^{1}_{35}
-4  c^{3}_{35}
-3  c^{4}_{35}
-c^{5}_{35}
-4  c^{6}_{35}
-2  c^{8}_{35}
-3  c^{10}_{35}
-c^{11}_{35}
$,
$ -5+6c^{1}_{35}
+2c^{2}_{35}
+3c^{4}_{35}
+c^{5}_{35}
+4c^{6}_{35}
-3  c^{7}_{35}
+4c^{8}_{35}
+c^{9}_{35}
-2  c^{10}_{35}
+3c^{11}_{35}
$,
$0$,
$ 4-c^{1}_{35}
+c^{2}_{35}
+c^{3}_{35}
+c^{4}_{35}
+c^{5}_{35}
+2c^{7}_{35}
-c^{8}_{35}
+c^{9}_{35}
+c^{10}_{35}
$;\ \ 
$ -1+c^{1}_{35}
-4  c^{2}_{35}
+2c^{3}_{35}
-c^{4}_{35}
-2  c^{5}_{35}
+c^{6}_{35}
-c^{7}_{35}
-2  c^{9}_{35}
+3c^{10}_{35}
-2  c^{11}_{35}
$,
$0$,
$ 4-c^{1}_{35}
+c^{2}_{35}
+c^{3}_{35}
+c^{4}_{35}
+c^{5}_{35}
+2c^{7}_{35}
-c^{8}_{35}
+c^{9}_{35}
+c^{10}_{35}
$;\ \ 
$ -9+3c^{1}_{35}
-3  c^{2}_{35}
-3  c^{3}_{35}
-2  c^{4}_{35}
-3  c^{5}_{35}
+c^{6}_{35}
-4  c^{7}_{35}
+3c^{8}_{35}
-2  c^{9}_{35}
-3  c^{10}_{35}
$,
$ 9-3  c^{1}_{35}
+3c^{2}_{35}
+3c^{3}_{35}
+2c^{4}_{35}
+3c^{5}_{35}
-c^{6}_{35}
+4c^{7}_{35}
-3  c^{8}_{35}
+2c^{9}_{35}
+3c^{10}_{35}
$;\ \ 
$ -6+2c^{1}_{35}
-2  c^{2}_{35}
-2  c^{3}_{35}
-c^{4}_{35}
-2  c^{5}_{35}
+c^{6}_{35}
-4  c^{7}_{35}
+2c^{8}_{35}
-c^{9}_{35}
-2  c^{10}_{35}
$)

Realization: unknown

\vskip 1ex 

}



\subsection{Rank 12}\label{ss:rank12}

{\small

\noindent1. $12_{1,40.}^{80,190}$ \irep{5194}:\ \ 
$d_i$ = ($1.0$,
$1.0$,
$1.0$,
$1.0$,
$2.0$,
$2.0$,
$2.0$,
$2.0$,
$2.236$,
$2.236$,
$2.236$,
$2.236$) 

\vskip 0.7ex
\hangindent=3em \hangafter=1
$D^2= 40.0 = 
40$

\vskip 0.7ex
\hangindent=3em \hangafter=1
$T = ( 0,
0,
\frac{1}{4},
\frac{1}{4},
\frac{1}{5},
\frac{4}{5},
\frac{1}{20},
\frac{9}{20},
\frac{1}{16},
\frac{1}{16},
\frac{9}{16},
\frac{9}{16} )
$,

\vskip 0.7ex
\hangindent=3em \hangafter=1
$S$ = ($ 1$,
$ 1$,
$ 1$,
$ 1$,
$ 2$,
$ 2$,
$ 2$,
$ 2$,
$ \sqrt{5}$,
$ \sqrt{5}$,
$ \sqrt{5}$,
$ \sqrt{5}$;\ \ 
$ 1$,
$ 1$,
$ 1$,
$ 2$,
$ 2$,
$ 2$,
$ 2$,
$ -\sqrt{5}$,
$ -\sqrt{5}$,
$ -\sqrt{5}$,
$ -\sqrt{5}$;\ \ 
$ -1$,
$ -1$,
$ 2$,
$ 2$,
$ -2$,
$ -2$,
$(-\sqrt{5})\mathrm{i}$,
$(\sqrt{5})\mathrm{i}$,
$(-\sqrt{5})\mathrm{i}$,
$(\sqrt{5})\mathrm{i}$;\ \ 
$ -1$,
$ 2$,
$ 2$,
$ -2$,
$ -2$,
$(\sqrt{5})\mathrm{i}$,
$(-\sqrt{5})\mathrm{i}$,
$(\sqrt{5})\mathrm{i}$,
$(-\sqrt{5})\mathrm{i}$;\ \ 
$  -1-\sqrt{5} $,
$  -1+\sqrt{5} $,
$  -1+\sqrt{5} $,
$  -1-\sqrt{5} $,
$0$,
$0$,
$0$,
$0$;\ \ 
$  -1-\sqrt{5} $,
$  -1-\sqrt{5} $,
$  -1+\sqrt{5} $,
$0$,
$0$,
$0$,
$0$;\ \ 
$  1+\sqrt{5} $,
$  1-\sqrt{5} $,
$0$,
$0$,
$0$,
$0$;\ \ 
$  1+\sqrt{5} $,
$0$,
$0$,
$0$,
$0$;\ \ 
$ -\sqrt{5}\zeta_{8}^{3}$,
$ \sqrt{5}\zeta_{8}^{1}$,
$ \sqrt{5}\zeta_{8}^{3}$,
$ -\sqrt{5}\zeta_{8}^{1}$;\ \ 
$ -\sqrt{5}\zeta_{8}^{3}$,
$ -\sqrt{5}\zeta_{8}^{1}$,
$ \sqrt{5}\zeta_{8}^{3}$;\ \ 
$ -\sqrt{5}\zeta_{8}^{3}$,
$ \sqrt{5}\zeta_{8}^{1}$;\ \ 
$ -\sqrt{5}\zeta_{8}^{3}$)

Connected to the orbit of $12_{5,40.}^{80,348}$ via a change of spherical
structure.

Realization:
$SO(10)_2$ (see \cite[Section 3]{BPR} for explicit data) or anyon condensation of $O_5$.

\vskip 1ex 

{\grey
\noindent2. $12_{5,40.}^{80,348}$ \irep{5194}:\ \ 
$d_i$ = ($1.0$,
$1.0$,
$1.0$,
$1.0$,
$2.0$,
$2.0$,
$2.0$,
$2.0$,
$2.236$,
$2.236$,
$2.236$,
$2.236$) 

\vskip 0.7ex
\hangindent=3em \hangafter=1
$D^2= 40.0 = 
40$

\vskip 0.7ex
\hangindent=3em \hangafter=1
$T = ( 0,
0,
\frac{1}{4},
\frac{1}{4},
\frac{2}{5},
\frac{3}{5},
\frac{13}{20},
\frac{17}{20},
\frac{1}{16},
\frac{1}{16},
\frac{9}{16},
\frac{9}{16} )
$,

\vskip 0.7ex
\hangindent=3em \hangafter=1
$S$ = ($ 1$,
$ 1$,
$ 1$,
$ 1$,
$ 2$,
$ 2$,
$ 2$,
$ 2$,
$ \sqrt{5}$,
$ \sqrt{5}$,
$ \sqrt{5}$,
$ \sqrt{5}$;\ \ 
$ 1$,
$ 1$,
$ 1$,
$ 2$,
$ 2$,
$ 2$,
$ 2$,
$ -\sqrt{5}$,
$ -\sqrt{5}$,
$ -\sqrt{5}$,
$ -\sqrt{5}$;\ \ 
$ -1$,
$ -1$,
$ 2$,
$ 2$,
$ -2$,
$ -2$,
$(-\sqrt{5})\mathrm{i}$,
$(\sqrt{5})\mathrm{i}$,
$(-\sqrt{5})\mathrm{i}$,
$(\sqrt{5})\mathrm{i}$;\ \ 
$ -1$,
$ 2$,
$ 2$,
$ -2$,
$ -2$,
$(\sqrt{5})\mathrm{i}$,
$(-\sqrt{5})\mathrm{i}$,
$(\sqrt{5})\mathrm{i}$,
$(-\sqrt{5})\mathrm{i}$;\ \ 
$  -1+\sqrt{5} $,
$  -1-\sqrt{5} $,
$  -1+\sqrt{5} $,
$  -1-\sqrt{5} $,
$0$,
$0$,
$0$,
$0$;\ \ 
$  -1+\sqrt{5} $,
$  -1-\sqrt{5} $,
$  -1+\sqrt{5} $,
$0$,
$0$,
$0$,
$0$;\ \ 
$  1-\sqrt{5} $,
$  1+\sqrt{5} $,
$0$,
$0$,
$0$,
$0$;\ \ 
$  1-\sqrt{5} $,
$0$,
$0$,
$0$,
$0$;\ \ 
$ \sqrt{5}\zeta_{8}^{3}$,
$ -\sqrt{5}\zeta_{8}^{1}$,
$ -\sqrt{5}\zeta_{8}^{3}$,
$ \sqrt{5}\zeta_{8}^{1}$;\ \ 
$ \sqrt{5}\zeta_{8}^{3}$,
$ \sqrt{5}\zeta_{8}^{1}$,
$ -\sqrt{5}\zeta_{8}^{3}$;\ \ 
$ \sqrt{5}\zeta_{8}^{3}$,
$ -\sqrt{5}\zeta_{8}^{1}$;\ \ 
$ \sqrt{5}\zeta_{8}^{3}$)

Connected to the orbit of $12_{1,40.}^{80,190}$ via a change of spherical
structure.

Realization:
Abelian anyon condensation of $SO(5)_2$.
}

\vskip 1ex

\noindent3. $12_{0,68.}^{34,116}$ \irep{4690}:\ \ 
$d_i$ = ($1.0$,
$1.0$,
$2.0$,
$2.0$,
$2.0$,
$2.0$,
$2.0$,
$2.0$,
$2.0$,
$2.0$,
$4.123$,
$4.123$) 

\vskip 0.7ex
\hangindent=3em \hangafter=1
$D^2= 68.0 = 
68$

\vskip 0.7ex
\hangindent=3em \hangafter=1
$T = ( 0,
0,
\frac{1}{17},
\frac{2}{17},
\frac{4}{17},
\frac{8}{17},
\frac{9}{17},
\frac{13}{17},
\frac{15}{17},
\frac{16}{17},
0,
\frac{1}{2} )
$,

\vskip 0.7ex
\hangindent=3em \hangafter=1
$S$ = ($ 1$,
$ 1$,
$ 2$,
$ 2$,
$ 2$,
$ 2$,
$ 2$,
$ 2$,
$ 2$,
$ 2$,
$ \sqrt{17}$,
$ \sqrt{17}$;\ \ 
$ 1$,
$ 2$,
$ 2$,
$ 2$,
$ 2$,
$ 2$,
$ 2$,
$ 2$,
$ 2$,
$ -\sqrt{17}$,
$ -\sqrt{17}$;\ \ 
$ 2c_{17}^{2}$,
$ 2c_{17}^{5}$,
$ 2c_{17}^{4}$,
$ 2c_{17}^{7}$,
$ 2c_{17}^{6}$,
$ 2c_{17}^{1}$,
$ 2c_{17}^{3}$,
$ 2c_{17}^{8}$,
$0$,
$0$;\ \ 
$ 2c_{17}^{4}$,
$ 2c_{17}^{7}$,
$ 2c_{17}^{8}$,
$ 2c_{17}^{2}$,
$ 2c_{17}^{6}$,
$ 2c_{17}^{1}$,
$ 2c_{17}^{3}$,
$0$,
$0$;\ \ 
$ 2c_{17}^{8}$,
$ 2c_{17}^{3}$,
$ 2c_{17}^{5}$,
$ 2c_{17}^{2}$,
$ 2c_{17}^{6}$,
$ 2c_{17}^{1}$,
$0$,
$0$;\ \ 
$ 2c_{17}^{1}$,
$ 2c_{17}^{4}$,
$ 2c_{17}^{5}$,
$ 2c_{17}^{2}$,
$ 2c_{17}^{6}$,
$0$,
$0$;\ \ 
$ 2c_{17}^{1}$,
$ 2c_{17}^{3}$,
$ 2c_{17}^{8}$,
$ 2c_{17}^{7}$,
$0$,
$0$;\ \ 
$ 2c_{17}^{8}$,
$ 2c_{17}^{7}$,
$ 2c_{17}^{4}$,
$0$,
$0$;\ \ 
$ 2c_{17}^{4}$,
$ 2c_{17}^{5}$,
$0$,
$0$;\ \ 
$ 2c_{17}^{2}$,
$0$,
$0$;\ \ 
$ \sqrt{17}$,
$ -\sqrt{17}$;\ \ 
$ \sqrt{17}$)

Connected to the orbit of $12_{4,68.}^{34,824}$ via a change of spherical
structure.

Realization:
$SO(17)_2$,
or Abelian anyon condensation of $\overline{SO(17)_2}$, $O_{17}$.

\vskip 1ex 

{\grey
\noindent4. $12_{4,68.}^{34,824}$ \irep{4690}:\ \ 
$d_i$ = ($1.0$,
$1.0$,
$2.0$,
$2.0$,
$2.0$,
$2.0$,
$2.0$,
$2.0$,
$2.0$,
$2.0$,
$4.123$,
$4.123$) 

\vskip 0.7ex
\hangindent=3em \hangafter=1
$D^2= 68.0 = 
68$

\vskip 0.7ex
\hangindent=3em \hangafter=1
$T = ( 0,
0,
\frac{3}{17},
\frac{5}{17},
\frac{6}{17},
\frac{7}{17},
\frac{10}{17},
\frac{11}{17},
\frac{12}{17},
\frac{14}{17},
0,
\frac{1}{2} )
$,

\vskip 0.7ex
\hangindent=3em \hangafter=1
$S$ = ($ 1$,
$ 1$,
$ 2$,
$ 2$,
$ 2$,
$ 2$,
$ 2$,
$ 2$,
$ 2$,
$ 2$,
$ \sqrt{17}$,
$ \sqrt{17}$;\ \ 
$ 1$,
$ 2$,
$ 2$,
$ 2$,
$ 2$,
$ 2$,
$ 2$,
$ 2$,
$ 2$,
$ -\sqrt{17}$,
$ -\sqrt{17}$;\ \ 
$ 2c_{17}^{6}$,
$ 2c_{17}^{3}$,
$ 2c_{17}^{2}$,
$ 2c_{17}^{4}$,
$ 2c_{17}^{1}$,
$ 2c_{17}^{8}$,
$ 2c_{17}^{5}$,
$ 2c_{17}^{7}$,
$0$,
$0$;\ \ 
$ 2c_{17}^{7}$,
$ 2c_{17}^{1}$,
$ 2c_{17}^{2}$,
$ 2c_{17}^{8}$,
$ 2c_{17}^{4}$,
$ 2c_{17}^{6}$,
$ 2c_{17}^{5}$,
$0$,
$0$;\ \ 
$ 2c_{17}^{5}$,
$ 2c_{17}^{7}$,
$ 2c_{17}^{6}$,
$ 2c_{17}^{3}$,
$ 2c_{17}^{4}$,
$ 2c_{17}^{8}$,
$0$,
$0$;\ \ 
$ 2c_{17}^{3}$,
$ 2c_{17}^{5}$,
$ 2c_{17}^{6}$,
$ 2c_{17}^{8}$,
$ 2c_{17}^{1}$,
$0$,
$0$;\ \ 
$ 2c_{17}^{3}$,
$ 2c_{17}^{7}$,
$ 2c_{17}^{2}$,
$ 2c_{17}^{4}$,
$0$,
$0$;\ \ 
$ 2c_{17}^{5}$,
$ 2c_{17}^{1}$,
$ 2c_{17}^{2}$,
$0$,
$0$;\ \ 
$ 2c_{17}^{7}$,
$ 2c_{17}^{3}$,
$0$,
$0$;\ \ 
$ 2c_{17}^{6}$,
$0$,
$0$;\ \ 
$ -\sqrt{17}$,
$ \sqrt{17}$;\ \ 
$ -\sqrt{17}$)

Connected to the orbit of $12_{0,68.}^{34,116}$ via a change of spherical
structure.

Realization: $SO(17)_2$ with a Galois conjugation and a change of spherical
structure.
}

\vskip 1ex

\noindent5. $12_{0,68.}^{68,166}$ \irep{5150}:\ \ 
$d_i$ = ($1.0$,
$1.0$,
$2.0$,
$2.0$,
$2.0$,
$2.0$,
$2.0$,
$2.0$,
$2.0$,
$2.0$,
$4.123$,
$4.123$) 

\vskip 0.7ex
\hangindent=3em \hangafter=1
$D^2= 68.0 = 
68$

\vskip 0.7ex
\hangindent=3em \hangafter=1
$T = ( 0,
0,
\frac{1}{17},
\frac{2}{17},
\frac{4}{17},
\frac{8}{17},
\frac{9}{17},
\frac{13}{17},
\frac{15}{17},
\frac{16}{17},
\frac{1}{4},
\frac{3}{4} )
$,

\vskip 0.7ex
\hangindent=3em \hangafter=1
$S$ = ($ 1$,
$ 1$,
$ 2$,
$ 2$,
$ 2$,
$ 2$,
$ 2$,
$ 2$,
$ 2$,
$ 2$,
$ \sqrt{17}$,
$ \sqrt{17}$;\ \ 
$ 1$,
$ 2$,
$ 2$,
$ 2$,
$ 2$,
$ 2$,
$ 2$,
$ 2$,
$ 2$,
$ -\sqrt{17}$,
$ -\sqrt{17}$;\ \ 
$ 2c_{17}^{2}$,
$ 2c_{17}^{5}$,
$ 2c_{17}^{4}$,
$ 2c_{17}^{7}$,
$ 2c_{17}^{6}$,
$ 2c_{17}^{1}$,
$ 2c_{17}^{3}$,
$ 2c_{17}^{8}$,
$0$,
$0$;\ \ 
$ 2c_{17}^{4}$,
$ 2c_{17}^{7}$,
$ 2c_{17}^{8}$,
$ 2c_{17}^{2}$,
$ 2c_{17}^{6}$,
$ 2c_{17}^{1}$,
$ 2c_{17}^{3}$,
$0$,
$0$;\ \ 
$ 2c_{17}^{8}$,
$ 2c_{17}^{3}$,
$ 2c_{17}^{5}$,
$ 2c_{17}^{2}$,
$ 2c_{17}^{6}$,
$ 2c_{17}^{1}$,
$0$,
$0$;\ \ 
$ 2c_{17}^{1}$,
$ 2c_{17}^{4}$,
$ 2c_{17}^{5}$,
$ 2c_{17}^{2}$,
$ 2c_{17}^{6}$,
$0$,
$0$;\ \ 
$ 2c_{17}^{1}$,
$ 2c_{17}^{3}$,
$ 2c_{17}^{8}$,
$ 2c_{17}^{7}$,
$0$,
$0$;\ \ 
$ 2c_{17}^{8}$,
$ 2c_{17}^{7}$,
$ 2c_{17}^{4}$,
$0$,
$0$;\ \ 
$ 2c_{17}^{4}$,
$ 2c_{17}^{5}$,
$0$,
$0$;\ \ 
$ 2c_{17}^{2}$,
$0$,
$0$;\ \ 
$ -\sqrt{17}$,
$ \sqrt{17}$;\ \ 
$ -\sqrt{17}$)

Connected to the orbit of $12_{4,68.}^{68,721}$ via a change of spherical
structure.

Realization:
Abelian anyon condensation of $SO(17)_2$, $\overline{SO(17)}_2$ or  $O_{17}$.

\vskip 1ex 

{\grey 
\noindent6. $12_{4,68.}^{68,721}$ \irep{5150}:\ \ 
$d_i$ = ($1.0$,
$1.0$,
$2.0$,
$2.0$,
$2.0$,
$2.0$,
$2.0$,
$2.0$,
$2.0$,
$2.0$,
$4.123$,
$4.123$) 

\vskip 0.7ex
\hangindent=3em \hangafter=1
$D^2= 68.0 = 
68$

\vskip 0.7ex
\hangindent=3em \hangafter=1
$T = ( 0,
0,
\frac{3}{17},
\frac{5}{17},
\frac{6}{17},
\frac{7}{17},
\frac{10}{17},
\frac{11}{17},
\frac{12}{17},
\frac{14}{17},
\frac{1}{4},
\frac{3}{4} )
$,

\vskip 0.7ex
\hangindent=3em \hangafter=1
$S$ = ($ 1$,
$ 1$,
$ 2$,
$ 2$,
$ 2$,
$ 2$,
$ 2$,
$ 2$,
$ 2$,
$ 2$,
$ \sqrt{17}$,
$ \sqrt{17}$;\ \ 
$ 1$,
$ 2$,
$ 2$,
$ 2$,
$ 2$,
$ 2$,
$ 2$,
$ 2$,
$ 2$,
$ -\sqrt{17}$,
$ -\sqrt{17}$;\ \ 
$ 2c_{17}^{6}$,
$ 2c_{17}^{3}$,
$ 2c_{17}^{2}$,
$ 2c_{17}^{4}$,
$ 2c_{17}^{1}$,
$ 2c_{17}^{8}$,
$ 2c_{17}^{5}$,
$ 2c_{17}^{7}$,
$0$,
$0$;\ \ 
$ 2c_{17}^{7}$,
$ 2c_{17}^{1}$,
$ 2c_{17}^{2}$,
$ 2c_{17}^{8}$,
$ 2c_{17}^{4}$,
$ 2c_{17}^{6}$,
$ 2c_{17}^{5}$,
$0$,
$0$;\ \ 
$ 2c_{17}^{5}$,
$ 2c_{17}^{7}$,
$ 2c_{17}^{6}$,
$ 2c_{17}^{3}$,
$ 2c_{17}^{4}$,
$ 2c_{17}^{8}$,
$0$,
$0$;\ \ 
$ 2c_{17}^{3}$,
$ 2c_{17}^{5}$,
$ 2c_{17}^{6}$,
$ 2c_{17}^{8}$,
$ 2c_{17}^{1}$,
$0$,
$0$;\ \ 
$ 2c_{17}^{3}$,
$ 2c_{17}^{7}$,
$ 2c_{17}^{2}$,
$ 2c_{17}^{4}$,
$0$,
$0$;\ \ 
$ 2c_{17}^{5}$,
$ 2c_{17}^{1}$,
$ 2c_{17}^{2}$,
$0$,
$0$;\ \ 
$ 2c_{17}^{7}$,
$ 2c_{17}^{3}$,
$0$,
$0$;\ \ 
$ 2c_{17}^{6}$,
$0$,
$0$;\ \ 
$ \sqrt{17}$,
$ -\sqrt{17}$;\ \ 
$ \sqrt{17}$)

Connected to the orbit of $12_{0,68.}^{68,166}$ via a change of spherical
structure.

Realization: anyon condensation of $SO(17)_2$ with a Galois conjugation and a change of spherical
structure.
}

\vskip 1ex

\noindent7. $12_{4,144.}^{48,120}$ \irep{4970}:\ \ 
$d_i$ = ($1.0$,
$1.0$,
$2.0$,
$3.0$,
$3.0$,
$4.0$,
$4.0$,
$4.0$,
$4.242$,
$4.242$,
$4.242$,
$4.242$) 

\vskip 0.7ex
\hangindent=3em \hangafter=1
$D^2= 144.0 = 
144$

\vskip 0.7ex
\hangindent=3em \hangafter=1
$T = ( 0,
0,
0,
\frac{1}{2},
\frac{1}{2},
0,
\frac{1}{3},
\frac{2}{3},
\frac{1}{16},
\frac{7}{16},
\frac{9}{16},
\frac{15}{16} )
$,

\vskip 0.7ex
\hangindent=3em \hangafter=1
$S$ = ($ 1$,
$ 1$,
$ 2$,
$ 3$,
$ 3$,
$ 4$,
$ 4$,
$ 4$,
$ 3\sqrt{2}$,
$ 3\sqrt{2}$,
$ 3\sqrt{2}$,
$ 3\sqrt{2}$;\ \ 
$ 1$,
$ 2$,
$ 3$,
$ 3$,
$ 4$,
$ 4$,
$ 4$,
$ -3\sqrt{2}$,
$ -3\sqrt{2}$,
$ -3\sqrt{2}$,
$ -3\sqrt{2}$;\ \ 
$ 4$,
$ 6$,
$ 6$,
$ -4$,
$ -4$,
$ -4$,
$0$,
$0$,
$0$,
$0$;\ \ 
$ -3$,
$ -3$,
$0$,
$0$,
$0$,
$ -3\sqrt{2}$,
$ 3\sqrt{2}$,
$ -3\sqrt{2}$,
$ 3\sqrt{2}$;\ \ 
$ -3$,
$0$,
$0$,
$0$,
$ 3\sqrt{2}$,
$ -3\sqrt{2}$,
$ 3\sqrt{2}$,
$ -3\sqrt{2}$;\ \ 
$ -8$,
$ 4$,
$ 4$,
$0$,
$0$,
$0$,
$0$;\ \ 
$ 4$,
$ -8$,
$0$,
$0$,
$0$,
$0$;\ \ 
$ 4$,
$0$,
$0$,
$0$,
$0$;\ \ 
$0$,
$ 6$,
$0$,
$ -6$;\ \ 
$0$,
$ -6$,
$0$;\ \ 
$0$,
$ 6$;\ \ 
$0$)

Connected to the orbit of $12_{4,144.}^{48,650}$ via a change of spherical
structure.

Realization: $S_3$-gauging of the rank-4 3-fermion MTC $4^{2,250}
_{4,4.}$ with a different minimal modular extension.
\vskip 1ex 

{\grey
\noindent8. $12_{4,144.}^{48,650}$ \irep{4970}:\ \ 
$d_i$ = ($1.0$,
$1.0$,
$2.0$,
$3.0$,
$3.0$,
$4.0$,
$4.0$,
$4.0$,
$4.242$,
$4.242$,
$4.242$,
$4.242$) 

\vskip 0.7ex
\hangindent=3em \hangafter=1
$D^2= 144.0 = 
144$

\vskip 0.7ex
\hangindent=3em \hangafter=1
$T = ( 0,
0,
0,
\frac{1}{2},
\frac{1}{2},
0,
\frac{1}{3},
\frac{2}{3},
\frac{3}{16},
\frac{5}{16},
\frac{11}{16},
\frac{13}{16} )
$,

\vskip 0.7ex
\hangindent=3em \hangafter=1
$S$ = ($ 1$,
$ 1$,
$ 2$,
$ 3$,
$ 3$,
$ 4$,
$ 4$,
$ 4$,
$ 3\sqrt{2}$,
$ 3\sqrt{2}$,
$ 3\sqrt{2}$,
$ 3\sqrt{2}$;\ \ 
$ 1$,
$ 2$,
$ 3$,
$ 3$,
$ 4$,
$ 4$,
$ 4$,
$ -3\sqrt{2}$,
$ -3\sqrt{2}$,
$ -3\sqrt{2}$,
$ -3\sqrt{2}$;\ \ 
$ 4$,
$ 6$,
$ 6$,
$ -4$,
$ -4$,
$ -4$,
$0$,
$0$,
$0$,
$0$;\ \ 
$ -3$,
$ -3$,
$0$,
$0$,
$0$,
$ -3\sqrt{2}$,
$ 3\sqrt{2}$,
$ -3\sqrt{2}$,
$ 3\sqrt{2}$;\ \ 
$ -3$,
$0$,
$0$,
$0$,
$ 3\sqrt{2}$,
$ -3\sqrt{2}$,
$ 3\sqrt{2}$,
$ -3\sqrt{2}$;\ \ 
$ -8$,
$ 4$,
$ 4$,
$0$,
$0$,
$0$,
$0$;\ \ 
$ 4$,
$ -8$,
$0$,
$0$,
$0$,
$0$;\ \ 
$ 4$,
$0$,
$0$,
$0$,
$0$;\ \ 
$0$,
$ 6$,
$0$,
$ -6$;\ \ 
$0$,
$ -6$,
$0$;\ \ 
$0$,
$ 6$;\ \ 
$0$)

Connected to the orbit of $12_{4,144.}^{48,120}$ via a change of spherical
structure.

Realization: $S_3$-gauging of the rank-4 3-fermion MTC $4^{2,250}
_{4,4.}$ with a different minimal modular extension. Also constructed by condensing the diagonal copy of $\Rep(S_3)$ in
$\cC\boxtimes \Rep(D^\omega S_3)$ where $\cC$ is $12_{4,144.}^{144,916}$.
}
\vskip 1ex

\noindent9. $12_{4,144.}^{144,916}$ \irep{5262}:\ \ 
$d_i$ = ($1.0$,
$1.0$,
$2.0$,
$3.0$,
$3.0$,
$4.0$,
$4.0$,
$4.0$,
$4.242$,
$4.242$,
$4.242$,
$4.242$) 

\vskip 0.7ex
\hangindent=3em \hangafter=1
$D^2= 144.0 = 
144$

\vskip 0.7ex
\hangindent=3em \hangafter=1
$T = ( 0,
0,
0,
\frac{1}{2},
\frac{1}{2},
\frac{1}{9},
\frac{4}{9},
\frac{7}{9},
\frac{1}{16},
\frac{7}{16},
\frac{9}{16},
\frac{15}{16} )
$,

\vskip 0.7ex
\hangindent=3em \hangafter=1
$S$ = ($ 1$,
$ 1$,
$ 2$,
$ 3$,
$ 3$,
$ 4$,
$ 4$,
$ 4$,
$ 3\sqrt{2}$,
$ 3\sqrt{2}$,
$ 3\sqrt{2}$,
$ 3\sqrt{2}$;\ \ 
$ 1$,
$ 2$,
$ 3$,
$ 3$,
$ 4$,
$ 4$,
$ 4$,
$ -3\sqrt{2}$,
$ -3\sqrt{2}$,
$ -3\sqrt{2}$,
$ -3\sqrt{2}$;\ \ 
$ 4$,
$ 6$,
$ 6$,
$ -4$,
$ -4$,
$ -4$,
$0$,
$0$,
$0$,
$0$;\ \ 
$ -3$,
$ -3$,
$0$,
$0$,
$0$,
$ -3\sqrt{2}$,
$ 3\sqrt{2}$,
$ -3\sqrt{2}$,
$ 3\sqrt{2}$;\ \ 
$ -3$,
$0$,
$0$,
$0$,
$ 3\sqrt{2}$,
$ -3\sqrt{2}$,
$ 3\sqrt{2}$,
$ -3\sqrt{2}$;\ \ 
$ -4c_{9}^{2}$,
$ -4c_{9}^{4}$,
$ -4c_{9}^{1}$,
$0$,
$0$,
$0$,
$0$;\ \ 
$ -4c_{9}^{1}$,
$ -4c_{9}^{2}$,
$0$,
$0$,
$0$,
$0$;\ \ 
$ -4c_{9}^{4}$,
$0$,
$0$,
$0$,
$0$;\ \ 
$0$,
$ 6$,
$0$,
$ -6$;\ \ 
$0$,
$ -6$,
$0$;\ \ 
$0$,
$ 6$;\ \ 
$0$)

Connected to the orbit of $12_{4,144.}^{144,386}$ via a change of spherical
structure.

Realization: $S_3$-gauging of the rank-4 3-fermion MTC $4^{2,250}
_{4,4.}$, see \cite{CuiGalindoPlavnikWang}.
\vskip 1ex 

{\grey
\noindent10. $12_{4,144.}^{144,386}$ \irep{5262}:\ \ 
$d_i$ = ($1.0$,
$1.0$,
$2.0$,
$3.0$,
$3.0$,
$4.0$,
$4.0$,
$4.0$,
$4.242$,
$4.242$,
$4.242$,
$4.242$) 

\vskip 0.7ex
\hangindent=3em \hangafter=1
$D^2= 144.0 = 
144$

\vskip 0.7ex
\hangindent=3em \hangafter=1
$T = ( 0,
0,
0,
\frac{1}{2},
\frac{1}{2},
\frac{1}{9},
\frac{4}{9},
\frac{7}{9},
\frac{3}{16},
\frac{5}{16},
\frac{11}{16},
\frac{13}{16} )
$,

\vskip 0.7ex
\hangindent=3em \hangafter=1
$S$ = ($ 1$,
$ 1$,
$ 2$,
$ 3$,
$ 3$,
$ 4$,
$ 4$,
$ 4$,
$ 3\sqrt{2}$,
$ 3\sqrt{2}$,
$ 3\sqrt{2}$,
$ 3\sqrt{2}$;\ \ 
$ 1$,
$ 2$,
$ 3$,
$ 3$,
$ 4$,
$ 4$,
$ 4$,
$ -3\sqrt{2}$,
$ -3\sqrt{2}$,
$ -3\sqrt{2}$,
$ -3\sqrt{2}$;\ \ 
$ 4$,
$ 6$,
$ 6$,
$ -4$,
$ -4$,
$ -4$,
$0$,
$0$,
$0$,
$0$;\ \ 
$ -3$,
$ -3$,
$0$,
$0$,
$0$,
$ -3\sqrt{2}$,
$ 3\sqrt{2}$,
$ -3\sqrt{2}$,
$ 3\sqrt{2}$;\ \ 
$ -3$,
$0$,
$0$,
$0$,
$ 3\sqrt{2}$,
$ -3\sqrt{2}$,
$ 3\sqrt{2}$,
$ -3\sqrt{2}$;\ \ 
$ -4c_{9}^{2}$,
$ -4c_{9}^{4}$,
$ -4c_{9}^{1}$,
$0$,
$0$,
$0$,
$0$;\ \ 
$ -4c_{9}^{1}$,
$ -4c_{9}^{2}$,
$0$,
$0$,
$0$,
$0$;\ \ 
$ -4c_{9}^{4}$,
$0$,
$0$,
$0$,
$0$;\ \ 
$0$,
$ 6$,
$0$,
$ -6$;\ \ 
$0$,
$ -6$,
$0$;\ \ 
$0$,
$ 6$;\ \ 
$0$)

Connected to the orbit of $12_{4,144.}^{144,916}$ via a change of spherical
structure.

Realization: $S_3$-gauging of the rank-4 3-fermion MTC $4^{2,250}
_{4,4.}$ with a different minimal modular extension. 
}

\vskip 1ex

\noindent11. $12_{\frac{94}{25},397.8}^{25,285}$ \irep{4446}:\ \ 
$d_i$ = ($1.0$,
$1.984$,
$2.937$,
$3.843$,
$4.689$,
$5.461$,
$6.147$,
$6.736$,
$7.219$,
$7.588$,
$7.837$,
$7.962$) 

\vskip 0.7ex
\hangindent=3em \hangafter=1
$D^2= 397.875 = 
 75+65c^{1}_{25}
+55c^{2}_{25}
+45c^{3}_{25}
+35c^{4}_{25}
+25c^{1}_{5}
+20c^{6}_{25}
+15c^{7}_{25}
+10c^{8}_{25}
+5c^{9}_{25}
$

\vskip 0.7ex
\hangindent=3em \hangafter=1
$T = ( 0,
\frac{7}{25},
\frac{2}{25},
\frac{2}{5},
\frac{6}{25},
\frac{3}{5},
\frac{12}{25},
\frac{22}{25},
\frac{4}{5},
\frac{6}{25},
\frac{1}{5},
\frac{17}{25} )
$,

\vskip 0.7ex
\hangindent=3em \hangafter=1
$S$ = ($ 1$,
$ -c_{25}^{12}$,
$ \xi_{25}^{3}$,
$ \xi_{25}^{21}$,
$ \xi_{25}^{5}$,
$ \xi_{25}^{19}$,
$ \xi_{25}^{7}$,
$ \xi_{25}^{17}$,
$ \xi_{25}^{9}$,
$ \frac{1+\sqrt{5}}{2}\xi_{25}^{5}$,
$ \xi_{25}^{11}$,
$ \xi_{25}^{13}$;\ \ 
$ -\xi_{25}^{21}$,
$ \xi_{25}^{19}$,
$ -\xi_{25}^{17}$,
$ \frac{1+\sqrt{5}}{2}\xi_{25}^{5}$,
$ -\xi_{25}^{13}$,
$ \xi_{25}^{11}$,
$ -\xi_{25}^{9}$,
$ \xi_{25}^{7}$,
$ -\xi_{25}^{5}$,
$ \xi_{25}^{3}$,
$ -1$;\ \ 
$ \xi_{25}^{9}$,
$ \xi_{25}^{13}$,
$ \frac{1+\sqrt{5}}{2}\xi_{25}^{5}$,
$ \xi_{25}^{7}$,
$ \xi_{25}^{21}$,
$ 1$,
$ c_{25}^{12}$,
$ -\xi_{25}^{5}$,
$ -\xi_{25}^{17}$,
$ -\xi_{25}^{11}$;\ \ 
$ -\xi_{25}^{9}$,
$ \xi_{25}^{5}$,
$ -1$,
$ -\xi_{25}^{3}$,
$ \xi_{25}^{7}$,
$ -\xi_{25}^{11}$,
$ \frac{1+\sqrt{5}}{2}\xi_{25}^{5}$,
$ -\xi_{25}^{19}$,
$ -c_{25}^{12}$;\ \ 
$0$,
$ -\xi_{25}^{5}$,
$ -\frac{1+\sqrt{5}}{2}\xi_{25}^{5}$,
$ -\frac{1+\sqrt{5}}{2}\xi_{25}^{5}$,
$ -\xi_{25}^{5}$,
$0$,
$ \xi_{25}^{5}$,
$ \frac{1+\sqrt{5}}{2}\xi_{25}^{5}$;\ \ 
$ \xi_{25}^{11}$,
$ -\xi_{25}^{17}$,
$ -c_{25}^{12}$,
$ \xi_{25}^{21}$,
$ -\frac{1+\sqrt{5}}{2}\xi_{25}^{5}$,
$ \xi_{25}^{9}$,
$ -\xi_{25}^{3}$;\ \ 
$ -1$,
$ \xi_{25}^{19}$,
$ \xi_{25}^{13}$,
$ \xi_{25}^{5}$,
$ c_{25}^{12}$,
$ -\xi_{25}^{9}$;\ \ 
$ -\xi_{25}^{11}$,
$ \xi_{25}^{3}$,
$ \xi_{25}^{5}$,
$ -\xi_{25}^{13}$,
$ \xi_{25}^{21}$;\ \ 
$ -\xi_{25}^{19}$,
$ -\frac{1+\sqrt{5}}{2}\xi_{25}^{5}$,
$ -1$,
$ \xi_{25}^{17}$;\ \ 
$0$,
$ \frac{1+\sqrt{5}}{2}\xi_{25}^{5}$,
$ -\xi_{25}^{5}$;\ \ 
$ \xi_{25}^{21}$,
$ -\xi_{25}^{7}$;\ \ 
$ \xi_{25}^{19}$)

Realization:
Abelian anyon condensation of $SU(2)_{23}$.

\vskip 1ex

\noindent12. $12_{\frac{22}{5},495.9}^{10,127}$ \irep{2073}:\ \ 
$d_i$ = ($1.0$,
$2.618$,
$2.618$,
$4.236$,
$4.236$,
$5.236$,
$5.236$,
$5.854$,
$8.472$,
$8.472$,
$9.472$,
$11.90$) 

\vskip 0.7ex
\hangindent=3em \hangafter=1
$D^2= 495.967 = 
 250+110\sqrt{5} $

\vskip 0.7ex
\hangindent=3em \hangafter=1
$T = ( 0,
\frac{1}{5},
\frac{1}{5},
0,
0,
\frac{2}{5},
\frac{2}{5},
\frac{7}{10},
\frac{4}{5},
\frac{4}{5},
\frac{1}{2},
\frac{1}{5} )
$,

\vskip 0.7ex
\hangindent=3em \hangafter=1
$S$ = ($ 1$,
$ \frac{3+\sqrt{5}}{2}$,
$ \frac{3+\sqrt{5}}{2}$,
$  2+\sqrt{5} $,
$  2+\sqrt{5} $,
$  3+\sqrt{5} $,
$  3+\sqrt{5} $,
$ \frac{5+3\sqrt{5}}{2}$,
$  4+2\sqrt{5} $,
$  4+2\sqrt{5} $,
$  5+2\sqrt{5} $,
$ \frac{11+5\sqrt{5}}{2}$;\ \ 
$  1+4\zeta^{1}_{5}
+4\zeta^{2}_{5}
+2\zeta^{3}_{5}
$,
$  -3-4\zeta^{1}_{5}
-2\zeta^{2}_{5}
$,
$  -3-4\zeta^{1}_{5}
-\zeta^{2}_{5}
+\zeta^{3}_{5}
$,
$  1+4\zeta^{1}_{5}
+5\zeta^{2}_{5}
+3\zeta^{3}_{5}
$,
$ -(3-\sqrt{5} )\zeta_{5}^{2}$,
$ (3+\sqrt{5} )\zeta_{10}^{1}$,
$  5+2\sqrt{5} $,
$ (4+2\sqrt{5} )\zeta_{10}^{3}$,
$ -(4-2\sqrt{5} )\zeta_{5}^{1}$,
$ -\frac{5+3\sqrt{5}}{2}$,
$  2+\sqrt{5} $;\ \ 
$  1+4\zeta^{1}_{5}
+4\zeta^{2}_{5}
+2\zeta^{3}_{5}
$,
$  1+4\zeta^{1}_{5}
+5\zeta^{2}_{5}
+3\zeta^{3}_{5}
$,
$  -3-4\zeta^{1}_{5}
-\zeta^{2}_{5}
+\zeta^{3}_{5}
$,
$ (3+\sqrt{5} )\zeta_{10}^{1}$,
$ -(3-\sqrt{5} )\zeta_{5}^{2}$,
$  5+2\sqrt{5} $,
$ -(4-2\sqrt{5} )\zeta_{5}^{1}$,
$ (4+2\sqrt{5} )\zeta_{10}^{3}$,
$ -\frac{5+3\sqrt{5}}{2}$,
$  2+\sqrt{5} $;\ \ 
$  3+4\zeta^{1}_{5}
+2\zeta^{2}_{5}
$,
$  -1-4\zeta^{1}_{5}
-4\zeta^{2}_{5}
-2\zeta^{3}_{5}
$,
$ -(4-2\sqrt{5} )\zeta_{5}^{1}$,
$ (4+2\sqrt{5} )\zeta_{10}^{3}$,
$ \frac{5+3\sqrt{5}}{2}$,
$ (3+\sqrt{5} )\zeta_{5}^{2}$,
$ -(3-\sqrt{5} )\zeta_{10}^{1}$,
$  5+2\sqrt{5} $,
$ -\frac{3+\sqrt{5}}{2}$;\ \ 
$  3+4\zeta^{1}_{5}
+2\zeta^{2}_{5}
$,
$ (4+2\sqrt{5} )\zeta_{10}^{3}$,
$ -(4-2\sqrt{5} )\zeta_{5}^{1}$,
$ \frac{5+3\sqrt{5}}{2}$,
$ -(3-\sqrt{5} )\zeta_{10}^{1}$,
$ (3+\sqrt{5} )\zeta_{5}^{2}$,
$  5+2\sqrt{5} $,
$ -\frac{3+\sqrt{5}}{2}$;\ \ 
$ -(4-2\sqrt{5} )\zeta_{10}^{1}$,
$ (4+2\sqrt{5} )\zeta_{5}^{2}$,
$0$,
$ -(3-\sqrt{5} )\zeta_{5}^{1}$,
$ (3+\sqrt{5} )\zeta_{10}^{3}$,
$0$,
$  4+2\sqrt{5} $;\ \ 
$ -(4-2\sqrt{5} )\zeta_{10}^{1}$,
$0$,
$ (3+\sqrt{5} )\zeta_{10}^{3}$,
$ -(3-\sqrt{5} )\zeta_{5}^{1}$,
$0$,
$  4+2\sqrt{5} $;\ \ 
$  5+2\sqrt{5} $,
$0$,
$0$,
$ -\frac{5+3\sqrt{5}}{2}$,
$  -5-2\sqrt{5} $;\ \ 
$ -(4-2\sqrt{5} )\zeta_{5}^{2}$,
$ (4+2\sqrt{5} )\zeta_{10}^{1}$,
$0$,
$  -3-\sqrt{5} $;\ \ 
$ -(4-2\sqrt{5} )\zeta_{5}^{2}$,
$0$,
$  -3-\sqrt{5} $;\ \ 
$  -5-2\sqrt{5} $,
$ \frac{5+3\sqrt{5}}{2}$;\ \ 
$ -1$)

Realization:
Abelian anyon condensation of $SU(7)_{3}$.

\vskip 1ex

\noindent13. $12_{\frac{20}{7},940.0}^{21,324}$ \irep{4003}:\ \ 
$d_i$ = ($1.0$,
$2.977$,
$4.888$,
$6.690$,
$6.690$,
$6.690$,
$8.343$,
$9.809$,
$11.56$,
$12.56$,
$12.786$,
$13.232$) 

\vskip 0.7ex
\hangindent=3em \hangafter=1
$D^2= 940.87 = 
 105+147c^{1}_{21}
+189c^{2}_{21}
+105c^{1}_{7}
+126c^{4}_{21}
+126c^{5}_{21}
$

\vskip 0.7ex
\hangindent=3em \hangafter=1
$T = ( 0,
\frac{1}{21},
\frac{1}{7},
\frac{2}{7},
\frac{13}{21},
\frac{13}{21},
\frac{10}{21},
\frac{5}{7},
0,
\frac{1}{3},
\frac{5}{7},
\frac{1}{7} )
$,

\vskip 0.7ex
\hangindent=3em \hangafter=1
$S$ = ($ 1$,
$ \xi_{42}^{3}$,
$ \xi_{42}^{5}$,
$ \xi_{21}^{11}$,
$ \xi_{21}^{11}$,
$ \xi_{21}^{11}$,
$ \xi_{42}^{9}$,
$ \xi_{42}^{11}$,
$ \frac{5+\sqrt{21}}{2}\xi_{21}^{16,2}$,
$ \xi_{42}^{15}$,
$ \xi_{42}^{17}$,
$ \xi_{42}^{19}$;\ \ 
$ \xi_{42}^{9}$,
$ \xi_{42}^{15}$,
$ 2\xi_{21}^{11}$,
$ -\xi_{21}^{11}$,
$ -\xi_{21}^{11}$,
$ \xi_{42}^{15}$,
$ \xi_{42}^{9}$,
$ \xi_{42}^{3}$,
$ -\xi_{42}^{3}$,
$ -\xi_{42}^{9}$,
$ -\xi_{42}^{15}$;\ \ 
$ \xi_{42}^{17}$,
$ \xi_{21}^{11}$,
$ \xi_{21}^{11}$,
$ \xi_{21}^{11}$,
$ -\xi_{42}^{3}$,
$ -\frac{5+\sqrt{21}}{2}\xi_{21}^{16,2}$,
$ -\xi_{42}^{19}$,
$ -\xi_{42}^{9}$,
$ 1$,
$ \xi_{42}^{11}$;\ \ 
$ -\xi_{21}^{11}$,
$ -\xi_{21}^{11}$,
$ -\xi_{21}^{11}$,
$ -2\xi_{21}^{11}$,
$ -\xi_{21}^{11}$,
$ \xi_{21}^{11}$,
$ 2\xi_{21}^{11}$,
$ \xi_{21}^{11}$,
$ -\xi_{21}^{11}$;\ \ 
$ \frac{1-\sqrt{21}}{2}\xi_{21}^{11}$,
$ \frac{1+\sqrt{21}}{2}\xi_{21}^{11}$,
$ \xi_{21}^{11}$,
$ -\xi_{21}^{11}$,
$ \xi_{21}^{11}$,
$ -\xi_{21}^{11}$,
$ \xi_{21}^{11}$,
$ -\xi_{21}^{11}$;\ \ 
$ \frac{1-\sqrt{21}}{2}\xi_{21}^{11}$,
$ \xi_{21}^{11}$,
$ -\xi_{21}^{11}$,
$ \xi_{21}^{11}$,
$ -\xi_{21}^{11}$,
$ \xi_{21}^{11}$,
$ -\xi_{21}^{11}$;\ \ 
$ -\xi_{42}^{3}$,
$ \xi_{42}^{15}$,
$ \xi_{42}^{9}$,
$ -\xi_{42}^{9}$,
$ -\xi_{42}^{15}$,
$ \xi_{42}^{3}$;\ \ 
$ \xi_{42}^{5}$,
$ -\xi_{42}^{17}$,
$ -\xi_{42}^{3}$,
$ \xi_{42}^{19}$,
$ 1$;\ \ 
$ 1$,
$ \xi_{42}^{15}$,
$ -\xi_{42}^{11}$,
$ -\xi_{42}^{5}$;\ \ 
$ -\xi_{42}^{15}$,
$ \xi_{42}^{3}$,
$ \xi_{42}^{9}$;\ \ 
$ \xi_{42}^{5}$,
$ -\frac{5+\sqrt{21}}{2}\xi_{21}^{16,2}$;\ \ 
$ \xi_{42}^{17}$)

Realization:
Abelian anyon condensation of $\overline{SO(20)}_{3}$.

\vskip 1ex

\noindent14. $12_{4,1276.}^{39,406}$ \irep{4852}:\ \ 
$d_i$ = ($1.0$,
$9.908$,
$9.908$,
$9.908$,
$9.908$,
$9.908$,
$9.908$,
$10.908$,
$11.908$,
$11.908$,
$11.908$,
$11.908$) 

\vskip 0.7ex
\hangindent=3em \hangafter=1
$D^2= 1276.274 = 
\frac{1287+351\sqrt{13}}{2}$

\vskip 0.7ex
\hangindent=3em \hangafter=1
$T = ( 0,
\frac{1}{13},
\frac{3}{13},
\frac{4}{13},
\frac{9}{13},
\frac{10}{13},
\frac{12}{13},
0,
\frac{1}{3},
\frac{1}{3},
\frac{2}{3},
\frac{2}{3} )
$,

\vskip 0.7ex
\hangindent=3em \hangafter=1
$S$ = ($ 1$,
$ \frac{9+3\sqrt{13}}{2}$,
$ \frac{9+3\sqrt{13}}{2}$,
$ \frac{9+3\sqrt{13}}{2}$,
$ \frac{9+3\sqrt{13}}{2}$,
$ \frac{9+3\sqrt{13}}{2}$,
$ \frac{9+3\sqrt{13}}{2}$,
$ \frac{11+3\sqrt{13}}{2}$,
$ \frac{13+3\sqrt{13}}{2}$,
$ \frac{13+3\sqrt{13}}{2}$,
$ \frac{13+3\sqrt{13}}{2}$,
$ \frac{13+3\sqrt{13}}{2}$;\ \ 
$ -\frac{9+3\sqrt{13}}{2}c_{13}^{2}$,
$ -\frac{9+3\sqrt{13}}{2}c_{13}^{5}$,
$ -\frac{9+3\sqrt{13}}{2}c_{13}^{4}$,
$ -\frac{9+3\sqrt{13}}{2}c_{13}^{6}$,
$ -\frac{9+3\sqrt{13}}{2}c_{13}^{1}$,
$ -\frac{9+3\sqrt{13}}{2}c_{13}^{3}$,
$ -\frac{9+3\sqrt{13}}{2}$,
$0$,
$0$,
$0$,
$0$;\ \ 
$ -\frac{9+3\sqrt{13}}{2}c_{13}^{6}$,
$ -\frac{9+3\sqrt{13}}{2}c_{13}^{3}$,
$ -\frac{9+3\sqrt{13}}{2}c_{13}^{2}$,
$ -\frac{9+3\sqrt{13}}{2}c_{13}^{4}$,
$ -\frac{9+3\sqrt{13}}{2}c_{13}^{1}$,
$ -\frac{9+3\sqrt{13}}{2}$,
$0$,
$0$,
$0$,
$0$;\ \ 
$ -\frac{9+3\sqrt{13}}{2}c_{13}^{5}$,
$ -\frac{9+3\sqrt{13}}{2}c_{13}^{1}$,
$ -\frac{9+3\sqrt{13}}{2}c_{13}^{2}$,
$ -\frac{9+3\sqrt{13}}{2}c_{13}^{6}$,
$ -\frac{9+3\sqrt{13}}{2}$,
$0$,
$0$,
$0$,
$0$;\ \ 
$ -\frac{9+3\sqrt{13}}{2}c_{13}^{5}$,
$ -\frac{9+3\sqrt{13}}{2}c_{13}^{3}$,
$ -\frac{9+3\sqrt{13}}{2}c_{13}^{4}$,
$ -\frac{9+3\sqrt{13}}{2}$,
$0$,
$0$,
$0$,
$0$;\ \ 
$ -\frac{9+3\sqrt{13}}{2}c_{13}^{6}$,
$ -\frac{9+3\sqrt{13}}{2}c_{13}^{5}$,
$ -\frac{9+3\sqrt{13}}{2}$,
$0$,
$0$,
$0$,
$0$;\ \ 
$ -\frac{9+3\sqrt{13}}{2}c_{13}^{2}$,
$ -\frac{9+3\sqrt{13}}{2}$,
$0$,
$0$,
$0$,
$0$;\ \ 
$ 1$,
$ \frac{13+3\sqrt{13}}{2}$,
$ \frac{13+3\sqrt{13}}{2}$,
$ \frac{13+3\sqrt{13}}{2}$,
$ \frac{13+3\sqrt{13}}{2}$;\ \ 
$ -\frac{13+3\sqrt{13}}{2}$,
$  13+3\sqrt{13} $,
$ -\frac{13+3\sqrt{13}}{2}$,
$ -\frac{13+3\sqrt{13}}{2}$;\ \ 
$ -\frac{13+3\sqrt{13}}{2}$,
$ -\frac{13+3\sqrt{13}}{2}$,
$ -\frac{13+3\sqrt{13}}{2}$;\ \ 
$ -\frac{13+3\sqrt{13}}{2}$,
$  13+3\sqrt{13} $;\ \ 
$ -\frac{13+3\sqrt{13}}{2}$)

Realization: May be related to condensation reductions of categories of the form $\eZ(\cNG(\Z_3\times \Z_3,9))$.
\vskip 1ex

\noindent15. $12_{0,1276.}^{39,560}$ \irep{4851}:\ \ 
$d_i$ = ($1.0$,
$9.908$,
$9.908$,
$9.908$,
$9.908$,
$9.908$,
$9.908$,
$10.908$,
$11.908$,
$11.908$,
$11.908$,
$11.908$) 

\vskip 0.7ex
\hangindent=3em \hangafter=1
$D^2= 1276.274 = 
\frac{1287+351\sqrt{13}}{2}$

\vskip 0.7ex
\hangindent=3em \hangafter=1
$T = ( 0,
\frac{2}{13},
\frac{5}{13},
\frac{6}{13},
\frac{7}{13},
\frac{8}{13},
\frac{11}{13},
0,
0,
0,
\frac{1}{3},
\frac{2}{3} )
$,

\vskip 0.7ex
\hangindent=3em \hangafter=1
$S$ = ($ 1$,
$ \frac{9+3\sqrt{13}}{2}$,
$ \frac{9+3\sqrt{13}}{2}$,
$ \frac{9+3\sqrt{13}}{2}$,
$ \frac{9+3\sqrt{13}}{2}$,
$ \frac{9+3\sqrt{13}}{2}$,
$ \frac{9+3\sqrt{13}}{2}$,
$ \frac{11+3\sqrt{13}}{2}$,
$ \frac{13+3\sqrt{13}}{2}$,
$ \frac{13+3\sqrt{13}}{2}$,
$ \frac{13+3\sqrt{13}}{2}$,
$ \frac{13+3\sqrt{13}}{2}$;\ \ 
$ -\frac{9+3\sqrt{13}}{2}c_{13}^{4}$,
$ -\frac{9+3\sqrt{13}}{2}c_{13}^{1}$,
$ -\frac{9+3\sqrt{13}}{2}c_{13}^{3}$,
$ -\frac{9+3\sqrt{13}}{2}c_{13}^{2}$,
$ -\frac{9+3\sqrt{13}}{2}c_{13}^{5}$,
$ -\frac{9+3\sqrt{13}}{2}c_{13}^{6}$,
$ -\frac{9+3\sqrt{13}}{2}$,
$0$,
$0$,
$0$,
$0$;\ \ 
$ -\frac{9+3\sqrt{13}}{2}c_{13}^{3}$,
$ -\frac{9+3\sqrt{13}}{2}c_{13}^{4}$,
$ -\frac{9+3\sqrt{13}}{2}c_{13}^{6}$,
$ -\frac{9+3\sqrt{13}}{2}c_{13}^{2}$,
$ -\frac{9+3\sqrt{13}}{2}c_{13}^{5}$,
$ -\frac{9+3\sqrt{13}}{2}$,
$0$,
$0$,
$0$,
$0$;\ \ 
$ -\frac{9+3\sqrt{13}}{2}c_{13}^{1}$,
$ -\frac{9+3\sqrt{13}}{2}c_{13}^{5}$,
$ -\frac{9+3\sqrt{13}}{2}c_{13}^{6}$,
$ -\frac{9+3\sqrt{13}}{2}c_{13}^{2}$,
$ -\frac{9+3\sqrt{13}}{2}$,
$0$,
$0$,
$0$,
$0$;\ \ 
$ -\frac{9+3\sqrt{13}}{2}c_{13}^{1}$,
$ -\frac{9+3\sqrt{13}}{2}c_{13}^{4}$,
$ -\frac{9+3\sqrt{13}}{2}c_{13}^{3}$,
$ -\frac{9+3\sqrt{13}}{2}$,
$0$,
$0$,
$0$,
$0$;\ \ 
$ -\frac{9+3\sqrt{13}}{2}c_{13}^{3}$,
$ -\frac{9+3\sqrt{13}}{2}c_{13}^{1}$,
$ -\frac{9+3\sqrt{13}}{2}$,
$0$,
$0$,
$0$,
$0$;\ \ 
$ -\frac{9+3\sqrt{13}}{2}c_{13}^{4}$,
$ -\frac{9+3\sqrt{13}}{2}$,
$0$,
$0$,
$0$,
$0$;\ \ 
$ 1$,
$ \frac{13+3\sqrt{13}}{2}$,
$ \frac{13+3\sqrt{13}}{2}$,
$ \frac{13+3\sqrt{13}}{2}$,
$ \frac{13+3\sqrt{13}}{2}$;\ \ 
$  13+3\sqrt{13} $,
$ -\frac{13+3\sqrt{13}}{2}$,
$ -\frac{13+3\sqrt{13}}{2}$,
$ -\frac{13+3\sqrt{13}}{2}$;\ \ 
$  13+3\sqrt{13} $,
$ -\frac{13+3\sqrt{13}}{2}$,
$ -\frac{13+3\sqrt{13}}{2}$;\ \ 
$ -\frac{13+3\sqrt{13}}{2}$,
$  13+3\sqrt{13} $;\ \ 
$ -\frac{13+3\sqrt{13}}{2}$)

Realization:
Haag$(1)_0$.
\vskip 1ex

\noindent16. $12_{0,1276.}^{117,251}$ \irep{5248}:\ \ 
$d_i$ = ($1.0$,
$9.908$,
$9.908$,
$9.908$,
$9.908$,
$9.908$,
$9.908$,
$10.908$,
$11.908$,
$11.908$,
$11.908$,
$11.908$) 

\vskip 0.7ex
\hangindent=3em \hangafter=1
$D^2= 1276.274 = 
\frac{1287+351\sqrt{13}}{2}$

\vskip 0.7ex
\hangindent=3em \hangafter=1
$T = ( 0,
\frac{2}{13},
\frac{5}{13},
\frac{6}{13},
\frac{7}{13},
\frac{8}{13},
\frac{11}{13},
0,
0,
\frac{1}{9},
\frac{4}{9},
\frac{7}{9} )
$,

\vskip 0.7ex
\hangindent=3em \hangafter=1
$S$ = ($ 1$,
$ \frac{9+3\sqrt{13}}{2}$,
$ \frac{9+3\sqrt{13}}{2}$,
$ \frac{9+3\sqrt{13}}{2}$,
$ \frac{9+3\sqrt{13}}{2}$,
$ \frac{9+3\sqrt{13}}{2}$,
$ \frac{9+3\sqrt{13}}{2}$,
$ \frac{11+3\sqrt{13}}{2}$,
$ \frac{13+3\sqrt{13}}{2}$,
$ \frac{13+3\sqrt{13}}{2}$,
$ \frac{13+3\sqrt{13}}{2}$,
$ \frac{13+3\sqrt{13}}{2}$;\ \ 
$ -\frac{9+3\sqrt{13}}{2}c_{13}^{4}$,
$ -\frac{9+3\sqrt{13}}{2}c_{13}^{1}$,
$ -\frac{9+3\sqrt{13}}{2}c_{13}^{3}$,
$ -\frac{9+3\sqrt{13}}{2}c_{13}^{2}$,
$ -\frac{9+3\sqrt{13}}{2}c_{13}^{5}$,
$ -\frac{9+3\sqrt{13}}{2}c_{13}^{6}$,
$ -\frac{9+3\sqrt{13}}{2}$,
$0$,
$0$,
$0$,
$0$;\ \ 
$ -\frac{9+3\sqrt{13}}{2}c_{13}^{3}$,
$ -\frac{9+3\sqrt{13}}{2}c_{13}^{4}$,
$ -\frac{9+3\sqrt{13}}{2}c_{13}^{6}$,
$ -\frac{9+3\sqrt{13}}{2}c_{13}^{2}$,
$ -\frac{9+3\sqrt{13}}{2}c_{13}^{5}$,
$ -\frac{9+3\sqrt{13}}{2}$,
$0$,
$0$,
$0$,
$0$;\ \ 
$ -\frac{9+3\sqrt{13}}{2}c_{13}^{1}$,
$ -\frac{9+3\sqrt{13}}{2}c_{13}^{5}$,
$ -\frac{9+3\sqrt{13}}{2}c_{13}^{6}$,
$ -\frac{9+3\sqrt{13}}{2}c_{13}^{2}$,
$ -\frac{9+3\sqrt{13}}{2}$,
$0$,
$0$,
$0$,
$0$;\ \ 
$ -\frac{9+3\sqrt{13}}{2}c_{13}^{1}$,
$ -\frac{9+3\sqrt{13}}{2}c_{13}^{4}$,
$ -\frac{9+3\sqrt{13}}{2}c_{13}^{3}$,
$ -\frac{9+3\sqrt{13}}{2}$,
$0$,
$0$,
$0$,
$0$;\ \ 
$ -\frac{9+3\sqrt{13}}{2}c_{13}^{3}$,
$ -\frac{9+3\sqrt{13}}{2}c_{13}^{1}$,
$ -\frac{9+3\sqrt{13}}{2}$,
$0$,
$0$,
$0$,
$0$;\ \ 
$ -\frac{9+3\sqrt{13}}{2}c_{13}^{4}$,
$ -\frac{9+3\sqrt{13}}{2}$,
$0$,
$0$,
$0$,
$0$;\ \ 
$ 1$,
$ \frac{13+3\sqrt{13}}{2}$,
$ \frac{13+3\sqrt{13}}{2}$,
$ \frac{13+3\sqrt{13}}{2}$,
$ \frac{13+3\sqrt{13}}{2}$;\ \ 
$  13+3\sqrt{13} $,
$ -\frac{13+3\sqrt{13}}{2}$,
$ -\frac{13+3\sqrt{13}}{2}$,
$ -\frac{13+3\sqrt{13}}{2}$;\ \ 
$ \frac{13+3\sqrt{13}}{2}c_{9}^{2}$,
$ \frac{13+3\sqrt{13}}{2}c_{9}^{4}$,
$ \frac{13+3\sqrt{13}}{2}c_{9}^{1}$;\ \ 
$ \frac{13+3\sqrt{13}}{2}c_{9}^{1}$,
$ \frac{13+3\sqrt{13}}{2}c_{9}^{2}$;\ \ 
$ \frac{13+3\sqrt{13}}{2}c_{9}^{4}$)

Realization:
Haag$(1)_1$.

\vskip 1ex

\noindent17. $12_{\frac{70}{9},1996.}^{27,115}$ \irep{4453}:\ \ 
$d_i$ = ($1.0$,
$4.932$,
$6.811$,
$7.562$,
$10.270$,
$11.585$,
$11.802$,
$14.369$,
$15.369$,
$16.734$,
$17.82$,
$21.773$) 

\vskip 0.7ex
\hangindent=3em \hangafter=1
$D^2= 1996.556 = 
 243+243c^{1}_{27}
+243c^{2}_{27}
+243c^{1}_{9}
+243c^{4}_{27}
+216c^{5}_{27}
+162c^{2}_{9}
+108c^{7}_{27}
+54c^{8}_{27}
$

\vskip 0.7ex
\hangindent=3em \hangafter=1
$T = ( 0,
\frac{2}{9},
\frac{4}{9},
\frac{23}{27},
\frac{1}{9},
\frac{14}{27},
\frac{5}{9},
\frac{2}{3},
\frac{1}{3},
\frac{8}{9},
\frac{7}{9},
\frac{5}{27} )
$,

\vskip 0.7ex
\hangindent=3em \hangafter=1
$S$ = ($ 1$,
$ \xi_{54}^{5}$,
$ \xi_{54}^{7}$,
$  1+c^{1}_{27}
+c^{2}_{27}
+c^{1}_{9}
+c^{4}_{27}
+c^{5}_{27}
+c^{7}_{27}
+c^{8}_{27}
$,
$ \xi_{54}^{11}$,
$  2+c^{1}_{27}
+c^{2}_{27}
+c^{1}_{9}
+2c^{4}_{27}
+2c^{5}_{27}
+c^{2}_{9}
$,
$ \xi_{54}^{13}$,
$ \xi_{54}^{17}$,
$ \xi_{54}^{19}$,
$ \xi_{54}^{23}$,
$ \xi_{54}^{25}$,
$  2+3c^{1}_{27}
+3c^{2}_{27}
+3c^{1}_{9}
+2c^{4}_{27}
+2c^{5}_{27}
+2c^{2}_{9}
+c^{7}_{27}
+c^{8}_{27}
$;\ \ 
$ \xi_{54}^{25}$,
$ \xi_{54}^{19}$,
$  -2-c^{1}_{27}
-c^{2}_{27}
-c^{1}_{9}
-2  c^{4}_{27}
-2  c^{5}_{27}
-c^{2}_{9}
$,
$ -1$,
$  2+3c^{1}_{27}
+3c^{2}_{27}
+3c^{1}_{9}
+2c^{4}_{27}
+2c^{5}_{27}
+2c^{2}_{9}
+c^{7}_{27}
+c^{8}_{27}
$,
$ -\xi_{54}^{11}$,
$ -\xi_{54}^{23}$,
$ -\xi_{54}^{13}$,
$ \xi_{54}^{7}$,
$ \xi_{54}^{17}$,
$  -1-c^{1}_{27}
-c^{2}_{27}
-c^{1}_{9}
-c^{4}_{27}
-c^{5}_{27}
-c^{7}_{27}
-c^{8}_{27}
$;\ \ 
$ \xi_{54}^{5}$,
$  2+3c^{1}_{27}
+3c^{2}_{27}
+3c^{1}_{9}
+2c^{4}_{27}
+2c^{5}_{27}
+2c^{2}_{9}
+c^{7}_{27}
+c^{8}_{27}
$,
$ -\xi_{54}^{23}$,
$  1+c^{1}_{27}
+c^{2}_{27}
+c^{1}_{9}
+c^{4}_{27}
+c^{5}_{27}
+c^{7}_{27}
+c^{8}_{27}
$,
$ -\xi_{54}^{17}$,
$ \xi_{54}^{11}$,
$ \xi_{54}^{25}$,
$ 1$,
$ -\xi_{54}^{13}$,
$  -2-c^{1}_{27}
-c^{2}_{27}
-c^{1}_{9}
-2  c^{4}_{27}
-2  c^{5}_{27}
-c^{2}_{9}
$;\ \ 
$0$,
$  2+3c^{1}_{27}
+3c^{2}_{27}
+3c^{1}_{9}
+2c^{4}_{27}
+2c^{5}_{27}
+2c^{2}_{9}
+c^{7}_{27}
+c^{8}_{27}
$,
$0$,
$  -2-c^{1}_{27}
-c^{2}_{27}
-c^{1}_{9}
-2  c^{4}_{27}
-2  c^{5}_{27}
-c^{2}_{9}
$,
$  1+c^{1}_{27}
+c^{2}_{27}
+c^{1}_{9}
+c^{4}_{27}
+c^{5}_{27}
+c^{7}_{27}
+c^{8}_{27}
$,
$  -1-c^{1}_{27}
-c^{2}_{27}
-c^{1}_{9}
-c^{4}_{27}
-c^{5}_{27}
-c^{7}_{27}
-c^{8}_{27}
$,
$  2+c^{1}_{27}
+c^{2}_{27}
+c^{1}_{9}
+2c^{4}_{27}
+2c^{5}_{27}
+c^{2}_{9}
$,
$  -2-3  c^{1}_{27}
-3  c^{2}_{27}
-3  c^{1}_{9}
-2  c^{4}_{27}
-2  c^{5}_{27}
-2  c^{2}_{9}
-c^{7}_{27}
-c^{8}_{27}
$,
$0$;\ \ 
$ \xi_{54}^{13}$,
$  1+c^{1}_{27}
+c^{2}_{27}
+c^{1}_{9}
+c^{4}_{27}
+c^{5}_{27}
+c^{7}_{27}
+c^{8}_{27}
$,
$ \xi_{54}^{19}$,
$ -\xi_{54}^{25}$,
$ -\xi_{54}^{7}$,
$ \xi_{54}^{17}$,
$ -\xi_{54}^{5}$,
$  -2-c^{1}_{27}
-c^{2}_{27}
-c^{1}_{9}
-2  c^{4}_{27}
-2  c^{5}_{27}
-c^{2}_{9}
$;\ \ 
$0$,
$  2+3c^{1}_{27}
+3c^{2}_{27}
+3c^{1}_{9}
+2c^{4}_{27}
+2c^{5}_{27}
+2c^{2}_{9}
+c^{7}_{27}
+c^{8}_{27}
$,
$  2+c^{1}_{27}
+c^{2}_{27}
+c^{1}_{9}
+2c^{4}_{27}
+2c^{5}_{27}
+c^{2}_{9}
$,
$  -2-c^{1}_{27}
-c^{2}_{27}
-c^{1}_{9}
-2  c^{4}_{27}
-2  c^{5}_{27}
-c^{2}_{9}
$,
$  -2-3  c^{1}_{27}
-3  c^{2}_{27}
-3  c^{1}_{9}
-2  c^{4}_{27}
-2  c^{5}_{27}
-2  c^{2}_{9}
-c^{7}_{27}
-c^{8}_{27}
$,
$  -1-c^{1}_{27}
-c^{2}_{27}
-c^{1}_{9}
-c^{4}_{27}
-c^{5}_{27}
-c^{7}_{27}
-c^{8}_{27}
$,
$0$;\ \ 
$ -\xi_{54}^{7}$,
$ \xi_{54}^{5}$,
$ \xi_{54}^{23}$,
$ -\xi_{54}^{25}$,
$ 1$,
$  -1-c^{1}_{27}
-c^{2}_{27}
-c^{1}_{9}
-c^{4}_{27}
-c^{5}_{27}
-c^{7}_{27}
-c^{8}_{27}
$;\ \ 
$ -\xi_{54}^{19}$,
$ -1$,
$ -\xi_{54}^{13}$,
$ -\xi_{54}^{7}$,
$  2+3c^{1}_{27}
+3c^{2}_{27}
+3c^{1}_{9}
+2c^{4}_{27}
+2c^{5}_{27}
+2c^{2}_{9}
+c^{7}_{27}
+c^{8}_{27}
$;\ \ 
$ \xi_{54}^{17}$,
$ \xi_{54}^{5}$,
$ \xi_{54}^{11}$,
$  -2-3  c^{1}_{27}
-3  c^{2}_{27}
-3  c^{1}_{9}
-2  c^{4}_{27}
-2  c^{5}_{27}
-2  c^{2}_{9}
-c^{7}_{27}
-c^{8}_{27}
$;\ \ 
$ -\xi_{54}^{11}$,
$ \xi_{54}^{19}$,
$  1+c^{1}_{27}
+c^{2}_{27}
+c^{1}_{9}
+c^{4}_{27}
+c^{5}_{27}
+c^{7}_{27}
+c^{8}_{27}
$;\ \ 
$ -\xi_{54}^{23}$,
$  2+c^{1}_{27}
+c^{2}_{27}
+c^{1}_{9}
+2c^{4}_{27}
+2c^{5}_{27}
+c^{2}_{9}
$;\ \ 
$0$)

Realization:
$G(2)_5$.

\vskip 1ex

\noindent18. $12_{0,3926.}^{21,464}$ \irep{4007}:\ \ 
$d_i$ = ($1.0$,
$12.146$,
$12.146$,
$13.146$,
$14.146$,
$15.146$,
$15.146$,
$20.887$,
$20.887$,
$20.887$,
$27.293$,
$27.293$) 

\vskip 0.7ex
\hangindent=3em \hangafter=1
$D^2= 3926.660 = 
 2142+1701c^{1}_{7}
+756c^{2}_{7}
$

\vskip 0.7ex
\hangindent=3em \hangafter=1
$T = ( 0,
\frac{3}{7},
\frac{4}{7},
0,
0,
\frac{1}{7},
\frac{6}{7},
0,
\frac{1}{3},
\frac{2}{3},
\frac{2}{7},
\frac{5}{7} )
$,

\vskip 0.7ex
\hangindent=3em \hangafter=1
$S$ = ($ 1$,
$ 3\xi_{14}^{5}$,
$ 3\xi_{14}^{5}$,
$  7+6c^{1}_{7}
+3c^{2}_{7}
$,
$  8+6c^{1}_{7}
+3c^{2}_{7}
$,
$  9+6c^{1}_{7}
+3c^{2}_{7}
$,
$  9+6c^{1}_{7}
+3c^{2}_{7}
$,
$  11+9c^{1}_{7}
+3c^{2}_{7}
$,
$  11+9c^{1}_{7}
+3c^{2}_{7}
$,
$  11+9c^{1}_{7}
+3c^{2}_{7}
$,
$  15+12c^{1}_{7}
+6c^{2}_{7}
$,
$  15+12c^{1}_{7}
+6c^{2}_{7}
$;\ \ 
$  15+12c^{1}_{7}
+6c^{2}_{7}
$,
$ 3\xi_{7}^{3}$,
$  -15-12  c^{1}_{7}
-6  c^{2}_{7}
$,
$  -9-6  c^{1}_{7}
-3  c^{2}_{7}
$,
$  18+15c^{1}_{7}
+6c^{2}_{7}
$,
$ -3\xi_{14}^{5}$,
$0$,
$0$,
$0$,
$  -12-9  c^{1}_{7}
-3  c^{2}_{7}
$,
$  9+6c^{1}_{7}
+3c^{2}_{7}
$;\ \ 
$  15+12c^{1}_{7}
+6c^{2}_{7}
$,
$  -15-12  c^{1}_{7}
-6  c^{2}_{7}
$,
$  -9-6  c^{1}_{7}
-3  c^{2}_{7}
$,
$ -3\xi_{14}^{5}$,
$  18+15c^{1}_{7}
+6c^{2}_{7}
$,
$0$,
$0$,
$0$,
$  9+6c^{1}_{7}
+3c^{2}_{7}
$,
$  -12-9  c^{1}_{7}
-3  c^{2}_{7}
$;\ \ 
$  -8-6  c^{1}_{7}
-3  c^{2}_{7}
$,
$ -1$,
$ 3\xi_{14}^{5}$,
$ 3\xi_{14}^{5}$,
$  11+9c^{1}_{7}
+3c^{2}_{7}
$,
$  11+9c^{1}_{7}
+3c^{2}_{7}
$,
$  11+9c^{1}_{7}
+3c^{2}_{7}
$,
$  -9-6  c^{1}_{7}
-3  c^{2}_{7}
$,
$  -9-6  c^{1}_{7}
-3  c^{2}_{7}
$;\ \ 
$  7+6c^{1}_{7}
+3c^{2}_{7}
$,
$  15+12c^{1}_{7}
+6c^{2}_{7}
$,
$  15+12c^{1}_{7}
+6c^{2}_{7}
$,
$  -11-9  c^{1}_{7}
-3  c^{2}_{7}
$,
$  -11-9  c^{1}_{7}
-3  c^{2}_{7}
$,
$  -11-9  c^{1}_{7}
-3  c^{2}_{7}
$,
$ 3\xi_{14}^{5}$,
$ 3\xi_{14}^{5}$;\ \ 
$  -9-6  c^{1}_{7}
-3  c^{2}_{7}
$,
$  12+9c^{1}_{7}
+3c^{2}_{7}
$,
$0$,
$0$,
$0$,
$  -15-12  c^{1}_{7}
-6  c^{2}_{7}
$,
$ -3\xi_{7}^{3}$;\ \ 
$  -9-6  c^{1}_{7}
-3  c^{2}_{7}
$,
$0$,
$0$,
$0$,
$ -3\xi_{7}^{3}$,
$  -15-12  c^{1}_{7}
-6  c^{2}_{7}
$;\ \ 
$  22+18c^{1}_{7}
+6c^{2}_{7}
$,
$  -11-9  c^{1}_{7}
-3  c^{2}_{7}
$,
$  -11-9  c^{1}_{7}
-3  c^{2}_{7}
$,
$0$,
$0$;\ \ 
$  -11-9  c^{1}_{7}
-3  c^{2}_{7}
$,
$  22+18c^{1}_{7}
+6c^{2}_{7}
$,
$0$,
$0$;\ \ 
$  -11-9  c^{1}_{7}
-3  c^{2}_{7}
$,
$0$,
$0$;\ \ 
$ -3\xi_{14}^{5}$,
$  18+15c^{1}_{7}
+6c^{2}_{7}
$;\ \ 
$ -3\xi_{14}^{5}$)

Realization: unknown

\vskip 1ex

\noindent19. $12_{0,3926.}^{63,774}$ \irep{5124}:\ \ 
$d_i$ = ($1.0$,
$12.146$,
$12.146$,
$13.146$,
$14.146$,
$15.146$,
$15.146$,
$20.887$,
$20.887$,
$20.887$,
$27.293$,
$27.293$) 

\vskip 0.7ex
\hangindent=3em \hangafter=1
$D^2= 3926.660 = 
 2142+1701c^{1}_{7}
+756c^{2}_{7}
$

\vskip 0.7ex
\hangindent=3em \hangafter=1
$T = ( 0,
\frac{3}{7},
\frac{4}{7},
0,
0,
\frac{1}{7},
\frac{6}{7},
\frac{1}{9},
\frac{4}{9},
\frac{7}{9},
\frac{2}{7},
\frac{5}{7} )
$,

\vskip 0.7ex
\hangindent=3em \hangafter=1
$S$ = ($ 1$,
$ 3\xi_{14}^{5}$,
$ 3\xi_{14}^{5}$,
$  7+6c^{1}_{7}
+3c^{2}_{7}
$,
$  8+6c^{1}_{7}
+3c^{2}_{7}
$,
$  9+6c^{1}_{7}
+3c^{2}_{7}
$,
$  9+6c^{1}_{7}
+3c^{2}_{7}
$,
$  11+9c^{1}_{7}
+3c^{2}_{7}
$,
$  11+9c^{1}_{7}
+3c^{2}_{7}
$,
$  11+9c^{1}_{7}
+3c^{2}_{7}
$,
$  15+12c^{1}_{7}
+6c^{2}_{7}
$,
$  15+12c^{1}_{7}
+6c^{2}_{7}
$;\ \ 
$  15+12c^{1}_{7}
+6c^{2}_{7}
$,
$ 3\xi_{7}^{3}$,
$  -15-12  c^{1}_{7}
-6  c^{2}_{7}
$,
$  -9-6  c^{1}_{7}
-3  c^{2}_{7}
$,
$  18+15c^{1}_{7}
+6c^{2}_{7}
$,
$ -3\xi_{14}^{5}$,
$0$,
$0$,
$0$,
$  -12-9  c^{1}_{7}
-3  c^{2}_{7}
$,
$  9+6c^{1}_{7}
+3c^{2}_{7}
$;\ \ 
$  15+12c^{1}_{7}
+6c^{2}_{7}
$,
$  -15-12  c^{1}_{7}
-6  c^{2}_{7}
$,
$  -9-6  c^{1}_{7}
-3  c^{2}_{7}
$,
$ -3\xi_{14}^{5}$,
$  18+15c^{1}_{7}
+6c^{2}_{7}
$,
$0$,
$0$,
$0$,
$  9+6c^{1}_{7}
+3c^{2}_{7}
$,
$  -12-9  c^{1}_{7}
-3  c^{2}_{7}
$;\ \ 
$  -8-6  c^{1}_{7}
-3  c^{2}_{7}
$,
$ -1$,
$ 3\xi_{14}^{5}$,
$ 3\xi_{14}^{5}$,
$  11+9c^{1}_{7}
+3c^{2}_{7}
$,
$  11+9c^{1}_{7}
+3c^{2}_{7}
$,
$  11+9c^{1}_{7}
+3c^{2}_{7}
$,
$  -9-6  c^{1}_{7}
-3  c^{2}_{7}
$,
$  -9-6  c^{1}_{7}
-3  c^{2}_{7}
$;\ \ 
$  7+6c^{1}_{7}
+3c^{2}_{7}
$,
$  15+12c^{1}_{7}
+6c^{2}_{7}
$,
$  15+12c^{1}_{7}
+6c^{2}_{7}
$,
$  -11-9  c^{1}_{7}
-3  c^{2}_{7}
$,
$  -11-9  c^{1}_{7}
-3  c^{2}_{7}
$,
$  -11-9  c^{1}_{7}
-3  c^{2}_{7}
$,
$ 3\xi_{14}^{5}$,
$ 3\xi_{14}^{5}$;\ \ 
$  -9-6  c^{1}_{7}
-3  c^{2}_{7}
$,
$  12+9c^{1}_{7}
+3c^{2}_{7}
$,
$0$,
$0$,
$0$,
$  -15-12  c^{1}_{7}
-6  c^{2}_{7}
$,
$ -3\xi_{7}^{3}$;\ \ 
$  -9-6  c^{1}_{7}
-3  c^{2}_{7}
$,
$0$,
$0$,
$0$,
$ -3\xi_{7}^{3}$,
$  -15-12  c^{1}_{7}
-6  c^{2}_{7}
$;\ \ 
$  9c^{1}_{63}
-9  c^{2}_{63}
+3c^{4}_{63}
-9  c^{1}_{9}
+9c^{8}_{63}
+6c^{10}_{63}
-3  c^{11}_{63}
+2c^{2}_{9}
-9  c^{16}_{63}
+9c^{17}_{63}
$,
$  -9  c^{1}_{63}
-2  c^{1}_{9}
-9  c^{8}_{63}
-6  c^{10}_{63}
-2  c^{2}_{9}
-6  c^{17}_{63}
$,
$  9c^{2}_{63}
-3  c^{4}_{63}
+11c^{1}_{9}
+3c^{11}_{63}
+9c^{16}_{63}
-3  c^{17}_{63}
$,
$0$,
$0$;\ \ 
$  9c^{2}_{63}
-3  c^{4}_{63}
+11c^{1}_{9}
+3c^{11}_{63}
+9c^{16}_{63}
-3  c^{17}_{63}
$,
$  9c^{1}_{63}
-9  c^{2}_{63}
+3c^{4}_{63}
-9  c^{1}_{9}
+9c^{8}_{63}
+6c^{10}_{63}
-3  c^{11}_{63}
+2c^{2}_{9}
-9  c^{16}_{63}
+9c^{17}_{63}
$,
$0$,
$0$;\ \ 
$  -9  c^{1}_{63}
-2  c^{1}_{9}
-9  c^{8}_{63}
-6  c^{10}_{63}
-2  c^{2}_{9}
-6  c^{17}_{63}
$,
$0$,
$0$;\ \ 
$ -3\xi_{14}^{5}$,
$  18+15c^{1}_{7}
+6c^{2}_{7}
$;\ \ 
$ -3\xi_{14}^{5}$)

Realization: unknown

\vskip 1ex 

}

\section{Realizations of exotic modular data}
\label{realizations}

We find some potential modular data that cannot be realized by modular tensor
categories from Kac-Moody algebra or twisted quantum doubles, nor from their
Abelian anyon condensations \cite{LW170107820}, their Galois conjugation, and
their change of spherical structure.  We refer to those potential modular data
as exotic potential modular data.  

All the exotic
potential modular data that we found are listed below:
\begin{enumerate}

\item
 two Galois orbits of rank-8 modular data 
represented by $8_{4,36.}^{6,102}$ and $8_{4,36.}^{12,972}$ 
with $D^2=36$, 

\item
 three Galois orbits of rank-8 modular data with
$D^2\approx 308.434$, 

\item
 one Galois orbit of rank-10 modular data 
represented by $ 10^{20,676}_{0,1435.}$
with $D^2\approx 1435.541$, 

\item
 two Galois orbits of rank-11 modular data with
$D^2\approx 1337.107$ and $D^2\approx 1964.590$, 

\item
four Galois orbits of rank-12 modular data with
$D^2= 144$ (which contain six unitary modular data).

\item
one Galois orbit of rank-12 modular data 
represented by $ 12^{39,406}_{4,1276.}$
with $D^2\approx 1276.274$.

\item
two Galois orbits of rank-12 modular data with
$D^2\approx 3926.660$.

\end{enumerate}

In the following, we will
discuss the realizations of those exotic potential modular data, to see if they
are actually exotic modular data that can realized by some modular tensor
categories.

\subsection{Near-Group fusion categories and their centers}\label{ss: near group}

The main references for near-group categories and their centers are
\cite{IzumiII,EvansGannon}.  We reproduce their results for the reader's
convenience.

Let $G$ be a finite group of order $|G|$ and $m$ a non-negative integer. A
\textbf{near-group category} of type $G+m$ is a rank $|G|+1$ fusion category
with simple objects labeled by elements $g\in G$ and an additional simple
object $\rho$ such that the fusion rules are given by the group operation in
$G$, $g\rho=\rho g=\rho$ for all $g\in G$, and $ \rho\otimes \rho =
m\rho+\sum_{g\in G}g$.  While the near-group fusion rule of type $G+m$ is
well-defined, not all are associated with fusion categories.   {We denote the (possibly empty) class of fusion categories of near-group type $G+m$ by $\cNG(G,m)$.}
In the literature, one typically finds results in the \emph{unitary} setting,
so we will focus on this situation.  It is known \cite[Theorem 2]{EvansGannon}
{ that, in order for $\cNG(G,m)$ to  {contain} an unitary fusion category with $H^2(G,\R/\Z)=0$, the only
possible values of $m$ are $|G|-1$ or $k|G|$ for some non-negative integer
$k$.}  The Tambara-Yamagami categories are the near-group fusion categories
with $m = 0$, which are classified in \cite{TY98}.  

To construct/classify unitary near-group fusion categories, one must solve a
system of non-linear equations \cite{IzumiII}, which is computationally
strenuous.  A precise statement for $G$ abelian can be found in \cite[Theorem
5.3]{IzumiII} and \cite[Corollary 5]{EvansGannon}.  In the abelian case one
first chooses a non-degenerate symmetric bicharacter $\langle\;,\;\rangle$ on
$G$, which facilitates the solution method found in \cite{IzumiII,EvansGannon}.
Realizations have been found for near-group fusion rules of the following
types, for example:
\begin{enumerate}
 \item $A+0$ for all abelian groups $A$ \cite{TY98}
     \item $A+|A|$ for $A$ abelian and $|A|\leq 13$ \cite{EvansGannon}
     \item $\Z_N+N$ for $N\leq 30$, except $N=19,29.$ \cite{Budthesis}.
     \item $\Z_3+6$ see \cite{Izumi15} who attributes this to Liu and Snyder, \cite{liupreprint}.
     \item $G+2^k$ with $G$ and extra-special $2$-group of order $|G|=2^{2k+1}$.
     \cite[Therem 6.1]{Izumi15}.
 \end{enumerate}

 If it is true that there are finitely many fusion categories of each rank, as is true in the braided setting \cite{JMNR}, there should be no near-group categories of type $G+m$ for $m$ above some bound.  For example there are no unitary near-group categories of the following types:
 \begin{enumerate}
     \item $\Z_2+m$ $m\geq 3$ \cite{Ostrikrank3},
     \item $(\Z_2)^k+m$ with $k\geq 3$ and $m\neq 0$ \cite{Schop22},
     \item $\Z_3+m$ $m\geq 7$ \cite{hannahlarson},
     \item $G+m$ where $G$ is non-abelian, except when $G$ is an extra-special $2$-group as described above \cite{Izumi15}.
 \end{enumerate}

The modular data of the center $\mathcal{Z}(\eC)$ of a near group-category $\eC$ can sometimes be obtained, but it is computationally difficult (employing tube algebra or Cuntz algebra methods, \cite{IzumiII,Izumi15}).

If $\eC$ is a unitary near-group fusion category of type $A+|A|$ for an abelian group $A$, formulae for the modular data are found in \cite{EvansGannon}, see \cite{RSZ} for some explicit computations.  In this case we have the following facts about $\eD:=\mathcal{Z}(\eC)$, where $\eC$ is of type $A+N$ where $|A|=N$:
\begin{enumerate}
    \item The rank of $\eD$ is $N(N+3)$
    \item The dimensions of simple objects in $\eD$ are $1,d,d+1,d+2$ where $d:=\frac{N+\sqrt{N^2+4N}}{2}$.
    \item There are $N$ invertible objects,
    \item $N$ simple objects of dimension $d+1$,
    \item $\binom{N}{2}$ simple objects of dimension $d+2$, and 
    \item $N(N+3)/2$ simple objects of dimension $d$.
    \item $\dim(\eD)=(N+d^2)^2$,
    \item the pointed part of $\eD$ is a ribbon fusion category of the form $\eC(A,q)$ where $q$ is a quadratic form given by $q(a)=\langle a,a\rangle$ with $\langle\/,\/\rangle$ a non-degenerate symmetric bicharacter on $A$.  In particular the $S$-matrix for the pointed subcategory has entries $\langle a,b\rangle^{-2}$ and the $T$-matrix has entries $\delta_{a,b}\langle a,b\rangle$.  
\end{enumerate}
The last point completely characterizes the pointed subcategory of $\eD$,
without the need of further computation.  Using the results of
\cite{EvansGannon} it can be shown (cf. \cite[Section 2.4]{RSZ}) that if $\eC$
is a near-group unitary fusion category of type $A+|A|$ for $|A|=N$ odd, then
the rank $N$ pointed subcategory $\eC(A,q)$ of $\eD$ is modular, and thus
$\eD\cong\eC(A,q)\boxtimes \eF$ where $\eF$ is a modular category. Moreover,
$\eF$ has rank $(N+3)$ and dimension $D^2 = {\frac {N \left( N+4 \right) }{2}
\left( N+\sqrt {N \left( N+4 \right) }+2 \right) }$$ $, with no non-trivial
invertible objects.  When $|A|$ is even, one often finds that $\eD$ contains invertible bosons, which
can be condensed.  

In the following we address the question of realizability of several of our modular data via centers of near-group categories.  In some cases this leads to
 definitive constructions, while in other we provide strong evidence of realizability.
\begin{enumerate}

    \item

One obtains rank $8$ modular categories as the non-pointed
Deligne factor of $\mathcal{Z}(\eC)$ for $\eC$ a near-group category of type
$\Z_3\times \Z_3+0$.  These provide \textbf{definitive} realizations for the two modular data of
the form $8_{4,36.}^{6,102}$ and $8_{4,36.}^{12,972}$.  

\item One also obtains rank
$8$ modular categories as the non-pointed Deligne factor of $\mathcal{Z}(\eC)$
for $\eC$ a near-group category of type $\Z_5+5$.  These provide \textbf{definitive} realizations
for the three modular data of the form $8_{y,308.4}^{15,x}$ (note that there
are exactly $3$ near-group categories of type $\Z_5+5$ in \cite{EvansGannon}).

\item
The case of $\eC$ being a near-group category of type $\Z_3+6$ is relevant to
us.  It exists by \cite{Izumi15}, but we do not know many details about
$\eD=\mathcal{Z}(\eC)$ other than its dimension: $144\cdot(7+4\sqrt{3})$. One
can show directly that $\eD$ contains a modular pointed subcategory with fusion
rules like $\Z_3$, so that $\eD\cong \eF\boxtimes \eC(\Z_3,q)$ where $\eF$ is a
modular category of dimension $48\cdot (7+4\sqrt{3})\approx 668.553$, and hence
a \textbf{candidate} for a realization of the potential modular data
$9_{6,668.5}^{12,567}$.  For our purposes it is sufficient to show that $\eD$
has rank $27$, which is a special case of \cite[Conj. 7.6]{IzumiGrossman} see
Example 7.1.7(1) in the same reference.  We have verified that our modular data
coincides with the conjectural data provided in \cite{IzumiGrossman}.  On the
other hand, it is shown in \cite{Edie-MichellCAMS} the existence of a category
with the same rank and dimensions as $9_{6,668.5}^{12,567}$ by condensing the
$\Z_3$-bosons in $SU(3)_9$, providing a \textbf{definitive} realization.  We
find that $SU(3)_{-9}\boxtimes 9_{6,668.5}^{12,567}$ has 9 potential Lagrangian
condensable algebras, supporting the above realization.

\item
Similarly as the rank $8$ case above, one obtains rank $10$ modular categories
as the non-pointed Deligne factor of $\mathcal{Z}(\eC)$ for $\eC$ a near-group
category of type $\Z_7+7$.  These provide \textbf{definitive} realizations for
the modular data $10_{6,684.3}^{77,298}$ (there is one fusion category of type
$\Z_7+7$ in \cite{EvansGannon} up to complex conjugation).

\item
Consider $\eD=\mathcal{Z}(\eC)$ where $\eC$ is a near-group fusion category of
type $A+|A|$ for $A=\Z_4\times\Z_4$ with symmetric bicharacter $\langle
(a,b),(c,d)\rangle=(\zeta_4)^{ac-bd}$ where $\zeta_4=e^{\pi i/2}$.  $\eD$ has
dimension $5\cdot2^{10}\cdot(9+4\sqrt{5})$ and rank $304$.  We find a Tannakian
pointed subcategory $\mathcal{P}\subset \eD$ with fusion rules like $\Z_4\times
\Z_2$, generated by  $(1,1)$ and $(2,0)$.  Note that
$\langle(1,0),(1,1)\rangle^{-2}=(\zeta_4)^{-2}=-1$, so that $(1,0)$ is not in
the centralizer of $\mathcal{P}$.  Thus the condensation
$[\eD_{\Z_4\times\Z_4}]_0$ has no invertible objects and has dimension
$\frac{5\cdot2^{10}\cdot(9+4\sqrt{5})}{8^2}\approx 1435.541$.  This method was
used to confirm that we have a realization of modular data
$10_{0,1435.}^{20,676}$ in \cite{YuZhang}, providing a \textbf{definitive}
construction for this modular data.  It is also found to be the center of a
rank $4$ fusion category, see \cite{EMIP} On the other hand, the modular data
$10_{4,1435.}^{10,168}$ looks very similar, but is realized by quantum groups:
a $\Z_5$ boson condensation of $SU(5)_5$.

\item 
Another example is the following: Consider a fusion category $\eC$ of
near-group type $\Z_{12}+12$ (which exists, by \cite{EvansGannon}). In this
case we find that $\eD=\mathcal{Z}(\eC)$ a modular category of rank $180$ and
dimension $48^2\cdot(7+4\sqrt{3})$. By example 2.5 of \cite{RSZ} We see that
$\eD\cong \eF\boxtimes\eC(\Z_3,q)$ for some quadratic form $q$ on $\Z_3$.  Now
$\eF$ contains an invertible boson $b$ with $b\otimes b\cong \one$ so that we
obtain a modular category $[\eF_{\Z_2}]_0$ by condensation.  Moreover,
computations as in \cite{RSZ} show that $[\eF_{\Z_2}]_0$ factors as
$\mathcal{B}\otimes \eC(\Z_2,q^\prime)$ where $\mathcal{B}$ is a modular
category of dimension $96\cdot(7+4\sqrt{3})\approx 1337.107510$.  A more
detailed analysis is necessary, but this is a strong evidence that this
provides a \textbf{candidate} realization of potential modular data
$11_{3,1337.}^{48,634}$.
\item The potential modular data $12^{39,406}_{4,1276.}$ is similar to that of $12^{39,560}_{0,1276.}$ which is realized by means of the center of the even part of the Haagerup subfactor.  As the latter may also be obtained from the center of the near-group category $\cNG(\Z_3\times\Z_3,9)$, this suggests a similar construction for this potential modular data, and therefore is a \textbf{candidate} realization.  Indeed, in \cite{EvansGannon} there are 2 distinct near-group categories associated with $\Z_3\times \Z_3$.
\end{enumerate}

\subsection{Realizations via minimal modular extensions}

The theory of minimal modular extensions could explain the multiplicity of a
collection of similar but not identical modular tensor categories. Let $\eC$ be
a modular tensor category which contains a fusion subcategory $\eF$ such that
the Frobenius-Perron dimension of $\eC$ is the product of $\FPdim(\eF)$ and the
FP-dimension of the M\"uger center $\eE$ of $\eF$. In this case, $\eC$ is
called a minimal modular extension of $\eF$ (cf. \cite{Mu2003}). Note that
$\eE$ is a symmetric fusion category, and so $\eE$ is equivalent as braided
fusion categories to the representation category of a finite group $G$, which
is  uniquely determined by $\eE$. According to the theorem of \cite{LKW2017},
the set of equivalent classes of the minimal modular extensions of $\eF$ is a
torsor of the 3rd cohomology group $H^3(G, U(1))$ of $G$. In fact, there is a
faithful transitive action of $H^3(G, U(1))$ on the set of minimal modular
extensions of $\eF$ by condensing $\eE$ of the Deligne product of $\eC$ and the
center $Z(Vec_G^\omega)$ of the pointed category $Vec_G^\omega$ for any
cohomology class $\omega$ in $H^3(G, U(1))$.  In particular, if a minimal
modular extension of $\eF$ shows up in the list, one should find exactly the
number of classes in the group  $H^3(G, U(1))$. 
 
 It is quite routine to check a given modular data is a minimal extension of
one of its fusion subcategory. We consider the six rank 12 modular data of
FPdim $D^2 = 144$ in our list as examples. Let $\eC$ be an MTC whose modular
data is any of these rank 12 data. One can see immediately that there are 3
bosons but one of them is nonabelian. It is can check from the fusion rules
computed from the $S$-matrix that these three objects form a Tannakian
subcategory $\eE$ of dimension 6. It has to be the representation category of
the symmetric group $S_3$. Using the modular data again, we see that the simple
objects of dimension 3 together with the simple object of $\eE$ form a fusion
subcategory $\eF$ of  dimension 24. By considering the twists of the simple
objects of $\eF$, we can confirm that $\eE$ is the M\"ger center of $\eF$ and
so $\eC$ is a minimal modular extension of $\eF$. Since $H^3(S_3, U(1))$ is a
cyclic group of order 6, there are exactly 6 minimal modular extensions of
$\eF$, and they are all there. 

The two rank 12 dimension $D^2 = 68$ modular data, $12_{0,68.}^{34,116}$ and
$12_{0,68.}^{68,166}$, are minimal modular extensions $\eC$ of its integral
fusion subcategory $\eF$, whose M\"uger center is obviously a fusion
subcategory generated by the abelian boson. A fusion category of FP-dimension 2
can only be equivalent to the representation category of  the cyclic group $G$
of order 2. Since $H^3(G, U(1))$ is also a cyclic group of order 2, there are
exactly two such minimal modular extensions.    \subsection{Gauging and
Zesting}

Some of the modular data realizations are related to other more familiar
categories by means of well-studied constructions such as gauging
\cite{CuiGalindoPlavnikWang} and zesting \cite{16fold,zesting}.  For example
$12_{4,144.}^{144,916}$ was first constructed by gauging, and
$10^{18,490}_{4,36.}$ was first constructed by zesting.  The condensed fiber
product \cite{CFP} is related to zesting and the minimal modular construction
described above.  We briefly describe gauging and zesting mathematically,
referring the interested reader to \emph{loc. cit.} for details.

\emph{Gauging} is a 2 step process.  One begins with a MTC $\cC$ with an action
of a finite group $G$ by braided tensor autoequivalences.  For example the rank
4 theory known as 3-fermions with fusion rules like $\Z_2\times\Z_2$ and
non-trivial objects having twist $-1$ admits an action of $S_3$, by permuting
the 3 nontrivial objects.  Then, assuming certain cohomological obstructions
vanish,  one constructs a $G$-crossed braided fusion category $\cD$ as a
$G$-extension of $\cC$.  There are typically choices to be made in this step,
which are parametrized by cohomological data.  Then one takes the
$G$-equivariantization $\cD^G$ of $\cD$. This will be a modular category of
dimension $\dim(\cC)\cdot |G|^2$ and will contain a symmetric subcategory
equivalent to $\Rep(G)$. It is often easy to see that a given modular category
is a gauging, by looking for these signatures.  Of course if one has modular
data with these signatures this helps in the search for a realization, which
must then be constructed explicitly by gauging.

\emph{Zesting} is similarly a process that uses cohomological data to build a
new category out of a given one.  In this case the set up is that one has an
$A$-graded MTC with a pointed (i.e., abelian) subcategory $\cB$ in the trivial
component (here $A$ is necessarily an abelian group).  For example $SU(3)_3$
i.e., $10^{6,152}_{4,36.}$ is $\Z_3$-graded and the trivial component contains
the pointed subcategory $\Rep(\Z_3)$.  Then one can twist the fusion rules
component-wise, using a $2$-cocycle on $A$ with values in $\cB$.  This requires
adjustment of the associativity constraints component-wise by a $3$-cochain.
In some cases the resulting category admits a braiding, which is obtained by
adjusting the braiding in $\cC$, again component-wise.  Finally, one may
typically adjust the twists to obtain a ribbon fusion category--which will be
modular under some mild conditions.  This all requires choices at each stage,
but is quite explicit. 

More recently, some instances of the zesting procedure has been explained in
more physically relevant terms, called the \emph{condensed fiber product}
\cite{CFP}.  This is also related to the minimal modular extension torsor
described above.  Here one takes two $A$-graded modular categories $\cD$ and $\cC$ with a common pointed symmetric
subcategory $\cB$ in the trivial components $\cD_0$ and $\cC_0$, and condenses the Tannakian diagonal subcategory $\Rep(B)\cong\Delta(\cB)\subset \cB\boxtimes\cB$ 
in $\cD\boxtimes \cC$, for some finite abelian group $B$.   The Tannakian subcategory $\Delta(\cB)$ centralizes the fiber product $\bigoplus_a \cD_a\boxtimes \cC_a$, so that $([\cD\boxtimes\cC]_B)_0$ is again a modular category.  If $\cD$ is pointed with $\dim(\cD)=|B|^2$ then the condensed fiber product coincides with zesting.  This is also related to the construction of \emph{anyon condensation} found in \cite{ZW240612068}.

\subsection{Potential modular data whose realizations are unknown}
\label{unknownrealization}

\begin{table}[tb] 
\caption{
A list of potential modular data (one for each Galois orbit)
whose realization are unknown or not definitively constructed.
}
\label{tab:exotic} 
\centerline{
\begin{tabular}{|l|l|p{3.7in}|}
\hline
$r_{c,D^2}^{\text{ord}T,\text{fp}}$ & $s_i$ & $d_i$ \\
\hline
$11_{3,1337.}^{48,634}$ &
$ 0,
0,
0,
\frac{1}{4},
\frac{1}{4},
\frac{3}{4},
\frac{3}{16},
\frac{11}{16},
0,
\frac{1}{3},
\frac{7}{12} $ &
$ 1$,
$ 3+2\sqrt{3}$,
$ 3+2\sqrt{3}$,
$ 4+2\sqrt{3}$,
$ 4+2\sqrt{3}$,
$ 6+4\sqrt{3}$,
$ 6+4\sqrt{3}$,
$ 6+4\sqrt{3}$,
$ 7+4\sqrt{3}$,
$ 8+4\sqrt{3}$,
$ 8+4\sqrt{3}$ (see Section \ref{ss: near group}) \\
\hline
$11_{\frac{32}{5},1964.}^{35,581}$ &
$0,
\frac{2}{35},
\frac{22}{35},
\frac{32}{35},
\frac{1}{5},
0,
\frac{3}{7},
\frac{5}{7},
\frac{6}{7},
\frac{3}{5},
\frac{1}{5} $ &
$1.0$,
$8.807$,
$8.807$,
$8.807$,
$11.632$,
$13.250$,
$14.250$,
$14.250$,
$14.250$,
$19.822$,
$20.440$  (for algebraic expressions, see Section \ref{sec:sum})
\\
\hline
$12_{4,1276.}^{39,406}$ &
$0,
\frac{1}{13},
\frac{3}{13},
\frac{4}{13},
\frac{9}{13},
\frac{10}{13},
\frac{12}{13},
0,
\frac{1}{3},
\frac{1}{3},
\frac{2}{3},
\frac{2}{3} $ &
$ 1$,
$ \frac{9+3\sqrt{13}}{2}$,
$ \frac{9+3\sqrt{13}}{2}$,
$ \frac{9+3\sqrt{13}}{2}$,
$ \frac{9+3\sqrt{13}}{2}$,
$ \frac{9+3\sqrt{13}}{2}$,
$ \frac{9+3\sqrt{13}}{2}$,
$ \frac{11+3\sqrt{13}}{2}$,
$ \frac{13+3\sqrt{13}}{2}$,
$ \frac{13+3\sqrt{13}}{2}$,
$ \frac{13+3\sqrt{13}}{2}$,
$ \frac{13+3\sqrt{13}}{2}$ 
\\
\hline
$12_{0,3926.}^{21,464}$ &
$0,
\frac{3}{7},
\frac{4}{7},
0,
0,
\frac{1}{7},
\frac{6}{7},
0,
\frac{1}{3},
\frac{2}{3},
\frac{2}{7},
\frac{5}{7} $&
$ 1$,
$ 3\xi_{14}^{5}$,
$ 3\xi_{14}^{5}$,
$  7+6c^{1}_{7} +3c^{2}_{7} $,
$  8+6c^{1}_{7} +3c^{2}_{7} $,
$  9+6c^{1}_{7} +3c^{2}_{7} $,
$  9+6c^{1}_{7} +3c^{2}_{7} $,
$  11+9c^{1}_{7} +3c^{2}_{7} $,
$  11+9c^{1}_{7} +3c^{2}_{7} $,
$  11+9c^{1}_{7} +3c^{2}_{7} $,
$  15+12c^{1}_{7} +6c^{2}_{7} $,
$  15+12c^{1}_{7} +6c^{2}_{7} $\\
\hline
$12_{0,3926.}^{63,774}$ &
$0,
\frac{3}{7},
\frac{4}{7},
0,
0,
\frac{1}{7},
\frac{6}{7},
\frac{1}{9},
\frac{4}{9},
\frac{7}{9},
\frac{2}{7},
\frac{5}{7} $ &
$ 1$,
$ 3\xi_{14}^{5}$,
$ 3\xi_{14}^{5}$,
$  7+6c^{1}_{7}
+3c^{2}_{7}
$,
$  8+6c^{1}_{7}
+3c^{2}_{7}
$,
$  9+6c^{1}_{7}
+3c^{2}_{7}
$,
$  9+6c^{1}_{7}
+3c^{2}_{7}
$,
$  11+9c^{1}_{7}
+3c^{2}_{7}
$,
$  11+9c^{1}_{7}
+3c^{2}_{7}
$,
$  11+9c^{1}_{7}
+3c^{2}_{7}
$,
$  15+12c^{1}_{7}
+6c^{2}_{7}
$,
$  15+12c^{1}_{7}
+6c^{2}_{7}
$ \\
\hline
\end{tabular}
}
\end{table}

Table \ref{tab:exotic} lists the potential modular data whose realizations
are still unknown or unsure.  Those potential modular data have different
fusion rings from the the modular tensor categories generated from Kac-Moody
algebra and twisted quantum doubles, plus Abelian anyon condensations.  Those
data may correspond to new modular tensor categories, or they are fake modular
data.  For the 
 rank-11 data
$11_{3,1337.}^{48,634}$, we have some evidence that they can be realized by
centers of near-group fusion categories, followed by some condensation
reductions.  

In order to gain some understanding, in this section, we are going compute  the
potential condensible algebra $\cA = \oplus_i A_i (d_i,s_i)$ (see Appendix
\ref{MMA}) for those potential modular data. Here $(d_i,s_i)$ is the
$i^\text{th}$ simple object, labeled by its quantum dimension $d_i$ and
topological spin $s_i$.

\begin{enumerate}

\item
\textbf{Potential modular data $11_{3,1337.}^{48,634}$} has no
Lagrangian condensible algebra since the central charge $c\neq 0$.  It has one
potential condensible algebra $\cA = (1,0)\oplus(7+4\sqrt 3,0)$.  The
condensation of anyon $(7+4\sqrt 3,0)$ reduces the modular data
$11_{3,1337.}^{48,634}$ to the modular data $6_{3,6}^{12,534} = 2_{1,2}^{4,437}
\boxtimes 3_{2,3}^{3,527}$ -- a pointed $\Z_2\boxtimes \Z_3$ MTC, since
$11_{3,1337.}^{48,634}\boxtimes \overline{6_{3,6}^{12,534}}$ has Lagrangian
condensible algebras.

\item
\textbf{Potential modular data $11_{\frac{32}{5},1964.}^{35,581}$} (see Section \ref{sec:sum} for a more detailed descrition) has no
Lagrangian condensible algebra.  It has one potential condensible algebra $\cA
= (1,0)\oplus(13.250,0)$.  The condensation of anyon $(13.250,0)$ should reduce
$11_{\frac{32}{5},1964.}^{35,581}$ to an unitary modular data with $D^2 =
\frac{35-7\sqrt{5}}{2} = 9.6737$.  But there are no unitary modular data with $D^2 =
\frac{35-7\sqrt{5}}{2} = 9.6737$ based on our classification.  Thus $\cA =
(1,0)\oplus(13.250,0)$ is a fake condensible algebra, and
$11_{\frac{32}{5},1964.}^{35,581}$ has no condensible algebra.

\item
\textbf{Potential modular data $12_{4,1276.}^{39,406}$} has no
Lagrangian condensible algebra.  It has one potential condensible algebra $\cA
= (1,0)\oplus(\frac{11+3\sqrt{13}}{2},0)$.  The condensation of anyon
$(\frac{11+3\sqrt{13}}{2},0)$ should reduce the modular data
$12_{4,1276.}^{39,406}$ to the modular data $9_{4,9}^{3,277} =
\overline{9_{4,9}^{3,277}} = 3_{2,3}^{3,527} \boxtimes 3_{2,3}^{3,527}$ -- a
pointed $\Z_3\boxtimes \Z_3$ MTC.  We remark that $12_{4,1276.}^{39,406}$ is
related to Haagerup-Izumi modular data. They share the same set of quantum
dimension $d_i$.  Also,  $12_{4,1276.}^{39,406}$ has a $\Z_2\times\Z_2$
automorphism, generated exchanging two simple objects with $(d,s) =
(\frac{13+3\sqrt{13}}{2},\frac13)$, and exchanging two simple objects with
$(d,s) = (\frac{13+3\sqrt{13}}{2},\frac23)$.

\item
\textbf{Potential modular data $12_{0,3926.}^{21,464}$} has no Lagrangian
condensible algebra.  It has three potential condensible algebra 
\begin{align}
\cA_1 &= (1,0)\oplus(8+6c^{1}_{7} +3c^{2}_{7},0), 
\nonumber\\
\cA_2 &= (1,0)\oplus(11+9c^{1}_{7} +3c^{2}_{7},0), 
\nonumber\\
\cA_3 &= (1,0)
\oplus(7+6c^{1}_{7} +3c^{2}_{7},0)
\oplus(8+6c^{1}_{7} +3c^{2}_{7},0)
\oplus(11+9c^{1}_{7} +3c^{2}_{7},0). 
\end{align}
The condensation of $\cA_1$ should reduce $12_{0,3926.}^{21,464}$ to an unitary modular
data with $D^2 = 14-7  c^{2}_{7} = 17.1152$.  The condensation of $\cA_2$
should reduce $12_{0,3926.}^{21,464}$ to an unitary modular data with $D^2 = 35-14
c^{1}_{7} +21c^{2}_{7} = 8.1964$.  The condensation of $\cA_3$ should reduce
$12_{0,3926.}^{21,464}$ to an unitary modular data with $D^2 = 49-28  c^{1}_{7}
+28c^{2}_{7} = 1.6233$.  There are some modular data with $D^2 = 14-7
c^{2}_{7} = 17.1152$ and $D^2 = 35-14 c^{1}_{7} +21c^{2}_{7} = 8.1964$ at rank
9.  But all those modular data are not unitary.  There is no unitary modular data with
$D^2 = 49-28  c^{1}_{7} +28c^{2}_{7} = 1.6233$.  Thus $\cA_{1,2,3}$ are fake
condensible algebras, and $12_{0,3926.}^{21,464}$ has no condensible algebra.

\item
\textbf{Potential modular data $12_{0,3926.}^{63,774}$} has no Lagrangian
condensible algebra.  It has one potential condensible algebra $\cA=
(1,0)\oplus(8+6c^{1}_{7} +3c^{2}_{7},0)$.  The condensation of anyon
$(8+6c^{1}_{7} +3c^{2}_{7},0)$ should reduce $12_{0,3926.}^{63,774}$ to an unitary
modular data with $D^2 = 14-7  c^{2}_{7} = 17.1152$.  There are some modular
data with $D^2 = 14-7 c^{2}_{7} = 17.1152$ at rank 9.  But all those modular
data are not unitary.  Thus $\cA= (1,0)\oplus(8+6c^{1}_{7} +3c^{2}_{7},0)$ is a
fake condensible algebra, and $12_{0,3926.}^{63,774}$ has no condensible
algebra.

\end{enumerate}

\section{Classify symmetries via $\text{symTOs}$ (\ie UMTCs in the trivial Witt class)}

\begin{figure}[t]
\begin{center}
\includegraphics[scale=0.8]{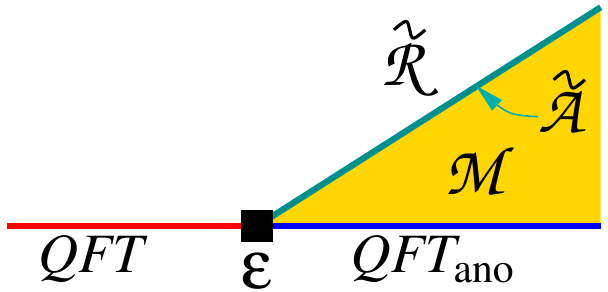}
\end{center}
\caption{The isomorphism $\veps$ (\ie the transparent domain wall in
space-time) between two anomaly-free (gapped or gapless) quantum field
theories, $QFT$ and $QFT_\mathrm{ano}\boxtimes_{\eM} \t \cR $, describes an low
energy equivalence of the two quantum field theories, below the energy gap of
the bulk topological order $\eM$ and its gapped boundary $\t\cR$
\cite{KZ150201690}.  Such an equivalence (called an isomorphic  holographic
decomposition) exposes the emergent symmetry in the quantum field theory $QFT$.
The emergent symmetry is described by the fusion higher category $\t\cR$ that
describes the excitations on the gapped boundary $\t\cR$, and will be referred
to as $\t\cR$-category symmetry.  The gapped boundary $\t\cR$ is induced by
Lagrangian condensable algebra $\t\cA$.  The holo-equivalent class of the
emergent symmetries, by definition, is described by a braided fusion higher
category $\eM$ that describes the excitations in the bulk topological order
$\eM$.  The bulk topological order $\eM$ will be referred to as the symTO,
which is a topological order with gappable boundary.} \label{QFTR} 
\end{figure}

We used to think symmetries are described by groups.  In recent years, we
realized that the low energy emergent symmetry in a quantum field theory $QFT$
can be a generalized symmetry beyond group and higher group.  It turns out that
finite generalized symmetries can all be described by higher fusion categories.

One way to obtain such a result is through the isomorphic holographic
decomposition in Fig. \ref{QFTR}, which was introduced to define homomorphism
between quantum field theories \cite{KZ150201690}.  If a quantum field theory
$QFT$ has an isomorphic holographic decomposition in Fig.  \ref{QFTR}, then we
say that the quantum field theory $QFT$ has a $\t\cR$-category symmetry.  Thus
a generalized symmetry is described by a fusion higher category $\t\cR$, which
describes the excitations on the gapped boundary $\t\cR$ of the bulk
topological order $\eM = \eZ(\t\cR)$ in Fig.  \ref{QFTR}.  Here $\eZ$ is the
generalized Drinfeld center that maps the fusion higher category $\t\cR$ to the
braided fusion higher category $\eM$ that describes the excitations in the bulk
topological order. 

The connection between boundary symmetry and bulk topological order was
observed in \Rf{LW0605}, where it was shown that topological entanglement
entropy arises from a boundary conservation law rooted in the bulk topological
order. This connection was later confirmed through numerical calculations
\cite{YS13094596}.  A systematic theory of symmetry topological-order (sym/TO)
correspondence was developed via holographic picture of emergent non-invertible
gravitational anomaly\cite{KZ150201690,KZ170501087,JW190513279}, holographic
picture of duality\cite{FT180600008,PV190906151}, which lead to holographic
picture of generalized
symmetry\cite{JW191213492,LB200304328,KZ200308898,KZ200514178,GK200805960,AS211202092,FT220907471}.
Symm/TO correspondence is also closely related to \emph{topological Wick
rotation} introduced in \Rf{KZ170501087,KZ190504924,KZ191201760}, which
summarizes a mathematical theory on how bulk can determine boundary.

A notion of holo-equivalence between symmetries was introduced in
\cite{TW191202817,KZ200308898,KZ200514178}.  Holo-equivalent symmetries impose
the same constraint on the dynamical properties of the associated systems
within the symmetric sub Hilbert spaces.  Thus holo-equivalent symmetries are
indistinguishable if we ignore the sectors with non-zero total symmetry charge.
The holo-equivalent class of the above $\t\cR$-category symmetry is described
by the bulk $\eM $ in Fig.  \ref{QFTR}, which will be referred to as the symTO
of $\t\cR$-category symmetry.  SymTOs are classified by topological orders $\eM
$ with gappable boundary.

For bosonic systems in 1-dimensional space, the holo-equivalent classes of the
generalized symmetries are classified by UMTCs in trivial Witt class.  The
symmetries in a holo-equivalent class described by an UMTC $\eM$ are classified
by the gapped boundaries (\ie fusion categories $\t\cR$) of $\eM$:
$\eM=\eZ(\t\cR)$.  In this section, we are going to use the classification of
UMTCs to classify symTOs via UMTCs in trivial Witt class, which in turn
classify (generalized) symmetries in 1-dimensional space.

How to determine if a UMTC is in the trivial Witt class or not?  We know some
necessary conditions for a UMTC to be in the trivial Witt class.  One set of
necessary conditions is the higher central charges \cite{NW181211234}.  Another
set of necessary conditions comes from Lagrangian condensable algebra
\cite{K13078244}, which is summarized in Appendix \ref{MMA}.  We apply these
two sets of necessary conditions and find the potential Witt-trivial UMTCs.  By
examining them one by one, we find that those potential Witt-trivial UMTCs are
actually in trivial Witt class.

Witt-trivial UMTCs do not exist for rank 2, 3, 5, 6, 7, 11.  In the following
we list Witt-trivial UMTCs (\ie symTOs) for rank 4, 8, 9, 10, 12.  We also list
the composite objects $\cA= \bigoplus_i A_i \, i$ (where $i$ are simple objects
and $A_i\in \N$) that give rise to the Lagrangian condensable algebras.
Actually, in next a few subsections, the composite object $\cA$ is listed by
its expansion coefficients $A_i$.  We note that a  Lagrangian condensable
algebra $\t\cA$ gives rise to a gapped boundary $\t\cR$ and a $\t\cR$-category
symmetry in the holo-equivalence class of the symTO.  

When $\t\cR$ is a local fusion category, it will describe an anomaly-free
symmetry\cite{KZ200514178}.  A $\t\cR$-category symmetry is anomaly-free if the
Lagrangian condensable algebra $\t\cA= \bigoplus_i \t A_i$ that induces $\t\cR$
has a dual Lagrangian condensable algebra $\cA  = \bigoplus_i A_i $ such that
\begin{align} 
\sum_i A_i \t A_i = 1.  
\end{align}
We will also indicate such an anomaly-free symmetry.

\subsection{Rank 4}

\noindent1. $4_{0,4.}^{2,750}$ \irep{0}:\ \ 
$d_i$ = ($1.0$,
$1.0$,
$1.0$,
$1.0$) 

\vskip 0.7ex
\hangindent=3em \hangafter=1
$D^2= 4.0 = 
4$

\vskip 0.7ex
\hangindent=3em \hangafter=1
$T = ( 0,
0,
0,
\frac{1}{2} )
$,

\vskip 0.7ex
\hangindent=3em \hangafter=1
$S$ = ($ 1$,
$ 1$,
$ 1$,
$ 1$;\ \ 
$ 1$,
$ -1$,
$ -1$;\ \ 
$ 1$,
$ -1$;\ \ 
$ 1$)

The holo-equivalence class of two symmetries.

 Lagrangian condensible algebra $A_i$:

$( 1,
1,
0,
0 )
$
$\to$ $\Vec_{\Z_2}$-category symmetry = $\Z_2$ symmetry

$( 1,
0,
1,
0 )
$
$\to$ $\Rep_{\Z_2}$-category symmetry.
Isomorphic to the above symmetry.

  \vskip 2ex

\noindent2. $4_{0,4.}^{4,375}$ \irep{0}:\ \ 
$d_i$ = ($1.0$,
$1.0$,
$1.0$,
$1.0$) 

\vskip 0.7ex
\hangindent=3em \hangafter=1
$D^2= 4.0 = 
4$

\vskip 0.7ex
\hangindent=3em \hangafter=1
$T = ( 0,
0,
\frac{1}{4},
\frac{3}{4} )
$,

\vskip 0.7ex
\hangindent=3em \hangafter=1
$S$ = ($ 1$,
$ 1$,
$ 1$,
$ 1$;\ \ 
$ 1$,
$ -1$,
$ -1$;\ \ 
$ -1$,
$ 1$;\ \ 
$ -1$)

Factors = $2_{1,2.}^{4,437}\boxtimes 2_{7,2.}^{4,625}$

The holo-equivalence class of one symmetry.

 Lagrangian condensible algebra $A_i$:

$( 1,
1,
0,
0 )
$
$\to$ $\t\cR_\text{semion}$-category symmetry = anomalous $\Z_2$ symmetry

  \vskip 2ex

\noindent3. $4_{0,13.09}^{5,872}$ \irep{5}:\ \ 
$d_i$ = ($1.0$,
$1.618$,
$1.618$,
$2.618$) 

\vskip 0.7ex
\hangindent=3em \hangafter=1
$D^2= 13.90 = 
\frac{15+5\sqrt{5}}{2}$

\vskip 0.7ex
\hangindent=3em \hangafter=1
$T = ( 0,
\frac{2}{5},
\frac{3}{5},
0 )
$,

\vskip 0.7ex
\hangindent=3em \hangafter=1
$S$ = ($ 1$,
$ \frac{1+\sqrt{5}}{2}$,
$ \frac{1+\sqrt{5}}{2}$,
$ \frac{3+\sqrt{5}}{2}$;\ \ 
$ -1$,
$ \frac{3+\sqrt{5}}{2}$,
$ -\frac{1+\sqrt{5}}{2}$;\ \ 
$ -1$,
$ -\frac{1+\sqrt{5}}{2}$;\ \ 
$ 1$)

Factors = $2_{\frac{14}{5},3.618}^{5,395}\boxtimes 2_{\frac{26}{5},3.618}^{5,720}$

The holo-equivalence class of one symmetry.

 Lagrangian condensible algebra $A_i$:

$( 1,
0,
0,
1 )
$
$\to$ $\t\cR_\text{Fib}$-category symmetry (beyond algebraic symmetry of \Rf{KZ200514178})

  \vskip 2ex

The holo-equivalent class of symmetries described by the symTO
$4_{0,4.}^{4,375}$ (the double-semion topological order)) contains one
symmetry, since the symTO has only one condensable algebra $A_i = (1,1,0,0)$.
The condensable algebra gives rise to a fusion category
$\t\cR=\t\cR_\mathrm{semion} \leftarrow 2_{1,2.}^{4,437}$, the fusion category
of a single semion\footnote{$\t\cR\leftarrow \eM$ means that the fusion
category $\t\cR$ is the fusion category formed by the objects in the braided
fusion category $\eM$, \ie $\leftarrow$ is the forgetful functor which ignores
the braiding. In this case the center of $\t\cR$ is given by $\eZ(\t\cR) = \eM
\boxtimes \overline{\eM}$.  }, whose Drinfeld center gives rise to the symTO
$4_{0,4.}^{4,375}$.  Such a symTO describes a single symmetry, the anomalous
$\Z_2$ symmetry \cite{CLW1141,LG1209}.

Although the symTO $4_{0,4.}^{2,750}$ (the $\Z_2$ topological order) has two
condensable algebra $A_i = (1,1,0,0)$ and $A_i = (1,0,1,0)$, the
holo-equivalent class of symmetries described by the symTO contains only one
symmetry.  The condensable algebra $A_i = (1,1,0,0)$ gives rise to a boundary
described by fusion category $\Vec_{\Z_2}$: $4_{0,4.}^{2,750} =
\eZ(\Vec_{\Z_2})$ \cite{KK11045047,CZ190312334}.  Thus, this condensable
algebra gives rise to the $\Vec_{\Z_2}$-category symmetry, which is nothing but
the usual $\Z_2$ symmetry.  The condensable algebra $A_i = (1,0,1,0)$ gives
rise to a boundary described by fusion category $\Rep_{\Z_2}$:
$4_{0,4.}^{2,750} = \eZ(\Rep_{\Z_2})$, and gives rise to
$\Rep_{\Z_2}$-category symmetry, which corresponds to the dual of the $\Z_2$
symmetry \cite{BT170402330,JW191209391}.  Since $\Z_2$ is an Abelian group, the
$\Z_2$ symmetry is isomorphic to the dual $\Z_2$ symmetry.  Therefore, the
holo-equivalent class $\eM = 4_{0,4.}^{2,750}$ contains only one symmetry.

The holo-equivalent class of symmetries described by the symTO
$4_{0,13.09}^{5,872}$ (the double-Fibonacci topological order) also contains
one symmetry.  The symTO is Drinfeld center of fusion category
$\t\cR=\t\cR_\mathrm{Fib} \leftarrow 2_{\frac{14}{5},3.618}^{5,395}$, the fusion category of a single Fibonacci anyon.
Such a symTO corresponds to an non-invertible symmetry in 1-dimensional space
that is beyond group theory description.  It is also beyond the algebraic
symmetry introduced in \cite{KZ200514178}, since it is anomalous.  This
anomalous non-invertible symmetry can appear as a low energy emergent symmetry
\cite{CW220506244}.

\subsection{Rank 8}

\noindent1. $8_{0,36.}^{6,213}$ \irep{0}:\ \ 
$d_i$ = ($1.0$,
$1.0$,
$2.0$,
$2.0$,
$2.0$,
$2.0$,
$3.0$,
$3.0$) 

\vskip 0.7ex
\hangindent=3em \hangafter=1
$D^2= 36.0 = 
36$

\vskip 0.7ex
\hangindent=3em \hangafter=1
$T = ( 0,
0,
0,
0,
\frac{1}{3},
\frac{2}{3},
0,
\frac{1}{2} )
$,

\vskip 0.7ex
\hangindent=3em \hangafter=1
$S$ = ($ 1$,
$ 1$,
$ 2$,
$ 2$,
$ 2$,
$ 2$,
$ 3$,
$ 3$;\ \ 
$ 1$,
$ 2$,
$ 2$,
$ 2$,
$ 2$,
$ -3$,
$ -3$;\ \ 
$ 4$,
$ -2$,
$ -2$,
$ -2$,
$0$,
$0$;\ \ 
$ 4$,
$ -2$,
$ -2$,
$0$,
$0$;\ \ 
$ -2$,
$ 4$,
$0$,
$0$;\ \ 
$ -2$,
$0$,
$0$;\ \ 
$ 3$,
$ -3$;\ \ 
$ 3$)

The holo-equivalence class of two symmetries.

 Lagrangian condensible algebra $A_i$:

$( 1, 1, 2, 0, 0, 0, 0, 0 ) $ $\to$ $\Vec_{S_3}$-category symmetry = $S_3$ symmetry 

$( 1, 1, 0, 2, 0, 0, 0, 0 ) $ $\to$ a symmetry isomorphic to the one above.

$( 1, 0, 0, 1, 0, 0, 1, 0 ) $ $\to$ $\Rep_{S_3}$-category symmetry = dual $S_3$ symmetry

$( 1, 0, 1, 0, 0, 0, 1, 0 ) $ $\to$ a symmetry isomorphic to the one above.

  \vskip 2ex

\noindent2. $8_{0,36.}^{12,101}$ \irep{0}:\ \ 
$d_i$ = ($1.0$,
$1.0$,
$2.0$,
$2.0$,
$2.0$,
$2.0$,
$3.0$,
$3.0$) 

\vskip 0.7ex
\hangindent=3em \hangafter=1
$D^2= 36.0 = 
36$

\vskip 0.7ex
\hangindent=3em \hangafter=1
$T = ( 0,
0,
0,
0,
\frac{1}{3},
\frac{2}{3},
\frac{1}{4},
\frac{3}{4} )
$,

\vskip 0.7ex
\hangindent=3em \hangafter=1
$S$ = ($ 1$,
$ 1$,
$ 2$,
$ 2$,
$ 2$,
$ 2$,
$ 3$,
$ 3$;\ \ 
$ 1$,
$ 2$,
$ 2$,
$ 2$,
$ 2$,
$ -3$,
$ -3$;\ \ 
$ 4$,
$ -2$,
$ -2$,
$ -2$,
$0$,
$0$;\ \ 
$ 4$,
$ -2$,
$ -2$,
$0$,
$0$;\ \ 
$ -2$,
$ 4$,
$0$,
$0$;\ \ 
$ -2$,
$0$,
$0$;\ \ 
$ -3$,
$ 3$;\ \ 
$ -3$)

The holo-equivalence class of one symmetry

 Lagrangian condensible algebra $A_i$:

$( 1, 1, 2, 0, 0, 0, 0, 0 ) $
$\to$ anomalous $S_3^{(3)}$ symmetry

$( 1, 1, 0, 2, 0, 0, 0, 0 ) $
$\to$ a symmetry isomorphic to the one above.

  \vskip 2ex

\noindent3. $8_{0,36.}^{18,162}$ \irep{0}:\ \ 
$d_i$ = ($1.0$,
$1.0$,
$2.0$,
$2.0$,
$2.0$,
$2.0$,
$3.0$,
$3.0$) 

\vskip 0.7ex
\hangindent=3em \hangafter=1
$D^2= 36.0 = 
36$

\vskip 0.7ex
\hangindent=3em \hangafter=1
$T = ( 0,
0,
0,
\frac{1}{9},
\frac{4}{9},
\frac{7}{9},
0,
\frac{1}{2} )
$,

\vskip 0.7ex
\hangindent=3em \hangafter=1
$S$ = ($ 1$,
$ 1$,
$ 2$,
$ 2$,
$ 2$,
$ 2$,
$ 3$,
$ 3$;\ \ 
$ 1$,
$ 2$,
$ 2$,
$ 2$,
$ 2$,
$ -3$,
$ -3$;\ \ 
$ 4$,
$ -2$,
$ -2$,
$ -2$,
$0$,
$0$;\ \ 
$ 2c_{9}^{2}$,
$ 2c_{9}^{4}$,
$ 2c_{9}^{1}$,
$0$,
$0$;\ \ 
$ 2c_{9}^{1}$,
$ 2c_{9}^{2}$,
$0$,
$0$;\ \ 
$ 2c_{9}^{4}$,
$0$,
$0$;\ \ 
$ 3$,
$ -3$;\ \ 
$ 3$)

The holo-equivalence class of two symmetries

 Lagrangian condensible algebra $A_i$:

$( 1, 1, 2, 0, 0, 0, 0, 0 ) $
$\to$ anomalous $S_3^{(4)}$ symmetry

$( 1, 0, 1, 0, 0, 0, 1, 0 ) $
$\to$ an anomalous non-invertible symmetry

  \vskip 2ex

\noindent4. $8_{0,36.}^{18,953}$ \irep{0}:\ \ 
$d_i$ = ($1.0$,
$1.0$,
$2.0$,
$2.0$,
$2.0$,
$2.0$,
$3.0$,
$3.0$) 

\vskip 0.7ex
\hangindent=3em \hangafter=1
$D^2= 36.0 = 
36$

\vskip 0.7ex
\hangindent=3em \hangafter=1
$T = ( 0,
0,
0,
\frac{2}{9},
\frac{5}{9},
\frac{8}{9},
0,
\frac{1}{2} )
$,

\vskip 0.7ex
\hangindent=3em \hangafter=1
$S$ = ($ 1$,
$ 1$,
$ 2$,
$ 2$,
$ 2$,
$ 2$,
$ 3$,
$ 3$;\ \ 
$ 1$,
$ 2$,
$ 2$,
$ 2$,
$ 2$,
$ -3$,
$ -3$;\ \ 
$ 4$,
$ -2$,
$ -2$,
$ -2$,
$0$,
$0$;\ \ 
$ 2c_{9}^{4}$,
$ 2c_{9}^{2}$,
$ 2c_{9}^{1}$,
$0$,
$0$;\ \ 
$ 2c_{9}^{1}$,
$ 2c_{9}^{4}$,
$0$,
$0$;\ \ 
$ 2c_{9}^{2}$,
$0$,
$0$;\ \ 
$ 3$,
$ -3$;\ \ 
$ 3$)

The holo-equivalence class of two symmetries

 Lagrangian condensible algebra $A_i$:

$( 1,
1,
2,
0,
0,
0,
0,
0 )
$
$\to$ anomalous $S_3^{(2)}$ symmetry

$( 1,
0,
1,
0,
0,
0,
1,
0 )
$
$\to$ an anomalous non-invertible symmetry

  \vskip 2ex

\noindent5. $8_{0,36.}^{36,495}$ \irep{0}:\ \ 
$d_i$ = ($1.0$,
$1.0$,
$2.0$,
$2.0$,
$2.0$,
$2.0$,
$3.0$,
$3.0$) 

\vskip 0.7ex
\hangindent=3em \hangafter=1
$D^2= 36.0 = 
36$

\vskip 0.7ex
\hangindent=3em \hangafter=1
$T = ( 0,
0,
0,
\frac{1}{9},
\frac{4}{9},
\frac{7}{9},
\frac{1}{4},
\frac{3}{4} )
$,

\vskip 0.7ex
\hangindent=3em \hangafter=1
$S$ = ($ 1$,
$ 1$,
$ 2$,
$ 2$,
$ 2$,
$ 2$,
$ 3$,
$ 3$;\ \ 
$ 1$,
$ 2$,
$ 2$,
$ 2$,
$ 2$,
$ -3$,
$ -3$;\ \ 
$ 4$,
$ -2$,
$ -2$,
$ -2$,
$0$,
$0$;\ \ 
$ 2c_{9}^{2}$,
$ 2c_{9}^{4}$,
$ 2c_{9}^{1}$,
$0$,
$0$;\ \ 
$ 2c_{9}^{1}$,
$ 2c_{9}^{2}$,
$0$,
$0$;\ \ 
$ 2c_{9}^{4}$,
$0$,
$0$;\ \ 
$ -3$,
$ 3$;\ \ 
$ -3$)

The holo-equivalence class of one symmetry

 Lagrangian condensible algebra $A_i$:

$( 1,
1,
2,
0,
0,
0,
0,
0 )
$
$\to$ anomalous $S_3^{(1)}$ symmetry

  \vskip 2ex

\noindent6. $8_{0,36.}^{36,171}$ \irep{0}:\ \ 
$d_i$ = ($1.0$,
$1.0$,
$2.0$,
$2.0$,
$2.0$,
$2.0$,
$3.0$,
$3.0$) 

\vskip 0.7ex
\hangindent=3em \hangafter=1
$D^2= 36.0 = 
36$

\vskip 0.7ex
\hangindent=3em \hangafter=1
$T = ( 0,
0,
0,
\frac{2}{9},
\frac{5}{9},
\frac{8}{9},
\frac{1}{4},
\frac{3}{4} )
$,

\vskip 0.7ex
\hangindent=3em \hangafter=1
$S$ = ($ 1$,
$ 1$,
$ 2$,
$ 2$,
$ 2$,
$ 2$,
$ 3$,
$ 3$;\ \ 
$ 1$,
$ 2$,
$ 2$,
$ 2$,
$ 2$,
$ -3$,
$ -3$;\ \ 
$ 4$,
$ -2$,
$ -2$,
$ -2$,
$0$,
$0$;\ \ 
$ 2c_{9}^{4}$,
$ 2c_{9}^{2}$,
$ 2c_{9}^{1}$,
$0$,
$0$;\ \ 
$ 2c_{9}^{1}$,
$ 2c_{9}^{4}$,
$0$,
$0$;\ \ 
$ 2c_{9}^{2}$,
$0$,
$0$;\ \ 
$ -3$,
$ 3$;\ \ 
$ -3$)

The holo-equivalence class of one symmetry

 Lagrangian condensible algebra $A_i$:

$( 1,
1,
2,
0,
0,
0,
0,
0 )
$
$\to$ anomalous $S_3^{(5)}$ symmetry

  \vskip 2ex

There are 6 symTOs at rank 8, giving rise to 6 holo-equivalent classes of
symmetries.  Each class contains an anomalous $S_3$-symmetry denoted as
$S_3^\om$.  In 1-dimensional space, the $S_3$-symmetry can have 6 types of
anomalies $\om \in H^3(S_3;\R/\Z) =\Z_6$ \cite{CGL1314} which correspond to
the 6 symTOs.

The first symTO $8_{0,36.}^{6,213}$ is the $S_3$-gauge theory in 2-dimensional
space.  One of its condensable algebra $A_i = (1,1,2,0,0,0,0,0)$ gives rise to
a fusion category $\Vec_{S_3}$.  Thus the condensable algebra $A_i =
(1,1,2,0,0,0,0,0)$ gives rise to a $\Vec_{S_3}$-category symmetry, which is
nothing but the usual $S_3$ symmetry.

The condensable algebra $A_i = (1,0,0,1,0,0,1,0)$ gives rise to a fusion
category $\Rep_{S_3}$.  Thus the condensable algebra $A_i = (1,0,0,1,0,0,1,0)$
gives rise to a holo-equivalent symmetry -- $\Rep_{S_3}$-category symmetry,
which was called dual $S_3$ symmetry in \Rf{JW191213492,KZ200514178}, and is a
non-invertible algebraic symmetry.  It is interesting to see an ordinary
group-like symmetry is holo-equivalent to a non-invertible algebraic symmetry.

We note that the symTO $8_{0,36.}^{6,213}$ has an automorphism of exchanging
simple objects $i=3$ and $i=4$.  The other two Lagrangian condensable algebras
are related to the previous two by this  automorphism.  Thus, the other two
condensable algebras also give rise to $S_3$ symmetry and dual $S_3$ symmetry,
respectively.  The symTO $8_{0,36.}^{12,101}$ also has an automorphism of
exchanging simple objects $i=3$ and $i=4$.  Its two Lagrangian condensable
algebras are related by this automorphism.  Such two Lagrangian condensable
algebras give rise to the same fusion category $\t\cR$ and the same
$\t\cR$-category symmetry.

Other rank-8 symTOs are Dijkgraaf-Witten $S_3$-gauge theories in 2-dimensional
space. The corresponding holo-equivalence classes contain anomalous $S_3^\om$
symmetry and anomalous non-invertible symmetry.  For more discussions, see
\Rf{CW220506244}.

\subsection{Rank 9 }

\noindent1. $9_{0,9.}^{3,113}$ \irep{0}:\ \ 
$d_i$ = ($1.0$,
$1.0$,
$1.0$,
$1.0$,
$1.0$,
$1.0$,
$1.0$,
$1.0$,
$1.0$) 

\vskip 0.7ex
\hangindent=3em \hangafter=1
$D^2= 9.0 = 
9$

\vskip 0.7ex
\hangindent=3em \hangafter=1
$T = ( 0,
0,
0,
0,
0,
\frac{1}{3},
\frac{1}{3},
\frac{2}{3},
\frac{2}{3} )
$,

\vskip 0.7ex
\hangindent=3em \hangafter=1
$S$ = ($ 1$,
$ 1$,
$ 1$,
$ 1$,
$ 1$,
$ 1$,
$ 1$,
$ 1$,
$ 1$;\ \ 
$ 1$,
$ 1$,
$ -\zeta_{6}^{1}$,
$ \zeta_{3}^{1}$,
$ -\zeta_{6}^{1}$,
$ \zeta_{3}^{1}$,
$ -\zeta_{6}^{1}$,
$ \zeta_{3}^{1}$;\ \ 
$ 1$,
$ \zeta_{3}^{1}$,
$ -\zeta_{6}^{1}$,
$ \zeta_{3}^{1}$,
$ -\zeta_{6}^{1}$,
$ \zeta_{3}^{1}$,
$ -\zeta_{6}^{1}$;\ \ 
$ 1$,
$ 1$,
$ -\zeta_{6}^{1}$,
$ \zeta_{3}^{1}$,
$ \zeta_{3}^{1}$,
$ -\zeta_{6}^{1}$;\ \ 
$ 1$,
$ \zeta_{3}^{1}$,
$ -\zeta_{6}^{1}$,
$ -\zeta_{6}^{1}$,
$ \zeta_{3}^{1}$;\ \ 
$ \zeta_{3}^{1}$,
$ -\zeta_{6}^{1}$,
$ 1$,
$ 1$;\ \ 
$ \zeta_{3}^{1}$,
$ 1$,
$ 1$;\ \ 
$ -\zeta_{6}^{1}$,
$ \zeta_{3}^{1}$;\ \ 
$ -\zeta_{6}^{1}$)

Factors = $3_{2,3.}^{3,527}\boxtimes 3_{6,3.}^{3,138}$

The holo-equivalence class of one symmetry

 Lagrangian condensible algebra $A_i$:

$( 1, 1, 1, 0, 0, 0, 0, 0, 0 ) $
$\to$ $\Vec_{\Z_3}$-category symmetry = $\Z_3$ symmetry 

$( 1, 0, 0, 1, 1, 0, 0, 0, 0 ) $
$\to$ $\Rep_{\Z_3}$-category symmetry.
Isomorphic to the above symmetry.

  \vskip 2ex

\noindent2. $9_{0,9.}^{9,620}$ \irep{0}:\ \ 
$d_i$ = ($1.0$,
$1.0$,
$1.0$,
$1.0$,
$1.0$,
$1.0$,
$1.0$,
$1.0$,
$1.0$) 

\vskip 0.7ex
\hangindent=3em \hangafter=1
$D^2= 9.0 = 
9$

\vskip 0.7ex
\hangindent=3em \hangafter=1
$T = ( 0,
0,
0,
\frac{1}{9},
\frac{1}{9},
\frac{4}{9},
\frac{4}{9},
\frac{7}{9},
\frac{7}{9} )
$,

\vskip 0.7ex
\hangindent=3em \hangafter=1
$S$ = ($ 1$,
$ 1$,
$ 1$,
$ 1$,
$ 1$,
$ 1$,
$ 1$,
$ 1$,
$ 1$;\ \ 
$ 1$,
$ 1$,
$ -\zeta_{6}^{1}$,
$ \zeta_{3}^{1}$,
$ -\zeta_{6}^{1}$,
$ \zeta_{3}^{1}$,
$ -\zeta_{6}^{1}$,
$ \zeta_{3}^{1}$;\ \ 
$ 1$,
$ \zeta_{3}^{1}$,
$ -\zeta_{6}^{1}$,
$ \zeta_{3}^{1}$,
$ -\zeta_{6}^{1}$,
$ \zeta_{3}^{1}$,
$ -\zeta_{6}^{1}$;\ \ 
$ -\zeta_{18}^{5}$,
$ \zeta_{9}^{2}$,
$ \zeta_{9}^{4}$,
$ -\zeta_{18}^{1}$,
$ \zeta_{9}^{1}$,
$ -\zeta_{18}^{7}$;\ \ 
$ -\zeta_{18}^{5}$,
$ -\zeta_{18}^{1}$,
$ \zeta_{9}^{4}$,
$ -\zeta_{18}^{7}$,
$ \zeta_{9}^{1}$;\ \ 
$ \zeta_{9}^{1}$,
$ -\zeta_{18}^{7}$,
$ -\zeta_{18}^{5}$,
$ \zeta_{9}^{2}$;\ \ 
$ \zeta_{9}^{1}$,
$ \zeta_{9}^{2}$,
$ -\zeta_{18}^{5}$;\ \ 
$ \zeta_{9}^{4}$,
$ -\zeta_{18}^{1}$;\ \ 
$ \zeta_{9}^{4}$)

The holo-equivalence class of one symmetry 

 Lagrangian condensible algebra $A_i$:

$( 1,
1,
1,
0,
0,
0,
0,
0,
0 )
$
$\to$ anomalous $\Z_3^{(1)}$ symmetry 

  \vskip 2ex

\noindent3. $9_{0,9.}^{9,462}$ \irep{0}:\ \ 
$d_i$ = ($1.0$,
$1.0$,
$1.0$,
$1.0$,
$1.0$,
$1.0$,
$1.0$,
$1.0$,
$1.0$) 

\vskip 0.7ex
\hangindent=3em \hangafter=1
$D^2= 9.0 = 
9$

\vskip 0.7ex
\hangindent=3em \hangafter=1
$T = ( 0,
0,
0,
\frac{2}{9},
\frac{2}{9},
\frac{5}{9},
\frac{5}{9},
\frac{8}{9},
\frac{8}{9} )
$,

\vskip 0.7ex
\hangindent=3em \hangafter=1
$S$ = ($ 1$,
$ 1$,
$ 1$,
$ 1$,
$ 1$,
$ 1$,
$ 1$,
$ 1$,
$ 1$;\ \ 
$ 1$,
$ 1$,
$ -\zeta_{6}^{1}$,
$ \zeta_{3}^{1}$,
$ -\zeta_{6}^{1}$,
$ \zeta_{3}^{1}$,
$ -\zeta_{6}^{1}$,
$ \zeta_{3}^{1}$;\ \ 
$ 1$,
$ \zeta_{3}^{1}$,
$ -\zeta_{6}^{1}$,
$ \zeta_{3}^{1}$,
$ -\zeta_{6}^{1}$,
$ \zeta_{3}^{1}$,
$ -\zeta_{6}^{1}$;\ \ 
$ -\zeta_{18}^{1}$,
$ \zeta_{9}^{4}$,
$ \zeta_{9}^{2}$,
$ -\zeta_{18}^{5}$,
$ -\zeta_{18}^{7}$,
$ \zeta_{9}^{1}$;\ \ 
$ -\zeta_{18}^{1}$,
$ -\zeta_{18}^{5}$,
$ \zeta_{9}^{2}$,
$ \zeta_{9}^{1}$,
$ -\zeta_{18}^{7}$;\ \ 
$ -\zeta_{18}^{7}$,
$ \zeta_{9}^{1}$,
$ -\zeta_{18}^{1}$,
$ \zeta_{9}^{4}$;\ \ 
$ -\zeta_{18}^{7}$,
$ \zeta_{9}^{4}$,
$ -\zeta_{18}^{1}$;\ \ 
$ \zeta_{9}^{2}$,
$ -\zeta_{18}^{5}$;\ \ 
$ \zeta_{9}^{2}$)

The holo-equivalence class of one symmetry 

 Lagrangian condensible algebra $A_i$:

$( 1,
1,
1,
0,
0,
0,
0,
0,
0 )
$
$\to$ anomalous $\Z_3^{(2)}$ symmetry 

  \vskip 2ex

\noindent4. $9_{0,16.}^{16,447}$ \irep{397}:\ \ 
$d_i$ = ($1.0$,
$1.0$,
$1.0$,
$1.0$,
$1.414$,
$1.414$,
$1.414$,
$1.414$,
$2.0$) 

\vskip 0.7ex
\hangindent=3em \hangafter=1
$D^2= 16.0 = 
16$

\vskip 0.7ex
\hangindent=3em \hangafter=1
$T = ( 0,
0,
\frac{1}{2},
\frac{1}{2},
\frac{1}{16},
\frac{7}{16},
\frac{9}{16},
\frac{15}{16},
0 )
$,

\vskip 0.7ex
\hangindent=3em \hangafter=1
$S$ = ($ 1$,
$ 1$,
$ 1$,
$ 1$,
$ \sqrt{2}$,
$ \sqrt{2}$,
$ \sqrt{2}$,
$ \sqrt{2}$,
$ 2$;\ \ 
$ 1$,
$ 1$,
$ 1$,
$ -\sqrt{2}$,
$ -\sqrt{2}$,
$ -\sqrt{2}$,
$ -\sqrt{2}$,
$ 2$;\ \ 
$ 1$,
$ 1$,
$ -\sqrt{2}$,
$ \sqrt{2}$,
$ -\sqrt{2}$,
$ \sqrt{2}$,
$ -2$;\ \ 
$ 1$,
$ \sqrt{2}$,
$ -\sqrt{2}$,
$ \sqrt{2}$,
$ -\sqrt{2}$,
$ -2$;\ \ 
$0$,
$ -2$,
$0$,
$ 2$,
$0$;\ \ 
$0$,
$ 2$,
$0$,
$0$;\ \ 
$0$,
$ -2$,
$0$;\ \ 
$0$,
$0$;\ \ 
$0$)

Factors = $3_{\frac{1}{2},4.}^{16,598}\boxtimes 3_{\frac{15}{2},4.}^{16,639}$

The holo-equivalence class of one symmetry 

 Lagrangian condensible algebra $A_i$:

$( 1,
1,
0,
0,
0,
0,
0,
0,
1 )
$
$\to$ anomalous non-invertible $\t\cR_\text{Ising}$-category symmetry  

  \vskip 2ex

\noindent5. $9_{0,16.}^{16,624}$ \irep{397}:\ \ 
$d_i$ = ($1.0$,
$1.0$,
$1.0$,
$1.0$,
$1.414$,
$1.414$,
$1.414$,
$1.414$,
$2.0$) 

\vskip 0.7ex
\hangindent=3em \hangafter=1
$D^2= 16.0 = 
16$

\vskip 0.7ex
\hangindent=3em \hangafter=1
$T = ( 0,
0,
\frac{1}{2},
\frac{1}{2},
\frac{3}{16},
\frac{5}{16},
\frac{11}{16},
\frac{13}{16},
0 )
$,

\vskip 0.7ex
\hangindent=3em \hangafter=1
$S$ = ($ 1$,
$ 1$,
$ 1$,
$ 1$,
$ \sqrt{2}$,
$ \sqrt{2}$,
$ \sqrt{2}$,
$ \sqrt{2}$,
$ 2$;\ \ 
$ 1$,
$ 1$,
$ 1$,
$ -\sqrt{2}$,
$ -\sqrt{2}$,
$ -\sqrt{2}$,
$ -\sqrt{2}$,
$ 2$;\ \ 
$ 1$,
$ 1$,
$ -\sqrt{2}$,
$ \sqrt{2}$,
$ -\sqrt{2}$,
$ \sqrt{2}$,
$ -2$;\ \ 
$ 1$,
$ \sqrt{2}$,
$ -\sqrt{2}$,
$ \sqrt{2}$,
$ -\sqrt{2}$,
$ -2$;\ \ 
$0$,
$ -2$,
$0$,
$ 2$,
$0$;\ \ 
$0$,
$ 2$,
$0$,
$0$;\ \ 
$0$,
$ -2$,
$0$;\ \ 
$0$,
$0$;\ \ 
$0$)

Factors = $3_{\frac{3}{2},4.}^{16,553}\boxtimes 3_{\frac{13}{2},4.}^{16,330}$

The holo-equivalence class of one symmetry 

 Lagrangian condensible algebra $A_i$:

$( 1,
1,
0,
0,
0,
0,
0,
0,
1 )
$
$\to$ anomalous non-invertible $\t\cR_\text{twIsing}$-category symmetry 

  \vskip 2ex

\noindent6. $9_{0,86.41}^{7,161}$ \irep{196}:\ \ 
$d_i$ = ($1.0$,
$1.801$,
$1.801$,
$2.246$,
$2.246$,
$3.246$,
$4.48$,
$4.48$,
$5.48$) 

\vskip 0.7ex
\hangindent=3em \hangafter=1
$D^2= 86.413 = 
49+35c^{1}_{7}
+14c^{2}_{7}
$

\vskip 0.7ex
\hangindent=3em \hangafter=1
$T = ( 0,
\frac{1}{7},
\frac{6}{7},
\frac{2}{7},
\frac{5}{7},
0,
\frac{3}{7},
\frac{4}{7},
0 )
$,

\vskip 0.7ex
\hangindent=3em \hangafter=1
$S$ = ($ 1$,
$ -c_{7}^{3}$,
$ -c_{7}^{3}$,
$ \xi_{7}^{3}$,
$ \xi_{7}^{3}$,
$ 2+c^{1}_{7}
$,
$ 2+2c^{1}_{7}
+c^{2}_{7}
$,
$ 2+2c^{1}_{7}
+c^{2}_{7}
$,
$ 3+2c^{1}_{7}
+c^{2}_{7}
$;\ \ 
$ -\xi_{7}^{3}$,
$ 2+c^{1}_{7}
$,
$ 2+2c^{1}_{7}
+c^{2}_{7}
$,
$ 1$,
$ -2-2  c^{1}_{7}
-c^{2}_{7}
$,
$ -3-2  c^{1}_{7}
-c^{2}_{7}
$,
$ -c_{7}^{3}$,
$ \xi_{7}^{3}$;\ \ 
$ -\xi_{7}^{3}$,
$ 1$,
$ 2+2c^{1}_{7}
+c^{2}_{7}
$,
$ -2-2  c^{1}_{7}
-c^{2}_{7}
$,
$ -c_{7}^{3}$,
$ -3-2  c^{1}_{7}
-c^{2}_{7}
$,
$ \xi_{7}^{3}$;\ \ 
$ c_{7}^{3}$,
$ 3+2c^{1}_{7}
+c^{2}_{7}
$,
$ -c_{7}^{3}$,
$ -2-c^{1}_{7}
$,
$ \xi_{7}^{3}$,
$ -2-2  c^{1}_{7}
-c^{2}_{7}
$;\ \ 
$ c_{7}^{3}$,
$ -c_{7}^{3}$,
$ \xi_{7}^{3}$,
$ -2-c^{1}_{7}
$,
$ -2-2  c^{1}_{7}
-c^{2}_{7}
$;\ \ 
$ 3+2c^{1}_{7}
+c^{2}_{7}
$,
$ -\xi_{7}^{3}$,
$ -\xi_{7}^{3}$,
$ 1$;\ \ 
$ 2+2c^{1}_{7}
+c^{2}_{7}
$,
$ 1$,
$ c_{7}^{3}$;\ \ 
$ 2+2c^{1}_{7}
+c^{2}_{7}
$,
$ c_{7}^{3}$;\ \ 
$ 2+c^{1}_{7}
$)

Factors = $3_{\frac{48}{7},9.295}^{7,790}\boxtimes 3_{\frac{8}{7},9.295}^{7,245}$

The holo-equivalence class of one symmetry 

 Lagrangian condensible algebra $A_i$:

$( 1,
0,
0,
0,
0,
1,
0,
0,
1 )
$
$\to$ anomalous non-invertible $\t\cR_{PSU(2)_5}$-category symmetry 

  \vskip 2ex

At rank 9, the first Abelian symTO $9_{0,9.}^{3,113}$ is the Drinfeld center of
$\Vec_{\Z_3}$. The corresponding holo-equivalent class of symmetries contains
$\Z_3$ symmetry.  The other two Abelian symTOs correspond to two
holo-equivalent classes of symmetries which contain anomalous $\Z_3$
symmetries.

The symTO $9_{0,16.}^{16,447}$ is the Drinfeld center of Ising fusion category
$\t\cR_\mathrm{Ising}\leftarrow 3_{\frac12,4.}^{16,598}$ (see Appendix
\ref{uni3}) of rank 3: $9_{0,16.}^{16,447} = \eZ(\t\cR_\mathrm{Ising})$.
Therefore, $9_{0,16.}^{16,447}$ describes the holo-equivalent class of
$\t\cR_\mathrm{Ising}$-category symmetry, which is a non-invertible symmetry
beyond the algebraic symmetry of \Rf{KZ200514178}.  Such a symmetry was
referred to as $\Z_2^e\vee \Z_2^m\vee \Z_2^{em}$ symmetry in \Rf{CW221214432}.
It contains a $\Z_2^m$, and a dual of $\Z_2^m$ symmetry which is the $\Z_2^e$
symmetry.  It also contains a $\Z_2^{em}$ duality symmetry that exchanges
$\Z_2^e$ and $\Z_2^m$.  This non-invertible $\t\cR_\mathrm{Ising}$-category
symmetry is important since it is the maximal emergent symmetry at the Ising
critical point.  In fact the self dual Ising model on 1-dimensional lattice 
\begin{align}
 H = - \sum_i Z_iZ_{i+1} + X_i
\end{align}
has this $\t\cR_\mathrm{Ising}$-category symmetry realized exactly, after
projecting into the $\Z_2$ symmetric sub-Hilbert space \cite{JW191213492}.  

In fact, the Ising fusion category $\t\cR_\mathrm{Ising}$ is the fusion
category that describes the fusion of the symmetry defects of the
$\t\cR_\mathrm{Ising}$-category symmetry \cite{FS09095013,KZ200514178}.  With
this understating, we can conclude that such a $\t\cR_\mathrm{Ising}$-category
symmetry is not anomaly-free.  This is because the fusion category that
describes the fusion of the symmetry defects of an anomaly-free symmetry must a
local fusion category \cite{TW191202817,KZ200514178}, whose objects must all
have integral quantum dimension.  Since, $\t\cR_\mathrm{Ising}$ has
non-integral quantum dimensions, it cannot be an anomaly-free symmetry.  In
fact, there is no local fusion category whose Drinfeld center is the symTO
$9_{0,16.}^{16,447}$.  Thus the holo-equivalent class of symmetries for this
symTO $9_{0,16.}^{16,447}$ contains no anomaly-free symmetries.

The symTOs $9_{0,16.}^{16,624}$ and $9_{0,86.41}^{7,161}$ also have a property
that the holo-equivalent class of symmetries for the symTOs contain no
anomaly-free symmetries.  The holo-equivalent class of  symTO
$9_{0,16.}^{16,624}$ contains a $\t\cR_\mathrm{twIsing}$-category symmetry,
where the twisted-Ising fusion category $\t\cR_\mathrm{twIsing}$ is given by
$\t\cR_\mathrm{twIsing} \leftarrow 3_{\frac32,4.}^{16,553}$.  The
holo-equivalent class of  symTO $9_{0,86.41}^{7,161}$ contains a
$\t\cR_{PSU(2)_5}$-category symmetry, where the fusion category
$\t\cR_{PSU(2)_5}$ is given by $\t\cR_{PSU(2)_5} \leftarrow
3_{\frac87,9.295}^{7,245}$.  Here $3_{\frac87,9.295}^{7,245}= PSU(2)_5$ is the
non-pointed Deligne factor of braided fusion category $SU(2)_5$: $SU(2)_5 =
PSU(2)_5\boxtimes \eM_\mathrm{pointed}$.

\subsection{Rank 10 }

\noindent1. $10_{0,89.56}^{12,155}$ \irep{587}:\ \ 
$d_i$ = ($1.0$,
$1.0$,
$2.732$,
$2.732$,
$2.732$,
$2.732$,
$2.732$,
$3.732$,
$3.732$,
$4.732$) 

\vskip 0.7ex
\hangindent=3em \hangafter=1
$D^2= 89.569 = 
48+24\sqrt{3}$

\vskip 0.7ex
\hangindent=3em \hangafter=1
$T = ( 0,
\frac{1}{2},
\frac{1}{3},
\frac{1}{4},
\frac{5}{6},
\frac{7}{12},
\frac{7}{12},
0,
\frac{1}{2},
0 )
$,

\vskip 0.7ex
\hangindent=3em \hangafter=1
$S$ = ($ 1$,
$ 1$,
$ 1+\sqrt{3}$,
$ 1+\sqrt{3}$,
$ 1+\sqrt{3}$,
$ 1+\sqrt{3}$,
$ 1+\sqrt{3}$,
$ 2+\sqrt{3}$,
$ 2+\sqrt{3}$,
$ 3+\sqrt{3}$;\ \ 
$ 1$,
$ 1+\sqrt{3}$,
$ -1-\sqrt{3}$,
$ 1+\sqrt{3}$,
$ -1-\sqrt{3}$,
$ -1-\sqrt{3}$,
$ 2+\sqrt{3}$,
$ 2+\sqrt{3}$,
$ -3-\sqrt{3}$;\ \ 
$ 1+\sqrt{3}$,
$ -2-2\sqrt{3}$,
$ 1+\sqrt{3}$,
$ 1+\sqrt{3}$,
$ 1+\sqrt{3}$,
$ -1-\sqrt{3}$,
$ -1-\sqrt{3}$,
$0$;\ \ 
$0$,
$ 2+2\sqrt{3}$,
$0$,
$0$,
$ -1-\sqrt{3}$,
$ 1+\sqrt{3}$,
$0$;\ \ 
$ 1+\sqrt{3}$,
$ -1-\sqrt{3}$,
$ -1-\sqrt{3}$,
$ -1-\sqrt{3}$,
$ -1-\sqrt{3}$,
$0$;\ \ 
$(3+\sqrt{3})\mathrm{i}$,
$(-3-\sqrt{3})\mathrm{i}$,
$ -1-\sqrt{3}$,
$ 1+\sqrt{3}$,
$0$;\ \ 
$(3+\sqrt{3})\mathrm{i}$,
$ -1-\sqrt{3}$,
$ 1+\sqrt{3}$,
$0$;\ \ 
$ 1$,
$ 1$,
$ 3+\sqrt{3}$;\ \ 
$ 1$,
$ -3-\sqrt{3}$;\ \ 
$0$)

The holo-equivalence class of $\t\cR_{\frac12 E(6)}$-category symmetry 

 Lagrangian condensible algebra $A_i$:

$( 1,
0,
0,
0,
0,
0,
0,
1,
0,
1 )
$
$\to$ $\t\cR_{\frac12 E(6)}$-category symmetry

  \vskip 2ex

\noindent2. $10_{0,89.56}^{12,200}$ \irep{587}:\ \ 
$d_i$ = ($1.0$,
$1.0$,
$2.732$,
$2.732$,
$2.732$,
$2.732$,
$2.732$,
$3.732$,
$3.732$,
$4.732$) 

\vskip 0.7ex
\hangindent=3em \hangafter=1
$D^2= 89.569 = 
48+24\sqrt{3}$

\vskip 0.7ex
\hangindent=3em \hangafter=1
$T = ( 0,
\frac{1}{2},
\frac{2}{3},
\frac{3}{4},
\frac{1}{6},
\frac{5}{12},
\frac{5}{12},
0,
\frac{1}{2},
0 )
$,

\vskip 0.7ex
\hangindent=3em \hangafter=1
$S$ = ($ 1$,
$ 1$,
$ 1+\sqrt{3}$,
$ 1+\sqrt{3}$,
$ 1+\sqrt{3}$,
$ 1+\sqrt{3}$,
$ 1+\sqrt{3}$,
$ 2+\sqrt{3}$,
$ 2+\sqrt{3}$,
$ 3+\sqrt{3}$;\ \ 
$ 1$,
$ 1+\sqrt{3}$,
$ -1-\sqrt{3}$,
$ 1+\sqrt{3}$,
$ -1-\sqrt{3}$,
$ -1-\sqrt{3}$,
$ 2+\sqrt{3}$,
$ 2+\sqrt{3}$,
$ -3-\sqrt{3}$;\ \ 
$ 1+\sqrt{3}$,
$ -2-2\sqrt{3}$,
$ 1+\sqrt{3}$,
$ 1+\sqrt{3}$,
$ 1+\sqrt{3}$,
$ -1-\sqrt{3}$,
$ -1-\sqrt{3}$,
$0$;\ \ 
$0$,
$ 2+2\sqrt{3}$,
$0$,
$0$,
$ -1-\sqrt{3}$,
$ 1+\sqrt{3}$,
$0$;\ \ 
$ 1+\sqrt{3}$,
$ -1-\sqrt{3}$,
$ -1-\sqrt{3}$,
$ -1-\sqrt{3}$,
$ -1-\sqrt{3}$,
$0$;\ \ 
$(-3-\sqrt{3})\mathrm{i}$,
$(3+\sqrt{3})\mathrm{i}$,
$ -1-\sqrt{3}$,
$ 1+\sqrt{3}$,
$0$;\ \ 
$(-3-\sqrt{3})\mathrm{i}$,
$ -1-\sqrt{3}$,
$ 1+\sqrt{3}$,
$0$;\ \ 
$ 1$,
$ 1$,
$ 3+\sqrt{3}$;\ \ 
$ 1$,
$ -3-\sqrt{3}$;\ \ 
$0$)

The holo-equivalence class of $\t\cR_{\frac12 \overline{E(6})}$-category symmetry 

 Lagrangian condensible algebra $A_i$:

$( 1,
0,
0,
0,
0,
0,
0,
1,
0,
1 )
$
$\to$ $\t\cR_{\frac12 \overline{E(6})}$-category symmetry 

  \vskip 2ex

\noindent3. $10_{0,1435.}^{20,676}$ \irep{948}:\ \ 
$d_i$ = ($1.0$,
$9.472$,
$9.472$,
$9.472$,
$9.472$,
$9.472$,
$9.472$,
$16.944$,
$16.944$,
$17.944$) 

\vskip 0.7ex
\hangindent=3em \hangafter=1
$D^2= 1435.541 = 
720+320\sqrt{5}$

\vskip 0.7ex
\hangindent=3em \hangafter=1
$T = ( 0,
0,
0,
\frac{1}{4},
\frac{1}{4},
\frac{3}{4},
\frac{3}{4},
\frac{2}{5},
\frac{3}{5},
0 )
$,

\vskip 0.7ex
\hangindent=3em \hangafter=1
$S$ = ($ 1$,
$ 5+2\sqrt{5}$,
$ 5+2\sqrt{5}$,
$ 5+2\sqrt{5}$,
$ 5+2\sqrt{5}$,
$ 5+2\sqrt{5}$,
$ 5+2\sqrt{5}$,
$ 8+4\sqrt{5}$,
$ 8+4\sqrt{5}$,
$ 9+4\sqrt{5}$;\ \ 
$ 15+6\sqrt{5}$,
$ -5-2\sqrt{5}$,
$ -5-2\sqrt{5}$,
$ -5-2\sqrt{5}$,
$ -5-2\sqrt{5}$,
$ -5-2\sqrt{5}$,
$0$,
$0$,
$ 5+2\sqrt{5}$;\ \ 
$ 15+6\sqrt{5}$,
$ -5-2\sqrt{5}$,
$ -5-2\sqrt{5}$,
$ -5-2\sqrt{5}$,
$ -5-2\sqrt{5}$,
$0$,
$0$,
$ 5+2\sqrt{5}$;\ \ 
$ -3-6  s^{1}_{20}
-4  c^{2}_{20}
+14s^{3}_{20}
$,
$ -3+6s^{1}_{20}
-4  c^{2}_{20}
-14  s^{3}_{20}
$,
$ 5+2\sqrt{5}$,
$ 5+2\sqrt{5}$,
$0$,
$0$,
$ 5+2\sqrt{5}$;\ \ 
$ -3-6  s^{1}_{20}
-4  c^{2}_{20}
+14s^{3}_{20}
$,
$ 5+2\sqrt{5}$,
$ 5+2\sqrt{5}$,
$0$,
$0$,
$ 5+2\sqrt{5}$;\ \ 
$ -3+6s^{1}_{20}
-4  c^{2}_{20}
-14  s^{3}_{20}
$,
$ -3-6  s^{1}_{20}
-4  c^{2}_{20}
+14s^{3}_{20}
$,
$0$,
$0$,
$ 5+2\sqrt{5}$;\ \ 
$ -3+6s^{1}_{20}
-4  c^{2}_{20}
-14  s^{3}_{20}
$,
$0$,
$0$,
$ 5+2\sqrt{5}$;\ \ 
$ -6-2\sqrt{5}$,
$ 14+6\sqrt{5}$,
$ -8-4\sqrt{5}$;\ \ 
$ -6-2\sqrt{5}$,
$ -8-4\sqrt{5}$;\ \ 
$ 1$)

The holo-equivalence class of two symmetries 

Lagrangian condensible algebra $A_i$:

$( 1, 2, 0, 0, 0, 0, 0, 0, 0, 1 ) $
$\to$ an anomalous non-invertible symmetry

$( 1, 0, 2, 0, 0, 0, 0, 0, 0, 1 ) $
$\to$ a symmetry isomorphic to the one above.

$( 1, 1, 1, 0, 0, 0, 0, 0, 0, 1 ) $
$\to$ an anomalous non-invertible symmetry

  \vskip 2ex

The first rank 10 symTO $10_{0,89.56}^{12,155}$ has only one condensable
algebra $A_i=(1,0,0,0,0,0,0,1,0,1)$.  Thus the symTO is the Drinfeld center of
only one fusion category, the so called $\frac12 E(6)$ fusion category
$\t\cR_{\frac12 E(6)}$\cite{HW07105761}.  The holo-equivalent class of
symmetries for this symTO contains only one symmetry which is the anomalous
non-invertible $\t\cR_{\frac12 E(6)}$-category symmetry and is beyond the
algebraic symmetry of \Rf{KZ200514178}.  The second rank 10 symTO
$10_{0,89.56}^{12,200}$ is the complex conjugation of the first, which contains
$\t\cR_{\frac12 \overline{E(6})}$-category symmetry.

The third rank 10 symTO $10_{0,1435.}^{20,676}$ is a condensation reduction of
$\eZ(\cC)$ {for some $\cC\in\cNG(\Z_4\times \Z_4, 16)$}, as described in
Section \ref{realizations}, see \cite{YuZhang}.  Such a symTO has a $\Z_2\times
\Z_2$ automorphism group generated by the following exchanges of simple objects
$i=4 \leftrightarrow i=5)$ and $i=6 \leftrightarrow i=7)$.  The symTO appears
to have three gapped boundaries, that give rise to only two non-invertible
$\t\cR$-category symmetries, due to the automorphisms.

\subsection{Rank 12 }

\noindent1. $12_{0,1276.}^{39,560}$ \irep{4851}:\ \ 
$d_i$ = ($1.0$,
$9.908$,
$9.908$,
$9.908$,
$9.908$,
$9.908$,
$9.908$,
$10.908$,
$11.908$,
$11.908$,
$11.908$,
$11.908$) 

\vskip 0.7ex
\hangindent=3em \hangafter=1
$D^2= 1276.274 = 
\frac{1287+351\sqrt{13}}{2}$

\vskip 0.7ex
\hangindent=3em \hangafter=1
$T = ( 0,
\frac{2}{13},
\frac{5}{13},
\frac{6}{13},
\frac{7}{13},
\frac{8}{13},
\frac{11}{13},
0,
0,
0,
\frac{1}{3},
\frac{2}{3} )
$,

\vskip 0.7ex
\hangindent=3em \hangafter=1
$S$ = ($ 1$,
$ \frac{9+3\sqrt{13}}{2}$,
$ \frac{9+3\sqrt{13}}{2}$,
$ \frac{9+3\sqrt{13}}{2}$,
$ \frac{9+3\sqrt{13}}{2}$,
$ \frac{9+3\sqrt{13}}{2}$,
$ \frac{9+3\sqrt{13}}{2}$,
$ \frac{11+3\sqrt{13}}{2}$,
$ \frac{13+3\sqrt{13}}{2}$,
$ \frac{13+3\sqrt{13}}{2}$,
$ \frac{13+3\sqrt{13}}{2}$,
$ \frac{13+3\sqrt{13}}{2}$;\ \ 
$ -\frac{9+3\sqrt{13}}{2}c_{13}^{4}$,
$ -\frac{9+3\sqrt{13}}{2}c_{13}^{1}$,
$ -\frac{9+3\sqrt{13}}{2}c_{13}^{3}$,
$ -\frac{9+3\sqrt{13}}{2}c_{13}^{2}$,
$ -\frac{9+3\sqrt{13}}{2}c_{13}^{5}$,
$ -\frac{9+3\sqrt{13}}{2}c_{13}^{6}$,
$ -\frac{9+3\sqrt{13}}{2}$,
$0$,
$0$,
$0$,
$0$;\ \ 
$ -\frac{9+3\sqrt{13}}{2}c_{13}^{3}$,
$ -\frac{9+3\sqrt{13}}{2}c_{13}^{4}$,
$ -\frac{9+3\sqrt{13}}{2}c_{13}^{6}$,
$ -\frac{9+3\sqrt{13}}{2}c_{13}^{2}$,
$ -\frac{9+3\sqrt{13}}{2}c_{13}^{5}$,
$ -\frac{9+3\sqrt{13}}{2}$,
$0$,
$0$,
$0$,
$0$;\ \ 
$ -\frac{9+3\sqrt{13}}{2}c_{13}^{1}$,
$ -\frac{9+3\sqrt{13}}{2}c_{13}^{5}$,
$ -\frac{9+3\sqrt{13}}{2}c_{13}^{6}$,
$ -\frac{9+3\sqrt{13}}{2}c_{13}^{2}$,
$ -\frac{9+3\sqrt{13}}{2}$,
$0$,
$0$,
$0$,
$0$;\ \ 
$ -\frac{9+3\sqrt{13}}{2}c_{13}^{1}$,
$ -\frac{9+3\sqrt{13}}{2}c_{13}^{4}$,
$ -\frac{9+3\sqrt{13}}{2}c_{13}^{3}$,
$ -\frac{9+3\sqrt{13}}{2}$,
$0$,
$0$,
$0$,
$0$;\ \ 
$ -\frac{9+3\sqrt{13}}{2}c_{13}^{3}$,
$ -\frac{9+3\sqrt{13}}{2}c_{13}^{1}$,
$ -\frac{9+3\sqrt{13}}{2}$,
$0$,
$0$,
$0$,
$0$;\ \ 
$ -\frac{9+3\sqrt{13}}{2}c_{13}^{4}$,
$ -\frac{9+3\sqrt{13}}{2}$,
$0$,
$0$,
$0$,
$0$;\ \ 
$ 1$,
$ \frac{13+3\sqrt{13}}{2}$,
$ \frac{13+3\sqrt{13}}{2}$,
$ \frac{13+3\sqrt{13}}{2}$,
$ \frac{13+3\sqrt{13}}{2}$;\ \ 
$  13+3\sqrt{13} $,
$ -\frac{13+3\sqrt{13}}{2}$,
$ -\frac{13+3\sqrt{13}}{2}$,
$ -\frac{13+3\sqrt{13}}{2}$;\ \ 
$  13+3\sqrt{13} $,
$ -\frac{13+3\sqrt{13}}{2}$,
$ -\frac{13+3\sqrt{13}}{2}$;\ \ 
$ -\frac{13+3\sqrt{13}}{2}$,
$  13+3\sqrt{13} $;\ \ 
$ -\frac{13+3\sqrt{13}}{2}$)

The holo-equivalent class of two symmetries.

 Lagrangian condensible algebra $A_i$:

$( 1, 0, 0, 0, 0, 0, 0, 1, 2, 0, 0, 0 ) $ 
$\to$ an anomalous non-invertible symmetry

$( 1, 0, 0, 0, 0, 0, 0, 1, 0, 2, 0, 0 ) $
$\to$ a symmetry isomorphic to the one above

$( 1, 0, 0, 0, 0, 0, 0, 1, 1, 1, 0, 0 ) $
$\to$ an anomalous non-invertible symmetry

  \vskip 2ex

\noindent2. $12_{0,1276.}^{117,251}$ \irep{5248}:\ \ 
$d_i$ = ($1.0$,
$9.908$,
$9.908$,
$9.908$,
$9.908$,
$9.908$,
$9.908$,
$10.908$,
$11.908$,
$11.908$,
$11.908$,
$11.908$) 

\vskip 0.7ex
\hangindent=3em \hangafter=1
$D^2= 1276.274 = 
\frac{1287+351\sqrt{13}}{2}$

\vskip 0.7ex
\hangindent=3em \hangafter=1
$T = ( 0,
\frac{2}{13},
\frac{5}{13},
\frac{6}{13},
\frac{7}{13},
\frac{8}{13},
\frac{11}{13},
0,
0,
\frac{1}{9},
\frac{4}{9},
\frac{7}{9} )
$,

\vskip 0.7ex
\hangindent=3em \hangafter=1
$S$ = ($ 1$,
$ \frac{9+3\sqrt{13}}{2}$,
$ \frac{9+3\sqrt{13}}{2}$,
$ \frac{9+3\sqrt{13}}{2}$,
$ \frac{9+3\sqrt{13}}{2}$,
$ \frac{9+3\sqrt{13}}{2}$,
$ \frac{9+3\sqrt{13}}{2}$,
$ \frac{11+3\sqrt{13}}{2}$,
$ \frac{13+3\sqrt{13}}{2}$,
$ \frac{13+3\sqrt{13}}{2}$,
$ \frac{13+3\sqrt{13}}{2}$,
$ \frac{13+3\sqrt{13}}{2}$;\ \ 
$ -\frac{9+3\sqrt{13}}{2}c_{13}^{4}$,
$ -\frac{9+3\sqrt{13}}{2}c_{13}^{1}$,
$ -\frac{9+3\sqrt{13}}{2}c_{13}^{3}$,
$ -\frac{9+3\sqrt{13}}{2}c_{13}^{2}$,
$ -\frac{9+3\sqrt{13}}{2}c_{13}^{5}$,
$ -\frac{9+3\sqrt{13}}{2}c_{13}^{6}$,
$ -\frac{9+3\sqrt{13}}{2}$,
$0$,
$0$,
$0$,
$0$;\ \ 
$ -\frac{9+3\sqrt{13}}{2}c_{13}^{3}$,
$ -\frac{9+3\sqrt{13}}{2}c_{13}^{4}$,
$ -\frac{9+3\sqrt{13}}{2}c_{13}^{6}$,
$ -\frac{9+3\sqrt{13}}{2}c_{13}^{2}$,
$ -\frac{9+3\sqrt{13}}{2}c_{13}^{5}$,
$ -\frac{9+3\sqrt{13}}{2}$,
$0$,
$0$,
$0$,
$0$;\ \ 
$ -\frac{9+3\sqrt{13}}{2}c_{13}^{1}$,
$ -\frac{9+3\sqrt{13}}{2}c_{13}^{5}$,
$ -\frac{9+3\sqrt{13}}{2}c_{13}^{6}$,
$ -\frac{9+3\sqrt{13}}{2}c_{13}^{2}$,
$ -\frac{9+3\sqrt{13}}{2}$,
$0$,
$0$,
$0$,
$0$;\ \ 
$ -\frac{9+3\sqrt{13}}{2}c_{13}^{1}$,
$ -\frac{9+3\sqrt{13}}{2}c_{13}^{4}$,
$ -\frac{9+3\sqrt{13}}{2}c_{13}^{3}$,
$ -\frac{9+3\sqrt{13}}{2}$,
$0$,
$0$,
$0$,
$0$;\ \ 
$ -\frac{9+3\sqrt{13}}{2}c_{13}^{3}$,
$ -\frac{9+3\sqrt{13}}{2}c_{13}^{1}$,
$ -\frac{9+3\sqrt{13}}{2}$,
$0$,
$0$,
$0$,
$0$;\ \ 
$ -\frac{9+3\sqrt{13}}{2}c_{13}^{4}$,
$ -\frac{9+3\sqrt{13}}{2}$,
$0$,
$0$,
$0$,
$0$;\ \ 
$ 1$,
$ \frac{13+3\sqrt{13}}{2}$,
$ \frac{13+3\sqrt{13}}{2}$,
$ \frac{13+3\sqrt{13}}{2}$,
$ \frac{13+3\sqrt{13}}{2}$;\ \ 
$  13+3\sqrt{13} $,
$ -\frac{13+3\sqrt{13}}{2}$,
$ -\frac{13+3\sqrt{13}}{2}$,
$ -\frac{13+3\sqrt{13}}{2}$;\ \ 
$ \frac{13+3\sqrt{13}}{2}c_{9}^{2}$,
$ \frac{13+3\sqrt{13}}{2}c_{9}^{4}$,
$ \frac{13+3\sqrt{13}}{2}c_{9}^{1}$;\ \ 
$ \frac{13+3\sqrt{13}}{2}c_{9}^{1}$,
$ \frac{13+3\sqrt{13}}{2}c_{9}^{2}$;\ \ 
$ \frac{13+3\sqrt{13}}{2}c_{9}^{4}$)

The holo-equivalent class of one symmetry.

 Lagrangian condensible algebra $A_i$:

$( 1,
0,
0,
0,
0,
0,
0,
1,
2,
0,
0,
0 )
$
$\to$ an anomalous non-invertible symmetry

  \vskip 2ex

\noindent3. $12_{0,1276.}^{117,145}$ \irep{5248}:\ \ 
$d_i$ = ($1.0$,
$9.908$,
$9.908$,
$9.908$,
$9.908$,
$9.908$,
$9.908$,
$10.908$,
$11.908$,
$11.908$,
$11.908$,
$11.908$) 

\vskip 0.7ex
\hangindent=3em \hangafter=1
$D^2= 1276.274 = 
\frac{1287+351\sqrt{13}}{2}$

\vskip 0.7ex
\hangindent=3em \hangafter=1
$T = ( 0,
\frac{2}{13},
\frac{5}{13},
\frac{6}{13},
\frac{7}{13},
\frac{8}{13},
\frac{11}{13},
0,
0,
\frac{2}{9},
\frac{5}{9},
\frac{8}{9} )
$,

\vskip 0.7ex
\hangindent=3em \hangafter=1
$S$ = ($ 1$,
$ \frac{9+3\sqrt{13}}{2}$,
$ \frac{9+3\sqrt{13}}{2}$,
$ \frac{9+3\sqrt{13}}{2}$,
$ \frac{9+3\sqrt{13}}{2}$,
$ \frac{9+3\sqrt{13}}{2}$,
$ \frac{9+3\sqrt{13}}{2}$,
$ \frac{11+3\sqrt{13}}{2}$,
$ \frac{13+3\sqrt{13}}{2}$,
$ \frac{13+3\sqrt{13}}{2}$,
$ \frac{13+3\sqrt{13}}{2}$,
$ \frac{13+3\sqrt{13}}{2}$;\ \ 
$ -\frac{9+3\sqrt{13}}{2}c_{13}^{4}$,
$ -\frac{9+3\sqrt{13}}{2}c_{13}^{1}$,
$ -\frac{9+3\sqrt{13}}{2}c_{13}^{3}$,
$ -\frac{9+3\sqrt{13}}{2}c_{13}^{2}$,
$ -\frac{9+3\sqrt{13}}{2}c_{13}^{5}$,
$ -\frac{9+3\sqrt{13}}{2}c_{13}^{6}$,
$ -\frac{9+3\sqrt{13}}{2}$,
$0$,
$0$,
$0$,
$0$;\ \ 
$ -\frac{9+3\sqrt{13}}{2}c_{13}^{3}$,
$ -\frac{9+3\sqrt{13}}{2}c_{13}^{4}$,
$ -\frac{9+3\sqrt{13}}{2}c_{13}^{6}$,
$ -\frac{9+3\sqrt{13}}{2}c_{13}^{2}$,
$ -\frac{9+3\sqrt{13}}{2}c_{13}^{5}$,
$ -\frac{9+3\sqrt{13}}{2}$,
$0$,
$0$,
$0$,
$0$;\ \ 
$ -\frac{9+3\sqrt{13}}{2}c_{13}^{1}$,
$ -\frac{9+3\sqrt{13}}{2}c_{13}^{5}$,
$ -\frac{9+3\sqrt{13}}{2}c_{13}^{6}$,
$ -\frac{9+3\sqrt{13}}{2}c_{13}^{2}$,
$ -\frac{9+3\sqrt{13}}{2}$,
$0$,
$0$,
$0$,
$0$;\ \ 
$ -\frac{9+3\sqrt{13}}{2}c_{13}^{1}$,
$ -\frac{9+3\sqrt{13}}{2}c_{13}^{4}$,
$ -\frac{9+3\sqrt{13}}{2}c_{13}^{3}$,
$ -\frac{9+3\sqrt{13}}{2}$,
$0$,
$0$,
$0$,
$0$;\ \ 
$ -\frac{9+3\sqrt{13}}{2}c_{13}^{3}$,
$ -\frac{9+3\sqrt{13}}{2}c_{13}^{1}$,
$ -\frac{9+3\sqrt{13}}{2}$,
$0$,
$0$,
$0$,
$0$;\ \ 
$ -\frac{9+3\sqrt{13}}{2}c_{13}^{4}$,
$ -\frac{9+3\sqrt{13}}{2}$,
$0$,
$0$,
$0$,
$0$;\ \ 
$ 1$,
$ \frac{13+3\sqrt{13}}{2}$,
$ \frac{13+3\sqrt{13}}{2}$,
$ \frac{13+3\sqrt{13}}{2}$,
$ \frac{13+3\sqrt{13}}{2}$;\ \ 
$  13+3\sqrt{13} $,
$ -\frac{13+3\sqrt{13}}{2}$,
$ -\frac{13+3\sqrt{13}}{2}$,
$ -\frac{13+3\sqrt{13}}{2}$;\ \ 
$ \frac{13+3\sqrt{13}}{2}c_{9}^{4}$,
$ \frac{13+3\sqrt{13}}{2}c_{9}^{2}$,
$ \frac{13+3\sqrt{13}}{2}c_{9}^{1}$;\ \ 
$ \frac{13+3\sqrt{13}}{2}c_{9}^{1}$,
$ \frac{13+3\sqrt{13}}{2}c_{9}^{4}$;\ \ 
$ \frac{13+3\sqrt{13}}{2}c_{9}^{2}$)

The holo-equivalent class of one symmetry.

 Lagrangian condensible algebra $A_i$:

$( 1,
0,
0,
0,
0,
0,
0,
1,
2,
0,
0,
0 )
$
$\to$ an anomalous non-invertible symmetry

  \vskip 2ex

The first rank 12 symTO $12_{0,1276.}^{39,560}$ is the Haagerup-Izumi MTC
Haag$(1)_0$, which has three condensable algebras.  Its holo-equivalent class
contains only two symmetries, due to an automorphism that exchange two simple
objects $i=9 \leftrightarrow i=10$.  The other two rank 12 symTOs are the
Haagerup-Izumi MTC Haag$(1)_1$ and Haag$(1)_{-1}$. Their corresponding
holo-equivalent classes each contains one anomalous non-invertible symmetry.

\

\

\noindent
{\bf Acknowledgements:}
We thank Noah Snyder, Cain Edie-Michell, Terry Gannon, Qing Zhang, Masaki Izumi and Pinhas Grossman for useful
communications.  S.-H. N. was partially supported by NSF grant  DMS-1664418 and the Simons Foundation MPS-TSM-00008039.
E.C.R. was partially supported by NSF grant DMS-2205962 and a UK Royal Society Wolfson Visiting Fellowship.  X.-G.W was partially
supported by NSF grant DMR-2022428 and by the Simons Collaboration on
Ultra-Quantum Matter, which is a grant from the Simons Foundation (651446,
XGW).  The authors have no relevant financial or non-financial interests to
disclose.

\

\

\

\appendix

\section{$d$-number and $1/\rho_\textrm{pMD}(\frs)_{ui}$}
\label{sec:dnum}

The unit row of $\rho_\textrm{pMD}(\frs)$ has an useful property that
$1/\rho_\textrm{pMD}(\frs)_{ui}$ are the so called $d$-numbers.  
\begin{defn}
Let $K/\Q$ be a Galois extension and $O_K$ the ring of algebraic integers in $K$. An element $x \in Q_K$ is called a $d$-number if for any $\si \in \Gal(K/\Q)$, $\si(x)=x u$ for some unit $u \in O_K$.  
\end{defn}

\begin{thm}\label{l1}
Let $K$ be a subfield of $\Q(\zeta_{p^\ell})$ such that $[K: \Q]=q$ is a prime number. Let $O_K$ be the ring of algebraic integers in $K$. Then $x \in O_K$ is a $d$-number if and only if $x = n u^a v$ for some integer $n$, unit $u \in O_K$, $a = 0,\dots, q-1$ where $v= Norm^{\Q(\zeta_{p^\ell})}_K(1-\zeta_{p^\ell})$. If $\si$ is a generator of $\Gal(\Q(\zeta_{p^\ell})/\Q)$, then 
$$
v = (1-\zeta_{p^\ell})\si^q(1-\zeta_{p^\ell})\si^{2q}(1-\zeta_{p^\ell})\cdots \si^{p^{\ell-1}(p-1)-q}(1-\zeta_{p^\ell}).
$$
\end{thm}
\begin{proof}
Note that $\si(1-\zeta_{p^\ell}) = (1-\zeta_{p^\ell}) u'$ for some unit $u'$ in $\Z[\zeta_{p^\ell}]$, and $K$ is the fixed field of $\langle \si^q \rangle$ by Galois correspondence.  So $v \in K$ and $\si(v) = v u$ for some unit $u \in O_K$. Direct verification shows that $x = n v^a u$ for some integer $n$, unit $u \in O_K$, $a = 0,\dots, q-1$ is a $d$-number. 

Conversely, let $x \in O_K$ be a $d$-number. Then $\si(x) = x u'$ for some unit $u' \in O_K$. Thus $x^q = m u''$ for some unit $u'' \in O_K$ where $m = |Norm(x)|$. Let $m = p_0^{k_0} p_1^{k_1} \cdots p_l^{k_l}$ be the prime  factorization of $m$, where $p_0=p$ and $k_0 \ge 0$. Since $O_K$ is a Dedekind domain, for each $i > 0$, 
$$
(p_i) = \mathfrak{P}_{i,1}^{e_i} \cdots \mathfrak{P}_{i,g_i}^{e_i}
$$
for some conjugate distinct prime ideals $\mathfrak{P}_{i,1}, \dots, \mathfrak{P}_{i,g_i}$ of $O_K$ and $e_i f_i g_i = q$ where $f_i = \log_{p_i} |O_K/\mathfrak{P}_{i,1}|$. Since $q$ is prime and $p$ is the only totally ramified prime, $e_i < q$ for $i >0$. Since $q$ is a prime number, $e_i = 1$ and hence $(p_i) = \mathfrak{P}_{i,1} \cdots \mathfrak{P}_{i,g_i}$ for $i > 0$. The prime ideal factorization of $(p)$ in $O_K$ is $(p) = \mathfrak{P}_0^q$ where $\mathfrak{P}_0=(u)$. Thus, we have
$$
(x^q) = (x)^q = (m) = \mathfrak{P}_0^{q k_0} \mathfrak{P}_{1,1}^{k_1} \cdots \mathfrak{P}_{1,g_1}^{k_1} \cdots \mathfrak{P}_{l,1}^{k_l} \cdots \mathfrak{P}_{l,g_l}^{k_l}
$$
and hence
$$
(x) = \mathfrak{P}_0^{\al_{0}} \mathfrak{P}_{1,1}^{\al_{1,1}} \cdots \mathfrak{P}_{1,g_1}^{\al_{1, g_1}} \cdots \mathfrak{P}_{l,1}^{\al_{l,1}} \cdots \mathfrak{P}_{l,g_l}^{\al_{l, g_l}}
$$
for some nonnegative integers $\al_{0}, \al_{i,j}$.  By the unique factorizations of primes ideals, we have
$$
q \al_0 = q k_0, \quad q \al_{i, j} = k_i \text{ for } i = 1, \dots, l. 
$$
Thus, $\al_0 = k_0$ and $\al_{i, j} = \al_{i,1}$  for $i > 0$. Hence, 
$$
(x) = \mathfrak{P}_0^{k_0} (p_1)^{\al_{1,1}} \cdots  (p_l)^{\al_{l,1}} 
$$
and so
$$
x = v^{k_0} p_1^{\al_{1,1}} \cdots p_l^{\al_{l,1}} u'''
$$
for some unit $u''' \in O_K$. Since $(v^q) = (p)$, $x= n v^a u$ for some integer $n$, $u$ an unit in $O_K$ and $0 \le a \le q-1$.
\end{proof}

\begin{cor}\label{c1}
Let $K=\Q(\sqrt{5})$. Then $x \in O_K$ is a $d$-number if and only if b $x =n u$ or $x = n\sqrt{5} u$ for some integer $n$ and an unit $u \in O_K$.
\end{cor}
\begin{proof}
One can apply the preceding theorem as $[K:\Q]=2$ and $K$ is a subfield of $\Q(\zeta_5)$. Let $\si \in \Gal(\Q(\zeta_5)/\Q)$. Then $\si^2$ is the complex conjugation and so
$$
v = (1-\zeta_5) \si^2(1-\zeta_5) = 2-2\cos(2 \pi/5) = \sqrt{5} u'
$$ 
for some unit $u' \in Q_K$. Now, the result follows immediately from Theorem \ref{l1}.
\end{proof}

\begin{cor}\label{c2}
Let $K=\Q(\cos(2 \pi/7))$. Then $x \in O_K$ is a $d$-number if and only if $x =n v^a u$ for some integer $n$,  unit $u \in O_K$, $a = 0,1,2$, where $v = 2-2\cos(2\pi/7)$.
\end{cor}
\begin{proof}
One can apply the preceding theorem as $[K:\Q]=3$ and $K$ is the real subfield of $\Q(\zeta_7)$. The assignment $\si : \zeta_7 \mapsto \zeta_7^3$ defines a generator $\Gal(\Q(\zeta_7)/\Q)$. Since $7-1 = 6$, 
$$
v = (1-\zeta_7)\si^3(1-\zeta_7)  = (1-\zeta_7)(1-\zeta_7^{-1}) = 2-2\cos(2 \pi/7).
$$
 Now, the result follows from Theorem \ref{l1}.
\end{proof}

Using the above results, we find that
a $d$-number has the following explicit expression for simple conductors.  When
conductor $cnd=1$, all the $d$-numbers are given by
\begin{align}
 \del &= n , \ \ \ n \in \Z.
\end{align}
When conductor $cnd=3$, all the $d$-numbers are given by
\begin{align}
 \del &= n \zeta_{3}^m , \ \ \ n \in \Z, \ \ \ \ m=0,1,2.
\end{align}
When conductor $cnd=4$, all the $d$-numbers are given by
\begin{align}
 \del &= n \zeta_{4}^m , \ \ \ n \in \Z, \ \ \ \ m=0,1,2,3.
\end{align}

When conductor $cnd=5$, all the $d$-numbers are given by
\begin{align}
 \del &= n  (\zeta_{10} -1)^k \zeta_{10}^m v^{a}, \ \ \ n,k \in \Z,\ \ \ \ m=0,1,2,3,4,5,6,7,8,9
\nonumber\\
 v &= \zeta_{10} - \zeta_{10}^{-1}, \ \ \ a = 0,1,2,3.
\end{align}
All real $d$-numbers of conductor 5 are given by
\begin{align}
 \del = n  (\frac{1+\sqrt 5}2)^k v^{2a}, \ \ \
n,k \in \Z, \ a=0,1
\end{align}

We note that the generic complex $\del$ has the following form
\begin{align}
  \del &= n  (\zeta_{10} -1)^k \zeta_{10}^m v^{a}
= n  (\zeta_{20} -\zeta_{20}^{-1})^k \zeta_{20}^{2m+ k} v^{a}
\end{align}
where $\zeta_{20} -\zeta_{20}^{-1}$ and $v = \zeta_{10} - \zeta_{10}^{-1}$ are
imaginary.  $\zeta_{20}^{2m+ k}$ is real if $2m+ k$ mod 10 = 0 and it is
imaginary if $2m+ k$ mod 10 = 5.  In order for $ \del$ to be imaginary,
$(\zeta_{20} -\zeta_{20}^{-1})^k \zeta_{20}^{2m+ k}$ must be real or imaginary,
since $v^{a}$ is real or imaginary.  When $(\zeta_{20} -\zeta_{20}^{-1})^k
\zeta_{20}^{2m+ k}$ is real or imaginary, it turns out that it can only be
real.  Thus $v^{a}$ must be imaginary and $a$ must be odd.  This leads to the
following result: all imaginary $d$-numbers of conductor 5 are given by
\begin{align}
\label{img5}
 \del = \begin{cases}
 n  (\zeta_{20} -\zeta_{20}^{-1})^{2k} v^{2a+1} \\
 n  \ii (\zeta_{20} -\zeta_{20}^{-1})^{2k+1} v^{2a+1} \\
\end{cases}
\end{align}
We can rewrite \eqref{img5} as
\begin{align}
\del =  n  (\ii \zeta_{20} -\ii \zeta_{20}^{-1})^{k} v^{2a+1}
\ \ \ \text{ or } \ \ \
\del =  n  (\frac{1+\sqrt 5}{2})^{k}  (\sqrt 5)^{a} v, \ \ \ a =0,1
\end{align}
Thus all the imaginary $d$-numbers can be obtain from the real one
by multiply a factor $v$.

When conductor $cnd=7$, all the $d$-numbers are given by
\begin{align}
 \del &= n  (\zeta_{14} -1)^k (\zeta_{14}^3 -1)^l \zeta_{14}^m v^a,\ \ \ \ n,k,l \in \Z, \ \ m=0,1,2,3,4,5,6,7,8,9,10,11,12,13 
\nonumber\\
v &=  \zeta_{14} - \zeta_{14}^{-1}, \ \ \ a=0,1,2,3,4,5.
\end{align}
If $\del$ is real, then
\begin{align}
\del &= n  (2\cos(4 \pi /7)^k (2\cos(4 \pi /7)+1)^l u^a,\ \ \ \ n,k,l \in \Z 
\nonumber\\
u &=  (1-\zeta_{7})(1-\zeta_7^{-1})= 2-2\cos(2 \pi /7), \ \ \ a=0,1,2,3,4,5.
\end{align}
We will use those results in next a few sections.

\section{Rank-12 representation-347}


Among 5288  SL$_2(\mathbb{Z})$ representations that we considered at rank 12,
the 347$^{th}$ representation
is a difficult one for our GAP code to handle.
The 347$^{th}$ representation is given by:\\
$\tilde s$ =$( 0, 0, \frac{2}{5}, \frac{2}{5}, \frac{3}{5}, \frac{3}{5}, \frac{1}{5}, \frac{1}{5}, \frac{1}{5}, \frac{4}{5}, \frac{4}{5}, \frac{4}{5} ) $.
$\tilde \rho(\mathfrak{s}) $ =\\
\begin{align}
 {\tiny \begin{pmatrix}
-\frac{1}{5}s_{5}^{1},
& \frac{1}{5}s_{5}^{2},
& \frac{1+\sqrt{5}}{10}s_{5}^{1},
& 0,
& \frac{3-\sqrt{5}}{10}s_{5}^{1},
& 0,
& 0,
& \frac{\sqrt{2}}{5}s_{5}^{1},
& 0,
& 0,
& \frac{\sqrt{2}}{5}s_{10}^{1},
& 0 \\ 
\frac{1}{5}s_{5}^{2},
& \frac{1}{5}s_{5}^{1},
& \frac{3-\sqrt{5}}{10}s_{5}^{1},
& 0,
& -\frac{1+\sqrt{5}}{10}s_{5}^{1},
& 0,
& 0,
& -\frac{\sqrt{2}}{5}s_{10}^{1},
& 0,
& 0,
& \frac{\sqrt{2}}{5}s_{5}^{1},
& 0 \\ 
\frac{1+\sqrt{5}}{10}s_{5}^{1},
& \frac{3-\sqrt{5}}{10}s_{5}^{1},
& \frac{1}{5}s_{5}^{1},
& 0,
& -\frac{1}{5}s_{5}^{2},
& 0,
& 0,
& \frac{\sqrt{2}}{5}s_{10}^{1},
& 0,
& 0,
& -\frac{\sqrt{2}}{5}s_{5}^{1},
& 0 \\ 
0,
& 0,
& 0,
& -\frac{3-\sqrt{5}}{10}s_{5}^{1},
& 0,
& \frac{1+\sqrt{5}}{10}s_{5}^{1},
& 0,
& 0,
& \frac{\sqrt{3}}{5}s_{5}^{1},
& \frac{\sqrt{3}}{5}s_{10}^{1},
& 0,
& 0 \\ 
\frac{3-\sqrt{5}}{10}s_{5}^{1},
& -\frac{1+\sqrt{5}}{10}s_{5}^{1},
& -\frac{1}{5}s_{5}^{2},
& 0,
& -\frac{1}{5}s_{5}^{1},
& 0,
& 0,
& \frac{\sqrt{2}}{5}s_{5}^{1},
& 0,
& 0,
& \frac{\sqrt{2}}{5}s_{10}^{1},
& 0 \\ 
0,
& 0,
& 0,
& \frac{1+\sqrt{5}}{10}s_{5}^{1},
& 0,
& \frac{3-\sqrt{5}}{10}s_{5}^{1},
& 0,
& 0,
& -\frac{\sqrt{3}}{5}s_{10}^{1},
& \frac{\sqrt{3}}{5}s_{5}^{1},
& 0,
& 0 \\ 
0,
& 0,
& 0,
& 0,
& 0,
& 0,
& -\frac{5-\sqrt{5}}{10}s_{5}^{1},
& 0,
& 0,
& 0,
& 0,
& \frac{\sqrt{5}}{5}s_{5}^{1} \\ 
\frac{\sqrt{2}}{5}s_{5}^{1},
& -\frac{\sqrt{2}}{5}s_{10}^{1},
& \frac{\sqrt{2}}{5}s_{10}^{1},
& 0,
& \frac{\sqrt{2}}{5}s_{5}^{1},
& 0,
& 0,
& -\frac{1}{5}s_{5}^{2},
& 0,
& 0,
& \frac{1}{5}s_{5}^{1},
& 0 \\ 
0,
& 0,
& 0,
& \frac{\sqrt{3}}{5}s_{5}^{1},
& 0,
& -\frac{\sqrt{3}}{5}s_{10}^{1},
& 0,
& 0,
& \frac{1+\sqrt{5}}{10}s_{5}^{1},
& -\frac{3-\sqrt{5}}{10}s_{5}^{1},
& 0,
& 0 \\ 
0,
& 0,
& 0,
& \frac{\sqrt{3}}{5}s_{10}^{1},
& 0,
& \frac{\sqrt{3}}{5}s_{5}^{1},
& 0,
& 0,
& -\frac{3-\sqrt{5}}{10}s_{5}^{1},
& -\frac{1+\sqrt{5}}{10}s_{5}^{1},
& 0,
& 0 \\ 
\frac{\sqrt{2}}{5}s_{10}^{1},
& \frac{\sqrt{2}}{5}s_{5}^{1},
& -\frac{\sqrt{2}}{5}s_{5}^{1},
& 0,
& \frac{\sqrt{2}}{5}s_{10}^{1},
& 0,
& 0,
& \frac{1}{5}s_{5}^{1},
& 0,
& 0,
& \frac{1}{5}s_{5}^{2},
& 0 \\ 
0, & 0, & 0, & 0, & 0, & 0, & \frac{\sqrt{5}}{5}s_{5}^{1}, & 0, & 0, & 0, & 0, & \frac{5-\sqrt{5}}{10}s_{5}^{1} \\ 
\end{pmatrix}
} 
\end{align}
We use an orthogonal matrix $U$ to transform the representation $\t\rho$ into
pMD representation $ \rho_\mathrm{pMD} = U \t\rho  U^\dag $.  We start by
finding all possible $D_{\rho_\mathrm{pMD}(\mathfrak{s})}(\sigma)$'s.  We find
that there is only one possible $D_{\rho_\mathrm{pMD}(\mathfrak{s})}(\sigma)$
(up to equivalence), generated by 
\begin{align}
D_{\rho_\mathrm{pMD}}(\sigma = 2) = 
{\scriptsize
\begin{pmatrix}
0, & -1, & 0, & 0, & 0, & 0, & 0, & 0, & 0, & 0, & 0, & 0 \\ 
1, & 0, & 0, & 0, & 0, & 0, & 0, & 0, & 0, & 0, & 0, & 0 \\ 
0, & 0, & 0, & 0, & 1, & 0, & 0, & 0, & 0, & 0, & 0, & 0 \\ 
0, & 0, & 0, & 0, & 0, & 1, & 0, & 0, & 0, & 0, & 0, & 0 \\ 
0, & 0, & -1, & 0, & 0, & 0, & 0, & 0, & 0, & 0, & 0, & 0 \\ 
0, & 0, & 0, & -1, & 0, & 0, & 0, & 0, & 0, & 0, & 0, & 0 \\ 
0, & 0, & 0, & 0, & 0, & 0, & 0, & 0, & 0, & 1, & 0, & 0 \\ 
0, & 0, & 0, & 0, & 0, & 0, & 0, & 0, & 0, & 0, & 1, & 0 \\ 
0, & 0, & 0, & 0, & 0, & 0, & 0, & 0, & 0, & 0, & 0, & 1 \\ 
0, & 0, & 0, & 0, & 0, & 0, & -1, & 0, & 0, & 0, & 0, & 0 \\ 
0, & 0, & 0, & 0, & 0, & 0, & 0, & -1, & 0, & 0, & 0, & 0 \\ 
0, & 0, & 0, & 0, & 0, & 0, & 0, & 0, & -1, & 0, & 0, & 0 \\ 
\end{pmatrix}
}
\end{align}

There are six possible choices of unit row $u=1,3,5,7,8,9$ (up to equivalence).
We will only describe two cases, $u=1$ and $u=9$ here.  The $u=3,5$ cases are
similar to $u=1$ case.  The $u=7,8$ cases are similar to $u=9$ case.

For the $u=1$ case and at $r$-stage, we have
$\rho_\mathrm{pMD}(\mathfrak{s}) = $
\begin{align}
{\tiny \begin{pmatrix}
{\bf   -\frac{3-\sqrt{5}}{4}s_{5}^{1}r_{168} 
-\frac{1+\sqrt{5}}{4}s_{5}^{1}r_{167} 
 },
& {\bf  \frac{1+\sqrt{5}}{4}s_{5}^{1}r_{168} 
-\frac{3-\sqrt{5}}{4}s_{5}^{1}r_{167} 
 },
& {\bf  \frac{3-\sqrt{5}}{4}s_{5}^{1}r_{136} 
+\frac{1+\sqrt{5}}{4}s_{5}^{1}r_{135} 
 },
& {\bf  \frac{3-\sqrt{5}}{4}s_{5}^{1}r_{128} 
+\frac{1+\sqrt{5}}{4}s_{5}^{1}r_{127} 
 },
&  \cdots\\
{\bf  \frac{1+\sqrt{5}}{4}s_{5}^{1}r_{168} 
-\frac{3-\sqrt{5}}{4}s_{5}^{1}r_{167} 
 },
&  \frac{3-\sqrt{5}}{4}s_{5}^{1}r_{168} 
+\frac{1+\sqrt{5}}{4}s_{5}^{1}r_{167} 
,
&  -\frac{1+\sqrt{5}}{4}s_{5}^{1}r_{136} 
+\frac{3-\sqrt{5}}{4}s_{5}^{1}r_{135} 
,
&  -\frac{1+\sqrt{5}}{4}s_{5}^{1}r_{128} 
+\frac{3-\sqrt{5}}{4}s_{5}^{1}r_{127} 
,
&  \cdots\\
{\bf  \frac{3-\sqrt{5}}{4}s_{5}^{1}r_{136} 
+\frac{1+\sqrt{5}}{4}s_{5}^{1}r_{135} 
 },
&  -\frac{1+\sqrt{5}}{4}s_{5}^{1}r_{136} 
+\frac{3-\sqrt{5}}{4}s_{5}^{1}r_{135} 
,
&  -\frac{2-\sqrt{5}}{15}s_{5}^{1}+\frac{5-\sqrt{5}}{6}s_{5}^{1}r_{120} 
,
& -\frac{5-\sqrt{5}}{2}s_{5}^{1}r_{104} 
,
&  \cdots\\
{\bf  \frac{3-\sqrt{5}}{4}s_{5}^{1}r_{128} 
+\frac{1+\sqrt{5}}{4}s_{5}^{1}r_{127} 
 },
&  -\frac{1+\sqrt{5}}{4}s_{5}^{1}r_{128} 
+\frac{3-\sqrt{5}}{4}s_{5}^{1}r_{127} 
,
& -\frac{5-\sqrt{5}}{2}s_{5}^{1}r_{104} 
,
&  \frac{1+\sqrt{5}}{30}s_{5}^{1}-\frac{5-\sqrt{5}}{6}s_{5}^{1}r_{120} 
,
&  \cdots\\
{\bf  \frac{1+\sqrt{5}}{4}s_{5}^{1}r_{136} 
-\frac{3-\sqrt{5}}{4}s_{5}^{1}r_{135} 
 },
&  \frac{3-\sqrt{5}}{4}s_{5}^{1}r_{136} 
+\frac{1+\sqrt{5}}{4}s_{5}^{1}r_{135} 
,
&  -\frac{3+\sqrt{5}}{30}s_{5}^{1}+\frac{\sqrt{5}}{3}s_{5}^{1}r_{120} 
,
& -\sqrt{5}s_{5}^{1}r_{104} 
,
&  \cdots\\
{\bf  \frac{1+\sqrt{5}}{4}s_{5}^{1}r_{128} 
-\frac{3-\sqrt{5}}{4}s_{5}^{1}r_{127} 
 },
&  \frac{3-\sqrt{5}}{4}s_{5}^{1}r_{128} 
+\frac{1+\sqrt{5}}{4}s_{5}^{1}r_{127} 
,
& -\sqrt{5}s_{5}^{1}r_{104} 
,
&  -\frac{3-\sqrt{5}}{30}s_{5}^{1}-\frac{\sqrt{5}}{3}s_{5}^{1}r_{120} 
,
&  \cdots\\
{\bf  \frac{3-\sqrt{5}}{4}s_{5}^{1}r_{144} 
+\frac{1+\sqrt{5}}{4}s_{5}^{1}r_{143} 
 },
&  -\frac{1+\sqrt{5}}{4}s_{5}^{1}r_{144} 
+\frac{3-\sqrt{5}}{4}s_{5}^{1}r_{143} 
,
&  -\frac{1+\sqrt{5}}{4}s_{5}^{1}r_{16} 
+\frac{3-\sqrt{5}}{4}s_{5}^{1}r_{15} 
,
&  -\frac{3-\sqrt{5}}{4}s_{5}^{1}r_{24} 
-\frac{1+\sqrt{5}}{4}s_{5}^{1}r_{23} 
,
&  \cdots\\
{\bf  \frac{3-\sqrt{5}}{4}s_{5}^{1}r_{152} 
+\frac{1+\sqrt{5}}{4}s_{5}^{1}r_{151} 
 },
&  -\frac{1+\sqrt{5}}{4}s_{5}^{1}r_{152} 
+\frac{3-\sqrt{5}}{4}s_{5}^{1}r_{151} 
,
&  -\frac{1+\sqrt{5}}{4}s_{5}^{1}r_{32} 
+\frac{3-\sqrt{5}}{4}s_{5}^{1}r_{31} 
,
&  -\frac{3-\sqrt{5}}{4}s_{5}^{1}r_{40} 
-\frac{1+\sqrt{5}}{4}s_{5}^{1}r_{39} 
,
&  \cdots\\
{\bf  \frac{3-\sqrt{5}}{4}s_{5}^{1}r_{160} 
+\frac{1+\sqrt{5}}{4}s_{5}^{1}r_{159} 
 },
&  -\frac{1+\sqrt{5}}{4}s_{5}^{1}r_{160} 
+\frac{3-\sqrt{5}}{4}s_{5}^{1}r_{159} 
,
&  -\frac{1+\sqrt{5}}{4}s_{5}^{1}r_{48} 
+\frac{3-\sqrt{5}}{4}s_{5}^{1}r_{47} 
,
&  -\frac{3-\sqrt{5}}{4}s_{5}^{1}r_{56} 
-\frac{1+\sqrt{5}}{4}s_{5}^{1}r_{55} 
,
&  \cdots\\
{\bf  \frac{1+\sqrt{5}}{4}s_{5}^{1}r_{144} 
-\frac{3-\sqrt{5}}{4}s_{5}^{1}r_{143} 
 },
&  \frac{3-\sqrt{5}}{4}s_{5}^{1}r_{144} 
+\frac{1+\sqrt{5}}{4}s_{5}^{1}r_{143} 
,
&  \frac{3-\sqrt{5}}{4}s_{5}^{1}r_{16} 
+\frac{1+\sqrt{5}}{4}s_{5}^{1}r_{15} 
,
&  -\frac{1+\sqrt{5}}{4}s_{5}^{1}r_{24} 
+\frac{3-\sqrt{5}}{4}s_{5}^{1}r_{23} 
,
&  \cdots\\
{\bf  \frac{1+\sqrt{5}}{4}s_{5}^{1}r_{152} 
-\frac{3-\sqrt{5}}{4}s_{5}^{1}r_{151} 
 },
&  \frac{3-\sqrt{5}}{4}s_{5}^{1}r_{152} 
+\frac{1+\sqrt{5}}{4}s_{5}^{1}r_{151} 
,
&  \frac{3-\sqrt{5}}{4}s_{5}^{1}r_{32} 
+\frac{1+\sqrt{5}}{4}s_{5}^{1}r_{31} 
,
&  -\frac{1+\sqrt{5}}{4}s_{5}^{1}r_{40} 
+\frac{3-\sqrt{5}}{4}s_{5}^{1}r_{39} 
,
&  \cdots\\
{\bf  \frac{1+\sqrt{5}}{4}s_{5}^{1}r_{160} 
-\frac{3-\sqrt{5}}{4}s_{5}^{1}r_{159} 
 },
&  \frac{3-\sqrt{5}}{4}s_{5}^{1}r_{160} 
+\frac{1+\sqrt{5}}{4}s_{5}^{1}r_{159} 
,
&  \frac{3-\sqrt{5}}{4}s_{5}^{1}r_{48} 
+\frac{1+\sqrt{5}}{4}s_{5}^{1}r_{47} 
,
&  -\frac{1+\sqrt{5}}{4}s_{5}^{1}r_{56} 
+\frac{3-\sqrt{5}}{4}s_{5}^{1}r_{55} 
,
&  \cdots\\
\end{pmatrix}
}  
\end{align} 
where only column $1,2,3,4$ are displayed.  The Galois conjugation and the
action of $D_{\rho_\mathrm{pMD}}(\sigma)$ on $\rho_\mathrm{pMD}(\mathfrak{s})$
simply $\rho_\mathrm{pMD}(\mathfrak{s})$ to have the above form.  The
SL$_2(\mathbb{Z})$ conditions, the simple SL$_2(\mathbb{Z})$ character
conditions \eqref{matchKar}, and the orthogonality conditions of
$\rho_\mathrm{pMD}(\mathfrak{s})$ lead to many zero conditions on the
$r$-variables. But those zero conditions are not enough to determine the
$r$-variables.

The make progress, we note that $1/\rho_\mathrm{pMD}(\mathfrak{s})_{uu}$ is a
$d$-number with conductor 5. The number of $d$-numbers of conductor 5 is
infinite. If $1/\rho_\mathrm{pMD}(\mathfrak{s})_{uu}$ contained only one
$r$-variable, the number of $d$-numbers of such a single-$r$ form could be
finite (see Section \ref{tSuucnd5}).  However, for our case,
$1/\rho_\mathrm{pMD}(\mathfrak{s})_{uu}$ contains two $r$-variables,  the
number of $d$-numbers of such a two-$r$ form turn out to be infinite.

To make the infinite possibility finite, we manage to isolate a zero condition
$ -\frac{2}{25}+r_{168}^{2} +r_{167}^{2} =0 $ involving the variables in
$\rho_\mathrm{pMD}(\mathfrak{s})_{uu}$.  With this additional zero condition,
we find, through an extensive search, that in order for the two variables
$r_{168}$ and $r_{167}$ to satisfy $ -\frac{2}{25}+r_{168}^{2} +r_{167}^{2} =0
$ and to make $1/\rho_\mathrm{pMD}(\mathfrak{s})_{uu}$ a $d$-number, they can
only take the following eight sets of possible values
\begin{align}
 (r_{168},r_{167})
=\ \ &  (   -7/25 ,\ \   1/25  ),
&& (   -1/5 ,\ \   -1/5  ),
&& (   -1/5 ,\ \   1/5  ), 
&& (   -1/25 ,\ \   -7/25  ),
\nonumber\\
& (   1/25 ,\ \   7/25  ),
&& (   1/5 ,\ \   -1/5  ), 
&& (   1/5 ,\ \   1/5  ),
&& (   7/25 ,\ \   -1/25  ).
\end{align}
This allows us to determine $\rho_\mathrm{pMD}(\mathfrak{s})_{uu}$.  We stress
that such a result is a conjecture at this stage, because our extensive search
is still a finte one.  To see how extensive is our search, a complicated
$(r_{168},r_{167})$ in our search is given by
\begin{align}
(r_{168},r_{167}) =
\Big( \frac{225851433717}{390625},  \frac{956722026041}{390625} \Big)
\end{align}

Once the value of $\rho_\mathrm{pMD}(\mathfrak{s})_{uu}$ is known, we can
construct many inverse-pairs of integer conditions as decribed in Section
\ref{invpair}.  We find that those inverse-pairs of integer conditions all lead
to contradictions.  So the above eight solutions are all rejected.

 
Next, we consider the $u=9$ case. We have
 $\rho_\mathrm{pMD}(\mathfrak{s})$ = 
\begin{align}
{\tiny \begin{pmatrix}
\cdots &
&  \frac{3-\sqrt{5}}{4}s_{5}^{1}r_{64} 
+\frac{1+\sqrt{5}}{4}s_{5}^{1}r_{63} 
,
&  \frac{3-\sqrt{5}}{4}s_{5}^{1}r_{56} 
+\frac{1+\sqrt{5}}{4}s_{5}^{1}r_{55} 
,
& \cdots &
& {\bf  \frac{3-\sqrt{5}}{4}s_{5}^{1}r_{160} 
+\frac{1+\sqrt{5}}{4}s_{5}^{1}r_{159} 
 },
&\cdots
 \\ 
\cdots &
&  -\frac{1+\sqrt{5}}{4}s_{5}^{1}r_{64} 
+\frac{3-\sqrt{5}}{4}s_{5}^{1}r_{63} 
,
&  -\frac{1+\sqrt{5}}{4}s_{5}^{1}r_{56} 
+\frac{3-\sqrt{5}}{4}s_{5}^{1}r_{55} 
,
& \cdots &
& {\bf  -\frac{1+\sqrt{5}}{4}s_{5}^{1}r_{160} 
+\frac{3-\sqrt{5}}{4}s_{5}^{1}r_{159} 
 },
&\cdots
 \\ 
\cdots &
&  -\frac{3-\sqrt{5}}{10}s_{5}^{1}-\frac{5-\sqrt{5}}{2}s_{5}^{1}r_{120} 
,
& \frac{5-\sqrt{5}}{6}s_{5}^{1}r_{104} 
,
& \cdots &
& {\bf  -\frac{1+\sqrt{5}}{4}s_{5}^{1}r_{128} 
+\frac{3-\sqrt{5}}{4}s_{5}^{1}r_{127} 
 },
&\cdots
 \\ 
\cdots &
& \frac{5-\sqrt{5}}{6}s_{5}^{1}r_{104} 
,
&  \frac{1}{5}s_{5}^{1}+\frac{5-\sqrt{5}}{2}s_{5}^{1}r_{120} 
,
& \cdots &
& {\bf  -\frac{3-\sqrt{5}}{4}s_{5}^{1}r_{136} 
-\frac{1+\sqrt{5}}{4}s_{5}^{1}r_{135} 
 },
&\cdots
 \\ 
\cdots &
&  -\frac{1+\sqrt{5}}{10}s_{5}^{1}-\sqrt{5}s_{5}^{1}r_{120} 
,
& \frac{\sqrt{5}}{3}s_{5}^{1}r_{104} 
,
& \cdots &
& {\bf  \frac{3-\sqrt{5}}{4}s_{5}^{1}r_{128} 
+\frac{1+\sqrt{5}}{4}s_{5}^{1}r_{127} 
 },
&\cdots
 \\ 
\cdots &
& \frac{\sqrt{5}}{3}s_{5}^{1}r_{104} 
,
&  \frac{1}{5}s_{5}^{2}+\sqrt{5}s_{5}^{1}r_{120} 
,
& \cdots &
& {\bf  -\frac{1+\sqrt{5}}{4}s_{5}^{1}r_{136} 
+\frac{3-\sqrt{5}}{4}s_{5}^{1}r_{135} 
 },
&\cdots
 \\ 
\cdots &
&  -\frac{1+\sqrt{5}}{4}s_{5}^{1}r_{24} 
+\frac{3-\sqrt{5}}{4}s_{5}^{1}r_{23} 
,
&  -\frac{1+\sqrt{5}}{4}s_{5}^{1}r_{28} 
+\frac{3-\sqrt{5}}{4}s_{5}^{1}r_{27} 
,
& \cdots &
& {\bf  -\frac{3-\sqrt{5}}{4}s_{5}^{1}r_{144} 
-\frac{1+\sqrt{5}}{4}s_{5}^{1}r_{143} 
 },
&\cdots
 \\ 
\cdots &
&  -\frac{1+\sqrt{5}}{4}s_{5}^{1}r_{40} 
+\frac{3-\sqrt{5}}{4}s_{5}^{1}r_{39} 
,
&  -\frac{3-\sqrt{5}}{4}s_{5}^{1}r_{48} 
-\frac{1+\sqrt{5}}{4}s_{5}^{1}r_{47} 
,
& \cdots &
& {\bf  -\frac{3-\sqrt{5}}{4}s_{5}^{1}r_{152} 
-\frac{1+\sqrt{5}}{4}s_{5}^{1}r_{151} 
 },
&\cdots
 \\ 
\cdots &
& {\bf  -\frac{1+\sqrt{5}}{4}s_{5}^{1}r_{128} 
+\frac{3-\sqrt{5}}{4}s_{5}^{1}r_{127} 
 },
& {\bf  -\frac{3-\sqrt{5}}{4}s_{5}^{1}r_{136} 
-\frac{1+\sqrt{5}}{4}s_{5}^{1}r_{135} 
 },
& \cdots &
& {\bf   -\frac{3-\sqrt{5}}{4}s_{5}^{1}r_{168} 
-\frac{1+\sqrt{5}}{4}s_{5}^{1}r_{167} 
 },
&\cdots
\\ 
\cdots &
&  \frac{3-\sqrt{5}}{4}s_{5}^{1}r_{24} 
+\frac{1+\sqrt{5}}{4}s_{5}^{1}r_{23} 
,
&  \frac{3-\sqrt{5}}{4}s_{5}^{1}r_{28} 
+\frac{1+\sqrt{5}}{4}s_{5}^{1}r_{27} 
,
& \cdots &
& {\bf  -\frac{1+\sqrt{5}}{4}s_{5}^{1}r_{144} 
+\frac{3-\sqrt{5}}{4}s_{5}^{1}r_{143} 
 },
&\cdots
 \\ 
\cdots &
&  \frac{3-\sqrt{5}}{4}s_{5}^{1}r_{40} 
+\frac{1+\sqrt{5}}{4}s_{5}^{1}r_{39} 
,
&  -\frac{1+\sqrt{5}}{4}s_{5}^{1}r_{48} 
+\frac{3-\sqrt{5}}{4}s_{5}^{1}r_{47} 
,
& \cdots &
& {\bf  -\frac{1+\sqrt{5}}{4}s_{5}^{1}r_{152} 
+\frac{3-\sqrt{5}}{4}s_{5}^{1}r_{151} 
 },
&\cdots
 \\ 
\cdots &
&  \frac{3-\sqrt{5}}{4}s_{5}^{1}r_{128} 
+\frac{1+\sqrt{5}}{4}s_{5}^{1}r_{127} 
,
&  -\frac{1+\sqrt{5}}{4}s_{5}^{1}r_{136} 
+\frac{3-\sqrt{5}}{4}s_{5}^{1}r_{135} 
,
& \cdots &
& {\bf  -\frac{1+\sqrt{5}}{4}s_{5}^{1}r_{168} 
+\frac{3-\sqrt{5}}{4}s_{5}^{1}r_{167} 
 },
&\cdots
 \\ 
\end{pmatrix}
}
\end{align}
where only column $3,4,9$ are displayed.
$1/\rho_\mathrm{pMD}(\mathfrak{s})_{ui}$ are $d$-number of conductor 5 for all
$i$'s.  The number of possible values of
$1/\rho_\mathrm{pMD}(\mathfrak{s})_{uu}$ is infinite.  To make the  number of
possible values finite, we need to isolate zero condition for the variables
$r_{167}$ and $r_{168}$ in $\rho_\mathrm{pMD}(\mathfrak{s})_{uu}$.
But we fail to isolate such zero condition.

To make progress, we consider variables in
$\rho_\mathrm{pMD}(\mathfrak{s})_{ij}$ for $i,j \in \{3,4,9\}$, which are $(
r_{104} $, $ r_{120} $, $ r_{127} $, $ r_{128} $, $ r_{135} $, $ r_{136} $, $
r_{167} $, $ r_{168} ) $.  We isolate the following zero conditions for those
variables:
\begin{align}
\label{zcASuuSuiSuj}
& \frac{9}{5}r_{120} +9r_{120}^{2} +r_{104}^{2} = 0 \text{\ \ and \ \ } 
 -\frac{2}{25}+\frac{1}{10}r_{168} 
+\frac{1}{10}r_{167} 
+r_{135}^{2} 
+r_{128}^{2} 
 = 0 \text{\ \ and \ \ } 
\nonumber\\
&
 -\frac{2}{25}+\frac{1}{10}r_{168} 
+\frac{1}{10}r_{167} 
+r_{136}^{2} 
+r_{127}^{2} 
 = 0 \text{\ \ and \ \ } 
 \frac{3}{10}r_{128} 
+\frac{3}{10}r_{127} 
+3r_{120} r_{127} 
+r_{104} r_{136} 
 = 0 \text{\ \ and \ \ } 
\nonumber\\
&
 -\frac{3}{10}r_{128} 
-\frac{3}{10}r_{127} 
-3r_{120} r_{128} 
+r_{104} r_{135} 
 = 0 \text{\ \ and \ \ } 
 \frac{3}{10}r_{136} 
+\frac{3}{10}r_{135} 
+3r_{120} r_{135} 
+r_{104} r_{128} 
 = 0 \text{\ \ and \ \ } 
\nonumber\\
&
 -\frac{3}{10}r_{136} 
-\frac{3}{10}r_{135} 
-3r_{120} r_{136} 
+r_{104} r_{127} 
 = 0 \text{\ \ and \ \ } 
 \frac{1}{10}r_{168} 
-\frac{3}{10}r_{167} 
-r_{135} r_{136} 
+r_{127} r_{128} 
\end{align} 
Also the following expressions of those variables must be cyclotomic integers:
\begin{align}
\label{icycA}
&\frac{1}{ -\frac{3-\sqrt{5}}{4}s_{5}^{1}r_{168} 
-\frac{1+\sqrt{5}}{4}s_{5}^{1}r_{167} 
}
 = \text{cyc-int},\ \ \ 
\frac{1}{ -\frac{1+\sqrt{5}}{4}s_{5}^{1}r_{128} 
+\frac{3-\sqrt{5}}{4}s_{5}^{1}r_{127} 
}
 = \text{cyc-int},
\nonumber\\
&\frac{ -\frac{1+\sqrt{5}}{4}s_{5}^{1}r_{128} 
+\frac{3-\sqrt{5}}{4}s_{5}^{1}r_{127} 
}{ -\frac{3-\sqrt{5}}{4}s_{5}^{1}r_{168} 
-\frac{1+\sqrt{5}}{4}s_{5}^{1}r_{167} 
}
 = \text{cyc-int},\ \ \ 
\frac{ -\frac{3-\sqrt{5}}{10}s_{5}^{1}-\frac{5-\sqrt{5}}{2}s_{5}^{1}r_{120} 
}{ -\frac{3-\sqrt{5}}{4}s_{5}^{1}r_{168} 
-\frac{1+\sqrt{5}}{4}s_{5}^{1}r_{167} 
}
 = \text{cyc-int},
\nonumber\\
&\frac{1}{ -\frac{3-\sqrt{5}}{4}s_{5}^{1}r_{136} 
-\frac{1+\sqrt{5}}{4}s_{5}^{1}r_{135} 
}
 = \text{cyc-int},\ \ \ 
\frac{ -\frac{3-\sqrt{5}}{4}s_{5}^{1}r_{136} 
-\frac{1+\sqrt{5}}{4}s_{5}^{1}r_{135} 
}{ -\frac{3-\sqrt{5}}{4}s_{5}^{1}r_{168} 
-\frac{1+\sqrt{5}}{4}s_{5}^{1}r_{167} 
}
 = \text{cyc-int},
\nonumber\\
&\frac{ \frac{1}{5}s_{5}^{1}+\frac{5-\sqrt{5}}{2}s_{5}^{1}r_{120} 
}{ -\frac{3-\sqrt{5}}{4}s_{5}^{1}r_{168} 
-\frac{1+\sqrt{5}}{4}s_{5}^{1}r_{167} 
}
 = \text{cyc-int},\ \ \ 
\frac{\frac{5-\sqrt{5}}{6}s_{5}^{1}r_{104} 
}{ -\frac{3-\sqrt{5}}{4}s_{5}^{1}r_{168} 
-\frac{1+\sqrt{5}}{4}s_{5}^{1}r_{167} 
}
 = \text{cyc-int},
\nonumber\\
&\frac{\frac{5-\sqrt{5}}{6}s_{5}^{1}r_{104} 
}{ -\frac{3-\sqrt{5}}{4}s_{5}^{1}r_{136} 
-\frac{1+\sqrt{5}}{4}s_{5}^{1}r_{135} 
}
 = \text{cyc-int},\ \ \ 
\frac{\frac{5-\sqrt{5}}{6}s_{5}^{1}r_{104} 
}{ -\frac{1+\sqrt{5}}{4}s_{5}^{1}r_{128} 
+\frac{3-\sqrt{5}}{4}s_{5}^{1}r_{127} 
}
= \text{cyc-int},
\end{align}

Although the number of possible values for the three inverse $d$-numbers,
$\rho_\mathrm{pMD}(\mathfrak{s})_{ui}$, $i=3,4,9$, is infinite, through an
extensive research, we conjecture that none of the three combined $d$-numbers
can satisfy the two sets of conditions \eqref{zcASuuSuiSuj} and \eqref{icycA}.
We note that, once we know $ r_{127} $, $ r_{128} $, $ r_{135} $, $ r_{136} $,
$ r_{167} $, $ r_{168}  $, we obtain a linear relation betwen $ r_{104} $ and $
r_{120} $, from the zero conditions \eqref{zcASuuSuiSuj}.  We find that such a
linear relation, together with the qudratic relation $\frac{9}{5}r_{120}
+9r_{120}^{2} +r_{104}^{2} = 0 $, only leads to a few solutions of $ r_{104} $
and $ r_{120} $. Those solutions are rejected by the cyclotomic-integer
conditions \eqref{icycA}.  To see how extensive is our search, a complicated
$(r_{168},r_{167})$ in our search is given by
\begin{align}
(r_{168},r_{167}) =
\Big( -\frac{3571}{78125},  -\frac{15127}{78125} \Big)
\end{align}

\section{Rank-12 representation-333}


The rank 12, 333$^{th}$ representation
is also hard for our GAP code.
The 333$^{th}$ representation is given by:
$\tilde s$ =$( \frac{1}{5}, \frac{1}{5}, \frac{2}{5}, \frac{2}{5}, \frac{3}{5}, \frac{3}{5}, \frac{4}{5}, \frac{4}{5}, 0, 0, 0, 0 ) $ and
 $\tilde \rho(\mathfrak{s}) $ = 
\begin{align}
{\tiny \begin{pmatrix}
\frac{3-\sqrt{5}}{10}, & 0, & -\frac{1+\sqrt{5}}{5}, & 0, & -\frac{1-\sqrt{5}}{5}, & 0, & -\frac{3+\sqrt{5}}{10}, & 0, & \frac{\sqrt{6}}{5}, & 0, & 0, & 0 \\ 
0, & -\frac{5+\sqrt{5}}{10}, & 0, & 0, & 0, & 0, & 0, & -\frac{5-\sqrt{5}}{10}, & 0, & 0, & 0, & -\frac{\sqrt{10}}{5} \\ 
-\frac{1+\sqrt{5}}{5}, & 0, & \frac{3+\sqrt{5}}{10}, & 0, & \frac{3-\sqrt{5}}{10}, & 0, & \frac{1-\sqrt{5}}{5}, & 0, & \frac{\sqrt{6}}{5}, & 0, & 0, & 0 \\ 
0, & 0, & 0, & -\frac{5-\sqrt{5}}{10}, & 0, & -\frac{5+\sqrt{5}}{10}, & 0, & 0, & 0, & 0, & -\frac{\sqrt{10}}{5}, & 0 \\ 
-\frac{1-\sqrt{5}}{5}, & 0, & \frac{3-\sqrt{5}}{10}, & 0, & \frac{3+\sqrt{5}}{10}, & 0, & \frac{1+\sqrt{5}}{5}, & 0, & \frac{\sqrt{6}}{5}, & 0, & 0, & 0 \\ 
0, & 0, & 0, & -\frac{5+\sqrt{5}}{10}, & 0, & -\frac{5-\sqrt{5}}{10}, & 0, & 0, & 0, & 0, & \frac{\sqrt{10}}{5}, & 0 \\ 
-\frac{3+\sqrt{5}}{10}, & 0, & \frac{1-\sqrt{5}}{5}, & 0, & \frac{1+\sqrt{5}}{5}, & 0, & \frac{3-\sqrt{5}}{10}, & 0, & -\frac{\sqrt{6}}{5}, & 0, & 0, & 0 \\ 
0, & -\frac{5-\sqrt{5}}{10}, & 0, & 0, & 0, & 0, & 0, & -\frac{5+\sqrt{5}}{10}, & 0, & 0, & 0, & \frac{\sqrt{10}}{5} \\ 
\frac{\sqrt{6}}{5}, & 0, & \frac{\sqrt{6}}{5}, & 0, & \frac{\sqrt{6}}{5}, & 0, & -\frac{\sqrt{6}}{5}, & 0, & -\frac{1}{5}, & 0, & 0, & 0 \\ 
0, & 0, & 0, & 0, & 0, & 0, & 0, & 0, & 0, & 1, & 0, & 0 \\ 
0, & 0, & 0, & -\frac{\sqrt{10}}{5}, & 0, & \frac{\sqrt{10}}{5}, & 0, & 0, & 0, & 0, & -\frac{\sqrt{5}}{5}, & 0 \\ 
0, & -\frac{\sqrt{10}}{5}, & 0, & 0, & 0, & 0, & 0, & \frac{\sqrt{10}}{5}, & 0, & 0, & 0, & \frac{\sqrt{5}}{5} \\ 
\end{pmatrix}
} 
\end{align}
At $r$-stage, there are 117 cases from different possible
$D_{\rho_{pMD}}(\si)$'s, different possible unit index $u$ and different
or-connected zero conditions.  For some cases,
$\rho_\mathrm{pMD}(\mathfrak{s})_{uu}$ contain only one $r$ variable.  In this
case, $1/\rho_\mathrm{pMD}(\mathfrak{s})_{uu}$ can only take a finite number of
possible values in order for it to be a $d$-number of conductor 5 (see Section
\ref{tSuucnd5}).

For some other cases, $\rho_\mathrm{pMD}(\mathfrak{s})_{uu}$ contain two $r$
variables.  We will discuss one such case here. Other two-variable cases are
similar.  One of the two-variable case has $\rho_\mathrm{pMD}(\mathfrak{s})$ =
\\
\begin{align}
{\tiny \begin{pmatrix}
 -\frac{3+\sqrt{5}}{10}+r_{109} 
,
& -r_{57} 
,
&\cdots
& {\bf  \frac{1-\sqrt{5}}{2}r_{144} 
-r_{141} 
 },
&  \frac{1+\sqrt{5}}{2}r_{144} 
-r_{141} 
,
&  \frac{1-\sqrt{5}}{2}r_{8} 
-r_{5} 
,
&  \frac{1+\sqrt{5}}{2}r_{8} 
-r_{5} 
 \\ 
-r_{57} 
,
&  \frac{1-\sqrt{5}}{10}-r_{109} 
,
&\cdots
& {\bf  \frac{1+\sqrt{5}}{2}r_{136} 
-r_{133} 
 },
&  \frac{1-\sqrt{5}}{2}r_{136} 
-r_{133} 
,
&  -\frac{1+\sqrt{5}}{2}r_{32} 
+r_{29} 
,
&  -\frac{1-\sqrt{5}}{2}r_{32} 
+r_{29} 
 \\ 
\frac{1+\sqrt{5}}{2}r_{56} 
,
& -\frac{1+\sqrt{5}}{2}r_{80} 
,
&\cdots
& {\bf  -\frac{1-\sqrt{5}}{2}r_{128} 
+r_{125} 
 },
&  -\frac{1+\sqrt{5}}{2}r_{128} 
+r_{125} 
,
&  \frac{1-\sqrt{5}}{2}r_{28} 
-r_{25} 
,
&  \frac{1+\sqrt{5}}{2}r_{28} 
-r_{25} 
 \\ 
-\frac{1+\sqrt{5}}{2}r_{68} 
,
& \frac{1+\sqrt{5}}{2}r_{76} 
,
&\cdots
& {\bf  \frac{1+\sqrt{5}}{2}r_{120} 
-r_{117} 
 },
&  \frac{1-\sqrt{5}}{2}r_{120} 
-r_{117} 
,
&  -\frac{1+\sqrt{5}}{2}r_{20} 
+r_{17} 
,
&  -\frac{1-\sqrt{5}}{2}r_{20} 
+r_{17} 
 \\ 
\frac{1-\sqrt{5}}{2}r_{56} 
,
& -\frac{1-\sqrt{5}}{2}r_{80} 
,
&\cdots
& {\bf  -\frac{1+\sqrt{5}}{2}r_{128} 
+r_{125} 
 },
&  -\frac{1-\sqrt{5}}{2}r_{128} 
+r_{125} 
,
&  \frac{1+\sqrt{5}}{2}r_{28} 
-r_{25} 
,
&  \frac{1-\sqrt{5}}{2}r_{28} 
-r_{25} 
 \\ 
-\frac{1-\sqrt{5}}{2}r_{68} 
,
& \frac{1-\sqrt{5}}{2}r_{76} 
,
&\cdots
& {\bf  \frac{1-\sqrt{5}}{2}r_{120} 
-r_{117} 
 },
&  \frac{1+\sqrt{5}}{2}r_{120} 
-r_{117} 
,
&  -\frac{1-\sqrt{5}}{2}r_{20} 
+r_{17} 
,
&  -\frac{1+\sqrt{5}}{2}r_{20} 
+r_{17} 
 \\ 
 -\frac{3-\sqrt{5}}{10}+r_{109} 
,
& -r_{57} 
,
&\cdots
& {\bf  \frac{1+\sqrt{5}}{2}r_{144} 
-r_{141} 
 },
&  \frac{1-\sqrt{5}}{2}r_{144} 
-r_{141} 
,
&  \frac{1+\sqrt{5}}{2}r_{8} 
-r_{5} 
,
&  \frac{1-\sqrt{5}}{2}r_{8} 
-r_{5} 
 \\ 
r_{57} 
,
&  -\frac{1+\sqrt{5}}{10}+r_{109} 
,
&\cdots
& {\bf  -\frac{1-\sqrt{5}}{2}r_{136} 
+r_{133} 
 },
&  -\frac{1+\sqrt{5}}{2}r_{136} 
+r_{133} 
,
&  \frac{1-\sqrt{5}}{2}r_{32} 
-r_{29} 
,
&  \frac{1+\sqrt{5}}{2}r_{32} 
-r_{29} 
 \\ 
{\bf  \frac{1-\sqrt{5}}{2}r_{144} 
-r_{141} 
 },
& {\bf  \frac{1+\sqrt{5}}{2}r_{136} 
-r_{133} 
 },
&\cdots
& {\bf  \frac{1+\sqrt{5}}{2}r_{160} 
-r_{157} 
 },
& {\bf  \frac{1-\sqrt{5}}{2}r_{160} 
-r_{157} 
 },
& {\bf  -\frac{1-\sqrt{5}}{2}r_{152} 
+r_{149} 
 },
& {\bf  -\frac{1+\sqrt{5}}{2}r_{152} 
+r_{149} 
 } \\ 
 \frac{1+\sqrt{5}}{2}r_{144} 
-r_{141} 
,
&  \frac{1-\sqrt{5}}{2}r_{136} 
-r_{133} 
,
&\cdots
& {\bf  \frac{1-\sqrt{5}}{2}r_{160} 
-r_{157} 
 },
&  \frac{1+\sqrt{5}}{2}r_{160} 
-r_{157} 
,
&  -\frac{1+\sqrt{5}}{2}r_{152} 
+r_{149} 
,
&  -\frac{1-\sqrt{5}}{2}r_{152} 
+r_{149} 
 \\ 
 \frac{1-\sqrt{5}}{2}r_{8} 
-r_{5} 
,
&  -\frac{1+\sqrt{5}}{2}r_{32} 
+r_{29} 
,
&\cdots
& {\bf  -\frac{1-\sqrt{5}}{2}r_{152} 
+r_{149} 
 },
&  -\frac{1+\sqrt{5}}{2}r_{152} 
+r_{149} 
,
&  \frac{2}{5}-\frac{1+\sqrt{5}}{2}r_{160} 
+r_{157} 
,
&  \frac{2}{5}-\frac{1-\sqrt{5}}{2}r_{160} 
+r_{157} 
 \\ 
 \frac{1+\sqrt{5}}{2}r_{8} 
-r_{5} 
,
&  -\frac{1-\sqrt{5}}{2}r_{32} 
+r_{29} 
,
&\cdots
& {\bf  -\frac{1+\sqrt{5}}{2}r_{152} 
+r_{149} 
 },
&  -\frac{1-\sqrt{5}}{2}r_{152} 
+r_{149} 
,
&  \frac{2}{5}-\frac{1-\sqrt{5}}{2}r_{160} 
+r_{157} 
,
&  \frac{2}{5}-\frac{1+\sqrt{5}}{2}r_{160} 
+r_{157} 
 \\ 
\end{pmatrix}
}  
\end{align}
where only column $1,2,9,10,11,12$ are displayed.
$1/\rho_\mathrm{pMD}(\mathfrak{s})_{ui}$ are $d$-number of conductor 5 for all
$i$'s.  The number of possible values of
$1/\rho_\mathrm{pMD}(\mathfrak{s})_{uu}$ is infinite.  To make the  number of
possible values finite, we need to isolate zero condition for the variables
$r_{157}$ and $r_{160}$ in $\rho_\mathrm{pMD}(\mathfrak{s})_{uu}$.
But we fail to isolate such zero condition.

To make progress, consider variables in $\rho_\mathrm{pMD}(\mathfrak{s})_{ij}$
for $i,j \in \{1,2,9\}$, which are $ r_{57} $, $ r_{109} $, $ r_{133} $, $
r_{136} $, $ r_{141} $, $ r_{144} $, $ r_{157} $, $ r_{160}  $. 
We isolate
the following zero conditions for those variables:
\begin{align}
\label{zcASuuSuiSuj1}
& -\frac{3}{25}-\frac{2}{5}r_{109} 
+r_{109}^{2} 
+r_{57}^{2} 
  = 0 \text{\ \ and \ \ } 
 -\frac{2}{25}-\frac{2}{5}r_{160} 
+r_{144}^{2} 
+r_{136}^{2} 
  = 0 \text{\ \ and \ \ } 
\nonumber\\ &
 -\frac{3}{25}-\frac{3}{5}r_{160} 
+\frac{1}{5}r_{109} 
+2r_{144}^{2} 
+r_{109} r_{160} 
  = 0 \text{\ \ and \ \ } 
 \frac{3}{5}r_{136} 
-r_{109} r_{136} 
+r_{57} r_{144} 
  = 0 \text{\ \ and \ \ } 
\nonumber\\ &
 \frac{1}{5}r_{144} 
+r_{109} r_{144} 
+r_{57} r_{136} 
  = 0 \text{\ \ and \ \ } 
 -\frac{3}{25}-\frac{1}{5}r_{157} 
+r_{144}^{2} 
-2r_{141} r_{144} 
+r_{141}^{2} 
+r_{133}^{2} 
  = 0 \text{\ \ and \ \ } 
\nonumber\\ &
 -\frac{3}{25}r_{136} 
-\frac{3}{5}r_{136} r_{160} 
+\frac{1}{5}r_{109} r_{136} 
+2r_{136} r_{144}^{2} 
+r_{109} r_{136} r_{160} 
  = 0 \text{\ \ and \ \ } 
  = 0 \text{\ \ and \ \ } 
\nonumber\\ &
 -\frac{1}{25}-\frac{1}{5}r_{160} 
+r_{144}^{2} 
-r_{141} r_{144} 
+r_{133} r_{136} 
  = 0 \text{\ \ and \ \ } 
\nonumber\\ &
 -\frac{1}{5}r_{144} 
+\frac{3}{5}r_{141} 
+r_{109} r_{144} 
-r_{109} r_{141} 
+r_{57} r_{133} 
  = 0 \text{\ \ and \ \ } 
\nonumber\\ &
 \frac{1}{5}r_{136} 
+\frac{1}{5}r_{133} 
-r_{109} r_{136} 
+r_{109} r_{133} 
+r_{57} r_{141} 
\end{align} 
Although the number of possible values for the three inverse $d$-numbers,
$\rho_\mathrm{pMD}(\mathfrak{s})_{ui}$, $i=1,2,9$, is infinite, through an
extensive research, we conjecture that there are only 12 sets of solutions can
satisfy the zero conditions \eqref{zcASuuSuiSuj1}:
\begin{align}
& (r[160], r[157]  , r[144]  , r[141]  , r[136] , && \hskip -1.7ex r[133]  , r[109] , r[57] ) 
\\ 
=\ \ &
(  -\frac15 ,  -\frac15 ,  0 ,  -\frac15 ,  0 ,  -\frac15 ,  \frac15 ,  -\frac25 ) ,\ \ 
&&
(  -\frac15 ,  -\frac15 ,  0 ,  -\frac15 ,  0 ,  \frac15 ,  \frac15 ,  \frac25 ) ,\ \ 
&&
(  -\frac15 ,  -\frac15 ,  0 ,  \frac15 ,  0 ,  -\frac15 ,  \frac15 ,  \frac25 ) ,\ \ 
\nonumber\\ &
(  -\frac15 ,  -\frac15 ,  0 ,  \frac15 ,  0 ,  \frac15 ,  \frac15 ,  -\frac25 ) ,\ \ 
&&
(  \frac15 ,  -\frac15 ,  -\frac25 ,  -\frac15 ,  0 ,  -\frac15 ,  0 ) ,\ \ 
&&
(  \frac15 ,  -\frac15 ,  -\frac25 ,  -\frac15 ,  \frac15 ,  -\frac15 ,  0 ) ,\ \ 
\nonumber\\ &
(  \frac15 ,  -\frac15 ,  0 ,  -\frac15 ,  -\frac25 ,  -\frac15 ,  \frac35 ,  0 ) ,\ \ 
&&
(  \frac15 ,  -\frac15 ,  0 ,  -\frac15 ,  \frac25 ,  \frac15 ,  \frac35 ,  0 ) ,\ \ 
&&
(  \frac15 ,  -\frac15 ,  0 ,  \frac15 ,  -\frac25 ,  -\frac15 ,  \frac35 ,  0 ) ,\ \ 
\nonumber\\ &
(  \frac15 ,  -\frac15 ,  0 ,  \frac15 ,  \frac25 ,  \frac15 ,  \frac35 ,  0 ) ,\ \ 
&&
(  \frac15 ,  -\frac15 ,  \frac25 ,  \frac15 ,  0 ,  -\frac15 ,  -\frac15 ,  0 ) ,\ \ 
&&
(  \frac15 ,  -\frac15 ,  \frac25 ,  \frac15 ,  0 ,  \frac15 ,  -\frac15 ,  0 ) 
\nonumber 
\end{align}
Once the value of $\rho_\mathrm{pMD}(\mathfrak{s})_{uu}$ is known, we can
construct many inverse-pairs of integer conditions as decribed in Section
\ref{invpair}.  We find that those inverse-pairs of integer conditions all lead
to contradictions.  So the above 12 solutions are all rejected.

\section{Rank-12 representation-3246}


The rank 12, 3246$^{th}$ representation has $\tilde s$ = $( 0, 0,
\frac{1}{2}, 0, \frac{1}{7}, \frac{2}{7}, \frac{3}{7}, \frac{4}{7},
\frac{5}{7}, \frac{6}{7}, 0, 0 ) $.
At $r$-stage, there is only one case with unit row $u=8$ and
$\rho_\mathrm{pMD}(\mathfrak{s})$ = 
\begin{align}
\hskip -8mm {\tiny \begin{pmatrix}
\cdots &
& {\bf r_{73} 
 },
& r_{73} 
,
& r_{73} 
,
& r_{19} 
,
& r_{25} 
 \\ 
\cdots &
& {\bf  -(\frac{1}{3}-c^{1}_{7}
)r_{72} 
-(\frac{1}{3}-c^{2}_{7}
)r_{71} 
 },
&  -(\frac{1}{3}-c^{2}_{7}
)r_{72} 
+(\frac{2}{3}+c^{1}_{7}
+c^{2}_{7}
)r_{71} 
,
&  (\frac{2}{3}+c^{1}_{7}
+c^{2}_{7}
)r_{72} 
-(\frac{1}{3}-c^{1}_{7}
)r_{71} 
,
& 0,
& 0 \\ 
\cdots &
& {\bf  (\frac{2}{3}+c^{1}_{7}
+c^{2}_{7}
)r_{72} 
-(\frac{1}{3}-c^{1}_{7}
)r_{71} 
 },
&  -(\frac{1}{3}-c^{1}_{7}
)r_{72} 
-(\frac{1}{3}-c^{2}_{7}
)r_{71} 
,
&  -(\frac{1}{3}-c^{2}_{7}
)r_{72} 
+(\frac{2}{3}+c^{1}_{7}
+c^{2}_{7}
)r_{71} 
,
& 0,
& 0 \\ 
\cdots &
& {\bf  (\frac{1}{3}+c^{2}_{7}
)r_{72} 
-(\frac{2}{3}-c^{1}_{7}
-c^{2}_{7}
)r_{71} 
 },
&  -(\frac{2}{3}-c^{1}_{7}
-c^{2}_{7}
)r_{72} 
+(\frac{1}{3}+c^{1}_{7}
)r_{71} 
,
&  (\frac{1}{3}+c^{1}_{7}
)r_{72} 
+(\frac{1}{3}+c^{2}_{7}
)r_{71} 
,
& 0,
& 0 \\ 
\cdots &
& {\bf  (\frac{1}{3}+c^{2}_{7}
)r_{72} 
-(\frac{2}{3}-c^{1}_{7}
-c^{2}_{7}
)r_{71} 
 },
&  -(\frac{2}{3}-c^{1}_{7}
-c^{2}_{7}
)r_{72} 
+(\frac{1}{3}+c^{1}_{7}
)r_{71} 
,
&  (\frac{1}{3}+c^{1}_{7}
)r_{72} 
+(\frac{1}{3}+c^{2}_{7}
)r_{71} 
,
& 0,
& 0 \\ 
\cdots &
& {\bf  -(\frac{2}{3}-c^{1}_{7}
-c^{2}_{7}
)r_{72} 
+(\frac{1}{3}+c^{1}_{7}
)r_{71} 
 },
&  (\frac{1}{3}+c^{1}_{7}
)r_{72} 
+(\frac{1}{3}+c^{2}_{7}
)r_{71} 
,
&  (\frac{1}{3}+c^{2}_{7}
)r_{72} 
-(\frac{2}{3}-c^{1}_{7}
-c^{2}_{7}
)r_{71} 
,
& 0,
& 0 \\ 
\cdots &
& {\bf  (\frac{1}{3}+c^{1}_{7}
)r_{72} 
+(\frac{1}{3}+c^{2}_{7}
)r_{71} 
 },
&  (\frac{1}{3}+c^{2}_{7}
)r_{72} 
-(\frac{2}{3}-c^{1}_{7}
-c^{2}_{7}
)r_{71} 
,
&  -(\frac{2}{3}-c^{1}_{7}
-c^{2}_{7}
)r_{72} 
+(\frac{1}{3}+c^{1}_{7}
)r_{71} 
,
& 0,
& 0 \\ 
\cdots &
& {\bf c_{7}^{2}r_{84} 
+c_{7}^{3}r_{83} 
+r_{79} 
 },
& {\bf  c_{7}^{3}r_{84} 
+c_{7}^{1}r_{83} 
+r_{79} 
 },
& {\bf  c_{7}^{1}r_{84} 
+c_{7}^{2}r_{83} 
+r_{79} 
 },
& {\bf -r_{43} 
 },
& {\bf -r_{49} 
 } \\ 
\cdots &
& {\bf  c_{7}^{3}r_{84} 
+c_{7}^{1}r_{83} 
+r_{79} 
 },
&  c_{7}^{1}r_{84} 
+c_{7}^{2}r_{83} 
+r_{79} 
,
&  c_{7}^{2}r_{84} 
+c_{7}^{3}r_{83} 
+r_{79} 
,
& -r_{43} 
,
& -r_{49} 
 \\ 
\cdots &
& {\bf  c_{7}^{1}r_{84} 
+c_{7}^{2}r_{83} 
+r_{79} 
 },
&  c_{7}^{2}r_{84} 
+c_{7}^{3}r_{83} 
+r_{79} 
,
&  c_{7}^{3}r_{84} 
+c_{7}^{1}r_{83} 
+r_{79} 
,
& -r_{43} 
,
& -r_{49} 
 \\ 
\cdots &
& {\bf -r_{43} 
 },
& -r_{43} 
,
& -r_{43} 
,
&  \frac{3}{2}+r_{84} 
+r_{83} 
-3r_{79} 
+r_{13} 
,
& -r_{1} 
 \\ 
\cdots &
& {\bf -r_{49} 
 },
& -r_{49} 
,
& -r_{49} 
,
& -r_{1} 
,
& -r_{13} 
 \\ 
\end{pmatrix}
}
\end{align}
where only column $8,9,10,11,12$ are displayed.
$1/\rho_\mathrm{pMD}(\mathfrak{s})_{ui}$ are $d$-number of conductor 7 for all
$i$'s.  The number of possible values of
$1/\rho_\mathrm{pMD}(\mathfrak{s})_{uu}$ is infinite.  To make the  number of
possible values finite, we need to isolate zero condition for the variables
$r[79]$, $r_{83}$, and $r_{84}$ in $\rho_\mathrm{pMD}(\mathfrak{s})_{uu}$.
However, the isolated zero conditions are not enough to make number of possible
values finite.

To make progress, we isolate zero condition for the variables $r_{79}$,
$r_{83}$, $r_{84}$, plus one more variable.  We obtain the following zero
conditions:
\begin{align}
\label{zcASuu}
& -\frac{1}{49}+r_{84}^{2} 
-r_{83} r_{84} 
+r_{83}^{2} 
 =0, \ \ \text{ and } \ \ 
 -\frac{1}{49}r_{84} 
-\frac{1}{49}r_{83} 
+r_{84}^{3} 
+r_{83}^{3} 
 =0, \ \ \text{ and } \ \ 
\nonumber\\ &
 -\frac{1}{49}r_{84} 
+r_{84}^{3} 
-r_{83} r_{84}^{2} 
+r_{83}^{2} r_{84} 
 =0, \ \ \text{ and } \ \ 
 -\frac{2}{49}-\frac{1}{7}r_{84} 
+\frac{2}{7}r_{83} 
+r_{72}^{2} 
 =0, \ \ \text{ and } \ \ 
\nonumber\\ &
 -\frac{2}{49}-\frac{1}{7}r_{84} 
-\frac{1}{7}r_{83} 
+r_{71}^{2} 
 =0, \ \ \text{ and } \ \ 
 -\frac{1}{6}-\frac{1}{6}r_{84} 
-\frac{1}{6}r_{83} 
+\frac{1}{2}r_{79} 
+r_{73}^{2} = 0.
\end{align} 
Although the number of possible values for the inverse $d$-number
$\rho_\mathrm{pMD}(\mathfrak{s})_{uu}$ is infinite, through an extensive
search, we conjecture that there are only three sets of solutions can satisfy
the zero conditions \eqref{zcASuu}:
\begin{align}
( r_{84}, r_{83}, r_{79})
& = 
(  -\frac17 ,  0 ,  -\frac{3}{14}  ),\ \ \
(  0 ,  -\frac17 ,  -\frac{3}{14}  ),\ \ \ 
(  \frac17 ,  \frac17 , -\frac{1}{14}  ).
\end{align}
To see how extensive is our search, an complicated triple $( r_{84}, r_{83},
r_{79})$ in our search is given by 
\begin{align}
( r_{84}, r_{83}, r_{79})
 = 
\Big( \frac{127603558175}{100352},
\frac{459867333397}{200704}, -\frac{16225407395}{28672} \Big).
\end{align}

We remark that we search for possible $d$-numbers
$1/\rho_\mathrm{pMD}(\mathfrak{s})_{uu}$, \ie search for different values of
$r_{84}, r_{83}, r_{79}$.  But the zero condition with an extra variable, say,
$-\frac{1}{6}-\frac{1}{6}r_{84} -\frac{1}{6}r_{83} +\frac{1}{2}r_{79}
+r_{73}^{2} = 0$ is still useful.  When $( r_{84}, r_{83}, r_{79}) = (  \frac17
, \frac17 , -\frac{1}{14}  )$, the zero condition
$-\frac{1}{6}-\frac{1}{6}r_{84} -\frac{1}{6}r_{83} +\frac{1}{2}r_{79}
+r_{73}^{2} = -1/4 + r_{73}^{2}$ has solutions for $r_{73} = \pm 1/2$.  But when
$( r_{84}, r_{83}, r_{79}) = (  \frac17 ,  \frac17 , \frac{1}{14}  )$, the zero
condition $-\frac{1}{6}-\frac{1}{6}r_{84} -\frac{1}{6}r_{83} +\frac{1}{2}r_{79}
+r_{73}^{2} = -5/28 + r_{73}^{2}$ has no solution for a rational $r_{73}$.
This is why the zero condition with one extra variable can still reject most
possible values of the $d$-numbers.

Once the value of $\rho_\mathrm{pMD}(\mathfrak{s})_{uu}$ is known, we can
construct many inverse-pairs of integer conditions as descried in Section
\ref{invpair}.  We find that those inverse-pairs of integer conditions all lead
to contradictions.  So the above three solutions are all rejected.

The rank 12 representation-1372 and representation-3251 can be
handled in a similar way.

\section{Condensable algebras, boundaries, and domain walls}
\label{MMA}

Let us use $a,b,c$ to label the simple objects in a modular tensor category
(MTC) $\eM$.  $\eM$ is characterized by modular matrices $\t S_\eM = (\t
S_\eM^{ab})$ and $\t T_\eM = (\t T_\eM^{ab})$, whose indices are labeled by the
simple objects.  $\t S_\eM, \t T_\eM$ are unitary matrices that generate a
representation of $SL(2,\Z_n)$, where $n$ is the smallest integer that satisfy
$\t T_\eM^n =\id$. We call $n$ as the order of $\t T_\eM$ and denote it as
$n=\ord(\t T_\eM)$.  $\t T_\eM$ is a diagonal matrix and $\t S_\eM$ is a
symmetric matrix.

From $\t S_\eM$ and $\t T_\eM$, we define
normalized $S,T$-matrices 
\begin{align}
  S_\eM &= \t S_\eM/\t S_\eM^{\one\one}, &
  T_\eM &= \t T_\eM/\t T_\eM^{\one\one}.
\end{align}
Let $d_a$ be the quantum dimension of simple object $a$, which is given by
$d_a=(S_\eM)^{a\one}$.  Let $s_a$ be the topological spin of simple
object $a$, which is given by $\ee^{\ii 2\pi s_a} = (T_\eM)^{aa}$.  The total
dimension of $\eM$ is defined as $D^2_{\eM} \equiv \sum_{a\in \eM} d_a^2$.
Also let $d_{\cA}$ be the quantum dimension of the condensable algebra $\cA$,
\ie if 
\begin{align}
\cA=\bigoplus_{a\in \eM} A^a a
\end{align}
then $d_\cA =\sum_a A^a d_a$.  
Similarly, we use $i,j,k$ to label the simple object in $\eM_{/\cA}$,
where $\eM_{/\cA}$ is the MTC obtained from $\eM$ by condensing $\cA$.  
Following the
above, we can define $\t S_{\eM_{/\cA}} = (\t S_{\eM_{/\cA}} ^{ij})$, $\t
T_{\eM_{/\cA}} = (\t T_{\eM_{/\cA}} ^{ij})$, $S_{\eM_{/\cA}}$,
$T_{\eM_{/\cA}}$, as well as
$d_i$, $s_i$, and $D^2_{\eM_{/\cA}}$.  Then we have the following properties
\begin{itemize}

\item The distinct $s_i$'s form a subset of $\{s_a \mid a \in \eM\}$.

\item $D_{\eM} = D_{\eM_{/\cA}} d_\cA$.  $D_{\eM}, D_{\eM_{/\cA}}, d_\cA$
are cyclotomic integers.

\item $A^a$ in $\cA$ are non-negative integers, $A^a=A^{\bar a}$, and
$A^\one =1$. 

\item For $a \in \cA$ (\ie for $A^a \neq 0$), the corresponding  $s_a  = 0$ mod
1. 

\item if $a,b \in \cA$, then at least one of the fusion products in $a\otimes b$ must be contained in $\cA$, \ie $ \exists  c\in \cA $ such that
$a\otimes b = c\oplus \cdots$.

\end{itemize}

Now, let us assume $\cA$ to be Lagrangian, then the
$\cA$-condensed boundary of $\eM$ is gapped. 
Let us use $x$ to label the (simple) excitations on the gapped boundary.
If we bring a bulk excitation $a$
to such a boundary and fuse it with a boundary excitation $y$, it will become a (composite) boundary excitation 
\begin{align}
a\otimes y = \bigoplus_{x} M^{ay}_x x, \ \ \ M^{ay}_x \in \N.
\end{align}
Then $A^a$ is given by $A^a =M^a_\one $, where $M^a_x \equiv  M^{a\one}_x$.  In other words,
$A^a\neq 0$ means that $a$ condenses on the boundary (\ie the bulk $a$ can
become the null excitation $\one$ on the boundary).  

The boundary excitations form a fusion category, denoted as $\cF$.  
We can fuse $a,b$ in bulk to $c$, and then fuse $c$ to the boundary.  Or we
can fuse $a,b$ in bulk to the boundary to get $x,y$, and then fuse $x,y$ on
boundary to $z$.  The two ways of fusion should give rise to the same result. This
requires $M^a_x$ to satisfy
\begin{align}
\label{NMMN}
\sum_c N^{ab}_{\eM,c} M^c_x = \sum_{y,z} M^a_y M^b_{z} N^{yz}_{B,x} 
\end{align}
where $N^{ab}_{\eM,c}$ describes the fusion ring of the bulk excitations in
$\eM$ and $N^{xy}_{B,z}$ describes the fusion ring of the boundary excitations.

Taking $x=\one$, \eqref{NMMN} reduces to
\begin{align}
\label{AAAMM}
\sum_c N^{ab}_{\eM,c} A^c 
&=
A^a A^b + \sum_{x\neq \one}  M^a_x M^b_{\bar x}
\end{align}
Since $M^a_x\geq 0$, we obtain an additional condition
on $A^a$
\begin{align}
 \label{NAAA}
\sum_c N^{ab}_{\eM,c} A^c & \geq A^a A^b.
\end{align}

From the conservation of quantum dimensions, we have
\begin{align}
 d_a = \sum_x M^a_x d_x 
= A^a + \sum_{x\neq \one}  M^a_x d_x.
\end{align}
This implies that, if $d_a\neq $integer, $\sum_{x\neq \one}  M^a_x d_x \geq 1$,
or
\begin{align}
 A^a \leq d_a -\del(d_a),
\end{align}
where $\del(d_a)$ is defined as
\begin{align}
\del(x) =
\begin{cases}
0 & \text{ if } x \in \Z \\
1 & \text{ if } x \not\in \Z \\
\end{cases}
\end{align}

By writing 
$\sum_y = \sum_{y=\one} +\sum_{y\neq \one}$ 
in \eqref{NMMN}, we
find
\begin{align}
\sum_c N^{ab}_{\eM,c} M^c_x & 
=  
A^a M^b_{x}  
+\sum_{y\neq  \one ,z} M^a_y M^b_{z} N^{yz}_{B,x} .
\end{align}
Noticing that $d_a-A^a = \sum_{x\neq \one}  M^a_x d_x$.
We multiply $d_x$ and sum over $x\neq \one$ in the above
\begin{align}
\sum_{c,x} N^{ab}_{\eM,c} (d_c-A^c) & 
=  
A^a (d_b-A^b) +\sum_{y\neq  \one ,z,x\neq \one} M^a_y M^b_{z} N^{yz}_{B,x} d_x
\nonumber\\
& =
A^a (d_b-A^b) 
+\sum_{x\neq \one} M^a_x A^b d_x
+\sum_{y\neq  \one ,z\neq \one,x\neq \one} M^a_y M^b_{z} N^{yz}_{B,x} d_x
\nonumber\\
&\geq A^a (d_b-A^b) + A^b (d_a-A^a) 
.
\end{align}
Taking $b=\bar a$ in
 \eqref{AAAMM}, we find
\begin{align}
\label{AaAa}
\sum_c N^{a\bar a}_{\eM,c} A^c 
&\geq
(A^a)^2  +\del(d_a).
\end{align}

Summarizing the above discussions and adding a modular covariant condition, we
have
\begin{align}
\label{Acond1}
  \v A &= D_\eM^{-1} S_\eM\v A, & \v A &= T_\eM\v A,
\nonumber\\
A^a & \leq  d_a -\del(d_a),  &
D_\eM &= d_{\cA} = \sum_a d_aA^a,
\nonumber\\
A^a A^b &\leq \sum_c N^{ab}_{\eM,c} A^c -\del_{a,\bar b}\del(d_a)
, &  
\sum_{c,x} N^{ab}_{\eM,c} (d_c-A^c) &\geq A^a (d_b-A^b) + A^b (d_a-A^a) 
,
\end{align}
where $\v A =(A^\one,\cdots,A^a,\cdots)^\top$.

Now, let us assume $\cA$ not to be Lagrangian.  In this case, the condensation
of $\cA$ will change $\eM$ to a non-trivial $\eM_{/\cA}$.  Let us consider the
domain wall between $\eM$ and $\eM_{/\cA}$.  Such a domain wall can be viewed
as a boundary of $\eM \boxtimes \overline \eM_{/\cA}$ topological order form by
stacking $\eM$ and the spatial reflection of $\eM_{/\cA}$.  Since the domain
wall, and hence the boundary, is gapped, there must be a Lagrangian condensable
algebra $\cA_{\eM \boxtimes \overline \eM_{/\cA}}$ in $\eM \boxtimes \overline
\eM_{/\cA}$, whose condensation gives rise to the gapped boundary.  Let 
\begin{align}
\cA_{\eM \boxtimes \overline \eM_{/\cA}} =
\bigoplus_{a\in \eM,\ i \in \eM_{/\cA}} A^{ai}\ a\otimes i
,
\end{align}
then the matrix $ A = (A^{ai})$ satisfies
\begin{align}
\label{STA}
& D_\eM^{-1}S_{\eM} A = A  D_{\eM_{/\cA}}^{-1}S_{\eM_{/\cA}}, \ \ \ \
   T_{\eM} A = A  T_{\eM_{/\cA}}, \ \ \ \
A^{ai} \leq d_ad_i -\del(d_a d_i),
 \\
&A^{ai} A^{bj}  \leq \sum_{c,k} 
N^{ab}_{\eM,c} 
N^{ij}_{\eM_{/\cA},k} 
A^{ck} -
\del_{a,\bar b} 
\del_{i,\bar j}\del(d_ad_i) 
\nonumber\\
&
\sum_{c,x} 
N^{ab}_{\eM,c} 
N^{ij}_{\eM_{/\cA},k} 
(d_cd_k-A^{ck}) \geq A^{ai} (d_bd_j-A^{bj}) + A^{bj} (d_ad_i-A^{ai}) 
\nonumber 
\end{align}
The above conditions only require
the domain wall between $\eM$ and $\eM_{/\cA}$ to be gapped.  However, since
$\eM$ and $\eM_{/\cA}$ are related by a condensation of $\cA$, there is a
special domain wall (called the canonical domain wall) that satisfies the
following condition: 
\begin{align}
\label{Aainzero}
\text{For any $i$, there exists an $a$
such that }\ A^{ai} \neq 0 
.
\end{align}
The above, together with $ T_{\eM} A = A  T_{\eM_{/\cA}}$, impies that the
eigenvalues of $T_{\eM_{/\cA}}$ are also  eigenvalues of $T_\eM$.

The canonical domain wall can be viewed as $\cA_{\eM\to \eM_{/\cA}}$-condensed
boundary of $\eM$ with 
\begin{align}
 \cA_{\eM\to \eM_{/\cA}} = \bigoplus A^{a\one} a .
\end{align}
We note that anyon $a$ in $\eM$ condenses on the canonical domain wall between
$\eM$ and $\eM_{/\cA}$, if and only if $A^{a\one}\neq 0$.  This implies that
\begin{align}
\label{AMA1}
\cA_{\eM\to \eM_{/\cA}}  = \cA ,\ \ \ \ A^a = A^{a\one}.
\end{align}
The  domain wall can also be viewed as $\cA_{\eM_{/\cA}\to\eM}$-condensed
boundary of $\eM_{/\cA}$ with 
\begin{align}
 \cA_{\eM_{/\cA}\to\eM} = \bigoplus A^{\one i} i .
\end{align}
Since $\eM_{/\cA}$ comes from a condensation of $\eM$, the canonical domain
wall must be an $\one$-condensed boundary of $\eM_{/\cA}$, \ie
\begin{align}
\label{AMA2}
 \cA_{\eM_{/\cA}\to\eM} = \one,\ \ \
A^{\one i} = \del_{\one,i}.
\end{align}

We can obtain more conditions on $A^a$. From \eqref{STA}, we find
\begin{align}
\label{ASSA}
\frac{D_{\eM_{/\cA}}}{D_\eM}\sum_{b\in \eM} (S_{\eM})^{ab} A^{bi} 
= 
\sum_{j\in \eM_{/\cA} } A^{aj} (S_{\eM_{/\cA}})^{ji} , 
\end{align}
Setting $a=\one$ in the above, we find that
\begin{align}
\label{DdAd}
\frac{D_{\eM_{/\cA}}}{D_\eM}
\sum_{b\in \eM} d_b A^{bi} 
=
d_i , 
\end{align}
or
\begin{align}
D_\eM =  D_{\eM_{/\cA}}  \sum_{b\in \eM} d_b A^{b} .
\end{align}
Eq. \eqref{ASSA} also
implies
\begin{align}
 \frac{\sum_{b\in \eM} 
(S_{\eM})^{ab} A^{bi} }{\sum_{b\in \eM} d_b A^b}
&= 
\sum_{j\in \eM_{/\cA} } A^{aj} (S_{\eM_{/\cA}})^{ji} 
=
\text{cyclotomic integer} ,
\ \ \ \
\text{for all }   a\in \eM ,\ \  i\in \eM_{/\cA}  
\end{align} 
Setting $i=\one$ in the above, we find that
$\cA  = \bigoplus_a A^{a} a$ must satisfies
\begin{align}
 \frac{\sum_{b\in \eM} (S_{\eM})^{ab} A^b}{\sum_{b\in \eM} d_b A^b}
= 
\sum_{j\in \eM_{/\cA} } A^{aj} d_j 
=
\text{cyclotomic integer for all } &  a\in \eM  
\end{align} 

Setting $i=j=\one$ in \eqref{STA}, we also obtain
\begin{align}
 A^{a} &\leq d_a -\del(d_a), 
\nonumber\\
A^a A^b &\leq \sum_c N^{ab}_{\eM,c} A^c -\del_{a,\bar a}\del(d_a)
\nonumber\\
\sum_{c} 
N^{ab}_{\eM,c} 
(d_c-A^{c}) &\geq A^{a} (d_b-A^{b}) + A^{b} (d_a-A^{a}) 
\end{align}
We also have
\begin{align}
 D_\eM^{-1}(S_{\eM})_{ab} A^b = 
A^{ai}  D_{\eM_{/\cA}}^{-1}(S_{\eM_{/\cA}})_{i\one}, \ \ \ \
   (T_{\eM})_{ab} A^b = A^{ai}  (T_{\eM_{/\cA}})_{i\one},
\end{align}
which only leads to a weak condition on the vector 
 $\v A =(A^\one,\cdots,A^a,\cdots)^\top$
\begin{align}
   T_{\eM} \v A = \v A . 
\end{align}
To summarize, a condensable algebra $A^a$, which may not a Lagrangian, must
satisfies
\begin{align}
\label{Aconds}
&   T_{\eM} \v A = \v A , \ \ \ \
\del_{a\one} \leq A^{a} \leq d_a -\del(d_a),\ \ \ \
\frac{ D_\eM }{\sum_{b\in \eM} d_b A^b} = D_{\cM_{/\cA}} =
d\text{-number} \geq 1, 
\nonumber\\
& 
\frac{\sum_{b\in \eM} 
(S_{\eM})^{ab} A^{b} }{\sum_{b\in \eM} d_b A^b}
= 
\text{cyclotomic integer} > 0,
\ \ \ \
\text{for all }   
a\in \eM ,
\\
& 
\sum_c N^{ab}_{\eM,c} A^c 
- A^a A^b 
-\del_{a,\bar b}\del(d_a)
\geq 0
,
\ \ \ \
(d_a-A^a)(d_b- A^b)
- (\sum_{c} N^{ab}_{\eM,c} A^{c} - A^a A^b)
\geq 0
.
\nonumber 
\end{align}
The last condition comes from $\sum_{c} 
N^{ab}_{\eM,c} (d_c-A^{c}) 
\geq 
A^{a} (d_b-A^{b}) + A^{b} (d_a-A^{a})$.
The last two conditions in \eqref{Aconds}
imply that
\begin{align}
A^{a} d_b + A^{b}d_a +2\del_{a,\bar b}\del(d_a)
\leq \sum_c N^{ab}_{\eM,c} A^c + d_ad_b \ \ \ &\text{linear in }A^a
\nonumber\\
(d_a-A^a)(d_b- A^b) 
-\del_{a,\bar b}\del(d_a) 
\geq 0
\ \ \ &\text{independent of }  N^{ab}_{\eM,c}
.
\end{align}


\section{Full lists of unitary modular data}
\label{UMD}

In this section, we list all the unitary modular data for rank 2 -- 12.  The
list is ordered by $D^2$.

The list also contain information on the factorization of modular data.  For
example, a rank-4 modular data $4_{0,4.}^{4,375}$ (the third entry in Appendix
\ref{uni4}), is the product of two rank-2 modular data $2_{1,2.}^{4,437}$ and $
2_{7,2.}^{4,625} $ (the first and the second entries in Appendix \ref{uni2}):
$4_{0,4.}^{4,375} = 2_{1,2.}^{4,437} \boxtimes 2_{7,2.}^{4,625}.  $

\subsection{Rank 2}
\label{uni2}

{\small

\noindent1. $2_{1,2.}^{4,437}$ \irep{1}:\ \ 
$d_i$ = ($1.0$,
$1.0$) 

\vskip 0.7ex
\hangindent=3em \hangafter=1
$D^2= 2.0 = 
2$

\vskip 0.7ex
\hangindent=3em \hangafter=1
$T = ( 0,
\frac{1}{4} )
$,

\vskip 0.7ex
\hangindent=3em \hangafter=1
$S$ = ($ 1$,
$ 1$;\ \ 
$ -1$)

  \vskip 2ex

\noindent2. $2_{7,2.}^{4,625}$ \irep{1}:\ \ 
$d_i$ = ($1.0$,
$1.0$) 

\vskip 0.7ex
\hangindent=3em \hangafter=1
$D^2= 2.0 = 
2$

\vskip 0.7ex
\hangindent=3em \hangafter=1
$T = ( 0,
\frac{3}{4} )
$,

\vskip 0.7ex
\hangindent=3em \hangafter=1
$S$ = ($ 1$,
$ 1$;\ \ 
$ -1$)

  \vskip 2ex

\noindent3. $2_{\frac{14}{5},3.618}^{5,395}$ \irep{2}:\ \ 
$d_i$ = ($1.0$,
$1.618$) 

\vskip 0.7ex
\hangindent=3em \hangafter=1
$D^2= 3.618 = 
\frac{5+\sqrt{5}}{2}$

\vskip 0.7ex
\hangindent=3em \hangafter=1
$T = ( 0,
\frac{2}{5} )
$,

\vskip 0.7ex
\hangindent=3em \hangafter=1
$S$ = ($ 1$,
$ \frac{1+\sqrt{5}}{2}$;\ \ 
$ -1$)

  \vskip 2ex

\noindent4. $2_{\frac{26}{5},3.618}^{5,720}$ \irep{2}:\ \ 
$d_i$ = ($1.0$,
$1.618$) 

\vskip 0.7ex
\hangindent=3em \hangafter=1
$D^2= 3.618 = 
\frac{5+\sqrt{5}}{2}$

\vskip 0.7ex
\hangindent=3em \hangafter=1
$T = ( 0,
\frac{3}{5} )
$,

\vskip 0.7ex
\hangindent=3em \hangafter=1
$S$ = ($ 1$,
$ \frac{1+\sqrt{5}}{2}$;\ \ 
$ -1$)

  \vskip 2ex 

}

\subsection{Rank 3}
\label{uni3}

{\small

\noindent1. $3_{2,3.}^{3,527}$ \irep{2}:\ \ 
$d_i$ = ($1.0$,
$1.0$,
$1.0$) 

\vskip 0.7ex
\hangindent=3em \hangafter=1
$D^2= 3.0 = 
3$

\vskip 0.7ex
\hangindent=3em \hangafter=1
$T = ( 0,
\frac{1}{3},
\frac{1}{3} )
$,

\vskip 0.7ex
\hangindent=3em \hangafter=1
$S$ = ($ 1$,
$ 1$,
$ 1$;\ \ 
$ \zeta_{3}^{1}$,
$ -\zeta_{6}^{1}$;\ \ 
$ \zeta_{3}^{1}$)

  \vskip 2ex

\noindent2. $3_{6,3.}^{3,138}$ \irep{2}:\ \ 
$d_i$ = ($1.0$,
$1.0$,
$1.0$) 

\vskip 0.7ex
\hangindent=3em \hangafter=1
$D^2= 3.0 = 
3$

\vskip 0.7ex
\hangindent=3em \hangafter=1
$T = ( 0,
\frac{2}{3},
\frac{2}{3} )
$,

\vskip 0.7ex
\hangindent=3em \hangafter=1
$S$ = ($ 1$,
$ 1$,
$ 1$;\ \ 
$ -\zeta_{6}^{1}$,
$ \zeta_{3}^{1}$;\ \ 
$ -\zeta_{6}^{1}$)

  \vskip 2ex

\noindent3. $3_{\frac{1}{2},4.}^{16,598}$ \irep{4}:\ \ 
$d_i$ = ($1.0$,
$1.0$,
$1.414$) 

\vskip 0.7ex
\hangindent=3em \hangafter=1
$D^2= 4.0 = 
4$

\vskip 0.7ex
\hangindent=3em \hangafter=1
$T = ( 0,
\frac{1}{2},
\frac{1}{16} )
$,

\vskip 0.7ex
\hangindent=3em \hangafter=1
$S$ = ($ 1$,
$ 1$,
$ \sqrt{2}$;\ \ 
$ 1$,
$ -\sqrt{2}$;\ \ 
$0$)

  \vskip 2ex

\noindent4. $3_{\frac{3}{2},4.}^{16,553}$ \irep{4}:\ \ 
$d_i$ = ($1.0$,
$1.0$,
$1.414$) 

\vskip 0.7ex
\hangindent=3em \hangafter=1
$D^2= 4.0 = 
4$

\vskip 0.7ex
\hangindent=3em \hangafter=1
$T = ( 0,
\frac{1}{2},
\frac{3}{16} )
$,

\vskip 0.7ex
\hangindent=3em \hangafter=1
$S$ = ($ 1$,
$ 1$,
$ \sqrt{2}$;\ \ 
$ 1$,
$ -\sqrt{2}$;\ \ 
$0$)

  \vskip 2ex

\noindent5. $3_{\frac{5}{2},4.}^{16,465}$ \irep{4}:\ \ 
$d_i$ = ($1.0$,
$1.0$,
$1.414$) 

\vskip 0.7ex
\hangindent=3em \hangafter=1
$D^2= 4.0 = 
4$

\vskip 0.7ex
\hangindent=3em \hangafter=1
$T = ( 0,
\frac{1}{2},
\frac{5}{16} )
$,

\vskip 0.7ex
\hangindent=3em \hangafter=1
$S$ = ($ 1$,
$ 1$,
$ \sqrt{2}$;\ \ 
$ 1$,
$ -\sqrt{2}$;\ \ 
$0$)

  \vskip 2ex

\noindent6. $3_{\frac{7}{2},4.}^{16,332}$ \irep{4}:\ \ 
$d_i$ = ($1.0$,
$1.0$,
$1.414$) 

\vskip 0.7ex
\hangindent=3em \hangafter=1
$D^2= 4.0 = 
4$

\vskip 0.7ex
\hangindent=3em \hangafter=1
$T = ( 0,
\frac{1}{2},
\frac{7}{16} )
$,

\vskip 0.7ex
\hangindent=3em \hangafter=1
$S$ = ($ 1$,
$ 1$,
$ \sqrt{2}$;\ \ 
$ 1$,
$ -\sqrt{2}$;\ \ 
$0$)

  \vskip 2ex

\noindent7. $3_{\frac{9}{2},4.}^{16,156}$ \irep{4}:\ \ 
$d_i$ = ($1.0$,
$1.0$,
$1.414$) 

\vskip 0.7ex
\hangindent=3em \hangafter=1
$D^2= 4.0 = 
4$

\vskip 0.7ex
\hangindent=3em \hangafter=1
$T = ( 0,
\frac{1}{2},
\frac{9}{16} )
$,

\vskip 0.7ex
\hangindent=3em \hangafter=1
$S$ = ($ 1$,
$ 1$,
$ \sqrt{2}$;\ \ 
$ 1$,
$ -\sqrt{2}$;\ \ 
$0$)

  \vskip 2ex

\noindent8. $3_{\frac{11}{2},4.}^{16,648}$ \irep{4}:\ \ 
$d_i$ = ($1.0$,
$1.0$,
$1.414$) 

\vskip 0.7ex
\hangindent=3em \hangafter=1
$D^2= 4.0 = 
4$

\vskip 0.7ex
\hangindent=3em \hangafter=1
$T = ( 0,
\frac{1}{2},
\frac{11}{16} )
$,

\vskip 0.7ex
\hangindent=3em \hangafter=1
$S$ = ($ 1$,
$ 1$,
$ \sqrt{2}$;\ \ 
$ 1$,
$ -\sqrt{2}$;\ \ 
$0$)

  \vskip 2ex

\noindent9. $3_{\frac{13}{2},4.}^{16,330}$ \irep{4}:\ \ 
$d_i$ = ($1.0$,
$1.0$,
$1.414$) 

\vskip 0.7ex
\hangindent=3em \hangafter=1
$D^2= 4.0 = 
4$

\vskip 0.7ex
\hangindent=3em \hangafter=1
$T = ( 0,
\frac{1}{2},
\frac{13}{16} )
$,

\vskip 0.7ex
\hangindent=3em \hangafter=1
$S$ = ($ 1$,
$ 1$,
$ \sqrt{2}$;\ \ 
$ 1$,
$ -\sqrt{2}$;\ \ 
$0$)

  \vskip 2ex

\noindent10. $3_{\frac{15}{2},4.}^{16,639}$ \irep{4}:\ \ 
$d_i$ = ($1.0$,
$1.0$,
$1.414$) 

\vskip 0.7ex
\hangindent=3em \hangafter=1
$D^2= 4.0 = 
4$

\vskip 0.7ex
\hangindent=3em \hangafter=1
$T = ( 0,
\frac{1}{2},
\frac{15}{16} )
$,

\vskip 0.7ex
\hangindent=3em \hangafter=1
$S$ = ($ 1$,
$ 1$,
$ \sqrt{2}$;\ \ 
$ 1$,
$ -\sqrt{2}$;\ \ 
$0$)

  \vskip 2ex

\noindent11. $3_{\frac{48}{7},9.295}^{7,790}$ \irep{3}:\ \ 
$d_i$ = ($1.0$,
$1.801$,
$2.246$) 

\vskip 0.7ex
\hangindent=3em \hangafter=1
$D^2= 9.295 = 
6+3c^{1}_{7}
+c^{2}_{7}
$

\vskip 0.7ex
\hangindent=3em \hangafter=1
$T = ( 0,
\frac{1}{7},
\frac{5}{7} )
$,

\vskip 0.7ex
\hangindent=3em \hangafter=1
$S$ = ($ 1$,
$ -c_{7}^{3}$,
$ \xi_{7}^{3}$;\ \ 
$ -\xi_{7}^{3}$,
$ 1$;\ \ 
$ c_{7}^{3}$)

  \vskip 2ex

\noindent12. $3_{\frac{8}{7},9.295}^{7,245}$ \irep{3}:\ \ 
$d_i$ = ($1.0$,
$1.801$,
$2.246$) 

\vskip 0.7ex
\hangindent=3em \hangafter=1
$D^2= 9.295 = 
6+3c^{1}_{7}
+c^{2}_{7}
$

\vskip 0.7ex
\hangindent=3em \hangafter=1
$T = ( 0,
\frac{6}{7},
\frac{2}{7} )
$,

\vskip 0.7ex
\hangindent=3em \hangafter=1
$S$ = ($ 1$,
$ -c_{7}^{3}$,
$ \xi_{7}^{3}$;\ \ 
$ -\xi_{7}^{3}$,
$ 1$;\ \ 
$ c_{7}^{3}$)

  \vskip 2ex 

}

\subsection{Rank 4}
\label{uni4}

{\small

\noindent1. $4_{0,4.}^{2,750}$ \irep{0}:\ \ 
$d_i$ = ($1.0$,
$1.0$,
$1.0$,
$1.0$) 

\vskip 0.7ex
\hangindent=3em \hangafter=1
$D^2= 4.0 = 
4$

\vskip 0.7ex
\hangindent=3em \hangafter=1
$T = ( 0,
0,
0,
\frac{1}{2} )
$,

\vskip 0.7ex
\hangindent=3em \hangafter=1
$S$ = ($ 1$,
$ 1$,
$ 1$,
$ 1$;\ \ 
$ 1$,
$ -1$,
$ -1$;\ \ 
$ 1$,
$ -1$;\ \ 
$ 1$)

  \vskip 2ex

\noindent2. $4_{4,4.}^{2,250}$ \irep{0}:\ \ 
$d_i$ = ($1.0$,
$1.0$,
$1.0$,
$1.0$) 

\vskip 0.7ex
\hangindent=3em \hangafter=1
$D^2= 4.0 = 
4$

\vskip 0.7ex
\hangindent=3em \hangafter=1
$T = ( 0,
\frac{1}{2},
\frac{1}{2},
\frac{1}{2} )
$,

\vskip 0.7ex
\hangindent=3em \hangafter=1
$S$ = ($ 1$,
$ 1$,
$ 1$,
$ 1$;\ \ 
$ 1$,
$ -1$,
$ -1$;\ \ 
$ 1$,
$ -1$;\ \ 
$ 1$)

  \vskip 2ex

\noindent3. $4_{0,4.}^{4,375}$ \irep{0}:\ \ 
$d_i$ = ($1.0$,
$1.0$,
$1.0$,
$1.0$) 

\vskip 0.7ex
\hangindent=3em \hangafter=1
$D^2= 4.0 = 
4$

\vskip 0.7ex
\hangindent=3em \hangafter=1
$T = ( 0,
0,
\frac{1}{4},
\frac{3}{4} )
$,

\vskip 0.7ex
\hangindent=3em \hangafter=1
$S$ = ($ 1$,
$ 1$,
$ 1$,
$ 1$;\ \ 
$ 1$,
$ -1$,
$ -1$;\ \ 
$ -1$,
$ 1$;\ \ 
$ -1$)

Factors = $2_{1,2.}^{4,437}\boxtimes 2_{7,2.}^{4,625}$

  \vskip 2ex

\noindent4. $4_{2,4.}^{4,625}$ \irep{0}:\ \ 
$d_i$ = ($1.0$,
$1.0$,
$1.0$,
$1.0$) 

\vskip 0.7ex
\hangindent=3em \hangafter=1
$D^2= 4.0 = 
4$

\vskip 0.7ex
\hangindent=3em \hangafter=1
$T = ( 0,
\frac{1}{2},
\frac{1}{4},
\frac{1}{4} )
$,

\vskip 0.7ex
\hangindent=3em \hangafter=1
$S$ = ($ 1$,
$ 1$,
$ 1$,
$ 1$;\ \ 
$ 1$,
$ -1$,
$ -1$;\ \ 
$ -1$,
$ 1$;\ \ 
$ -1$)

Factors = $2_{1,2.}^{4,437}\boxtimes 2_{1,2.}^{4,437}$

  \vskip 2ex

\noindent5. $4_{6,4.}^{4,375}$ \irep{0}:\ \ 
$d_i$ = ($1.0$,
$1.0$,
$1.0$,
$1.0$) 

\vskip 0.7ex
\hangindent=3em \hangafter=1
$D^2= 4.0 = 
4$

\vskip 0.7ex
\hangindent=3em \hangafter=1
$T = ( 0,
\frac{1}{2},
\frac{3}{4},
\frac{3}{4} )
$,

\vskip 0.7ex
\hangindent=3em \hangafter=1
$S$ = ($ 1$,
$ 1$,
$ 1$,
$ 1$;\ \ 
$ 1$,
$ -1$,
$ -1$;\ \ 
$ -1$,
$ 1$;\ \ 
$ -1$)

Factors = $2_{7,2.}^{4,625}\boxtimes 2_{7,2.}^{4,625}$

  \vskip 2ex

\noindent6. $4_{1,4.}^{8,718}$ \irep{6}:\ \ 
$d_i$ = ($1.0$,
$1.0$,
$1.0$,
$1.0$) 

\vskip 0.7ex
\hangindent=3em \hangafter=1
$D^2= 4.0 = 
4$

\vskip 0.7ex
\hangindent=3em \hangafter=1
$T = ( 0,
\frac{1}{2},
\frac{1}{8},
\frac{1}{8} )
$,

\vskip 0.7ex
\hangindent=3em \hangafter=1
$S$ = ($ 1$,
$ 1$,
$ 1$,
$ 1$;\ \ 
$ 1$,
$ -1$,
$ -1$;\ \ 
$-\mathrm{i}$,
$\mathrm{i}$;\ \ 
$-\mathrm{i}$)

  \vskip 2ex

\noindent7. $4_{3,4.}^{8,468}$ \irep{6}:\ \ 
$d_i$ = ($1.0$,
$1.0$,
$1.0$,
$1.0$) 

\vskip 0.7ex
\hangindent=3em \hangafter=1
$D^2= 4.0 = 
4$

\vskip 0.7ex
\hangindent=3em \hangafter=1
$T = ( 0,
\frac{1}{2},
\frac{3}{8},
\frac{3}{8} )
$,

\vskip 0.7ex
\hangindent=3em \hangafter=1
$S$ = ($ 1$,
$ 1$,
$ 1$,
$ 1$;\ \ 
$ 1$,
$ -1$,
$ -1$;\ \ 
$\mathrm{i}$,
$-\mathrm{i}$;\ \ 
$\mathrm{i}$)

  \vskip 2ex

\noindent8. $4_{5,4.}^{8,312}$ \irep{6}:\ \ 
$d_i$ = ($1.0$,
$1.0$,
$1.0$,
$1.0$) 

\vskip 0.7ex
\hangindent=3em \hangafter=1
$D^2= 4.0 = 
4$

\vskip 0.7ex
\hangindent=3em \hangafter=1
$T = ( 0,
\frac{1}{2},
\frac{5}{8},
\frac{5}{8} )
$,

\vskip 0.7ex
\hangindent=3em \hangafter=1
$S$ = ($ 1$,
$ 1$,
$ 1$,
$ 1$;\ \ 
$ 1$,
$ -1$,
$ -1$;\ \ 
$-\mathrm{i}$,
$\mathrm{i}$;\ \ 
$-\mathrm{i}$)

  \vskip 2ex

\noindent9. $4_{7,4.}^{8,781}$ \irep{6}:\ \ 
$d_i$ = ($1.0$,
$1.0$,
$1.0$,
$1.0$) 

\vskip 0.7ex
\hangindent=3em \hangafter=1
$D^2= 4.0 = 
4$

\vskip 0.7ex
\hangindent=3em \hangafter=1
$T = ( 0,
\frac{1}{2},
\frac{7}{8},
\frac{7}{8} )
$,

\vskip 0.7ex
\hangindent=3em \hangafter=1
$S$ = ($ 1$,
$ 1$,
$ 1$,
$ 1$;\ \ 
$ 1$,
$ -1$,
$ -1$;\ \ 
$\mathrm{i}$,
$-\mathrm{i}$;\ \ 
$\mathrm{i}$)

  \vskip 2ex

\noindent10. $4_{\frac{19}{5},7.236}^{20,304}$ \irep{8}:\ \ 
$d_i$ = ($1.0$,
$1.0$,
$1.618$,
$1.618$) 

\vskip 0.7ex
\hangindent=3em \hangafter=1
$D^2= 7.236 = 
5+\sqrt{5}$

\vskip 0.7ex
\hangindent=3em \hangafter=1
$T = ( 0,
\frac{1}{4},
\frac{2}{5},
\frac{13}{20} )
$,

\vskip 0.7ex
\hangindent=3em \hangafter=1
$S$ = ($ 1$,
$ 1$,
$ \frac{1+\sqrt{5}}{2}$,
$ \frac{1+\sqrt{5}}{2}$;\ \ 
$ -1$,
$ \frac{1+\sqrt{5}}{2}$,
$ -\frac{1+\sqrt{5}}{2}$;\ \ 
$ -1$,
$ -1$;\ \ 
$ 1$)

Factors = $2_{1,2.}^{4,437}\boxtimes 2_{\frac{14}{5},3.618}^{5,395}$

  \vskip 2ex

\noindent11. $4_{\frac{31}{5},7.236}^{20,505}$ \irep{8}:\ \ 
$d_i$ = ($1.0$,
$1.0$,
$1.618$,
$1.618$) 

\vskip 0.7ex
\hangindent=3em \hangafter=1
$D^2= 7.236 = 
5+\sqrt{5}$

\vskip 0.7ex
\hangindent=3em \hangafter=1
$T = ( 0,
\frac{1}{4},
\frac{3}{5},
\frac{17}{20} )
$,

\vskip 0.7ex
\hangindent=3em \hangafter=1
$S$ = ($ 1$,
$ 1$,
$ \frac{1+\sqrt{5}}{2}$,
$ \frac{1+\sqrt{5}}{2}$;\ \ 
$ -1$,
$ \frac{1+\sqrt{5}}{2}$,
$ -\frac{1+\sqrt{5}}{2}$;\ \ 
$ -1$,
$ -1$;\ \ 
$ 1$)

Factors = $2_{1,2.}^{4,437}\boxtimes 2_{\frac{26}{5},3.618}^{5,720}$

  \vskip 2ex

\noindent12. $4_{\frac{9}{5},7.236}^{20,451}$ \irep{8}:\ \ 
$d_i$ = ($1.0$,
$1.0$,
$1.618$,
$1.618$) 

\vskip 0.7ex
\hangindent=3em \hangafter=1
$D^2= 7.236 = 
5+\sqrt{5}$

\vskip 0.7ex
\hangindent=3em \hangafter=1
$T = ( 0,
\frac{3}{4},
\frac{2}{5},
\frac{3}{20} )
$,

\vskip 0.7ex
\hangindent=3em \hangafter=1
$S$ = ($ 1$,
$ 1$,
$ \frac{1+\sqrt{5}}{2}$,
$ \frac{1+\sqrt{5}}{2}$;\ \ 
$ -1$,
$ \frac{1+\sqrt{5}}{2}$,
$ -\frac{1+\sqrt{5}}{2}$;\ \ 
$ -1$,
$ -1$;\ \ 
$ 1$)

Factors = $2_{7,2.}^{4,625}\boxtimes 2_{\frac{14}{5},3.618}^{5,395}$

  \vskip 2ex

\noindent13. $4_{\frac{21}{5},7.236}^{20,341}$ \irep{8}:\ \ 
$d_i$ = ($1.0$,
$1.0$,
$1.618$,
$1.618$) 

\vskip 0.7ex
\hangindent=3em \hangafter=1
$D^2= 7.236 = 
5+\sqrt{5}$

\vskip 0.7ex
\hangindent=3em \hangafter=1
$T = ( 0,
\frac{3}{4},
\frac{3}{5},
\frac{7}{20} )
$,

\vskip 0.7ex
\hangindent=3em \hangafter=1
$S$ = ($ 1$,
$ 1$,
$ \frac{1+\sqrt{5}}{2}$,
$ \frac{1+\sqrt{5}}{2}$;\ \ 
$ -1$,
$ \frac{1+\sqrt{5}}{2}$,
$ -\frac{1+\sqrt{5}}{2}$;\ \ 
$ -1$,
$ -1$;\ \ 
$ 1$)

Factors = $2_{7,2.}^{4,625}\boxtimes 2_{\frac{26}{5},3.618}^{5,720}$

  \vskip 2ex

\noindent14. $4_{\frac{28}{5},13.09}^{5,479}$ \irep{5}:\ \ 
$d_i$ = ($1.0$,
$1.618$,
$1.618$,
$2.618$) 

\vskip 0.7ex
\hangindent=3em \hangafter=1
$D^2= 13.90 = 
\frac{15+5\sqrt{5}}{2}$

\vskip 0.7ex
\hangindent=3em \hangafter=1
$T = ( 0,
\frac{2}{5},
\frac{2}{5},
\frac{4}{5} )
$,

\vskip 0.7ex
\hangindent=3em \hangafter=1
$S$ = ($ 1$,
$ \frac{1+\sqrt{5}}{2}$,
$ \frac{1+\sqrt{5}}{2}$,
$ \frac{3+\sqrt{5}}{2}$;\ \ 
$ -1$,
$ \frac{3+\sqrt{5}}{2}$,
$ -\frac{1+\sqrt{5}}{2}$;\ \ 
$ -1$,
$ -\frac{1+\sqrt{5}}{2}$;\ \ 
$ 1$)

Factors = $2_{\frac{14}{5},3.618}^{5,395}\boxtimes 2_{\frac{14}{5},3.618}^{5,395}$

  \vskip 2ex

\noindent15. $4_{0,13.09}^{5,872}$ \irep{5}:\ \ 
$d_i$ = ($1.0$,
$1.618$,
$1.618$,
$2.618$) 

\vskip 0.7ex
\hangindent=3em \hangafter=1
$D^2= 13.90 = 
\frac{15+5\sqrt{5}}{2}$

\vskip 0.7ex
\hangindent=3em \hangafter=1
$T = ( 0,
\frac{2}{5},
\frac{3}{5},
0 )
$,

\vskip 0.7ex
\hangindent=3em \hangafter=1
$S$ = ($ 1$,
$ \frac{1+\sqrt{5}}{2}$,
$ \frac{1+\sqrt{5}}{2}$,
$ \frac{3+\sqrt{5}}{2}$;\ \ 
$ -1$,
$ \frac{3+\sqrt{5}}{2}$,
$ -\frac{1+\sqrt{5}}{2}$;\ \ 
$ -1$,
$ -\frac{1+\sqrt{5}}{2}$;\ \ 
$ 1$)

Factors = $2_{\frac{14}{5},3.618}^{5,395}\boxtimes 2_{\frac{26}{5},3.618}^{5,720}$

  \vskip 2ex

\noindent16. $4_{\frac{12}{5},13.09}^{5,443}$ \irep{5}:\ \ 
$d_i$ = ($1.0$,
$1.618$,
$1.618$,
$2.618$) 

\vskip 0.7ex
\hangindent=3em \hangafter=1
$D^2= 13.90 = 
\frac{15+5\sqrt{5}}{2}$

\vskip 0.7ex
\hangindent=3em \hangafter=1
$T = ( 0,
\frac{3}{5},
\frac{3}{5},
\frac{1}{5} )
$,

\vskip 0.7ex
\hangindent=3em \hangafter=1
$S$ = ($ 1$,
$ \frac{1+\sqrt{5}}{2}$,
$ \frac{1+\sqrt{5}}{2}$,
$ \frac{3+\sqrt{5}}{2}$;\ \ 
$ -1$,
$ \frac{3+\sqrt{5}}{2}$,
$ -\frac{1+\sqrt{5}}{2}$;\ \ 
$ -1$,
$ -\frac{1+\sqrt{5}}{2}$;\ \ 
$ 1$)

Factors = $2_{\frac{26}{5},3.618}^{5,720}\boxtimes 2_{\frac{26}{5},3.618}^{5,720}$

  \vskip 2ex

\noindent17. $4_{\frac{10}{3},19.23}^{9,459}$ \irep{7}:\ \ 
$d_i$ = ($1.0$,
$1.879$,
$2.532$,
$2.879$) 

\vskip 0.7ex
\hangindent=3em \hangafter=1
$D^2= 19.234 = 
9+6c^{1}_{9}
+3c^{2}_{9}
$

\vskip 0.7ex
\hangindent=3em \hangafter=1
$T = ( 0,
\frac{1}{3},
\frac{2}{9},
\frac{2}{3} )
$,

\vskip 0.7ex
\hangindent=3em \hangafter=1
$S$ = ($ 1$,
$ -c_{9}^{4}$,
$ \xi_{9}^{3}$,
$ \xi_{9}^{5}$;\ \ 
$ -\xi_{9}^{5}$,
$ \xi_{9}^{3}$,
$ -1$;\ \ 
$0$,
$ -\xi_{9}^{3}$;\ \ 
$ -c_{9}^{4}$)

  \vskip 2ex

\noindent18. $4_{\frac{14}{3},19.23}^{9,614}$ \irep{7}:\ \ 
$d_i$ = ($1.0$,
$1.879$,
$2.532$,
$2.879$) 

\vskip 0.7ex
\hangindent=3em \hangafter=1
$D^2= 19.234 = 
9+6c^{1}_{9}
+3c^{2}_{9}
$

\vskip 0.7ex
\hangindent=3em \hangafter=1
$T = ( 0,
\frac{2}{3},
\frac{7}{9},
\frac{1}{3} )
$,

\vskip 0.7ex
\hangindent=3em \hangafter=1
$S$ = ($ 1$,
$ -c_{9}^{4}$,
$ \xi_{9}^{3}$,
$ \xi_{9}^{5}$;\ \ 
$ -\xi_{9}^{5}$,
$ \xi_{9}^{3}$,
$ -1$;\ \ 
$0$,
$ -\xi_{9}^{3}$;\ \ 
$ -c_{9}^{4}$)

  \vskip 2ex 

}

\subsection{Rank 5}
\label{uni5}

{\small

\noindent1. $5_{0,5.}^{5,110}$ \irep{9}:\ \ 
$d_i$ = ($1.0$,
$1.0$,
$1.0$,
$1.0$,
$1.0$) 

\vskip 0.7ex
\hangindent=3em \hangafter=1
$D^2= 5.0 = 
5$

\vskip 0.7ex
\hangindent=3em \hangafter=1
$T = ( 0,
\frac{1}{5},
\frac{1}{5},
\frac{4}{5},
\frac{4}{5} )
$,

\vskip 0.7ex
\hangindent=3em \hangafter=1
$S$ = ($ 1$,
$ 1$,
$ 1$,
$ 1$,
$ 1$;\ \ 
$ -\zeta_{10}^{1}$,
$ \zeta_{5}^{2}$,
$ -\zeta_{10}^{3}$,
$ \zeta_{5}^{1}$;\ \ 
$ -\zeta_{10}^{1}$,
$ \zeta_{5}^{1}$,
$ -\zeta_{10}^{3}$;\ \ 
$ \zeta_{5}^{2}$,
$ -\zeta_{10}^{1}$;\ \ 
$ \zeta_{5}^{2}$)

  \vskip 2ex

\noindent2. $5_{4,5.}^{5,210}$ \irep{9}:\ \ 
$d_i$ = ($1.0$,
$1.0$,
$1.0$,
$1.0$,
$1.0$) 

\vskip 0.7ex
\hangindent=3em \hangafter=1
$D^2= 5.0 = 
5$

\vskip 0.7ex
\hangindent=3em \hangafter=1
$T = ( 0,
\frac{2}{5},
\frac{2}{5},
\frac{3}{5},
\frac{3}{5} )
$,

\vskip 0.7ex
\hangindent=3em \hangafter=1
$S$ = ($ 1$,
$ 1$,
$ 1$,
$ 1$,
$ 1$;\ \ 
$ \zeta_{5}^{1}$,
$ -\zeta_{10}^{3}$,
$ -\zeta_{10}^{1}$,
$ \zeta_{5}^{2}$;\ \ 
$ \zeta_{5}^{1}$,
$ \zeta_{5}^{2}$,
$ -\zeta_{10}^{1}$;\ \ 
$ -\zeta_{10}^{3}$,
$ \zeta_{5}^{1}$;\ \ 
$ -\zeta_{10}^{3}$)

  \vskip 2ex

\noindent3. $5_{2,12.}^{24,940}$ \irep{17}:\ \ 
$d_i$ = ($1.0$,
$1.0$,
$1.732$,
$1.732$,
$2.0$) 

\vskip 0.7ex
\hangindent=3em \hangafter=1
$D^2= 12.0 = 
12$

\vskip 0.7ex
\hangindent=3em \hangafter=1
$T = ( 0,
0,
\frac{1}{8},
\frac{5}{8},
\frac{1}{3} )
$,

\vskip 0.7ex
\hangindent=3em \hangafter=1
$S$ = ($ 1$,
$ 1$,
$ \sqrt{3}$,
$ \sqrt{3}$,
$ 2$;\ \ 
$ 1$,
$ -\sqrt{3}$,
$ -\sqrt{3}$,
$ 2$;\ \ 
$ \sqrt{3}$,
$ -\sqrt{3}$,
$0$;\ \ 
$ \sqrt{3}$,
$0$;\ \ 
$ -2$)

  \vskip 2ex

\noindent4. $5_{6,12.}^{24,273}$ \irep{17}:\ \ 
$d_i$ = ($1.0$,
$1.0$,
$1.732$,
$1.732$,
$2.0$) 

\vskip 0.7ex
\hangindent=3em \hangafter=1
$D^2= 12.0 = 
12$

\vskip 0.7ex
\hangindent=3em \hangafter=1
$T = ( 0,
0,
\frac{1}{8},
\frac{5}{8},
\frac{2}{3} )
$,

\vskip 0.7ex
\hangindent=3em \hangafter=1
$S$ = ($ 1$,
$ 1$,
$ \sqrt{3}$,
$ \sqrt{3}$,
$ 2$;\ \ 
$ 1$,
$ -\sqrt{3}$,
$ -\sqrt{3}$,
$ 2$;\ \ 
$ -\sqrt{3}$,
$ \sqrt{3}$,
$0$;\ \ 
$ -\sqrt{3}$,
$0$;\ \ 
$ -2$)

  \vskip 2ex

\noindent5. $5_{2,12.}^{24,741}$ \irep{17}:\ \ 
$d_i$ = ($1.0$,
$1.0$,
$1.732$,
$1.732$,
$2.0$) 

\vskip 0.7ex
\hangindent=3em \hangafter=1
$D^2= 12.0 = 
12$

\vskip 0.7ex
\hangindent=3em \hangafter=1
$T = ( 0,
0,
\frac{3}{8},
\frac{7}{8},
\frac{1}{3} )
$,

\vskip 0.7ex
\hangindent=3em \hangafter=1
$S$ = ($ 1$,
$ 1$,
$ \sqrt{3}$,
$ \sqrt{3}$,
$ 2$;\ \ 
$ 1$,
$ -\sqrt{3}$,
$ -\sqrt{3}$,
$ 2$;\ \ 
$ -\sqrt{3}$,
$ \sqrt{3}$,
$0$;\ \ 
$ -\sqrt{3}$,
$0$;\ \ 
$ -2$)

  \vskip 2ex

\noindent6. $5_{6,12.}^{24,592}$ \irep{17}:\ \ 
$d_i$ = ($1.0$,
$1.0$,
$1.732$,
$1.732$,
$2.0$) 

\vskip 0.7ex
\hangindent=3em \hangafter=1
$D^2= 12.0 = 
12$

\vskip 0.7ex
\hangindent=3em \hangafter=1
$T = ( 0,
0,
\frac{3}{8},
\frac{7}{8},
\frac{2}{3} )
$,

\vskip 0.7ex
\hangindent=3em \hangafter=1
$S$ = ($ 1$,
$ 1$,
$ \sqrt{3}$,
$ \sqrt{3}$,
$ 2$;\ \ 
$ 1$,
$ -\sqrt{3}$,
$ -\sqrt{3}$,
$ 2$;\ \ 
$ \sqrt{3}$,
$ -\sqrt{3}$,
$0$;\ \ 
$ \sqrt{3}$,
$0$;\ \ 
$ -2$)

  \vskip 2ex

\noindent7. $5_{\frac{72}{11},34.64}^{11,216}$ \irep{15}:\ \ 
$d_i$ = ($1.0$,
$1.918$,
$2.682$,
$3.228$,
$3.513$) 

\vskip 0.7ex
\hangindent=3em \hangafter=1
$D^2= 34.646 = 
15+10c^{1}_{11}
+6c^{2}_{11}
+3c^{3}_{11}
+c^{4}_{11}
$

\vskip 0.7ex
\hangindent=3em \hangafter=1
$T = ( 0,
\frac{2}{11},
\frac{9}{11},
\frac{10}{11},
\frac{5}{11} )
$,

\vskip 0.7ex
\hangindent=3em \hangafter=1
$S$ = ($ 1$,
$ -c_{11}^{5}$,
$ \xi_{11}^{3}$,
$ \xi_{11}^{7}$,
$ \xi_{11}^{5}$;\ \ 
$ -\xi_{11}^{7}$,
$ \xi_{11}^{5}$,
$ -\xi_{11}^{3}$,
$ 1$;\ \ 
$ -c_{11}^{5}$,
$ -1$,
$ -\xi_{11}^{7}$;\ \ 
$ \xi_{11}^{5}$,
$ c_{11}^{5}$;\ \ 
$ \xi_{11}^{3}$)

  \vskip 2ex

\noindent8. $5_{\frac{16}{11},34.64}^{11,640}$ \irep{15}:\ \ 
$d_i$ = ($1.0$,
$1.918$,
$2.682$,
$3.228$,
$3.513$) 

\vskip 0.7ex
\hangindent=3em \hangafter=1
$D^2= 34.646 = 
15+10c^{1}_{11}
+6c^{2}_{11}
+3c^{3}_{11}
+c^{4}_{11}
$

\vskip 0.7ex
\hangindent=3em \hangafter=1
$T = ( 0,
\frac{9}{11},
\frac{2}{11},
\frac{1}{11},
\frac{6}{11} )
$,

\vskip 0.7ex
\hangindent=3em \hangafter=1
$S$ = ($ 1$,
$ -c_{11}^{5}$,
$ \xi_{11}^{3}$,
$ \xi_{11}^{7}$,
$ \xi_{11}^{5}$;\ \ 
$ -\xi_{11}^{7}$,
$ \xi_{11}^{5}$,
$ -\xi_{11}^{3}$,
$ 1$;\ \ 
$ -c_{11}^{5}$,
$ -1$,
$ -\xi_{11}^{7}$;\ \ 
$ \xi_{11}^{5}$,
$ c_{11}^{5}$;\ \ 
$ \xi_{11}^{3}$)

  \vskip 2ex

\noindent9. $5_{\frac{38}{7},35.34}^{7,386}$ \irep{11}:\ \ 
$d_i$ = ($1.0$,
$2.246$,
$2.246$,
$2.801$,
$4.48$) 

\vskip 0.7ex
\hangindent=3em \hangafter=1
$D^2= 35.342 = 
21+14c^{1}_{7}
+7c^{2}_{7}
$

\vskip 0.7ex
\hangindent=3em \hangafter=1
$T = ( 0,
\frac{1}{7},
\frac{1}{7},
\frac{6}{7},
\frac{4}{7} )
$,

\vskip 0.7ex
\hangindent=3em \hangafter=1
$S$ = ($ 1$,
$ \xi_{7}^{3}$,
$ \xi_{7}^{3}$,
$ 2+c^{1}_{7}
+c^{2}_{7}
$,
$ 2+2c^{1}_{7}
+c^{2}_{7}
$;\ \ 
$ s^{1}_{7}
+\zeta^{2}_{7}
+\zeta^{3}_{7}
$,
$ -1-2  \zeta^{1}_{7}
-\zeta^{2}_{7}
-\zeta^{3}_{7}
$,
$ -\xi_{7}^{3}$,
$ \xi_{7}^{3}$;\ \ 
$ s^{1}_{7}
+\zeta^{2}_{7}
+\zeta^{3}_{7}
$,
$ -\xi_{7}^{3}$,
$ \xi_{7}^{3}$;\ \ 
$ 2+2c^{1}_{7}
+c^{2}_{7}
$,
$ -1$;\ \ 
$ -2-c^{1}_{7}
-c^{2}_{7}
$)

  \vskip 2ex

\noindent10. $5_{\frac{18}{7},35.34}^{7,101}$ \irep{11}:\ \ 
$d_i$ = ($1.0$,
$2.246$,
$2.246$,
$2.801$,
$4.48$) 

\vskip 0.7ex
\hangindent=3em \hangafter=1
$D^2= 35.342 = 
21+14c^{1}_{7}
+7c^{2}_{7}
$

\vskip 0.7ex
\hangindent=3em \hangafter=1
$T = ( 0,
\frac{6}{7},
\frac{6}{7},
\frac{1}{7},
\frac{3}{7} )
$,

\vskip 0.7ex
\hangindent=3em \hangafter=1
$S$ = ($ 1$,
$ \xi_{7}^{3}$,
$ \xi_{7}^{3}$,
$ 2+c^{1}_{7}
+c^{2}_{7}
$,
$ 2+2c^{1}_{7}
+c^{2}_{7}
$;\ \ 
$ -1-2  \zeta^{1}_{7}
-\zeta^{2}_{7}
-\zeta^{3}_{7}
$,
$ s^{1}_{7}
+\zeta^{2}_{7}
+\zeta^{3}_{7}
$,
$ -\xi_{7}^{3}$,
$ \xi_{7}^{3}$;\ \ 
$ -1-2  \zeta^{1}_{7}
-\zeta^{2}_{7}
-\zeta^{3}_{7}
$,
$ -\xi_{7}^{3}$,
$ \xi_{7}^{3}$;\ \ 
$ 2+2c^{1}_{7}
+c^{2}_{7}
$,
$ -1$;\ \ 
$ -2-c^{1}_{7}
-c^{2}_{7}
$)

  \vskip 2ex 

}

\subsection{Rank 6}
\label{uni6}

{\small

\noindent1. $6_{3,6.}^{12,534}$ \irep{34}:\ \ 
$d_i$ = ($1.0$,
$1.0$,
$1.0$,
$1.0$,
$1.0$,
$1.0$) 

\vskip 0.7ex
\hangindent=3em \hangafter=1
$D^2= 6.0 = 
6$

\vskip 0.7ex
\hangindent=3em \hangafter=1
$T = ( 0,
\frac{1}{3},
\frac{1}{3},
\frac{1}{4},
\frac{7}{12},
\frac{7}{12} )
$,

\vskip 0.7ex
\hangindent=3em \hangafter=1
$S$ = ($ 1$,
$ 1$,
$ 1$,
$ 1$,
$ 1$,
$ 1$;\ \ 
$ \zeta_{3}^{1}$,
$ -\zeta_{6}^{1}$,
$ 1$,
$ -\zeta_{6}^{1}$,
$ \zeta_{3}^{1}$;\ \ 
$ \zeta_{3}^{1}$,
$ 1$,
$ \zeta_{3}^{1}$,
$ -\zeta_{6}^{1}$;\ \ 
$ -1$,
$ -1$,
$ -1$;\ \ 
$ -\zeta_{3}^{1}$,
$ \zeta_{6}^{1}$;\ \ 
$ -\zeta_{3}^{1}$)

Factors = $2_{1,2.}^{4,437}\boxtimes 3_{2,3.}^{3,527}$

  \vskip 2ex

\noindent2. $6_{1,6.}^{12,701}$ \irep{34}:\ \ 
$d_i$ = ($1.0$,
$1.0$,
$1.0$,
$1.0$,
$1.0$,
$1.0$) 

\vskip 0.7ex
\hangindent=3em \hangafter=1
$D^2= 6.0 = 
6$

\vskip 0.7ex
\hangindent=3em \hangafter=1
$T = ( 0,
\frac{1}{3},
\frac{1}{3},
\frac{3}{4},
\frac{1}{12},
\frac{1}{12} )
$,

\vskip 0.7ex
\hangindent=3em \hangafter=1
$S$ = ($ 1$,
$ 1$,
$ 1$,
$ 1$,
$ 1$,
$ 1$;\ \ 
$ \zeta_{3}^{1}$,
$ -\zeta_{6}^{1}$,
$ 1$,
$ -\zeta_{6}^{1}$,
$ \zeta_{3}^{1}$;\ \ 
$ \zeta_{3}^{1}$,
$ 1$,
$ \zeta_{3}^{1}$,
$ -\zeta_{6}^{1}$;\ \ 
$ -1$,
$ -1$,
$ -1$;\ \ 
$ -\zeta_{3}^{1}$,
$ \zeta_{6}^{1}$;\ \ 
$ -\zeta_{3}^{1}$)

Factors = $2_{7,2.}^{4,625}\boxtimes 3_{2,3.}^{3,527}$

  \vskip 2ex

\noindent3. $6_{7,6.}^{12,113}$ \irep{34}:\ \ 
$d_i$ = ($1.0$,
$1.0$,
$1.0$,
$1.0$,
$1.0$,
$1.0$) 

\vskip 0.7ex
\hangindent=3em \hangafter=1
$D^2= 6.0 = 
6$

\vskip 0.7ex
\hangindent=3em \hangafter=1
$T = ( 0,
\frac{2}{3},
\frac{2}{3},
\frac{1}{4},
\frac{11}{12},
\frac{11}{12} )
$,

\vskip 0.7ex
\hangindent=3em \hangafter=1
$S$ = ($ 1$,
$ 1$,
$ 1$,
$ 1$,
$ 1$,
$ 1$;\ \ 
$ -\zeta_{6}^{1}$,
$ \zeta_{3}^{1}$,
$ 1$,
$ -\zeta_{6}^{1}$,
$ \zeta_{3}^{1}$;\ \ 
$ -\zeta_{6}^{1}$,
$ 1$,
$ \zeta_{3}^{1}$,
$ -\zeta_{6}^{1}$;\ \ 
$ -1$,
$ -1$,
$ -1$;\ \ 
$ \zeta_{6}^{1}$,
$ -\zeta_{3}^{1}$;\ \ 
$ \zeta_{6}^{1}$)

Factors = $2_{1,2.}^{4,437}\boxtimes 3_{6,3.}^{3,138}$

  \vskip 2ex

\noindent4. $6_{5,6.}^{12,298}$ \irep{34}:\ \ 
$d_i$ = ($1.0$,
$1.0$,
$1.0$,
$1.0$,
$1.0$,
$1.0$) 

\vskip 0.7ex
\hangindent=3em \hangafter=1
$D^2= 6.0 = 
6$

\vskip 0.7ex
\hangindent=3em \hangafter=1
$T = ( 0,
\frac{2}{3},
\frac{2}{3},
\frac{3}{4},
\frac{5}{12},
\frac{5}{12} )
$,

\vskip 0.7ex
\hangindent=3em \hangafter=1
$S$ = ($ 1$,
$ 1$,
$ 1$,
$ 1$,
$ 1$,
$ 1$;\ \ 
$ -\zeta_{6}^{1}$,
$ \zeta_{3}^{1}$,
$ 1$,
$ -\zeta_{6}^{1}$,
$ \zeta_{3}^{1}$;\ \ 
$ -\zeta_{6}^{1}$,
$ 1$,
$ \zeta_{3}^{1}$,
$ -\zeta_{6}^{1}$;\ \ 
$ -1$,
$ -1$,
$ -1$;\ \ 
$ \zeta_{6}^{1}$,
$ -\zeta_{3}^{1}$;\ \ 
$ \zeta_{6}^{1}$)

Factors = $2_{7,2.}^{4,625}\boxtimes 3_{6,3.}^{3,138}$

  \vskip 2ex

\noindent5. $6_{\frac{3}{2},8.}^{16,688}$ \irep{39}:\ \ 
$d_i$ = ($1.0$,
$1.0$,
$1.0$,
$1.0$,
$1.414$,
$1.414$) 

\vskip 0.7ex
\hangindent=3em \hangafter=1
$D^2= 8.0 = 
8$

\vskip 0.7ex
\hangindent=3em \hangafter=1
$T = ( 0,
\frac{1}{2},
\frac{1}{4},
\frac{3}{4},
\frac{1}{16},
\frac{5}{16} )
$,

\vskip 0.7ex
\hangindent=3em \hangafter=1
$S$ = ($ 1$,
$ 1$,
$ 1$,
$ 1$,
$ \sqrt{2}$,
$ \sqrt{2}$;\ \ 
$ 1$,
$ 1$,
$ 1$,
$ -\sqrt{2}$,
$ -\sqrt{2}$;\ \ 
$ -1$,
$ -1$,
$ \sqrt{2}$,
$ -\sqrt{2}$;\ \ 
$ -1$,
$ -\sqrt{2}$,
$ \sqrt{2}$;\ \ 
$0$,
$0$;\ \ 
$0$)

Factors = $2_{1,2.}^{4,437}\boxtimes 3_{\frac{1}{2},4.}^{16,598}$

  \vskip 2ex

\noindent6. $6_{\frac{15}{2},8.}^{16,107}$ \irep{39}:\ \ 
$d_i$ = ($1.0$,
$1.0$,
$1.0$,
$1.0$,
$1.414$,
$1.414$) 

\vskip 0.7ex
\hangindent=3em \hangafter=1
$D^2= 8.0 = 
8$

\vskip 0.7ex
\hangindent=3em \hangafter=1
$T = ( 0,
\frac{1}{2},
\frac{1}{4},
\frac{3}{4},
\frac{1}{16},
\frac{13}{16} )
$,

\vskip 0.7ex
\hangindent=3em \hangafter=1
$S$ = ($ 1$,
$ 1$,
$ 1$,
$ 1$,
$ \sqrt{2}$,
$ \sqrt{2}$;\ \ 
$ 1$,
$ 1$,
$ 1$,
$ -\sqrt{2}$,
$ -\sqrt{2}$;\ \ 
$ -1$,
$ -1$,
$ -\sqrt{2}$,
$ \sqrt{2}$;\ \ 
$ -1$,
$ \sqrt{2}$,
$ -\sqrt{2}$;\ \ 
$0$,
$0$;\ \ 
$0$)

Factors = $2_{7,2.}^{4,625}\boxtimes 3_{\frac{1}{2},4.}^{16,598}$

  \vskip 2ex

\noindent7. $6_{\frac{5}{2},8.}^{16,511}$ \irep{39}:\ \ 
$d_i$ = ($1.0$,
$1.0$,
$1.0$,
$1.0$,
$1.414$,
$1.414$) 

\vskip 0.7ex
\hangindent=3em \hangafter=1
$D^2= 8.0 = 
8$

\vskip 0.7ex
\hangindent=3em \hangafter=1
$T = ( 0,
\frac{1}{2},
\frac{1}{4},
\frac{3}{4},
\frac{3}{16},
\frac{7}{16} )
$,

\vskip 0.7ex
\hangindent=3em \hangafter=1
$S$ = ($ 1$,
$ 1$,
$ 1$,
$ 1$,
$ \sqrt{2}$,
$ \sqrt{2}$;\ \ 
$ 1$,
$ 1$,
$ 1$,
$ -\sqrt{2}$,
$ -\sqrt{2}$;\ \ 
$ -1$,
$ -1$,
$ \sqrt{2}$,
$ -\sqrt{2}$;\ \ 
$ -1$,
$ -\sqrt{2}$,
$ \sqrt{2}$;\ \ 
$0$,
$0$;\ \ 
$0$)

Factors = $2_{7,2.}^{4,625}\boxtimes 3_{\frac{7}{2},4.}^{16,332}$

  \vskip 2ex

\noindent8. $6_{\frac{1}{2},8.}^{16,460}$ \irep{39}:\ \ 
$d_i$ = ($1.0$,
$1.0$,
$1.0$,
$1.0$,
$1.414$,
$1.414$) 

\vskip 0.7ex
\hangindent=3em \hangafter=1
$D^2= 8.0 = 
8$

\vskip 0.7ex
\hangindent=3em \hangafter=1
$T = ( 0,
\frac{1}{2},
\frac{1}{4},
\frac{3}{4},
\frac{3}{16},
\frac{15}{16} )
$,

\vskip 0.7ex
\hangindent=3em \hangafter=1
$S$ = ($ 1$,
$ 1$,
$ 1$,
$ 1$,
$ \sqrt{2}$,
$ \sqrt{2}$;\ \ 
$ 1$,
$ 1$,
$ 1$,
$ -\sqrt{2}$,
$ -\sqrt{2}$;\ \ 
$ -1$,
$ -1$,
$ -\sqrt{2}$,
$ \sqrt{2}$;\ \ 
$ -1$,
$ \sqrt{2}$,
$ -\sqrt{2}$;\ \ 
$0$,
$0$;\ \ 
$0$)

Factors = $2_{1,2.}^{4,437}\boxtimes 3_{\frac{15}{2},4.}^{16,639}$

  \vskip 2ex

\noindent9. $6_{\frac{7}{2},8.}^{16,246}$ \irep{39}:\ \ 
$d_i$ = ($1.0$,
$1.0$,
$1.0$,
$1.0$,
$1.414$,
$1.414$) 

\vskip 0.7ex
\hangindent=3em \hangafter=1
$D^2= 8.0 = 
8$

\vskip 0.7ex
\hangindent=3em \hangafter=1
$T = ( 0,
\frac{1}{2},
\frac{1}{4},
\frac{3}{4},
\frac{5}{16},
\frac{9}{16} )
$,

\vskip 0.7ex
\hangindent=3em \hangafter=1
$S$ = ($ 1$,
$ 1$,
$ 1$,
$ 1$,
$ \sqrt{2}$,
$ \sqrt{2}$;\ \ 
$ 1$,
$ 1$,
$ 1$,
$ -\sqrt{2}$,
$ -\sqrt{2}$;\ \ 
$ -1$,
$ -1$,
$ \sqrt{2}$,
$ -\sqrt{2}$;\ \ 
$ -1$,
$ -\sqrt{2}$,
$ \sqrt{2}$;\ \ 
$0$,
$0$;\ \ 
$0$)

Factors = $2_{7,2.}^{4,625}\boxtimes 3_{\frac{9}{2},4.}^{16,156}$

  \vskip 2ex

\noindent10. $6_{\frac{9}{2},8.}^{16,107}$ \irep{39}:\ \ 
$d_i$ = ($1.0$,
$1.0$,
$1.0$,
$1.0$,
$1.414$,
$1.414$) 

\vskip 0.7ex
\hangindent=3em \hangafter=1
$D^2= 8.0 = 
8$

\vskip 0.7ex
\hangindent=3em \hangafter=1
$T = ( 0,
\frac{1}{2},
\frac{1}{4},
\frac{3}{4},
\frac{7}{16},
\frac{11}{16} )
$,

\vskip 0.7ex
\hangindent=3em \hangafter=1
$S$ = ($ 1$,
$ 1$,
$ 1$,
$ 1$,
$ \sqrt{2}$,
$ \sqrt{2}$;\ \ 
$ 1$,
$ 1$,
$ 1$,
$ -\sqrt{2}$,
$ -\sqrt{2}$;\ \ 
$ -1$,
$ -1$,
$ \sqrt{2}$,
$ -\sqrt{2}$;\ \ 
$ -1$,
$ -\sqrt{2}$,
$ \sqrt{2}$;\ \ 
$0$,
$0$;\ \ 
$0$)

Factors = $2_{1,2.}^{4,437}\boxtimes 3_{\frac{7}{2},4.}^{16,332}$

  \vskip 2ex

\noindent11. $6_{\frac{11}{2},8.}^{16,548}$ \irep{39}:\ \ 
$d_i$ = ($1.0$,
$1.0$,
$1.0$,
$1.0$,
$1.414$,
$1.414$) 

\vskip 0.7ex
\hangindent=3em \hangafter=1
$D^2= 8.0 = 
8$

\vskip 0.7ex
\hangindent=3em \hangafter=1
$T = ( 0,
\frac{1}{2},
\frac{1}{4},
\frac{3}{4},
\frac{9}{16},
\frac{13}{16} )
$,

\vskip 0.7ex
\hangindent=3em \hangafter=1
$S$ = ($ 1$,
$ 1$,
$ 1$,
$ 1$,
$ \sqrt{2}$,
$ \sqrt{2}$;\ \ 
$ 1$,
$ 1$,
$ 1$,
$ -\sqrt{2}$,
$ -\sqrt{2}$;\ \ 
$ -1$,
$ -1$,
$ \sqrt{2}$,
$ -\sqrt{2}$;\ \ 
$ -1$,
$ -\sqrt{2}$,
$ \sqrt{2}$;\ \ 
$0$,
$0$;\ \ 
$0$)

Factors = $2_{1,2.}^{4,437}\boxtimes 3_{\frac{9}{2},4.}^{16,156}$

  \vskip 2ex

\noindent12. $6_{\frac{13}{2},8.}^{16,107}$ \irep{39}:\ \ 
$d_i$ = ($1.0$,
$1.0$,
$1.0$,
$1.0$,
$1.414$,
$1.414$) 

\vskip 0.7ex
\hangindent=3em \hangafter=1
$D^2= 8.0 = 
8$

\vskip 0.7ex
\hangindent=3em \hangafter=1
$T = ( 0,
\frac{1}{2},
\frac{1}{4},
\frac{3}{4},
\frac{11}{16},
\frac{15}{16} )
$,

\vskip 0.7ex
\hangindent=3em \hangafter=1
$S$ = ($ 1$,
$ 1$,
$ 1$,
$ 1$,
$ \sqrt{2}$,
$ \sqrt{2}$;\ \ 
$ 1$,
$ 1$,
$ 1$,
$ -\sqrt{2}$,
$ -\sqrt{2}$;\ \ 
$ -1$,
$ -1$,
$ \sqrt{2}$,
$ -\sqrt{2}$;\ \ 
$ -1$,
$ -\sqrt{2}$,
$ \sqrt{2}$;\ \ 
$0$,
$0$;\ \ 
$0$)

Factors = $2_{7,2.}^{4,625}\boxtimes 3_{\frac{15}{2},4.}^{16,639}$

  \vskip 2ex

\noindent13. $6_{\frac{24}{5},10.85}^{15,257}$ \irep{38}:\ \ 
$d_i$ = ($1.0$,
$1.0$,
$1.0$,
$1.618$,
$1.618$,
$1.618$) 

\vskip 0.7ex
\hangindent=3em \hangafter=1
$D^2= 10.854 = 
\frac{15+3\sqrt{5}}{2}$

\vskip 0.7ex
\hangindent=3em \hangafter=1
$T = ( 0,
\frac{1}{3},
\frac{1}{3},
\frac{2}{5},
\frac{11}{15},
\frac{11}{15} )
$,

\vskip 0.7ex
\hangindent=3em \hangafter=1
$S$ = ($ 1$,
$ 1$,
$ 1$,
$ \frac{1+\sqrt{5}}{2}$,
$ \frac{1+\sqrt{5}}{2}$,
$ \frac{1+\sqrt{5}}{2}$;\ \ 
$ \zeta_{3}^{1}$,
$ -\zeta_{6}^{1}$,
$ \frac{1+\sqrt{5}}{2}$,
$ -\frac{1+\sqrt{5}}{2}\zeta_{6}^{1}$,
$ \frac{1+\sqrt{5}}{2}\zeta_{3}^{1}$;\ \ 
$ \zeta_{3}^{1}$,
$ \frac{1+\sqrt{5}}{2}$,
$ \frac{1+\sqrt{5}}{2}\zeta_{3}^{1}$,
$ -\frac{1+\sqrt{5}}{2}\zeta_{6}^{1}$;\ \ 
$ -1$,
$ -1$,
$ -1$;\ \ 
$ -\zeta_{3}^{1}$,
$ \zeta_{6}^{1}$;\ \ 
$ -\zeta_{3}^{1}$)

Factors = $2_{\frac{14}{5},3.618}^{5,395}\boxtimes 3_{2,3.}^{3,527}$

  \vskip 2ex

\noindent14. $6_{\frac{36}{5},10.85}^{15,166}$ \irep{38}:\ \ 
$d_i$ = ($1.0$,
$1.0$,
$1.0$,
$1.618$,
$1.618$,
$1.618$) 

\vskip 0.7ex
\hangindent=3em \hangafter=1
$D^2= 10.854 = 
\frac{15+3\sqrt{5}}{2}$

\vskip 0.7ex
\hangindent=3em \hangafter=1
$T = ( 0,
\frac{1}{3},
\frac{1}{3},
\frac{3}{5},
\frac{14}{15},
\frac{14}{15} )
$,

\vskip 0.7ex
\hangindent=3em \hangafter=1
$S$ = ($ 1$,
$ 1$,
$ 1$,
$ \frac{1+\sqrt{5}}{2}$,
$ \frac{1+\sqrt{5}}{2}$,
$ \frac{1+\sqrt{5}}{2}$;\ \ 
$ \zeta_{3}^{1}$,
$ -\zeta_{6}^{1}$,
$ \frac{1+\sqrt{5}}{2}$,
$ -\frac{1+\sqrt{5}}{2}\zeta_{6}^{1}$,
$ \frac{1+\sqrt{5}}{2}\zeta_{3}^{1}$;\ \ 
$ \zeta_{3}^{1}$,
$ \frac{1+\sqrt{5}}{2}$,
$ \frac{1+\sqrt{5}}{2}\zeta_{3}^{1}$,
$ -\frac{1+\sqrt{5}}{2}\zeta_{6}^{1}$;\ \ 
$ -1$,
$ -1$,
$ -1$;\ \ 
$ -\zeta_{3}^{1}$,
$ \zeta_{6}^{1}$;\ \ 
$ -\zeta_{3}^{1}$)

Factors = $2_{\frac{26}{5},3.618}^{5,720}\boxtimes 3_{2,3.}^{3,527}$

  \vskip 2ex

\noindent15. $6_{\frac{4}{5},10.85}^{15,801}$ \irep{38}:\ \ 
$d_i$ = ($1.0$,
$1.0$,
$1.0$,
$1.618$,
$1.618$,
$1.618$) 

\vskip 0.7ex
\hangindent=3em \hangafter=1
$D^2= 10.854 = 
\frac{15+3\sqrt{5}}{2}$

\vskip 0.7ex
\hangindent=3em \hangafter=1
$T = ( 0,
\frac{2}{3},
\frac{2}{3},
\frac{2}{5},
\frac{1}{15},
\frac{1}{15} )
$,

\vskip 0.7ex
\hangindent=3em \hangafter=1
$S$ = ($ 1$,
$ 1$,
$ 1$,
$ \frac{1+\sqrt{5}}{2}$,
$ \frac{1+\sqrt{5}}{2}$,
$ \frac{1+\sqrt{5}}{2}$;\ \ 
$ -\zeta_{6}^{1}$,
$ \zeta_{3}^{1}$,
$ \frac{1+\sqrt{5}}{2}$,
$ -\frac{1+\sqrt{5}}{2}\zeta_{6}^{1}$,
$ \frac{1+\sqrt{5}}{2}\zeta_{3}^{1}$;\ \ 
$ -\zeta_{6}^{1}$,
$ \frac{1+\sqrt{5}}{2}$,
$ \frac{1+\sqrt{5}}{2}\zeta_{3}^{1}$,
$ -\frac{1+\sqrt{5}}{2}\zeta_{6}^{1}$;\ \ 
$ -1$,
$ -1$,
$ -1$;\ \ 
$ \zeta_{6}^{1}$,
$ -\zeta_{3}^{1}$;\ \ 
$ \zeta_{6}^{1}$)

Factors = $2_{\frac{14}{5},3.618}^{5,395}\boxtimes 3_{6,3.}^{3,138}$

  \vskip 2ex

\noindent16. $6_{\frac{16}{5},10.85}^{15,262}$ \irep{38}:\ \ 
$d_i$ = ($1.0$,
$1.0$,
$1.0$,
$1.618$,
$1.618$,
$1.618$) 

\vskip 0.7ex
\hangindent=3em \hangafter=1
$D^2= 10.854 = 
\frac{15+3\sqrt{5}}{2}$

\vskip 0.7ex
\hangindent=3em \hangafter=1
$T = ( 0,
\frac{2}{3},
\frac{2}{3},
\frac{3}{5},
\frac{4}{15},
\frac{4}{15} )
$,

\vskip 0.7ex
\hangindent=3em \hangafter=1
$S$ = ($ 1$,
$ 1$,
$ 1$,
$ \frac{1+\sqrt{5}}{2}$,
$ \frac{1+\sqrt{5}}{2}$,
$ \frac{1+\sqrt{5}}{2}$;\ \ 
$ -\zeta_{6}^{1}$,
$ \zeta_{3}^{1}$,
$ \frac{1+\sqrt{5}}{2}$,
$ -\frac{1+\sqrt{5}}{2}\zeta_{6}^{1}$,
$ \frac{1+\sqrt{5}}{2}\zeta_{3}^{1}$;\ \ 
$ -\zeta_{6}^{1}$,
$ \frac{1+\sqrt{5}}{2}$,
$ \frac{1+\sqrt{5}}{2}\zeta_{3}^{1}$,
$ -\frac{1+\sqrt{5}}{2}\zeta_{6}^{1}$;\ \ 
$ -1$,
$ -1$,
$ -1$;\ \ 
$ \zeta_{6}^{1}$,
$ -\zeta_{3}^{1}$;\ \ 
$ \zeta_{6}^{1}$)

Factors = $2_{\frac{26}{5},3.618}^{5,720}\boxtimes 3_{6,3.}^{3,138}$

  \vskip 2ex

\noindent17. $6_{\frac{33}{10},14.47}^{80,798}$ \irep{48}:\ \ 
$d_i$ = ($1.0$,
$1.0$,
$1.414$,
$1.618$,
$1.618$,
$2.288$) 

\vskip 0.7ex
\hangindent=3em \hangafter=1
$D^2= 14.472 = 
10+2\sqrt{5}$

\vskip 0.7ex
\hangindent=3em \hangafter=1
$T = ( 0,
\frac{1}{2},
\frac{1}{16},
\frac{2}{5},
\frac{9}{10},
\frac{37}{80} )
$,

\vskip 0.7ex
\hangindent=3em \hangafter=1
$S$ = ($ 1$,
$ 1$,
$ \sqrt{2}$,
$ \frac{1+\sqrt{5}}{2}$,
$ \frac{1+\sqrt{5}}{2}$,
$ \frac{5+\sqrt{5}}{\sqrt{10}}$;\ \ 
$ 1$,
$ -\sqrt{2}$,
$ \frac{1+\sqrt{5}}{2}$,
$ \frac{1+\sqrt{5}}{2}$,
$ \frac{-5-\sqrt{5}}{\sqrt{10}}$;\ \ 
$0$,
$ \frac{5+\sqrt{5}}{\sqrt{10}}$,
$ \frac{-5-\sqrt{5}}{\sqrt{10}}$,
$0$;\ \ 
$ -1$,
$ -1$,
$ -\sqrt{2}$;\ \ 
$ -1$,
$ \sqrt{2}$;\ \ 
$0$)

Factors = $2_{\frac{14}{5},3.618}^{5,395}\boxtimes 3_{\frac{1}{2},4.}^{16,598}$

  \vskip 2ex

\noindent18. $6_{\frac{57}{10},14.47}^{80,376}$ \irep{48}:\ \ 
$d_i$ = ($1.0$,
$1.0$,
$1.414$,
$1.618$,
$1.618$,
$2.288$) 

\vskip 0.7ex
\hangindent=3em \hangafter=1
$D^2= 14.472 = 
10+2\sqrt{5}$

\vskip 0.7ex
\hangindent=3em \hangafter=1
$T = ( 0,
\frac{1}{2},
\frac{1}{16},
\frac{3}{5},
\frac{1}{10},
\frac{53}{80} )
$,

\vskip 0.7ex
\hangindent=3em \hangafter=1
$S$ = ($ 1$,
$ 1$,
$ \sqrt{2}$,
$ \frac{1+\sqrt{5}}{2}$,
$ \frac{1+\sqrt{5}}{2}$,
$ \frac{5+\sqrt{5}}{\sqrt{10}}$;\ \ 
$ 1$,
$ -\sqrt{2}$,
$ \frac{1+\sqrt{5}}{2}$,
$ \frac{1+\sqrt{5}}{2}$,
$ \frac{-5-\sqrt{5}}{\sqrt{10}}$;\ \ 
$0$,
$ \frac{5+\sqrt{5}}{\sqrt{10}}$,
$ \frac{-5-\sqrt{5}}{\sqrt{10}}$,
$0$;\ \ 
$ -1$,
$ -1$,
$ -\sqrt{2}$;\ \ 
$ -1$,
$ \sqrt{2}$;\ \ 
$0$)

Factors = $2_{\frac{26}{5},3.618}^{5,720}\boxtimes 3_{\frac{1}{2},4.}^{16,598}$

  \vskip 2ex

\noindent19. $6_{\frac{43}{10},14.47}^{80,424}$ \irep{48}:\ \ 
$d_i$ = ($1.0$,
$1.0$,
$1.414$,
$1.618$,
$1.618$,
$2.288$) 

\vskip 0.7ex
\hangindent=3em \hangafter=1
$D^2= 14.472 = 
10+2\sqrt{5}$

\vskip 0.7ex
\hangindent=3em \hangafter=1
$T = ( 0,
\frac{1}{2},
\frac{3}{16},
\frac{2}{5},
\frac{9}{10},
\frac{47}{80} )
$,

\vskip 0.7ex
\hangindent=3em \hangafter=1
$S$ = ($ 1$,
$ 1$,
$ \sqrt{2}$,
$ \frac{1+\sqrt{5}}{2}$,
$ \frac{1+\sqrt{5}}{2}$,
$ \frac{5+\sqrt{5}}{\sqrt{10}}$;\ \ 
$ 1$,
$ -\sqrt{2}$,
$ \frac{1+\sqrt{5}}{2}$,
$ \frac{1+\sqrt{5}}{2}$,
$ \frac{-5-\sqrt{5}}{\sqrt{10}}$;\ \ 
$0$,
$ \frac{5+\sqrt{5}}{\sqrt{10}}$,
$ \frac{-5-\sqrt{5}}{\sqrt{10}}$,
$0$;\ \ 
$ -1$,
$ -1$,
$ -\sqrt{2}$;\ \ 
$ -1$,
$ \sqrt{2}$;\ \ 
$0$)

Factors = $2_{\frac{14}{5},3.618}^{5,395}\boxtimes 3_{\frac{3}{2},4.}^{16,553}$

  \vskip 2ex

\noindent20. $6_{\frac{67}{10},14.47}^{80,828}$ \irep{48}:\ \ 
$d_i$ = ($1.0$,
$1.0$,
$1.414$,
$1.618$,
$1.618$,
$2.288$) 

\vskip 0.7ex
\hangindent=3em \hangafter=1
$D^2= 14.472 = 
10+2\sqrt{5}$

\vskip 0.7ex
\hangindent=3em \hangafter=1
$T = ( 0,
\frac{1}{2},
\frac{3}{16},
\frac{3}{5},
\frac{1}{10},
\frac{63}{80} )
$,

\vskip 0.7ex
\hangindent=3em \hangafter=1
$S$ = ($ 1$,
$ 1$,
$ \sqrt{2}$,
$ \frac{1+\sqrt{5}}{2}$,
$ \frac{1+\sqrt{5}}{2}$,
$ \frac{5+\sqrt{5}}{\sqrt{10}}$;\ \ 
$ 1$,
$ -\sqrt{2}$,
$ \frac{1+\sqrt{5}}{2}$,
$ \frac{1+\sqrt{5}}{2}$,
$ \frac{-5-\sqrt{5}}{\sqrt{10}}$;\ \ 
$0$,
$ \frac{5+\sqrt{5}}{\sqrt{10}}$,
$ \frac{-5-\sqrt{5}}{\sqrt{10}}$,
$0$;\ \ 
$ -1$,
$ -1$,
$ -\sqrt{2}$;\ \ 
$ -1$,
$ \sqrt{2}$;\ \ 
$0$)

Factors = $2_{\frac{26}{5},3.618}^{5,720}\boxtimes 3_{\frac{3}{2},4.}^{16,553}$

  \vskip 2ex

\noindent21. $6_{\frac{53}{10},14.47}^{80,884}$ \irep{48}:\ \ 
$d_i$ = ($1.0$,
$1.0$,
$1.414$,
$1.618$,
$1.618$,
$2.288$) 

\vskip 0.7ex
\hangindent=3em \hangafter=1
$D^2= 14.472 = 
10+2\sqrt{5}$

\vskip 0.7ex
\hangindent=3em \hangafter=1
$T = ( 0,
\frac{1}{2},
\frac{5}{16},
\frac{2}{5},
\frac{9}{10},
\frac{57}{80} )
$,

\vskip 0.7ex
\hangindent=3em \hangafter=1
$S$ = ($ 1$,
$ 1$,
$ \sqrt{2}$,
$ \frac{1+\sqrt{5}}{2}$,
$ \frac{1+\sqrt{5}}{2}$,
$ \frac{5+\sqrt{5}}{\sqrt{10}}$;\ \ 
$ 1$,
$ -\sqrt{2}$,
$ \frac{1+\sqrt{5}}{2}$,
$ \frac{1+\sqrt{5}}{2}$,
$ \frac{-5-\sqrt{5}}{\sqrt{10}}$;\ \ 
$0$,
$ \frac{5+\sqrt{5}}{\sqrt{10}}$,
$ \frac{-5-\sqrt{5}}{\sqrt{10}}$,
$0$;\ \ 
$ -1$,
$ -1$,
$ -\sqrt{2}$;\ \ 
$ -1$,
$ \sqrt{2}$;\ \ 
$0$)

Factors = $2_{\frac{14}{5},3.618}^{5,395}\boxtimes 3_{\frac{5}{2},4.}^{16,465}$

  \vskip 2ex

\noindent22. $6_{\frac{77}{10},14.47}^{80,657}$ \irep{48}:\ \ 
$d_i$ = ($1.0$,
$1.0$,
$1.414$,
$1.618$,
$1.618$,
$2.288$) 

\vskip 0.7ex
\hangindent=3em \hangafter=1
$D^2= 14.472 = 
10+2\sqrt{5}$

\vskip 0.7ex
\hangindent=3em \hangafter=1
$T = ( 0,
\frac{1}{2},
\frac{5}{16},
\frac{3}{5},
\frac{1}{10},
\frac{73}{80} )
$,

\vskip 0.7ex
\hangindent=3em \hangafter=1
$S$ = ($ 1$,
$ 1$,
$ \sqrt{2}$,
$ \frac{1+\sqrt{5}}{2}$,
$ \frac{1+\sqrt{5}}{2}$,
$ \frac{5+\sqrt{5}}{\sqrt{10}}$;\ \ 
$ 1$,
$ -\sqrt{2}$,
$ \frac{1+\sqrt{5}}{2}$,
$ \frac{1+\sqrt{5}}{2}$,
$ \frac{-5-\sqrt{5}}{\sqrt{10}}$;\ \ 
$0$,
$ \frac{5+\sqrt{5}}{\sqrt{10}}$,
$ \frac{-5-\sqrt{5}}{\sqrt{10}}$,
$0$;\ \ 
$ -1$,
$ -1$,
$ -\sqrt{2}$;\ \ 
$ -1$,
$ \sqrt{2}$;\ \ 
$0$)

Factors = $2_{\frac{26}{5},3.618}^{5,720}\boxtimes 3_{\frac{5}{2},4.}^{16,465}$

  \vskip 2ex

\noindent23. $6_{\frac{63}{10},14.47}^{80,146}$ \irep{48}:\ \ 
$d_i$ = ($1.0$,
$1.0$,
$1.414$,
$1.618$,
$1.618$,
$2.288$) 

\vskip 0.7ex
\hangindent=3em \hangafter=1
$D^2= 14.472 = 
10+2\sqrt{5}$

\vskip 0.7ex
\hangindent=3em \hangafter=1
$T = ( 0,
\frac{1}{2},
\frac{7}{16},
\frac{2}{5},
\frac{9}{10},
\frac{67}{80} )
$,

\vskip 0.7ex
\hangindent=3em \hangafter=1
$S$ = ($ 1$,
$ 1$,
$ \sqrt{2}$,
$ \frac{1+\sqrt{5}}{2}$,
$ \frac{1+\sqrt{5}}{2}$,
$ \frac{5+\sqrt{5}}{\sqrt{10}}$;\ \ 
$ 1$,
$ -\sqrt{2}$,
$ \frac{1+\sqrt{5}}{2}$,
$ \frac{1+\sqrt{5}}{2}$,
$ \frac{-5-\sqrt{5}}{\sqrt{10}}$;\ \ 
$0$,
$ \frac{5+\sqrt{5}}{\sqrt{10}}$,
$ \frac{-5-\sqrt{5}}{\sqrt{10}}$,
$0$;\ \ 
$ -1$,
$ -1$,
$ -\sqrt{2}$;\ \ 
$ -1$,
$ \sqrt{2}$;\ \ 
$0$)

Factors = $2_{\frac{14}{5},3.618}^{5,395}\boxtimes 3_{\frac{7}{2},4.}^{16,332}$

  \vskip 2ex

\noindent24. $6_{\frac{7}{10},14.47}^{80,111}$ \irep{48}:\ \ 
$d_i$ = ($1.0$,
$1.0$,
$1.414$,
$1.618$,
$1.618$,
$2.288$) 

\vskip 0.7ex
\hangindent=3em \hangafter=1
$D^2= 14.472 = 
10+2\sqrt{5}$

\vskip 0.7ex
\hangindent=3em \hangafter=1
$T = ( 0,
\frac{1}{2},
\frac{7}{16},
\frac{3}{5},
\frac{1}{10},
\frac{3}{80} )
$,

\vskip 0.7ex
\hangindent=3em \hangafter=1
$S$ = ($ 1$,
$ 1$,
$ \sqrt{2}$,
$ \frac{1+\sqrt{5}}{2}$,
$ \frac{1+\sqrt{5}}{2}$,
$ \frac{5+\sqrt{5}}{\sqrt{10}}$;\ \ 
$ 1$,
$ -\sqrt{2}$,
$ \frac{1+\sqrt{5}}{2}$,
$ \frac{1+\sqrt{5}}{2}$,
$ \frac{-5-\sqrt{5}}{\sqrt{10}}$;\ \ 
$0$,
$ \frac{5+\sqrt{5}}{\sqrt{10}}$,
$ \frac{-5-\sqrt{5}}{\sqrt{10}}$,
$0$;\ \ 
$ -1$,
$ -1$,
$ -\sqrt{2}$;\ \ 
$ -1$,
$ \sqrt{2}$;\ \ 
$0$)

Factors = $2_{\frac{26}{5},3.618}^{5,720}\boxtimes 3_{\frac{7}{2},4.}^{16,332}$

  \vskip 2ex

\noindent25. $6_{\frac{73}{10},14.47}^{80,215}$ \irep{48}:\ \ 
$d_i$ = ($1.0$,
$1.0$,
$1.414$,
$1.618$,
$1.618$,
$2.288$) 

\vskip 0.7ex
\hangindent=3em \hangafter=1
$D^2= 14.472 = 
10+2\sqrt{5}$

\vskip 0.7ex
\hangindent=3em \hangafter=1
$T = ( 0,
\frac{1}{2},
\frac{9}{16},
\frac{2}{5},
\frac{9}{10},
\frac{77}{80} )
$,

\vskip 0.7ex
\hangindent=3em \hangafter=1
$S$ = ($ 1$,
$ 1$,
$ \sqrt{2}$,
$ \frac{1+\sqrt{5}}{2}$,
$ \frac{1+\sqrt{5}}{2}$,
$ \frac{5+\sqrt{5}}{\sqrt{10}}$;\ \ 
$ 1$,
$ -\sqrt{2}$,
$ \frac{1+\sqrt{5}}{2}$,
$ \frac{1+\sqrt{5}}{2}$,
$ \frac{-5-\sqrt{5}}{\sqrt{10}}$;\ \ 
$0$,
$ \frac{5+\sqrt{5}}{\sqrt{10}}$,
$ \frac{-5-\sqrt{5}}{\sqrt{10}}$,
$0$;\ \ 
$ -1$,
$ -1$,
$ -\sqrt{2}$;\ \ 
$ -1$,
$ \sqrt{2}$;\ \ 
$0$)

Factors = $2_{\frac{14}{5},3.618}^{5,395}\boxtimes 3_{\frac{9}{2},4.}^{16,156}$

  \vskip 2ex

\noindent26. $6_{\frac{17}{10},14.47}^{80,878}$ \irep{48}:\ \ 
$d_i$ = ($1.0$,
$1.0$,
$1.414$,
$1.618$,
$1.618$,
$2.288$) 

\vskip 0.7ex
\hangindent=3em \hangafter=1
$D^2= 14.472 = 
10+2\sqrt{5}$

\vskip 0.7ex
\hangindent=3em \hangafter=1
$T = ( 0,
\frac{1}{2},
\frac{9}{16},
\frac{3}{5},
\frac{1}{10},
\frac{13}{80} )
$,

\vskip 0.7ex
\hangindent=3em \hangafter=1
$S$ = ($ 1$,
$ 1$,
$ \sqrt{2}$,
$ \frac{1+\sqrt{5}}{2}$,
$ \frac{1+\sqrt{5}}{2}$,
$ \frac{5+\sqrt{5}}{\sqrt{10}}$;\ \ 
$ 1$,
$ -\sqrt{2}$,
$ \frac{1+\sqrt{5}}{2}$,
$ \frac{1+\sqrt{5}}{2}$,
$ \frac{-5-\sqrt{5}}{\sqrt{10}}$;\ \ 
$0$,
$ \frac{5+\sqrt{5}}{\sqrt{10}}$,
$ \frac{-5-\sqrt{5}}{\sqrt{10}}$,
$0$;\ \ 
$ -1$,
$ -1$,
$ -\sqrt{2}$;\ \ 
$ -1$,
$ \sqrt{2}$;\ \ 
$0$)

Factors = $2_{\frac{26}{5},3.618}^{5,720}\boxtimes 3_{\frac{9}{2},4.}^{16,156}$

  \vskip 2ex

\noindent27. $6_{\frac{3}{10},14.47}^{80,270}$ \irep{48}:\ \ 
$d_i$ = ($1.0$,
$1.0$,
$1.414$,
$1.618$,
$1.618$,
$2.288$) 

\vskip 0.7ex
\hangindent=3em \hangafter=1
$D^2= 14.472 = 
10+2\sqrt{5}$

\vskip 0.7ex
\hangindent=3em \hangafter=1
$T = ( 0,
\frac{1}{2},
\frac{11}{16},
\frac{2}{5},
\frac{9}{10},
\frac{7}{80} )
$,

\vskip 0.7ex
\hangindent=3em \hangafter=1
$S$ = ($ 1$,
$ 1$,
$ \sqrt{2}$,
$ \frac{1+\sqrt{5}}{2}$,
$ \frac{1+\sqrt{5}}{2}$,
$ \frac{5+\sqrt{5}}{\sqrt{10}}$;\ \ 
$ 1$,
$ -\sqrt{2}$,
$ \frac{1+\sqrt{5}}{2}$,
$ \frac{1+\sqrt{5}}{2}$,
$ \frac{-5-\sqrt{5}}{\sqrt{10}}$;\ \ 
$0$,
$ \frac{5+\sqrt{5}}{\sqrt{10}}$,
$ \frac{-5-\sqrt{5}}{\sqrt{10}}$,
$0$;\ \ 
$ -1$,
$ -1$,
$ -\sqrt{2}$;\ \ 
$ -1$,
$ \sqrt{2}$;\ \ 
$0$)

Factors = $2_{\frac{14}{5},3.618}^{5,395}\boxtimes 3_{\frac{11}{2},4.}^{16,648}$

  \vskip 2ex

\noindent28. $6_{\frac{27}{10},14.47}^{80,528}$ \irep{48}:\ \ 
$d_i$ = ($1.0$,
$1.0$,
$1.414$,
$1.618$,
$1.618$,
$2.288$) 

\vskip 0.7ex
\hangindent=3em \hangafter=1
$D^2= 14.472 = 
10+2\sqrt{5}$

\vskip 0.7ex
\hangindent=3em \hangafter=1
$T = ( 0,
\frac{1}{2},
\frac{11}{16},
\frac{3}{5},
\frac{1}{10},
\frac{23}{80} )
$,

\vskip 0.7ex
\hangindent=3em \hangafter=1
$S$ = ($ 1$,
$ 1$,
$ \sqrt{2}$,
$ \frac{1+\sqrt{5}}{2}$,
$ \frac{1+\sqrt{5}}{2}$,
$ \frac{5+\sqrt{5}}{\sqrt{10}}$;\ \ 
$ 1$,
$ -\sqrt{2}$,
$ \frac{1+\sqrt{5}}{2}$,
$ \frac{1+\sqrt{5}}{2}$,
$ \frac{-5-\sqrt{5}}{\sqrt{10}}$;\ \ 
$0$,
$ \frac{5+\sqrt{5}}{\sqrt{10}}$,
$ \frac{-5-\sqrt{5}}{\sqrt{10}}$,
$0$;\ \ 
$ -1$,
$ -1$,
$ -\sqrt{2}$;\ \ 
$ -1$,
$ \sqrt{2}$;\ \ 
$0$)

Factors = $2_{\frac{26}{5},3.618}^{5,720}\boxtimes 3_{\frac{11}{2},4.}^{16,648}$

  \vskip 2ex

\noindent29. $6_{\frac{13}{10},14.47}^{80,621}$ \irep{48}:\ \ 
$d_i$ = ($1.0$,
$1.0$,
$1.414$,
$1.618$,
$1.618$,
$2.288$) 

\vskip 0.7ex
\hangindent=3em \hangafter=1
$D^2= 14.472 = 
10+2\sqrt{5}$

\vskip 0.7ex
\hangindent=3em \hangafter=1
$T = ( 0,
\frac{1}{2},
\frac{13}{16},
\frac{2}{5},
\frac{9}{10},
\frac{17}{80} )
$,

\vskip 0.7ex
\hangindent=3em \hangafter=1
$S$ = ($ 1$,
$ 1$,
$ \sqrt{2}$,
$ \frac{1+\sqrt{5}}{2}$,
$ \frac{1+\sqrt{5}}{2}$,
$ \frac{5+\sqrt{5}}{\sqrt{10}}$;\ \ 
$ 1$,
$ -\sqrt{2}$,
$ \frac{1+\sqrt{5}}{2}$,
$ \frac{1+\sqrt{5}}{2}$,
$ \frac{-5-\sqrt{5}}{\sqrt{10}}$;\ \ 
$0$,
$ \frac{5+\sqrt{5}}{\sqrt{10}}$,
$ \frac{-5-\sqrt{5}}{\sqrt{10}}$,
$0$;\ \ 
$ -1$,
$ -1$,
$ -\sqrt{2}$;\ \ 
$ -1$,
$ \sqrt{2}$;\ \ 
$0$)

Factors = $2_{\frac{14}{5},3.618}^{5,395}\boxtimes 3_{\frac{13}{2},4.}^{16,330}$

  \vskip 2ex

\noindent30. $6_{\frac{37}{10},14.47}^{80,629}$ \irep{48}:\ \ 
$d_i$ = ($1.0$,
$1.0$,
$1.414$,
$1.618$,
$1.618$,
$2.288$) 

\vskip 0.7ex
\hangindent=3em \hangafter=1
$D^2= 14.472 = 
10+2\sqrt{5}$

\vskip 0.7ex
\hangindent=3em \hangafter=1
$T = ( 0,
\frac{1}{2},
\frac{13}{16},
\frac{3}{5},
\frac{1}{10},
\frac{33}{80} )
$,

\vskip 0.7ex
\hangindent=3em \hangafter=1
$S$ = ($ 1$,
$ 1$,
$ \sqrt{2}$,
$ \frac{1+\sqrt{5}}{2}$,
$ \frac{1+\sqrt{5}}{2}$,
$ \frac{5+\sqrt{5}}{\sqrt{10}}$;\ \ 
$ 1$,
$ -\sqrt{2}$,
$ \frac{1+\sqrt{5}}{2}$,
$ \frac{1+\sqrt{5}}{2}$,
$ \frac{-5-\sqrt{5}}{\sqrt{10}}$;\ \ 
$0$,
$ \frac{5+\sqrt{5}}{\sqrt{10}}$,
$ \frac{-5-\sqrt{5}}{\sqrt{10}}$,
$0$;\ \ 
$ -1$,
$ -1$,
$ -\sqrt{2}$;\ \ 
$ -1$,
$ \sqrt{2}$;\ \ 
$0$)

Factors = $2_{\frac{26}{5},3.618}^{5,720}\boxtimes 3_{\frac{13}{2},4.}^{16,330}$

  \vskip 2ex

\noindent31. $6_{\frac{23}{10},14.47}^{80,108}$ \irep{48}:\ \ 
$d_i$ = ($1.0$,
$1.0$,
$1.414$,
$1.618$,
$1.618$,
$2.288$) 

\vskip 0.7ex
\hangindent=3em \hangafter=1
$D^2= 14.472 = 
10+2\sqrt{5}$

\vskip 0.7ex
\hangindent=3em \hangafter=1
$T = ( 0,
\frac{1}{2},
\frac{15}{16},
\frac{2}{5},
\frac{9}{10},
\frac{27}{80} )
$,

\vskip 0.7ex
\hangindent=3em \hangafter=1
$S$ = ($ 1$,
$ 1$,
$ \sqrt{2}$,
$ \frac{1+\sqrt{5}}{2}$,
$ \frac{1+\sqrt{5}}{2}$,
$ \frac{5+\sqrt{5}}{\sqrt{10}}$;\ \ 
$ 1$,
$ -\sqrt{2}$,
$ \frac{1+\sqrt{5}}{2}$,
$ \frac{1+\sqrt{5}}{2}$,
$ \frac{-5-\sqrt{5}}{\sqrt{10}}$;\ \ 
$0$,
$ \frac{5+\sqrt{5}}{\sqrt{10}}$,
$ \frac{-5-\sqrt{5}}{\sqrt{10}}$,
$0$;\ \ 
$ -1$,
$ -1$,
$ -\sqrt{2}$;\ \ 
$ -1$,
$ \sqrt{2}$;\ \ 
$0$)

Factors = $2_{\frac{14}{5},3.618}^{5,395}\boxtimes 3_{\frac{15}{2},4.}^{16,639}$

  \vskip 2ex

\noindent32. $6_{\frac{47}{10},14.47}^{80,518}$ \irep{48}:\ \ 
$d_i$ = ($1.0$,
$1.0$,
$1.414$,
$1.618$,
$1.618$,
$2.288$) 

\vskip 0.7ex
\hangindent=3em \hangafter=1
$D^2= 14.472 = 
10+2\sqrt{5}$

\vskip 0.7ex
\hangindent=3em \hangafter=1
$T = ( 0,
\frac{1}{2},
\frac{15}{16},
\frac{3}{5},
\frac{1}{10},
\frac{43}{80} )
$,

\vskip 0.7ex
\hangindent=3em \hangafter=1
$S$ = ($ 1$,
$ 1$,
$ \sqrt{2}$,
$ \frac{1+\sqrt{5}}{2}$,
$ \frac{1+\sqrt{5}}{2}$,
$ \frac{5+\sqrt{5}}{\sqrt{10}}$;\ \ 
$ 1$,
$ -\sqrt{2}$,
$ \frac{1+\sqrt{5}}{2}$,
$ \frac{1+\sqrt{5}}{2}$,
$ \frac{-5-\sqrt{5}}{\sqrt{10}}$;\ \ 
$0$,
$ \frac{5+\sqrt{5}}{\sqrt{10}}$,
$ \frac{-5-\sqrt{5}}{\sqrt{10}}$,
$0$;\ \ 
$ -1$,
$ -1$,
$ -\sqrt{2}$;\ \ 
$ -1$,
$ \sqrt{2}$;\ \ 
$0$)

Factors = $2_{\frac{26}{5},3.618}^{5,720}\boxtimes 3_{\frac{15}{2},4.}^{16,639}$

  \vskip 2ex

\noindent33. $6_{\frac{55}{7},18.59}^{28,108}$ \irep{46}:\ \ 
$d_i$ = ($1.0$,
$1.0$,
$1.801$,
$1.801$,
$2.246$,
$2.246$) 

\vskip 0.7ex
\hangindent=3em \hangafter=1
$D^2= 18.591 = 
12+6c^{1}_{7}
+2c^{2}_{7}
$

\vskip 0.7ex
\hangindent=3em \hangafter=1
$T = ( 0,
\frac{1}{4},
\frac{1}{7},
\frac{11}{28},
\frac{5}{7},
\frac{27}{28} )
$,

\vskip 0.7ex
\hangindent=3em \hangafter=1
$S$ = ($ 1$,
$ 1$,
$ -c_{7}^{3}$,
$ -c_{7}^{3}$,
$ \xi_{7}^{3}$,
$ \xi_{7}^{3}$;\ \ 
$ -1$,
$ -c_{7}^{3}$,
$ c_{7}^{3}$,
$ \xi_{7}^{3}$,
$ -\xi_{7}^{3}$;\ \ 
$ -\xi_{7}^{3}$,
$ -\xi_{7}^{3}$,
$ 1$,
$ 1$;\ \ 
$ \xi_{7}^{3}$,
$ 1$,
$ -1$;\ \ 
$ c_{7}^{3}$,
$ c_{7}^{3}$;\ \ 
$ -c_{7}^{3}$)

Factors = $2_{1,2.}^{4,437}\boxtimes 3_{\frac{48}{7},9.295}^{7,790}$

  \vskip 2ex

\noindent34. $6_{\frac{15}{7},18.59}^{28,289}$ \irep{46}:\ \ 
$d_i$ = ($1.0$,
$1.0$,
$1.801$,
$1.801$,
$2.246$,
$2.246$) 

\vskip 0.7ex
\hangindent=3em \hangafter=1
$D^2= 18.591 = 
12+6c^{1}_{7}
+2c^{2}_{7}
$

\vskip 0.7ex
\hangindent=3em \hangafter=1
$T = ( 0,
\frac{1}{4},
\frac{6}{7},
\frac{3}{28},
\frac{2}{7},
\frac{15}{28} )
$,

\vskip 0.7ex
\hangindent=3em \hangafter=1
$S$ = ($ 1$,
$ 1$,
$ -c_{7}^{3}$,
$ -c_{7}^{3}$,
$ \xi_{7}^{3}$,
$ \xi_{7}^{3}$;\ \ 
$ -1$,
$ -c_{7}^{3}$,
$ c_{7}^{3}$,
$ \xi_{7}^{3}$,
$ -\xi_{7}^{3}$;\ \ 
$ -\xi_{7}^{3}$,
$ -\xi_{7}^{3}$,
$ 1$,
$ 1$;\ \ 
$ \xi_{7}^{3}$,
$ 1$,
$ -1$;\ \ 
$ c_{7}^{3}$,
$ c_{7}^{3}$;\ \ 
$ -c_{7}^{3}$)

Factors = $2_{1,2.}^{4,437}\boxtimes 3_{\frac{8}{7},9.295}^{7,245}$

  \vskip 2ex

\noindent35. $6_{\frac{41}{7},18.59}^{28,114}$ \irep{46}:\ \ 
$d_i$ = ($1.0$,
$1.0$,
$1.801$,
$1.801$,
$2.246$,
$2.246$) 

\vskip 0.7ex
\hangindent=3em \hangafter=1
$D^2= 18.591 = 
12+6c^{1}_{7}
+2c^{2}_{7}
$

\vskip 0.7ex
\hangindent=3em \hangafter=1
$T = ( 0,
\frac{3}{4},
\frac{1}{7},
\frac{25}{28},
\frac{5}{7},
\frac{13}{28} )
$,

\vskip 0.7ex
\hangindent=3em \hangafter=1
$S$ = ($ 1$,
$ 1$,
$ -c_{7}^{3}$,
$ -c_{7}^{3}$,
$ \xi_{7}^{3}$,
$ \xi_{7}^{3}$;\ \ 
$ -1$,
$ -c_{7}^{3}$,
$ c_{7}^{3}$,
$ \xi_{7}^{3}$,
$ -\xi_{7}^{3}$;\ \ 
$ -\xi_{7}^{3}$,
$ -\xi_{7}^{3}$,
$ 1$,
$ 1$;\ \ 
$ \xi_{7}^{3}$,
$ 1$,
$ -1$;\ \ 
$ c_{7}^{3}$,
$ c_{7}^{3}$;\ \ 
$ -c_{7}^{3}$)

Factors = $2_{7,2.}^{4,625}\boxtimes 3_{\frac{48}{7},9.295}^{7,790}$

  \vskip 2ex

\noindent36. $6_{\frac{1}{7},18.59}^{28,212}$ \irep{46}:\ \ 
$d_i$ = ($1.0$,
$1.0$,
$1.801$,
$1.801$,
$2.246$,
$2.246$) 

\vskip 0.7ex
\hangindent=3em \hangafter=1
$D^2= 18.591 = 
12+6c^{1}_{7}
+2c^{2}_{7}
$

\vskip 0.7ex
\hangindent=3em \hangafter=1
$T = ( 0,
\frac{3}{4},
\frac{6}{7},
\frac{17}{28},
\frac{2}{7},
\frac{1}{28} )
$,

\vskip 0.7ex
\hangindent=3em \hangafter=1
$S$ = ($ 1$,
$ 1$,
$ -c_{7}^{3}$,
$ -c_{7}^{3}$,
$ \xi_{7}^{3}$,
$ \xi_{7}^{3}$;\ \ 
$ -1$,
$ -c_{7}^{3}$,
$ c_{7}^{3}$,
$ \xi_{7}^{3}$,
$ -\xi_{7}^{3}$;\ \ 
$ -\xi_{7}^{3}$,
$ -\xi_{7}^{3}$,
$ 1$,
$ 1$;\ \ 
$ \xi_{7}^{3}$,
$ 1$,
$ -1$;\ \ 
$ c_{7}^{3}$,
$ c_{7}^{3}$;\ \ 
$ -c_{7}^{3}$)

Factors = $2_{7,2.}^{4,625}\boxtimes 3_{\frac{8}{7},9.295}^{7,245}$

  \vskip 2ex

\noindent37. $6_{0,20.}^{10,699}$ \irep{27}:\ \ 
$d_i$ = ($1.0$,
$1.0$,
$2.0$,
$2.0$,
$2.236$,
$2.236$) 

\vskip 0.7ex
\hangindent=3em \hangafter=1
$D^2= 20.0 = 
20$

\vskip 0.7ex
\hangindent=3em \hangafter=1
$T = ( 0,
0,
\frac{1}{5},
\frac{4}{5},
0,
\frac{1}{2} )
$,

\vskip 0.7ex
\hangindent=3em \hangafter=1
$S$ = ($ 1$,
$ 1$,
$ 2$,
$ 2$,
$ \sqrt{5}$,
$ \sqrt{5}$;\ \ 
$ 1$,
$ 2$,
$ 2$,
$ -\sqrt{5}$,
$ -\sqrt{5}$;\ \ 
$ -1-\sqrt{5}$,
$ -1+\sqrt{5}$,
$0$,
$0$;\ \ 
$ -1-\sqrt{5}$,
$0$,
$0$;\ \ 
$ \sqrt{5}$,
$ -\sqrt{5}$;\ \ 
$ \sqrt{5}$)

  \vskip 2ex

\noindent38. $6_{4,20.}^{10,101}$ \irep{27}:\ \ 
$d_i$ = ($1.0$,
$1.0$,
$2.0$,
$2.0$,
$2.236$,
$2.236$) 

\vskip 0.7ex
\hangindent=3em \hangafter=1
$D^2= 20.0 = 
20$

\vskip 0.7ex
\hangindent=3em \hangafter=1
$T = ( 0,
0,
\frac{2}{5},
\frac{3}{5},
0,
\frac{1}{2} )
$,

\vskip 0.7ex
\hangindent=3em \hangafter=1
$S$ = ($ 1$,
$ 1$,
$ 2$,
$ 2$,
$ \sqrt{5}$,
$ \sqrt{5}$;\ \ 
$ 1$,
$ 2$,
$ 2$,
$ -\sqrt{5}$,
$ -\sqrt{5}$;\ \ 
$ -1+\sqrt{5}$,
$ -1-\sqrt{5}$,
$0$,
$0$;\ \ 
$ -1+\sqrt{5}$,
$0$,
$0$;\ \ 
$ -\sqrt{5}$,
$ \sqrt{5}$;\ \ 
$ -\sqrt{5}$)

  \vskip 2ex

\noindent39. $6_{0,20.}^{20,139}$ \irep{42}:\ \ 
$d_i$ = ($1.0$,
$1.0$,
$2.0$,
$2.0$,
$2.236$,
$2.236$) 

\vskip 0.7ex
\hangindent=3em \hangafter=1
$D^2= 20.0 = 
20$

\vskip 0.7ex
\hangindent=3em \hangafter=1
$T = ( 0,
0,
\frac{1}{5},
\frac{4}{5},
\frac{1}{4},
\frac{3}{4} )
$,

\vskip 0.7ex
\hangindent=3em \hangafter=1
$S$ = ($ 1$,
$ 1$,
$ 2$,
$ 2$,
$ \sqrt{5}$,
$ \sqrt{5}$;\ \ 
$ 1$,
$ 2$,
$ 2$,
$ -\sqrt{5}$,
$ -\sqrt{5}$;\ \ 
$ -1-\sqrt{5}$,
$ -1+\sqrt{5}$,
$0$,
$0$;\ \ 
$ -1-\sqrt{5}$,
$0$,
$0$;\ \ 
$ -\sqrt{5}$,
$ \sqrt{5}$;\ \ 
$ -\sqrt{5}$)

  \vskip 2ex

\noindent40. $6_{4,20.}^{20,180}$ \irep{42}:\ \ 
$d_i$ = ($1.0$,
$1.0$,
$2.0$,
$2.0$,
$2.236$,
$2.236$) 

\vskip 0.7ex
\hangindent=3em \hangafter=1
$D^2= 20.0 = 
20$

\vskip 0.7ex
\hangindent=3em \hangafter=1
$T = ( 0,
0,
\frac{2}{5},
\frac{3}{5},
\frac{1}{4},
\frac{3}{4} )
$,

\vskip 0.7ex
\hangindent=3em \hangafter=1
$S$ = ($ 1$,
$ 1$,
$ 2$,
$ 2$,
$ \sqrt{5}$,
$ \sqrt{5}$;\ \ 
$ 1$,
$ 2$,
$ 2$,
$ -\sqrt{5}$,
$ -\sqrt{5}$;\ \ 
$ -1+\sqrt{5}$,
$ -1-\sqrt{5}$,
$0$,
$0$;\ \ 
$ -1+\sqrt{5}$,
$0$,
$0$;\ \ 
$ \sqrt{5}$,
$ -\sqrt{5}$;\ \ 
$ \sqrt{5}$)

  \vskip 2ex

\noindent41. $6_{\frac{58}{35},33.63}^{35,955}$ \irep{47}:\ \ 
$d_i$ = ($1.0$,
$1.618$,
$1.801$,
$2.246$,
$2.915$,
$3.635$) 

\vskip 0.7ex
\hangindent=3em \hangafter=1
$D^2= 33.632 = 
15+3c^{1}_{35}
+2c^{4}_{35}
+6c^{5}_{35}
+3c^{6}_{35}
+3c^{7}_{35}
+2c^{10}_{35}
+2c^{11}_{35}
$

\vskip 0.7ex
\hangindent=3em \hangafter=1
$T = ( 0,
\frac{2}{5},
\frac{1}{7},
\frac{5}{7},
\frac{19}{35},
\frac{4}{35} )
$,

\vskip 0.7ex
\hangindent=3em \hangafter=1
$S$ = ($ 1$,
$ \frac{1+\sqrt{5}}{2}$,
$ -c_{7}^{3}$,
$ \xi_{7}^{3}$,
$ c^{1}_{35}
+c^{6}_{35}
$,
$ c^{1}_{35}
+c^{4}_{35}
+c^{6}_{35}
+c^{11}_{35}
$;\ \ 
$ -1$,
$ c^{1}_{35}
+c^{6}_{35}
$,
$ c^{1}_{35}
+c^{4}_{35}
+c^{6}_{35}
+c^{11}_{35}
$,
$ c_{7}^{3}$,
$ -\xi_{7}^{3}$;\ \ 
$ -\xi_{7}^{3}$,
$ 1$,
$ -c^{1}_{35}
-c^{4}_{35}
-c^{6}_{35}
-c^{11}_{35}
$,
$ \frac{1+\sqrt{5}}{2}$;\ \ 
$ c_{7}^{3}$,
$ \frac{1+\sqrt{5}}{2}$,
$ -c^{1}_{35}
-c^{6}_{35}
$;\ \ 
$ \xi_{7}^{3}$,
$ -1$;\ \ 
$ -c_{7}^{3}$)

Factors = $2_{\frac{14}{5},3.618}^{5,395}\boxtimes 3_{\frac{48}{7},9.295}^{7,790}$

  \vskip 2ex

\noindent42. $6_{\frac{138}{35},33.63}^{35,363}$ \irep{47}:\ \ 
$d_i$ = ($1.0$,
$1.618$,
$1.801$,
$2.246$,
$2.915$,
$3.635$) 

\vskip 0.7ex
\hangindent=3em \hangafter=1
$D^2= 33.632 = 
15+3c^{1}_{35}
+2c^{4}_{35}
+6c^{5}_{35}
+3c^{6}_{35}
+3c^{7}_{35}
+2c^{10}_{35}
+2c^{11}_{35}
$

\vskip 0.7ex
\hangindent=3em \hangafter=1
$T = ( 0,
\frac{2}{5},
\frac{6}{7},
\frac{2}{7},
\frac{9}{35},
\frac{24}{35} )
$,

\vskip 0.7ex
\hangindent=3em \hangafter=1
$S$ = ($ 1$,
$ \frac{1+\sqrt{5}}{2}$,
$ -c_{7}^{3}$,
$ \xi_{7}^{3}$,
$ c^{1}_{35}
+c^{6}_{35}
$,
$ c^{1}_{35}
+c^{4}_{35}
+c^{6}_{35}
+c^{11}_{35}
$;\ \ 
$ -1$,
$ c^{1}_{35}
+c^{6}_{35}
$,
$ c^{1}_{35}
+c^{4}_{35}
+c^{6}_{35}
+c^{11}_{35}
$,
$ c_{7}^{3}$,
$ -\xi_{7}^{3}$;\ \ 
$ -\xi_{7}^{3}$,
$ 1$,
$ -c^{1}_{35}
-c^{4}_{35}
-c^{6}_{35}
-c^{11}_{35}
$,
$ \frac{1+\sqrt{5}}{2}$;\ \ 
$ c_{7}^{3}$,
$ \frac{1+\sqrt{5}}{2}$,
$ -c^{1}_{35}
-c^{6}_{35}
$;\ \ 
$ \xi_{7}^{3}$,
$ -1$;\ \ 
$ -c_{7}^{3}$)

Factors = $2_{\frac{14}{5},3.618}^{5,395}\boxtimes 3_{\frac{8}{7},9.295}^{7,245}$

  \vskip 2ex

\noindent43. $6_{\frac{142}{35},33.63}^{35,429}$ \irep{47}:\ \ 
$d_i$ = ($1.0$,
$1.618$,
$1.801$,
$2.246$,
$2.915$,
$3.635$) 

\vskip 0.7ex
\hangindent=3em \hangafter=1
$D^2= 33.632 = 
15+3c^{1}_{35}
+2c^{4}_{35}
+6c^{5}_{35}
+3c^{6}_{35}
+3c^{7}_{35}
+2c^{10}_{35}
+2c^{11}_{35}
$

\vskip 0.7ex
\hangindent=3em \hangafter=1
$T = ( 0,
\frac{3}{5},
\frac{1}{7},
\frac{5}{7},
\frac{26}{35},
\frac{11}{35} )
$,

\vskip 0.7ex
\hangindent=3em \hangafter=1
$S$ = ($ 1$,
$ \frac{1+\sqrt{5}}{2}$,
$ -c_{7}^{3}$,
$ \xi_{7}^{3}$,
$ c^{1}_{35}
+c^{6}_{35}
$,
$ c^{1}_{35}
+c^{4}_{35}
+c^{6}_{35}
+c^{11}_{35}
$;\ \ 
$ -1$,
$ c^{1}_{35}
+c^{6}_{35}
$,
$ c^{1}_{35}
+c^{4}_{35}
+c^{6}_{35}
+c^{11}_{35}
$,
$ c_{7}^{3}$,
$ -\xi_{7}^{3}$;\ \ 
$ -\xi_{7}^{3}$,
$ 1$,
$ -c^{1}_{35}
-c^{4}_{35}
-c^{6}_{35}
-c^{11}_{35}
$,
$ \frac{1+\sqrt{5}}{2}$;\ \ 
$ c_{7}^{3}$,
$ \frac{1+\sqrt{5}}{2}$,
$ -c^{1}_{35}
-c^{6}_{35}
$;\ \ 
$ \xi_{7}^{3}$,
$ -1$;\ \ 
$ -c_{7}^{3}$)

Factors = $2_{\frac{26}{5},3.618}^{5,720}\boxtimes 3_{\frac{48}{7},9.295}^{7,790}$

  \vskip 2ex

\noindent44. $6_{\frac{222}{35},33.63}^{35,224}$ \irep{47}:\ \ 
$d_i$ = ($1.0$,
$1.618$,
$1.801$,
$2.246$,
$2.915$,
$3.635$) 

\vskip 0.7ex
\hangindent=3em \hangafter=1
$D^2= 33.632 = 
15+3c^{1}_{35}
+2c^{4}_{35}
+6c^{5}_{35}
+3c^{6}_{35}
+3c^{7}_{35}
+2c^{10}_{35}
+2c^{11}_{35}
$

\vskip 0.7ex
\hangindent=3em \hangafter=1
$T = ( 0,
\frac{3}{5},
\frac{6}{7},
\frac{2}{7},
\frac{16}{35},
\frac{31}{35} )
$,

\vskip 0.7ex
\hangindent=3em \hangafter=1
$S$ = ($ 1$,
$ \frac{1+\sqrt{5}}{2}$,
$ -c_{7}^{3}$,
$ \xi_{7}^{3}$,
$ c^{1}_{35}
+c^{6}_{35}
$,
$ c^{1}_{35}
+c^{4}_{35}
+c^{6}_{35}
+c^{11}_{35}
$;\ \ 
$ -1$,
$ c^{1}_{35}
+c^{6}_{35}
$,
$ c^{1}_{35}
+c^{4}_{35}
+c^{6}_{35}
+c^{11}_{35}
$,
$ c_{7}^{3}$,
$ -\xi_{7}^{3}$;\ \ 
$ -\xi_{7}^{3}$,
$ 1$,
$ -c^{1}_{35}
-c^{4}_{35}
-c^{6}_{35}
-c^{11}_{35}
$,
$ \frac{1+\sqrt{5}}{2}$;\ \ 
$ c_{7}^{3}$,
$ \frac{1+\sqrt{5}}{2}$,
$ -c^{1}_{35}
-c^{6}_{35}
$;\ \ 
$ \xi_{7}^{3}$,
$ -1$;\ \ 
$ -c_{7}^{3}$)

Factors = $2_{\frac{26}{5},3.618}^{5,720}\boxtimes 3_{\frac{8}{7},9.295}^{7,245}$

  \vskip 2ex

\noindent45. $6_{\frac{46}{13},56.74}^{13,131}$ \irep{35}:\ \ 
$d_i$ = ($1.0$,
$1.941$,
$2.770$,
$3.438$,
$3.907$,
$4.148$) 

\vskip 0.7ex
\hangindent=3em \hangafter=1
$D^2= 56.746 = 
21+15c^{1}_{13}
+10c^{2}_{13}
+6c^{3}_{13}
+3c^{4}_{13}
+c^{5}_{13}
$

\vskip 0.7ex
\hangindent=3em \hangafter=1
$T = ( 0,
\frac{4}{13},
\frac{2}{13},
\frac{7}{13},
\frac{6}{13},
\frac{12}{13} )
$,

\vskip 0.7ex
\hangindent=3em \hangafter=1
$S$ = ($ 1$,
$ -c_{13}^{6}$,
$ \xi_{13}^{3}$,
$ \xi_{13}^{9}$,
$ \xi_{13}^{5}$,
$ \xi_{13}^{7}$;\ \ 
$ -\xi_{13}^{9}$,
$ \xi_{13}^{7}$,
$ -\xi_{13}^{5}$,
$ \xi_{13}^{3}$,
$ -1$;\ \ 
$ \xi_{13}^{9}$,
$ 1$,
$ c_{13}^{6}$,
$ -\xi_{13}^{5}$;\ \ 
$ \xi_{13}^{3}$,
$ -\xi_{13}^{7}$,
$ -c_{13}^{6}$;\ \ 
$ -1$,
$ \xi_{13}^{9}$;\ \ 
$ -\xi_{13}^{3}$)

  \vskip 2ex

\noindent46. $6_{\frac{58}{13},56.74}^{13,502}$ \irep{35}:\ \ 
$d_i$ = ($1.0$,
$1.941$,
$2.770$,
$3.438$,
$3.907$,
$4.148$) 

\vskip 0.7ex
\hangindent=3em \hangafter=1
$D^2= 56.746 = 
21+15c^{1}_{13}
+10c^{2}_{13}
+6c^{3}_{13}
+3c^{4}_{13}
+c^{5}_{13}
$

\vskip 0.7ex
\hangindent=3em \hangafter=1
$T = ( 0,
\frac{9}{13},
\frac{11}{13},
\frac{6}{13},
\frac{7}{13},
\frac{1}{13} )
$,

\vskip 0.7ex
\hangindent=3em \hangafter=1
$S$ = ($ 1$,
$ -c_{13}^{6}$,
$ \xi_{13}^{3}$,
$ \xi_{13}^{9}$,
$ \xi_{13}^{5}$,
$ \xi_{13}^{7}$;\ \ 
$ -\xi_{13}^{9}$,
$ \xi_{13}^{7}$,
$ -\xi_{13}^{5}$,
$ \xi_{13}^{3}$,
$ -1$;\ \ 
$ \xi_{13}^{9}$,
$ 1$,
$ c_{13}^{6}$,
$ -\xi_{13}^{5}$;\ \ 
$ \xi_{13}^{3}$,
$ -\xi_{13}^{7}$,
$ -c_{13}^{6}$;\ \ 
$ -1$,
$ \xi_{13}^{9}$;\ \ 
$ -\xi_{13}^{3}$)

  \vskip 2ex

\noindent47. $6_{\frac{8}{3},74.61}^{9,186}$ \irep{26}:\ \ 
$d_i$ = ($1.0$,
$2.879$,
$2.879$,
$2.879$,
$4.411$,
$5.411$) 

\vskip 0.7ex
\hangindent=3em \hangafter=1
$D^2= 74.617 = 
27+27c^{1}_{9}
+18c^{2}_{9}
$

\vskip 0.7ex
\hangindent=3em \hangafter=1
$T = ( 0,
\frac{1}{9},
\frac{1}{9},
\frac{1}{9},
\frac{1}{3},
\frac{2}{3} )
$,

\vskip 0.7ex
\hangindent=3em \hangafter=1
$S$ = ($ 1$,
$ \xi_{9}^{5}$,
$ \xi_{9}^{5}$,
$ \xi_{9}^{5}$,
$ 1+2c^{1}_{9}
+c^{2}_{9}
$,
$ 2+2c^{1}_{9}
+c^{2}_{9}
$;\ \ 
$ 2\xi_{9}^{5}$,
$ -\xi_{9}^{5}$,
$ -\xi_{9}^{5}$,
$ \xi_{9}^{5}$,
$ -\xi_{9}^{5}$;\ \ 
$ 2\xi_{9}^{5}$,
$ -\xi_{9}^{5}$,
$ \xi_{9}^{5}$,
$ -\xi_{9}^{5}$;\ \ 
$ 2\xi_{9}^{5}$,
$ \xi_{9}^{5}$,
$ -\xi_{9}^{5}$;\ \ 
$ -2-2  c^{1}_{9}
-c^{2}_{9}
$,
$ -1$;\ \ 
$ 1+2c^{1}_{9}
+c^{2}_{9}
$)

  \vskip 2ex

\noindent48. $6_{\frac{16}{3},74.61}^{9,452}$ \irep{26}:\ \ 
$d_i$ = ($1.0$,
$2.879$,
$2.879$,
$2.879$,
$4.411$,
$5.411$) 

\vskip 0.7ex
\hangindent=3em \hangafter=1
$D^2= 74.617 = 
27+27c^{1}_{9}
+18c^{2}_{9}
$

\vskip 0.7ex
\hangindent=3em \hangafter=1
$T = ( 0,
\frac{8}{9},
\frac{8}{9},
\frac{8}{9},
\frac{2}{3},
\frac{1}{3} )
$,

\vskip 0.7ex
\hangindent=3em \hangafter=1
$S$ = ($ 1$,
$ \xi_{9}^{5}$,
$ \xi_{9}^{5}$,
$ \xi_{9}^{5}$,
$ 1+2c^{1}_{9}
+c^{2}_{9}
$,
$ 2+2c^{1}_{9}
+c^{2}_{9}
$;\ \ 
$ 2\xi_{9}^{5}$,
$ -\xi_{9}^{5}$,
$ -\xi_{9}^{5}$,
$ \xi_{9}^{5}$,
$ -\xi_{9}^{5}$;\ \ 
$ 2\xi_{9}^{5}$,
$ -\xi_{9}^{5}$,
$ \xi_{9}^{5}$,
$ -\xi_{9}^{5}$;\ \ 
$ 2\xi_{9}^{5}$,
$ \xi_{9}^{5}$,
$ -\xi_{9}^{5}$;\ \ 
$ -2-2  c^{1}_{9}
-c^{2}_{9}
$,
$ -1$;\ \ 
$ 1+2c^{1}_{9}
+c^{2}_{9}
$)

  \vskip 2ex

\noindent49. $6_{6,100.6}^{21,154}$ \irep{43}:\ \ 
$d_i$ = ($1.0$,
$3.791$,
$3.791$,
$3.791$,
$4.791$,
$5.791$) 

\vskip 0.7ex
\hangindent=3em \hangafter=1
$D^2= 100.617 = 
\frac{105+21\sqrt{21}}{2}$

\vskip 0.7ex
\hangindent=3em \hangafter=1
$T = ( 0,
\frac{1}{7},
\frac{2}{7},
\frac{4}{7},
0,
\frac{2}{3} )
$,

\vskip 0.7ex
\hangindent=3em \hangafter=1
$S$ = ($ 1$,
$ \frac{3+\sqrt{21}}{2}$,
$ \frac{3+\sqrt{21}}{2}$,
$ \frac{3+\sqrt{21}}{2}$,
$ \frac{5+\sqrt{21}}{2}$,
$ \frac{7+\sqrt{21}}{2}$;\ \ 
$ 2-c^{1}_{21}
-2  c^{2}_{21}
+3c^{3}_{21}
+2c^{4}_{21}
-2  c^{5}_{21}
$,
$ -c^{2}_{21}
-2  c^{3}_{21}
-c^{4}_{21}
+c^{5}_{21}
$,
$ -1+2c^{1}_{21}
+3c^{2}_{21}
-c^{3}_{21}
+2c^{5}_{21}
$,
$ -\frac{3+\sqrt{21}}{2}$,
$0$;\ \ 
$ -1+2c^{1}_{21}
+3c^{2}_{21}
-c^{3}_{21}
+2c^{5}_{21}
$,
$ 2-c^{1}_{21}
-2  c^{2}_{21}
+3c^{3}_{21}
+2c^{4}_{21}
-2  c^{5}_{21}
$,
$ -\frac{3+\sqrt{21}}{2}$,
$0$;\ \ 
$ -c^{2}_{21}
-2  c^{3}_{21}
-c^{4}_{21}
+c^{5}_{21}
$,
$ -\frac{3+\sqrt{21}}{2}$,
$0$;\ \ 
$ 1$,
$ \frac{7+\sqrt{21}}{2}$;\ \ 
$ -\frac{7+\sqrt{21}}{2}$)

  \vskip 2ex

\noindent50. $6_{2,100.6}^{21,320}$ \irep{43}:\ \ 
$d_i$ = ($1.0$,
$3.791$,
$3.791$,
$3.791$,
$4.791$,
$5.791$) 

\vskip 0.7ex
\hangindent=3em \hangafter=1
$D^2= 100.617 = 
\frac{105+21\sqrt{21}}{2}$

\vskip 0.7ex
\hangindent=3em \hangafter=1
$T = ( 0,
\frac{3}{7},
\frac{5}{7},
\frac{6}{7},
0,
\frac{1}{3} )
$,

\vskip 0.7ex
\hangindent=3em \hangafter=1
$S$ = ($ 1$,
$ \frac{3+\sqrt{21}}{2}$,
$ \frac{3+\sqrt{21}}{2}$,
$ \frac{3+\sqrt{21}}{2}$,
$ \frac{5+\sqrt{21}}{2}$,
$ \frac{7+\sqrt{21}}{2}$;\ \ 
$ -c^{2}_{21}
-2  c^{3}_{21}
-c^{4}_{21}
+c^{5}_{21}
$,
$ 2-c^{1}_{21}
-2  c^{2}_{21}
+3c^{3}_{21}
+2c^{4}_{21}
-2  c^{5}_{21}
$,
$ -1+2c^{1}_{21}
+3c^{2}_{21}
-c^{3}_{21}
+2c^{5}_{21}
$,
$ -\frac{3+\sqrt{21}}{2}$,
$0$;\ \ 
$ -1+2c^{1}_{21}
+3c^{2}_{21}
-c^{3}_{21}
+2c^{5}_{21}
$,
$ -c^{2}_{21}
-2  c^{3}_{21}
-c^{4}_{21}
+c^{5}_{21}
$,
$ -\frac{3+\sqrt{21}}{2}$,
$0$;\ \ 
$ 2-c^{1}_{21}
-2  c^{2}_{21}
+3c^{3}_{21}
+2c^{4}_{21}
-2  c^{5}_{21}
$,
$ -\frac{3+\sqrt{21}}{2}$,
$0$;\ \ 
$ 1$,
$ \frac{7+\sqrt{21}}{2}$;\ \ 
$ -\frac{7+\sqrt{21}}{2}$)

  \vskip 2ex 

}

\subsection{Rank 7}
\label{uni7}

{\small

\noindent1. $7_{2,7.}^{7,892}$ \irep{48}:\ \ 
$d_i$ = ($1.0$,
$1.0$,
$1.0$,
$1.0$,
$1.0$,
$1.0$,
$1.0$) 

\vskip 0.7ex
\hangindent=3em \hangafter=1
$D^2= 7.0 = 
7$

\vskip 0.7ex
\hangindent=3em \hangafter=1
$T = ( 0,
\frac{1}{7},
\frac{1}{7},
\frac{2}{7},
\frac{2}{7},
\frac{4}{7},
\frac{4}{7} )
$,

\vskip 0.7ex
\hangindent=3em \hangafter=1
$S$ = ($ 1$,
$ 1$,
$ 1$,
$ 1$,
$ 1$,
$ 1$,
$ 1$;\ \ 
$ -\zeta_{14}^{3}$,
$ \zeta_{7}^{2}$,
$ -\zeta_{14}^{5}$,
$ \zeta_{7}^{1}$,
$ -\zeta_{14}^{1}$,
$ \zeta_{7}^{3}$;\ \ 
$ -\zeta_{14}^{3}$,
$ \zeta_{7}^{1}$,
$ -\zeta_{14}^{5}$,
$ \zeta_{7}^{3}$,
$ -\zeta_{14}^{1}$;\ \ 
$ \zeta_{7}^{3}$,
$ -\zeta_{14}^{1}$,
$ \zeta_{7}^{2}$,
$ -\zeta_{14}^{3}$;\ \ 
$ \zeta_{7}^{3}$,
$ -\zeta_{14}^{3}$,
$ \zeta_{7}^{2}$;\ \ 
$ -\zeta_{14}^{5}$,
$ \zeta_{7}^{1}$;\ \ 
$ -\zeta_{14}^{5}$)

  \vskip 2ex

\noindent2. $7_{6,7.}^{7,110}$ \irep{48}:\ \ 
$d_i$ = ($1.0$,
$1.0$,
$1.0$,
$1.0$,
$1.0$,
$1.0$,
$1.0$) 

\vskip 0.7ex
\hangindent=3em \hangafter=1
$D^2= 7.0 = 
7$

\vskip 0.7ex
\hangindent=3em \hangafter=1
$T = ( 0,
\frac{3}{7},
\frac{3}{7},
\frac{5}{7},
\frac{5}{7},
\frac{6}{7},
\frac{6}{7} )
$,

\vskip 0.7ex
\hangindent=3em \hangafter=1
$S$ = ($ 1$,
$ 1$,
$ 1$,
$ 1$,
$ 1$,
$ 1$,
$ 1$;\ \ 
$ \zeta_{7}^{1}$,
$ -\zeta_{14}^{5}$,
$ -\zeta_{14}^{3}$,
$ \zeta_{7}^{2}$,
$ -\zeta_{14}^{1}$,
$ \zeta_{7}^{3}$;\ \ 
$ \zeta_{7}^{1}$,
$ \zeta_{7}^{2}$,
$ -\zeta_{14}^{3}$,
$ \zeta_{7}^{3}$,
$ -\zeta_{14}^{1}$;\ \ 
$ -\zeta_{14}^{1}$,
$ \zeta_{7}^{3}$,
$ -\zeta_{14}^{5}$,
$ \zeta_{7}^{1}$;\ \ 
$ -\zeta_{14}^{1}$,
$ \zeta_{7}^{1}$,
$ -\zeta_{14}^{5}$;\ \ 
$ \zeta_{7}^{2}$,
$ -\zeta_{14}^{3}$;\ \ 
$ \zeta_{7}^{2}$)

  \vskip 2ex

\noindent3. $7_{\frac{27}{4},27.31}^{32,396}$ \irep{96}:\ \ 
$d_i$ = ($1.0$,
$1.0$,
$1.847$,
$1.847$,
$2.414$,
$2.414$,
$2.613$) 

\vskip 0.7ex
\hangindent=3em \hangafter=1
$D^2= 27.313 = 
16+8\sqrt{2}$

\vskip 0.7ex
\hangindent=3em \hangafter=1
$T = ( 0,
\frac{1}{2},
\frac{1}{32},
\frac{1}{32},
\frac{1}{4},
\frac{3}{4},
\frac{21}{32} )
$,

\vskip 0.7ex
\hangindent=3em \hangafter=1
$S$ = ($ 1$,
$ 1$,
$ c_{16}^{1}$,
$ c_{16}^{1}$,
$ 1+\sqrt{2}$,
$ 1+\sqrt{2}$,
$ c^{1}_{16}
+c^{3}_{16}
$;\ \ 
$ 1$,
$ -c_{16}^{1}$,
$ -c_{16}^{1}$,
$ 1+\sqrt{2}$,
$ 1+\sqrt{2}$,
$ -c^{1}_{16}
-c^{3}_{16}
$;\ \ 
$(-c^{1}_{16}
-c^{3}_{16}
)\mathrm{i}$,
$(c^{1}_{16}
+c^{3}_{16}
)\mathrm{i}$,
$ -c_{16}^{1}$,
$ c_{16}^{1}$,
$0$;\ \ 
$(-c^{1}_{16}
-c^{3}_{16}
)\mathrm{i}$,
$ -c_{16}^{1}$,
$ c_{16}^{1}$,
$0$;\ \ 
$ -1$,
$ -1$,
$ c^{1}_{16}
+c^{3}_{16}
$;\ \ 
$ -1$,
$ -c^{1}_{16}
-c^{3}_{16}
$;\ \ 
$0$)

  \vskip 2ex

\noindent4. $7_{\frac{9}{4},27.31}^{32,918}$ \irep{97}:\ \ 
$d_i$ = ($1.0$,
$1.0$,
$1.847$,
$1.847$,
$2.414$,
$2.414$,
$2.613$) 

\vskip 0.7ex
\hangindent=3em \hangafter=1
$D^2= 27.313 = 
16+8\sqrt{2}$

\vskip 0.7ex
\hangindent=3em \hangafter=1
$T = ( 0,
\frac{1}{2},
\frac{3}{32},
\frac{3}{32},
\frac{1}{4},
\frac{3}{4},
\frac{15}{32} )
$,

\vskip 0.7ex
\hangindent=3em \hangafter=1
$S$ = ($ 1$,
$ 1$,
$ c_{16}^{1}$,
$ c_{16}^{1}$,
$ 1+\sqrt{2}$,
$ 1+\sqrt{2}$,
$ c^{1}_{16}
+c^{3}_{16}
$;\ \ 
$ 1$,
$ -c_{16}^{1}$,
$ -c_{16}^{1}$,
$ 1+\sqrt{2}$,
$ 1+\sqrt{2}$,
$ -c^{1}_{16}
-c^{3}_{16}
$;\ \ 
$ c^{1}_{16}
+c^{3}_{16}
$,
$ -c^{1}_{16}
-c^{3}_{16}
$,
$ c_{16}^{1}$,
$ -c_{16}^{1}$,
$0$;\ \ 
$ c^{1}_{16}
+c^{3}_{16}
$,
$ c_{16}^{1}$,
$ -c_{16}^{1}$,
$0$;\ \ 
$ -1$,
$ -1$,
$ -c^{1}_{16}
-c^{3}_{16}
$;\ \ 
$ -1$,
$ c^{1}_{16}
+c^{3}_{16}
$;\ \ 
$0$)

  \vskip 2ex

\noindent5. $7_{\frac{31}{4},27.31}^{32,159}$ \irep{97}:\ \ 
$d_i$ = ($1.0$,
$1.0$,
$1.847$,
$1.847$,
$2.414$,
$2.414$,
$2.613$) 

\vskip 0.7ex
\hangindent=3em \hangafter=1
$D^2= 27.313 = 
16+8\sqrt{2}$

\vskip 0.7ex
\hangindent=3em \hangafter=1
$T = ( 0,
\frac{1}{2},
\frac{5}{32},
\frac{5}{32},
\frac{1}{4},
\frac{3}{4},
\frac{25}{32} )
$,

\vskip 0.7ex
\hangindent=3em \hangafter=1
$S$ = ($ 1$,
$ 1$,
$ c_{16}^{1}$,
$ c_{16}^{1}$,
$ 1+\sqrt{2}$,
$ 1+\sqrt{2}$,
$ c^{1}_{16}
+c^{3}_{16}
$;\ \ 
$ 1$,
$ -c_{16}^{1}$,
$ -c_{16}^{1}$,
$ 1+\sqrt{2}$,
$ 1+\sqrt{2}$,
$ -c^{1}_{16}
-c^{3}_{16}
$;\ \ 
$ -c^{1}_{16}
-c^{3}_{16}
$,
$ c^{1}_{16}
+c^{3}_{16}
$,
$ -c_{16}^{1}$,
$ c_{16}^{1}$,
$0$;\ \ 
$ -c^{1}_{16}
-c^{3}_{16}
$,
$ -c_{16}^{1}$,
$ c_{16}^{1}$,
$0$;\ \ 
$ -1$,
$ -1$,
$ c^{1}_{16}
+c^{3}_{16}
$;\ \ 
$ -1$,
$ -c^{1}_{16}
-c^{3}_{16}
$;\ \ 
$0$)

  \vskip 2ex

\noindent6. $7_{\frac{13}{4},27.31}^{32,427}$ \irep{96}:\ \ 
$d_i$ = ($1.0$,
$1.0$,
$1.847$,
$1.847$,
$2.414$,
$2.414$,
$2.613$) 

\vskip 0.7ex
\hangindent=3em \hangafter=1
$D^2= 27.313 = 
16+8\sqrt{2}$

\vskip 0.7ex
\hangindent=3em \hangafter=1
$T = ( 0,
\frac{1}{2},
\frac{7}{32},
\frac{7}{32},
\frac{1}{4},
\frac{3}{4},
\frac{19}{32} )
$,

\vskip 0.7ex
\hangindent=3em \hangafter=1
$S$ = ($ 1$,
$ 1$,
$ c_{16}^{1}$,
$ c_{16}^{1}$,
$ 1+\sqrt{2}$,
$ 1+\sqrt{2}$,
$ c^{1}_{16}
+c^{3}_{16}
$;\ \ 
$ 1$,
$ -c_{16}^{1}$,
$ -c_{16}^{1}$,
$ 1+\sqrt{2}$,
$ 1+\sqrt{2}$,
$ -c^{1}_{16}
-c^{3}_{16}
$;\ \ 
$(-c^{1}_{16}
-c^{3}_{16}
)\mathrm{i}$,
$(c^{1}_{16}
+c^{3}_{16}
)\mathrm{i}$,
$ c_{16}^{1}$,
$ -c_{16}^{1}$,
$0$;\ \ 
$(-c^{1}_{16}
-c^{3}_{16}
)\mathrm{i}$,
$ c_{16}^{1}$,
$ -c_{16}^{1}$,
$0$;\ \ 
$ -1$,
$ -1$,
$ -c^{1}_{16}
-c^{3}_{16}
$;\ \ 
$ -1$,
$ c^{1}_{16}
+c^{3}_{16}
$;\ \ 
$0$)

  \vskip 2ex

\noindent7. $7_{\frac{3}{4},27.31}^{32,913}$ \irep{96}:\ \ 
$d_i$ = ($1.0$,
$1.0$,
$1.847$,
$1.847$,
$2.414$,
$2.414$,
$2.613$) 

\vskip 0.7ex
\hangindent=3em \hangafter=1
$D^2= 27.313 = 
16+8\sqrt{2}$

\vskip 0.7ex
\hangindent=3em \hangafter=1
$T = ( 0,
\frac{1}{2},
\frac{9}{32},
\frac{9}{32},
\frac{1}{4},
\frac{3}{4},
\frac{29}{32} )
$,

\vskip 0.7ex
\hangindent=3em \hangafter=1
$S$ = ($ 1$,
$ 1$,
$ c_{16}^{1}$,
$ c_{16}^{1}$,
$ 1+\sqrt{2}$,
$ 1+\sqrt{2}$,
$ c^{1}_{16}
+c^{3}_{16}
$;\ \ 
$ 1$,
$ -c_{16}^{1}$,
$ -c_{16}^{1}$,
$ 1+\sqrt{2}$,
$ 1+\sqrt{2}$,
$ -c^{1}_{16}
-c^{3}_{16}
$;\ \ 
$(c^{1}_{16}
+c^{3}_{16}
)\mathrm{i}$,
$(-c^{1}_{16}
-c^{3}_{16}
)\mathrm{i}$,
$ -c_{16}^{1}$,
$ c_{16}^{1}$,
$0$;\ \ 
$(c^{1}_{16}
+c^{3}_{16}
)\mathrm{i}$,
$ -c_{16}^{1}$,
$ c_{16}^{1}$,
$0$;\ \ 
$ -1$,
$ -1$,
$ c^{1}_{16}
+c^{3}_{16}
$;\ \ 
$ -1$,
$ -c^{1}_{16}
-c^{3}_{16}
$;\ \ 
$0$)

  \vskip 2ex

\noindent8. $7_{\frac{17}{4},27.31}^{32,261}$ \irep{97}:\ \ 
$d_i$ = ($1.0$,
$1.0$,
$1.847$,
$1.847$,
$2.414$,
$2.414$,
$2.613$) 

\vskip 0.7ex
\hangindent=3em \hangafter=1
$D^2= 27.313 = 
16+8\sqrt{2}$

\vskip 0.7ex
\hangindent=3em \hangafter=1
$T = ( 0,
\frac{1}{2},
\frac{11}{32},
\frac{11}{32},
\frac{1}{4},
\frac{3}{4},
\frac{23}{32} )
$,

\vskip 0.7ex
\hangindent=3em \hangafter=1
$S$ = ($ 1$,
$ 1$,
$ c_{16}^{1}$,
$ c_{16}^{1}$,
$ 1+\sqrt{2}$,
$ 1+\sqrt{2}$,
$ c^{1}_{16}
+c^{3}_{16}
$;\ \ 
$ 1$,
$ -c_{16}^{1}$,
$ -c_{16}^{1}$,
$ 1+\sqrt{2}$,
$ 1+\sqrt{2}$,
$ -c^{1}_{16}
-c^{3}_{16}
$;\ \ 
$ -c^{1}_{16}
-c^{3}_{16}
$,
$ c^{1}_{16}
+c^{3}_{16}
$,
$ c_{16}^{1}$,
$ -c_{16}^{1}$,
$0$;\ \ 
$ -c^{1}_{16}
-c^{3}_{16}
$,
$ c_{16}^{1}$,
$ -c_{16}^{1}$,
$0$;\ \ 
$ -1$,
$ -1$,
$ -c^{1}_{16}
-c^{3}_{16}
$;\ \ 
$ -1$,
$ c^{1}_{16}
+c^{3}_{16}
$;\ \ 
$0$)

  \vskip 2ex

\noindent9. $7_{\frac{7}{4},27.31}^{32,912}$ \irep{97}:\ \ 
$d_i$ = ($1.0$,
$1.0$,
$1.847$,
$1.847$,
$2.414$,
$2.414$,
$2.613$) 

\vskip 0.7ex
\hangindent=3em \hangafter=1
$D^2= 27.313 = 
16+8\sqrt{2}$

\vskip 0.7ex
\hangindent=3em \hangafter=1
$T = ( 0,
\frac{1}{2},
\frac{13}{32},
\frac{13}{32},
\frac{1}{4},
\frac{3}{4},
\frac{1}{32} )
$,

\vskip 0.7ex
\hangindent=3em \hangafter=1
$S$ = ($ 1$,
$ 1$,
$ c_{16}^{1}$,
$ c_{16}^{1}$,
$ 1+\sqrt{2}$,
$ 1+\sqrt{2}$,
$ c^{1}_{16}
+c^{3}_{16}
$;\ \ 
$ 1$,
$ -c_{16}^{1}$,
$ -c_{16}^{1}$,
$ 1+\sqrt{2}$,
$ 1+\sqrt{2}$,
$ -c^{1}_{16}
-c^{3}_{16}
$;\ \ 
$ c^{1}_{16}
+c^{3}_{16}
$,
$ -c^{1}_{16}
-c^{3}_{16}
$,
$ -c_{16}^{1}$,
$ c_{16}^{1}$,
$0$;\ \ 
$ c^{1}_{16}
+c^{3}_{16}
$,
$ -c_{16}^{1}$,
$ c_{16}^{1}$,
$0$;\ \ 
$ -1$,
$ -1$,
$ c^{1}_{16}
+c^{3}_{16}
$;\ \ 
$ -1$,
$ -c^{1}_{16}
-c^{3}_{16}
$;\ \ 
$0$)

  \vskip 2ex

\noindent10. $7_{\frac{21}{4},27.31}^{32,114}$ \irep{96}:\ \ 
$d_i$ = ($1.0$,
$1.0$,
$1.847$,
$1.847$,
$2.414$,
$2.414$,
$2.613$) 

\vskip 0.7ex
\hangindent=3em \hangafter=1
$D^2= 27.313 = 
16+8\sqrt{2}$

\vskip 0.7ex
\hangindent=3em \hangafter=1
$T = ( 0,
\frac{1}{2},
\frac{15}{32},
\frac{15}{32},
\frac{1}{4},
\frac{3}{4},
\frac{27}{32} )
$,

\vskip 0.7ex
\hangindent=3em \hangafter=1
$S$ = ($ 1$,
$ 1$,
$ c_{16}^{1}$,
$ c_{16}^{1}$,
$ 1+\sqrt{2}$,
$ 1+\sqrt{2}$,
$ c^{1}_{16}
+c^{3}_{16}
$;\ \ 
$ 1$,
$ -c_{16}^{1}$,
$ -c_{16}^{1}$,
$ 1+\sqrt{2}$,
$ 1+\sqrt{2}$,
$ -c^{1}_{16}
-c^{3}_{16}
$;\ \ 
$(c^{1}_{16}
+c^{3}_{16}
)\mathrm{i}$,
$(-c^{1}_{16}
-c^{3}_{16}
)\mathrm{i}$,
$ c_{16}^{1}$,
$ -c_{16}^{1}$,
$0$;\ \ 
$(c^{1}_{16}
+c^{3}_{16}
)\mathrm{i}$,
$ c_{16}^{1}$,
$ -c_{16}^{1}$,
$0$;\ \ 
$ -1$,
$ -1$,
$ -c^{1}_{16}
-c^{3}_{16}
$;\ \ 
$ -1$,
$ c^{1}_{16}
+c^{3}_{16}
$;\ \ 
$0$)

  \vskip 2ex

\noindent11. $7_{\frac{11}{4},27.31}^{32,418}$ \irep{96}:\ \ 
$d_i$ = ($1.0$,
$1.0$,
$1.847$,
$1.847$,
$2.414$,
$2.414$,
$2.613$) 

\vskip 0.7ex
\hangindent=3em \hangafter=1
$D^2= 27.313 = 
16+8\sqrt{2}$

\vskip 0.7ex
\hangindent=3em \hangafter=1
$T = ( 0,
\frac{1}{2},
\frac{17}{32},
\frac{17}{32},
\frac{1}{4},
\frac{3}{4},
\frac{5}{32} )
$,

\vskip 0.7ex
\hangindent=3em \hangafter=1
$S$ = ($ 1$,
$ 1$,
$ c_{16}^{1}$,
$ c_{16}^{1}$,
$ 1+\sqrt{2}$,
$ 1+\sqrt{2}$,
$ c^{1}_{16}
+c^{3}_{16}
$;\ \ 
$ 1$,
$ -c_{16}^{1}$,
$ -c_{16}^{1}$,
$ 1+\sqrt{2}$,
$ 1+\sqrt{2}$,
$ -c^{1}_{16}
-c^{3}_{16}
$;\ \ 
$(-c^{1}_{16}
-c^{3}_{16}
)\mathrm{i}$,
$(c^{1}_{16}
+c^{3}_{16}
)\mathrm{i}$,
$ -c_{16}^{1}$,
$ c_{16}^{1}$,
$0$;\ \ 
$(-c^{1}_{16}
-c^{3}_{16}
)\mathrm{i}$,
$ -c_{16}^{1}$,
$ c_{16}^{1}$,
$0$;\ \ 
$ -1$,
$ -1$,
$ c^{1}_{16}
+c^{3}_{16}
$;\ \ 
$ -1$,
$ -c^{1}_{16}
-c^{3}_{16}
$;\ \ 
$0$)

  \vskip 2ex

\noindent12. $7_{\frac{25}{4},27.31}^{32,222}$ \irep{97}:\ \ 
$d_i$ = ($1.0$,
$1.0$,
$1.847$,
$1.847$,
$2.414$,
$2.414$,
$2.613$) 

\vskip 0.7ex
\hangindent=3em \hangafter=1
$D^2= 27.313 = 
16+8\sqrt{2}$

\vskip 0.7ex
\hangindent=3em \hangafter=1
$T = ( 0,
\frac{1}{2},
\frac{19}{32},
\frac{19}{32},
\frac{1}{4},
\frac{3}{4},
\frac{31}{32} )
$,

\vskip 0.7ex
\hangindent=3em \hangafter=1
$S$ = ($ 1$,
$ 1$,
$ c_{16}^{1}$,
$ c_{16}^{1}$,
$ 1+\sqrt{2}$,
$ 1+\sqrt{2}$,
$ c^{1}_{16}
+c^{3}_{16}
$;\ \ 
$ 1$,
$ -c_{16}^{1}$,
$ -c_{16}^{1}$,
$ 1+\sqrt{2}$,
$ 1+\sqrt{2}$,
$ -c^{1}_{16}
-c^{3}_{16}
$;\ \ 
$ c^{1}_{16}
+c^{3}_{16}
$,
$ -c^{1}_{16}
-c^{3}_{16}
$,
$ c_{16}^{1}$,
$ -c_{16}^{1}$,
$0$;\ \ 
$ c^{1}_{16}
+c^{3}_{16}
$,
$ c_{16}^{1}$,
$ -c_{16}^{1}$,
$0$;\ \ 
$ -1$,
$ -1$,
$ -c^{1}_{16}
-c^{3}_{16}
$;\ \ 
$ -1$,
$ c^{1}_{16}
+c^{3}_{16}
$;\ \ 
$0$)

  \vskip 2ex

\noindent13. $7_{\frac{15}{4},27.31}^{32,272}$ \irep{97}:\ \ 
$d_i$ = ($1.0$,
$1.0$,
$1.847$,
$1.847$,
$2.414$,
$2.414$,
$2.613$) 

\vskip 0.7ex
\hangindent=3em \hangafter=1
$D^2= 27.313 = 
16+8\sqrt{2}$

\vskip 0.7ex
\hangindent=3em \hangafter=1
$T = ( 0,
\frac{1}{2},
\frac{21}{32},
\frac{21}{32},
\frac{1}{4},
\frac{3}{4},
\frac{9}{32} )
$,

\vskip 0.7ex
\hangindent=3em \hangafter=1
$S$ = ($ 1$,
$ 1$,
$ c_{16}^{1}$,
$ c_{16}^{1}$,
$ 1+\sqrt{2}$,
$ 1+\sqrt{2}$,
$ c^{1}_{16}
+c^{3}_{16}
$;\ \ 
$ 1$,
$ -c_{16}^{1}$,
$ -c_{16}^{1}$,
$ 1+\sqrt{2}$,
$ 1+\sqrt{2}$,
$ -c^{1}_{16}
-c^{3}_{16}
$;\ \ 
$ -c^{1}_{16}
-c^{3}_{16}
$,
$ c^{1}_{16}
+c^{3}_{16}
$,
$ -c_{16}^{1}$,
$ c_{16}^{1}$,
$0$;\ \ 
$ -c^{1}_{16}
-c^{3}_{16}
$,
$ -c_{16}^{1}$,
$ c_{16}^{1}$,
$0$;\ \ 
$ -1$,
$ -1$,
$ c^{1}_{16}
+c^{3}_{16}
$;\ \ 
$ -1$,
$ -c^{1}_{16}
-c^{3}_{16}
$;\ \ 
$0$)

  \vskip 2ex

\noindent14. $7_{\frac{29}{4},27.31}^{32,406}$ \irep{96}:\ \ 
$d_i$ = ($1.0$,
$1.0$,
$1.847$,
$1.847$,
$2.414$,
$2.414$,
$2.613$) 

\vskip 0.7ex
\hangindent=3em \hangafter=1
$D^2= 27.313 = 
16+8\sqrt{2}$

\vskip 0.7ex
\hangindent=3em \hangafter=1
$T = ( 0,
\frac{1}{2},
\frac{23}{32},
\frac{23}{32},
\frac{1}{4},
\frac{3}{4},
\frac{3}{32} )
$,

\vskip 0.7ex
\hangindent=3em \hangafter=1
$S$ = ($ 1$,
$ 1$,
$ c_{16}^{1}$,
$ c_{16}^{1}$,
$ 1+\sqrt{2}$,
$ 1+\sqrt{2}$,
$ c^{1}_{16}
+c^{3}_{16}
$;\ \ 
$ 1$,
$ -c_{16}^{1}$,
$ -c_{16}^{1}$,
$ 1+\sqrt{2}$,
$ 1+\sqrt{2}$,
$ -c^{1}_{16}
-c^{3}_{16}
$;\ \ 
$(-c^{1}_{16}
-c^{3}_{16}
)\mathrm{i}$,
$(c^{1}_{16}
+c^{3}_{16}
)\mathrm{i}$,
$ c_{16}^{1}$,
$ -c_{16}^{1}$,
$0$;\ \ 
$(-c^{1}_{16}
-c^{3}_{16}
)\mathrm{i}$,
$ c_{16}^{1}$,
$ -c_{16}^{1}$,
$0$;\ \ 
$ -1$,
$ -1$,
$ -c^{1}_{16}
-c^{3}_{16}
$;\ \ 
$ -1$,
$ c^{1}_{16}
+c^{3}_{16}
$;\ \ 
$0$)

  \vskip 2ex

\noindent15. $7_{\frac{19}{4},27.31}^{32,116}$ \irep{96}:\ \ 
$d_i$ = ($1.0$,
$1.0$,
$1.847$,
$1.847$,
$2.414$,
$2.414$,
$2.613$) 

\vskip 0.7ex
\hangindent=3em \hangafter=1
$D^2= 27.313 = 
16+8\sqrt{2}$

\vskip 0.7ex
\hangindent=3em \hangafter=1
$T = ( 0,
\frac{1}{2},
\frac{25}{32},
\frac{25}{32},
\frac{1}{4},
\frac{3}{4},
\frac{13}{32} )
$,

\vskip 0.7ex
\hangindent=3em \hangafter=1
$S$ = ($ 1$,
$ 1$,
$ c_{16}^{1}$,
$ c_{16}^{1}$,
$ 1+\sqrt{2}$,
$ 1+\sqrt{2}$,
$ c^{1}_{16}
+c^{3}_{16}
$;\ \ 
$ 1$,
$ -c_{16}^{1}$,
$ -c_{16}^{1}$,
$ 1+\sqrt{2}$,
$ 1+\sqrt{2}$,
$ -c^{1}_{16}
-c^{3}_{16}
$;\ \ 
$(c^{1}_{16}
+c^{3}_{16}
)\mathrm{i}$,
$(-c^{1}_{16}
-c^{3}_{16}
)\mathrm{i}$,
$ -c_{16}^{1}$,
$ c_{16}^{1}$,
$0$;\ \ 
$(c^{1}_{16}
+c^{3}_{16}
)\mathrm{i}$,
$ -c_{16}^{1}$,
$ c_{16}^{1}$,
$0$;\ \ 
$ -1$,
$ -1$,
$ c^{1}_{16}
+c^{3}_{16}
$;\ \ 
$ -1$,
$ -c^{1}_{16}
-c^{3}_{16}
$;\ \ 
$0$)

  \vskip 2ex

\noindent16. $7_{\frac{1}{4},27.31}^{32,123}$ \irep{97}:\ \ 
$d_i$ = ($1.0$,
$1.0$,
$1.847$,
$1.847$,
$2.414$,
$2.414$,
$2.613$) 

\vskip 0.7ex
\hangindent=3em \hangafter=1
$D^2= 27.313 = 
16+8\sqrt{2}$

\vskip 0.7ex
\hangindent=3em \hangafter=1
$T = ( 0,
\frac{1}{2},
\frac{27}{32},
\frac{27}{32},
\frac{1}{4},
\frac{3}{4},
\frac{7}{32} )
$,

\vskip 0.7ex
\hangindent=3em \hangafter=1
$S$ = ($ 1$,
$ 1$,
$ c_{16}^{1}$,
$ c_{16}^{1}$,
$ 1+\sqrt{2}$,
$ 1+\sqrt{2}$,
$ c^{1}_{16}
+c^{3}_{16}
$;\ \ 
$ 1$,
$ -c_{16}^{1}$,
$ -c_{16}^{1}$,
$ 1+\sqrt{2}$,
$ 1+\sqrt{2}$,
$ -c^{1}_{16}
-c^{3}_{16}
$;\ \ 
$ -c^{1}_{16}
-c^{3}_{16}
$,
$ c^{1}_{16}
+c^{3}_{16}
$,
$ c_{16}^{1}$,
$ -c_{16}^{1}$,
$0$;\ \ 
$ -c^{1}_{16}
-c^{3}_{16}
$,
$ c_{16}^{1}$,
$ -c_{16}^{1}$,
$0$;\ \ 
$ -1$,
$ -1$,
$ -c^{1}_{16}
-c^{3}_{16}
$;\ \ 
$ -1$,
$ c^{1}_{16}
+c^{3}_{16}
$;\ \ 
$0$)

  \vskip 2ex

\noindent17. $7_{\frac{23}{4},27.31}^{32,224}$ \irep{97}:\ \ 
$d_i$ = ($1.0$,
$1.0$,
$1.847$,
$1.847$,
$2.414$,
$2.414$,
$2.613$) 

\vskip 0.7ex
\hangindent=3em \hangafter=1
$D^2= 27.313 = 
16+8\sqrt{2}$

\vskip 0.7ex
\hangindent=3em \hangafter=1
$T = ( 0,
\frac{1}{2},
\frac{29}{32},
\frac{29}{32},
\frac{1}{4},
\frac{3}{4},
\frac{17}{32} )
$,

\vskip 0.7ex
\hangindent=3em \hangafter=1
$S$ = ($ 1$,
$ 1$,
$ c_{16}^{1}$,
$ c_{16}^{1}$,
$ 1+\sqrt{2}$,
$ 1+\sqrt{2}$,
$ c^{1}_{16}
+c^{3}_{16}
$;\ \ 
$ 1$,
$ -c_{16}^{1}$,
$ -c_{16}^{1}$,
$ 1+\sqrt{2}$,
$ 1+\sqrt{2}$,
$ -c^{1}_{16}
-c^{3}_{16}
$;\ \ 
$ c^{1}_{16}
+c^{3}_{16}
$,
$ -c^{1}_{16}
-c^{3}_{16}
$,
$ -c_{16}^{1}$,
$ c_{16}^{1}$,
$0$;\ \ 
$ c^{1}_{16}
+c^{3}_{16}
$,
$ -c_{16}^{1}$,
$ c_{16}^{1}$,
$0$;\ \ 
$ -1$,
$ -1$,
$ c^{1}_{16}
+c^{3}_{16}
$;\ \ 
$ -1$,
$ -c^{1}_{16}
-c^{3}_{16}
$;\ \ 
$0$)

  \vskip 2ex

\noindent18. $7_{\frac{5}{4},27.31}^{32,225}$ \irep{96}:\ \ 
$d_i$ = ($1.0$,
$1.0$,
$1.847$,
$1.847$,
$2.414$,
$2.414$,
$2.613$) 

\vskip 0.7ex
\hangindent=3em \hangafter=1
$D^2= 27.313 = 
16+8\sqrt{2}$

\vskip 0.7ex
\hangindent=3em \hangafter=1
$T = ( 0,
\frac{1}{2},
\frac{31}{32},
\frac{31}{32},
\frac{1}{4},
\frac{3}{4},
\frac{11}{32} )
$,

\vskip 0.7ex
\hangindent=3em \hangafter=1
$S$ = ($ 1$,
$ 1$,
$ c_{16}^{1}$,
$ c_{16}^{1}$,
$ 1+\sqrt{2}$,
$ 1+\sqrt{2}$,
$ c^{1}_{16}
+c^{3}_{16}
$;\ \ 
$ 1$,
$ -c_{16}^{1}$,
$ -c_{16}^{1}$,
$ 1+\sqrt{2}$,
$ 1+\sqrt{2}$,
$ -c^{1}_{16}
-c^{3}_{16}
$;\ \ 
$(c^{1}_{16}
+c^{3}_{16}
)\mathrm{i}$,
$(-c^{1}_{16}
-c^{3}_{16}
)\mathrm{i}$,
$ c_{16}^{1}$,
$ -c_{16}^{1}$,
$0$;\ \ 
$(c^{1}_{16}
+c^{3}_{16}
)\mathrm{i}$,
$ c_{16}^{1}$,
$ -c_{16}^{1}$,
$0$;\ \ 
$ -1$,
$ -1$,
$ -c^{1}_{16}
-c^{3}_{16}
$;\ \ 
$ -1$,
$ c^{1}_{16}
+c^{3}_{16}
$;\ \ 
$0$)

  \vskip 2ex

\noindent19. $7_{2,28.}^{56,139}$ \irep{99}:\ \ 
$d_i$ = ($1.0$,
$1.0$,
$2.0$,
$2.0$,
$2.0$,
$2.645$,
$2.645$) 

\vskip 0.7ex
\hangindent=3em \hangafter=1
$D^2= 28.0 = 
28$

\vskip 0.7ex
\hangindent=3em \hangafter=1
$T = ( 0,
0,
\frac{1}{7},
\frac{2}{7},
\frac{4}{7},
\frac{1}{8},
\frac{5}{8} )
$,

\vskip 0.7ex
\hangindent=3em \hangafter=1
$S$ = ($ 1$,
$ 1$,
$ 2$,
$ 2$,
$ 2$,
$ \sqrt{7}$,
$ \sqrt{7}$;\ \ 
$ 1$,
$ 2$,
$ 2$,
$ 2$,
$ -\sqrt{7}$,
$ -\sqrt{7}$;\ \ 
$ 2c_{7}^{2}$,
$ 2c_{7}^{1}$,
$ 2c_{7}^{3}$,
$0$,
$0$;\ \ 
$ 2c_{7}^{3}$,
$ 2c_{7}^{2}$,
$0$,
$0$;\ \ 
$ 2c_{7}^{1}$,
$0$,
$0$;\ \ 
$ \sqrt{7}$,
$ -\sqrt{7}$;\ \ 
$ \sqrt{7}$)

  \vskip 2ex

\noindent20. $7_{2,28.}^{56,680}$ \irep{99}:\ \ 
$d_i$ = ($1.0$,
$1.0$,
$2.0$,
$2.0$,
$2.0$,
$2.645$,
$2.645$) 

\vskip 0.7ex
\hangindent=3em \hangafter=1
$D^2= 28.0 = 
28$

\vskip 0.7ex
\hangindent=3em \hangafter=1
$T = ( 0,
0,
\frac{1}{7},
\frac{2}{7},
\frac{4}{7},
\frac{3}{8},
\frac{7}{8} )
$,

\vskip 0.7ex
\hangindent=3em \hangafter=1
$S$ = ($ 1$,
$ 1$,
$ 2$,
$ 2$,
$ 2$,
$ \sqrt{7}$,
$ \sqrt{7}$;\ \ 
$ 1$,
$ 2$,
$ 2$,
$ 2$,
$ -\sqrt{7}$,
$ -\sqrt{7}$;\ \ 
$ 2c_{7}^{2}$,
$ 2c_{7}^{1}$,
$ 2c_{7}^{3}$,
$0$,
$0$;\ \ 
$ 2c_{7}^{3}$,
$ 2c_{7}^{2}$,
$0$,
$0$;\ \ 
$ 2c_{7}^{1}$,
$0$,
$0$;\ \ 
$ -\sqrt{7}$,
$ \sqrt{7}$;\ \ 
$ -\sqrt{7}$)

  \vskip 2ex

\noindent21. $7_{6,28.}^{56,609}$ \irep{99}:\ \ 
$d_i$ = ($1.0$,
$1.0$,
$2.0$,
$2.0$,
$2.0$,
$2.645$,
$2.645$) 

\vskip 0.7ex
\hangindent=3em \hangafter=1
$D^2= 28.0 = 
28$

\vskip 0.7ex
\hangindent=3em \hangafter=1
$T = ( 0,
0,
\frac{3}{7},
\frac{5}{7},
\frac{6}{7},
\frac{1}{8},
\frac{5}{8} )
$,

\vskip 0.7ex
\hangindent=3em \hangafter=1
$S$ = ($ 1$,
$ 1$,
$ 2$,
$ 2$,
$ 2$,
$ \sqrt{7}$,
$ \sqrt{7}$;\ \ 
$ 1$,
$ 2$,
$ 2$,
$ 2$,
$ -\sqrt{7}$,
$ -\sqrt{7}$;\ \ 
$ 2c_{7}^{1}$,
$ 2c_{7}^{2}$,
$ 2c_{7}^{3}$,
$0$,
$0$;\ \ 
$ 2c_{7}^{3}$,
$ 2c_{7}^{1}$,
$0$,
$0$;\ \ 
$ 2c_{7}^{2}$,
$0$,
$0$;\ \ 
$ -\sqrt{7}$,
$ \sqrt{7}$;\ \ 
$ -\sqrt{7}$)

  \vskip 2ex

\noindent22. $7_{6,28.}^{56,193}$ \irep{99}:\ \ 
$d_i$ = ($1.0$,
$1.0$,
$2.0$,
$2.0$,
$2.0$,
$2.645$,
$2.645$) 

\vskip 0.7ex
\hangindent=3em \hangafter=1
$D^2= 28.0 = 
28$

\vskip 0.7ex
\hangindent=3em \hangafter=1
$T = ( 0,
0,
\frac{3}{7},
\frac{5}{7},
\frac{6}{7},
\frac{3}{8},
\frac{7}{8} )
$,

\vskip 0.7ex
\hangindent=3em \hangafter=1
$S$ = ($ 1$,
$ 1$,
$ 2$,
$ 2$,
$ 2$,
$ \sqrt{7}$,
$ \sqrt{7}$;\ \ 
$ 1$,
$ 2$,
$ 2$,
$ 2$,
$ -\sqrt{7}$,
$ -\sqrt{7}$;\ \ 
$ 2c_{7}^{1}$,
$ 2c_{7}^{2}$,
$ 2c_{7}^{3}$,
$0$,
$0$;\ \ 
$ 2c_{7}^{3}$,
$ 2c_{7}^{1}$,
$0$,
$0$;\ \ 
$ 2c_{7}^{2}$,
$0$,
$0$;\ \ 
$ \sqrt{7}$,
$ -\sqrt{7}$;\ \ 
$ \sqrt{7}$)

  \vskip 2ex

\noindent23. $7_{\frac{32}{5},86.75}^{15,205}$ \irep{80}:\ \ 
$d_i$ = ($1.0$,
$1.956$,
$2.827$,
$3.574$,
$4.165$,
$4.574$,
$4.783$) 

\vskip 0.7ex
\hangindent=3em \hangafter=1
$D^2= 86.750 = 
30+15c^{1}_{15}
+15c^{2}_{15}
+15c^{3}_{15}
$

\vskip 0.7ex
\hangindent=3em \hangafter=1
$T = ( 0,
\frac{1}{5},
\frac{13}{15},
0,
\frac{3}{5},
\frac{2}{3},
\frac{1}{5} )
$,

\vskip 0.7ex
\hangindent=3em \hangafter=1
$S$ = ($ 1$,
$ -c_{15}^{7}$,
$ \xi_{15}^{3}$,
$ \xi_{15}^{11}$,
$ \xi_{15}^{5}$,
$ \xi_{15}^{9}$,
$ \xi_{15}^{7}$;\ \ 
$ -\xi_{15}^{11}$,
$ \xi_{15}^{9}$,
$ -\xi_{15}^{7}$,
$ \xi_{15}^{5}$,
$ -\xi_{15}^{3}$,
$ 1$;\ \ 
$ \xi_{15}^{9}$,
$ \xi_{15}^{3}$,
$0$,
$ -\xi_{15}^{3}$,
$ -\xi_{15}^{9}$;\ \ 
$ 1$,
$ -\xi_{15}^{5}$,
$ \xi_{15}^{9}$,
$ c_{15}^{7}$;\ \ 
$ -\xi_{15}^{5}$,
$0$,
$ \xi_{15}^{5}$;\ \ 
$ -\xi_{15}^{9}$,
$ \xi_{15}^{3}$;\ \ 
$ -\xi_{15}^{11}$)

  \vskip 2ex

\noindent24. $7_{\frac{8}{5},86.75}^{15,181}$ \irep{80}:\ \ 
$d_i$ = ($1.0$,
$1.956$,
$2.827$,
$3.574$,
$4.165$,
$4.574$,
$4.783$) 

\vskip 0.7ex
\hangindent=3em \hangafter=1
$D^2= 86.750 = 
30+15c^{1}_{15}
+15c^{2}_{15}
+15c^{3}_{15}
$

\vskip 0.7ex
\hangindent=3em \hangafter=1
$T = ( 0,
\frac{4}{5},
\frac{2}{15},
0,
\frac{2}{5},
\frac{1}{3},
\frac{4}{5} )
$,

\vskip 0.7ex
\hangindent=3em \hangafter=1
$S$ = ($ 1$,
$ -c_{15}^{7}$,
$ \xi_{15}^{3}$,
$ \xi_{15}^{11}$,
$ \xi_{15}^{5}$,
$ \xi_{15}^{9}$,
$ \xi_{15}^{7}$;\ \ 
$ -\xi_{15}^{11}$,
$ \xi_{15}^{9}$,
$ -\xi_{15}^{7}$,
$ \xi_{15}^{5}$,
$ -\xi_{15}^{3}$,
$ 1$;\ \ 
$ \xi_{15}^{9}$,
$ \xi_{15}^{3}$,
$0$,
$ -\xi_{15}^{3}$,
$ -\xi_{15}^{9}$;\ \ 
$ 1$,
$ -\xi_{15}^{5}$,
$ \xi_{15}^{9}$,
$ c_{15}^{7}$;\ \ 
$ -\xi_{15}^{5}$,
$0$,
$ \xi_{15}^{5}$;\ \ 
$ -\xi_{15}^{9}$,
$ \xi_{15}^{3}$;\ \ 
$ -\xi_{15}^{11}$)

  \vskip 2ex

\noindent25. $7_{1,93.25}^{8,230}$ \irep{54}:\ \ 
$d_i$ = ($1.0$,
$2.414$,
$2.414$,
$3.414$,
$3.414$,
$4.828$,
$5.828$) 

\vskip 0.7ex
\hangindent=3em \hangafter=1
$D^2= 93.254 = 
48+32\sqrt{2}$

\vskip 0.7ex
\hangindent=3em \hangafter=1
$T = ( 0,
\frac{1}{2},
\frac{1}{2},
\frac{1}{4},
\frac{1}{4},
\frac{5}{8},
0 )
$,

\vskip 0.7ex
\hangindent=3em \hangafter=1
$S$ = ($ 1$,
$ 1+\sqrt{2}$,
$ 1+\sqrt{2}$,
$ 2+\sqrt{2}$,
$ 2+\sqrt{2}$,
$ 2+2\sqrt{2}$,
$ 3+2\sqrt{2}$;\ \ 
$ -1-2  \zeta^{1}_{8}
-2  \zeta^{2}_{8}
$,
$ -1-2  \zeta^{-1}_{8}
+2\zeta^{2}_{8}
$,
$(-2-\sqrt{2})\mathrm{i}$,
$(2+\sqrt{2})\mathrm{i}$,
$ 2+2\sqrt{2}$,
$ -1-\sqrt{2}$;\ \ 
$ -1-2  \zeta^{1}_{8}
-2  \zeta^{2}_{8}
$,
$(2+\sqrt{2})\mathrm{i}$,
$(-2-\sqrt{2})\mathrm{i}$,
$ 2+2\sqrt{2}$,
$ -1-\sqrt{2}$;\ \ 
$ (2+2\sqrt{2})\zeta_{8}^{3}$,
$ (-2-2\sqrt{2})\zeta_{8}^{1}$,
$0$,
$ 2+\sqrt{2}$;\ \ 
$ (2+2\sqrt{2})\zeta_{8}^{3}$,
$0$,
$ 2+\sqrt{2}$;\ \ 
$0$,
$ -2-2\sqrt{2}$;\ \ 
$ 1$)

  \vskip 2ex

\noindent26. $7_{7,93.25}^{8,101}$ \irep{54}:\ \ 
$d_i$ = ($1.0$,
$2.414$,
$2.414$,
$3.414$,
$3.414$,
$4.828$,
$5.828$) 

\vskip 0.7ex
\hangindent=3em \hangafter=1
$D^2= 93.254 = 
48+32\sqrt{2}$

\vskip 0.7ex
\hangindent=3em \hangafter=1
$T = ( 0,
\frac{1}{2},
\frac{1}{2},
\frac{3}{4},
\frac{3}{4},
\frac{3}{8},
0 )
$,

\vskip 0.7ex
\hangindent=3em \hangafter=1
$S$ = ($ 1$,
$ 1+\sqrt{2}$,
$ 1+\sqrt{2}$,
$ 2+\sqrt{2}$,
$ 2+\sqrt{2}$,
$ 2+2\sqrt{2}$,
$ 3+2\sqrt{2}$;\ \ 
$ -1-2  \zeta^{-1}_{8}
+2\zeta^{2}_{8}
$,
$ -1-2  \zeta^{1}_{8}
-2  \zeta^{2}_{8}
$,
$(-2-\sqrt{2})\mathrm{i}$,
$(2+\sqrt{2})\mathrm{i}$,
$ 2+2\sqrt{2}$,
$ -1-\sqrt{2}$;\ \ 
$ -1-2  \zeta^{-1}_{8}
+2\zeta^{2}_{8}
$,
$(2+\sqrt{2})\mathrm{i}$,
$(-2-\sqrt{2})\mathrm{i}$,
$ 2+2\sqrt{2}$,
$ -1-\sqrt{2}$;\ \ 
$ (-2-2\sqrt{2})\zeta_{8}^{1}$,
$ (2+2\sqrt{2})\zeta_{8}^{3}$,
$0$,
$ 2+\sqrt{2}$;\ \ 
$ (-2-2\sqrt{2})\zeta_{8}^{1}$,
$0$,
$ 2+\sqrt{2}$;\ \ 
$0$,
$ -2-2\sqrt{2}$;\ \ 
$ 1$)

  \vskip 2ex

\noindent27. $7_{\frac{30}{11},135.7}^{11,157}$ \irep{66}:\ \ 
$d_i$ = ($1.0$,
$2.918$,
$3.513$,
$3.513$,
$4.601$,
$5.911$,
$6.742$) 

\vskip 0.7ex
\hangindent=3em \hangafter=1
$D^2= 135.778 = 
55+44c^{1}_{11}
+33c^{2}_{11}
+22c^{3}_{11}
+11c^{4}_{11}
$

\vskip 0.7ex
\hangindent=3em \hangafter=1
$T = ( 0,
\frac{1}{11},
\frac{4}{11},
\frac{4}{11},
\frac{3}{11},
\frac{6}{11},
\frac{10}{11} )
$,

\vskip 0.7ex
\hangindent=3em \hangafter=1
$S$ = ($ 1$,
$ 2+c^{1}_{11}
+c^{2}_{11}
+c^{3}_{11}
+c^{4}_{11}
$,
$ \xi_{11}^{5}$,
$ \xi_{11}^{5}$,
$ 2+2c^{1}_{11}
+c^{2}_{11}
+c^{3}_{11}
+c^{4}_{11}
$,
$ 2+2c^{1}_{11}
+c^{2}_{11}
+c^{3}_{11}
$,
$ 2+2c^{1}_{11}
+2c^{2}_{11}
+c^{3}_{11}
$;\ \ 
$ 2+2c^{1}_{11}
+2c^{2}_{11}
+c^{3}_{11}
$,
$ -\xi_{11}^{5}$,
$ -\xi_{11}^{5}$,
$ 2+2c^{1}_{11}
+c^{2}_{11}
+c^{3}_{11}
$,
$ 1$,
$ -2-2  c^{1}_{11}
-c^{2}_{11}
-c^{3}_{11}
-c^{4}_{11}
$;\ \ 
$ s^{2}_{11}
+2\zeta^{3}_{11}
-\zeta^{-3}_{11}
+\zeta^{4}_{11}
+\zeta^{5}_{11}
$,
$ -1-c^{1}_{11}
-2  \zeta^{2}_{11}
-2  \zeta^{3}_{11}
+\zeta^{-3}_{11}
-\zeta^{4}_{11}
-\zeta^{5}_{11}
$,
$ \xi_{11}^{5}$,
$ -\xi_{11}^{5}$,
$ \xi_{11}^{5}$;\ \ 
$ s^{2}_{11}
+2\zeta^{3}_{11}
-\zeta^{-3}_{11}
+\zeta^{4}_{11}
+\zeta^{5}_{11}
$,
$ \xi_{11}^{5}$,
$ -\xi_{11}^{5}$,
$ \xi_{11}^{5}$;\ \ 
$ -2-c^{1}_{11}
-c^{2}_{11}
-c^{3}_{11}
-c^{4}_{11}
$,
$ -2-2  c^{1}_{11}
-2  c^{2}_{11}
-c^{3}_{11}
$,
$ 1$;\ \ 
$ 2+2c^{1}_{11}
+c^{2}_{11}
+c^{3}_{11}
+c^{4}_{11}
$,
$ 2+c^{1}_{11}
+c^{2}_{11}
+c^{3}_{11}
+c^{4}_{11}
$;\ \ 
$ -2-2  c^{1}_{11}
-c^{2}_{11}
-c^{3}_{11}
$)

  \vskip 2ex

\noindent28. $7_{\frac{58}{11},135.7}^{11,191}$ \irep{66}:\ \ 
$d_i$ = ($1.0$,
$2.918$,
$3.513$,
$3.513$,
$4.601$,
$5.911$,
$6.742$) 

\vskip 0.7ex
\hangindent=3em \hangafter=1
$D^2= 135.778 = 
55+44c^{1}_{11}
+33c^{2}_{11}
+22c^{3}_{11}
+11c^{4}_{11}
$

\vskip 0.7ex
\hangindent=3em \hangafter=1
$T = ( 0,
\frac{10}{11},
\frac{7}{11},
\frac{7}{11},
\frac{8}{11},
\frac{5}{11},
\frac{1}{11} )
$,

\vskip 0.7ex
\hangindent=3em \hangafter=1
$S$ = ($ 1$,
$ 2+c^{1}_{11}
+c^{2}_{11}
+c^{3}_{11}
+c^{4}_{11}
$,
$ \xi_{11}^{5}$,
$ \xi_{11}^{5}$,
$ 2+2c^{1}_{11}
+c^{2}_{11}
+c^{3}_{11}
+c^{4}_{11}
$,
$ 2+2c^{1}_{11}
+c^{2}_{11}
+c^{3}_{11}
$,
$ 2+2c^{1}_{11}
+2c^{2}_{11}
+c^{3}_{11}
$;\ \ 
$ 2+2c^{1}_{11}
+2c^{2}_{11}
+c^{3}_{11}
$,
$ -\xi_{11}^{5}$,
$ -\xi_{11}^{5}$,
$ 2+2c^{1}_{11}
+c^{2}_{11}
+c^{3}_{11}
$,
$ 1$,
$ -2-2  c^{1}_{11}
-c^{2}_{11}
-c^{3}_{11}
-c^{4}_{11}
$;\ \ 
$ -1-c^{1}_{11}
-2  \zeta^{2}_{11}
-2  \zeta^{3}_{11}
+\zeta^{-3}_{11}
-\zeta^{4}_{11}
-\zeta^{5}_{11}
$,
$ s^{2}_{11}
+2\zeta^{3}_{11}
-\zeta^{-3}_{11}
+\zeta^{4}_{11}
+\zeta^{5}_{11}
$,
$ \xi_{11}^{5}$,
$ -\xi_{11}^{5}$,
$ \xi_{11}^{5}$;\ \ 
$ -1-c^{1}_{11}
-2  \zeta^{2}_{11}
-2  \zeta^{3}_{11}
+\zeta^{-3}_{11}
-\zeta^{4}_{11}
-\zeta^{5}_{11}
$,
$ \xi_{11}^{5}$,
$ -\xi_{11}^{5}$,
$ \xi_{11}^{5}$;\ \ 
$ -2-c^{1}_{11}
-c^{2}_{11}
-c^{3}_{11}
-c^{4}_{11}
$,
$ -2-2  c^{1}_{11}
-2  c^{2}_{11}
-c^{3}_{11}
$,
$ 1$;\ \ 
$ 2+2c^{1}_{11}
+c^{2}_{11}
+c^{3}_{11}
+c^{4}_{11}
$,
$ 2+c^{1}_{11}
+c^{2}_{11}
+c^{3}_{11}
+c^{4}_{11}
$;\ \ 
$ -2-2  c^{1}_{11}
-c^{2}_{11}
-c^{3}_{11}
$)

  \vskip 2ex 

}

\subsection{Rank 8}
\label{uni8}

{\small

\noindent1. $8_{1,8.}^{4,100}$ \irep{0}:\ \ 
$d_i$ = ($1.0$,
$1.0$,
$1.0$,
$1.0$,
$1.0$,
$1.0$,
$1.0$,
$1.0$) 

\vskip 0.7ex
\hangindent=3em \hangafter=1
$D^2= 8.0 = 
8$

\vskip 0.7ex
\hangindent=3em \hangafter=1
$T = ( 0,
0,
0,
\frac{1}{2},
\frac{1}{4},
\frac{1}{4},
\frac{1}{4},
\frac{3}{4} )
$,

\vskip 0.7ex
\hangindent=3em \hangafter=1
$S$ = ($ 1$,
$ 1$,
$ 1$,
$ 1$,
$ 1$,
$ 1$,
$ 1$,
$ 1$;\ \ 
$ 1$,
$ -1$,
$ -1$,
$ -1$,
$ 1$,
$ 1$,
$ -1$;\ \ 
$ 1$,
$ -1$,
$ 1$,
$ 1$,
$ -1$,
$ -1$;\ \ 
$ 1$,
$ -1$,
$ 1$,
$ -1$,
$ 1$;\ \ 
$ -1$,
$ -1$,
$ 1$,
$ 1$;\ \ 
$ -1$,
$ -1$,
$ -1$;\ \ 
$ -1$,
$ 1$;\ \ 
$ -1$)

Factors = $2_{1,2.}^{4,437}\boxtimes 4_{0,4.}^{2,750}$

  \vskip 2ex

\noindent2. $8_{7,8.}^{4,000}$ \irep{0}:\ \ 
$d_i$ = ($1.0$,
$1.0$,
$1.0$,
$1.0$,
$1.0$,
$1.0$,
$1.0$,
$1.0$) 

\vskip 0.7ex
\hangindent=3em \hangafter=1
$D^2= 8.0 = 
8$

\vskip 0.7ex
\hangindent=3em \hangafter=1
$T = ( 0,
0,
0,
\frac{1}{2},
\frac{1}{4},
\frac{3}{4},
\frac{3}{4},
\frac{3}{4} )
$,

\vskip 0.7ex
\hangindent=3em \hangafter=1
$S$ = ($ 1$,
$ 1$,
$ 1$,
$ 1$,
$ 1$,
$ 1$,
$ 1$,
$ 1$;\ \ 
$ 1$,
$ -1$,
$ -1$,
$ -1$,
$ 1$,
$ -1$,
$ 1$;\ \ 
$ 1$,
$ -1$,
$ -1$,
$ -1$,
$ 1$,
$ 1$;\ \ 
$ 1$,
$ 1$,
$ -1$,
$ -1$,
$ 1$;\ \ 
$ -1$,
$ 1$,
$ 1$,
$ -1$;\ \ 
$ -1$,
$ 1$,
$ -1$;\ \ 
$ -1$,
$ -1$;\ \ 
$ -1$)

Factors = $2_{7,2.}^{4,625}\boxtimes 4_{0,4.}^{2,750}$

  \vskip 2ex

\noindent3. $8_{3,8.}^{4,500}$ \irep{0}:\ \ 
$d_i$ = ($1.0$,
$1.0$,
$1.0$,
$1.0$,
$1.0$,
$1.0$,
$1.0$,
$1.0$) 

\vskip 0.7ex
\hangindent=3em \hangafter=1
$D^2= 8.0 = 
8$

\vskip 0.7ex
\hangindent=3em \hangafter=1
$T = ( 0,
\frac{1}{2},
\frac{1}{2},
\frac{1}{2},
\frac{1}{4},
\frac{1}{4},
\frac{1}{4},
\frac{3}{4} )
$,

\vskip 0.7ex
\hangindent=3em \hangafter=1
$S$ = ($ 1$,
$ 1$,
$ 1$,
$ 1$,
$ 1$,
$ 1$,
$ 1$,
$ 1$;\ \ 
$ 1$,
$ -1$,
$ -1$,
$ -1$,
$ -1$,
$ 1$,
$ 1$;\ \ 
$ 1$,
$ -1$,
$ -1$,
$ 1$,
$ -1$,
$ 1$;\ \ 
$ 1$,
$ 1$,
$ -1$,
$ -1$,
$ 1$;\ \ 
$ -1$,
$ 1$,
$ 1$,
$ -1$;\ \ 
$ -1$,
$ 1$,
$ -1$;\ \ 
$ -1$,
$ -1$;\ \ 
$ -1$)

Factors = $2_{7,2.}^{4,625}\boxtimes 4_{4,4.}^{2,250}$

  \vskip 2ex

\noindent4. $8_{5,8.}^{4,500}$ \irep{0}:\ \ 
$d_i$ = ($1.0$,
$1.0$,
$1.0$,
$1.0$,
$1.0$,
$1.0$,
$1.0$,
$1.0$) 

\vskip 0.7ex
\hangindent=3em \hangafter=1
$D^2= 8.0 = 
8$

\vskip 0.7ex
\hangindent=3em \hangafter=1
$T = ( 0,
\frac{1}{2},
\frac{1}{2},
\frac{1}{2},
\frac{1}{4},
\frac{3}{4},
\frac{3}{4},
\frac{3}{4} )
$,

\vskip 0.7ex
\hangindent=3em \hangafter=1
$S$ = ($ 1$,
$ 1$,
$ 1$,
$ 1$,
$ 1$,
$ 1$,
$ 1$,
$ 1$;\ \ 
$ 1$,
$ -1$,
$ -1$,
$ 1$,
$ -1$,
$ -1$,
$ 1$;\ \ 
$ 1$,
$ -1$,
$ 1$,
$ -1$,
$ 1$,
$ -1$;\ \ 
$ 1$,
$ 1$,
$ 1$,
$ -1$,
$ -1$;\ \ 
$ -1$,
$ -1$,
$ -1$,
$ -1$;\ \ 
$ -1$,
$ 1$,
$ 1$;\ \ 
$ -1$,
$ 1$;\ \ 
$ -1$)

Factors = $2_{1,2.}^{4,437}\boxtimes 4_{4,4.}^{2,250}$

  \vskip 2ex

\noindent5. $8_{2,8.}^{8,812}$ \irep{103}:\ \ 
$d_i$ = ($1.0$,
$1.0$,
$1.0$,
$1.0$,
$1.0$,
$1.0$,
$1.0$,
$1.0$) 

\vskip 0.7ex
\hangindent=3em \hangafter=1
$D^2= 8.0 = 
8$

\vskip 0.7ex
\hangindent=3em \hangafter=1
$T = ( 0,
\frac{1}{2},
\frac{1}{4},
\frac{3}{4},
\frac{1}{8},
\frac{1}{8},
\frac{3}{8},
\frac{3}{8} )
$,

\vskip 0.7ex
\hangindent=3em \hangafter=1
$S$ = ($ 1$,
$ 1$,
$ 1$,
$ 1$,
$ 1$,
$ 1$,
$ 1$,
$ 1$;\ \ 
$ 1$,
$ 1$,
$ 1$,
$ -1$,
$ -1$,
$ -1$,
$ -1$;\ \ 
$ -1$,
$ -1$,
$ 1$,
$ 1$,
$ -1$,
$ -1$;\ \ 
$ -1$,
$ -1$,
$ -1$,
$ 1$,
$ 1$;\ \ 
$-\mathrm{i}$,
$\mathrm{i}$,
$-\mathrm{i}$,
$\mathrm{i}$;\ \ 
$-\mathrm{i}$,
$\mathrm{i}$,
$-\mathrm{i}$;\ \ 
$\mathrm{i}$,
$-\mathrm{i}$;\ \ 
$\mathrm{i}$)

Factors = $2_{1,2.}^{4,437}\boxtimes 4_{1,4.}^{8,718}$

  \vskip 2ex

\noindent6. $8_{0,8.}^{8,437}$ \irep{103}:\ \ 
$d_i$ = ($1.0$,
$1.0$,
$1.0$,
$1.0$,
$1.0$,
$1.0$,
$1.0$,
$1.0$) 

\vskip 0.7ex
\hangindent=3em \hangafter=1
$D^2= 8.0 = 
8$

\vskip 0.7ex
\hangindent=3em \hangafter=1
$T = ( 0,
\frac{1}{2},
\frac{1}{4},
\frac{3}{4},
\frac{1}{8},
\frac{1}{8},
\frac{7}{8},
\frac{7}{8} )
$,

\vskip 0.7ex
\hangindent=3em \hangafter=1
$S$ = ($ 1$,
$ 1$,
$ 1$,
$ 1$,
$ 1$,
$ 1$,
$ 1$,
$ 1$;\ \ 
$ 1$,
$ 1$,
$ 1$,
$ -1$,
$ -1$,
$ -1$,
$ -1$;\ \ 
$ -1$,
$ -1$,
$ -1$,
$ -1$,
$ 1$,
$ 1$;\ \ 
$ -1$,
$ 1$,
$ 1$,
$ -1$,
$ -1$;\ \ 
$-\mathrm{i}$,
$\mathrm{i}$,
$-\mathrm{i}$,
$\mathrm{i}$;\ \ 
$-\mathrm{i}$,
$\mathrm{i}$,
$-\mathrm{i}$;\ \ 
$\mathrm{i}$,
$-\mathrm{i}$;\ \ 
$\mathrm{i}$)

Factors = $2_{1,2.}^{4,437}\boxtimes 4_{7,4.}^{8,781}$

  \vskip 2ex

\noindent7. $8_{4,8.}^{8,625}$ \irep{103}:\ \ 
$d_i$ = ($1.0$,
$1.0$,
$1.0$,
$1.0$,
$1.0$,
$1.0$,
$1.0$,
$1.0$) 

\vskip 0.7ex
\hangindent=3em \hangafter=1
$D^2= 8.0 = 
8$

\vskip 0.7ex
\hangindent=3em \hangafter=1
$T = ( 0,
\frac{1}{2},
\frac{1}{4},
\frac{3}{4},
\frac{3}{8},
\frac{3}{8},
\frac{5}{8},
\frac{5}{8} )
$,

\vskip 0.7ex
\hangindent=3em \hangafter=1
$S$ = ($ 1$,
$ 1$,
$ 1$,
$ 1$,
$ 1$,
$ 1$,
$ 1$,
$ 1$;\ \ 
$ 1$,
$ 1$,
$ 1$,
$ -1$,
$ -1$,
$ -1$,
$ -1$;\ \ 
$ -1$,
$ -1$,
$ 1$,
$ 1$,
$ -1$,
$ -1$;\ \ 
$ -1$,
$ -1$,
$ -1$,
$ 1$,
$ 1$;\ \ 
$\mathrm{i}$,
$-\mathrm{i}$,
$-\mathrm{i}$,
$\mathrm{i}$;\ \ 
$\mathrm{i}$,
$\mathrm{i}$,
$-\mathrm{i}$;\ \ 
$-\mathrm{i}$,
$\mathrm{i}$;\ \ 
$-\mathrm{i}$)

Factors = $2_{1,2.}^{4,437}\boxtimes 4_{3,4.}^{8,468}$

  \vskip 2ex

\noindent8. $8_{6,8.}^{8,118}$ \irep{103}:\ \ 
$d_i$ = ($1.0$,
$1.0$,
$1.0$,
$1.0$,
$1.0$,
$1.0$,
$1.0$,
$1.0$) 

\vskip 0.7ex
\hangindent=3em \hangafter=1
$D^2= 8.0 = 
8$

\vskip 0.7ex
\hangindent=3em \hangafter=1
$T = ( 0,
\frac{1}{2},
\frac{1}{4},
\frac{3}{4},
\frac{5}{8},
\frac{5}{8},
\frac{7}{8},
\frac{7}{8} )
$,

\vskip 0.7ex
\hangindent=3em \hangafter=1
$S$ = ($ 1$,
$ 1$,
$ 1$,
$ 1$,
$ 1$,
$ 1$,
$ 1$,
$ 1$;\ \ 
$ 1$,
$ 1$,
$ 1$,
$ -1$,
$ -1$,
$ -1$,
$ -1$;\ \ 
$ -1$,
$ -1$,
$ 1$,
$ 1$,
$ -1$,
$ -1$;\ \ 
$ -1$,
$ -1$,
$ -1$,
$ 1$,
$ 1$;\ \ 
$-\mathrm{i}$,
$\mathrm{i}$,
$-\mathrm{i}$,
$\mathrm{i}$;\ \ 
$-\mathrm{i}$,
$\mathrm{i}$,
$-\mathrm{i}$;\ \ 
$\mathrm{i}$,
$-\mathrm{i}$;\ \ 
$\mathrm{i}$)

Factors = $2_{1,2.}^{4,437}\boxtimes 4_{5,4.}^{8,312}$

  \vskip 2ex

\noindent9. $8_{1,8.}^{16,123}$ \irep{0}:\ \ 
$d_i$ = ($1.0$,
$1.0$,
$1.0$,
$1.0$,
$1.0$,
$1.0$,
$1.0$,
$1.0$) 

\vskip 0.7ex
\hangindent=3em \hangafter=1
$D^2= 8.0 = 
8$

\vskip 0.7ex
\hangindent=3em \hangafter=1
$T = ( 0,
0,
\frac{1}{4},
\frac{1}{4},
\frac{1}{16},
\frac{1}{16},
\frac{9}{16},
\frac{9}{16} )
$,

\vskip 0.7ex
\hangindent=3em \hangafter=1
$S$ = ($ 1$,
$ 1$,
$ 1$,
$ 1$,
$ 1$,
$ 1$,
$ 1$,
$ 1$;\ \ 
$ 1$,
$ 1$,
$ 1$,
$ -1$,
$ -1$,
$ -1$,
$ -1$;\ \ 
$ -1$,
$ -1$,
$-\mathrm{i}$,
$\mathrm{i}$,
$-\mathrm{i}$,
$\mathrm{i}$;\ \ 
$ -1$,
$\mathrm{i}$,
$-\mathrm{i}$,
$\mathrm{i}$,
$-\mathrm{i}$;\ \ 
$ -\zeta_{8}^{3}$,
$ \zeta_{8}^{1}$,
$ \zeta_{8}^{3}$,
$ -\zeta_{8}^{1}$;\ \ 
$ -\zeta_{8}^{3}$,
$ -\zeta_{8}^{1}$,
$ \zeta_{8}^{3}$;\ \ 
$ -\zeta_{8}^{3}$,
$ \zeta_{8}^{1}$;\ \ 
$ -\zeta_{8}^{3}$)

  \vskip 2ex

\noindent10. $8_{1,8.}^{16,359}$ \irep{0}:\ \ 
$d_i$ = ($1.0$,
$1.0$,
$1.0$,
$1.0$,
$1.0$,
$1.0$,
$1.0$,
$1.0$) 

\vskip 0.7ex
\hangindent=3em \hangafter=1
$D^2= 8.0 = 
8$

\vskip 0.7ex
\hangindent=3em \hangafter=1
$T = ( 0,
0,
\frac{1}{4},
\frac{1}{4},
\frac{5}{16},
\frac{5}{16},
\frac{13}{16},
\frac{13}{16} )
$,

\vskip 0.7ex
\hangindent=3em \hangafter=1
$S$ = ($ 1$,
$ 1$,
$ 1$,
$ 1$,
$ 1$,
$ 1$,
$ 1$,
$ 1$;\ \ 
$ 1$,
$ 1$,
$ 1$,
$ -1$,
$ -1$,
$ -1$,
$ -1$;\ \ 
$ -1$,
$ -1$,
$-\mathrm{i}$,
$\mathrm{i}$,
$-\mathrm{i}$,
$\mathrm{i}$;\ \ 
$ -1$,
$\mathrm{i}$,
$-\mathrm{i}$,
$\mathrm{i}$,
$-\mathrm{i}$;\ \ 
$ \zeta_{8}^{3}$,
$ -\zeta_{8}^{1}$,
$ -\zeta_{8}^{3}$,
$ \zeta_{8}^{1}$;\ \ 
$ \zeta_{8}^{3}$,
$ \zeta_{8}^{1}$,
$ -\zeta_{8}^{3}$;\ \ 
$ \zeta_{8}^{3}$,
$ -\zeta_{8}^{1}$;\ \ 
$ \zeta_{8}^{3}$)

  \vskip 2ex

\noindent11. $8_{7,8.}^{16,140}$ \irep{0}:\ \ 
$d_i$ = ($1.0$,
$1.0$,
$1.0$,
$1.0$,
$1.0$,
$1.0$,
$1.0$,
$1.0$) 

\vskip 0.7ex
\hangindent=3em \hangafter=1
$D^2= 8.0 = 
8$

\vskip 0.7ex
\hangindent=3em \hangafter=1
$T = ( 0,
0,
\frac{3}{4},
\frac{3}{4},
\frac{3}{16},
\frac{3}{16},
\frac{11}{16},
\frac{11}{16} )
$,

\vskip 0.7ex
\hangindent=3em \hangafter=1
$S$ = ($ 1$,
$ 1$,
$ 1$,
$ 1$,
$ 1$,
$ 1$,
$ 1$,
$ 1$;\ \ 
$ 1$,
$ 1$,
$ 1$,
$ -1$,
$ -1$,
$ -1$,
$ -1$;\ \ 
$ -1$,
$ -1$,
$-\mathrm{i}$,
$\mathrm{i}$,
$-\mathrm{i}$,
$\mathrm{i}$;\ \ 
$ -1$,
$\mathrm{i}$,
$-\mathrm{i}$,
$\mathrm{i}$,
$-\mathrm{i}$;\ \ 
$ -\zeta_{8}^{1}$,
$ \zeta_{8}^{3}$,
$ \zeta_{8}^{1}$,
$ -\zeta_{8}^{3}$;\ \ 
$ -\zeta_{8}^{1}$,
$ -\zeta_{8}^{3}$,
$ \zeta_{8}^{1}$;\ \ 
$ -\zeta_{8}^{1}$,
$ \zeta_{8}^{3}$;\ \ 
$ -\zeta_{8}^{1}$)

  \vskip 2ex

\noindent12. $8_{7,8.}^{16,126}$ \irep{0}:\ \ 
$d_i$ = ($1.0$,
$1.0$,
$1.0$,
$1.0$,
$1.0$,
$1.0$,
$1.0$,
$1.0$) 

\vskip 0.7ex
\hangindent=3em \hangafter=1
$D^2= 8.0 = 
8$

\vskip 0.7ex
\hangindent=3em \hangafter=1
$T = ( 0,
0,
\frac{3}{4},
\frac{3}{4},
\frac{7}{16},
\frac{7}{16},
\frac{15}{16},
\frac{15}{16} )
$,

\vskip 0.7ex
\hangindent=3em \hangafter=1
$S$ = ($ 1$,
$ 1$,
$ 1$,
$ 1$,
$ 1$,
$ 1$,
$ 1$,
$ 1$;\ \ 
$ 1$,
$ 1$,
$ 1$,
$ -1$,
$ -1$,
$ -1$,
$ -1$;\ \ 
$ -1$,
$ -1$,
$-\mathrm{i}$,
$\mathrm{i}$,
$-\mathrm{i}$,
$\mathrm{i}$;\ \ 
$ -1$,
$\mathrm{i}$,
$-\mathrm{i}$,
$\mathrm{i}$,
$-\mathrm{i}$;\ \ 
$ \zeta_{8}^{1}$,
$ -\zeta_{8}^{3}$,
$ -\zeta_{8}^{1}$,
$ \zeta_{8}^{3}$;\ \ 
$ \zeta_{8}^{1}$,
$ \zeta_{8}^{3}$,
$ -\zeta_{8}^{1}$;\ \ 
$ \zeta_{8}^{1}$,
$ -\zeta_{8}^{3}$;\ \ 
$ \zeta_{8}^{1}$)

  \vskip 2ex

\noindent13. $8_{\frac{14}{5},14.47}^{10,280}$ \irep{138}:\ \ 
$d_i$ = ($1.0$,
$1.0$,
$1.0$,
$1.0$,
$1.618$,
$1.618$,
$1.618$,
$1.618$) 

\vskip 0.7ex
\hangindent=3em \hangafter=1
$D^2= 14.472 = 
10+2\sqrt{5}$

\vskip 0.7ex
\hangindent=3em \hangafter=1
$T = ( 0,
0,
0,
\frac{1}{2},
\frac{2}{5},
\frac{2}{5},
\frac{2}{5},
\frac{9}{10} )
$,

\vskip 0.7ex
\hangindent=3em \hangafter=1
$S$ = ($ 1$,
$ 1$,
$ 1$,
$ 1$,
$ \frac{1+\sqrt{5}}{2}$,
$ \frac{1+\sqrt{5}}{2}$,
$ \frac{1+\sqrt{5}}{2}$,
$ \frac{1+\sqrt{5}}{2}$;\ \ 
$ 1$,
$ -1$,
$ -1$,
$ \frac{1+\sqrt{5}}{2}$,
$ \frac{1+\sqrt{5}}{2}$,
$ -\frac{1+\sqrt{5}}{2}$,
$ -\frac{1+\sqrt{5}}{2}$;\ \ 
$ 1$,
$ -1$,
$ \frac{1+\sqrt{5}}{2}$,
$ -\frac{1+\sqrt{5}}{2}$,
$ \frac{1+\sqrt{5}}{2}$,
$ -\frac{1+\sqrt{5}}{2}$;\ \ 
$ 1$,
$ \frac{1+\sqrt{5}}{2}$,
$ -\frac{1+\sqrt{5}}{2}$,
$ -\frac{1+\sqrt{5}}{2}$,
$ \frac{1+\sqrt{5}}{2}$;\ \ 
$ -1$,
$ -1$,
$ -1$,
$ -1$;\ \ 
$ -1$,
$ 1$,
$ 1$;\ \ 
$ -1$,
$ 1$;\ \ 
$ -1$)

Factors = $2_{\frac{14}{5},3.618}^{5,395}\boxtimes 4_{0,4.}^{2,750}$

  \vskip 2ex

\noindent14. $8_{\frac{26}{5},14.47}^{10,604}$ \irep{138}:\ \ 
$d_i$ = ($1.0$,
$1.0$,
$1.0$,
$1.0$,
$1.618$,
$1.618$,
$1.618$,
$1.618$) 

\vskip 0.7ex
\hangindent=3em \hangafter=1
$D^2= 14.472 = 
10+2\sqrt{5}$

\vskip 0.7ex
\hangindent=3em \hangafter=1
$T = ( 0,
0,
0,
\frac{1}{2},
\frac{3}{5},
\frac{3}{5},
\frac{3}{5},
\frac{1}{10} )
$,

\vskip 0.7ex
\hangindent=3em \hangafter=1
$S$ = ($ 1$,
$ 1$,
$ 1$,
$ 1$,
$ \frac{1+\sqrt{5}}{2}$,
$ \frac{1+\sqrt{5}}{2}$,
$ \frac{1+\sqrt{5}}{2}$,
$ \frac{1+\sqrt{5}}{2}$;\ \ 
$ 1$,
$ -1$,
$ -1$,
$ \frac{1+\sqrt{5}}{2}$,
$ \frac{1+\sqrt{5}}{2}$,
$ -\frac{1+\sqrt{5}}{2}$,
$ -\frac{1+\sqrt{5}}{2}$;\ \ 
$ 1$,
$ -1$,
$ \frac{1+\sqrt{5}}{2}$,
$ -\frac{1+\sqrt{5}}{2}$,
$ \frac{1+\sqrt{5}}{2}$,
$ -\frac{1+\sqrt{5}}{2}$;\ \ 
$ 1$,
$ \frac{1+\sqrt{5}}{2}$,
$ -\frac{1+\sqrt{5}}{2}$,
$ -\frac{1+\sqrt{5}}{2}$,
$ \frac{1+\sqrt{5}}{2}$;\ \ 
$ -1$,
$ -1$,
$ -1$,
$ -1$;\ \ 
$ -1$,
$ 1$,
$ 1$;\ \ 
$ -1$,
$ 1$;\ \ 
$ -1$)

Factors = $2_{\frac{26}{5},3.618}^{5,720}\boxtimes 4_{0,4.}^{2,750}$

  \vskip 2ex

\noindent15. $8_{\frac{34}{5},14.47}^{10,232}$ \irep{138}:\ \ 
$d_i$ = ($1.0$,
$1.0$,
$1.0$,
$1.0$,
$1.618$,
$1.618$,
$1.618$,
$1.618$) 

\vskip 0.7ex
\hangindent=3em \hangafter=1
$D^2= 14.472 = 
10+2\sqrt{5}$

\vskip 0.7ex
\hangindent=3em \hangafter=1
$T = ( 0,
\frac{1}{2},
\frac{1}{2},
\frac{1}{2},
\frac{2}{5},
\frac{9}{10},
\frac{9}{10},
\frac{9}{10} )
$,

\vskip 0.7ex
\hangindent=3em \hangafter=1
$S$ = ($ 1$,
$ 1$,
$ 1$,
$ 1$,
$ \frac{1+\sqrt{5}}{2}$,
$ \frac{1+\sqrt{5}}{2}$,
$ \frac{1+\sqrt{5}}{2}$,
$ \frac{1+\sqrt{5}}{2}$;\ \ 
$ 1$,
$ -1$,
$ -1$,
$ \frac{1+\sqrt{5}}{2}$,
$ \frac{1+\sqrt{5}}{2}$,
$ -\frac{1+\sqrt{5}}{2}$,
$ -\frac{1+\sqrt{5}}{2}$;\ \ 
$ 1$,
$ -1$,
$ \frac{1+\sqrt{5}}{2}$,
$ -\frac{1+\sqrt{5}}{2}$,
$ \frac{1+\sqrt{5}}{2}$,
$ -\frac{1+\sqrt{5}}{2}$;\ \ 
$ 1$,
$ \frac{1+\sqrt{5}}{2}$,
$ -\frac{1+\sqrt{5}}{2}$,
$ -\frac{1+\sqrt{5}}{2}$,
$ \frac{1+\sqrt{5}}{2}$;\ \ 
$ -1$,
$ -1$,
$ -1$,
$ -1$;\ \ 
$ -1$,
$ 1$,
$ 1$;\ \ 
$ -1$,
$ 1$;\ \ 
$ -1$)

Factors = $2_{\frac{14}{5},3.618}^{5,395}\boxtimes 4_{4,4.}^{2,250}$

  \vskip 2ex

\noindent16. $8_{\frac{6}{5},14.47}^{10,123}$ \irep{138}:\ \ 
$d_i$ = ($1.0$,
$1.0$,
$1.0$,
$1.0$,
$1.618$,
$1.618$,
$1.618$,
$1.618$) 

\vskip 0.7ex
\hangindent=3em \hangafter=1
$D^2= 14.472 = 
10+2\sqrt{5}$

\vskip 0.7ex
\hangindent=3em \hangafter=1
$T = ( 0,
\frac{1}{2},
\frac{1}{2},
\frac{1}{2},
\frac{3}{5},
\frac{1}{10},
\frac{1}{10},
\frac{1}{10} )
$,

\vskip 0.7ex
\hangindent=3em \hangafter=1
$S$ = ($ 1$,
$ 1$,
$ 1$,
$ 1$,
$ \frac{1+\sqrt{5}}{2}$,
$ \frac{1+\sqrt{5}}{2}$,
$ \frac{1+\sqrt{5}}{2}$,
$ \frac{1+\sqrt{5}}{2}$;\ \ 
$ 1$,
$ -1$,
$ -1$,
$ \frac{1+\sqrt{5}}{2}$,
$ \frac{1+\sqrt{5}}{2}$,
$ -\frac{1+\sqrt{5}}{2}$,
$ -\frac{1+\sqrt{5}}{2}$;\ \ 
$ 1$,
$ -1$,
$ \frac{1+\sqrt{5}}{2}$,
$ -\frac{1+\sqrt{5}}{2}$,
$ \frac{1+\sqrt{5}}{2}$,
$ -\frac{1+\sqrt{5}}{2}$;\ \ 
$ 1$,
$ \frac{1+\sqrt{5}}{2}$,
$ -\frac{1+\sqrt{5}}{2}$,
$ -\frac{1+\sqrt{5}}{2}$,
$ \frac{1+\sqrt{5}}{2}$;\ \ 
$ -1$,
$ -1$,
$ -1$,
$ -1$;\ \ 
$ -1$,
$ 1$,
$ 1$;\ \ 
$ -1$,
$ 1$;\ \ 
$ -1$)

Factors = $2_{\frac{26}{5},3.618}^{5,720}\boxtimes 4_{4,4.}^{2,250}$

  \vskip 2ex

\noindent17. $8_{\frac{14}{5},14.47}^{20,755}$ \irep{207}:\ \ 
$d_i$ = ($1.0$,
$1.0$,
$1.0$,
$1.0$,
$1.618$,
$1.618$,
$1.618$,
$1.618$) 

\vskip 0.7ex
\hangindent=3em \hangafter=1
$D^2= 14.472 = 
10+2\sqrt{5}$

\vskip 0.7ex
\hangindent=3em \hangafter=1
$T = ( 0,
0,
\frac{1}{4},
\frac{3}{4},
\frac{2}{5},
\frac{2}{5},
\frac{3}{20},
\frac{13}{20} )
$,

\vskip 0.7ex
\hangindent=3em \hangafter=1
$S$ = ($ 1$,
$ 1$,
$ 1$,
$ 1$,
$ \frac{1+\sqrt{5}}{2}$,
$ \frac{1+\sqrt{5}}{2}$,
$ \frac{1+\sqrt{5}}{2}$,
$ \frac{1+\sqrt{5}}{2}$;\ \ 
$ 1$,
$ -1$,
$ -1$,
$ \frac{1+\sqrt{5}}{2}$,
$ \frac{1+\sqrt{5}}{2}$,
$ -\frac{1+\sqrt{5}}{2}$,
$ -\frac{1+\sqrt{5}}{2}$;\ \ 
$ -1$,
$ 1$,
$ \frac{1+\sqrt{5}}{2}$,
$ -\frac{1+\sqrt{5}}{2}$,
$ \frac{1+\sqrt{5}}{2}$,
$ -\frac{1+\sqrt{5}}{2}$;\ \ 
$ -1$,
$ \frac{1+\sqrt{5}}{2}$,
$ -\frac{1+\sqrt{5}}{2}$,
$ -\frac{1+\sqrt{5}}{2}$,
$ \frac{1+\sqrt{5}}{2}$;\ \ 
$ -1$,
$ -1$,
$ -1$,
$ -1$;\ \ 
$ -1$,
$ 1$,
$ 1$;\ \ 
$ 1$,
$ -1$;\ \ 
$ 1$)

Factors = $2_{1,2.}^{4,437}\boxtimes 4_{\frac{9}{5},7.236}^{20,451}$

  \vskip 2ex

\noindent18. $8_{\frac{26}{5},14.47}^{20,539}$ \irep{207}:\ \ 
$d_i$ = ($1.0$,
$1.0$,
$1.0$,
$1.0$,
$1.618$,
$1.618$,
$1.618$,
$1.618$) 

\vskip 0.7ex
\hangindent=3em \hangafter=1
$D^2= 14.472 = 
10+2\sqrt{5}$

\vskip 0.7ex
\hangindent=3em \hangafter=1
$T = ( 0,
0,
\frac{1}{4},
\frac{3}{4},
\frac{3}{5},
\frac{3}{5},
\frac{7}{20},
\frac{17}{20} )
$,

\vskip 0.7ex
\hangindent=3em \hangafter=1
$S$ = ($ 1$,
$ 1$,
$ 1$,
$ 1$,
$ \frac{1+\sqrt{5}}{2}$,
$ \frac{1+\sqrt{5}}{2}$,
$ \frac{1+\sqrt{5}}{2}$,
$ \frac{1+\sqrt{5}}{2}$;\ \ 
$ 1$,
$ -1$,
$ -1$,
$ \frac{1+\sqrt{5}}{2}$,
$ \frac{1+\sqrt{5}}{2}$,
$ -\frac{1+\sqrt{5}}{2}$,
$ -\frac{1+\sqrt{5}}{2}$;\ \ 
$ -1$,
$ 1$,
$ \frac{1+\sqrt{5}}{2}$,
$ -\frac{1+\sqrt{5}}{2}$,
$ \frac{1+\sqrt{5}}{2}$,
$ -\frac{1+\sqrt{5}}{2}$;\ \ 
$ -1$,
$ \frac{1+\sqrt{5}}{2}$,
$ -\frac{1+\sqrt{5}}{2}$,
$ -\frac{1+\sqrt{5}}{2}$,
$ \frac{1+\sqrt{5}}{2}$;\ \ 
$ -1$,
$ -1$,
$ -1$,
$ -1$;\ \ 
$ -1$,
$ 1$,
$ 1$;\ \ 
$ 1$,
$ -1$;\ \ 
$ 1$)

Factors = $2_{1,2.}^{4,437}\boxtimes 4_{\frac{21}{5},7.236}^{20,341}$

  \vskip 2ex

\noindent19. $8_{\frac{24}{5},14.47}^{20,693}$ \irep{207}:\ \ 
$d_i$ = ($1.0$,
$1.0$,
$1.0$,
$1.0$,
$1.618$,
$1.618$,
$1.618$,
$1.618$) 

\vskip 0.7ex
\hangindent=3em \hangafter=1
$D^2= 14.472 = 
10+2\sqrt{5}$

\vskip 0.7ex
\hangindent=3em \hangafter=1
$T = ( 0,
\frac{1}{2},
\frac{1}{4},
\frac{1}{4},
\frac{2}{5},
\frac{9}{10},
\frac{13}{20},
\frac{13}{20} )
$,

\vskip 0.7ex
\hangindent=3em \hangafter=1
$S$ = ($ 1$,
$ 1$,
$ 1$,
$ 1$,
$ \frac{1+\sqrt{5}}{2}$,
$ \frac{1+\sqrt{5}}{2}$,
$ \frac{1+\sqrt{5}}{2}$,
$ \frac{1+\sqrt{5}}{2}$;\ \ 
$ 1$,
$ -1$,
$ -1$,
$ \frac{1+\sqrt{5}}{2}$,
$ \frac{1+\sqrt{5}}{2}$,
$ -\frac{1+\sqrt{5}}{2}$,
$ -\frac{1+\sqrt{5}}{2}$;\ \ 
$ -1$,
$ 1$,
$ \frac{1+\sqrt{5}}{2}$,
$ -\frac{1+\sqrt{5}}{2}$,
$ \frac{1+\sqrt{5}}{2}$,
$ -\frac{1+\sqrt{5}}{2}$;\ \ 
$ -1$,
$ \frac{1+\sqrt{5}}{2}$,
$ -\frac{1+\sqrt{5}}{2}$,
$ -\frac{1+\sqrt{5}}{2}$,
$ \frac{1+\sqrt{5}}{2}$;\ \ 
$ -1$,
$ -1$,
$ -1$,
$ -1$;\ \ 
$ -1$,
$ 1$,
$ 1$;\ \ 
$ 1$,
$ -1$;\ \ 
$ 1$)

Factors = $2_{1,2.}^{4,437}\boxtimes 4_{\frac{19}{5},7.236}^{20,304}$

  \vskip 2ex

\noindent20. $8_{\frac{36}{5},14.47}^{20,693}$ \irep{207}:\ \ 
$d_i$ = ($1.0$,
$1.0$,
$1.0$,
$1.0$,
$1.618$,
$1.618$,
$1.618$,
$1.618$) 

\vskip 0.7ex
\hangindent=3em \hangafter=1
$D^2= 14.472 = 
10+2\sqrt{5}$

\vskip 0.7ex
\hangindent=3em \hangafter=1
$T = ( 0,
\frac{1}{2},
\frac{1}{4},
\frac{1}{4},
\frac{3}{5},
\frac{1}{10},
\frac{17}{20},
\frac{17}{20} )
$,

\vskip 0.7ex
\hangindent=3em \hangafter=1
$S$ = ($ 1$,
$ 1$,
$ 1$,
$ 1$,
$ \frac{1+\sqrt{5}}{2}$,
$ \frac{1+\sqrt{5}}{2}$,
$ \frac{1+\sqrt{5}}{2}$,
$ \frac{1+\sqrt{5}}{2}$;\ \ 
$ 1$,
$ -1$,
$ -1$,
$ \frac{1+\sqrt{5}}{2}$,
$ \frac{1+\sqrt{5}}{2}$,
$ -\frac{1+\sqrt{5}}{2}$,
$ -\frac{1+\sqrt{5}}{2}$;\ \ 
$ -1$,
$ 1$,
$ \frac{1+\sqrt{5}}{2}$,
$ -\frac{1+\sqrt{5}}{2}$,
$ \frac{1+\sqrt{5}}{2}$,
$ -\frac{1+\sqrt{5}}{2}$;\ \ 
$ -1$,
$ \frac{1+\sqrt{5}}{2}$,
$ -\frac{1+\sqrt{5}}{2}$,
$ -\frac{1+\sqrt{5}}{2}$,
$ \frac{1+\sqrt{5}}{2}$;\ \ 
$ -1$,
$ -1$,
$ -1$,
$ -1$;\ \ 
$ -1$,
$ 1$,
$ 1$;\ \ 
$ 1$,
$ -1$;\ \ 
$ 1$)

Factors = $2_{1,2.}^{4,437}\boxtimes 4_{\frac{31}{5},7.236}^{20,505}$

  \vskip 2ex

\noindent21. $8_{\frac{4}{5},14.47}^{20,399}$ \irep{207}:\ \ 
$d_i$ = ($1.0$,
$1.0$,
$1.0$,
$1.0$,
$1.618$,
$1.618$,
$1.618$,
$1.618$) 

\vskip 0.7ex
\hangindent=3em \hangafter=1
$D^2= 14.472 = 
10+2\sqrt{5}$

\vskip 0.7ex
\hangindent=3em \hangafter=1
$T = ( 0,
\frac{1}{2},
\frac{3}{4},
\frac{3}{4},
\frac{2}{5},
\frac{9}{10},
\frac{3}{20},
\frac{3}{20} )
$,

\vskip 0.7ex
\hangindent=3em \hangafter=1
$S$ = ($ 1$,
$ 1$,
$ 1$,
$ 1$,
$ \frac{1+\sqrt{5}}{2}$,
$ \frac{1+\sqrt{5}}{2}$,
$ \frac{1+\sqrt{5}}{2}$,
$ \frac{1+\sqrt{5}}{2}$;\ \ 
$ 1$,
$ -1$,
$ -1$,
$ \frac{1+\sqrt{5}}{2}$,
$ \frac{1+\sqrt{5}}{2}$,
$ -\frac{1+\sqrt{5}}{2}$,
$ -\frac{1+\sqrt{5}}{2}$;\ \ 
$ -1$,
$ 1$,
$ \frac{1+\sqrt{5}}{2}$,
$ -\frac{1+\sqrt{5}}{2}$,
$ \frac{1+\sqrt{5}}{2}$,
$ -\frac{1+\sqrt{5}}{2}$;\ \ 
$ -1$,
$ \frac{1+\sqrt{5}}{2}$,
$ -\frac{1+\sqrt{5}}{2}$,
$ -\frac{1+\sqrt{5}}{2}$,
$ \frac{1+\sqrt{5}}{2}$;\ \ 
$ -1$,
$ -1$,
$ -1$,
$ -1$;\ \ 
$ -1$,
$ 1$,
$ 1$;\ \ 
$ 1$,
$ -1$;\ \ 
$ 1$)

Factors = $2_{7,2.}^{4,625}\boxtimes 4_{\frac{9}{5},7.236}^{20,451}$

  \vskip 2ex

\noindent22. $8_{\frac{16}{5},14.47}^{20,247}$ \irep{207}:\ \ 
$d_i$ = ($1.0$,
$1.0$,
$1.0$,
$1.0$,
$1.618$,
$1.618$,
$1.618$,
$1.618$) 

\vskip 0.7ex
\hangindent=3em \hangafter=1
$D^2= 14.472 = 
10+2\sqrt{5}$

\vskip 0.7ex
\hangindent=3em \hangafter=1
$T = ( 0,
\frac{1}{2},
\frac{3}{4},
\frac{3}{4},
\frac{3}{5},
\frac{1}{10},
\frac{7}{20},
\frac{7}{20} )
$,

\vskip 0.7ex
\hangindent=3em \hangafter=1
$S$ = ($ 1$,
$ 1$,
$ 1$,
$ 1$,
$ \frac{1+\sqrt{5}}{2}$,
$ \frac{1+\sqrt{5}}{2}$,
$ \frac{1+\sqrt{5}}{2}$,
$ \frac{1+\sqrt{5}}{2}$;\ \ 
$ 1$,
$ -1$,
$ -1$,
$ \frac{1+\sqrt{5}}{2}$,
$ \frac{1+\sqrt{5}}{2}$,
$ -\frac{1+\sqrt{5}}{2}$,
$ -\frac{1+\sqrt{5}}{2}$;\ \ 
$ -1$,
$ 1$,
$ \frac{1+\sqrt{5}}{2}$,
$ -\frac{1+\sqrt{5}}{2}$,
$ \frac{1+\sqrt{5}}{2}$,
$ -\frac{1+\sqrt{5}}{2}$;\ \ 
$ -1$,
$ \frac{1+\sqrt{5}}{2}$,
$ -\frac{1+\sqrt{5}}{2}$,
$ -\frac{1+\sqrt{5}}{2}$,
$ \frac{1+\sqrt{5}}{2}$;\ \ 
$ -1$,
$ -1$,
$ -1$,
$ -1$;\ \ 
$ -1$,
$ 1$,
$ 1$;\ \ 
$ 1$,
$ -1$;\ \ 
$ 1$)

Factors = $2_{7,2.}^{4,625}\boxtimes 4_{\frac{21}{5},7.236}^{20,341}$

  \vskip 2ex

\noindent23. $8_{\frac{19}{5},14.47}^{40,124}$ \irep{233}:\ \ 
$d_i$ = ($1.0$,
$1.0$,
$1.0$,
$1.0$,
$1.618$,
$1.618$,
$1.618$,
$1.618$) 

\vskip 0.7ex
\hangindent=3em \hangafter=1
$D^2= 14.472 = 
10+2\sqrt{5}$

\vskip 0.7ex
\hangindent=3em \hangafter=1
$T = ( 0,
\frac{1}{2},
\frac{1}{8},
\frac{1}{8},
\frac{2}{5},
\frac{9}{10},
\frac{21}{40},
\frac{21}{40} )
$,

\vskip 0.7ex
\hangindent=3em \hangafter=1
$S$ = ($ 1$,
$ 1$,
$ 1$,
$ 1$,
$ \frac{1+\sqrt{5}}{2}$,
$ \frac{1+\sqrt{5}}{2}$,
$ \frac{1+\sqrt{5}}{2}$,
$ \frac{1+\sqrt{5}}{2}$;\ \ 
$ 1$,
$ -1$,
$ -1$,
$ \frac{1+\sqrt{5}}{2}$,
$ \frac{1+\sqrt{5}}{2}$,
$ -\frac{1+\sqrt{5}}{2}$,
$ -\frac{1+\sqrt{5}}{2}$;\ \ 
$-\mathrm{i}$,
$\mathrm{i}$,
$ \frac{1+\sqrt{5}}{2}$,
$ -\frac{1+\sqrt{5}}{2}$,
$(\frac{1+\sqrt{5}}{2})\mathrm{i}$,
$(-\frac{1+\sqrt{5}}{2})\mathrm{i}$;\ \ 
$-\mathrm{i}$,
$ \frac{1+\sqrt{5}}{2}$,
$ -\frac{1+\sqrt{5}}{2}$,
$(-\frac{1+\sqrt{5}}{2})\mathrm{i}$,
$(\frac{1+\sqrt{5}}{2})\mathrm{i}$;\ \ 
$ -1$,
$ -1$,
$ -1$,
$ -1$;\ \ 
$ -1$,
$ 1$,
$ 1$;\ \ 
$\mathrm{i}$,
$-\mathrm{i}$;\ \ 
$\mathrm{i}$)

Factors = $2_{\frac{14}{5},3.618}^{5,395}\boxtimes 4_{1,4.}^{8,718}$

  \vskip 2ex

\noindent24. $8_{\frac{31}{5},14.47}^{40,371}$ \irep{233}:\ \ 
$d_i$ = ($1.0$,
$1.0$,
$1.0$,
$1.0$,
$1.618$,
$1.618$,
$1.618$,
$1.618$) 

\vskip 0.7ex
\hangindent=3em \hangafter=1
$D^2= 14.472 = 
10+2\sqrt{5}$

\vskip 0.7ex
\hangindent=3em \hangafter=1
$T = ( 0,
\frac{1}{2},
\frac{1}{8},
\frac{1}{8},
\frac{3}{5},
\frac{1}{10},
\frac{29}{40},
\frac{29}{40} )
$,

\vskip 0.7ex
\hangindent=3em \hangafter=1
$S$ = ($ 1$,
$ 1$,
$ 1$,
$ 1$,
$ \frac{1+\sqrt{5}}{2}$,
$ \frac{1+\sqrt{5}}{2}$,
$ \frac{1+\sqrt{5}}{2}$,
$ \frac{1+\sqrt{5}}{2}$;\ \ 
$ 1$,
$ -1$,
$ -1$,
$ \frac{1+\sqrt{5}}{2}$,
$ \frac{1+\sqrt{5}}{2}$,
$ -\frac{1+\sqrt{5}}{2}$,
$ -\frac{1+\sqrt{5}}{2}$;\ \ 
$-\mathrm{i}$,
$\mathrm{i}$,
$ \frac{1+\sqrt{5}}{2}$,
$ -\frac{1+\sqrt{5}}{2}$,
$(\frac{1+\sqrt{5}}{2})\mathrm{i}$,
$(-\frac{1+\sqrt{5}}{2})\mathrm{i}$;\ \ 
$-\mathrm{i}$,
$ \frac{1+\sqrt{5}}{2}$,
$ -\frac{1+\sqrt{5}}{2}$,
$(-\frac{1+\sqrt{5}}{2})\mathrm{i}$,
$(\frac{1+\sqrt{5}}{2})\mathrm{i}$;\ \ 
$ -1$,
$ -1$,
$ -1$,
$ -1$;\ \ 
$ -1$,
$ 1$,
$ 1$;\ \ 
$\mathrm{i}$,
$-\mathrm{i}$;\ \ 
$\mathrm{i}$)

Factors = $2_{\frac{26}{5},3.618}^{5,720}\boxtimes 4_{1,4.}^{8,718}$

  \vskip 2ex

\noindent25. $8_{\frac{29}{5},14.47}^{40,142}$ \irep{233}:\ \ 
$d_i$ = ($1.0$,
$1.0$,
$1.0$,
$1.0$,
$1.618$,
$1.618$,
$1.618$,
$1.618$) 

\vskip 0.7ex
\hangindent=3em \hangafter=1
$D^2= 14.472 = 
10+2\sqrt{5}$

\vskip 0.7ex
\hangindent=3em \hangafter=1
$T = ( 0,
\frac{1}{2},
\frac{3}{8},
\frac{3}{8},
\frac{2}{5},
\frac{9}{10},
\frac{31}{40},
\frac{31}{40} )
$,

\vskip 0.7ex
\hangindent=3em \hangafter=1
$S$ = ($ 1$,
$ 1$,
$ 1$,
$ 1$,
$ \frac{1+\sqrt{5}}{2}$,
$ \frac{1+\sqrt{5}}{2}$,
$ \frac{1+\sqrt{5}}{2}$,
$ \frac{1+\sqrt{5}}{2}$;\ \ 
$ 1$,
$ -1$,
$ -1$,
$ \frac{1+\sqrt{5}}{2}$,
$ \frac{1+\sqrt{5}}{2}$,
$ -\frac{1+\sqrt{5}}{2}$,
$ -\frac{1+\sqrt{5}}{2}$;\ \ 
$\mathrm{i}$,
$-\mathrm{i}$,
$ \frac{1+\sqrt{5}}{2}$,
$ -\frac{1+\sqrt{5}}{2}$,
$(\frac{1+\sqrt{5}}{2})\mathrm{i}$,
$(-\frac{1+\sqrt{5}}{2})\mathrm{i}$;\ \ 
$\mathrm{i}$,
$ \frac{1+\sqrt{5}}{2}$,
$ -\frac{1+\sqrt{5}}{2}$,
$(-\frac{1+\sqrt{5}}{2})\mathrm{i}$,
$(\frac{1+\sqrt{5}}{2})\mathrm{i}$;\ \ 
$ -1$,
$ -1$,
$ -1$,
$ -1$;\ \ 
$ -1$,
$ 1$,
$ 1$;\ \ 
$-\mathrm{i}$,
$\mathrm{i}$;\ \ 
$-\mathrm{i}$)

Factors = $2_{\frac{14}{5},3.618}^{5,395}\boxtimes 4_{3,4.}^{8,468}$

  \vskip 2ex

\noindent26. $8_{\frac{1}{5},14.47}^{40,158}$ \irep{233}:\ \ 
$d_i$ = ($1.0$,
$1.0$,
$1.0$,
$1.0$,
$1.618$,
$1.618$,
$1.618$,
$1.618$) 

\vskip 0.7ex
\hangindent=3em \hangafter=1
$D^2= 14.472 = 
10+2\sqrt{5}$

\vskip 0.7ex
\hangindent=3em \hangafter=1
$T = ( 0,
\frac{1}{2},
\frac{3}{8},
\frac{3}{8},
\frac{3}{5},
\frac{1}{10},
\frac{39}{40},
\frac{39}{40} )
$,

\vskip 0.7ex
\hangindent=3em \hangafter=1
$S$ = ($ 1$,
$ 1$,
$ 1$,
$ 1$,
$ \frac{1+\sqrt{5}}{2}$,
$ \frac{1+\sqrt{5}}{2}$,
$ \frac{1+\sqrt{5}}{2}$,
$ \frac{1+\sqrt{5}}{2}$;\ \ 
$ 1$,
$ -1$,
$ -1$,
$ \frac{1+\sqrt{5}}{2}$,
$ \frac{1+\sqrt{5}}{2}$,
$ -\frac{1+\sqrt{5}}{2}$,
$ -\frac{1+\sqrt{5}}{2}$;\ \ 
$\mathrm{i}$,
$-\mathrm{i}$,
$ \frac{1+\sqrt{5}}{2}$,
$ -\frac{1+\sqrt{5}}{2}$,
$(\frac{1+\sqrt{5}}{2})\mathrm{i}$,
$(-\frac{1+\sqrt{5}}{2})\mathrm{i}$;\ \ 
$\mathrm{i}$,
$ \frac{1+\sqrt{5}}{2}$,
$ -\frac{1+\sqrt{5}}{2}$,
$(-\frac{1+\sqrt{5}}{2})\mathrm{i}$,
$(\frac{1+\sqrt{5}}{2})\mathrm{i}$;\ \ 
$ -1$,
$ -1$,
$ -1$,
$ -1$;\ \ 
$ -1$,
$ 1$,
$ 1$;\ \ 
$-\mathrm{i}$,
$\mathrm{i}$;\ \ 
$-\mathrm{i}$)

Factors = $2_{\frac{26}{5},3.618}^{5,720}\boxtimes 4_{3,4.}^{8,468}$

  \vskip 2ex

\noindent27. $8_{\frac{39}{5},14.47}^{40,152}$ \irep{233}:\ \ 
$d_i$ = ($1.0$,
$1.0$,
$1.0$,
$1.0$,
$1.618$,
$1.618$,
$1.618$,
$1.618$) 

\vskip 0.7ex
\hangindent=3em \hangafter=1
$D^2= 14.472 = 
10+2\sqrt{5}$

\vskip 0.7ex
\hangindent=3em \hangafter=1
$T = ( 0,
\frac{1}{2},
\frac{5}{8},
\frac{5}{8},
\frac{2}{5},
\frac{9}{10},
\frac{1}{40},
\frac{1}{40} )
$,

\vskip 0.7ex
\hangindent=3em \hangafter=1
$S$ = ($ 1$,
$ 1$,
$ 1$,
$ 1$,
$ \frac{1+\sqrt{5}}{2}$,
$ \frac{1+\sqrt{5}}{2}$,
$ \frac{1+\sqrt{5}}{2}$,
$ \frac{1+\sqrt{5}}{2}$;\ \ 
$ 1$,
$ -1$,
$ -1$,
$ \frac{1+\sqrt{5}}{2}$,
$ \frac{1+\sqrt{5}}{2}$,
$ -\frac{1+\sqrt{5}}{2}$,
$ -\frac{1+\sqrt{5}}{2}$;\ \ 
$-\mathrm{i}$,
$\mathrm{i}$,
$ \frac{1+\sqrt{5}}{2}$,
$ -\frac{1+\sqrt{5}}{2}$,
$(\frac{1+\sqrt{5}}{2})\mathrm{i}$,
$(-\frac{1+\sqrt{5}}{2})\mathrm{i}$;\ \ 
$-\mathrm{i}$,
$ \frac{1+\sqrt{5}}{2}$,
$ -\frac{1+\sqrt{5}}{2}$,
$(-\frac{1+\sqrt{5}}{2})\mathrm{i}$,
$(\frac{1+\sqrt{5}}{2})\mathrm{i}$;\ \ 
$ -1$,
$ -1$,
$ -1$,
$ -1$;\ \ 
$ -1$,
$ 1$,
$ 1$;\ \ 
$\mathrm{i}$,
$-\mathrm{i}$;\ \ 
$\mathrm{i}$)

Factors = $2_{\frac{14}{5},3.618}^{5,395}\boxtimes 4_{5,4.}^{8,312}$

  \vskip 2ex

\noindent28. $8_{\frac{11}{5},14.47}^{40,824}$ \irep{233}:\ \ 
$d_i$ = ($1.0$,
$1.0$,
$1.0$,
$1.0$,
$1.618$,
$1.618$,
$1.618$,
$1.618$) 

\vskip 0.7ex
\hangindent=3em \hangafter=1
$D^2= 14.472 = 
10+2\sqrt{5}$

\vskip 0.7ex
\hangindent=3em \hangafter=1
$T = ( 0,
\frac{1}{2},
\frac{5}{8},
\frac{5}{8},
\frac{3}{5},
\frac{1}{10},
\frac{9}{40},
\frac{9}{40} )
$,

\vskip 0.7ex
\hangindent=3em \hangafter=1
$S$ = ($ 1$,
$ 1$,
$ 1$,
$ 1$,
$ \frac{1+\sqrt{5}}{2}$,
$ \frac{1+\sqrt{5}}{2}$,
$ \frac{1+\sqrt{5}}{2}$,
$ \frac{1+\sqrt{5}}{2}$;\ \ 
$ 1$,
$ -1$,
$ -1$,
$ \frac{1+\sqrt{5}}{2}$,
$ \frac{1+\sqrt{5}}{2}$,
$ -\frac{1+\sqrt{5}}{2}$,
$ -\frac{1+\sqrt{5}}{2}$;\ \ 
$-\mathrm{i}$,
$\mathrm{i}$,
$ \frac{1+\sqrt{5}}{2}$,
$ -\frac{1+\sqrt{5}}{2}$,
$(\frac{1+\sqrt{5}}{2})\mathrm{i}$,
$(-\frac{1+\sqrt{5}}{2})\mathrm{i}$;\ \ 
$-\mathrm{i}$,
$ \frac{1+\sqrt{5}}{2}$,
$ -\frac{1+\sqrt{5}}{2}$,
$(-\frac{1+\sqrt{5}}{2})\mathrm{i}$,
$(\frac{1+\sqrt{5}}{2})\mathrm{i}$;\ \ 
$ -1$,
$ -1$,
$ -1$,
$ -1$;\ \ 
$ -1$,
$ 1$,
$ 1$;\ \ 
$\mathrm{i}$,
$-\mathrm{i}$;\ \ 
$\mathrm{i}$)

Factors = $2_{\frac{26}{5},3.618}^{5,720}\boxtimes 4_{5,4.}^{8,312}$

  \vskip 2ex

\noindent29. $8_{\frac{9}{5},14.47}^{40,977}$ \irep{233}:\ \ 
$d_i$ = ($1.0$,
$1.0$,
$1.0$,
$1.0$,
$1.618$,
$1.618$,
$1.618$,
$1.618$) 

\vskip 0.7ex
\hangindent=3em \hangafter=1
$D^2= 14.472 = 
10+2\sqrt{5}$

\vskip 0.7ex
\hangindent=3em \hangafter=1
$T = ( 0,
\frac{1}{2},
\frac{7}{8},
\frac{7}{8},
\frac{2}{5},
\frac{9}{10},
\frac{11}{40},
\frac{11}{40} )
$,

\vskip 0.7ex
\hangindent=3em \hangafter=1
$S$ = ($ 1$,
$ 1$,
$ 1$,
$ 1$,
$ \frac{1+\sqrt{5}}{2}$,
$ \frac{1+\sqrt{5}}{2}$,
$ \frac{1+\sqrt{5}}{2}$,
$ \frac{1+\sqrt{5}}{2}$;\ \ 
$ 1$,
$ -1$,
$ -1$,
$ \frac{1+\sqrt{5}}{2}$,
$ \frac{1+\sqrt{5}}{2}$,
$ -\frac{1+\sqrt{5}}{2}$,
$ -\frac{1+\sqrt{5}}{2}$;\ \ 
$\mathrm{i}$,
$-\mathrm{i}$,
$ \frac{1+\sqrt{5}}{2}$,
$ -\frac{1+\sqrt{5}}{2}$,
$(\frac{1+\sqrt{5}}{2})\mathrm{i}$,
$(-\frac{1+\sqrt{5}}{2})\mathrm{i}$;\ \ 
$\mathrm{i}$,
$ \frac{1+\sqrt{5}}{2}$,
$ -\frac{1+\sqrt{5}}{2}$,
$(-\frac{1+\sqrt{5}}{2})\mathrm{i}$,
$(\frac{1+\sqrt{5}}{2})\mathrm{i}$;\ \ 
$ -1$,
$ -1$,
$ -1$,
$ -1$;\ \ 
$ -1$,
$ 1$,
$ 1$;\ \ 
$-\mathrm{i}$,
$\mathrm{i}$;\ \ 
$-\mathrm{i}$)

Factors = $2_{\frac{14}{5},3.618}^{5,395}\boxtimes 4_{7,4.}^{8,781}$

  \vskip 2ex

\noindent30. $8_{\frac{21}{5},14.47}^{40,492}$ \irep{233}:\ \ 
$d_i$ = ($1.0$,
$1.0$,
$1.0$,
$1.0$,
$1.618$,
$1.618$,
$1.618$,
$1.618$) 

\vskip 0.7ex
\hangindent=3em \hangafter=1
$D^2= 14.472 = 
10+2\sqrt{5}$

\vskip 0.7ex
\hangindent=3em \hangafter=1
$T = ( 0,
\frac{1}{2},
\frac{7}{8},
\frac{7}{8},
\frac{3}{5},
\frac{1}{10},
\frac{19}{40},
\frac{19}{40} )
$,

\vskip 0.7ex
\hangindent=3em \hangafter=1
$S$ = ($ 1$,
$ 1$,
$ 1$,
$ 1$,
$ \frac{1+\sqrt{5}}{2}$,
$ \frac{1+\sqrt{5}}{2}$,
$ \frac{1+\sqrt{5}}{2}$,
$ \frac{1+\sqrt{5}}{2}$;\ \ 
$ 1$,
$ -1$,
$ -1$,
$ \frac{1+\sqrt{5}}{2}$,
$ \frac{1+\sqrt{5}}{2}$,
$ -\frac{1+\sqrt{5}}{2}$,
$ -\frac{1+\sqrt{5}}{2}$;\ \ 
$\mathrm{i}$,
$-\mathrm{i}$,
$ \frac{1+\sqrt{5}}{2}$,
$ -\frac{1+\sqrt{5}}{2}$,
$(\frac{1+\sqrt{5}}{2})\mathrm{i}$,
$(-\frac{1+\sqrt{5}}{2})\mathrm{i}$;\ \ 
$\mathrm{i}$,
$ \frac{1+\sqrt{5}}{2}$,
$ -\frac{1+\sqrt{5}}{2}$,
$(-\frac{1+\sqrt{5}}{2})\mathrm{i}$,
$(\frac{1+\sqrt{5}}{2})\mathrm{i}$;\ \ 
$ -1$,
$ -1$,
$ -1$,
$ -1$;\ \ 
$ -1$,
$ 1$,
$ 1$;\ \ 
$-\mathrm{i}$,
$\mathrm{i}$;\ \ 
$-\mathrm{i}$)

Factors = $2_{\frac{26}{5},3.618}^{5,720}\boxtimes 4_{7,4.}^{8,781}$

  \vskip 2ex

\noindent31. $8_{\frac{33}{5},26.18}^{20,202}$ \irep{208}:\ \ 
$d_i$ = ($1.0$,
$1.0$,
$1.618$,
$1.618$,
$1.618$,
$1.618$,
$2.618$,
$2.618$) 

\vskip 0.7ex
\hangindent=3em \hangafter=1
$D^2= 26.180 = 
15+5\sqrt{5}$

\vskip 0.7ex
\hangindent=3em \hangafter=1
$T = ( 0,
\frac{1}{4},
\frac{2}{5},
\frac{2}{5},
\frac{13}{20},
\frac{13}{20},
\frac{4}{5},
\frac{1}{20} )
$,

\vskip 0.7ex
\hangindent=3em \hangafter=1
$S$ = ($ 1$,
$ 1$,
$ \frac{1+\sqrt{5}}{2}$,
$ \frac{1+\sqrt{5}}{2}$,
$ \frac{1+\sqrt{5}}{2}$,
$ \frac{1+\sqrt{5}}{2}$,
$ \frac{3+\sqrt{5}}{2}$,
$ \frac{3+\sqrt{5}}{2}$;\ \ 
$ -1$,
$ \frac{1+\sqrt{5}}{2}$,
$ \frac{1+\sqrt{5}}{2}$,
$ -\frac{1+\sqrt{5}}{2}$,
$ -\frac{1+\sqrt{5}}{2}$,
$ \frac{3+\sqrt{5}}{2}$,
$ -\frac{3+\sqrt{5}}{2}$;\ \ 
$ -1$,
$ \frac{3+\sqrt{5}}{2}$,
$ -1$,
$ \frac{3+\sqrt{5}}{2}$,
$ -\frac{1+\sqrt{5}}{2}$,
$ -\frac{1+\sqrt{5}}{2}$;\ \ 
$ -1$,
$ \frac{3+\sqrt{5}}{2}$,
$ -1$,
$ -\frac{1+\sqrt{5}}{2}$,
$ -\frac{1+\sqrt{5}}{2}$;\ \ 
$ 1$,
$ -\frac{3+\sqrt{5}}{2}$,
$ -\frac{1+\sqrt{5}}{2}$,
$ \frac{1+\sqrt{5}}{2}$;\ \ 
$ 1$,
$ -\frac{1+\sqrt{5}}{2}$,
$ \frac{1+\sqrt{5}}{2}$;\ \ 
$ 1$,
$ 1$;\ \ 
$ -1$)

Factors = $2_{1,2.}^{4,437}\boxtimes 4_{\frac{28}{5},13.09}^{5,479}$

  \vskip 2ex

\noindent32. $8_{1,26.18}^{20,506}$ \irep{208}:\ \ 
$d_i$ = ($1.0$,
$1.0$,
$1.618$,
$1.618$,
$1.618$,
$1.618$,
$2.618$,
$2.618$) 

\vskip 0.7ex
\hangindent=3em \hangafter=1
$D^2= 26.180 = 
15+5\sqrt{5}$

\vskip 0.7ex
\hangindent=3em \hangafter=1
$T = ( 0,
\frac{1}{4},
\frac{2}{5},
\frac{3}{5},
\frac{13}{20},
\frac{17}{20},
0,
\frac{1}{4} )
$,

\vskip 0.7ex
\hangindent=3em \hangafter=1
$S$ = ($ 1$,
$ 1$,
$ \frac{1+\sqrt{5}}{2}$,
$ \frac{1+\sqrt{5}}{2}$,
$ \frac{1+\sqrt{5}}{2}$,
$ \frac{1+\sqrt{5}}{2}$,
$ \frac{3+\sqrt{5}}{2}$,
$ \frac{3+\sqrt{5}}{2}$;\ \ 
$ -1$,
$ \frac{1+\sqrt{5}}{2}$,
$ \frac{1+\sqrt{5}}{2}$,
$ -\frac{1+\sqrt{5}}{2}$,
$ -\frac{1+\sqrt{5}}{2}$,
$ \frac{3+\sqrt{5}}{2}$,
$ -\frac{3+\sqrt{5}}{2}$;\ \ 
$ -1$,
$ \frac{3+\sqrt{5}}{2}$,
$ -1$,
$ \frac{3+\sqrt{5}}{2}$,
$ -\frac{1+\sqrt{5}}{2}$,
$ -\frac{1+\sqrt{5}}{2}$;\ \ 
$ -1$,
$ \frac{3+\sqrt{5}}{2}$,
$ -1$,
$ -\frac{1+\sqrt{5}}{2}$,
$ -\frac{1+\sqrt{5}}{2}$;\ \ 
$ 1$,
$ -\frac{3+\sqrt{5}}{2}$,
$ -\frac{1+\sqrt{5}}{2}$,
$ \frac{1+\sqrt{5}}{2}$;\ \ 
$ 1$,
$ -\frac{1+\sqrt{5}}{2}$,
$ \frac{1+\sqrt{5}}{2}$;\ \ 
$ 1$,
$ 1$;\ \ 
$ -1$)

Factors = $2_{1,2.}^{4,437}\boxtimes 4_{0,13.09}^{5,872}$

  \vskip 2ex

\noindent33. $8_{\frac{17}{5},26.18}^{20,773}$ \irep{208}:\ \ 
$d_i$ = ($1.0$,
$1.0$,
$1.618$,
$1.618$,
$1.618$,
$1.618$,
$2.618$,
$2.618$) 

\vskip 0.7ex
\hangindent=3em \hangafter=1
$D^2= 26.180 = 
15+5\sqrt{5}$

\vskip 0.7ex
\hangindent=3em \hangafter=1
$T = ( 0,
\frac{1}{4},
\frac{3}{5},
\frac{3}{5},
\frac{17}{20},
\frac{17}{20},
\frac{1}{5},
\frac{9}{20} )
$,

\vskip 0.7ex
\hangindent=3em \hangafter=1
$S$ = ($ 1$,
$ 1$,
$ \frac{1+\sqrt{5}}{2}$,
$ \frac{1+\sqrt{5}}{2}$,
$ \frac{1+\sqrt{5}}{2}$,
$ \frac{1+\sqrt{5}}{2}$,
$ \frac{3+\sqrt{5}}{2}$,
$ \frac{3+\sqrt{5}}{2}$;\ \ 
$ -1$,
$ \frac{1+\sqrt{5}}{2}$,
$ \frac{1+\sqrt{5}}{2}$,
$ -\frac{1+\sqrt{5}}{2}$,
$ -\frac{1+\sqrt{5}}{2}$,
$ \frac{3+\sqrt{5}}{2}$,
$ -\frac{3+\sqrt{5}}{2}$;\ \ 
$ -1$,
$ \frac{3+\sqrt{5}}{2}$,
$ -1$,
$ \frac{3+\sqrt{5}}{2}$,
$ -\frac{1+\sqrt{5}}{2}$,
$ -\frac{1+\sqrt{5}}{2}$;\ \ 
$ -1$,
$ \frac{3+\sqrt{5}}{2}$,
$ -1$,
$ -\frac{1+\sqrt{5}}{2}$,
$ -\frac{1+\sqrt{5}}{2}$;\ \ 
$ 1$,
$ -\frac{3+\sqrt{5}}{2}$,
$ -\frac{1+\sqrt{5}}{2}$,
$ \frac{1+\sqrt{5}}{2}$;\ \ 
$ 1$,
$ -\frac{1+\sqrt{5}}{2}$,
$ \frac{1+\sqrt{5}}{2}$;\ \ 
$ 1$,
$ 1$;\ \ 
$ -1$)

Factors = $2_{1,2.}^{4,437}\boxtimes 4_{\frac{12}{5},13.09}^{5,443}$

  \vskip 2ex

\noindent34. $8_{\frac{23}{5},26.18}^{20,193}$ \irep{208}:\ \ 
$d_i$ = ($1.0$,
$1.0$,
$1.618$,
$1.618$,
$1.618$,
$1.618$,
$2.618$,
$2.618$) 

\vskip 0.7ex
\hangindent=3em \hangafter=1
$D^2= 26.180 = 
15+5\sqrt{5}$

\vskip 0.7ex
\hangindent=3em \hangafter=1
$T = ( 0,
\frac{3}{4},
\frac{2}{5},
\frac{2}{5},
\frac{3}{20},
\frac{3}{20},
\frac{4}{5},
\frac{11}{20} )
$,

\vskip 0.7ex
\hangindent=3em \hangafter=1
$S$ = ($ 1$,
$ 1$,
$ \frac{1+\sqrt{5}}{2}$,
$ \frac{1+\sqrt{5}}{2}$,
$ \frac{1+\sqrt{5}}{2}$,
$ \frac{1+\sqrt{5}}{2}$,
$ \frac{3+\sqrt{5}}{2}$,
$ \frac{3+\sqrt{5}}{2}$;\ \ 
$ -1$,
$ \frac{1+\sqrt{5}}{2}$,
$ \frac{1+\sqrt{5}}{2}$,
$ -\frac{1+\sqrt{5}}{2}$,
$ -\frac{1+\sqrt{5}}{2}$,
$ \frac{3+\sqrt{5}}{2}$,
$ -\frac{3+\sqrt{5}}{2}$;\ \ 
$ -1$,
$ \frac{3+\sqrt{5}}{2}$,
$ -1$,
$ \frac{3+\sqrt{5}}{2}$,
$ -\frac{1+\sqrt{5}}{2}$,
$ -\frac{1+\sqrt{5}}{2}$;\ \ 
$ -1$,
$ \frac{3+\sqrt{5}}{2}$,
$ -1$,
$ -\frac{1+\sqrt{5}}{2}$,
$ -\frac{1+\sqrt{5}}{2}$;\ \ 
$ 1$,
$ -\frac{3+\sqrt{5}}{2}$,
$ -\frac{1+\sqrt{5}}{2}$,
$ \frac{1+\sqrt{5}}{2}$;\ \ 
$ 1$,
$ -\frac{1+\sqrt{5}}{2}$,
$ \frac{1+\sqrt{5}}{2}$;\ \ 
$ 1$,
$ 1$;\ \ 
$ -1$)

Factors = $2_{7,2.}^{4,625}\boxtimes 4_{\frac{28}{5},13.09}^{5,479}$

  \vskip 2ex

\noindent35. $8_{7,26.18}^{20,315}$ \irep{208}:\ \ 
$d_i$ = ($1.0$,
$1.0$,
$1.618$,
$1.618$,
$1.618$,
$1.618$,
$2.618$,
$2.618$) 

\vskip 0.7ex
\hangindent=3em \hangafter=1
$D^2= 26.180 = 
15+5\sqrt{5}$

\vskip 0.7ex
\hangindent=3em \hangafter=1
$T = ( 0,
\frac{3}{4},
\frac{2}{5},
\frac{3}{5},
\frac{3}{20},
\frac{7}{20},
0,
\frac{3}{4} )
$,

\vskip 0.7ex
\hangindent=3em \hangafter=1
$S$ = ($ 1$,
$ 1$,
$ \frac{1+\sqrt{5}}{2}$,
$ \frac{1+\sqrt{5}}{2}$,
$ \frac{1+\sqrt{5}}{2}$,
$ \frac{1+\sqrt{5}}{2}$,
$ \frac{3+\sqrt{5}}{2}$,
$ \frac{3+\sqrt{5}}{2}$;\ \ 
$ -1$,
$ \frac{1+\sqrt{5}}{2}$,
$ \frac{1+\sqrt{5}}{2}$,
$ -\frac{1+\sqrt{5}}{2}$,
$ -\frac{1+\sqrt{5}}{2}$,
$ \frac{3+\sqrt{5}}{2}$,
$ -\frac{3+\sqrt{5}}{2}$;\ \ 
$ -1$,
$ \frac{3+\sqrt{5}}{2}$,
$ -1$,
$ \frac{3+\sqrt{5}}{2}$,
$ -\frac{1+\sqrt{5}}{2}$,
$ -\frac{1+\sqrt{5}}{2}$;\ \ 
$ -1$,
$ \frac{3+\sqrt{5}}{2}$,
$ -1$,
$ -\frac{1+\sqrt{5}}{2}$,
$ -\frac{1+\sqrt{5}}{2}$;\ \ 
$ 1$,
$ -\frac{3+\sqrt{5}}{2}$,
$ -\frac{1+\sqrt{5}}{2}$,
$ \frac{1+\sqrt{5}}{2}$;\ \ 
$ 1$,
$ -\frac{1+\sqrt{5}}{2}$,
$ \frac{1+\sqrt{5}}{2}$;\ \ 
$ 1$,
$ 1$;\ \ 
$ -1$)

Factors = $2_{7,2.}^{4,625}\boxtimes 4_{0,13.09}^{5,872}$

  \vskip 2ex

\noindent36. $8_{\frac{7}{5},26.18}^{20,116}$ \irep{208}:\ \ 
$d_i$ = ($1.0$,
$1.0$,
$1.618$,
$1.618$,
$1.618$,
$1.618$,
$2.618$,
$2.618$) 

\vskip 0.7ex
\hangindent=3em \hangafter=1
$D^2= 26.180 = 
15+5\sqrt{5}$

\vskip 0.7ex
\hangindent=3em \hangafter=1
$T = ( 0,
\frac{3}{4},
\frac{3}{5},
\frac{3}{5},
\frac{7}{20},
\frac{7}{20},
\frac{1}{5},
\frac{19}{20} )
$,

\vskip 0.7ex
\hangindent=3em \hangafter=1
$S$ = ($ 1$,
$ 1$,
$ \frac{1+\sqrt{5}}{2}$,
$ \frac{1+\sqrt{5}}{2}$,
$ \frac{1+\sqrt{5}}{2}$,
$ \frac{1+\sqrt{5}}{2}$,
$ \frac{3+\sqrt{5}}{2}$,
$ \frac{3+\sqrt{5}}{2}$;\ \ 
$ -1$,
$ \frac{1+\sqrt{5}}{2}$,
$ \frac{1+\sqrt{5}}{2}$,
$ -\frac{1+\sqrt{5}}{2}$,
$ -\frac{1+\sqrt{5}}{2}$,
$ \frac{3+\sqrt{5}}{2}$,
$ -\frac{3+\sqrt{5}}{2}$;\ \ 
$ -1$,
$ \frac{3+\sqrt{5}}{2}$,
$ -1$,
$ \frac{3+\sqrt{5}}{2}$,
$ -\frac{1+\sqrt{5}}{2}$,
$ -\frac{1+\sqrt{5}}{2}$;\ \ 
$ -1$,
$ \frac{3+\sqrt{5}}{2}$,
$ -1$,
$ -\frac{1+\sqrt{5}}{2}$,
$ -\frac{1+\sqrt{5}}{2}$;\ \ 
$ 1$,
$ -\frac{3+\sqrt{5}}{2}$,
$ -\frac{1+\sqrt{5}}{2}$,
$ \frac{1+\sqrt{5}}{2}$;\ \ 
$ 1$,
$ -\frac{1+\sqrt{5}}{2}$,
$ \frac{1+\sqrt{5}}{2}$;\ \ 
$ 1$,
$ 1$;\ \ 
$ -1$)

Factors = $2_{7,2.}^{4,625}\boxtimes 4_{\frac{12}{5},13.09}^{5,443}$

  \vskip 2ex

\noindent37. $8_{0,36.}^{6,213}$ \irep{0}:\ \ 
$d_i$ = ($1.0$,
$1.0$,
$2.0$,
$2.0$,
$2.0$,
$2.0$,
$3.0$,
$3.0$) 

\vskip 0.7ex
\hangindent=3em \hangafter=1
$D^2= 36.0 = 
36$

\vskip 0.7ex
\hangindent=3em \hangafter=1
$T = ( 0,
0,
0,
0,
\frac{1}{3},
\frac{2}{3},
0,
\frac{1}{2} )
$,

\vskip 0.7ex
\hangindent=3em \hangafter=1
$S$ = ($ 1$,
$ 1$,
$ 2$,
$ 2$,
$ 2$,
$ 2$,
$ 3$,
$ 3$;\ \ 
$ 1$,
$ 2$,
$ 2$,
$ 2$,
$ 2$,
$ -3$,
$ -3$;\ \ 
$ 4$,
$ -2$,
$ -2$,
$ -2$,
$0$,
$0$;\ \ 
$ 4$,
$ -2$,
$ -2$,
$0$,
$0$;\ \ 
$ -2$,
$ 4$,
$0$,
$0$;\ \ 
$ -2$,
$0$,
$0$;\ \ 
$ 3$,
$ -3$;\ \ 
$ 3$)

  \vskip 2ex

\noindent38. $8_{4,36.}^{6,102}$ \irep{0}:\ \ 
$d_i$ = ($1.0$,
$1.0$,
$2.0$,
$2.0$,
$2.0$,
$2.0$,
$3.0$,
$3.0$) 

\vskip 0.7ex
\hangindent=3em \hangafter=1
$D^2= 36.0 = 
36$

\vskip 0.7ex
\hangindent=3em \hangafter=1
$T = ( 0,
0,
\frac{1}{3},
\frac{1}{3},
\frac{2}{3},
\frac{2}{3},
0,
\frac{1}{2} )
$,

\vskip 0.7ex
\hangindent=3em \hangafter=1
$S$ = ($ 1$,
$ 1$,
$ 2$,
$ 2$,
$ 2$,
$ 2$,
$ 3$,
$ 3$;\ \ 
$ 1$,
$ 2$,
$ 2$,
$ 2$,
$ 2$,
$ -3$,
$ -3$;\ \ 
$ -2$,
$ 4$,
$ -2$,
$ -2$,
$0$,
$0$;\ \ 
$ -2$,
$ -2$,
$ -2$,
$0$,
$0$;\ \ 
$ -2$,
$ 4$,
$0$,
$0$;\ \ 
$ -2$,
$0$,
$0$;\ \ 
$ -3$,
$ 3$;\ \ 
$ -3$)

  \vskip 2ex

\noindent39. $8_{0,36.}^{12,101}$ \irep{0}:\ \ 
$d_i$ = ($1.0$,
$1.0$,
$2.0$,
$2.0$,
$2.0$,
$2.0$,
$3.0$,
$3.0$) 

\vskip 0.7ex
\hangindent=3em \hangafter=1
$D^2= 36.0 = 
36$

\vskip 0.7ex
\hangindent=3em \hangafter=1
$T = ( 0,
0,
0,
0,
\frac{1}{3},
\frac{2}{3},
\frac{1}{4},
\frac{3}{4} )
$,

\vskip 0.7ex
\hangindent=3em \hangafter=1
$S$ = ($ 1$,
$ 1$,
$ 2$,
$ 2$,
$ 2$,
$ 2$,
$ 3$,
$ 3$;\ \ 
$ 1$,
$ 2$,
$ 2$,
$ 2$,
$ 2$,
$ -3$,
$ -3$;\ \ 
$ 4$,
$ -2$,
$ -2$,
$ -2$,
$0$,
$0$;\ \ 
$ 4$,
$ -2$,
$ -2$,
$0$,
$0$;\ \ 
$ -2$,
$ 4$,
$0$,
$0$;\ \ 
$ -2$,
$0$,
$0$;\ \ 
$ -3$,
$ 3$;\ \ 
$ -3$)

  \vskip 2ex

\noindent40. $8_{4,36.}^{12,972}$ \irep{0}:\ \ 
$d_i$ = ($1.0$,
$1.0$,
$2.0$,
$2.0$,
$2.0$,
$2.0$,
$3.0$,
$3.0$) 

\vskip 0.7ex
\hangindent=3em \hangafter=1
$D^2= 36.0 = 
36$

\vskip 0.7ex
\hangindent=3em \hangafter=1
$T = ( 0,
0,
\frac{1}{3},
\frac{1}{3},
\frac{2}{3},
\frac{2}{3},
\frac{1}{4},
\frac{3}{4} )
$,

\vskip 0.7ex
\hangindent=3em \hangafter=1
$S$ = ($ 1$,
$ 1$,
$ 2$,
$ 2$,
$ 2$,
$ 2$,
$ 3$,
$ 3$;\ \ 
$ 1$,
$ 2$,
$ 2$,
$ 2$,
$ 2$,
$ -3$,
$ -3$;\ \ 
$ -2$,
$ 4$,
$ -2$,
$ -2$,
$0$,
$0$;\ \ 
$ -2$,
$ -2$,
$ -2$,
$0$,
$0$;\ \ 
$ -2$,
$ 4$,
$0$,
$0$;\ \ 
$ -2$,
$0$,
$0$;\ \ 
$ 3$,
$ -3$;\ \ 
$ 3$)

  \vskip 2ex

\noindent41. $8_{0,36.}^{18,162}$ \irep{0}:\ \ 
$d_i$ = ($1.0$,
$1.0$,
$2.0$,
$2.0$,
$2.0$,
$2.0$,
$3.0$,
$3.0$) 

\vskip 0.7ex
\hangindent=3em \hangafter=1
$D^2= 36.0 = 
36$

\vskip 0.7ex
\hangindent=3em \hangafter=1
$T = ( 0,
0,
0,
\frac{1}{9},
\frac{4}{9},
\frac{7}{9},
0,
\frac{1}{2} )
$,

\vskip 0.7ex
\hangindent=3em \hangafter=1
$S$ = ($ 1$,
$ 1$,
$ 2$,
$ 2$,
$ 2$,
$ 2$,
$ 3$,
$ 3$;\ \ 
$ 1$,
$ 2$,
$ 2$,
$ 2$,
$ 2$,
$ -3$,
$ -3$;\ \ 
$ 4$,
$ -2$,
$ -2$,
$ -2$,
$0$,
$0$;\ \ 
$ 2c_{9}^{2}$,
$ 2c_{9}^{4}$,
$ 2c_{9}^{1}$,
$0$,
$0$;\ \ 
$ 2c_{9}^{1}$,
$ 2c_{9}^{2}$,
$0$,
$0$;\ \ 
$ 2c_{9}^{4}$,
$0$,
$0$;\ \ 
$ 3$,
$ -3$;\ \ 
$ 3$)

  \vskip 2ex

\noindent42. $8_{0,36.}^{18,953}$ \irep{0}:\ \ 
$d_i$ = ($1.0$,
$1.0$,
$2.0$,
$2.0$,
$2.0$,
$2.0$,
$3.0$,
$3.0$) 

\vskip 0.7ex
\hangindent=3em \hangafter=1
$D^2= 36.0 = 
36$

\vskip 0.7ex
\hangindent=3em \hangafter=1
$T = ( 0,
0,
0,
\frac{2}{9},
\frac{5}{9},
\frac{8}{9},
0,
\frac{1}{2} )
$,

\vskip 0.7ex
\hangindent=3em \hangafter=1
$S$ = ($ 1$,
$ 1$,
$ 2$,
$ 2$,
$ 2$,
$ 2$,
$ 3$,
$ 3$;\ \ 
$ 1$,
$ 2$,
$ 2$,
$ 2$,
$ 2$,
$ -3$,
$ -3$;\ \ 
$ 4$,
$ -2$,
$ -2$,
$ -2$,
$0$,
$0$;\ \ 
$ 2c_{9}^{4}$,
$ 2c_{9}^{2}$,
$ 2c_{9}^{1}$,
$0$,
$0$;\ \ 
$ 2c_{9}^{1}$,
$ 2c_{9}^{4}$,
$0$,
$0$;\ \ 
$ 2c_{9}^{2}$,
$0$,
$0$;\ \ 
$ 3$,
$ -3$;\ \ 
$ 3$)

  \vskip 2ex

\noindent43. $8_{0,36.}^{36,495}$ \irep{0}:\ \ 
$d_i$ = ($1.0$,
$1.0$,
$2.0$,
$2.0$,
$2.0$,
$2.0$,
$3.0$,
$3.0$) 

\vskip 0.7ex
\hangindent=3em \hangafter=1
$D^2= 36.0 = 
36$

\vskip 0.7ex
\hangindent=3em \hangafter=1
$T = ( 0,
0,
0,
\frac{1}{9},
\frac{4}{9},
\frac{7}{9},
\frac{1}{4},
\frac{3}{4} )
$,

\vskip 0.7ex
\hangindent=3em \hangafter=1
$S$ = ($ 1$,
$ 1$,
$ 2$,
$ 2$,
$ 2$,
$ 2$,
$ 3$,
$ 3$;\ \ 
$ 1$,
$ 2$,
$ 2$,
$ 2$,
$ 2$,
$ -3$,
$ -3$;\ \ 
$ 4$,
$ -2$,
$ -2$,
$ -2$,
$0$,
$0$;\ \ 
$ 2c_{9}^{2}$,
$ 2c_{9}^{4}$,
$ 2c_{9}^{1}$,
$0$,
$0$;\ \ 
$ 2c_{9}^{1}$,
$ 2c_{9}^{2}$,
$0$,
$0$;\ \ 
$ 2c_{9}^{4}$,
$0$,
$0$;\ \ 
$ -3$,
$ 3$;\ \ 
$ -3$)

  \vskip 2ex

\noindent44. $8_{0,36.}^{36,171}$ \irep{0}:\ \ 
$d_i$ = ($1.0$,
$1.0$,
$2.0$,
$2.0$,
$2.0$,
$2.0$,
$3.0$,
$3.0$) 

\vskip 0.7ex
\hangindent=3em \hangafter=1
$D^2= 36.0 = 
36$

\vskip 0.7ex
\hangindent=3em \hangafter=1
$T = ( 0,
0,
0,
\frac{2}{9},
\frac{5}{9},
\frac{8}{9},
\frac{1}{4},
\frac{3}{4} )
$,

\vskip 0.7ex
\hangindent=3em \hangafter=1
$S$ = ($ 1$,
$ 1$,
$ 2$,
$ 2$,
$ 2$,
$ 2$,
$ 3$,
$ 3$;\ \ 
$ 1$,
$ 2$,
$ 2$,
$ 2$,
$ 2$,
$ -3$,
$ -3$;\ \ 
$ 4$,
$ -2$,
$ -2$,
$ -2$,
$0$,
$0$;\ \ 
$ 2c_{9}^{4}$,
$ 2c_{9}^{2}$,
$ 2c_{9}^{1}$,
$0$,
$0$;\ \ 
$ 2c_{9}^{1}$,
$ 2c_{9}^{4}$,
$0$,
$0$;\ \ 
$ 2c_{9}^{2}$,
$0$,
$0$;\ \ 
$ -3$,
$ 3$;\ \ 
$ -3$)

  \vskip 2ex

\noindent45. $8_{\frac{13}{3},38.46}^{36,115}$ \irep{229}:\ \ 
$d_i$ = ($1.0$,
$1.0$,
$1.879$,
$1.879$,
$2.532$,
$2.532$,
$2.879$,
$2.879$) 

\vskip 0.7ex
\hangindent=3em \hangafter=1
$D^2= 38.468 = 
18+12c^{1}_{9}
+6c^{2}_{9}
$

\vskip 0.7ex
\hangindent=3em \hangafter=1
$T = ( 0,
\frac{1}{4},
\frac{1}{3},
\frac{7}{12},
\frac{2}{9},
\frac{17}{36},
\frac{2}{3},
\frac{11}{12} )
$,

\vskip 0.7ex
\hangindent=3em \hangafter=1
$S$ = ($ 1$,
$ 1$,
$ -c_{9}^{4}$,
$ -c_{9}^{4}$,
$ \xi_{9}^{3}$,
$ \xi_{9}^{3}$,
$ \xi_{9}^{5}$,
$ \xi_{9}^{5}$;\ \ 
$ -1$,
$ -c_{9}^{4}$,
$ c_{9}^{4}$,
$ \xi_{9}^{3}$,
$ -\xi_{9}^{3}$,
$ \xi_{9}^{5}$,
$ -\xi_{9}^{5}$;\ \ 
$ -\xi_{9}^{5}$,
$ -\xi_{9}^{5}$,
$ \xi_{9}^{3}$,
$ \xi_{9}^{3}$,
$ -1$,
$ -1$;\ \ 
$ \xi_{9}^{5}$,
$ \xi_{9}^{3}$,
$ -\xi_{9}^{3}$,
$ -1$,
$ 1$;\ \ 
$0$,
$0$,
$ -\xi_{9}^{3}$,
$ -\xi_{9}^{3}$;\ \ 
$0$,
$ -\xi_{9}^{3}$,
$ \xi_{9}^{3}$;\ \ 
$ -c_{9}^{4}$,
$ -c_{9}^{4}$;\ \ 
$ c_{9}^{4}$)

Factors = $2_{1,2.}^{4,437}\boxtimes 4_{\frac{10}{3},19.23}^{9,459}$

  \vskip 2ex

\noindent46. $8_{\frac{17}{3},38.46}^{36,116}$ \irep{229}:\ \ 
$d_i$ = ($1.0$,
$1.0$,
$1.879$,
$1.879$,
$2.532$,
$2.532$,
$2.879$,
$2.879$) 

\vskip 0.7ex
\hangindent=3em \hangafter=1
$D^2= 38.468 = 
18+12c^{1}_{9}
+6c^{2}_{9}
$

\vskip 0.7ex
\hangindent=3em \hangafter=1
$T = ( 0,
\frac{1}{4},
\frac{2}{3},
\frac{11}{12},
\frac{7}{9},
\frac{1}{36},
\frac{1}{3},
\frac{7}{12} )
$,

\vskip 0.7ex
\hangindent=3em \hangafter=1
$S$ = ($ 1$,
$ 1$,
$ -c_{9}^{4}$,
$ -c_{9}^{4}$,
$ \xi_{9}^{3}$,
$ \xi_{9}^{3}$,
$ \xi_{9}^{5}$,
$ \xi_{9}^{5}$;\ \ 
$ -1$,
$ -c_{9}^{4}$,
$ c_{9}^{4}$,
$ \xi_{9}^{3}$,
$ -\xi_{9}^{3}$,
$ \xi_{9}^{5}$,
$ -\xi_{9}^{5}$;\ \ 
$ -\xi_{9}^{5}$,
$ -\xi_{9}^{5}$,
$ \xi_{9}^{3}$,
$ \xi_{9}^{3}$,
$ -1$,
$ -1$;\ \ 
$ \xi_{9}^{5}$,
$ \xi_{9}^{3}$,
$ -\xi_{9}^{3}$,
$ -1$,
$ 1$;\ \ 
$0$,
$0$,
$ -\xi_{9}^{3}$,
$ -\xi_{9}^{3}$;\ \ 
$0$,
$ -\xi_{9}^{3}$,
$ \xi_{9}^{3}$;\ \ 
$ -c_{9}^{4}$,
$ -c_{9}^{4}$;\ \ 
$ c_{9}^{4}$)

Factors = $2_{1,2.}^{4,437}\boxtimes 4_{\frac{14}{3},19.23}^{9,614}$

  \vskip 2ex

\noindent47. $8_{\frac{7}{3},38.46}^{36,936}$ \irep{229}:\ \ 
$d_i$ = ($1.0$,
$1.0$,
$1.879$,
$1.879$,
$2.532$,
$2.532$,
$2.879$,
$2.879$) 

\vskip 0.7ex
\hangindent=3em \hangafter=1
$D^2= 38.468 = 
18+12c^{1}_{9}
+6c^{2}_{9}
$

\vskip 0.7ex
\hangindent=3em \hangafter=1
$T = ( 0,
\frac{3}{4},
\frac{1}{3},
\frac{1}{12},
\frac{2}{9},
\frac{35}{36},
\frac{2}{3},
\frac{5}{12} )
$,

\vskip 0.7ex
\hangindent=3em \hangafter=1
$S$ = ($ 1$,
$ 1$,
$ -c_{9}^{4}$,
$ -c_{9}^{4}$,
$ \xi_{9}^{3}$,
$ \xi_{9}^{3}$,
$ \xi_{9}^{5}$,
$ \xi_{9}^{5}$;\ \ 
$ -1$,
$ -c_{9}^{4}$,
$ c_{9}^{4}$,
$ \xi_{9}^{3}$,
$ -\xi_{9}^{3}$,
$ \xi_{9}^{5}$,
$ -\xi_{9}^{5}$;\ \ 
$ -\xi_{9}^{5}$,
$ -\xi_{9}^{5}$,
$ \xi_{9}^{3}$,
$ \xi_{9}^{3}$,
$ -1$,
$ -1$;\ \ 
$ \xi_{9}^{5}$,
$ \xi_{9}^{3}$,
$ -\xi_{9}^{3}$,
$ -1$,
$ 1$;\ \ 
$0$,
$0$,
$ -\xi_{9}^{3}$,
$ -\xi_{9}^{3}$;\ \ 
$0$,
$ -\xi_{9}^{3}$,
$ \xi_{9}^{3}$;\ \ 
$ -c_{9}^{4}$,
$ -c_{9}^{4}$;\ \ 
$ c_{9}^{4}$)

Factors = $2_{7,2.}^{4,625}\boxtimes 4_{\frac{10}{3},19.23}^{9,459}$

  \vskip 2ex

\noindent48. $8_{\frac{11}{3},38.46}^{36,155}$ \irep{229}:\ \ 
$d_i$ = ($1.0$,
$1.0$,
$1.879$,
$1.879$,
$2.532$,
$2.532$,
$2.879$,
$2.879$) 

\vskip 0.7ex
\hangindent=3em \hangafter=1
$D^2= 38.468 = 
18+12c^{1}_{9}
+6c^{2}_{9}
$

\vskip 0.7ex
\hangindent=3em \hangafter=1
$T = ( 0,
\frac{3}{4},
\frac{2}{3},
\frac{5}{12},
\frac{7}{9},
\frac{19}{36},
\frac{1}{3},
\frac{1}{12} )
$,

\vskip 0.7ex
\hangindent=3em \hangafter=1
$S$ = ($ 1$,
$ 1$,
$ -c_{9}^{4}$,
$ -c_{9}^{4}$,
$ \xi_{9}^{3}$,
$ \xi_{9}^{3}$,
$ \xi_{9}^{5}$,
$ \xi_{9}^{5}$;\ \ 
$ -1$,
$ -c_{9}^{4}$,
$ c_{9}^{4}$,
$ \xi_{9}^{3}$,
$ -\xi_{9}^{3}$,
$ \xi_{9}^{5}$,
$ -\xi_{9}^{5}$;\ \ 
$ -\xi_{9}^{5}$,
$ -\xi_{9}^{5}$,
$ \xi_{9}^{3}$,
$ \xi_{9}^{3}$,
$ -1$,
$ -1$;\ \ 
$ \xi_{9}^{5}$,
$ \xi_{9}^{3}$,
$ -\xi_{9}^{3}$,
$ -1$,
$ 1$;\ \ 
$0$,
$0$,
$ -\xi_{9}^{3}$,
$ -\xi_{9}^{3}$;\ \ 
$0$,
$ -\xi_{9}^{3}$,
$ \xi_{9}^{3}$;\ \ 
$ -c_{9}^{4}$,
$ -c_{9}^{4}$;\ \ 
$ c_{9}^{4}$)

Factors = $2_{7,2.}^{4,625}\boxtimes 4_{\frac{14}{3},19.23}^{9,614}$

  \vskip 2ex

\noindent49. $8_{\frac{2}{5},47.36}^{5,148}$ \irep{52}:\ \ 
$d_i$ = ($1.0$,
$1.618$,
$1.618$,
$1.618$,
$2.618$,
$2.618$,
$2.618$,
$4.236$) 

\vskip 0.7ex
\hangindent=3em \hangafter=1
$D^2= 47.360 = 
25+10\sqrt{5}$

\vskip 0.7ex
\hangindent=3em \hangafter=1
$T = ( 0,
\frac{2}{5},
\frac{2}{5},
\frac{2}{5},
\frac{4}{5},
\frac{4}{5},
\frac{4}{5},
\frac{1}{5} )
$,

\vskip 0.7ex
\hangindent=3em \hangafter=1
$S$ = ($ 1$,
$ \frac{1+\sqrt{5}}{2}$,
$ \frac{1+\sqrt{5}}{2}$,
$ \frac{1+\sqrt{5}}{2}$,
$ \frac{3+\sqrt{5}}{2}$,
$ \frac{3+\sqrt{5}}{2}$,
$ \frac{3+\sqrt{5}}{2}$,
$ 2+\sqrt{5}$;\ \ 
$ -1$,
$ \frac{3+\sqrt{5}}{2}$,
$ \frac{3+\sqrt{5}}{2}$,
$ 2+\sqrt{5}$,
$ -\frac{1+\sqrt{5}}{2}$,
$ -\frac{1+\sqrt{5}}{2}$,
$ -\frac{3+\sqrt{5}}{2}$;\ \ 
$ -1$,
$ \frac{3+\sqrt{5}}{2}$,
$ -\frac{1+\sqrt{5}}{2}$,
$ 2+\sqrt{5}$,
$ -\frac{1+\sqrt{5}}{2}$,
$ -\frac{3+\sqrt{5}}{2}$;\ \ 
$ -1$,
$ -\frac{1+\sqrt{5}}{2}$,
$ -\frac{1+\sqrt{5}}{2}$,
$ 2+\sqrt{5}$,
$ -\frac{3+\sqrt{5}}{2}$;\ \ 
$ 1$,
$ -\frac{3+\sqrt{5}}{2}$,
$ -\frac{3+\sqrt{5}}{2}$,
$ \frac{1+\sqrt{5}}{2}$;\ \ 
$ 1$,
$ -\frac{3+\sqrt{5}}{2}$,
$ \frac{1+\sqrt{5}}{2}$;\ \ 
$ 1$,
$ \frac{1+\sqrt{5}}{2}$;\ \ 
$ -1$)

Factors = $2_{\frac{14}{5},3.618}^{5,395}\boxtimes 4_{\frac{28}{5},13.09}^{5,479}$

  \vskip 2ex

\noindent50. $8_{\frac{14}{5},47.36}^{5,103}$ \irep{52}:\ \ 
$d_i$ = ($1.0$,
$1.618$,
$1.618$,
$1.618$,
$2.618$,
$2.618$,
$2.618$,
$4.236$) 

\vskip 0.7ex
\hangindent=3em \hangafter=1
$D^2= 47.360 = 
25+10\sqrt{5}$

\vskip 0.7ex
\hangindent=3em \hangafter=1
$T = ( 0,
\frac{2}{5},
\frac{2}{5},
\frac{3}{5},
0,
0,
\frac{4}{5},
\frac{2}{5} )
$,

\vskip 0.7ex
\hangindent=3em \hangafter=1
$S$ = ($ 1$,
$ \frac{1+\sqrt{5}}{2}$,
$ \frac{1+\sqrt{5}}{2}$,
$ \frac{1+\sqrt{5}}{2}$,
$ \frac{3+\sqrt{5}}{2}$,
$ \frac{3+\sqrt{5}}{2}$,
$ \frac{3+\sqrt{5}}{2}$,
$ 2+\sqrt{5}$;\ \ 
$ -1$,
$ \frac{3+\sqrt{5}}{2}$,
$ \frac{3+\sqrt{5}}{2}$,
$ 2+\sqrt{5}$,
$ -\frac{1+\sqrt{5}}{2}$,
$ -\frac{1+\sqrt{5}}{2}$,
$ -\frac{3+\sqrt{5}}{2}$;\ \ 
$ -1$,
$ \frac{3+\sqrt{5}}{2}$,
$ -\frac{1+\sqrt{5}}{2}$,
$ 2+\sqrt{5}$,
$ -\frac{1+\sqrt{5}}{2}$,
$ -\frac{3+\sqrt{5}}{2}$;\ \ 
$ -1$,
$ -\frac{1+\sqrt{5}}{2}$,
$ -\frac{1+\sqrt{5}}{2}$,
$ 2+\sqrt{5}$,
$ -\frac{3+\sqrt{5}}{2}$;\ \ 
$ 1$,
$ -\frac{3+\sqrt{5}}{2}$,
$ -\frac{3+\sqrt{5}}{2}$,
$ \frac{1+\sqrt{5}}{2}$;\ \ 
$ 1$,
$ -\frac{3+\sqrt{5}}{2}$,
$ \frac{1+\sqrt{5}}{2}$;\ \ 
$ 1$,
$ \frac{1+\sqrt{5}}{2}$;\ \ 
$ -1$)

Factors = $2_{\frac{26}{5},3.618}^{5,720}\boxtimes 4_{\frac{28}{5},13.09}^{5,479}$

  \vskip 2ex

\noindent51. $8_{\frac{26}{5},47.36}^{5,143}$ \irep{52}:\ \ 
$d_i$ = ($1.0$,
$1.618$,
$1.618$,
$1.618$,
$2.618$,
$2.618$,
$2.618$,
$4.236$) 

\vskip 0.7ex
\hangindent=3em \hangafter=1
$D^2= 47.360 = 
25+10\sqrt{5}$

\vskip 0.7ex
\hangindent=3em \hangafter=1
$T = ( 0,
\frac{2}{5},
\frac{3}{5},
\frac{3}{5},
0,
0,
\frac{1}{5},
\frac{3}{5} )
$,

\vskip 0.7ex
\hangindent=3em \hangafter=1
$S$ = ($ 1$,
$ \frac{1+\sqrt{5}}{2}$,
$ \frac{1+\sqrt{5}}{2}$,
$ \frac{1+\sqrt{5}}{2}$,
$ \frac{3+\sqrt{5}}{2}$,
$ \frac{3+\sqrt{5}}{2}$,
$ \frac{3+\sqrt{5}}{2}$,
$ 2+\sqrt{5}$;\ \ 
$ -1$,
$ \frac{3+\sqrt{5}}{2}$,
$ \frac{3+\sqrt{5}}{2}$,
$ -\frac{1+\sqrt{5}}{2}$,
$ -\frac{1+\sqrt{5}}{2}$,
$ 2+\sqrt{5}$,
$ -\frac{3+\sqrt{5}}{2}$;\ \ 
$ -1$,
$ \frac{3+\sqrt{5}}{2}$,
$ 2+\sqrt{5}$,
$ -\frac{1+\sqrt{5}}{2}$,
$ -\frac{1+\sqrt{5}}{2}$,
$ -\frac{3+\sqrt{5}}{2}$;\ \ 
$ -1$,
$ -\frac{1+\sqrt{5}}{2}$,
$ 2+\sqrt{5}$,
$ -\frac{1+\sqrt{5}}{2}$,
$ -\frac{3+\sqrt{5}}{2}$;\ \ 
$ 1$,
$ -\frac{3+\sqrt{5}}{2}$,
$ -\frac{3+\sqrt{5}}{2}$,
$ \frac{1+\sqrt{5}}{2}$;\ \ 
$ 1$,
$ -\frac{3+\sqrt{5}}{2}$,
$ \frac{1+\sqrt{5}}{2}$;\ \ 
$ 1$,
$ \frac{1+\sqrt{5}}{2}$;\ \ 
$ -1$)

Factors = $2_{\frac{14}{5},3.618}^{5,395}\boxtimes 4_{\frac{12}{5},13.09}^{5,443}$

  \vskip 2ex

\noindent52. $8_{\frac{38}{5},47.36}^{5,286}$ \irep{52}:\ \ 
$d_i$ = ($1.0$,
$1.618$,
$1.618$,
$1.618$,
$2.618$,
$2.618$,
$2.618$,
$4.236$) 

\vskip 0.7ex
\hangindent=3em \hangafter=1
$D^2= 47.360 = 
25+10\sqrt{5}$

\vskip 0.7ex
\hangindent=3em \hangafter=1
$T = ( 0,
\frac{3}{5},
\frac{3}{5},
\frac{3}{5},
\frac{1}{5},
\frac{1}{5},
\frac{1}{5},
\frac{4}{5} )
$,

\vskip 0.7ex
\hangindent=3em \hangafter=1
$S$ = ($ 1$,
$ \frac{1+\sqrt{5}}{2}$,
$ \frac{1+\sqrt{5}}{2}$,
$ \frac{1+\sqrt{5}}{2}$,
$ \frac{3+\sqrt{5}}{2}$,
$ \frac{3+\sqrt{5}}{2}$,
$ \frac{3+\sqrt{5}}{2}$,
$ 2+\sqrt{5}$;\ \ 
$ -1$,
$ \frac{3+\sqrt{5}}{2}$,
$ \frac{3+\sqrt{5}}{2}$,
$ 2+\sqrt{5}$,
$ -\frac{1+\sqrt{5}}{2}$,
$ -\frac{1+\sqrt{5}}{2}$,
$ -\frac{3+\sqrt{5}}{2}$;\ \ 
$ -1$,
$ \frac{3+\sqrt{5}}{2}$,
$ -\frac{1+\sqrt{5}}{2}$,
$ 2+\sqrt{5}$,
$ -\frac{1+\sqrt{5}}{2}$,
$ -\frac{3+\sqrt{5}}{2}$;\ \ 
$ -1$,
$ -\frac{1+\sqrt{5}}{2}$,
$ -\frac{1+\sqrt{5}}{2}$,
$ 2+\sqrt{5}$,
$ -\frac{3+\sqrt{5}}{2}$;\ \ 
$ 1$,
$ -\frac{3+\sqrt{5}}{2}$,
$ -\frac{3+\sqrt{5}}{2}$,
$ \frac{1+\sqrt{5}}{2}$;\ \ 
$ 1$,
$ -\frac{3+\sqrt{5}}{2}$,
$ \frac{1+\sqrt{5}}{2}$;\ \ 
$ 1$,
$ \frac{1+\sqrt{5}}{2}$;\ \ 
$ -1$)

Factors = $2_{\frac{26}{5},3.618}^{5,720}\boxtimes 4_{\frac{12}{5},13.09}^{5,443}$

  \vskip 2ex

\noindent53. $8_{\frac{92}{15},69.59}^{45,311}$ \irep{234}:\ \ 
$d_i$ = ($1.0$,
$1.618$,
$1.879$,
$2.532$,
$2.879$,
$3.40$,
$4.97$,
$4.658$) 

\vskip 0.7ex
\hangindent=3em \hangafter=1
$D^2= 69.590 = 
27-3  c^{1}_{45}
+6c^{2}_{45}
+3c^{4}_{45}
+12c^{5}_{45}
+6c^{7}_{45}
+9c^{9}_{45}
+3c^{10}_{45}
+3c^{11}_{45}
$

\vskip 0.7ex
\hangindent=3em \hangafter=1
$T = ( 0,
\frac{2}{5},
\frac{1}{3},
\frac{2}{9},
\frac{2}{3},
\frac{11}{15},
\frac{28}{45},
\frac{1}{15} )
$,

\vskip 0.7ex
\hangindent=3em \hangafter=1
$S$ = ($ 1$,
$ \frac{1+\sqrt{5}}{2}$,
$ -c_{9}^{4}$,
$ \xi_{9}^{3}$,
$ \xi_{9}^{5}$,
$ c^{2}_{45}
+c^{7}_{45}
$,
$ 1-c^{1}_{45}
+c^{2}_{45}
+c^{4}_{45}
+c^{7}_{45}
+c^{9}_{45}
-c^{10}_{45}
+c^{11}_{45}
$,
$ 1+c^{2}_{45}
+c^{7}_{45}
+c^{9}_{45}
$;\ \ 
$ -1$,
$ c^{2}_{45}
+c^{7}_{45}
$,
$ 1-c^{1}_{45}
+c^{2}_{45}
+c^{4}_{45}
+c^{7}_{45}
+c^{9}_{45}
-c^{10}_{45}
+c^{11}_{45}
$,
$ 1+c^{2}_{45}
+c^{7}_{45}
+c^{9}_{45}
$,
$ c_{9}^{4}$,
$ -\xi_{9}^{3}$,
$ -\xi_{9}^{5}$;\ \ 
$ -\xi_{9}^{5}$,
$ \xi_{9}^{3}$,
$ -1$,
$ -1-c^{2}_{45}
-c^{7}_{45}
-c^{9}_{45}
$,
$ 1-c^{1}_{45}
+c^{2}_{45}
+c^{4}_{45}
+c^{7}_{45}
+c^{9}_{45}
-c^{10}_{45}
+c^{11}_{45}
$,
$ -\frac{1+\sqrt{5}}{2}$;\ \ 
$0$,
$ -\xi_{9}^{3}$,
$ 1-c^{1}_{45}
+c^{2}_{45}
+c^{4}_{45}
+c^{7}_{45}
+c^{9}_{45}
-c^{10}_{45}
+c^{11}_{45}
$,
$0$,
$ -1+c^{1}_{45}
-c^{2}_{45}
-c^{4}_{45}
-c^{7}_{45}
-c^{9}_{45}
+c^{10}_{45}
-c^{11}_{45}
$;\ \ 
$ -c_{9}^{4}$,
$ -\frac{1+\sqrt{5}}{2}$,
$ -1+c^{1}_{45}
-c^{2}_{45}
-c^{4}_{45}
-c^{7}_{45}
-c^{9}_{45}
+c^{10}_{45}
-c^{11}_{45}
$,
$ c^{2}_{45}
+c^{7}_{45}
$;\ \ 
$ \xi_{9}^{5}$,
$ -\xi_{9}^{3}$,
$ 1$;\ \ 
$0$,
$ \xi_{9}^{3}$;\ \ 
$ c_{9}^{4}$)

Factors = $2_{\frac{14}{5},3.618}^{5,395}\boxtimes 4_{\frac{10}{3},19.23}^{9,459}$

  \vskip 2ex

\noindent54. $8_{\frac{112}{15},69.59}^{45,167}$ \irep{234}:\ \ 
$d_i$ = ($1.0$,
$1.618$,
$1.879$,
$2.532$,
$2.879$,
$3.40$,
$4.97$,
$4.658$) 

\vskip 0.7ex
\hangindent=3em \hangafter=1
$D^2= 69.590 = 
27-3  c^{1}_{45}
+6c^{2}_{45}
+3c^{4}_{45}
+12c^{5}_{45}
+6c^{7}_{45}
+9c^{9}_{45}
+3c^{10}_{45}
+3c^{11}_{45}
$

\vskip 0.7ex
\hangindent=3em \hangafter=1
$T = ( 0,
\frac{2}{5},
\frac{2}{3},
\frac{7}{9},
\frac{1}{3},
\frac{1}{15},
\frac{8}{45},
\frac{11}{15} )
$,

\vskip 0.7ex
\hangindent=3em \hangafter=1
$S$ = ($ 1$,
$ \frac{1+\sqrt{5}}{2}$,
$ -c_{9}^{4}$,
$ \xi_{9}^{3}$,
$ \xi_{9}^{5}$,
$ c^{2}_{45}
+c^{7}_{45}
$,
$ 1-c^{1}_{45}
+c^{2}_{45}
+c^{4}_{45}
+c^{7}_{45}
+c^{9}_{45}
-c^{10}_{45}
+c^{11}_{45}
$,
$ 1+c^{2}_{45}
+c^{7}_{45}
+c^{9}_{45}
$;\ \ 
$ -1$,
$ c^{2}_{45}
+c^{7}_{45}
$,
$ 1-c^{1}_{45}
+c^{2}_{45}
+c^{4}_{45}
+c^{7}_{45}
+c^{9}_{45}
-c^{10}_{45}
+c^{11}_{45}
$,
$ 1+c^{2}_{45}
+c^{7}_{45}
+c^{9}_{45}
$,
$ c_{9}^{4}$,
$ -\xi_{9}^{3}$,
$ -\xi_{9}^{5}$;\ \ 
$ -\xi_{9}^{5}$,
$ \xi_{9}^{3}$,
$ -1$,
$ -1-c^{2}_{45}
-c^{7}_{45}
-c^{9}_{45}
$,
$ 1-c^{1}_{45}
+c^{2}_{45}
+c^{4}_{45}
+c^{7}_{45}
+c^{9}_{45}
-c^{10}_{45}
+c^{11}_{45}
$,
$ -\frac{1+\sqrt{5}}{2}$;\ \ 
$0$,
$ -\xi_{9}^{3}$,
$ 1-c^{1}_{45}
+c^{2}_{45}
+c^{4}_{45}
+c^{7}_{45}
+c^{9}_{45}
-c^{10}_{45}
+c^{11}_{45}
$,
$0$,
$ -1+c^{1}_{45}
-c^{2}_{45}
-c^{4}_{45}
-c^{7}_{45}
-c^{9}_{45}
+c^{10}_{45}
-c^{11}_{45}
$;\ \ 
$ -c_{9}^{4}$,
$ -\frac{1+\sqrt{5}}{2}$,
$ -1+c^{1}_{45}
-c^{2}_{45}
-c^{4}_{45}
-c^{7}_{45}
-c^{9}_{45}
+c^{10}_{45}
-c^{11}_{45}
$,
$ c^{2}_{45}
+c^{7}_{45}
$;\ \ 
$ \xi_{9}^{5}$,
$ -\xi_{9}^{3}$,
$ 1$;\ \ 
$0$,
$ \xi_{9}^{3}$;\ \ 
$ c_{9}^{4}$)

Factors = $2_{\frac{14}{5},3.618}^{5,395}\boxtimes 4_{\frac{14}{3},19.23}^{9,614}$

  \vskip 2ex

\noindent55. $8_{\frac{8}{15},69.59}^{45,251}$ \irep{234}:\ \ 
$d_i$ = ($1.0$,
$1.618$,
$1.879$,
$2.532$,
$2.879$,
$3.40$,
$4.97$,
$4.658$) 

\vskip 0.7ex
\hangindent=3em \hangafter=1
$D^2= 69.590 = 
27-3  c^{1}_{45}
+6c^{2}_{45}
+3c^{4}_{45}
+12c^{5}_{45}
+6c^{7}_{45}
+9c^{9}_{45}
+3c^{10}_{45}
+3c^{11}_{45}
$

\vskip 0.7ex
\hangindent=3em \hangafter=1
$T = ( 0,
\frac{3}{5},
\frac{1}{3},
\frac{2}{9},
\frac{2}{3},
\frac{14}{15},
\frac{37}{45},
\frac{4}{15} )
$,

\vskip 0.7ex
\hangindent=3em \hangafter=1
$S$ = ($ 1$,
$ \frac{1+\sqrt{5}}{2}$,
$ -c_{9}^{4}$,
$ \xi_{9}^{3}$,
$ \xi_{9}^{5}$,
$ c^{2}_{45}
+c^{7}_{45}
$,
$ 1-c^{1}_{45}
+c^{2}_{45}
+c^{4}_{45}
+c^{7}_{45}
+c^{9}_{45}
-c^{10}_{45}
+c^{11}_{45}
$,
$ 1+c^{2}_{45}
+c^{7}_{45}
+c^{9}_{45}
$;\ \ 
$ -1$,
$ c^{2}_{45}
+c^{7}_{45}
$,
$ 1-c^{1}_{45}
+c^{2}_{45}
+c^{4}_{45}
+c^{7}_{45}
+c^{9}_{45}
-c^{10}_{45}
+c^{11}_{45}
$,
$ 1+c^{2}_{45}
+c^{7}_{45}
+c^{9}_{45}
$,
$ c_{9}^{4}$,
$ -\xi_{9}^{3}$,
$ -\xi_{9}^{5}$;\ \ 
$ -\xi_{9}^{5}$,
$ \xi_{9}^{3}$,
$ -1$,
$ -1-c^{2}_{45}
-c^{7}_{45}
-c^{9}_{45}
$,
$ 1-c^{1}_{45}
+c^{2}_{45}
+c^{4}_{45}
+c^{7}_{45}
+c^{9}_{45}
-c^{10}_{45}
+c^{11}_{45}
$,
$ -\frac{1+\sqrt{5}}{2}$;\ \ 
$0$,
$ -\xi_{9}^{3}$,
$ 1-c^{1}_{45}
+c^{2}_{45}
+c^{4}_{45}
+c^{7}_{45}
+c^{9}_{45}
-c^{10}_{45}
+c^{11}_{45}
$,
$0$,
$ -1+c^{1}_{45}
-c^{2}_{45}
-c^{4}_{45}
-c^{7}_{45}
-c^{9}_{45}
+c^{10}_{45}
-c^{11}_{45}
$;\ \ 
$ -c_{9}^{4}$,
$ -\frac{1+\sqrt{5}}{2}$,
$ -1+c^{1}_{45}
-c^{2}_{45}
-c^{4}_{45}
-c^{7}_{45}
-c^{9}_{45}
+c^{10}_{45}
-c^{11}_{45}
$,
$ c^{2}_{45}
+c^{7}_{45}
$;\ \ 
$ \xi_{9}^{5}$,
$ -\xi_{9}^{3}$,
$ 1$;\ \ 
$0$,
$ \xi_{9}^{3}$;\ \ 
$ c_{9}^{4}$)

Factors = $2_{\frac{26}{5},3.618}^{5,720}\boxtimes 4_{\frac{10}{3},19.23}^{9,459}$

  \vskip 2ex

\noindent56. $8_{\frac{28}{15},69.59}^{45,270}$ \irep{234}:\ \ 
$d_i$ = ($1.0$,
$1.618$,
$1.879$,
$2.532$,
$2.879$,
$3.40$,
$4.97$,
$4.658$) 

\vskip 0.7ex
\hangindent=3em \hangafter=1
$D^2= 69.590 = 
27-3  c^{1}_{45}
+6c^{2}_{45}
+3c^{4}_{45}
+12c^{5}_{45}
+6c^{7}_{45}
+9c^{9}_{45}
+3c^{10}_{45}
+3c^{11}_{45}
$

\vskip 0.7ex
\hangindent=3em \hangafter=1
$T = ( 0,
\frac{3}{5},
\frac{2}{3},
\frac{7}{9},
\frac{1}{3},
\frac{4}{15},
\frac{17}{45},
\frac{14}{15} )
$,

\vskip 0.7ex
\hangindent=3em \hangafter=1
$S$ = ($ 1$,
$ \frac{1+\sqrt{5}}{2}$,
$ -c_{9}^{4}$,
$ \xi_{9}^{3}$,
$ \xi_{9}^{5}$,
$ c^{2}_{45}
+c^{7}_{45}
$,
$ 1-c^{1}_{45}
+c^{2}_{45}
+c^{4}_{45}
+c^{7}_{45}
+c^{9}_{45}
-c^{10}_{45}
+c^{11}_{45}
$,
$ 1+c^{2}_{45}
+c^{7}_{45}
+c^{9}_{45}
$;\ \ 
$ -1$,
$ c^{2}_{45}
+c^{7}_{45}
$,
$ 1-c^{1}_{45}
+c^{2}_{45}
+c^{4}_{45}
+c^{7}_{45}
+c^{9}_{45}
-c^{10}_{45}
+c^{11}_{45}
$,
$ 1+c^{2}_{45}
+c^{7}_{45}
+c^{9}_{45}
$,
$ c_{9}^{4}$,
$ -\xi_{9}^{3}$,
$ -\xi_{9}^{5}$;\ \ 
$ -\xi_{9}^{5}$,
$ \xi_{9}^{3}$,
$ -1$,
$ -1-c^{2}_{45}
-c^{7}_{45}
-c^{9}_{45}
$,
$ 1-c^{1}_{45}
+c^{2}_{45}
+c^{4}_{45}
+c^{7}_{45}
+c^{9}_{45}
-c^{10}_{45}
+c^{11}_{45}
$,
$ -\frac{1+\sqrt{5}}{2}$;\ \ 
$0$,
$ -\xi_{9}^{3}$,
$ 1-c^{1}_{45}
+c^{2}_{45}
+c^{4}_{45}
+c^{7}_{45}
+c^{9}_{45}
-c^{10}_{45}
+c^{11}_{45}
$,
$0$,
$ -1+c^{1}_{45}
-c^{2}_{45}
-c^{4}_{45}
-c^{7}_{45}
-c^{9}_{45}
+c^{10}_{45}
-c^{11}_{45}
$;\ \ 
$ -c_{9}^{4}$,
$ -\frac{1+\sqrt{5}}{2}$,
$ -1+c^{1}_{45}
-c^{2}_{45}
-c^{4}_{45}
-c^{7}_{45}
-c^{9}_{45}
+c^{10}_{45}
-c^{11}_{45}
$,
$ c^{2}_{45}
+c^{7}_{45}
$;\ \ 
$ \xi_{9}^{5}$,
$ -\xi_{9}^{3}$,
$ 1$;\ \ 
$0$,
$ \xi_{9}^{3}$;\ \ 
$ c_{9}^{4}$)

Factors = $2_{\frac{26}{5},3.618}^{5,720}\boxtimes 4_{\frac{14}{3},19.23}^{9,614}$

  \vskip 2ex

\noindent57. $8_{\frac{62}{17},125.8}^{17,152}$ \irep{199}:\ \ 
$d_i$ = ($1.0$,
$1.965$,
$2.864$,
$3.666$,
$4.342$,
$4.871$,
$5.234$,
$5.418$) 

\vskip 0.7ex
\hangindent=3em \hangafter=1
$D^2= 125.874 = 
36+28c^{1}_{17}
+21c^{2}_{17}
+15c^{3}_{17}
+10c^{4}_{17}
+6c^{5}_{17}
+3c^{6}_{17}
+c^{7}_{17}
$

\vskip 0.7ex
\hangindent=3em \hangafter=1
$T = ( 0,
\frac{5}{17},
\frac{2}{17},
\frac{8}{17},
\frac{6}{17},
\frac{13}{17},
\frac{12}{17},
\frac{3}{17} )
$,

\vskip 0.7ex
\hangindent=3em \hangafter=1
$S$ = ($ 1$,
$ -c_{17}^{8}$,
$ \xi_{17}^{3}$,
$ \xi_{17}^{13}$,
$ \xi_{17}^{5}$,
$ \xi_{17}^{11}$,
$ \xi_{17}^{7}$,
$ \xi_{17}^{9}$;\ \ 
$ -\xi_{17}^{13}$,
$ \xi_{17}^{11}$,
$ -\xi_{17}^{9}$,
$ \xi_{17}^{7}$,
$ -\xi_{17}^{5}$,
$ \xi_{17}^{3}$,
$ -1$;\ \ 
$ \xi_{17}^{9}$,
$ \xi_{17}^{5}$,
$ -c_{17}^{8}$,
$ -1$,
$ -\xi_{17}^{13}$,
$ -\xi_{17}^{7}$;\ \ 
$ -1$,
$ -\xi_{17}^{3}$,
$ \xi_{17}^{7}$,
$ -\xi_{17}^{11}$,
$ -c_{17}^{8}$;\ \ 
$ -\xi_{17}^{9}$,
$ -\xi_{17}^{13}$,
$ 1$,
$ \xi_{17}^{11}$;\ \ 
$ c_{17}^{8}$,
$ \xi_{17}^{9}$,
$ -\xi_{17}^{3}$;\ \ 
$ -c_{17}^{8}$,
$ -\xi_{17}^{5}$;\ \ 
$ \xi_{17}^{13}$)

  \vskip 2ex

\noindent58. $8_{\frac{74}{17},125.8}^{17,311}$ \irep{199}:\ \ 
$d_i$ = ($1.0$,
$1.965$,
$2.864$,
$3.666$,
$4.342$,
$4.871$,
$5.234$,
$5.418$) 

\vskip 0.7ex
\hangindent=3em \hangafter=1
$D^2= 125.874 = 
36+28c^{1}_{17}
+21c^{2}_{17}
+15c^{3}_{17}
+10c^{4}_{17}
+6c^{5}_{17}
+3c^{6}_{17}
+c^{7}_{17}
$

\vskip 0.7ex
\hangindent=3em \hangafter=1
$T = ( 0,
\frac{12}{17},
\frac{15}{17},
\frac{9}{17},
\frac{11}{17},
\frac{4}{17},
\frac{5}{17},
\frac{14}{17} )
$,

\vskip 0.7ex
\hangindent=3em \hangafter=1
$S$ = ($ 1$,
$ -c_{17}^{8}$,
$ \xi_{17}^{3}$,
$ \xi_{17}^{13}$,
$ \xi_{17}^{5}$,
$ \xi_{17}^{11}$,
$ \xi_{17}^{7}$,
$ \xi_{17}^{9}$;\ \ 
$ -\xi_{17}^{13}$,
$ \xi_{17}^{11}$,
$ -\xi_{17}^{9}$,
$ \xi_{17}^{7}$,
$ -\xi_{17}^{5}$,
$ \xi_{17}^{3}$,
$ -1$;\ \ 
$ \xi_{17}^{9}$,
$ \xi_{17}^{5}$,
$ -c_{17}^{8}$,
$ -1$,
$ -\xi_{17}^{13}$,
$ -\xi_{17}^{7}$;\ \ 
$ -1$,
$ -\xi_{17}^{3}$,
$ \xi_{17}^{7}$,
$ -\xi_{17}^{11}$,
$ -c_{17}^{8}$;\ \ 
$ -\xi_{17}^{9}$,
$ -\xi_{17}^{13}$,
$ 1$,
$ \xi_{17}^{11}$;\ \ 
$ c_{17}^{8}$,
$ \xi_{17}^{9}$,
$ -\xi_{17}^{3}$;\ \ 
$ -c_{17}^{8}$,
$ -\xi_{17}^{5}$;\ \ 
$ \xi_{17}^{13}$)

  \vskip 2ex

\noindent59. $8_{\frac{36}{13},223.6}^{13,370}$ \irep{179}:\ \ 
$d_i$ = ($1.0$,
$2.941$,
$4.148$,
$4.148$,
$4.712$,
$6.209$,
$7.345$,
$8.55$) 

\vskip 0.7ex
\hangindent=3em \hangafter=1
$D^2= 223.689 = 
78+65c^{1}_{13}
+52c^{2}_{13}
+39c^{3}_{13}
+26c^{4}_{13}
+13c^{5}_{13}
$

\vskip 0.7ex
\hangindent=3em \hangafter=1
$T = ( 0,
\frac{1}{13},
\frac{8}{13},
\frac{8}{13},
\frac{3}{13},
\frac{6}{13},
\frac{10}{13},
\frac{2}{13} )
$,

\vskip 0.7ex
\hangindent=3em \hangafter=1
$S$ = ($ 1$,
$ 2+c^{1}_{13}
+c^{2}_{13}
+c^{3}_{13}
+c^{4}_{13}
+c^{5}_{13}
$,
$ \xi_{13}^{7}$,
$ \xi_{13}^{7}$,
$ 2+2c^{1}_{13}
+c^{2}_{13}
+c^{3}_{13}
+c^{4}_{13}
+c^{5}_{13}
$,
$ 2+2c^{1}_{13}
+c^{2}_{13}
+c^{3}_{13}
+c^{4}_{13}
$,
$ 2+2c^{1}_{13}
+2c^{2}_{13}
+c^{3}_{13}
+c^{4}_{13}
$,
$ 2+2c^{1}_{13}
+2c^{2}_{13}
+c^{3}_{13}
$;\ \ 
$ 2+2c^{1}_{13}
+2c^{2}_{13}
+c^{3}_{13}
+c^{4}_{13}
$,
$ -\xi_{13}^{7}$,
$ -\xi_{13}^{7}$,
$ 2+2c^{1}_{13}
+2c^{2}_{13}
+c^{3}_{13}
$,
$ 2+2c^{1}_{13}
+c^{2}_{13}
+c^{3}_{13}
+c^{4}_{13}
+c^{5}_{13}
$,
$ -1$,
$ -2-2  c^{1}_{13}
-c^{2}_{13}
-c^{3}_{13}
-c^{4}_{13}
$;\ \ 
$ -1-c^{1}_{13}
-c^{2}_{13}
+c^{5}_{13}
$,
$ 2+2c^{1}_{13}
+2c^{2}_{13}
+c^{3}_{13}
-c^{5}_{13}
$,
$ \xi_{13}^{7}$,
$ -\xi_{13}^{7}$,
$ \xi_{13}^{7}$,
$ -\xi_{13}^{7}$;\ \ 
$ -1-c^{1}_{13}
-c^{2}_{13}
+c^{5}_{13}
$,
$ \xi_{13}^{7}$,
$ -\xi_{13}^{7}$,
$ \xi_{13}^{7}$,
$ -\xi_{13}^{7}$;\ \ 
$ 1$,
$ -2-2  c^{1}_{13}
-2  c^{2}_{13}
-c^{3}_{13}
-c^{4}_{13}
$,
$ -2-2  c^{1}_{13}
-c^{2}_{13}
-c^{3}_{13}
-c^{4}_{13}
$,
$ 2+c^{1}_{13}
+c^{2}_{13}
+c^{3}_{13}
+c^{4}_{13}
+c^{5}_{13}
$;\ \ 
$ -2-c^{1}_{13}
-c^{2}_{13}
-c^{3}_{13}
-c^{4}_{13}
-c^{5}_{13}
$,
$ 2+2c^{1}_{13}
+2c^{2}_{13}
+c^{3}_{13}
$,
$ 1$;\ \ 
$ -2-c^{1}_{13}
-c^{2}_{13}
-c^{3}_{13}
-c^{4}_{13}
-c^{5}_{13}
$,
$ -2-2  c^{1}_{13}
-c^{2}_{13}
-c^{3}_{13}
-c^{4}_{13}
-c^{5}_{13}
$;\ \ 
$ 2+2c^{1}_{13}
+2c^{2}_{13}
+c^{3}_{13}
+c^{4}_{13}
$)

  \vskip 2ex

\noindent60. $8_{\frac{68}{13},223.6}^{13,484}$ \irep{179}:\ \ 
$d_i$ = ($1.0$,
$2.941$,
$4.148$,
$4.148$,
$4.712$,
$6.209$,
$7.345$,
$8.55$) 

\vskip 0.7ex
\hangindent=3em \hangafter=1
$D^2= 223.689 = 
78+65c^{1}_{13}
+52c^{2}_{13}
+39c^{3}_{13}
+26c^{4}_{13}
+13c^{5}_{13}
$

\vskip 0.7ex
\hangindent=3em \hangafter=1
$T = ( 0,
\frac{12}{13},
\frac{5}{13},
\frac{5}{13},
\frac{10}{13},
\frac{7}{13},
\frac{3}{13},
\frac{11}{13} )
$,

\vskip 0.7ex
\hangindent=3em \hangafter=1
$S$ = ($ 1$,
$ 2+c^{1}_{13}
+c^{2}_{13}
+c^{3}_{13}
+c^{4}_{13}
+c^{5}_{13}
$,
$ \xi_{13}^{7}$,
$ \xi_{13}^{7}$,
$ 2+2c^{1}_{13}
+c^{2}_{13}
+c^{3}_{13}
+c^{4}_{13}
+c^{5}_{13}
$,
$ 2+2c^{1}_{13}
+c^{2}_{13}
+c^{3}_{13}
+c^{4}_{13}
$,
$ 2+2c^{1}_{13}
+2c^{2}_{13}
+c^{3}_{13}
+c^{4}_{13}
$,
$ 2+2c^{1}_{13}
+2c^{2}_{13}
+c^{3}_{13}
$;\ \ 
$ 2+2c^{1}_{13}
+2c^{2}_{13}
+c^{3}_{13}
+c^{4}_{13}
$,
$ -\xi_{13}^{7}$,
$ -\xi_{13}^{7}$,
$ 2+2c^{1}_{13}
+2c^{2}_{13}
+c^{3}_{13}
$,
$ 2+2c^{1}_{13}
+c^{2}_{13}
+c^{3}_{13}
+c^{4}_{13}
+c^{5}_{13}
$,
$ -1$,
$ -2-2  c^{1}_{13}
-c^{2}_{13}
-c^{3}_{13}
-c^{4}_{13}
$;\ \ 
$ -1-c^{1}_{13}
-c^{2}_{13}
+c^{5}_{13}
$,
$ 2+2c^{1}_{13}
+2c^{2}_{13}
+c^{3}_{13}
-c^{5}_{13}
$,
$ \xi_{13}^{7}$,
$ -\xi_{13}^{7}$,
$ \xi_{13}^{7}$,
$ -\xi_{13}^{7}$;\ \ 
$ -1-c^{1}_{13}
-c^{2}_{13}
+c^{5}_{13}
$,
$ \xi_{13}^{7}$,
$ -\xi_{13}^{7}$,
$ \xi_{13}^{7}$,
$ -\xi_{13}^{7}$;\ \ 
$ 1$,
$ -2-2  c^{1}_{13}
-2  c^{2}_{13}
-c^{3}_{13}
-c^{4}_{13}
$,
$ -2-2  c^{1}_{13}
-c^{2}_{13}
-c^{3}_{13}
-c^{4}_{13}
$,
$ 2+c^{1}_{13}
+c^{2}_{13}
+c^{3}_{13}
+c^{4}_{13}
+c^{5}_{13}
$;\ \ 
$ -2-c^{1}_{13}
-c^{2}_{13}
-c^{3}_{13}
-c^{4}_{13}
-c^{5}_{13}
$,
$ 2+2c^{1}_{13}
+2c^{2}_{13}
+c^{3}_{13}
$,
$ 1$;\ \ 
$ -2-c^{1}_{13}
-c^{2}_{13}
-c^{3}_{13}
-c^{4}_{13}
-c^{5}_{13}
$,
$ -2-2  c^{1}_{13}
-c^{2}_{13}
-c^{3}_{13}
-c^{4}_{13}
-c^{5}_{13}
$;\ \ 
$ 2+2c^{1}_{13}
+2c^{2}_{13}
+c^{3}_{13}
+c^{4}_{13}
$)

  \vskip 2ex

\noindent61. $8_{4,308.4}^{15,440}$ \irep{183}:\ \ 
$d_i$ = ($1.0$,
$5.854$,
$5.854$,
$5.854$,
$5.854$,
$6.854$,
$7.854$,
$7.854$) 

\vskip 0.7ex
\hangindent=3em \hangafter=1
$D^2= 308.434 = 
\frac{315+135\sqrt{5}}{2}$

\vskip 0.7ex
\hangindent=3em \hangafter=1
$T = ( 0,
0,
0,
\frac{1}{3},
\frac{2}{3},
0,
\frac{2}{5},
\frac{3}{5} )
$,

\vskip 0.7ex
\hangindent=3em \hangafter=1
$S$ = ($ 1$,
$ \frac{5+3\sqrt{5}}{2}$,
$ \frac{5+3\sqrt{5}}{2}$,
$ \frac{5+3\sqrt{5}}{2}$,
$ \frac{5+3\sqrt{5}}{2}$,
$ \frac{7+3\sqrt{5}}{2}$,
$ \frac{9+3\sqrt{5}}{2}$,
$ \frac{9+3\sqrt{5}}{2}$;\ \ 
$ -5-3\sqrt{5}$,
$ \frac{5+3\sqrt{5}}{2}$,
$ \frac{5+3\sqrt{5}}{2}$,
$ \frac{5+3\sqrt{5}}{2}$,
$ -\frac{5+3\sqrt{5}}{2}$,
$0$,
$0$;\ \ 
$ -5-3\sqrt{5}$,
$ \frac{5+3\sqrt{5}}{2}$,
$ \frac{5+3\sqrt{5}}{2}$,
$ -\frac{5+3\sqrt{5}}{2}$,
$0$,
$0$;\ \ 
$ \frac{5+3\sqrt{5}}{2}$,
$ -5-3\sqrt{5}$,
$ -\frac{5+3\sqrt{5}}{2}$,
$0$,
$0$;\ \ 
$ \frac{5+3\sqrt{5}}{2}$,
$ -\frac{5+3\sqrt{5}}{2}$,
$0$,
$0$;\ \ 
$ 1$,
$ \frac{9+3\sqrt{5}}{2}$,
$ \frac{9+3\sqrt{5}}{2}$;\ \ 
$ \frac{3+3\sqrt{5}}{2}$,
$ -6-3\sqrt{5}$;\ \ 
$ \frac{3+3\sqrt{5}}{2}$)

  \vskip 2ex

\noindent62. $8_{0,308.4}^{15,100}$ \irep{188}:\ \ 
$d_i$ = ($1.0$,
$5.854$,
$5.854$,
$5.854$,
$5.854$,
$6.854$,
$7.854$,
$7.854$) 

\vskip 0.7ex
\hangindent=3em \hangafter=1
$D^2= 308.434 = 
\frac{315+135\sqrt{5}}{2}$

\vskip 0.7ex
\hangindent=3em \hangafter=1
$T = ( 0,
\frac{1}{3},
\frac{1}{3},
\frac{2}{3},
\frac{2}{3},
0,
\frac{1}{5},
\frac{4}{5} )
$,

\vskip 0.7ex
\hangindent=3em \hangafter=1
$S$ = ($ 1$,
$ \frac{5+3\sqrt{5}}{2}$,
$ \frac{5+3\sqrt{5}}{2}$,
$ \frac{5+3\sqrt{5}}{2}$,
$ \frac{5+3\sqrt{5}}{2}$,
$ \frac{7+3\sqrt{5}}{2}$,
$ \frac{9+3\sqrt{5}}{2}$,
$ \frac{9+3\sqrt{5}}{2}$;\ \ 
$ \frac{5+3\sqrt{5}}{2}$,
$ -5-3\sqrt{5}$,
$ \frac{5+3\sqrt{5}}{2}$,
$ \frac{5+3\sqrt{5}}{2}$,
$ -\frac{5+3\sqrt{5}}{2}$,
$0$,
$0$;\ \ 
$ \frac{5+3\sqrt{5}}{2}$,
$ \frac{5+3\sqrt{5}}{2}$,
$ \frac{5+3\sqrt{5}}{2}$,
$ -\frac{5+3\sqrt{5}}{2}$,
$0$,
$0$;\ \ 
$ \frac{5+3\sqrt{5}}{2}$,
$ -5-3\sqrt{5}$,
$ -\frac{5+3\sqrt{5}}{2}$,
$0$,
$0$;\ \ 
$ \frac{5+3\sqrt{5}}{2}$,
$ -\frac{5+3\sqrt{5}}{2}$,
$0$,
$0$;\ \ 
$ 1$,
$ \frac{9+3\sqrt{5}}{2}$,
$ \frac{9+3\sqrt{5}}{2}$;\ \ 
$ -6-3\sqrt{5}$,
$ \frac{3+3\sqrt{5}}{2}$;\ \ 
$ -6-3\sqrt{5}$)

  \vskip 2ex

\noindent63. $8_{4,308.4}^{45,289}$ \irep{235}:\ \ 
$d_i$ = ($1.0$,
$5.854$,
$5.854$,
$5.854$,
$5.854$,
$6.854$,
$7.854$,
$7.854$) 

\vskip 0.7ex
\hangindent=3em \hangafter=1
$D^2= 308.434 = 
\frac{315+135\sqrt{5}}{2}$

\vskip 0.7ex
\hangindent=3em \hangafter=1
$T = ( 0,
0,
\frac{1}{9},
\frac{4}{9},
\frac{7}{9},
0,
\frac{2}{5},
\frac{3}{5} )
$,

\vskip 0.7ex
\hangindent=3em \hangafter=1
$S$ = ($ 1$,
$ \frac{5+3\sqrt{5}}{2}$,
$ \frac{5+3\sqrt{5}}{2}$,
$ \frac{5+3\sqrt{5}}{2}$,
$ \frac{5+3\sqrt{5}}{2}$,
$ \frac{7+3\sqrt{5}}{2}$,
$ \frac{9+3\sqrt{5}}{2}$,
$ \frac{9+3\sqrt{5}}{2}$;\ \ 
$ -5-3\sqrt{5}$,
$ \frac{5+3\sqrt{5}}{2}$,
$ \frac{5+3\sqrt{5}}{2}$,
$ \frac{5+3\sqrt{5}}{2}$,
$ -\frac{5+3\sqrt{5}}{2}$,
$0$,
$0$;\ \ 
$ -3  c^{1}_{45}
+3c^{4}_{45}
-4  c^{10}_{45}
+3c^{11}_{45}
$,
$ 3c^{2}_{45}
+c^{5}_{45}
+3c^{7}_{45}
+c^{10}_{45}
$,
$ 3c^{1}_{45}
-3  c^{2}_{45}
-3  c^{4}_{45}
-c^{5}_{45}
-3  c^{7}_{45}
+3c^{10}_{45}
-3  c^{11}_{45}
$,
$ -\frac{5+3\sqrt{5}}{2}$,
$0$,
$0$;\ \ 
$ 3c^{1}_{45}
-3  c^{2}_{45}
-3  c^{4}_{45}
-c^{5}_{45}
-3  c^{7}_{45}
+3c^{10}_{45}
-3  c^{11}_{45}
$,
$ -3  c^{1}_{45}
+3c^{4}_{45}
-4  c^{10}_{45}
+3c^{11}_{45}
$,
$ -\frac{5+3\sqrt{5}}{2}$,
$0$,
$0$;\ \ 
$ 3c^{2}_{45}
+c^{5}_{45}
+3c^{7}_{45}
+c^{10}_{45}
$,
$ -\frac{5+3\sqrt{5}}{2}$,
$0$,
$0$;\ \ 
$ 1$,
$ \frac{9+3\sqrt{5}}{2}$,
$ \frac{9+3\sqrt{5}}{2}$;\ \ 
$ \frac{3+3\sqrt{5}}{2}$,
$ -6-3\sqrt{5}$;\ \ 
$ \frac{3+3\sqrt{5}}{2}$)

  \vskip 2ex

\noindent64. $8_{4,308.4}^{45,939}$ \irep{235}:\ \ 
$d_i$ = ($1.0$,
$5.854$,
$5.854$,
$5.854$,
$5.854$,
$6.854$,
$7.854$,
$7.854$) 

\vskip 0.7ex
\hangindent=3em \hangafter=1
$D^2= 308.434 = 
\frac{315+135\sqrt{5}}{2}$

\vskip 0.7ex
\hangindent=3em \hangafter=1
$T = ( 0,
0,
\frac{2}{9},
\frac{5}{9},
\frac{8}{9},
0,
\frac{2}{5},
\frac{3}{5} )
$,

\vskip 0.7ex
\hangindent=3em \hangafter=1
$S$ = ($ 1$,
$ \frac{5+3\sqrt{5}}{2}$,
$ \frac{5+3\sqrt{5}}{2}$,
$ \frac{5+3\sqrt{5}}{2}$,
$ \frac{5+3\sqrt{5}}{2}$,
$ \frac{7+3\sqrt{5}}{2}$,
$ \frac{9+3\sqrt{5}}{2}$,
$ \frac{9+3\sqrt{5}}{2}$;\ \ 
$ -5-3\sqrt{5}$,
$ \frac{5+3\sqrt{5}}{2}$,
$ \frac{5+3\sqrt{5}}{2}$,
$ \frac{5+3\sqrt{5}}{2}$,
$ -\frac{5+3\sqrt{5}}{2}$,
$0$,
$0$;\ \ 
$ 3c^{2}_{45}
+c^{5}_{45}
+3c^{7}_{45}
+c^{10}_{45}
$,
$ -3  c^{1}_{45}
+3c^{4}_{45}
-4  c^{10}_{45}
+3c^{11}_{45}
$,
$ 3c^{1}_{45}
-3  c^{2}_{45}
-3  c^{4}_{45}
-c^{5}_{45}
-3  c^{7}_{45}
+3c^{10}_{45}
-3  c^{11}_{45}
$,
$ -\frac{5+3\sqrt{5}}{2}$,
$0$,
$0$;\ \ 
$ 3c^{1}_{45}
-3  c^{2}_{45}
-3  c^{4}_{45}
-c^{5}_{45}
-3  c^{7}_{45}
+3c^{10}_{45}
-3  c^{11}_{45}
$,
$ 3c^{2}_{45}
+c^{5}_{45}
+3c^{7}_{45}
+c^{10}_{45}
$,
$ -\frac{5+3\sqrt{5}}{2}$,
$0$,
$0$;\ \ 
$ -3  c^{1}_{45}
+3c^{4}_{45}
-4  c^{10}_{45}
+3c^{11}_{45}
$,
$ -\frac{5+3\sqrt{5}}{2}$,
$0$,
$0$;\ \ 
$ 1$,
$ \frac{9+3\sqrt{5}}{2}$,
$ \frac{9+3\sqrt{5}}{2}$;\ \ 
$ \frac{3+3\sqrt{5}}{2}$,
$ -6-3\sqrt{5}$;\ \ 
$ \frac{3+3\sqrt{5}}{2}$)

  \vskip 2ex 

}

\subsection{Rank 9 }
\label{uni9}

{\small

\noindent1. $9_{0,9.}^{3,113}$ \irep{0}:\ \ 
$d_i$ = ($1.0$,
$1.0$,
$1.0$,
$1.0$,
$1.0$,
$1.0$,
$1.0$,
$1.0$,
$1.0$) 

\vskip 0.7ex
\hangindent=3em \hangafter=1
$D^2= 9.0 = 
9$

\vskip 0.7ex
\hangindent=3em \hangafter=1
$T = ( 0,
0,
0,
0,
0,
\frac{1}{3},
\frac{1}{3},
\frac{2}{3},
\frac{2}{3} )
$,

\vskip 0.7ex
\hangindent=3em \hangafter=1
$S$ = ($ 1$,
$ 1$,
$ 1$,
$ 1$,
$ 1$,
$ 1$,
$ 1$,
$ 1$,
$ 1$;\ \ 
$ 1$,
$ 1$,
$ -\zeta_{6}^{1}$,
$ -\zeta_{6}^{5}$,
$ -\zeta_{6}^{1}$,
$ -\zeta_{6}^{5}$,
$ -\zeta_{6}^{1}$,
$ -\zeta_{6}^{5}$;\ \ 
$ 1$,
$ -\zeta_{6}^{5}$,
$ -\zeta_{6}^{1}$,
$ -\zeta_{6}^{5}$,
$ -\zeta_{6}^{1}$,
$ -\zeta_{6}^{5}$,
$ -\zeta_{6}^{1}$;\ \ 
$ 1$,
$ 1$,
$ -\zeta_{6}^{1}$,
$ -\zeta_{6}^{5}$,
$ -\zeta_{6}^{5}$,
$ -\zeta_{6}^{1}$;\ \ 
$ 1$,
$ -\zeta_{6}^{5}$,
$ -\zeta_{6}^{1}$,
$ -\zeta_{6}^{1}$,
$ -\zeta_{6}^{5}$;\ \ 
$ -\zeta_{6}^{5}$,
$ -\zeta_{6}^{1}$,
$ 1$,
$ 1$;\ \ 
$ -\zeta_{6}^{5}$,
$ 1$,
$ 1$;\ \ 
$ -\zeta_{6}^{1}$,
$ -\zeta_{6}^{5}$;\ \ 
$ -\zeta_{6}^{1}$)

Factors = $3_{2,3.}^{3,527}\boxtimes 3_{6,3.}^{3,138}$

  \vskip 2ex

\noindent2. $9_{4,9.}^{3,277}$ \irep{0}:\ \ 
$d_i$ = ($1.0$,
$1.0$,
$1.0$,
$1.0$,
$1.0$,
$1.0$,
$1.0$,
$1.0$,
$1.0$) 

\vskip 0.7ex
\hangindent=3em \hangafter=1
$D^2= 9.0 = 
9$

\vskip 0.7ex
\hangindent=3em \hangafter=1
$T = ( 0,
\frac{1}{3},
\frac{1}{3},
\frac{1}{3},
\frac{1}{3},
\frac{2}{3},
\frac{2}{3},
\frac{2}{3},
\frac{2}{3} )
$,

\vskip 0.7ex
\hangindent=3em \hangafter=1
$S$ = ($ 1$,
$ 1$,
$ 1$,
$ 1$,
$ 1$,
$ 1$,
$ 1$,
$ 1$,
$ 1$;\ \ 
$ -\zeta_{6}^{5}$,
$ 1$,
$ 1$,
$ -\zeta_{6}^{1}$,
$ -\zeta_{6}^{1}$,
$ -\zeta_{6}^{5}$,
$ -\zeta_{6}^{1}$,
$ -\zeta_{6}^{5}$;\ \ 
$ -\zeta_{6}^{5}$,
$ -\zeta_{6}^{1}$,
$ 1$,
$ -\zeta_{6}^{5}$,
$ -\zeta_{6}^{1}$,
$ -\zeta_{6}^{1}$,
$ -\zeta_{6}^{5}$;\ \ 
$ -\zeta_{6}^{5}$,
$ 1$,
$ -\zeta_{6}^{1}$,
$ -\zeta_{6}^{5}$,
$ -\zeta_{6}^{5}$,
$ -\zeta_{6}^{1}$;\ \ 
$ -\zeta_{6}^{5}$,
$ -\zeta_{6}^{5}$,
$ -\zeta_{6}^{1}$,
$ -\zeta_{6}^{5}$,
$ -\zeta_{6}^{1}$;\ \ 
$ -\zeta_{6}^{1}$,
$ -\zeta_{6}^{5}$,
$ 1$,
$ 1$;\ \ 
$ -\zeta_{6}^{1}$,
$ 1$,
$ 1$;\ \ 
$ -\zeta_{6}^{1}$,
$ -\zeta_{6}^{5}$;\ \ 
$ -\zeta_{6}^{1}$)

Factors = $3_{2,3.}^{3,527}\boxtimes 3_{2,3.}^{3,527}$

  \vskip 2ex

\noindent3. $9_{0,9.}^{9,620}$ \irep{0}:\ \ 
$d_i$ = ($1.0$,
$1.0$,
$1.0$,
$1.0$,
$1.0$,
$1.0$,
$1.0$,
$1.0$,
$1.0$) 

\vskip 0.7ex
\hangindent=3em \hangafter=1
$D^2= 9.0 = 
9$

\vskip 0.7ex
\hangindent=3em \hangafter=1
$T = ( 0,
0,
0,
\frac{1}{9},
\frac{1}{9},
\frac{4}{9},
\frac{4}{9},
\frac{7}{9},
\frac{7}{9} )
$,

\vskip 0.7ex
\hangindent=3em \hangafter=1
$S$ = ($ 1$,
$ 1$,
$ 1$,
$ 1$,
$ 1$,
$ 1$,
$ 1$,
$ 1$,
$ 1$;\ \ 
$ 1$,
$ 1$,
$ -\zeta_{6}^{1}$,
$ -\zeta_{6}^{5}$,
$ -\zeta_{6}^{1}$,
$ -\zeta_{6}^{5}$,
$ -\zeta_{6}^{1}$,
$ -\zeta_{6}^{5}$;\ \ 
$ 1$,
$ -\zeta_{6}^{5}$,
$ -\zeta_{6}^{1}$,
$ -\zeta_{6}^{5}$,
$ -\zeta_{6}^{1}$,
$ -\zeta_{6}^{5}$,
$ -\zeta_{6}^{1}$;\ \ 
$ -\zeta_{18}^{5}$,
$ -\zeta_{18}^{13}$,
$ -\zeta_{18}^{17}$,
$ -\zeta_{18}^{1}$,
$ -\zeta_{18}^{11}$,
$ -\zeta_{18}^{7}$;\ \ 
$ -\zeta_{18}^{5}$,
$ -\zeta_{18}^{1}$,
$ -\zeta_{18}^{17}$,
$ -\zeta_{18}^{7}$,
$ -\zeta_{18}^{11}$;\ \ 
$ -\zeta_{18}^{11}$,
$ -\zeta_{18}^{7}$,
$ -\zeta_{18}^{5}$,
$ -\zeta_{18}^{13}$;\ \ 
$ -\zeta_{18}^{11}$,
$ -\zeta_{18}^{13}$,
$ -\zeta_{18}^{5}$;\ \ 
$ -\zeta_{18}^{17}$,
$ -\zeta_{18}^{1}$;\ \ 
$ -\zeta_{18}^{17}$)

  \vskip 2ex

\noindent4. $9_{0,9.}^{9,462}$ \irep{0}:\ \ 
$d_i$ = ($1.0$,
$1.0$,
$1.0$,
$1.0$,
$1.0$,
$1.0$,
$1.0$,
$1.0$,
$1.0$) 

\vskip 0.7ex
\hangindent=3em \hangafter=1
$D^2= 9.0 = 
9$

\vskip 0.7ex
\hangindent=3em \hangafter=1
$T = ( 0,
0,
0,
\frac{2}{9},
\frac{2}{9},
\frac{5}{9},
\frac{5}{9},
\frac{8}{9},
\frac{8}{9} )
$,

\vskip 0.7ex
\hangindent=3em \hangafter=1
$S$ = ($ 1$,
$ 1$,
$ 1$,
$ 1$,
$ 1$,
$ 1$,
$ 1$,
$ 1$,
$ 1$;\ \ 
$ 1$,
$ 1$,
$ -\zeta_{6}^{1}$,
$ -\zeta_{6}^{5}$,
$ -\zeta_{6}^{1}$,
$ -\zeta_{6}^{5}$,
$ -\zeta_{6}^{1}$,
$ -\zeta_{6}^{5}$;\ \ 
$ 1$,
$ -\zeta_{6}^{5}$,
$ -\zeta_{6}^{1}$,
$ -\zeta_{6}^{5}$,
$ -\zeta_{6}^{1}$,
$ -\zeta_{6}^{5}$,
$ -\zeta_{6}^{1}$;\ \ 
$ -\zeta_{18}^{1}$,
$ -\zeta_{18}^{17}$,
$ -\zeta_{18}^{13}$,
$ -\zeta_{18}^{5}$,
$ -\zeta_{18}^{7}$,
$ -\zeta_{18}^{11}$;\ \ 
$ -\zeta_{18}^{1}$,
$ -\zeta_{18}^{5}$,
$ -\zeta_{18}^{13}$,
$ -\zeta_{18}^{11}$,
$ -\zeta_{18}^{7}$;\ \ 
$ -\zeta_{18}^{7}$,
$ -\zeta_{18}^{11}$,
$ -\zeta_{18}^{1}$,
$ -\zeta_{18}^{17}$;\ \ 
$ -\zeta_{18}^{7}$,
$ -\zeta_{18}^{17}$,
$ -\zeta_{18}^{1}$;\ \ 
$ -\zeta_{18}^{13}$,
$ -\zeta_{18}^{5}$;\ \ 
$ -\zeta_{18}^{13}$)

  \vskip 2ex

\noindent5. $9_{\frac{5}{2},12.}^{48,250}$ \irep{524}:\ \ 
$d_i$ = ($1.0$,
$1.0$,
$1.0$,
$1.0$,
$1.0$,
$1.0$,
$1.414$,
$1.414$,
$1.414$) 

\vskip 0.7ex
\hangindent=3em \hangafter=1
$D^2= 12.0 = 
12$

\vskip 0.7ex
\hangindent=3em \hangafter=1
$T = ( 0,
\frac{1}{2},
\frac{1}{3},
\frac{1}{3},
\frac{5}{6},
\frac{5}{6},
\frac{1}{16},
\frac{19}{48},
\frac{19}{48} )
$,

\vskip 0.7ex
\hangindent=3em \hangafter=1
$S$ = ($ 1$,
$ 1$,
$ 1$,
$ 1$,
$ 1$,
$ 1$,
$ \sqrt{2}$,
$ \sqrt{2}$,
$ \sqrt{2}$;\ \ 
$ 1$,
$ 1$,
$ 1$,
$ 1$,
$ 1$,
$ -\sqrt{2}$,
$ -\sqrt{2}$,
$ -\sqrt{2}$;\ \ 
$ -\zeta_{6}^{5}$,
$ -\zeta_{6}^{1}$,
$ -\zeta_{6}^{1}$,
$ -\zeta_{6}^{5}$,
$ \sqrt{2}$,
$ -\sqrt{2}\zeta_{6}^{1}$,
$ \sqrt{2}\zeta_{3}^{1}$;\ \ 
$ -\zeta_{6}^{5}$,
$ -\zeta_{6}^{5}$,
$ -\zeta_{6}^{1}$,
$ \sqrt{2}$,
$ \sqrt{2}\zeta_{3}^{1}$,
$ -\sqrt{2}\zeta_{6}^{1}$;\ \ 
$ -\zeta_{6}^{5}$,
$ -\zeta_{6}^{1}$,
$ -\sqrt{2}$,
$ -\sqrt{2}\zeta_{3}^{1}$,
$ \sqrt{2}\zeta_{6}^{1}$;\ \ 
$ -\zeta_{6}^{5}$,
$ -\sqrt{2}$,
$ \sqrt{2}\zeta_{6}^{1}$,
$ -\sqrt{2}\zeta_{3}^{1}$;\ \ 
$0$,
$0$,
$0$;\ \ 
$0$,
$0$;\ \ 
$0$)

Factors = $3_{2,3.}^{3,527}\boxtimes 3_{\frac{1}{2},4.}^{16,598}$

  \vskip 2ex

\noindent6. $9_{\frac{7}{2},12.}^{48,117}$ \irep{524}:\ \ 
$d_i$ = ($1.0$,
$1.0$,
$1.0$,
$1.0$,
$1.0$,
$1.0$,
$1.414$,
$1.414$,
$1.414$) 

\vskip 0.7ex
\hangindent=3em \hangafter=1
$D^2= 12.0 = 
12$

\vskip 0.7ex
\hangindent=3em \hangafter=1
$T = ( 0,
\frac{1}{2},
\frac{1}{3},
\frac{1}{3},
\frac{5}{6},
\frac{5}{6},
\frac{3}{16},
\frac{25}{48},
\frac{25}{48} )
$,

\vskip 0.7ex
\hangindent=3em \hangafter=1
$S$ = ($ 1$,
$ 1$,
$ 1$,
$ 1$,
$ 1$,
$ 1$,
$ \sqrt{2}$,
$ \sqrt{2}$,
$ \sqrt{2}$;\ \ 
$ 1$,
$ 1$,
$ 1$,
$ 1$,
$ 1$,
$ -\sqrt{2}$,
$ -\sqrt{2}$,
$ -\sqrt{2}$;\ \ 
$ -\zeta_{6}^{5}$,
$ -\zeta_{6}^{1}$,
$ -\zeta_{6}^{1}$,
$ -\zeta_{6}^{5}$,
$ \sqrt{2}$,
$ -\sqrt{2}\zeta_{6}^{1}$,
$ \sqrt{2}\zeta_{3}^{1}$;\ \ 
$ -\zeta_{6}^{5}$,
$ -\zeta_{6}^{5}$,
$ -\zeta_{6}^{1}$,
$ \sqrt{2}$,
$ \sqrt{2}\zeta_{3}^{1}$,
$ -\sqrt{2}\zeta_{6}^{1}$;\ \ 
$ -\zeta_{6}^{5}$,
$ -\zeta_{6}^{1}$,
$ -\sqrt{2}$,
$ -\sqrt{2}\zeta_{3}^{1}$,
$ \sqrt{2}\zeta_{6}^{1}$;\ \ 
$ -\zeta_{6}^{5}$,
$ -\sqrt{2}$,
$ \sqrt{2}\zeta_{6}^{1}$,
$ -\sqrt{2}\zeta_{3}^{1}$;\ \ 
$0$,
$0$,
$0$;\ \ 
$0$,
$0$;\ \ 
$0$)

Factors = $3_{2,3.}^{3,527}\boxtimes 3_{\frac{3}{2},4.}^{16,553}$

  \vskip 2ex

\noindent7. $9_{\frac{9}{2},12.}^{48,618}$ \irep{524}:\ \ 
$d_i$ = ($1.0$,
$1.0$,
$1.0$,
$1.0$,
$1.0$,
$1.0$,
$1.414$,
$1.414$,
$1.414$) 

\vskip 0.7ex
\hangindent=3em \hangafter=1
$D^2= 12.0 = 
12$

\vskip 0.7ex
\hangindent=3em \hangafter=1
$T = ( 0,
\frac{1}{2},
\frac{1}{3},
\frac{1}{3},
\frac{5}{6},
\frac{5}{6},
\frac{5}{16},
\frac{31}{48},
\frac{31}{48} )
$,

\vskip 0.7ex
\hangindent=3em \hangafter=1
$S$ = ($ 1$,
$ 1$,
$ 1$,
$ 1$,
$ 1$,
$ 1$,
$ \sqrt{2}$,
$ \sqrt{2}$,
$ \sqrt{2}$;\ \ 
$ 1$,
$ 1$,
$ 1$,
$ 1$,
$ 1$,
$ -\sqrt{2}$,
$ -\sqrt{2}$,
$ -\sqrt{2}$;\ \ 
$ -\zeta_{6}^{5}$,
$ -\zeta_{6}^{1}$,
$ -\zeta_{6}^{1}$,
$ -\zeta_{6}^{5}$,
$ \sqrt{2}$,
$ -\sqrt{2}\zeta_{6}^{1}$,
$ \sqrt{2}\zeta_{3}^{1}$;\ \ 
$ -\zeta_{6}^{5}$,
$ -\zeta_{6}^{5}$,
$ -\zeta_{6}^{1}$,
$ \sqrt{2}$,
$ \sqrt{2}\zeta_{3}^{1}$,
$ -\sqrt{2}\zeta_{6}^{1}$;\ \ 
$ -\zeta_{6}^{5}$,
$ -\zeta_{6}^{1}$,
$ -\sqrt{2}$,
$ -\sqrt{2}\zeta_{3}^{1}$,
$ \sqrt{2}\zeta_{6}^{1}$;\ \ 
$ -\zeta_{6}^{5}$,
$ -\sqrt{2}$,
$ \sqrt{2}\zeta_{6}^{1}$,
$ -\sqrt{2}\zeta_{3}^{1}$;\ \ 
$0$,
$0$,
$0$;\ \ 
$0$,
$0$;\ \ 
$0$)

Factors = $3_{2,3.}^{3,527}\boxtimes 3_{\frac{5}{2},4.}^{16,465}$

  \vskip 2ex

\noindent8. $9_{\frac{11}{2},12.}^{48,125}$ \irep{524}:\ \ 
$d_i$ = ($1.0$,
$1.0$,
$1.0$,
$1.0$,
$1.0$,
$1.0$,
$1.414$,
$1.414$,
$1.414$) 

\vskip 0.7ex
\hangindent=3em \hangafter=1
$D^2= 12.0 = 
12$

\vskip 0.7ex
\hangindent=3em \hangafter=1
$T = ( 0,
\frac{1}{2},
\frac{1}{3},
\frac{1}{3},
\frac{5}{6},
\frac{5}{6},
\frac{7}{16},
\frac{37}{48},
\frac{37}{48} )
$,

\vskip 0.7ex
\hangindent=3em \hangafter=1
$S$ = ($ 1$,
$ 1$,
$ 1$,
$ 1$,
$ 1$,
$ 1$,
$ \sqrt{2}$,
$ \sqrt{2}$,
$ \sqrt{2}$;\ \ 
$ 1$,
$ 1$,
$ 1$,
$ 1$,
$ 1$,
$ -\sqrt{2}$,
$ -\sqrt{2}$,
$ -\sqrt{2}$;\ \ 
$ -\zeta_{6}^{5}$,
$ -\zeta_{6}^{1}$,
$ -\zeta_{6}^{1}$,
$ -\zeta_{6}^{5}$,
$ \sqrt{2}$,
$ -\sqrt{2}\zeta_{6}^{1}$,
$ \sqrt{2}\zeta_{3}^{1}$;\ \ 
$ -\zeta_{6}^{5}$,
$ -\zeta_{6}^{5}$,
$ -\zeta_{6}^{1}$,
$ \sqrt{2}$,
$ \sqrt{2}\zeta_{3}^{1}$,
$ -\sqrt{2}\zeta_{6}^{1}$;\ \ 
$ -\zeta_{6}^{5}$,
$ -\zeta_{6}^{1}$,
$ -\sqrt{2}$,
$ -\sqrt{2}\zeta_{3}^{1}$,
$ \sqrt{2}\zeta_{6}^{1}$;\ \ 
$ -\zeta_{6}^{5}$,
$ -\sqrt{2}$,
$ \sqrt{2}\zeta_{6}^{1}$,
$ -\sqrt{2}\zeta_{3}^{1}$;\ \ 
$0$,
$0$,
$0$;\ \ 
$0$,
$0$;\ \ 
$0$)

Factors = $3_{2,3.}^{3,527}\boxtimes 3_{\frac{7}{2},4.}^{16,332}$

  \vskip 2ex

\noindent9. $9_{\frac{13}{2},12.}^{48,201}$ \irep{524}:\ \ 
$d_i$ = ($1.0$,
$1.0$,
$1.0$,
$1.0$,
$1.0$,
$1.0$,
$1.414$,
$1.414$,
$1.414$) 

\vskip 0.7ex
\hangindent=3em \hangafter=1
$D^2= 12.0 = 
12$

\vskip 0.7ex
\hangindent=3em \hangafter=1
$T = ( 0,
\frac{1}{2},
\frac{1}{3},
\frac{1}{3},
\frac{5}{6},
\frac{5}{6},
\frac{9}{16},
\frac{43}{48},
\frac{43}{48} )
$,

\vskip 0.7ex
\hangindent=3em \hangafter=1
$S$ = ($ 1$,
$ 1$,
$ 1$,
$ 1$,
$ 1$,
$ 1$,
$ \sqrt{2}$,
$ \sqrt{2}$,
$ \sqrt{2}$;\ \ 
$ 1$,
$ 1$,
$ 1$,
$ 1$,
$ 1$,
$ -\sqrt{2}$,
$ -\sqrt{2}$,
$ -\sqrt{2}$;\ \ 
$ -\zeta_{6}^{5}$,
$ -\zeta_{6}^{1}$,
$ -\zeta_{6}^{1}$,
$ -\zeta_{6}^{5}$,
$ \sqrt{2}$,
$ -\sqrt{2}\zeta_{6}^{1}$,
$ \sqrt{2}\zeta_{3}^{1}$;\ \ 
$ -\zeta_{6}^{5}$,
$ -\zeta_{6}^{5}$,
$ -\zeta_{6}^{1}$,
$ \sqrt{2}$,
$ \sqrt{2}\zeta_{3}^{1}$,
$ -\sqrt{2}\zeta_{6}^{1}$;\ \ 
$ -\zeta_{6}^{5}$,
$ -\zeta_{6}^{1}$,
$ -\sqrt{2}$,
$ -\sqrt{2}\zeta_{3}^{1}$,
$ \sqrt{2}\zeta_{6}^{1}$;\ \ 
$ -\zeta_{6}^{5}$,
$ -\sqrt{2}$,
$ \sqrt{2}\zeta_{6}^{1}$,
$ -\sqrt{2}\zeta_{3}^{1}$;\ \ 
$0$,
$0$,
$0$;\ \ 
$0$,
$0$;\ \ 
$0$)

Factors = $3_{2,3.}^{3,527}\boxtimes 3_{\frac{9}{2},4.}^{16,156}$

  \vskip 2ex

\noindent10. $9_{\frac{15}{2},12.}^{48,298}$ \irep{524}:\ \ 
$d_i$ = ($1.0$,
$1.0$,
$1.0$,
$1.0$,
$1.0$,
$1.0$,
$1.414$,
$1.414$,
$1.414$) 

\vskip 0.7ex
\hangindent=3em \hangafter=1
$D^2= 12.0 = 
12$

\vskip 0.7ex
\hangindent=3em \hangafter=1
$T = ( 0,
\frac{1}{2},
\frac{1}{3},
\frac{1}{3},
\frac{5}{6},
\frac{5}{6},
\frac{11}{16},
\frac{1}{48},
\frac{1}{48} )
$,

\vskip 0.7ex
\hangindent=3em \hangafter=1
$S$ = ($ 1$,
$ 1$,
$ 1$,
$ 1$,
$ 1$,
$ 1$,
$ \sqrt{2}$,
$ \sqrt{2}$,
$ \sqrt{2}$;\ \ 
$ 1$,
$ 1$,
$ 1$,
$ 1$,
$ 1$,
$ -\sqrt{2}$,
$ -\sqrt{2}$,
$ -\sqrt{2}$;\ \ 
$ -\zeta_{6}^{5}$,
$ -\zeta_{6}^{1}$,
$ -\zeta_{6}^{1}$,
$ -\zeta_{6}^{5}$,
$ \sqrt{2}$,
$ -\sqrt{2}\zeta_{6}^{1}$,
$ \sqrt{2}\zeta_{3}^{1}$;\ \ 
$ -\zeta_{6}^{5}$,
$ -\zeta_{6}^{5}$,
$ -\zeta_{6}^{1}$,
$ \sqrt{2}$,
$ \sqrt{2}\zeta_{3}^{1}$,
$ -\sqrt{2}\zeta_{6}^{1}$;\ \ 
$ -\zeta_{6}^{5}$,
$ -\zeta_{6}^{1}$,
$ -\sqrt{2}$,
$ -\sqrt{2}\zeta_{3}^{1}$,
$ \sqrt{2}\zeta_{6}^{1}$;\ \ 
$ -\zeta_{6}^{5}$,
$ -\sqrt{2}$,
$ \sqrt{2}\zeta_{6}^{1}$,
$ -\sqrt{2}\zeta_{3}^{1}$;\ \ 
$0$,
$0$,
$0$;\ \ 
$0$,
$0$;\ \ 
$0$)

Factors = $3_{2,3.}^{3,527}\boxtimes 3_{\frac{11}{2},4.}^{16,648}$

  \vskip 2ex

\noindent11. $9_{\frac{1}{2},12.}^{48,294}$ \irep{524}:\ \ 
$d_i$ = ($1.0$,
$1.0$,
$1.0$,
$1.0$,
$1.0$,
$1.0$,
$1.414$,
$1.414$,
$1.414$) 

\vskip 0.7ex
\hangindent=3em \hangafter=1
$D^2= 12.0 = 
12$

\vskip 0.7ex
\hangindent=3em \hangafter=1
$T = ( 0,
\frac{1}{2},
\frac{1}{3},
\frac{1}{3},
\frac{5}{6},
\frac{5}{6},
\frac{13}{16},
\frac{7}{48},
\frac{7}{48} )
$,

\vskip 0.7ex
\hangindent=3em \hangafter=1
$S$ = ($ 1$,
$ 1$,
$ 1$,
$ 1$,
$ 1$,
$ 1$,
$ \sqrt{2}$,
$ \sqrt{2}$,
$ \sqrt{2}$;\ \ 
$ 1$,
$ 1$,
$ 1$,
$ 1$,
$ 1$,
$ -\sqrt{2}$,
$ -\sqrt{2}$,
$ -\sqrt{2}$;\ \ 
$ -\zeta_{6}^{5}$,
$ -\zeta_{6}^{1}$,
$ -\zeta_{6}^{1}$,
$ -\zeta_{6}^{5}$,
$ \sqrt{2}$,
$ -\sqrt{2}\zeta_{6}^{1}$,
$ \sqrt{2}\zeta_{3}^{1}$;\ \ 
$ -\zeta_{6}^{5}$,
$ -\zeta_{6}^{5}$,
$ -\zeta_{6}^{1}$,
$ \sqrt{2}$,
$ \sqrt{2}\zeta_{3}^{1}$,
$ -\sqrt{2}\zeta_{6}^{1}$;\ \ 
$ -\zeta_{6}^{5}$,
$ -\zeta_{6}^{1}$,
$ -\sqrt{2}$,
$ -\sqrt{2}\zeta_{3}^{1}$,
$ \sqrt{2}\zeta_{6}^{1}$;\ \ 
$ -\zeta_{6}^{5}$,
$ -\sqrt{2}$,
$ \sqrt{2}\zeta_{6}^{1}$,
$ -\sqrt{2}\zeta_{3}^{1}$;\ \ 
$0$,
$0$,
$0$;\ \ 
$0$,
$0$;\ \ 
$0$)

Factors = $3_{2,3.}^{3,527}\boxtimes 3_{\frac{13}{2},4.}^{16,330}$

  \vskip 2ex

\noindent12. $9_{\frac{3}{2},12.}^{48,750}$ \irep{524}:\ \ 
$d_i$ = ($1.0$,
$1.0$,
$1.0$,
$1.0$,
$1.0$,
$1.0$,
$1.414$,
$1.414$,
$1.414$) 

\vskip 0.7ex
\hangindent=3em \hangafter=1
$D^2= 12.0 = 
12$

\vskip 0.7ex
\hangindent=3em \hangafter=1
$T = ( 0,
\frac{1}{2},
\frac{1}{3},
\frac{1}{3},
\frac{5}{6},
\frac{5}{6},
\frac{15}{16},
\frac{13}{48},
\frac{13}{48} )
$,

\vskip 0.7ex
\hangindent=3em \hangafter=1
$S$ = ($ 1$,
$ 1$,
$ 1$,
$ 1$,
$ 1$,
$ 1$,
$ \sqrt{2}$,
$ \sqrt{2}$,
$ \sqrt{2}$;\ \ 
$ 1$,
$ 1$,
$ 1$,
$ 1$,
$ 1$,
$ -\sqrt{2}$,
$ -\sqrt{2}$,
$ -\sqrt{2}$;\ \ 
$ -\zeta_{6}^{5}$,
$ -\zeta_{6}^{1}$,
$ -\zeta_{6}^{1}$,
$ -\zeta_{6}^{5}$,
$ \sqrt{2}$,
$ -\sqrt{2}\zeta_{6}^{1}$,
$ \sqrt{2}\zeta_{3}^{1}$;\ \ 
$ -\zeta_{6}^{5}$,
$ -\zeta_{6}^{5}$,
$ -\zeta_{6}^{1}$,
$ \sqrt{2}$,
$ \sqrt{2}\zeta_{3}^{1}$,
$ -\sqrt{2}\zeta_{6}^{1}$;\ \ 
$ -\zeta_{6}^{5}$,
$ -\zeta_{6}^{1}$,
$ -\sqrt{2}$,
$ -\sqrt{2}\zeta_{3}^{1}$,
$ \sqrt{2}\zeta_{6}^{1}$;\ \ 
$ -\zeta_{6}^{5}$,
$ -\sqrt{2}$,
$ \sqrt{2}\zeta_{6}^{1}$,
$ -\sqrt{2}\zeta_{3}^{1}$;\ \ 
$0$,
$0$,
$0$;\ \ 
$0$,
$0$;\ \ 
$0$)

Factors = $3_{2,3.}^{3,527}\boxtimes 3_{\frac{15}{2},4.}^{16,639}$

  \vskip 2ex

\noindent13. $9_{\frac{13}{2},12.}^{48,143}$ \irep{524}:\ \ 
$d_i$ = ($1.0$,
$1.0$,
$1.0$,
$1.0$,
$1.0$,
$1.0$,
$1.414$,
$1.414$,
$1.414$) 

\vskip 0.7ex
\hangindent=3em \hangafter=1
$D^2= 12.0 = 
12$

\vskip 0.7ex
\hangindent=3em \hangafter=1
$T = ( 0,
\frac{1}{2},
\frac{2}{3},
\frac{2}{3},
\frac{1}{6},
\frac{1}{6},
\frac{1}{16},
\frac{35}{48},
\frac{35}{48} )
$,

\vskip 0.7ex
\hangindent=3em \hangafter=1
$S$ = ($ 1$,
$ 1$,
$ 1$,
$ 1$,
$ 1$,
$ 1$,
$ \sqrt{2}$,
$ \sqrt{2}$,
$ \sqrt{2}$;\ \ 
$ 1$,
$ 1$,
$ 1$,
$ 1$,
$ 1$,
$ -\sqrt{2}$,
$ -\sqrt{2}$,
$ -\sqrt{2}$;\ \ 
$ -\zeta_{6}^{1}$,
$ -\zeta_{6}^{5}$,
$ -\zeta_{6}^{1}$,
$ -\zeta_{6}^{5}$,
$ \sqrt{2}$,
$ -\sqrt{2}\zeta_{6}^{1}$,
$ \sqrt{2}\zeta_{3}^{1}$;\ \ 
$ -\zeta_{6}^{1}$,
$ -\zeta_{6}^{5}$,
$ -\zeta_{6}^{1}$,
$ \sqrt{2}$,
$ \sqrt{2}\zeta_{3}^{1}$,
$ -\sqrt{2}\zeta_{6}^{1}$;\ \ 
$ -\zeta_{6}^{1}$,
$ -\zeta_{6}^{5}$,
$ -\sqrt{2}$,
$ \sqrt{2}\zeta_{6}^{1}$,
$ -\sqrt{2}\zeta_{3}^{1}$;\ \ 
$ -\zeta_{6}^{1}$,
$ -\sqrt{2}$,
$ -\sqrt{2}\zeta_{3}^{1}$,
$ \sqrt{2}\zeta_{6}^{1}$;\ \ 
$0$,
$0$,
$0$;\ \ 
$0$,
$0$;\ \ 
$0$)

Factors = $3_{6,3.}^{3,138}\boxtimes 3_{\frac{1}{2},4.}^{16,598}$

  \vskip 2ex

\noindent14. $9_{\frac{15}{2},12.}^{48,747}$ \irep{524}:\ \ 
$d_i$ = ($1.0$,
$1.0$,
$1.0$,
$1.0$,
$1.0$,
$1.0$,
$1.414$,
$1.414$,
$1.414$) 

\vskip 0.7ex
\hangindent=3em \hangafter=1
$D^2= 12.0 = 
12$

\vskip 0.7ex
\hangindent=3em \hangafter=1
$T = ( 0,
\frac{1}{2},
\frac{2}{3},
\frac{2}{3},
\frac{1}{6},
\frac{1}{6},
\frac{3}{16},
\frac{41}{48},
\frac{41}{48} )
$,

\vskip 0.7ex
\hangindent=3em \hangafter=1
$S$ = ($ 1$,
$ 1$,
$ 1$,
$ 1$,
$ 1$,
$ 1$,
$ \sqrt{2}$,
$ \sqrt{2}$,
$ \sqrt{2}$;\ \ 
$ 1$,
$ 1$,
$ 1$,
$ 1$,
$ 1$,
$ -\sqrt{2}$,
$ -\sqrt{2}$,
$ -\sqrt{2}$;\ \ 
$ -\zeta_{6}^{1}$,
$ -\zeta_{6}^{5}$,
$ -\zeta_{6}^{1}$,
$ -\zeta_{6}^{5}$,
$ \sqrt{2}$,
$ -\sqrt{2}\zeta_{6}^{1}$,
$ \sqrt{2}\zeta_{3}^{1}$;\ \ 
$ -\zeta_{6}^{1}$,
$ -\zeta_{6}^{5}$,
$ -\zeta_{6}^{1}$,
$ \sqrt{2}$,
$ \sqrt{2}\zeta_{3}^{1}$,
$ -\sqrt{2}\zeta_{6}^{1}$;\ \ 
$ -\zeta_{6}^{1}$,
$ -\zeta_{6}^{5}$,
$ -\sqrt{2}$,
$ \sqrt{2}\zeta_{6}^{1}$,
$ -\sqrt{2}\zeta_{3}^{1}$;\ \ 
$ -\zeta_{6}^{1}$,
$ -\sqrt{2}$,
$ -\sqrt{2}\zeta_{3}^{1}$,
$ \sqrt{2}\zeta_{6}^{1}$;\ \ 
$0$,
$0$,
$0$;\ \ 
$0$,
$0$;\ \ 
$0$)

Factors = $3_{6,3.}^{3,138}\boxtimes 3_{\frac{3}{2},4.}^{16,553}$

  \vskip 2ex

\noindent15. $9_{\frac{1}{2},12.}^{48,148}$ \irep{524}:\ \ 
$d_i$ = ($1.0$,
$1.0$,
$1.0$,
$1.0$,
$1.0$,
$1.0$,
$1.414$,
$1.414$,
$1.414$) 

\vskip 0.7ex
\hangindent=3em \hangafter=1
$D^2= 12.0 = 
12$

\vskip 0.7ex
\hangindent=3em \hangafter=1
$T = ( 0,
\frac{1}{2},
\frac{2}{3},
\frac{2}{3},
\frac{1}{6},
\frac{1}{6},
\frac{5}{16},
\frac{47}{48},
\frac{47}{48} )
$,

\vskip 0.7ex
\hangindent=3em \hangafter=1
$S$ = ($ 1$,
$ 1$,
$ 1$,
$ 1$,
$ 1$,
$ 1$,
$ \sqrt{2}$,
$ \sqrt{2}$,
$ \sqrt{2}$;\ \ 
$ 1$,
$ 1$,
$ 1$,
$ 1$,
$ 1$,
$ -\sqrt{2}$,
$ -\sqrt{2}$,
$ -\sqrt{2}$;\ \ 
$ -\zeta_{6}^{1}$,
$ -\zeta_{6}^{5}$,
$ -\zeta_{6}^{1}$,
$ -\zeta_{6}^{5}$,
$ \sqrt{2}$,
$ -\sqrt{2}\zeta_{6}^{1}$,
$ \sqrt{2}\zeta_{3}^{1}$;\ \ 
$ -\zeta_{6}^{1}$,
$ -\zeta_{6}^{5}$,
$ -\zeta_{6}^{1}$,
$ \sqrt{2}$,
$ \sqrt{2}\zeta_{3}^{1}$,
$ -\sqrt{2}\zeta_{6}^{1}$;\ \ 
$ -\zeta_{6}^{1}$,
$ -\zeta_{6}^{5}$,
$ -\sqrt{2}$,
$ \sqrt{2}\zeta_{6}^{1}$,
$ -\sqrt{2}\zeta_{3}^{1}$;\ \ 
$ -\zeta_{6}^{1}$,
$ -\sqrt{2}$,
$ -\sqrt{2}\zeta_{3}^{1}$,
$ \sqrt{2}\zeta_{6}^{1}$;\ \ 
$0$,
$0$,
$0$;\ \ 
$0$,
$0$;\ \ 
$0$)

Factors = $3_{6,3.}^{3,138}\boxtimes 3_{\frac{5}{2},4.}^{16,465}$

  \vskip 2ex

\noindent16. $9_{\frac{3}{2},12.}^{48,106}$ \irep{524}:\ \ 
$d_i$ = ($1.0$,
$1.0$,
$1.0$,
$1.0$,
$1.0$,
$1.0$,
$1.414$,
$1.414$,
$1.414$) 

\vskip 0.7ex
\hangindent=3em \hangafter=1
$D^2= 12.0 = 
12$

\vskip 0.7ex
\hangindent=3em \hangafter=1
$T = ( 0,
\frac{1}{2},
\frac{2}{3},
\frac{2}{3},
\frac{1}{6},
\frac{1}{6},
\frac{7}{16},
\frac{5}{48},
\frac{5}{48} )
$,

\vskip 0.7ex
\hangindent=3em \hangafter=1
$S$ = ($ 1$,
$ 1$,
$ 1$,
$ 1$,
$ 1$,
$ 1$,
$ \sqrt{2}$,
$ \sqrt{2}$,
$ \sqrt{2}$;\ \ 
$ 1$,
$ 1$,
$ 1$,
$ 1$,
$ 1$,
$ -\sqrt{2}$,
$ -\sqrt{2}$,
$ -\sqrt{2}$;\ \ 
$ -\zeta_{6}^{1}$,
$ -\zeta_{6}^{5}$,
$ -\zeta_{6}^{1}$,
$ -\zeta_{6}^{5}$,
$ \sqrt{2}$,
$ -\sqrt{2}\zeta_{6}^{1}$,
$ \sqrt{2}\zeta_{3}^{1}$;\ \ 
$ -\zeta_{6}^{1}$,
$ -\zeta_{6}^{5}$,
$ -\zeta_{6}^{1}$,
$ \sqrt{2}$,
$ \sqrt{2}\zeta_{3}^{1}$,
$ -\sqrt{2}\zeta_{6}^{1}$;\ \ 
$ -\zeta_{6}^{1}$,
$ -\zeta_{6}^{5}$,
$ -\sqrt{2}$,
$ \sqrt{2}\zeta_{6}^{1}$,
$ -\sqrt{2}\zeta_{3}^{1}$;\ \ 
$ -\zeta_{6}^{1}$,
$ -\sqrt{2}$,
$ -\sqrt{2}\zeta_{3}^{1}$,
$ \sqrt{2}\zeta_{6}^{1}$;\ \ 
$0$,
$0$,
$0$;\ \ 
$0$,
$0$;\ \ 
$0$)

Factors = $3_{6,3.}^{3,138}\boxtimes 3_{\frac{7}{2},4.}^{16,332}$

  \vskip 2ex

\noindent17. $9_{\frac{5}{2},12.}^{48,770}$ \irep{524}:\ \ 
$d_i$ = ($1.0$,
$1.0$,
$1.0$,
$1.0$,
$1.0$,
$1.0$,
$1.414$,
$1.414$,
$1.414$) 

\vskip 0.7ex
\hangindent=3em \hangafter=1
$D^2= 12.0 = 
12$

\vskip 0.7ex
\hangindent=3em \hangafter=1
$T = ( 0,
\frac{1}{2},
\frac{2}{3},
\frac{2}{3},
\frac{1}{6},
\frac{1}{6},
\frac{9}{16},
\frac{11}{48},
\frac{11}{48} )
$,

\vskip 0.7ex
\hangindent=3em \hangafter=1
$S$ = ($ 1$,
$ 1$,
$ 1$,
$ 1$,
$ 1$,
$ 1$,
$ \sqrt{2}$,
$ \sqrt{2}$,
$ \sqrt{2}$;\ \ 
$ 1$,
$ 1$,
$ 1$,
$ 1$,
$ 1$,
$ -\sqrt{2}$,
$ -\sqrt{2}$,
$ -\sqrt{2}$;\ \ 
$ -\zeta_{6}^{1}$,
$ -\zeta_{6}^{5}$,
$ -\zeta_{6}^{1}$,
$ -\zeta_{6}^{5}$,
$ \sqrt{2}$,
$ -\sqrt{2}\zeta_{6}^{1}$,
$ \sqrt{2}\zeta_{3}^{1}$;\ \ 
$ -\zeta_{6}^{1}$,
$ -\zeta_{6}^{5}$,
$ -\zeta_{6}^{1}$,
$ \sqrt{2}$,
$ \sqrt{2}\zeta_{3}^{1}$,
$ -\sqrt{2}\zeta_{6}^{1}$;\ \ 
$ -\zeta_{6}^{1}$,
$ -\zeta_{6}^{5}$,
$ -\sqrt{2}$,
$ \sqrt{2}\zeta_{6}^{1}$,
$ -\sqrt{2}\zeta_{3}^{1}$;\ \ 
$ -\zeta_{6}^{1}$,
$ -\sqrt{2}$,
$ -\sqrt{2}\zeta_{3}^{1}$,
$ \sqrt{2}\zeta_{6}^{1}$;\ \ 
$0$,
$0$,
$0$;\ \ 
$0$,
$0$;\ \ 
$0$)

Factors = $3_{6,3.}^{3,138}\boxtimes 3_{\frac{9}{2},4.}^{16,156}$

  \vskip 2ex

\noindent18. $9_{\frac{7}{2},12.}^{48,342}$ \irep{524}:\ \ 
$d_i$ = ($1.0$,
$1.0$,
$1.0$,
$1.0$,
$1.0$,
$1.0$,
$1.414$,
$1.414$,
$1.414$) 

\vskip 0.7ex
\hangindent=3em \hangafter=1
$D^2= 12.0 = 
12$

\vskip 0.7ex
\hangindent=3em \hangafter=1
$T = ( 0,
\frac{1}{2},
\frac{2}{3},
\frac{2}{3},
\frac{1}{6},
\frac{1}{6},
\frac{11}{16},
\frac{17}{48},
\frac{17}{48} )
$,

\vskip 0.7ex
\hangindent=3em \hangafter=1
$S$ = ($ 1$,
$ 1$,
$ 1$,
$ 1$,
$ 1$,
$ 1$,
$ \sqrt{2}$,
$ \sqrt{2}$,
$ \sqrt{2}$;\ \ 
$ 1$,
$ 1$,
$ 1$,
$ 1$,
$ 1$,
$ -\sqrt{2}$,
$ -\sqrt{2}$,
$ -\sqrt{2}$;\ \ 
$ -\zeta_{6}^{1}$,
$ -\zeta_{6}^{5}$,
$ -\zeta_{6}^{1}$,
$ -\zeta_{6}^{5}$,
$ \sqrt{2}$,
$ -\sqrt{2}\zeta_{6}^{1}$,
$ \sqrt{2}\zeta_{3}^{1}$;\ \ 
$ -\zeta_{6}^{1}$,
$ -\zeta_{6}^{5}$,
$ -\zeta_{6}^{1}$,
$ \sqrt{2}$,
$ \sqrt{2}\zeta_{3}^{1}$,
$ -\sqrt{2}\zeta_{6}^{1}$;\ \ 
$ -\zeta_{6}^{1}$,
$ -\zeta_{6}^{5}$,
$ -\sqrt{2}$,
$ \sqrt{2}\zeta_{6}^{1}$,
$ -\sqrt{2}\zeta_{3}^{1}$;\ \ 
$ -\zeta_{6}^{1}$,
$ -\sqrt{2}$,
$ -\sqrt{2}\zeta_{3}^{1}$,
$ \sqrt{2}\zeta_{6}^{1}$;\ \ 
$0$,
$0$,
$0$;\ \ 
$0$,
$0$;\ \ 
$0$)

Factors = $3_{6,3.}^{3,138}\boxtimes 3_{\frac{11}{2},4.}^{16,648}$

  \vskip 2ex

\noindent19. $9_{\frac{9}{2},12.}^{48,216}$ \irep{524}:\ \ 
$d_i$ = ($1.0$,
$1.0$,
$1.0$,
$1.0$,
$1.0$,
$1.0$,
$1.414$,
$1.414$,
$1.414$) 

\vskip 0.7ex
\hangindent=3em \hangafter=1
$D^2= 12.0 = 
12$

\vskip 0.7ex
\hangindent=3em \hangafter=1
$T = ( 0,
\frac{1}{2},
\frac{2}{3},
\frac{2}{3},
\frac{1}{6},
\frac{1}{6},
\frac{13}{16},
\frac{23}{48},
\frac{23}{48} )
$,

\vskip 0.7ex
\hangindent=3em \hangafter=1
$S$ = ($ 1$,
$ 1$,
$ 1$,
$ 1$,
$ 1$,
$ 1$,
$ \sqrt{2}$,
$ \sqrt{2}$,
$ \sqrt{2}$;\ \ 
$ 1$,
$ 1$,
$ 1$,
$ 1$,
$ 1$,
$ -\sqrt{2}$,
$ -\sqrt{2}$,
$ -\sqrt{2}$;\ \ 
$ -\zeta_{6}^{1}$,
$ -\zeta_{6}^{5}$,
$ -\zeta_{6}^{1}$,
$ -\zeta_{6}^{5}$,
$ \sqrt{2}$,
$ -\sqrt{2}\zeta_{6}^{1}$,
$ \sqrt{2}\zeta_{3}^{1}$;\ \ 
$ -\zeta_{6}^{1}$,
$ -\zeta_{6}^{5}$,
$ -\zeta_{6}^{1}$,
$ \sqrt{2}$,
$ \sqrt{2}\zeta_{3}^{1}$,
$ -\sqrt{2}\zeta_{6}^{1}$;\ \ 
$ -\zeta_{6}^{1}$,
$ -\zeta_{6}^{5}$,
$ -\sqrt{2}$,
$ \sqrt{2}\zeta_{6}^{1}$,
$ -\sqrt{2}\zeta_{3}^{1}$;\ \ 
$ -\zeta_{6}^{1}$,
$ -\sqrt{2}$,
$ -\sqrt{2}\zeta_{3}^{1}$,
$ \sqrt{2}\zeta_{6}^{1}$;\ \ 
$0$,
$0$,
$0$;\ \ 
$0$,
$0$;\ \ 
$0$)

Factors = $3_{6,3.}^{3,138}\boxtimes 3_{\frac{13}{2},4.}^{16,330}$

  \vskip 2ex

\noindent20. $9_{\frac{11}{2},12.}^{48,909}$ \irep{524}:\ \ 
$d_i$ = ($1.0$,
$1.0$,
$1.0$,
$1.0$,
$1.0$,
$1.0$,
$1.414$,
$1.414$,
$1.414$) 

\vskip 0.7ex
\hangindent=3em \hangafter=1
$D^2= 12.0 = 
12$

\vskip 0.7ex
\hangindent=3em \hangafter=1
$T = ( 0,
\frac{1}{2},
\frac{2}{3},
\frac{2}{3},
\frac{1}{6},
\frac{1}{6},
\frac{15}{16},
\frac{29}{48},
\frac{29}{48} )
$,

\vskip 0.7ex
\hangindent=3em \hangafter=1
$S$ = ($ 1$,
$ 1$,
$ 1$,
$ 1$,
$ 1$,
$ 1$,
$ \sqrt{2}$,
$ \sqrt{2}$,
$ \sqrt{2}$;\ \ 
$ 1$,
$ 1$,
$ 1$,
$ 1$,
$ 1$,
$ -\sqrt{2}$,
$ -\sqrt{2}$,
$ -\sqrt{2}$;\ \ 
$ -\zeta_{6}^{1}$,
$ -\zeta_{6}^{5}$,
$ -\zeta_{6}^{1}$,
$ -\zeta_{6}^{5}$,
$ \sqrt{2}$,
$ -\sqrt{2}\zeta_{6}^{1}$,
$ \sqrt{2}\zeta_{3}^{1}$;\ \ 
$ -\zeta_{6}^{1}$,
$ -\zeta_{6}^{5}$,
$ -\zeta_{6}^{1}$,
$ \sqrt{2}$,
$ \sqrt{2}\zeta_{3}^{1}$,
$ -\sqrt{2}\zeta_{6}^{1}$;\ \ 
$ -\zeta_{6}^{1}$,
$ -\zeta_{6}^{5}$,
$ -\sqrt{2}$,
$ \sqrt{2}\zeta_{6}^{1}$,
$ -\sqrt{2}\zeta_{3}^{1}$;\ \ 
$ -\zeta_{6}^{1}$,
$ -\sqrt{2}$,
$ -\sqrt{2}\zeta_{3}^{1}$,
$ \sqrt{2}\zeta_{6}^{1}$;\ \ 
$0$,
$0$,
$0$;\ \ 
$0$,
$0$;\ \ 
$0$)

Factors = $3_{6,3.}^{3,138}\boxtimes 3_{\frac{15}{2},4.}^{16,639}$

  \vskip 2ex

\noindent21. $9_{1,16.}^{16,147}$ \irep{429}:\ \ 
$d_i$ = ($1.0$,
$1.0$,
$1.0$,
$1.0$,
$1.414$,
$1.414$,
$1.414$,
$1.414$,
$2.0$) 

\vskip 0.7ex
\hangindent=3em \hangafter=1
$D^2= 16.0 = 
16$

\vskip 0.7ex
\hangindent=3em \hangafter=1
$T = ( 0,
0,
\frac{1}{2},
\frac{1}{2},
\frac{1}{16},
\frac{1}{16},
\frac{9}{16},
\frac{9}{16},
\frac{1}{8} )
$,

\vskip 0.7ex
\hangindent=3em \hangafter=1
$S$ = ($ 1$,
$ 1$,
$ 1$,
$ 1$,
$ \sqrt{2}$,
$ \sqrt{2}$,
$ \sqrt{2}$,
$ \sqrt{2}$,
$ 2$;\ \ 
$ 1$,
$ 1$,
$ 1$,
$ -\sqrt{2}$,
$ -\sqrt{2}$,
$ -\sqrt{2}$,
$ -\sqrt{2}$,
$ 2$;\ \ 
$ 1$,
$ 1$,
$ -\sqrt{2}$,
$ \sqrt{2}$,
$ -\sqrt{2}$,
$ \sqrt{2}$,
$ -2$;\ \ 
$ 1$,
$ \sqrt{2}$,
$ -\sqrt{2}$,
$ \sqrt{2}$,
$ -\sqrt{2}$,
$ -2$;\ \ 
$0$,
$ 2$,
$0$,
$ -2$,
$0$;\ \ 
$0$,
$ -2$,
$0$,
$0$;\ \ 
$0$,
$ 2$,
$0$;\ \ 
$0$,
$0$;\ \ 
$0$)

Factors = $3_{\frac{1}{2},4.}^{16,598}\boxtimes 3_{\frac{1}{2},4.}^{16,598}$

  \vskip 2ex

\noindent22. $9_{5,16.}^{16,726}$ \irep{429}:\ \ 
$d_i$ = ($1.0$,
$1.0$,
$1.0$,
$1.0$,
$1.414$,
$1.414$,
$1.414$,
$1.414$,
$2.0$) 

\vskip 0.7ex
\hangindent=3em \hangafter=1
$D^2= 16.0 = 
16$

\vskip 0.7ex
\hangindent=3em \hangafter=1
$T = ( 0,
0,
\frac{1}{2},
\frac{1}{2},
\frac{1}{16},
\frac{1}{16},
\frac{9}{16},
\frac{9}{16},
\frac{5}{8} )
$,

\vskip 0.7ex
\hangindent=3em \hangafter=1
$S$ = ($ 1$,
$ 1$,
$ 1$,
$ 1$,
$ \sqrt{2}$,
$ \sqrt{2}$,
$ \sqrt{2}$,
$ \sqrt{2}$,
$ 2$;\ \ 
$ 1$,
$ 1$,
$ 1$,
$ -\sqrt{2}$,
$ -\sqrt{2}$,
$ -\sqrt{2}$,
$ -\sqrt{2}$,
$ 2$;\ \ 
$ 1$,
$ 1$,
$ -\sqrt{2}$,
$ \sqrt{2}$,
$ -\sqrt{2}$,
$ \sqrt{2}$,
$ -2$;\ \ 
$ 1$,
$ \sqrt{2}$,
$ -\sqrt{2}$,
$ \sqrt{2}$,
$ -\sqrt{2}$,
$ -2$;\ \ 
$0$,
$ -2$,
$0$,
$ 2$,
$0$;\ \ 
$0$,
$ 2$,
$0$,
$0$;\ \ 
$0$,
$ -2$,
$0$;\ \ 
$0$,
$0$;\ \ 
$0$)

Factors = $3_{\frac{1}{2},4.}^{16,598}\boxtimes 3_{\frac{9}{2},4.}^{16,156}$

  \vskip 2ex

\noindent23. $9_{2,16.}^{16,111}$ \irep{419}:\ \ 
$d_i$ = ($1.0$,
$1.0$,
$1.0$,
$1.0$,
$1.414$,
$1.414$,
$1.414$,
$1.414$,
$2.0$) 

\vskip 0.7ex
\hangindent=3em \hangafter=1
$D^2= 16.0 = 
16$

\vskip 0.7ex
\hangindent=3em \hangafter=1
$T = ( 0,
0,
\frac{1}{2},
\frac{1}{2},
\frac{1}{16},
\frac{3}{16},
\frac{9}{16},
\frac{11}{16},
\frac{1}{4} )
$,

\vskip 0.7ex
\hangindent=3em \hangafter=1
$S$ = ($ 1$,
$ 1$,
$ 1$,
$ 1$,
$ \sqrt{2}$,
$ \sqrt{2}$,
$ \sqrt{2}$,
$ \sqrt{2}$,
$ 2$;\ \ 
$ 1$,
$ 1$,
$ 1$,
$ -\sqrt{2}$,
$ -\sqrt{2}$,
$ -\sqrt{2}$,
$ -\sqrt{2}$,
$ 2$;\ \ 
$ 1$,
$ 1$,
$ -\sqrt{2}$,
$ \sqrt{2}$,
$ -\sqrt{2}$,
$ \sqrt{2}$,
$ -2$;\ \ 
$ 1$,
$ \sqrt{2}$,
$ -\sqrt{2}$,
$ \sqrt{2}$,
$ -\sqrt{2}$,
$ -2$;\ \ 
$0$,
$ 2$,
$0$,
$ -2$,
$0$;\ \ 
$0$,
$ -2$,
$0$,
$0$;\ \ 
$0$,
$ 2$,
$0$;\ \ 
$0$,
$0$;\ \ 
$0$)

Factors = $3_{\frac{1}{2},4.}^{16,598}\boxtimes 3_{\frac{3}{2},4.}^{16,553}$

  \vskip 2ex

\noindent24. $9_{6,16.}^{16,118}$ \irep{419}:\ \ 
$d_i$ = ($1.0$,
$1.0$,
$1.0$,
$1.0$,
$1.414$,
$1.414$,
$1.414$,
$1.414$,
$2.0$) 

\vskip 0.7ex
\hangindent=3em \hangafter=1
$D^2= 16.0 = 
16$

\vskip 0.7ex
\hangindent=3em \hangafter=1
$T = ( 0,
0,
\frac{1}{2},
\frac{1}{2},
\frac{1}{16},
\frac{3}{16},
\frac{9}{16},
\frac{11}{16},
\frac{3}{4} )
$,

\vskip 0.7ex
\hangindent=3em \hangafter=1
$S$ = ($ 1$,
$ 1$,
$ 1$,
$ 1$,
$ \sqrt{2}$,
$ \sqrt{2}$,
$ \sqrt{2}$,
$ \sqrt{2}$,
$ 2$;\ \ 
$ 1$,
$ 1$,
$ 1$,
$ -\sqrt{2}$,
$ -\sqrt{2}$,
$ -\sqrt{2}$,
$ -\sqrt{2}$,
$ 2$;\ \ 
$ 1$,
$ 1$,
$ -\sqrt{2}$,
$ \sqrt{2}$,
$ -\sqrt{2}$,
$ \sqrt{2}$,
$ -2$;\ \ 
$ 1$,
$ \sqrt{2}$,
$ -\sqrt{2}$,
$ \sqrt{2}$,
$ -\sqrt{2}$,
$ -2$;\ \ 
$0$,
$ -2$,
$0$,
$ 2$,
$0$;\ \ 
$0$,
$ 2$,
$0$,
$0$;\ \ 
$0$,
$ -2$,
$0$;\ \ 
$0$,
$0$;\ \ 
$0$)

Factors = $3_{\frac{1}{2},4.}^{16,598}\boxtimes 3_{\frac{11}{2},4.}^{16,648}$

  \vskip 2ex

\noindent25. $9_{3,16.}^{16,608}$ \irep{425}:\ \ 
$d_i$ = ($1.0$,
$1.0$,
$1.0$,
$1.0$,
$1.414$,
$1.414$,
$1.414$,
$1.414$,
$2.0$) 

\vskip 0.7ex
\hangindent=3em \hangafter=1
$D^2= 16.0 = 
16$

\vskip 0.7ex
\hangindent=3em \hangafter=1
$T = ( 0,
0,
\frac{1}{2},
\frac{1}{2},
\frac{1}{16},
\frac{5}{16},
\frac{9}{16},
\frac{13}{16},
\frac{3}{8} )
$,

\vskip 0.7ex
\hangindent=3em \hangafter=1
$S$ = ($ 1$,
$ 1$,
$ 1$,
$ 1$,
$ \sqrt{2}$,
$ \sqrt{2}$,
$ \sqrt{2}$,
$ \sqrt{2}$,
$ 2$;\ \ 
$ 1$,
$ 1$,
$ 1$,
$ -\sqrt{2}$,
$ -\sqrt{2}$,
$ -\sqrt{2}$,
$ -\sqrt{2}$,
$ 2$;\ \ 
$ 1$,
$ 1$,
$ -\sqrt{2}$,
$ \sqrt{2}$,
$ -\sqrt{2}$,
$ \sqrt{2}$,
$ -2$;\ \ 
$ 1$,
$ \sqrt{2}$,
$ -\sqrt{2}$,
$ \sqrt{2}$,
$ -\sqrt{2}$,
$ -2$;\ \ 
$0$,
$ 2$,
$0$,
$ -2$,
$0$;\ \ 
$0$,
$ -2$,
$0$,
$0$;\ \ 
$0$,
$ 2$,
$0$;\ \ 
$0$,
$0$;\ \ 
$0$)

Factors = $3_{\frac{1}{2},4.}^{16,598}\boxtimes 3_{\frac{5}{2},4.}^{16,465}$

  \vskip 2ex

\noindent26. $9_{7,16.}^{16,641}$ \irep{425}:\ \ 
$d_i$ = ($1.0$,
$1.0$,
$1.0$,
$1.0$,
$1.414$,
$1.414$,
$1.414$,
$1.414$,
$2.0$) 

\vskip 0.7ex
\hangindent=3em \hangafter=1
$D^2= 16.0 = 
16$

\vskip 0.7ex
\hangindent=3em \hangafter=1
$T = ( 0,
0,
\frac{1}{2},
\frac{1}{2},
\frac{1}{16},
\frac{5}{16},
\frac{9}{16},
\frac{13}{16},
\frac{7}{8} )
$,

\vskip 0.7ex
\hangindent=3em \hangafter=1
$S$ = ($ 1$,
$ 1$,
$ 1$,
$ 1$,
$ \sqrt{2}$,
$ \sqrt{2}$,
$ \sqrt{2}$,
$ \sqrt{2}$,
$ 2$;\ \ 
$ 1$,
$ 1$,
$ 1$,
$ -\sqrt{2}$,
$ -\sqrt{2}$,
$ -\sqrt{2}$,
$ -\sqrt{2}$,
$ 2$;\ \ 
$ 1$,
$ 1$,
$ -\sqrt{2}$,
$ \sqrt{2}$,
$ -\sqrt{2}$,
$ \sqrt{2}$,
$ -2$;\ \ 
$ 1$,
$ \sqrt{2}$,
$ -\sqrt{2}$,
$ \sqrt{2}$,
$ -\sqrt{2}$,
$ -2$;\ \ 
$0$,
$ -2$,
$0$,
$ 2$,
$0$;\ \ 
$0$,
$ 2$,
$0$,
$0$;\ \ 
$0$,
$ -2$,
$0$;\ \ 
$0$,
$0$;\ \ 
$0$)

Factors = $3_{\frac{1}{2},4.}^{16,598}\boxtimes 3_{\frac{13}{2},4.}^{16,330}$

  \vskip 2ex

\noindent27. $9_{0,16.}^{16,447}$ \irep{424}:\ \ 
$d_i$ = ($1.0$,
$1.0$,
$1.0$,
$1.0$,
$1.414$,
$1.414$,
$1.414$,
$1.414$,
$2.0$) 

\vskip 0.7ex
\hangindent=3em \hangafter=1
$D^2= 16.0 = 
16$

\vskip 0.7ex
\hangindent=3em \hangafter=1
$T = ( 0,
0,
\frac{1}{2},
\frac{1}{2},
\frac{1}{16},
\frac{7}{16},
\frac{9}{16},
\frac{15}{16},
0 )
$,

\vskip 0.7ex
\hangindent=3em \hangafter=1
$S$ = ($ 1$,
$ 1$,
$ 1$,
$ 1$,
$ \sqrt{2}$,
$ \sqrt{2}$,
$ \sqrt{2}$,
$ \sqrt{2}$,
$ 2$;\ \ 
$ 1$,
$ 1$,
$ 1$,
$ -\sqrt{2}$,
$ -\sqrt{2}$,
$ -\sqrt{2}$,
$ -\sqrt{2}$,
$ 2$;\ \ 
$ 1$,
$ 1$,
$ -\sqrt{2}$,
$ \sqrt{2}$,
$ -\sqrt{2}$,
$ \sqrt{2}$,
$ -2$;\ \ 
$ 1$,
$ \sqrt{2}$,
$ -\sqrt{2}$,
$ \sqrt{2}$,
$ -\sqrt{2}$,
$ -2$;\ \ 
$0$,
$ -2$,
$0$,
$ 2$,
$0$;\ \ 
$0$,
$ 2$,
$0$,
$0$;\ \ 
$0$,
$ -2$,
$0$;\ \ 
$0$,
$0$;\ \ 
$0$)

Factors = $3_{\frac{1}{2},4.}^{16,598}\boxtimes 3_{\frac{15}{2},4.}^{16,639}$

  \vskip 2ex

\noindent28. $9_{4,16.}^{16,524}$ \irep{424}:\ \ 
$d_i$ = ($1.0$,
$1.0$,
$1.0$,
$1.0$,
$1.414$,
$1.414$,
$1.414$,
$1.414$,
$2.0$) 

\vskip 0.7ex
\hangindent=3em \hangafter=1
$D^2= 16.0 = 
16$

\vskip 0.7ex
\hangindent=3em \hangafter=1
$T = ( 0,
0,
\frac{1}{2},
\frac{1}{2},
\frac{1}{16},
\frac{7}{16},
\frac{9}{16},
\frac{15}{16},
\frac{1}{2} )
$,

\vskip 0.7ex
\hangindent=3em \hangafter=1
$S$ = ($ 1$,
$ 1$,
$ 1$,
$ 1$,
$ \sqrt{2}$,
$ \sqrt{2}$,
$ \sqrt{2}$,
$ \sqrt{2}$,
$ 2$;\ \ 
$ 1$,
$ 1$,
$ 1$,
$ -\sqrt{2}$,
$ -\sqrt{2}$,
$ -\sqrt{2}$,
$ -\sqrt{2}$,
$ 2$;\ \ 
$ 1$,
$ 1$,
$ -\sqrt{2}$,
$ \sqrt{2}$,
$ -\sqrt{2}$,
$ \sqrt{2}$,
$ -2$;\ \ 
$ 1$,
$ \sqrt{2}$,
$ -\sqrt{2}$,
$ \sqrt{2}$,
$ -\sqrt{2}$,
$ -2$;\ \ 
$0$,
$ 2$,
$0$,
$ -2$,
$0$;\ \ 
$0$,
$ -2$,
$0$,
$0$;\ \ 
$0$,
$ 2$,
$0$;\ \ 
$0$,
$0$;\ \ 
$0$)

Factors = $3_{\frac{1}{2},4.}^{16,598}\boxtimes 3_{\frac{7}{2},4.}^{16,332}$

  \vskip 2ex

\noindent29. $9_{3,16.}^{16,696}$ \irep{429}:\ \ 
$d_i$ = ($1.0$,
$1.0$,
$1.0$,
$1.0$,
$1.414$,
$1.414$,
$1.414$,
$1.414$,
$2.0$) 

\vskip 0.7ex
\hangindent=3em \hangafter=1
$D^2= 16.0 = 
16$

\vskip 0.7ex
\hangindent=3em \hangafter=1
$T = ( 0,
0,
\frac{1}{2},
\frac{1}{2},
\frac{3}{16},
\frac{3}{16},
\frac{11}{16},
\frac{11}{16},
\frac{3}{8} )
$,

\vskip 0.7ex
\hangindent=3em \hangafter=1
$S$ = ($ 1$,
$ 1$,
$ 1$,
$ 1$,
$ \sqrt{2}$,
$ \sqrt{2}$,
$ \sqrt{2}$,
$ \sqrt{2}$,
$ 2$;\ \ 
$ 1$,
$ 1$,
$ 1$,
$ -\sqrt{2}$,
$ -\sqrt{2}$,
$ -\sqrt{2}$,
$ -\sqrt{2}$,
$ 2$;\ \ 
$ 1$,
$ 1$,
$ -\sqrt{2}$,
$ \sqrt{2}$,
$ -\sqrt{2}$,
$ \sqrt{2}$,
$ -2$;\ \ 
$ 1$,
$ \sqrt{2}$,
$ -\sqrt{2}$,
$ \sqrt{2}$,
$ -\sqrt{2}$,
$ -2$;\ \ 
$0$,
$ 2$,
$0$,
$ -2$,
$0$;\ \ 
$0$,
$ -2$,
$0$,
$0$;\ \ 
$0$,
$ 2$,
$0$;\ \ 
$0$,
$0$;\ \ 
$0$)

Factors = $3_{\frac{3}{2},4.}^{16,553}\boxtimes 3_{\frac{3}{2},4.}^{16,553}$

  \vskip 2ex

\noindent30. $9_{7,16.}^{16,553}$ \irep{429}:\ \ 
$d_i$ = ($1.0$,
$1.0$,
$1.0$,
$1.0$,
$1.414$,
$1.414$,
$1.414$,
$1.414$,
$2.0$) 

\vskip 0.7ex
\hangindent=3em \hangafter=1
$D^2= 16.0 = 
16$

\vskip 0.7ex
\hangindent=3em \hangafter=1
$T = ( 0,
0,
\frac{1}{2},
\frac{1}{2},
\frac{3}{16},
\frac{3}{16},
\frac{11}{16},
\frac{11}{16},
\frac{7}{8} )
$,

\vskip 0.7ex
\hangindent=3em \hangafter=1
$S$ = ($ 1$,
$ 1$,
$ 1$,
$ 1$,
$ \sqrt{2}$,
$ \sqrt{2}$,
$ \sqrt{2}$,
$ \sqrt{2}$,
$ 2$;\ \ 
$ 1$,
$ 1$,
$ 1$,
$ -\sqrt{2}$,
$ -\sqrt{2}$,
$ -\sqrt{2}$,
$ -\sqrt{2}$,
$ 2$;\ \ 
$ 1$,
$ 1$,
$ -\sqrt{2}$,
$ \sqrt{2}$,
$ -\sqrt{2}$,
$ \sqrt{2}$,
$ -2$;\ \ 
$ 1$,
$ \sqrt{2}$,
$ -\sqrt{2}$,
$ \sqrt{2}$,
$ -\sqrt{2}$,
$ -2$;\ \ 
$0$,
$ -2$,
$0$,
$ 2$,
$0$;\ \ 
$0$,
$ 2$,
$0$,
$0$;\ \ 
$0$,
$ -2$,
$0$;\ \ 
$0$,
$0$;\ \ 
$0$)

Factors = $3_{\frac{3}{2},4.}^{16,553}\boxtimes 3_{\frac{11}{2},4.}^{16,648}$

  \vskip 2ex

\noindent31. $9_{0,16.}^{16,624}$ \irep{424}:\ \ 
$d_i$ = ($1.0$,
$1.0$,
$1.0$,
$1.0$,
$1.414$,
$1.414$,
$1.414$,
$1.414$,
$2.0$) 

\vskip 0.7ex
\hangindent=3em \hangafter=1
$D^2= 16.0 = 
16$

\vskip 0.7ex
\hangindent=3em \hangafter=1
$T = ( 0,
0,
\frac{1}{2},
\frac{1}{2},
\frac{3}{16},
\frac{5}{16},
\frac{11}{16},
\frac{13}{16},
0 )
$,

\vskip 0.7ex
\hangindent=3em \hangafter=1
$S$ = ($ 1$,
$ 1$,
$ 1$,
$ 1$,
$ \sqrt{2}$,
$ \sqrt{2}$,
$ \sqrt{2}$,
$ \sqrt{2}$,
$ 2$;\ \ 
$ 1$,
$ 1$,
$ 1$,
$ -\sqrt{2}$,
$ -\sqrt{2}$,
$ -\sqrt{2}$,
$ -\sqrt{2}$,
$ 2$;\ \ 
$ 1$,
$ 1$,
$ -\sqrt{2}$,
$ \sqrt{2}$,
$ -\sqrt{2}$,
$ \sqrt{2}$,
$ -2$;\ \ 
$ 1$,
$ \sqrt{2}$,
$ -\sqrt{2}$,
$ \sqrt{2}$,
$ -\sqrt{2}$,
$ -2$;\ \ 
$0$,
$ -2$,
$0$,
$ 2$,
$0$;\ \ 
$0$,
$ 2$,
$0$,
$0$;\ \ 
$0$,
$ -2$,
$0$;\ \ 
$0$,
$0$;\ \ 
$0$)

Factors = $3_{\frac{3}{2},4.}^{16,553}\boxtimes 3_{\frac{13}{2},4.}^{16,330}$

  \vskip 2ex

\noindent32. $9_{4,16.}^{16,124}$ \irep{424}:\ \ 
$d_i$ = ($1.0$,
$1.0$,
$1.0$,
$1.0$,
$1.414$,
$1.414$,
$1.414$,
$1.414$,
$2.0$) 

\vskip 0.7ex
\hangindent=3em \hangafter=1
$D^2= 16.0 = 
16$

\vskip 0.7ex
\hangindent=3em \hangafter=1
$T = ( 0,
0,
\frac{1}{2},
\frac{1}{2},
\frac{3}{16},
\frac{5}{16},
\frac{11}{16},
\frac{13}{16},
\frac{1}{2} )
$,

\vskip 0.7ex
\hangindent=3em \hangafter=1
$S$ = ($ 1$,
$ 1$,
$ 1$,
$ 1$,
$ \sqrt{2}$,
$ \sqrt{2}$,
$ \sqrt{2}$,
$ \sqrt{2}$,
$ 2$;\ \ 
$ 1$,
$ 1$,
$ 1$,
$ -\sqrt{2}$,
$ -\sqrt{2}$,
$ -\sqrt{2}$,
$ -\sqrt{2}$,
$ 2$;\ \ 
$ 1$,
$ 1$,
$ -\sqrt{2}$,
$ \sqrt{2}$,
$ -\sqrt{2}$,
$ \sqrt{2}$,
$ -2$;\ \ 
$ 1$,
$ \sqrt{2}$,
$ -\sqrt{2}$,
$ \sqrt{2}$,
$ -\sqrt{2}$,
$ -2$;\ \ 
$0$,
$ 2$,
$0$,
$ -2$,
$0$;\ \ 
$0$,
$ -2$,
$0$,
$0$;\ \ 
$0$,
$ 2$,
$0$;\ \ 
$0$,
$0$;\ \ 
$0$)

Factors = $3_{\frac{3}{2},4.}^{16,553}\boxtimes 3_{\frac{5}{2},4.}^{16,465}$

  \vskip 2ex

\noindent33. $9_{1,16.}^{16,151}$ \irep{425}:\ \ 
$d_i$ = ($1.0$,
$1.0$,
$1.0$,
$1.0$,
$1.414$,
$1.414$,
$1.414$,
$1.414$,
$2.0$) 

\vskip 0.7ex
\hangindent=3em \hangafter=1
$D^2= 16.0 = 
16$

\vskip 0.7ex
\hangindent=3em \hangafter=1
$T = ( 0,
0,
\frac{1}{2},
\frac{1}{2},
\frac{3}{16},
\frac{7}{16},
\frac{11}{16},
\frac{15}{16},
\frac{1}{8} )
$,

\vskip 0.7ex
\hangindent=3em \hangafter=1
$S$ = ($ 1$,
$ 1$,
$ 1$,
$ 1$,
$ \sqrt{2}$,
$ \sqrt{2}$,
$ \sqrt{2}$,
$ \sqrt{2}$,
$ 2$;\ \ 
$ 1$,
$ 1$,
$ 1$,
$ -\sqrt{2}$,
$ -\sqrt{2}$,
$ -\sqrt{2}$,
$ -\sqrt{2}$,
$ 2$;\ \ 
$ 1$,
$ 1$,
$ -\sqrt{2}$,
$ \sqrt{2}$,
$ -\sqrt{2}$,
$ \sqrt{2}$,
$ -2$;\ \ 
$ 1$,
$ \sqrt{2}$,
$ -\sqrt{2}$,
$ \sqrt{2}$,
$ -\sqrt{2}$,
$ -2$;\ \ 
$0$,
$ -2$,
$0$,
$ 2$,
$0$;\ \ 
$0$,
$ 2$,
$0$,
$0$;\ \ 
$0$,
$ -2$,
$0$;\ \ 
$0$,
$0$;\ \ 
$0$)

Factors = $3_{\frac{7}{2},4.}^{16,332}\boxtimes 3_{\frac{11}{2},4.}^{16,648}$

  \vskip 2ex

\noindent34. $9_{5,16.}^{16,598}$ \irep{425}:\ \ 
$d_i$ = ($1.0$,
$1.0$,
$1.0$,
$1.0$,
$1.414$,
$1.414$,
$1.414$,
$1.414$,
$2.0$) 

\vskip 0.7ex
\hangindent=3em \hangafter=1
$D^2= 16.0 = 
16$

\vskip 0.7ex
\hangindent=3em \hangafter=1
$T = ( 0,
0,
\frac{1}{2},
\frac{1}{2},
\frac{3}{16},
\frac{7}{16},
\frac{11}{16},
\frac{15}{16},
\frac{5}{8} )
$,

\vskip 0.7ex
\hangindent=3em \hangafter=1
$S$ = ($ 1$,
$ 1$,
$ 1$,
$ 1$,
$ \sqrt{2}$,
$ \sqrt{2}$,
$ \sqrt{2}$,
$ \sqrt{2}$,
$ 2$;\ \ 
$ 1$,
$ 1$,
$ 1$,
$ -\sqrt{2}$,
$ -\sqrt{2}$,
$ -\sqrt{2}$,
$ -\sqrt{2}$,
$ 2$;\ \ 
$ 1$,
$ 1$,
$ -\sqrt{2}$,
$ \sqrt{2}$,
$ -\sqrt{2}$,
$ \sqrt{2}$,
$ -2$;\ \ 
$ 1$,
$ \sqrt{2}$,
$ -\sqrt{2}$,
$ \sqrt{2}$,
$ -\sqrt{2}$,
$ -2$;\ \ 
$0$,
$ 2$,
$0$,
$ -2$,
$0$;\ \ 
$0$,
$ -2$,
$0$,
$0$;\ \ 
$0$,
$ 2$,
$0$;\ \ 
$0$,
$0$;\ \ 
$0$)

Factors = $3_{\frac{7}{2},4.}^{16,332}\boxtimes 3_{\frac{3}{2},4.}^{16,553}$

  \vskip 2ex

\noindent35. $9_{1,16.}^{16,239}$ \irep{429}:\ \ 
$d_i$ = ($1.0$,
$1.0$,
$1.0$,
$1.0$,
$1.414$,
$1.414$,
$1.414$,
$1.414$,
$2.0$) 

\vskip 0.7ex
\hangindent=3em \hangafter=1
$D^2= 16.0 = 
16$

\vskip 0.7ex
\hangindent=3em \hangafter=1
$T = ( 0,
0,
\frac{1}{2},
\frac{1}{2},
\frac{5}{16},
\frac{5}{16},
\frac{13}{16},
\frac{13}{16},
\frac{1}{8} )
$,

\vskip 0.7ex
\hangindent=3em \hangafter=1
$S$ = ($ 1$,
$ 1$,
$ 1$,
$ 1$,
$ \sqrt{2}$,
$ \sqrt{2}$,
$ \sqrt{2}$,
$ \sqrt{2}$,
$ 2$;\ \ 
$ 1$,
$ 1$,
$ 1$,
$ -\sqrt{2}$,
$ -\sqrt{2}$,
$ -\sqrt{2}$,
$ -\sqrt{2}$,
$ 2$;\ \ 
$ 1$,
$ 1$,
$ -\sqrt{2}$,
$ \sqrt{2}$,
$ -\sqrt{2}$,
$ \sqrt{2}$,
$ -2$;\ \ 
$ 1$,
$ \sqrt{2}$,
$ -\sqrt{2}$,
$ \sqrt{2}$,
$ -\sqrt{2}$,
$ -2$;\ \ 
$0$,
$ -2$,
$0$,
$ 2$,
$0$;\ \ 
$0$,
$ 2$,
$0$,
$0$;\ \ 
$0$,
$ -2$,
$0$;\ \ 
$0$,
$0$;\ \ 
$0$)

Factors = $3_{\frac{5}{2},4.}^{16,465}\boxtimes 3_{\frac{13}{2},4.}^{16,330}$

  \vskip 2ex

\noindent36. $9_{5,16.}^{16,510}$ \irep{429}:\ \ 
$d_i$ = ($1.0$,
$1.0$,
$1.0$,
$1.0$,
$1.414$,
$1.414$,
$1.414$,
$1.414$,
$2.0$) 

\vskip 0.7ex
\hangindent=3em \hangafter=1
$D^2= 16.0 = 
16$

\vskip 0.7ex
\hangindent=3em \hangafter=1
$T = ( 0,
0,
\frac{1}{2},
\frac{1}{2},
\frac{5}{16},
\frac{5}{16},
\frac{13}{16},
\frac{13}{16},
\frac{5}{8} )
$,

\vskip 0.7ex
\hangindent=3em \hangafter=1
$S$ = ($ 1$,
$ 1$,
$ 1$,
$ 1$,
$ \sqrt{2}$,
$ \sqrt{2}$,
$ \sqrt{2}$,
$ \sqrt{2}$,
$ 2$;\ \ 
$ 1$,
$ 1$,
$ 1$,
$ -\sqrt{2}$,
$ -\sqrt{2}$,
$ -\sqrt{2}$,
$ -\sqrt{2}$,
$ 2$;\ \ 
$ 1$,
$ 1$,
$ -\sqrt{2}$,
$ \sqrt{2}$,
$ -\sqrt{2}$,
$ \sqrt{2}$,
$ -2$;\ \ 
$ 1$,
$ \sqrt{2}$,
$ -\sqrt{2}$,
$ \sqrt{2}$,
$ -\sqrt{2}$,
$ -2$;\ \ 
$0$,
$ 2$,
$0$,
$ -2$,
$0$;\ \ 
$0$,
$ -2$,
$0$,
$0$;\ \ 
$0$,
$ 2$,
$0$;\ \ 
$0$,
$0$;\ \ 
$0$)

Factors = $3_{\frac{5}{2},4.}^{16,465}\boxtimes 3_{\frac{5}{2},4.}^{16,465}$

  \vskip 2ex

\noindent37. $9_{2,16.}^{16,296}$ \irep{419}:\ \ 
$d_i$ = ($1.0$,
$1.0$,
$1.0$,
$1.0$,
$1.414$,
$1.414$,
$1.414$,
$1.414$,
$2.0$) 

\vskip 0.7ex
\hangindent=3em \hangafter=1
$D^2= 16.0 = 
16$

\vskip 0.7ex
\hangindent=3em \hangafter=1
$T = ( 0,
0,
\frac{1}{2},
\frac{1}{2},
\frac{5}{16},
\frac{7}{16},
\frac{13}{16},
\frac{15}{16},
\frac{1}{4} )
$,

\vskip 0.7ex
\hangindent=3em \hangafter=1
$S$ = ($ 1$,
$ 1$,
$ 1$,
$ 1$,
$ \sqrt{2}$,
$ \sqrt{2}$,
$ \sqrt{2}$,
$ \sqrt{2}$,
$ 2$;\ \ 
$ 1$,
$ 1$,
$ 1$,
$ -\sqrt{2}$,
$ -\sqrt{2}$,
$ -\sqrt{2}$,
$ -\sqrt{2}$,
$ 2$;\ \ 
$ 1$,
$ 1$,
$ -\sqrt{2}$,
$ \sqrt{2}$,
$ -\sqrt{2}$,
$ \sqrt{2}$,
$ -2$;\ \ 
$ 1$,
$ \sqrt{2}$,
$ -\sqrt{2}$,
$ \sqrt{2}$,
$ -\sqrt{2}$,
$ -2$;\ \ 
$0$,
$ -2$,
$0$,
$ 2$,
$0$;\ \ 
$0$,
$ 2$,
$0$,
$0$;\ \ 
$0$,
$ -2$,
$0$;\ \ 
$0$,
$0$;\ \ 
$0$)

Factors = $3_{\frac{7}{2},4.}^{16,332}\boxtimes 3_{\frac{13}{2},4.}^{16,330}$

  \vskip 2ex

\noindent38. $9_{6,16.}^{16,129}$ \irep{419}:\ \ 
$d_i$ = ($1.0$,
$1.0$,
$1.0$,
$1.0$,
$1.414$,
$1.414$,
$1.414$,
$1.414$,
$2.0$) 

\vskip 0.7ex
\hangindent=3em \hangafter=1
$D^2= 16.0 = 
16$

\vskip 0.7ex
\hangindent=3em \hangafter=1
$T = ( 0,
0,
\frac{1}{2},
\frac{1}{2},
\frac{5}{16},
\frac{7}{16},
\frac{13}{16},
\frac{15}{16},
\frac{3}{4} )
$,

\vskip 0.7ex
\hangindent=3em \hangafter=1
$S$ = ($ 1$,
$ 1$,
$ 1$,
$ 1$,
$ \sqrt{2}$,
$ \sqrt{2}$,
$ \sqrt{2}$,
$ \sqrt{2}$,
$ 2$;\ \ 
$ 1$,
$ 1$,
$ 1$,
$ -\sqrt{2}$,
$ -\sqrt{2}$,
$ -\sqrt{2}$,
$ -\sqrt{2}$,
$ 2$;\ \ 
$ 1$,
$ 1$,
$ -\sqrt{2}$,
$ \sqrt{2}$,
$ -\sqrt{2}$,
$ \sqrt{2}$,
$ -2$;\ \ 
$ 1$,
$ \sqrt{2}$,
$ -\sqrt{2}$,
$ \sqrt{2}$,
$ -\sqrt{2}$,
$ -2$;\ \ 
$0$,
$ 2$,
$0$,
$ -2$,
$0$;\ \ 
$0$,
$ -2$,
$0$,
$0$;\ \ 
$0$,
$ 2$,
$0$;\ \ 
$0$,
$0$;\ \ 
$0$)

Factors = $3_{\frac{7}{2},4.}^{16,332}\boxtimes 3_{\frac{5}{2},4.}^{16,465}$

  \vskip 2ex

\noindent39. $9_{3,16.}^{16,894}$ \irep{429}:\ \ 
$d_i$ = ($1.0$,
$1.0$,
$1.0$,
$1.0$,
$1.414$,
$1.414$,
$1.414$,
$1.414$,
$2.0$) 

\vskip 0.7ex
\hangindent=3em \hangafter=1
$D^2= 16.0 = 
16$

\vskip 0.7ex
\hangindent=3em \hangafter=1
$T = ( 0,
0,
\frac{1}{2},
\frac{1}{2},
\frac{7}{16},
\frac{7}{16},
\frac{15}{16},
\frac{15}{16},
\frac{3}{8} )
$,

\vskip 0.7ex
\hangindent=3em \hangafter=1
$S$ = ($ 1$,
$ 1$,
$ 1$,
$ 1$,
$ \sqrt{2}$,
$ \sqrt{2}$,
$ \sqrt{2}$,
$ \sqrt{2}$,
$ 2$;\ \ 
$ 1$,
$ 1$,
$ 1$,
$ -\sqrt{2}$,
$ -\sqrt{2}$,
$ -\sqrt{2}$,
$ -\sqrt{2}$,
$ 2$;\ \ 
$ 1$,
$ 1$,
$ -\sqrt{2}$,
$ \sqrt{2}$,
$ -\sqrt{2}$,
$ \sqrt{2}$,
$ -2$;\ \ 
$ 1$,
$ \sqrt{2}$,
$ -\sqrt{2}$,
$ \sqrt{2}$,
$ -\sqrt{2}$,
$ -2$;\ \ 
$0$,
$ -2$,
$0$,
$ 2$,
$0$;\ \ 
$0$,
$ 2$,
$0$,
$0$;\ \ 
$0$,
$ -2$,
$0$;\ \ 
$0$,
$0$;\ \ 
$0$)

Factors = $3_{\frac{7}{2},4.}^{16,332}\boxtimes 3_{\frac{15}{2},4.}^{16,639}$

  \vskip 2ex

\noindent40. $9_{7,16.}^{16,214}$ \irep{429}:\ \ 
$d_i$ = ($1.0$,
$1.0$,
$1.0$,
$1.0$,
$1.414$,
$1.414$,
$1.414$,
$1.414$,
$2.0$) 

\vskip 0.7ex
\hangindent=3em \hangafter=1
$D^2= 16.0 = 
16$

\vskip 0.7ex
\hangindent=3em \hangafter=1
$T = ( 0,
0,
\frac{1}{2},
\frac{1}{2},
\frac{7}{16},
\frac{7}{16},
\frac{15}{16},
\frac{15}{16},
\frac{7}{8} )
$,

\vskip 0.7ex
\hangindent=3em \hangafter=1
$S$ = ($ 1$,
$ 1$,
$ 1$,
$ 1$,
$ \sqrt{2}$,
$ \sqrt{2}$,
$ \sqrt{2}$,
$ \sqrt{2}$,
$ 2$;\ \ 
$ 1$,
$ 1$,
$ 1$,
$ -\sqrt{2}$,
$ -\sqrt{2}$,
$ -\sqrt{2}$,
$ -\sqrt{2}$,
$ 2$;\ \ 
$ 1$,
$ 1$,
$ -\sqrt{2}$,
$ \sqrt{2}$,
$ -\sqrt{2}$,
$ \sqrt{2}$,
$ -2$;\ \ 
$ 1$,
$ \sqrt{2}$,
$ -\sqrt{2}$,
$ \sqrt{2}$,
$ -\sqrt{2}$,
$ -2$;\ \ 
$0$,
$ 2$,
$0$,
$ -2$,
$0$;\ \ 
$0$,
$ -2$,
$0$,
$0$;\ \ 
$0$,
$ 2$,
$0$;\ \ 
$0$,
$0$;\ \ 
$0$)

Factors = $3_{\frac{7}{2},4.}^{16,332}\boxtimes 3_{\frac{7}{2},4.}^{16,332}$

  \vskip 2ex

\noindent41. $9_{\frac{6}{7},27.88}^{21,155}$ \irep{467}:\ \ 
$d_i$ = ($1.0$,
$1.0$,
$1.0$,
$1.801$,
$1.801$,
$1.801$,
$2.246$,
$2.246$,
$2.246$) 

\vskip 0.7ex
\hangindent=3em \hangafter=1
$D^2= 27.887 = 
 18+9c^{1}_{7}
+3c^{2}_{7}
$

\vskip 0.7ex
\hangindent=3em \hangafter=1
$T = ( 0,
\frac{1}{3},
\frac{1}{3},
\frac{1}{7},
\frac{10}{21},
\frac{10}{21},
\frac{5}{7},
\frac{1}{21},
\frac{1}{21} )
$,

\vskip 0.7ex
\hangindent=3em \hangafter=1
$S$ = ($ 1$,
$ 1$,
$ 1$,
$ -c_{7}^{3}$,
$ -c_{7}^{3}$,
$ -c_{7}^{3}$,
$ \xi_{7}^{3}$,
$ \xi_{7}^{3}$,
$ \xi_{7}^{3}$;\ \ 
$ -\zeta_{6}^{5}$,
$ -\zeta_{6}^{1}$,
$ -c_{7}^{3}$,
$  \zeta^{3}_{7}
+\zeta^{4}_{7}
+\zeta^{16}_{21}
+\zeta^{19}_{21}
$,
$  -\zeta^{16}_{21}
-\zeta^{19}_{21}
$,
$ \xi_{7}^{3}$,
$  -\zeta^{1}_{21}
-\zeta^{13}_{21}
-\zeta^{16}_{21}
-\zeta^{19}_{21}
$,
$  \zeta^{1}_{21}
+\zeta^{2}_{7}
+\zeta^{3}_{7}
+\zeta^{4}_{7}
+\zeta^{13}_{21}
+\zeta^{5}_{7}
+\zeta^{16}_{21}
+\zeta^{19}_{21}
$;\ \ 
$ -\zeta_{6}^{5}$,
$ -c_{7}^{3}$,
$  -\zeta^{16}_{21}
-\zeta^{19}_{21}
$,
$  \zeta^{3}_{7}
+\zeta^{4}_{7}
+\zeta^{16}_{21}
+\zeta^{19}_{21}
$,
$ \xi_{7}^{3}$,
$  \zeta^{1}_{21}
+\zeta^{2}_{7}
+\zeta^{3}_{7}
+\zeta^{4}_{7}
+\zeta^{13}_{21}
+\zeta^{5}_{7}
+\zeta^{16}_{21}
+\zeta^{19}_{21}
$,
$  -\zeta^{1}_{21}
-\zeta^{13}_{21}
-\zeta^{16}_{21}
-\zeta^{19}_{21}
$;\ \ 
$ -\xi_{7}^{3}$,
$ -\xi_{7}^{3}$,
$ -\xi_{7}^{3}$,
$ 1$,
$ 1$,
$ 1$;\ \ 
$  \zeta^{1}_{21}
+\zeta^{13}_{21}
+\zeta^{16}_{21}
+\zeta^{19}_{21}
$,
$  -\zeta^{1}_{21}
-\zeta^{2}_{7}
-\zeta^{3}_{7}
-\zeta^{4}_{7}
-\zeta^{13}_{21}
-\zeta^{5}_{7}
-\zeta^{16}_{21}
-\zeta^{19}_{21}
$,
$ 1$,
$ -\zeta_{6}^{1}$,
$ -\zeta_{6}^{5}$;\ \ 
$  \zeta^{1}_{21}
+\zeta^{13}_{21}
+\zeta^{16}_{21}
+\zeta^{19}_{21}
$,
$ 1$,
$ -\zeta_{6}^{5}$,
$ -\zeta_{6}^{1}$;\ \ 
$ c_{7}^{3}$,
$ c_{7}^{3}$,
$ c_{7}^{3}$;\ \ 
$  \zeta^{16}_{21}
+\zeta^{19}_{21}
$,
$  -\zeta^{3}_{7}
-\zeta^{4}_{7}
-\zeta^{16}_{21}
-\zeta^{19}_{21}
$;\ \ 
$  \zeta^{16}_{21}
+\zeta^{19}_{21}
$)

Factors = $3_{2,3.}^{3,527}\boxtimes 3_{\frac{48}{7},9.295}^{7,790}$

  \vskip 2ex

\noindent42. $9_{\frac{22}{7},27.88}^{21,204}$ \irep{467}:\ \ 
$d_i$ = ($1.0$,
$1.0$,
$1.0$,
$1.801$,
$1.801$,
$1.801$,
$2.246$,
$2.246$,
$2.246$) 

\vskip 0.7ex
\hangindent=3em \hangafter=1
$D^2= 27.887 = 
 18+9c^{1}_{7}
+3c^{2}_{7}
$

\vskip 0.7ex
\hangindent=3em \hangafter=1
$T = ( 0,
\frac{1}{3},
\frac{1}{3},
\frac{6}{7},
\frac{4}{21},
\frac{4}{21},
\frac{2}{7},
\frac{13}{21},
\frac{13}{21} )
$,

\vskip 0.7ex
\hangindent=3em \hangafter=1
$S$ = ($ 1$,
$ 1$,
$ 1$,
$ -c_{7}^{3}$,
$ -c_{7}^{3}$,
$ -c_{7}^{3}$,
$ \xi_{7}^{3}$,
$ \xi_{7}^{3}$,
$ \xi_{7}^{3}$;\ \ 
$ -\zeta_{6}^{5}$,
$ -\zeta_{6}^{1}$,
$ -c_{7}^{3}$,
$  \zeta^{3}_{7}
+\zeta^{4}_{7}
+\zeta^{16}_{21}
+\zeta^{19}_{21}
$,
$  -\zeta^{16}_{21}
-\zeta^{19}_{21}
$,
$ \xi_{7}^{3}$,
$  -\zeta^{1}_{21}
-\zeta^{13}_{21}
-\zeta^{16}_{21}
-\zeta^{19}_{21}
$,
$  \zeta^{1}_{21}
+\zeta^{2}_{7}
+\zeta^{3}_{7}
+\zeta^{4}_{7}
+\zeta^{13}_{21}
+\zeta^{5}_{7}
+\zeta^{16}_{21}
+\zeta^{19}_{21}
$;\ \ 
$ -\zeta_{6}^{5}$,
$ -c_{7}^{3}$,
$  -\zeta^{16}_{21}
-\zeta^{19}_{21}
$,
$  \zeta^{3}_{7}
+\zeta^{4}_{7}
+\zeta^{16}_{21}
+\zeta^{19}_{21}
$,
$ \xi_{7}^{3}$,
$  \zeta^{1}_{21}
+\zeta^{2}_{7}
+\zeta^{3}_{7}
+\zeta^{4}_{7}
+\zeta^{13}_{21}
+\zeta^{5}_{7}
+\zeta^{16}_{21}
+\zeta^{19}_{21}
$,
$  -\zeta^{1}_{21}
-\zeta^{13}_{21}
-\zeta^{16}_{21}
-\zeta^{19}_{21}
$;\ \ 
$ -\xi_{7}^{3}$,
$ -\xi_{7}^{3}$,
$ -\xi_{7}^{3}$,
$ 1$,
$ 1$,
$ 1$;\ \ 
$  \zeta^{1}_{21}
+\zeta^{13}_{21}
+\zeta^{16}_{21}
+\zeta^{19}_{21}
$,
$  -\zeta^{1}_{21}
-\zeta^{2}_{7}
-\zeta^{3}_{7}
-\zeta^{4}_{7}
-\zeta^{13}_{21}
-\zeta^{5}_{7}
-\zeta^{16}_{21}
-\zeta^{19}_{21}
$,
$ 1$,
$ -\zeta_{6}^{1}$,
$ -\zeta_{6}^{5}$;\ \ 
$  \zeta^{1}_{21}
+\zeta^{13}_{21}
+\zeta^{16}_{21}
+\zeta^{19}_{21}
$,
$ 1$,
$ -\zeta_{6}^{5}$,
$ -\zeta_{6}^{1}$;\ \ 
$ c_{7}^{3}$,
$ c_{7}^{3}$,
$ c_{7}^{3}$;\ \ 
$  \zeta^{16}_{21}
+\zeta^{19}_{21}
$,
$  -\zeta^{3}_{7}
-\zeta^{4}_{7}
-\zeta^{16}_{21}
-\zeta^{19}_{21}
$;\ \ 
$  \zeta^{16}_{21}
+\zeta^{19}_{21}
$)

Factors = $3_{2,3.}^{3,527}\boxtimes 3_{\frac{8}{7},9.295}^{7,245}$

  \vskip 2ex

\noindent43. $9_{\frac{34}{7},27.88}^{21,129}$ \irep{467}:\ \ 
$d_i$ = ($1.0$,
$1.0$,
$1.0$,
$1.801$,
$1.801$,
$1.801$,
$2.246$,
$2.246$,
$2.246$) 

\vskip 0.7ex
\hangindent=3em \hangafter=1
$D^2= 27.887 = 
 18+9c^{1}_{7}
+3c^{2}_{7}
$

\vskip 0.7ex
\hangindent=3em \hangafter=1
$T = ( 0,
\frac{2}{3},
\frac{2}{3},
\frac{1}{7},
\frac{17}{21},
\frac{17}{21},
\frac{5}{7},
\frac{8}{21},
\frac{8}{21} )
$,

\vskip 0.7ex
\hangindent=3em \hangafter=1
$S$ = ($ 1$,
$ 1$,
$ 1$,
$ -c_{7}^{3}$,
$ -c_{7}^{3}$,
$ -c_{7}^{3}$,
$ \xi_{7}^{3}$,
$ \xi_{7}^{3}$,
$ \xi_{7}^{3}$;\ \ 
$ -\zeta_{6}^{1}$,
$ -\zeta_{6}^{5}$,
$ -c_{7}^{3}$,
$  \zeta^{3}_{7}
+\zeta^{4}_{7}
+\zeta^{16}_{21}
+\zeta^{19}_{21}
$,
$  -\zeta^{16}_{21}
-\zeta^{19}_{21}
$,
$ \xi_{7}^{3}$,
$  -\zeta^{1}_{21}
-\zeta^{13}_{21}
-\zeta^{16}_{21}
-\zeta^{19}_{21}
$,
$  \zeta^{1}_{21}
+\zeta^{2}_{7}
+\zeta^{3}_{7}
+\zeta^{4}_{7}
+\zeta^{13}_{21}
+\zeta^{5}_{7}
+\zeta^{16}_{21}
+\zeta^{19}_{21}
$;\ \ 
$ -\zeta_{6}^{1}$,
$ -c_{7}^{3}$,
$  -\zeta^{16}_{21}
-\zeta^{19}_{21}
$,
$  \zeta^{3}_{7}
+\zeta^{4}_{7}
+\zeta^{16}_{21}
+\zeta^{19}_{21}
$,
$ \xi_{7}^{3}$,
$  \zeta^{1}_{21}
+\zeta^{2}_{7}
+\zeta^{3}_{7}
+\zeta^{4}_{7}
+\zeta^{13}_{21}
+\zeta^{5}_{7}
+\zeta^{16}_{21}
+\zeta^{19}_{21}
$,
$  -\zeta^{1}_{21}
-\zeta^{13}_{21}
-\zeta^{16}_{21}
-\zeta^{19}_{21}
$;\ \ 
$ -\xi_{7}^{3}$,
$ -\xi_{7}^{3}$,
$ -\xi_{7}^{3}$,
$ 1$,
$ 1$,
$ 1$;\ \ 
$  -\zeta^{1}_{21}
-\zeta^{2}_{7}
-\zeta^{3}_{7}
-\zeta^{4}_{7}
-\zeta^{13}_{21}
-\zeta^{5}_{7}
-\zeta^{16}_{21}
-\zeta^{19}_{21}
$,
$  \zeta^{1}_{21}
+\zeta^{13}_{21}
+\zeta^{16}_{21}
+\zeta^{19}_{21}
$,
$ 1$,
$ -\zeta_{6}^{5}$,
$ -\zeta_{6}^{1}$;\ \ 
$  -\zeta^{1}_{21}
-\zeta^{2}_{7}
-\zeta^{3}_{7}
-\zeta^{4}_{7}
-\zeta^{13}_{21}
-\zeta^{5}_{7}
-\zeta^{16}_{21}
-\zeta^{19}_{21}
$,
$ 1$,
$ -\zeta_{6}^{1}$,
$ -\zeta_{6}^{5}$;\ \ 
$ c_{7}^{3}$,
$ c_{7}^{3}$,
$ c_{7}^{3}$;\ \ 
$  -\zeta^{3}_{7}
-\zeta^{4}_{7}
-\zeta^{16}_{21}
-\zeta^{19}_{21}
$,
$  \zeta^{16}_{21}
+\zeta^{19}_{21}
$;\ \ 
$  -\zeta^{3}_{7}
-\zeta^{4}_{7}
-\zeta^{16}_{21}
-\zeta^{19}_{21}
$)

Factors = $3_{6,3.}^{3,138}\boxtimes 3_{\frac{48}{7},9.295}^{7,790}$

  \vskip 2ex

\noindent44. $9_{\frac{50}{7},27.88}^{21,367}$ \irep{467}:\ \ 
$d_i$ = ($1.0$,
$1.0$,
$1.0$,
$1.801$,
$1.801$,
$1.801$,
$2.246$,
$2.246$,
$2.246$) 

\vskip 0.7ex
\hangindent=3em \hangafter=1
$D^2= 27.887 = 
 18+9c^{1}_{7}
+3c^{2}_{7}
$

\vskip 0.7ex
\hangindent=3em \hangafter=1
$T = ( 0,
\frac{2}{3},
\frac{2}{3},
\frac{6}{7},
\frac{11}{21},
\frac{11}{21},
\frac{2}{7},
\frac{20}{21},
\frac{20}{21} )
$,

\vskip 0.7ex
\hangindent=3em \hangafter=1
$S$ = ($ 1$,
$ 1$,
$ 1$,
$ -c_{7}^{3}$,
$ -c_{7}^{3}$,
$ -c_{7}^{3}$,
$ \xi_{7}^{3}$,
$ \xi_{7}^{3}$,
$ \xi_{7}^{3}$;\ \ 
$ -\zeta_{6}^{1}$,
$ -\zeta_{6}^{5}$,
$ -c_{7}^{3}$,
$  \zeta^{3}_{7}
+\zeta^{4}_{7}
+\zeta^{16}_{21}
+\zeta^{19}_{21}
$,
$  -\zeta^{16}_{21}
-\zeta^{19}_{21}
$,
$ \xi_{7}^{3}$,
$  -\zeta^{1}_{21}
-\zeta^{13}_{21}
-\zeta^{16}_{21}
-\zeta^{19}_{21}
$,
$  \zeta^{1}_{21}
+\zeta^{2}_{7}
+\zeta^{3}_{7}
+\zeta^{4}_{7}
+\zeta^{13}_{21}
+\zeta^{5}_{7}
+\zeta^{16}_{21}
+\zeta^{19}_{21}
$;\ \ 
$ -\zeta_{6}^{1}$,
$ -c_{7}^{3}$,
$  -\zeta^{16}_{21}
-\zeta^{19}_{21}
$,
$  \zeta^{3}_{7}
+\zeta^{4}_{7}
+\zeta^{16}_{21}
+\zeta^{19}_{21}
$,
$ \xi_{7}^{3}$,
$  \zeta^{1}_{21}
+\zeta^{2}_{7}
+\zeta^{3}_{7}
+\zeta^{4}_{7}
+\zeta^{13}_{21}
+\zeta^{5}_{7}
+\zeta^{16}_{21}
+\zeta^{19}_{21}
$,
$  -\zeta^{1}_{21}
-\zeta^{13}_{21}
-\zeta^{16}_{21}
-\zeta^{19}_{21}
$;\ \ 
$ -\xi_{7}^{3}$,
$ -\xi_{7}^{3}$,
$ -\xi_{7}^{3}$,
$ 1$,
$ 1$,
$ 1$;\ \ 
$  -\zeta^{1}_{21}
-\zeta^{2}_{7}
-\zeta^{3}_{7}
-\zeta^{4}_{7}
-\zeta^{13}_{21}
-\zeta^{5}_{7}
-\zeta^{16}_{21}
-\zeta^{19}_{21}
$,
$  \zeta^{1}_{21}
+\zeta^{13}_{21}
+\zeta^{16}_{21}
+\zeta^{19}_{21}
$,
$ 1$,
$ -\zeta_{6}^{5}$,
$ -\zeta_{6}^{1}$;\ \ 
$  -\zeta^{1}_{21}
-\zeta^{2}_{7}
-\zeta^{3}_{7}
-\zeta^{4}_{7}
-\zeta^{13}_{21}
-\zeta^{5}_{7}
-\zeta^{16}_{21}
-\zeta^{19}_{21}
$,
$ 1$,
$ -\zeta_{6}^{1}$,
$ -\zeta_{6}^{5}$;\ \ 
$ c_{7}^{3}$,
$ c_{7}^{3}$,
$ c_{7}^{3}$;\ \ 
$  -\zeta^{3}_{7}
-\zeta^{4}_{7}
-\zeta^{16}_{21}
-\zeta^{19}_{21}
$,
$  \zeta^{16}_{21}
+\zeta^{19}_{21}
$;\ \ 
$  -\zeta^{3}_{7}
-\zeta^{4}_{7}
-\zeta^{16}_{21}
-\zeta^{19}_{21}
$)

Factors = $3_{6,3.}^{3,138}\boxtimes 3_{\frac{8}{7},9.295}^{7,245}$

  \vskip 2ex

\noindent45. $9_{\frac{103}{14},37.18}^{112,193}$ \irep{532}:\ \ 
$d_i$ = ($1.0$,
$1.0$,
$1.414$,
$1.801$,
$1.801$,
$2.246$,
$2.246$,
$2.548$,
$3.177$) 

\vskip 0.7ex
\hangindent=3em \hangafter=1
$D^2= 37.183 = 
 24+12c^{1}_{7}
+4c^{2}_{7}
$

\vskip 0.7ex
\hangindent=3em \hangafter=1
$T = ( 0,
\frac{1}{2},
\frac{1}{16},
\frac{1}{7},
\frac{9}{14},
\frac{5}{7},
\frac{3}{14},
\frac{23}{112},
\frac{87}{112} )
$,

\vskip 0.7ex
\hangindent=3em \hangafter=1
$S$ = ($ 1$,
$ 1$,
$ \sqrt{2}$,
$ -c_{7}^{3}$,
$ -c_{7}^{3}$,
$ \xi_{7}^{3}$,
$ \xi_{7}^{3}$,
$ \sqrt{2}c_{14}^{1}$,
$  c^{3}_{56}
+c^{5}_{56}
-c^{9}_{56}
+c^{11}_{56}
$;\ \ 
$ 1$,
$ -\sqrt{2}$,
$ -c_{7}^{3}$,
$ -c_{7}^{3}$,
$ \xi_{7}^{3}$,
$ \xi_{7}^{3}$,
$ -\sqrt{2}c_{14}^{1}$,
$  -c^{3}_{56}
-c^{5}_{56}
+c^{9}_{56}
-c^{11}_{56}
$;\ \ 
$0$,
$ \sqrt{2}c_{14}^{1}$,
$ -\sqrt{2}c_{14}^{1}$,
$  c^{3}_{56}
+c^{5}_{56}
-c^{9}_{56}
+c^{11}_{56}
$,
$  -c^{3}_{56}
-c^{5}_{56}
+c^{9}_{56}
-c^{11}_{56}
$,
$0$,
$0$;\ \ 
$ -\xi_{7}^{3}$,
$ -\xi_{7}^{3}$,
$ 1$,
$ 1$,
$  -c^{3}_{56}
-c^{5}_{56}
+c^{9}_{56}
-c^{11}_{56}
$,
$ \sqrt{2}$;\ \ 
$ -\xi_{7}^{3}$,
$ 1$,
$ 1$,
$  c^{3}_{56}
+c^{5}_{56}
-c^{9}_{56}
+c^{11}_{56}
$,
$ -\sqrt{2}$;\ \ 
$ c_{7}^{3}$,
$ c_{7}^{3}$,
$ \sqrt{2}$,
$ -\sqrt{2}c_{14}^{1}$;\ \ 
$ c_{7}^{3}$,
$ -\sqrt{2}$,
$ \sqrt{2}c_{14}^{1}$;\ \ 
$0$,
$0$;\ \ 
$0$)

Factors = $3_{\frac{1}{2},4.}^{16,598}\boxtimes 3_{\frac{48}{7},9.295}^{7,790}$

  \vskip 2ex

\noindent46. $9_{\frac{23}{14},37.18}^{112,161}$ \irep{532}:\ \ 
$d_i$ = ($1.0$,
$1.0$,
$1.414$,
$1.801$,
$1.801$,
$2.246$,
$2.246$,
$2.548$,
$3.177$) 

\vskip 0.7ex
\hangindent=3em \hangafter=1
$D^2= 37.183 = 
 24+12c^{1}_{7}
+4c^{2}_{7}
$

\vskip 0.7ex
\hangindent=3em \hangafter=1
$T = ( 0,
\frac{1}{2},
\frac{1}{16},
\frac{6}{7},
\frac{5}{14},
\frac{2}{7},
\frac{11}{14},
\frac{103}{112},
\frac{39}{112} )
$,

\vskip 0.7ex
\hangindent=3em \hangafter=1
$S$ = ($ 1$,
$ 1$,
$ \sqrt{2}$,
$ -c_{7}^{3}$,
$ -c_{7}^{3}$,
$ \xi_{7}^{3}$,
$ \xi_{7}^{3}$,
$ \sqrt{2}c_{14}^{1}$,
$  c^{3}_{56}
+c^{5}_{56}
-c^{9}_{56}
+c^{11}_{56}
$;\ \ 
$ 1$,
$ -\sqrt{2}$,
$ -c_{7}^{3}$,
$ -c_{7}^{3}$,
$ \xi_{7}^{3}$,
$ \xi_{7}^{3}$,
$ -\sqrt{2}c_{14}^{1}$,
$  -c^{3}_{56}
-c^{5}_{56}
+c^{9}_{56}
-c^{11}_{56}
$;\ \ 
$0$,
$ \sqrt{2}c_{14}^{1}$,
$ -\sqrt{2}c_{14}^{1}$,
$  c^{3}_{56}
+c^{5}_{56}
-c^{9}_{56}
+c^{11}_{56}
$,
$  -c^{3}_{56}
-c^{5}_{56}
+c^{9}_{56}
-c^{11}_{56}
$,
$0$,
$0$;\ \ 
$ -\xi_{7}^{3}$,
$ -\xi_{7}^{3}$,
$ 1$,
$ 1$,
$  -c^{3}_{56}
-c^{5}_{56}
+c^{9}_{56}
-c^{11}_{56}
$,
$ \sqrt{2}$;\ \ 
$ -\xi_{7}^{3}$,
$ 1$,
$ 1$,
$  c^{3}_{56}
+c^{5}_{56}
-c^{9}_{56}
+c^{11}_{56}
$,
$ -\sqrt{2}$;\ \ 
$ c_{7}^{3}$,
$ c_{7}^{3}$,
$ \sqrt{2}$,
$ -\sqrt{2}c_{14}^{1}$;\ \ 
$ c_{7}^{3}$,
$ -\sqrt{2}$,
$ \sqrt{2}c_{14}^{1}$;\ \ 
$0$,
$0$;\ \ 
$0$)

Factors = $3_{\frac{1}{2},4.}^{16,598}\boxtimes 3_{\frac{8}{7},9.295}^{7,245}$

  \vskip 2ex

\noindent47. $9_{\frac{5}{14},37.18}^{112,883}$ \irep{532}:\ \ 
$d_i$ = ($1.0$,
$1.0$,
$1.414$,
$1.801$,
$1.801$,
$2.246$,
$2.246$,
$2.548$,
$3.177$) 

\vskip 0.7ex
\hangindent=3em \hangafter=1
$D^2= 37.183 = 
 24+12c^{1}_{7}
+4c^{2}_{7}
$

\vskip 0.7ex
\hangindent=3em \hangafter=1
$T = ( 0,
\frac{1}{2},
\frac{3}{16},
\frac{1}{7},
\frac{9}{14},
\frac{5}{7},
\frac{3}{14},
\frac{37}{112},
\frac{101}{112} )
$,

\vskip 0.7ex
\hangindent=3em \hangafter=1
$S$ = ($ 1$,
$ 1$,
$ \sqrt{2}$,
$ -c_{7}^{3}$,
$ -c_{7}^{3}$,
$ \xi_{7}^{3}$,
$ \xi_{7}^{3}$,
$ \sqrt{2}c_{14}^{1}$,
$  c^{3}_{56}
+c^{5}_{56}
-c^{9}_{56}
+c^{11}_{56}
$;\ \ 
$ 1$,
$ -\sqrt{2}$,
$ -c_{7}^{3}$,
$ -c_{7}^{3}$,
$ \xi_{7}^{3}$,
$ \xi_{7}^{3}$,
$ -\sqrt{2}c_{14}^{1}$,
$  -c^{3}_{56}
-c^{5}_{56}
+c^{9}_{56}
-c^{11}_{56}
$;\ \ 
$0$,
$ \sqrt{2}c_{14}^{1}$,
$ -\sqrt{2}c_{14}^{1}$,
$  c^{3}_{56}
+c^{5}_{56}
-c^{9}_{56}
+c^{11}_{56}
$,
$  -c^{3}_{56}
-c^{5}_{56}
+c^{9}_{56}
-c^{11}_{56}
$,
$0$,
$0$;\ \ 
$ -\xi_{7}^{3}$,
$ -\xi_{7}^{3}$,
$ 1$,
$ 1$,
$  -c^{3}_{56}
-c^{5}_{56}
+c^{9}_{56}
-c^{11}_{56}
$,
$ \sqrt{2}$;\ \ 
$ -\xi_{7}^{3}$,
$ 1$,
$ 1$,
$  c^{3}_{56}
+c^{5}_{56}
-c^{9}_{56}
+c^{11}_{56}
$,
$ -\sqrt{2}$;\ \ 
$ c_{7}^{3}$,
$ c_{7}^{3}$,
$ \sqrt{2}$,
$ -\sqrt{2}c_{14}^{1}$;\ \ 
$ c_{7}^{3}$,
$ -\sqrt{2}$,
$ \sqrt{2}c_{14}^{1}$;\ \ 
$0$,
$0$;\ \ 
$0$)

Factors = $3_{\frac{3}{2},4.}^{16,553}\boxtimes 3_{\frac{48}{7},9.295}^{7,790}$

  \vskip 2ex

\noindent48. $9_{\frac{37}{14},37.18}^{112,168}$ \irep{532}:\ \ 
$d_i$ = ($1.0$,
$1.0$,
$1.414$,
$1.801$,
$1.801$,
$2.246$,
$2.246$,
$2.548$,
$3.177$) 

\vskip 0.7ex
\hangindent=3em \hangafter=1
$D^2= 37.183 = 
 24+12c^{1}_{7}
+4c^{2}_{7}
$

\vskip 0.7ex
\hangindent=3em \hangafter=1
$T = ( 0,
\frac{1}{2},
\frac{3}{16},
\frac{6}{7},
\frac{5}{14},
\frac{2}{7},
\frac{11}{14},
\frac{5}{112},
\frac{53}{112} )
$,

\vskip 0.7ex
\hangindent=3em \hangafter=1
$S$ = ($ 1$,
$ 1$,
$ \sqrt{2}$,
$ -c_{7}^{3}$,
$ -c_{7}^{3}$,
$ \xi_{7}^{3}$,
$ \xi_{7}^{3}$,
$ \sqrt{2}c_{14}^{1}$,
$  c^{3}_{56}
+c^{5}_{56}
-c^{9}_{56}
+c^{11}_{56}
$;\ \ 
$ 1$,
$ -\sqrt{2}$,
$ -c_{7}^{3}$,
$ -c_{7}^{3}$,
$ \xi_{7}^{3}$,
$ \xi_{7}^{3}$,
$ -\sqrt{2}c_{14}^{1}$,
$  -c^{3}_{56}
-c^{5}_{56}
+c^{9}_{56}
-c^{11}_{56}
$;\ \ 
$0$,
$ \sqrt{2}c_{14}^{1}$,
$ -\sqrt{2}c_{14}^{1}$,
$  c^{3}_{56}
+c^{5}_{56}
-c^{9}_{56}
+c^{11}_{56}
$,
$  -c^{3}_{56}
-c^{5}_{56}
+c^{9}_{56}
-c^{11}_{56}
$,
$0$,
$0$;\ \ 
$ -\xi_{7}^{3}$,
$ -\xi_{7}^{3}$,
$ 1$,
$ 1$,
$  -c^{3}_{56}
-c^{5}_{56}
+c^{9}_{56}
-c^{11}_{56}
$,
$ \sqrt{2}$;\ \ 
$ -\xi_{7}^{3}$,
$ 1$,
$ 1$,
$  c^{3}_{56}
+c^{5}_{56}
-c^{9}_{56}
+c^{11}_{56}
$,
$ -\sqrt{2}$;\ \ 
$ c_{7}^{3}$,
$ c_{7}^{3}$,
$ \sqrt{2}$,
$ -\sqrt{2}c_{14}^{1}$;\ \ 
$ c_{7}^{3}$,
$ -\sqrt{2}$,
$ \sqrt{2}c_{14}^{1}$;\ \ 
$0$,
$0$;\ \ 
$0$)

Factors = $3_{\frac{3}{2},4.}^{16,553}\boxtimes 3_{\frac{8}{7},9.295}^{7,245}$

  \vskip 2ex

\noindent49. $9_{\frac{19}{14},37.18}^{112,135}$ \irep{532}:\ \ 
$d_i$ = ($1.0$,
$1.0$,
$1.414$,
$1.801$,
$1.801$,
$2.246$,
$2.246$,
$2.548$,
$3.177$) 

\vskip 0.7ex
\hangindent=3em \hangafter=1
$D^2= 37.183 = 
 24+12c^{1}_{7}
+4c^{2}_{7}
$

\vskip 0.7ex
\hangindent=3em \hangafter=1
$T = ( 0,
\frac{1}{2},
\frac{5}{16},
\frac{1}{7},
\frac{9}{14},
\frac{5}{7},
\frac{3}{14},
\frac{51}{112},
\frac{3}{112} )
$,

\vskip 0.7ex
\hangindent=3em \hangafter=1
$S$ = ($ 1$,
$ 1$,
$ \sqrt{2}$,
$ -c_{7}^{3}$,
$ -c_{7}^{3}$,
$ \xi_{7}^{3}$,
$ \xi_{7}^{3}$,
$ \sqrt{2}c_{14}^{1}$,
$  c^{3}_{56}
+c^{5}_{56}
-c^{9}_{56}
+c^{11}_{56}
$;\ \ 
$ 1$,
$ -\sqrt{2}$,
$ -c_{7}^{3}$,
$ -c_{7}^{3}$,
$ \xi_{7}^{3}$,
$ \xi_{7}^{3}$,
$ -\sqrt{2}c_{14}^{1}$,
$  -c^{3}_{56}
-c^{5}_{56}
+c^{9}_{56}
-c^{11}_{56}
$;\ \ 
$0$,
$ \sqrt{2}c_{14}^{1}$,
$ -\sqrt{2}c_{14}^{1}$,
$  c^{3}_{56}
+c^{5}_{56}
-c^{9}_{56}
+c^{11}_{56}
$,
$  -c^{3}_{56}
-c^{5}_{56}
+c^{9}_{56}
-c^{11}_{56}
$,
$0$,
$0$;\ \ 
$ -\xi_{7}^{3}$,
$ -\xi_{7}^{3}$,
$ 1$,
$ 1$,
$  -c^{3}_{56}
-c^{5}_{56}
+c^{9}_{56}
-c^{11}_{56}
$,
$ \sqrt{2}$;\ \ 
$ -\xi_{7}^{3}$,
$ 1$,
$ 1$,
$  c^{3}_{56}
+c^{5}_{56}
-c^{9}_{56}
+c^{11}_{56}
$,
$ -\sqrt{2}$;\ \ 
$ c_{7}^{3}$,
$ c_{7}^{3}$,
$ \sqrt{2}$,
$ -\sqrt{2}c_{14}^{1}$;\ \ 
$ c_{7}^{3}$,
$ -\sqrt{2}$,
$ \sqrt{2}c_{14}^{1}$;\ \ 
$0$,
$0$;\ \ 
$0$)

Factors = $3_{\frac{5}{2},4.}^{16,465}\boxtimes 3_{\frac{48}{7},9.295}^{7,790}$

  \vskip 2ex

\noindent50. $9_{\frac{51}{14},37.18}^{112,413}$ \irep{532}:\ \ 
$d_i$ = ($1.0$,
$1.0$,
$1.414$,
$1.801$,
$1.801$,
$2.246$,
$2.246$,
$2.548$,
$3.177$) 

\vskip 0.7ex
\hangindent=3em \hangafter=1
$D^2= 37.183 = 
 24+12c^{1}_{7}
+4c^{2}_{7}
$

\vskip 0.7ex
\hangindent=3em \hangafter=1
$T = ( 0,
\frac{1}{2},
\frac{5}{16},
\frac{6}{7},
\frac{5}{14},
\frac{2}{7},
\frac{11}{14},
\frac{19}{112},
\frac{67}{112} )
$,

\vskip 0.7ex
\hangindent=3em \hangafter=1
$S$ = ($ 1$,
$ 1$,
$ \sqrt{2}$,
$ -c_{7}^{3}$,
$ -c_{7}^{3}$,
$ \xi_{7}^{3}$,
$ \xi_{7}^{3}$,
$ \sqrt{2}c_{14}^{1}$,
$  c^{3}_{56}
+c^{5}_{56}
-c^{9}_{56}
+c^{11}_{56}
$;\ \ 
$ 1$,
$ -\sqrt{2}$,
$ -c_{7}^{3}$,
$ -c_{7}^{3}$,
$ \xi_{7}^{3}$,
$ \xi_{7}^{3}$,
$ -\sqrt{2}c_{14}^{1}$,
$  -c^{3}_{56}
-c^{5}_{56}
+c^{9}_{56}
-c^{11}_{56}
$;\ \ 
$0$,
$ \sqrt{2}c_{14}^{1}$,
$ -\sqrt{2}c_{14}^{1}$,
$  c^{3}_{56}
+c^{5}_{56}
-c^{9}_{56}
+c^{11}_{56}
$,
$  -c^{3}_{56}
-c^{5}_{56}
+c^{9}_{56}
-c^{11}_{56}
$,
$0$,
$0$;\ \ 
$ -\xi_{7}^{3}$,
$ -\xi_{7}^{3}$,
$ 1$,
$ 1$,
$  -c^{3}_{56}
-c^{5}_{56}
+c^{9}_{56}
-c^{11}_{56}
$,
$ \sqrt{2}$;\ \ 
$ -\xi_{7}^{3}$,
$ 1$,
$ 1$,
$  c^{3}_{56}
+c^{5}_{56}
-c^{9}_{56}
+c^{11}_{56}
$,
$ -\sqrt{2}$;\ \ 
$ c_{7}^{3}$,
$ c_{7}^{3}$,
$ \sqrt{2}$,
$ -\sqrt{2}c_{14}^{1}$;\ \ 
$ c_{7}^{3}$,
$ -\sqrt{2}$,
$ \sqrt{2}c_{14}^{1}$;\ \ 
$0$,
$0$;\ \ 
$0$)

Factors = $3_{\frac{5}{2},4.}^{16,465}\boxtimes 3_{\frac{8}{7},9.295}^{7,245}$

  \vskip 2ex

\noindent51. $9_{\frac{33}{14},37.18}^{112,826}$ \irep{532}:\ \ 
$d_i$ = ($1.0$,
$1.0$,
$1.414$,
$1.801$,
$1.801$,
$2.246$,
$2.246$,
$2.548$,
$3.177$) 

\vskip 0.7ex
\hangindent=3em \hangafter=1
$D^2= 37.183 = 
 24+12c^{1}_{7}
+4c^{2}_{7}
$

\vskip 0.7ex
\hangindent=3em \hangafter=1
$T = ( 0,
\frac{1}{2},
\frac{7}{16},
\frac{1}{7},
\frac{9}{14},
\frac{5}{7},
\frac{3}{14},
\frac{65}{112},
\frac{17}{112} )
$,

\vskip 0.7ex
\hangindent=3em \hangafter=1
$S$ = ($ 1$,
$ 1$,
$ \sqrt{2}$,
$ -c_{7}^{3}$,
$ -c_{7}^{3}$,
$ \xi_{7}^{3}$,
$ \xi_{7}^{3}$,
$ \sqrt{2}c_{14}^{1}$,
$  c^{3}_{56}
+c^{5}_{56}
-c^{9}_{56}
+c^{11}_{56}
$;\ \ 
$ 1$,
$ -\sqrt{2}$,
$ -c_{7}^{3}$,
$ -c_{7}^{3}$,
$ \xi_{7}^{3}$,
$ \xi_{7}^{3}$,
$ -\sqrt{2}c_{14}^{1}$,
$  -c^{3}_{56}
-c^{5}_{56}
+c^{9}_{56}
-c^{11}_{56}
$;\ \ 
$0$,
$ \sqrt{2}c_{14}^{1}$,
$ -\sqrt{2}c_{14}^{1}$,
$  c^{3}_{56}
+c^{5}_{56}
-c^{9}_{56}
+c^{11}_{56}
$,
$  -c^{3}_{56}
-c^{5}_{56}
+c^{9}_{56}
-c^{11}_{56}
$,
$0$,
$0$;\ \ 
$ -\xi_{7}^{3}$,
$ -\xi_{7}^{3}$,
$ 1$,
$ 1$,
$  -c^{3}_{56}
-c^{5}_{56}
+c^{9}_{56}
-c^{11}_{56}
$,
$ \sqrt{2}$;\ \ 
$ -\xi_{7}^{3}$,
$ 1$,
$ 1$,
$  c^{3}_{56}
+c^{5}_{56}
-c^{9}_{56}
+c^{11}_{56}
$,
$ -\sqrt{2}$;\ \ 
$ c_{7}^{3}$,
$ c_{7}^{3}$,
$ \sqrt{2}$,
$ -\sqrt{2}c_{14}^{1}$;\ \ 
$ c_{7}^{3}$,
$ -\sqrt{2}$,
$ \sqrt{2}c_{14}^{1}$;\ \ 
$0$,
$0$;\ \ 
$0$)

Factors = $3_{\frac{7}{2},4.}^{16,332}\boxtimes 3_{\frac{48}{7},9.295}^{7,790}$

  \vskip 2ex

\noindent52. $9_{\frac{65}{14},37.18}^{112,121}$ \irep{532}:\ \ 
$d_i$ = ($1.0$,
$1.0$,
$1.414$,
$1.801$,
$1.801$,
$2.246$,
$2.246$,
$2.548$,
$3.177$) 

\vskip 0.7ex
\hangindent=3em \hangafter=1
$D^2= 37.183 = 
 24+12c^{1}_{7}
+4c^{2}_{7}
$

\vskip 0.7ex
\hangindent=3em \hangafter=1
$T = ( 0,
\frac{1}{2},
\frac{7}{16},
\frac{6}{7},
\frac{5}{14},
\frac{2}{7},
\frac{11}{14},
\frac{33}{112},
\frac{81}{112} )
$,

\vskip 0.7ex
\hangindent=3em \hangafter=1
$S$ = ($ 1$,
$ 1$,
$ \sqrt{2}$,
$ -c_{7}^{3}$,
$ -c_{7}^{3}$,
$ \xi_{7}^{3}$,
$ \xi_{7}^{3}$,
$ \sqrt{2}c_{14}^{1}$,
$  c^{3}_{56}
+c^{5}_{56}
-c^{9}_{56}
+c^{11}_{56}
$;\ \ 
$ 1$,
$ -\sqrt{2}$,
$ -c_{7}^{3}$,
$ -c_{7}^{3}$,
$ \xi_{7}^{3}$,
$ \xi_{7}^{3}$,
$ -\sqrt{2}c_{14}^{1}$,
$  -c^{3}_{56}
-c^{5}_{56}
+c^{9}_{56}
-c^{11}_{56}
$;\ \ 
$0$,
$ \sqrt{2}c_{14}^{1}$,
$ -\sqrt{2}c_{14}^{1}$,
$  c^{3}_{56}
+c^{5}_{56}
-c^{9}_{56}
+c^{11}_{56}
$,
$  -c^{3}_{56}
-c^{5}_{56}
+c^{9}_{56}
-c^{11}_{56}
$,
$0$,
$0$;\ \ 
$ -\xi_{7}^{3}$,
$ -\xi_{7}^{3}$,
$ 1$,
$ 1$,
$  -c^{3}_{56}
-c^{5}_{56}
+c^{9}_{56}
-c^{11}_{56}
$,
$ \sqrt{2}$;\ \ 
$ -\xi_{7}^{3}$,
$ 1$,
$ 1$,
$  c^{3}_{56}
+c^{5}_{56}
-c^{9}_{56}
+c^{11}_{56}
$,
$ -\sqrt{2}$;\ \ 
$ c_{7}^{3}$,
$ c_{7}^{3}$,
$ \sqrt{2}$,
$ -\sqrt{2}c_{14}^{1}$;\ \ 
$ c_{7}^{3}$,
$ -\sqrt{2}$,
$ \sqrt{2}c_{14}^{1}$;\ \ 
$0$,
$0$;\ \ 
$0$)

Factors = $3_{\frac{7}{2},4.}^{16,332}\boxtimes 3_{\frac{8}{7},9.295}^{7,245}$

  \vskip 2ex

\noindent53. $9_{\frac{47}{14},37.18}^{112,696}$ \irep{532}:\ \ 
$d_i$ = ($1.0$,
$1.0$,
$1.414$,
$1.801$,
$1.801$,
$2.246$,
$2.246$,
$2.548$,
$3.177$) 

\vskip 0.7ex
\hangindent=3em \hangafter=1
$D^2= 37.183 = 
 24+12c^{1}_{7}
+4c^{2}_{7}
$

\vskip 0.7ex
\hangindent=3em \hangafter=1
$T = ( 0,
\frac{1}{2},
\frac{9}{16},
\frac{1}{7},
\frac{9}{14},
\frac{5}{7},
\frac{3}{14},
\frac{79}{112},
\frac{31}{112} )
$,

\vskip 0.7ex
\hangindent=3em \hangafter=1
$S$ = ($ 1$,
$ 1$,
$ \sqrt{2}$,
$ -c_{7}^{3}$,
$ -c_{7}^{3}$,
$ \xi_{7}^{3}$,
$ \xi_{7}^{3}$,
$ \sqrt{2}c_{14}^{1}$,
$  c^{3}_{56}
+c^{5}_{56}
-c^{9}_{56}
+c^{11}_{56}
$;\ \ 
$ 1$,
$ -\sqrt{2}$,
$ -c_{7}^{3}$,
$ -c_{7}^{3}$,
$ \xi_{7}^{3}$,
$ \xi_{7}^{3}$,
$ -\sqrt{2}c_{14}^{1}$,
$  -c^{3}_{56}
-c^{5}_{56}
+c^{9}_{56}
-c^{11}_{56}
$;\ \ 
$0$,
$ \sqrt{2}c_{14}^{1}$,
$ -\sqrt{2}c_{14}^{1}$,
$  c^{3}_{56}
+c^{5}_{56}
-c^{9}_{56}
+c^{11}_{56}
$,
$  -c^{3}_{56}
-c^{5}_{56}
+c^{9}_{56}
-c^{11}_{56}
$,
$0$,
$0$;\ \ 
$ -\xi_{7}^{3}$,
$ -\xi_{7}^{3}$,
$ 1$,
$ 1$,
$  -c^{3}_{56}
-c^{5}_{56}
+c^{9}_{56}
-c^{11}_{56}
$,
$ \sqrt{2}$;\ \ 
$ -\xi_{7}^{3}$,
$ 1$,
$ 1$,
$  c^{3}_{56}
+c^{5}_{56}
-c^{9}_{56}
+c^{11}_{56}
$,
$ -\sqrt{2}$;\ \ 
$ c_{7}^{3}$,
$ c_{7}^{3}$,
$ \sqrt{2}$,
$ -\sqrt{2}c_{14}^{1}$;\ \ 
$ c_{7}^{3}$,
$ -\sqrt{2}$,
$ \sqrt{2}c_{14}^{1}$;\ \ 
$0$,
$0$;\ \ 
$0$)

Factors = $3_{\frac{9}{2},4.}^{16,156}\boxtimes 3_{\frac{48}{7},9.295}^{7,790}$

  \vskip 2ex

\noindent54. $9_{\frac{79}{14},37.18}^{112,224}$ \irep{532}:\ \ 
$d_i$ = ($1.0$,
$1.0$,
$1.414$,
$1.801$,
$1.801$,
$2.246$,
$2.246$,
$2.548$,
$3.177$) 

\vskip 0.7ex
\hangindent=3em \hangafter=1
$D^2= 37.183 = 
 24+12c^{1}_{7}
+4c^{2}_{7}
$

\vskip 0.7ex
\hangindent=3em \hangafter=1
$T = ( 0,
\frac{1}{2},
\frac{9}{16},
\frac{6}{7},
\frac{5}{14},
\frac{2}{7},
\frac{11}{14},
\frac{47}{112},
\frac{95}{112} )
$,

\vskip 0.7ex
\hangindent=3em \hangafter=1
$S$ = ($ 1$,
$ 1$,
$ \sqrt{2}$,
$ -c_{7}^{3}$,
$ -c_{7}^{3}$,
$ \xi_{7}^{3}$,
$ \xi_{7}^{3}$,
$ \sqrt{2}c_{14}^{1}$,
$  c^{3}_{56}
+c^{5}_{56}
-c^{9}_{56}
+c^{11}_{56}
$;\ \ 
$ 1$,
$ -\sqrt{2}$,
$ -c_{7}^{3}$,
$ -c_{7}^{3}$,
$ \xi_{7}^{3}$,
$ \xi_{7}^{3}$,
$ -\sqrt{2}c_{14}^{1}$,
$  -c^{3}_{56}
-c^{5}_{56}
+c^{9}_{56}
-c^{11}_{56}
$;\ \ 
$0$,
$ \sqrt{2}c_{14}^{1}$,
$ -\sqrt{2}c_{14}^{1}$,
$  c^{3}_{56}
+c^{5}_{56}
-c^{9}_{56}
+c^{11}_{56}
$,
$  -c^{3}_{56}
-c^{5}_{56}
+c^{9}_{56}
-c^{11}_{56}
$,
$0$,
$0$;\ \ 
$ -\xi_{7}^{3}$,
$ -\xi_{7}^{3}$,
$ 1$,
$ 1$,
$  -c^{3}_{56}
-c^{5}_{56}
+c^{9}_{56}
-c^{11}_{56}
$,
$ \sqrt{2}$;\ \ 
$ -\xi_{7}^{3}$,
$ 1$,
$ 1$,
$  c^{3}_{56}
+c^{5}_{56}
-c^{9}_{56}
+c^{11}_{56}
$,
$ -\sqrt{2}$;\ \ 
$ c_{7}^{3}$,
$ c_{7}^{3}$,
$ \sqrt{2}$,
$ -\sqrt{2}c_{14}^{1}$;\ \ 
$ c_{7}^{3}$,
$ -\sqrt{2}$,
$ \sqrt{2}c_{14}^{1}$;\ \ 
$0$,
$0$;\ \ 
$0$)

Factors = $3_{\frac{9}{2},4.}^{16,156}\boxtimes 3_{\frac{8}{7},9.295}^{7,245}$

  \vskip 2ex

\noindent55. $9_{\frac{61}{14},37.18}^{112,909}$ \irep{532}:\ \ 
$d_i$ = ($1.0$,
$1.0$,
$1.414$,
$1.801$,
$1.801$,
$2.246$,
$2.246$,
$2.548$,
$3.177$) 

\vskip 0.7ex
\hangindent=3em \hangafter=1
$D^2= 37.183 = 
 24+12c^{1}_{7}
+4c^{2}_{7}
$

\vskip 0.7ex
\hangindent=3em \hangafter=1
$T = ( 0,
\frac{1}{2},
\frac{11}{16},
\frac{1}{7},
\frac{9}{14},
\frac{5}{7},
\frac{3}{14},
\frac{93}{112},
\frac{45}{112} )
$,

\vskip 0.7ex
\hangindent=3em \hangafter=1
$S$ = ($ 1$,
$ 1$,
$ \sqrt{2}$,
$ -c_{7}^{3}$,
$ -c_{7}^{3}$,
$ \xi_{7}^{3}$,
$ \xi_{7}^{3}$,
$ \sqrt{2}c_{14}^{1}$,
$  c^{3}_{56}
+c^{5}_{56}
-c^{9}_{56}
+c^{11}_{56}
$;\ \ 
$ 1$,
$ -\sqrt{2}$,
$ -c_{7}^{3}$,
$ -c_{7}^{3}$,
$ \xi_{7}^{3}$,
$ \xi_{7}^{3}$,
$ -\sqrt{2}c_{14}^{1}$,
$  -c^{3}_{56}
-c^{5}_{56}
+c^{9}_{56}
-c^{11}_{56}
$;\ \ 
$0$,
$ \sqrt{2}c_{14}^{1}$,
$ -\sqrt{2}c_{14}^{1}$,
$  c^{3}_{56}
+c^{5}_{56}
-c^{9}_{56}
+c^{11}_{56}
$,
$  -c^{3}_{56}
-c^{5}_{56}
+c^{9}_{56}
-c^{11}_{56}
$,
$0$,
$0$;\ \ 
$ -\xi_{7}^{3}$,
$ -\xi_{7}^{3}$,
$ 1$,
$ 1$,
$  -c^{3}_{56}
-c^{5}_{56}
+c^{9}_{56}
-c^{11}_{56}
$,
$ \sqrt{2}$;\ \ 
$ -\xi_{7}^{3}$,
$ 1$,
$ 1$,
$  c^{3}_{56}
+c^{5}_{56}
-c^{9}_{56}
+c^{11}_{56}
$,
$ -\sqrt{2}$;\ \ 
$ c_{7}^{3}$,
$ c_{7}^{3}$,
$ \sqrt{2}$,
$ -\sqrt{2}c_{14}^{1}$;\ \ 
$ c_{7}^{3}$,
$ -\sqrt{2}$,
$ \sqrt{2}c_{14}^{1}$;\ \ 
$0$,
$0$;\ \ 
$0$)

Factors = $3_{\frac{11}{2},4.}^{16,648}\boxtimes 3_{\frac{48}{7},9.295}^{7,790}$

  \vskip 2ex

\noindent56. $9_{\frac{93}{14},37.18}^{112,349}$ \irep{532}:\ \ 
$d_i$ = ($1.0$,
$1.0$,
$1.414$,
$1.801$,
$1.801$,
$2.246$,
$2.246$,
$2.548$,
$3.177$) 

\vskip 0.7ex
\hangindent=3em \hangafter=1
$D^2= 37.183 = 
 24+12c^{1}_{7}
+4c^{2}_{7}
$

\vskip 0.7ex
\hangindent=3em \hangafter=1
$T = ( 0,
\frac{1}{2},
\frac{11}{16},
\frac{6}{7},
\frac{5}{14},
\frac{2}{7},
\frac{11}{14},
\frac{61}{112},
\frac{109}{112} )
$,

\vskip 0.7ex
\hangindent=3em \hangafter=1
$S$ = ($ 1$,
$ 1$,
$ \sqrt{2}$,
$ -c_{7}^{3}$,
$ -c_{7}^{3}$,
$ \xi_{7}^{3}$,
$ \xi_{7}^{3}$,
$ \sqrt{2}c_{14}^{1}$,
$  c^{3}_{56}
+c^{5}_{56}
-c^{9}_{56}
+c^{11}_{56}
$;\ \ 
$ 1$,
$ -\sqrt{2}$,
$ -c_{7}^{3}$,
$ -c_{7}^{3}$,
$ \xi_{7}^{3}$,
$ \xi_{7}^{3}$,
$ -\sqrt{2}c_{14}^{1}$,
$  -c^{3}_{56}
-c^{5}_{56}
+c^{9}_{56}
-c^{11}_{56}
$;\ \ 
$0$,
$ \sqrt{2}c_{14}^{1}$,
$ -\sqrt{2}c_{14}^{1}$,
$  c^{3}_{56}
+c^{5}_{56}
-c^{9}_{56}
+c^{11}_{56}
$,
$  -c^{3}_{56}
-c^{5}_{56}
+c^{9}_{56}
-c^{11}_{56}
$,
$0$,
$0$;\ \ 
$ -\xi_{7}^{3}$,
$ -\xi_{7}^{3}$,
$ 1$,
$ 1$,
$  -c^{3}_{56}
-c^{5}_{56}
+c^{9}_{56}
-c^{11}_{56}
$,
$ \sqrt{2}$;\ \ 
$ -\xi_{7}^{3}$,
$ 1$,
$ 1$,
$  c^{3}_{56}
+c^{5}_{56}
-c^{9}_{56}
+c^{11}_{56}
$,
$ -\sqrt{2}$;\ \ 
$ c_{7}^{3}$,
$ c_{7}^{3}$,
$ \sqrt{2}$,
$ -\sqrt{2}c_{14}^{1}$;\ \ 
$ c_{7}^{3}$,
$ -\sqrt{2}$,
$ \sqrt{2}c_{14}^{1}$;\ \ 
$0$,
$0$;\ \ 
$0$)

Factors = $3_{\frac{11}{2},4.}^{16,648}\boxtimes 3_{\frac{8}{7},9.295}^{7,245}$

  \vskip 2ex

\noindent57. $9_{\frac{75}{14},37.18}^{112,211}$ \irep{532}:\ \ 
$d_i$ = ($1.0$,
$1.0$,
$1.414$,
$1.801$,
$1.801$,
$2.246$,
$2.246$,
$2.548$,
$3.177$) 

\vskip 0.7ex
\hangindent=3em \hangafter=1
$D^2= 37.183 = 
 24+12c^{1}_{7}
+4c^{2}_{7}
$

\vskip 0.7ex
\hangindent=3em \hangafter=1
$T = ( 0,
\frac{1}{2},
\frac{13}{16},
\frac{1}{7},
\frac{9}{14},
\frac{5}{7},
\frac{3}{14},
\frac{107}{112},
\frac{59}{112} )
$,

\vskip 0.7ex
\hangindent=3em \hangafter=1
$S$ = ($ 1$,
$ 1$,
$ \sqrt{2}$,
$ -c_{7}^{3}$,
$ -c_{7}^{3}$,
$ \xi_{7}^{3}$,
$ \xi_{7}^{3}$,
$ \sqrt{2}c_{14}^{1}$,
$  c^{3}_{56}
+c^{5}_{56}
-c^{9}_{56}
+c^{11}_{56}
$;\ \ 
$ 1$,
$ -\sqrt{2}$,
$ -c_{7}^{3}$,
$ -c_{7}^{3}$,
$ \xi_{7}^{3}$,
$ \xi_{7}^{3}$,
$ -\sqrt{2}c_{14}^{1}$,
$  -c^{3}_{56}
-c^{5}_{56}
+c^{9}_{56}
-c^{11}_{56}
$;\ \ 
$0$,
$ \sqrt{2}c_{14}^{1}$,
$ -\sqrt{2}c_{14}^{1}$,
$  c^{3}_{56}
+c^{5}_{56}
-c^{9}_{56}
+c^{11}_{56}
$,
$  -c^{3}_{56}
-c^{5}_{56}
+c^{9}_{56}
-c^{11}_{56}
$,
$0$,
$0$;\ \ 
$ -\xi_{7}^{3}$,
$ -\xi_{7}^{3}$,
$ 1$,
$ 1$,
$  -c^{3}_{56}
-c^{5}_{56}
+c^{9}_{56}
-c^{11}_{56}
$,
$ \sqrt{2}$;\ \ 
$ -\xi_{7}^{3}$,
$ 1$,
$ 1$,
$  c^{3}_{56}
+c^{5}_{56}
-c^{9}_{56}
+c^{11}_{56}
$,
$ -\sqrt{2}$;\ \ 
$ c_{7}^{3}$,
$ c_{7}^{3}$,
$ \sqrt{2}$,
$ -\sqrt{2}c_{14}^{1}$;\ \ 
$ c_{7}^{3}$,
$ -\sqrt{2}$,
$ \sqrt{2}c_{14}^{1}$;\ \ 
$0$,
$0$;\ \ 
$0$)

Factors = $3_{\frac{13}{2},4.}^{16,330}\boxtimes 3_{\frac{48}{7},9.295}^{7,790}$

  \vskip 2ex

\noindent58. $9_{\frac{107}{14},37.18}^{112,117}$ \irep{532}:\ \ 
$d_i$ = ($1.0$,
$1.0$,
$1.414$,
$1.801$,
$1.801$,
$2.246$,
$2.246$,
$2.548$,
$3.177$) 

\vskip 0.7ex
\hangindent=3em \hangafter=1
$D^2= 37.183 = 
 24+12c^{1}_{7}
+4c^{2}_{7}
$

\vskip 0.7ex
\hangindent=3em \hangafter=1
$T = ( 0,
\frac{1}{2},
\frac{13}{16},
\frac{6}{7},
\frac{5}{14},
\frac{2}{7},
\frac{11}{14},
\frac{75}{112},
\frac{11}{112} )
$,

\vskip 0.7ex
\hangindent=3em \hangafter=1
$S$ = ($ 1$,
$ 1$,
$ \sqrt{2}$,
$ -c_{7}^{3}$,
$ -c_{7}^{3}$,
$ \xi_{7}^{3}$,
$ \xi_{7}^{3}$,
$ \sqrt{2}c_{14}^{1}$,
$  c^{3}_{56}
+c^{5}_{56}
-c^{9}_{56}
+c^{11}_{56}
$;\ \ 
$ 1$,
$ -\sqrt{2}$,
$ -c_{7}^{3}$,
$ -c_{7}^{3}$,
$ \xi_{7}^{3}$,
$ \xi_{7}^{3}$,
$ -\sqrt{2}c_{14}^{1}$,
$  -c^{3}_{56}
-c^{5}_{56}
+c^{9}_{56}
-c^{11}_{56}
$;\ \ 
$0$,
$ \sqrt{2}c_{14}^{1}$,
$ -\sqrt{2}c_{14}^{1}$,
$  c^{3}_{56}
+c^{5}_{56}
-c^{9}_{56}
+c^{11}_{56}
$,
$  -c^{3}_{56}
-c^{5}_{56}
+c^{9}_{56}
-c^{11}_{56}
$,
$0$,
$0$;\ \ 
$ -\xi_{7}^{3}$,
$ -\xi_{7}^{3}$,
$ 1$,
$ 1$,
$  -c^{3}_{56}
-c^{5}_{56}
+c^{9}_{56}
-c^{11}_{56}
$,
$ \sqrt{2}$;\ \ 
$ -\xi_{7}^{3}$,
$ 1$,
$ 1$,
$  c^{3}_{56}
+c^{5}_{56}
-c^{9}_{56}
+c^{11}_{56}
$,
$ -\sqrt{2}$;\ \ 
$ c_{7}^{3}$,
$ c_{7}^{3}$,
$ \sqrt{2}$,
$ -\sqrt{2}c_{14}^{1}$;\ \ 
$ c_{7}^{3}$,
$ -\sqrt{2}$,
$ \sqrt{2}c_{14}^{1}$;\ \ 
$0$,
$0$;\ \ 
$0$)

Factors = $3_{\frac{13}{2},4.}^{16,330}\boxtimes 3_{\frac{8}{7},9.295}^{7,245}$

  \vskip 2ex

\noindent59. $9_{\frac{89}{14},37.18}^{112,580}$ \irep{532}:\ \ 
$d_i$ = ($1.0$,
$1.0$,
$1.414$,
$1.801$,
$1.801$,
$2.246$,
$2.246$,
$2.548$,
$3.177$) 

\vskip 0.7ex
\hangindent=3em \hangafter=1
$D^2= 37.183 = 
 24+12c^{1}_{7}
+4c^{2}_{7}
$

\vskip 0.7ex
\hangindent=3em \hangafter=1
$T = ( 0,
\frac{1}{2},
\frac{15}{16},
\frac{1}{7},
\frac{9}{14},
\frac{5}{7},
\frac{3}{14},
\frac{9}{112},
\frac{73}{112} )
$,

\vskip 0.7ex
\hangindent=3em \hangafter=1
$S$ = ($ 1$,
$ 1$,
$ \sqrt{2}$,
$ -c_{7}^{3}$,
$ -c_{7}^{3}$,
$ \xi_{7}^{3}$,
$ \xi_{7}^{3}$,
$ \sqrt{2}c_{14}^{1}$,
$  c^{3}_{56}
+c^{5}_{56}
-c^{9}_{56}
+c^{11}_{56}
$;\ \ 
$ 1$,
$ -\sqrt{2}$,
$ -c_{7}^{3}$,
$ -c_{7}^{3}$,
$ \xi_{7}^{3}$,
$ \xi_{7}^{3}$,
$ -\sqrt{2}c_{14}^{1}$,
$  -c^{3}_{56}
-c^{5}_{56}
+c^{9}_{56}
-c^{11}_{56}
$;\ \ 
$0$,
$ \sqrt{2}c_{14}^{1}$,
$ -\sqrt{2}c_{14}^{1}$,
$  c^{3}_{56}
+c^{5}_{56}
-c^{9}_{56}
+c^{11}_{56}
$,
$  -c^{3}_{56}
-c^{5}_{56}
+c^{9}_{56}
-c^{11}_{56}
$,
$0$,
$0$;\ \ 
$ -\xi_{7}^{3}$,
$ -\xi_{7}^{3}$,
$ 1$,
$ 1$,
$  -c^{3}_{56}
-c^{5}_{56}
+c^{9}_{56}
-c^{11}_{56}
$,
$ \sqrt{2}$;\ \ 
$ -\xi_{7}^{3}$,
$ 1$,
$ 1$,
$  c^{3}_{56}
+c^{5}_{56}
-c^{9}_{56}
+c^{11}_{56}
$,
$ -\sqrt{2}$;\ \ 
$ c_{7}^{3}$,
$ c_{7}^{3}$,
$ \sqrt{2}$,
$ -\sqrt{2}c_{14}^{1}$;\ \ 
$ c_{7}^{3}$,
$ -\sqrt{2}$,
$ \sqrt{2}c_{14}^{1}$;\ \ 
$0$,
$0$;\ \ 
$0$)

Factors = $3_{\frac{15}{2},4.}^{16,639}\boxtimes 3_{\frac{48}{7},9.295}^{7,790}$

  \vskip 2ex

\noindent60. $9_{\frac{9}{14},37.18}^{112,207}$ \irep{532}:\ \ 
$d_i$ = ($1.0$,
$1.0$,
$1.414$,
$1.801$,
$1.801$,
$2.246$,
$2.246$,
$2.548$,
$3.177$) 

\vskip 0.7ex
\hangindent=3em \hangafter=1
$D^2= 37.183 = 
 24+12c^{1}_{7}
+4c^{2}_{7}
$

\vskip 0.7ex
\hangindent=3em \hangafter=1
$T = ( 0,
\frac{1}{2},
\frac{15}{16},
\frac{6}{7},
\frac{5}{14},
\frac{2}{7},
\frac{11}{14},
\frac{89}{112},
\frac{25}{112} )
$,

\vskip 0.7ex
\hangindent=3em \hangafter=1
$S$ = ($ 1$,
$ 1$,
$ \sqrt{2}$,
$ -c_{7}^{3}$,
$ -c_{7}^{3}$,
$ \xi_{7}^{3}$,
$ \xi_{7}^{3}$,
$ \sqrt{2}c_{14}^{1}$,
$  c^{3}_{56}
+c^{5}_{56}
-c^{9}_{56}
+c^{11}_{56}
$;\ \ 
$ 1$,
$ -\sqrt{2}$,
$ -c_{7}^{3}$,
$ -c_{7}^{3}$,
$ \xi_{7}^{3}$,
$ \xi_{7}^{3}$,
$ -\sqrt{2}c_{14}^{1}$,
$  -c^{3}_{56}
-c^{5}_{56}
+c^{9}_{56}
-c^{11}_{56}
$;\ \ 
$0$,
$ \sqrt{2}c_{14}^{1}$,
$ -\sqrt{2}c_{14}^{1}$,
$  c^{3}_{56}
+c^{5}_{56}
-c^{9}_{56}
+c^{11}_{56}
$,
$  -c^{3}_{56}
-c^{5}_{56}
+c^{9}_{56}
-c^{11}_{56}
$,
$0$,
$0$;\ \ 
$ -\xi_{7}^{3}$,
$ -\xi_{7}^{3}$,
$ 1$,
$ 1$,
$  -c^{3}_{56}
-c^{5}_{56}
+c^{9}_{56}
-c^{11}_{56}
$,
$ \sqrt{2}$;\ \ 
$ -\xi_{7}^{3}$,
$ 1$,
$ 1$,
$  c^{3}_{56}
+c^{5}_{56}
-c^{9}_{56}
+c^{11}_{56}
$,
$ -\sqrt{2}$;\ \ 
$ c_{7}^{3}$,
$ c_{7}^{3}$,
$ \sqrt{2}$,
$ -\sqrt{2}c_{14}^{1}$;\ \ 
$ c_{7}^{3}$,
$ -\sqrt{2}$,
$ \sqrt{2}c_{14}^{1}$;\ \ 
$0$,
$0$;\ \ 
$0$)

Factors = $3_{\frac{15}{2},4.}^{16,639}\boxtimes 3_{\frac{8}{7},9.295}^{7,245}$

  \vskip 2ex

\noindent61. $9_{2,44.}^{88,112}$ \irep{530}:\ \ 
$d_i$ = ($1.0$,
$1.0$,
$2.0$,
$2.0$,
$2.0$,
$2.0$,
$2.0$,
$3.316$,
$3.316$) 

\vskip 0.7ex
\hangindent=3em \hangafter=1
$D^2= 44.0 = 
44$

\vskip 0.7ex
\hangindent=3em \hangafter=1
$T = ( 0,
0,
\frac{1}{11},
\frac{3}{11},
\frac{4}{11},
\frac{5}{11},
\frac{9}{11},
\frac{1}{8},
\frac{5}{8} )
$,

\vskip 0.7ex
\hangindent=3em \hangafter=1
$S$ = ($ 1$,
$ 1$,
$ 2$,
$ 2$,
$ 2$,
$ 2$,
$ 2$,
$ \sqrt{11}$,
$ \sqrt{11}$;\ \ 
$ 1$,
$ 2$,
$ 2$,
$ 2$,
$ 2$,
$ 2$,
$ -\sqrt{11}$,
$ -\sqrt{11}$;\ \ 
$ 2c_{11}^{2}$,
$ 2c_{11}^{1}$,
$ 2c_{11}^{4}$,
$ 2c_{11}^{3}$,
$ 2c_{11}^{5}$,
$0$,
$0$;\ \ 
$ 2c_{11}^{5}$,
$ 2c_{11}^{2}$,
$ 2c_{11}^{4}$,
$ 2c_{11}^{3}$,
$0$,
$0$;\ \ 
$ 2c_{11}^{3}$,
$ 2c_{11}^{5}$,
$ 2c_{11}^{1}$,
$0$,
$0$;\ \ 
$ 2c_{11}^{1}$,
$ 2c_{11}^{2}$,
$0$,
$0$;\ \ 
$ 2c_{11}^{4}$,
$0$,
$0$;\ \ 
$ \sqrt{11}$,
$ -\sqrt{11}$;\ \ 
$ \sqrt{11}$)

  \vskip 2ex

\noindent62. $9_{2,44.}^{88,529}$ \irep{530}:\ \ 
$d_i$ = ($1.0$,
$1.0$,
$2.0$,
$2.0$,
$2.0$,
$2.0$,
$2.0$,
$3.316$,
$3.316$) 

\vskip 0.7ex
\hangindent=3em \hangafter=1
$D^2= 44.0 = 
44$

\vskip 0.7ex
\hangindent=3em \hangafter=1
$T = ( 0,
0,
\frac{1}{11},
\frac{3}{11},
\frac{4}{11},
\frac{5}{11},
\frac{9}{11},
\frac{3}{8},
\frac{7}{8} )
$,

\vskip 0.7ex
\hangindent=3em \hangafter=1
$S$ = ($ 1$,
$ 1$,
$ 2$,
$ 2$,
$ 2$,
$ 2$,
$ 2$,
$ \sqrt{11}$,
$ \sqrt{11}$;\ \ 
$ 1$,
$ 2$,
$ 2$,
$ 2$,
$ 2$,
$ 2$,
$ -\sqrt{11}$,
$ -\sqrt{11}$;\ \ 
$ 2c_{11}^{2}$,
$ 2c_{11}^{1}$,
$ 2c_{11}^{4}$,
$ 2c_{11}^{3}$,
$ 2c_{11}^{5}$,
$0$,
$0$;\ \ 
$ 2c_{11}^{5}$,
$ 2c_{11}^{2}$,
$ 2c_{11}^{4}$,
$ 2c_{11}^{3}$,
$0$,
$0$;\ \ 
$ 2c_{11}^{3}$,
$ 2c_{11}^{5}$,
$ 2c_{11}^{1}$,
$0$,
$0$;\ \ 
$ 2c_{11}^{1}$,
$ 2c_{11}^{2}$,
$0$,
$0$;\ \ 
$ 2c_{11}^{4}$,
$0$,
$0$;\ \ 
$ -\sqrt{11}$,
$ \sqrt{11}$;\ \ 
$ -\sqrt{11}$)

  \vskip 2ex

\noindent63. $9_{6,44.}^{88,870}$ \irep{530}:\ \ 
$d_i$ = ($1.0$,
$1.0$,
$2.0$,
$2.0$,
$2.0$,
$2.0$,
$2.0$,
$3.316$,
$3.316$) 

\vskip 0.7ex
\hangindent=3em \hangafter=1
$D^2= 44.0 = 
44$

\vskip 0.7ex
\hangindent=3em \hangafter=1
$T = ( 0,
0,
\frac{2}{11},
\frac{6}{11},
\frac{7}{11},
\frac{8}{11},
\frac{10}{11},
\frac{1}{8},
\frac{5}{8} )
$,

\vskip 0.7ex
\hangindent=3em \hangafter=1
$S$ = ($ 1$,
$ 1$,
$ 2$,
$ 2$,
$ 2$,
$ 2$,
$ 2$,
$ \sqrt{11}$,
$ \sqrt{11}$;\ \ 
$ 1$,
$ 2$,
$ 2$,
$ 2$,
$ 2$,
$ 2$,
$ -\sqrt{11}$,
$ -\sqrt{11}$;\ \ 
$ 2c_{11}^{4}$,
$ 2c_{11}^{2}$,
$ 2c_{11}^{1}$,
$ 2c_{11}^{3}$,
$ 2c_{11}^{5}$,
$0$,
$0$;\ \ 
$ 2c_{11}^{1}$,
$ 2c_{11}^{5}$,
$ 2c_{11}^{4}$,
$ 2c_{11}^{3}$,
$0$,
$0$;\ \ 
$ 2c_{11}^{3}$,
$ 2c_{11}^{2}$,
$ 2c_{11}^{4}$,
$0$,
$0$;\ \ 
$ 2c_{11}^{5}$,
$ 2c_{11}^{1}$,
$0$,
$0$;\ \ 
$ 2c_{11}^{2}$,
$0$,
$0$;\ \ 
$ -\sqrt{11}$,
$ \sqrt{11}$;\ \ 
$ -\sqrt{11}$)

  \vskip 2ex

\noindent64. $9_{6,44.}^{88,252}$ \irep{530}:\ \ 
$d_i$ = ($1.0$,
$1.0$,
$2.0$,
$2.0$,
$2.0$,
$2.0$,
$2.0$,
$3.316$,
$3.316$) 

\vskip 0.7ex
\hangindent=3em \hangafter=1
$D^2= 44.0 = 
44$

\vskip 0.7ex
\hangindent=3em \hangafter=1
$T = ( 0,
0,
\frac{2}{11},
\frac{6}{11},
\frac{7}{11},
\frac{8}{11},
\frac{10}{11},
\frac{3}{8},
\frac{7}{8} )
$,

\vskip 0.7ex
\hangindent=3em \hangafter=1
$S$ = ($ 1$,
$ 1$,
$ 2$,
$ 2$,
$ 2$,
$ 2$,
$ 2$,
$ \sqrt{11}$,
$ \sqrt{11}$;\ \ 
$ 1$,
$ 2$,
$ 2$,
$ 2$,
$ 2$,
$ 2$,
$ -\sqrt{11}$,
$ -\sqrt{11}$;\ \ 
$ 2c_{11}^{4}$,
$ 2c_{11}^{2}$,
$ 2c_{11}^{1}$,
$ 2c_{11}^{3}$,
$ 2c_{11}^{5}$,
$0$,
$0$;\ \ 
$ 2c_{11}^{1}$,
$ 2c_{11}^{5}$,
$ 2c_{11}^{4}$,
$ 2c_{11}^{3}$,
$0$,
$0$;\ \ 
$ 2c_{11}^{3}$,
$ 2c_{11}^{2}$,
$ 2c_{11}^{4}$,
$0$,
$0$;\ \ 
$ 2c_{11}^{5}$,
$ 2c_{11}^{1}$,
$0$,
$0$;\ \ 
$ 2c_{11}^{2}$,
$0$,
$0$;\ \ 
$ \sqrt{11}$,
$ -\sqrt{11}$;\ \ 
$ \sqrt{11}$)

  \vskip 2ex

\noindent65. $9_{\frac{12}{5},52.36}^{40,304}$ \irep{516}:\ \ 
$d_i$ = ($1.0$,
$1.0$,
$1.902$,
$1.902$,
$2.618$,
$2.618$,
$3.77$,
$3.77$,
$3.236$) 

\vskip 0.7ex
\hangindent=3em \hangafter=1
$D^2= 52.360 = 
 30+10\sqrt{5} $

\vskip 0.7ex
\hangindent=3em \hangafter=1
$T = ( 0,
0,
\frac{3}{40},
\frac{23}{40},
\frac{1}{5},
\frac{1}{5},
\frac{3}{8},
\frac{7}{8},
\frac{3}{5} )
$,

\vskip 0.7ex
\hangindent=3em \hangafter=1
$S$ = ($ 1$,
$ 1$,
$ c_{20}^{1}$,
$ c_{20}^{1}$,
$ \frac{3+\sqrt{5}}{2}$,
$ \frac{3+\sqrt{5}}{2}$,
$ \frac{1+\sqrt{5}}{2}c_{20}^{1}$,
$ \frac{1+\sqrt{5}}{2}c_{20}^{1}$,
$  1+\sqrt{5} $;\ \ 
$ 1$,
$ -c_{20}^{1}$,
$ -c_{20}^{1}$,
$ \frac{3+\sqrt{5}}{2}$,
$ \frac{3+\sqrt{5}}{2}$,
$ -\frac{1+\sqrt{5}}{2}c_{20}^{1}$,
$ -\frac{1+\sqrt{5}}{2}c_{20}^{1}$,
$  1+\sqrt{5} $;\ \ 
$ \frac{1+\sqrt{5}}{2}c_{20}^{1}$,
$ -\frac{1+\sqrt{5}}{2}c_{20}^{1}$,
$ -\frac{1+\sqrt{5}}{2}c_{20}^{1}$,
$ \frac{1+\sqrt{5}}{2}c_{20}^{1}$,
$ c_{20}^{1}$,
$ -c_{20}^{1}$,
$0$;\ \ 
$ \frac{1+\sqrt{5}}{2}c_{20}^{1}$,
$ -\frac{1+\sqrt{5}}{2}c_{20}^{1}$,
$ \frac{1+\sqrt{5}}{2}c_{20}^{1}$,
$ -c_{20}^{1}$,
$ c_{20}^{1}$,
$0$;\ \ 
$ 1$,
$ 1$,
$ c_{20}^{1}$,
$ c_{20}^{1}$,
$  -1-\sqrt{5} $;\ \ 
$ 1$,
$ -c_{20}^{1}$,
$ -c_{20}^{1}$,
$  -1-\sqrt{5} $;\ \ 
$ -\frac{1+\sqrt{5}}{2}c_{20}^{1}$,
$ \frac{1+\sqrt{5}}{2}c_{20}^{1}$,
$0$;\ \ 
$ -\frac{1+\sqrt{5}}{2}c_{20}^{1}$,
$0$;\ \ 
$  1+\sqrt{5} $)

  \vskip 2ex

\noindent66. $9_{\frac{28}{5},52.36}^{40,247}$ \irep{516}:\ \ 
$d_i$ = ($1.0$,
$1.0$,
$1.902$,
$1.902$,
$2.618$,
$2.618$,
$3.77$,
$3.77$,
$3.236$) 

\vskip 0.7ex
\hangindent=3em \hangafter=1
$D^2= 52.360 = 
 30+10\sqrt{5} $

\vskip 0.7ex
\hangindent=3em \hangafter=1
$T = ( 0,
0,
\frac{7}{40},
\frac{27}{40},
\frac{4}{5},
\frac{4}{5},
\frac{3}{8},
\frac{7}{8},
\frac{2}{5} )
$,

\vskip 0.7ex
\hangindent=3em \hangafter=1
$S$ = ($ 1$,
$ 1$,
$ c_{20}^{1}$,
$ c_{20}^{1}$,
$ \frac{3+\sqrt{5}}{2}$,
$ \frac{3+\sqrt{5}}{2}$,
$ \frac{1+\sqrt{5}}{2}c_{20}^{1}$,
$ \frac{1+\sqrt{5}}{2}c_{20}^{1}$,
$  1+\sqrt{5} $;\ \ 
$ 1$,
$ -c_{20}^{1}$,
$ -c_{20}^{1}$,
$ \frac{3+\sqrt{5}}{2}$,
$ \frac{3+\sqrt{5}}{2}$,
$ -\frac{1+\sqrt{5}}{2}c_{20}^{1}$,
$ -\frac{1+\sqrt{5}}{2}c_{20}^{1}$,
$  1+\sqrt{5} $;\ \ 
$ -\frac{1+\sqrt{5}}{2}c_{20}^{1}$,
$ \frac{1+\sqrt{5}}{2}c_{20}^{1}$,
$ -\frac{1+\sqrt{5}}{2}c_{20}^{1}$,
$ \frac{1+\sqrt{5}}{2}c_{20}^{1}$,
$ c_{20}^{1}$,
$ -c_{20}^{1}$,
$0$;\ \ 
$ -\frac{1+\sqrt{5}}{2}c_{20}^{1}$,
$ -\frac{1+\sqrt{5}}{2}c_{20}^{1}$,
$ \frac{1+\sqrt{5}}{2}c_{20}^{1}$,
$ -c_{20}^{1}$,
$ c_{20}^{1}$,
$0$;\ \ 
$ 1$,
$ 1$,
$ c_{20}^{1}$,
$ c_{20}^{1}$,
$  -1-\sqrt{5} $;\ \ 
$ 1$,
$ -c_{20}^{1}$,
$ -c_{20}^{1}$,
$  -1-\sqrt{5} $;\ \ 
$ \frac{1+\sqrt{5}}{2}c_{20}^{1}$,
$ -\frac{1+\sqrt{5}}{2}c_{20}^{1}$,
$0$;\ \ 
$ \frac{1+\sqrt{5}}{2}c_{20}^{1}$,
$0$;\ \ 
$  1+\sqrt{5} $)

  \vskip 2ex

\noindent67. $9_{\frac{12}{5},52.36}^{40,987}$ \irep{516}:\ \ 
$d_i$ = ($1.0$,
$1.0$,
$1.902$,
$1.902$,
$2.618$,
$2.618$,
$3.77$,
$3.77$,
$3.236$) 

\vskip 0.7ex
\hangindent=3em \hangafter=1
$D^2= 52.360 = 
 30+10\sqrt{5} $

\vskip 0.7ex
\hangindent=3em \hangafter=1
$T = ( 0,
0,
\frac{13}{40},
\frac{33}{40},
\frac{1}{5},
\frac{1}{5},
\frac{1}{8},
\frac{5}{8},
\frac{3}{5} )
$,

\vskip 0.7ex
\hangindent=3em \hangafter=1
$S$ = ($ 1$,
$ 1$,
$ c_{20}^{1}$,
$ c_{20}^{1}$,
$ \frac{3+\sqrt{5}}{2}$,
$ \frac{3+\sqrt{5}}{2}$,
$ \frac{1+\sqrt{5}}{2}c_{20}^{1}$,
$ \frac{1+\sqrt{5}}{2}c_{20}^{1}$,
$  1+\sqrt{5} $;\ \ 
$ 1$,
$ -c_{20}^{1}$,
$ -c_{20}^{1}$,
$ \frac{3+\sqrt{5}}{2}$,
$ \frac{3+\sqrt{5}}{2}$,
$ -\frac{1+\sqrt{5}}{2}c_{20}^{1}$,
$ -\frac{1+\sqrt{5}}{2}c_{20}^{1}$,
$  1+\sqrt{5} $;\ \ 
$ -\frac{1+\sqrt{5}}{2}c_{20}^{1}$,
$ \frac{1+\sqrt{5}}{2}c_{20}^{1}$,
$ -\frac{1+\sqrt{5}}{2}c_{20}^{1}$,
$ \frac{1+\sqrt{5}}{2}c_{20}^{1}$,
$ c_{20}^{1}$,
$ -c_{20}^{1}$,
$0$;\ \ 
$ -\frac{1+\sqrt{5}}{2}c_{20}^{1}$,
$ -\frac{1+\sqrt{5}}{2}c_{20}^{1}$,
$ \frac{1+\sqrt{5}}{2}c_{20}^{1}$,
$ -c_{20}^{1}$,
$ c_{20}^{1}$,
$0$;\ \ 
$ 1$,
$ 1$,
$ c_{20}^{1}$,
$ c_{20}^{1}$,
$  -1-\sqrt{5} $;\ \ 
$ 1$,
$ -c_{20}^{1}$,
$ -c_{20}^{1}$,
$  -1-\sqrt{5} $;\ \ 
$ \frac{1+\sqrt{5}}{2}c_{20}^{1}$,
$ -\frac{1+\sqrt{5}}{2}c_{20}^{1}$,
$0$;\ \ 
$ \frac{1+\sqrt{5}}{2}c_{20}^{1}$,
$0$;\ \ 
$  1+\sqrt{5} $)

  \vskip 2ex

\noindent68. $9_{\frac{28}{5},52.36}^{40,198}$ \irep{516}:\ \ 
$d_i$ = ($1.0$,
$1.0$,
$1.902$,
$1.902$,
$2.618$,
$2.618$,
$3.77$,
$3.77$,
$3.236$) 

\vskip 0.7ex
\hangindent=3em \hangafter=1
$D^2= 52.360 = 
 30+10\sqrt{5} $

\vskip 0.7ex
\hangindent=3em \hangafter=1
$T = ( 0,
0,
\frac{17}{40},
\frac{37}{40},
\frac{4}{5},
\frac{4}{5},
\frac{1}{8},
\frac{5}{8},
\frac{2}{5} )
$,

\vskip 0.7ex
\hangindent=3em \hangafter=1
$S$ = ($ 1$,
$ 1$,
$ c_{20}^{1}$,
$ c_{20}^{1}$,
$ \frac{3+\sqrt{5}}{2}$,
$ \frac{3+\sqrt{5}}{2}$,
$ \frac{1+\sqrt{5}}{2}c_{20}^{1}$,
$ \frac{1+\sqrt{5}}{2}c_{20}^{1}$,
$  1+\sqrt{5} $;\ \ 
$ 1$,
$ -c_{20}^{1}$,
$ -c_{20}^{1}$,
$ \frac{3+\sqrt{5}}{2}$,
$ \frac{3+\sqrt{5}}{2}$,
$ -\frac{1+\sqrt{5}}{2}c_{20}^{1}$,
$ -\frac{1+\sqrt{5}}{2}c_{20}^{1}$,
$  1+\sqrt{5} $;\ \ 
$ \frac{1+\sqrt{5}}{2}c_{20}^{1}$,
$ -\frac{1+\sqrt{5}}{2}c_{20}^{1}$,
$ -\frac{1+\sqrt{5}}{2}c_{20}^{1}$,
$ \frac{1+\sqrt{5}}{2}c_{20}^{1}$,
$ c_{20}^{1}$,
$ -c_{20}^{1}$,
$0$;\ \ 
$ \frac{1+\sqrt{5}}{2}c_{20}^{1}$,
$ -\frac{1+\sqrt{5}}{2}c_{20}^{1}$,
$ \frac{1+\sqrt{5}}{2}c_{20}^{1}$,
$ -c_{20}^{1}$,
$ c_{20}^{1}$,
$0$;\ \ 
$ 1$,
$ 1$,
$ c_{20}^{1}$,
$ c_{20}^{1}$,
$  -1-\sqrt{5} $;\ \ 
$ 1$,
$ -c_{20}^{1}$,
$ -c_{20}^{1}$,
$  -1-\sqrt{5} $;\ \ 
$ -\frac{1+\sqrt{5}}{2}c_{20}^{1}$,
$ \frac{1+\sqrt{5}}{2}c_{20}^{1}$,
$0$;\ \ 
$ -\frac{1+\sqrt{5}}{2}c_{20}^{1}$,
$0$;\ \ 
$  1+\sqrt{5} $)

  \vskip 2ex

\noindent69. $9_{\frac{40}{7},86.41}^{7,313}$ \irep{215}:\ \ 
$d_i$ = ($1.0$,
$1.801$,
$1.801$,
$2.246$,
$2.246$,
$3.246$,
$4.48$,
$4.48$,
$5.48$) 

\vskip 0.7ex
\hangindent=3em \hangafter=1
$D^2= 86.413 = 
 49+35c^{1}_{7}
+14c^{2}_{7}
$

\vskip 0.7ex
\hangindent=3em \hangafter=1
$T = ( 0,
\frac{1}{7},
\frac{1}{7},
\frac{5}{7},
\frac{5}{7},
\frac{2}{7},
\frac{6}{7},
\frac{6}{7},
\frac{3}{7} )
$,

\vskip 0.7ex
\hangindent=3em \hangafter=1
$S$ = ($ 1$,
$ -c_{7}^{3}$,
$ -c_{7}^{3}$,
$ \xi_{7}^{3}$,
$ \xi_{7}^{3}$,
$  2+c^{1}_{7}
$,
$ \xi_{14}^{5}$,
$ \xi_{14}^{5}$,
$  3+2c^{1}_{7}
+c^{2}_{7}
$;\ \ 
$ -\xi_{7}^{3}$,
$  2+c^{1}_{7}
$,
$ \xi_{14}^{5}$,
$ 1$,
$ -\xi_{14}^{5}$,
$  -3-2  c^{1}_{7}
-c^{2}_{7}
$,
$ -c_{7}^{3}$,
$ \xi_{7}^{3}$;\ \ 
$ -\xi_{7}^{3}$,
$ 1$,
$ \xi_{14}^{5}$,
$ -\xi_{14}^{5}$,
$ -c_{7}^{3}$,
$  -3-2  c^{1}_{7}
-c^{2}_{7}
$,
$ \xi_{7}^{3}$;\ \ 
$ c_{7}^{3}$,
$  3+2c^{1}_{7}
+c^{2}_{7}
$,
$ -c_{7}^{3}$,
$  -2-c^{1}_{7}
$,
$ \xi_{7}^{3}$,
$ -\xi_{14}^{5}$;\ \ 
$ c_{7}^{3}$,
$ -c_{7}^{3}$,
$ \xi_{7}^{3}$,
$  -2-c^{1}_{7}
$,
$ -\xi_{14}^{5}$;\ \ 
$  3+2c^{1}_{7}
+c^{2}_{7}
$,
$ -\xi_{7}^{3}$,
$ -\xi_{7}^{3}$,
$ 1$;\ \ 
$ \xi_{14}^{5}$,
$ 1$,
$ c_{7}^{3}$;\ \ 
$ \xi_{14}^{5}$,
$ c_{7}^{3}$;\ \ 
$  2+c^{1}_{7}
$)

Factors = $3_{\frac{48}{7},9.295}^{7,790}\boxtimes 3_{\frac{48}{7},9.295}^{7,790}$

  \vskip 2ex

\noindent70. $9_{0,86.41}^{7,161}$ \irep{213}:\ \ 
$d_i$ = ($1.0$,
$1.801$,
$1.801$,
$2.246$,
$2.246$,
$3.246$,
$4.48$,
$4.48$,
$5.48$) 

\vskip 0.7ex
\hangindent=3em \hangafter=1
$D^2= 86.413 = 
 49+35c^{1}_{7}
+14c^{2}_{7}
$

\vskip 0.7ex
\hangindent=3em \hangafter=1
$T = ( 0,
\frac{1}{7},
\frac{6}{7},
\frac{2}{7},
\frac{5}{7},
0,
\frac{3}{7},
\frac{4}{7},
0 )
$,

\vskip 0.7ex
\hangindent=3em \hangafter=1
$S$ = ($ 1$,
$ -c_{7}^{3}$,
$ -c_{7}^{3}$,
$ \xi_{7}^{3}$,
$ \xi_{7}^{3}$,
$  2+c^{1}_{7}
$,
$ \xi_{14}^{5}$,
$ \xi_{14}^{5}$,
$  3+2c^{1}_{7}
+c^{2}_{7}
$;\ \ 
$ -\xi_{7}^{3}$,
$  2+c^{1}_{7}
$,
$ \xi_{14}^{5}$,
$ 1$,
$ -\xi_{14}^{5}$,
$  -3-2  c^{1}_{7}
-c^{2}_{7}
$,
$ -c_{7}^{3}$,
$ \xi_{7}^{3}$;\ \ 
$ -\xi_{7}^{3}$,
$ 1$,
$ \xi_{14}^{5}$,
$ -\xi_{14}^{5}$,
$ -c_{7}^{3}$,
$  -3-2  c^{1}_{7}
-c^{2}_{7}
$,
$ \xi_{7}^{3}$;\ \ 
$ c_{7}^{3}$,
$  3+2c^{1}_{7}
+c^{2}_{7}
$,
$ -c_{7}^{3}$,
$  -2-c^{1}_{7}
$,
$ \xi_{7}^{3}$,
$ -\xi_{14}^{5}$;\ \ 
$ c_{7}^{3}$,
$ -c_{7}^{3}$,
$ \xi_{7}^{3}$,
$  -2-c^{1}_{7}
$,
$ -\xi_{14}^{5}$;\ \ 
$  3+2c^{1}_{7}
+c^{2}_{7}
$,
$ -\xi_{7}^{3}$,
$ -\xi_{7}^{3}$,
$ 1$;\ \ 
$ \xi_{14}^{5}$,
$ 1$,
$ c_{7}^{3}$;\ \ 
$ \xi_{14}^{5}$,
$ c_{7}^{3}$;\ \ 
$  2+c^{1}_{7}
$)

Factors = $3_{\frac{48}{7},9.295}^{7,790}\boxtimes 3_{\frac{8}{7},9.295}^{7,245}$

  \vskip 2ex

\noindent71. $9_{\frac{16}{7},86.41}^{7,112}$ \irep{215}:\ \ 
$d_i$ = ($1.0$,
$1.801$,
$1.801$,
$2.246$,
$2.246$,
$3.246$,
$4.48$,
$4.48$,
$5.48$) 

\vskip 0.7ex
\hangindent=3em \hangafter=1
$D^2= 86.413 = 
 49+35c^{1}_{7}
+14c^{2}_{7}
$

\vskip 0.7ex
\hangindent=3em \hangafter=1
$T = ( 0,
\frac{6}{7},
\frac{6}{7},
\frac{2}{7},
\frac{2}{7},
\frac{5}{7},
\frac{1}{7},
\frac{1}{7},
\frac{4}{7} )
$,

\vskip 0.7ex
\hangindent=3em \hangafter=1
$S$ = ($ 1$,
$ -c_{7}^{3}$,
$ -c_{7}^{3}$,
$ \xi_{7}^{3}$,
$ \xi_{7}^{3}$,
$  2+c^{1}_{7}
$,
$ \xi_{14}^{5}$,
$ \xi_{14}^{5}$,
$  3+2c^{1}_{7}
+c^{2}_{7}
$;\ \ 
$ -\xi_{7}^{3}$,
$  2+c^{1}_{7}
$,
$ \xi_{14}^{5}$,
$ 1$,
$ -\xi_{14}^{5}$,
$  -3-2  c^{1}_{7}
-c^{2}_{7}
$,
$ -c_{7}^{3}$,
$ \xi_{7}^{3}$;\ \ 
$ -\xi_{7}^{3}$,
$ 1$,
$ \xi_{14}^{5}$,
$ -\xi_{14}^{5}$,
$ -c_{7}^{3}$,
$  -3-2  c^{1}_{7}
-c^{2}_{7}
$,
$ \xi_{7}^{3}$;\ \ 
$ c_{7}^{3}$,
$  3+2c^{1}_{7}
+c^{2}_{7}
$,
$ -c_{7}^{3}$,
$  -2-c^{1}_{7}
$,
$ \xi_{7}^{3}$,
$ -\xi_{14}^{5}$;\ \ 
$ c_{7}^{3}$,
$ -c_{7}^{3}$,
$ \xi_{7}^{3}$,
$  -2-c^{1}_{7}
$,
$ -\xi_{14}^{5}$;\ \ 
$  3+2c^{1}_{7}
+c^{2}_{7}
$,
$ -\xi_{7}^{3}$,
$ -\xi_{7}^{3}$,
$ 1$;\ \ 
$ \xi_{14}^{5}$,
$ 1$,
$ c_{7}^{3}$;\ \ 
$ \xi_{14}^{5}$,
$ c_{7}^{3}$;\ \ 
$  2+c^{1}_{7}
$)

Factors = $3_{\frac{8}{7},9.295}^{7,245}\boxtimes 3_{\frac{8}{7},9.295}^{7,245}$

  \vskip 2ex

\noindent72. $9_{\frac{120}{19},175.3}^{19,574}$ \irep{447}:\ \ 
$d_i$ = ($1.0$,
$1.972$,
$2.891$,
$3.731$,
$4.469$,
$5.86$,
$5.563$,
$5.889$,
$6.54$) 

\vskip 0.7ex
\hangindent=3em \hangafter=1
$D^2= 175.332 = 
 \# $

\vskip 0.7ex
\hangindent=3em \hangafter=1
$T = ( 0,
\frac{4}{19},
\frac{17}{19},
\frac{1}{19},
\frac{13}{19},
\frac{15}{19},
\frac{7}{19},
\frac{8}{19},
\frac{18}{19} )
$,

\vskip 0.7ex
\hangindent=3em \hangafter=1
$S$ = ($ 1$,
$ -c_{19}^{9}$,
$ \xi_{19}^{3}$,
$ \xi_{19}^{15}$,
$ \xi_{19}^{5}$,
$ \xi_{19}^{13}$,
$ \xi_{19}^{7}$,
$ \xi_{19}^{11}$,
$ \xi_{19}^{9}$;\ \ 
$ -\xi_{19}^{15}$,
$ \xi_{19}^{13}$,
$ -\xi_{19}^{11}$,
$ \xi_{19}^{9}$,
$ -\xi_{19}^{7}$,
$ \xi_{19}^{5}$,
$ -\xi_{19}^{3}$,
$ 1$;\ \ 
$ \xi_{19}^{9}$,
$ \xi_{19}^{7}$,
$ \xi_{19}^{15}$,
$ 1$,
$ c_{19}^{9}$,
$ -\xi_{19}^{5}$,
$ -\xi_{19}^{11}$;\ \ 
$ -\xi_{19}^{3}$,
$ -1$,
$ \xi_{19}^{5}$,
$ -\xi_{19}^{9}$,
$ \xi_{19}^{13}$,
$ c_{19}^{9}$;\ \ 
$ -\xi_{19}^{13}$,
$ -\xi_{19}^{11}$,
$ -\xi_{19}^{3}$,
$ -c_{19}^{9}$,
$ \xi_{19}^{7}$;\ \ 
$ -c_{19}^{9}$,
$ \xi_{19}^{15}$,
$ -\xi_{19}^{9}$,
$ \xi_{19}^{3}$;\ \ 
$ \xi_{19}^{11}$,
$ 1$,
$ -\xi_{19}^{13}$;\ \ 
$ \xi_{19}^{7}$,
$ -\xi_{19}^{15}$;\ \ 
$ \xi_{19}^{5}$)

  \vskip 2ex

\noindent73. $9_{\frac{32}{19},175.3}^{19,327}$ \irep{447}:\ \ 
$d_i$ = ($1.0$,
$1.972$,
$2.891$,
$3.731$,
$4.469$,
$5.86$,
$5.563$,
$5.889$,
$6.54$) 

\vskip 0.7ex
\hangindent=3em \hangafter=1
$D^2= 175.332 = 
 \# $

\vskip 0.7ex
\hangindent=3em \hangafter=1
$T = ( 0,
\frac{15}{19},
\frac{2}{19},
\frac{18}{19},
\frac{6}{19},
\frac{4}{19},
\frac{12}{19},
\frac{11}{19},
\frac{1}{19} )
$,

\vskip 0.7ex
\hangindent=3em \hangafter=1
$S$ = ($ 1$,
$ -c_{19}^{9}$,
$ \xi_{19}^{3}$,
$ \xi_{19}^{15}$,
$ \xi_{19}^{5}$,
$ \xi_{19}^{13}$,
$ \xi_{19}^{7}$,
$ \xi_{19}^{11}$,
$ \xi_{19}^{9}$;\ \ 
$ -\xi_{19}^{15}$,
$ \xi_{19}^{13}$,
$ -\xi_{19}^{11}$,
$ \xi_{19}^{9}$,
$ -\xi_{19}^{7}$,
$ \xi_{19}^{5}$,
$ -\xi_{19}^{3}$,
$ 1$;\ \ 
$ \xi_{19}^{9}$,
$ \xi_{19}^{7}$,
$ \xi_{19}^{15}$,
$ 1$,
$ c_{19}^{9}$,
$ -\xi_{19}^{5}$,
$ -\xi_{19}^{11}$;\ \ 
$ -\xi_{19}^{3}$,
$ -1$,
$ \xi_{19}^{5}$,
$ -\xi_{19}^{9}$,
$ \xi_{19}^{13}$,
$ c_{19}^{9}$;\ \ 
$ -\xi_{19}^{13}$,
$ -\xi_{19}^{11}$,
$ -\xi_{19}^{3}$,
$ -c_{19}^{9}$,
$ \xi_{19}^{7}$;\ \ 
$ -c_{19}^{9}$,
$ \xi_{19}^{15}$,
$ -\xi_{19}^{9}$,
$ \xi_{19}^{3}$;\ \ 
$ \xi_{19}^{11}$,
$ 1$,
$ -\xi_{19}^{13}$;\ \ 
$ \xi_{19}^{7}$,
$ -\xi_{19}^{15}$;\ \ 
$ \xi_{19}^{5}$)

  \vskip 2ex

\noindent74. $9_{\frac{14}{5},343.2}^{15,715}$ \irep{400}:\ \ 
$d_i$ = ($1.0$,
$2.956$,
$4.783$,
$4.783$,
$4.783$,
$6.401$,
$7.739$,
$8.739$,
$9.357$) 

\vskip 0.7ex
\hangindent=3em \hangafter=1
$D^2= 343.211 = 
\frac{45+15\sqrt{5}}{2}\xi_{30}^{11}$

\vskip 0.7ex
\hangindent=3em \hangafter=1
$T = ( 0,
\frac{1}{15},
\frac{1}{5},
\frac{13}{15},
\frac{13}{15},
\frac{2}{5},
\frac{2}{3},
0,
\frac{2}{5} )
$,

\vskip 0.7ex
\hangindent=3em \hangafter=1
$S$ = ($ 1$,
$ \frac{1+\sqrt{5}}{2}c_{15}^{1}$,
$ \frac{3+\sqrt{5}}{2}c_{15}^{1}$,
$ \frac{3+\sqrt{5}}{2}c_{15}^{1}$,
$ \frac{3+\sqrt{5}}{2}c_{15}^{1}$,
$ \frac{3+\sqrt{5}}{2}\xi_{15}^{8,2}$,
$ (2+\sqrt{5} )c_{15}^{1}$,
$ \xi_{30}^{11}$,
$ \frac{3+\sqrt{5}}{2}\xi_{15}^{11}$;\ \ 
$ (2+\sqrt{5} )c_{15}^{1}$,
$ (3+\sqrt{5} )c_{15}^{1}$,
$ -\frac{3+\sqrt{5}}{2}c_{15}^{1}$,
$ -\frac{3+\sqrt{5}}{2}c_{15}^{1}$,
$ (2+\sqrt{5} )c_{15}^{1}$,
$ \frac{1+\sqrt{5}}{2}c_{15}^{1}$,
$ -\frac{1+\sqrt{5}}{2}c_{15}^{1}$,
$ -(2-\sqrt{5} )c_{15}^{1}$;\ \ 
$ \frac{3+\sqrt{5}}{2}c_{15}^{1}$,
$ \frac{3+\sqrt{5}}{2}c_{15}^{1}$,
$ \frac{3+\sqrt{5}}{2}c_{15}^{1}$,
$ -\frac{3+\sqrt{5}}{2}c_{15}^{1}$,
$ -(3-\sqrt{5} )c_{15}^{1}$,
$ -\frac{3+\sqrt{5}}{2}c_{15}^{1}$,
$ \frac{3+\sqrt{5}}{2}c_{15}^{1}$;\ \ 
$  -3-4\zeta^{1}_{5}
-\zeta^{1}_{3}
-\zeta^{2}_{5}
+\zeta^{8}_{15}
+3\zeta^{11}_{15}
+2\zeta^{14}_{15}
$,
$  3+5\zeta^{1}_{5}
+2\zeta^{1}_{3}
+3\zeta^{2}_{5}
+\zeta^{8}_{15}
+\zeta^{3}_{5}
-\zeta^{11}_{15}
-2\zeta^{14}_{15}
$,
$ -\frac{3+\sqrt{5}}{2}c_{15}^{1}$,
$ \frac{3+\sqrt{5}}{2}c_{15}^{1}$,
$ -\frac{3+\sqrt{5}}{2}c_{15}^{1}$,
$ \frac{3+\sqrt{5}}{2}c_{15}^{1}$;\ \ 
$  -3-4\zeta^{1}_{5}
-\zeta^{1}_{3}
-\zeta^{2}_{5}
+\zeta^{8}_{15}
+3\zeta^{11}_{15}
+2\zeta^{14}_{15}
$,
$ -\frac{3+\sqrt{5}}{2}c_{15}^{1}$,
$ \frac{3+\sqrt{5}}{2}c_{15}^{1}$,
$ -\frac{3+\sqrt{5}}{2}c_{15}^{1}$,
$ \frac{3+\sqrt{5}}{2}c_{15}^{1}$;\ \ 
$ -\xi_{30}^{11}$,
$ \frac{1+\sqrt{5}}{2}c_{15}^{1}$,
$ \frac{3+\sqrt{5}}{2}\xi_{15}^{11}$,
$ -1$;\ \ 
$ (2+\sqrt{5} )c_{15}^{1}$,
$ -(2-\sqrt{5} )c_{15}^{1}$,
$ -\frac{1+\sqrt{5}}{2}c_{15}^{1}$;\ \ 
$ 1$,
$ \frac{3+\sqrt{5}}{2}\xi_{15}^{8,2}$;\ \ 
$ -\xi_{30}^{11}$)

  \vskip 2ex

\noindent75. $9_{\frac{26}{5},343.2}^{15,296}$ \irep{400}:\ \ 
$d_i$ = ($1.0$,
$2.956$,
$4.783$,
$4.783$,
$4.783$,
$6.401$,
$7.739$,
$8.739$,
$9.357$) 

\vskip 0.7ex
\hangindent=3em \hangafter=1
$D^2= 343.211 = 
\frac{45+15\sqrt{5}}{2}\xi_{30}^{11}$

\vskip 0.7ex
\hangindent=3em \hangafter=1
$T = ( 0,
\frac{14}{15},
\frac{4}{5},
\frac{2}{15},
\frac{2}{15},
\frac{3}{5},
\frac{1}{3},
0,
\frac{3}{5} )
$,

\vskip 0.7ex
\hangindent=3em \hangafter=1
$S$ = ($ 1$,
$ \frac{1+\sqrt{5}}{2}c_{15}^{1}$,
$ \frac{3+\sqrt{5}}{2}c_{15}^{1}$,
$ \frac{3+\sqrt{5}}{2}c_{15}^{1}$,
$ \frac{3+\sqrt{5}}{2}c_{15}^{1}$,
$ \frac{3+\sqrt{5}}{2}\xi_{15}^{8,2}$,
$ (2+\sqrt{5} )c_{15}^{1}$,
$ \xi_{30}^{11}$,
$ \frac{3+\sqrt{5}}{2}\xi_{15}^{11}$;\ \ 
$ (2+\sqrt{5} )c_{15}^{1}$,
$ (3+\sqrt{5} )c_{15}^{1}$,
$ -\frac{3+\sqrt{5}}{2}c_{15}^{1}$,
$ -\frac{3+\sqrt{5}}{2}c_{15}^{1}$,
$ (2+\sqrt{5} )c_{15}^{1}$,
$ \frac{1+\sqrt{5}}{2}c_{15}^{1}$,
$ -\frac{1+\sqrt{5}}{2}c_{15}^{1}$,
$ -(2-\sqrt{5} )c_{15}^{1}$;\ \ 
$ \frac{3+\sqrt{5}}{2}c_{15}^{1}$,
$ \frac{3+\sqrt{5}}{2}c_{15}^{1}$,
$ \frac{3+\sqrt{5}}{2}c_{15}^{1}$,
$ -\frac{3+\sqrt{5}}{2}c_{15}^{1}$,
$ -(3-\sqrt{5} )c_{15}^{1}$,
$ -\frac{3+\sqrt{5}}{2}c_{15}^{1}$,
$ \frac{3+\sqrt{5}}{2}c_{15}^{1}$;\ \ 
$  3+5\zeta^{1}_{5}
+2\zeta^{1}_{3}
+3\zeta^{2}_{5}
+\zeta^{8}_{15}
+\zeta^{3}_{5}
-\zeta^{11}_{15}
-2\zeta^{14}_{15}
$,
$  -3-4\zeta^{1}_{5}
-\zeta^{1}_{3}
-\zeta^{2}_{5}
+\zeta^{8}_{15}
+3\zeta^{11}_{15}
+2\zeta^{14}_{15}
$,
$ -\frac{3+\sqrt{5}}{2}c_{15}^{1}$,
$ \frac{3+\sqrt{5}}{2}c_{15}^{1}$,
$ -\frac{3+\sqrt{5}}{2}c_{15}^{1}$,
$ \frac{3+\sqrt{5}}{2}c_{15}^{1}$;\ \ 
$  3+5\zeta^{1}_{5}
+2\zeta^{1}_{3}
+3\zeta^{2}_{5}
+\zeta^{8}_{15}
+\zeta^{3}_{5}
-\zeta^{11}_{15}
-2\zeta^{14}_{15}
$,
$ -\frac{3+\sqrt{5}}{2}c_{15}^{1}$,
$ \frac{3+\sqrt{5}}{2}c_{15}^{1}$,
$ -\frac{3+\sqrt{5}}{2}c_{15}^{1}$,
$ \frac{3+\sqrt{5}}{2}c_{15}^{1}$;\ \ 
$ -\xi_{30}^{11}$,
$ \frac{1+\sqrt{5}}{2}c_{15}^{1}$,
$ \frac{3+\sqrt{5}}{2}\xi_{15}^{11}$,
$ -1$;\ \ 
$ (2+\sqrt{5} )c_{15}^{1}$,
$ -(2-\sqrt{5} )c_{15}^{1}$,
$ -\frac{1+\sqrt{5}}{2}c_{15}^{1}$;\ \ 
$ 1$,
$ \frac{3+\sqrt{5}}{2}\xi_{15}^{8,2}$;\ \ 
$ -\xi_{30}^{11}$)

  \vskip 2ex

\noindent76. $9_{7,475.1}^{24,793}$ \irep{471}:\ \ 
$d_i$ = ($1.0$,
$4.449$,
$4.449$,
$5.449$,
$5.449$,
$8.898$,
$8.898$,
$9.898$,
$10.898$) 

\vskip 0.7ex
\hangindent=3em \hangafter=1
$D^2= 475.151 = 
 240+96\sqrt{6} $

\vskip 0.7ex
\hangindent=3em \hangafter=1
$T = ( 0,
\frac{1}{4},
\frac{1}{4},
\frac{1}{2},
\frac{1}{2},
\frac{1}{3},
\frac{7}{12},
0,
\frac{7}{8} )
$,

\vskip 0.7ex
\hangindent=3em \hangafter=1
$S$ = ($ 1$,
$  2+\sqrt{6} $,
$  2+\sqrt{6} $,
$  3+\sqrt{6} $,
$  3+\sqrt{6} $,
$  4+2\sqrt{6} $,
$  4+2\sqrt{6} $,
$  5+2\sqrt{6} $,
$  6+2\sqrt{6} $;\ \ 
$ 2\xi_{24}^{7}$,
$ -(4+2\sqrt{3} )\xi_{24}^{7}$,
$ -(2-\sqrt{6} )c_{12}^{1}$,
$ (2+\sqrt{6} )c_{12}^{1}$,
$  -4-2\sqrt{6} $,
$  4+2\sqrt{6} $,
$  -2-\sqrt{6} $,
$0$;\ \ 
$ 2\xi_{24}^{7}$,
$ (2+\sqrt{6} )c_{12}^{1}$,
$ -(2-\sqrt{6} )c_{12}^{1}$,
$  -4-2\sqrt{6} $,
$  4+2\sqrt{6} $,
$  -2-\sqrt{6} $,
$0$;\ \ 
$ -\sqrt{3}\xi_{24}^{7,5}$,
$ \sqrt{3}\xi_{24}^{11}$,
$0$,
$0$,
$  3+\sqrt{6} $,
$  -6-2\sqrt{6} $;\ \ 
$ -\sqrt{3}\xi_{24}^{7,5}$,
$0$,
$0$,
$  3+\sqrt{6} $,
$  -6-2\sqrt{6} $;\ \ 
$  4+2\sqrt{6} $,
$  4+2\sqrt{6} $,
$  -4-2\sqrt{6} $,
$0$;\ \ 
$  -4-2\sqrt{6} $,
$  -4-2\sqrt{6} $,
$0$;\ \ 
$ 1$,
$  6+2\sqrt{6} $;\ \ 
$0$)

  \vskip 2ex

\noindent77. $9_{3,475.1}^{24,139}$ \irep{471}:\ \ 
$d_i$ = ($1.0$,
$4.449$,
$4.449$,
$5.449$,
$5.449$,
$8.898$,
$8.898$,
$9.898$,
$10.898$) 

\vskip 0.7ex
\hangindent=3em \hangafter=1
$D^2= 475.151 = 
 240+96\sqrt{6} $

\vskip 0.7ex
\hangindent=3em \hangafter=1
$T = ( 0,
\frac{1}{4},
\frac{1}{4},
\frac{1}{2},
\frac{1}{2},
\frac{2}{3},
\frac{11}{12},
0,
\frac{3}{8} )
$,

\vskip 0.7ex
\hangindent=3em \hangafter=1
$S$ = ($ 1$,
$  2+\sqrt{6} $,
$  2+\sqrt{6} $,
$  3+\sqrt{6} $,
$  3+\sqrt{6} $,
$  4+2\sqrt{6} $,
$  4+2\sqrt{6} $,
$  5+2\sqrt{6} $,
$  6+2\sqrt{6} $;\ \ 
$ -(4+2\sqrt{3} )\xi_{24}^{7}$,
$ 2\xi_{24}^{7}$,
$ -(2-\sqrt{6} )c_{12}^{1}$,
$ (2+\sqrt{6} )c_{12}^{1}$,
$  -4-2\sqrt{6} $,
$  4+2\sqrt{6} $,
$  -2-\sqrt{6} $,
$0$;\ \ 
$ -(4+2\sqrt{3} )\xi_{24}^{7}$,
$ (2+\sqrt{6} )c_{12}^{1}$,
$ -(2-\sqrt{6} )c_{12}^{1}$,
$  -4-2\sqrt{6} $,
$  4+2\sqrt{6} $,
$  -2-\sqrt{6} $,
$0$;\ \ 
$ \sqrt{3}\xi_{24}^{11}$,
$ -\sqrt{3}\xi_{24}^{7,5}$,
$0$,
$0$,
$  3+\sqrt{6} $,
$  -6-2\sqrt{6} $;\ \ 
$ \sqrt{3}\xi_{24}^{11}$,
$0$,
$0$,
$  3+\sqrt{6} $,
$  -6-2\sqrt{6} $;\ \ 
$  4+2\sqrt{6} $,
$  4+2\sqrt{6} $,
$  -4-2\sqrt{6} $,
$0$;\ \ 
$  -4-2\sqrt{6} $,
$  -4-2\sqrt{6} $,
$0$;\ \ 
$ 1$,
$  6+2\sqrt{6} $;\ \ 
$0$)

  \vskip 2ex

\noindent78. $9_{5,475.1}^{24,181}$ \irep{471}:\ \ 
$d_i$ = ($1.0$,
$4.449$,
$4.449$,
$5.449$,
$5.449$,
$8.898$,
$8.898$,
$9.898$,
$10.898$) 

\vskip 0.7ex
\hangindent=3em \hangafter=1
$D^2= 475.151 = 
 240+96\sqrt{6} $

\vskip 0.7ex
\hangindent=3em \hangafter=1
$T = ( 0,
\frac{3}{4},
\frac{3}{4},
\frac{1}{2},
\frac{1}{2},
\frac{1}{3},
\frac{1}{12},
0,
\frac{5}{8} )
$,

\vskip 0.7ex
\hangindent=3em \hangafter=1
$S$ = ($ 1$,
$  2+\sqrt{6} $,
$  2+\sqrt{6} $,
$  3+\sqrt{6} $,
$  3+\sqrt{6} $,
$  4+2\sqrt{6} $,
$  4+2\sqrt{6} $,
$  5+2\sqrt{6} $,
$  6+2\sqrt{6} $;\ \ 
$ -(4+2\sqrt{3} )\xi_{24}^{7}$,
$ 2\xi_{24}^{7}$,
$ -(2-\sqrt{6} )c_{12}^{1}$,
$ (2+\sqrt{6} )c_{12}^{1}$,
$  -4-2\sqrt{6} $,
$  4+2\sqrt{6} $,
$  -2-\sqrt{6} $,
$0$;\ \ 
$ -(4+2\sqrt{3} )\xi_{24}^{7}$,
$ (2+\sqrt{6} )c_{12}^{1}$,
$ -(2-\sqrt{6} )c_{12}^{1}$,
$  -4-2\sqrt{6} $,
$  4+2\sqrt{6} $,
$  -2-\sqrt{6} $,
$0$;\ \ 
$ \sqrt{3}\xi_{24}^{11}$,
$ -\sqrt{3}\xi_{24}^{7,5}$,
$0$,
$0$,
$  3+\sqrt{6} $,
$  -6-2\sqrt{6} $;\ \ 
$ \sqrt{3}\xi_{24}^{11}$,
$0$,
$0$,
$  3+\sqrt{6} $,
$  -6-2\sqrt{6} $;\ \ 
$  4+2\sqrt{6} $,
$  4+2\sqrt{6} $,
$  -4-2\sqrt{6} $,
$0$;\ \ 
$  -4-2\sqrt{6} $,
$  -4-2\sqrt{6} $,
$0$;\ \ 
$ 1$,
$  6+2\sqrt{6} $;\ \ 
$0$)

  \vskip 2ex

\noindent79. $9_{1,475.1}^{24,144}$ \irep{471}:\ \ 
$d_i$ = ($1.0$,
$4.449$,
$4.449$,
$5.449$,
$5.449$,
$8.898$,
$8.898$,
$9.898$,
$10.898$) 

\vskip 0.7ex
\hangindent=3em \hangafter=1
$D^2= 475.151 = 
 240+96\sqrt{6} $

\vskip 0.7ex
\hangindent=3em \hangafter=1
$T = ( 0,
\frac{3}{4},
\frac{3}{4},
\frac{1}{2},
\frac{1}{2},
\frac{2}{3},
\frac{5}{12},
0,
\frac{1}{8} )
$,

\vskip 0.7ex
\hangindent=3em \hangafter=1
$S$ = ($ 1$,
$  2+\sqrt{6} $,
$  2+\sqrt{6} $,
$  3+\sqrt{6} $,
$  3+\sqrt{6} $,
$  4+2\sqrt{6} $,
$  4+2\sqrt{6} $,
$  5+2\sqrt{6} $,
$  6+2\sqrt{6} $;\ \ 
$ 2\xi_{24}^{7}$,
$ -(4+2\sqrt{3} )\xi_{24}^{7}$,
$ -(2-\sqrt{6} )c_{12}^{1}$,
$ (2+\sqrt{6} )c_{12}^{1}$,
$  -4-2\sqrt{6} $,
$  4+2\sqrt{6} $,
$  -2-\sqrt{6} $,
$0$;\ \ 
$ 2\xi_{24}^{7}$,
$ (2+\sqrt{6} )c_{12}^{1}$,
$ -(2-\sqrt{6} )c_{12}^{1}$,
$  -4-2\sqrt{6} $,
$  4+2\sqrt{6} $,
$  -2-\sqrt{6} $,
$0$;\ \ 
$ -\sqrt{3}\xi_{24}^{7,5}$,
$ \sqrt{3}\xi_{24}^{11}$,
$0$,
$0$,
$  3+\sqrt{6} $,
$  -6-2\sqrt{6} $;\ \ 
$ -\sqrt{3}\xi_{24}^{7,5}$,
$0$,
$0$,
$  3+\sqrt{6} $,
$  -6-2\sqrt{6} $;\ \ 
$  4+2\sqrt{6} $,
$  4+2\sqrt{6} $,
$  -4-2\sqrt{6} $,
$0$;\ \ 
$  -4-2\sqrt{6} $,
$  -4-2\sqrt{6} $,
$0$;\ \ 
$ 1$,
$  6+2\sqrt{6} $;\ \ 
$0$)

  \vskip 2ex

\noindent80. $9_{6,668.5}^{12,567}$ \irep{333}:\ \ 
$d_i$ = ($1.0$,
$6.464$,
$6.464$,
$6.464$,
$6.464$,
$6.464$,
$6.464$,
$13.928$,
$14.928$) 

\vskip 0.7ex
\hangindent=3em \hangafter=1
$D^2= 668.553 = 
 336+192\sqrt{3} $

\vskip 0.7ex
\hangindent=3em \hangafter=1
$T = ( 0,
\frac{1}{2},
\frac{1}{2},
\frac{1}{4},
\frac{1}{4},
\frac{1}{4},
\frac{1}{4},
0,
\frac{2}{3} )
$,

\vskip 0.7ex
\hangindent=3em \hangafter=1
$S$ = ($ 1$,
$  3+2\sqrt{3} $,
$  3+2\sqrt{3} $,
$  3+2\sqrt{3} $,
$  3+2\sqrt{3} $,
$  3+2\sqrt{3} $,
$  3+2\sqrt{3} $,
$  7+4\sqrt{3} $,
$  8+4\sqrt{3} $;\ \ 
$  1+8\zeta^{1}_{4}
+8\zeta^{1}_{3}
+4\zeta^{7}_{12}
$,
$  -7-4\zeta^{1}_{4}
-8\zeta^{1}_{3}
+4\zeta^{7}_{12}
$,
$  3+2\sqrt{3} $,
$  3+2\sqrt{3} $,
$  3+2\sqrt{3} $,
$  3+2\sqrt{3} $,
$  -3-2\sqrt{3} $,
$0$;\ \ 
$  1+8\zeta^{1}_{4}
+8\zeta^{1}_{3}
+4\zeta^{7}_{12}
$,
$  3+2\sqrt{3} $,
$  3+2\sqrt{3} $,
$  3+2\sqrt{3} $,
$  3+2\sqrt{3} $,
$  -3-2\sqrt{3} $,
$0$;\ \ 
$  9+6\sqrt{3} $,
$  -3-2\sqrt{3} $,
$  -3-2\sqrt{3} $,
$  -3-2\sqrt{3} $,
$  -3-2\sqrt{3} $,
$0$;\ \ 
$  9+6\sqrt{3} $,
$  -3-2\sqrt{3} $,
$  -3-2\sqrt{3} $,
$  -3-2\sqrt{3} $,
$0$;\ \ 
$  9+6\sqrt{3} $,
$  -3-2\sqrt{3} $,
$  -3-2\sqrt{3} $,
$0$;\ \ 
$  9+6\sqrt{3} $,
$  -3-2\sqrt{3} $,
$0$;\ \ 
$ 1$,
$  8+4\sqrt{3} $;\ \ 
$  -8-4\sqrt{3} $)

  \vskip 2ex

\noindent81. $9_{2,668.5}^{12,227}$ \irep{333}:\ \ 
$d_i$ = ($1.0$,
$6.464$,
$6.464$,
$6.464$,
$6.464$,
$6.464$,
$6.464$,
$13.928$,
$14.928$) 

\vskip 0.7ex
\hangindent=3em \hangafter=1
$D^2= 668.553 = 
 336+192\sqrt{3} $

\vskip 0.7ex
\hangindent=3em \hangafter=1
$T = ( 0,
\frac{1}{2},
\frac{1}{2},
\frac{3}{4},
\frac{3}{4},
\frac{3}{4},
\frac{3}{4},
0,
\frac{1}{3} )
$,

\vskip 0.7ex
\hangindent=3em \hangafter=1
$S$ = ($ 1$,
$  3+2\sqrt{3} $,
$  3+2\sqrt{3} $,
$  3+2\sqrt{3} $,
$  3+2\sqrt{3} $,
$  3+2\sqrt{3} $,
$  3+2\sqrt{3} $,
$  7+4\sqrt{3} $,
$  8+4\sqrt{3} $;\ \ 
$  -7-4\zeta^{1}_{4}
-8\zeta^{1}_{3}
+4\zeta^{7}_{12}
$,
$  1+8\zeta^{1}_{4}
+8\zeta^{1}_{3}
+4\zeta^{7}_{12}
$,
$  3+2\sqrt{3} $,
$  3+2\sqrt{3} $,
$  3+2\sqrt{3} $,
$  3+2\sqrt{3} $,
$  -3-2\sqrt{3} $,
$0$;\ \ 
$  -7-4\zeta^{1}_{4}
-8\zeta^{1}_{3}
+4\zeta^{7}_{12}
$,
$  3+2\sqrt{3} $,
$  3+2\sqrt{3} $,
$  3+2\sqrt{3} $,
$  3+2\sqrt{3} $,
$  -3-2\sqrt{3} $,
$0$;\ \ 
$  9+6\sqrt{3} $,
$  -3-2\sqrt{3} $,
$  -3-2\sqrt{3} $,
$  -3-2\sqrt{3} $,
$  -3-2\sqrt{3} $,
$0$;\ \ 
$  9+6\sqrt{3} $,
$  -3-2\sqrt{3} $,
$  -3-2\sqrt{3} $,
$  -3-2\sqrt{3} $,
$0$;\ \ 
$  9+6\sqrt{3} $,
$  -3-2\sqrt{3} $,
$  -3-2\sqrt{3} $,
$0$;\ \ 
$  9+6\sqrt{3} $,
$  -3-2\sqrt{3} $,
$0$;\ \ 
$ 1$,
$  8+4\sqrt{3} $;\ \ 
$  -8-4\sqrt{3} $)

  \vskip 2ex 

}

\subsection{Rank 10 }
\label{uni10}

{\small

\noindent1. $10_{1,10.}^{20,667}$ \irep{947}:\ \ 
$d_i$ = ($1.0$,
$1.0$,
$1.0$,
$1.0$,
$1.0$,
$1.0$,
$1.0$,
$1.0$,
$1.0$,
$1.0$) 

\vskip 0.7ex
\hangindent=3em \hangafter=1
$D^2= 10.0 = 
10$

\vskip 0.7ex
\hangindent=3em \hangafter=1
$T = ( 0,
\frac{1}{4},
\frac{1}{5},
\frac{1}{5},
\frac{4}{5},
\frac{4}{5},
\frac{1}{20},
\frac{1}{20},
\frac{9}{20},
\frac{9}{20} )
$,

\vskip 0.7ex
\hangindent=3em \hangafter=1
$S$ = ($ 1$,
$ 1$,
$ 1$,
$ 1$,
$ 1$,
$ 1$,
$ 1$,
$ 1$,
$ 1$,
$ 1$;\ \ 
$ -1$,
$ 1$,
$ 1$,
$ 1$,
$ 1$,
$ -1$,
$ -1$,
$ -1$,
$ -1$;\ \ 
$ -\zeta_{10}^{1}$,
$ \zeta_{5}^{2}$,
$ -\zeta_{10}^{3}$,
$ \zeta_{5}^{1}$,
$ -\zeta_{10}^{3}$,
$ \zeta_{5}^{1}$,
$ -\zeta_{10}^{1}$,
$ \zeta_{5}^{2}$;\ \ 
$ -\zeta_{10}^{1}$,
$ \zeta_{5}^{1}$,
$ -\zeta_{10}^{3}$,
$ \zeta_{5}^{1}$,
$ -\zeta_{10}^{3}$,
$ \zeta_{5}^{2}$,
$ -\zeta_{10}^{1}$;\ \ 
$ \zeta_{5}^{2}$,
$ -\zeta_{10}^{1}$,
$ \zeta_{5}^{2}$,
$ -\zeta_{10}^{1}$,
$ -\zeta_{10}^{3}$,
$ \zeta_{5}^{1}$;\ \ 
$ \zeta_{5}^{2}$,
$ -\zeta_{10}^{1}$,
$ \zeta_{5}^{2}$,
$ \zeta_{5}^{1}$,
$ -\zeta_{10}^{3}$;\ \ 
$ -\zeta_{5}^{2}$,
$ \zeta_{10}^{1}$,
$ \zeta_{10}^{3}$,
$ -\zeta_{5}^{1}$;\ \ 
$ -\zeta_{5}^{2}$,
$ -\zeta_{5}^{1}$,
$ \zeta_{10}^{3}$;\ \ 
$ \zeta_{10}^{1}$,
$ -\zeta_{5}^{2}$;\ \ 
$ \zeta_{10}^{1}$)

Factors = $2_{1,2.}^{4,437}\boxtimes 5_{0,5.}^{5,110}$

  \vskip 2ex

\noindent2. $10_{5,10.}^{20,892}$ \irep{947}:\ \ 
$d_i$ = ($1.0$,
$1.0$,
$1.0$,
$1.0$,
$1.0$,
$1.0$,
$1.0$,
$1.0$,
$1.0$,
$1.0$) 

\vskip 0.7ex
\hangindent=3em \hangafter=1
$D^2= 10.0 = 
10$

\vskip 0.7ex
\hangindent=3em \hangafter=1
$T = ( 0,
\frac{1}{4},
\frac{2}{5},
\frac{2}{5},
\frac{3}{5},
\frac{3}{5},
\frac{13}{20},
\frac{13}{20},
\frac{17}{20},
\frac{17}{20} )
$,

\vskip 0.7ex
\hangindent=3em \hangafter=1
$S$ = ($ 1$,
$ 1$,
$ 1$,
$ 1$,
$ 1$,
$ 1$,
$ 1$,
$ 1$,
$ 1$,
$ 1$;\ \ 
$ -1$,
$ 1$,
$ 1$,
$ 1$,
$ 1$,
$ -1$,
$ -1$,
$ -1$,
$ -1$;\ \ 
$ \zeta_{5}^{1}$,
$ -\zeta_{10}^{3}$,
$ -\zeta_{10}^{1}$,
$ \zeta_{5}^{2}$,
$ -\zeta_{10}^{3}$,
$ \zeta_{5}^{1}$,
$ \zeta_{5}^{2}$,
$ -\zeta_{10}^{1}$;\ \ 
$ \zeta_{5}^{1}$,
$ \zeta_{5}^{2}$,
$ -\zeta_{10}^{1}$,
$ \zeta_{5}^{1}$,
$ -\zeta_{10}^{3}$,
$ -\zeta_{10}^{1}$,
$ \zeta_{5}^{2}$;\ \ 
$ -\zeta_{10}^{3}$,
$ \zeta_{5}^{1}$,
$ \zeta_{5}^{2}$,
$ -\zeta_{10}^{1}$,
$ \zeta_{5}^{1}$,
$ -\zeta_{10}^{3}$;\ \ 
$ -\zeta_{10}^{3}$,
$ -\zeta_{10}^{1}$,
$ \zeta_{5}^{2}$,
$ -\zeta_{10}^{3}$,
$ \zeta_{5}^{1}$;\ \ 
$ -\zeta_{5}^{1}$,
$ \zeta_{10}^{3}$,
$ \zeta_{10}^{1}$,
$ -\zeta_{5}^{2}$;\ \ 
$ -\zeta_{5}^{1}$,
$ -\zeta_{5}^{2}$,
$ \zeta_{10}^{1}$;\ \ 
$ \zeta_{10}^{3}$,
$ -\zeta_{5}^{1}$;\ \ 
$ \zeta_{10}^{3}$)

Factors = $2_{1,2.}^{4,437}\boxtimes 5_{4,5.}^{5,210}$

  \vskip 2ex

\noindent3. $10_{7,10.}^{20,183}$ \irep{947}:\ \ 
$d_i$ = ($1.0$,
$1.0$,
$1.0$,
$1.0$,
$1.0$,
$1.0$,
$1.0$,
$1.0$,
$1.0$,
$1.0$) 

\vskip 0.7ex
\hangindent=3em \hangafter=1
$D^2= 10.0 = 
10$

\vskip 0.7ex
\hangindent=3em \hangafter=1
$T = ( 0,
\frac{3}{4},
\frac{1}{5},
\frac{1}{5},
\frac{4}{5},
\frac{4}{5},
\frac{11}{20},
\frac{11}{20},
\frac{19}{20},
\frac{19}{20} )
$,

\vskip 0.7ex
\hangindent=3em \hangafter=1
$S$ = ($ 1$,
$ 1$,
$ 1$,
$ 1$,
$ 1$,
$ 1$,
$ 1$,
$ 1$,
$ 1$,
$ 1$;\ \ 
$ -1$,
$ 1$,
$ 1$,
$ 1$,
$ 1$,
$ -1$,
$ -1$,
$ -1$,
$ -1$;\ \ 
$ -\zeta_{10}^{1}$,
$ \zeta_{5}^{2}$,
$ -\zeta_{10}^{3}$,
$ \zeta_{5}^{1}$,
$ -\zeta_{10}^{3}$,
$ \zeta_{5}^{1}$,
$ -\zeta_{10}^{1}$,
$ \zeta_{5}^{2}$;\ \ 
$ -\zeta_{10}^{1}$,
$ \zeta_{5}^{1}$,
$ -\zeta_{10}^{3}$,
$ \zeta_{5}^{1}$,
$ -\zeta_{10}^{3}$,
$ \zeta_{5}^{2}$,
$ -\zeta_{10}^{1}$;\ \ 
$ \zeta_{5}^{2}$,
$ -\zeta_{10}^{1}$,
$ \zeta_{5}^{2}$,
$ -\zeta_{10}^{1}$,
$ -\zeta_{10}^{3}$,
$ \zeta_{5}^{1}$;\ \ 
$ \zeta_{5}^{2}$,
$ -\zeta_{10}^{1}$,
$ \zeta_{5}^{2}$,
$ \zeta_{5}^{1}$,
$ -\zeta_{10}^{3}$;\ \ 
$ -\zeta_{5}^{2}$,
$ \zeta_{10}^{1}$,
$ \zeta_{10}^{3}$,
$ -\zeta_{5}^{1}$;\ \ 
$ -\zeta_{5}^{2}$,
$ -\zeta_{5}^{1}$,
$ \zeta_{10}^{3}$;\ \ 
$ \zeta_{10}^{1}$,
$ -\zeta_{5}^{2}$;\ \ 
$ \zeta_{10}^{1}$)

Factors = $2_{7,2.}^{4,625}\boxtimes 5_{0,5.}^{5,110}$

  \vskip 2ex

\noindent4. $10_{3,10.}^{20,607}$ \irep{947}:\ \ 
$d_i$ = ($1.0$,
$1.0$,
$1.0$,
$1.0$,
$1.0$,
$1.0$,
$1.0$,
$1.0$,
$1.0$,
$1.0$) 

\vskip 0.7ex
\hangindent=3em \hangafter=1
$D^2= 10.0 = 
10$

\vskip 0.7ex
\hangindent=3em \hangafter=1
$T = ( 0,
\frac{3}{4},
\frac{2}{5},
\frac{2}{5},
\frac{3}{5},
\frac{3}{5},
\frac{3}{20},
\frac{3}{20},
\frac{7}{20},
\frac{7}{20} )
$,

\vskip 0.7ex
\hangindent=3em \hangafter=1
$S$ = ($ 1$,
$ 1$,
$ 1$,
$ 1$,
$ 1$,
$ 1$,
$ 1$,
$ 1$,
$ 1$,
$ 1$;\ \ 
$ -1$,
$ 1$,
$ 1$,
$ 1$,
$ 1$,
$ -1$,
$ -1$,
$ -1$,
$ -1$;\ \ 
$ \zeta_{5}^{1}$,
$ -\zeta_{10}^{3}$,
$ -\zeta_{10}^{1}$,
$ \zeta_{5}^{2}$,
$ -\zeta_{10}^{3}$,
$ \zeta_{5}^{1}$,
$ \zeta_{5}^{2}$,
$ -\zeta_{10}^{1}$;\ \ 
$ \zeta_{5}^{1}$,
$ \zeta_{5}^{2}$,
$ -\zeta_{10}^{1}$,
$ \zeta_{5}^{1}$,
$ -\zeta_{10}^{3}$,
$ -\zeta_{10}^{1}$,
$ \zeta_{5}^{2}$;\ \ 
$ -\zeta_{10}^{3}$,
$ \zeta_{5}^{1}$,
$ \zeta_{5}^{2}$,
$ -\zeta_{10}^{1}$,
$ \zeta_{5}^{1}$,
$ -\zeta_{10}^{3}$;\ \ 
$ -\zeta_{10}^{3}$,
$ -\zeta_{10}^{1}$,
$ \zeta_{5}^{2}$,
$ -\zeta_{10}^{3}$,
$ \zeta_{5}^{1}$;\ \ 
$ -\zeta_{5}^{1}$,
$ \zeta_{10}^{3}$,
$ \zeta_{10}^{1}$,
$ -\zeta_{5}^{2}$;\ \ 
$ -\zeta_{5}^{1}$,
$ -\zeta_{5}^{2}$,
$ \zeta_{10}^{1}$;\ \ 
$ \zeta_{10}^{3}$,
$ -\zeta_{5}^{1}$;\ \ 
$ \zeta_{10}^{3}$)

Factors = $2_{7,2.}^{4,625}\boxtimes 5_{4,5.}^{5,210}$

  \vskip 2ex

\noindent5. $10_{\frac{14}{5},18.09}^{5,359}$ \irep{125}:\ \ 
$d_i$ = ($1.0$,
$1.0$,
$1.0$,
$1.0$,
$1.0$,
$1.618$,
$1.618$,
$1.618$,
$1.618$,
$1.618$) 

\vskip 0.7ex
\hangindent=3em \hangafter=1
$D^2= 18.90 = 
\frac{25+5\sqrt{5}}{2}$

\vskip 0.7ex
\hangindent=3em \hangafter=1
$T = ( 0,
\frac{1}{5},
\frac{1}{5},
\frac{4}{5},
\frac{4}{5},
\frac{1}{5},
\frac{1}{5},
\frac{2}{5},
\frac{3}{5},
\frac{3}{5} )
$,

\vskip 0.7ex
\hangindent=3em \hangafter=1
$S$ = ($ 1$,
$ 1$,
$ 1$,
$ 1$,
$ 1$,
$ \frac{1+\sqrt{5}}{2}$,
$ \frac{1+\sqrt{5}}{2}$,
$ \frac{1+\sqrt{5}}{2}$,
$ \frac{1+\sqrt{5}}{2}$,
$ \frac{1+\sqrt{5}}{2}$;\ \ 
$ -\zeta_{10}^{1}$,
$ \zeta_{5}^{2}$,
$ -\zeta_{10}^{3}$,
$ \zeta_{5}^{1}$,
$ -\frac{1+\sqrt{5}}{2}\zeta_{10}^{3}$,
$ \frac{1+\sqrt{5}}{2}\zeta_{5}^{1}$,
$ \frac{1+\sqrt{5}}{2}$,
$ -\frac{1+\sqrt{5}}{2}\zeta_{10}^{1}$,
$ \frac{1+\sqrt{5}}{2}\zeta_{5}^{2}$;\ \ 
$ -\zeta_{10}^{1}$,
$ \zeta_{5}^{1}$,
$ -\zeta_{10}^{3}$,
$ \frac{1+\sqrt{5}}{2}\zeta_{5}^{1}$,
$ -\frac{1+\sqrt{5}}{2}\zeta_{10}^{3}$,
$ \frac{1+\sqrt{5}}{2}$,
$ \frac{1+\sqrt{5}}{2}\zeta_{5}^{2}$,
$ -\frac{1+\sqrt{5}}{2}\zeta_{10}^{1}$;\ \ 
$ \zeta_{5}^{2}$,
$ -\zeta_{10}^{1}$,
$ \frac{1+\sqrt{5}}{2}\zeta_{5}^{2}$,
$ -\frac{1+\sqrt{5}}{2}\zeta_{10}^{1}$,
$ \frac{1+\sqrt{5}}{2}$,
$ -\frac{1+\sqrt{5}}{2}\zeta_{10}^{3}$,
$ \frac{1+\sqrt{5}}{2}\zeta_{5}^{1}$;\ \ 
$ \zeta_{5}^{2}$,
$ -\frac{1+\sqrt{5}}{2}\zeta_{10}^{1}$,
$ \frac{1+\sqrt{5}}{2}\zeta_{5}^{2}$,
$ \frac{1+\sqrt{5}}{2}$,
$ \frac{1+\sqrt{5}}{2}\zeta_{5}^{1}$,
$ -\frac{1+\sqrt{5}}{2}\zeta_{10}^{3}$;\ \ 
$ -\zeta_{5}^{2}$,
$ \zeta_{10}^{1}$,
$ -1$,
$ \zeta_{10}^{3}$,
$ -\zeta_{5}^{1}$;\ \ 
$ -\zeta_{5}^{2}$,
$ -1$,
$ -\zeta_{5}^{1}$,
$ \zeta_{10}^{3}$;\ \ 
$ -1$,
$ -1$,
$ -1$;\ \ 
$ \zeta_{10}^{1}$,
$ -\zeta_{5}^{2}$;\ \ 
$ \zeta_{10}^{1}$)

Factors = $2_{\frac{14}{5},3.618}^{5,395}\boxtimes 5_{0,5.}^{5,110}$

  \vskip 2ex

\noindent6. $10_{\frac{26}{5},18.09}^{5,125}$ \irep{125}:\ \ 
$d_i$ = ($1.0$,
$1.0$,
$1.0$,
$1.0$,
$1.0$,
$1.618$,
$1.618$,
$1.618$,
$1.618$,
$1.618$) 

\vskip 0.7ex
\hangindent=3em \hangafter=1
$D^2= 18.90 = 
\frac{25+5\sqrt{5}}{2}$

\vskip 0.7ex
\hangindent=3em \hangafter=1
$T = ( 0,
\frac{1}{5},
\frac{1}{5},
\frac{4}{5},
\frac{4}{5},
\frac{2}{5},
\frac{2}{5},
\frac{3}{5},
\frac{4}{5},
\frac{4}{5} )
$,

\vskip 0.7ex
\hangindent=3em \hangafter=1
$S$ = ($ 1$,
$ 1$,
$ 1$,
$ 1$,
$ 1$,
$ \frac{1+\sqrt{5}}{2}$,
$ \frac{1+\sqrt{5}}{2}$,
$ \frac{1+\sqrt{5}}{2}$,
$ \frac{1+\sqrt{5}}{2}$,
$ \frac{1+\sqrt{5}}{2}$;\ \ 
$ -\zeta_{10}^{1}$,
$ \zeta_{5}^{2}$,
$ -\zeta_{10}^{3}$,
$ \zeta_{5}^{1}$,
$ -\frac{1+\sqrt{5}}{2}\zeta_{10}^{3}$,
$ \frac{1+\sqrt{5}}{2}\zeta_{5}^{1}$,
$ \frac{1+\sqrt{5}}{2}$,
$ -\frac{1+\sqrt{5}}{2}\zeta_{10}^{1}$,
$ \frac{1+\sqrt{5}}{2}\zeta_{5}^{2}$;\ \ 
$ -\zeta_{10}^{1}$,
$ \zeta_{5}^{1}$,
$ -\zeta_{10}^{3}$,
$ \frac{1+\sqrt{5}}{2}\zeta_{5}^{1}$,
$ -\frac{1+\sqrt{5}}{2}\zeta_{10}^{3}$,
$ \frac{1+\sqrt{5}}{2}$,
$ \frac{1+\sqrt{5}}{2}\zeta_{5}^{2}$,
$ -\frac{1+\sqrt{5}}{2}\zeta_{10}^{1}$;\ \ 
$ \zeta_{5}^{2}$,
$ -\zeta_{10}^{1}$,
$ \frac{1+\sqrt{5}}{2}\zeta_{5}^{2}$,
$ -\frac{1+\sqrt{5}}{2}\zeta_{10}^{1}$,
$ \frac{1+\sqrt{5}}{2}$,
$ -\frac{1+\sqrt{5}}{2}\zeta_{10}^{3}$,
$ \frac{1+\sqrt{5}}{2}\zeta_{5}^{1}$;\ \ 
$ \zeta_{5}^{2}$,
$ -\frac{1+\sqrt{5}}{2}\zeta_{10}^{1}$,
$ \frac{1+\sqrt{5}}{2}\zeta_{5}^{2}$,
$ \frac{1+\sqrt{5}}{2}$,
$ \frac{1+\sqrt{5}}{2}\zeta_{5}^{1}$,
$ -\frac{1+\sqrt{5}}{2}\zeta_{10}^{3}$;\ \ 
$ -\zeta_{5}^{2}$,
$ \zeta_{10}^{1}$,
$ -1$,
$ \zeta_{10}^{3}$,
$ -\zeta_{5}^{1}$;\ \ 
$ -\zeta_{5}^{2}$,
$ -1$,
$ -\zeta_{5}^{1}$,
$ \zeta_{10}^{3}$;\ \ 
$ -1$,
$ -1$,
$ -1$;\ \ 
$ \zeta_{10}^{1}$,
$ -\zeta_{5}^{2}$;\ \ 
$ \zeta_{10}^{1}$)

Factors = $2_{\frac{26}{5},3.618}^{5,720}\boxtimes 5_{0,5.}^{5,110}$

  \vskip 2ex

\noindent7. $10_{\frac{6}{5},18.09}^{5,152}$ \irep{145}:\ \ 
$d_i$ = ($1.0$,
$1.0$,
$1.0$,
$1.0$,
$1.0$,
$1.618$,
$1.618$,
$1.618$,
$1.618$,
$1.618$) 

\vskip 0.7ex
\hangindent=3em \hangafter=1
$D^2= 18.90 = 
\frac{25+5\sqrt{5}}{2}$

\vskip 0.7ex
\hangindent=3em \hangafter=1
$T = ( 0,
\frac{2}{5},
\frac{2}{5},
\frac{3}{5},
\frac{3}{5},
0,
0,
\frac{1}{5},
\frac{1}{5},
\frac{3}{5} )
$,

\vskip 0.7ex
\hangindent=3em \hangafter=1
$S$ = ($ 1$,
$ 1$,
$ 1$,
$ 1$,
$ 1$,
$ \frac{1+\sqrt{5}}{2}$,
$ \frac{1+\sqrt{5}}{2}$,
$ \frac{1+\sqrt{5}}{2}$,
$ \frac{1+\sqrt{5}}{2}$,
$ \frac{1+\sqrt{5}}{2}$;\ \ 
$ \zeta_{5}^{1}$,
$ -\zeta_{10}^{3}$,
$ -\zeta_{10}^{1}$,
$ \zeta_{5}^{2}$,
$ -\frac{1+\sqrt{5}}{2}\zeta_{10}^{3}$,
$ \frac{1+\sqrt{5}}{2}\zeta_{5}^{1}$,
$ \frac{1+\sqrt{5}}{2}\zeta_{5}^{2}$,
$ -\frac{1+\sqrt{5}}{2}\zeta_{10}^{1}$,
$ \frac{1+\sqrt{5}}{2}$;\ \ 
$ \zeta_{5}^{1}$,
$ \zeta_{5}^{2}$,
$ -\zeta_{10}^{1}$,
$ \frac{1+\sqrt{5}}{2}\zeta_{5}^{1}$,
$ -\frac{1+\sqrt{5}}{2}\zeta_{10}^{3}$,
$ -\frac{1+\sqrt{5}}{2}\zeta_{10}^{1}$,
$ \frac{1+\sqrt{5}}{2}\zeta_{5}^{2}$,
$ \frac{1+\sqrt{5}}{2}$;\ \ 
$ -\zeta_{10}^{3}$,
$ \zeta_{5}^{1}$,
$ \frac{1+\sqrt{5}}{2}\zeta_{5}^{2}$,
$ -\frac{1+\sqrt{5}}{2}\zeta_{10}^{1}$,
$ \frac{1+\sqrt{5}}{2}\zeta_{5}^{1}$,
$ -\frac{1+\sqrt{5}}{2}\zeta_{10}^{3}$,
$ \frac{1+\sqrt{5}}{2}$;\ \ 
$ -\zeta_{10}^{3}$,
$ -\frac{1+\sqrt{5}}{2}\zeta_{10}^{1}$,
$ \frac{1+\sqrt{5}}{2}\zeta_{5}^{2}$,
$ -\frac{1+\sqrt{5}}{2}\zeta_{10}^{3}$,
$ \frac{1+\sqrt{5}}{2}\zeta_{5}^{1}$,
$ \frac{1+\sqrt{5}}{2}$;\ \ 
$ -\zeta_{5}^{1}$,
$ \zeta_{10}^{3}$,
$ \zeta_{10}^{1}$,
$ -\zeta_{5}^{2}$,
$ -1$;\ \ 
$ -\zeta_{5}^{1}$,
$ -\zeta_{5}^{2}$,
$ \zeta_{10}^{1}$,
$ -1$;\ \ 
$ \zeta_{10}^{3}$,
$ -\zeta_{5}^{1}$,
$ -1$;\ \ 
$ \zeta_{10}^{3}$,
$ -1$;\ \ 
$ -1$)

Factors = $2_{\frac{26}{5},3.618}^{5,720}\boxtimes 5_{4,5.}^{5,210}$

  \vskip 2ex

\noindent8. $10_{\frac{34}{5},18.09}^{5,974}$ \irep{145}:\ \ 
$d_i$ = ($1.0$,
$1.0$,
$1.0$,
$1.0$,
$1.0$,
$1.618$,
$1.618$,
$1.618$,
$1.618$,
$1.618$) 

\vskip 0.7ex
\hangindent=3em \hangafter=1
$D^2= 18.90 = 
\frac{25+5\sqrt{5}}{2}$

\vskip 0.7ex
\hangindent=3em \hangafter=1
$T = ( 0,
\frac{2}{5},
\frac{2}{5},
\frac{3}{5},
\frac{3}{5},
0,
0,
\frac{2}{5},
\frac{4}{5},
\frac{4}{5} )
$,

\vskip 0.7ex
\hangindent=3em \hangafter=1
$S$ = ($ 1$,
$ 1$,
$ 1$,
$ 1$,
$ 1$,
$ \frac{1+\sqrt{5}}{2}$,
$ \frac{1+\sqrt{5}}{2}$,
$ \frac{1+\sqrt{5}}{2}$,
$ \frac{1+\sqrt{5}}{2}$,
$ \frac{1+\sqrt{5}}{2}$;\ \ 
$ \zeta_{5}^{1}$,
$ -\zeta_{10}^{3}$,
$ -\zeta_{10}^{1}$,
$ \zeta_{5}^{2}$,
$ -\frac{1+\sqrt{5}}{2}\zeta_{10}^{1}$,
$ \frac{1+\sqrt{5}}{2}\zeta_{5}^{2}$,
$ \frac{1+\sqrt{5}}{2}$,
$ \frac{1+\sqrt{5}}{2}\zeta_{5}^{1}$,
$ -\frac{1+\sqrt{5}}{2}\zeta_{10}^{3}$;\ \ 
$ \zeta_{5}^{1}$,
$ \zeta_{5}^{2}$,
$ -\zeta_{10}^{1}$,
$ \frac{1+\sqrt{5}}{2}\zeta_{5}^{2}$,
$ -\frac{1+\sqrt{5}}{2}\zeta_{10}^{1}$,
$ \frac{1+\sqrt{5}}{2}$,
$ -\frac{1+\sqrt{5}}{2}\zeta_{10}^{3}$,
$ \frac{1+\sqrt{5}}{2}\zeta_{5}^{1}$;\ \ 
$ -\zeta_{10}^{3}$,
$ \zeta_{5}^{1}$,
$ -\frac{1+\sqrt{5}}{2}\zeta_{10}^{3}$,
$ \frac{1+\sqrt{5}}{2}\zeta_{5}^{1}$,
$ \frac{1+\sqrt{5}}{2}$,
$ -\frac{1+\sqrt{5}}{2}\zeta_{10}^{1}$,
$ \frac{1+\sqrt{5}}{2}\zeta_{5}^{2}$;\ \ 
$ -\zeta_{10}^{3}$,
$ \frac{1+\sqrt{5}}{2}\zeta_{5}^{1}$,
$ -\frac{1+\sqrt{5}}{2}\zeta_{10}^{3}$,
$ \frac{1+\sqrt{5}}{2}$,
$ \frac{1+\sqrt{5}}{2}\zeta_{5}^{2}$,
$ -\frac{1+\sqrt{5}}{2}\zeta_{10}^{1}$;\ \ 
$ \zeta_{10}^{3}$,
$ -\zeta_{5}^{1}$,
$ -1$,
$ \zeta_{10}^{1}$,
$ -\zeta_{5}^{2}$;\ \ 
$ \zeta_{10}^{3}$,
$ -1$,
$ -\zeta_{5}^{2}$,
$ \zeta_{10}^{1}$;\ \ 
$ -1$,
$ -1$,
$ -1$;\ \ 
$ -\zeta_{5}^{1}$,
$ \zeta_{10}^{3}$;\ \ 
$ -\zeta_{5}^{1}$)

Factors = $2_{\frac{14}{5},3.618}^{5,395}\boxtimes 5_{4,5.}^{5,210}$

  \vskip 2ex

\noindent9. $10_{3,24.}^{24,430}$ \irep{972}:\ \ 
$d_i$ = ($1.0$,
$1.0$,
$1.0$,
$1.0$,
$1.732$,
$1.732$,
$1.732$,
$1.732$,
$2.0$,
$2.0$) 

\vskip 0.7ex
\hangindent=3em \hangafter=1
$D^2= 24.0 = 
24$

\vskip 0.7ex
\hangindent=3em \hangafter=1
$T = ( 0,
0,
\frac{1}{4},
\frac{1}{4},
\frac{1}{8},
\frac{3}{8},
\frac{5}{8},
\frac{7}{8},
\frac{1}{3},
\frac{7}{12} )
$,

\vskip 0.7ex
\hangindent=3em \hangafter=1
$S$ = ($ 1$,
$ 1$,
$ 1$,
$ 1$,
$ \sqrt{3}$,
$ \sqrt{3}$,
$ \sqrt{3}$,
$ \sqrt{3}$,
$ 2$,
$ 2$;\ \ 
$ 1$,
$ 1$,
$ 1$,
$ -\sqrt{3}$,
$ -\sqrt{3}$,
$ -\sqrt{3}$,
$ -\sqrt{3}$,
$ 2$,
$ 2$;\ \ 
$ -1$,
$ -1$,
$ \sqrt{3}$,
$ -\sqrt{3}$,
$ \sqrt{3}$,
$ -\sqrt{3}$,
$ 2$,
$ -2$;\ \ 
$ -1$,
$ -\sqrt{3}$,
$ \sqrt{3}$,
$ -\sqrt{3}$,
$ \sqrt{3}$,
$ 2$,
$ -2$;\ \ 
$ \sqrt{3}$,
$ \sqrt{3}$,
$ -\sqrt{3}$,
$ -\sqrt{3}$,
$0$,
$0$;\ \ 
$ -\sqrt{3}$,
$ -\sqrt{3}$,
$ \sqrt{3}$,
$0$,
$0$;\ \ 
$ \sqrt{3}$,
$ \sqrt{3}$,
$0$,
$0$;\ \ 
$ -\sqrt{3}$,
$0$,
$0$;\ \ 
$ -2$,
$ -2$;\ \ 
$ 2$)

Factors = $2_{1,2.}^{4,437}\boxtimes 5_{2,12.}^{24,940}$

  \vskip 2ex

\noindent10. $10_{7,24.}^{24,123}$ \irep{972}:\ \ 
$d_i$ = ($1.0$,
$1.0$,
$1.0$,
$1.0$,
$1.732$,
$1.732$,
$1.732$,
$1.732$,
$2.0$,
$2.0$) 

\vskip 0.7ex
\hangindent=3em \hangafter=1
$D^2= 24.0 = 
24$

\vskip 0.7ex
\hangindent=3em \hangafter=1
$T = ( 0,
0,
\frac{1}{4},
\frac{1}{4},
\frac{1}{8},
\frac{3}{8},
\frac{5}{8},
\frac{7}{8},
\frac{2}{3},
\frac{11}{12} )
$,

\vskip 0.7ex
\hangindent=3em \hangafter=1
$S$ = ($ 1$,
$ 1$,
$ 1$,
$ 1$,
$ \sqrt{3}$,
$ \sqrt{3}$,
$ \sqrt{3}$,
$ \sqrt{3}$,
$ 2$,
$ 2$;\ \ 
$ 1$,
$ 1$,
$ 1$,
$ -\sqrt{3}$,
$ -\sqrt{3}$,
$ -\sqrt{3}$,
$ -\sqrt{3}$,
$ 2$,
$ 2$;\ \ 
$ -1$,
$ -1$,
$ \sqrt{3}$,
$ -\sqrt{3}$,
$ \sqrt{3}$,
$ -\sqrt{3}$,
$ 2$,
$ -2$;\ \ 
$ -1$,
$ -\sqrt{3}$,
$ \sqrt{3}$,
$ -\sqrt{3}$,
$ \sqrt{3}$,
$ 2$,
$ -2$;\ \ 
$ -\sqrt{3}$,
$ -\sqrt{3}$,
$ \sqrt{3}$,
$ \sqrt{3}$,
$0$,
$0$;\ \ 
$ \sqrt{3}$,
$ \sqrt{3}$,
$ -\sqrt{3}$,
$0$,
$0$;\ \ 
$ -\sqrt{3}$,
$ -\sqrt{3}$,
$0$,
$0$;\ \ 
$ \sqrt{3}$,
$0$,
$0$;\ \ 
$ -2$,
$ -2$;\ \ 
$ 2$)

Factors = $2_{1,2.}^{4,437}\boxtimes 5_{6,12.}^{24,592}$

  \vskip 2ex

\noindent11. $10_{1,24.}^{24,976}$ \irep{972}:\ \ 
$d_i$ = ($1.0$,
$1.0$,
$1.0$,
$1.0$,
$1.732$,
$1.732$,
$1.732$,
$1.732$,
$2.0$,
$2.0$) 

\vskip 0.7ex
\hangindent=3em \hangafter=1
$D^2= 24.0 = 
24$

\vskip 0.7ex
\hangindent=3em \hangafter=1
$T = ( 0,
0,
\frac{3}{4},
\frac{3}{4},
\frac{1}{8},
\frac{3}{8},
\frac{5}{8},
\frac{7}{8},
\frac{1}{3},
\frac{1}{12} )
$,

\vskip 0.7ex
\hangindent=3em \hangafter=1
$S$ = ($ 1$,
$ 1$,
$ 1$,
$ 1$,
$ \sqrt{3}$,
$ \sqrt{3}$,
$ \sqrt{3}$,
$ \sqrt{3}$,
$ 2$,
$ 2$;\ \ 
$ 1$,
$ 1$,
$ 1$,
$ -\sqrt{3}$,
$ -\sqrt{3}$,
$ -\sqrt{3}$,
$ -\sqrt{3}$,
$ 2$,
$ 2$;\ \ 
$ -1$,
$ -1$,
$ \sqrt{3}$,
$ -\sqrt{3}$,
$ \sqrt{3}$,
$ -\sqrt{3}$,
$ 2$,
$ -2$;\ \ 
$ -1$,
$ -\sqrt{3}$,
$ \sqrt{3}$,
$ -\sqrt{3}$,
$ \sqrt{3}$,
$ 2$,
$ -2$;\ \ 
$ \sqrt{3}$,
$ -\sqrt{3}$,
$ -\sqrt{3}$,
$ \sqrt{3}$,
$0$,
$0$;\ \ 
$ -\sqrt{3}$,
$ \sqrt{3}$,
$ \sqrt{3}$,
$0$,
$0$;\ \ 
$ \sqrt{3}$,
$ -\sqrt{3}$,
$0$,
$0$;\ \ 
$ -\sqrt{3}$,
$0$,
$0$;\ \ 
$ -2$,
$ -2$;\ \ 
$ 2$)

Factors = $2_{7,2.}^{4,625}\boxtimes 5_{2,12.}^{24,940}$

  \vskip 2ex

\noindent12. $10_{5,24.}^{24,902}$ \irep{972}:\ \ 
$d_i$ = ($1.0$,
$1.0$,
$1.0$,
$1.0$,
$1.732$,
$1.732$,
$1.732$,
$1.732$,
$2.0$,
$2.0$) 

\vskip 0.7ex
\hangindent=3em \hangafter=1
$D^2= 24.0 = 
24$

\vskip 0.7ex
\hangindent=3em \hangafter=1
$T = ( 0,
0,
\frac{3}{4},
\frac{3}{4},
\frac{1}{8},
\frac{3}{8},
\frac{5}{8},
\frac{7}{8},
\frac{2}{3},
\frac{5}{12} )
$,

\vskip 0.7ex
\hangindent=3em \hangafter=1
$S$ = ($ 1$,
$ 1$,
$ 1$,
$ 1$,
$ \sqrt{3}$,
$ \sqrt{3}$,
$ \sqrt{3}$,
$ \sqrt{3}$,
$ 2$,
$ 2$;\ \ 
$ 1$,
$ 1$,
$ 1$,
$ -\sqrt{3}$,
$ -\sqrt{3}$,
$ -\sqrt{3}$,
$ -\sqrt{3}$,
$ 2$,
$ 2$;\ \ 
$ -1$,
$ -1$,
$ \sqrt{3}$,
$ -\sqrt{3}$,
$ \sqrt{3}$,
$ -\sqrt{3}$,
$ 2$,
$ -2$;\ \ 
$ -1$,
$ -\sqrt{3}$,
$ \sqrt{3}$,
$ -\sqrt{3}$,
$ \sqrt{3}$,
$ 2$,
$ -2$;\ \ 
$ -\sqrt{3}$,
$ \sqrt{3}$,
$ \sqrt{3}$,
$ -\sqrt{3}$,
$0$,
$0$;\ \ 
$ \sqrt{3}$,
$ -\sqrt{3}$,
$ -\sqrt{3}$,
$0$,
$0$;\ \ 
$ -\sqrt{3}$,
$ \sqrt{3}$,
$0$,
$0$;\ \ 
$ \sqrt{3}$,
$0$,
$0$;\ \ 
$ -2$,
$ -2$;\ \ 
$ 2$)

Factors = $2_{7,2.}^{4,625}\boxtimes 5_{6,12.}^{24,592}$

  \vskip 2ex

\noindent13. $10_{3,24.}^{48,945}$ \irep{1119}:\ \ 
$d_i$ = ($1.0$,
$1.0$,
$1.0$,
$1.0$,
$1.732$,
$1.732$,
$1.732$,
$1.732$,
$2.0$,
$2.0$) 

\vskip 0.7ex
\hangindent=3em \hangafter=1
$D^2= 24.0 = 
24$

\vskip 0.7ex
\hangindent=3em \hangafter=1
$T = ( 0,
0,
\frac{1}{4},
\frac{1}{4},
\frac{3}{16},
\frac{3}{16},
\frac{11}{16},
\frac{11}{16},
\frac{1}{3},
\frac{7}{12} )
$,

\vskip 0.7ex
\hangindent=3em \hangafter=1
$S$ = ($ 1$,
$ 1$,
$ 1$,
$ 1$,
$ \sqrt{3}$,
$ \sqrt{3}$,
$ \sqrt{3}$,
$ \sqrt{3}$,
$ 2$,
$ 2$;\ \ 
$ 1$,
$ 1$,
$ 1$,
$ -\sqrt{3}$,
$ -\sqrt{3}$,
$ -\sqrt{3}$,
$ -\sqrt{3}$,
$ 2$,
$ 2$;\ \ 
$ -1$,
$ -1$,
$(-\sqrt{3})\mathrm{i}$,
$(\sqrt{3})\mathrm{i}$,
$(-\sqrt{3})\mathrm{i}$,
$(\sqrt{3})\mathrm{i}$,
$ 2$,
$ -2$;\ \ 
$ -1$,
$(\sqrt{3})\mathrm{i}$,
$(-\sqrt{3})\mathrm{i}$,
$(\sqrt{3})\mathrm{i}$,
$(-\sqrt{3})\mathrm{i}$,
$ 2$,
$ -2$;\ \ 
$ -\sqrt{3}\zeta_{8}^{3}$,
$ \sqrt{3}\zeta_{8}^{1}$,
$ \sqrt{3}\zeta_{8}^{3}$,
$ -\sqrt{3}\zeta_{8}^{1}$,
$0$,
$0$;\ \ 
$ -\sqrt{3}\zeta_{8}^{3}$,
$ -\sqrt{3}\zeta_{8}^{1}$,
$ \sqrt{3}\zeta_{8}^{3}$,
$0$,
$0$;\ \ 
$ -\sqrt{3}\zeta_{8}^{3}$,
$ \sqrt{3}\zeta_{8}^{1}$,
$0$,
$0$;\ \ 
$ -\sqrt{3}\zeta_{8}^{3}$,
$0$,
$0$;\ \ 
$ -2$,
$ -2$;\ \ 
$ 2$)

  \vskip 2ex

\noindent14. $10_{7,24.}^{48,721}$ \irep{1119}:\ \ 
$d_i$ = ($1.0$,
$1.0$,
$1.0$,
$1.0$,
$1.732$,
$1.732$,
$1.732$,
$1.732$,
$2.0$,
$2.0$) 

\vskip 0.7ex
\hangindent=3em \hangafter=1
$D^2= 24.0 = 
24$

\vskip 0.7ex
\hangindent=3em \hangafter=1
$T = ( 0,
0,
\frac{1}{4},
\frac{1}{4},
\frac{3}{16},
\frac{3}{16},
\frac{11}{16},
\frac{11}{16},
\frac{2}{3},
\frac{11}{12} )
$,

\vskip 0.7ex
\hangindent=3em \hangafter=1
$S$ = ($ 1$,
$ 1$,
$ 1$,
$ 1$,
$ \sqrt{3}$,
$ \sqrt{3}$,
$ \sqrt{3}$,
$ \sqrt{3}$,
$ 2$,
$ 2$;\ \ 
$ 1$,
$ 1$,
$ 1$,
$ -\sqrt{3}$,
$ -\sqrt{3}$,
$ -\sqrt{3}$,
$ -\sqrt{3}$,
$ 2$,
$ 2$;\ \ 
$ -1$,
$ -1$,
$(-\sqrt{3})\mathrm{i}$,
$(\sqrt{3})\mathrm{i}$,
$(-\sqrt{3})\mathrm{i}$,
$(\sqrt{3})\mathrm{i}$,
$ 2$,
$ -2$;\ \ 
$ -1$,
$(\sqrt{3})\mathrm{i}$,
$(-\sqrt{3})\mathrm{i}$,
$(\sqrt{3})\mathrm{i}$,
$(-\sqrt{3})\mathrm{i}$,
$ 2$,
$ -2$;\ \ 
$ \sqrt{3}\zeta_{8}^{3}$,
$ -\sqrt{3}\zeta_{8}^{1}$,
$ -\sqrt{3}\zeta_{8}^{3}$,
$ \sqrt{3}\zeta_{8}^{1}$,
$0$,
$0$;\ \ 
$ \sqrt{3}\zeta_{8}^{3}$,
$ \sqrt{3}\zeta_{8}^{1}$,
$ -\sqrt{3}\zeta_{8}^{3}$,
$0$,
$0$;\ \ 
$ \sqrt{3}\zeta_{8}^{3}$,
$ -\sqrt{3}\zeta_{8}^{1}$,
$0$,
$0$;\ \ 
$ \sqrt{3}\zeta_{8}^{3}$,
$0$,
$0$;\ \ 
$ -2$,
$ -2$;\ \ 
$ 2$)

  \vskip 2ex

\noindent15. $10_{3,24.}^{48,100}$ \irep{1119}:\ \ 
$d_i$ = ($1.0$,
$1.0$,
$1.0$,
$1.0$,
$1.732$,
$1.732$,
$1.732$,
$1.732$,
$2.0$,
$2.0$) 

\vskip 0.7ex
\hangindent=3em \hangafter=1
$D^2= 24.0 = 
24$

\vskip 0.7ex
\hangindent=3em \hangafter=1
$T = ( 0,
0,
\frac{1}{4},
\frac{1}{4},
\frac{7}{16},
\frac{7}{16},
\frac{15}{16},
\frac{15}{16},
\frac{1}{3},
\frac{7}{12} )
$,

\vskip 0.7ex
\hangindent=3em \hangafter=1
$S$ = ($ 1$,
$ 1$,
$ 1$,
$ 1$,
$ \sqrt{3}$,
$ \sqrt{3}$,
$ \sqrt{3}$,
$ \sqrt{3}$,
$ 2$,
$ 2$;\ \ 
$ 1$,
$ 1$,
$ 1$,
$ -\sqrt{3}$,
$ -\sqrt{3}$,
$ -\sqrt{3}$,
$ -\sqrt{3}$,
$ 2$,
$ 2$;\ \ 
$ -1$,
$ -1$,
$(-\sqrt{3})\mathrm{i}$,
$(\sqrt{3})\mathrm{i}$,
$(-\sqrt{3})\mathrm{i}$,
$(\sqrt{3})\mathrm{i}$,
$ 2$,
$ -2$;\ \ 
$ -1$,
$(\sqrt{3})\mathrm{i}$,
$(-\sqrt{3})\mathrm{i}$,
$(\sqrt{3})\mathrm{i}$,
$(-\sqrt{3})\mathrm{i}$,
$ 2$,
$ -2$;\ \ 
$ \sqrt{3}\zeta_{8}^{3}$,
$ -\sqrt{3}\zeta_{8}^{1}$,
$ -\sqrt{3}\zeta_{8}^{3}$,
$ \sqrt{3}\zeta_{8}^{1}$,
$0$,
$0$;\ \ 
$ \sqrt{3}\zeta_{8}^{3}$,
$ \sqrt{3}\zeta_{8}^{1}$,
$ -\sqrt{3}\zeta_{8}^{3}$,
$0$,
$0$;\ \ 
$ \sqrt{3}\zeta_{8}^{3}$,
$ -\sqrt{3}\zeta_{8}^{1}$,
$0$,
$0$;\ \ 
$ \sqrt{3}\zeta_{8}^{3}$,
$0$,
$0$;\ \ 
$ -2$,
$ -2$;\ \ 
$ 2$)

  \vskip 2ex

\noindent16. $10_{7,24.}^{48,267}$ \irep{1119}:\ \ 
$d_i$ = ($1.0$,
$1.0$,
$1.0$,
$1.0$,
$1.732$,
$1.732$,
$1.732$,
$1.732$,
$2.0$,
$2.0$) 

\vskip 0.7ex
\hangindent=3em \hangafter=1
$D^2= 24.0 = 
24$

\vskip 0.7ex
\hangindent=3em \hangafter=1
$T = ( 0,
0,
\frac{1}{4},
\frac{1}{4},
\frac{7}{16},
\frac{7}{16},
\frac{15}{16},
\frac{15}{16},
\frac{2}{3},
\frac{11}{12} )
$,

\vskip 0.7ex
\hangindent=3em \hangafter=1
$S$ = ($ 1$,
$ 1$,
$ 1$,
$ 1$,
$ \sqrt{3}$,
$ \sqrt{3}$,
$ \sqrt{3}$,
$ \sqrt{3}$,
$ 2$,
$ 2$;\ \ 
$ 1$,
$ 1$,
$ 1$,
$ -\sqrt{3}$,
$ -\sqrt{3}$,
$ -\sqrt{3}$,
$ -\sqrt{3}$,
$ 2$,
$ 2$;\ \ 
$ -1$,
$ -1$,
$(-\sqrt{3})\mathrm{i}$,
$(\sqrt{3})\mathrm{i}$,
$(-\sqrt{3})\mathrm{i}$,
$(\sqrt{3})\mathrm{i}$,
$ 2$,
$ -2$;\ \ 
$ -1$,
$(\sqrt{3})\mathrm{i}$,
$(-\sqrt{3})\mathrm{i}$,
$(\sqrt{3})\mathrm{i}$,
$(-\sqrt{3})\mathrm{i}$,
$ 2$,
$ -2$;\ \ 
$ -\sqrt{3}\zeta_{8}^{3}$,
$ \sqrt{3}\zeta_{8}^{1}$,
$ \sqrt{3}\zeta_{8}^{3}$,
$ -\sqrt{3}\zeta_{8}^{1}$,
$0$,
$0$;\ \ 
$ -\sqrt{3}\zeta_{8}^{3}$,
$ -\sqrt{3}\zeta_{8}^{1}$,
$ \sqrt{3}\zeta_{8}^{3}$,
$0$,
$0$;\ \ 
$ -\sqrt{3}\zeta_{8}^{3}$,
$ \sqrt{3}\zeta_{8}^{1}$,
$0$,
$0$;\ \ 
$ -\sqrt{3}\zeta_{8}^{3}$,
$0$,
$0$;\ \ 
$ -2$,
$ -2$;\ \ 
$ 2$)

  \vskip 2ex

\noindent17. $10_{1,24.}^{48,126}$ \irep{1119}:\ \ 
$d_i$ = ($1.0$,
$1.0$,
$1.0$,
$1.0$,
$1.732$,
$1.732$,
$1.732$,
$1.732$,
$2.0$,
$2.0$) 

\vskip 0.7ex
\hangindent=3em \hangafter=1
$D^2= 24.0 = 
24$

\vskip 0.7ex
\hangindent=3em \hangafter=1
$T = ( 0,
0,
\frac{3}{4},
\frac{3}{4},
\frac{1}{16},
\frac{1}{16},
\frac{9}{16},
\frac{9}{16},
\frac{1}{3},
\frac{1}{12} )
$,

\vskip 0.7ex
\hangindent=3em \hangafter=1
$S$ = ($ 1$,
$ 1$,
$ 1$,
$ 1$,
$ \sqrt{3}$,
$ \sqrt{3}$,
$ \sqrt{3}$,
$ \sqrt{3}$,
$ 2$,
$ 2$;\ \ 
$ 1$,
$ 1$,
$ 1$,
$ -\sqrt{3}$,
$ -\sqrt{3}$,
$ -\sqrt{3}$,
$ -\sqrt{3}$,
$ 2$,
$ 2$;\ \ 
$ -1$,
$ -1$,
$(-\sqrt{3})\mathrm{i}$,
$(\sqrt{3})\mathrm{i}$,
$(-\sqrt{3})\mathrm{i}$,
$(\sqrt{3})\mathrm{i}$,
$ 2$,
$ -2$;\ \ 
$ -1$,
$(\sqrt{3})\mathrm{i}$,
$(-\sqrt{3})\mathrm{i}$,
$(\sqrt{3})\mathrm{i}$,
$(-\sqrt{3})\mathrm{i}$,
$ 2$,
$ -2$;\ \ 
$ \sqrt{3}\zeta_{8}^{1}$,
$ -\sqrt{3}\zeta_{8}^{3}$,
$ -\sqrt{3}\zeta_{8}^{1}$,
$ \sqrt{3}\zeta_{8}^{3}$,
$0$,
$0$;\ \ 
$ \sqrt{3}\zeta_{8}^{1}$,
$ \sqrt{3}\zeta_{8}^{3}$,
$ -\sqrt{3}\zeta_{8}^{1}$,
$0$,
$0$;\ \ 
$ \sqrt{3}\zeta_{8}^{1}$,
$ -\sqrt{3}\zeta_{8}^{3}$,
$0$,
$0$;\ \ 
$ \sqrt{3}\zeta_{8}^{1}$,
$0$,
$0$;\ \ 
$ -2$,
$ -2$;\ \ 
$ 2$)

  \vskip 2ex

\noindent18. $10_{5,24.}^{48,261}$ \irep{1119}:\ \ 
$d_i$ = ($1.0$,
$1.0$,
$1.0$,
$1.0$,
$1.732$,
$1.732$,
$1.732$,
$1.732$,
$2.0$,
$2.0$) 

\vskip 0.7ex
\hangindent=3em \hangafter=1
$D^2= 24.0 = 
24$

\vskip 0.7ex
\hangindent=3em \hangafter=1
$T = ( 0,
0,
\frac{3}{4},
\frac{3}{4},
\frac{1}{16},
\frac{1}{16},
\frac{9}{16},
\frac{9}{16},
\frac{2}{3},
\frac{5}{12} )
$,

\vskip 0.7ex
\hangindent=3em \hangafter=1
$S$ = ($ 1$,
$ 1$,
$ 1$,
$ 1$,
$ \sqrt{3}$,
$ \sqrt{3}$,
$ \sqrt{3}$,
$ \sqrt{3}$,
$ 2$,
$ 2$;\ \ 
$ 1$,
$ 1$,
$ 1$,
$ -\sqrt{3}$,
$ -\sqrt{3}$,
$ -\sqrt{3}$,
$ -\sqrt{3}$,
$ 2$,
$ 2$;\ \ 
$ -1$,
$ -1$,
$(-\sqrt{3})\mathrm{i}$,
$(\sqrt{3})\mathrm{i}$,
$(-\sqrt{3})\mathrm{i}$,
$(\sqrt{3})\mathrm{i}$,
$ 2$,
$ -2$;\ \ 
$ -1$,
$(\sqrt{3})\mathrm{i}$,
$(-\sqrt{3})\mathrm{i}$,
$(\sqrt{3})\mathrm{i}$,
$(-\sqrt{3})\mathrm{i}$,
$ 2$,
$ -2$;\ \ 
$ -\sqrt{3}\zeta_{8}^{1}$,
$ \sqrt{3}\zeta_{8}^{3}$,
$ \sqrt{3}\zeta_{8}^{1}$,
$ -\sqrt{3}\zeta_{8}^{3}$,
$0$,
$0$;\ \ 
$ -\sqrt{3}\zeta_{8}^{1}$,
$ -\sqrt{3}\zeta_{8}^{3}$,
$ \sqrt{3}\zeta_{8}^{1}$,
$0$,
$0$;\ \ 
$ -\sqrt{3}\zeta_{8}^{1}$,
$ \sqrt{3}\zeta_{8}^{3}$,
$0$,
$0$;\ \ 
$ -\sqrt{3}\zeta_{8}^{1}$,
$0$,
$0$;\ \ 
$ -2$,
$ -2$;\ \ 
$ 2$)

  \vskip 2ex

\noindent19. $10_{1,24.}^{48,254}$ \irep{1119}:\ \ 
$d_i$ = ($1.0$,
$1.0$,
$1.0$,
$1.0$,
$1.732$,
$1.732$,
$1.732$,
$1.732$,
$2.0$,
$2.0$) 

\vskip 0.7ex
\hangindent=3em \hangafter=1
$D^2= 24.0 = 
24$

\vskip 0.7ex
\hangindent=3em \hangafter=1
$T = ( 0,
0,
\frac{3}{4},
\frac{3}{4},
\frac{5}{16},
\frac{5}{16},
\frac{13}{16},
\frac{13}{16},
\frac{1}{3},
\frac{1}{12} )
$,

\vskip 0.7ex
\hangindent=3em \hangafter=1
$S$ = ($ 1$,
$ 1$,
$ 1$,
$ 1$,
$ \sqrt{3}$,
$ \sqrt{3}$,
$ \sqrt{3}$,
$ \sqrt{3}$,
$ 2$,
$ 2$;\ \ 
$ 1$,
$ 1$,
$ 1$,
$ -\sqrt{3}$,
$ -\sqrt{3}$,
$ -\sqrt{3}$,
$ -\sqrt{3}$,
$ 2$,
$ 2$;\ \ 
$ -1$,
$ -1$,
$(-\sqrt{3})\mathrm{i}$,
$(\sqrt{3})\mathrm{i}$,
$(-\sqrt{3})\mathrm{i}$,
$(\sqrt{3})\mathrm{i}$,
$ 2$,
$ -2$;\ \ 
$ -1$,
$(\sqrt{3})\mathrm{i}$,
$(-\sqrt{3})\mathrm{i}$,
$(\sqrt{3})\mathrm{i}$,
$(-\sqrt{3})\mathrm{i}$,
$ 2$,
$ -2$;\ \ 
$ -\sqrt{3}\zeta_{8}^{1}$,
$ \sqrt{3}\zeta_{8}^{3}$,
$ \sqrt{3}\zeta_{8}^{1}$,
$ -\sqrt{3}\zeta_{8}^{3}$,
$0$,
$0$;\ \ 
$ -\sqrt{3}\zeta_{8}^{1}$,
$ -\sqrt{3}\zeta_{8}^{3}$,
$ \sqrt{3}\zeta_{8}^{1}$,
$0$,
$0$;\ \ 
$ -\sqrt{3}\zeta_{8}^{1}$,
$ \sqrt{3}\zeta_{8}^{3}$,
$0$,
$0$;\ \ 
$ -\sqrt{3}\zeta_{8}^{1}$,
$0$,
$0$;\ \ 
$ -2$,
$ -2$;\ \ 
$ 2$)

  \vskip 2ex

\noindent20. $10_{5,24.}^{48,125}$ \irep{1119}:\ \ 
$d_i$ = ($1.0$,
$1.0$,
$1.0$,
$1.0$,
$1.732$,
$1.732$,
$1.732$,
$1.732$,
$2.0$,
$2.0$) 

\vskip 0.7ex
\hangindent=3em \hangafter=1
$D^2= 24.0 = 
24$

\vskip 0.7ex
\hangindent=3em \hangafter=1
$T = ( 0,
0,
\frac{3}{4},
\frac{3}{4},
\frac{5}{16},
\frac{5}{16},
\frac{13}{16},
\frac{13}{16},
\frac{2}{3},
\frac{5}{12} )
$,

\vskip 0.7ex
\hangindent=3em \hangafter=1
$S$ = ($ 1$,
$ 1$,
$ 1$,
$ 1$,
$ \sqrt{3}$,
$ \sqrt{3}$,
$ \sqrt{3}$,
$ \sqrt{3}$,
$ 2$,
$ 2$;\ \ 
$ 1$,
$ 1$,
$ 1$,
$ -\sqrt{3}$,
$ -\sqrt{3}$,
$ -\sqrt{3}$,
$ -\sqrt{3}$,
$ 2$,
$ 2$;\ \ 
$ -1$,
$ -1$,
$(-\sqrt{3})\mathrm{i}$,
$(\sqrt{3})\mathrm{i}$,
$(-\sqrt{3})\mathrm{i}$,
$(\sqrt{3})\mathrm{i}$,
$ 2$,
$ -2$;\ \ 
$ -1$,
$(\sqrt{3})\mathrm{i}$,
$(-\sqrt{3})\mathrm{i}$,
$(\sqrt{3})\mathrm{i}$,
$(-\sqrt{3})\mathrm{i}$,
$ 2$,
$ -2$;\ \ 
$ \sqrt{3}\zeta_{8}^{1}$,
$ -\sqrt{3}\zeta_{8}^{3}$,
$ -\sqrt{3}\zeta_{8}^{1}$,
$ \sqrt{3}\zeta_{8}^{3}$,
$0$,
$0$;\ \ 
$ \sqrt{3}\zeta_{8}^{1}$,
$ \sqrt{3}\zeta_{8}^{3}$,
$ -\sqrt{3}\zeta_{8}^{1}$,
$0$,
$0$;\ \ 
$ \sqrt{3}\zeta_{8}^{1}$,
$ -\sqrt{3}\zeta_{8}^{3}$,
$0$,
$0$;\ \ 
$ \sqrt{3}\zeta_{8}^{1}$,
$0$,
$0$;\ \ 
$ -2$,
$ -2$;\ \ 
$ 2$)

  \vskip 2ex

\noindent21. $10_{4,36.}^{6,152}$ \irep{0}:\ \ 
$d_i$ = ($1.0$,
$1.0$,
$1.0$,
$2.0$,
$2.0$,
$2.0$,
$2.0$,
$2.0$,
$2.0$,
$3.0$) 

\vskip 0.7ex
\hangindent=3em \hangafter=1
$D^2= 36.0 = 
36$

\vskip 0.7ex
\hangindent=3em \hangafter=1
$T = ( 0,
0,
0,
0,
0,
\frac{1}{3},
\frac{1}{3},
\frac{2}{3},
\frac{2}{3},
\frac{1}{2} )
$,

\vskip 0.7ex
\hangindent=3em \hangafter=1
$S$ = ($ 1$,
$ 1$,
$ 1$,
$ 2$,
$ 2$,
$ 2$,
$ 2$,
$ 2$,
$ 2$,
$ 3$;\ \ 
$ 1$,
$ 1$,
$ -2\zeta_{6}^{1}$,
$ 2\zeta_{3}^{1}$,
$ -2\zeta_{6}^{1}$,
$ 2\zeta_{3}^{1}$,
$ -2\zeta_{6}^{1}$,
$ 2\zeta_{3}^{1}$,
$ 3$;\ \ 
$ 1$,
$ 2\zeta_{3}^{1}$,
$ -2\zeta_{6}^{1}$,
$ 2\zeta_{3}^{1}$,
$ -2\zeta_{6}^{1}$,
$ 2\zeta_{3}^{1}$,
$ -2\zeta_{6}^{1}$,
$ 3$;\ \ 
$ -2$,
$ -2$,
$ 2\zeta_{6}^{1}$,
$ -2\zeta_{3}^{1}$,
$ -2\zeta_{3}^{1}$,
$ 2\zeta_{6}^{1}$,
$0$;\ \ 
$ -2$,
$ -2\zeta_{3}^{1}$,
$ 2\zeta_{6}^{1}$,
$ 2\zeta_{6}^{1}$,
$ -2\zeta_{3}^{1}$,
$0$;\ \ 
$ -2\zeta_{3}^{1}$,
$ 2\zeta_{6}^{1}$,
$ -2$,
$ -2$,
$0$;\ \ 
$ -2\zeta_{3}^{1}$,
$ -2$,
$ -2$,
$0$;\ \ 
$ 2\zeta_{6}^{1}$,
$ -2\zeta_{3}^{1}$,
$0$;\ \ 
$ 2\zeta_{6}^{1}$,
$0$;\ \ 
$ -3$)

  \vskip 2ex

\noindent22. $10_{4,36.}^{18,490}$ \irep{0}:\ \ 
$d_i$ = ($1.0$,
$1.0$,
$1.0$,
$2.0$,
$2.0$,
$2.0$,
$2.0$,
$2.0$,
$2.0$,
$3.0$) 

\vskip 0.7ex
\hangindent=3em \hangafter=1
$D^2= 36.0 = 
36$

\vskip 0.7ex
\hangindent=3em \hangafter=1
$T = ( 0,
0,
0,
\frac{1}{9},
\frac{1}{9},
\frac{4}{9},
\frac{4}{9},
\frac{7}{9},
\frac{7}{9},
\frac{1}{2} )
$,

\vskip 0.7ex
\hangindent=3em \hangafter=1
$S$ = ($ 1$,
$ 1$,
$ 1$,
$ 2$,
$ 2$,
$ 2$,
$ 2$,
$ 2$,
$ 2$,
$ 3$;\ \ 
$ 1$,
$ 1$,
$ -2\zeta_{6}^{1}$,
$ 2\zeta_{3}^{1}$,
$ -2\zeta_{6}^{1}$,
$ 2\zeta_{3}^{1}$,
$ -2\zeta_{6}^{1}$,
$ 2\zeta_{3}^{1}$,
$ 3$;\ \ 
$ 1$,
$ 2\zeta_{3}^{1}$,
$ -2\zeta_{6}^{1}$,
$ 2\zeta_{3}^{1}$,
$ -2\zeta_{6}^{1}$,
$ 2\zeta_{3}^{1}$,
$ -2\zeta_{6}^{1}$,
$ 3$;\ \ 
$ 2\zeta_{18}^{5}$,
$ -2\zeta_{9}^{2}$,
$ -2\zeta_{9}^{4}$,
$ 2\zeta_{18}^{1}$,
$ -2\zeta_{9}^{1}$,
$ 2\zeta_{18}^{7}$,
$0$;\ \ 
$ 2\zeta_{18}^{5}$,
$ 2\zeta_{18}^{1}$,
$ -2\zeta_{9}^{4}$,
$ 2\zeta_{18}^{7}$,
$ -2\zeta_{9}^{1}$,
$0$;\ \ 
$ -2\zeta_{9}^{1}$,
$ 2\zeta_{18}^{7}$,
$ 2\zeta_{18}^{5}$,
$ -2\zeta_{9}^{2}$,
$0$;\ \ 
$ -2\zeta_{9}^{1}$,
$ -2\zeta_{9}^{2}$,
$ 2\zeta_{18}^{5}$,
$0$;\ \ 
$ -2\zeta_{9}^{4}$,
$ 2\zeta_{18}^{1}$,
$0$;\ \ 
$ -2\zeta_{9}^{4}$,
$0$;\ \ 
$ -3$)

  \vskip 2ex

\noindent23. $10_{4,36.}^{18,842}$ \irep{0}:\ \ 
$d_i$ = ($1.0$,
$1.0$,
$1.0$,
$2.0$,
$2.0$,
$2.0$,
$2.0$,
$2.0$,
$2.0$,
$3.0$) 

\vskip 0.7ex
\hangindent=3em \hangafter=1
$D^2= 36.0 = 
36$

\vskip 0.7ex
\hangindent=3em \hangafter=1
$T = ( 0,
0,
0,
\frac{2}{9},
\frac{2}{9},
\frac{5}{9},
\frac{5}{9},
\frac{8}{9},
\frac{8}{9},
\frac{1}{2} )
$,

\vskip 0.7ex
\hangindent=3em \hangafter=1
$S$ = ($ 1$,
$ 1$,
$ 1$,
$ 2$,
$ 2$,
$ 2$,
$ 2$,
$ 2$,
$ 2$,
$ 3$;\ \ 
$ 1$,
$ 1$,
$ -2\zeta_{6}^{1}$,
$ 2\zeta_{3}^{1}$,
$ -2\zeta_{6}^{1}$,
$ 2\zeta_{3}^{1}$,
$ -2\zeta_{6}^{1}$,
$ 2\zeta_{3}^{1}$,
$ 3$;\ \ 
$ 1$,
$ 2\zeta_{3}^{1}$,
$ -2\zeta_{6}^{1}$,
$ 2\zeta_{3}^{1}$,
$ -2\zeta_{6}^{1}$,
$ 2\zeta_{3}^{1}$,
$ -2\zeta_{6}^{1}$,
$ 3$;\ \ 
$ 2\zeta_{18}^{1}$,
$ -2\zeta_{9}^{4}$,
$ -2\zeta_{9}^{2}$,
$ 2\zeta_{18}^{5}$,
$ 2\zeta_{18}^{7}$,
$ -2\zeta_{9}^{1}$,
$0$;\ \ 
$ 2\zeta_{18}^{1}$,
$ 2\zeta_{18}^{5}$,
$ -2\zeta_{9}^{2}$,
$ -2\zeta_{9}^{1}$,
$ 2\zeta_{18}^{7}$,
$0$;\ \ 
$ 2\zeta_{18}^{7}$,
$ -2\zeta_{9}^{1}$,
$ 2\zeta_{18}^{1}$,
$ -2\zeta_{9}^{4}$,
$0$;\ \ 
$ 2\zeta_{18}^{7}$,
$ -2\zeta_{9}^{4}$,
$ 2\zeta_{18}^{1}$,
$0$;\ \ 
$ -2\zeta_{9}^{2}$,
$ 2\zeta_{18}^{5}$,
$0$;\ \ 
$ -2\zeta_{9}^{2}$,
$0$;\ \ 
$ -3$)

  \vskip 2ex

\noindent24. $10_{\frac{24}{5},43.41}^{120,927}$ \irep{1146}:\ \ 
$d_i$ = ($1.0$,
$1.0$,
$1.618$,
$1.618$,
$1.732$,
$1.732$,
$2.0$,
$2.802$,
$2.802$,
$3.236$) 

\vskip 0.7ex
\hangindent=3em \hangafter=1
$D^2= 43.416 = 
30+6\sqrt{5}$

\vskip 0.7ex
\hangindent=3em \hangafter=1
$T = ( 0,
0,
\frac{2}{5},
\frac{2}{5},
\frac{1}{8},
\frac{5}{8},
\frac{1}{3},
\frac{1}{40},
\frac{21}{40},
\frac{11}{15} )
$,

\vskip 0.7ex
\hangindent=3em \hangafter=1
$S$ = ($ 1$,
$ 1$,
$ \frac{1+\sqrt{5}}{2}$,
$ \frac{1+\sqrt{5}}{2}$,
$ \sqrt{3}$,
$ \sqrt{3}$,
$ 2$,
$ \frac{15+3\sqrt{5}}{2\sqrt{15}}$,
$ \frac{15+3\sqrt{5}}{2\sqrt{15}}$,
$ 1+\sqrt{5}$;\ \ 
$ 1$,
$ \frac{1+\sqrt{5}}{2}$,
$ \frac{1+\sqrt{5}}{2}$,
$ -\sqrt{3}$,
$ -\sqrt{3}$,
$ 2$,
$ -\frac{15+3\sqrt{5}}{2\sqrt{15}}$,
$ -\frac{15+3\sqrt{5}}{2\sqrt{15}}$,
$ 1+\sqrt{5}$;\ \ 
$ -1$,
$ -1$,
$ -\frac{15+3\sqrt{5}}{2\sqrt{15}}$,
$ -\frac{15+3\sqrt{5}}{2\sqrt{15}}$,
$ 1+\sqrt{5}$,
$ \sqrt{3}$,
$ \sqrt{3}$,
$ -2$;\ \ 
$ -1$,
$ \frac{15+3\sqrt{5}}{2\sqrt{15}}$,
$ \frac{15+3\sqrt{5}}{2\sqrt{15}}$,
$ 1+\sqrt{5}$,
$ -\sqrt{3}$,
$ -\sqrt{3}$,
$ -2$;\ \ 
$ \sqrt{3}$,
$ -\sqrt{3}$,
$0$,
$ -\frac{15+3\sqrt{5}}{2\sqrt{15}}$,
$ \frac{15+3\sqrt{5}}{2\sqrt{15}}$,
$0$;\ \ 
$ \sqrt{3}$,
$0$,
$ \frac{15+3\sqrt{5}}{2\sqrt{15}}$,
$ -\frac{15+3\sqrt{5}}{2\sqrt{15}}$,
$0$;\ \ 
$ -2$,
$0$,
$0$,
$ -1-\sqrt{5}$;\ \ 
$ -\sqrt{3}$,
$ \sqrt{3}$,
$0$;\ \ 
$ -\sqrt{3}$,
$0$;\ \ 
$ 2$)

Factors = $2_{\frac{14}{5},3.618}^{5,395}\boxtimes 5_{2,12.}^{24,940}$

  \vskip 2ex

\noindent25. $10_{\frac{4}{5},43.41}^{120,198}$ \irep{1146}:\ \ 
$d_i$ = ($1.0$,
$1.0$,
$1.618$,
$1.618$,
$1.732$,
$1.732$,
$2.0$,
$2.802$,
$2.802$,
$3.236$) 

\vskip 0.7ex
\hangindent=3em \hangafter=1
$D^2= 43.416 = 
30+6\sqrt{5}$

\vskip 0.7ex
\hangindent=3em \hangafter=1
$T = ( 0,
0,
\frac{2}{5},
\frac{2}{5},
\frac{1}{8},
\frac{5}{8},
\frac{2}{3},
\frac{1}{40},
\frac{21}{40},
\frac{1}{15} )
$,

\vskip 0.7ex
\hangindent=3em \hangafter=1
$S$ = ($ 1$,
$ 1$,
$ \frac{1+\sqrt{5}}{2}$,
$ \frac{1+\sqrt{5}}{2}$,
$ \sqrt{3}$,
$ \sqrt{3}$,
$ 2$,
$ \frac{15+3\sqrt{5}}{2\sqrt{15}}$,
$ \frac{15+3\sqrt{5}}{2\sqrt{15}}$,
$ 1+\sqrt{5}$;\ \ 
$ 1$,
$ \frac{1+\sqrt{5}}{2}$,
$ \frac{1+\sqrt{5}}{2}$,
$ -\sqrt{3}$,
$ -\sqrt{3}$,
$ 2$,
$ -\frac{15+3\sqrt{5}}{2\sqrt{15}}$,
$ -\frac{15+3\sqrt{5}}{2\sqrt{15}}$,
$ 1+\sqrt{5}$;\ \ 
$ -1$,
$ -1$,
$ -\frac{15+3\sqrt{5}}{2\sqrt{15}}$,
$ -\frac{15+3\sqrt{5}}{2\sqrt{15}}$,
$ 1+\sqrt{5}$,
$ \sqrt{3}$,
$ \sqrt{3}$,
$ -2$;\ \ 
$ -1$,
$ \frac{15+3\sqrt{5}}{2\sqrt{15}}$,
$ \frac{15+3\sqrt{5}}{2\sqrt{15}}$,
$ 1+\sqrt{5}$,
$ -\sqrt{3}$,
$ -\sqrt{3}$,
$ -2$;\ \ 
$ -\sqrt{3}$,
$ \sqrt{3}$,
$0$,
$ \frac{15+3\sqrt{5}}{2\sqrt{15}}$,
$ -\frac{15+3\sqrt{5}}{2\sqrt{15}}$,
$0$;\ \ 
$ -\sqrt{3}$,
$0$,
$ -\frac{15+3\sqrt{5}}{2\sqrt{15}}$,
$ \frac{15+3\sqrt{5}}{2\sqrt{15}}$,
$0$;\ \ 
$ -2$,
$0$,
$0$,
$ -1-\sqrt{5}$;\ \ 
$ \sqrt{3}$,
$ -\sqrt{3}$,
$0$;\ \ 
$ \sqrt{3}$,
$0$;\ \ 
$ 2$)

Factors = $2_{\frac{14}{5},3.618}^{5,395}\boxtimes 5_{6,12.}^{24,273}$

  \vskip 2ex

\noindent26. $10_{\frac{24}{5},43.41}^{120,105}$ \irep{1146}:\ \ 
$d_i$ = ($1.0$,
$1.0$,
$1.618$,
$1.618$,
$1.732$,
$1.732$,
$2.0$,
$2.802$,
$2.802$,
$3.236$) 

\vskip 0.7ex
\hangindent=3em \hangafter=1
$D^2= 43.416 = 
30+6\sqrt{5}$

\vskip 0.7ex
\hangindent=3em \hangafter=1
$T = ( 0,
0,
\frac{2}{5},
\frac{2}{5},
\frac{3}{8},
\frac{7}{8},
\frac{1}{3},
\frac{11}{40},
\frac{31}{40},
\frac{11}{15} )
$,

\vskip 0.7ex
\hangindent=3em \hangafter=1
$S$ = ($ 1$,
$ 1$,
$ \frac{1+\sqrt{5}}{2}$,
$ \frac{1+\sqrt{5}}{2}$,
$ \sqrt{3}$,
$ \sqrt{3}$,
$ 2$,
$ \frac{15+3\sqrt{5}}{2\sqrt{15}}$,
$ \frac{15+3\sqrt{5}}{2\sqrt{15}}$,
$ 1+\sqrt{5}$;\ \ 
$ 1$,
$ \frac{1+\sqrt{5}}{2}$,
$ \frac{1+\sqrt{5}}{2}$,
$ -\sqrt{3}$,
$ -\sqrt{3}$,
$ 2$,
$ -\frac{15+3\sqrt{5}}{2\sqrt{15}}$,
$ -\frac{15+3\sqrt{5}}{2\sqrt{15}}$,
$ 1+\sqrt{5}$;\ \ 
$ -1$,
$ -1$,
$ -\frac{15+3\sqrt{5}}{2\sqrt{15}}$,
$ -\frac{15+3\sqrt{5}}{2\sqrt{15}}$,
$ 1+\sqrt{5}$,
$ \sqrt{3}$,
$ \sqrt{3}$,
$ -2$;\ \ 
$ -1$,
$ \frac{15+3\sqrt{5}}{2\sqrt{15}}$,
$ \frac{15+3\sqrt{5}}{2\sqrt{15}}$,
$ 1+\sqrt{5}$,
$ -\sqrt{3}$,
$ -\sqrt{3}$,
$ -2$;\ \ 
$ -\sqrt{3}$,
$ \sqrt{3}$,
$0$,
$ \frac{15+3\sqrt{5}}{2\sqrt{15}}$,
$ -\frac{15+3\sqrt{5}}{2\sqrt{15}}$,
$0$;\ \ 
$ -\sqrt{3}$,
$0$,
$ -\frac{15+3\sqrt{5}}{2\sqrt{15}}$,
$ \frac{15+3\sqrt{5}}{2\sqrt{15}}$,
$0$;\ \ 
$ -2$,
$0$,
$0$,
$ -1-\sqrt{5}$;\ \ 
$ \sqrt{3}$,
$ -\sqrt{3}$,
$0$;\ \ 
$ \sqrt{3}$,
$0$;\ \ 
$ 2$)

Factors = $2_{\frac{14}{5},3.618}^{5,395}\boxtimes 5_{2,12.}^{24,741}$

  \vskip 2ex

\noindent27. $10_{\frac{4}{5},43.41}^{120,597}$ \irep{1146}:\ \ 
$d_i$ = ($1.0$,
$1.0$,
$1.618$,
$1.618$,
$1.732$,
$1.732$,
$2.0$,
$2.802$,
$2.802$,
$3.236$) 

\vskip 0.7ex
\hangindent=3em \hangafter=1
$D^2= 43.416 = 
30+6\sqrt{5}$

\vskip 0.7ex
\hangindent=3em \hangafter=1
$T = ( 0,
0,
\frac{2}{5},
\frac{2}{5},
\frac{3}{8},
\frac{7}{8},
\frac{2}{3},
\frac{11}{40},
\frac{31}{40},
\frac{1}{15} )
$,

\vskip 0.7ex
\hangindent=3em \hangafter=1
$S$ = ($ 1$,
$ 1$,
$ \frac{1+\sqrt{5}}{2}$,
$ \frac{1+\sqrt{5}}{2}$,
$ \sqrt{3}$,
$ \sqrt{3}$,
$ 2$,
$ \frac{15+3\sqrt{5}}{2\sqrt{15}}$,
$ \frac{15+3\sqrt{5}}{2\sqrt{15}}$,
$ 1+\sqrt{5}$;\ \ 
$ 1$,
$ \frac{1+\sqrt{5}}{2}$,
$ \frac{1+\sqrt{5}}{2}$,
$ -\sqrt{3}$,
$ -\sqrt{3}$,
$ 2$,
$ -\frac{15+3\sqrt{5}}{2\sqrt{15}}$,
$ -\frac{15+3\sqrt{5}}{2\sqrt{15}}$,
$ 1+\sqrt{5}$;\ \ 
$ -1$,
$ -1$,
$ -\frac{15+3\sqrt{5}}{2\sqrt{15}}$,
$ -\frac{15+3\sqrt{5}}{2\sqrt{15}}$,
$ 1+\sqrt{5}$,
$ \sqrt{3}$,
$ \sqrt{3}$,
$ -2$;\ \ 
$ -1$,
$ \frac{15+3\sqrt{5}}{2\sqrt{15}}$,
$ \frac{15+3\sqrt{5}}{2\sqrt{15}}$,
$ 1+\sqrt{5}$,
$ -\sqrt{3}$,
$ -\sqrt{3}$,
$ -2$;\ \ 
$ \sqrt{3}$,
$ -\sqrt{3}$,
$0$,
$ -\frac{15+3\sqrt{5}}{2\sqrt{15}}$,
$ \frac{15+3\sqrt{5}}{2\sqrt{15}}$,
$0$;\ \ 
$ \sqrt{3}$,
$0$,
$ \frac{15+3\sqrt{5}}{2\sqrt{15}}$,
$ -\frac{15+3\sqrt{5}}{2\sqrt{15}}$,
$0$;\ \ 
$ -2$,
$0$,
$0$,
$ -1-\sqrt{5}$;\ \ 
$ -\sqrt{3}$,
$ \sqrt{3}$,
$0$;\ \ 
$ -\sqrt{3}$,
$0$;\ \ 
$ 2$)

Factors = $2_{\frac{14}{5},3.618}^{5,395}\boxtimes 5_{6,12.}^{24,592}$

  \vskip 2ex

\noindent28. $10_{\frac{36}{5},43.41}^{120,163}$ \irep{1146}:\ \ 
$d_i$ = ($1.0$,
$1.0$,
$1.618$,
$1.618$,
$1.732$,
$1.732$,
$2.0$,
$2.802$,
$2.802$,
$3.236$) 

\vskip 0.7ex
\hangindent=3em \hangafter=1
$D^2= 43.416 = 
30+6\sqrt{5}$

\vskip 0.7ex
\hangindent=3em \hangafter=1
$T = ( 0,
0,
\frac{3}{5},
\frac{3}{5},
\frac{1}{8},
\frac{5}{8},
\frac{1}{3},
\frac{9}{40},
\frac{29}{40},
\frac{14}{15} )
$,

\vskip 0.7ex
\hangindent=3em \hangafter=1
$S$ = ($ 1$,
$ 1$,
$ \frac{1+\sqrt{5}}{2}$,
$ \frac{1+\sqrt{5}}{2}$,
$ \sqrt{3}$,
$ \sqrt{3}$,
$ 2$,
$ \frac{15+3\sqrt{5}}{2\sqrt{15}}$,
$ \frac{15+3\sqrt{5}}{2\sqrt{15}}$,
$ 1+\sqrt{5}$;\ \ 
$ 1$,
$ \frac{1+\sqrt{5}}{2}$,
$ \frac{1+\sqrt{5}}{2}$,
$ -\sqrt{3}$,
$ -\sqrt{3}$,
$ 2$,
$ -\frac{15+3\sqrt{5}}{2\sqrt{15}}$,
$ -\frac{15+3\sqrt{5}}{2\sqrt{15}}$,
$ 1+\sqrt{5}$;\ \ 
$ -1$,
$ -1$,
$ -\frac{15+3\sqrt{5}}{2\sqrt{15}}$,
$ -\frac{15+3\sqrt{5}}{2\sqrt{15}}$,
$ 1+\sqrt{5}$,
$ \sqrt{3}$,
$ \sqrt{3}$,
$ -2$;\ \ 
$ -1$,
$ \frac{15+3\sqrt{5}}{2\sqrt{15}}$,
$ \frac{15+3\sqrt{5}}{2\sqrt{15}}$,
$ 1+\sqrt{5}$,
$ -\sqrt{3}$,
$ -\sqrt{3}$,
$ -2$;\ \ 
$ \sqrt{3}$,
$ -\sqrt{3}$,
$0$,
$ -\frac{15+3\sqrt{5}}{2\sqrt{15}}$,
$ \frac{15+3\sqrt{5}}{2\sqrt{15}}$,
$0$;\ \ 
$ \sqrt{3}$,
$0$,
$ \frac{15+3\sqrt{5}}{2\sqrt{15}}$,
$ -\frac{15+3\sqrt{5}}{2\sqrt{15}}$,
$0$;\ \ 
$ -2$,
$0$,
$0$,
$ -1-\sqrt{5}$;\ \ 
$ -\sqrt{3}$,
$ \sqrt{3}$,
$0$;\ \ 
$ -\sqrt{3}$,
$0$;\ \ 
$ 2$)

Factors = $2_{\frac{26}{5},3.618}^{5,720}\boxtimes 5_{2,12.}^{24,940}$

  \vskip 2ex

\noindent29. $10_{\frac{16}{5},43.41}^{120,282}$ \irep{1146}:\ \ 
$d_i$ = ($1.0$,
$1.0$,
$1.618$,
$1.618$,
$1.732$,
$1.732$,
$2.0$,
$2.802$,
$2.802$,
$3.236$) 

\vskip 0.7ex
\hangindent=3em \hangafter=1
$D^2= 43.416 = 
30+6\sqrt{5}$

\vskip 0.7ex
\hangindent=3em \hangafter=1
$T = ( 0,
0,
\frac{3}{5},
\frac{3}{5},
\frac{1}{8},
\frac{5}{8},
\frac{2}{3},
\frac{9}{40},
\frac{29}{40},
\frac{4}{15} )
$,

\vskip 0.7ex
\hangindent=3em \hangafter=1
$S$ = ($ 1$,
$ 1$,
$ \frac{1+\sqrt{5}}{2}$,
$ \frac{1+\sqrt{5}}{2}$,
$ \sqrt{3}$,
$ \sqrt{3}$,
$ 2$,
$ \frac{15+3\sqrt{5}}{2\sqrt{15}}$,
$ \frac{15+3\sqrt{5}}{2\sqrt{15}}$,
$ 1+\sqrt{5}$;\ \ 
$ 1$,
$ \frac{1+\sqrt{5}}{2}$,
$ \frac{1+\sqrt{5}}{2}$,
$ -\sqrt{3}$,
$ -\sqrt{3}$,
$ 2$,
$ -\frac{15+3\sqrt{5}}{2\sqrt{15}}$,
$ -\frac{15+3\sqrt{5}}{2\sqrt{15}}$,
$ 1+\sqrt{5}$;\ \ 
$ -1$,
$ -1$,
$ -\frac{15+3\sqrt{5}}{2\sqrt{15}}$,
$ -\frac{15+3\sqrt{5}}{2\sqrt{15}}$,
$ 1+\sqrt{5}$,
$ \sqrt{3}$,
$ \sqrt{3}$,
$ -2$;\ \ 
$ -1$,
$ \frac{15+3\sqrt{5}}{2\sqrt{15}}$,
$ \frac{15+3\sqrt{5}}{2\sqrt{15}}$,
$ 1+\sqrt{5}$,
$ -\sqrt{3}$,
$ -\sqrt{3}$,
$ -2$;\ \ 
$ -\sqrt{3}$,
$ \sqrt{3}$,
$0$,
$ \frac{15+3\sqrt{5}}{2\sqrt{15}}$,
$ -\frac{15+3\sqrt{5}}{2\sqrt{15}}$,
$0$;\ \ 
$ -\sqrt{3}$,
$0$,
$ -\frac{15+3\sqrt{5}}{2\sqrt{15}}$,
$ \frac{15+3\sqrt{5}}{2\sqrt{15}}$,
$0$;\ \ 
$ -2$,
$0$,
$0$,
$ -1-\sqrt{5}$;\ \ 
$ \sqrt{3}$,
$ -\sqrt{3}$,
$0$;\ \ 
$ \sqrt{3}$,
$0$;\ \ 
$ 2$)

Factors = $2_{\frac{26}{5},3.618}^{5,720}\boxtimes 5_{6,12.}^{24,273}$

  \vskip 2ex

\noindent30. $10_{\frac{36}{5},43.41}^{120,418}$ \irep{1146}:\ \ 
$d_i$ = ($1.0$,
$1.0$,
$1.618$,
$1.618$,
$1.732$,
$1.732$,
$2.0$,
$2.802$,
$2.802$,
$3.236$) 

\vskip 0.7ex
\hangindent=3em \hangafter=1
$D^2= 43.416 = 
30+6\sqrt{5}$

\vskip 0.7ex
\hangindent=3em \hangafter=1
$T = ( 0,
0,
\frac{3}{5},
\frac{3}{5},
\frac{3}{8},
\frac{7}{8},
\frac{1}{3},
\frac{19}{40},
\frac{39}{40},
\frac{14}{15} )
$,

\vskip 0.7ex
\hangindent=3em \hangafter=1
$S$ = ($ 1$,
$ 1$,
$ \frac{1+\sqrt{5}}{2}$,
$ \frac{1+\sqrt{5}}{2}$,
$ \sqrt{3}$,
$ \sqrt{3}$,
$ 2$,
$ \frac{15+3\sqrt{5}}{2\sqrt{15}}$,
$ \frac{15+3\sqrt{5}}{2\sqrt{15}}$,
$ 1+\sqrt{5}$;\ \ 
$ 1$,
$ \frac{1+\sqrt{5}}{2}$,
$ \frac{1+\sqrt{5}}{2}$,
$ -\sqrt{3}$,
$ -\sqrt{3}$,
$ 2$,
$ -\frac{15+3\sqrt{5}}{2\sqrt{15}}$,
$ -\frac{15+3\sqrt{5}}{2\sqrt{15}}$,
$ 1+\sqrt{5}$;\ \ 
$ -1$,
$ -1$,
$ -\frac{15+3\sqrt{5}}{2\sqrt{15}}$,
$ -\frac{15+3\sqrt{5}}{2\sqrt{15}}$,
$ 1+\sqrt{5}$,
$ \sqrt{3}$,
$ \sqrt{3}$,
$ -2$;\ \ 
$ -1$,
$ \frac{15+3\sqrt{5}}{2\sqrt{15}}$,
$ \frac{15+3\sqrt{5}}{2\sqrt{15}}$,
$ 1+\sqrt{5}$,
$ -\sqrt{3}$,
$ -\sqrt{3}$,
$ -2$;\ \ 
$ -\sqrt{3}$,
$ \sqrt{3}$,
$0$,
$ \frac{15+3\sqrt{5}}{2\sqrt{15}}$,
$ -\frac{15+3\sqrt{5}}{2\sqrt{15}}$,
$0$;\ \ 
$ -\sqrt{3}$,
$0$,
$ -\frac{15+3\sqrt{5}}{2\sqrt{15}}$,
$ \frac{15+3\sqrt{5}}{2\sqrt{15}}$,
$0$;\ \ 
$ -2$,
$0$,
$0$,
$ -1-\sqrt{5}$;\ \ 
$ \sqrt{3}$,
$ -\sqrt{3}$,
$0$;\ \ 
$ \sqrt{3}$,
$0$;\ \ 
$ 2$)

Factors = $2_{\frac{26}{5},3.618}^{5,720}\boxtimes 5_{2,12.}^{24,741}$

  \vskip 2ex

\noindent31. $10_{\frac{16}{5},43.41}^{120,226}$ \irep{1146}:\ \ 
$d_i$ = ($1.0$,
$1.0$,
$1.618$,
$1.618$,
$1.732$,
$1.732$,
$2.0$,
$2.802$,
$2.802$,
$3.236$) 

\vskip 0.7ex
\hangindent=3em \hangafter=1
$D^2= 43.416 = 
30+6\sqrt{5}$

\vskip 0.7ex
\hangindent=3em \hangafter=1
$T = ( 0,
0,
\frac{3}{5},
\frac{3}{5},
\frac{3}{8},
\frac{7}{8},
\frac{2}{3},
\frac{19}{40},
\frac{39}{40},
\frac{4}{15} )
$,

\vskip 0.7ex
\hangindent=3em \hangafter=1
$S$ = ($ 1$,
$ 1$,
$ \frac{1+\sqrt{5}}{2}$,
$ \frac{1+\sqrt{5}}{2}$,
$ \sqrt{3}$,
$ \sqrt{3}$,
$ 2$,
$ \frac{15+3\sqrt{5}}{2\sqrt{15}}$,
$ \frac{15+3\sqrt{5}}{2\sqrt{15}}$,
$ 1+\sqrt{5}$;\ \ 
$ 1$,
$ \frac{1+\sqrt{5}}{2}$,
$ \frac{1+\sqrt{5}}{2}$,
$ -\sqrt{3}$,
$ -\sqrt{3}$,
$ 2$,
$ -\frac{15+3\sqrt{5}}{2\sqrt{15}}$,
$ -\frac{15+3\sqrt{5}}{2\sqrt{15}}$,
$ 1+\sqrt{5}$;\ \ 
$ -1$,
$ -1$,
$ -\frac{15+3\sqrt{5}}{2\sqrt{15}}$,
$ -\frac{15+3\sqrt{5}}{2\sqrt{15}}$,
$ 1+\sqrt{5}$,
$ \sqrt{3}$,
$ \sqrt{3}$,
$ -2$;\ \ 
$ -1$,
$ \frac{15+3\sqrt{5}}{2\sqrt{15}}$,
$ \frac{15+3\sqrt{5}}{2\sqrt{15}}$,
$ 1+\sqrt{5}$,
$ -\sqrt{3}$,
$ -\sqrt{3}$,
$ -2$;\ \ 
$ \sqrt{3}$,
$ -\sqrt{3}$,
$0$,
$ -\frac{15+3\sqrt{5}}{2\sqrt{15}}$,
$ \frac{15+3\sqrt{5}}{2\sqrt{15}}$,
$0$;\ \ 
$ \sqrt{3}$,
$0$,
$ \frac{15+3\sqrt{5}}{2\sqrt{15}}$,
$ -\frac{15+3\sqrt{5}}{2\sqrt{15}}$,
$0$;\ \ 
$ -2$,
$0$,
$0$,
$ -1-\sqrt{5}$;\ \ 
$ -\sqrt{3}$,
$ \sqrt{3}$,
$0$;\ \ 
$ -\sqrt{3}$,
$0$;\ \ 
$ 2$)

Factors = $2_{\frac{26}{5},3.618}^{5,720}\boxtimes 5_{6,12.}^{24,592}$

  \vskip 2ex

\noindent32. $10_{0,52.}^{26,247}$ \irep{1035}:\ \ 
$d_i$ = ($1.0$,
$1.0$,
$2.0$,
$2.0$,
$2.0$,
$2.0$,
$2.0$,
$2.0$,
$3.605$,
$3.605$) 

\vskip 0.7ex
\hangindent=3em \hangafter=1
$D^2= 52.0 = 
52$

\vskip 0.7ex
\hangindent=3em \hangafter=1
$T = ( 0,
0,
\frac{1}{13},
\frac{3}{13},
\frac{4}{13},
\frac{9}{13},
\frac{10}{13},
\frac{12}{13},
0,
\frac{1}{2} )
$,

\vskip 0.7ex
\hangindent=3em \hangafter=1
$S$ = ($ 1$,
$ 1$,
$ 2$,
$ 2$,
$ 2$,
$ 2$,
$ 2$,
$ 2$,
$ \sqrt{13}$,
$ \sqrt{13}$;\ \ 
$ 1$,
$ 2$,
$ 2$,
$ 2$,
$ 2$,
$ 2$,
$ 2$,
$ -\sqrt{13}$,
$ -\sqrt{13}$;\ \ 
$ 2c_{13}^{2}$,
$ 2c_{13}^{5}$,
$ 2c_{13}^{4}$,
$ 2c_{13}^{6}$,
$ 2c_{13}^{1}$,
$ 2c_{13}^{3}$,
$0$,
$0$;\ \ 
$ 2c_{13}^{6}$,
$ 2c_{13}^{3}$,
$ 2c_{13}^{2}$,
$ 2c_{13}^{4}$,
$ 2c_{13}^{1}$,
$0$,
$0$;\ \ 
$ 2c_{13}^{5}$,
$ 2c_{13}^{1}$,
$ 2c_{13}^{2}$,
$ 2c_{13}^{6}$,
$0$,
$0$;\ \ 
$ 2c_{13}^{5}$,
$ 2c_{13}^{3}$,
$ 2c_{13}^{4}$,
$0$,
$0$;\ \ 
$ 2c_{13}^{6}$,
$ 2c_{13}^{5}$,
$0$,
$0$;\ \ 
$ 2c_{13}^{2}$,
$0$,
$0$;\ \ 
$ \sqrt{13}$,
$ -\sqrt{13}$;\ \ 
$ \sqrt{13}$)

  \vskip 2ex

\noindent33. $10_{4,52.}^{26,862}$ \irep{1035}:\ \ 
$d_i$ = ($1.0$,
$1.0$,
$2.0$,
$2.0$,
$2.0$,
$2.0$,
$2.0$,
$2.0$,
$3.605$,
$3.605$) 

\vskip 0.7ex
\hangindent=3em \hangafter=1
$D^2= 52.0 = 
52$

\vskip 0.7ex
\hangindent=3em \hangafter=1
$T = ( 0,
0,
\frac{2}{13},
\frac{5}{13},
\frac{6}{13},
\frac{7}{13},
\frac{8}{13},
\frac{11}{13},
0,
\frac{1}{2} )
$,

\vskip 0.7ex
\hangindent=3em \hangafter=1
$S$ = ($ 1$,
$ 1$,
$ 2$,
$ 2$,
$ 2$,
$ 2$,
$ 2$,
$ 2$,
$ \sqrt{13}$,
$ \sqrt{13}$;\ \ 
$ 1$,
$ 2$,
$ 2$,
$ 2$,
$ 2$,
$ 2$,
$ 2$,
$ -\sqrt{13}$,
$ -\sqrt{13}$;\ \ 
$ 2c_{13}^{4}$,
$ 2c_{13}^{1}$,
$ 2c_{13}^{3}$,
$ 2c_{13}^{2}$,
$ 2c_{13}^{5}$,
$ 2c_{13}^{6}$,
$0$,
$0$;\ \ 
$ 2c_{13}^{3}$,
$ 2c_{13}^{4}$,
$ 2c_{13}^{6}$,
$ 2c_{13}^{2}$,
$ 2c_{13}^{5}$,
$0$,
$0$;\ \ 
$ 2c_{13}^{1}$,
$ 2c_{13}^{5}$,
$ 2c_{13}^{6}$,
$ 2c_{13}^{2}$,
$0$,
$0$;\ \ 
$ 2c_{13}^{1}$,
$ 2c_{13}^{4}$,
$ 2c_{13}^{3}$,
$0$,
$0$;\ \ 
$ 2c_{13}^{3}$,
$ 2c_{13}^{1}$,
$0$,
$0$;\ \ 
$ 2c_{13}^{4}$,
$0$,
$0$;\ \ 
$ -\sqrt{13}$,
$ \sqrt{13}$;\ \ 
$ -\sqrt{13}$)

  \vskip 2ex

\noindent34. $10_{0,52.}^{52,110}$ \irep{1120}:\ \ 
$d_i$ = ($1.0$,
$1.0$,
$2.0$,
$2.0$,
$2.0$,
$2.0$,
$2.0$,
$2.0$,
$3.605$,
$3.605$) 

\vskip 0.7ex
\hangindent=3em \hangafter=1
$D^2= 52.0 = 
52$

\vskip 0.7ex
\hangindent=3em \hangafter=1
$T = ( 0,
0,
\frac{1}{13},
\frac{3}{13},
\frac{4}{13},
\frac{9}{13},
\frac{10}{13},
\frac{12}{13},
\frac{1}{4},
\frac{3}{4} )
$,

\vskip 0.7ex
\hangindent=3em \hangafter=1
$S$ = ($ 1$,
$ 1$,
$ 2$,
$ 2$,
$ 2$,
$ 2$,
$ 2$,
$ 2$,
$ \sqrt{13}$,
$ \sqrt{13}$;\ \ 
$ 1$,
$ 2$,
$ 2$,
$ 2$,
$ 2$,
$ 2$,
$ 2$,
$ -\sqrt{13}$,
$ -\sqrt{13}$;\ \ 
$ 2c_{13}^{2}$,
$ 2c_{13}^{5}$,
$ 2c_{13}^{4}$,
$ 2c_{13}^{6}$,
$ 2c_{13}^{1}$,
$ 2c_{13}^{3}$,
$0$,
$0$;\ \ 
$ 2c_{13}^{6}$,
$ 2c_{13}^{3}$,
$ 2c_{13}^{2}$,
$ 2c_{13}^{4}$,
$ 2c_{13}^{1}$,
$0$,
$0$;\ \ 
$ 2c_{13}^{5}$,
$ 2c_{13}^{1}$,
$ 2c_{13}^{2}$,
$ 2c_{13}^{6}$,
$0$,
$0$;\ \ 
$ 2c_{13}^{5}$,
$ 2c_{13}^{3}$,
$ 2c_{13}^{4}$,
$0$,
$0$;\ \ 
$ 2c_{13}^{6}$,
$ 2c_{13}^{5}$,
$0$,
$0$;\ \ 
$ 2c_{13}^{2}$,
$0$,
$0$;\ \ 
$ -\sqrt{13}$,
$ \sqrt{13}$;\ \ 
$ -\sqrt{13}$)

  \vskip 2ex

\noindent35. $10_{4,52.}^{52,489}$ \irep{1120}:\ \ 
$d_i$ = ($1.0$,
$1.0$,
$2.0$,
$2.0$,
$2.0$,
$2.0$,
$2.0$,
$2.0$,
$3.605$,
$3.605$) 

\vskip 0.7ex
\hangindent=3em \hangafter=1
$D^2= 52.0 = 
52$

\vskip 0.7ex
\hangindent=3em \hangafter=1
$T = ( 0,
0,
\frac{2}{13},
\frac{5}{13},
\frac{6}{13},
\frac{7}{13},
\frac{8}{13},
\frac{11}{13},
\frac{1}{4},
\frac{3}{4} )
$,

\vskip 0.7ex
\hangindent=3em \hangafter=1
$S$ = ($ 1$,
$ 1$,
$ 2$,
$ 2$,
$ 2$,
$ 2$,
$ 2$,
$ 2$,
$ \sqrt{13}$,
$ \sqrt{13}$;\ \ 
$ 1$,
$ 2$,
$ 2$,
$ 2$,
$ 2$,
$ 2$,
$ 2$,
$ -\sqrt{13}$,
$ -\sqrt{13}$;\ \ 
$ 2c_{13}^{4}$,
$ 2c_{13}^{1}$,
$ 2c_{13}^{3}$,
$ 2c_{13}^{2}$,
$ 2c_{13}^{5}$,
$ 2c_{13}^{6}$,
$0$,
$0$;\ \ 
$ 2c_{13}^{3}$,
$ 2c_{13}^{4}$,
$ 2c_{13}^{6}$,
$ 2c_{13}^{2}$,
$ 2c_{13}^{5}$,
$0$,
$0$;\ \ 
$ 2c_{13}^{1}$,
$ 2c_{13}^{5}$,
$ 2c_{13}^{6}$,
$ 2c_{13}^{2}$,
$0$,
$0$;\ \ 
$ 2c_{13}^{1}$,
$ 2c_{13}^{4}$,
$ 2c_{13}^{3}$,
$0$,
$0$;\ \ 
$ 2c_{13}^{3}$,
$ 2c_{13}^{1}$,
$0$,
$0$;\ \ 
$ 2c_{13}^{4}$,
$0$,
$0$;\ \ 
$ \sqrt{13}$,
$ -\sqrt{13}$;\ \ 
$ \sqrt{13}$)

  \vskip 2ex

\noindent36. $10_{\frac{83}{11},69.29}^{44,134}$ \irep{1106}:\ \ 
$d_i$ = ($1.0$,
$1.0$,
$1.918$,
$1.918$,
$2.682$,
$2.682$,
$3.228$,
$3.228$,
$3.513$,
$3.513$) 

\vskip 0.7ex
\hangindent=3em \hangafter=1
$D^2= 69.292 = 
30+20c^{1}_{11}
+12c^{2}_{11}
+6c^{3}_{11}
+2c^{4}_{11}
$

\vskip 0.7ex
\hangindent=3em \hangafter=1
$T = ( 0,
\frac{1}{4},
\frac{2}{11},
\frac{19}{44},
\frac{9}{11},
\frac{3}{44},
\frac{10}{11},
\frac{7}{44},
\frac{5}{11},
\frac{31}{44} )
$,

\vskip 0.7ex
\hangindent=3em \hangafter=1
$S$ = ($ 1$,
$ 1$,
$ -c_{11}^{5}$,
$ -c_{11}^{5}$,
$ \xi_{11}^{3}$,
$ \xi_{11}^{3}$,
$ \xi_{11}^{7}$,
$ \xi_{11}^{7}$,
$ \xi_{11}^{5}$,
$ \xi_{11}^{5}$;\ \ 
$ -1$,
$ -c_{11}^{5}$,
$ c_{11}^{5}$,
$ \xi_{11}^{3}$,
$ -\xi_{11}^{3}$,
$ \xi_{11}^{7}$,
$ -\xi_{11}^{7}$,
$ \xi_{11}^{5}$,
$ -\xi_{11}^{5}$;\ \ 
$ -\xi_{11}^{7}$,
$ -\xi_{11}^{7}$,
$ \xi_{11}^{5}$,
$ \xi_{11}^{5}$,
$ -\xi_{11}^{3}$,
$ -\xi_{11}^{3}$,
$ 1$,
$ 1$;\ \ 
$ \xi_{11}^{7}$,
$ \xi_{11}^{5}$,
$ -\xi_{11}^{5}$,
$ -\xi_{11}^{3}$,
$ \xi_{11}^{3}$,
$ 1$,
$ -1$;\ \ 
$ -c_{11}^{5}$,
$ -c_{11}^{5}$,
$ -1$,
$ -1$,
$ -\xi_{11}^{7}$,
$ -\xi_{11}^{7}$;\ \ 
$ c_{11}^{5}$,
$ -1$,
$ 1$,
$ -\xi_{11}^{7}$,
$ \xi_{11}^{7}$;\ \ 
$ \xi_{11}^{5}$,
$ \xi_{11}^{5}$,
$ c_{11}^{5}$,
$ c_{11}^{5}$;\ \ 
$ -\xi_{11}^{5}$,
$ c_{11}^{5}$,
$ -c_{11}^{5}$;\ \ 
$ \xi_{11}^{3}$,
$ \xi_{11}^{3}$;\ \ 
$ -\xi_{11}^{3}$)

Factors = $2_{1,2.}^{4,437}\boxtimes 5_{\frac{72}{11},34.64}^{11,216}$

  \vskip 2ex

\noindent37. $10_{\frac{27}{11},69.29}^{44,556}$ \irep{1106}:\ \ 
$d_i$ = ($1.0$,
$1.0$,
$1.918$,
$1.918$,
$2.682$,
$2.682$,
$3.228$,
$3.228$,
$3.513$,
$3.513$) 

\vskip 0.7ex
\hangindent=3em \hangafter=1
$D^2= 69.292 = 
30+20c^{1}_{11}
+12c^{2}_{11}
+6c^{3}_{11}
+2c^{4}_{11}
$

\vskip 0.7ex
\hangindent=3em \hangafter=1
$T = ( 0,
\frac{1}{4},
\frac{9}{11},
\frac{3}{44},
\frac{2}{11},
\frac{19}{44},
\frac{1}{11},
\frac{15}{44},
\frac{6}{11},
\frac{35}{44} )
$,

\vskip 0.7ex
\hangindent=3em \hangafter=1
$S$ = ($ 1$,
$ 1$,
$ -c_{11}^{5}$,
$ -c_{11}^{5}$,
$ \xi_{11}^{3}$,
$ \xi_{11}^{3}$,
$ \xi_{11}^{7}$,
$ \xi_{11}^{7}$,
$ \xi_{11}^{5}$,
$ \xi_{11}^{5}$;\ \ 
$ -1$,
$ -c_{11}^{5}$,
$ c_{11}^{5}$,
$ \xi_{11}^{3}$,
$ -\xi_{11}^{3}$,
$ \xi_{11}^{7}$,
$ -\xi_{11}^{7}$,
$ \xi_{11}^{5}$,
$ -\xi_{11}^{5}$;\ \ 
$ -\xi_{11}^{7}$,
$ -\xi_{11}^{7}$,
$ \xi_{11}^{5}$,
$ \xi_{11}^{5}$,
$ -\xi_{11}^{3}$,
$ -\xi_{11}^{3}$,
$ 1$,
$ 1$;\ \ 
$ \xi_{11}^{7}$,
$ \xi_{11}^{5}$,
$ -\xi_{11}^{5}$,
$ -\xi_{11}^{3}$,
$ \xi_{11}^{3}$,
$ 1$,
$ -1$;\ \ 
$ -c_{11}^{5}$,
$ -c_{11}^{5}$,
$ -1$,
$ -1$,
$ -\xi_{11}^{7}$,
$ -\xi_{11}^{7}$;\ \ 
$ c_{11}^{5}$,
$ -1$,
$ 1$,
$ -\xi_{11}^{7}$,
$ \xi_{11}^{7}$;\ \ 
$ \xi_{11}^{5}$,
$ \xi_{11}^{5}$,
$ c_{11}^{5}$,
$ c_{11}^{5}$;\ \ 
$ -\xi_{11}^{5}$,
$ c_{11}^{5}$,
$ -c_{11}^{5}$;\ \ 
$ \xi_{11}^{3}$,
$ \xi_{11}^{3}$;\ \ 
$ -\xi_{11}^{3}$)

Factors = $2_{1,2.}^{4,437}\boxtimes 5_{\frac{16}{11},34.64}^{11,640}$

  \vskip 2ex

\noindent38. $10_{\frac{61}{11},69.29}^{44,372}$ \irep{1106}:\ \ 
$d_i$ = ($1.0$,
$1.0$,
$1.918$,
$1.918$,
$2.682$,
$2.682$,
$3.228$,
$3.228$,
$3.513$,
$3.513$) 

\vskip 0.7ex
\hangindent=3em \hangafter=1
$D^2= 69.292 = 
30+20c^{1}_{11}
+12c^{2}_{11}
+6c^{3}_{11}
+2c^{4}_{11}
$

\vskip 0.7ex
\hangindent=3em \hangafter=1
$T = ( 0,
\frac{3}{4},
\frac{2}{11},
\frac{41}{44},
\frac{9}{11},
\frac{25}{44},
\frac{10}{11},
\frac{29}{44},
\frac{5}{11},
\frac{9}{44} )
$,

\vskip 0.7ex
\hangindent=3em \hangafter=1
$S$ = ($ 1$,
$ 1$,
$ -c_{11}^{5}$,
$ -c_{11}^{5}$,
$ \xi_{11}^{3}$,
$ \xi_{11}^{3}$,
$ \xi_{11}^{7}$,
$ \xi_{11}^{7}$,
$ \xi_{11}^{5}$,
$ \xi_{11}^{5}$;\ \ 
$ -1$,
$ -c_{11}^{5}$,
$ c_{11}^{5}$,
$ \xi_{11}^{3}$,
$ -\xi_{11}^{3}$,
$ \xi_{11}^{7}$,
$ -\xi_{11}^{7}$,
$ \xi_{11}^{5}$,
$ -\xi_{11}^{5}$;\ \ 
$ -\xi_{11}^{7}$,
$ -\xi_{11}^{7}$,
$ \xi_{11}^{5}$,
$ \xi_{11}^{5}$,
$ -\xi_{11}^{3}$,
$ -\xi_{11}^{3}$,
$ 1$,
$ 1$;\ \ 
$ \xi_{11}^{7}$,
$ \xi_{11}^{5}$,
$ -\xi_{11}^{5}$,
$ -\xi_{11}^{3}$,
$ \xi_{11}^{3}$,
$ 1$,
$ -1$;\ \ 
$ -c_{11}^{5}$,
$ -c_{11}^{5}$,
$ -1$,
$ -1$,
$ -\xi_{11}^{7}$,
$ -\xi_{11}^{7}$;\ \ 
$ c_{11}^{5}$,
$ -1$,
$ 1$,
$ -\xi_{11}^{7}$,
$ \xi_{11}^{7}$;\ \ 
$ \xi_{11}^{5}$,
$ \xi_{11}^{5}$,
$ c_{11}^{5}$,
$ c_{11}^{5}$;\ \ 
$ -\xi_{11}^{5}$,
$ c_{11}^{5}$,
$ -c_{11}^{5}$;\ \ 
$ \xi_{11}^{3}$,
$ \xi_{11}^{3}$;\ \ 
$ -\xi_{11}^{3}$)

Factors = $2_{7,2.}^{4,625}\boxtimes 5_{\frac{72}{11},34.64}^{11,216}$

  \vskip 2ex

\noindent39. $10_{\frac{5}{11},69.29}^{44,237}$ \irep{1106}:\ \ 
$d_i$ = ($1.0$,
$1.0$,
$1.918$,
$1.918$,
$2.682$,
$2.682$,
$3.228$,
$3.228$,
$3.513$,
$3.513$) 

\vskip 0.7ex
\hangindent=3em \hangafter=1
$D^2= 69.292 = 
30+20c^{1}_{11}
+12c^{2}_{11}
+6c^{3}_{11}
+2c^{4}_{11}
$

\vskip 0.7ex
\hangindent=3em \hangafter=1
$T = ( 0,
\frac{3}{4},
\frac{9}{11},
\frac{25}{44},
\frac{2}{11},
\frac{41}{44},
\frac{1}{11},
\frac{37}{44},
\frac{6}{11},
\frac{13}{44} )
$,

\vskip 0.7ex
\hangindent=3em \hangafter=1
$S$ = ($ 1$,
$ 1$,
$ -c_{11}^{5}$,
$ -c_{11}^{5}$,
$ \xi_{11}^{3}$,
$ \xi_{11}^{3}$,
$ \xi_{11}^{7}$,
$ \xi_{11}^{7}$,
$ \xi_{11}^{5}$,
$ \xi_{11}^{5}$;\ \ 
$ -1$,
$ -c_{11}^{5}$,
$ c_{11}^{5}$,
$ \xi_{11}^{3}$,
$ -\xi_{11}^{3}$,
$ \xi_{11}^{7}$,
$ -\xi_{11}^{7}$,
$ \xi_{11}^{5}$,
$ -\xi_{11}^{5}$;\ \ 
$ -\xi_{11}^{7}$,
$ -\xi_{11}^{7}$,
$ \xi_{11}^{5}$,
$ \xi_{11}^{5}$,
$ -\xi_{11}^{3}$,
$ -\xi_{11}^{3}$,
$ 1$,
$ 1$;\ \ 
$ \xi_{11}^{7}$,
$ \xi_{11}^{5}$,
$ -\xi_{11}^{5}$,
$ -\xi_{11}^{3}$,
$ \xi_{11}^{3}$,
$ 1$,
$ -1$;\ \ 
$ -c_{11}^{5}$,
$ -c_{11}^{5}$,
$ -1$,
$ -1$,
$ -\xi_{11}^{7}$,
$ -\xi_{11}^{7}$;\ \ 
$ c_{11}^{5}$,
$ -1$,
$ 1$,
$ -\xi_{11}^{7}$,
$ \xi_{11}^{7}$;\ \ 
$ \xi_{11}^{5}$,
$ \xi_{11}^{5}$,
$ c_{11}^{5}$,
$ c_{11}^{5}$;\ \ 
$ -\xi_{11}^{5}$,
$ c_{11}^{5}$,
$ -c_{11}^{5}$;\ \ 
$ \xi_{11}^{3}$,
$ \xi_{11}^{3}$;\ \ 
$ -\xi_{11}^{3}$)

Factors = $2_{7,2.}^{4,625}\boxtimes 5_{\frac{16}{11},34.64}^{11,640}$

  \vskip 2ex

\noindent40. $10_{\frac{45}{7},70.68}^{28,820}$ \irep{1038}:\ \ 
$d_i$ = ($1.0$,
$1.0$,
$2.246$,
$2.246$,
$2.246$,
$2.246$,
$2.801$,
$2.801$,
$4.48$,
$4.48$) 

\vskip 0.7ex
\hangindent=3em \hangafter=1
$D^2= 70.684 = 
42+28c^{1}_{7}
+14c^{2}_{7}
$

\vskip 0.7ex
\hangindent=3em \hangafter=1
$T = ( 0,
\frac{1}{4},
\frac{1}{7},
\frac{1}{7},
\frac{11}{28},
\frac{11}{28},
\frac{6}{7},
\frac{3}{28},
\frac{4}{7},
\frac{23}{28} )
$,

\vskip 0.7ex
\hangindent=3em \hangafter=1
$S$ = ($ 1$,
$ 1$,
$ \xi_{7}^{3}$,
$ \xi_{7}^{3}$,
$ \xi_{7}^{3}$,
$ \xi_{7}^{3}$,
$ 2+c^{1}_{7}
+c^{2}_{7}
$,
$ 2+c^{1}_{7}
+c^{2}_{7}
$,
$ 2+2c^{1}_{7}
+c^{2}_{7}
$,
$ 2+2c^{1}_{7}
+c^{2}_{7}
$;\ \ 
$ -1$,
$ \xi_{7}^{3}$,
$ \xi_{7}^{3}$,
$ -\xi_{7}^{3}$,
$ -\xi_{7}^{3}$,
$ 2+c^{1}_{7}
+c^{2}_{7}
$,
$ -2-c^{1}_{7}
-c^{2}_{7}
$,
$ 2+2c^{1}_{7}
+c^{2}_{7}
$,
$ -2-2  c^{1}_{7}
-c^{2}_{7}
$;\ \ 
$ s^{1}_{7}
+\zeta^{2}_{7}
+\zeta^{3}_{7}
$,
$ -1-2  \zeta^{1}_{7}
-\zeta^{2}_{7}
-\zeta^{3}_{7}
$,
$ -1-2  \zeta^{1}_{7}
-\zeta^{2}_{7}
-\zeta^{3}_{7}
$,
$ s^{1}_{7}
+\zeta^{2}_{7}
+\zeta^{3}_{7}
$,
$ -\xi_{7}^{3}$,
$ -\xi_{7}^{3}$,
$ \xi_{7}^{3}$,
$ \xi_{7}^{3}$;\ \ 
$ s^{1}_{7}
+\zeta^{2}_{7}
+\zeta^{3}_{7}
$,
$ s^{1}_{7}
+\zeta^{2}_{7}
+\zeta^{3}_{7}
$,
$ -1-2  \zeta^{1}_{7}
-\zeta^{2}_{7}
-\zeta^{3}_{7}
$,
$ -\xi_{7}^{3}$,
$ -\xi_{7}^{3}$,
$ \xi_{7}^{3}$,
$ \xi_{7}^{3}$;\ \ 
$ -s^{1}_{7}
-\zeta^{2}_{7}
-\zeta^{3}_{7}
$,
$ 1+2\zeta^{1}_{7}
+\zeta^{2}_{7}
+\zeta^{3}_{7}
$,
$ -\xi_{7}^{3}$,
$ \xi_{7}^{3}$,
$ \xi_{7}^{3}$,
$ -\xi_{7}^{3}$;\ \ 
$ -s^{1}_{7}
-\zeta^{2}_{7}
-\zeta^{3}_{7}
$,
$ -\xi_{7}^{3}$,
$ \xi_{7}^{3}$,
$ \xi_{7}^{3}$,
$ -\xi_{7}^{3}$;\ \ 
$ 2+2c^{1}_{7}
+c^{2}_{7}
$,
$ 2+2c^{1}_{7}
+c^{2}_{7}
$,
$ -1$,
$ -1$;\ \ 
$ -2-2  c^{1}_{7}
-c^{2}_{7}
$,
$ -1$,
$ 1$;\ \ 
$ -2-c^{1}_{7}
-c^{2}_{7}
$,
$ -2-c^{1}_{7}
-c^{2}_{7}
$;\ \ 
$ 2+c^{1}_{7}
+c^{2}_{7}
$)

Factors = $2_{1,2.}^{4,437}\boxtimes 5_{\frac{38}{7},35.34}^{7,386}$

  \vskip 2ex

\noindent41. $10_{\frac{25}{7},70.68}^{28,341}$ \irep{1038}:\ \ 
$d_i$ = ($1.0$,
$1.0$,
$2.246$,
$2.246$,
$2.246$,
$2.246$,
$2.801$,
$2.801$,
$4.48$,
$4.48$) 

\vskip 0.7ex
\hangindent=3em \hangafter=1
$D^2= 70.684 = 
42+28c^{1}_{7}
+14c^{2}_{7}
$

\vskip 0.7ex
\hangindent=3em \hangafter=1
$T = ( 0,
\frac{1}{4},
\frac{6}{7},
\frac{6}{7},
\frac{3}{28},
\frac{3}{28},
\frac{1}{7},
\frac{11}{28},
\frac{3}{7},
\frac{19}{28} )
$,

\vskip 0.7ex
\hangindent=3em \hangafter=1
$S$ = ($ 1$,
$ 1$,
$ \xi_{7}^{3}$,
$ \xi_{7}^{3}$,
$ \xi_{7}^{3}$,
$ \xi_{7}^{3}$,
$ 2+c^{1}_{7}
+c^{2}_{7}
$,
$ 2+c^{1}_{7}
+c^{2}_{7}
$,
$ 2+2c^{1}_{7}
+c^{2}_{7}
$,
$ 2+2c^{1}_{7}
+c^{2}_{7}
$;\ \ 
$ -1$,
$ \xi_{7}^{3}$,
$ \xi_{7}^{3}$,
$ -\xi_{7}^{3}$,
$ -\xi_{7}^{3}$,
$ 2+c^{1}_{7}
+c^{2}_{7}
$,
$ -2-c^{1}_{7}
-c^{2}_{7}
$,
$ 2+2c^{1}_{7}
+c^{2}_{7}
$,
$ -2-2  c^{1}_{7}
-c^{2}_{7}
$;\ \ 
$ -1-2  \zeta^{1}_{7}
-\zeta^{2}_{7}
-\zeta^{3}_{7}
$,
$ s^{1}_{7}
+\zeta^{2}_{7}
+\zeta^{3}_{7}
$,
$ -1-2  \zeta^{1}_{7}
-\zeta^{2}_{7}
-\zeta^{3}_{7}
$,
$ s^{1}_{7}
+\zeta^{2}_{7}
+\zeta^{3}_{7}
$,
$ -\xi_{7}^{3}$,
$ -\xi_{7}^{3}$,
$ \xi_{7}^{3}$,
$ \xi_{7}^{3}$;\ \ 
$ -1-2  \zeta^{1}_{7}
-\zeta^{2}_{7}
-\zeta^{3}_{7}
$,
$ s^{1}_{7}
+\zeta^{2}_{7}
+\zeta^{3}_{7}
$,
$ -1-2  \zeta^{1}_{7}
-\zeta^{2}_{7}
-\zeta^{3}_{7}
$,
$ -\xi_{7}^{3}$,
$ -\xi_{7}^{3}$,
$ \xi_{7}^{3}$,
$ \xi_{7}^{3}$;\ \ 
$ 1+2\zeta^{1}_{7}
+\zeta^{2}_{7}
+\zeta^{3}_{7}
$,
$ -s^{1}_{7}
-\zeta^{2}_{7}
-\zeta^{3}_{7}
$,
$ -\xi_{7}^{3}$,
$ \xi_{7}^{3}$,
$ \xi_{7}^{3}$,
$ -\xi_{7}^{3}$;\ \ 
$ 1+2\zeta^{1}_{7}
+\zeta^{2}_{7}
+\zeta^{3}_{7}
$,
$ -\xi_{7}^{3}$,
$ \xi_{7}^{3}$,
$ \xi_{7}^{3}$,
$ -\xi_{7}^{3}$;\ \ 
$ 2+2c^{1}_{7}
+c^{2}_{7}
$,
$ 2+2c^{1}_{7}
+c^{2}_{7}
$,
$ -1$,
$ -1$;\ \ 
$ -2-2  c^{1}_{7}
-c^{2}_{7}
$,
$ -1$,
$ 1$;\ \ 
$ -2-c^{1}_{7}
-c^{2}_{7}
$,
$ -2-c^{1}_{7}
-c^{2}_{7}
$;\ \ 
$ 2+c^{1}_{7}
+c^{2}_{7}
$)

Factors = $2_{1,2.}^{4,437}\boxtimes 5_{\frac{18}{7},35.34}^{7,101}$

  \vskip 2ex

\noindent42. $10_{\frac{31}{7},70.68}^{28,289}$ \irep{1038}:\ \ 
$d_i$ = ($1.0$,
$1.0$,
$2.246$,
$2.246$,
$2.246$,
$2.246$,
$2.801$,
$2.801$,
$4.48$,
$4.48$) 

\vskip 0.7ex
\hangindent=3em \hangafter=1
$D^2= 70.684 = 
42+28c^{1}_{7}
+14c^{2}_{7}
$

\vskip 0.7ex
\hangindent=3em \hangafter=1
$T = ( 0,
\frac{3}{4},
\frac{1}{7},
\frac{1}{7},
\frac{25}{28},
\frac{25}{28},
\frac{6}{7},
\frac{17}{28},
\frac{4}{7},
\frac{9}{28} )
$,

\vskip 0.7ex
\hangindent=3em \hangafter=1
$S$ = ($ 1$,
$ 1$,
$ \xi_{7}^{3}$,
$ \xi_{7}^{3}$,
$ \xi_{7}^{3}$,
$ \xi_{7}^{3}$,
$ 2+c^{1}_{7}
+c^{2}_{7}
$,
$ 2+c^{1}_{7}
+c^{2}_{7}
$,
$ 2+2c^{1}_{7}
+c^{2}_{7}
$,
$ 2+2c^{1}_{7}
+c^{2}_{7}
$;\ \ 
$ -1$,
$ \xi_{7}^{3}$,
$ \xi_{7}^{3}$,
$ -\xi_{7}^{3}$,
$ -\xi_{7}^{3}$,
$ 2+c^{1}_{7}
+c^{2}_{7}
$,
$ -2-c^{1}_{7}
-c^{2}_{7}
$,
$ 2+2c^{1}_{7}
+c^{2}_{7}
$,
$ -2-2  c^{1}_{7}
-c^{2}_{7}
$;\ \ 
$ s^{1}_{7}
+\zeta^{2}_{7}
+\zeta^{3}_{7}
$,
$ -1-2  \zeta^{1}_{7}
-\zeta^{2}_{7}
-\zeta^{3}_{7}
$,
$ -1-2  \zeta^{1}_{7}
-\zeta^{2}_{7}
-\zeta^{3}_{7}
$,
$ s^{1}_{7}
+\zeta^{2}_{7}
+\zeta^{3}_{7}
$,
$ -\xi_{7}^{3}$,
$ -\xi_{7}^{3}$,
$ \xi_{7}^{3}$,
$ \xi_{7}^{3}$;\ \ 
$ s^{1}_{7}
+\zeta^{2}_{7}
+\zeta^{3}_{7}
$,
$ s^{1}_{7}
+\zeta^{2}_{7}
+\zeta^{3}_{7}
$,
$ -1-2  \zeta^{1}_{7}
-\zeta^{2}_{7}
-\zeta^{3}_{7}
$,
$ -\xi_{7}^{3}$,
$ -\xi_{7}^{3}$,
$ \xi_{7}^{3}$,
$ \xi_{7}^{3}$;\ \ 
$ -s^{1}_{7}
-\zeta^{2}_{7}
-\zeta^{3}_{7}
$,
$ 1+2\zeta^{1}_{7}
+\zeta^{2}_{7}
+\zeta^{3}_{7}
$,
$ -\xi_{7}^{3}$,
$ \xi_{7}^{3}$,
$ \xi_{7}^{3}$,
$ -\xi_{7}^{3}$;\ \ 
$ -s^{1}_{7}
-\zeta^{2}_{7}
-\zeta^{3}_{7}
$,
$ -\xi_{7}^{3}$,
$ \xi_{7}^{3}$,
$ \xi_{7}^{3}$,
$ -\xi_{7}^{3}$;\ \ 
$ 2+2c^{1}_{7}
+c^{2}_{7}
$,
$ 2+2c^{1}_{7}
+c^{2}_{7}
$,
$ -1$,
$ -1$;\ \ 
$ -2-2  c^{1}_{7}
-c^{2}_{7}
$,
$ -1$,
$ 1$;\ \ 
$ -2-c^{1}_{7}
-c^{2}_{7}
$,
$ -2-c^{1}_{7}
-c^{2}_{7}
$;\ \ 
$ 2+c^{1}_{7}
+c^{2}_{7}
$)

Factors = $2_{7,2.}^{4,625}\boxtimes 5_{\frac{38}{7},35.34}^{7,386}$

  \vskip 2ex

\noindent43. $10_{\frac{11}{7},70.68}^{28,251}$ \irep{1038}:\ \ 
$d_i$ = ($1.0$,
$1.0$,
$2.246$,
$2.246$,
$2.246$,
$2.246$,
$2.801$,
$2.801$,
$4.48$,
$4.48$) 

\vskip 0.7ex
\hangindent=3em \hangafter=1
$D^2= 70.684 = 
42+28c^{1}_{7}
+14c^{2}_{7}
$

\vskip 0.7ex
\hangindent=3em \hangafter=1
$T = ( 0,
\frac{3}{4},
\frac{6}{7},
\frac{6}{7},
\frac{17}{28},
\frac{17}{28},
\frac{1}{7},
\frac{25}{28},
\frac{3}{7},
\frac{5}{28} )
$,

\vskip 0.7ex
\hangindent=3em \hangafter=1
$S$ = ($ 1$,
$ 1$,
$ \xi_{7}^{3}$,
$ \xi_{7}^{3}$,
$ \xi_{7}^{3}$,
$ \xi_{7}^{3}$,
$ 2+c^{1}_{7}
+c^{2}_{7}
$,
$ 2+c^{1}_{7}
+c^{2}_{7}
$,
$ 2+2c^{1}_{7}
+c^{2}_{7}
$,
$ 2+2c^{1}_{7}
+c^{2}_{7}
$;\ \ 
$ -1$,
$ \xi_{7}^{3}$,
$ \xi_{7}^{3}$,
$ -\xi_{7}^{3}$,
$ -\xi_{7}^{3}$,
$ 2+c^{1}_{7}
+c^{2}_{7}
$,
$ -2-c^{1}_{7}
-c^{2}_{7}
$,
$ 2+2c^{1}_{7}
+c^{2}_{7}
$,
$ -2-2  c^{1}_{7}
-c^{2}_{7}
$;\ \ 
$ -1-2  \zeta^{1}_{7}
-\zeta^{2}_{7}
-\zeta^{3}_{7}
$,
$ s^{1}_{7}
+\zeta^{2}_{7}
+\zeta^{3}_{7}
$,
$ -1-2  \zeta^{1}_{7}
-\zeta^{2}_{7}
-\zeta^{3}_{7}
$,
$ s^{1}_{7}
+\zeta^{2}_{7}
+\zeta^{3}_{7}
$,
$ -\xi_{7}^{3}$,
$ -\xi_{7}^{3}$,
$ \xi_{7}^{3}$,
$ \xi_{7}^{3}$;\ \ 
$ -1-2  \zeta^{1}_{7}
-\zeta^{2}_{7}
-\zeta^{3}_{7}
$,
$ s^{1}_{7}
+\zeta^{2}_{7}
+\zeta^{3}_{7}
$,
$ -1-2  \zeta^{1}_{7}
-\zeta^{2}_{7}
-\zeta^{3}_{7}
$,
$ -\xi_{7}^{3}$,
$ -\xi_{7}^{3}$,
$ \xi_{7}^{3}$,
$ \xi_{7}^{3}$;\ \ 
$ 1+2\zeta^{1}_{7}
+\zeta^{2}_{7}
+\zeta^{3}_{7}
$,
$ -s^{1}_{7}
-\zeta^{2}_{7}
-\zeta^{3}_{7}
$,
$ -\xi_{7}^{3}$,
$ \xi_{7}^{3}$,
$ \xi_{7}^{3}$,
$ -\xi_{7}^{3}$;\ \ 
$ 1+2\zeta^{1}_{7}
+\zeta^{2}_{7}
+\zeta^{3}_{7}
$,
$ -\xi_{7}^{3}$,
$ \xi_{7}^{3}$,
$ \xi_{7}^{3}$,
$ -\xi_{7}^{3}$;\ \ 
$ 2+2c^{1}_{7}
+c^{2}_{7}
$,
$ 2+2c^{1}_{7}
+c^{2}_{7}
$,
$ -1$,
$ -1$;\ \ 
$ -2-2  c^{1}_{7}
-c^{2}_{7}
$,
$ -1$,
$ 1$;\ \ 
$ -2-c^{1}_{7}
-c^{2}_{7}
$,
$ -2-c^{1}_{7}
-c^{2}_{7}
$;\ \ 
$ 2+c^{1}_{7}
+c^{2}_{7}
$)

Factors = $2_{7,2.}^{4,625}\boxtimes 5_{\frac{18}{7},35.34}^{7,101}$

  \vskip 2ex

\noindent44. $10_{6,89.56}^{12,311}$ \irep{680}:\ \ 
$d_i$ = ($1.0$,
$1.0$,
$2.732$,
$2.732$,
$2.732$,
$2.732$,
$2.732$,
$3.732$,
$3.732$,
$4.732$) 

\vskip 0.7ex
\hangindent=3em \hangafter=1
$D^2= 89.569 = 
48+24\sqrt{3}$

\vskip 0.7ex
\hangindent=3em \hangafter=1
$T = ( 0,
\frac{1}{2},
0,
\frac{1}{3},
\frac{1}{3},
\frac{1}{3},
\frac{5}{6},
0,
\frac{1}{2},
\frac{3}{4} )
$,

\vskip 0.7ex
\hangindent=3em \hangafter=1
$S$ = ($ 1$,
$ 1$,
$ 1+\sqrt{3}$,
$ 1+\sqrt{3}$,
$ 1+\sqrt{3}$,
$ 1+\sqrt{3}$,
$ 1+\sqrt{3}$,
$ 2+\sqrt{3}$,
$ 2+\sqrt{3}$,
$ 3+\sqrt{3}$;\ \ 
$ 1$,
$ -1-\sqrt{3}$,
$ 1+\sqrt{3}$,
$ -1-\sqrt{3}$,
$ -1-\sqrt{3}$,
$ 1+\sqrt{3}$,
$ 2+\sqrt{3}$,
$ 2+\sqrt{3}$,
$ -3-\sqrt{3}$;\ \ 
$0$,
$ -2-2\sqrt{3}$,
$0$,
$0$,
$ 2+2\sqrt{3}$,
$ -1-\sqrt{3}$,
$ 1+\sqrt{3}$,
$0$;\ \ 
$ 1+\sqrt{3}$,
$ 1+\sqrt{3}$,
$ 1+\sqrt{3}$,
$ 1+\sqrt{3}$,
$ -1-\sqrt{3}$,
$ -1-\sqrt{3}$,
$0$;\ \ 
$(-3-\sqrt{3})\mathrm{i}$,
$(3+\sqrt{3})\mathrm{i}$,
$ -1-\sqrt{3}$,
$ -1-\sqrt{3}$,
$ 1+\sqrt{3}$,
$0$;\ \ 
$(-3-\sqrt{3})\mathrm{i}$,
$ -1-\sqrt{3}$,
$ -1-\sqrt{3}$,
$ 1+\sqrt{3}$,
$0$;\ \ 
$ 1+\sqrt{3}$,
$ -1-\sqrt{3}$,
$ -1-\sqrt{3}$,
$0$;\ \ 
$ 1$,
$ 1$,
$ 3+\sqrt{3}$;\ \ 
$ 1$,
$ -3-\sqrt{3}$;\ \ 
$0$)

  \vskip 2ex

\noindent45. $10_{2,89.56}^{12,176}$ \irep{680}:\ \ 
$d_i$ = ($1.0$,
$1.0$,
$2.732$,
$2.732$,
$2.732$,
$2.732$,
$2.732$,
$3.732$,
$3.732$,
$4.732$) 

\vskip 0.7ex
\hangindent=3em \hangafter=1
$D^2= 89.569 = 
48+24\sqrt{3}$

\vskip 0.7ex
\hangindent=3em \hangafter=1
$T = ( 0,
\frac{1}{2},
0,
\frac{2}{3},
\frac{2}{3},
\frac{2}{3},
\frac{1}{6},
0,
\frac{1}{2},
\frac{1}{4} )
$,

\vskip 0.7ex
\hangindent=3em \hangafter=1
$S$ = ($ 1$,
$ 1$,
$ 1+\sqrt{3}$,
$ 1+\sqrt{3}$,
$ 1+\sqrt{3}$,
$ 1+\sqrt{3}$,
$ 1+\sqrt{3}$,
$ 2+\sqrt{3}$,
$ 2+\sqrt{3}$,
$ 3+\sqrt{3}$;\ \ 
$ 1$,
$ -1-\sqrt{3}$,
$ 1+\sqrt{3}$,
$ -1-\sqrt{3}$,
$ -1-\sqrt{3}$,
$ 1+\sqrt{3}$,
$ 2+\sqrt{3}$,
$ 2+\sqrt{3}$,
$ -3-\sqrt{3}$;\ \ 
$0$,
$ -2-2\sqrt{3}$,
$0$,
$0$,
$ 2+2\sqrt{3}$,
$ -1-\sqrt{3}$,
$ 1+\sqrt{3}$,
$0$;\ \ 
$ 1+\sqrt{3}$,
$ 1+\sqrt{3}$,
$ 1+\sqrt{3}$,
$ 1+\sqrt{3}$,
$ -1-\sqrt{3}$,
$ -1-\sqrt{3}$,
$0$;\ \ 
$(3+\sqrt{3})\mathrm{i}$,
$(-3-\sqrt{3})\mathrm{i}$,
$ -1-\sqrt{3}$,
$ -1-\sqrt{3}$,
$ 1+\sqrt{3}$,
$0$;\ \ 
$(3+\sqrt{3})\mathrm{i}$,
$ -1-\sqrt{3}$,
$ -1-\sqrt{3}$,
$ 1+\sqrt{3}$,
$0$;\ \ 
$ 1+\sqrt{3}$,
$ -1-\sqrt{3}$,
$ -1-\sqrt{3}$,
$0$;\ \ 
$ 1$,
$ 1$,
$ 3+\sqrt{3}$;\ \ 
$ 1$,
$ -3-\sqrt{3}$;\ \ 
$0$)

  \vskip 2ex

\noindent46. $10_{0,89.56}^{12,155}$ \irep{587}:\ \ 
$d_i$ = ($1.0$,
$1.0$,
$2.732$,
$2.732$,
$2.732$,
$2.732$,
$2.732$,
$3.732$,
$3.732$,
$4.732$) 

\vskip 0.7ex
\hangindent=3em \hangafter=1
$D^2= 89.569 = 
48+24\sqrt{3}$

\vskip 0.7ex
\hangindent=3em \hangafter=1
$T = ( 0,
\frac{1}{2},
\frac{1}{3},
\frac{1}{4},
\frac{5}{6},
\frac{7}{12},
\frac{7}{12},
0,
\frac{1}{2},
0 )
$,

\vskip 0.7ex
\hangindent=3em \hangafter=1
$S$ = ($ 1$,
$ 1$,
$ 1+\sqrt{3}$,
$ 1+\sqrt{3}$,
$ 1+\sqrt{3}$,
$ 1+\sqrt{3}$,
$ 1+\sqrt{3}$,
$ 2+\sqrt{3}$,
$ 2+\sqrt{3}$,
$ 3+\sqrt{3}$;\ \ 
$ 1$,
$ 1+\sqrt{3}$,
$ -1-\sqrt{3}$,
$ 1+\sqrt{3}$,
$ -1-\sqrt{3}$,
$ -1-\sqrt{3}$,
$ 2+\sqrt{3}$,
$ 2+\sqrt{3}$,
$ -3-\sqrt{3}$;\ \ 
$ 1+\sqrt{3}$,
$ -2-2\sqrt{3}$,
$ 1+\sqrt{3}$,
$ 1+\sqrt{3}$,
$ 1+\sqrt{3}$,
$ -1-\sqrt{3}$,
$ -1-\sqrt{3}$,
$0$;\ \ 
$0$,
$ 2+2\sqrt{3}$,
$0$,
$0$,
$ -1-\sqrt{3}$,
$ 1+\sqrt{3}$,
$0$;\ \ 
$ 1+\sqrt{3}$,
$ -1-\sqrt{3}$,
$ -1-\sqrt{3}$,
$ -1-\sqrt{3}$,
$ -1-\sqrt{3}$,
$0$;\ \ 
$(3+\sqrt{3})\mathrm{i}$,
$(-3-\sqrt{3})\mathrm{i}$,
$ -1-\sqrt{3}$,
$ 1+\sqrt{3}$,
$0$;\ \ 
$(3+\sqrt{3})\mathrm{i}$,
$ -1-\sqrt{3}$,
$ 1+\sqrt{3}$,
$0$;\ \ 
$ 1$,
$ 1$,
$ 3+\sqrt{3}$;\ \ 
$ 1$,
$ -3-\sqrt{3}$;\ \ 
$0$)

  \vskip 2ex

\noindent47. $10_{4,89.56}^{12,822}$ \irep{587}:\ \ 
$d_i$ = ($1.0$,
$1.0$,
$2.732$,
$2.732$,
$2.732$,
$2.732$,
$2.732$,
$3.732$,
$3.732$,
$4.732$) 

\vskip 0.7ex
\hangindent=3em \hangafter=1
$D^2= 89.569 = 
48+24\sqrt{3}$

\vskip 0.7ex
\hangindent=3em \hangafter=1
$T = ( 0,
\frac{1}{2},
\frac{1}{3},
\frac{3}{4},
\frac{5}{6},
\frac{1}{12},
\frac{1}{12},
0,
\frac{1}{2},
\frac{1}{2} )
$,

\vskip 0.7ex
\hangindent=3em \hangafter=1
$S$ = ($ 1$,
$ 1$,
$ 1+\sqrt{3}$,
$ 1+\sqrt{3}$,
$ 1+\sqrt{3}$,
$ 1+\sqrt{3}$,
$ 1+\sqrt{3}$,
$ 2+\sqrt{3}$,
$ 2+\sqrt{3}$,
$ 3+\sqrt{3}$;\ \ 
$ 1$,
$ 1+\sqrt{3}$,
$ -1-\sqrt{3}$,
$ 1+\sqrt{3}$,
$ -1-\sqrt{3}$,
$ -1-\sqrt{3}$,
$ 2+\sqrt{3}$,
$ 2+\sqrt{3}$,
$ -3-\sqrt{3}$;\ \ 
$ 1+\sqrt{3}$,
$ -2-2\sqrt{3}$,
$ 1+\sqrt{3}$,
$ 1+\sqrt{3}$,
$ 1+\sqrt{3}$,
$ -1-\sqrt{3}$,
$ -1-\sqrt{3}$,
$0$;\ \ 
$0$,
$ 2+2\sqrt{3}$,
$0$,
$0$,
$ -1-\sqrt{3}$,
$ 1+\sqrt{3}$,
$0$;\ \ 
$ 1+\sqrt{3}$,
$ -1-\sqrt{3}$,
$ -1-\sqrt{3}$,
$ -1-\sqrt{3}$,
$ -1-\sqrt{3}$,
$0$;\ \ 
$(3+\sqrt{3})\mathrm{i}$,
$(-3-\sqrt{3})\mathrm{i}$,
$ -1-\sqrt{3}$,
$ 1+\sqrt{3}$,
$0$;\ \ 
$(3+\sqrt{3})\mathrm{i}$,
$ -1-\sqrt{3}$,
$ 1+\sqrt{3}$,
$0$;\ \ 
$ 1$,
$ 1$,
$ 3+\sqrt{3}$;\ \ 
$ 1$,
$ -3-\sqrt{3}$;\ \ 
$0$)

  \vskip 2ex

\noindent48. $10_{2,89.56}^{12,119}$ \irep{680}:\ \ 
$d_i$ = ($1.0$,
$1.0$,
$2.732$,
$2.732$,
$2.732$,
$2.732$,
$2.732$,
$3.732$,
$3.732$,
$4.732$) 

\vskip 0.7ex
\hangindent=3em \hangafter=1
$D^2= 89.569 = 
48+24\sqrt{3}$

\vskip 0.7ex
\hangindent=3em \hangafter=1
$T = ( 0,
\frac{1}{2},
\frac{1}{2},
\frac{1}{3},
\frac{5}{6},
\frac{5}{6},
\frac{5}{6},
0,
\frac{1}{2},
\frac{1}{4} )
$,

\vskip 0.7ex
\hangindent=3em \hangafter=1
$S$ = ($ 1$,
$ 1$,
$ 1+\sqrt{3}$,
$ 1+\sqrt{3}$,
$ 1+\sqrt{3}$,
$ 1+\sqrt{3}$,
$ 1+\sqrt{3}$,
$ 2+\sqrt{3}$,
$ 2+\sqrt{3}$,
$ 3+\sqrt{3}$;\ \ 
$ 1$,
$ -1-\sqrt{3}$,
$ 1+\sqrt{3}$,
$ 1+\sqrt{3}$,
$ -1-\sqrt{3}$,
$ -1-\sqrt{3}$,
$ 2+\sqrt{3}$,
$ 2+\sqrt{3}$,
$ -3-\sqrt{3}$;\ \ 
$0$,
$ -2-2\sqrt{3}$,
$ 2+2\sqrt{3}$,
$0$,
$0$,
$ -1-\sqrt{3}$,
$ 1+\sqrt{3}$,
$0$;\ \ 
$ 1+\sqrt{3}$,
$ 1+\sqrt{3}$,
$ 1+\sqrt{3}$,
$ 1+\sqrt{3}$,
$ -1-\sqrt{3}$,
$ -1-\sqrt{3}$,
$0$;\ \ 
$ 1+\sqrt{3}$,
$ -1-\sqrt{3}$,
$ -1-\sqrt{3}$,
$ -1-\sqrt{3}$,
$ -1-\sqrt{3}$,
$0$;\ \ 
$(-3-\sqrt{3})\mathrm{i}$,
$(3+\sqrt{3})\mathrm{i}$,
$ -1-\sqrt{3}$,
$ 1+\sqrt{3}$,
$0$;\ \ 
$(-3-\sqrt{3})\mathrm{i}$,
$ -1-\sqrt{3}$,
$ 1+\sqrt{3}$,
$0$;\ \ 
$ 1$,
$ 1$,
$ 3+\sqrt{3}$;\ \ 
$ 1$,
$ -3-\sqrt{3}$;\ \ 
$0$)

  \vskip 2ex

\noindent49. $10_{6,89.56}^{12,994}$ \irep{680}:\ \ 
$d_i$ = ($1.0$,
$1.0$,
$2.732$,
$2.732$,
$2.732$,
$2.732$,
$2.732$,
$3.732$,
$3.732$,
$4.732$) 

\vskip 0.7ex
\hangindent=3em \hangafter=1
$D^2= 89.569 = 
48+24\sqrt{3}$

\vskip 0.7ex
\hangindent=3em \hangafter=1
$T = ( 0,
\frac{1}{2},
\frac{1}{2},
\frac{2}{3},
\frac{1}{6},
\frac{1}{6},
\frac{1}{6},
0,
\frac{1}{2},
\frac{3}{4} )
$,

\vskip 0.7ex
\hangindent=3em \hangafter=1
$S$ = ($ 1$,
$ 1$,
$ 1+\sqrt{3}$,
$ 1+\sqrt{3}$,
$ 1+\sqrt{3}$,
$ 1+\sqrt{3}$,
$ 1+\sqrt{3}$,
$ 2+\sqrt{3}$,
$ 2+\sqrt{3}$,
$ 3+\sqrt{3}$;\ \ 
$ 1$,
$ -1-\sqrt{3}$,
$ 1+\sqrt{3}$,
$ 1+\sqrt{3}$,
$ -1-\sqrt{3}$,
$ -1-\sqrt{3}$,
$ 2+\sqrt{3}$,
$ 2+\sqrt{3}$,
$ -3-\sqrt{3}$;\ \ 
$0$,
$ -2-2\sqrt{3}$,
$ 2+2\sqrt{3}$,
$0$,
$0$,
$ -1-\sqrt{3}$,
$ 1+\sqrt{3}$,
$0$;\ \ 
$ 1+\sqrt{3}$,
$ 1+\sqrt{3}$,
$ 1+\sqrt{3}$,
$ 1+\sqrt{3}$,
$ -1-\sqrt{3}$,
$ -1-\sqrt{3}$,
$0$;\ \ 
$ 1+\sqrt{3}$,
$ -1-\sqrt{3}$,
$ -1-\sqrt{3}$,
$ -1-\sqrt{3}$,
$ -1-\sqrt{3}$,
$0$;\ \ 
$(3+\sqrt{3})\mathrm{i}$,
$(-3-\sqrt{3})\mathrm{i}$,
$ -1-\sqrt{3}$,
$ 1+\sqrt{3}$,
$0$;\ \ 
$(3+\sqrt{3})\mathrm{i}$,
$ -1-\sqrt{3}$,
$ 1+\sqrt{3}$,
$0$;\ \ 
$ 1$,
$ 1$,
$ 3+\sqrt{3}$;\ \ 
$ 1$,
$ -3-\sqrt{3}$;\ \ 
$0$)

  \vskip 2ex

\noindent50. $10_{4,89.56}^{12,145}$ \irep{587}:\ \ 
$d_i$ = ($1.0$,
$1.0$,
$2.732$,
$2.732$,
$2.732$,
$2.732$,
$2.732$,
$3.732$,
$3.732$,
$4.732$) 

\vskip 0.7ex
\hangindent=3em \hangafter=1
$D^2= 89.569 = 
48+24\sqrt{3}$

\vskip 0.7ex
\hangindent=3em \hangafter=1
$T = ( 0,
\frac{1}{2},
\frac{2}{3},
\frac{1}{4},
\frac{1}{6},
\frac{11}{12},
\frac{11}{12},
0,
\frac{1}{2},
\frac{1}{2} )
$,

\vskip 0.7ex
\hangindent=3em \hangafter=1
$S$ = ($ 1$,
$ 1$,
$ 1+\sqrt{3}$,
$ 1+\sqrt{3}$,
$ 1+\sqrt{3}$,
$ 1+\sqrt{3}$,
$ 1+\sqrt{3}$,
$ 2+\sqrt{3}$,
$ 2+\sqrt{3}$,
$ 3+\sqrt{3}$;\ \ 
$ 1$,
$ 1+\sqrt{3}$,
$ -1-\sqrt{3}$,
$ 1+\sqrt{3}$,
$ -1-\sqrt{3}$,
$ -1-\sqrt{3}$,
$ 2+\sqrt{3}$,
$ 2+\sqrt{3}$,
$ -3-\sqrt{3}$;\ \ 
$ 1+\sqrt{3}$,
$ -2-2\sqrt{3}$,
$ 1+\sqrt{3}$,
$ 1+\sqrt{3}$,
$ 1+\sqrt{3}$,
$ -1-\sqrt{3}$,
$ -1-\sqrt{3}$,
$0$;\ \ 
$0$,
$ 2+2\sqrt{3}$,
$0$,
$0$,
$ -1-\sqrt{3}$,
$ 1+\sqrt{3}$,
$0$;\ \ 
$ 1+\sqrt{3}$,
$ -1-\sqrt{3}$,
$ -1-\sqrt{3}$,
$ -1-\sqrt{3}$,
$ -1-\sqrt{3}$,
$0$;\ \ 
$(-3-\sqrt{3})\mathrm{i}$,
$(3+\sqrt{3})\mathrm{i}$,
$ -1-\sqrt{3}$,
$ 1+\sqrt{3}$,
$0$;\ \ 
$(-3-\sqrt{3})\mathrm{i}$,
$ -1-\sqrt{3}$,
$ 1+\sqrt{3}$,
$0$;\ \ 
$ 1$,
$ 1$,
$ 3+\sqrt{3}$;\ \ 
$ 1$,
$ -3-\sqrt{3}$;\ \ 
$0$)

  \vskip 2ex

\noindent51. $10_{0,89.56}^{12,200}$ \irep{587}:\ \ 
$d_i$ = ($1.0$,
$1.0$,
$2.732$,
$2.732$,
$2.732$,
$2.732$,
$2.732$,
$3.732$,
$3.732$,
$4.732$) 

\vskip 0.7ex
\hangindent=3em \hangafter=1
$D^2= 89.569 = 
48+24\sqrt{3}$

\vskip 0.7ex
\hangindent=3em \hangafter=1
$T = ( 0,
\frac{1}{2},
\frac{2}{3},
\frac{3}{4},
\frac{1}{6},
\frac{5}{12},
\frac{5}{12},
0,
\frac{1}{2},
0 )
$,

\vskip 0.7ex
\hangindent=3em \hangafter=1
$S$ = ($ 1$,
$ 1$,
$ 1+\sqrt{3}$,
$ 1+\sqrt{3}$,
$ 1+\sqrt{3}$,
$ 1+\sqrt{3}$,
$ 1+\sqrt{3}$,
$ 2+\sqrt{3}$,
$ 2+\sqrt{3}$,
$ 3+\sqrt{3}$;\ \ 
$ 1$,
$ 1+\sqrt{3}$,
$ -1-\sqrt{3}$,
$ 1+\sqrt{3}$,
$ -1-\sqrt{3}$,
$ -1-\sqrt{3}$,
$ 2+\sqrt{3}$,
$ 2+\sqrt{3}$,
$ -3-\sqrt{3}$;\ \ 
$ 1+\sqrt{3}$,
$ -2-2\sqrt{3}$,
$ 1+\sqrt{3}$,
$ 1+\sqrt{3}$,
$ 1+\sqrt{3}$,
$ -1-\sqrt{3}$,
$ -1-\sqrt{3}$,
$0$;\ \ 
$0$,
$ 2+2\sqrt{3}$,
$0$,
$0$,
$ -1-\sqrt{3}$,
$ 1+\sqrt{3}$,
$0$;\ \ 
$ 1+\sqrt{3}$,
$ -1-\sqrt{3}$,
$ -1-\sqrt{3}$,
$ -1-\sqrt{3}$,
$ -1-\sqrt{3}$,
$0$;\ \ 
$(-3-\sqrt{3})\mathrm{i}$,
$(3+\sqrt{3})\mathrm{i}$,
$ -1-\sqrt{3}$,
$ 1+\sqrt{3}$,
$0$;\ \ 
$(-3-\sqrt{3})\mathrm{i}$,
$ -1-\sqrt{3}$,
$ 1+\sqrt{3}$,
$0$;\ \ 
$ 1$,
$ 1$,
$ 3+\sqrt{3}$;\ \ 
$ 1$,
$ -3-\sqrt{3}$;\ \ 
$0$)

  \vskip 2ex

\noindent52. $10_{7,89.56}^{24,123}$ \irep{978}:\ \ 
$d_i$ = ($1.0$,
$1.0$,
$2.732$,
$2.732$,
$2.732$,
$2.732$,
$2.732$,
$3.732$,
$3.732$,
$4.732$) 

\vskip 0.7ex
\hangindent=3em \hangafter=1
$D^2= 89.569 = 
48+24\sqrt{3}$

\vskip 0.7ex
\hangindent=3em \hangafter=1
$T = ( 0,
\frac{1}{2},
\frac{1}{3},
\frac{5}{6},
\frac{1}{8},
\frac{11}{24},
\frac{11}{24},
0,
\frac{1}{2},
\frac{7}{8} )
$,

\vskip 0.7ex
\hangindent=3em \hangafter=1
$S$ = ($ 1$,
$ 1$,
$ 1+\sqrt{3}$,
$ 1+\sqrt{3}$,
$ 1+\sqrt{3}$,
$ 1+\sqrt{3}$,
$ 1+\sqrt{3}$,
$ 2+\sqrt{3}$,
$ 2+\sqrt{3}$,
$ 3+\sqrt{3}$;\ \ 
$ 1$,
$ 1+\sqrt{3}$,
$ 1+\sqrt{3}$,
$ -1-\sqrt{3}$,
$ -1-\sqrt{3}$,
$ -1-\sqrt{3}$,
$ 2+\sqrt{3}$,
$ 2+\sqrt{3}$,
$ -3-\sqrt{3}$;\ \ 
$ 1+\sqrt{3}$,
$ 1+\sqrt{3}$,
$ -2-2\sqrt{3}$,
$ 1+\sqrt{3}$,
$ 1+\sqrt{3}$,
$ -1-\sqrt{3}$,
$ -1-\sqrt{3}$,
$0$;\ \ 
$ 1+\sqrt{3}$,
$ 2+2\sqrt{3}$,
$ -1-\sqrt{3}$,
$ -1-\sqrt{3}$,
$ -1-\sqrt{3}$,
$ -1-\sqrt{3}$,
$0$;\ \ 
$0$,
$0$,
$0$,
$ -1-\sqrt{3}$,
$ 1+\sqrt{3}$,
$0$;\ \ 
$ -3-\sqrt{3}$,
$ 3+\sqrt{3}$,
$ -1-\sqrt{3}$,
$ 1+\sqrt{3}$,
$0$;\ \ 
$ -3-\sqrt{3}$,
$ -1-\sqrt{3}$,
$ 1+\sqrt{3}$,
$0$;\ \ 
$ 1$,
$ 1$,
$ 3+\sqrt{3}$;\ \ 
$ 1$,
$ -3-\sqrt{3}$;\ \ 
$0$)

  \vskip 2ex

\noindent53. $10_{1,89.56}^{24,380}$ \irep{977}:\ \ 
$d_i$ = ($1.0$,
$1.0$,
$2.732$,
$2.732$,
$2.732$,
$2.732$,
$2.732$,
$3.732$,
$3.732$,
$4.732$) 

\vskip 0.7ex
\hangindent=3em \hangafter=1
$D^2= 89.569 = 
48+24\sqrt{3}$

\vskip 0.7ex
\hangindent=3em \hangafter=1
$T = ( 0,
\frac{1}{2},
\frac{1}{3},
\frac{5}{6},
\frac{3}{8},
\frac{17}{24},
\frac{17}{24},
0,
\frac{1}{2},
\frac{1}{8} )
$,

\vskip 0.7ex
\hangindent=3em \hangafter=1
$S$ = ($ 1$,
$ 1$,
$ 1+\sqrt{3}$,
$ 1+\sqrt{3}$,
$ 1+\sqrt{3}$,
$ 1+\sqrt{3}$,
$ 1+\sqrt{3}$,
$ 2+\sqrt{3}$,
$ 2+\sqrt{3}$,
$ 3+\sqrt{3}$;\ \ 
$ 1$,
$ 1+\sqrt{3}$,
$ 1+\sqrt{3}$,
$ -1-\sqrt{3}$,
$ -1-\sqrt{3}$,
$ -1-\sqrt{3}$,
$ 2+\sqrt{3}$,
$ 2+\sqrt{3}$,
$ -3-\sqrt{3}$;\ \ 
$ 1+\sqrt{3}$,
$ 1+\sqrt{3}$,
$ -2-2\sqrt{3}$,
$ 1+\sqrt{3}$,
$ 1+\sqrt{3}$,
$ -1-\sqrt{3}$,
$ -1-\sqrt{3}$,
$0$;\ \ 
$ 1+\sqrt{3}$,
$ 2+2\sqrt{3}$,
$ -1-\sqrt{3}$,
$ -1-\sqrt{3}$,
$ -1-\sqrt{3}$,
$ -1-\sqrt{3}$,
$0$;\ \ 
$0$,
$0$,
$0$,
$ -1-\sqrt{3}$,
$ 1+\sqrt{3}$,
$0$;\ \ 
$ 3+\sqrt{3}$,
$ -3-\sqrt{3}$,
$ -1-\sqrt{3}$,
$ 1+\sqrt{3}$,
$0$;\ \ 
$ 3+\sqrt{3}$,
$ -1-\sqrt{3}$,
$ 1+\sqrt{3}$,
$0$;\ \ 
$ 1$,
$ 1$,
$ 3+\sqrt{3}$;\ \ 
$ 1$,
$ -3-\sqrt{3}$;\ \ 
$0$)

  \vskip 2ex

\noindent54. $10_{3,89.56}^{24,317}$ \irep{978}:\ \ 
$d_i$ = ($1.0$,
$1.0$,
$2.732$,
$2.732$,
$2.732$,
$2.732$,
$2.732$,
$3.732$,
$3.732$,
$4.732$) 

\vskip 0.7ex
\hangindent=3em \hangafter=1
$D^2= 89.569 = 
48+24\sqrt{3}$

\vskip 0.7ex
\hangindent=3em \hangafter=1
$T = ( 0,
\frac{1}{2},
\frac{1}{3},
\frac{5}{6},
\frac{5}{8},
\frac{23}{24},
\frac{23}{24},
0,
\frac{1}{2},
\frac{3}{8} )
$,

\vskip 0.7ex
\hangindent=3em \hangafter=1
$S$ = ($ 1$,
$ 1$,
$ 1+\sqrt{3}$,
$ 1+\sqrt{3}$,
$ 1+\sqrt{3}$,
$ 1+\sqrt{3}$,
$ 1+\sqrt{3}$,
$ 2+\sqrt{3}$,
$ 2+\sqrt{3}$,
$ 3+\sqrt{3}$;\ \ 
$ 1$,
$ 1+\sqrt{3}$,
$ 1+\sqrt{3}$,
$ -1-\sqrt{3}$,
$ -1-\sqrt{3}$,
$ -1-\sqrt{3}$,
$ 2+\sqrt{3}$,
$ 2+\sqrt{3}$,
$ -3-\sqrt{3}$;\ \ 
$ 1+\sqrt{3}$,
$ 1+\sqrt{3}$,
$ -2-2\sqrt{3}$,
$ 1+\sqrt{3}$,
$ 1+\sqrt{3}$,
$ -1-\sqrt{3}$,
$ -1-\sqrt{3}$,
$0$;\ \ 
$ 1+\sqrt{3}$,
$ 2+2\sqrt{3}$,
$ -1-\sqrt{3}$,
$ -1-\sqrt{3}$,
$ -1-\sqrt{3}$,
$ -1-\sqrt{3}$,
$0$;\ \ 
$0$,
$0$,
$0$,
$ -1-\sqrt{3}$,
$ 1+\sqrt{3}$,
$0$;\ \ 
$ -3-\sqrt{3}$,
$ 3+\sqrt{3}$,
$ -1-\sqrt{3}$,
$ 1+\sqrt{3}$,
$0$;\ \ 
$ -3-\sqrt{3}$,
$ -1-\sqrt{3}$,
$ 1+\sqrt{3}$,
$0$;\ \ 
$ 1$,
$ 1$,
$ 3+\sqrt{3}$;\ \ 
$ 1$,
$ -3-\sqrt{3}$;\ \ 
$0$)

  \vskip 2ex

\noindent55. $10_{5,89.56}^{24,597}$ \irep{977}:\ \ 
$d_i$ = ($1.0$,
$1.0$,
$2.732$,
$2.732$,
$2.732$,
$2.732$,
$2.732$,
$3.732$,
$3.732$,
$4.732$) 

\vskip 0.7ex
\hangindent=3em \hangafter=1
$D^2= 89.569 = 
48+24\sqrt{3}$

\vskip 0.7ex
\hangindent=3em \hangafter=1
$T = ( 0,
\frac{1}{2},
\frac{1}{3},
\frac{5}{6},
\frac{7}{8},
\frac{5}{24},
\frac{5}{24},
0,
\frac{1}{2},
\frac{5}{8} )
$,

\vskip 0.7ex
\hangindent=3em \hangafter=1
$S$ = ($ 1$,
$ 1$,
$ 1+\sqrt{3}$,
$ 1+\sqrt{3}$,
$ 1+\sqrt{3}$,
$ 1+\sqrt{3}$,
$ 1+\sqrt{3}$,
$ 2+\sqrt{3}$,
$ 2+\sqrt{3}$,
$ 3+\sqrt{3}$;\ \ 
$ 1$,
$ 1+\sqrt{3}$,
$ 1+\sqrt{3}$,
$ -1-\sqrt{3}$,
$ -1-\sqrt{3}$,
$ -1-\sqrt{3}$,
$ 2+\sqrt{3}$,
$ 2+\sqrt{3}$,
$ -3-\sqrt{3}$;\ \ 
$ 1+\sqrt{3}$,
$ 1+\sqrt{3}$,
$ -2-2\sqrt{3}$,
$ 1+\sqrt{3}$,
$ 1+\sqrt{3}$,
$ -1-\sqrt{3}$,
$ -1-\sqrt{3}$,
$0$;\ \ 
$ 1+\sqrt{3}$,
$ 2+2\sqrt{3}$,
$ -1-\sqrt{3}$,
$ -1-\sqrt{3}$,
$ -1-\sqrt{3}$,
$ -1-\sqrt{3}$,
$0$;\ \ 
$0$,
$0$,
$0$,
$ -1-\sqrt{3}$,
$ 1+\sqrt{3}$,
$0$;\ \ 
$ 3+\sqrt{3}$,
$ -3-\sqrt{3}$,
$ -1-\sqrt{3}$,
$ 1+\sqrt{3}$,
$0$;\ \ 
$ 3+\sqrt{3}$,
$ -1-\sqrt{3}$,
$ 1+\sqrt{3}$,
$0$;\ \ 
$ 1$,
$ 1$,
$ 3+\sqrt{3}$;\ \ 
$ 1$,
$ -3-\sqrt{3}$;\ \ 
$0$)

  \vskip 2ex

\noindent56. $10_{3,89.56}^{24,358}$ \irep{977}:\ \ 
$d_i$ = ($1.0$,
$1.0$,
$2.732$,
$2.732$,
$2.732$,
$2.732$,
$2.732$,
$3.732$,
$3.732$,
$4.732$) 

\vskip 0.7ex
\hangindent=3em \hangafter=1
$D^2= 89.569 = 
48+24\sqrt{3}$

\vskip 0.7ex
\hangindent=3em \hangafter=1
$T = ( 0,
\frac{1}{2},
\frac{2}{3},
\frac{1}{6},
\frac{1}{8},
\frac{19}{24},
\frac{19}{24},
0,
\frac{1}{2},
\frac{3}{8} )
$,

\vskip 0.7ex
\hangindent=3em \hangafter=1
$S$ = ($ 1$,
$ 1$,
$ 1+\sqrt{3}$,
$ 1+\sqrt{3}$,
$ 1+\sqrt{3}$,
$ 1+\sqrt{3}$,
$ 1+\sqrt{3}$,
$ 2+\sqrt{3}$,
$ 2+\sqrt{3}$,
$ 3+\sqrt{3}$;\ \ 
$ 1$,
$ 1+\sqrt{3}$,
$ 1+\sqrt{3}$,
$ -1-\sqrt{3}$,
$ -1-\sqrt{3}$,
$ -1-\sqrt{3}$,
$ 2+\sqrt{3}$,
$ 2+\sqrt{3}$,
$ -3-\sqrt{3}$;\ \ 
$ 1+\sqrt{3}$,
$ 1+\sqrt{3}$,
$ -2-2\sqrt{3}$,
$ 1+\sqrt{3}$,
$ 1+\sqrt{3}$,
$ -1-\sqrt{3}$,
$ -1-\sqrt{3}$,
$0$;\ \ 
$ 1+\sqrt{3}$,
$ 2+2\sqrt{3}$,
$ -1-\sqrt{3}$,
$ -1-\sqrt{3}$,
$ -1-\sqrt{3}$,
$ -1-\sqrt{3}$,
$0$;\ \ 
$0$,
$0$,
$0$,
$ -1-\sqrt{3}$,
$ 1+\sqrt{3}$,
$0$;\ \ 
$ 3+\sqrt{3}$,
$ -3-\sqrt{3}$,
$ -1-\sqrt{3}$,
$ 1+\sqrt{3}$,
$0$;\ \ 
$ 3+\sqrt{3}$,
$ -1-\sqrt{3}$,
$ 1+\sqrt{3}$,
$0$;\ \ 
$ 1$,
$ 1$,
$ 3+\sqrt{3}$;\ \ 
$ 1$,
$ -3-\sqrt{3}$;\ \ 
$0$)

  \vskip 2ex

\noindent57. $10_{5,89.56}^{24,224}$ \irep{978}:\ \ 
$d_i$ = ($1.0$,
$1.0$,
$2.732$,
$2.732$,
$2.732$,
$2.732$,
$2.732$,
$3.732$,
$3.732$,
$4.732$) 

\vskip 0.7ex
\hangindent=3em \hangafter=1
$D^2= 89.569 = 
48+24\sqrt{3}$

\vskip 0.7ex
\hangindent=3em \hangafter=1
$T = ( 0,
\frac{1}{2},
\frac{2}{3},
\frac{1}{6},
\frac{3}{8},
\frac{1}{24},
\frac{1}{24},
0,
\frac{1}{2},
\frac{5}{8} )
$,

\vskip 0.7ex
\hangindent=3em \hangafter=1
$S$ = ($ 1$,
$ 1$,
$ 1+\sqrt{3}$,
$ 1+\sqrt{3}$,
$ 1+\sqrt{3}$,
$ 1+\sqrt{3}$,
$ 1+\sqrt{3}$,
$ 2+\sqrt{3}$,
$ 2+\sqrt{3}$,
$ 3+\sqrt{3}$;\ \ 
$ 1$,
$ 1+\sqrt{3}$,
$ 1+\sqrt{3}$,
$ -1-\sqrt{3}$,
$ -1-\sqrt{3}$,
$ -1-\sqrt{3}$,
$ 2+\sqrt{3}$,
$ 2+\sqrt{3}$,
$ -3-\sqrt{3}$;\ \ 
$ 1+\sqrt{3}$,
$ 1+\sqrt{3}$,
$ -2-2\sqrt{3}$,
$ 1+\sqrt{3}$,
$ 1+\sqrt{3}$,
$ -1-\sqrt{3}$,
$ -1-\sqrt{3}$,
$0$;\ \ 
$ 1+\sqrt{3}$,
$ 2+2\sqrt{3}$,
$ -1-\sqrt{3}$,
$ -1-\sqrt{3}$,
$ -1-\sqrt{3}$,
$ -1-\sqrt{3}$,
$0$;\ \ 
$0$,
$0$,
$0$,
$ -1-\sqrt{3}$,
$ 1+\sqrt{3}$,
$0$;\ \ 
$ -3-\sqrt{3}$,
$ 3+\sqrt{3}$,
$ -1-\sqrt{3}$,
$ 1+\sqrt{3}$,
$0$;\ \ 
$ -3-\sqrt{3}$,
$ -1-\sqrt{3}$,
$ 1+\sqrt{3}$,
$0$;\ \ 
$ 1$,
$ 1$,
$ 3+\sqrt{3}$;\ \ 
$ 1$,
$ -3-\sqrt{3}$;\ \ 
$0$)

  \vskip 2ex

\noindent58. $10_{7,89.56}^{24,664}$ \irep{977}:\ \ 
$d_i$ = ($1.0$,
$1.0$,
$2.732$,
$2.732$,
$2.732$,
$2.732$,
$2.732$,
$3.732$,
$3.732$,
$4.732$) 

\vskip 0.7ex
\hangindent=3em \hangafter=1
$D^2= 89.569 = 
48+24\sqrt{3}$

\vskip 0.7ex
\hangindent=3em \hangafter=1
$T = ( 0,
\frac{1}{2},
\frac{2}{3},
\frac{1}{6},
\frac{5}{8},
\frac{7}{24},
\frac{7}{24},
0,
\frac{1}{2},
\frac{7}{8} )
$,

\vskip 0.7ex
\hangindent=3em \hangafter=1
$S$ = ($ 1$,
$ 1$,
$ 1+\sqrt{3}$,
$ 1+\sqrt{3}$,
$ 1+\sqrt{3}$,
$ 1+\sqrt{3}$,
$ 1+\sqrt{3}$,
$ 2+\sqrt{3}$,
$ 2+\sqrt{3}$,
$ 3+\sqrt{3}$;\ \ 
$ 1$,
$ 1+\sqrt{3}$,
$ 1+\sqrt{3}$,
$ -1-\sqrt{3}$,
$ -1-\sqrt{3}$,
$ -1-\sqrt{3}$,
$ 2+\sqrt{3}$,
$ 2+\sqrt{3}$,
$ -3-\sqrt{3}$;\ \ 
$ 1+\sqrt{3}$,
$ 1+\sqrt{3}$,
$ -2-2\sqrt{3}$,
$ 1+\sqrt{3}$,
$ 1+\sqrt{3}$,
$ -1-\sqrt{3}$,
$ -1-\sqrt{3}$,
$0$;\ \ 
$ 1+\sqrt{3}$,
$ 2+2\sqrt{3}$,
$ -1-\sqrt{3}$,
$ -1-\sqrt{3}$,
$ -1-\sqrt{3}$,
$ -1-\sqrt{3}$,
$0$;\ \ 
$0$,
$0$,
$0$,
$ -1-\sqrt{3}$,
$ 1+\sqrt{3}$,
$0$;\ \ 
$ 3+\sqrt{3}$,
$ -3-\sqrt{3}$,
$ -1-\sqrt{3}$,
$ 1+\sqrt{3}$,
$0$;\ \ 
$ 3+\sqrt{3}$,
$ -1-\sqrt{3}$,
$ 1+\sqrt{3}$,
$0$;\ \ 
$ 1$,
$ 1$,
$ 3+\sqrt{3}$;\ \ 
$ 1$,
$ -3-\sqrt{3}$;\ \ 
$0$)

  \vskip 2ex

\noindent59. $10_{1,89.56}^{24,722}$ \irep{978}:\ \ 
$d_i$ = ($1.0$,
$1.0$,
$2.732$,
$2.732$,
$2.732$,
$2.732$,
$2.732$,
$3.732$,
$3.732$,
$4.732$) 

\vskip 0.7ex
\hangindent=3em \hangafter=1
$D^2= 89.569 = 
48+24\sqrt{3}$

\vskip 0.7ex
\hangindent=3em \hangafter=1
$T = ( 0,
\frac{1}{2},
\frac{2}{3},
\frac{1}{6},
\frac{7}{8},
\frac{13}{24},
\frac{13}{24},
0,
\frac{1}{2},
\frac{1}{8} )
$,

\vskip 0.7ex
\hangindent=3em \hangafter=1
$S$ = ($ 1$,
$ 1$,
$ 1+\sqrt{3}$,
$ 1+\sqrt{3}$,
$ 1+\sqrt{3}$,
$ 1+\sqrt{3}$,
$ 1+\sqrt{3}$,
$ 2+\sqrt{3}$,
$ 2+\sqrt{3}$,
$ 3+\sqrt{3}$;\ \ 
$ 1$,
$ 1+\sqrt{3}$,
$ 1+\sqrt{3}$,
$ -1-\sqrt{3}$,
$ -1-\sqrt{3}$,
$ -1-\sqrt{3}$,
$ 2+\sqrt{3}$,
$ 2+\sqrt{3}$,
$ -3-\sqrt{3}$;\ \ 
$ 1+\sqrt{3}$,
$ 1+\sqrt{3}$,
$ -2-2\sqrt{3}$,
$ 1+\sqrt{3}$,
$ 1+\sqrt{3}$,
$ -1-\sqrt{3}$,
$ -1-\sqrt{3}$,
$0$;\ \ 
$ 1+\sqrt{3}$,
$ 2+2\sqrt{3}$,
$ -1-\sqrt{3}$,
$ -1-\sqrt{3}$,
$ -1-\sqrt{3}$,
$ -1-\sqrt{3}$,
$0$;\ \ 
$0$,
$0$,
$0$,
$ -1-\sqrt{3}$,
$ 1+\sqrt{3}$,
$0$;\ \ 
$ -3-\sqrt{3}$,
$ 3+\sqrt{3}$,
$ -1-\sqrt{3}$,
$ 1+\sqrt{3}$,
$0$;\ \ 
$ -3-\sqrt{3}$,
$ -1-\sqrt{3}$,
$ 1+\sqrt{3}$,
$0$;\ \ 
$ 1$,
$ 1$,
$ 3+\sqrt{3}$;\ \ 
$ 1$,
$ -3-\sqrt{3}$;\ \ 
$0$)

  \vskip 2ex

\noindent60. $10_{\frac{74}{55},125.3}^{55,334}$ \irep{1121}:\ \ 
$d_i$ = ($1.0$,
$1.618$,
$1.918$,
$2.682$,
$3.104$,
$3.228$,
$3.513$,
$4.340$,
$5.224$,
$5.684$) 

\vskip 0.7ex
\hangindent=3em \hangafter=1
$D^2= 125.351 = 
37+5c^{2}_{55}
+8c^{3}_{55}
+22c^{5}_{55}
+4c^{6}_{55}
+3c^{7}_{55}
+8c^{8}_{55}
+10c^{10}_{55}
+9c^{11}_{55}
+5c^{13}_{55}
+2c^{14}_{55}
+4c^{15}_{55}
+4c^{16}_{55}
+3c^{18}_{55}
+2c^{19}_{55}
$

\vskip 0.7ex
\hangindent=3em \hangafter=1
$T = ( 0,
\frac{2}{5},
\frac{2}{11},
\frac{9}{11},
\frac{32}{55},
\frac{10}{11},
\frac{5}{11},
\frac{12}{55},
\frac{17}{55},
\frac{47}{55} )
$,

\vskip 0.7ex
\hangindent=3em \hangafter=1
$S$ = ($ 1$,
$ \frac{1+\sqrt{5}}{2}$,
$ -c_{11}^{5}$,
$ \xi_{11}^{3}$,
$ c^{3}_{55}
+c^{8}_{55}
$,
$ \xi_{11}^{7}$,
$ \xi_{11}^{5}$,
$ 1+c^{5}_{55}
+c^{6}_{55}
+c^{11}_{55}
+c^{16}_{55}
$,
$ c^{2}_{55}
+c^{3}_{55}
+c^{8}_{55}
+c^{13}_{55}
$,
$ c^{2}_{55}
+c^{3}_{55}
+c^{7}_{55}
+c^{8}_{55}
+c^{13}_{55}
+c^{18}_{55}
$;\ \ 
$ -1$,
$ c^{3}_{55}
+c^{8}_{55}
$,
$ 1+c^{5}_{55}
+c^{6}_{55}
+c^{11}_{55}
+c^{16}_{55}
$,
$ c_{11}^{5}$,
$ c^{2}_{55}
+c^{3}_{55}
+c^{8}_{55}
+c^{13}_{55}
$,
$ c^{2}_{55}
+c^{3}_{55}
+c^{7}_{55}
+c^{8}_{55}
+c^{13}_{55}
+c^{18}_{55}
$,
$ -\xi_{11}^{3}$,
$ -\xi_{11}^{7}$,
$ -\xi_{11}^{5}$;\ \ 
$ -\xi_{11}^{7}$,
$ \xi_{11}^{5}$,
$ -c^{2}_{55}
-c^{3}_{55}
-c^{8}_{55}
-c^{13}_{55}
$,
$ -\xi_{11}^{3}$,
$ 1$,
$ c^{2}_{55}
+c^{3}_{55}
+c^{7}_{55}
+c^{8}_{55}
+c^{13}_{55}
+c^{18}_{55}
$,
$ -1-c^{5}_{55}
-c^{6}_{55}
-c^{11}_{55}
-c^{16}_{55}
$,
$ \frac{1+\sqrt{5}}{2}$;\ \ 
$ -c_{11}^{5}$,
$ c^{2}_{55}
+c^{3}_{55}
+c^{7}_{55}
+c^{8}_{55}
+c^{13}_{55}
+c^{18}_{55}
$,
$ -1$,
$ -\xi_{11}^{7}$,
$ c^{3}_{55}
+c^{8}_{55}
$,
$ -\frac{1+\sqrt{5}}{2}$,
$ -c^{2}_{55}
-c^{3}_{55}
-c^{8}_{55}
-c^{13}_{55}
$;\ \ 
$ \xi_{11}^{7}$,
$ -1-c^{5}_{55}
-c^{6}_{55}
-c^{11}_{55}
-c^{16}_{55}
$,
$ \frac{1+\sqrt{5}}{2}$,
$ -\xi_{11}^{5}$,
$ \xi_{11}^{3}$,
$ -1$;\ \ 
$ \xi_{11}^{5}$,
$ c_{11}^{5}$,
$ -\frac{1+\sqrt{5}}{2}$,
$ c^{2}_{55}
+c^{3}_{55}
+c^{7}_{55}
+c^{8}_{55}
+c^{13}_{55}
+c^{18}_{55}
$,
$ -c^{3}_{55}
-c^{8}_{55}
$;\ \ 
$ \xi_{11}^{3}$,
$ -c^{2}_{55}
-c^{3}_{55}
-c^{8}_{55}
-c^{13}_{55}
$,
$ -c^{3}_{55}
-c^{8}_{55}
$,
$ 1+c^{5}_{55}
+c^{6}_{55}
+c^{11}_{55}
+c^{16}_{55}
$;\ \ 
$ c_{11}^{5}$,
$ 1$,
$ \xi_{11}^{7}$;\ \ 
$ -\xi_{11}^{5}$,
$ -c_{11}^{5}$;\ \ 
$ -\xi_{11}^{3}$)

Factors = $2_{\frac{14}{5},3.618}^{5,395}\boxtimes 5_{\frac{72}{11},34.64}^{11,216}$

  \vskip 2ex

\noindent61. $10_{\frac{234}{55},125.3}^{55,258}$ \irep{1121}:\ \ 
$d_i$ = ($1.0$,
$1.618$,
$1.918$,
$2.682$,
$3.104$,
$3.228$,
$3.513$,
$4.340$,
$5.224$,
$5.684$) 

\vskip 0.7ex
\hangindent=3em \hangafter=1
$D^2= 125.351 = 
37+5c^{2}_{55}
+8c^{3}_{55}
+22c^{5}_{55}
+4c^{6}_{55}
+3c^{7}_{55}
+8c^{8}_{55}
+10c^{10}_{55}
+9c^{11}_{55}
+5c^{13}_{55}
+2c^{14}_{55}
+4c^{15}_{55}
+4c^{16}_{55}
+3c^{18}_{55}
+2c^{19}_{55}
$

\vskip 0.7ex
\hangindent=3em \hangafter=1
$T = ( 0,
\frac{2}{5},
\frac{9}{11},
\frac{2}{11},
\frac{12}{55},
\frac{1}{11},
\frac{6}{11},
\frac{32}{55},
\frac{27}{55},
\frac{52}{55} )
$,

\vskip 0.7ex
\hangindent=3em \hangafter=1
$S$ = ($ 1$,
$ \frac{1+\sqrt{5}}{2}$,
$ -c_{11}^{5}$,
$ \xi_{11}^{3}$,
$ c^{3}_{55}
+c^{8}_{55}
$,
$ \xi_{11}^{7}$,
$ \xi_{11}^{5}$,
$ 1+c^{5}_{55}
+c^{6}_{55}
+c^{11}_{55}
+c^{16}_{55}
$,
$ c^{2}_{55}
+c^{3}_{55}
+c^{8}_{55}
+c^{13}_{55}
$,
$ c^{2}_{55}
+c^{3}_{55}
+c^{7}_{55}
+c^{8}_{55}
+c^{13}_{55}
+c^{18}_{55}
$;\ \ 
$ -1$,
$ c^{3}_{55}
+c^{8}_{55}
$,
$ 1+c^{5}_{55}
+c^{6}_{55}
+c^{11}_{55}
+c^{16}_{55}
$,
$ c_{11}^{5}$,
$ c^{2}_{55}
+c^{3}_{55}
+c^{8}_{55}
+c^{13}_{55}
$,
$ c^{2}_{55}
+c^{3}_{55}
+c^{7}_{55}
+c^{8}_{55}
+c^{13}_{55}
+c^{18}_{55}
$,
$ -\xi_{11}^{3}$,
$ -\xi_{11}^{7}$,
$ -\xi_{11}^{5}$;\ \ 
$ -\xi_{11}^{7}$,
$ \xi_{11}^{5}$,
$ -c^{2}_{55}
-c^{3}_{55}
-c^{8}_{55}
-c^{13}_{55}
$,
$ -\xi_{11}^{3}$,
$ 1$,
$ c^{2}_{55}
+c^{3}_{55}
+c^{7}_{55}
+c^{8}_{55}
+c^{13}_{55}
+c^{18}_{55}
$,
$ -1-c^{5}_{55}
-c^{6}_{55}
-c^{11}_{55}
-c^{16}_{55}
$,
$ \frac{1+\sqrt{5}}{2}$;\ \ 
$ -c_{11}^{5}$,
$ c^{2}_{55}
+c^{3}_{55}
+c^{7}_{55}
+c^{8}_{55}
+c^{13}_{55}
+c^{18}_{55}
$,
$ -1$,
$ -\xi_{11}^{7}$,
$ c^{3}_{55}
+c^{8}_{55}
$,
$ -\frac{1+\sqrt{5}}{2}$,
$ -c^{2}_{55}
-c^{3}_{55}
-c^{8}_{55}
-c^{13}_{55}
$;\ \ 
$ \xi_{11}^{7}$,
$ -1-c^{5}_{55}
-c^{6}_{55}
-c^{11}_{55}
-c^{16}_{55}
$,
$ \frac{1+\sqrt{5}}{2}$,
$ -\xi_{11}^{5}$,
$ \xi_{11}^{3}$,
$ -1$;\ \ 
$ \xi_{11}^{5}$,
$ c_{11}^{5}$,
$ -\frac{1+\sqrt{5}}{2}$,
$ c^{2}_{55}
+c^{3}_{55}
+c^{7}_{55}
+c^{8}_{55}
+c^{13}_{55}
+c^{18}_{55}
$,
$ -c^{3}_{55}
-c^{8}_{55}
$;\ \ 
$ \xi_{11}^{3}$,
$ -c^{2}_{55}
-c^{3}_{55}
-c^{8}_{55}
-c^{13}_{55}
$,
$ -c^{3}_{55}
-c^{8}_{55}
$,
$ 1+c^{5}_{55}
+c^{6}_{55}
+c^{11}_{55}
+c^{16}_{55}
$;\ \ 
$ c_{11}^{5}$,
$ 1$,
$ \xi_{11}^{7}$;\ \ 
$ -\xi_{11}^{5}$,
$ -c_{11}^{5}$;\ \ 
$ -\xi_{11}^{3}$)

Factors = $2_{\frac{14}{5},3.618}^{5,395}\boxtimes 5_{\frac{16}{11},34.64}^{11,640}$

  \vskip 2ex

\noindent62. $10_{\frac{206}{55},125.3}^{55,178}$ \irep{1121}:\ \ 
$d_i$ = ($1.0$,
$1.618$,
$1.918$,
$2.682$,
$3.104$,
$3.228$,
$3.513$,
$4.340$,
$5.224$,
$5.684$) 

\vskip 0.7ex
\hangindent=3em \hangafter=1
$D^2= 125.351 = 
37+5c^{2}_{55}
+8c^{3}_{55}
+22c^{5}_{55}
+4c^{6}_{55}
+3c^{7}_{55}
+8c^{8}_{55}
+10c^{10}_{55}
+9c^{11}_{55}
+5c^{13}_{55}
+2c^{14}_{55}
+4c^{15}_{55}
+4c^{16}_{55}
+3c^{18}_{55}
+2c^{19}_{55}
$

\vskip 0.7ex
\hangindent=3em \hangafter=1
$T = ( 0,
\frac{3}{5},
\frac{2}{11},
\frac{9}{11},
\frac{43}{55},
\frac{10}{11},
\frac{5}{11},
\frac{23}{55},
\frac{28}{55},
\frac{3}{55} )
$,

\vskip 0.7ex
\hangindent=3em \hangafter=1
$S$ = ($ 1$,
$ \frac{1+\sqrt{5}}{2}$,
$ -c_{11}^{5}$,
$ \xi_{11}^{3}$,
$ c^{3}_{55}
+c^{8}_{55}
$,
$ \xi_{11}^{7}$,
$ \xi_{11}^{5}$,
$ 1+c^{5}_{55}
+c^{6}_{55}
+c^{11}_{55}
+c^{16}_{55}
$,
$ c^{2}_{55}
+c^{3}_{55}
+c^{8}_{55}
+c^{13}_{55}
$,
$ c^{2}_{55}
+c^{3}_{55}
+c^{7}_{55}
+c^{8}_{55}
+c^{13}_{55}
+c^{18}_{55}
$;\ \ 
$ -1$,
$ c^{3}_{55}
+c^{8}_{55}
$,
$ 1+c^{5}_{55}
+c^{6}_{55}
+c^{11}_{55}
+c^{16}_{55}
$,
$ c_{11}^{5}$,
$ c^{2}_{55}
+c^{3}_{55}
+c^{8}_{55}
+c^{13}_{55}
$,
$ c^{2}_{55}
+c^{3}_{55}
+c^{7}_{55}
+c^{8}_{55}
+c^{13}_{55}
+c^{18}_{55}
$,
$ -\xi_{11}^{3}$,
$ -\xi_{11}^{7}$,
$ -\xi_{11}^{5}$;\ \ 
$ -\xi_{11}^{7}$,
$ \xi_{11}^{5}$,
$ -c^{2}_{55}
-c^{3}_{55}
-c^{8}_{55}
-c^{13}_{55}
$,
$ -\xi_{11}^{3}$,
$ 1$,
$ c^{2}_{55}
+c^{3}_{55}
+c^{7}_{55}
+c^{8}_{55}
+c^{13}_{55}
+c^{18}_{55}
$,
$ -1-c^{5}_{55}
-c^{6}_{55}
-c^{11}_{55}
-c^{16}_{55}
$,
$ \frac{1+\sqrt{5}}{2}$;\ \ 
$ -c_{11}^{5}$,
$ c^{2}_{55}
+c^{3}_{55}
+c^{7}_{55}
+c^{8}_{55}
+c^{13}_{55}
+c^{18}_{55}
$,
$ -1$,
$ -\xi_{11}^{7}$,
$ c^{3}_{55}
+c^{8}_{55}
$,
$ -\frac{1+\sqrt{5}}{2}$,
$ -c^{2}_{55}
-c^{3}_{55}
-c^{8}_{55}
-c^{13}_{55}
$;\ \ 
$ \xi_{11}^{7}$,
$ -1-c^{5}_{55}
-c^{6}_{55}
-c^{11}_{55}
-c^{16}_{55}
$,
$ \frac{1+\sqrt{5}}{2}$,
$ -\xi_{11}^{5}$,
$ \xi_{11}^{3}$,
$ -1$;\ \ 
$ \xi_{11}^{5}$,
$ c_{11}^{5}$,
$ -\frac{1+\sqrt{5}}{2}$,
$ c^{2}_{55}
+c^{3}_{55}
+c^{7}_{55}
+c^{8}_{55}
+c^{13}_{55}
+c^{18}_{55}
$,
$ -c^{3}_{55}
-c^{8}_{55}
$;\ \ 
$ \xi_{11}^{3}$,
$ -c^{2}_{55}
-c^{3}_{55}
-c^{8}_{55}
-c^{13}_{55}
$,
$ -c^{3}_{55}
-c^{8}_{55}
$,
$ 1+c^{5}_{55}
+c^{6}_{55}
+c^{11}_{55}
+c^{16}_{55}
$;\ \ 
$ c_{11}^{5}$,
$ 1$,
$ \xi_{11}^{7}$;\ \ 
$ -\xi_{11}^{5}$,
$ -c_{11}^{5}$;\ \ 
$ -\xi_{11}^{3}$)

Factors = $2_{\frac{26}{5},3.618}^{5,720}\boxtimes 5_{\frac{72}{11},34.64}^{11,216}$

  \vskip 2ex

\noindent63. $10_{\frac{366}{55},125.3}^{55,758}$ \irep{1121}:\ \ 
$d_i$ = ($1.0$,
$1.618$,
$1.918$,
$2.682$,
$3.104$,
$3.228$,
$3.513$,
$4.340$,
$5.224$,
$5.684$) 

\vskip 0.7ex
\hangindent=3em \hangafter=1
$D^2= 125.351 = 
37+5c^{2}_{55}
+8c^{3}_{55}
+22c^{5}_{55}
+4c^{6}_{55}
+3c^{7}_{55}
+8c^{8}_{55}
+10c^{10}_{55}
+9c^{11}_{55}
+5c^{13}_{55}
+2c^{14}_{55}
+4c^{15}_{55}
+4c^{16}_{55}
+3c^{18}_{55}
+2c^{19}_{55}
$

\vskip 0.7ex
\hangindent=3em \hangafter=1
$T = ( 0,
\frac{3}{5},
\frac{9}{11},
\frac{2}{11},
\frac{23}{55},
\frac{1}{11},
\frac{6}{11},
\frac{43}{55},
\frac{38}{55},
\frac{8}{55} )
$,

\vskip 0.7ex
\hangindent=3em \hangafter=1
$S$ = ($ 1$,
$ \frac{1+\sqrt{5}}{2}$,
$ -c_{11}^{5}$,
$ \xi_{11}^{3}$,
$ c^{3}_{55}
+c^{8}_{55}
$,
$ \xi_{11}^{7}$,
$ \xi_{11}^{5}$,
$ 1+c^{5}_{55}
+c^{6}_{55}
+c^{11}_{55}
+c^{16}_{55}
$,
$ c^{2}_{55}
+c^{3}_{55}
+c^{8}_{55}
+c^{13}_{55}
$,
$ c^{2}_{55}
+c^{3}_{55}
+c^{7}_{55}
+c^{8}_{55}
+c^{13}_{55}
+c^{18}_{55}
$;\ \ 
$ -1$,
$ c^{3}_{55}
+c^{8}_{55}
$,
$ 1+c^{5}_{55}
+c^{6}_{55}
+c^{11}_{55}
+c^{16}_{55}
$,
$ c_{11}^{5}$,
$ c^{2}_{55}
+c^{3}_{55}
+c^{8}_{55}
+c^{13}_{55}
$,
$ c^{2}_{55}
+c^{3}_{55}
+c^{7}_{55}
+c^{8}_{55}
+c^{13}_{55}
+c^{18}_{55}
$,
$ -\xi_{11}^{3}$,
$ -\xi_{11}^{7}$,
$ -\xi_{11}^{5}$;\ \ 
$ -\xi_{11}^{7}$,
$ \xi_{11}^{5}$,
$ -c^{2}_{55}
-c^{3}_{55}
-c^{8}_{55}
-c^{13}_{55}
$,
$ -\xi_{11}^{3}$,
$ 1$,
$ c^{2}_{55}
+c^{3}_{55}
+c^{7}_{55}
+c^{8}_{55}
+c^{13}_{55}
+c^{18}_{55}
$,
$ -1-c^{5}_{55}
-c^{6}_{55}
-c^{11}_{55}
-c^{16}_{55}
$,
$ \frac{1+\sqrt{5}}{2}$;\ \ 
$ -c_{11}^{5}$,
$ c^{2}_{55}
+c^{3}_{55}
+c^{7}_{55}
+c^{8}_{55}
+c^{13}_{55}
+c^{18}_{55}
$,
$ -1$,
$ -\xi_{11}^{7}$,
$ c^{3}_{55}
+c^{8}_{55}
$,
$ -\frac{1+\sqrt{5}}{2}$,
$ -c^{2}_{55}
-c^{3}_{55}
-c^{8}_{55}
-c^{13}_{55}
$;\ \ 
$ \xi_{11}^{7}$,
$ -1-c^{5}_{55}
-c^{6}_{55}
-c^{11}_{55}
-c^{16}_{55}
$,
$ \frac{1+\sqrt{5}}{2}$,
$ -\xi_{11}^{5}$,
$ \xi_{11}^{3}$,
$ -1$;\ \ 
$ \xi_{11}^{5}$,
$ c_{11}^{5}$,
$ -\frac{1+\sqrt{5}}{2}$,
$ c^{2}_{55}
+c^{3}_{55}
+c^{7}_{55}
+c^{8}_{55}
+c^{13}_{55}
+c^{18}_{55}
$,
$ -c^{3}_{55}
-c^{8}_{55}
$;\ \ 
$ \xi_{11}^{3}$,
$ -c^{2}_{55}
-c^{3}_{55}
-c^{8}_{55}
-c^{13}_{55}
$,
$ -c^{3}_{55}
-c^{8}_{55}
$,
$ 1+c^{5}_{55}
+c^{6}_{55}
+c^{11}_{55}
+c^{16}_{55}
$;\ \ 
$ c_{11}^{5}$,
$ 1$,
$ \xi_{11}^{7}$;\ \ 
$ -\xi_{11}^{5}$,
$ -c_{11}^{5}$;\ \ 
$ -\xi_{11}^{3}$)

Factors = $2_{\frac{26}{5},3.618}^{5,720}\boxtimes 5_{\frac{16}{11},34.64}^{11,640}$

  \vskip 2ex

\noindent64. $10_{\frac{8}{35},127.8}^{35,427}$ \irep{1074}:\ \ 
$d_i$ = ($1.0$,
$1.618$,
$2.246$,
$2.246$,
$2.801$,
$3.635$,
$3.635$,
$4.48$,
$4.533$,
$6.551$) 

\vskip 0.7ex
\hangindent=3em \hangafter=1
$D^2= 127.870 = 
49+14c^{1}_{35}
+7c^{4}_{35}
+28c^{5}_{35}
+14c^{6}_{35}
+7c^{7}_{35}
+14c^{10}_{35}
+7c^{11}_{35}
$

\vskip 0.7ex
\hangindent=3em \hangafter=1
$T = ( 0,
\frac{2}{5},
\frac{1}{7},
\frac{1}{7},
\frac{6}{7},
\frac{19}{35},
\frac{19}{35},
\frac{4}{7},
\frac{9}{35},
\frac{34}{35} )
$,

\vskip 0.7ex
\hangindent=3em \hangafter=1
$S$ = ($ 1$,
$ \frac{1+\sqrt{5}}{2}$,
$ \xi_{7}^{3}$,
$ \xi_{7}^{3}$,
$ 2+c^{1}_{7}
+c^{2}_{7}
$,
$ c^{1}_{35}
+c^{4}_{35}
+c^{6}_{35}
+c^{11}_{35}
$,
$ c^{1}_{35}
+c^{4}_{35}
+c^{6}_{35}
+c^{11}_{35}
$,
$ 2+2c^{1}_{7}
+c^{2}_{7}
$,
$ 1+c^{1}_{35}
+c^{6}_{35}
+c^{7}_{35}
$,
$ 2c^{1}_{35}
+c^{4}_{35}
+2c^{6}_{35}
+c^{11}_{35}
$;\ \ 
$ -1$,
$ c^{1}_{35}
+c^{4}_{35}
+c^{6}_{35}
+c^{11}_{35}
$,
$ c^{1}_{35}
+c^{4}_{35}
+c^{6}_{35}
+c^{11}_{35}
$,
$ 1+c^{1}_{35}
+c^{6}_{35}
+c^{7}_{35}
$,
$ -\xi_{7}^{3}$,
$ -\xi_{7}^{3}$,
$ 2c^{1}_{35}
+c^{4}_{35}
+2c^{6}_{35}
+c^{11}_{35}
$,
$ -2-c^{1}_{7}
-c^{2}_{7}
$,
$ -2-2  c^{1}_{7}
-c^{2}_{7}
$;\ \ 
$ s^{1}_{7}
+\zeta^{2}_{7}
+\zeta^{3}_{7}
$,
$ -1-2  \zeta^{1}_{7}
-\zeta^{2}_{7}
-\zeta^{3}_{7}
$,
$ -\xi_{7}^{3}$,
$ 1-2  \zeta^{1}_{35}
-\zeta^{-1}_{35}
+2\zeta^{-2}_{35}
-\zeta^{4}_{35}
-2  \zeta^{-4}_{35}
+2\zeta^{5}_{35}
-\zeta^{6}_{35}
-2  \zeta^{-6}_{35}
+c^{7}_{35}
-\zeta^{11}_{35}
-2  \zeta^{-11}_{35}
+2\zeta^{12}_{35}
$,
$ -1+\zeta^{1}_{35}
-2  \zeta^{-2}_{35}
+\zeta^{-4}_{35}
-2  \zeta^{5}_{35}
+\zeta^{-6}_{35}
-c^{7}_{35}
+\zeta^{-11}_{35}
-2  \zeta^{12}_{35}
$,
$ \xi_{7}^{3}$,
$ -c^{1}_{35}
-c^{4}_{35}
-c^{6}_{35}
-c^{11}_{35}
$,
$ c^{1}_{35}
+c^{4}_{35}
+c^{6}_{35}
+c^{11}_{35}
$;\ \ 
$ s^{1}_{7}
+\zeta^{2}_{7}
+\zeta^{3}_{7}
$,
$ -\xi_{7}^{3}$,
$ -1+\zeta^{1}_{35}
-2  \zeta^{-2}_{35}
+\zeta^{-4}_{35}
-2  \zeta^{5}_{35}
+\zeta^{-6}_{35}
-c^{7}_{35}
+\zeta^{-11}_{35}
-2  \zeta^{12}_{35}
$,
$ 1-2  \zeta^{1}_{35}
-\zeta^{-1}_{35}
+2\zeta^{-2}_{35}
-\zeta^{4}_{35}
-2  \zeta^{-4}_{35}
+2\zeta^{5}_{35}
-\zeta^{6}_{35}
-2  \zeta^{-6}_{35}
+c^{7}_{35}
-\zeta^{11}_{35}
-2  \zeta^{-11}_{35}
+2\zeta^{12}_{35}
$,
$ \xi_{7}^{3}$,
$ -c^{1}_{35}
-c^{4}_{35}
-c^{6}_{35}
-c^{11}_{35}
$,
$ c^{1}_{35}
+c^{4}_{35}
+c^{6}_{35}
+c^{11}_{35}
$;\ \ 
$ 2+2c^{1}_{7}
+c^{2}_{7}
$,
$ -c^{1}_{35}
-c^{4}_{35}
-c^{6}_{35}
-c^{11}_{35}
$,
$ -c^{1}_{35}
-c^{4}_{35}
-c^{6}_{35}
-c^{11}_{35}
$,
$ -1$,
$ 2c^{1}_{35}
+c^{4}_{35}
+2c^{6}_{35}
+c^{11}_{35}
$,
$ -\frac{1+\sqrt{5}}{2}$;\ \ 
$ -s^{1}_{7}
-\zeta^{2}_{7}
-\zeta^{3}_{7}
$,
$ 1+2\zeta^{1}_{7}
+\zeta^{2}_{7}
+\zeta^{3}_{7}
$,
$ c^{1}_{35}
+c^{4}_{35}
+c^{6}_{35}
+c^{11}_{35}
$,
$ \xi_{7}^{3}$,
$ -\xi_{7}^{3}$;\ \ 
$ -s^{1}_{7}
-\zeta^{2}_{7}
-\zeta^{3}_{7}
$,
$ c^{1}_{35}
+c^{4}_{35}
+c^{6}_{35}
+c^{11}_{35}
$,
$ \xi_{7}^{3}$,
$ -\xi_{7}^{3}$;\ \ 
$ -2-c^{1}_{7}
-c^{2}_{7}
$,
$ -\frac{1+\sqrt{5}}{2}$,
$ -1-c^{1}_{35}
-c^{6}_{35}
-c^{7}_{35}
$;\ \ 
$ -2-2  c^{1}_{7}
-c^{2}_{7}
$,
$ 1$;\ \ 
$ 2+c^{1}_{7}
+c^{2}_{7}
$)

Factors = $2_{\frac{14}{5},3.618}^{5,395}\boxtimes 5_{\frac{38}{7},35.34}^{7,386}$

  \vskip 2ex

\noindent65. $10_{\frac{188}{35},127.8}^{35,259}$ \irep{1074}:\ \ 
$d_i$ = ($1.0$,
$1.618$,
$2.246$,
$2.246$,
$2.801$,
$3.635$,
$3.635$,
$4.48$,
$4.533$,
$6.551$) 

\vskip 0.7ex
\hangindent=3em \hangafter=1
$D^2= 127.870 = 
49+14c^{1}_{35}
+7c^{4}_{35}
+28c^{5}_{35}
+14c^{6}_{35}
+7c^{7}_{35}
+14c^{10}_{35}
+7c^{11}_{35}
$

\vskip 0.7ex
\hangindent=3em \hangafter=1
$T = ( 0,
\frac{2}{5},
\frac{6}{7},
\frac{6}{7},
\frac{1}{7},
\frac{9}{35},
\frac{9}{35},
\frac{3}{7},
\frac{19}{35},
\frac{29}{35} )
$,

\vskip 0.7ex
\hangindent=3em \hangafter=1
$S$ = ($ 1$,
$ \frac{1+\sqrt{5}}{2}$,
$ \xi_{7}^{3}$,
$ \xi_{7}^{3}$,
$ 2+c^{1}_{7}
+c^{2}_{7}
$,
$ c^{1}_{35}
+c^{4}_{35}
+c^{6}_{35}
+c^{11}_{35}
$,
$ c^{1}_{35}
+c^{4}_{35}
+c^{6}_{35}
+c^{11}_{35}
$,
$ 2+2c^{1}_{7}
+c^{2}_{7}
$,
$ 1+c^{1}_{35}
+c^{6}_{35}
+c^{7}_{35}
$,
$ 2c^{1}_{35}
+c^{4}_{35}
+2c^{6}_{35}
+c^{11}_{35}
$;\ \ 
$ -1$,
$ c^{1}_{35}
+c^{4}_{35}
+c^{6}_{35}
+c^{11}_{35}
$,
$ c^{1}_{35}
+c^{4}_{35}
+c^{6}_{35}
+c^{11}_{35}
$,
$ 1+c^{1}_{35}
+c^{6}_{35}
+c^{7}_{35}
$,
$ -\xi_{7}^{3}$,
$ -\xi_{7}^{3}$,
$ 2c^{1}_{35}
+c^{4}_{35}
+2c^{6}_{35}
+c^{11}_{35}
$,
$ -2-c^{1}_{7}
-c^{2}_{7}
$,
$ -2-2  c^{1}_{7}
-c^{2}_{7}
$;\ \ 
$ -1-2  \zeta^{1}_{7}
-\zeta^{2}_{7}
-\zeta^{3}_{7}
$,
$ s^{1}_{7}
+\zeta^{2}_{7}
+\zeta^{3}_{7}
$,
$ -\xi_{7}^{3}$,
$ 1-2  \zeta^{1}_{35}
-\zeta^{-1}_{35}
+2\zeta^{-2}_{35}
-\zeta^{4}_{35}
-2  \zeta^{-4}_{35}
+2\zeta^{5}_{35}
-\zeta^{6}_{35}
-2  \zeta^{-6}_{35}
+c^{7}_{35}
-\zeta^{11}_{35}
-2  \zeta^{-11}_{35}
+2\zeta^{12}_{35}
$,
$ -1+\zeta^{1}_{35}
-2  \zeta^{-2}_{35}
+\zeta^{-4}_{35}
-2  \zeta^{5}_{35}
+\zeta^{-6}_{35}
-c^{7}_{35}
+\zeta^{-11}_{35}
-2  \zeta^{12}_{35}
$,
$ \xi_{7}^{3}$,
$ -c^{1}_{35}
-c^{4}_{35}
-c^{6}_{35}
-c^{11}_{35}
$,
$ c^{1}_{35}
+c^{4}_{35}
+c^{6}_{35}
+c^{11}_{35}
$;\ \ 
$ -1-2  \zeta^{1}_{7}
-\zeta^{2}_{7}
-\zeta^{3}_{7}
$,
$ -\xi_{7}^{3}$,
$ -1+\zeta^{1}_{35}
-2  \zeta^{-2}_{35}
+\zeta^{-4}_{35}
-2  \zeta^{5}_{35}
+\zeta^{-6}_{35}
-c^{7}_{35}
+\zeta^{-11}_{35}
-2  \zeta^{12}_{35}
$,
$ 1-2  \zeta^{1}_{35}
-\zeta^{-1}_{35}
+2\zeta^{-2}_{35}
-\zeta^{4}_{35}
-2  \zeta^{-4}_{35}
+2\zeta^{5}_{35}
-\zeta^{6}_{35}
-2  \zeta^{-6}_{35}
+c^{7}_{35}
-\zeta^{11}_{35}
-2  \zeta^{-11}_{35}
+2\zeta^{12}_{35}
$,
$ \xi_{7}^{3}$,
$ -c^{1}_{35}
-c^{4}_{35}
-c^{6}_{35}
-c^{11}_{35}
$,
$ c^{1}_{35}
+c^{4}_{35}
+c^{6}_{35}
+c^{11}_{35}
$;\ \ 
$ 2+2c^{1}_{7}
+c^{2}_{7}
$,
$ -c^{1}_{35}
-c^{4}_{35}
-c^{6}_{35}
-c^{11}_{35}
$,
$ -c^{1}_{35}
-c^{4}_{35}
-c^{6}_{35}
-c^{11}_{35}
$,
$ -1$,
$ 2c^{1}_{35}
+c^{4}_{35}
+2c^{6}_{35}
+c^{11}_{35}
$,
$ -\frac{1+\sqrt{5}}{2}$;\ \ 
$ 1+2\zeta^{1}_{7}
+\zeta^{2}_{7}
+\zeta^{3}_{7}
$,
$ -s^{1}_{7}
-\zeta^{2}_{7}
-\zeta^{3}_{7}
$,
$ c^{1}_{35}
+c^{4}_{35}
+c^{6}_{35}
+c^{11}_{35}
$,
$ \xi_{7}^{3}$,
$ -\xi_{7}^{3}$;\ \ 
$ 1+2\zeta^{1}_{7}
+\zeta^{2}_{7}
+\zeta^{3}_{7}
$,
$ c^{1}_{35}
+c^{4}_{35}
+c^{6}_{35}
+c^{11}_{35}
$,
$ \xi_{7}^{3}$,
$ -\xi_{7}^{3}$;\ \ 
$ -2-c^{1}_{7}
-c^{2}_{7}
$,
$ -\frac{1+\sqrt{5}}{2}$,
$ -1-c^{1}_{35}
-c^{6}_{35}
-c^{7}_{35}
$;\ \ 
$ -2-2  c^{1}_{7}
-c^{2}_{7}
$,
$ 1$;\ \ 
$ 2+c^{1}_{7}
+c^{2}_{7}
$)

Factors = $2_{\frac{14}{5},3.618}^{5,395}\boxtimes 5_{\frac{18}{7},35.34}^{7,101}$

  \vskip 2ex

\noindent66. $10_{\frac{92}{35},127.8}^{35,112}$ \irep{1074}:\ \ 
$d_i$ = ($1.0$,
$1.618$,
$2.246$,
$2.246$,
$2.801$,
$3.635$,
$3.635$,
$4.48$,
$4.533$,
$6.551$) 

\vskip 0.7ex
\hangindent=3em \hangafter=1
$D^2= 127.870 = 
49+14c^{1}_{35}
+7c^{4}_{35}
+28c^{5}_{35}
+14c^{6}_{35}
+7c^{7}_{35}
+14c^{10}_{35}
+7c^{11}_{35}
$

\vskip 0.7ex
\hangindent=3em \hangafter=1
$T = ( 0,
\frac{3}{5},
\frac{1}{7},
\frac{1}{7},
\frac{6}{7},
\frac{26}{35},
\frac{26}{35},
\frac{4}{7},
\frac{16}{35},
\frac{6}{35} )
$,

\vskip 0.7ex
\hangindent=3em \hangafter=1
$S$ = ($ 1$,
$ \frac{1+\sqrt{5}}{2}$,
$ \xi_{7}^{3}$,
$ \xi_{7}^{3}$,
$ 2+c^{1}_{7}
+c^{2}_{7}
$,
$ c^{1}_{35}
+c^{4}_{35}
+c^{6}_{35}
+c^{11}_{35}
$,
$ c^{1}_{35}
+c^{4}_{35}
+c^{6}_{35}
+c^{11}_{35}
$,
$ 2+2c^{1}_{7}
+c^{2}_{7}
$,
$ 1+c^{1}_{35}
+c^{6}_{35}
+c^{7}_{35}
$,
$ 2c^{1}_{35}
+c^{4}_{35}
+2c^{6}_{35}
+c^{11}_{35}
$;\ \ 
$ -1$,
$ c^{1}_{35}
+c^{4}_{35}
+c^{6}_{35}
+c^{11}_{35}
$,
$ c^{1}_{35}
+c^{4}_{35}
+c^{6}_{35}
+c^{11}_{35}
$,
$ 1+c^{1}_{35}
+c^{6}_{35}
+c^{7}_{35}
$,
$ -\xi_{7}^{3}$,
$ -\xi_{7}^{3}$,
$ 2c^{1}_{35}
+c^{4}_{35}
+2c^{6}_{35}
+c^{11}_{35}
$,
$ -2-c^{1}_{7}
-c^{2}_{7}
$,
$ -2-2  c^{1}_{7}
-c^{2}_{7}
$;\ \ 
$ s^{1}_{7}
+\zeta^{2}_{7}
+\zeta^{3}_{7}
$,
$ -1-2  \zeta^{1}_{7}
-\zeta^{2}_{7}
-\zeta^{3}_{7}
$,
$ -\xi_{7}^{3}$,
$ 1-2  \zeta^{1}_{35}
-\zeta^{-1}_{35}
+2\zeta^{-2}_{35}
-\zeta^{4}_{35}
-2  \zeta^{-4}_{35}
+2\zeta^{5}_{35}
-\zeta^{6}_{35}
-2  \zeta^{-6}_{35}
+c^{7}_{35}
-\zeta^{11}_{35}
-2  \zeta^{-11}_{35}
+2\zeta^{12}_{35}
$,
$ -1+\zeta^{1}_{35}
-2  \zeta^{-2}_{35}
+\zeta^{-4}_{35}
-2  \zeta^{5}_{35}
+\zeta^{-6}_{35}
-c^{7}_{35}
+\zeta^{-11}_{35}
-2  \zeta^{12}_{35}
$,
$ \xi_{7}^{3}$,
$ -c^{1}_{35}
-c^{4}_{35}
-c^{6}_{35}
-c^{11}_{35}
$,
$ c^{1}_{35}
+c^{4}_{35}
+c^{6}_{35}
+c^{11}_{35}
$;\ \ 
$ s^{1}_{7}
+\zeta^{2}_{7}
+\zeta^{3}_{7}
$,
$ -\xi_{7}^{3}$,
$ -1+\zeta^{1}_{35}
-2  \zeta^{-2}_{35}
+\zeta^{-4}_{35}
-2  \zeta^{5}_{35}
+\zeta^{-6}_{35}
-c^{7}_{35}
+\zeta^{-11}_{35}
-2  \zeta^{12}_{35}
$,
$ 1-2  \zeta^{1}_{35}
-\zeta^{-1}_{35}
+2\zeta^{-2}_{35}
-\zeta^{4}_{35}
-2  \zeta^{-4}_{35}
+2\zeta^{5}_{35}
-\zeta^{6}_{35}
-2  \zeta^{-6}_{35}
+c^{7}_{35}
-\zeta^{11}_{35}
-2  \zeta^{-11}_{35}
+2\zeta^{12}_{35}
$,
$ \xi_{7}^{3}$,
$ -c^{1}_{35}
-c^{4}_{35}
-c^{6}_{35}
-c^{11}_{35}
$,
$ c^{1}_{35}
+c^{4}_{35}
+c^{6}_{35}
+c^{11}_{35}
$;\ \ 
$ 2+2c^{1}_{7}
+c^{2}_{7}
$,
$ -c^{1}_{35}
-c^{4}_{35}
-c^{6}_{35}
-c^{11}_{35}
$,
$ -c^{1}_{35}
-c^{4}_{35}
-c^{6}_{35}
-c^{11}_{35}
$,
$ -1$,
$ 2c^{1}_{35}
+c^{4}_{35}
+2c^{6}_{35}
+c^{11}_{35}
$,
$ -\frac{1+\sqrt{5}}{2}$;\ \ 
$ -s^{1}_{7}
-\zeta^{2}_{7}
-\zeta^{3}_{7}
$,
$ 1+2\zeta^{1}_{7}
+\zeta^{2}_{7}
+\zeta^{3}_{7}
$,
$ c^{1}_{35}
+c^{4}_{35}
+c^{6}_{35}
+c^{11}_{35}
$,
$ \xi_{7}^{3}$,
$ -\xi_{7}^{3}$;\ \ 
$ -s^{1}_{7}
-\zeta^{2}_{7}
-\zeta^{3}_{7}
$,
$ c^{1}_{35}
+c^{4}_{35}
+c^{6}_{35}
+c^{11}_{35}
$,
$ \xi_{7}^{3}$,
$ -\xi_{7}^{3}$;\ \ 
$ -2-c^{1}_{7}
-c^{2}_{7}
$,
$ -\frac{1+\sqrt{5}}{2}$,
$ -1-c^{1}_{35}
-c^{6}_{35}
-c^{7}_{35}
$;\ \ 
$ -2-2  c^{1}_{7}
-c^{2}_{7}
$,
$ 1$;\ \ 
$ 2+c^{1}_{7}
+c^{2}_{7}
$)

Factors = $2_{\frac{26}{5},3.618}^{5,720}\boxtimes 5_{\frac{38}{7},35.34}^{7,386}$

  \vskip 2ex

\noindent67. $10_{\frac{272}{35},127.8}^{35,631}$ \irep{1074}:\ \ 
$d_i$ = ($1.0$,
$1.618$,
$2.246$,
$2.246$,
$2.801$,
$3.635$,
$3.635$,
$4.48$,
$4.533$,
$6.551$) 

\vskip 0.7ex
\hangindent=3em \hangafter=1
$D^2= 127.870 = 
49+14c^{1}_{35}
+7c^{4}_{35}
+28c^{5}_{35}
+14c^{6}_{35}
+7c^{7}_{35}
+14c^{10}_{35}
+7c^{11}_{35}
$

\vskip 0.7ex
\hangindent=3em \hangafter=1
$T = ( 0,
\frac{3}{5},
\frac{6}{7},
\frac{6}{7},
\frac{1}{7},
\frac{16}{35},
\frac{16}{35},
\frac{3}{7},
\frac{26}{35},
\frac{1}{35} )
$,

\vskip 0.7ex
\hangindent=3em \hangafter=1
$S$ = ($ 1$,
$ \frac{1+\sqrt{5}}{2}$,
$ \xi_{7}^{3}$,
$ \xi_{7}^{3}$,
$ 2+c^{1}_{7}
+c^{2}_{7}
$,
$ c^{1}_{35}
+c^{4}_{35}
+c^{6}_{35}
+c^{11}_{35}
$,
$ c^{1}_{35}
+c^{4}_{35}
+c^{6}_{35}
+c^{11}_{35}
$,
$ 2+2c^{1}_{7}
+c^{2}_{7}
$,
$ 1+c^{1}_{35}
+c^{6}_{35}
+c^{7}_{35}
$,
$ 2c^{1}_{35}
+c^{4}_{35}
+2c^{6}_{35}
+c^{11}_{35}
$;\ \ 
$ -1$,
$ c^{1}_{35}
+c^{4}_{35}
+c^{6}_{35}
+c^{11}_{35}
$,
$ c^{1}_{35}
+c^{4}_{35}
+c^{6}_{35}
+c^{11}_{35}
$,
$ 1+c^{1}_{35}
+c^{6}_{35}
+c^{7}_{35}
$,
$ -\xi_{7}^{3}$,
$ -\xi_{7}^{3}$,
$ 2c^{1}_{35}
+c^{4}_{35}
+2c^{6}_{35}
+c^{11}_{35}
$,
$ -2-c^{1}_{7}
-c^{2}_{7}
$,
$ -2-2  c^{1}_{7}
-c^{2}_{7}
$;\ \ 
$ -1-2  \zeta^{1}_{7}
-\zeta^{2}_{7}
-\zeta^{3}_{7}
$,
$ s^{1}_{7}
+\zeta^{2}_{7}
+\zeta^{3}_{7}
$,
$ -\xi_{7}^{3}$,
$ 1-2  \zeta^{1}_{35}
-\zeta^{-1}_{35}
+2\zeta^{-2}_{35}
-\zeta^{4}_{35}
-2  \zeta^{-4}_{35}
+2\zeta^{5}_{35}
-\zeta^{6}_{35}
-2  \zeta^{-6}_{35}
+c^{7}_{35}
-\zeta^{11}_{35}
-2  \zeta^{-11}_{35}
+2\zeta^{12}_{35}
$,
$ -1+\zeta^{1}_{35}
-2  \zeta^{-2}_{35}
+\zeta^{-4}_{35}
-2  \zeta^{5}_{35}
+\zeta^{-6}_{35}
-c^{7}_{35}
+\zeta^{-11}_{35}
-2  \zeta^{12}_{35}
$,
$ \xi_{7}^{3}$,
$ -c^{1}_{35}
-c^{4}_{35}
-c^{6}_{35}
-c^{11}_{35}
$,
$ c^{1}_{35}
+c^{4}_{35}
+c^{6}_{35}
+c^{11}_{35}
$;\ \ 
$ -1-2  \zeta^{1}_{7}
-\zeta^{2}_{7}
-\zeta^{3}_{7}
$,
$ -\xi_{7}^{3}$,
$ -1+\zeta^{1}_{35}
-2  \zeta^{-2}_{35}
+\zeta^{-4}_{35}
-2  \zeta^{5}_{35}
+\zeta^{-6}_{35}
-c^{7}_{35}
+\zeta^{-11}_{35}
-2  \zeta^{12}_{35}
$,
$ 1-2  \zeta^{1}_{35}
-\zeta^{-1}_{35}
+2\zeta^{-2}_{35}
-\zeta^{4}_{35}
-2  \zeta^{-4}_{35}
+2\zeta^{5}_{35}
-\zeta^{6}_{35}
-2  \zeta^{-6}_{35}
+c^{7}_{35}
-\zeta^{11}_{35}
-2  \zeta^{-11}_{35}
+2\zeta^{12}_{35}
$,
$ \xi_{7}^{3}$,
$ -c^{1}_{35}
-c^{4}_{35}
-c^{6}_{35}
-c^{11}_{35}
$,
$ c^{1}_{35}
+c^{4}_{35}
+c^{6}_{35}
+c^{11}_{35}
$;\ \ 
$ 2+2c^{1}_{7}
+c^{2}_{7}
$,
$ -c^{1}_{35}
-c^{4}_{35}
-c^{6}_{35}
-c^{11}_{35}
$,
$ -c^{1}_{35}
-c^{4}_{35}
-c^{6}_{35}
-c^{11}_{35}
$,
$ -1$,
$ 2c^{1}_{35}
+c^{4}_{35}
+2c^{6}_{35}
+c^{11}_{35}
$,
$ -\frac{1+\sqrt{5}}{2}$;\ \ 
$ 1+2\zeta^{1}_{7}
+\zeta^{2}_{7}
+\zeta^{3}_{7}
$,
$ -s^{1}_{7}
-\zeta^{2}_{7}
-\zeta^{3}_{7}
$,
$ c^{1}_{35}
+c^{4}_{35}
+c^{6}_{35}
+c^{11}_{35}
$,
$ \xi_{7}^{3}$,
$ -\xi_{7}^{3}$;\ \ 
$ 1+2\zeta^{1}_{7}
+\zeta^{2}_{7}
+\zeta^{3}_{7}
$,
$ c^{1}_{35}
+c^{4}_{35}
+c^{6}_{35}
+c^{11}_{35}
$,
$ \xi_{7}^{3}$,
$ -\xi_{7}^{3}$;\ \ 
$ -2-c^{1}_{7}
-c^{2}_{7}
$,
$ -\frac{1+\sqrt{5}}{2}$,
$ -1-c^{1}_{35}
-c^{6}_{35}
-c^{7}_{35}
$;\ \ 
$ -2-2  c^{1}_{7}
-c^{2}_{7}
$,
$ 1$;\ \ 
$ 2+c^{1}_{7}
+c^{2}_{7}
$)

Factors = $2_{\frac{26}{5},3.618}^{5,720}\boxtimes 5_{\frac{18}{7},35.34}^{7,101}$

  \vskip 2ex

\noindent68. $10_{\frac{26}{7},236.3}^{21,145}$ \irep{960}:\ \ 
$d_i$ = ($1.0$,
$1.977$,
$2.911$,
$3.779$,
$4.563$,
$5.245$,
$5.810$,
$6.245$,
$6.541$,
$6.690$) 

\vskip 0.7ex
\hangindent=3em \hangafter=1
$D^2= 236.341 = 
42+42c^{1}_{21}
+42c^{2}_{21}
+21c^{3}_{21}
+21c^{4}_{21}
+21c^{5}_{21}
$

\vskip 0.7ex
\hangindent=3em \hangafter=1
$T = ( 0,
\frac{2}{7},
\frac{2}{21},
\frac{3}{7},
\frac{2}{7},
\frac{2}{3},
\frac{4}{7},
0,
\frac{20}{21},
\frac{3}{7} )
$,

\vskip 0.7ex
\hangindent=3em \hangafter=1
$S$ = ($ 1$,
$ -c_{21}^{10}$,
$ \xi_{21}^{3}$,
$ \xi_{21}^{17}$,
$ \xi_{21}^{5}$,
$ \xi_{21}^{15}$,
$ \xi_{21}^{7}$,
$ \xi_{21}^{13}$,
$ \xi_{21}^{9}$,
$ \xi_{21}^{11}$;\ \ 
$ -\xi_{21}^{17}$,
$ \xi_{21}^{15}$,
$ -\xi_{21}^{13}$,
$ \xi_{21}^{11}$,
$ -\xi_{21}^{9}$,
$ \xi_{21}^{7}$,
$ -\xi_{21}^{5}$,
$ \xi_{21}^{3}$,
$ -1$;\ \ 
$ \xi_{21}^{9}$,
$ \xi_{21}^{9}$,
$ \xi_{21}^{15}$,
$ \xi_{21}^{3}$,
$0$,
$ -\xi_{21}^{3}$,
$ -\xi_{21}^{15}$,
$ -\xi_{21}^{9}$;\ \ 
$ -\xi_{21}^{5}$,
$ 1$,
$ \xi_{21}^{3}$,
$ -\xi_{21}^{7}$,
$ \xi_{21}^{11}$,
$ -\xi_{21}^{15}$,
$ -c_{21}^{10}$;\ \ 
$ -\xi_{21}^{17}$,
$ -\xi_{21}^{9}$,
$ -\xi_{21}^{7}$,
$ c_{21}^{10}$,
$ \xi_{21}^{3}$,
$ \xi_{21}^{13}$;\ \ 
$ \xi_{21}^{15}$,
$0$,
$ -\xi_{21}^{15}$,
$ \xi_{21}^{9}$,
$ -\xi_{21}^{3}$;\ \ 
$ \xi_{21}^{7}$,
$ \xi_{21}^{7}$,
$0$,
$ -\xi_{21}^{7}$;\ \ 
$ 1$,
$ -\xi_{21}^{9}$,
$ \xi_{21}^{17}$;\ \ 
$ -\xi_{21}^{3}$,
$ \xi_{21}^{15}$;\ \ 
$ -\xi_{21}^{5}$)

  \vskip 2ex

\noindent69. $10_{\frac{30}{7},236.3}^{21,387}$ \irep{960}:\ \ 
$d_i$ = ($1.0$,
$1.977$,
$2.911$,
$3.779$,
$4.563$,
$5.245$,
$5.810$,
$6.245$,
$6.541$,
$6.690$) 

\vskip 0.7ex
\hangindent=3em \hangafter=1
$D^2= 236.341 = 
42+42c^{1}_{21}
+42c^{2}_{21}
+21c^{3}_{21}
+21c^{4}_{21}
+21c^{5}_{21}
$

\vskip 0.7ex
\hangindent=3em \hangafter=1
$T = ( 0,
\frac{5}{7},
\frac{19}{21},
\frac{4}{7},
\frac{5}{7},
\frac{1}{3},
\frac{3}{7},
0,
\frac{1}{21},
\frac{4}{7} )
$,

\vskip 0.7ex
\hangindent=3em \hangafter=1
$S$ = ($ 1$,
$ -c_{21}^{10}$,
$ \xi_{21}^{3}$,
$ \xi_{21}^{17}$,
$ \xi_{21}^{5}$,
$ \xi_{21}^{15}$,
$ \xi_{21}^{7}$,
$ \xi_{21}^{13}$,
$ \xi_{21}^{9}$,
$ \xi_{21}^{11}$;\ \ 
$ -\xi_{21}^{17}$,
$ \xi_{21}^{15}$,
$ -\xi_{21}^{13}$,
$ \xi_{21}^{11}$,
$ -\xi_{21}^{9}$,
$ \xi_{21}^{7}$,
$ -\xi_{21}^{5}$,
$ \xi_{21}^{3}$,
$ -1$;\ \ 
$ \xi_{21}^{9}$,
$ \xi_{21}^{9}$,
$ \xi_{21}^{15}$,
$ \xi_{21}^{3}$,
$0$,
$ -\xi_{21}^{3}$,
$ -\xi_{21}^{15}$,
$ -\xi_{21}^{9}$;\ \ 
$ -\xi_{21}^{5}$,
$ 1$,
$ \xi_{21}^{3}$,
$ -\xi_{21}^{7}$,
$ \xi_{21}^{11}$,
$ -\xi_{21}^{15}$,
$ -c_{21}^{10}$;\ \ 
$ -\xi_{21}^{17}$,
$ -\xi_{21}^{9}$,
$ -\xi_{21}^{7}$,
$ c_{21}^{10}$,
$ \xi_{21}^{3}$,
$ \xi_{21}^{13}$;\ \ 
$ \xi_{21}^{15}$,
$0$,
$ -\xi_{21}^{15}$,
$ \xi_{21}^{9}$,
$ -\xi_{21}^{3}$;\ \ 
$ \xi_{21}^{7}$,
$ \xi_{21}^{7}$,
$0$,
$ -\xi_{21}^{7}$;\ \ 
$ 1$,
$ -\xi_{21}^{9}$,
$ \xi_{21}^{17}$;\ \ 
$ -\xi_{21}^{3}$,
$ \xi_{21}^{15}$;\ \ 
$ -\xi_{21}^{5}$)

  \vskip 2ex

\noindent70. $10_{\frac{48}{17},499.2}^{17,522}$ \irep{888}:\ \ 
$d_i$ = ($1.0$,
$2.965$,
$4.830$,
$5.418$,
$5.418$,
$6.531$,
$8.9$,
$9.214$,
$10.106$,
$10.653$) 

\vskip 0.7ex
\hangindent=3em \hangafter=1
$D^2= 499.210 = 
136+119c^{1}_{17}
+102c^{2}_{17}
+85c^{3}_{17}
+68c^{4}_{17}
+51c^{5}_{17}
+34c^{6}_{17}
+17c^{7}_{17}
$

\vskip 0.7ex
\hangindent=3em \hangafter=1
$T = ( 0,
\frac{1}{17},
\frac{3}{17},
\frac{2}{17},
\frac{2}{17},
\frac{6}{17},
\frac{10}{17},
\frac{15}{17},
\frac{4}{17},
\frac{11}{17} )
$,

\vskip 0.7ex
\hangindent=3em \hangafter=1
$S$ = ($ 1$,
$ 2+c^{1}_{17}
+c^{2}_{17}
+c^{3}_{17}
+c^{4}_{17}
+c^{5}_{17}
+c^{6}_{17}
+c^{7}_{17}
$,
$ 2+2c^{1}_{17}
+c^{2}_{17}
+c^{3}_{17}
+c^{4}_{17}
+c^{5}_{17}
+c^{6}_{17}
+c^{7}_{17}
$,
$ \xi_{17}^{9}$,
$ \xi_{17}^{9}$,
$ 2+2c^{1}_{17}
+c^{2}_{17}
+c^{3}_{17}
+c^{4}_{17}
+c^{5}_{17}
+c^{6}_{17}
$,
$ 2+2c^{1}_{17}
+2c^{2}_{17}
+c^{3}_{17}
+c^{4}_{17}
+c^{5}_{17}
+c^{6}_{17}
$,
$ 2+2c^{1}_{17}
+2c^{2}_{17}
+c^{3}_{17}
+c^{4}_{17}
+c^{5}_{17}
$,
$ 2+2c^{1}_{17}
+2c^{2}_{17}
+2c^{3}_{17}
+c^{4}_{17}
+c^{5}_{17}
$,
$ 2+2c^{1}_{17}
+2c^{2}_{17}
+2c^{3}_{17}
+c^{4}_{17}
$;\ \ 
$ 2+2c^{1}_{17}
+2c^{2}_{17}
+c^{3}_{17}
+c^{4}_{17}
+c^{5}_{17}
+c^{6}_{17}
$,
$ 2+2c^{1}_{17}
+2c^{2}_{17}
+2c^{3}_{17}
+c^{4}_{17}
$,
$ -\xi_{17}^{9}$,
$ -\xi_{17}^{9}$,
$ 2+2c^{1}_{17}
+2c^{2}_{17}
+2c^{3}_{17}
+c^{4}_{17}
+c^{5}_{17}
$,
$ 2+2c^{1}_{17}
+c^{2}_{17}
+c^{3}_{17}
+c^{4}_{17}
+c^{5}_{17}
+c^{6}_{17}
$,
$ 1$,
$ -2-2  c^{1}_{17}
-c^{2}_{17}
-c^{3}_{17}
-c^{4}_{17}
-c^{5}_{17}
-c^{6}_{17}
-c^{7}_{17}
$,
$ -2-2  c^{1}_{17}
-2  c^{2}_{17}
-c^{3}_{17}
-c^{4}_{17}
-c^{5}_{17}
$;\ \ 
$ 2+2c^{1}_{17}
+2c^{2}_{17}
+c^{3}_{17}
+c^{4}_{17}
+c^{5}_{17}
+c^{6}_{17}
$,
$ \xi_{17}^{9}$,
$ \xi_{17}^{9}$,
$ -1$,
$ -2-2  c^{1}_{17}
-2  c^{2}_{17}
-c^{3}_{17}
-c^{4}_{17}
-c^{5}_{17}
$,
$ -2-2  c^{1}_{17}
-2  c^{2}_{17}
-2  c^{3}_{17}
-c^{4}_{17}
-c^{5}_{17}
$,
$ -2-c^{1}_{17}
-c^{2}_{17}
-c^{3}_{17}
-c^{4}_{17}
-c^{5}_{17}
-c^{6}_{17}
-c^{7}_{17}
$,
$ 2+2c^{1}_{17}
+c^{2}_{17}
+c^{3}_{17}
+c^{4}_{17}
+c^{5}_{17}
+c^{6}_{17}
$;\ \ 
$ 4+3c^{1}_{17}
+3c^{2}_{17}
+2c^{3}_{17}
+2c^{4}_{17}
+2c^{5}_{17}
+c^{6}_{17}
$,
$ -3-2  c^{1}_{17}
-2  c^{2}_{17}
-c^{3}_{17}
-c^{4}_{17}
-2  c^{5}_{17}
-c^{6}_{17}
$,
$ -\xi_{17}^{9}$,
$ \xi_{17}^{9}$,
$ -\xi_{17}^{9}$,
$ \xi_{17}^{9}$,
$ -\xi_{17}^{9}$;\ \ 
$ 4+3c^{1}_{17}
+3c^{2}_{17}
+2c^{3}_{17}
+2c^{4}_{17}
+2c^{5}_{17}
+c^{6}_{17}
$,
$ -\xi_{17}^{9}$,
$ \xi_{17}^{9}$,
$ -\xi_{17}^{9}$,
$ \xi_{17}^{9}$,
$ -\xi_{17}^{9}$;\ \ 
$ -2-2  c^{1}_{17}
-2  c^{2}_{17}
-2  c^{3}_{17}
-c^{4}_{17}
$,
$ -2-2  c^{1}_{17}
-c^{2}_{17}
-c^{3}_{17}
-c^{4}_{17}
-c^{5}_{17}
-c^{6}_{17}
-c^{7}_{17}
$,
$ 2+2c^{1}_{17}
+2c^{2}_{17}
+c^{3}_{17}
+c^{4}_{17}
+c^{5}_{17}
+c^{6}_{17}
$,
$ 2+2c^{1}_{17}
+2c^{2}_{17}
+c^{3}_{17}
+c^{4}_{17}
+c^{5}_{17}
$,
$ -2-c^{1}_{17}
-c^{2}_{17}
-c^{3}_{17}
-c^{4}_{17}
-c^{5}_{17}
-c^{6}_{17}
-c^{7}_{17}
$;\ \ 
$ 2+2c^{1}_{17}
+2c^{2}_{17}
+2c^{3}_{17}
+c^{4}_{17}
+c^{5}_{17}
$,
$ 2+c^{1}_{17}
+c^{2}_{17}
+c^{3}_{17}
+c^{4}_{17}
+c^{5}_{17}
+c^{6}_{17}
+c^{7}_{17}
$,
$ -2-2  c^{1}_{17}
-2  c^{2}_{17}
-2  c^{3}_{17}
-c^{4}_{17}
$,
$ -1$;\ \ 
$ -2-2  c^{1}_{17}
-2  c^{2}_{17}
-2  c^{3}_{17}
-c^{4}_{17}
$,
$ 2+2c^{1}_{17}
+c^{2}_{17}
+c^{3}_{17}
+c^{4}_{17}
+c^{5}_{17}
+c^{6}_{17}
$,
$ 2+2c^{1}_{17}
+c^{2}_{17}
+c^{3}_{17}
+c^{4}_{17}
+c^{5}_{17}
+c^{6}_{17}
+c^{7}_{17}
$;\ \ 
$ 1$,
$ -2-2  c^{1}_{17}
-2  c^{2}_{17}
-c^{3}_{17}
-c^{4}_{17}
-c^{5}_{17}
-c^{6}_{17}
$;\ \ 
$ 2+2c^{1}_{17}
+2c^{2}_{17}
+2c^{3}_{17}
+c^{4}_{17}
+c^{5}_{17}
$)

  \vskip 2ex

\noindent71. $10_{\frac{88}{17},499.2}^{17,976}$ \irep{888}:\ \ 
$d_i$ = ($1.0$,
$2.965$,
$4.830$,
$5.418$,
$5.418$,
$6.531$,
$8.9$,
$9.214$,
$10.106$,
$10.653$) 

\vskip 0.7ex
\hangindent=3em \hangafter=1
$D^2= 499.210 = 
136+119c^{1}_{17}
+102c^{2}_{17}
+85c^{3}_{17}
+68c^{4}_{17}
+51c^{5}_{17}
+34c^{6}_{17}
+17c^{7}_{17}
$

\vskip 0.7ex
\hangindent=3em \hangafter=1
$T = ( 0,
\frac{16}{17},
\frac{14}{17},
\frac{15}{17},
\frac{15}{17},
\frac{11}{17},
\frac{7}{17},
\frac{2}{17},
\frac{13}{17},
\frac{6}{17} )
$,

\vskip 0.7ex
\hangindent=3em \hangafter=1
$S$ = ($ 1$,
$ 2+c^{1}_{17}
+c^{2}_{17}
+c^{3}_{17}
+c^{4}_{17}
+c^{5}_{17}
+c^{6}_{17}
+c^{7}_{17}
$,
$ 2+2c^{1}_{17}
+c^{2}_{17}
+c^{3}_{17}
+c^{4}_{17}
+c^{5}_{17}
+c^{6}_{17}
+c^{7}_{17}
$,
$ \xi_{17}^{9}$,
$ \xi_{17}^{9}$,
$ 2+2c^{1}_{17}
+c^{2}_{17}
+c^{3}_{17}
+c^{4}_{17}
+c^{5}_{17}
+c^{6}_{17}
$,
$ 2+2c^{1}_{17}
+2c^{2}_{17}
+c^{3}_{17}
+c^{4}_{17}
+c^{5}_{17}
+c^{6}_{17}
$,
$ 2+2c^{1}_{17}
+2c^{2}_{17}
+c^{3}_{17}
+c^{4}_{17}
+c^{5}_{17}
$,
$ 2+2c^{1}_{17}
+2c^{2}_{17}
+2c^{3}_{17}
+c^{4}_{17}
+c^{5}_{17}
$,
$ 2+2c^{1}_{17}
+2c^{2}_{17}
+2c^{3}_{17}
+c^{4}_{17}
$;\ \ 
$ 2+2c^{1}_{17}
+2c^{2}_{17}
+c^{3}_{17}
+c^{4}_{17}
+c^{5}_{17}
+c^{6}_{17}
$,
$ 2+2c^{1}_{17}
+2c^{2}_{17}
+2c^{3}_{17}
+c^{4}_{17}
$,
$ -\xi_{17}^{9}$,
$ -\xi_{17}^{9}$,
$ 2+2c^{1}_{17}
+2c^{2}_{17}
+2c^{3}_{17}
+c^{4}_{17}
+c^{5}_{17}
$,
$ 2+2c^{1}_{17}
+c^{2}_{17}
+c^{3}_{17}
+c^{4}_{17}
+c^{5}_{17}
+c^{6}_{17}
$,
$ 1$,
$ -2-2  c^{1}_{17}
-c^{2}_{17}
-c^{3}_{17}
-c^{4}_{17}
-c^{5}_{17}
-c^{6}_{17}
-c^{7}_{17}
$,
$ -2-2  c^{1}_{17}
-2  c^{2}_{17}
-c^{3}_{17}
-c^{4}_{17}
-c^{5}_{17}
$;\ \ 
$ 2+2c^{1}_{17}
+2c^{2}_{17}
+c^{3}_{17}
+c^{4}_{17}
+c^{5}_{17}
+c^{6}_{17}
$,
$ \xi_{17}^{9}$,
$ \xi_{17}^{9}$,
$ -1$,
$ -2-2  c^{1}_{17}
-2  c^{2}_{17}
-c^{3}_{17}
-c^{4}_{17}
-c^{5}_{17}
$,
$ -2-2  c^{1}_{17}
-2  c^{2}_{17}
-2  c^{3}_{17}
-c^{4}_{17}
-c^{5}_{17}
$,
$ -2-c^{1}_{17}
-c^{2}_{17}
-c^{3}_{17}
-c^{4}_{17}
-c^{5}_{17}
-c^{6}_{17}
-c^{7}_{17}
$,
$ 2+2c^{1}_{17}
+c^{2}_{17}
+c^{3}_{17}
+c^{4}_{17}
+c^{5}_{17}
+c^{6}_{17}
$;\ \ 
$ 4+3c^{1}_{17}
+3c^{2}_{17}
+2c^{3}_{17}
+2c^{4}_{17}
+2c^{5}_{17}
+c^{6}_{17}
$,
$ -3-2  c^{1}_{17}
-2  c^{2}_{17}
-c^{3}_{17}
-c^{4}_{17}
-2  c^{5}_{17}
-c^{6}_{17}
$,
$ -\xi_{17}^{9}$,
$ \xi_{17}^{9}$,
$ -\xi_{17}^{9}$,
$ \xi_{17}^{9}$,
$ -\xi_{17}^{9}$;\ \ 
$ 4+3c^{1}_{17}
+3c^{2}_{17}
+2c^{3}_{17}
+2c^{4}_{17}
+2c^{5}_{17}
+c^{6}_{17}
$,
$ -\xi_{17}^{9}$,
$ \xi_{17}^{9}$,
$ -\xi_{17}^{9}$,
$ \xi_{17}^{9}$,
$ -\xi_{17}^{9}$;\ \ 
$ -2-2  c^{1}_{17}
-2  c^{2}_{17}
-2  c^{3}_{17}
-c^{4}_{17}
$,
$ -2-2  c^{1}_{17}
-c^{2}_{17}
-c^{3}_{17}
-c^{4}_{17}
-c^{5}_{17}
-c^{6}_{17}
-c^{7}_{17}
$,
$ 2+2c^{1}_{17}
+2c^{2}_{17}
+c^{3}_{17}
+c^{4}_{17}
+c^{5}_{17}
+c^{6}_{17}
$,
$ 2+2c^{1}_{17}
+2c^{2}_{17}
+c^{3}_{17}
+c^{4}_{17}
+c^{5}_{17}
$,
$ -2-c^{1}_{17}
-c^{2}_{17}
-c^{3}_{17}
-c^{4}_{17}
-c^{5}_{17}
-c^{6}_{17}
-c^{7}_{17}
$;\ \ 
$ 2+2c^{1}_{17}
+2c^{2}_{17}
+2c^{3}_{17}
+c^{4}_{17}
+c^{5}_{17}
$,
$ 2+c^{1}_{17}
+c^{2}_{17}
+c^{3}_{17}
+c^{4}_{17}
+c^{5}_{17}
+c^{6}_{17}
+c^{7}_{17}
$,
$ -2-2  c^{1}_{17}
-2  c^{2}_{17}
-2  c^{3}_{17}
-c^{4}_{17}
$,
$ -1$;\ \ 
$ -2-2  c^{1}_{17}
-2  c^{2}_{17}
-2  c^{3}_{17}
-c^{4}_{17}
$,
$ 2+2c^{1}_{17}
+c^{2}_{17}
+c^{3}_{17}
+c^{4}_{17}
+c^{5}_{17}
+c^{6}_{17}
$,
$ 2+2c^{1}_{17}
+c^{2}_{17}
+c^{3}_{17}
+c^{4}_{17}
+c^{5}_{17}
+c^{6}_{17}
+c^{7}_{17}
$;\ \ 
$ 1$,
$ -2-2  c^{1}_{17}
-2  c^{2}_{17}
-c^{3}_{17}
-c^{4}_{17}
-c^{5}_{17}
-c^{6}_{17}
$;\ \ 
$ 2+2c^{1}_{17}
+2c^{2}_{17}
+2c^{3}_{17}
+c^{4}_{17}
+c^{5}_{17}
$)

  \vskip 2ex

\noindent72. $10_{0,537.4}^{14,352}$ \irep{783}:\ \ 
$d_i$ = ($1.0$,
$3.493$,
$4.493$,
$4.493$,
$5.603$,
$5.603$,
$9.97$,
$10.97$,
$10.97$,
$11.591$) 

\vskip 0.7ex
\hangindent=3em \hangafter=1
$D^2= 537.478 = 
308+224c^{1}_{7}
+112c^{2}_{7}
$

\vskip 0.7ex
\hangindent=3em \hangafter=1
$T = ( 0,
0,
\frac{2}{7},
\frac{5}{7},
\frac{3}{7},
\frac{4}{7},
0,
\frac{1}{7},
\frac{6}{7},
\frac{1}{2} )
$,

\vskip 0.7ex
\hangindent=3em \hangafter=1
$S$ = ($ 1$,
$ 1+2c^{1}_{7}
$,
$ 2\xi_{7}^{3}$,
$ 2\xi_{7}^{3}$,
$ 4+2c^{1}_{7}
+2c^{2}_{7}
$,
$ 4+2c^{1}_{7}
+2c^{2}_{7}
$,
$ 5+4c^{1}_{7}
+2c^{2}_{7}
$,
$ 6+4c^{1}_{7}
+2c^{2}_{7}
$,
$ 6+4c^{1}_{7}
+2c^{2}_{7}
$,
$ 5+6c^{1}_{7}
+2c^{2}_{7}
$;\ \ 
$ -5-4  c^{1}_{7}
-2  c^{2}_{7}
$,
$ -4-2  c^{1}_{7}
-2  c^{2}_{7}
$,
$ -4-2  c^{1}_{7}
-2  c^{2}_{7}
$,
$ 6+4c^{1}_{7}
+2c^{2}_{7}
$,
$ 6+4c^{1}_{7}
+2c^{2}_{7}
$,
$ 1$,
$ 2\xi_{7}^{3}$,
$ 2\xi_{7}^{3}$,
$ -5-6  c^{1}_{7}
-2  c^{2}_{7}
$;\ \ 
$ -2\xi_{7}^{3}$,
$ 6+6c^{1}_{7}
+2c^{2}_{7}
$,
$ 4+4c^{1}_{7}
+2c^{2}_{7}
$,
$ -4-2  c^{1}_{7}
-2  c^{2}_{7}
$,
$ 6+4c^{1}_{7}
+2c^{2}_{7}
$,
$ -6-4  c^{1}_{7}
-2  c^{2}_{7}
$,
$ -2c_{7}^{1}$,
$0$;\ \ 
$ -2\xi_{7}^{3}$,
$ -4-2  c^{1}_{7}
-2  c^{2}_{7}
$,
$ 4+4c^{1}_{7}
+2c^{2}_{7}
$,
$ 6+4c^{1}_{7}
+2c^{2}_{7}
$,
$ -2c_{7}^{1}$,
$ -6-4  c^{1}_{7}
-2  c^{2}_{7}
$,
$0$;\ \ 
$ 6+4c^{1}_{7}
+2c^{2}_{7}
$,
$ 2c_{7}^{1}$,
$ -2\xi_{7}^{3}$,
$ -6-6  c^{1}_{7}
-2  c^{2}_{7}
$,
$ 2\xi_{7}^{3}$,
$0$;\ \ 
$ 6+4c^{1}_{7}
+2c^{2}_{7}
$,
$ -2\xi_{7}^{3}$,
$ 2\xi_{7}^{3}$,
$ -6-6  c^{1}_{7}
-2  c^{2}_{7}
$,
$0$;\ \ 
$ -1-2  c^{1}_{7}
$,
$ 4+2c^{1}_{7}
+2c^{2}_{7}
$,
$ 4+2c^{1}_{7}
+2c^{2}_{7}
$,
$ -5-6  c^{1}_{7}
-2  c^{2}_{7}
$;\ \ 
$ -4-2  c^{1}_{7}
-2  c^{2}_{7}
$,
$ 4+4c^{1}_{7}
+2c^{2}_{7}
$,
$0$;\ \ 
$ -4-2  c^{1}_{7}
-2  c^{2}_{7}
$,
$0$;\ \ 
$ 5+6c^{1}_{7}
+2c^{2}_{7}
$)

  \vskip 2ex

\noindent73. $10_{6,684.3}^{77,298}$ \irep{1138}:\ \ 
$d_i$ = ($1.0$,
$7.887$,
$7.887$,
$7.887$,
$7.887$,
$7.887$,
$8.887$,
$9.887$,
$9.887$,
$9.887$) 

\vskip 0.7ex
\hangindent=3em \hangafter=1
$D^2= 684.336 = 
\frac{693+77\sqrt{77}}{2}$

\vskip 0.7ex
\hangindent=3em \hangafter=1
$T = ( 0,
\frac{1}{11},
\frac{3}{11},
\frac{4}{11},
\frac{5}{11},
\frac{9}{11},
0,
\frac{3}{7},
\frac{5}{7},
\frac{6}{7} )
$,

\vskip 0.7ex
\hangindent=3em \hangafter=1
$S$ = ($ 1$,
$ \frac{7+\sqrt{77}}{2}$,
$ \frac{7+\sqrt{77}}{2}$,
$ \frac{7+\sqrt{77}}{2}$,
$ \frac{7+\sqrt{77}}{2}$,
$ \frac{7+\sqrt{77}}{2}$,
$ \frac{9+\sqrt{77}}{2}$,
$ \frac{11+\sqrt{77}}{2}$,
$ \frac{11+\sqrt{77}}{2}$,
$ \frac{11+\sqrt{77}}{2}$;\ \ 
$ -1-2  c^{1}_{77}
+c^{2}_{77}
-c^{3}_{77}
-c^{4}_{77}
+c^{5}_{77}
+2c^{6}_{77}
-c^{7}_{77}
-c^{8}_{77}
+c^{9}_{77}
-2  c^{10}_{77}
-c^{11}_{77}
+c^{12}_{77}
-4  c^{14}_{77}
-c^{15}_{77}
+3c^{16}_{77}
+2c^{17}_{77}
-c^{18}_{77}
+c^{19}_{77}
-c^{21}_{77}
-c^{22}_{77}
-c^{23}_{77}
+c^{28}_{77}
-c^{29}_{77}
$,
$ -2+2c^{2}_{77}
-2  c^{3}_{77}
-2  c^{4}_{77}
+2c^{5}_{77}
-2  c^{6}_{77}
-6  c^{7}_{77}
-2  c^{8}_{77}
+2c^{9}_{77}
-2  c^{11}_{77}
+2c^{12}_{77}
-c^{14}_{77}
-2  c^{15}_{77}
-2  c^{17}_{77}
-3  c^{18}_{77}
+2c^{19}_{77}
-2  c^{22}_{77}
+2c^{23}_{77}
-c^{26}_{77}
-c^{28}_{77}
-3  c^{29}_{77}
$,
$ 1+2c^{1}_{77}
+c^{6}_{77}
+c^{7}_{77}
-2  c^{9}_{77}
+2c^{10}_{77}
-2  c^{13}_{77}
+c^{14}_{77}
+c^{16}_{77}
+c^{17}_{77}
+2c^{21}_{77}
+2c^{23}_{77}
-2  c^{24}_{77}
-2  c^{28}_{77}
$,
$ -1+c^{1}_{77}
+2c^{9}_{77}
+c^{10}_{77}
+2c^{13}_{77}
-c^{14}_{77}
+2c^{18}_{77}
-4  c^{21}_{77}
+c^{23}_{77}
+2c^{24}_{77}
+2c^{26}_{77}
-c^{28}_{77}
+2c^{29}_{77}
$,
$ 5-2  c^{2}_{77}
+2c^{3}_{77}
+2c^{4}_{77}
-2  c^{5}_{77}
+4c^{7}_{77}
+2c^{8}_{77}
-c^{9}_{77}
+2c^{11}_{77}
-2  c^{12}_{77}
+c^{13}_{77}
+4c^{14}_{77}
+2c^{15}_{77}
-2  c^{16}_{77}
-2  c^{19}_{77}
+3c^{21}_{77}
+2c^{22}_{77}
-2  c^{23}_{77}
+c^{24}_{77}
-2  c^{26}_{77}
+3c^{28}_{77}
$,
$ -\frac{7+\sqrt{77}}{2}$,
$0$,
$0$,
$0$;\ \ 
$ 5-2  c^{2}_{77}
+2c^{3}_{77}
+2c^{4}_{77}
-2  c^{5}_{77}
+4c^{7}_{77}
+2c^{8}_{77}
-c^{9}_{77}
+2c^{11}_{77}
-2  c^{12}_{77}
+c^{13}_{77}
+4c^{14}_{77}
+2c^{15}_{77}
-2  c^{16}_{77}
-2  c^{19}_{77}
+3c^{21}_{77}
+2c^{22}_{77}
-2  c^{23}_{77}
+c^{24}_{77}
-2  c^{26}_{77}
+3c^{28}_{77}
$,
$ -1-2  c^{1}_{77}
+c^{2}_{77}
-c^{3}_{77}
-c^{4}_{77}
+c^{5}_{77}
+2c^{6}_{77}
-c^{7}_{77}
-c^{8}_{77}
+c^{9}_{77}
-2  c^{10}_{77}
-c^{11}_{77}
+c^{12}_{77}
-4  c^{14}_{77}
-c^{15}_{77}
+3c^{16}_{77}
+2c^{17}_{77}
-c^{18}_{77}
+c^{19}_{77}
-c^{21}_{77}
-c^{22}_{77}
-c^{23}_{77}
+c^{28}_{77}
-c^{29}_{77}
$,
$ 1+2c^{1}_{77}
+c^{6}_{77}
+c^{7}_{77}
-2  c^{9}_{77}
+2c^{10}_{77}
-2  c^{13}_{77}
+c^{14}_{77}
+c^{16}_{77}
+c^{17}_{77}
+2c^{21}_{77}
+2c^{23}_{77}
-2  c^{24}_{77}
-2  c^{28}_{77}
$,
$ -1+c^{1}_{77}
+2c^{9}_{77}
+c^{10}_{77}
+2c^{13}_{77}
-c^{14}_{77}
+2c^{18}_{77}
-4  c^{21}_{77}
+c^{23}_{77}
+2c^{24}_{77}
+2c^{26}_{77}
-c^{28}_{77}
+2c^{29}_{77}
$,
$ -\frac{7+\sqrt{77}}{2}$,
$0$,
$0$,
$0$;\ \ 
$ -1+c^{1}_{77}
+2c^{9}_{77}
+c^{10}_{77}
+2c^{13}_{77}
-c^{14}_{77}
+2c^{18}_{77}
-4  c^{21}_{77}
+c^{23}_{77}
+2c^{24}_{77}
+2c^{26}_{77}
-c^{28}_{77}
+2c^{29}_{77}
$,
$ 5-2  c^{2}_{77}
+2c^{3}_{77}
+2c^{4}_{77}
-2  c^{5}_{77}
+4c^{7}_{77}
+2c^{8}_{77}
-c^{9}_{77}
+2c^{11}_{77}
-2  c^{12}_{77}
+c^{13}_{77}
+4c^{14}_{77}
+2c^{15}_{77}
-2  c^{16}_{77}
-2  c^{19}_{77}
+3c^{21}_{77}
+2c^{22}_{77}
-2  c^{23}_{77}
+c^{24}_{77}
-2  c^{26}_{77}
+3c^{28}_{77}
$,
$ -2+2c^{2}_{77}
-2  c^{3}_{77}
-2  c^{4}_{77}
+2c^{5}_{77}
-2  c^{6}_{77}
-6  c^{7}_{77}
-2  c^{8}_{77}
+2c^{9}_{77}
-2  c^{11}_{77}
+2c^{12}_{77}
-c^{14}_{77}
-2  c^{15}_{77}
-2  c^{17}_{77}
-3  c^{18}_{77}
+2c^{19}_{77}
-2  c^{22}_{77}
+2c^{23}_{77}
-c^{26}_{77}
-c^{28}_{77}
-3  c^{29}_{77}
$,
$ -\frac{7+\sqrt{77}}{2}$,
$0$,
$0$,
$0$;\ \ 
$ -2+2c^{2}_{77}
-2  c^{3}_{77}
-2  c^{4}_{77}
+2c^{5}_{77}
-2  c^{6}_{77}
-6  c^{7}_{77}
-2  c^{8}_{77}
+2c^{9}_{77}
-2  c^{11}_{77}
+2c^{12}_{77}
-c^{14}_{77}
-2  c^{15}_{77}
-2  c^{17}_{77}
-3  c^{18}_{77}
+2c^{19}_{77}
-2  c^{22}_{77}
+2c^{23}_{77}
-c^{26}_{77}
-c^{28}_{77}
-3  c^{29}_{77}
$,
$ -1-2  c^{1}_{77}
+c^{2}_{77}
-c^{3}_{77}
-c^{4}_{77}
+c^{5}_{77}
+2c^{6}_{77}
-c^{7}_{77}
-c^{8}_{77}
+c^{9}_{77}
-2  c^{10}_{77}
-c^{11}_{77}
+c^{12}_{77}
-4  c^{14}_{77}
-c^{15}_{77}
+3c^{16}_{77}
+2c^{17}_{77}
-c^{18}_{77}
+c^{19}_{77}
-c^{21}_{77}
-c^{22}_{77}
-c^{23}_{77}
+c^{28}_{77}
-c^{29}_{77}
$,
$ -\frac{7+\sqrt{77}}{2}$,
$0$,
$0$,
$0$;\ \ 
$ 1+2c^{1}_{77}
+c^{6}_{77}
+c^{7}_{77}
-2  c^{9}_{77}
+2c^{10}_{77}
-2  c^{13}_{77}
+c^{14}_{77}
+c^{16}_{77}
+c^{17}_{77}
+2c^{21}_{77}
+2c^{23}_{77}
-2  c^{24}_{77}
-2  c^{28}_{77}
$,
$ -\frac{7+\sqrt{77}}{2}$,
$0$,
$0$,
$0$;\ \ 
$ 1$,
$ \frac{11+\sqrt{77}}{2}$,
$ \frac{11+\sqrt{77}}{2}$,
$ \frac{11+\sqrt{77}}{2}$;\ \ 
$ 1+3c^{4}_{77}
+2c^{7}_{77}
-2  c^{9}_{77}
+c^{10}_{77}
+7c^{11}_{77}
+2c^{15}_{77}
-2  c^{16}_{77}
+c^{17}_{77}
+2c^{18}_{77}
-2  c^{19}_{77}
+c^{22}_{77}
-2  c^{23}_{77}
+c^{24}_{77}
+c^{25}_{77}
+2c^{26}_{77}
+2c^{29}_{77}
$,
$ -1+c^{1}_{77}
+c^{2}_{77}
-c^{3}_{77}
-3  c^{4}_{77}
+c^{5}_{77}
+c^{6}_{77}
-c^{7}_{77}
-c^{8}_{77}
+c^{9}_{77}
-2  c^{10}_{77}
-2  c^{11}_{77}
+c^{12}_{77}
+c^{13}_{77}
-c^{14}_{77}
+c^{16}_{77}
-2  c^{17}_{77}
-c^{18}_{77}
+5c^{22}_{77}
+c^{23}_{77}
-2  c^{24}_{77}
-3  c^{25}_{77}
-c^{29}_{77}
$,
$ -4-2  c^{1}_{77}
-2  c^{2}_{77}
+2c^{3}_{77}
+c^{4}_{77}
-2  c^{5}_{77}
-2  c^{6}_{77}
+c^{7}_{77}
+2c^{8}_{77}
-c^{9}_{77}
-4  c^{11}_{77}
-2  c^{12}_{77}
-2  c^{13}_{77}
+2c^{14}_{77}
-c^{15}_{77}
-c^{16}_{77}
+c^{18}_{77}
+c^{19}_{77}
-5  c^{22}_{77}
-c^{23}_{77}
+2c^{25}_{77}
-c^{26}_{77}
+c^{29}_{77}
$;\ \ 
$ -4-2  c^{1}_{77}
-2  c^{2}_{77}
+2c^{3}_{77}
+c^{4}_{77}
-2  c^{5}_{77}
-2  c^{6}_{77}
+c^{7}_{77}
+2c^{8}_{77}
-c^{9}_{77}
-4  c^{11}_{77}
-2  c^{12}_{77}
-2  c^{13}_{77}
+2c^{14}_{77}
-c^{15}_{77}
-c^{16}_{77}
+c^{18}_{77}
+c^{19}_{77}
-5  c^{22}_{77}
-c^{23}_{77}
+2c^{25}_{77}
-c^{26}_{77}
+c^{29}_{77}
$,
$ 1+3c^{4}_{77}
+2c^{7}_{77}
-2  c^{9}_{77}
+c^{10}_{77}
+7c^{11}_{77}
+2c^{15}_{77}
-2  c^{16}_{77}
+c^{17}_{77}
+2c^{18}_{77}
-2  c^{19}_{77}
+c^{22}_{77}
-2  c^{23}_{77}
+c^{24}_{77}
+c^{25}_{77}
+2c^{26}_{77}
+2c^{29}_{77}
$;\ \ 
$ -1+c^{1}_{77}
+c^{2}_{77}
-c^{3}_{77}
-3  c^{4}_{77}
+c^{5}_{77}
+c^{6}_{77}
-c^{7}_{77}
-c^{8}_{77}
+c^{9}_{77}
-2  c^{10}_{77}
-2  c^{11}_{77}
+c^{12}_{77}
+c^{13}_{77}
-c^{14}_{77}
+c^{16}_{77}
-2  c^{17}_{77}
-c^{18}_{77}
+5c^{22}_{77}
+c^{23}_{77}
-2  c^{24}_{77}
-3  c^{25}_{77}
-c^{29}_{77}
$)

  \vskip 2ex

\noindent74. $10_{2,684.3}^{77,982}$ \irep{1138}:\ \ 
$d_i$ = ($1.0$,
$7.887$,
$7.887$,
$7.887$,
$7.887$,
$7.887$,
$8.887$,
$9.887$,
$9.887$,
$9.887$) 

\vskip 0.7ex
\hangindent=3em \hangafter=1
$D^2= 684.336 = 
\frac{693+77\sqrt{77}}{2}$

\vskip 0.7ex
\hangindent=3em \hangafter=1
$T = ( 0,
\frac{2}{11},
\frac{6}{11},
\frac{7}{11},
\frac{8}{11},
\frac{10}{11},
0,
\frac{1}{7},
\frac{2}{7},
\frac{4}{7} )
$,

\vskip 0.7ex
\hangindent=3em \hangafter=1
$S$ = ($ 1$,
$ \frac{7+\sqrt{77}}{2}$,
$ \frac{7+\sqrt{77}}{2}$,
$ \frac{7+\sqrt{77}}{2}$,
$ \frac{7+\sqrt{77}}{2}$,
$ \frac{7+\sqrt{77}}{2}$,
$ \frac{9+\sqrt{77}}{2}$,
$ \frac{11+\sqrt{77}}{2}$,
$ \frac{11+\sqrt{77}}{2}$,
$ \frac{11+\sqrt{77}}{2}$;\ \ 
$ 1+2c^{1}_{77}
+c^{6}_{77}
+c^{7}_{77}
-2  c^{9}_{77}
+2c^{10}_{77}
-2  c^{13}_{77}
+c^{14}_{77}
+c^{16}_{77}
+c^{17}_{77}
+2c^{21}_{77}
+2c^{23}_{77}
-2  c^{24}_{77}
-2  c^{28}_{77}
$,
$ -1-2  c^{1}_{77}
+c^{2}_{77}
-c^{3}_{77}
-c^{4}_{77}
+c^{5}_{77}
+2c^{6}_{77}
-c^{7}_{77}
-c^{8}_{77}
+c^{9}_{77}
-2  c^{10}_{77}
-c^{11}_{77}
+c^{12}_{77}
-4  c^{14}_{77}
-c^{15}_{77}
+3c^{16}_{77}
+2c^{17}_{77}
-c^{18}_{77}
+c^{19}_{77}
-c^{21}_{77}
-c^{22}_{77}
-c^{23}_{77}
+c^{28}_{77}
-c^{29}_{77}
$,
$ -2+2c^{2}_{77}
-2  c^{3}_{77}
-2  c^{4}_{77}
+2c^{5}_{77}
-2  c^{6}_{77}
-6  c^{7}_{77}
-2  c^{8}_{77}
+2c^{9}_{77}
-2  c^{11}_{77}
+2c^{12}_{77}
-c^{14}_{77}
-2  c^{15}_{77}
-2  c^{17}_{77}
-3  c^{18}_{77}
+2c^{19}_{77}
-2  c^{22}_{77}
+2c^{23}_{77}
-c^{26}_{77}
-c^{28}_{77}
-3  c^{29}_{77}
$,
$ -1+c^{1}_{77}
+2c^{9}_{77}
+c^{10}_{77}
+2c^{13}_{77}
-c^{14}_{77}
+2c^{18}_{77}
-4  c^{21}_{77}
+c^{23}_{77}
+2c^{24}_{77}
+2c^{26}_{77}
-c^{28}_{77}
+2c^{29}_{77}
$,
$ 5-2  c^{2}_{77}
+2c^{3}_{77}
+2c^{4}_{77}
-2  c^{5}_{77}
+4c^{7}_{77}
+2c^{8}_{77}
-c^{9}_{77}
+2c^{11}_{77}
-2  c^{12}_{77}
+c^{13}_{77}
+4c^{14}_{77}
+2c^{15}_{77}
-2  c^{16}_{77}
-2  c^{19}_{77}
+3c^{21}_{77}
+2c^{22}_{77}
-2  c^{23}_{77}
+c^{24}_{77}
-2  c^{26}_{77}
+3c^{28}_{77}
$,
$ -\frac{7+\sqrt{77}}{2}$,
$0$,
$0$,
$0$;\ \ 
$ -2+2c^{2}_{77}
-2  c^{3}_{77}
-2  c^{4}_{77}
+2c^{5}_{77}
-2  c^{6}_{77}
-6  c^{7}_{77}
-2  c^{8}_{77}
+2c^{9}_{77}
-2  c^{11}_{77}
+2c^{12}_{77}
-c^{14}_{77}
-2  c^{15}_{77}
-2  c^{17}_{77}
-3  c^{18}_{77}
+2c^{19}_{77}
-2  c^{22}_{77}
+2c^{23}_{77}
-c^{26}_{77}
-c^{28}_{77}
-3  c^{29}_{77}
$,
$ 5-2  c^{2}_{77}
+2c^{3}_{77}
+2c^{4}_{77}
-2  c^{5}_{77}
+4c^{7}_{77}
+2c^{8}_{77}
-c^{9}_{77}
+2c^{11}_{77}
-2  c^{12}_{77}
+c^{13}_{77}
+4c^{14}_{77}
+2c^{15}_{77}
-2  c^{16}_{77}
-2  c^{19}_{77}
+3c^{21}_{77}
+2c^{22}_{77}
-2  c^{23}_{77}
+c^{24}_{77}
-2  c^{26}_{77}
+3c^{28}_{77}
$,
$ 1+2c^{1}_{77}
+c^{6}_{77}
+c^{7}_{77}
-2  c^{9}_{77}
+2c^{10}_{77}
-2  c^{13}_{77}
+c^{14}_{77}
+c^{16}_{77}
+c^{17}_{77}
+2c^{21}_{77}
+2c^{23}_{77}
-2  c^{24}_{77}
-2  c^{28}_{77}
$,
$ -1+c^{1}_{77}
+2c^{9}_{77}
+c^{10}_{77}
+2c^{13}_{77}
-c^{14}_{77}
+2c^{18}_{77}
-4  c^{21}_{77}
+c^{23}_{77}
+2c^{24}_{77}
+2c^{26}_{77}
-c^{28}_{77}
+2c^{29}_{77}
$,
$ -\frac{7+\sqrt{77}}{2}$,
$0$,
$0$,
$0$;\ \ 
$ -1+c^{1}_{77}
+2c^{9}_{77}
+c^{10}_{77}
+2c^{13}_{77}
-c^{14}_{77}
+2c^{18}_{77}
-4  c^{21}_{77}
+c^{23}_{77}
+2c^{24}_{77}
+2c^{26}_{77}
-c^{28}_{77}
+2c^{29}_{77}
$,
$ -1-2  c^{1}_{77}
+c^{2}_{77}
-c^{3}_{77}
-c^{4}_{77}
+c^{5}_{77}
+2c^{6}_{77}
-c^{7}_{77}
-c^{8}_{77}
+c^{9}_{77}
-2  c^{10}_{77}
-c^{11}_{77}
+c^{12}_{77}
-4  c^{14}_{77}
-c^{15}_{77}
+3c^{16}_{77}
+2c^{17}_{77}
-c^{18}_{77}
+c^{19}_{77}
-c^{21}_{77}
-c^{22}_{77}
-c^{23}_{77}
+c^{28}_{77}
-c^{29}_{77}
$,
$ 1+2c^{1}_{77}
+c^{6}_{77}
+c^{7}_{77}
-2  c^{9}_{77}
+2c^{10}_{77}
-2  c^{13}_{77}
+c^{14}_{77}
+c^{16}_{77}
+c^{17}_{77}
+2c^{21}_{77}
+2c^{23}_{77}
-2  c^{24}_{77}
-2  c^{28}_{77}
$,
$ -\frac{7+\sqrt{77}}{2}$,
$0$,
$0$,
$0$;\ \ 
$ 5-2  c^{2}_{77}
+2c^{3}_{77}
+2c^{4}_{77}
-2  c^{5}_{77}
+4c^{7}_{77}
+2c^{8}_{77}
-c^{9}_{77}
+2c^{11}_{77}
-2  c^{12}_{77}
+c^{13}_{77}
+4c^{14}_{77}
+2c^{15}_{77}
-2  c^{16}_{77}
-2  c^{19}_{77}
+3c^{21}_{77}
+2c^{22}_{77}
-2  c^{23}_{77}
+c^{24}_{77}
-2  c^{26}_{77}
+3c^{28}_{77}
$,
$ -2+2c^{2}_{77}
-2  c^{3}_{77}
-2  c^{4}_{77}
+2c^{5}_{77}
-2  c^{6}_{77}
-6  c^{7}_{77}
-2  c^{8}_{77}
+2c^{9}_{77}
-2  c^{11}_{77}
+2c^{12}_{77}
-c^{14}_{77}
-2  c^{15}_{77}
-2  c^{17}_{77}
-3  c^{18}_{77}
+2c^{19}_{77}
-2  c^{22}_{77}
+2c^{23}_{77}
-c^{26}_{77}
-c^{28}_{77}
-3  c^{29}_{77}
$,
$ -\frac{7+\sqrt{77}}{2}$,
$0$,
$0$,
$0$;\ \ 
$ -1-2  c^{1}_{77}
+c^{2}_{77}
-c^{3}_{77}
-c^{4}_{77}
+c^{5}_{77}
+2c^{6}_{77}
-c^{7}_{77}
-c^{8}_{77}
+c^{9}_{77}
-2  c^{10}_{77}
-c^{11}_{77}
+c^{12}_{77}
-4  c^{14}_{77}
-c^{15}_{77}
+3c^{16}_{77}
+2c^{17}_{77}
-c^{18}_{77}
+c^{19}_{77}
-c^{21}_{77}
-c^{22}_{77}
-c^{23}_{77}
+c^{28}_{77}
-c^{29}_{77}
$,
$ -\frac{7+\sqrt{77}}{2}$,
$0$,
$0$,
$0$;\ \ 
$ 1$,
$ \frac{11+\sqrt{77}}{2}$,
$ \frac{11+\sqrt{77}}{2}$,
$ \frac{11+\sqrt{77}}{2}$;\ \ 
$ -1+c^{1}_{77}
+c^{2}_{77}
-c^{3}_{77}
-3  c^{4}_{77}
+c^{5}_{77}
+c^{6}_{77}
-c^{7}_{77}
-c^{8}_{77}
+c^{9}_{77}
-2  c^{10}_{77}
-2  c^{11}_{77}
+c^{12}_{77}
+c^{13}_{77}
-c^{14}_{77}
+c^{16}_{77}
-2  c^{17}_{77}
-c^{18}_{77}
+5c^{22}_{77}
+c^{23}_{77}
-2  c^{24}_{77}
-3  c^{25}_{77}
-c^{29}_{77}
$,
$ 1+3c^{4}_{77}
+2c^{7}_{77}
-2  c^{9}_{77}
+c^{10}_{77}
+7c^{11}_{77}
+2c^{15}_{77}
-2  c^{16}_{77}
+c^{17}_{77}
+2c^{18}_{77}
-2  c^{19}_{77}
+c^{22}_{77}
-2  c^{23}_{77}
+c^{24}_{77}
+c^{25}_{77}
+2c^{26}_{77}
+2c^{29}_{77}
$,
$ -4-2  c^{1}_{77}
-2  c^{2}_{77}
+2c^{3}_{77}
+c^{4}_{77}
-2  c^{5}_{77}
-2  c^{6}_{77}
+c^{7}_{77}
+2c^{8}_{77}
-c^{9}_{77}
-4  c^{11}_{77}
-2  c^{12}_{77}
-2  c^{13}_{77}
+2c^{14}_{77}
-c^{15}_{77}
-c^{16}_{77}
+c^{18}_{77}
+c^{19}_{77}
-5  c^{22}_{77}
-c^{23}_{77}
+2c^{25}_{77}
-c^{26}_{77}
+c^{29}_{77}
$;\ \ 
$ -4-2  c^{1}_{77}
-2  c^{2}_{77}
+2c^{3}_{77}
+c^{4}_{77}
-2  c^{5}_{77}
-2  c^{6}_{77}
+c^{7}_{77}
+2c^{8}_{77}
-c^{9}_{77}
-4  c^{11}_{77}
-2  c^{12}_{77}
-2  c^{13}_{77}
+2c^{14}_{77}
-c^{15}_{77}
-c^{16}_{77}
+c^{18}_{77}
+c^{19}_{77}
-5  c^{22}_{77}
-c^{23}_{77}
+2c^{25}_{77}
-c^{26}_{77}
+c^{29}_{77}
$,
$ -1+c^{1}_{77}
+c^{2}_{77}
-c^{3}_{77}
-3  c^{4}_{77}
+c^{5}_{77}
+c^{6}_{77}
-c^{7}_{77}
-c^{8}_{77}
+c^{9}_{77}
-2  c^{10}_{77}
-2  c^{11}_{77}
+c^{12}_{77}
+c^{13}_{77}
-c^{14}_{77}
+c^{16}_{77}
-2  c^{17}_{77}
-c^{18}_{77}
+5c^{22}_{77}
+c^{23}_{77}
-2  c^{24}_{77}
-3  c^{25}_{77}
-c^{29}_{77}
$;\ \ 
$ 1+3c^{4}_{77}
+2c^{7}_{77}
-2  c^{9}_{77}
+c^{10}_{77}
+7c^{11}_{77}
+2c^{15}_{77}
-2  c^{16}_{77}
+c^{17}_{77}
+2c^{18}_{77}
-2  c^{19}_{77}
+c^{22}_{77}
-2  c^{23}_{77}
+c^{24}_{77}
+c^{25}_{77}
+2c^{26}_{77}
+2c^{29}_{77}
$)

  \vskip 2ex

\noindent75. $10_{4,1435.}^{10,168}$ \irep{538}:\ \ 
$d_i$ = ($1.0$,
$9.472$,
$9.472$,
$9.472$,
$9.472$,
$9.472$,
$9.472$,
$16.944$,
$16.944$,
$17.944$) 

\vskip 0.7ex
\hangindent=3em \hangafter=1
$D^2= 1435.541 = 
720+320\sqrt{5}$

\vskip 0.7ex
\hangindent=3em \hangafter=1
$T = ( 0,
\frac{1}{2},
\frac{1}{2},
\frac{1}{2},
\frac{1}{2},
\frac{1}{2},
\frac{1}{2},
\frac{1}{5},
\frac{4}{5},
0 )
$,

\vskip 0.7ex
\hangindent=3em \hangafter=1
$S$ = ($ 1$,
$ 5+2\sqrt{5}$,
$ 5+2\sqrt{5}$,
$ 5+2\sqrt{5}$,
$ 5+2\sqrt{5}$,
$ 5+2\sqrt{5}$,
$ 5+2\sqrt{5}$,
$ 8+4\sqrt{5}$,
$ 8+4\sqrt{5}$,
$ 9+4\sqrt{5}$;\ \ 
$ 15+6\sqrt{5}$,
$ -5-2\sqrt{5}$,
$ -5-2\sqrt{5}$,
$ -5-2\sqrt{5}$,
$ -5-2\sqrt{5}$,
$ -5-2\sqrt{5}$,
$0$,
$0$,
$ 5+2\sqrt{5}$;\ \ 
$ 15+6\sqrt{5}$,
$ -5-2\sqrt{5}$,
$ -5-2\sqrt{5}$,
$ -5-2\sqrt{5}$,
$ -5-2\sqrt{5}$,
$0$,
$0$,
$ 5+2\sqrt{5}$;\ \ 
$ 15+6\sqrt{5}$,
$ -5-2\sqrt{5}$,
$ -5-2\sqrt{5}$,
$ -5-2\sqrt{5}$,
$0$,
$0$,
$ 5+2\sqrt{5}$;\ \ 
$ 15+6\sqrt{5}$,
$ -5-2\sqrt{5}$,
$ -5-2\sqrt{5}$,
$0$,
$0$,
$ 5+2\sqrt{5}$;\ \ 
$ 15+6\sqrt{5}$,
$ -5-2\sqrt{5}$,
$0$,
$0$,
$ 5+2\sqrt{5}$;\ \ 
$ 15+6\sqrt{5}$,
$0$,
$0$,
$ 5+2\sqrt{5}$;\ \ 
$ 14+6\sqrt{5}$,
$ -6-2\sqrt{5}$,
$ -8-4\sqrt{5}$;\ \ 
$ 14+6\sqrt{5}$,
$ -8-4\sqrt{5}$;\ \ 
$ 1$)

  \vskip 2ex

\noindent76. $10_{0,1435.}^{20,676}$ \irep{948}:\ \ 
$d_i$ = ($1.0$,
$9.472$,
$9.472$,
$9.472$,
$9.472$,
$9.472$,
$9.472$,
$16.944$,
$16.944$,
$17.944$) 

\vskip 0.7ex
\hangindent=3em \hangafter=1
$D^2= 1435.541 = 
720+320\sqrt{5}$

\vskip 0.7ex
\hangindent=3em \hangafter=1
$T = ( 0,
0,
0,
\frac{1}{4},
\frac{1}{4},
\frac{3}{4},
\frac{3}{4},
\frac{2}{5},
\frac{3}{5},
0 )
$,

\vskip 0.7ex
\hangindent=3em \hangafter=1
$S$ = ($ 1$,
$ 5+2\sqrt{5}$,
$ 5+2\sqrt{5}$,
$ 5+2\sqrt{5}$,
$ 5+2\sqrt{5}$,
$ 5+2\sqrt{5}$,
$ 5+2\sqrt{5}$,
$ 8+4\sqrt{5}$,
$ 8+4\sqrt{5}$,
$ 9+4\sqrt{5}$;\ \ 
$ 15+6\sqrt{5}$,
$ -5-2\sqrt{5}$,
$ -5-2\sqrt{5}$,
$ -5-2\sqrt{5}$,
$ -5-2\sqrt{5}$,
$ -5-2\sqrt{5}$,
$0$,
$0$,
$ 5+2\sqrt{5}$;\ \ 
$ 15+6\sqrt{5}$,
$ -5-2\sqrt{5}$,
$ -5-2\sqrt{5}$,
$ -5-2\sqrt{5}$,
$ -5-2\sqrt{5}$,
$0$,
$0$,
$ 5+2\sqrt{5}$;\ \ 
$ -3-6  s^{1}_{20}
-4  c^{2}_{20}
+14s^{3}_{20}
$,
$ -3+6s^{1}_{20}
-4  c^{2}_{20}
-14  s^{3}_{20}
$,
$ 5+2\sqrt{5}$,
$ 5+2\sqrt{5}$,
$0$,
$0$,
$ 5+2\sqrt{5}$;\ \ 
$ -3-6  s^{1}_{20}
-4  c^{2}_{20}
+14s^{3}_{20}
$,
$ 5+2\sqrt{5}$,
$ 5+2\sqrt{5}$,
$0$,
$0$,
$ 5+2\sqrt{5}$;\ \ 
$ -3+6s^{1}_{20}
-4  c^{2}_{20}
-14  s^{3}_{20}
$,
$ -3-6  s^{1}_{20}
-4  c^{2}_{20}
+14s^{3}_{20}
$,
$0$,
$0$,
$ 5+2\sqrt{5}$;\ \ 
$ -3+6s^{1}_{20}
-4  c^{2}_{20}
-14  s^{3}_{20}
$,
$0$,
$0$,
$ 5+2\sqrt{5}$;\ \ 
$ -6-2\sqrt{5}$,
$ 14+6\sqrt{5}$,
$ -8-4\sqrt{5}$;\ \ 
$ -6-2\sqrt{5}$,
$ -8-4\sqrt{5}$;\ \ 
$ 1$)

  \vskip 2ex 

}

\subsection{Rank 11}
\label{uni11}

{\small

\noindent1. $11_{2,11.}^{11,568}$ \irep{983}:\ \ 
$d_i$ = ($1.0$,
$1.0$,
$1.0$,
$1.0$,
$1.0$,
$1.0$,
$1.0$,
$1.0$,
$1.0$,
$1.0$,
$1.0$) 

\vskip 0.7ex
\hangindent=3em \hangafter=1
$D^2= 11.0 = 
11$

\vskip 0.7ex
\hangindent=3em \hangafter=1
$T = ( 0,
\frac{1}{11},
\frac{1}{11},
\frac{3}{11},
\frac{3}{11},
\frac{4}{11},
\frac{4}{11},
\frac{5}{11},
\frac{5}{11},
\frac{9}{11},
\frac{9}{11} )
$,

\vskip 0.7ex
\hangindent=3em \hangafter=1
$S$ = ($ 1$,
$ 1$,
$ 1$,
$ 1$,
$ 1$,
$ 1$,
$ 1$,
$ 1$,
$ 1$,
$ 1$,
$ 1$;\ \ 
$ -\zeta_{22}^{7}$,
$ \zeta_{11}^{2}$,
$ -\zeta_{22}^{9}$,
$ \zeta_{11}^{1}$,
$ -\zeta_{22}^{3}$,
$ \zeta_{11}^{4}$,
$ -\zeta_{22}^{5}$,
$ \zeta_{11}^{3}$,
$ -\zeta_{22}^{1}$,
$ \zeta_{11}^{5}$;\ \ 
$ -\zeta_{22}^{7}$,
$ \zeta_{11}^{1}$,
$ -\zeta_{22}^{9}$,
$ \zeta_{11}^{4}$,
$ -\zeta_{22}^{3}$,
$ \zeta_{11}^{3}$,
$ -\zeta_{22}^{5}$,
$ \zeta_{11}^{5}$,
$ -\zeta_{22}^{1}$;\ \ 
$ \zeta_{11}^{5}$,
$ -\zeta_{22}^{1}$,
$ -\zeta_{22}^{7}$,
$ \zeta_{11}^{2}$,
$ \zeta_{11}^{4}$,
$ -\zeta_{22}^{3}$,
$ \zeta_{11}^{3}$,
$ -\zeta_{22}^{5}$;\ \ 
$ \zeta_{11}^{5}$,
$ \zeta_{11}^{2}$,
$ -\zeta_{22}^{7}$,
$ -\zeta_{22}^{3}$,
$ \zeta_{11}^{4}$,
$ -\zeta_{22}^{5}$,
$ \zeta_{11}^{3}$;\ \ 
$ \zeta_{11}^{3}$,
$ -\zeta_{22}^{5}$,
$ \zeta_{11}^{5}$,
$ -\zeta_{22}^{1}$,
$ \zeta_{11}^{1}$,
$ -\zeta_{22}^{9}$;\ \ 
$ \zeta_{11}^{3}$,
$ -\zeta_{22}^{1}$,
$ \zeta_{11}^{5}$,
$ -\zeta_{22}^{9}$,
$ \zeta_{11}^{1}$;\ \ 
$ \zeta_{11}^{1}$,
$ -\zeta_{22}^{9}$,
$ -\zeta_{22}^{7}$,
$ \zeta_{11}^{2}$;\ \ 
$ \zeta_{11}^{1}$,
$ \zeta_{11}^{2}$,
$ -\zeta_{22}^{7}$;\ \ 
$ \zeta_{11}^{4}$,
$ -\zeta_{22}^{3}$;\ \ 
$ \zeta_{11}^{4}$)

  \vskip 2ex

\noindent2. $11_{6,11.}^{11,143}$ \irep{983}:\ \ 
$d_i$ = ($1.0$,
$1.0$,
$1.0$,
$1.0$,
$1.0$,
$1.0$,
$1.0$,
$1.0$,
$1.0$,
$1.0$,
$1.0$) 

\vskip 0.7ex
\hangindent=3em \hangafter=1
$D^2= 11.0 = 
11$

\vskip 0.7ex
\hangindent=3em \hangafter=1
$T = ( 0,
\frac{2}{11},
\frac{2}{11},
\frac{6}{11},
\frac{6}{11},
\frac{7}{11},
\frac{7}{11},
\frac{8}{11},
\frac{8}{11},
\frac{10}{11},
\frac{10}{11} )
$,

\vskip 0.7ex
\hangindent=3em \hangafter=1
$S$ = ($ 1$,
$ 1$,
$ 1$,
$ 1$,
$ 1$,
$ 1$,
$ 1$,
$ 1$,
$ 1$,
$ 1$,
$ 1$;\ \ 
$ -\zeta_{22}^{3}$,
$ \zeta_{11}^{4}$,
$ -\zeta_{22}^{7}$,
$ \zeta_{11}^{2}$,
$ -\zeta_{22}^{9}$,
$ \zeta_{11}^{1}$,
$ -\zeta_{22}^{5}$,
$ \zeta_{11}^{3}$,
$ -\zeta_{22}^{1}$,
$ \zeta_{11}^{5}$;\ \ 
$ -\zeta_{22}^{3}$,
$ \zeta_{11}^{2}$,
$ -\zeta_{22}^{7}$,
$ \zeta_{11}^{1}$,
$ -\zeta_{22}^{9}$,
$ \zeta_{11}^{3}$,
$ -\zeta_{22}^{5}$,
$ \zeta_{11}^{5}$,
$ -\zeta_{22}^{1}$;\ \ 
$ -\zeta_{22}^{9}$,
$ \zeta_{11}^{1}$,
$ \zeta_{11}^{5}$,
$ -\zeta_{22}^{1}$,
$ \zeta_{11}^{4}$,
$ -\zeta_{22}^{3}$,
$ \zeta_{11}^{3}$,
$ -\zeta_{22}^{5}$;\ \ 
$ -\zeta_{22}^{9}$,
$ -\zeta_{22}^{1}$,
$ \zeta_{11}^{5}$,
$ -\zeta_{22}^{3}$,
$ \zeta_{11}^{4}$,
$ -\zeta_{22}^{5}$,
$ \zeta_{11}^{3}$;\ \ 
$ -\zeta_{22}^{5}$,
$ \zeta_{11}^{3}$,
$ \zeta_{11}^{2}$,
$ -\zeta_{22}^{7}$,
$ -\zeta_{22}^{3}$,
$ \zeta_{11}^{4}$;\ \ 
$ -\zeta_{22}^{5}$,
$ -\zeta_{22}^{7}$,
$ \zeta_{11}^{2}$,
$ \zeta_{11}^{4}$,
$ -\zeta_{22}^{3}$;\ \ 
$ -\zeta_{22}^{1}$,
$ \zeta_{11}^{5}$,
$ -\zeta_{22}^{9}$,
$ \zeta_{11}^{1}$;\ \ 
$ -\zeta_{22}^{1}$,
$ \zeta_{11}^{1}$,
$ -\zeta_{22}^{9}$;\ \ 
$ \zeta_{11}^{2}$,
$ -\zeta_{22}^{7}$;\ \ 
$ \zeta_{11}^{2}$)

  \vskip 2ex

\noindent3. $11_{1,32.}^{16,245}$ \irep{0}:\ \ 
$d_i$ = ($1.0$,
$1.0$,
$1.0$,
$1.0$,
$2.0$,
$2.0$,
$2.0$,
$2.0$,
$2.0$,
$2.0$,
$2.0$) 

\vskip 0.7ex
\hangindent=3em \hangafter=1
$D^2= 32.0 = 
32$

\vskip 0.7ex
\hangindent=3em \hangafter=1
$T = ( 0,
0,
0,
0,
\frac{1}{4},
\frac{1}{16},
\frac{1}{16},
\frac{1}{16},
\frac{9}{16},
\frac{9}{16},
\frac{9}{16} )
$,

\vskip 0.7ex
\hangindent=3em \hangafter=1
$S$ = ($ 1$,
$ 1$,
$ 1$,
$ 1$,
$ 2$,
$ 2$,
$ 2$,
$ 2$,
$ 2$,
$ 2$,
$ 2$;\ \ 
$ 1$,
$ 1$,
$ 1$,
$ 2$,
$ -2$,
$ -2$,
$ 2$,
$ -2$,
$ -2$,
$ 2$;\ \ 
$ 1$,
$ 1$,
$ 2$,
$ -2$,
$ 2$,
$ -2$,
$ -2$,
$ 2$,
$ -2$;\ \ 
$ 1$,
$ 2$,
$ 2$,
$ -2$,
$ -2$,
$ 2$,
$ -2$,
$ -2$;\ \ 
$ -4$,
$0$,
$0$,
$0$,
$0$,
$0$,
$0$;\ \ 
$ 2\sqrt{2}$,
$0$,
$0$,
$ -2\sqrt{2}$,
$0$,
$0$;\ \ 
$ 2\sqrt{2}$,
$0$,
$0$,
$ -2\sqrt{2}$,
$0$;\ \ 
$ 2\sqrt{2}$,
$0$,
$0$,
$ -2\sqrt{2}$;\ \ 
$ 2\sqrt{2}$,
$0$,
$0$;\ \ 
$ 2\sqrt{2}$,
$0$;\ \ 
$ 2\sqrt{2}$)

  \vskip 2ex

\noindent4. $11_{1,32.}^{16,157}$ \irep{0}:\ \ 
$d_i$ = ($1.0$,
$1.0$,
$1.0$,
$1.0$,
$2.0$,
$2.0$,
$2.0$,
$2.0$,
$2.0$,
$2.0$,
$2.0$) 

\vskip 0.7ex
\hangindent=3em \hangafter=1
$D^2= 32.0 = 
32$

\vskip 0.7ex
\hangindent=3em \hangafter=1
$T = ( 0,
0,
0,
0,
\frac{1}{4},
\frac{1}{16},
\frac{1}{16},
\frac{5}{16},
\frac{9}{16},
\frac{9}{16},
\frac{13}{16} )
$,

\vskip 0.7ex
\hangindent=3em \hangafter=1
$S$ = ($ 1$,
$ 1$,
$ 1$,
$ 1$,
$ 2$,
$ 2$,
$ 2$,
$ 2$,
$ 2$,
$ 2$,
$ 2$;\ \ 
$ 1$,
$ 1$,
$ 1$,
$ 2$,
$ -2$,
$ -2$,
$ 2$,
$ -2$,
$ -2$,
$ 2$;\ \ 
$ 1$,
$ 1$,
$ 2$,
$ -2$,
$ 2$,
$ -2$,
$ -2$,
$ 2$,
$ -2$;\ \ 
$ 1$,
$ 2$,
$ 2$,
$ -2$,
$ -2$,
$ 2$,
$ -2$,
$ -2$;\ \ 
$ -4$,
$0$,
$0$,
$0$,
$0$,
$0$,
$0$;\ \ 
$ 2\sqrt{2}$,
$0$,
$0$,
$ -2\sqrt{2}$,
$0$,
$0$;\ \ 
$ 2\sqrt{2}$,
$0$,
$0$,
$ -2\sqrt{2}$,
$0$;\ \ 
$ -2\sqrt{2}$,
$0$,
$0$,
$ 2\sqrt{2}$;\ \ 
$ 2\sqrt{2}$,
$0$,
$0$;\ \ 
$ 2\sqrt{2}$,
$0$;\ \ 
$ -2\sqrt{2}$)

  \vskip 2ex

\noindent5. $11_{1,32.}^{16,703}$ \irep{0}:\ \ 
$d_i$ = ($1.0$,
$1.0$,
$1.0$,
$1.0$,
$2.0$,
$2.0$,
$2.0$,
$2.0$,
$2.0$,
$2.0$,
$2.0$) 

\vskip 0.7ex
\hangindent=3em \hangafter=1
$D^2= 32.0 = 
32$

\vskip 0.7ex
\hangindent=3em \hangafter=1
$T = ( 0,
0,
0,
0,
\frac{1}{4},
\frac{1}{16},
\frac{5}{16},
\frac{5}{16},
\frac{9}{16},
\frac{13}{16},
\frac{13}{16} )
$,

\vskip 0.7ex
\hangindent=3em \hangafter=1
$S$ = ($ 1$,
$ 1$,
$ 1$,
$ 1$,
$ 2$,
$ 2$,
$ 2$,
$ 2$,
$ 2$,
$ 2$,
$ 2$;\ \ 
$ 1$,
$ 1$,
$ 1$,
$ 2$,
$ -2$,
$ -2$,
$ 2$,
$ -2$,
$ -2$,
$ 2$;\ \ 
$ 1$,
$ 1$,
$ 2$,
$ -2$,
$ 2$,
$ -2$,
$ -2$,
$ 2$,
$ -2$;\ \ 
$ 1$,
$ 2$,
$ 2$,
$ -2$,
$ -2$,
$ 2$,
$ -2$,
$ -2$;\ \ 
$ -4$,
$0$,
$0$,
$0$,
$0$,
$0$,
$0$;\ \ 
$ 2\sqrt{2}$,
$0$,
$0$,
$ -2\sqrt{2}$,
$0$,
$0$;\ \ 
$ -2\sqrt{2}$,
$0$,
$0$,
$ 2\sqrt{2}$,
$0$;\ \ 
$ -2\sqrt{2}$,
$0$,
$0$,
$ 2\sqrt{2}$;\ \ 
$ 2\sqrt{2}$,
$0$,
$0$;\ \ 
$ -2\sqrt{2}$,
$0$;\ \ 
$ -2\sqrt{2}$)

  \vskip 2ex

\noindent6. $11_{1,32.}^{16,171}$ \irep{0}:\ \ 
$d_i$ = ($1.0$,
$1.0$,
$1.0$,
$1.0$,
$2.0$,
$2.0$,
$2.0$,
$2.0$,
$2.0$,
$2.0$,
$2.0$) 

\vskip 0.7ex
\hangindent=3em \hangafter=1
$D^2= 32.0 = 
32$

\vskip 0.7ex
\hangindent=3em \hangafter=1
$T = ( 0,
0,
0,
0,
\frac{1}{4},
\frac{5}{16},
\frac{5}{16},
\frac{5}{16},
\frac{13}{16},
\frac{13}{16},
\frac{13}{16} )
$,

\vskip 0.7ex
\hangindent=3em \hangafter=1
$S$ = ($ 1$,
$ 1$,
$ 1$,
$ 1$,
$ 2$,
$ 2$,
$ 2$,
$ 2$,
$ 2$,
$ 2$,
$ 2$;\ \ 
$ 1$,
$ 1$,
$ 1$,
$ 2$,
$ -2$,
$ -2$,
$ 2$,
$ -2$,
$ -2$,
$ 2$;\ \ 
$ 1$,
$ 1$,
$ 2$,
$ -2$,
$ 2$,
$ -2$,
$ -2$,
$ 2$,
$ -2$;\ \ 
$ 1$,
$ 2$,
$ 2$,
$ -2$,
$ -2$,
$ 2$,
$ -2$,
$ -2$;\ \ 
$ -4$,
$0$,
$0$,
$0$,
$0$,
$0$,
$0$;\ \ 
$ -2\sqrt{2}$,
$0$,
$0$,
$ 2\sqrt{2}$,
$0$,
$0$;\ \ 
$ -2\sqrt{2}$,
$0$,
$0$,
$ 2\sqrt{2}$,
$0$;\ \ 
$ -2\sqrt{2}$,
$0$,
$0$,
$ 2\sqrt{2}$;\ \ 
$ -2\sqrt{2}$,
$0$,
$0$;\ \ 
$ -2\sqrt{2}$,
$0$;\ \ 
$ -2\sqrt{2}$)

  \vskip 2ex

\noindent7. $11_{7,32.}^{16,328}$ \irep{0}:\ \ 
$d_i$ = ($1.0$,
$1.0$,
$1.0$,
$1.0$,
$2.0$,
$2.0$,
$2.0$,
$2.0$,
$2.0$,
$2.0$,
$2.0$) 

\vskip 0.7ex
\hangindent=3em \hangafter=1
$D^2= 32.0 = 
32$

\vskip 0.7ex
\hangindent=3em \hangafter=1
$T = ( 0,
0,
0,
0,
\frac{3}{4},
\frac{3}{16},
\frac{3}{16},
\frac{3}{16},
\frac{11}{16},
\frac{11}{16},
\frac{11}{16} )
$,

\vskip 0.7ex
\hangindent=3em \hangafter=1
$S$ = ($ 1$,
$ 1$,
$ 1$,
$ 1$,
$ 2$,
$ 2$,
$ 2$,
$ 2$,
$ 2$,
$ 2$,
$ 2$;\ \ 
$ 1$,
$ 1$,
$ 1$,
$ 2$,
$ -2$,
$ -2$,
$ 2$,
$ -2$,
$ -2$,
$ 2$;\ \ 
$ 1$,
$ 1$,
$ 2$,
$ -2$,
$ 2$,
$ -2$,
$ -2$,
$ 2$,
$ -2$;\ \ 
$ 1$,
$ 2$,
$ 2$,
$ -2$,
$ -2$,
$ 2$,
$ -2$,
$ -2$;\ \ 
$ -4$,
$0$,
$0$,
$0$,
$0$,
$0$,
$0$;\ \ 
$ -2\sqrt{2}$,
$0$,
$0$,
$ 2\sqrt{2}$,
$0$,
$0$;\ \ 
$ -2\sqrt{2}$,
$0$,
$0$,
$ 2\sqrt{2}$,
$0$;\ \ 
$ -2\sqrt{2}$,
$0$,
$0$,
$ 2\sqrt{2}$;\ \ 
$ -2\sqrt{2}$,
$0$,
$0$;\ \ 
$ -2\sqrt{2}$,
$0$;\ \ 
$ -2\sqrt{2}$)

  \vskip 2ex

\noindent8. $11_{7,32.}^{16,796}$ \irep{0}:\ \ 
$d_i$ = ($1.0$,
$1.0$,
$1.0$,
$1.0$,
$2.0$,
$2.0$,
$2.0$,
$2.0$,
$2.0$,
$2.0$,
$2.0$) 

\vskip 0.7ex
\hangindent=3em \hangafter=1
$D^2= 32.0 = 
32$

\vskip 0.7ex
\hangindent=3em \hangafter=1
$T = ( 0,
0,
0,
0,
\frac{3}{4},
\frac{3}{16},
\frac{3}{16},
\frac{7}{16},
\frac{11}{16},
\frac{11}{16},
\frac{15}{16} )
$,

\vskip 0.7ex
\hangindent=3em \hangafter=1
$S$ = ($ 1$,
$ 1$,
$ 1$,
$ 1$,
$ 2$,
$ 2$,
$ 2$,
$ 2$,
$ 2$,
$ 2$,
$ 2$;\ \ 
$ 1$,
$ 1$,
$ 1$,
$ 2$,
$ -2$,
$ -2$,
$ 2$,
$ -2$,
$ -2$,
$ 2$;\ \ 
$ 1$,
$ 1$,
$ 2$,
$ -2$,
$ 2$,
$ -2$,
$ -2$,
$ 2$,
$ -2$;\ \ 
$ 1$,
$ 2$,
$ 2$,
$ -2$,
$ -2$,
$ 2$,
$ -2$,
$ -2$;\ \ 
$ -4$,
$0$,
$0$,
$0$,
$0$,
$0$,
$0$;\ \ 
$ -2\sqrt{2}$,
$0$,
$0$,
$ 2\sqrt{2}$,
$0$,
$0$;\ \ 
$ -2\sqrt{2}$,
$0$,
$0$,
$ 2\sqrt{2}$,
$0$;\ \ 
$ 2\sqrt{2}$,
$0$,
$0$,
$ -2\sqrt{2}$;\ \ 
$ -2\sqrt{2}$,
$0$,
$0$;\ \ 
$ -2\sqrt{2}$,
$0$;\ \ 
$ 2\sqrt{2}$)

  \vskip 2ex

\noindent9. $11_{7,32.}^{16,192}$ \irep{0}:\ \ 
$d_i$ = ($1.0$,
$1.0$,
$1.0$,
$1.0$,
$2.0$,
$2.0$,
$2.0$,
$2.0$,
$2.0$,
$2.0$,
$2.0$) 

\vskip 0.7ex
\hangindent=3em \hangafter=1
$D^2= 32.0 = 
32$

\vskip 0.7ex
\hangindent=3em \hangafter=1
$T = ( 0,
0,
0,
0,
\frac{3}{4},
\frac{3}{16},
\frac{7}{16},
\frac{7}{16},
\frac{11}{16},
\frac{15}{16},
\frac{15}{16} )
$,

\vskip 0.7ex
\hangindent=3em \hangafter=1
$S$ = ($ 1$,
$ 1$,
$ 1$,
$ 1$,
$ 2$,
$ 2$,
$ 2$,
$ 2$,
$ 2$,
$ 2$,
$ 2$;\ \ 
$ 1$,
$ 1$,
$ 1$,
$ 2$,
$ -2$,
$ -2$,
$ 2$,
$ -2$,
$ -2$,
$ 2$;\ \ 
$ 1$,
$ 1$,
$ 2$,
$ -2$,
$ 2$,
$ -2$,
$ -2$,
$ 2$,
$ -2$;\ \ 
$ 1$,
$ 2$,
$ 2$,
$ -2$,
$ -2$,
$ 2$,
$ -2$,
$ -2$;\ \ 
$ -4$,
$0$,
$0$,
$0$,
$0$,
$0$,
$0$;\ \ 
$ -2\sqrt{2}$,
$0$,
$0$,
$ 2\sqrt{2}$,
$0$,
$0$;\ \ 
$ 2\sqrt{2}$,
$0$,
$0$,
$ -2\sqrt{2}$,
$0$;\ \ 
$ 2\sqrt{2}$,
$0$,
$0$,
$ -2\sqrt{2}$;\ \ 
$ -2\sqrt{2}$,
$0$,
$0$;\ \ 
$ 2\sqrt{2}$,
$0$;\ \ 
$ 2\sqrt{2}$)

  \vskip 2ex

\noindent10. $11_{7,32.}^{16,304}$ \irep{0}:\ \ 
$d_i$ = ($1.0$,
$1.0$,
$1.0$,
$1.0$,
$2.0$,
$2.0$,
$2.0$,
$2.0$,
$2.0$,
$2.0$,
$2.0$) 

\vskip 0.7ex
\hangindent=3em \hangafter=1
$D^2= 32.0 = 
32$

\vskip 0.7ex
\hangindent=3em \hangafter=1
$T = ( 0,
0,
0,
0,
\frac{3}{4},
\frac{7}{16},
\frac{7}{16},
\frac{7}{16},
\frac{15}{16},
\frac{15}{16},
\frac{15}{16} )
$,

\vskip 0.7ex
\hangindent=3em \hangafter=1
$S$ = ($ 1$,
$ 1$,
$ 1$,
$ 1$,
$ 2$,
$ 2$,
$ 2$,
$ 2$,
$ 2$,
$ 2$,
$ 2$;\ \ 
$ 1$,
$ 1$,
$ 1$,
$ 2$,
$ -2$,
$ -2$,
$ 2$,
$ -2$,
$ -2$,
$ 2$;\ \ 
$ 1$,
$ 1$,
$ 2$,
$ -2$,
$ 2$,
$ -2$,
$ -2$,
$ 2$,
$ -2$;\ \ 
$ 1$,
$ 2$,
$ 2$,
$ -2$,
$ -2$,
$ 2$,
$ -2$,
$ -2$;\ \ 
$ -4$,
$0$,
$0$,
$0$,
$0$,
$0$,
$0$;\ \ 
$ 2\sqrt{2}$,
$0$,
$0$,
$ -2\sqrt{2}$,
$0$,
$0$;\ \ 
$ 2\sqrt{2}$,
$0$,
$0$,
$ -2\sqrt{2}$,
$0$;\ \ 
$ 2\sqrt{2}$,
$0$,
$0$,
$ -2\sqrt{2}$;\ \ 
$ 2\sqrt{2}$,
$0$,
$0$;\ \ 
$ 2\sqrt{2}$,
$0$;\ \ 
$ 2\sqrt{2}$)

  \vskip 2ex

\noindent11. $11_{2,60.}^{120,157}$ \irep{2274}:\ \ 
$d_i$ = ($1.0$,
$1.0$,
$2.0$,
$2.0$,
$2.0$,
$2.0$,
$2.0$,
$2.0$,
$2.0$,
$3.872$,
$3.872$) 

\vskip 0.7ex
\hangindent=3em \hangafter=1
$D^2= 60.0 = 
60$

\vskip 0.7ex
\hangindent=3em \hangafter=1
$T = ( 0,
0,
\frac{1}{3},
\frac{1}{5},
\frac{4}{5},
\frac{2}{15},
\frac{2}{15},
\frac{8}{15},
\frac{8}{15},
\frac{1}{8},
\frac{5}{8} )
$,

\vskip 0.7ex
\hangindent=3em \hangafter=1
$S$ = ($ 1$,
$ 1$,
$ 2$,
$ 2$,
$ 2$,
$ 2$,
$ 2$,
$ 2$,
$ 2$,
$ \sqrt{15}$,
$ \sqrt{15}$;\ \ 
$ 1$,
$ 2$,
$ 2$,
$ 2$,
$ 2$,
$ 2$,
$ 2$,
$ 2$,
$ -\sqrt{15}$,
$ -\sqrt{15}$;\ \ 
$ -2$,
$ 4$,
$ 4$,
$ -2$,
$ -2$,
$ -2$,
$ -2$,
$0$,
$0$;\ \ 
$ -1-\sqrt{5}$,
$ -1+\sqrt{5}$,
$ -1+\sqrt{5}$,
$ -1+\sqrt{5}$,
$ -1-\sqrt{5}$,
$ -1-\sqrt{5}$,
$0$,
$0$;\ \ 
$ -1-\sqrt{5}$,
$ -1-\sqrt{5}$,
$ -1-\sqrt{5}$,
$ -1+\sqrt{5}$,
$ -1+\sqrt{5}$,
$0$,
$0$;\ \ 
$ 2c_{15}^{4}$,
$ 2c_{15}^{1}$,
$ 2c_{15}^{7}$,
$ 2c_{15}^{2}$,
$0$,
$0$;\ \ 
$ 2c_{15}^{4}$,
$ 2c_{15}^{2}$,
$ 2c_{15}^{7}$,
$0$,
$0$;\ \ 
$ 2c_{15}^{1}$,
$ 2c_{15}^{4}$,
$0$,
$0$;\ \ 
$ 2c_{15}^{1}$,
$0$,
$0$;\ \ 
$ \sqrt{15}$,
$ -\sqrt{15}$;\ \ 
$ \sqrt{15}$)

  \vskip 2ex

\noindent12. $11_{2,60.}^{120,364}$ \irep{2274}:\ \ 
$d_i$ = ($1.0$,
$1.0$,
$2.0$,
$2.0$,
$2.0$,
$2.0$,
$2.0$,
$2.0$,
$2.0$,
$3.872$,
$3.872$) 

\vskip 0.7ex
\hangindent=3em \hangafter=1
$D^2= 60.0 = 
60$

\vskip 0.7ex
\hangindent=3em \hangafter=1
$T = ( 0,
0,
\frac{1}{3},
\frac{1}{5},
\frac{4}{5},
\frac{2}{15},
\frac{2}{15},
\frac{8}{15},
\frac{8}{15},
\frac{3}{8},
\frac{7}{8} )
$,

\vskip 0.7ex
\hangindent=3em \hangafter=1
$S$ = ($ 1$,
$ 1$,
$ 2$,
$ 2$,
$ 2$,
$ 2$,
$ 2$,
$ 2$,
$ 2$,
$ \sqrt{15}$,
$ \sqrt{15}$;\ \ 
$ 1$,
$ 2$,
$ 2$,
$ 2$,
$ 2$,
$ 2$,
$ 2$,
$ 2$,
$ -\sqrt{15}$,
$ -\sqrt{15}$;\ \ 
$ -2$,
$ 4$,
$ 4$,
$ -2$,
$ -2$,
$ -2$,
$ -2$,
$0$,
$0$;\ \ 
$ -1-\sqrt{5}$,
$ -1+\sqrt{5}$,
$ -1+\sqrt{5}$,
$ -1+\sqrt{5}$,
$ -1-\sqrt{5}$,
$ -1-\sqrt{5}$,
$0$,
$0$;\ \ 
$ -1-\sqrt{5}$,
$ -1-\sqrt{5}$,
$ -1-\sqrt{5}$,
$ -1+\sqrt{5}$,
$ -1+\sqrt{5}$,
$0$,
$0$;\ \ 
$ 2c_{15}^{4}$,
$ 2c_{15}^{1}$,
$ 2c_{15}^{7}$,
$ 2c_{15}^{2}$,
$0$,
$0$;\ \ 
$ 2c_{15}^{4}$,
$ 2c_{15}^{2}$,
$ 2c_{15}^{7}$,
$0$,
$0$;\ \ 
$ 2c_{15}^{1}$,
$ 2c_{15}^{4}$,
$0$,
$0$;\ \ 
$ 2c_{15}^{1}$,
$0$,
$0$;\ \ 
$ -\sqrt{15}$,
$ \sqrt{15}$;\ \ 
$ -\sqrt{15}$)

  \vskip 2ex

\noindent13. $11_{6,60.}^{120,253}$ \irep{2274}:\ \ 
$d_i$ = ($1.0$,
$1.0$,
$2.0$,
$2.0$,
$2.0$,
$2.0$,
$2.0$,
$2.0$,
$2.0$,
$3.872$,
$3.872$) 

\vskip 0.7ex
\hangindent=3em \hangafter=1
$D^2= 60.0 = 
60$

\vskip 0.7ex
\hangindent=3em \hangafter=1
$T = ( 0,
0,
\frac{1}{3},
\frac{2}{5},
\frac{3}{5},
\frac{11}{15},
\frac{11}{15},
\frac{14}{15},
\frac{14}{15},
\frac{1}{8},
\frac{5}{8} )
$,

\vskip 0.7ex
\hangindent=3em \hangafter=1
$S$ = ($ 1$,
$ 1$,
$ 2$,
$ 2$,
$ 2$,
$ 2$,
$ 2$,
$ 2$,
$ 2$,
$ \sqrt{15}$,
$ \sqrt{15}$;\ \ 
$ 1$,
$ 2$,
$ 2$,
$ 2$,
$ 2$,
$ 2$,
$ 2$,
$ 2$,
$ -\sqrt{15}$,
$ -\sqrt{15}$;\ \ 
$ -2$,
$ 4$,
$ 4$,
$ -2$,
$ -2$,
$ -2$,
$ -2$,
$0$,
$0$;\ \ 
$ -1+\sqrt{5}$,
$ -1-\sqrt{5}$,
$ -1+\sqrt{5}$,
$ -1+\sqrt{5}$,
$ -1-\sqrt{5}$,
$ -1-\sqrt{5}$,
$0$,
$0$;\ \ 
$ -1+\sqrt{5}$,
$ -1-\sqrt{5}$,
$ -1-\sqrt{5}$,
$ -1+\sqrt{5}$,
$ -1+\sqrt{5}$,
$0$,
$0$;\ \ 
$ 2c_{15}^{7}$,
$ 2c_{15}^{2}$,
$ 2c_{15}^{4}$,
$ 2c_{15}^{1}$,
$0$,
$0$;\ \ 
$ 2c_{15}^{7}$,
$ 2c_{15}^{1}$,
$ 2c_{15}^{4}$,
$0$,
$0$;\ \ 
$ 2c_{15}^{2}$,
$ 2c_{15}^{7}$,
$0$,
$0$;\ \ 
$ 2c_{15}^{2}$,
$0$,
$0$;\ \ 
$ -\sqrt{15}$,
$ \sqrt{15}$;\ \ 
$ -\sqrt{15}$)

  \vskip 2ex

\noindent14. $11_{6,60.}^{120,447}$ \irep{2274}:\ \ 
$d_i$ = ($1.0$,
$1.0$,
$2.0$,
$2.0$,
$2.0$,
$2.0$,
$2.0$,
$2.0$,
$2.0$,
$3.872$,
$3.872$) 

\vskip 0.7ex
\hangindent=3em \hangafter=1
$D^2= 60.0 = 
60$

\vskip 0.7ex
\hangindent=3em \hangafter=1
$T = ( 0,
0,
\frac{1}{3},
\frac{2}{5},
\frac{3}{5},
\frac{11}{15},
\frac{11}{15},
\frac{14}{15},
\frac{14}{15},
\frac{3}{8},
\frac{7}{8} )
$,

\vskip 0.7ex
\hangindent=3em \hangafter=1
$S$ = ($ 1$,
$ 1$,
$ 2$,
$ 2$,
$ 2$,
$ 2$,
$ 2$,
$ 2$,
$ 2$,
$ \sqrt{15}$,
$ \sqrt{15}$;\ \ 
$ 1$,
$ 2$,
$ 2$,
$ 2$,
$ 2$,
$ 2$,
$ 2$,
$ 2$,
$ -\sqrt{15}$,
$ -\sqrt{15}$;\ \ 
$ -2$,
$ 4$,
$ 4$,
$ -2$,
$ -2$,
$ -2$,
$ -2$,
$0$,
$0$;\ \ 
$ -1+\sqrt{5}$,
$ -1-\sqrt{5}$,
$ -1+\sqrt{5}$,
$ -1+\sqrt{5}$,
$ -1-\sqrt{5}$,
$ -1-\sqrt{5}$,
$0$,
$0$;\ \ 
$ -1+\sqrt{5}$,
$ -1-\sqrt{5}$,
$ -1-\sqrt{5}$,
$ -1+\sqrt{5}$,
$ -1+\sqrt{5}$,
$0$,
$0$;\ \ 
$ 2c_{15}^{7}$,
$ 2c_{15}^{2}$,
$ 2c_{15}^{4}$,
$ 2c_{15}^{1}$,
$0$,
$0$;\ \ 
$ 2c_{15}^{7}$,
$ 2c_{15}^{1}$,
$ 2c_{15}^{4}$,
$0$,
$0$;\ \ 
$ 2c_{15}^{2}$,
$ 2c_{15}^{7}$,
$0$,
$0$;\ \ 
$ 2c_{15}^{2}$,
$0$,
$0$;\ \ 
$ \sqrt{15}$,
$ -\sqrt{15}$;\ \ 
$ \sqrt{15}$)

  \vskip 2ex

\noindent15. $11_{6,60.}^{120,176}$ \irep{2274}:\ \ 
$d_i$ = ($1.0$,
$1.0$,
$2.0$,
$2.0$,
$2.0$,
$2.0$,
$2.0$,
$2.0$,
$2.0$,
$3.872$,
$3.872$) 

\vskip 0.7ex
\hangindent=3em \hangafter=1
$D^2= 60.0 = 
60$

\vskip 0.7ex
\hangindent=3em \hangafter=1
$T = ( 0,
0,
\frac{2}{3},
\frac{1}{5},
\frac{4}{5},
\frac{7}{15},
\frac{7}{15},
\frac{13}{15},
\frac{13}{15},
\frac{1}{8},
\frac{5}{8} )
$,

\vskip 0.7ex
\hangindent=3em \hangafter=1
$S$ = ($ 1$,
$ 1$,
$ 2$,
$ 2$,
$ 2$,
$ 2$,
$ 2$,
$ 2$,
$ 2$,
$ \sqrt{15}$,
$ \sqrt{15}$;\ \ 
$ 1$,
$ 2$,
$ 2$,
$ 2$,
$ 2$,
$ 2$,
$ 2$,
$ 2$,
$ -\sqrt{15}$,
$ -\sqrt{15}$;\ \ 
$ -2$,
$ 4$,
$ 4$,
$ -2$,
$ -2$,
$ -2$,
$ -2$,
$0$,
$0$;\ \ 
$ -1-\sqrt{5}$,
$ -1+\sqrt{5}$,
$ -1+\sqrt{5}$,
$ -1+\sqrt{5}$,
$ -1-\sqrt{5}$,
$ -1-\sqrt{5}$,
$0$,
$0$;\ \ 
$ -1-\sqrt{5}$,
$ -1-\sqrt{5}$,
$ -1-\sqrt{5}$,
$ -1+\sqrt{5}$,
$ -1+\sqrt{5}$,
$0$,
$0$;\ \ 
$ 2c_{15}^{1}$,
$ 2c_{15}^{4}$,
$ 2c_{15}^{7}$,
$ 2c_{15}^{2}$,
$0$,
$0$;\ \ 
$ 2c_{15}^{1}$,
$ 2c_{15}^{2}$,
$ 2c_{15}^{7}$,
$0$,
$0$;\ \ 
$ 2c_{15}^{4}$,
$ 2c_{15}^{1}$,
$0$,
$0$;\ \ 
$ 2c_{15}^{4}$,
$0$,
$0$;\ \ 
$ -\sqrt{15}$,
$ \sqrt{15}$;\ \ 
$ -\sqrt{15}$)

  \vskip 2ex

\noindent16. $11_{6,60.}^{120,369}$ \irep{2274}:\ \ 
$d_i$ = ($1.0$,
$1.0$,
$2.0$,
$2.0$,
$2.0$,
$2.0$,
$2.0$,
$2.0$,
$2.0$,
$3.872$,
$3.872$) 

\vskip 0.7ex
\hangindent=3em \hangafter=1
$D^2= 60.0 = 
60$

\vskip 0.7ex
\hangindent=3em \hangafter=1
$T = ( 0,
0,
\frac{2}{3},
\frac{1}{5},
\frac{4}{5},
\frac{7}{15},
\frac{7}{15},
\frac{13}{15},
\frac{13}{15},
\frac{3}{8},
\frac{7}{8} )
$,

\vskip 0.7ex
\hangindent=3em \hangafter=1
$S$ = ($ 1$,
$ 1$,
$ 2$,
$ 2$,
$ 2$,
$ 2$,
$ 2$,
$ 2$,
$ 2$,
$ \sqrt{15}$,
$ \sqrt{15}$;\ \ 
$ 1$,
$ 2$,
$ 2$,
$ 2$,
$ 2$,
$ 2$,
$ 2$,
$ 2$,
$ -\sqrt{15}$,
$ -\sqrt{15}$;\ \ 
$ -2$,
$ 4$,
$ 4$,
$ -2$,
$ -2$,
$ -2$,
$ -2$,
$0$,
$0$;\ \ 
$ -1-\sqrt{5}$,
$ -1+\sqrt{5}$,
$ -1+\sqrt{5}$,
$ -1+\sqrt{5}$,
$ -1-\sqrt{5}$,
$ -1-\sqrt{5}$,
$0$,
$0$;\ \ 
$ -1-\sqrt{5}$,
$ -1-\sqrt{5}$,
$ -1-\sqrt{5}$,
$ -1+\sqrt{5}$,
$ -1+\sqrt{5}$,
$0$,
$0$;\ \ 
$ 2c_{15}^{1}$,
$ 2c_{15}^{4}$,
$ 2c_{15}^{7}$,
$ 2c_{15}^{2}$,
$0$,
$0$;\ \ 
$ 2c_{15}^{1}$,
$ 2c_{15}^{2}$,
$ 2c_{15}^{7}$,
$0$,
$0$;\ \ 
$ 2c_{15}^{4}$,
$ 2c_{15}^{1}$,
$0$,
$0$;\ \ 
$ 2c_{15}^{4}$,
$0$,
$0$;\ \ 
$ \sqrt{15}$,
$ -\sqrt{15}$;\ \ 
$ \sqrt{15}$)

  \vskip 2ex

\noindent17. $11_{2,60.}^{120,213}$ \irep{2274}:\ \ 
$d_i$ = ($1.0$,
$1.0$,
$2.0$,
$2.0$,
$2.0$,
$2.0$,
$2.0$,
$2.0$,
$2.0$,
$3.872$,
$3.872$) 

\vskip 0.7ex
\hangindent=3em \hangafter=1
$D^2= 60.0 = 
60$

\vskip 0.7ex
\hangindent=3em \hangafter=1
$T = ( 0,
0,
\frac{2}{3},
\frac{2}{5},
\frac{3}{5},
\frac{1}{15},
\frac{1}{15},
\frac{4}{15},
\frac{4}{15},
\frac{1}{8},
\frac{5}{8} )
$,

\vskip 0.7ex
\hangindent=3em \hangafter=1
$S$ = ($ 1$,
$ 1$,
$ 2$,
$ 2$,
$ 2$,
$ 2$,
$ 2$,
$ 2$,
$ 2$,
$ \sqrt{15}$,
$ \sqrt{15}$;\ \ 
$ 1$,
$ 2$,
$ 2$,
$ 2$,
$ 2$,
$ 2$,
$ 2$,
$ 2$,
$ -\sqrt{15}$,
$ -\sqrt{15}$;\ \ 
$ -2$,
$ 4$,
$ 4$,
$ -2$,
$ -2$,
$ -2$,
$ -2$,
$0$,
$0$;\ \ 
$ -1+\sqrt{5}$,
$ -1-\sqrt{5}$,
$ -1+\sqrt{5}$,
$ -1+\sqrt{5}$,
$ -1-\sqrt{5}$,
$ -1-\sqrt{5}$,
$0$,
$0$;\ \ 
$ -1+\sqrt{5}$,
$ -1-\sqrt{5}$,
$ -1-\sqrt{5}$,
$ -1+\sqrt{5}$,
$ -1+\sqrt{5}$,
$0$,
$0$;\ \ 
$ 2c_{15}^{2}$,
$ 2c_{15}^{7}$,
$ 2c_{15}^{4}$,
$ 2c_{15}^{1}$,
$0$,
$0$;\ \ 
$ 2c_{15}^{2}$,
$ 2c_{15}^{1}$,
$ 2c_{15}^{4}$,
$0$,
$0$;\ \ 
$ 2c_{15}^{7}$,
$ 2c_{15}^{2}$,
$0$,
$0$;\ \ 
$ 2c_{15}^{7}$,
$0$,
$0$;\ \ 
$ \sqrt{15}$,
$ -\sqrt{15}$;\ \ 
$ \sqrt{15}$)

  \vskip 2ex

\noindent18. $11_{2,60.}^{120,195}$ \irep{2274}:\ \ 
$d_i$ = ($1.0$,
$1.0$,
$2.0$,
$2.0$,
$2.0$,
$2.0$,
$2.0$,
$2.0$,
$2.0$,
$3.872$,
$3.872$) 

\vskip 0.7ex
\hangindent=3em \hangafter=1
$D^2= 60.0 = 
60$

\vskip 0.7ex
\hangindent=3em \hangafter=1
$T = ( 0,
0,
\frac{2}{3},
\frac{2}{5},
\frac{3}{5},
\frac{1}{15},
\frac{1}{15},
\frac{4}{15},
\frac{4}{15},
\frac{3}{8},
\frac{7}{8} )
$,

\vskip 0.7ex
\hangindent=3em \hangafter=1
$S$ = ($ 1$,
$ 1$,
$ 2$,
$ 2$,
$ 2$,
$ 2$,
$ 2$,
$ 2$,
$ 2$,
$ \sqrt{15}$,
$ \sqrt{15}$;\ \ 
$ 1$,
$ 2$,
$ 2$,
$ 2$,
$ 2$,
$ 2$,
$ 2$,
$ 2$,
$ -\sqrt{15}$,
$ -\sqrt{15}$;\ \ 
$ -2$,
$ 4$,
$ 4$,
$ -2$,
$ -2$,
$ -2$,
$ -2$,
$0$,
$0$;\ \ 
$ -1+\sqrt{5}$,
$ -1-\sqrt{5}$,
$ -1+\sqrt{5}$,
$ -1+\sqrt{5}$,
$ -1-\sqrt{5}$,
$ -1-\sqrt{5}$,
$0$,
$0$;\ \ 
$ -1+\sqrt{5}$,
$ -1-\sqrt{5}$,
$ -1-\sqrt{5}$,
$ -1+\sqrt{5}$,
$ -1+\sqrt{5}$,
$0$,
$0$;\ \ 
$ 2c_{15}^{2}$,
$ 2c_{15}^{7}$,
$ 2c_{15}^{4}$,
$ 2c_{15}^{1}$,
$0$,
$0$;\ \ 
$ 2c_{15}^{2}$,
$ 2c_{15}^{1}$,
$ 2c_{15}^{4}$,
$0$,
$0$;\ \ 
$ 2c_{15}^{7}$,
$ 2c_{15}^{2}$,
$0$,
$0$;\ \ 
$ 2c_{15}^{7}$,
$0$,
$0$;\ \ 
$ -\sqrt{15}$,
$ \sqrt{15}$;\ \ 
$ -\sqrt{15}$)

  \vskip 2ex

\noindent19. $11_{\frac{13}{2},89.56}^{48,108}$ \irep{2191}:\ \ 
$d_i$ = ($1.0$,
$1.0$,
$1.931$,
$1.931$,
$2.732$,
$2.732$,
$3.346$,
$3.346$,
$3.732$,
$3.732$,
$3.863$) 

\vskip 0.7ex
\hangindent=3em \hangafter=1
$D^2= 89.569 = 
48+24\sqrt{3}$

\vskip 0.7ex
\hangindent=3em \hangafter=1
$T = ( 0,
\frac{1}{2},
\frac{1}{16},
\frac{1}{16},
\frac{1}{3},
\frac{5}{6},
\frac{13}{16},
\frac{13}{16},
0,
\frac{1}{2},
\frac{19}{48} )
$,

\vskip 0.7ex
\hangindent=3em \hangafter=1
$S$ = ($ 1$,
$ 1$,
$ c_{24}^{1}$,
$ c_{24}^{1}$,
$ 1+\sqrt{3}$,
$ 1+\sqrt{3}$,
$ \frac{3+3\sqrt{3}}{\sqrt{6}}$,
$ \frac{3+3\sqrt{3}}{\sqrt{6}}$,
$ 2+\sqrt{3}$,
$ 2+\sqrt{3}$,
$ 2c_{24}^{1}$;\ \ 
$ 1$,
$ -c_{24}^{1}$,
$ -c_{24}^{1}$,
$ 1+\sqrt{3}$,
$ 1+\sqrt{3}$,
$ \frac{-3-3\sqrt{3}}{\sqrt{6}}$,
$ \frac{-3-3\sqrt{3}}{\sqrt{6}}$,
$ 2+\sqrt{3}$,
$ 2+\sqrt{3}$,
$ -2c_{24}^{1}$;\ \ 
$(\frac{-3-3\sqrt{3}}{\sqrt{6}})\mathrm{i}$,
$(\frac{3+3\sqrt{3}}{\sqrt{6}})\mathrm{i}$,
$ -2c_{24}^{1}$,
$ 2c_{24}^{1}$,
$(\frac{3+3\sqrt{3}}{\sqrt{6}})\mathrm{i}$,
$(\frac{-3-3\sqrt{3}}{\sqrt{6}})\mathrm{i}$,
$ -c_{24}^{1}$,
$ c_{24}^{1}$,
$0$;\ \ 
$(\frac{-3-3\sqrt{3}}{\sqrt{6}})\mathrm{i}$,
$ -2c_{24}^{1}$,
$ 2c_{24}^{1}$,
$(\frac{-3-3\sqrt{3}}{\sqrt{6}})\mathrm{i}$,
$(\frac{3+3\sqrt{3}}{\sqrt{6}})\mathrm{i}$,
$ -c_{24}^{1}$,
$ c_{24}^{1}$,
$0$;\ \ 
$ 1+\sqrt{3}$,
$ 1+\sqrt{3}$,
$0$,
$0$,
$ -1-\sqrt{3}$,
$ -1-\sqrt{3}$,
$ 2c_{24}^{1}$;\ \ 
$ 1+\sqrt{3}$,
$0$,
$0$,
$ -1-\sqrt{3}$,
$ -1-\sqrt{3}$,
$ -2c_{24}^{1}$;\ \ 
$(\frac{3+3\sqrt{3}}{\sqrt{6}})\mathrm{i}$,
$(\frac{-3-3\sqrt{3}}{\sqrt{6}})\mathrm{i}$,
$ \frac{3+3\sqrt{3}}{\sqrt{6}}$,
$ \frac{-3-3\sqrt{3}}{\sqrt{6}}$,
$0$;\ \ 
$(\frac{3+3\sqrt{3}}{\sqrt{6}})\mathrm{i}$,
$ \frac{3+3\sqrt{3}}{\sqrt{6}}$,
$ \frac{-3-3\sqrt{3}}{\sqrt{6}}$,
$0$;\ \ 
$ 1$,
$ 1$,
$ -2c_{24}^{1}$;\ \ 
$ 1$,
$ 2c_{24}^{1}$;\ \ 
$0$)

  \vskip 2ex

\noindent20. $11_{\frac{5}{2},89.56}^{48,214}$ \irep{2192}:\ \ 
$d_i$ = ($1.0$,
$1.0$,
$1.931$,
$1.931$,
$2.732$,
$2.732$,
$3.346$,
$3.346$,
$3.732$,
$3.732$,
$3.863$) 

\vskip 0.7ex
\hangindent=3em \hangafter=1
$D^2= 89.569 = 
48+24\sqrt{3}$

\vskip 0.7ex
\hangindent=3em \hangafter=1
$T = ( 0,
\frac{1}{2},
\frac{1}{16},
\frac{1}{16},
\frac{2}{3},
\frac{1}{6},
\frac{5}{16},
\frac{5}{16},
0,
\frac{1}{2},
\frac{35}{48} )
$,

\vskip 0.7ex
\hangindent=3em \hangafter=1
$S$ = ($ 1$,
$ 1$,
$ c_{24}^{1}$,
$ c_{24}^{1}$,
$ 1+\sqrt{3}$,
$ 1+\sqrt{3}$,
$ \frac{3+3\sqrt{3}}{\sqrt{6}}$,
$ \frac{3+3\sqrt{3}}{\sqrt{6}}$,
$ 2+\sqrt{3}$,
$ 2+\sqrt{3}$,
$ 2c_{24}^{1}$;\ \ 
$ 1$,
$ -c_{24}^{1}$,
$ -c_{24}^{1}$,
$ 1+\sqrt{3}$,
$ 1+\sqrt{3}$,
$ \frac{-3-3\sqrt{3}}{\sqrt{6}}$,
$ \frac{-3-3\sqrt{3}}{\sqrt{6}}$,
$ 2+\sqrt{3}$,
$ 2+\sqrt{3}$,
$ -2c_{24}^{1}$;\ \ 
$ \frac{3+3\sqrt{3}}{\sqrt{6}}$,
$ \frac{-3-3\sqrt{3}}{\sqrt{6}}$,
$ -2c_{24}^{1}$,
$ 2c_{24}^{1}$,
$ \frac{-3-3\sqrt{3}}{\sqrt{6}}$,
$ \frac{3+3\sqrt{3}}{\sqrt{6}}$,
$ -c_{24}^{1}$,
$ c_{24}^{1}$,
$0$;\ \ 
$ \frac{3+3\sqrt{3}}{\sqrt{6}}$,
$ -2c_{24}^{1}$,
$ 2c_{24}^{1}$,
$ \frac{3+3\sqrt{3}}{\sqrt{6}}$,
$ \frac{-3-3\sqrt{3}}{\sqrt{6}}$,
$ -c_{24}^{1}$,
$ c_{24}^{1}$,
$0$;\ \ 
$ 1+\sqrt{3}$,
$ 1+\sqrt{3}$,
$0$,
$0$,
$ -1-\sqrt{3}$,
$ -1-\sqrt{3}$,
$ 2c_{24}^{1}$;\ \ 
$ 1+\sqrt{3}$,
$0$,
$0$,
$ -1-\sqrt{3}$,
$ -1-\sqrt{3}$,
$ -2c_{24}^{1}$;\ \ 
$ \frac{-3-3\sqrt{3}}{\sqrt{6}}$,
$ \frac{3+3\sqrt{3}}{\sqrt{6}}$,
$ \frac{3+3\sqrt{3}}{\sqrt{6}}$,
$ \frac{-3-3\sqrt{3}}{\sqrt{6}}$,
$0$;\ \ 
$ \frac{-3-3\sqrt{3}}{\sqrt{6}}$,
$ \frac{3+3\sqrt{3}}{\sqrt{6}}$,
$ \frac{-3-3\sqrt{3}}{\sqrt{6}}$,
$0$;\ \ 
$ 1$,
$ 1$,
$ -2c_{24}^{1}$;\ \ 
$ 1$,
$ 2c_{24}^{1}$;\ \ 
$0$)

  \vskip 2ex

\noindent21. $11_{\frac{15}{2},89.56}^{48,311}$ \irep{2192}:\ \ 
$d_i$ = ($1.0$,
$1.0$,
$1.931$,
$1.931$,
$2.732$,
$2.732$,
$3.346$,
$3.346$,
$3.732$,
$3.732$,
$3.863$) 

\vskip 0.7ex
\hangindent=3em \hangafter=1
$D^2= 89.569 = 
48+24\sqrt{3}$

\vskip 0.7ex
\hangindent=3em \hangafter=1
$T = ( 0,
\frac{1}{2},
\frac{3}{16},
\frac{3}{16},
\frac{1}{3},
\frac{5}{6},
\frac{15}{16},
\frac{15}{16},
0,
\frac{1}{2},
\frac{25}{48} )
$,

\vskip 0.7ex
\hangindent=3em \hangafter=1
$S$ = ($ 1$,
$ 1$,
$ c_{24}^{1}$,
$ c_{24}^{1}$,
$ 1+\sqrt{3}$,
$ 1+\sqrt{3}$,
$ \frac{3+3\sqrt{3}}{\sqrt{6}}$,
$ \frac{3+3\sqrt{3}}{\sqrt{6}}$,
$ 2+\sqrt{3}$,
$ 2+\sqrt{3}$,
$ 2c_{24}^{1}$;\ \ 
$ 1$,
$ -c_{24}^{1}$,
$ -c_{24}^{1}$,
$ 1+\sqrt{3}$,
$ 1+\sqrt{3}$,
$ \frac{-3-3\sqrt{3}}{\sqrt{6}}$,
$ \frac{-3-3\sqrt{3}}{\sqrt{6}}$,
$ 2+\sqrt{3}$,
$ 2+\sqrt{3}$,
$ -2c_{24}^{1}$;\ \ 
$ \frac{-3-3\sqrt{3}}{\sqrt{6}}$,
$ \frac{3+3\sqrt{3}}{\sqrt{6}}$,
$ -2c_{24}^{1}$,
$ 2c_{24}^{1}$,
$ \frac{-3-3\sqrt{3}}{\sqrt{6}}$,
$ \frac{3+3\sqrt{3}}{\sqrt{6}}$,
$ -c_{24}^{1}$,
$ c_{24}^{1}$,
$0$;\ \ 
$ \frac{-3-3\sqrt{3}}{\sqrt{6}}$,
$ -2c_{24}^{1}$,
$ 2c_{24}^{1}$,
$ \frac{3+3\sqrt{3}}{\sqrt{6}}$,
$ \frac{-3-3\sqrt{3}}{\sqrt{6}}$,
$ -c_{24}^{1}$,
$ c_{24}^{1}$,
$0$;\ \ 
$ 1+\sqrt{3}$,
$ 1+\sqrt{3}$,
$0$,
$0$,
$ -1-\sqrt{3}$,
$ -1-\sqrt{3}$,
$ 2c_{24}^{1}$;\ \ 
$ 1+\sqrt{3}$,
$0$,
$0$,
$ -1-\sqrt{3}$,
$ -1-\sqrt{3}$,
$ -2c_{24}^{1}$;\ \ 
$ \frac{3+3\sqrt{3}}{\sqrt{6}}$,
$ \frac{-3-3\sqrt{3}}{\sqrt{6}}$,
$ \frac{3+3\sqrt{3}}{\sqrt{6}}$,
$ \frac{-3-3\sqrt{3}}{\sqrt{6}}$,
$0$;\ \ 
$ \frac{3+3\sqrt{3}}{\sqrt{6}}$,
$ \frac{3+3\sqrt{3}}{\sqrt{6}}$,
$ \frac{-3-3\sqrt{3}}{\sqrt{6}}$,
$0$;\ \ 
$ 1$,
$ 1$,
$ -2c_{24}^{1}$;\ \ 
$ 1$,
$ 2c_{24}^{1}$;\ \ 
$0$)

  \vskip 2ex

\noindent22. $11_{\frac{7}{2},89.56}^{48,628}$ \irep{2191}:\ \ 
$d_i$ = ($1.0$,
$1.0$,
$1.931$,
$1.931$,
$2.732$,
$2.732$,
$3.346$,
$3.346$,
$3.732$,
$3.732$,
$3.863$) 

\vskip 0.7ex
\hangindent=3em \hangafter=1
$D^2= 89.569 = 
48+24\sqrt{3}$

\vskip 0.7ex
\hangindent=3em \hangafter=1
$T = ( 0,
\frac{1}{2},
\frac{3}{16},
\frac{3}{16},
\frac{2}{3},
\frac{1}{6},
\frac{7}{16},
\frac{7}{16},
0,
\frac{1}{2},
\frac{41}{48} )
$,

\vskip 0.7ex
\hangindent=3em \hangafter=1
$S$ = ($ 1$,
$ 1$,
$ c_{24}^{1}$,
$ c_{24}^{1}$,
$ 1+\sqrt{3}$,
$ 1+\sqrt{3}$,
$ \frac{3+3\sqrt{3}}{\sqrt{6}}$,
$ \frac{3+3\sqrt{3}}{\sqrt{6}}$,
$ 2+\sqrt{3}$,
$ 2+\sqrt{3}$,
$ 2c_{24}^{1}$;\ \ 
$ 1$,
$ -c_{24}^{1}$,
$ -c_{24}^{1}$,
$ 1+\sqrt{3}$,
$ 1+\sqrt{3}$,
$ \frac{-3-3\sqrt{3}}{\sqrt{6}}$,
$ \frac{-3-3\sqrt{3}}{\sqrt{6}}$,
$ 2+\sqrt{3}$,
$ 2+\sqrt{3}$,
$ -2c_{24}^{1}$;\ \ 
$(\frac{-3-3\sqrt{3}}{\sqrt{6}})\mathrm{i}$,
$(\frac{3+3\sqrt{3}}{\sqrt{6}})\mathrm{i}$,
$ -2c_{24}^{1}$,
$ 2c_{24}^{1}$,
$(\frac{3+3\sqrt{3}}{\sqrt{6}})\mathrm{i}$,
$(\frac{-3-3\sqrt{3}}{\sqrt{6}})\mathrm{i}$,
$ -c_{24}^{1}$,
$ c_{24}^{1}$,
$0$;\ \ 
$(\frac{-3-3\sqrt{3}}{\sqrt{6}})\mathrm{i}$,
$ -2c_{24}^{1}$,
$ 2c_{24}^{1}$,
$(\frac{-3-3\sqrt{3}}{\sqrt{6}})\mathrm{i}$,
$(\frac{3+3\sqrt{3}}{\sqrt{6}})\mathrm{i}$,
$ -c_{24}^{1}$,
$ c_{24}^{1}$,
$0$;\ \ 
$ 1+\sqrt{3}$,
$ 1+\sqrt{3}$,
$0$,
$0$,
$ -1-\sqrt{3}$,
$ -1-\sqrt{3}$,
$ 2c_{24}^{1}$;\ \ 
$ 1+\sqrt{3}$,
$0$,
$0$,
$ -1-\sqrt{3}$,
$ -1-\sqrt{3}$,
$ -2c_{24}^{1}$;\ \ 
$(\frac{3+3\sqrt{3}}{\sqrt{6}})\mathrm{i}$,
$(\frac{-3-3\sqrt{3}}{\sqrt{6}})\mathrm{i}$,
$ \frac{3+3\sqrt{3}}{\sqrt{6}}$,
$ \frac{-3-3\sqrt{3}}{\sqrt{6}}$,
$0$;\ \ 
$(\frac{3+3\sqrt{3}}{\sqrt{6}})\mathrm{i}$,
$ \frac{3+3\sqrt{3}}{\sqrt{6}}$,
$ \frac{-3-3\sqrt{3}}{\sqrt{6}}$,
$0$;\ \ 
$ 1$,
$ 1$,
$ -2c_{24}^{1}$;\ \ 
$ 1$,
$ 2c_{24}^{1}$;\ \ 
$0$)

  \vskip 2ex

\noindent23. $11_{\frac{1}{2},89.56}^{48,193}$ \irep{2191}:\ \ 
$d_i$ = ($1.0$,
$1.0$,
$1.931$,
$1.931$,
$2.732$,
$2.732$,
$3.346$,
$3.346$,
$3.732$,
$3.732$,
$3.863$) 

\vskip 0.7ex
\hangindent=3em \hangafter=1
$D^2= 89.569 = 
48+24\sqrt{3}$

\vskip 0.7ex
\hangindent=3em \hangafter=1
$T = ( 0,
\frac{1}{2},
\frac{5}{16},
\frac{5}{16},
\frac{1}{3},
\frac{5}{6},
\frac{1}{16},
\frac{1}{16},
0,
\frac{1}{2},
\frac{31}{48} )
$,

\vskip 0.7ex
\hangindent=3em \hangafter=1
$S$ = ($ 1$,
$ 1$,
$ c_{24}^{1}$,
$ c_{24}^{1}$,
$ 1+\sqrt{3}$,
$ 1+\sqrt{3}$,
$ \frac{3+3\sqrt{3}}{\sqrt{6}}$,
$ \frac{3+3\sqrt{3}}{\sqrt{6}}$,
$ 2+\sqrt{3}$,
$ 2+\sqrt{3}$,
$ 2c_{24}^{1}$;\ \ 
$ 1$,
$ -c_{24}^{1}$,
$ -c_{24}^{1}$,
$ 1+\sqrt{3}$,
$ 1+\sqrt{3}$,
$ \frac{-3-3\sqrt{3}}{\sqrt{6}}$,
$ \frac{-3-3\sqrt{3}}{\sqrt{6}}$,
$ 2+\sqrt{3}$,
$ 2+\sqrt{3}$,
$ -2c_{24}^{1}$;\ \ 
$(\frac{3+3\sqrt{3}}{\sqrt{6}})\mathrm{i}$,
$(\frac{-3-3\sqrt{3}}{\sqrt{6}})\mathrm{i}$,
$ -2c_{24}^{1}$,
$ 2c_{24}^{1}$,
$(\frac{3+3\sqrt{3}}{\sqrt{6}})\mathrm{i}$,
$(\frac{-3-3\sqrt{3}}{\sqrt{6}})\mathrm{i}$,
$ -c_{24}^{1}$,
$ c_{24}^{1}$,
$0$;\ \ 
$(\frac{3+3\sqrt{3}}{\sqrt{6}})\mathrm{i}$,
$ -2c_{24}^{1}$,
$ 2c_{24}^{1}$,
$(\frac{-3-3\sqrt{3}}{\sqrt{6}})\mathrm{i}$,
$(\frac{3+3\sqrt{3}}{\sqrt{6}})\mathrm{i}$,
$ -c_{24}^{1}$,
$ c_{24}^{1}$,
$0$;\ \ 
$ 1+\sqrt{3}$,
$ 1+\sqrt{3}$,
$0$,
$0$,
$ -1-\sqrt{3}$,
$ -1-\sqrt{3}$,
$ 2c_{24}^{1}$;\ \ 
$ 1+\sqrt{3}$,
$0$,
$0$,
$ -1-\sqrt{3}$,
$ -1-\sqrt{3}$,
$ -2c_{24}^{1}$;\ \ 
$(\frac{-3-3\sqrt{3}}{\sqrt{6}})\mathrm{i}$,
$(\frac{3+3\sqrt{3}}{\sqrt{6}})\mathrm{i}$,
$ \frac{3+3\sqrt{3}}{\sqrt{6}}$,
$ \frac{-3-3\sqrt{3}}{\sqrt{6}}$,
$0$;\ \ 
$(\frac{-3-3\sqrt{3}}{\sqrt{6}})\mathrm{i}$,
$ \frac{3+3\sqrt{3}}{\sqrt{6}}$,
$ \frac{-3-3\sqrt{3}}{\sqrt{6}}$,
$0$;\ \ 
$ 1$,
$ 1$,
$ -2c_{24}^{1}$;\ \ 
$ 1$,
$ 2c_{24}^{1}$;\ \ 
$0$)

  \vskip 2ex

\noindent24. $11_{\frac{9}{2},89.56}^{48,133}$ \irep{2192}:\ \ 
$d_i$ = ($1.0$,
$1.0$,
$1.931$,
$1.931$,
$2.732$,
$2.732$,
$3.346$,
$3.346$,
$3.732$,
$3.732$,
$3.863$) 

\vskip 0.7ex
\hangindent=3em \hangafter=1
$D^2= 89.569 = 
48+24\sqrt{3}$

\vskip 0.7ex
\hangindent=3em \hangafter=1
$T = ( 0,
\frac{1}{2},
\frac{5}{16},
\frac{5}{16},
\frac{2}{3},
\frac{1}{6},
\frac{9}{16},
\frac{9}{16},
0,
\frac{1}{2},
\frac{47}{48} )
$,

\vskip 0.7ex
\hangindent=3em \hangafter=1
$S$ = ($ 1$,
$ 1$,
$ c_{24}^{1}$,
$ c_{24}^{1}$,
$ 1+\sqrt{3}$,
$ 1+\sqrt{3}$,
$ \frac{3+3\sqrt{3}}{\sqrt{6}}$,
$ \frac{3+3\sqrt{3}}{\sqrt{6}}$,
$ 2+\sqrt{3}$,
$ 2+\sqrt{3}$,
$ 2c_{24}^{1}$;\ \ 
$ 1$,
$ -c_{24}^{1}$,
$ -c_{24}^{1}$,
$ 1+\sqrt{3}$,
$ 1+\sqrt{3}$,
$ \frac{-3-3\sqrt{3}}{\sqrt{6}}$,
$ \frac{-3-3\sqrt{3}}{\sqrt{6}}$,
$ 2+\sqrt{3}$,
$ 2+\sqrt{3}$,
$ -2c_{24}^{1}$;\ \ 
$ \frac{-3-3\sqrt{3}}{\sqrt{6}}$,
$ \frac{3+3\sqrt{3}}{\sqrt{6}}$,
$ -2c_{24}^{1}$,
$ 2c_{24}^{1}$,
$ \frac{-3-3\sqrt{3}}{\sqrt{6}}$,
$ \frac{3+3\sqrt{3}}{\sqrt{6}}$,
$ -c_{24}^{1}$,
$ c_{24}^{1}$,
$0$;\ \ 
$ \frac{-3-3\sqrt{3}}{\sqrt{6}}$,
$ -2c_{24}^{1}$,
$ 2c_{24}^{1}$,
$ \frac{3+3\sqrt{3}}{\sqrt{6}}$,
$ \frac{-3-3\sqrt{3}}{\sqrt{6}}$,
$ -c_{24}^{1}$,
$ c_{24}^{1}$,
$0$;\ \ 
$ 1+\sqrt{3}$,
$ 1+\sqrt{3}$,
$0$,
$0$,
$ -1-\sqrt{3}$,
$ -1-\sqrt{3}$,
$ 2c_{24}^{1}$;\ \ 
$ 1+\sqrt{3}$,
$0$,
$0$,
$ -1-\sqrt{3}$,
$ -1-\sqrt{3}$,
$ -2c_{24}^{1}$;\ \ 
$ \frac{3+3\sqrt{3}}{\sqrt{6}}$,
$ \frac{-3-3\sqrt{3}}{\sqrt{6}}$,
$ \frac{3+3\sqrt{3}}{\sqrt{6}}$,
$ \frac{-3-3\sqrt{3}}{\sqrt{6}}$,
$0$;\ \ 
$ \frac{3+3\sqrt{3}}{\sqrt{6}}$,
$ \frac{3+3\sqrt{3}}{\sqrt{6}}$,
$ \frac{-3-3\sqrt{3}}{\sqrt{6}}$,
$0$;\ \ 
$ 1$,
$ 1$,
$ -2c_{24}^{1}$;\ \ 
$ 1$,
$ 2c_{24}^{1}$;\ \ 
$0$)

  \vskip 2ex

\noindent25. $11_{\frac{3}{2},89.56}^{48,682}$ \irep{2192}:\ \ 
$d_i$ = ($1.0$,
$1.0$,
$1.931$,
$1.931$,
$2.732$,
$2.732$,
$3.346$,
$3.346$,
$3.732$,
$3.732$,
$3.863$) 

\vskip 0.7ex
\hangindent=3em \hangafter=1
$D^2= 89.569 = 
48+24\sqrt{3}$

\vskip 0.7ex
\hangindent=3em \hangafter=1
$T = ( 0,
\frac{1}{2},
\frac{7}{16},
\frac{7}{16},
\frac{1}{3},
\frac{5}{6},
\frac{3}{16},
\frac{3}{16},
0,
\frac{1}{2},
\frac{37}{48} )
$,

\vskip 0.7ex
\hangindent=3em \hangafter=1
$S$ = ($ 1$,
$ 1$,
$ c_{24}^{1}$,
$ c_{24}^{1}$,
$ 1+\sqrt{3}$,
$ 1+\sqrt{3}$,
$ \frac{3+3\sqrt{3}}{\sqrt{6}}$,
$ \frac{3+3\sqrt{3}}{\sqrt{6}}$,
$ 2+\sqrt{3}$,
$ 2+\sqrt{3}$,
$ 2c_{24}^{1}$;\ \ 
$ 1$,
$ -c_{24}^{1}$,
$ -c_{24}^{1}$,
$ 1+\sqrt{3}$,
$ 1+\sqrt{3}$,
$ \frac{-3-3\sqrt{3}}{\sqrt{6}}$,
$ \frac{-3-3\sqrt{3}}{\sqrt{6}}$,
$ 2+\sqrt{3}$,
$ 2+\sqrt{3}$,
$ -2c_{24}^{1}$;\ \ 
$ \frac{3+3\sqrt{3}}{\sqrt{6}}$,
$ \frac{-3-3\sqrt{3}}{\sqrt{6}}$,
$ -2c_{24}^{1}$,
$ 2c_{24}^{1}$,
$ \frac{-3-3\sqrt{3}}{\sqrt{6}}$,
$ \frac{3+3\sqrt{3}}{\sqrt{6}}$,
$ -c_{24}^{1}$,
$ c_{24}^{1}$,
$0$;\ \ 
$ \frac{3+3\sqrt{3}}{\sqrt{6}}$,
$ -2c_{24}^{1}$,
$ 2c_{24}^{1}$,
$ \frac{3+3\sqrt{3}}{\sqrt{6}}$,
$ \frac{-3-3\sqrt{3}}{\sqrt{6}}$,
$ -c_{24}^{1}$,
$ c_{24}^{1}$,
$0$;\ \ 
$ 1+\sqrt{3}$,
$ 1+\sqrt{3}$,
$0$,
$0$,
$ -1-\sqrt{3}$,
$ -1-\sqrt{3}$,
$ 2c_{24}^{1}$;\ \ 
$ 1+\sqrt{3}$,
$0$,
$0$,
$ -1-\sqrt{3}$,
$ -1-\sqrt{3}$,
$ -2c_{24}^{1}$;\ \ 
$ \frac{-3-3\sqrt{3}}{\sqrt{6}}$,
$ \frac{3+3\sqrt{3}}{\sqrt{6}}$,
$ \frac{3+3\sqrt{3}}{\sqrt{6}}$,
$ \frac{-3-3\sqrt{3}}{\sqrt{6}}$,
$0$;\ \ 
$ \frac{-3-3\sqrt{3}}{\sqrt{6}}$,
$ \frac{3+3\sqrt{3}}{\sqrt{6}}$,
$ \frac{-3-3\sqrt{3}}{\sqrt{6}}$,
$0$;\ \ 
$ 1$,
$ 1$,
$ -2c_{24}^{1}$;\ \ 
$ 1$,
$ 2c_{24}^{1}$;\ \ 
$0$)

  \vskip 2ex

\noindent26. $11_{\frac{11}{2},89.56}^{48,919}$ \irep{2191}:\ \ 
$d_i$ = ($1.0$,
$1.0$,
$1.931$,
$1.931$,
$2.732$,
$2.732$,
$3.346$,
$3.346$,
$3.732$,
$3.732$,
$3.863$) 

\vskip 0.7ex
\hangindent=3em \hangafter=1
$D^2= 89.569 = 
48+24\sqrt{3}$

\vskip 0.7ex
\hangindent=3em \hangafter=1
$T = ( 0,
\frac{1}{2},
\frac{7}{16},
\frac{7}{16},
\frac{2}{3},
\frac{1}{6},
\frac{11}{16},
\frac{11}{16},
0,
\frac{1}{2},
\frac{5}{48} )
$,

\vskip 0.7ex
\hangindent=3em \hangafter=1
$S$ = ($ 1$,
$ 1$,
$ c_{24}^{1}$,
$ c_{24}^{1}$,
$ 1+\sqrt{3}$,
$ 1+\sqrt{3}$,
$ \frac{3+3\sqrt{3}}{\sqrt{6}}$,
$ \frac{3+3\sqrt{3}}{\sqrt{6}}$,
$ 2+\sqrt{3}$,
$ 2+\sqrt{3}$,
$ 2c_{24}^{1}$;\ \ 
$ 1$,
$ -c_{24}^{1}$,
$ -c_{24}^{1}$,
$ 1+\sqrt{3}$,
$ 1+\sqrt{3}$,
$ \frac{-3-3\sqrt{3}}{\sqrt{6}}$,
$ \frac{-3-3\sqrt{3}}{\sqrt{6}}$,
$ 2+\sqrt{3}$,
$ 2+\sqrt{3}$,
$ -2c_{24}^{1}$;\ \ 
$(\frac{3+3\sqrt{3}}{\sqrt{6}})\mathrm{i}$,
$(\frac{-3-3\sqrt{3}}{\sqrt{6}})\mathrm{i}$,
$ -2c_{24}^{1}$,
$ 2c_{24}^{1}$,
$(\frac{3+3\sqrt{3}}{\sqrt{6}})\mathrm{i}$,
$(\frac{-3-3\sqrt{3}}{\sqrt{6}})\mathrm{i}$,
$ -c_{24}^{1}$,
$ c_{24}^{1}$,
$0$;\ \ 
$(\frac{3+3\sqrt{3}}{\sqrt{6}})\mathrm{i}$,
$ -2c_{24}^{1}$,
$ 2c_{24}^{1}$,
$(\frac{-3-3\sqrt{3}}{\sqrt{6}})\mathrm{i}$,
$(\frac{3+3\sqrt{3}}{\sqrt{6}})\mathrm{i}$,
$ -c_{24}^{1}$,
$ c_{24}^{1}$,
$0$;\ \ 
$ 1+\sqrt{3}$,
$ 1+\sqrt{3}$,
$0$,
$0$,
$ -1-\sqrt{3}$,
$ -1-\sqrt{3}$,
$ 2c_{24}^{1}$;\ \ 
$ 1+\sqrt{3}$,
$0$,
$0$,
$ -1-\sqrt{3}$,
$ -1-\sqrt{3}$,
$ -2c_{24}^{1}$;\ \ 
$(\frac{-3-3\sqrt{3}}{\sqrt{6}})\mathrm{i}$,
$(\frac{3+3\sqrt{3}}{\sqrt{6}})\mathrm{i}$,
$ \frac{3+3\sqrt{3}}{\sqrt{6}}$,
$ \frac{-3-3\sqrt{3}}{\sqrt{6}}$,
$0$;\ \ 
$(\frac{-3-3\sqrt{3}}{\sqrt{6}})\mathrm{i}$,
$ \frac{3+3\sqrt{3}}{\sqrt{6}}$,
$ \frac{-3-3\sqrt{3}}{\sqrt{6}}$,
$0$;\ \ 
$ 1$,
$ 1$,
$ -2c_{24}^{1}$;\ \ 
$ 1$,
$ 2c_{24}^{1}$;\ \ 
$0$)

  \vskip 2ex

\noindent27. $11_{\frac{5}{2},89.56}^{48,102}$ \irep{2191}:\ \ 
$d_i$ = ($1.0$,
$1.0$,
$1.931$,
$1.931$,
$2.732$,
$2.732$,
$3.346$,
$3.346$,
$3.732$,
$3.732$,
$3.863$) 

\vskip 0.7ex
\hangindent=3em \hangafter=1
$D^2= 89.569 = 
48+24\sqrt{3}$

\vskip 0.7ex
\hangindent=3em \hangafter=1
$T = ( 0,
\frac{1}{2},
\frac{9}{16},
\frac{9}{16},
\frac{1}{3},
\frac{5}{6},
\frac{5}{16},
\frac{5}{16},
0,
\frac{1}{2},
\frac{43}{48} )
$,

\vskip 0.7ex
\hangindent=3em \hangafter=1
$S$ = ($ 1$,
$ 1$,
$ c_{24}^{1}$,
$ c_{24}^{1}$,
$ 1+\sqrt{3}$,
$ 1+\sqrt{3}$,
$ \frac{3+3\sqrt{3}}{\sqrt{6}}$,
$ \frac{3+3\sqrt{3}}{\sqrt{6}}$,
$ 2+\sqrt{3}$,
$ 2+\sqrt{3}$,
$ 2c_{24}^{1}$;\ \ 
$ 1$,
$ -c_{24}^{1}$,
$ -c_{24}^{1}$,
$ 1+\sqrt{3}$,
$ 1+\sqrt{3}$,
$ \frac{-3-3\sqrt{3}}{\sqrt{6}}$,
$ \frac{-3-3\sqrt{3}}{\sqrt{6}}$,
$ 2+\sqrt{3}$,
$ 2+\sqrt{3}$,
$ -2c_{24}^{1}$;\ \ 
$(\frac{-3-3\sqrt{3}}{\sqrt{6}})\mathrm{i}$,
$(\frac{3+3\sqrt{3}}{\sqrt{6}})\mathrm{i}$,
$ -2c_{24}^{1}$,
$ 2c_{24}^{1}$,
$(\frac{3+3\sqrt{3}}{\sqrt{6}})\mathrm{i}$,
$(\frac{-3-3\sqrt{3}}{\sqrt{6}})\mathrm{i}$,
$ -c_{24}^{1}$,
$ c_{24}^{1}$,
$0$;\ \ 
$(\frac{-3-3\sqrt{3}}{\sqrt{6}})\mathrm{i}$,
$ -2c_{24}^{1}$,
$ 2c_{24}^{1}$,
$(\frac{-3-3\sqrt{3}}{\sqrt{6}})\mathrm{i}$,
$(\frac{3+3\sqrt{3}}{\sqrt{6}})\mathrm{i}$,
$ -c_{24}^{1}$,
$ c_{24}^{1}$,
$0$;\ \ 
$ 1+\sqrt{3}$,
$ 1+\sqrt{3}$,
$0$,
$0$,
$ -1-\sqrt{3}$,
$ -1-\sqrt{3}$,
$ 2c_{24}^{1}$;\ \ 
$ 1+\sqrt{3}$,
$0$,
$0$,
$ -1-\sqrt{3}$,
$ -1-\sqrt{3}$,
$ -2c_{24}^{1}$;\ \ 
$(\frac{3+3\sqrt{3}}{\sqrt{6}})\mathrm{i}$,
$(\frac{-3-3\sqrt{3}}{\sqrt{6}})\mathrm{i}$,
$ \frac{3+3\sqrt{3}}{\sqrt{6}}$,
$ \frac{-3-3\sqrt{3}}{\sqrt{6}}$,
$0$;\ \ 
$(\frac{3+3\sqrt{3}}{\sqrt{6}})\mathrm{i}$,
$ \frac{3+3\sqrt{3}}{\sqrt{6}}$,
$ \frac{-3-3\sqrt{3}}{\sqrt{6}}$,
$0$;\ \ 
$ 1$,
$ 1$,
$ -2c_{24}^{1}$;\ \ 
$ 1$,
$ 2c_{24}^{1}$;\ \ 
$0$)

  \vskip 2ex

\noindent28. $11_{\frac{13}{2},89.56}^{48,979}$ \irep{2192}:\ \ 
$d_i$ = ($1.0$,
$1.0$,
$1.931$,
$1.931$,
$2.732$,
$2.732$,
$3.346$,
$3.346$,
$3.732$,
$3.732$,
$3.863$) 

\vskip 0.7ex
\hangindent=3em \hangafter=1
$D^2= 89.569 = 
48+24\sqrt{3}$

\vskip 0.7ex
\hangindent=3em \hangafter=1
$T = ( 0,
\frac{1}{2},
\frac{9}{16},
\frac{9}{16},
\frac{2}{3},
\frac{1}{6},
\frac{13}{16},
\frac{13}{16},
0,
\frac{1}{2},
\frac{11}{48} )
$,

\vskip 0.7ex
\hangindent=3em \hangafter=1
$S$ = ($ 1$,
$ 1$,
$ c_{24}^{1}$,
$ c_{24}^{1}$,
$ 1+\sqrt{3}$,
$ 1+\sqrt{3}$,
$ \frac{3+3\sqrt{3}}{\sqrt{6}}$,
$ \frac{3+3\sqrt{3}}{\sqrt{6}}$,
$ 2+\sqrt{3}$,
$ 2+\sqrt{3}$,
$ 2c_{24}^{1}$;\ \ 
$ 1$,
$ -c_{24}^{1}$,
$ -c_{24}^{1}$,
$ 1+\sqrt{3}$,
$ 1+\sqrt{3}$,
$ \frac{-3-3\sqrt{3}}{\sqrt{6}}$,
$ \frac{-3-3\sqrt{3}}{\sqrt{6}}$,
$ 2+\sqrt{3}$,
$ 2+\sqrt{3}$,
$ -2c_{24}^{1}$;\ \ 
$ \frac{3+3\sqrt{3}}{\sqrt{6}}$,
$ \frac{-3-3\sqrt{3}}{\sqrt{6}}$,
$ -2c_{24}^{1}$,
$ 2c_{24}^{1}$,
$ \frac{-3-3\sqrt{3}}{\sqrt{6}}$,
$ \frac{3+3\sqrt{3}}{\sqrt{6}}$,
$ -c_{24}^{1}$,
$ c_{24}^{1}$,
$0$;\ \ 
$ \frac{3+3\sqrt{3}}{\sqrt{6}}$,
$ -2c_{24}^{1}$,
$ 2c_{24}^{1}$,
$ \frac{3+3\sqrt{3}}{\sqrt{6}}$,
$ \frac{-3-3\sqrt{3}}{\sqrt{6}}$,
$ -c_{24}^{1}$,
$ c_{24}^{1}$,
$0$;\ \ 
$ 1+\sqrt{3}$,
$ 1+\sqrt{3}$,
$0$,
$0$,
$ -1-\sqrt{3}$,
$ -1-\sqrt{3}$,
$ 2c_{24}^{1}$;\ \ 
$ 1+\sqrt{3}$,
$0$,
$0$,
$ -1-\sqrt{3}$,
$ -1-\sqrt{3}$,
$ -2c_{24}^{1}$;\ \ 
$ \frac{-3-3\sqrt{3}}{\sqrt{6}}$,
$ \frac{3+3\sqrt{3}}{\sqrt{6}}$,
$ \frac{3+3\sqrt{3}}{\sqrt{6}}$,
$ \frac{-3-3\sqrt{3}}{\sqrt{6}}$,
$0$;\ \ 
$ \frac{-3-3\sqrt{3}}{\sqrt{6}}$,
$ \frac{3+3\sqrt{3}}{\sqrt{6}}$,
$ \frac{-3-3\sqrt{3}}{\sqrt{6}}$,
$0$;\ \ 
$ 1$,
$ 1$,
$ -2c_{24}^{1}$;\ \ 
$ 1$,
$ 2c_{24}^{1}$;\ \ 
$0$)

  \vskip 2ex

\noindent29. $11_{\frac{7}{2},89.56}^{48,844}$ \irep{2192}:\ \ 
$d_i$ = ($1.0$,
$1.0$,
$1.931$,
$1.931$,
$2.732$,
$2.732$,
$3.346$,
$3.346$,
$3.732$,
$3.732$,
$3.863$) 

\vskip 0.7ex
\hangindent=3em \hangafter=1
$D^2= 89.569 = 
48+24\sqrt{3}$

\vskip 0.7ex
\hangindent=3em \hangafter=1
$T = ( 0,
\frac{1}{2},
\frac{11}{16},
\frac{11}{16},
\frac{1}{3},
\frac{5}{6},
\frac{7}{16},
\frac{7}{16},
0,
\frac{1}{2},
\frac{1}{48} )
$,

\vskip 0.7ex
\hangindent=3em \hangafter=1
$S$ = ($ 1$,
$ 1$,
$ c_{24}^{1}$,
$ c_{24}^{1}$,
$ 1+\sqrt{3}$,
$ 1+\sqrt{3}$,
$ \frac{3+3\sqrt{3}}{\sqrt{6}}$,
$ \frac{3+3\sqrt{3}}{\sqrt{6}}$,
$ 2+\sqrt{3}$,
$ 2+\sqrt{3}$,
$ 2c_{24}^{1}$;\ \ 
$ 1$,
$ -c_{24}^{1}$,
$ -c_{24}^{1}$,
$ 1+\sqrt{3}$,
$ 1+\sqrt{3}$,
$ \frac{-3-3\sqrt{3}}{\sqrt{6}}$,
$ \frac{-3-3\sqrt{3}}{\sqrt{6}}$,
$ 2+\sqrt{3}$,
$ 2+\sqrt{3}$,
$ -2c_{24}^{1}$;\ \ 
$ \frac{-3-3\sqrt{3}}{\sqrt{6}}$,
$ \frac{3+3\sqrt{3}}{\sqrt{6}}$,
$ -2c_{24}^{1}$,
$ 2c_{24}^{1}$,
$ \frac{-3-3\sqrt{3}}{\sqrt{6}}$,
$ \frac{3+3\sqrt{3}}{\sqrt{6}}$,
$ -c_{24}^{1}$,
$ c_{24}^{1}$,
$0$;\ \ 
$ \frac{-3-3\sqrt{3}}{\sqrt{6}}$,
$ -2c_{24}^{1}$,
$ 2c_{24}^{1}$,
$ \frac{3+3\sqrt{3}}{\sqrt{6}}$,
$ \frac{-3-3\sqrt{3}}{\sqrt{6}}$,
$ -c_{24}^{1}$,
$ c_{24}^{1}$,
$0$;\ \ 
$ 1+\sqrt{3}$,
$ 1+\sqrt{3}$,
$0$,
$0$,
$ -1-\sqrt{3}$,
$ -1-\sqrt{3}$,
$ 2c_{24}^{1}$;\ \ 
$ 1+\sqrt{3}$,
$0$,
$0$,
$ -1-\sqrt{3}$,
$ -1-\sqrt{3}$,
$ -2c_{24}^{1}$;\ \ 
$ \frac{3+3\sqrt{3}}{\sqrt{6}}$,
$ \frac{-3-3\sqrt{3}}{\sqrt{6}}$,
$ \frac{3+3\sqrt{3}}{\sqrt{6}}$,
$ \frac{-3-3\sqrt{3}}{\sqrt{6}}$,
$0$;\ \ 
$ \frac{3+3\sqrt{3}}{\sqrt{6}}$,
$ \frac{3+3\sqrt{3}}{\sqrt{6}}$,
$ \frac{-3-3\sqrt{3}}{\sqrt{6}}$,
$0$;\ \ 
$ 1$,
$ 1$,
$ -2c_{24}^{1}$;\ \ 
$ 1$,
$ 2c_{24}^{1}$;\ \ 
$0$)

  \vskip 2ex

\noindent30. $11_{\frac{15}{2},89.56}^{48,332}$ \irep{2191}:\ \ 
$d_i$ = ($1.0$,
$1.0$,
$1.931$,
$1.931$,
$2.732$,
$2.732$,
$3.346$,
$3.346$,
$3.732$,
$3.732$,
$3.863$) 

\vskip 0.7ex
\hangindent=3em \hangafter=1
$D^2= 89.569 = 
48+24\sqrt{3}$

\vskip 0.7ex
\hangindent=3em \hangafter=1
$T = ( 0,
\frac{1}{2},
\frac{11}{16},
\frac{11}{16},
\frac{2}{3},
\frac{1}{6},
\frac{15}{16},
\frac{15}{16},
0,
\frac{1}{2},
\frac{17}{48} )
$,

\vskip 0.7ex
\hangindent=3em \hangafter=1
$S$ = ($ 1$,
$ 1$,
$ c_{24}^{1}$,
$ c_{24}^{1}$,
$ 1+\sqrt{3}$,
$ 1+\sqrt{3}$,
$ \frac{3+3\sqrt{3}}{\sqrt{6}}$,
$ \frac{3+3\sqrt{3}}{\sqrt{6}}$,
$ 2+\sqrt{3}$,
$ 2+\sqrt{3}$,
$ 2c_{24}^{1}$;\ \ 
$ 1$,
$ -c_{24}^{1}$,
$ -c_{24}^{1}$,
$ 1+\sqrt{3}$,
$ 1+\sqrt{3}$,
$ \frac{-3-3\sqrt{3}}{\sqrt{6}}$,
$ \frac{-3-3\sqrt{3}}{\sqrt{6}}$,
$ 2+\sqrt{3}$,
$ 2+\sqrt{3}$,
$ -2c_{24}^{1}$;\ \ 
$(\frac{-3-3\sqrt{3}}{\sqrt{6}})\mathrm{i}$,
$(\frac{3+3\sqrt{3}}{\sqrt{6}})\mathrm{i}$,
$ -2c_{24}^{1}$,
$ 2c_{24}^{1}$,
$(\frac{3+3\sqrt{3}}{\sqrt{6}})\mathrm{i}$,
$(\frac{-3-3\sqrt{3}}{\sqrt{6}})\mathrm{i}$,
$ -c_{24}^{1}$,
$ c_{24}^{1}$,
$0$;\ \ 
$(\frac{-3-3\sqrt{3}}{\sqrt{6}})\mathrm{i}$,
$ -2c_{24}^{1}$,
$ 2c_{24}^{1}$,
$(\frac{-3-3\sqrt{3}}{\sqrt{6}})\mathrm{i}$,
$(\frac{3+3\sqrt{3}}{\sqrt{6}})\mathrm{i}$,
$ -c_{24}^{1}$,
$ c_{24}^{1}$,
$0$;\ \ 
$ 1+\sqrt{3}$,
$ 1+\sqrt{3}$,
$0$,
$0$,
$ -1-\sqrt{3}$,
$ -1-\sqrt{3}$,
$ 2c_{24}^{1}$;\ \ 
$ 1+\sqrt{3}$,
$0$,
$0$,
$ -1-\sqrt{3}$,
$ -1-\sqrt{3}$,
$ -2c_{24}^{1}$;\ \ 
$(\frac{3+3\sqrt{3}}{\sqrt{6}})\mathrm{i}$,
$(\frac{-3-3\sqrt{3}}{\sqrt{6}})\mathrm{i}$,
$ \frac{3+3\sqrt{3}}{\sqrt{6}}$,
$ \frac{-3-3\sqrt{3}}{\sqrt{6}}$,
$0$;\ \ 
$(\frac{3+3\sqrt{3}}{\sqrt{6}})\mathrm{i}$,
$ \frac{3+3\sqrt{3}}{\sqrt{6}}$,
$ \frac{-3-3\sqrt{3}}{\sqrt{6}}$,
$0$;\ \ 
$ 1$,
$ 1$,
$ -2c_{24}^{1}$;\ \ 
$ 1$,
$ 2c_{24}^{1}$;\ \ 
$0$)

  \vskip 2ex

\noindent31. $11_{\frac{9}{2},89.56}^{48,797}$ \irep{2191}:\ \ 
$d_i$ = ($1.0$,
$1.0$,
$1.931$,
$1.931$,
$2.732$,
$2.732$,
$3.346$,
$3.346$,
$3.732$,
$3.732$,
$3.863$) 

\vskip 0.7ex
\hangindent=3em \hangafter=1
$D^2= 89.569 = 
48+24\sqrt{3}$

\vskip 0.7ex
\hangindent=3em \hangafter=1
$T = ( 0,
\frac{1}{2},
\frac{13}{16},
\frac{13}{16},
\frac{1}{3},
\frac{5}{6},
\frac{9}{16},
\frac{9}{16},
0,
\frac{1}{2},
\frac{7}{48} )
$,

\vskip 0.7ex
\hangindent=3em \hangafter=1
$S$ = ($ 1$,
$ 1$,
$ c_{24}^{1}$,
$ c_{24}^{1}$,
$ 1+\sqrt{3}$,
$ 1+\sqrt{3}$,
$ \frac{3+3\sqrt{3}}{\sqrt{6}}$,
$ \frac{3+3\sqrt{3}}{\sqrt{6}}$,
$ 2+\sqrt{3}$,
$ 2+\sqrt{3}$,
$ 2c_{24}^{1}$;\ \ 
$ 1$,
$ -c_{24}^{1}$,
$ -c_{24}^{1}$,
$ 1+\sqrt{3}$,
$ 1+\sqrt{3}$,
$ \frac{-3-3\sqrt{3}}{\sqrt{6}}$,
$ \frac{-3-3\sqrt{3}}{\sqrt{6}}$,
$ 2+\sqrt{3}$,
$ 2+\sqrt{3}$,
$ -2c_{24}^{1}$;\ \ 
$(\frac{3+3\sqrt{3}}{\sqrt{6}})\mathrm{i}$,
$(\frac{-3-3\sqrt{3}}{\sqrt{6}})\mathrm{i}$,
$ -2c_{24}^{1}$,
$ 2c_{24}^{1}$,
$(\frac{3+3\sqrt{3}}{\sqrt{6}})\mathrm{i}$,
$(\frac{-3-3\sqrt{3}}{\sqrt{6}})\mathrm{i}$,
$ -c_{24}^{1}$,
$ c_{24}^{1}$,
$0$;\ \ 
$(\frac{3+3\sqrt{3}}{\sqrt{6}})\mathrm{i}$,
$ -2c_{24}^{1}$,
$ 2c_{24}^{1}$,
$(\frac{-3-3\sqrt{3}}{\sqrt{6}})\mathrm{i}$,
$(\frac{3+3\sqrt{3}}{\sqrt{6}})\mathrm{i}$,
$ -c_{24}^{1}$,
$ c_{24}^{1}$,
$0$;\ \ 
$ 1+\sqrt{3}$,
$ 1+\sqrt{3}$,
$0$,
$0$,
$ -1-\sqrt{3}$,
$ -1-\sqrt{3}$,
$ 2c_{24}^{1}$;\ \ 
$ 1+\sqrt{3}$,
$0$,
$0$,
$ -1-\sqrt{3}$,
$ -1-\sqrt{3}$,
$ -2c_{24}^{1}$;\ \ 
$(\frac{-3-3\sqrt{3}}{\sqrt{6}})\mathrm{i}$,
$(\frac{3+3\sqrt{3}}{\sqrt{6}})\mathrm{i}$,
$ \frac{3+3\sqrt{3}}{\sqrt{6}}$,
$ \frac{-3-3\sqrt{3}}{\sqrt{6}}$,
$0$;\ \ 
$(\frac{-3-3\sqrt{3}}{\sqrt{6}})\mathrm{i}$,
$ \frac{3+3\sqrt{3}}{\sqrt{6}}$,
$ \frac{-3-3\sqrt{3}}{\sqrt{6}}$,
$0$;\ \ 
$ 1$,
$ 1$,
$ -2c_{24}^{1}$;\ \ 
$ 1$,
$ 2c_{24}^{1}$;\ \ 
$0$)

  \vskip 2ex

\noindent32. $11_{\frac{1}{2},89.56}^{48,139}$ \irep{2192}:\ \ 
$d_i$ = ($1.0$,
$1.0$,
$1.931$,
$1.931$,
$2.732$,
$2.732$,
$3.346$,
$3.346$,
$3.732$,
$3.732$,
$3.863$) 

\vskip 0.7ex
\hangindent=3em \hangafter=1
$D^2= 89.569 = 
48+24\sqrt{3}$

\vskip 0.7ex
\hangindent=3em \hangafter=1
$T = ( 0,
\frac{1}{2},
\frac{13}{16},
\frac{13}{16},
\frac{2}{3},
\frac{1}{6},
\frac{1}{16},
\frac{1}{16},
0,
\frac{1}{2},
\frac{23}{48} )
$,

\vskip 0.7ex
\hangindent=3em \hangafter=1
$S$ = ($ 1$,
$ 1$,
$ c_{24}^{1}$,
$ c_{24}^{1}$,
$ 1+\sqrt{3}$,
$ 1+\sqrt{3}$,
$ \frac{3+3\sqrt{3}}{\sqrt{6}}$,
$ \frac{3+3\sqrt{3}}{\sqrt{6}}$,
$ 2+\sqrt{3}$,
$ 2+\sqrt{3}$,
$ 2c_{24}^{1}$;\ \ 
$ 1$,
$ -c_{24}^{1}$,
$ -c_{24}^{1}$,
$ 1+\sqrt{3}$,
$ 1+\sqrt{3}$,
$ \frac{-3-3\sqrt{3}}{\sqrt{6}}$,
$ \frac{-3-3\sqrt{3}}{\sqrt{6}}$,
$ 2+\sqrt{3}$,
$ 2+\sqrt{3}$,
$ -2c_{24}^{1}$;\ \ 
$ \frac{-3-3\sqrt{3}}{\sqrt{6}}$,
$ \frac{3+3\sqrt{3}}{\sqrt{6}}$,
$ -2c_{24}^{1}$,
$ 2c_{24}^{1}$,
$ \frac{-3-3\sqrt{3}}{\sqrt{6}}$,
$ \frac{3+3\sqrt{3}}{\sqrt{6}}$,
$ -c_{24}^{1}$,
$ c_{24}^{1}$,
$0$;\ \ 
$ \frac{-3-3\sqrt{3}}{\sqrt{6}}$,
$ -2c_{24}^{1}$,
$ 2c_{24}^{1}$,
$ \frac{3+3\sqrt{3}}{\sqrt{6}}$,
$ \frac{-3-3\sqrt{3}}{\sqrt{6}}$,
$ -c_{24}^{1}$,
$ c_{24}^{1}$,
$0$;\ \ 
$ 1+\sqrt{3}$,
$ 1+\sqrt{3}$,
$0$,
$0$,
$ -1-\sqrt{3}$,
$ -1-\sqrt{3}$,
$ 2c_{24}^{1}$;\ \ 
$ 1+\sqrt{3}$,
$0$,
$0$,
$ -1-\sqrt{3}$,
$ -1-\sqrt{3}$,
$ -2c_{24}^{1}$;\ \ 
$ \frac{3+3\sqrt{3}}{\sqrt{6}}$,
$ \frac{-3-3\sqrt{3}}{\sqrt{6}}$,
$ \frac{3+3\sqrt{3}}{\sqrt{6}}$,
$ \frac{-3-3\sqrt{3}}{\sqrt{6}}$,
$0$;\ \ 
$ \frac{3+3\sqrt{3}}{\sqrt{6}}$,
$ \frac{3+3\sqrt{3}}{\sqrt{6}}$,
$ \frac{-3-3\sqrt{3}}{\sqrt{6}}$,
$0$;\ \ 
$ 1$,
$ 1$,
$ -2c_{24}^{1}$;\ \ 
$ 1$,
$ 2c_{24}^{1}$;\ \ 
$0$)

  \vskip 2ex

\noindent33. $11_{\frac{11}{2},89.56}^{48,288}$ \irep{2192}:\ \ 
$d_i$ = ($1.0$,
$1.0$,
$1.931$,
$1.931$,
$2.732$,
$2.732$,
$3.346$,
$3.346$,
$3.732$,
$3.732$,
$3.863$) 

\vskip 0.7ex
\hangindent=3em \hangafter=1
$D^2= 89.569 = 
48+24\sqrt{3}$

\vskip 0.7ex
\hangindent=3em \hangafter=1
$T = ( 0,
\frac{1}{2},
\frac{15}{16},
\frac{15}{16},
\frac{1}{3},
\frac{5}{6},
\frac{11}{16},
\frac{11}{16},
0,
\frac{1}{2},
\frac{13}{48} )
$,

\vskip 0.7ex
\hangindent=3em \hangafter=1
$S$ = ($ 1$,
$ 1$,
$ c_{24}^{1}$,
$ c_{24}^{1}$,
$ 1+\sqrt{3}$,
$ 1+\sqrt{3}$,
$ \frac{3+3\sqrt{3}}{\sqrt{6}}$,
$ \frac{3+3\sqrt{3}}{\sqrt{6}}$,
$ 2+\sqrt{3}$,
$ 2+\sqrt{3}$,
$ 2c_{24}^{1}$;\ \ 
$ 1$,
$ -c_{24}^{1}$,
$ -c_{24}^{1}$,
$ 1+\sqrt{3}$,
$ 1+\sqrt{3}$,
$ \frac{-3-3\sqrt{3}}{\sqrt{6}}$,
$ \frac{-3-3\sqrt{3}}{\sqrt{6}}$,
$ 2+\sqrt{3}$,
$ 2+\sqrt{3}$,
$ -2c_{24}^{1}$;\ \ 
$ \frac{3+3\sqrt{3}}{\sqrt{6}}$,
$ \frac{-3-3\sqrt{3}}{\sqrt{6}}$,
$ -2c_{24}^{1}$,
$ 2c_{24}^{1}$,
$ \frac{-3-3\sqrt{3}}{\sqrt{6}}$,
$ \frac{3+3\sqrt{3}}{\sqrt{6}}$,
$ -c_{24}^{1}$,
$ c_{24}^{1}$,
$0$;\ \ 
$ \frac{3+3\sqrt{3}}{\sqrt{6}}$,
$ -2c_{24}^{1}$,
$ 2c_{24}^{1}$,
$ \frac{3+3\sqrt{3}}{\sqrt{6}}$,
$ \frac{-3-3\sqrt{3}}{\sqrt{6}}$,
$ -c_{24}^{1}$,
$ c_{24}^{1}$,
$0$;\ \ 
$ 1+\sqrt{3}$,
$ 1+\sqrt{3}$,
$0$,
$0$,
$ -1-\sqrt{3}$,
$ -1-\sqrt{3}$,
$ 2c_{24}^{1}$;\ \ 
$ 1+\sqrt{3}$,
$0$,
$0$,
$ -1-\sqrt{3}$,
$ -1-\sqrt{3}$,
$ -2c_{24}^{1}$;\ \ 
$ \frac{-3-3\sqrt{3}}{\sqrt{6}}$,
$ \frac{3+3\sqrt{3}}{\sqrt{6}}$,
$ \frac{3+3\sqrt{3}}{\sqrt{6}}$,
$ \frac{-3-3\sqrt{3}}{\sqrt{6}}$,
$0$;\ \ 
$ \frac{-3-3\sqrt{3}}{\sqrt{6}}$,
$ \frac{3+3\sqrt{3}}{\sqrt{6}}$,
$ \frac{-3-3\sqrt{3}}{\sqrt{6}}$,
$0$;\ \ 
$ 1$,
$ 1$,
$ -2c_{24}^{1}$;\ \ 
$ 1$,
$ 2c_{24}^{1}$;\ \ 
$0$)

  \vskip 2ex

\noindent34. $11_{\frac{3}{2},89.56}^{48,177}$ \irep{2191}:\ \ 
$d_i$ = ($1.0$,
$1.0$,
$1.931$,
$1.931$,
$2.732$,
$2.732$,
$3.346$,
$3.346$,
$3.732$,
$3.732$,
$3.863$) 

\vskip 0.7ex
\hangindent=3em \hangafter=1
$D^2= 89.569 = 
48+24\sqrt{3}$

\vskip 0.7ex
\hangindent=3em \hangafter=1
$T = ( 0,
\frac{1}{2},
\frac{15}{16},
\frac{15}{16},
\frac{2}{3},
\frac{1}{6},
\frac{3}{16},
\frac{3}{16},
0,
\frac{1}{2},
\frac{29}{48} )
$,

\vskip 0.7ex
\hangindent=3em \hangafter=1
$S$ = ($ 1$,
$ 1$,
$ c_{24}^{1}$,
$ c_{24}^{1}$,
$ 1+\sqrt{3}$,
$ 1+\sqrt{3}$,
$ \frac{3+3\sqrt{3}}{\sqrt{6}}$,
$ \frac{3+3\sqrt{3}}{\sqrt{6}}$,
$ 2+\sqrt{3}$,
$ 2+\sqrt{3}$,
$ 2c_{24}^{1}$;\ \ 
$ 1$,
$ -c_{24}^{1}$,
$ -c_{24}^{1}$,
$ 1+\sqrt{3}$,
$ 1+\sqrt{3}$,
$ \frac{-3-3\sqrt{3}}{\sqrt{6}}$,
$ \frac{-3-3\sqrt{3}}{\sqrt{6}}$,
$ 2+\sqrt{3}$,
$ 2+\sqrt{3}$,
$ -2c_{24}^{1}$;\ \ 
$(\frac{3+3\sqrt{3}}{\sqrt{6}})\mathrm{i}$,
$(\frac{-3-3\sqrt{3}}{\sqrt{6}})\mathrm{i}$,
$ -2c_{24}^{1}$,
$ 2c_{24}^{1}$,
$(\frac{3+3\sqrt{3}}{\sqrt{6}})\mathrm{i}$,
$(\frac{-3-3\sqrt{3}}{\sqrt{6}})\mathrm{i}$,
$ -c_{24}^{1}$,
$ c_{24}^{1}$,
$0$;\ \ 
$(\frac{3+3\sqrt{3}}{\sqrt{6}})\mathrm{i}$,
$ -2c_{24}^{1}$,
$ 2c_{24}^{1}$,
$(\frac{-3-3\sqrt{3}}{\sqrt{6}})\mathrm{i}$,
$(\frac{3+3\sqrt{3}}{\sqrt{6}})\mathrm{i}$,
$ -c_{24}^{1}$,
$ c_{24}^{1}$,
$0$;\ \ 
$ 1+\sqrt{3}$,
$ 1+\sqrt{3}$,
$0$,
$0$,
$ -1-\sqrt{3}$,
$ -1-\sqrt{3}$,
$ 2c_{24}^{1}$;\ \ 
$ 1+\sqrt{3}$,
$0$,
$0$,
$ -1-\sqrt{3}$,
$ -1-\sqrt{3}$,
$ -2c_{24}^{1}$;\ \ 
$(\frac{-3-3\sqrt{3}}{\sqrt{6}})\mathrm{i}$,
$(\frac{3+3\sqrt{3}}{\sqrt{6}})\mathrm{i}$,
$ \frac{3+3\sqrt{3}}{\sqrt{6}}$,
$ \frac{-3-3\sqrt{3}}{\sqrt{6}}$,
$0$;\ \ 
$(\frac{-3-3\sqrt{3}}{\sqrt{6}})\mathrm{i}$,
$ \frac{3+3\sqrt{3}}{\sqrt{6}}$,
$ \frac{-3-3\sqrt{3}}{\sqrt{6}}$,
$0$;\ \ 
$ 1$,
$ 1$,
$ -2c_{24}^{1}$;\ \ 
$ 1$,
$ 2c_{24}^{1}$;\ \ 
$0$)

  \vskip 2ex

\noindent35. $11_{\frac{144}{23},310.1}^{23,306}$ \irep{1836}:\ \ 
$d_i$ = ($1.0$,
$1.981$,
$2.925$,
$3.815$,
$4.634$,
$5.367$,
$5.999$,
$6.520$,
$6.919$,
$7.190$,
$7.326$) 

\vskip 0.7ex
\hangindent=3em \hangafter=1
$D^2= 310.117 = 
66+55c^{1}_{23}
+45c^{2}_{23}
+36c^{3}_{23}
+28c^{4}_{23}
+21c^{5}_{23}
+15c^{6}_{23}
+10c^{7}_{23}
+6c^{8}_{23}
+3c^{9}_{23}
+c^{10}_{23}
$

\vskip 0.7ex
\hangindent=3em \hangafter=1
$T = ( 0,
\frac{5}{23},
\frac{21}{23},
\frac{2}{23},
\frac{17}{23},
\frac{20}{23},
\frac{11}{23},
\frac{13}{23},
\frac{3}{23},
\frac{4}{23},
\frac{16}{23} )
$,

\vskip 0.7ex
\hangindent=3em \hangafter=1
$S$ = ($ 1$,
$ -c_{23}^{11}$,
$ \xi_{23}^{3}$,
$ \xi_{23}^{19}$,
$ \xi_{23}^{5}$,
$ \xi_{23}^{17}$,
$ \xi_{23}^{7}$,
$ \xi_{23}^{15}$,
$ \xi_{23}^{9}$,
$ \xi_{23}^{13}$,
$ \xi_{23}^{11}$;\ \ 
$ -\xi_{23}^{19}$,
$ \xi_{23}^{17}$,
$ -\xi_{23}^{15}$,
$ \xi_{23}^{13}$,
$ -\xi_{23}^{11}$,
$ \xi_{23}^{9}$,
$ -\xi_{23}^{7}$,
$ \xi_{23}^{5}$,
$ -\xi_{23}^{3}$,
$ 1$;\ \ 
$ \xi_{23}^{9}$,
$ \xi_{23}^{11}$,
$ \xi_{23}^{15}$,
$ \xi_{23}^{5}$,
$ -c_{23}^{11}$,
$ -1$,
$ -\xi_{23}^{19}$,
$ -\xi_{23}^{7}$,
$ -\xi_{23}^{13}$;\ \ 
$ -\xi_{23}^{7}$,
$ \xi_{23}^{3}$,
$ 1$,
$ -\xi_{23}^{5}$,
$ \xi_{23}^{9}$,
$ -\xi_{23}^{13}$,
$ \xi_{23}^{17}$,
$ c_{23}^{11}$;\ \ 
$ c_{23}^{11}$,
$ -\xi_{23}^{7}$,
$ -\xi_{23}^{11}$,
$ -\xi_{23}^{17}$,
$ -1$,
$ \xi_{23}^{19}$,
$ \xi_{23}^{9}$;\ \ 
$ \xi_{23}^{13}$,
$ -\xi_{23}^{19}$,
$ c_{23}^{11}$,
$ \xi_{23}^{15}$,
$ -\xi_{23}^{9}$,
$ \xi_{23}^{3}$;\ \ 
$ \xi_{23}^{3}$,
$ \xi_{23}^{13}$,
$ \xi_{23}^{17}$,
$ -1$,
$ -\xi_{23}^{15}$;\ \ 
$ -\xi_{23}^{5}$,
$ -\xi_{23}^{3}$,
$ \xi_{23}^{11}$,
$ -\xi_{23}^{19}$;\ \ 
$ -\xi_{23}^{11}$,
$ c_{23}^{11}$,
$ \xi_{23}^{7}$;\ \ 
$ -\xi_{23}^{15}$,
$ \xi_{23}^{5}$;\ \ 
$ -\xi_{23}^{17}$)

  \vskip 2ex

\noindent36. $11_{\frac{40}{23},310.1}^{23,508}$ \irep{1836}:\ \ 
$d_i$ = ($1.0$,
$1.981$,
$2.925$,
$3.815$,
$4.634$,
$5.367$,
$5.999$,
$6.520$,
$6.919$,
$7.190$,
$7.326$) 

\vskip 0.7ex
\hangindent=3em \hangafter=1
$D^2= 310.117 = 
66+55c^{1}_{23}
+45c^{2}_{23}
+36c^{3}_{23}
+28c^{4}_{23}
+21c^{5}_{23}
+15c^{6}_{23}
+10c^{7}_{23}
+6c^{8}_{23}
+3c^{9}_{23}
+c^{10}_{23}
$

\vskip 0.7ex
\hangindent=3em \hangafter=1
$T = ( 0,
\frac{18}{23},
\frac{2}{23},
\frac{21}{23},
\frac{6}{23},
\frac{3}{23},
\frac{12}{23},
\frac{10}{23},
\frac{20}{23},
\frac{19}{23},
\frac{7}{23} )
$,

\vskip 0.7ex
\hangindent=3em \hangafter=1
$S$ = ($ 1$,
$ -c_{23}^{11}$,
$ \xi_{23}^{3}$,
$ \xi_{23}^{19}$,
$ \xi_{23}^{5}$,
$ \xi_{23}^{17}$,
$ \xi_{23}^{7}$,
$ \xi_{23}^{15}$,
$ \xi_{23}^{9}$,
$ \xi_{23}^{13}$,
$ \xi_{23}^{11}$;\ \ 
$ -\xi_{23}^{19}$,
$ \xi_{23}^{17}$,
$ -\xi_{23}^{15}$,
$ \xi_{23}^{13}$,
$ -\xi_{23}^{11}$,
$ \xi_{23}^{9}$,
$ -\xi_{23}^{7}$,
$ \xi_{23}^{5}$,
$ -\xi_{23}^{3}$,
$ 1$;\ \ 
$ \xi_{23}^{9}$,
$ \xi_{23}^{11}$,
$ \xi_{23}^{15}$,
$ \xi_{23}^{5}$,
$ -c_{23}^{11}$,
$ -1$,
$ -\xi_{23}^{19}$,
$ -\xi_{23}^{7}$,
$ -\xi_{23}^{13}$;\ \ 
$ -\xi_{23}^{7}$,
$ \xi_{23}^{3}$,
$ 1$,
$ -\xi_{23}^{5}$,
$ \xi_{23}^{9}$,
$ -\xi_{23}^{13}$,
$ \xi_{23}^{17}$,
$ c_{23}^{11}$;\ \ 
$ c_{23}^{11}$,
$ -\xi_{23}^{7}$,
$ -\xi_{23}^{11}$,
$ -\xi_{23}^{17}$,
$ -1$,
$ \xi_{23}^{19}$,
$ \xi_{23}^{9}$;\ \ 
$ \xi_{23}^{13}$,
$ -\xi_{23}^{19}$,
$ c_{23}^{11}$,
$ \xi_{23}^{15}$,
$ -\xi_{23}^{9}$,
$ \xi_{23}^{3}$;\ \ 
$ \xi_{23}^{3}$,
$ \xi_{23}^{13}$,
$ \xi_{23}^{17}$,
$ -1$,
$ -\xi_{23}^{15}$;\ \ 
$ -\xi_{23}^{5}$,
$ -\xi_{23}^{3}$,
$ \xi_{23}^{11}$,
$ -\xi_{23}^{19}$;\ \ 
$ -\xi_{23}^{11}$,
$ c_{23}^{11}$,
$ \xi_{23}^{7}$;\ \ 
$ -\xi_{23}^{15}$,
$ \xi_{23}^{5}$;\ \ 
$ -\xi_{23}^{17}$)

  \vskip 2ex

\noindent37. $11_{\frac{54}{19},696.5}^{19,306}$ \irep{1746}:\ \ 
$d_i$ = ($1.0$,
$2.972$,
$4.864$,
$6.54$,
$6.54$,
$6.623$,
$8.201$,
$9.556$,
$10.650$,
$11.453$,
$11.944$) 

\vskip 0.7ex
\hangindent=3em \hangafter=1
$D^2= 696.547 = 
171+152c^{1}_{19}
+133c^{2}_{19}
+114c^{3}_{19}
+95c^{4}_{19}
+76c^{5}_{19}
+57c^{6}_{19}
+38c^{7}_{19}
+19c^{8}_{19}
$

\vskip 0.7ex
\hangindent=3em \hangafter=1
$T = ( 0,
\frac{1}{19},
\frac{3}{19},
\frac{7}{19},
\frac{7}{19},
\frac{6}{19},
\frac{10}{19},
\frac{15}{19},
\frac{2}{19},
\frac{9}{19},
\frac{17}{19} )
$,

\vskip 0.7ex
\hangindent=3em \hangafter=1
$S$ = ($ 1$,
$ 2+c^{1}_{19}
+c^{2}_{19}
+c^{3}_{19}
+c^{4}_{19}
+c^{5}_{19}
+c^{6}_{19}
+c^{7}_{19}
+c^{8}_{19}
$,
$ 2+2c^{1}_{19}
+c^{2}_{19}
+c^{3}_{19}
+c^{4}_{19}
+c^{5}_{19}
+c^{6}_{19}
+c^{7}_{19}
+c^{8}_{19}
$,
$ \xi_{19}^{9}$,
$ \xi_{19}^{9}$,
$ 2+2c^{1}_{19}
+c^{2}_{19}
+c^{3}_{19}
+c^{4}_{19}
+c^{5}_{19}
+c^{6}_{19}
+c^{7}_{19}
$,
$ 2+2c^{1}_{19}
+2c^{2}_{19}
+c^{3}_{19}
+c^{4}_{19}
+c^{5}_{19}
+c^{6}_{19}
+c^{7}_{19}
$,
$ 2+2c^{1}_{19}
+2c^{2}_{19}
+c^{3}_{19}
+c^{4}_{19}
+c^{5}_{19}
+c^{6}_{19}
$,
$ 2+2c^{1}_{19}
+2c^{2}_{19}
+2c^{3}_{19}
+c^{4}_{19}
+c^{5}_{19}
+c^{6}_{19}
$,
$ 2+2c^{1}_{19}
+2c^{2}_{19}
+2c^{3}_{19}
+c^{4}_{19}
+c^{5}_{19}
$,
$ 2+2c^{1}_{19}
+2c^{2}_{19}
+2c^{3}_{19}
+2c^{4}_{19}
+c^{5}_{19}
$;\ \ 
$ 2+2c^{1}_{19}
+2c^{2}_{19}
+c^{3}_{19}
+c^{4}_{19}
+c^{5}_{19}
+c^{6}_{19}
+c^{7}_{19}
$,
$ 2+2c^{1}_{19}
+2c^{2}_{19}
+2c^{3}_{19}
+c^{4}_{19}
+c^{5}_{19}
$,
$ -\xi_{19}^{9}$,
$ -\xi_{19}^{9}$,
$ 2+2c^{1}_{19}
+2c^{2}_{19}
+2c^{3}_{19}
+2c^{4}_{19}
+c^{5}_{19}
$,
$ 2+2c^{1}_{19}
+2c^{2}_{19}
+c^{3}_{19}
+c^{4}_{19}
+c^{5}_{19}
+c^{6}_{19}
$,
$ 2+2c^{1}_{19}
+c^{2}_{19}
+c^{3}_{19}
+c^{4}_{19}
+c^{5}_{19}
+c^{6}_{19}
+c^{7}_{19}
+c^{8}_{19}
$,
$ -1$,
$ -2-2  c^{1}_{19}
-c^{2}_{19}
-c^{3}_{19}
-c^{4}_{19}
-c^{5}_{19}
-c^{6}_{19}
-c^{7}_{19}
$,
$ -2-2  c^{1}_{19}
-2  c^{2}_{19}
-2  c^{3}_{19}
-c^{4}_{19}
-c^{5}_{19}
-c^{6}_{19}
$;\ \ 
$ 2+2c^{1}_{19}
+2c^{2}_{19}
+2c^{3}_{19}
+c^{4}_{19}
+c^{5}_{19}
+c^{6}_{19}
$,
$ \xi_{19}^{9}$,
$ \xi_{19}^{9}$,
$ 2+c^{1}_{19}
+c^{2}_{19}
+c^{3}_{19}
+c^{4}_{19}
+c^{5}_{19}
+c^{6}_{19}
+c^{7}_{19}
+c^{8}_{19}
$,
$ -2-2  c^{1}_{19}
-c^{2}_{19}
-c^{3}_{19}
-c^{4}_{19}
-c^{5}_{19}
-c^{6}_{19}
-c^{7}_{19}
$,
$ -2-2  c^{1}_{19}
-2  c^{2}_{19}
-2  c^{3}_{19}
-2  c^{4}_{19}
-c^{5}_{19}
$,
$ -2-2  c^{1}_{19}
-2  c^{2}_{19}
-c^{3}_{19}
-c^{4}_{19}
-c^{5}_{19}
-c^{6}_{19}
$,
$ -1$,
$ 2+2c^{1}_{19}
+2c^{2}_{19}
+c^{3}_{19}
+c^{4}_{19}
+c^{5}_{19}
+c^{6}_{19}
+c^{7}_{19}
$;\ \ 
$ s^{1}_{19}
+s^{2}_{19}
+s^{3}_{19}
+s^{4}_{19}
+2\zeta^{5}_{19}
-\zeta^{-5}_{19}
+\zeta^{6}_{19}
+2\zeta^{7}_{19}
-\zeta^{-7}_{19}
+2\zeta^{8}_{19}
-\zeta^{-8}_{19}
+\zeta^{9}_{19}
$,
$ -1-2  \zeta^{1}_{19}
-2  \zeta^{2}_{19}
-2  \zeta^{3}_{19}
-2  \zeta^{4}_{19}
-2  \zeta^{5}_{19}
+\zeta^{-5}_{19}
-\zeta^{6}_{19}
-2  \zeta^{7}_{19}
+\zeta^{-7}_{19}
-2  \zeta^{8}_{19}
+\zeta^{-8}_{19}
-\zeta^{9}_{19}
$,
$ -\xi_{19}^{9}$,
$ \xi_{19}^{9}$,
$ -\xi_{19}^{9}$,
$ \xi_{19}^{9}$,
$ -\xi_{19}^{9}$,
$ \xi_{19}^{9}$;\ \ 
$ s^{1}_{19}
+s^{2}_{19}
+s^{3}_{19}
+s^{4}_{19}
+2\zeta^{5}_{19}
-\zeta^{-5}_{19}
+\zeta^{6}_{19}
+2\zeta^{7}_{19}
-\zeta^{-7}_{19}
+2\zeta^{8}_{19}
-\zeta^{-8}_{19}
+\zeta^{9}_{19}
$,
$ -\xi_{19}^{9}$,
$ \xi_{19}^{9}$,
$ -\xi_{19}^{9}$,
$ \xi_{19}^{9}$,
$ -\xi_{19}^{9}$,
$ \xi_{19}^{9}$;\ \ 
$ -2-2  c^{1}_{19}
-2  c^{2}_{19}
-c^{3}_{19}
-c^{4}_{19}
-c^{5}_{19}
-c^{6}_{19}
$,
$ -2-2  c^{1}_{19}
-2  c^{2}_{19}
-2  c^{3}_{19}
-c^{4}_{19}
-c^{5}_{19}
-c^{6}_{19}
$,
$ 1$,
$ 2+2c^{1}_{19}
+2c^{2}_{19}
+2c^{3}_{19}
+c^{4}_{19}
+c^{5}_{19}
$,
$ 2+2c^{1}_{19}
+2c^{2}_{19}
+c^{3}_{19}
+c^{4}_{19}
+c^{5}_{19}
+c^{6}_{19}
+c^{7}_{19}
$,
$ -2-2  c^{1}_{19}
-c^{2}_{19}
-c^{3}_{19}
-c^{4}_{19}
-c^{5}_{19}
-c^{6}_{19}
-c^{7}_{19}
-c^{8}_{19}
$;\ \ 
$ 2+2c^{1}_{19}
+c^{2}_{19}
+c^{3}_{19}
+c^{4}_{19}
+c^{5}_{19}
+c^{6}_{19}
+c^{7}_{19}
+c^{8}_{19}
$,
$ 2+2c^{1}_{19}
+2c^{2}_{19}
+2c^{3}_{19}
+c^{4}_{19}
+c^{5}_{19}
$,
$ -2-c^{1}_{19}
-c^{2}_{19}
-c^{3}_{19}
-c^{4}_{19}
-c^{5}_{19}
-c^{6}_{19}
-c^{7}_{19}
-c^{8}_{19}
$,
$ -2-2  c^{1}_{19}
-2  c^{2}_{19}
-2  c^{3}_{19}
-2  c^{4}_{19}
-c^{5}_{19}
$,
$ 1$;\ \ 
$ -2-2  c^{1}_{19}
-c^{2}_{19}
-c^{3}_{19}
-c^{4}_{19}
-c^{5}_{19}
-c^{6}_{19}
-c^{7}_{19}
$,
$ -2-2  c^{1}_{19}
-2  c^{2}_{19}
-c^{3}_{19}
-c^{4}_{19}
-c^{5}_{19}
-c^{6}_{19}
-c^{7}_{19}
$,
$ 2+2c^{1}_{19}
+2c^{2}_{19}
+2c^{3}_{19}
+c^{4}_{19}
+c^{5}_{19}
+c^{6}_{19}
$,
$ 2+c^{1}_{19}
+c^{2}_{19}
+c^{3}_{19}
+c^{4}_{19}
+c^{5}_{19}
+c^{6}_{19}
+c^{7}_{19}
+c^{8}_{19}
$;\ \ 
$ 2+2c^{1}_{19}
+2c^{2}_{19}
+2c^{3}_{19}
+2c^{4}_{19}
+c^{5}_{19}
$,
$ -2-2  c^{1}_{19}
-c^{2}_{19}
-c^{3}_{19}
-c^{4}_{19}
-c^{5}_{19}
-c^{6}_{19}
-c^{7}_{19}
-c^{8}_{19}
$,
$ -2-2  c^{1}_{19}
-c^{2}_{19}
-c^{3}_{19}
-c^{4}_{19}
-c^{5}_{19}
-c^{6}_{19}
-c^{7}_{19}
$;\ \ 
$ -2-c^{1}_{19}
-c^{2}_{19}
-c^{3}_{19}
-c^{4}_{19}
-c^{5}_{19}
-c^{6}_{19}
-c^{7}_{19}
-c^{8}_{19}
$,
$ 2+2c^{1}_{19}
+2c^{2}_{19}
+c^{3}_{19}
+c^{4}_{19}
+c^{5}_{19}
+c^{6}_{19}
$;\ \ 
$ -2-2  c^{1}_{19}
-2  c^{2}_{19}
-2  c^{3}_{19}
-c^{4}_{19}
-c^{5}_{19}
$)

  \vskip 2ex

\noindent38. $11_{\frac{98}{19},696.5}^{19,829}$ \irep{1746}:\ \ 
$d_i$ = ($1.0$,
$2.972$,
$4.864$,
$6.54$,
$6.54$,
$6.623$,
$8.201$,
$9.556$,
$10.650$,
$11.453$,
$11.944$) 

\vskip 0.7ex
\hangindent=3em \hangafter=1
$D^2= 696.547 = 
171+152c^{1}_{19}
+133c^{2}_{19}
+114c^{3}_{19}
+95c^{4}_{19}
+76c^{5}_{19}
+57c^{6}_{19}
+38c^{7}_{19}
+19c^{8}_{19}
$

\vskip 0.7ex
\hangindent=3em \hangafter=1
$T = ( 0,
\frac{18}{19},
\frac{16}{19},
\frac{12}{19},
\frac{12}{19},
\frac{13}{19},
\frac{9}{19},
\frac{4}{19},
\frac{17}{19},
\frac{10}{19},
\frac{2}{19} )
$,

\vskip 0.7ex
\hangindent=3em \hangafter=1
$S$ = ($ 1$,
$ 2+c^{1}_{19}
+c^{2}_{19}
+c^{3}_{19}
+c^{4}_{19}
+c^{5}_{19}
+c^{6}_{19}
+c^{7}_{19}
+c^{8}_{19}
$,
$ 2+2c^{1}_{19}
+c^{2}_{19}
+c^{3}_{19}
+c^{4}_{19}
+c^{5}_{19}
+c^{6}_{19}
+c^{7}_{19}
+c^{8}_{19}
$,
$ \xi_{19}^{9}$,
$ \xi_{19}^{9}$,
$ 2+2c^{1}_{19}
+c^{2}_{19}
+c^{3}_{19}
+c^{4}_{19}
+c^{5}_{19}
+c^{6}_{19}
+c^{7}_{19}
$,
$ 2+2c^{1}_{19}
+2c^{2}_{19}
+c^{3}_{19}
+c^{4}_{19}
+c^{5}_{19}
+c^{6}_{19}
+c^{7}_{19}
$,
$ 2+2c^{1}_{19}
+2c^{2}_{19}
+c^{3}_{19}
+c^{4}_{19}
+c^{5}_{19}
+c^{6}_{19}
$,
$ 2+2c^{1}_{19}
+2c^{2}_{19}
+2c^{3}_{19}
+c^{4}_{19}
+c^{5}_{19}
+c^{6}_{19}
$,
$ 2+2c^{1}_{19}
+2c^{2}_{19}
+2c^{3}_{19}
+c^{4}_{19}
+c^{5}_{19}
$,
$ 2+2c^{1}_{19}
+2c^{2}_{19}
+2c^{3}_{19}
+2c^{4}_{19}
+c^{5}_{19}
$;\ \ 
$ 2+2c^{1}_{19}
+2c^{2}_{19}
+c^{3}_{19}
+c^{4}_{19}
+c^{5}_{19}
+c^{6}_{19}
+c^{7}_{19}
$,
$ 2+2c^{1}_{19}
+2c^{2}_{19}
+2c^{3}_{19}
+c^{4}_{19}
+c^{5}_{19}
$,
$ -\xi_{19}^{9}$,
$ -\xi_{19}^{9}$,
$ 2+2c^{1}_{19}
+2c^{2}_{19}
+2c^{3}_{19}
+2c^{4}_{19}
+c^{5}_{19}
$,
$ 2+2c^{1}_{19}
+2c^{2}_{19}
+c^{3}_{19}
+c^{4}_{19}
+c^{5}_{19}
+c^{6}_{19}
$,
$ 2+2c^{1}_{19}
+c^{2}_{19}
+c^{3}_{19}
+c^{4}_{19}
+c^{5}_{19}
+c^{6}_{19}
+c^{7}_{19}
+c^{8}_{19}
$,
$ -1$,
$ -2-2  c^{1}_{19}
-c^{2}_{19}
-c^{3}_{19}
-c^{4}_{19}
-c^{5}_{19}
-c^{6}_{19}
-c^{7}_{19}
$,
$ -2-2  c^{1}_{19}
-2  c^{2}_{19}
-2  c^{3}_{19}
-c^{4}_{19}
-c^{5}_{19}
-c^{6}_{19}
$;\ \ 
$ 2+2c^{1}_{19}
+2c^{2}_{19}
+2c^{3}_{19}
+c^{4}_{19}
+c^{5}_{19}
+c^{6}_{19}
$,
$ \xi_{19}^{9}$,
$ \xi_{19}^{9}$,
$ 2+c^{1}_{19}
+c^{2}_{19}
+c^{3}_{19}
+c^{4}_{19}
+c^{5}_{19}
+c^{6}_{19}
+c^{7}_{19}
+c^{8}_{19}
$,
$ -2-2  c^{1}_{19}
-c^{2}_{19}
-c^{3}_{19}
-c^{4}_{19}
-c^{5}_{19}
-c^{6}_{19}
-c^{7}_{19}
$,
$ -2-2  c^{1}_{19}
-2  c^{2}_{19}
-2  c^{3}_{19}
-2  c^{4}_{19}
-c^{5}_{19}
$,
$ -2-2  c^{1}_{19}
-2  c^{2}_{19}
-c^{3}_{19}
-c^{4}_{19}
-c^{5}_{19}
-c^{6}_{19}
$,
$ -1$,
$ 2+2c^{1}_{19}
+2c^{2}_{19}
+c^{3}_{19}
+c^{4}_{19}
+c^{5}_{19}
+c^{6}_{19}
+c^{7}_{19}
$;\ \ 
$ -1-2  \zeta^{1}_{19}
-2  \zeta^{2}_{19}
-2  \zeta^{3}_{19}
-2  \zeta^{4}_{19}
-2  \zeta^{5}_{19}
+\zeta^{-5}_{19}
-\zeta^{6}_{19}
-2  \zeta^{7}_{19}
+\zeta^{-7}_{19}
-2  \zeta^{8}_{19}
+\zeta^{-8}_{19}
-\zeta^{9}_{19}
$,
$ s^{1}_{19}
+s^{2}_{19}
+s^{3}_{19}
+s^{4}_{19}
+2\zeta^{5}_{19}
-\zeta^{-5}_{19}
+\zeta^{6}_{19}
+2\zeta^{7}_{19}
-\zeta^{-7}_{19}
+2\zeta^{8}_{19}
-\zeta^{-8}_{19}
+\zeta^{9}_{19}
$,
$ -\xi_{19}^{9}$,
$ \xi_{19}^{9}$,
$ -\xi_{19}^{9}$,
$ \xi_{19}^{9}$,
$ -\xi_{19}^{9}$,
$ \xi_{19}^{9}$;\ \ 
$ -1-2  \zeta^{1}_{19}
-2  \zeta^{2}_{19}
-2  \zeta^{3}_{19}
-2  \zeta^{4}_{19}
-2  \zeta^{5}_{19}
+\zeta^{-5}_{19}
-\zeta^{6}_{19}
-2  \zeta^{7}_{19}
+\zeta^{-7}_{19}
-2  \zeta^{8}_{19}
+\zeta^{-8}_{19}
-\zeta^{9}_{19}
$,
$ -\xi_{19}^{9}$,
$ \xi_{19}^{9}$,
$ -\xi_{19}^{9}$,
$ \xi_{19}^{9}$,
$ -\xi_{19}^{9}$,
$ \xi_{19}^{9}$;\ \ 
$ -2-2  c^{1}_{19}
-2  c^{2}_{19}
-c^{3}_{19}
-c^{4}_{19}
-c^{5}_{19}
-c^{6}_{19}
$,
$ -2-2  c^{1}_{19}
-2  c^{2}_{19}
-2  c^{3}_{19}
-c^{4}_{19}
-c^{5}_{19}
-c^{6}_{19}
$,
$ 1$,
$ 2+2c^{1}_{19}
+2c^{2}_{19}
+2c^{3}_{19}
+c^{4}_{19}
+c^{5}_{19}
$,
$ 2+2c^{1}_{19}
+2c^{2}_{19}
+c^{3}_{19}
+c^{4}_{19}
+c^{5}_{19}
+c^{6}_{19}
+c^{7}_{19}
$,
$ -2-2  c^{1}_{19}
-c^{2}_{19}
-c^{3}_{19}
-c^{4}_{19}
-c^{5}_{19}
-c^{6}_{19}
-c^{7}_{19}
-c^{8}_{19}
$;\ \ 
$ 2+2c^{1}_{19}
+c^{2}_{19}
+c^{3}_{19}
+c^{4}_{19}
+c^{5}_{19}
+c^{6}_{19}
+c^{7}_{19}
+c^{8}_{19}
$,
$ 2+2c^{1}_{19}
+2c^{2}_{19}
+2c^{3}_{19}
+c^{4}_{19}
+c^{5}_{19}
$,
$ -2-c^{1}_{19}
-c^{2}_{19}
-c^{3}_{19}
-c^{4}_{19}
-c^{5}_{19}
-c^{6}_{19}
-c^{7}_{19}
-c^{8}_{19}
$,
$ -2-2  c^{1}_{19}
-2  c^{2}_{19}
-2  c^{3}_{19}
-2  c^{4}_{19}
-c^{5}_{19}
$,
$ 1$;\ \ 
$ -2-2  c^{1}_{19}
-c^{2}_{19}
-c^{3}_{19}
-c^{4}_{19}
-c^{5}_{19}
-c^{6}_{19}
-c^{7}_{19}
$,
$ -2-2  c^{1}_{19}
-2  c^{2}_{19}
-c^{3}_{19}
-c^{4}_{19}
-c^{5}_{19}
-c^{6}_{19}
-c^{7}_{19}
$,
$ 2+2c^{1}_{19}
+2c^{2}_{19}
+2c^{3}_{19}
+c^{4}_{19}
+c^{5}_{19}
+c^{6}_{19}
$,
$ 2+c^{1}_{19}
+c^{2}_{19}
+c^{3}_{19}
+c^{4}_{19}
+c^{5}_{19}
+c^{6}_{19}
+c^{7}_{19}
+c^{8}_{19}
$;\ \ 
$ 2+2c^{1}_{19}
+2c^{2}_{19}
+2c^{3}_{19}
+2c^{4}_{19}
+c^{5}_{19}
$,
$ -2-2  c^{1}_{19}
-c^{2}_{19}
-c^{3}_{19}
-c^{4}_{19}
-c^{5}_{19}
-c^{6}_{19}
-c^{7}_{19}
-c^{8}_{19}
$,
$ -2-2  c^{1}_{19}
-c^{2}_{19}
-c^{3}_{19}
-c^{4}_{19}
-c^{5}_{19}
-c^{6}_{19}
-c^{7}_{19}
$;\ \ 
$ -2-c^{1}_{19}
-c^{2}_{19}
-c^{3}_{19}
-c^{4}_{19}
-c^{5}_{19}
-c^{6}_{19}
-c^{7}_{19}
-c^{8}_{19}
$,
$ 2+2c^{1}_{19}
+2c^{2}_{19}
+c^{3}_{19}
+c^{4}_{19}
+c^{5}_{19}
+c^{6}_{19}
$;\ \ 
$ -2-2  c^{1}_{19}
-2  c^{2}_{19}
-2  c^{3}_{19}
-c^{4}_{19}
-c^{5}_{19}
$)

  \vskip 2ex

\noindent39. $11_{3,1337.}^{48,634}$ \irep{2187}:\ \ 
$d_i$ = ($1.0$,
$6.464$,
$6.464$,
$7.464$,
$7.464$,
$12.928$,
$12.928$,
$12.928$,
$13.928$,
$14.928$,
$14.928$) 

\vskip 0.7ex
\hangindent=3em \hangafter=1
$D^2= 1337.107 = 
672+384\sqrt{3}$

\vskip 0.7ex
\hangindent=3em \hangafter=1
$T = ( 0,
0,
0,
\frac{1}{4},
\frac{1}{4},
\frac{3}{4},
\frac{3}{16},
\frac{11}{16},
0,
\frac{1}{3},
\frac{7}{12} )
$,

\vskip 0.7ex
\hangindent=3em \hangafter=1
$S$ = ($ 1$,
$ 3+2\sqrt{3}$,
$ 3+2\sqrt{3}$,
$ 4+2\sqrt{3}$,
$ 4+2\sqrt{3}$,
$ 6+4\sqrt{3}$,
$ 6+4\sqrt{3}$,
$ 6+4\sqrt{3}$,
$ 7+4\sqrt{3}$,
$ 8+4\sqrt{3}$,
$ 8+4\sqrt{3}$;\ \ 
$ -7+4\zeta^{1}_{12}
-8  \zeta^{-1}_{12}
+8\zeta^{2}_{12}
$,
$ 1-8  \zeta^{1}_{12}
+4\zeta^{-1}_{12}
-8  \zeta^{2}_{12}
$,
$(-6-4\sqrt{3})\mathrm{i}$,
$(6+4\sqrt{3})\mathrm{i}$,
$ -6-4\sqrt{3}$,
$ 6+4\sqrt{3}$,
$ 6+4\sqrt{3}$,
$ -3-2\sqrt{3}$,
$0$,
$0$;\ \ 
$ -7+4\zeta^{1}_{12}
-8  \zeta^{-1}_{12}
+8\zeta^{2}_{12}
$,
$(6+4\sqrt{3})\mathrm{i}$,
$(-6-4\sqrt{3})\mathrm{i}$,
$ -6-4\sqrt{3}$,
$ 6+4\sqrt{3}$,
$ 6+4\sqrt{3}$,
$ -3-2\sqrt{3}$,
$0$,
$0$;\ \ 
$ (-8-4\sqrt{3})\zeta_{6}^{1}$,
$ (8+4\sqrt{3})\zeta_{3}^{1}$,
$0$,
$0$,
$0$,
$ 4+2\sqrt{3}$,
$ 8+4\sqrt{3}$,
$ -8-4\sqrt{3}$;\ \ 
$ (-8-4\sqrt{3})\zeta_{6}^{1}$,
$0$,
$0$,
$0$,
$ 4+2\sqrt{3}$,
$ 8+4\sqrt{3}$,
$ -8-4\sqrt{3}$;\ \ 
$ 12+8\sqrt{3}$,
$0$,
$0$,
$ -6-4\sqrt{3}$,
$0$,
$0$;\ \ 
$ \frac{24+12\sqrt{3}}{\sqrt{6}}$,
$ \frac{-24-12\sqrt{3}}{\sqrt{6}}$,
$ -6-4\sqrt{3}$,
$0$,
$0$;\ \ 
$ \frac{24+12\sqrt{3}}{\sqrt{6}}$,
$ -6-4\sqrt{3}$,
$0$,
$0$;\ \ 
$ 1$,
$ 8+4\sqrt{3}$,
$ 8+4\sqrt{3}$;\ \ 
$ -8-4\sqrt{3}$,
$ -8-4\sqrt{3}$;\ \ 
$ 8+4\sqrt{3}$)

  \vskip 2ex

\noindent40. $11_{3,1337.}^{48,924}$ \irep{2187}:\ \ 
$d_i$ = ($1.0$,
$6.464$,
$6.464$,
$7.464$,
$7.464$,
$12.928$,
$12.928$,
$12.928$,
$13.928$,
$14.928$,
$14.928$) 

\vskip 0.7ex
\hangindent=3em \hangafter=1
$D^2= 1337.107 = 
672+384\sqrt{3}$

\vskip 0.7ex
\hangindent=3em \hangafter=1
$T = ( 0,
0,
0,
\frac{1}{4},
\frac{1}{4},
\frac{3}{4},
\frac{7}{16},
\frac{15}{16},
0,
\frac{1}{3},
\frac{7}{12} )
$,

\vskip 0.7ex
\hangindent=3em \hangafter=1
$S$ = ($ 1$,
$ 3+2\sqrt{3}$,
$ 3+2\sqrt{3}$,
$ 4+2\sqrt{3}$,
$ 4+2\sqrt{3}$,
$ 6+4\sqrt{3}$,
$ 6+4\sqrt{3}$,
$ 6+4\sqrt{3}$,
$ 7+4\sqrt{3}$,
$ 8+4\sqrt{3}$,
$ 8+4\sqrt{3}$;\ \ 
$ -7+4\zeta^{1}_{12}
-8  \zeta^{-1}_{12}
+8\zeta^{2}_{12}
$,
$ 1-8  \zeta^{1}_{12}
+4\zeta^{-1}_{12}
-8  \zeta^{2}_{12}
$,
$(-6-4\sqrt{3})\mathrm{i}$,
$(6+4\sqrt{3})\mathrm{i}$,
$ -6-4\sqrt{3}$,
$ 6+4\sqrt{3}$,
$ 6+4\sqrt{3}$,
$ -3-2\sqrt{3}$,
$0$,
$0$;\ \ 
$ -7+4\zeta^{1}_{12}
-8  \zeta^{-1}_{12}
+8\zeta^{2}_{12}
$,
$(6+4\sqrt{3})\mathrm{i}$,
$(-6-4\sqrt{3})\mathrm{i}$,
$ -6-4\sqrt{3}$,
$ 6+4\sqrt{3}$,
$ 6+4\sqrt{3}$,
$ -3-2\sqrt{3}$,
$0$,
$0$;\ \ 
$ (-8-4\sqrt{3})\zeta_{6}^{1}$,
$ (8+4\sqrt{3})\zeta_{3}^{1}$,
$0$,
$0$,
$0$,
$ 4+2\sqrt{3}$,
$ 8+4\sqrt{3}$,
$ -8-4\sqrt{3}$;\ \ 
$ (-8-4\sqrt{3})\zeta_{6}^{1}$,
$0$,
$0$,
$0$,
$ 4+2\sqrt{3}$,
$ 8+4\sqrt{3}$,
$ -8-4\sqrt{3}$;\ \ 
$ 12+8\sqrt{3}$,
$0$,
$0$,
$ -6-4\sqrt{3}$,
$0$,
$0$;\ \ 
$ \frac{-24-12\sqrt{3}}{\sqrt{6}}$,
$ \frac{24+12\sqrt{3}}{\sqrt{6}}$,
$ -6-4\sqrt{3}$,
$0$,
$0$;\ \ 
$ \frac{-24-12\sqrt{3}}{\sqrt{6}}$,
$ -6-4\sqrt{3}$,
$0$,
$0$;\ \ 
$ 1$,
$ 8+4\sqrt{3}$,
$ 8+4\sqrt{3}$;\ \ 
$ -8-4\sqrt{3}$,
$ -8-4\sqrt{3}$;\ \ 
$ 8+4\sqrt{3}$)

  \vskip 2ex

\noindent41. $11_{5,1337.}^{48,528}$ \irep{2187}:\ \ 
$d_i$ = ($1.0$,
$6.464$,
$6.464$,
$7.464$,
$7.464$,
$12.928$,
$12.928$,
$12.928$,
$13.928$,
$14.928$,
$14.928$) 

\vskip 0.7ex
\hangindent=3em \hangafter=1
$D^2= 1337.107 = 
672+384\sqrt{3}$

\vskip 0.7ex
\hangindent=3em \hangafter=1
$T = ( 0,
0,
0,
\frac{3}{4},
\frac{3}{4},
\frac{1}{4},
\frac{1}{16},
\frac{9}{16},
0,
\frac{2}{3},
\frac{5}{12} )
$,

\vskip 0.7ex
\hangindent=3em \hangafter=1
$S$ = ($ 1$,
$ 3+2\sqrt{3}$,
$ 3+2\sqrt{3}$,
$ 4+2\sqrt{3}$,
$ 4+2\sqrt{3}$,
$ 6+4\sqrt{3}$,
$ 6+4\sqrt{3}$,
$ 6+4\sqrt{3}$,
$ 7+4\sqrt{3}$,
$ 8+4\sqrt{3}$,
$ 8+4\sqrt{3}$;\ \ 
$ 1-8  \zeta^{1}_{12}
+4\zeta^{-1}_{12}
-8  \zeta^{2}_{12}
$,
$ -7+4\zeta^{1}_{12}
-8  \zeta^{-1}_{12}
+8\zeta^{2}_{12}
$,
$(-6-4\sqrt{3})\mathrm{i}$,
$(6+4\sqrt{3})\mathrm{i}$,
$ -6-4\sqrt{3}$,
$ 6+4\sqrt{3}$,
$ 6+4\sqrt{3}$,
$ -3-2\sqrt{3}$,
$0$,
$0$;\ \ 
$ 1-8  \zeta^{1}_{12}
+4\zeta^{-1}_{12}
-8  \zeta^{2}_{12}
$,
$(6+4\sqrt{3})\mathrm{i}$,
$(-6-4\sqrt{3})\mathrm{i}$,
$ -6-4\sqrt{3}$,
$ 6+4\sqrt{3}$,
$ 6+4\sqrt{3}$,
$ -3-2\sqrt{3}$,
$0$,
$0$;\ \ 
$ (8+4\sqrt{3})\zeta_{3}^{1}$,
$ (-8-4\sqrt{3})\zeta_{6}^{1}$,
$0$,
$0$,
$0$,
$ 4+2\sqrt{3}$,
$ 8+4\sqrt{3}$,
$ -8-4\sqrt{3}$;\ \ 
$ (8+4\sqrt{3})\zeta_{3}^{1}$,
$0$,
$0$,
$0$,
$ 4+2\sqrt{3}$,
$ 8+4\sqrt{3}$,
$ -8-4\sqrt{3}$;\ \ 
$ 12+8\sqrt{3}$,
$0$,
$0$,
$ -6-4\sqrt{3}$,
$0$,
$0$;\ \ 
$ \frac{-24-12\sqrt{3}}{\sqrt{6}}$,
$ \frac{24+12\sqrt{3}}{\sqrt{6}}$,
$ -6-4\sqrt{3}$,
$0$,
$0$;\ \ 
$ \frac{-24-12\sqrt{3}}{\sqrt{6}}$,
$ -6-4\sqrt{3}$,
$0$,
$0$;\ \ 
$ 1$,
$ 8+4\sqrt{3}$,
$ 8+4\sqrt{3}$;\ \ 
$ -8-4\sqrt{3}$,
$ -8-4\sqrt{3}$;\ \ 
$ 8+4\sqrt{3}$)

  \vskip 2ex

\noindent42. $11_{5,1337.}^{48,372}$ \irep{2187}:\ \ 
$d_i$ = ($1.0$,
$6.464$,
$6.464$,
$7.464$,
$7.464$,
$12.928$,
$12.928$,
$12.928$,
$13.928$,
$14.928$,
$14.928$) 

\vskip 0.7ex
\hangindent=3em \hangafter=1
$D^2= 1337.107 = 
672+384\sqrt{3}$

\vskip 0.7ex
\hangindent=3em \hangafter=1
$T = ( 0,
0,
0,
\frac{3}{4},
\frac{3}{4},
\frac{1}{4},
\frac{5}{16},
\frac{13}{16},
0,
\frac{2}{3},
\frac{5}{12} )
$,

\vskip 0.7ex
\hangindent=3em \hangafter=1
$S$ = ($ 1$,
$ 3+2\sqrt{3}$,
$ 3+2\sqrt{3}$,
$ 4+2\sqrt{3}$,
$ 4+2\sqrt{3}$,
$ 6+4\sqrt{3}$,
$ 6+4\sqrt{3}$,
$ 6+4\sqrt{3}$,
$ 7+4\sqrt{3}$,
$ 8+4\sqrt{3}$,
$ 8+4\sqrt{3}$;\ \ 
$ 1-8  \zeta^{1}_{12}
+4\zeta^{-1}_{12}
-8  \zeta^{2}_{12}
$,
$ -7+4\zeta^{1}_{12}
-8  \zeta^{-1}_{12}
+8\zeta^{2}_{12}
$,
$(-6-4\sqrt{3})\mathrm{i}$,
$(6+4\sqrt{3})\mathrm{i}$,
$ -6-4\sqrt{3}$,
$ 6+4\sqrt{3}$,
$ 6+4\sqrt{3}$,
$ -3-2\sqrt{3}$,
$0$,
$0$;\ \ 
$ 1-8  \zeta^{1}_{12}
+4\zeta^{-1}_{12}
-8  \zeta^{2}_{12}
$,
$(6+4\sqrt{3})\mathrm{i}$,
$(-6-4\sqrt{3})\mathrm{i}$,
$ -6-4\sqrt{3}$,
$ 6+4\sqrt{3}$,
$ 6+4\sqrt{3}$,
$ -3-2\sqrt{3}$,
$0$,
$0$;\ \ 
$ (8+4\sqrt{3})\zeta_{3}^{1}$,
$ (-8-4\sqrt{3})\zeta_{6}^{1}$,
$0$,
$0$,
$0$,
$ 4+2\sqrt{3}$,
$ 8+4\sqrt{3}$,
$ -8-4\sqrt{3}$;\ \ 
$ (8+4\sqrt{3})\zeta_{3}^{1}$,
$0$,
$0$,
$0$,
$ 4+2\sqrt{3}$,
$ 8+4\sqrt{3}$,
$ -8-4\sqrt{3}$;\ \ 
$ 12+8\sqrt{3}$,
$0$,
$0$,
$ -6-4\sqrt{3}$,
$0$,
$0$;\ \ 
$ \frac{24+12\sqrt{3}}{\sqrt{6}}$,
$ \frac{-24-12\sqrt{3}}{\sqrt{6}}$,
$ -6-4\sqrt{3}$,
$0$,
$0$;\ \ 
$ \frac{24+12\sqrt{3}}{\sqrt{6}}$,
$ -6-4\sqrt{3}$,
$0$,
$0$;\ \ 
$ 1$,
$ 8+4\sqrt{3}$,
$ 8+4\sqrt{3}$;\ \ 
$ -8-4\sqrt{3}$,
$ -8-4\sqrt{3}$;\ \ 
$ 8+4\sqrt{3}$)

  \vskip 2ex

\noindent43. $11_{\frac{32}{5},1964.}^{35,581}$ \irep{2077}:\ \ 
$d_i$ = ($1.0$,
$8.807$,
$8.807$,
$8.807$,
$11.632$,
$13.250$,
$14.250$,
$14.250$,
$14.250$,
$19.822$,
$20.440$) 

\vskip 0.7ex
\hangindent=3em \hangafter=1
$D^2= 1964.590 = 
910-280  c^{1}_{35}
+280c^{2}_{35}
+280c^{3}_{35}
+175c^{4}_{35}
+280c^{5}_{35}
-105  c^{6}_{35}
+490c^{7}_{35}
-280  c^{8}_{35}
+175c^{9}_{35}
+280c^{10}_{35}
$

\vskip 0.7ex
\hangindent=3em \hangafter=1
$T = ( 0,
\frac{2}{35},
\frac{22}{35},
\frac{32}{35},
\frac{1}{5},
0,
\frac{3}{7},
\frac{5}{7},
\frac{6}{7},
\frac{3}{5},
\frac{1}{5} )
$,

\vskip 0.7ex
\hangindent=3em \hangafter=1
$S$ = ($ 1$,
$ 4-c^{1}_{35}
+c^{2}_{35}
+c^{3}_{35}
+c^{4}_{35}
+c^{5}_{35}
+2c^{7}_{35}
-c^{8}_{35}
+c^{9}_{35}
+c^{10}_{35}
$,
$ 4-c^{1}_{35}
+c^{2}_{35}
+c^{3}_{35}
+c^{4}_{35}
+c^{5}_{35}
+2c^{7}_{35}
-c^{8}_{35}
+c^{9}_{35}
+c^{10}_{35}
$,
$ 4-c^{1}_{35}
+c^{2}_{35}
+c^{3}_{35}
+c^{4}_{35}
+c^{5}_{35}
+2c^{7}_{35}
-c^{8}_{35}
+c^{9}_{35}
+c^{10}_{35}
$,
$ 5-2  c^{1}_{35}
+2c^{2}_{35}
+2c^{3}_{35}
+c^{4}_{35}
+2c^{5}_{35}
-c^{6}_{35}
+3c^{7}_{35}
-2  c^{8}_{35}
+c^{9}_{35}
+2c^{10}_{35}
$,
$ 6-2  c^{1}_{35}
+2c^{2}_{35}
+2c^{3}_{35}
+c^{4}_{35}
+2c^{5}_{35}
-c^{6}_{35}
+4c^{7}_{35}
-2  c^{8}_{35}
+c^{9}_{35}
+2c^{10}_{35}
$,
$ 7-2  c^{1}_{35}
+2c^{2}_{35}
+2c^{3}_{35}
+c^{4}_{35}
+2c^{5}_{35}
-c^{6}_{35}
+4c^{7}_{35}
-2  c^{8}_{35}
+c^{9}_{35}
+2c^{10}_{35}
$,
$ 7-2  c^{1}_{35}
+2c^{2}_{35}
+2c^{3}_{35}
+c^{4}_{35}
+2c^{5}_{35}
-c^{6}_{35}
+4c^{7}_{35}
-2  c^{8}_{35}
+c^{9}_{35}
+2c^{10}_{35}
$,
$ 7-2  c^{1}_{35}
+2c^{2}_{35}
+2c^{3}_{35}
+c^{4}_{35}
+2c^{5}_{35}
-c^{6}_{35}
+4c^{7}_{35}
-2  c^{8}_{35}
+c^{9}_{35}
+2c^{10}_{35}
$,
$ 9-3  c^{1}_{35}
+3c^{2}_{35}
+3c^{3}_{35}
+2c^{4}_{35}
+3c^{5}_{35}
-c^{6}_{35}
+4c^{7}_{35}
-3  c^{8}_{35}
+2c^{9}_{35}
+3c^{10}_{35}
$,
$ 9-3  c^{1}_{35}
+3c^{2}_{35}
+3c^{3}_{35}
+2c^{4}_{35}
+3c^{5}_{35}
-c^{6}_{35}
+5c^{7}_{35}
-3  c^{8}_{35}
+2c^{9}_{35}
+3c^{10}_{35}
$;\ \ 
$ 1-c^{1}_{35}
+4c^{2}_{35}
-2  c^{3}_{35}
+c^{4}_{35}
+2c^{5}_{35}
-c^{6}_{35}
+c^{7}_{35}
+2c^{9}_{35}
-3  c^{10}_{35}
+2c^{11}_{35}
$,
$ 1+5c^{1}_{35}
+4c^{3}_{35}
+3c^{4}_{35}
+c^{5}_{35}
+4c^{6}_{35}
+2c^{8}_{35}
+3c^{10}_{35}
+c^{11}_{35}
$,
$ 5-6  c^{1}_{35}
-2  c^{2}_{35}
-3  c^{4}_{35}
-c^{5}_{35}
-4  c^{6}_{35}
+3c^{7}_{35}
-4  c^{8}_{35}
-c^{9}_{35}
+2c^{10}_{35}
-3  c^{11}_{35}
$,
$ 7-2  c^{1}_{35}
+2c^{2}_{35}
+2c^{3}_{35}
+c^{4}_{35}
+2c^{5}_{35}
-c^{6}_{35}
+4c^{7}_{35}
-2  c^{8}_{35}
+c^{9}_{35}
+2c^{10}_{35}
$,
$ 4-c^{1}_{35}
+c^{2}_{35}
+c^{3}_{35}
+c^{4}_{35}
+c^{5}_{35}
+2c^{7}_{35}
-c^{8}_{35}
+c^{9}_{35}
+c^{10}_{35}
$,
$ -1-3  c^{1}_{35}
-2  c^{3}_{35}
-2  c^{4}_{35}
-c^{5}_{35}
-2  c^{6}_{35}
-c^{8}_{35}
-2  c^{10}_{35}
$,
$ -3+4c^{1}_{35}
+c^{2}_{35}
+2c^{4}_{35}
+c^{5}_{35}
+2c^{6}_{35}
-2  c^{7}_{35}
+2c^{8}_{35}
+c^{9}_{35}
-c^{10}_{35}
+2c^{11}_{35}
$,
$ -2  c^{2}_{35}
+c^{3}_{35}
-c^{4}_{35}
-c^{5}_{35}
-2  c^{9}_{35}
+2c^{10}_{35}
-2  c^{11}_{35}
$,
$0$,
$ -7+2c^{1}_{35}
-2  c^{2}_{35}
-2  c^{3}_{35}
-c^{4}_{35}
-2  c^{5}_{35}
+c^{6}_{35}
-4  c^{7}_{35}
+2c^{8}_{35}
-c^{9}_{35}
-2  c^{10}_{35}
$;\ \ 
$ 5-6  c^{1}_{35}
-2  c^{2}_{35}
-3  c^{4}_{35}
-c^{5}_{35}
-4  c^{6}_{35}
+3c^{7}_{35}
-4  c^{8}_{35}
-c^{9}_{35}
+2c^{10}_{35}
-3  c^{11}_{35}
$,
$ 1-c^{1}_{35}
+4c^{2}_{35}
-2  c^{3}_{35}
+c^{4}_{35}
+2c^{5}_{35}
-c^{6}_{35}
+c^{7}_{35}
+2c^{9}_{35}
-3  c^{10}_{35}
+2c^{11}_{35}
$,
$ 7-2  c^{1}_{35}
+2c^{2}_{35}
+2c^{3}_{35}
+c^{4}_{35}
+2c^{5}_{35}
-c^{6}_{35}
+4c^{7}_{35}
-2  c^{8}_{35}
+c^{9}_{35}
+2c^{10}_{35}
$,
$ 4-c^{1}_{35}
+c^{2}_{35}
+c^{3}_{35}
+c^{4}_{35}
+c^{5}_{35}
+2c^{7}_{35}
-c^{8}_{35}
+c^{9}_{35}
+c^{10}_{35}
$,
$ -3+4c^{1}_{35}
+c^{2}_{35}
+2c^{4}_{35}
+c^{5}_{35}
+2c^{6}_{35}
-2  c^{7}_{35}
+2c^{8}_{35}
+c^{9}_{35}
-c^{10}_{35}
+2c^{11}_{35}
$,
$ -2  c^{2}_{35}
+c^{3}_{35}
-c^{4}_{35}
-c^{5}_{35}
-2  c^{9}_{35}
+2c^{10}_{35}
-2  c^{11}_{35}
$,
$ -1-3  c^{1}_{35}
-2  c^{3}_{35}
-2  c^{4}_{35}
-c^{5}_{35}
-2  c^{6}_{35}
-c^{8}_{35}
-2  c^{10}_{35}
$,
$0$,
$ -7+2c^{1}_{35}
-2  c^{2}_{35}
-2  c^{3}_{35}
-c^{4}_{35}
-2  c^{5}_{35}
+c^{6}_{35}
-4  c^{7}_{35}
+2c^{8}_{35}
-c^{9}_{35}
-2  c^{10}_{35}
$;\ \ 
$ 1+5c^{1}_{35}
+4c^{3}_{35}
+3c^{4}_{35}
+c^{5}_{35}
+4c^{6}_{35}
+2c^{8}_{35}
+3c^{10}_{35}
+c^{11}_{35}
$,
$ 7-2  c^{1}_{35}
+2c^{2}_{35}
+2c^{3}_{35}
+c^{4}_{35}
+2c^{5}_{35}
-c^{6}_{35}
+4c^{7}_{35}
-2  c^{8}_{35}
+c^{9}_{35}
+2c^{10}_{35}
$,
$ 4-c^{1}_{35}
+c^{2}_{35}
+c^{3}_{35}
+c^{4}_{35}
+c^{5}_{35}
+2c^{7}_{35}
-c^{8}_{35}
+c^{9}_{35}
+c^{10}_{35}
$,
$ -2  c^{2}_{35}
+c^{3}_{35}
-c^{4}_{35}
-c^{5}_{35}
-2  c^{9}_{35}
+2c^{10}_{35}
-2  c^{11}_{35}
$,
$ -1-3  c^{1}_{35}
-2  c^{3}_{35}
-2  c^{4}_{35}
-c^{5}_{35}
-2  c^{6}_{35}
-c^{8}_{35}
-2  c^{10}_{35}
$,
$ -3+4c^{1}_{35}
+c^{2}_{35}
+2c^{4}_{35}
+c^{5}_{35}
+2c^{6}_{35}
-2  c^{7}_{35}
+2c^{8}_{35}
+c^{9}_{35}
-c^{10}_{35}
+2c^{11}_{35}
$,
$0$,
$ -7+2c^{1}_{35}
-2  c^{2}_{35}
-2  c^{3}_{35}
-c^{4}_{35}
-2  c^{5}_{35}
+c^{6}_{35}
-4  c^{7}_{35}
+2c^{8}_{35}
-c^{9}_{35}
-2  c^{10}_{35}
$;\ \ 
$ -6+2c^{1}_{35}
-2  c^{2}_{35}
-2  c^{3}_{35}
-c^{4}_{35}
-2  c^{5}_{35}
+c^{6}_{35}
-4  c^{7}_{35}
+2c^{8}_{35}
-c^{9}_{35}
-2  c^{10}_{35}
$,
$ -9+3c^{1}_{35}
-3  c^{2}_{35}
-3  c^{3}_{35}
-2  c^{4}_{35}
-3  c^{5}_{35}
+c^{6}_{35}
-5  c^{7}_{35}
+3c^{8}_{35}
-2  c^{9}_{35}
-3  c^{10}_{35}
$,
$ -4+c^{1}_{35}
-c^{2}_{35}
-c^{3}_{35}
-c^{4}_{35}
-c^{5}_{35}
-2  c^{7}_{35}
+c^{8}_{35}
-c^{9}_{35}
-c^{10}_{35}
$,
$ -4+c^{1}_{35}
-c^{2}_{35}
-c^{3}_{35}
-c^{4}_{35}
-c^{5}_{35}
-2  c^{7}_{35}
+c^{8}_{35}
-c^{9}_{35}
-c^{10}_{35}
$,
$ -4+c^{1}_{35}
-c^{2}_{35}
-c^{3}_{35}
-c^{4}_{35}
-c^{5}_{35}
-2  c^{7}_{35}
+c^{8}_{35}
-c^{9}_{35}
-c^{10}_{35}
$,
$ 9-3  c^{1}_{35}
+3c^{2}_{35}
+3c^{3}_{35}
+2c^{4}_{35}
+3c^{5}_{35}
-c^{6}_{35}
+4c^{7}_{35}
-3  c^{8}_{35}
+2c^{9}_{35}
+3c^{10}_{35}
$,
$ 1$;\ \ 
$ 1$,
$ 7-2  c^{1}_{35}
+2c^{2}_{35}
+2c^{3}_{35}
+c^{4}_{35}
+2c^{5}_{35}
-c^{6}_{35}
+4c^{7}_{35}
-2  c^{8}_{35}
+c^{9}_{35}
+2c^{10}_{35}
$,
$ 7-2  c^{1}_{35}
+2c^{2}_{35}
+2c^{3}_{35}
+c^{4}_{35}
+2c^{5}_{35}
-c^{6}_{35}
+4c^{7}_{35}
-2  c^{8}_{35}
+c^{9}_{35}
+2c^{10}_{35}
$,
$ 7-2  c^{1}_{35}
+2c^{2}_{35}
+2c^{3}_{35}
+c^{4}_{35}
+2c^{5}_{35}
-c^{6}_{35}
+4c^{7}_{35}
-2  c^{8}_{35}
+c^{9}_{35}
+2c^{10}_{35}
$,
$ -9+3c^{1}_{35}
-3  c^{2}_{35}
-3  c^{3}_{35}
-2  c^{4}_{35}
-3  c^{5}_{35}
+c^{6}_{35}
-4  c^{7}_{35}
+3c^{8}_{35}
-2  c^{9}_{35}
-3  c^{10}_{35}
$,
$ -5+2c^{1}_{35}
-2  c^{2}_{35}
-2  c^{3}_{35}
-c^{4}_{35}
-2  c^{5}_{35}
+c^{6}_{35}
-3  c^{7}_{35}
+2c^{8}_{35}
-c^{9}_{35}
-2  c^{10}_{35}
$;\ \ 
$ -5+6c^{1}_{35}
+2c^{2}_{35}
+3c^{4}_{35}
+c^{5}_{35}
+4c^{6}_{35}
-3  c^{7}_{35}
+4c^{8}_{35}
+c^{9}_{35}
-2  c^{10}_{35}
+3c^{11}_{35}
$,
$ -1+c^{1}_{35}
-4  c^{2}_{35}
+2c^{3}_{35}
-c^{4}_{35}
-2  c^{5}_{35}
+c^{6}_{35}
-c^{7}_{35}
-2  c^{9}_{35}
+3c^{10}_{35}
-2  c^{11}_{35}
$,
$ -1-5  c^{1}_{35}
-4  c^{3}_{35}
-3  c^{4}_{35}
-c^{5}_{35}
-4  c^{6}_{35}
-2  c^{8}_{35}
-3  c^{10}_{35}
-c^{11}_{35}
$,
$0$,
$ 4-c^{1}_{35}
+c^{2}_{35}
+c^{3}_{35}
+c^{4}_{35}
+c^{5}_{35}
+2c^{7}_{35}
-c^{8}_{35}
+c^{9}_{35}
+c^{10}_{35}
$;\ \ 
$ -1-5  c^{1}_{35}
-4  c^{3}_{35}
-3  c^{4}_{35}
-c^{5}_{35}
-4  c^{6}_{35}
-2  c^{8}_{35}
-3  c^{10}_{35}
-c^{11}_{35}
$,
$ -5+6c^{1}_{35}
+2c^{2}_{35}
+3c^{4}_{35}
+c^{5}_{35}
+4c^{6}_{35}
-3  c^{7}_{35}
+4c^{8}_{35}
+c^{9}_{35}
-2  c^{10}_{35}
+3c^{11}_{35}
$,
$0$,
$ 4-c^{1}_{35}
+c^{2}_{35}
+c^{3}_{35}
+c^{4}_{35}
+c^{5}_{35}
+2c^{7}_{35}
-c^{8}_{35}
+c^{9}_{35}
+c^{10}_{35}
$;\ \ 
$ -1+c^{1}_{35}
-4  c^{2}_{35}
+2c^{3}_{35}
-c^{4}_{35}
-2  c^{5}_{35}
+c^{6}_{35}
-c^{7}_{35}
-2  c^{9}_{35}
+3c^{10}_{35}
-2  c^{11}_{35}
$,
$0$,
$ 4-c^{1}_{35}
+c^{2}_{35}
+c^{3}_{35}
+c^{4}_{35}
+c^{5}_{35}
+2c^{7}_{35}
-c^{8}_{35}
+c^{9}_{35}
+c^{10}_{35}
$;\ \ 
$ -9+3c^{1}_{35}
-3  c^{2}_{35}
-3  c^{3}_{35}
-2  c^{4}_{35}
-3  c^{5}_{35}
+c^{6}_{35}
-4  c^{7}_{35}
+3c^{8}_{35}
-2  c^{9}_{35}
-3  c^{10}_{35}
$,
$ 9-3  c^{1}_{35}
+3c^{2}_{35}
+3c^{3}_{35}
+2c^{4}_{35}
+3c^{5}_{35}
-c^{6}_{35}
+4c^{7}_{35}
-3  c^{8}_{35}
+2c^{9}_{35}
+3c^{10}_{35}
$;\ \ 
$ -6+2c^{1}_{35}
-2  c^{2}_{35}
-2  c^{3}_{35}
-c^{4}_{35}
-2  c^{5}_{35}
+c^{6}_{35}
-4  c^{7}_{35}
+2c^{8}_{35}
-c^{9}_{35}
-2  c^{10}_{35}
$)

  \vskip 2ex

\noindent44. $11_{\frac{8}{5},1964.}^{35,508}$ \irep{2077}:\ \ 
$d_i$ = ($1.0$,
$8.807$,
$8.807$,
$8.807$,
$11.632$,
$13.250$,
$14.250$,
$14.250$,
$14.250$,
$19.822$,
$20.440$) 

\vskip 0.7ex
\hangindent=3em \hangafter=1
$D^2= 1964.590 = 
910-280  c^{1}_{35}
+280c^{2}_{35}
+280c^{3}_{35}
+175c^{4}_{35}
+280c^{5}_{35}
-105  c^{6}_{35}
+490c^{7}_{35}
-280  c^{8}_{35}
+175c^{9}_{35}
+280c^{10}_{35}
$

\vskip 0.7ex
\hangindent=3em \hangafter=1
$T = ( 0,
\frac{3}{35},
\frac{13}{35},
\frac{33}{35},
\frac{4}{5},
0,
\frac{1}{7},
\frac{2}{7},
\frac{4}{7},
\frac{2}{5},
\frac{4}{5} )
$,

\vskip 0.7ex
\hangindent=3em \hangafter=1
$S$ = ($ 1$,
$ 4-c^{1}_{35}
+c^{2}_{35}
+c^{3}_{35}
+c^{4}_{35}
+c^{5}_{35}
+2c^{7}_{35}
-c^{8}_{35}
+c^{9}_{35}
+c^{10}_{35}
$,
$ 4-c^{1}_{35}
+c^{2}_{35}
+c^{3}_{35}
+c^{4}_{35}
+c^{5}_{35}
+2c^{7}_{35}
-c^{8}_{35}
+c^{9}_{35}
+c^{10}_{35}
$,
$ 4-c^{1}_{35}
+c^{2}_{35}
+c^{3}_{35}
+c^{4}_{35}
+c^{5}_{35}
+2c^{7}_{35}
-c^{8}_{35}
+c^{9}_{35}
+c^{10}_{35}
$,
$ 5-2  c^{1}_{35}
+2c^{2}_{35}
+2c^{3}_{35}
+c^{4}_{35}
+2c^{5}_{35}
-c^{6}_{35}
+3c^{7}_{35}
-2  c^{8}_{35}
+c^{9}_{35}
+2c^{10}_{35}
$,
$ 6-2  c^{1}_{35}
+2c^{2}_{35}
+2c^{3}_{35}
+c^{4}_{35}
+2c^{5}_{35}
-c^{6}_{35}
+4c^{7}_{35}
-2  c^{8}_{35}
+c^{9}_{35}
+2c^{10}_{35}
$,
$ 7-2  c^{1}_{35}
+2c^{2}_{35}
+2c^{3}_{35}
+c^{4}_{35}
+2c^{5}_{35}
-c^{6}_{35}
+4c^{7}_{35}
-2  c^{8}_{35}
+c^{9}_{35}
+2c^{10}_{35}
$,
$ 7-2  c^{1}_{35}
+2c^{2}_{35}
+2c^{3}_{35}
+c^{4}_{35}
+2c^{5}_{35}
-c^{6}_{35}
+4c^{7}_{35}
-2  c^{8}_{35}
+c^{9}_{35}
+2c^{10}_{35}
$,
$ 7-2  c^{1}_{35}
+2c^{2}_{35}
+2c^{3}_{35}
+c^{4}_{35}
+2c^{5}_{35}
-c^{6}_{35}
+4c^{7}_{35}
-2  c^{8}_{35}
+c^{9}_{35}
+2c^{10}_{35}
$,
$ 9-3  c^{1}_{35}
+3c^{2}_{35}
+3c^{3}_{35}
+2c^{4}_{35}
+3c^{5}_{35}
-c^{6}_{35}
+4c^{7}_{35}
-3  c^{8}_{35}
+2c^{9}_{35}
+3c^{10}_{35}
$,
$ 9-3  c^{1}_{35}
+3c^{2}_{35}
+3c^{3}_{35}
+2c^{4}_{35}
+3c^{5}_{35}
-c^{6}_{35}
+5c^{7}_{35}
-3  c^{8}_{35}
+2c^{9}_{35}
+3c^{10}_{35}
$;\ \ 
$ 1+5c^{1}_{35}
+4c^{3}_{35}
+3c^{4}_{35}
+c^{5}_{35}
+4c^{6}_{35}
+2c^{8}_{35}
+3c^{10}_{35}
+c^{11}_{35}
$,
$ 1-c^{1}_{35}
+4c^{2}_{35}
-2  c^{3}_{35}
+c^{4}_{35}
+2c^{5}_{35}
-c^{6}_{35}
+c^{7}_{35}
+2c^{9}_{35}
-3  c^{10}_{35}
+2c^{11}_{35}
$,
$ 5-6  c^{1}_{35}
-2  c^{2}_{35}
-3  c^{4}_{35}
-c^{5}_{35}
-4  c^{6}_{35}
+3c^{7}_{35}
-4  c^{8}_{35}
-c^{9}_{35}
+2c^{10}_{35}
-3  c^{11}_{35}
$,
$ 7-2  c^{1}_{35}
+2c^{2}_{35}
+2c^{3}_{35}
+c^{4}_{35}
+2c^{5}_{35}
-c^{6}_{35}
+4c^{7}_{35}
-2  c^{8}_{35}
+c^{9}_{35}
+2c^{10}_{35}
$,
$ 4-c^{1}_{35}
+c^{2}_{35}
+c^{3}_{35}
+c^{4}_{35}
+c^{5}_{35}
+2c^{7}_{35}
-c^{8}_{35}
+c^{9}_{35}
+c^{10}_{35}
$,
$ -3+4c^{1}_{35}
+c^{2}_{35}
+2c^{4}_{35}
+c^{5}_{35}
+2c^{6}_{35}
-2  c^{7}_{35}
+2c^{8}_{35}
+c^{9}_{35}
-c^{10}_{35}
+2c^{11}_{35}
$,
$ -1-3  c^{1}_{35}
-2  c^{3}_{35}
-2  c^{4}_{35}
-c^{5}_{35}
-2  c^{6}_{35}
-c^{8}_{35}
-2  c^{10}_{35}
$,
$ -2  c^{2}_{35}
+c^{3}_{35}
-c^{4}_{35}
-c^{5}_{35}
-2  c^{9}_{35}
+2c^{10}_{35}
-2  c^{11}_{35}
$,
$0$,
$ -7+2c^{1}_{35}
-2  c^{2}_{35}
-2  c^{3}_{35}
-c^{4}_{35}
-2  c^{5}_{35}
+c^{6}_{35}
-4  c^{7}_{35}
+2c^{8}_{35}
-c^{9}_{35}
-2  c^{10}_{35}
$;\ \ 
$ 5-6  c^{1}_{35}
-2  c^{2}_{35}
-3  c^{4}_{35}
-c^{5}_{35}
-4  c^{6}_{35}
+3c^{7}_{35}
-4  c^{8}_{35}
-c^{9}_{35}
+2c^{10}_{35}
-3  c^{11}_{35}
$,
$ 1+5c^{1}_{35}
+4c^{3}_{35}
+3c^{4}_{35}
+c^{5}_{35}
+4c^{6}_{35}
+2c^{8}_{35}
+3c^{10}_{35}
+c^{11}_{35}
$,
$ 7-2  c^{1}_{35}
+2c^{2}_{35}
+2c^{3}_{35}
+c^{4}_{35}
+2c^{5}_{35}
-c^{6}_{35}
+4c^{7}_{35}
-2  c^{8}_{35}
+c^{9}_{35}
+2c^{10}_{35}
$,
$ 4-c^{1}_{35}
+c^{2}_{35}
+c^{3}_{35}
+c^{4}_{35}
+c^{5}_{35}
+2c^{7}_{35}
-c^{8}_{35}
+c^{9}_{35}
+c^{10}_{35}
$,
$ -1-3  c^{1}_{35}
-2  c^{3}_{35}
-2  c^{4}_{35}
-c^{5}_{35}
-2  c^{6}_{35}
-c^{8}_{35}
-2  c^{10}_{35}
$,
$ -2  c^{2}_{35}
+c^{3}_{35}
-c^{4}_{35}
-c^{5}_{35}
-2  c^{9}_{35}
+2c^{10}_{35}
-2  c^{11}_{35}
$,
$ -3+4c^{1}_{35}
+c^{2}_{35}
+2c^{4}_{35}
+c^{5}_{35}
+2c^{6}_{35}
-2  c^{7}_{35}
+2c^{8}_{35}
+c^{9}_{35}
-c^{10}_{35}
+2c^{11}_{35}
$,
$0$,
$ -7+2c^{1}_{35}
-2  c^{2}_{35}
-2  c^{3}_{35}
-c^{4}_{35}
-2  c^{5}_{35}
+c^{6}_{35}
-4  c^{7}_{35}
+2c^{8}_{35}
-c^{9}_{35}
-2  c^{10}_{35}
$;\ \ 
$ 1-c^{1}_{35}
+4c^{2}_{35}
-2  c^{3}_{35}
+c^{4}_{35}
+2c^{5}_{35}
-c^{6}_{35}
+c^{7}_{35}
+2c^{9}_{35}
-3  c^{10}_{35}
+2c^{11}_{35}
$,
$ 7-2  c^{1}_{35}
+2c^{2}_{35}
+2c^{3}_{35}
+c^{4}_{35}
+2c^{5}_{35}
-c^{6}_{35}
+4c^{7}_{35}
-2  c^{8}_{35}
+c^{9}_{35}
+2c^{10}_{35}
$,
$ 4-c^{1}_{35}
+c^{2}_{35}
+c^{3}_{35}
+c^{4}_{35}
+c^{5}_{35}
+2c^{7}_{35}
-c^{8}_{35}
+c^{9}_{35}
+c^{10}_{35}
$,
$ -2  c^{2}_{35}
+c^{3}_{35}
-c^{4}_{35}
-c^{5}_{35}
-2  c^{9}_{35}
+2c^{10}_{35}
-2  c^{11}_{35}
$,
$ -3+4c^{1}_{35}
+c^{2}_{35}
+2c^{4}_{35}
+c^{5}_{35}
+2c^{6}_{35}
-2  c^{7}_{35}
+2c^{8}_{35}
+c^{9}_{35}
-c^{10}_{35}
+2c^{11}_{35}
$,
$ -1-3  c^{1}_{35}
-2  c^{3}_{35}
-2  c^{4}_{35}
-c^{5}_{35}
-2  c^{6}_{35}
-c^{8}_{35}
-2  c^{10}_{35}
$,
$0$,
$ -7+2c^{1}_{35}
-2  c^{2}_{35}
-2  c^{3}_{35}
-c^{4}_{35}
-2  c^{5}_{35}
+c^{6}_{35}
-4  c^{7}_{35}
+2c^{8}_{35}
-c^{9}_{35}
-2  c^{10}_{35}
$;\ \ 
$ -6+2c^{1}_{35}
-2  c^{2}_{35}
-2  c^{3}_{35}
-c^{4}_{35}
-2  c^{5}_{35}
+c^{6}_{35}
-4  c^{7}_{35}
+2c^{8}_{35}
-c^{9}_{35}
-2  c^{10}_{35}
$,
$ -9+3c^{1}_{35}
-3  c^{2}_{35}
-3  c^{3}_{35}
-2  c^{4}_{35}
-3  c^{5}_{35}
+c^{6}_{35}
-5  c^{7}_{35}
+3c^{8}_{35}
-2  c^{9}_{35}
-3  c^{10}_{35}
$,
$ -4+c^{1}_{35}
-c^{2}_{35}
-c^{3}_{35}
-c^{4}_{35}
-c^{5}_{35}
-2  c^{7}_{35}
+c^{8}_{35}
-c^{9}_{35}
-c^{10}_{35}
$,
$ -4+c^{1}_{35}
-c^{2}_{35}
-c^{3}_{35}
-c^{4}_{35}
-c^{5}_{35}
-2  c^{7}_{35}
+c^{8}_{35}
-c^{9}_{35}
-c^{10}_{35}
$,
$ -4+c^{1}_{35}
-c^{2}_{35}
-c^{3}_{35}
-c^{4}_{35}
-c^{5}_{35}
-2  c^{7}_{35}
+c^{8}_{35}
-c^{9}_{35}
-c^{10}_{35}
$,
$ 9-3  c^{1}_{35}
+3c^{2}_{35}
+3c^{3}_{35}
+2c^{4}_{35}
+3c^{5}_{35}
-c^{6}_{35}
+4c^{7}_{35}
-3  c^{8}_{35}
+2c^{9}_{35}
+3c^{10}_{35}
$,
$ 1$;\ \ 
$ 1$,
$ 7-2  c^{1}_{35}
+2c^{2}_{35}
+2c^{3}_{35}
+c^{4}_{35}
+2c^{5}_{35}
-c^{6}_{35}
+4c^{7}_{35}
-2  c^{8}_{35}
+c^{9}_{35}
+2c^{10}_{35}
$,
$ 7-2  c^{1}_{35}
+2c^{2}_{35}
+2c^{3}_{35}
+c^{4}_{35}
+2c^{5}_{35}
-c^{6}_{35}
+4c^{7}_{35}
-2  c^{8}_{35}
+c^{9}_{35}
+2c^{10}_{35}
$,
$ 7-2  c^{1}_{35}
+2c^{2}_{35}
+2c^{3}_{35}
+c^{4}_{35}
+2c^{5}_{35}
-c^{6}_{35}
+4c^{7}_{35}
-2  c^{8}_{35}
+c^{9}_{35}
+2c^{10}_{35}
$,
$ -9+3c^{1}_{35}
-3  c^{2}_{35}
-3  c^{3}_{35}
-2  c^{4}_{35}
-3  c^{5}_{35}
+c^{6}_{35}
-4  c^{7}_{35}
+3c^{8}_{35}
-2  c^{9}_{35}
-3  c^{10}_{35}
$,
$ -5+2c^{1}_{35}
-2  c^{2}_{35}
-2  c^{3}_{35}
-c^{4}_{35}
-2  c^{5}_{35}
+c^{6}_{35}
-3  c^{7}_{35}
+2c^{8}_{35}
-c^{9}_{35}
-2  c^{10}_{35}
$;\ \ 
$ -1+c^{1}_{35}
-4  c^{2}_{35}
+2c^{3}_{35}
-c^{4}_{35}
-2  c^{5}_{35}
+c^{6}_{35}
-c^{7}_{35}
-2  c^{9}_{35}
+3c^{10}_{35}
-2  c^{11}_{35}
$,
$ -5+6c^{1}_{35}
+2c^{2}_{35}
+3c^{4}_{35}
+c^{5}_{35}
+4c^{6}_{35}
-3  c^{7}_{35}
+4c^{8}_{35}
+c^{9}_{35}
-2  c^{10}_{35}
+3c^{11}_{35}
$,
$ -1-5  c^{1}_{35}
-4  c^{3}_{35}
-3  c^{4}_{35}
-c^{5}_{35}
-4  c^{6}_{35}
-2  c^{8}_{35}
-3  c^{10}_{35}
-c^{11}_{35}
$,
$0$,
$ 4-c^{1}_{35}
+c^{2}_{35}
+c^{3}_{35}
+c^{4}_{35}
+c^{5}_{35}
+2c^{7}_{35}
-c^{8}_{35}
+c^{9}_{35}
+c^{10}_{35}
$;\ \ 
$ -1-5  c^{1}_{35}
-4  c^{3}_{35}
-3  c^{4}_{35}
-c^{5}_{35}
-4  c^{6}_{35}
-2  c^{8}_{35}
-3  c^{10}_{35}
-c^{11}_{35}
$,
$ -1+c^{1}_{35}
-4  c^{2}_{35}
+2c^{3}_{35}
-c^{4}_{35}
-2  c^{5}_{35}
+c^{6}_{35}
-c^{7}_{35}
-2  c^{9}_{35}
+3c^{10}_{35}
-2  c^{11}_{35}
$,
$0$,
$ 4-c^{1}_{35}
+c^{2}_{35}
+c^{3}_{35}
+c^{4}_{35}
+c^{5}_{35}
+2c^{7}_{35}
-c^{8}_{35}
+c^{9}_{35}
+c^{10}_{35}
$;\ \ 
$ -5+6c^{1}_{35}
+2c^{2}_{35}
+3c^{4}_{35}
+c^{5}_{35}
+4c^{6}_{35}
-3  c^{7}_{35}
+4c^{8}_{35}
+c^{9}_{35}
-2  c^{10}_{35}
+3c^{11}_{35}
$,
$0$,
$ 4-c^{1}_{35}
+c^{2}_{35}
+c^{3}_{35}
+c^{4}_{35}
+c^{5}_{35}
+2c^{7}_{35}
-c^{8}_{35}
+c^{9}_{35}
+c^{10}_{35}
$;\ \ 
$ -9+3c^{1}_{35}
-3  c^{2}_{35}
-3  c^{3}_{35}
-2  c^{4}_{35}
-3  c^{5}_{35}
+c^{6}_{35}
-4  c^{7}_{35}
+3c^{8}_{35}
-2  c^{9}_{35}
-3  c^{10}_{35}
$,
$ 9-3  c^{1}_{35}
+3c^{2}_{35}
+3c^{3}_{35}
+2c^{4}_{35}
+3c^{5}_{35}
-c^{6}_{35}
+4c^{7}_{35}
-3  c^{8}_{35}
+2c^{9}_{35}
+3c^{10}_{35}
$;\ \ 
$ -6+2c^{1}_{35}
-2  c^{2}_{35}
-2  c^{3}_{35}
-c^{4}_{35}
-2  c^{5}_{35}
+c^{6}_{35}
-4  c^{7}_{35}
+2c^{8}_{35}
-c^{9}_{35}
-2  c^{10}_{35}
$)

  \vskip 2ex 

}

\subsection{Rank 12}
\label{uni12}

{\small
\input{SsL12U_}
}

\section{Full lists of modular data}
\label{allMD}

In this section, we list all the modular data for rank 2 -- 12.  The blue and
black entries, plus the following grey entries, correspond to the set of the
modular data that forms a Galois orbit.  The list is ordered by total quantum
dimension $D^2$.  Each modular data is labeled by two integers (iGO, iMD),
where iGO is the index for Galois orbit, and iMD is the index for modular data
in each Galois orbit.  We also provide GAP/human readable files for those
modular data.  The file NsdGOL\#.g defines a GAP record, NsdGOL, for rank-\#
modular data.  The  modular data labeled by (iGO, iMD) in the following list is
given by NsdGOL[iGO][iMD]:\\
\noindent
NsdGOL[iGO][iMD].Nij\_k is the fusion coefficient $N^{ij}_k =
\text{NsdGOL[iGO][iMD].Nij\_k}[i][j][k] $.\\
  NsdGOL[iGO][iMD].s is a list of the
topological spins $s_i = \text{NsdGOL[iGO][iMD]}.s[i]$.\\
  NsdGOL[iGO][iMD].d is
a list of the quantum dimensions $d_i = \text{NsdGOL[iGO][iMD]}.d[i]$.\\
NsdGOL[iGO][iMD].df is a list of the quantum dimensions in floating point
format.

\subsection{Rank 2}
\label{long2}

{\small
\black

\noindent(1,1). $2_{1,2.}^{4,437}$ \irep{0}:\ \ 
$d_i$ = ($1.0$,
$1.0$) 

\vskip 0.7ex
\hangindent=3em \hangafter=1
$D^2= 2.0 = 
2$

\vskip 0.7ex
\hangindent=3em \hangafter=1
$T = ( 0,
\frac{1}{4} )
$,

\vskip 0.7ex
\hangindent=3em \hangafter=1
$S$ = ($ 1$,
$ 1$;\ \ 
$ -1$)

Prime. 

\vskip 1ex 
\color{grey}

\noindent(1,2). $2_{7,2.}^{4,625}$ \irep{0}:\ \ 
$d_i$ = ($1.0$,
$1.0$) 

\vskip 0.7ex
\hangindent=3em \hangafter=1
$D^2= 2.0 = 
2$

\vskip 0.7ex
\hangindent=3em \hangafter=1
$T = ( 0,
\frac{3}{4} )
$,

\vskip 0.7ex
\hangindent=3em \hangafter=1
$S$ = ($ 1$,
$ 1$;\ \ 
$ -1$)

Prime. 

\vskip 1ex 
\black

\noindent(2,1). $2_{\frac{14}{5},3.618}^{5,395}$ \irep{2}:\ \ 
$d_i$ = ($1.0$,
$1.618$) 

\vskip 0.7ex
\hangindent=3em \hangafter=1
$D^2= 3.618 = 
\frac{5+\sqrt{5}}{2}$

\vskip 0.7ex
\hangindent=3em \hangafter=1
$T = ( 0,
\frac{2}{5} )
$,

\vskip 0.7ex
\hangindent=3em \hangafter=1
$S$ = ($ 1$,
$ \frac{1+\sqrt{5}}{2}$;\ \ 
$ -1$)

Prime. 

\vskip 1ex 
\color{grey}

\noindent(2,2). $2_{\frac{26}{5},3.618}^{5,720}$ \irep{2}:\ \ 
$d_i$ = ($1.0$,
$1.618$) 

\vskip 0.7ex
\hangindent=3em \hangafter=1
$D^2= 3.618 = 
\frac{5+\sqrt{5}}{2}$

\vskip 0.7ex
\hangindent=3em \hangafter=1
$T = ( 0,
\frac{3}{5} )
$,

\vskip 0.7ex
\hangindent=3em \hangafter=1
$S$ = ($ 1$,
$ \frac{1+\sqrt{5}}{2}$;\ \ 
$ -1$)

Prime. 

\vskip 1ex 
\color{grey}

\noindent(2,3). $2_{\frac{2}{5},1.381}^{5,120}$ \irep{2}:\ \ 
$d_i$ = ($1.0$,
$-0.618$) 

\vskip 0.7ex
\hangindent=3em \hangafter=1
$D^2= 1.381 = 
\frac{5-\sqrt{5}}{2}$

\vskip 0.7ex
\hangindent=3em \hangafter=1
$T = ( 0,
\frac{1}{5} )
$,

\vskip 0.7ex
\hangindent=3em \hangafter=1
$S$ = ($ 1$,
$ \frac{1-\sqrt{5}}{2}$;\ \ 
$ -1$)

Prime. 

Not pseudo-unitary. 

\vskip 1ex 
\color{grey}

\noindent(2,4). $2_{\frac{38}{5},1.381}^{5,491}$ \irep{2}:\ \ 
$d_i$ = ($1.0$,
$-0.618$) 

\vskip 0.7ex
\hangindent=3em \hangafter=1
$D^2= 1.381 = 
\frac{5-\sqrt{5}}{2}$

\vskip 0.7ex
\hangindent=3em \hangafter=1
$T = ( 0,
\frac{4}{5} )
$,

\vskip 0.7ex
\hangindent=3em \hangafter=1
$S$ = ($ 1$,
$ \frac{1-\sqrt{5}}{2}$;\ \ 
$ -1$)

Prime. 

Not pseudo-unitary. 

\vskip 1ex 
\color{blue}

\noindent(3,1). $2_{1,2.}^{4,625}$ \irep{0}:\ \ 
$d_i$ = ($1.0$,
$-1.0$) 

\vskip 0.7ex
\hangindent=3em \hangafter=1
$D^2= 2.0 = 
2$

\vskip 0.7ex
\hangindent=3em \hangafter=1
$T = ( 0,
\frac{1}{4} )
$,

\vskip 0.7ex
\hangindent=3em \hangafter=1
$S$ = ($ 1$,
$ -1$;\ \ 
$ -1$)

Prime. 

Pseudo-unitary $\sim$  
$2_{7,2.}^{4,625}$

\vskip 1ex 
\color{grey}

\noindent(3,2). $2_{7,2.}^{4,562}$ \irep{0}:\ \ 
$d_i$ = ($1.0$,
$-1.0$) 

\vskip 0.7ex
\hangindent=3em \hangafter=1
$D^2= 2.0 = 
2$

\vskip 0.7ex
\hangindent=3em \hangafter=1
$T = ( 0,
\frac{3}{4} )
$,

\vskip 0.7ex
\hangindent=3em \hangafter=1
$S$ = ($ 1$,
$ -1$;\ \ 
$ -1$)

Prime. 

Pseudo-unitary $\sim$  
$2_{1,2.}^{4,437}$

\vskip 1ex 

}

\subsection{Rank 3}

{\small
\black

\noindent(1,1). $3_{2,3.}^{3,527}$ \irep{0}:\ \ 
$d_i$ = ($1.0$,
$1.0$,
$1.0$) 

\vskip 0.7ex
\hangindent=3em \hangafter=1
$D^2= 3.0 = 
3$

\vskip 0.7ex
\hangindent=3em \hangafter=1
$T = ( 0,
\frac{1}{3},
\frac{1}{3} )
$,

\vskip 0.7ex
\hangindent=3em \hangafter=1
$S$ = ($ 1$,
$ 1$,
$ 1$;\ \ 
$ \zeta_{3}^{1}$,
$ -\zeta_{6}^{1}$;\ \ 
$ \zeta_{3}^{1}$)

Prime. 

\vskip 1ex 
\color{grey}

\noindent(1,2). $3_{6,3.}^{3,138}$ \irep{0}:\ \ 
$d_i$ = ($1.0$,
$1.0$,
$1.0$) 

\vskip 0.7ex
\hangindent=3em \hangafter=1
$D^2= 3.0 = 
3$

\vskip 0.7ex
\hangindent=3em \hangafter=1
$T = ( 0,
\frac{2}{3},
\frac{2}{3} )
$,

\vskip 0.7ex
\hangindent=3em \hangafter=1
$S$ = ($ 1$,
$ 1$,
$ 1$;\ \ 
$ -\zeta_{6}^{1}$,
$ \zeta_{3}^{1}$;\ \ 
$ -\zeta_{6}^{1}$)

Prime. 

\vskip 1ex 
\black

\noindent(2,1). $3_{\frac{1}{2},4.}^{16,598}$ \irep{4}:\ \ 
$d_i$ = ($1.0$,
$1.0$,
$1.414$) 

\vskip 0.7ex
\hangindent=3em \hangafter=1
$D^2= 4.0 = 
4$

\vskip 0.7ex
\hangindent=3em \hangafter=1
$T = ( 0,
\frac{1}{2},
\frac{1}{16} )
$,

\vskip 0.7ex
\hangindent=3em \hangafter=1
$S$ = ($ 1$,
$ 1$,
$ \sqrt{2}$;\ \ 
$ 1$,
$ -\sqrt{2}$;\ \ 
$0$)

Prime. 

\vskip 1ex 
\color{grey}

\noindent(2,2). $3_{\frac{7}{2},4.}^{16,332}$ \irep{4}:\ \ 
$d_i$ = ($1.0$,
$1.0$,
$1.414$) 

\vskip 0.7ex
\hangindent=3em \hangafter=1
$D^2= 4.0 = 
4$

\vskip 0.7ex
\hangindent=3em \hangafter=1
$T = ( 0,
\frac{1}{2},
\frac{7}{16} )
$,

\vskip 0.7ex
\hangindent=3em \hangafter=1
$S$ = ($ 1$,
$ 1$,
$ \sqrt{2}$;\ \ 
$ 1$,
$ -\sqrt{2}$;\ \ 
$0$)

Prime. 

\vskip 1ex 
\color{grey}

\noindent(2,3). $3_{\frac{9}{2},4.}^{16,156}$ \irep{4}:\ \ 
$d_i$ = ($1.0$,
$1.0$,
$1.414$) 

\vskip 0.7ex
\hangindent=3em \hangafter=1
$D^2= 4.0 = 
4$

\vskip 0.7ex
\hangindent=3em \hangafter=1
$T = ( 0,
\frac{1}{2},
\frac{9}{16} )
$,

\vskip 0.7ex
\hangindent=3em \hangafter=1
$S$ = ($ 1$,
$ 1$,
$ \sqrt{2}$;\ \ 
$ 1$,
$ -\sqrt{2}$;\ \ 
$0$)

Prime. 

\vskip 1ex 
\color{grey}

\noindent(2,4). $3_{\frac{15}{2},4.}^{16,639}$ \irep{4}:\ \ 
$d_i$ = ($1.0$,
$1.0$,
$1.414$) 

\vskip 0.7ex
\hangindent=3em \hangafter=1
$D^2= 4.0 = 
4$

\vskip 0.7ex
\hangindent=3em \hangafter=1
$T = ( 0,
\frac{1}{2},
\frac{15}{16} )
$,

\vskip 0.7ex
\hangindent=3em \hangafter=1
$S$ = ($ 1$,
$ 1$,
$ \sqrt{2}$;\ \ 
$ 1$,
$ -\sqrt{2}$;\ \ 
$0$)

Prime. 

\vskip 1ex 
\color{grey}

\noindent(2,5). $3_{\frac{3}{2},4.}^{16,538}$ \irep{4}:\ \ 
$d_i$ = ($1.0$,
$1.0$,
$-1.414$) 

\vskip 0.7ex
\hangindent=3em \hangafter=1
$D^2= 4.0 = 
4$

\vskip 0.7ex
\hangindent=3em \hangafter=1
$T = ( 0,
\frac{1}{2},
\frac{3}{16} )
$,

\vskip 0.7ex
\hangindent=3em \hangafter=1
$S$ = ($ 1$,
$ 1$,
$ -\sqrt{2}$;\ \ 
$ 1$,
$ \sqrt{2}$;\ \ 
$0$)

Prime. 

Pseudo-unitary $\sim$  
$3_{\frac{11}{2},4.}^{16,648}$

\vskip 1ex 
\color{grey}

\noindent(2,6). $3_{\frac{5}{2},4.}^{16,345}$ \irep{4}:\ \ 
$d_i$ = ($1.0$,
$1.0$,
$-1.414$) 

\vskip 0.7ex
\hangindent=3em \hangafter=1
$D^2= 4.0 = 
4$

\vskip 0.7ex
\hangindent=3em \hangafter=1
$T = ( 0,
\frac{1}{2},
\frac{5}{16} )
$,

\vskip 0.7ex
\hangindent=3em \hangafter=1
$S$ = ($ 1$,
$ 1$,
$ -\sqrt{2}$;\ \ 
$ 1$,
$ \sqrt{2}$;\ \ 
$0$)

Prime. 

Pseudo-unitary $\sim$  
$3_{\frac{13}{2},4.}^{16,330}$

\vskip 1ex 
\color{grey}

\noindent(2,7). $3_{\frac{11}{2},4.}^{16,564}$ \irep{4}:\ \ 
$d_i$ = ($1.0$,
$1.0$,
$-1.414$) 

\vskip 0.7ex
\hangindent=3em \hangafter=1
$D^2= 4.0 = 
4$

\vskip 0.7ex
\hangindent=3em \hangafter=1
$T = ( 0,
\frac{1}{2},
\frac{11}{16} )
$,

\vskip 0.7ex
\hangindent=3em \hangafter=1
$S$ = ($ 1$,
$ 1$,
$ -\sqrt{2}$;\ \ 
$ 1$,
$ \sqrt{2}$;\ \ 
$0$)

Prime. 

Pseudo-unitary $\sim$  
$3_{\frac{3}{2},4.}^{16,553}$

\vskip 1ex 
\color{grey}

\noindent(2,8). $3_{\frac{13}{2},4.}^{16,830}$ \irep{4}:\ \ 
$d_i$ = ($1.0$,
$1.0$,
$-1.414$) 

\vskip 0.7ex
\hangindent=3em \hangafter=1
$D^2= 4.0 = 
4$

\vskip 0.7ex
\hangindent=3em \hangafter=1
$T = ( 0,
\frac{1}{2},
\frac{13}{16} )
$,

\vskip 0.7ex
\hangindent=3em \hangafter=1
$S$ = ($ 1$,
$ 1$,
$ -\sqrt{2}$;\ \ 
$ 1$,
$ \sqrt{2}$;\ \ 
$0$)

Prime. 

Pseudo-unitary $\sim$  
$3_{\frac{5}{2},4.}^{16,465}$

\vskip 1ex 
\black

\noindent(3,1). $3_{\frac{3}{2},4.}^{16,553}$ \irep{4}:\ \ 
$d_i$ = ($1.0$,
$1.0$,
$1.414$) 

\vskip 0.7ex
\hangindent=3em \hangafter=1
$D^2= 4.0 = 
4$

\vskip 0.7ex
\hangindent=3em \hangafter=1
$T = ( 0,
\frac{1}{2},
\frac{3}{16} )
$,

\vskip 0.7ex
\hangindent=3em \hangafter=1
$S$ = ($ 1$,
$ 1$,
$ \sqrt{2}$;\ \ 
$ 1$,
$ -\sqrt{2}$;\ \ 
$0$)

Prime. 

\vskip 1ex 
\color{grey}

\noindent(3,2). $3_{\frac{5}{2},4.}^{16,465}$ \irep{4}:\ \ 
$d_i$ = ($1.0$,
$1.0$,
$1.414$) 

\vskip 0.7ex
\hangindent=3em \hangafter=1
$D^2= 4.0 = 
4$

\vskip 0.7ex
\hangindent=3em \hangafter=1
$T = ( 0,
\frac{1}{2},
\frac{5}{16} )
$,

\vskip 0.7ex
\hangindent=3em \hangafter=1
$S$ = ($ 1$,
$ 1$,
$ \sqrt{2}$;\ \ 
$ 1$,
$ -\sqrt{2}$;\ \ 
$0$)

Prime. 

\vskip 1ex 
\color{grey}

\noindent(3,3). $3_{\frac{11}{2},4.}^{16,648}$ \irep{4}:\ \ 
$d_i$ = ($1.0$,
$1.0$,
$1.414$) 

\vskip 0.7ex
\hangindent=3em \hangafter=1
$D^2= 4.0 = 
4$

\vskip 0.7ex
\hangindent=3em \hangafter=1
$T = ( 0,
\frac{1}{2},
\frac{11}{16} )
$,

\vskip 0.7ex
\hangindent=3em \hangafter=1
$S$ = ($ 1$,
$ 1$,
$ \sqrt{2}$;\ \ 
$ 1$,
$ -\sqrt{2}$;\ \ 
$0$)

Prime. 

\vskip 1ex 
\color{grey}

\noindent(3,4). $3_{\frac{13}{2},4.}^{16,330}$ \irep{4}:\ \ 
$d_i$ = ($1.0$,
$1.0$,
$1.414$) 

\vskip 0.7ex
\hangindent=3em \hangafter=1
$D^2= 4.0 = 
4$

\vskip 0.7ex
\hangindent=3em \hangafter=1
$T = ( 0,
\frac{1}{2},
\frac{13}{16} )
$,

\vskip 0.7ex
\hangindent=3em \hangafter=1
$S$ = ($ 1$,
$ 1$,
$ \sqrt{2}$;\ \ 
$ 1$,
$ -\sqrt{2}$;\ \ 
$0$)

Prime. 

\vskip 1ex 
\color{grey}

\noindent(3,5). $3_{\frac{1}{2},4.}^{16,980}$ \irep{4}:\ \ 
$d_i$ = ($1.0$,
$1.0$,
$-1.414$) 

\vskip 0.7ex
\hangindent=3em \hangafter=1
$D^2= 4.0 = 
4$

\vskip 0.7ex
\hangindent=3em \hangafter=1
$T = ( 0,
\frac{1}{2},
\frac{1}{16} )
$,

\vskip 0.7ex
\hangindent=3em \hangafter=1
$S$ = ($ 1$,
$ 1$,
$ -\sqrt{2}$;\ \ 
$ 1$,
$ \sqrt{2}$;\ \ 
$0$)

Prime. 

Pseudo-unitary $\sim$  
$3_{\frac{9}{2},4.}^{16,156}$

\vskip 1ex 
\color{grey}

\noindent(3,6). $3_{\frac{7}{2},4.}^{16,167}$ \irep{4}:\ \ 
$d_i$ = ($1.0$,
$1.0$,
$-1.414$) 

\vskip 0.7ex
\hangindent=3em \hangafter=1
$D^2= 4.0 = 
4$

\vskip 0.7ex
\hangindent=3em \hangafter=1
$T = ( 0,
\frac{1}{2},
\frac{7}{16} )
$,

\vskip 0.7ex
\hangindent=3em \hangafter=1
$S$ = ($ 1$,
$ 1$,
$ -\sqrt{2}$;\ \ 
$ 1$,
$ \sqrt{2}$;\ \ 
$0$)

Prime. 

Pseudo-unitary $\sim$  
$3_{\frac{15}{2},4.}^{16,639}$

\vskip 1ex 
\color{grey}

\noindent(3,7). $3_{\frac{9}{2},4.}^{16,343}$ \irep{4}:\ \ 
$d_i$ = ($1.0$,
$1.0$,
$-1.414$) 

\vskip 0.7ex
\hangindent=3em \hangafter=1
$D^2= 4.0 = 
4$

\vskip 0.7ex
\hangindent=3em \hangafter=1
$T = ( 0,
\frac{1}{2},
\frac{9}{16} )
$,

\vskip 0.7ex
\hangindent=3em \hangafter=1
$S$ = ($ 1$,
$ 1$,
$ -\sqrt{2}$;\ \ 
$ 1$,
$ \sqrt{2}$;\ \ 
$0$)

Prime. 

Pseudo-unitary $\sim$  
$3_{\frac{1}{2},4.}^{16,598}$

\vskip 1ex 
\color{grey}

\noindent(3,8). $3_{\frac{15}{2},4.}^{16,113}$ \irep{4}:\ \ 
$d_i$ = ($1.0$,
$1.0$,
$-1.414$) 

\vskip 0.7ex
\hangindent=3em \hangafter=1
$D^2= 4.0 = 
4$

\vskip 0.7ex
\hangindent=3em \hangafter=1
$T = ( 0,
\frac{1}{2},
\frac{15}{16} )
$,

\vskip 0.7ex
\hangindent=3em \hangafter=1
$S$ = ($ 1$,
$ 1$,
$ -\sqrt{2}$;\ \ 
$ 1$,
$ \sqrt{2}$;\ \ 
$0$)

Prime. 

Pseudo-unitary $\sim$  
$3_{\frac{7}{2},4.}^{16,332}$

\vskip 1ex 
\black

\noindent(4,1). $3_{\frac{48}{7},9.295}^{7,790}$ \irep{3}:\ \ 
$d_i$ = ($1.0$,
$1.801$,
$2.246$) 

\vskip 0.7ex
\hangindent=3em \hangafter=1
$D^2= 9.295 = 
6+3c^{1}_{7}
+c^{2}_{7}
$

\vskip 0.7ex
\hangindent=3em \hangafter=1
$T = ( 0,
\frac{1}{7},
\frac{5}{7} )
$,

\vskip 0.7ex
\hangindent=3em \hangafter=1
$S$ = ($ 1$,
$ -c_{7}^{3}$,
$ \xi_{7}^{3}$;\ \ 
$ -\xi_{7}^{3}$,
$ 1$;\ \ 
$ c_{7}^{3}$)

Prime. 

\vskip 1ex 
\color{grey}

\noindent(4,2). $3_{\frac{8}{7},9.295}^{7,245}$ \irep{3}:\ \ 
$d_i$ = ($1.0$,
$1.801$,
$2.246$) 

\vskip 0.7ex
\hangindent=3em \hangafter=1
$D^2= 9.295 = 
6+3c^{1}_{7}
+c^{2}_{7}
$

\vskip 0.7ex
\hangindent=3em \hangafter=1
$T = ( 0,
\frac{6}{7},
\frac{2}{7} )
$,

\vskip 0.7ex
\hangindent=3em \hangafter=1
$S$ = ($ 1$,
$ -c_{7}^{3}$,
$ \xi_{7}^{3}$;\ \ 
$ -\xi_{7}^{3}$,
$ 1$;\ \ 
$ c_{7}^{3}$)

Prime. 

\vskip 1ex 
\color{grey}

\noindent(4,3). $3_{\frac{12}{7},2.862}^{7,768}$ \irep{3}:\ \ 
$d_i$ = ($1.0$,
$0.554$,
$-1.246$) 

\vskip 0.7ex
\hangindent=3em \hangafter=1
$D^2= 2.862 = 
5-c^{1}_{7}
+2c^{2}_{7}
$

\vskip 0.7ex
\hangindent=3em \hangafter=1
$T = ( 0,
\frac{3}{7},
\frac{2}{7} )
$,

\vskip 0.7ex
\hangindent=3em \hangafter=1
$S$ = ($ 1$,
$ 1+c^{2}_{7}
$,
$ -c_{7}^{1}$;\ \ 
$ c_{7}^{1}$,
$ 1$;\ \ 
$ -1-c^{2}_{7}
$)

Prime. 

Not pseudo-unitary. 

\vskip 1ex 
\color{grey}

\noindent(4,4). $3_{\frac{44}{7},2.862}^{7,531}$ \irep{3}:\ \ 
$d_i$ = ($1.0$,
$0.554$,
$-1.246$) 

\vskip 0.7ex
\hangindent=3em \hangafter=1
$D^2= 2.862 = 
5-c^{1}_{7}
+2c^{2}_{7}
$

\vskip 0.7ex
\hangindent=3em \hangafter=1
$T = ( 0,
\frac{4}{7},
\frac{5}{7} )
$,

\vskip 0.7ex
\hangindent=3em \hangafter=1
$S$ = ($ 1$,
$ 1+c^{2}_{7}
$,
$ -c_{7}^{1}$;\ \ 
$ c_{7}^{1}$,
$ 1$;\ \ 
$ -1-c^{2}_{7}
$)

Prime. 

Not pseudo-unitary. 

\vskip 1ex 
\color{grey}

\noindent(4,5). $3_{\frac{4}{7},1.841}^{7,953}$ \irep{3}:\ \ 
$d_i$ = ($1.0$,
$0.445$,
$-0.801$) 

\vskip 0.7ex
\hangindent=3em \hangafter=1
$D^2= 1.841 = 
3-2  c^{1}_{7}
-3  c^{2}_{7}
$

\vskip 0.7ex
\hangindent=3em \hangafter=1
$T = ( 0,
\frac{3}{7},
\frac{1}{7} )
$,

\vskip 0.7ex
\hangindent=3em \hangafter=1
$S$ = ($ 1$,
$ -c_{7}^{2}$,
$ -c^{1}_{7}
-c^{2}_{7}
$;\ \ 
$ c^{1}_{7}
+c^{2}_{7}
$,
$ 1$;\ \ 
$ c_{7}^{2}$)

Prime. 

Not pseudo-unitary. 

\vskip 1ex 
\color{grey}

\noindent(4,6). $3_{\frac{52}{7},1.841}^{7,604}$ \irep{3}:\ \ 
$d_i$ = ($1.0$,
$0.445$,
$-0.801$) 

\vskip 0.7ex
\hangindent=3em \hangafter=1
$D^2= 1.841 = 
3-2  c^{1}_{7}
-3  c^{2}_{7}
$

\vskip 0.7ex
\hangindent=3em \hangafter=1
$T = ( 0,
\frac{4}{7},
\frac{6}{7} )
$,

\vskip 0.7ex
\hangindent=3em \hangafter=1
$S$ = ($ 1$,
$ -c_{7}^{2}$,
$ -c^{1}_{7}
-c^{2}_{7}
$;\ \ 
$ c^{1}_{7}
+c^{2}_{7}
$,
$ 1$;\ \ 
$ c_{7}^{2}$)

Prime. 

Not pseudo-unitary. 

\vskip 1ex 

}

\subsection{Rank 4}
\label{long4}

{\small
\black

\noindent(1,1). $4_{0,4.}^{2,750}$ \irep{0}:\ \ 
$d_i$ = ($1.0$,
$1.0$,
$1.0$,
$1.0$) 

\vskip 0.7ex
\hangindent=3em \hangafter=1
$D^2= 4.0 = 
4$

\vskip 0.7ex
\hangindent=3em \hangafter=1
$T = ( 0,
0,
0,
\frac{1}{2} )
$,

\vskip 0.7ex
\hangindent=3em \hangafter=1
$S$ = ($ 1$,
$ 1$,
$ 1$,
$ 1$;\ \ 
$ 1$,
$ -1$,
$ -1$;\ \ 
$ 1$,
$ -1$;\ \ 
$ 1$)

Prime. 

\vskip 1ex 
\black

\noindent(2,1). $4_{4,4.}^{2,250}$ \irep{0}:\ \ 
$d_i$ = ($1.0$,
$1.0$,
$1.0$,
$1.0$) 

\vskip 0.7ex
\hangindent=3em \hangafter=1
$D^2= 4.0 = 
4$

\vskip 0.7ex
\hangindent=3em \hangafter=1
$T = ( 0,
\frac{1}{2},
\frac{1}{2},
\frac{1}{2} )
$,

\vskip 0.7ex
\hangindent=3em \hangafter=1
$S$ = ($ 1$,
$ 1$,
$ 1$,
$ 1$;\ \ 
$ 1$,
$ -1$,
$ -1$;\ \ 
$ 1$,
$ -1$;\ \ 
$ 1$)

Prime. 

\vskip 1ex 
\black

\noindent(3,1). $4_{0,4.}^{4,375}$ \irep{0}:\ \ 
$d_i$ = ($1.0$,
$1.0$,
$1.0$,
$1.0$) 

\vskip 0.7ex
\hangindent=3em \hangafter=1
$D^2= 4.0 = 
4$

\vskip 0.7ex
\hangindent=3em \hangafter=1
$T = ( 0,
0,
\frac{1}{4},
\frac{3}{4} )
$,

\vskip 0.7ex
\hangindent=3em \hangafter=1
$S$ = ($ 1$,
$ 1$,
$ 1$,
$ 1$;\ \ 
$ 1$,
$ -1$,
$ -1$;\ \ 
$ -1$,
$ 1$;\ \ 
$ -1$)

Factors = $2_{1,2.}^{4,437}\boxtimes 2_{7,2.}^{4,625} $

\vskip 1ex 
\black

\noindent(4,1). $4_{2,4.}^{4,625}$ \irep{0}:\ \ 
$d_i$ = ($1.0$,
$1.0$,
$1.0$,
$1.0$) 

\vskip 0.7ex
\hangindent=3em \hangafter=1
$D^2= 4.0 = 
4$

\vskip 0.7ex
\hangindent=3em \hangafter=1
$T = ( 0,
\frac{1}{2},
\frac{1}{4},
\frac{1}{4} )
$,

\vskip 0.7ex
\hangindent=3em \hangafter=1
$S$ = ($ 1$,
$ 1$,
$ 1$,
$ 1$;\ \ 
$ 1$,
$ -1$,
$ -1$;\ \ 
$ -1$,
$ 1$;\ \ 
$ -1$)

Factors = $2_{1,2.}^{4,437}\boxtimes 2_{1,2.}^{4,437} $

\vskip 1ex 
\color{grey}

\noindent(4,2). $4_{6,4.}^{4,375}$ \irep{0}:\ \ 
$d_i$ = ($1.0$,
$1.0$,
$1.0$,
$1.0$) 

\vskip 0.7ex
\hangindent=3em \hangafter=1
$D^2= 4.0 = 
4$

\vskip 0.7ex
\hangindent=3em \hangafter=1
$T = ( 0,
\frac{1}{2},
\frac{3}{4},
\frac{3}{4} )
$,

\vskip 0.7ex
\hangindent=3em \hangafter=1
$S$ = ($ 1$,
$ 1$,
$ 1$,
$ 1$;\ \ 
$ 1$,
$ -1$,
$ -1$;\ \ 
$ -1$,
$ 1$;\ \ 
$ -1$)

Factors = $2_{7,2.}^{4,625}\boxtimes 2_{7,2.}^{4,625} $

\vskip 1ex 
\black

\noindent(5,1). $4_{1,4.}^{8,718}$ \irep{6}:\ \ 
$d_i$ = ($1.0$,
$1.0$,
$1.0$,
$1.0$) 

\vskip 0.7ex
\hangindent=3em \hangafter=1
$D^2= 4.0 = 
4$

\vskip 0.7ex
\hangindent=3em \hangafter=1
$T = ( 0,
\frac{1}{2},
\frac{1}{8},
\frac{1}{8} )
$,

\vskip 0.7ex
\hangindent=3em \hangafter=1
$S$ = ($ 1$,
$ 1$,
$ 1$,
$ 1$;\ \ 
$ 1$,
$ -1$,
$ -1$;\ \ 
$-\mathrm{i}$,
$\mathrm{i}$;\ \ 
$-\mathrm{i}$)

Prime. 

\vskip 1ex 
\color{grey}

\noindent(5,2). $4_{3,4.}^{8,468}$ \irep{6}:\ \ 
$d_i$ = ($1.0$,
$1.0$,
$1.0$,
$1.0$) 

\vskip 0.7ex
\hangindent=3em \hangafter=1
$D^2= 4.0 = 
4$

\vskip 0.7ex
\hangindent=3em \hangafter=1
$T = ( 0,
\frac{1}{2},
\frac{3}{8},
\frac{3}{8} )
$,

\vskip 0.7ex
\hangindent=3em \hangafter=1
$S$ = ($ 1$,
$ 1$,
$ 1$,
$ 1$;\ \ 
$ 1$,
$ -1$,
$ -1$;\ \ 
$\mathrm{i}$,
$-\mathrm{i}$;\ \ 
$\mathrm{i}$)

Prime. 

\vskip 1ex 
\color{grey}

\noindent(5,3). $4_{5,4.}^{8,312}$ \irep{6}:\ \ 
$d_i$ = ($1.0$,
$1.0$,
$1.0$,
$1.0$) 

\vskip 0.7ex
\hangindent=3em \hangafter=1
$D^2= 4.0 = 
4$

\vskip 0.7ex
\hangindent=3em \hangafter=1
$T = ( 0,
\frac{1}{2},
\frac{5}{8},
\frac{5}{8} )
$,

\vskip 0.7ex
\hangindent=3em \hangafter=1
$S$ = ($ 1$,
$ 1$,
$ 1$,
$ 1$;\ \ 
$ 1$,
$ -1$,
$ -1$;\ \ 
$-\mathrm{i}$,
$\mathrm{i}$;\ \ 
$-\mathrm{i}$)

Prime. 

\vskip 1ex 
\color{grey}

\noindent(5,4). $4_{7,4.}^{8,781}$ \irep{6}:\ \ 
$d_i$ = ($1.0$,
$1.0$,
$1.0$,
$1.0$) 

\vskip 0.7ex
\hangindent=3em \hangafter=1
$D^2= 4.0 = 
4$

\vskip 0.7ex
\hangindent=3em \hangafter=1
$T = ( 0,
\frac{1}{2},
\frac{7}{8},
\frac{7}{8} )
$,

\vskip 0.7ex
\hangindent=3em \hangafter=1
$S$ = ($ 1$,
$ 1$,
$ 1$,
$ 1$;\ \ 
$ 1$,
$ -1$,
$ -1$;\ \ 
$\mathrm{i}$,
$-\mathrm{i}$;\ \ 
$\mathrm{i}$)

Prime. 

\vskip 1ex 
\black

\noindent(6,1). $4_{\frac{19}{5},7.236}^{20,304}$ \irep{8}:\ \ 
$d_i$ = ($1.0$,
$1.0$,
$1.618$,
$1.618$) 

\vskip 0.7ex
\hangindent=3em \hangafter=1
$D^2= 7.236 = 
5+\sqrt{5}$

\vskip 0.7ex
\hangindent=3em \hangafter=1
$T = ( 0,
\frac{1}{4},
\frac{2}{5},
\frac{13}{20} )
$,

\vskip 0.7ex
\hangindent=3em \hangafter=1
$S$ = ($ 1$,
$ 1$,
$ \frac{1+\sqrt{5}}{2}$,
$ \frac{1+\sqrt{5}}{2}$;\ \ 
$ -1$,
$ \frac{1+\sqrt{5}}{2}$,
$ -\frac{1+\sqrt{5}}{2}$;\ \ 
$ -1$,
$ -1$;\ \ 
$ 1$)

Factors = $2_{1,2.}^{4,437}\boxtimes 2_{\frac{14}{5},3.618}^{5,395} $

\vskip 1ex 
\color{grey}

\noindent(6,2). $4_{\frac{31}{5},7.236}^{20,505}$ \irep{8}:\ \ 
$d_i$ = ($1.0$,
$1.0$,
$1.618$,
$1.618$) 

\vskip 0.7ex
\hangindent=3em \hangafter=1
$D^2= 7.236 = 
5+\sqrt{5}$

\vskip 0.7ex
\hangindent=3em \hangafter=1
$T = ( 0,
\frac{1}{4},
\frac{3}{5},
\frac{17}{20} )
$,

\vskip 0.7ex
\hangindent=3em \hangafter=1
$S$ = ($ 1$,
$ 1$,
$ \frac{1+\sqrt{5}}{2}$,
$ \frac{1+\sqrt{5}}{2}$;\ \ 
$ -1$,
$ \frac{1+\sqrt{5}}{2}$,
$ -\frac{1+\sqrt{5}}{2}$;\ \ 
$ -1$,
$ -1$;\ \ 
$ 1$)

Factors = $2_{1,2.}^{4,437}\boxtimes 2_{\frac{26}{5},3.618}^{5,720} $

\vskip 1ex 
\color{grey}

\noindent(6,3). $4_{\frac{9}{5},7.236}^{20,451}$ \irep{8}:\ \ 
$d_i$ = ($1.0$,
$1.0$,
$1.618$,
$1.618$) 

\vskip 0.7ex
\hangindent=3em \hangafter=1
$D^2= 7.236 = 
5+\sqrt{5}$

\vskip 0.7ex
\hangindent=3em \hangafter=1
$T = ( 0,
\frac{3}{4},
\frac{2}{5},
\frac{3}{20} )
$,

\vskip 0.7ex
\hangindent=3em \hangafter=1
$S$ = ($ 1$,
$ 1$,
$ \frac{1+\sqrt{5}}{2}$,
$ \frac{1+\sqrt{5}}{2}$;\ \ 
$ -1$,
$ \frac{1+\sqrt{5}}{2}$,
$ -\frac{1+\sqrt{5}}{2}$;\ \ 
$ -1$,
$ -1$;\ \ 
$ 1$)

Factors = $2_{7,2.}^{4,625}\boxtimes 2_{\frac{14}{5},3.618}^{5,395} $

\vskip 1ex 
\color{grey}

\noindent(6,4). $4_{\frac{21}{5},7.236}^{20,341}$ \irep{8}:\ \ 
$d_i$ = ($1.0$,
$1.0$,
$1.618$,
$1.618$) 

\vskip 0.7ex
\hangindent=3em \hangafter=1
$D^2= 7.236 = 
5+\sqrt{5}$

\vskip 0.7ex
\hangindent=3em \hangafter=1
$T = ( 0,
\frac{3}{4},
\frac{3}{5},
\frac{7}{20} )
$,

\vskip 0.7ex
\hangindent=3em \hangafter=1
$S$ = ($ 1$,
$ 1$,
$ \frac{1+\sqrt{5}}{2}$,
$ \frac{1+\sqrt{5}}{2}$;\ \ 
$ -1$,
$ \frac{1+\sqrt{5}}{2}$,
$ -\frac{1+\sqrt{5}}{2}$;\ \ 
$ -1$,
$ -1$;\ \ 
$ 1$)

Factors = $2_{7,2.}^{4,625}\boxtimes 2_{\frac{26}{5},3.618}^{5,720} $

\vskip 1ex 
\color{grey}

\noindent(6,5). $4_{\frac{7}{5},2.763}^{20,278}$ \irep{8}:\ \ 
$d_i$ = ($1.0$,
$1.0$,
$-0.618$,
$-0.618$) 

\vskip 0.7ex
\hangindent=3em \hangafter=1
$D^2= 2.763 = 
5-\sqrt{5}$

\vskip 0.7ex
\hangindent=3em \hangafter=1
$T = ( 0,
\frac{1}{4},
\frac{1}{5},
\frac{9}{20} )
$,

\vskip 0.7ex
\hangindent=3em \hangafter=1
$S$ = ($ 1$,
$ 1$,
$ \frac{1-\sqrt{5}}{2}$,
$ \frac{1-\sqrt{5}}{2}$;\ \ 
$ -1$,
$ \frac{1-\sqrt{5}}{2}$,
$ -\frac{1-\sqrt{5}}{2}$;\ \ 
$ -1$,
$ -1$;\ \ 
$ 1$)

Factors = $2_{1,2.}^{4,437}\boxtimes 2_{\frac{2}{5},1.381}^{5,120} $

Not pseudo-unitary. 

\vskip 1ex 
\color{grey}

\noindent(6,6). $4_{\frac{3}{5},2.763}^{20,525}$ \irep{8}:\ \ 
$d_i$ = ($1.0$,
$1.0$,
$-0.618$,
$-0.618$) 

\vskip 0.7ex
\hangindent=3em \hangafter=1
$D^2= 2.763 = 
5-\sqrt{5}$

\vskip 0.7ex
\hangindent=3em \hangafter=1
$T = ( 0,
\frac{1}{4},
\frac{4}{5},
\frac{1}{20} )
$,

\vskip 0.7ex
\hangindent=3em \hangafter=1
$S$ = ($ 1$,
$ 1$,
$ \frac{1-\sqrt{5}}{2}$,
$ \frac{1-\sqrt{5}}{2}$;\ \ 
$ -1$,
$ \frac{1-\sqrt{5}}{2}$,
$ -\frac{1-\sqrt{5}}{2}$;\ \ 
$ -1$,
$ -1$;\ \ 
$ 1$)

Factors = $2_{1,2.}^{4,437}\boxtimes 2_{\frac{38}{5},1.381}^{5,491} $

Not pseudo-unitary. 

\vskip 1ex 
\color{grey}

\noindent(6,7). $4_{\frac{37}{5},2.763}^{20,210}$ \irep{8}:\ \ 
$d_i$ = ($1.0$,
$1.0$,
$-0.618$,
$-0.618$) 

\vskip 0.7ex
\hangindent=3em \hangafter=1
$D^2= 2.763 = 
5-\sqrt{5}$

\vskip 0.7ex
\hangindent=3em \hangafter=1
$T = ( 0,
\frac{3}{4},
\frac{1}{5},
\frac{19}{20} )
$,

\vskip 0.7ex
\hangindent=3em \hangafter=1
$S$ = ($ 1$,
$ 1$,
$ \frac{1-\sqrt{5}}{2}$,
$ \frac{1-\sqrt{5}}{2}$;\ \ 
$ -1$,
$ \frac{1-\sqrt{5}}{2}$,
$ -\frac{1-\sqrt{5}}{2}$;\ \ 
$ -1$,
$ -1$;\ \ 
$ 1$)

Factors = $2_{7,2.}^{4,625}\boxtimes 2_{\frac{2}{5},1.381}^{5,120} $

Not pseudo-unitary. 

\vskip 1ex 
\color{grey}

\noindent(6,8). $4_{\frac{33}{5},2.763}^{20,210}$ \irep{8}:\ \ 
$d_i$ = ($1.0$,
$1.0$,
$-0.618$,
$-0.618$) 

\vskip 0.7ex
\hangindent=3em \hangafter=1
$D^2= 2.763 = 
5-\sqrt{5}$

\vskip 0.7ex
\hangindent=3em \hangafter=1
$T = ( 0,
\frac{3}{4},
\frac{4}{5},
\frac{11}{20} )
$,

\vskip 0.7ex
\hangindent=3em \hangafter=1
$S$ = ($ 1$,
$ 1$,
$ \frac{1-\sqrt{5}}{2}$,
$ \frac{1-\sqrt{5}}{2}$;\ \ 
$ -1$,
$ \frac{1-\sqrt{5}}{2}$,
$ -\frac{1-\sqrt{5}}{2}$;\ \ 
$ -1$,
$ -1$;\ \ 
$ 1$)

Factors = $2_{7,2.}^{4,625}\boxtimes 2_{\frac{38}{5},1.381}^{5,491} $

Not pseudo-unitary. 

\vskip 1ex 
\black

\noindent(7,1). $4_{\frac{28}{5},13.09}^{5,479}$ \irep{5}:\ \ 
$d_i$ = ($1.0$,
$1.618$,
$1.618$,
$2.618$) 

\vskip 0.7ex
\hangindent=3em \hangafter=1
$D^2= 13.90 = 
\frac{15+5\sqrt{5}}{2}$

\vskip 0.7ex
\hangindent=3em \hangafter=1
$T = ( 0,
\frac{2}{5},
\frac{2}{5},
\frac{4}{5} )
$,

\vskip 0.7ex
\hangindent=3em \hangafter=1
$S$ = ($ 1$,
$ \frac{1+\sqrt{5}}{2}$,
$ \frac{1+\sqrt{5}}{2}$,
$ \frac{3+\sqrt{5}}{2}$;\ \ 
$ -1$,
$ \frac{3+\sqrt{5}}{2}$,
$ -\frac{1+\sqrt{5}}{2}$;\ \ 
$ -1$,
$ -\frac{1+\sqrt{5}}{2}$;\ \ 
$ 1$)

Factors = $2_{\frac{14}{5},3.618}^{5,395}\boxtimes 2_{\frac{14}{5},3.618}^{5,395} $

\vskip 1ex 
\color{grey}

\noindent(7,2). $4_{\frac{12}{5},13.09}^{5,443}$ \irep{5}:\ \ 
$d_i$ = ($1.0$,
$1.618$,
$1.618$,
$2.618$) 

\vskip 0.7ex
\hangindent=3em \hangafter=1
$D^2= 13.90 = 
\frac{15+5\sqrt{5}}{2}$

\vskip 0.7ex
\hangindent=3em \hangafter=1
$T = ( 0,
\frac{3}{5},
\frac{3}{5},
\frac{1}{5} )
$,

\vskip 0.7ex
\hangindent=3em \hangafter=1
$S$ = ($ 1$,
$ \frac{1+\sqrt{5}}{2}$,
$ \frac{1+\sqrt{5}}{2}$,
$ \frac{3+\sqrt{5}}{2}$;\ \ 
$ -1$,
$ \frac{3+\sqrt{5}}{2}$,
$ -\frac{1+\sqrt{5}}{2}$;\ \ 
$ -1$,
$ -\frac{1+\sqrt{5}}{2}$;\ \ 
$ 1$)

Factors = $2_{\frac{26}{5},3.618}^{5,720}\boxtimes 2_{\frac{26}{5},3.618}^{5,720} $

\vskip 1ex 
\color{grey}

\noindent(7,3). $4_{\frac{4}{5},1.909}^{5,248}$ \irep{5}:\ \ 
$d_i$ = ($1.0$,
$0.381$,
$-0.618$,
$-0.618$) 

\vskip 0.7ex
\hangindent=3em \hangafter=1
$D^2= 1.909 = 
\frac{15-5\sqrt{5}}{2}$

\vskip 0.7ex
\hangindent=3em \hangafter=1
$T = ( 0,
\frac{2}{5},
\frac{1}{5},
\frac{1}{5} )
$,

\vskip 0.7ex
\hangindent=3em \hangafter=1
$S$ = ($ 1$,
$ \frac{3-\sqrt{5}}{2}$,
$ \frac{1-\sqrt{5}}{2}$,
$ \frac{1-\sqrt{5}}{2}$;\ \ 
$ 1$,
$ -\frac{1-\sqrt{5}}{2}$,
$ -\frac{1-\sqrt{5}}{2}$;\ \ 
$ -1$,
$ \frac{3-\sqrt{5}}{2}$;\ \ 
$ -1$)

Factors = $2_{\frac{2}{5},1.381}^{5,120}\boxtimes 2_{\frac{2}{5},1.381}^{5,120} $

Not pseudo-unitary. 

\vskip 1ex 
\color{grey}

\noindent(7,4). $4_{\frac{36}{5},1.909}^{5,690}$ \irep{5}:\ \ 
$d_i$ = ($1.0$,
$0.381$,
$-0.618$,
$-0.618$) 

\vskip 0.7ex
\hangindent=3em \hangafter=1
$D^2= 1.909 = 
\frac{15-5\sqrt{5}}{2}$

\vskip 0.7ex
\hangindent=3em \hangafter=1
$T = ( 0,
\frac{3}{5},
\frac{4}{5},
\frac{4}{5} )
$,

\vskip 0.7ex
\hangindent=3em \hangafter=1
$S$ = ($ 1$,
$ \frac{3-\sqrt{5}}{2}$,
$ \frac{1-\sqrt{5}}{2}$,
$ \frac{1-\sqrt{5}}{2}$;\ \ 
$ 1$,
$ -\frac{1-\sqrt{5}}{2}$,
$ -\frac{1-\sqrt{5}}{2}$;\ \ 
$ -1$,
$ \frac{3-\sqrt{5}}{2}$;\ \ 
$ -1$)

Factors = $2_{\frac{38}{5},1.381}^{5,491}\boxtimes 2_{\frac{38}{5},1.381}^{5,491} $

Not pseudo-unitary. 

\vskip 1ex 
\black

\noindent(8,1). $4_{0,13.09}^{5,872}$ \irep{5}:\ \ 
$d_i$ = ($1.0$,
$1.618$,
$1.618$,
$2.618$) 

\vskip 0.7ex
\hangindent=3em \hangafter=1
$D^2= 13.90 = 
\frac{15+5\sqrt{5}}{2}$

\vskip 0.7ex
\hangindent=3em \hangafter=1
$T = ( 0,
\frac{2}{5},
\frac{3}{5},
0 )
$,

\vskip 0.7ex
\hangindent=3em \hangafter=1
$S$ = ($ 1$,
$ \frac{1+\sqrt{5}}{2}$,
$ \frac{1+\sqrt{5}}{2}$,
$ \frac{3+\sqrt{5}}{2}$;\ \ 
$ -1$,
$ \frac{3+\sqrt{5}}{2}$,
$ -\frac{1+\sqrt{5}}{2}$;\ \ 
$ -1$,
$ -\frac{1+\sqrt{5}}{2}$;\ \ 
$ 1$)

Factors = $2_{\frac{14}{5},3.618}^{5,395}\boxtimes 2_{\frac{26}{5},3.618}^{5,720} $

\vskip 1ex 
\color{grey}

\noindent(8,2). $4_{0,1.909}^{5,456}$ \irep{5}:\ \ 
$d_i$ = ($1.0$,
$0.381$,
$-0.618$,
$-0.618$) 

\vskip 0.7ex
\hangindent=3em \hangafter=1
$D^2= 1.909 = 
\frac{15-5\sqrt{5}}{2}$

\vskip 0.7ex
\hangindent=3em \hangafter=1
$T = ( 0,
0,
\frac{1}{5},
\frac{4}{5} )
$,

\vskip 0.7ex
\hangindent=3em \hangafter=1
$S$ = ($ 1$,
$ \frac{3-\sqrt{5}}{2}$,
$ \frac{1-\sqrt{5}}{2}$,
$ \frac{1-\sqrt{5}}{2}$;\ \ 
$ 1$,
$ -\frac{1-\sqrt{5}}{2}$,
$ -\frac{1-\sqrt{5}}{2}$;\ \ 
$ -1$,
$ \frac{3-\sqrt{5}}{2}$;\ \ 
$ -1$)

Factors = $2_{\frac{2}{5},1.381}^{5,120}\boxtimes 2_{\frac{38}{5},1.381}^{5,491} $

Not pseudo-unitary. 

\vskip 1ex 
\black

\noindent(9,1). $4_{\frac{10}{3},19.23}^{9,459}$ \irep{7}:\ \ 
$d_i$ = ($1.0$,
$1.879$,
$2.532$,
$2.879$) 

\vskip 0.7ex
\hangindent=3em \hangafter=1
$D^2= 19.234 = 
9+6c^{1}_{9}
+3c^{2}_{9}
$

\vskip 0.7ex
\hangindent=3em \hangafter=1
$T = ( 0,
\frac{1}{3},
\frac{2}{9},
\frac{2}{3} )
$,

\vskip 0.7ex
\hangindent=3em \hangafter=1
$S$ = ($ 1$,
$ -c_{9}^{4}$,
$ \xi_{9}^{3}$,
$ \xi_{9}^{5}$;\ \ 
$ -\xi_{9}^{5}$,
$ \xi_{9}^{3}$,
$ -1$;\ \ 
$0$,
$ -\xi_{9}^{3}$;\ \ 
$ -c_{9}^{4}$)

Prime. 

\vskip 1ex 
\color{grey}

\noindent(9,2). $4_{\frac{14}{3},19.23}^{9,614}$ \irep{7}:\ \ 
$d_i$ = ($1.0$,
$1.879$,
$2.532$,
$2.879$) 

\vskip 0.7ex
\hangindent=3em \hangafter=1
$D^2= 19.234 = 
9+6c^{1}_{9}
+3c^{2}_{9}
$

\vskip 0.7ex
\hangindent=3em \hangafter=1
$T = ( 0,
\frac{2}{3},
\frac{7}{9},
\frac{1}{3} )
$,

\vskip 0.7ex
\hangindent=3em \hangafter=1
$S$ = ($ 1$,
$ -c_{9}^{4}$,
$ \xi_{9}^{3}$,
$ \xi_{9}^{5}$;\ \ 
$ -\xi_{9}^{5}$,
$ \xi_{9}^{3}$,
$ -1$;\ \ 
$0$,
$ -\xi_{9}^{3}$;\ \ 
$ -c_{9}^{4}$)

Prime. 

\vskip 1ex 
\color{grey}

\noindent(9,3). $4_{\frac{14}{3},5.445}^{9,544}$ \irep{7}:\ \ 
$d_i$ = ($1.0$,
$1.347$,
$-0.532$,
$-1.532$) 

\vskip 0.7ex
\hangindent=3em \hangafter=1
$D^2= 5.445 = 
9-3  c^{1}_{9}
+3c^{2}_{9}
$

\vskip 0.7ex
\hangindent=3em \hangafter=1
$T = ( 0,
\frac{4}{9},
\frac{1}{3},
\frac{2}{3} )
$,

\vskip 0.7ex
\hangindent=3em \hangafter=1
$S$ = ($ 1$,
$ 1+c^{2}_{9}
$,
$ 1-c^{1}_{9}
$,
$ -c_{9}^{1}$;\ \ 
$0$,
$ -1-c^{2}_{9}
$,
$ 1+c^{2}_{9}
$;\ \ 
$ -c_{9}^{1}$,
$ -1$;\ \ 
$ -1+c^{1}_{9}
$)

Prime. 

Not pseudo-unitary. 

\vskip 1ex 
\color{grey}

\noindent(9,4). $4_{\frac{10}{3},5.445}^{9,616}$ \irep{7}:\ \ 
$d_i$ = ($1.0$,
$1.347$,
$-0.532$,
$-1.532$) 

\vskip 0.7ex
\hangindent=3em \hangafter=1
$D^2= 5.445 = 
9-3  c^{1}_{9}
+3c^{2}_{9}
$

\vskip 0.7ex
\hangindent=3em \hangafter=1
$T = ( 0,
\frac{5}{9},
\frac{2}{3},
\frac{1}{3} )
$,

\vskip 0.7ex
\hangindent=3em \hangafter=1
$S$ = ($ 1$,
$ 1+c^{2}_{9}
$,
$ 1-c^{1}_{9}
$,
$ -c_{9}^{1}$;\ \ 
$0$,
$ -1-c^{2}_{9}
$,
$ 1+c^{2}_{9}
$;\ \ 
$ -c_{9}^{1}$,
$ -1$;\ \ 
$ -1+c^{1}_{9}
$)

Prime. 

Not pseudo-unitary. 

\vskip 1ex 
\color{grey}

\noindent(9,5). $4_{\frac{2}{3},2.319}^{9,199}$ \irep{7}:\ \ 
$d_i$ = ($1.0$,
$0.652$,
$-0.347$,
$-0.879$) 

\vskip 0.7ex
\hangindent=3em \hangafter=1
$D^2= 2.319 = 
9-3  c^{1}_{9}
-6  c^{2}_{9}
$

\vskip 0.7ex
\hangindent=3em \hangafter=1
$T = ( 0,
\frac{1}{3},
\frac{2}{3},
\frac{1}{9} )
$,

\vskip 0.7ex
\hangindent=3em \hangafter=1
$S$ = ($ 1$,
$ 1-c^{2}_{9}
$,
$ -c_{9}^{2}$,
$ 1-c^{1}_{9}
-c^{2}_{9}
$;\ \ 
$ -c_{9}^{2}$,
$ -1$,
$ -1+c^{1}_{9}
+c^{2}_{9}
$;\ \ 
$ -1+c^{2}_{9}
$,
$ 1-c^{1}_{9}
-c^{2}_{9}
$;\ \ 
$0$)

Prime. 

Not pseudo-unitary. 

\vskip 1ex 
\color{grey}

\noindent(9,6). $4_{\frac{22}{3},2.319}^{9,549}$ \irep{7}:\ \ 
$d_i$ = ($1.0$,
$0.652$,
$-0.347$,
$-0.879$) 

\vskip 0.7ex
\hangindent=3em \hangafter=1
$D^2= 2.319 = 
9-3  c^{1}_{9}
-6  c^{2}_{9}
$

\vskip 0.7ex
\hangindent=3em \hangafter=1
$T = ( 0,
\frac{2}{3},
\frac{1}{3},
\frac{8}{9} )
$,

\vskip 0.7ex
\hangindent=3em \hangafter=1
$S$ = ($ 1$,
$ 1-c^{2}_{9}
$,
$ -c_{9}^{2}$,
$ 1-c^{1}_{9}
-c^{2}_{9}
$;\ \ 
$ -c_{9}^{2}$,
$ -1$,
$ -1+c^{1}_{9}
+c^{2}_{9}
$;\ \ 
$ -1+c^{2}_{9}
$,
$ 1-c^{1}_{9}
-c^{2}_{9}
$;\ \ 
$0$)

Prime. 

Not pseudo-unitary. 

\vskip 1ex 
\color{blue}

\noindent(10,1). $4_{\frac{19}{5},7.236}^{20,487}$ \irep{8}:\ \ 
$d_i$ = ($1.0$,
$1.618$,
$-1.0$,
$-1.618$) 

\vskip 0.7ex
\hangindent=3em \hangafter=1
$D^2= 7.236 = 
5+\sqrt{5}$

\vskip 0.7ex
\hangindent=3em \hangafter=1
$T = ( 0,
\frac{2}{5},
\frac{1}{4},
\frac{13}{20} )
$,

\vskip 0.7ex
\hangindent=3em \hangafter=1
$S$ = ($ 1$,
$ \frac{1+\sqrt{5}}{2}$,
$ -1$,
$ -\frac{1+\sqrt{5}}{2}$;\ \ 
$ -1$,
$ -\frac{1+\sqrt{5}}{2}$,
$ 1$;\ \ 
$ -1$,
$ -\frac{1+\sqrt{5}}{2}$;\ \ 
$ 1$)

Factors = $2_{\frac{14}{5},3.618}^{5,395}\boxtimes 2_{1,2.}^{4,625} $

Pseudo-unitary $\sim$  
$4_{\frac{9}{5},7.236}^{20,451}$

\vskip 1ex 
\color{grey}

\noindent(10,2). $4_{\frac{9}{5},7.236}^{20,340}$ \irep{8}:\ \ 
$d_i$ = ($1.0$,
$1.618$,
$-1.0$,
$-1.618$) 

\vskip 0.7ex
\hangindent=3em \hangafter=1
$D^2= 7.236 = 
5+\sqrt{5}$

\vskip 0.7ex
\hangindent=3em \hangafter=1
$T = ( 0,
\frac{2}{5},
\frac{3}{4},
\frac{3}{20} )
$,

\vskip 0.7ex
\hangindent=3em \hangafter=1
$S$ = ($ 1$,
$ \frac{1+\sqrt{5}}{2}$,
$ -1$,
$ -\frac{1+\sqrt{5}}{2}$;\ \ 
$ -1$,
$ -\frac{1+\sqrt{5}}{2}$,
$ 1$;\ \ 
$ -1$,
$ -\frac{1+\sqrt{5}}{2}$;\ \ 
$ 1$)

Factors = $2_{\frac{14}{5},3.618}^{5,395}\boxtimes 2_{7,2.}^{4,562} $

Pseudo-unitary $\sim$  
$4_{\frac{19}{5},7.236}^{20,304}$

\vskip 1ex 
\color{grey}

\noindent(10,3). $4_{\frac{31}{5},7.236}^{20,649}$ \irep{8}:\ \ 
$d_i$ = ($1.0$,
$1.618$,
$-1.0$,
$-1.618$) 

\vskip 0.7ex
\hangindent=3em \hangafter=1
$D^2= 7.236 = 
5+\sqrt{5}$

\vskip 0.7ex
\hangindent=3em \hangafter=1
$T = ( 0,
\frac{3}{5},
\frac{1}{4},
\frac{17}{20} )
$,

\vskip 0.7ex
\hangindent=3em \hangafter=1
$S$ = ($ 1$,
$ \frac{1+\sqrt{5}}{2}$,
$ -1$,
$ -\frac{1+\sqrt{5}}{2}$;\ \ 
$ -1$,
$ -\frac{1+\sqrt{5}}{2}$,
$ 1$;\ \ 
$ -1$,
$ -\frac{1+\sqrt{5}}{2}$;\ \ 
$ 1$)

Factors = $2_{\frac{26}{5},3.618}^{5,720}\boxtimes 2_{1,2.}^{4,625} $

Pseudo-unitary $\sim$  
$4_{\frac{21}{5},7.236}^{20,341}$

\vskip 1ex 
\color{grey}

\noindent(10,4). $4_{\frac{21}{5},7.236}^{20,178}$ \irep{8}:\ \ 
$d_i$ = ($1.0$,
$1.618$,
$-1.0$,
$-1.618$) 

\vskip 0.7ex
\hangindent=3em \hangafter=1
$D^2= 7.236 = 
5+\sqrt{5}$

\vskip 0.7ex
\hangindent=3em \hangafter=1
$T = ( 0,
\frac{3}{5},
\frac{3}{4},
\frac{7}{20} )
$,

\vskip 0.7ex
\hangindent=3em \hangafter=1
$S$ = ($ 1$,
$ \frac{1+\sqrt{5}}{2}$,
$ -1$,
$ -\frac{1+\sqrt{5}}{2}$;\ \ 
$ -1$,
$ -\frac{1+\sqrt{5}}{2}$,
$ 1$;\ \ 
$ -1$,
$ -\frac{1+\sqrt{5}}{2}$;\ \ 
$ 1$)

Factors = $2_{\frac{26}{5},3.618}^{5,720}\boxtimes 2_{7,2.}^{4,562} $

Pseudo-unitary $\sim$  
$4_{\frac{31}{5},7.236}^{20,505}$

\vskip 1ex 
\color{grey}

\noindent(10,5). $4_{\frac{3}{5},2.763}^{20,456}$ \irep{8}:\ \ 
$d_i$ = ($1.0$,
$0.618$,
$-0.618$,
$-1.0$) 

\vskip 0.7ex
\hangindent=3em \hangafter=1
$D^2= 2.763 = 
5-\sqrt{5}$

\vskip 0.7ex
\hangindent=3em \hangafter=1
$T = ( 0,
\frac{1}{20},
\frac{4}{5},
\frac{1}{4} )
$,

\vskip 0.7ex
\hangindent=3em \hangafter=1
$S$ = ($ 1$,
$ -\frac{1-\sqrt{5}}{2}$,
$ \frac{1-\sqrt{5}}{2}$,
$ -1$;\ \ 
$ 1$,
$ 1$,
$ -\frac{1-\sqrt{5}}{2}$;\ \ 
$ -1$,
$ -\frac{1-\sqrt{5}}{2}$;\ \ 
$ -1$)

Factors = $2_{\frac{38}{5},1.381}^{5,491}\boxtimes 2_{1,2.}^{4,625} $

Not pseudo-unitary. 

\vskip 1ex 
\color{grey}

\noindent(10,6). $4_{\frac{7}{5},2.763}^{20,379}$ \irep{8}:\ \ 
$d_i$ = ($1.0$,
$0.618$,
$-0.618$,
$-1.0$) 

\vskip 0.7ex
\hangindent=3em \hangafter=1
$D^2= 2.763 = 
5-\sqrt{5}$

\vskip 0.7ex
\hangindent=3em \hangafter=1
$T = ( 0,
\frac{9}{20},
\frac{1}{5},
\frac{1}{4} )
$,

\vskip 0.7ex
\hangindent=3em \hangafter=1
$S$ = ($ 1$,
$ -\frac{1-\sqrt{5}}{2}$,
$ \frac{1-\sqrt{5}}{2}$,
$ -1$;\ \ 
$ 1$,
$ 1$,
$ -\frac{1-\sqrt{5}}{2}$;\ \ 
$ -1$,
$ -\frac{1-\sqrt{5}}{2}$;\ \ 
$ -1$)

Factors = $2_{\frac{2}{5},1.381}^{5,120}\boxtimes 2_{1,2.}^{4,625} $

Not pseudo-unitary. 

\vskip 1ex 
\color{grey}

\noindent(10,7). $4_{\frac{33}{5},2.763}^{20,771}$ \irep{8}:\ \ 
$d_i$ = ($1.0$,
$0.618$,
$-0.618$,
$-1.0$) 

\vskip 0.7ex
\hangindent=3em \hangafter=1
$D^2= 2.763 = 
5-\sqrt{5}$

\vskip 0.7ex
\hangindent=3em \hangafter=1
$T = ( 0,
\frac{11}{20},
\frac{4}{5},
\frac{3}{4} )
$,

\vskip 0.7ex
\hangindent=3em \hangafter=1
$S$ = ($ 1$,
$ -\frac{1-\sqrt{5}}{2}$,
$ \frac{1-\sqrt{5}}{2}$,
$ -1$;\ \ 
$ 1$,
$ 1$,
$ -\frac{1-\sqrt{5}}{2}$;\ \ 
$ -1$,
$ -\frac{1-\sqrt{5}}{2}$;\ \ 
$ -1$)

Factors = $2_{\frac{38}{5},1.381}^{5,491}\boxtimes 2_{7,2.}^{4,562} $

Not pseudo-unitary. 

\vskip 1ex 
\color{grey}

\noindent(10,8). $4_{\frac{37}{5},2.763}^{20,294}$ \irep{8}:\ \ 
$d_i$ = ($1.0$,
$0.618$,
$-0.618$,
$-1.0$) 

\vskip 0.7ex
\hangindent=3em \hangafter=1
$D^2= 2.763 = 
5-\sqrt{5}$

\vskip 0.7ex
\hangindent=3em \hangafter=1
$T = ( 0,
\frac{19}{20},
\frac{1}{5},
\frac{3}{4} )
$,

\vskip 0.7ex
\hangindent=3em \hangafter=1
$S$ = ($ 1$,
$ -\frac{1-\sqrt{5}}{2}$,
$ \frac{1-\sqrt{5}}{2}$,
$ -1$;\ \ 
$ 1$,
$ 1$,
$ -\frac{1-\sqrt{5}}{2}$;\ \ 
$ -1$,
$ -\frac{1-\sqrt{5}}{2}$;\ \ 
$ -1$)

Factors = $2_{\frac{2}{5},1.381}^{5,120}\boxtimes 2_{7,2.}^{4,562} $

Not pseudo-unitary. 

\vskip 1ex 
\color{blue}

\noindent(11,1). $4_{\frac{16}{5},5.}^{5,375}$ \irep{4}:\ \ 
$d_i$ = ($1.0$,
$1.618$,
$-0.618$,
$-1.0$) 

\vskip 0.7ex
\hangindent=3em \hangafter=1
$D^2= 5.0 = 
5$

\vskip 0.7ex
\hangindent=3em \hangafter=1
$T = ( 0,
\frac{2}{5},
\frac{1}{5},
\frac{3}{5} )
$,

\vskip 0.7ex
\hangindent=3em \hangafter=1
$S$ = ($ 1$,
$ \frac{1+\sqrt{5}}{2}$,
$ \frac{1-\sqrt{5}}{2}$,
$ -1$;\ \ 
$ -1$,
$ -1$,
$ -\frac{1-\sqrt{5}}{2}$;\ \ 
$ -1$,
$ -\frac{1+\sqrt{5}}{2}$;\ \ 
$ 1$)

Factors = $2_{\frac{14}{5},3.618}^{5,395}\boxtimes 2_{\frac{2}{5},1.381}^{5,120} $

Not pseudo-unitary. 

\vskip 1ex 
\color{grey}

\noindent(11,2). $4_{\frac{12}{5},5.}^{5,426}$ \irep{4}:\ \ 
$d_i$ = ($1.0$,
$1.618$,
$-0.618$,
$-1.0$) 

\vskip 0.7ex
\hangindent=3em \hangafter=1
$D^2= 5.0 = 
5$

\vskip 0.7ex
\hangindent=3em \hangafter=1
$T = ( 0,
\frac{2}{5},
\frac{4}{5},
\frac{1}{5} )
$,

\vskip 0.7ex
\hangindent=3em \hangafter=1
$S$ = ($ 1$,
$ \frac{1+\sqrt{5}}{2}$,
$ \frac{1-\sqrt{5}}{2}$,
$ -1$;\ \ 
$ -1$,
$ -1$,
$ -\frac{1-\sqrt{5}}{2}$;\ \ 
$ -1$,
$ -\frac{1+\sqrt{5}}{2}$;\ \ 
$ 1$)

Factors = $2_{\frac{14}{5},3.618}^{5,395}\boxtimes 2_{\frac{38}{5},1.381}^{5,491} $

Not pseudo-unitary. 

\vskip 1ex 
\color{grey}

\noindent(11,3). $4_{\frac{28}{5},5.}^{5,332}$ \irep{4}:\ \ 
$d_i$ = ($1.0$,
$1.618$,
$-0.618$,
$-1.0$) 

\vskip 0.7ex
\hangindent=3em \hangafter=1
$D^2= 5.0 = 
5$

\vskip 0.7ex
\hangindent=3em \hangafter=1
$T = ( 0,
\frac{3}{5},
\frac{1}{5},
\frac{4}{5} )
$,

\vskip 0.7ex
\hangindent=3em \hangafter=1
$S$ = ($ 1$,
$ \frac{1+\sqrt{5}}{2}$,
$ \frac{1-\sqrt{5}}{2}$,
$ -1$;\ \ 
$ -1$,
$ -1$,
$ -\frac{1-\sqrt{5}}{2}$;\ \ 
$ -1$,
$ -\frac{1+\sqrt{5}}{2}$;\ \ 
$ 1$)

Factors = $2_{\frac{26}{5},3.618}^{5,720}\boxtimes 2_{\frac{2}{5},1.381}^{5,120} $

Not pseudo-unitary. 

\vskip 1ex 
\color{grey}

\noindent(11,4). $4_{\frac{24}{5},5.}^{5,223}$ \irep{4}:\ \ 
$d_i$ = ($1.0$,
$1.618$,
$-0.618$,
$-1.0$) 

\vskip 0.7ex
\hangindent=3em \hangafter=1
$D^2= 5.0 = 
5$

\vskip 0.7ex
\hangindent=3em \hangafter=1
$T = ( 0,
\frac{3}{5},
\frac{4}{5},
\frac{2}{5} )
$,

\vskip 0.7ex
\hangindent=3em \hangafter=1
$S$ = ($ 1$,
$ \frac{1+\sqrt{5}}{2}$,
$ \frac{1-\sqrt{5}}{2}$,
$ -1$;\ \ 
$ -1$,
$ -1$,
$ -\frac{1-\sqrt{5}}{2}$;\ \ 
$ -1$,
$ -\frac{1+\sqrt{5}}{2}$;\ \ 
$ 1$)

Factors = $2_{\frac{26}{5},3.618}^{5,720}\boxtimes 2_{\frac{38}{5},1.381}^{5,491} $

Not pseudo-unitary. 

\vskip 1ex 
\color{blue}

\noindent(12,1). $4_{0,4.}^{2,250}$ \irep{0}:\ \ 
$d_i$ = ($1.0$,
$1.0$,
$-1.0$,
$-1.0$) 

\vskip 0.7ex
\hangindent=3em \hangafter=1
$D^2= 4.0 = 
4$

\vskip 0.7ex
\hangindent=3em \hangafter=1
$T = ( 0,
0,
0,
\frac{1}{2} )
$,

\vskip 0.7ex
\hangindent=3em \hangafter=1
$S$ = ($ 1$,
$ 1$,
$ -1$,
$ -1$;\ \ 
$ 1$,
$ 1$,
$ 1$;\ \ 
$ 1$,
$ -1$;\ \ 
$ 1$)

Prime. 

Pseudo-unitary $\sim$  
$4_{0,4.}^{2,750}$

\vskip 1ex 
\color{blue}

\noindent(13,1). $4_{0,4.}^{2,250}$ \irep{0}:\ \ 
$d_i$ = ($1.0$,
$1.0$,
$-1.0$,
$-1.0$) 

\vskip 0.7ex
\hangindent=3em \hangafter=1
$D^2= 4.0 = 
4$

\vskip 0.7ex
\hangindent=3em \hangafter=1
$T = ( 0,
\frac{1}{2},
0,
0 )
$,

\vskip 0.7ex
\hangindent=3em \hangafter=1
$S$ = ($ 1$,
$ 1$,
$ -1$,
$ -1$;\ \ 
$ 1$,
$ 1$,
$ 1$;\ \ 
$ 1$,
$ -1$;\ \ 
$ 1$)

Prime. 

Pseudo-unitary $\sim$  
$4_{4,4.}^{2,250}$

\vskip 1ex 
\color{blue}

\noindent(14,1). $4_{4,4.}^{2,250}$ \irep{0}:\ \ 
$d_i$ = ($1.0$,
$1.0$,
$-1.0$,
$-1.0$) 

\vskip 0.7ex
\hangindent=3em \hangafter=1
$D^2= 4.0 = 
4$

\vskip 0.7ex
\hangindent=3em \hangafter=1
$T = ( 0,
\frac{1}{2},
\frac{1}{2},
\frac{1}{2} )
$,

\vskip 0.7ex
\hangindent=3em \hangafter=1
$S$ = ($ 1$,
$ 1$,
$ -1$,
$ -1$;\ \ 
$ 1$,
$ 1$,
$ 1$;\ \ 
$ 1$,
$ -1$;\ \ 
$ 1$)

Prime. 

Pseudo-unitary $\sim$  
$4_{0,4.}^{2,750}$

\vskip 1ex 
\color{blue}

\noindent(15,1). $4_{0,4.}^{4,625}$ \irep{0}:\ \ 
$d_i$ = ($1.0$,
$1.0$,
$-1.0$,
$-1.0$) 

\vskip 0.7ex
\hangindent=3em \hangafter=1
$D^2= 4.0 = 
4$

\vskip 0.7ex
\hangindent=3em \hangafter=1
$T = ( 0,
0,
\frac{1}{4},
\frac{3}{4} )
$,

\vskip 0.7ex
\hangindent=3em \hangafter=1
$S$ = ($ 1$,
$ 1$,
$ -1$,
$ -1$;\ \ 
$ 1$,
$ 1$,
$ 1$;\ \ 
$ -1$,
$ 1$;\ \ 
$ -1$)

Factors = $2_{1,2.}^{4,625}\boxtimes 2_{7,2.}^{4,562} $

Pseudo-unitary $\sim$  
$4_{0,4.}^{4,375}$

\vskip 1ex 
\color{blue}

\noindent(16,1). $4_{0,4.}^{4,500}$ \irep{0}:\ \ 
$d_i$ = ($1.0$,
$1.0$,
$-1.0$,
$-1.0$) 

\vskip 0.7ex
\hangindent=3em \hangafter=1
$D^2= 4.0 = 
4$

\vskip 0.7ex
\hangindent=3em \hangafter=1
$T = ( 0,
\frac{1}{4},
0,
\frac{3}{4} )
$,

\vskip 0.7ex
\hangindent=3em \hangafter=1
$S$ = ($ 1$,
$ 1$,
$ -1$,
$ -1$;\ \ 
$ -1$,
$ 1$,
$ -1$;\ \ 
$ 1$,
$ -1$;\ \ 
$ -1$)

Factors = $2_{1,2.}^{4,437}\boxtimes 2_{7,2.}^{4,562} $

Pseudo-unitary $\sim$  
$4_{2,4.}^{4,625}$

\vskip 1ex 
\color{grey}

\noindent(16,2). $4_{0,4.}^{4,500}$ \irep{0}:\ \ 
$d_i$ = ($1.0$,
$1.0$,
$-1.0$,
$-1.0$) 

\vskip 0.7ex
\hangindent=3em \hangafter=1
$D^2= 4.0 = 
4$

\vskip 0.7ex
\hangindent=3em \hangafter=1
$T = ( 0,
\frac{3}{4},
0,
\frac{1}{4} )
$,

\vskip 0.7ex
\hangindent=3em \hangafter=1
$S$ = ($ 1$,
$ 1$,
$ -1$,
$ -1$;\ \ 
$ -1$,
$ 1$,
$ -1$;\ \ 
$ 1$,
$ -1$;\ \ 
$ -1$)

Factors = $2_{7,2.}^{4,625}\boxtimes 2_{1,2.}^{4,625} $

Pseudo-unitary $\sim$  
$4_{6,4.}^{4,375}$

\vskip 1ex 
\color{blue}

\noindent(17,1). $4_{2,4.}^{4,250}$ \irep{0}:\ \ 
$d_i$ = ($1.0$,
$1.0$,
$-1.0$,
$-1.0$) 

\vskip 0.7ex
\hangindent=3em \hangafter=1
$D^2= 4.0 = 
4$

\vskip 0.7ex
\hangindent=3em \hangafter=1
$T = ( 0,
\frac{1}{4},
\frac{1}{2},
\frac{1}{4} )
$,

\vskip 0.7ex
\hangindent=3em \hangafter=1
$S$ = ($ 1$,
$ 1$,
$ -1$,
$ -1$;\ \ 
$ -1$,
$ 1$,
$ -1$;\ \ 
$ 1$,
$ -1$;\ \ 
$ -1$)

Factors = $2_{1,2.}^{4,437}\boxtimes 2_{1,2.}^{4,625} $

Pseudo-unitary $\sim$  
$4_{0,4.}^{4,375}$

\vskip 1ex 
\color{grey}

\noindent(17,2). $4_{6,4.}^{4,250}$ \irep{0}:\ \ 
$d_i$ = ($1.0$,
$1.0$,
$-1.0$,
$-1.0$) 

\vskip 0.7ex
\hangindent=3em \hangafter=1
$D^2= 4.0 = 
4$

\vskip 0.7ex
\hangindent=3em \hangafter=1
$T = ( 0,
\frac{3}{4},
\frac{1}{2},
\frac{3}{4} )
$,

\vskip 0.7ex
\hangindent=3em \hangafter=1
$S$ = ($ 1$,
$ 1$,
$ -1$,
$ -1$;\ \ 
$ -1$,
$ 1$,
$ -1$;\ \ 
$ 1$,
$ -1$;\ \ 
$ -1$)

Factors = $2_{7,2.}^{4,625}\boxtimes 2_{7,2.}^{4,562} $

Pseudo-unitary $\sim$  
$4_{0,4.}^{4,375}$

\vskip 1ex 
\color{blue}

\noindent(18,1). $4_{2,4.}^{4,125}$ \irep{0}:\ \ 
$d_i$ = ($1.0$,
$1.0$,
$-1.0$,
$-1.0$) 

\vskip 0.7ex
\hangindent=3em \hangafter=1
$D^2= 4.0 = 
4$

\vskip 0.7ex
\hangindent=3em \hangafter=1
$T = ( 0,
\frac{1}{2},
\frac{1}{4},
\frac{1}{4} )
$,

\vskip 0.7ex
\hangindent=3em \hangafter=1
$S$ = ($ 1$,
$ 1$,
$ -1$,
$ -1$;\ \ 
$ 1$,
$ 1$,
$ 1$;\ \ 
$ -1$,
$ 1$;\ \ 
$ -1$)

Factors = $2_{1,2.}^{4,625}\boxtimes 2_{1,2.}^{4,625} $

Pseudo-unitary $\sim$  
$4_{6,4.}^{4,375}$

\vskip 1ex 
\color{grey}

\noindent(18,2). $4_{6,4.}^{4,875}$ \irep{0}:\ \ 
$d_i$ = ($1.0$,
$1.0$,
$-1.0$,
$-1.0$) 

\vskip 0.7ex
\hangindent=3em \hangafter=1
$D^2= 4.0 = 
4$

\vskip 0.7ex
\hangindent=3em \hangafter=1
$T = ( 0,
\frac{1}{2},
\frac{3}{4},
\frac{3}{4} )
$,

\vskip 0.7ex
\hangindent=3em \hangafter=1
$S$ = ($ 1$,
$ 1$,
$ -1$,
$ -1$;\ \ 
$ 1$,
$ 1$,
$ 1$;\ \ 
$ -1$,
$ 1$;\ \ 
$ -1$)

Factors = $2_{7,2.}^{4,562}\boxtimes 2_{7,2.}^{4,562} $

Pseudo-unitary $\sim$  
$4_{2,4.}^{4,625}$

\vskip 1ex 
\color{blue}

\noindent(19,1). $4_{1,4.}^{8,218}$ \irep{6}:\ \ 
$d_i$ = ($1.0$,
$1.0$,
$-1.0$,
$-1.0$) 

\vskip 0.7ex
\hangindent=3em \hangafter=1
$D^2= 4.0 = 
4$

\vskip 0.7ex
\hangindent=3em \hangafter=1
$T = ( 0,
\frac{1}{2},
\frac{1}{8},
\frac{1}{8} )
$,

\vskip 0.7ex
\hangindent=3em \hangafter=1
$S$ = ($ 1$,
$ 1$,
$ -1$,
$ -1$;\ \ 
$ 1$,
$ 1$,
$ 1$;\ \ 
$-\mathrm{i}$,
$\mathrm{i}$;\ \ 
$-\mathrm{i}$)

Prime. 

Pseudo-unitary $\sim$  
$4_{5,4.}^{8,312}$

\vskip 1ex 
\color{grey}

\noindent(19,2). $4_{3,4.}^{8,312}$ \irep{6}:\ \ 
$d_i$ = ($1.0$,
$1.0$,
$-1.0$,
$-1.0$) 

\vskip 0.7ex
\hangindent=3em \hangafter=1
$D^2= 4.0 = 
4$

\vskip 0.7ex
\hangindent=3em \hangafter=1
$T = ( 0,
\frac{1}{2},
\frac{3}{8},
\frac{3}{8} )
$,

\vskip 0.7ex
\hangindent=3em \hangafter=1
$S$ = ($ 1$,
$ 1$,
$ -1$,
$ -1$;\ \ 
$ 1$,
$ 1$,
$ 1$;\ \ 
$\mathrm{i}$,
$-\mathrm{i}$;\ \ 
$\mathrm{i}$)

Prime. 

Pseudo-unitary $\sim$  
$4_{7,4.}^{8,781}$

\vskip 1ex 
\color{grey}

\noindent(19,3). $4_{5,4.}^{8,531}$ \irep{6}:\ \ 
$d_i$ = ($1.0$,
$1.0$,
$-1.0$,
$-1.0$) 

\vskip 0.7ex
\hangindent=3em \hangafter=1
$D^2= 4.0 = 
4$

\vskip 0.7ex
\hangindent=3em \hangafter=1
$T = ( 0,
\frac{1}{2},
\frac{5}{8},
\frac{5}{8} )
$,

\vskip 0.7ex
\hangindent=3em \hangafter=1
$S$ = ($ 1$,
$ 1$,
$ -1$,
$ -1$;\ \ 
$ 1$,
$ 1$,
$ 1$;\ \ 
$-\mathrm{i}$,
$\mathrm{i}$;\ \ 
$-\mathrm{i}$)

Prime. 

Pseudo-unitary $\sim$  
$4_{1,4.}^{8,718}$

\vskip 1ex 
\color{grey}

\noindent(19,4). $4_{7,4.}^{8,128}$ \irep{6}:\ \ 
$d_i$ = ($1.0$,
$1.0$,
$-1.0$,
$-1.0$) 

\vskip 0.7ex
\hangindent=3em \hangafter=1
$D^2= 4.0 = 
4$

\vskip 0.7ex
\hangindent=3em \hangafter=1
$T = ( 0,
\frac{1}{2},
\frac{7}{8},
\frac{7}{8} )
$,

\vskip 0.7ex
\hangindent=3em \hangafter=1
$S$ = ($ 1$,
$ 1$,
$ -1$,
$ -1$;\ \ 
$ 1$,
$ 1$,
$ 1$;\ \ 
$\mathrm{i}$,
$-\mathrm{i}$;\ \ 
$\mathrm{i}$)

Prime. 

Pseudo-unitary $\sim$  
$4_{3,4.}^{8,468}$

\vskip 1ex 

}

\subsection{Rank 5}

{\small
\black

\noindent(1,1). $5_{0,5.}^{5,110}$ \irep{9}:\ \ 
$d_i$ = ($1.0$,
$1.0$,
$1.0$,
$1.0$,
$1.0$) 

\vskip 0.7ex
\hangindent=3em \hangafter=1
$D^2= 5.0 = 
5$

\vskip 0.7ex
\hangindent=3em \hangafter=1
$T = ( 0,
\frac{1}{5},
\frac{1}{5},
\frac{4}{5},
\frac{4}{5} )
$,

\vskip 0.7ex
\hangindent=3em \hangafter=1
$S$ = ($ 1$,
$ 1$,
$ 1$,
$ 1$,
$ 1$;\ \ 
$ -\zeta_{10}^{1}$,
$ \zeta_{5}^{2}$,
$ -\zeta_{10}^{3}$,
$ \zeta_{5}^{1}$;\ \ 
$ -\zeta_{10}^{1}$,
$ \zeta_{5}^{1}$,
$ -\zeta_{10}^{3}$;\ \ 
$ \zeta_{5}^{2}$,
$ -\zeta_{10}^{1}$;\ \ 
$ \zeta_{5}^{2}$)

Prime. 

\vskip 1ex 
\color{grey}

\noindent(1,2). $5_{4,5.}^{5,210}$ \irep{9}:\ \ 
$d_i$ = ($1.0$,
$1.0$,
$1.0$,
$1.0$,
$1.0$) 

\vskip 0.7ex
\hangindent=3em \hangafter=1
$D^2= 5.0 = 
5$

\vskip 0.7ex
\hangindent=3em \hangafter=1
$T = ( 0,
\frac{2}{5},
\frac{2}{5},
\frac{3}{5},
\frac{3}{5} )
$,

\vskip 0.7ex
\hangindent=3em \hangafter=1
$S$ = ($ 1$,
$ 1$,
$ 1$,
$ 1$,
$ 1$;\ \ 
$ \zeta_{5}^{1}$,
$ -\zeta_{10}^{3}$,
$ -\zeta_{10}^{1}$,
$ \zeta_{5}^{2}$;\ \ 
$ \zeta_{5}^{1}$,
$ \zeta_{5}^{2}$,
$ -\zeta_{10}^{1}$;\ \ 
$ -\zeta_{10}^{3}$,
$ \zeta_{5}^{1}$;\ \ 
$ -\zeta_{10}^{3}$)

Prime. 

\vskip 1ex 
\black

\noindent(2,1). $5_{2,12.}^{24,940}$ \irep{17}:\ \ 
$d_i$ = ($1.0$,
$1.0$,
$1.732$,
$1.732$,
$2.0$) 

\vskip 0.7ex
\hangindent=3em \hangafter=1
$D^2= 12.0 = 
12$

\vskip 0.7ex
\hangindent=3em \hangafter=1
$T = ( 0,
0,
\frac{1}{8},
\frac{5}{8},
\frac{1}{3} )
$,

\vskip 0.7ex
\hangindent=3em \hangafter=1
$S$ = ($ 1$,
$ 1$,
$ \sqrt{3}$,
$ \sqrt{3}$,
$ 2$;\ \ 
$ 1$,
$ -\sqrt{3}$,
$ -\sqrt{3}$,
$ 2$;\ \ 
$ \sqrt{3}$,
$ -\sqrt{3}$,
$0$;\ \ 
$ \sqrt{3}$,
$0$;\ \ 
$ -2$)

Prime. 

\vskip 1ex 
\color{grey}

\noindent(2,2). $5_{6,12.}^{24,592}$ \irep{17}:\ \ 
$d_i$ = ($1.0$,
$1.0$,
$1.732$,
$1.732$,
$2.0$) 

\vskip 0.7ex
\hangindent=3em \hangafter=1
$D^2= 12.0 = 
12$

\vskip 0.7ex
\hangindent=3em \hangafter=1
$T = ( 0,
0,
\frac{3}{8},
\frac{7}{8},
\frac{2}{3} )
$,

\vskip 0.7ex
\hangindent=3em \hangafter=1
$S$ = ($ 1$,
$ 1$,
$ \sqrt{3}$,
$ \sqrt{3}$,
$ 2$;\ \ 
$ 1$,
$ -\sqrt{3}$,
$ -\sqrt{3}$,
$ 2$;\ \ 
$ \sqrt{3}$,
$ -\sqrt{3}$,
$0$;\ \ 
$ \sqrt{3}$,
$0$;\ \ 
$ -2$)

Prime. 

\vskip 1ex 
\color{grey}

\noindent(2,3). $5_{2,12.}^{24,148}$ \irep{17}:\ \ 
$d_i$ = ($1.0$,
$1.0$,
$2.0$,
$-1.732$,
$-1.732$) 

\vskip 0.7ex
\hangindent=3em \hangafter=1
$D^2= 12.0 = 
12$

\vskip 0.7ex
\hangindent=3em \hangafter=1
$T = ( 0,
0,
\frac{1}{3},
\frac{3}{8},
\frac{7}{8} )
$,

\vskip 0.7ex
\hangindent=3em \hangafter=1
$S$ = ($ 1$,
$ 1$,
$ 2$,
$ -\sqrt{3}$,
$ -\sqrt{3}$;\ \ 
$ 1$,
$ 2$,
$ \sqrt{3}$,
$ \sqrt{3}$;\ \ 
$ -2$,
$0$,
$0$;\ \ 
$ -\sqrt{3}$,
$ \sqrt{3}$;\ \ 
$ -\sqrt{3}$)

Prime. 

Pseudo-unitary $\sim$  
$5_{2,12.}^{24,741}$

\vskip 1ex 
\color{grey}

\noindent(2,4). $5_{6,12.}^{24,512}$ \irep{17}:\ \ 
$d_i$ = ($1.0$,
$1.0$,
$2.0$,
$-1.732$,
$-1.732$) 

\vskip 0.7ex
\hangindent=3em \hangafter=1
$D^2= 12.0 = 
12$

\vskip 0.7ex
\hangindent=3em \hangafter=1
$T = ( 0,
0,
\frac{2}{3},
\frac{1}{8},
\frac{5}{8} )
$,

\vskip 0.7ex
\hangindent=3em \hangafter=1
$S$ = ($ 1$,
$ 1$,
$ 2$,
$ -\sqrt{3}$,
$ -\sqrt{3}$;\ \ 
$ 1$,
$ 2$,
$ \sqrt{3}$,
$ \sqrt{3}$;\ \ 
$ -2$,
$0$,
$0$;\ \ 
$ -\sqrt{3}$,
$ \sqrt{3}$;\ \ 
$ -\sqrt{3}$)

Prime. 

Pseudo-unitary $\sim$  
$5_{6,12.}^{24,273}$

\vskip 1ex 
\black

\noindent(3,1). $5_{6,12.}^{24,273}$ \irep{17}:\ \ 
$d_i$ = ($1.0$,
$1.0$,
$1.732$,
$1.732$,
$2.0$) 

\vskip 0.7ex
\hangindent=3em \hangafter=1
$D^2= 12.0 = 
12$

\vskip 0.7ex
\hangindent=3em \hangafter=1
$T = ( 0,
0,
\frac{1}{8},
\frac{5}{8},
\frac{2}{3} )
$,

\vskip 0.7ex
\hangindent=3em \hangafter=1
$S$ = ($ 1$,
$ 1$,
$ \sqrt{3}$,
$ \sqrt{3}$,
$ 2$;\ \ 
$ 1$,
$ -\sqrt{3}$,
$ -\sqrt{3}$,
$ 2$;\ \ 
$ -\sqrt{3}$,
$ \sqrt{3}$,
$0$;\ \ 
$ -\sqrt{3}$,
$0$;\ \ 
$ -2$)

Prime. 

\vskip 1ex 
\color{grey}

\noindent(3,2). $5_{2,12.}^{24,741}$ \irep{17}:\ \ 
$d_i$ = ($1.0$,
$1.0$,
$1.732$,
$1.732$,
$2.0$) 

\vskip 0.7ex
\hangindent=3em \hangafter=1
$D^2= 12.0 = 
12$

\vskip 0.7ex
\hangindent=3em \hangafter=1
$T = ( 0,
0,
\frac{3}{8},
\frac{7}{8},
\frac{1}{3} )
$,

\vskip 0.7ex
\hangindent=3em \hangafter=1
$S$ = ($ 1$,
$ 1$,
$ \sqrt{3}$,
$ \sqrt{3}$,
$ 2$;\ \ 
$ 1$,
$ -\sqrt{3}$,
$ -\sqrt{3}$,
$ 2$;\ \ 
$ -\sqrt{3}$,
$ \sqrt{3}$,
$0$;\ \ 
$ -\sqrt{3}$,
$0$;\ \ 
$ -2$)

Prime. 

\vskip 1ex 
\color{grey}

\noindent(3,3). $5_{2,12.}^{24,615}$ \irep{17}:\ \ 
$d_i$ = ($1.0$,
$1.0$,
$2.0$,
$-1.732$,
$-1.732$) 

\vskip 0.7ex
\hangindent=3em \hangafter=1
$D^2= 12.0 = 
12$

\vskip 0.7ex
\hangindent=3em \hangafter=1
$T = ( 0,
0,
\frac{1}{3},
\frac{1}{8},
\frac{5}{8} )
$,

\vskip 0.7ex
\hangindent=3em \hangafter=1
$S$ = ($ 1$,
$ 1$,
$ 2$,
$ -\sqrt{3}$,
$ -\sqrt{3}$;\ \ 
$ 1$,
$ 2$,
$ \sqrt{3}$,
$ \sqrt{3}$;\ \ 
$ -2$,
$0$,
$0$;\ \ 
$ \sqrt{3}$,
$ -\sqrt{3}$;\ \ 
$ \sqrt{3}$)

Prime. 

Pseudo-unitary $\sim$  
$5_{2,12.}^{24,940}$

\vskip 1ex 
\color{grey}

\noindent(3,4). $5_{6,12.}^{24,814}$ \irep{17}:\ \ 
$d_i$ = ($1.0$,
$1.0$,
$2.0$,
$-1.732$,
$-1.732$) 

\vskip 0.7ex
\hangindent=3em \hangafter=1
$D^2= 12.0 = 
12$

\vskip 0.7ex
\hangindent=3em \hangafter=1
$T = ( 0,
0,
\frac{2}{3},
\frac{3}{8},
\frac{7}{8} )
$,

\vskip 0.7ex
\hangindent=3em \hangafter=1
$S$ = ($ 1$,
$ 1$,
$ 2$,
$ -\sqrt{3}$,
$ -\sqrt{3}$;\ \ 
$ 1$,
$ 2$,
$ \sqrt{3}$,
$ \sqrt{3}$;\ \ 
$ -2$,
$0$,
$0$;\ \ 
$ \sqrt{3}$,
$ -\sqrt{3}$;\ \ 
$ \sqrt{3}$)

Prime. 

Pseudo-unitary $\sim$  
$5_{6,12.}^{24,592}$

\vskip 1ex 
\black

\noindent(4,1). $5_{\frac{72}{11},34.64}^{11,216}$ \irep{15}:\ \ 
$d_i$ = ($1.0$,
$1.918$,
$2.682$,
$3.228$,
$3.513$) 

\vskip 0.7ex
\hangindent=3em \hangafter=1
$D^2= 34.646 = 
15+10c^{1}_{11}
+6c^{2}_{11}
+3c^{3}_{11}
+c^{4}_{11}
$

\vskip 0.7ex
\hangindent=3em \hangafter=1
$T = ( 0,
\frac{2}{11},
\frac{9}{11},
\frac{10}{11},
\frac{5}{11} )
$,

\vskip 0.7ex
\hangindent=3em \hangafter=1
$S$ = ($ 1$,
$ -c_{11}^{5}$,
$ \xi_{11}^{3}$,
$ \xi_{11}^{7}$,
$ \xi_{11}^{5}$;\ \ 
$ -\xi_{11}^{7}$,
$ \xi_{11}^{5}$,
$ -\xi_{11}^{3}$,
$ 1$;\ \ 
$ -c_{11}^{5}$,
$ -1$,
$ -\xi_{11}^{7}$;\ \ 
$ \xi_{11}^{5}$,
$ c_{11}^{5}$;\ \ 
$ \xi_{11}^{3}$)

Prime. 

\vskip 1ex 
\color{grey}

\noindent(4,2). $5_{\frac{16}{11},34.64}^{11,640}$ \irep{15}:\ \ 
$d_i$ = ($1.0$,
$1.918$,
$2.682$,
$3.228$,
$3.513$) 

\vskip 0.7ex
\hangindent=3em \hangafter=1
$D^2= 34.646 = 
15+10c^{1}_{11}
+6c^{2}_{11}
+3c^{3}_{11}
+c^{4}_{11}
$

\vskip 0.7ex
\hangindent=3em \hangafter=1
$T = ( 0,
\frac{9}{11},
\frac{2}{11},
\frac{1}{11},
\frac{6}{11} )
$,

\vskip 0.7ex
\hangindent=3em \hangafter=1
$S$ = ($ 1$,
$ -c_{11}^{5}$,
$ \xi_{11}^{3}$,
$ \xi_{11}^{7}$,
$ \xi_{11}^{5}$;\ \ 
$ -\xi_{11}^{7}$,
$ \xi_{11}^{5}$,
$ -\xi_{11}^{3}$,
$ 1$;\ \ 
$ -c_{11}^{5}$,
$ -1$,
$ -\xi_{11}^{7}$;\ \ 
$ \xi_{11}^{5}$,
$ c_{11}^{5}$;\ \ 
$ \xi_{11}^{3}$)

Prime. 

\vskip 1ex 
\color{grey}

\noindent(4,3). $5_{\frac{32}{11},9.408}^{11,549}$ \irep{15}:\ \ 
$d_i$ = ($1.0$,
$0.521$,
$1.830$,
$-1.397$,
$-1.682$) 

\vskip 0.7ex
\hangindent=3em \hangafter=1
$D^2= 9.408 = 
12-3  c^{1}_{11}
+7c^{2}_{11}
-2  c^{3}_{11}
+3c^{4}_{11}
$

\vskip 0.7ex
\hangindent=3em \hangafter=1
$T = ( 0,
\frac{1}{11},
\frac{4}{11},
\frac{2}{11},
\frac{7}{11} )
$,

\vskip 0.7ex
\hangindent=3em \hangafter=1
$S$ = ($ 1$,
$ 1+c^{2}_{11}
+c^{4}_{11}
$,
$ 1+c^{2}_{11}
$,
$ -c^{1}_{11}
-c^{3}_{11}
$,
$ -c_{11}^{1}$;\ \ 
$ 1+c^{2}_{11}
$,
$ c^{1}_{11}
+c^{3}_{11}
$,
$ c_{11}^{1}$,
$ 1$;\ \ 
$ -c_{11}^{1}$,
$ -1$,
$ 1+c^{2}_{11}
+c^{4}_{11}
$;\ \ 
$ 1+c^{2}_{11}
+c^{4}_{11}
$,
$ -1-c^{2}_{11}
$;\ \ 
$ c^{1}_{11}
+c^{3}_{11}
$)

Prime. 

Not pseudo-unitary. 

\vskip 1ex 
\color{grey}

\noindent(4,4). $5_{\frac{56}{11},9.408}^{11,540}$ \irep{15}:\ \ 
$d_i$ = ($1.0$,
$0.521$,
$1.830$,
$-1.397$,
$-1.682$) 

\vskip 0.7ex
\hangindent=3em \hangafter=1
$D^2= 9.408 = 
12-3  c^{1}_{11}
+7c^{2}_{11}
-2  c^{3}_{11}
+3c^{4}_{11}
$

\vskip 0.7ex
\hangindent=3em \hangafter=1
$T = ( 0,
\frac{10}{11},
\frac{7}{11},
\frac{9}{11},
\frac{4}{11} )
$,

\vskip 0.7ex
\hangindent=3em \hangafter=1
$S$ = ($ 1$,
$ 1+c^{2}_{11}
+c^{4}_{11}
$,
$ 1+c^{2}_{11}
$,
$ -c^{1}_{11}
-c^{3}_{11}
$,
$ -c_{11}^{1}$;\ \ 
$ 1+c^{2}_{11}
$,
$ c^{1}_{11}
+c^{3}_{11}
$,
$ c_{11}^{1}$,
$ 1$;\ \ 
$ -c_{11}^{1}$,
$ -1$,
$ 1+c^{2}_{11}
+c^{4}_{11}
$;\ \ 
$ 1+c^{2}_{11}
+c^{4}_{11}
$,
$ -1-c^{2}_{11}
$;\ \ 
$ c^{1}_{11}
+c^{3}_{11}
$)

Prime. 

Not pseudo-unitary. 

\vskip 1ex 
\color{grey}

\noindent(4,5). $5_{\frac{40}{11},4.814}^{11,181}$ \irep{15}:\ \ 
$d_i$ = ($1.0$,
$0.715$,
$1.309$,
$-0.372$,
$-1.203$) 

\vskip 0.7ex
\hangindent=3em \hangafter=1
$D^2= 4.814 = 
9-5  c^{1}_{11}
-3  c^{2}_{11}
+4c^{3}_{11}
-6  c^{4}_{11}
$

\vskip 0.7ex
\hangindent=3em \hangafter=1
$T = ( 0,
\frac{5}{11},
\frac{6}{11},
\frac{8}{11},
\frac{4}{11} )
$,

\vskip 0.7ex
\hangindent=3em \hangafter=1
$S$ = ($ 1$,
$ 1+c^{3}_{11}
$,
$ -c_{11}^{4}$,
$ -c^{1}_{11}
-c^{4}_{11}
$,
$ -c^{1}_{11}
-c^{2}_{11}
-c^{4}_{11}
$;\ \ 
$ -c_{11}^{4}$,
$ -c^{1}_{11}
-c^{2}_{11}
-c^{4}_{11}
$,
$ -1$,
$ c^{1}_{11}
+c^{4}_{11}
$;\ \ 
$ c^{1}_{11}
+c^{4}_{11}
$,
$ -1-c^{3}_{11}
$,
$ 1$;\ \ 
$ -c^{1}_{11}
-c^{2}_{11}
-c^{4}_{11}
$,
$ c_{11}^{4}$;\ \ 
$ 1+c^{3}_{11}
$)

Prime. 

Not pseudo-unitary. 

\vskip 1ex 
\color{grey}

\noindent(4,6). $5_{\frac{48}{11},4.814}^{11,393}$ \irep{15}:\ \ 
$d_i$ = ($1.0$,
$0.715$,
$1.309$,
$-0.372$,
$-1.203$) 

\vskip 0.7ex
\hangindent=3em \hangafter=1
$D^2= 4.814 = 
9-5  c^{1}_{11}
-3  c^{2}_{11}
+4c^{3}_{11}
-6  c^{4}_{11}
$

\vskip 0.7ex
\hangindent=3em \hangafter=1
$T = ( 0,
\frac{6}{11},
\frac{5}{11},
\frac{3}{11},
\frac{7}{11} )
$,

\vskip 0.7ex
\hangindent=3em \hangafter=1
$S$ = ($ 1$,
$ 1+c^{3}_{11}
$,
$ -c_{11}^{4}$,
$ -c^{1}_{11}
-c^{4}_{11}
$,
$ -c^{1}_{11}
-c^{2}_{11}
-c^{4}_{11}
$;\ \ 
$ -c_{11}^{4}$,
$ -c^{1}_{11}
-c^{2}_{11}
-c^{4}_{11}
$,
$ -1$,
$ c^{1}_{11}
+c^{4}_{11}
$;\ \ 
$ c^{1}_{11}
+c^{4}_{11}
$,
$ -1-c^{3}_{11}
$,
$ 1$;\ \ 
$ -c^{1}_{11}
-c^{2}_{11}
-c^{4}_{11}
$,
$ c_{11}^{4}$;\ \ 
$ 1+c^{3}_{11}
$)

Prime. 

Not pseudo-unitary. 

\vskip 1ex 
\color{grey}

\noindent(4,7). $5_{\frac{20}{11},3.323}^{11,189}$ \irep{15}:\ \ 
$d_i$ = ($1.0$,
$1.88$,
$-0.309$,
$-0.594$,
$-0.830$) 

\vskip 0.7ex
\hangindent=3em \hangafter=1
$D^2= 3.323 = 
14+2c^{1}_{11}
-c^{2}_{11}
+5c^{3}_{11}
+9c^{4}_{11}
$

\vskip 0.7ex
\hangindent=3em \hangafter=1
$T = ( 0,
\frac{4}{11},
\frac{8}{11},
\frac{2}{11},
\frac{3}{11} )
$,

\vskip 0.7ex
\hangindent=3em \hangafter=1
$S$ = ($ 1$,
$ 1+c^{1}_{11}
+c^{3}_{11}
+c^{4}_{11}
$,
$ 1+c^{4}_{11}
$,
$ 1+c^{3}_{11}
+c^{4}_{11}
$,
$ -c_{11}^{2}$;\ \ 
$ 1+c^{3}_{11}
+c^{4}_{11}
$,
$ -1$,
$ c_{11}^{2}$,
$ -1-c^{4}_{11}
$;\ \ 
$ -c_{11}^{2}$,
$ -1-c^{1}_{11}
-c^{3}_{11}
-c^{4}_{11}
$,
$ 1+c^{3}_{11}
+c^{4}_{11}
$;\ \ 
$ 1+c^{4}_{11}
$,
$ 1$;\ \ 
$ -1-c^{1}_{11}
-c^{3}_{11}
-c^{4}_{11}
$)

Prime. 

Not pseudo-unitary. 

\vskip 1ex 
\color{grey}

\noindent(4,8). $5_{\frac{68}{11},3.323}^{11,508}$ \irep{15}:\ \ 
$d_i$ = ($1.0$,
$1.88$,
$-0.309$,
$-0.594$,
$-0.830$) 

\vskip 0.7ex
\hangindent=3em \hangafter=1
$D^2= 3.323 = 
14+2c^{1}_{11}
-c^{2}_{11}
+5c^{3}_{11}
+9c^{4}_{11}
$

\vskip 0.7ex
\hangindent=3em \hangafter=1
$T = ( 0,
\frac{7}{11},
\frac{3}{11},
\frac{9}{11},
\frac{8}{11} )
$,

\vskip 0.7ex
\hangindent=3em \hangafter=1
$S$ = ($ 1$,
$ 1+c^{1}_{11}
+c^{3}_{11}
+c^{4}_{11}
$,
$ 1+c^{4}_{11}
$,
$ 1+c^{3}_{11}
+c^{4}_{11}
$,
$ -c_{11}^{2}$;\ \ 
$ 1+c^{3}_{11}
+c^{4}_{11}
$,
$ -1$,
$ c_{11}^{2}$,
$ -1-c^{4}_{11}
$;\ \ 
$ -c_{11}^{2}$,
$ -1-c^{1}_{11}
-c^{3}_{11}
-c^{4}_{11}
$,
$ 1+c^{3}_{11}
+c^{4}_{11}
$;\ \ 
$ 1+c^{4}_{11}
$,
$ 1$;\ \ 
$ -1-c^{1}_{11}
-c^{3}_{11}
-c^{4}_{11}
$)

Prime. 

Not pseudo-unitary. 

\vskip 1ex 
\color{grey}

\noindent(4,9). $5_{\frac{80}{11},2.806}^{11,611}$ \irep{15}:\ \ 
$d_i$ = ($1.0$,
$0.284$,
$0.763$,
$-0.546$,
$-0.918$) 

\vskip 0.7ex
\hangindent=3em \hangafter=1
$D^2= 2.806 = 
5-4  c^{1}_{11}
-9  c^{2}_{11}
-10  c^{3}_{11}
-7  c^{4}_{11}
$

\vskip 0.7ex
\hangindent=3em \hangafter=1
$T = ( 0,
\frac{1}{11},
\frac{8}{11},
\frac{5}{11},
\frac{10}{11} )
$,

\vskip 0.7ex
\hangindent=3em \hangafter=1
$S$ = ($ 1$,
$ -c_{11}^{3}$,
$ -c^{2}_{11}
-c^{3}_{11}
-c^{4}_{11}
$,
$ -c^{2}_{11}
-c^{3}_{11}
$,
$ -c^{1}_{11}
-c^{2}_{11}
-c^{3}_{11}
-c^{4}_{11}
$;\ \ 
$ c^{2}_{11}
+c^{3}_{11}
$,
$ 1$,
$ c^{1}_{11}
+c^{2}_{11}
+c^{3}_{11}
+c^{4}_{11}
$,
$ -c^{2}_{11}
-c^{3}_{11}
-c^{4}_{11}
$;\ \ 
$ -c^{1}_{11}
-c^{2}_{11}
-c^{3}_{11}
-c^{4}_{11}
$,
$ c_{11}^{3}$,
$ c^{2}_{11}
+c^{3}_{11}
$;\ \ 
$ -c^{2}_{11}
-c^{3}_{11}
-c^{4}_{11}
$,
$ -1$;\ \ 
$ -c_{11}^{3}$)

Prime. 

Not pseudo-unitary. 

\vskip 1ex 
\color{grey}

\noindent(4,10). $5_{\frac{8}{11},2.806}^{11,238}$ \irep{15}:\ \ 
$d_i$ = ($1.0$,
$0.284$,
$0.763$,
$-0.546$,
$-0.918$) 

\vskip 0.7ex
\hangindent=3em \hangafter=1
$D^2= 2.806 = 
5-4  c^{1}_{11}
-9  c^{2}_{11}
-10  c^{3}_{11}
-7  c^{4}_{11}
$

\vskip 0.7ex
\hangindent=3em \hangafter=1
$T = ( 0,
\frac{10}{11},
\frac{3}{11},
\frac{6}{11},
\frac{1}{11} )
$,

\vskip 0.7ex
\hangindent=3em \hangafter=1
$S$ = ($ 1$,
$ -c_{11}^{3}$,
$ -c^{2}_{11}
-c^{3}_{11}
-c^{4}_{11}
$,
$ -c^{2}_{11}
-c^{3}_{11}
$,
$ -c^{1}_{11}
-c^{2}_{11}
-c^{3}_{11}
-c^{4}_{11}
$;\ \ 
$ c^{2}_{11}
+c^{3}_{11}
$,
$ 1$,
$ c^{1}_{11}
+c^{2}_{11}
+c^{3}_{11}
+c^{4}_{11}
$,
$ -c^{2}_{11}
-c^{3}_{11}
-c^{4}_{11}
$;\ \ 
$ -c^{1}_{11}
-c^{2}_{11}
-c^{3}_{11}
-c^{4}_{11}
$,
$ c_{11}^{3}$,
$ c^{2}_{11}
+c^{3}_{11}
$;\ \ 
$ -c^{2}_{11}
-c^{3}_{11}
-c^{4}_{11}
$,
$ -1$;\ \ 
$ -c_{11}^{3}$)

Prime. 

Not pseudo-unitary. 

\vskip 1ex 
\black

\noindent(5,1). $5_{\frac{38}{7},35.34}^{7,386}$ \irep{11}:\ \ 
$d_i$ = ($1.0$,
$2.246$,
$2.246$,
$2.801$,
$4.48$) 

\vskip 0.7ex
\hangindent=3em \hangafter=1
$D^2= 35.342 = 
21+14c^{1}_{7}
+7c^{2}_{7}
$

\vskip 0.7ex
\hangindent=3em \hangafter=1
$T = ( 0,
\frac{1}{7},
\frac{1}{7},
\frac{6}{7},
\frac{4}{7} )
$,

\vskip 0.7ex
\hangindent=3em \hangafter=1
$S$ = ($ 1$,
$ \xi_{7}^{3}$,
$ \xi_{7}^{3}$,
$ 2+c^{1}_{7}
+c^{2}_{7}
$,
$ 2+2c^{1}_{7}
+c^{2}_{7}
$;\ \ 
$ s^{1}_{7}
+\zeta^{2}_{7}
+\zeta^{3}_{7}
$,
$ -1-2  \zeta^{1}_{7}
-\zeta^{2}_{7}
-\zeta^{3}_{7}
$,
$ -\xi_{7}^{3}$,
$ \xi_{7}^{3}$;\ \ 
$ s^{1}_{7}
+\zeta^{2}_{7}
+\zeta^{3}_{7}
$,
$ -\xi_{7}^{3}$,
$ \xi_{7}^{3}$;\ \ 
$ 2+2c^{1}_{7}
+c^{2}_{7}
$,
$ -1$;\ \ 
$ -2-c^{1}_{7}
-c^{2}_{7}
$)

Prime. 

\vskip 1ex 
\color{grey}

\noindent(5,2). $5_{\frac{18}{7},35.34}^{7,101}$ \irep{11}:\ \ 
$d_i$ = ($1.0$,
$2.246$,
$2.246$,
$2.801$,
$4.48$) 

\vskip 0.7ex
\hangindent=3em \hangafter=1
$D^2= 35.342 = 
21+14c^{1}_{7}
+7c^{2}_{7}
$

\vskip 0.7ex
\hangindent=3em \hangafter=1
$T = ( 0,
\frac{6}{7},
\frac{6}{7},
\frac{1}{7},
\frac{3}{7} )
$,

\vskip 0.7ex
\hangindent=3em \hangafter=1
$S$ = ($ 1$,
$ \xi_{7}^{3}$,
$ \xi_{7}^{3}$,
$ 2+c^{1}_{7}
+c^{2}_{7}
$,
$ 2+2c^{1}_{7}
+c^{2}_{7}
$;\ \ 
$ -1-2  \zeta^{1}_{7}
-\zeta^{2}_{7}
-\zeta^{3}_{7}
$,
$ s^{1}_{7}
+\zeta^{2}_{7}
+\zeta^{3}_{7}
$,
$ -\xi_{7}^{3}$,
$ \xi_{7}^{3}$;\ \ 
$ -1-2  \zeta^{1}_{7}
-\zeta^{2}_{7}
-\zeta^{3}_{7}
$,
$ -\xi_{7}^{3}$,
$ \xi_{7}^{3}$;\ \ 
$ 2+2c^{1}_{7}
+c^{2}_{7}
$,
$ -1$;\ \ 
$ -2-c^{1}_{7}
-c^{2}_{7}
$)

Prime. 

\vskip 1ex 
\color{grey}

\noindent(5,3). $5_{\frac{26}{7},4.501}^{7,408}$ \irep{11}:\ \ 
$d_i$ = ($1.0$,
$1.445$,
$-0.356$,
$-0.801$,
$-0.801$) 

\vskip 0.7ex
\hangindent=3em \hangafter=1
$D^2= 4.501 = 
7-7  c^{1}_{7}
-14  c^{2}_{7}
$

\vskip 0.7ex
\hangindent=3em \hangafter=1
$T = ( 0,
\frac{3}{7},
\frac{2}{7},
\frac{4}{7},
\frac{4}{7} )
$,

\vskip 0.7ex
\hangindent=3em \hangafter=1
$S$ = ($ 1$,
$ 1-c^{2}_{7}
$,
$ -c^{1}_{7}
-2  c^{2}_{7}
$,
$ -c^{1}_{7}
-c^{2}_{7}
$,
$ -c^{1}_{7}
-c^{2}_{7}
$;\ \ 
$ -c^{1}_{7}
-2  c^{2}_{7}
$,
$ -1$,
$ c^{1}_{7}
+c^{2}_{7}
$,
$ c^{1}_{7}
+c^{2}_{7}
$;\ \ 
$ -1+c^{2}_{7}
$,
$ -c^{1}_{7}
-c^{2}_{7}
$,
$ -c^{1}_{7}
-c^{2}_{7}
$;\ \ 
$ -1-\zeta^{-1}_{7}
-\zeta^{2}_{7}
-2  \zeta^{3}_{7}
$,
$ 1+\zeta^{1}_{7}
+2\zeta^{-1}_{7}
+2\zeta^{2}_{7}
+\zeta^{-2}_{7}
+2\zeta^{3}_{7}
$;\ \ 
$ -1-\zeta^{-1}_{7}
-\zeta^{2}_{7}
-2  \zeta^{3}_{7}
$)

Prime. 

Not pseudo-unitary. 

\vskip 1ex 
\color{grey}

\noindent(5,4). $5_{\frac{30}{7},4.501}^{7,125}$ \irep{11}:\ \ 
$d_i$ = ($1.0$,
$1.445$,
$-0.356$,
$-0.801$,
$-0.801$) 

\vskip 0.7ex
\hangindent=3em \hangafter=1
$D^2= 4.501 = 
7-7  c^{1}_{7}
-14  c^{2}_{7}
$

\vskip 0.7ex
\hangindent=3em \hangafter=1
$T = ( 0,
\frac{4}{7},
\frac{5}{7},
\frac{3}{7},
\frac{3}{7} )
$,

\vskip 0.7ex
\hangindent=3em \hangafter=1
$S$ = ($ 1$,
$ 1-c^{2}_{7}
$,
$ -c^{1}_{7}
-2  c^{2}_{7}
$,
$ -c^{1}_{7}
-c^{2}_{7}
$,
$ -c^{1}_{7}
-c^{2}_{7}
$;\ \ 
$ -c^{1}_{7}
-2  c^{2}_{7}
$,
$ -1$,
$ c^{1}_{7}
+c^{2}_{7}
$,
$ c^{1}_{7}
+c^{2}_{7}
$;\ \ 
$ -1+c^{2}_{7}
$,
$ -c^{1}_{7}
-c^{2}_{7}
$,
$ -c^{1}_{7}
-c^{2}_{7}
$;\ \ 
$ 1+\zeta^{1}_{7}
+2\zeta^{-1}_{7}
+2\zeta^{2}_{7}
+\zeta^{-2}_{7}
+2\zeta^{3}_{7}
$,
$ -1-\zeta^{-1}_{7}
-\zeta^{2}_{7}
-2  \zeta^{3}_{7}
$;\ \ 
$ 1+\zeta^{1}_{7}
+2\zeta^{-1}_{7}
+2\zeta^{2}_{7}
+\zeta^{-2}_{7}
+2\zeta^{3}_{7}
$)

Prime. 

Not pseudo-unitary. 

\vskip 1ex 
\color{grey}

\noindent(5,5). $5_{\frac{6}{7},2.155}^{7,342}$ \irep{11}:\ \ 
$d_i$ = ($1.0$,
$0.554$,
$0.554$,
$-0.246$,
$-0.692$) 

\vskip 0.7ex
\hangindent=3em \hangafter=1
$D^2= 2.155 = 
14-7  c^{1}_{7}
+7c^{2}_{7}
$

\vskip 0.7ex
\hangindent=3em \hangafter=1
$T = ( 0,
\frac{2}{7},
\frac{2}{7},
\frac{5}{7},
\frac{1}{7} )
$,

\vskip 0.7ex
\hangindent=3em \hangafter=1
$S$ = ($ 1$,
$ 1+c^{2}_{7}
$,
$ 1+c^{2}_{7}
$,
$ 1-c^{1}_{7}
$,
$ 1-c^{1}_{7}
+c^{2}_{7}
$;\ \ 
$ -1-\zeta^{1}_{7}
-2  \zeta^{-2}_{7}
-\zeta^{3}_{7}
$,
$ \zeta^{1}_{7}
-s^{2}_{7}
+\zeta^{3}_{7}
$,
$ -1-c^{2}_{7}
$,
$ 1+c^{2}_{7}
$;\ \ 
$ -1-\zeta^{1}_{7}
-2  \zeta^{-2}_{7}
-\zeta^{3}_{7}
$,
$ -1-c^{2}_{7}
$,
$ 1+c^{2}_{7}
$;\ \ 
$ 1-c^{1}_{7}
+c^{2}_{7}
$,
$ -1$;\ \ 
$ -1+c^{1}_{7}
$)

Prime. 

Not pseudo-unitary. 

\vskip 1ex 
\color{grey}

\noindent(5,6). $5_{\frac{50}{7},2.155}^{7,255}$ \irep{11}:\ \ 
$d_i$ = ($1.0$,
$0.554$,
$0.554$,
$-0.246$,
$-0.692$) 

\vskip 0.7ex
\hangindent=3em \hangafter=1
$D^2= 2.155 = 
14-7  c^{1}_{7}
+7c^{2}_{7}
$

\vskip 0.7ex
\hangindent=3em \hangafter=1
$T = ( 0,
\frac{5}{7},
\frac{5}{7},
\frac{2}{7},
\frac{6}{7} )
$,

\vskip 0.7ex
\hangindent=3em \hangafter=1
$S$ = ($ 1$,
$ 1+c^{2}_{7}
$,
$ 1+c^{2}_{7}
$,
$ 1-c^{1}_{7}
$,
$ 1-c^{1}_{7}
+c^{2}_{7}
$;\ \ 
$ \zeta^{1}_{7}
-s^{2}_{7}
+\zeta^{3}_{7}
$,
$ -1-\zeta^{1}_{7}
-2  \zeta^{-2}_{7}
-\zeta^{3}_{7}
$,
$ -1-c^{2}_{7}
$,
$ 1+c^{2}_{7}
$;\ \ 
$ \zeta^{1}_{7}
-s^{2}_{7}
+\zeta^{3}_{7}
$,
$ -1-c^{2}_{7}
$,
$ 1+c^{2}_{7}
$;\ \ 
$ 1-c^{1}_{7}
+c^{2}_{7}
$,
$ -1$;\ \ 
$ -1+c^{1}_{7}
$)

Prime. 

Not pseudo-unitary. 

\vskip 1ex 

}

\subsection{Rank 6}

{\small
\black

\noindent(1,1). $6_{3,6.}^{12,534}$ \irep{34}:\ \ 
$d_i$ = ($1.0$,
$1.0$,
$1.0$,
$1.0$,
$1.0$,
$1.0$) 

\vskip 0.7ex
\hangindent=3em \hangafter=1
$D^2= 6.0 = 
6$

\vskip 0.7ex
\hangindent=3em \hangafter=1
$T = ( 0,
\frac{1}{3},
\frac{1}{3},
\frac{1}{4},
\frac{7}{12},
\frac{7}{12} )
$,

\vskip 0.7ex
\hangindent=3em \hangafter=1
$S$ = ($ 1$,
$ 1$,
$ 1$,
$ 1$,
$ 1$,
$ 1$;\ \ 
$ \zeta_{3}^{1}$,
$ -\zeta_{6}^{1}$,
$ 1$,
$ -\zeta_{6}^{1}$,
$ \zeta_{3}^{1}$;\ \ 
$ \zeta_{3}^{1}$,
$ 1$,
$ \zeta_{3}^{1}$,
$ -\zeta_{6}^{1}$;\ \ 
$ -1$,
$ -1$,
$ -1$;\ \ 
$ -\zeta_{3}^{1}$,
$ \zeta_{6}^{1}$;\ \ 
$ -\zeta_{3}^{1}$)

Factors = $2_{1,2.}^{4,437}\boxtimes 3_{2,3.}^{3,527} $

\vskip 1ex 
\color{grey}

\noindent(1,2). $6_{1,6.}^{12,701}$ \irep{34}:\ \ 
$d_i$ = ($1.0$,
$1.0$,
$1.0$,
$1.0$,
$1.0$,
$1.0$) 

\vskip 0.7ex
\hangindent=3em \hangafter=1
$D^2= 6.0 = 
6$

\vskip 0.7ex
\hangindent=3em \hangafter=1
$T = ( 0,
\frac{1}{3},
\frac{1}{3},
\frac{3}{4},
\frac{1}{12},
\frac{1}{12} )
$,

\vskip 0.7ex
\hangindent=3em \hangafter=1
$S$ = ($ 1$,
$ 1$,
$ 1$,
$ 1$,
$ 1$,
$ 1$;\ \ 
$ \zeta_{3}^{1}$,
$ -\zeta_{6}^{1}$,
$ 1$,
$ -\zeta_{6}^{1}$,
$ \zeta_{3}^{1}$;\ \ 
$ \zeta_{3}^{1}$,
$ 1$,
$ \zeta_{3}^{1}$,
$ -\zeta_{6}^{1}$;\ \ 
$ -1$,
$ -1$,
$ -1$;\ \ 
$ -\zeta_{3}^{1}$,
$ \zeta_{6}^{1}$;\ \ 
$ -\zeta_{3}^{1}$)

Factors = $2_{7,2.}^{4,625}\boxtimes 3_{2,3.}^{3,527} $

\vskip 1ex 
\color{grey}

\noindent(1,3). $6_{7,6.}^{12,113}$ \irep{34}:\ \ 
$d_i$ = ($1.0$,
$1.0$,
$1.0$,
$1.0$,
$1.0$,
$1.0$) 

\vskip 0.7ex
\hangindent=3em \hangafter=1
$D^2= 6.0 = 
6$

\vskip 0.7ex
\hangindent=3em \hangafter=1
$T = ( 0,
\frac{2}{3},
\frac{2}{3},
\frac{1}{4},
\frac{11}{12},
\frac{11}{12} )
$,

\vskip 0.7ex
\hangindent=3em \hangafter=1
$S$ = ($ 1$,
$ 1$,
$ 1$,
$ 1$,
$ 1$,
$ 1$;\ \ 
$ -\zeta_{6}^{1}$,
$ \zeta_{3}^{1}$,
$ 1$,
$ -\zeta_{6}^{1}$,
$ \zeta_{3}^{1}$;\ \ 
$ -\zeta_{6}^{1}$,
$ 1$,
$ \zeta_{3}^{1}$,
$ -\zeta_{6}^{1}$;\ \ 
$ -1$,
$ -1$,
$ -1$;\ \ 
$ \zeta_{6}^{1}$,
$ -\zeta_{3}^{1}$;\ \ 
$ \zeta_{6}^{1}$)

Factors = $2_{1,2.}^{4,437}\boxtimes 3_{6,3.}^{3,138} $

\vskip 1ex 
\color{grey}

\noindent(1,4). $6_{5,6.}^{12,298}$ \irep{34}:\ \ 
$d_i$ = ($1.0$,
$1.0$,
$1.0$,
$1.0$,
$1.0$,
$1.0$) 

\vskip 0.7ex
\hangindent=3em \hangafter=1
$D^2= 6.0 = 
6$

\vskip 0.7ex
\hangindent=3em \hangafter=1
$T = ( 0,
\frac{2}{3},
\frac{2}{3},
\frac{3}{4},
\frac{5}{12},
\frac{5}{12} )
$,

\vskip 0.7ex
\hangindent=3em \hangafter=1
$S$ = ($ 1$,
$ 1$,
$ 1$,
$ 1$,
$ 1$,
$ 1$;\ \ 
$ -\zeta_{6}^{1}$,
$ \zeta_{3}^{1}$,
$ 1$,
$ -\zeta_{6}^{1}$,
$ \zeta_{3}^{1}$;\ \ 
$ -\zeta_{6}^{1}$,
$ 1$,
$ \zeta_{3}^{1}$,
$ -\zeta_{6}^{1}$;\ \ 
$ -1$,
$ -1$,
$ -1$;\ \ 
$ \zeta_{6}^{1}$,
$ -\zeta_{3}^{1}$;\ \ 
$ \zeta_{6}^{1}$)

Factors = $2_{7,2.}^{4,625}\boxtimes 3_{6,3.}^{3,138} $

\vskip 1ex 
\black

\noindent(2,1). $6_{\frac{3}{2},8.}^{16,688}$ \irep{39}:\ \ 
$d_i$ = ($1.0$,
$1.0$,
$1.0$,
$1.0$,
$1.414$,
$1.414$) 

\vskip 0.7ex
\hangindent=3em \hangafter=1
$D^2= 8.0 = 
8$

\vskip 0.7ex
\hangindent=3em \hangafter=1
$T = ( 0,
\frac{1}{2},
\frac{1}{4},
\frac{3}{4},
\frac{1}{16},
\frac{5}{16} )
$,

\vskip 0.7ex
\hangindent=3em \hangafter=1
$S$ = ($ 1$,
$ 1$,
$ 1$,
$ 1$,
$ \sqrt{2}$,
$ \sqrt{2}$;\ \ 
$ 1$,
$ 1$,
$ 1$,
$ -\sqrt{2}$,
$ -\sqrt{2}$;\ \ 
$ -1$,
$ -1$,
$ \sqrt{2}$,
$ -\sqrt{2}$;\ \ 
$ -1$,
$ -\sqrt{2}$,
$ \sqrt{2}$;\ \ 
$0$,
$0$;\ \ 
$0$)

Factors = $2_{1,2.}^{4,437}\boxtimes 3_{\frac{1}{2},4.}^{16,598} $

\vskip 1ex 
\color{grey}

\noindent(2,2). $6_{\frac{5}{2},8.}^{16,511}$ \irep{39}:\ \ 
$d_i$ = ($1.0$,
$1.0$,
$1.0$,
$1.0$,
$1.414$,
$1.414$) 

\vskip 0.7ex
\hangindent=3em \hangafter=1
$D^2= 8.0 = 
8$

\vskip 0.7ex
\hangindent=3em \hangafter=1
$T = ( 0,
\frac{1}{2},
\frac{1}{4},
\frac{3}{4},
\frac{3}{16},
\frac{7}{16} )
$,

\vskip 0.7ex
\hangindent=3em \hangafter=1
$S$ = ($ 1$,
$ 1$,
$ 1$,
$ 1$,
$ \sqrt{2}$,
$ \sqrt{2}$;\ \ 
$ 1$,
$ 1$,
$ 1$,
$ -\sqrt{2}$,
$ -\sqrt{2}$;\ \ 
$ -1$,
$ -1$,
$ \sqrt{2}$,
$ -\sqrt{2}$;\ \ 
$ -1$,
$ -\sqrt{2}$,
$ \sqrt{2}$;\ \ 
$0$,
$0$;\ \ 
$0$)

Factors = $2_{7,2.}^{4,625}\boxtimes 3_{\frac{7}{2},4.}^{16,332} $

\vskip 1ex 
\color{grey}

\noindent(2,3). $6_{\frac{11}{2},8.}^{16,548}$ \irep{39}:\ \ 
$d_i$ = ($1.0$,
$1.0$,
$1.0$,
$1.0$,
$1.414$,
$1.414$) 

\vskip 0.7ex
\hangindent=3em \hangafter=1
$D^2= 8.0 = 
8$

\vskip 0.7ex
\hangindent=3em \hangafter=1
$T = ( 0,
\frac{1}{2},
\frac{1}{4},
\frac{3}{4},
\frac{9}{16},
\frac{13}{16} )
$,

\vskip 0.7ex
\hangindent=3em \hangafter=1
$S$ = ($ 1$,
$ 1$,
$ 1$,
$ 1$,
$ \sqrt{2}$,
$ \sqrt{2}$;\ \ 
$ 1$,
$ 1$,
$ 1$,
$ -\sqrt{2}$,
$ -\sqrt{2}$;\ \ 
$ -1$,
$ -1$,
$ \sqrt{2}$,
$ -\sqrt{2}$;\ \ 
$ -1$,
$ -\sqrt{2}$,
$ \sqrt{2}$;\ \ 
$0$,
$0$;\ \ 
$0$)

Factors = $2_{1,2.}^{4,437}\boxtimes 3_{\frac{9}{2},4.}^{16,156} $

\vskip 1ex 
\color{grey}

\noindent(2,4). $6_{\frac{13}{2},8.}^{16,107}$ \irep{39}:\ \ 
$d_i$ = ($1.0$,
$1.0$,
$1.0$,
$1.0$,
$1.414$,
$1.414$) 

\vskip 0.7ex
\hangindent=3em \hangafter=1
$D^2= 8.0 = 
8$

\vskip 0.7ex
\hangindent=3em \hangafter=1
$T = ( 0,
\frac{1}{2},
\frac{1}{4},
\frac{3}{4},
\frac{11}{16},
\frac{15}{16} )
$,

\vskip 0.7ex
\hangindent=3em \hangafter=1
$S$ = ($ 1$,
$ 1$,
$ 1$,
$ 1$,
$ \sqrt{2}$,
$ \sqrt{2}$;\ \ 
$ 1$,
$ 1$,
$ 1$,
$ -\sqrt{2}$,
$ -\sqrt{2}$;\ \ 
$ -1$,
$ -1$,
$ \sqrt{2}$,
$ -\sqrt{2}$;\ \ 
$ -1$,
$ -\sqrt{2}$,
$ \sqrt{2}$;\ \ 
$0$,
$0$;\ \ 
$0$)

Factors = $2_{7,2.}^{4,625}\boxtimes 3_{\frac{15}{2},4.}^{16,639} $

\vskip 1ex 
\color{grey}

\noindent(2,5). $6_{\frac{15}{2},8.}^{16,357}$ \irep{39}:\ \ 
$d_i$ = ($1.0$,
$1.0$,
$1.0$,
$1.0$,
$-1.414$,
$-1.414$) 

\vskip 0.7ex
\hangindent=3em \hangafter=1
$D^2= 8.0 = 
8$

\vskip 0.7ex
\hangindent=3em \hangafter=1
$T = ( 0,
\frac{1}{2},
\frac{1}{4},
\frac{3}{4},
\frac{1}{16},
\frac{13}{16} )
$,

\vskip 0.7ex
\hangindent=3em \hangafter=1
$S$ = ($ 1$,
$ 1$,
$ 1$,
$ 1$,
$ -\sqrt{2}$,
$ -\sqrt{2}$;\ \ 
$ 1$,
$ 1$,
$ 1$,
$ \sqrt{2}$,
$ \sqrt{2}$;\ \ 
$ -1$,
$ -1$,
$ \sqrt{2}$,
$ -\sqrt{2}$;\ \ 
$ -1$,
$ -\sqrt{2}$,
$ \sqrt{2}$;\ \ 
$0$,
$0$;\ \ 
$0$)

Factors = $2_{1,2.}^{4,437}\boxtimes 3_{\frac{13}{2},4.}^{16,830} $

Pseudo-unitary $\sim$  
$6_{\frac{7}{2},8.}^{16,246}$

\vskip 1ex 
\color{grey}

\noindent(2,6). $6_{\frac{1}{2},8.}^{16,710}$ \irep{39}:\ \ 
$d_i$ = ($1.0$,
$1.0$,
$1.0$,
$1.0$,
$-1.414$,
$-1.414$) 

\vskip 0.7ex
\hangindent=3em \hangafter=1
$D^2= 8.0 = 
8$

\vskip 0.7ex
\hangindent=3em \hangafter=1
$T = ( 0,
\frac{1}{2},
\frac{1}{4},
\frac{3}{4},
\frac{3}{16},
\frac{15}{16} )
$,

\vskip 0.7ex
\hangindent=3em \hangafter=1
$S$ = ($ 1$,
$ 1$,
$ 1$,
$ 1$,
$ -\sqrt{2}$,
$ -\sqrt{2}$;\ \ 
$ 1$,
$ 1$,
$ 1$,
$ \sqrt{2}$,
$ \sqrt{2}$;\ \ 
$ -1$,
$ -1$,
$ \sqrt{2}$,
$ -\sqrt{2}$;\ \ 
$ -1$,
$ -\sqrt{2}$,
$ \sqrt{2}$;\ \ 
$0$,
$0$;\ \ 
$0$)

Factors = $2_{7,2.}^{4,625}\boxtimes 3_{\frac{3}{2},4.}^{16,538} $

Pseudo-unitary $\sim$  
$6_{\frac{9}{2},8.}^{16,107}$

\vskip 1ex 
\color{grey}

\noindent(2,7). $6_{\frac{7}{2},8.}^{16,346}$ \irep{39}:\ \ 
$d_i$ = ($1.0$,
$1.0$,
$1.0$,
$1.0$,
$-1.414$,
$-1.414$) 

\vskip 0.7ex
\hangindent=3em \hangafter=1
$D^2= 8.0 = 
8$

\vskip 0.7ex
\hangindent=3em \hangafter=1
$T = ( 0,
\frac{1}{2},
\frac{1}{4},
\frac{3}{4},
\frac{5}{16},
\frac{9}{16} )
$,

\vskip 0.7ex
\hangindent=3em \hangafter=1
$S$ = ($ 1$,
$ 1$,
$ 1$,
$ 1$,
$ -\sqrt{2}$,
$ -\sqrt{2}$;\ \ 
$ 1$,
$ 1$,
$ 1$,
$ \sqrt{2}$,
$ \sqrt{2}$;\ \ 
$ -1$,
$ -1$,
$ -\sqrt{2}$,
$ \sqrt{2}$;\ \ 
$ -1$,
$ \sqrt{2}$,
$ -\sqrt{2}$;\ \ 
$0$,
$0$;\ \ 
$0$)

Factors = $2_{1,2.}^{4,437}\boxtimes 3_{\frac{5}{2},4.}^{16,345} $

Pseudo-unitary $\sim$  
$6_{\frac{15}{2},8.}^{16,107}$

\vskip 1ex 
\color{grey}

\noindent(2,8). $6_{\frac{9}{2},8.}^{16,357}$ \irep{39}:\ \ 
$d_i$ = ($1.0$,
$1.0$,
$1.0$,
$1.0$,
$-1.414$,
$-1.414$) 

\vskip 0.7ex
\hangindent=3em \hangafter=1
$D^2= 8.0 = 
8$

\vskip 0.7ex
\hangindent=3em \hangafter=1
$T = ( 0,
\frac{1}{2},
\frac{1}{4},
\frac{3}{4},
\frac{7}{16},
\frac{11}{16} )
$,

\vskip 0.7ex
\hangindent=3em \hangafter=1
$S$ = ($ 1$,
$ 1$,
$ 1$,
$ 1$,
$ -\sqrt{2}$,
$ -\sqrt{2}$;\ \ 
$ 1$,
$ 1$,
$ 1$,
$ \sqrt{2}$,
$ \sqrt{2}$;\ \ 
$ -1$,
$ -1$,
$ -\sqrt{2}$,
$ \sqrt{2}$;\ \ 
$ -1$,
$ \sqrt{2}$,
$ -\sqrt{2}$;\ \ 
$0$,
$0$;\ \ 
$0$)

Factors = $2_{7,2.}^{4,625}\boxtimes 3_{\frac{11}{2},4.}^{16,564} $

Pseudo-unitary $\sim$  
$6_{\frac{1}{2},8.}^{16,460}$

\vskip 1ex 
\black

\noindent(3,1). $6_{\frac{15}{2},8.}^{16,107}$ \irep{39}:\ \ 
$d_i$ = ($1.0$,
$1.0$,
$1.0$,
$1.0$,
$1.414$,
$1.414$) 

\vskip 0.7ex
\hangindent=3em \hangafter=1
$D^2= 8.0 = 
8$

\vskip 0.7ex
\hangindent=3em \hangafter=1
$T = ( 0,
\frac{1}{2},
\frac{1}{4},
\frac{3}{4},
\frac{1}{16},
\frac{13}{16} )
$,

\vskip 0.7ex
\hangindent=3em \hangafter=1
$S$ = ($ 1$,
$ 1$,
$ 1$,
$ 1$,
$ \sqrt{2}$,
$ \sqrt{2}$;\ \ 
$ 1$,
$ 1$,
$ 1$,
$ -\sqrt{2}$,
$ -\sqrt{2}$;\ \ 
$ -1$,
$ -1$,
$ -\sqrt{2}$,
$ \sqrt{2}$;\ \ 
$ -1$,
$ \sqrt{2}$,
$ -\sqrt{2}$;\ \ 
$0$,
$0$;\ \ 
$0$)

Factors = $2_{7,2.}^{4,625}\boxtimes 3_{\frac{1}{2},4.}^{16,598} $

\vskip 1ex 
\color{grey}

\noindent(3,2). $6_{\frac{1}{2},8.}^{16,460}$ \irep{39}:\ \ 
$d_i$ = ($1.0$,
$1.0$,
$1.0$,
$1.0$,
$1.414$,
$1.414$) 

\vskip 0.7ex
\hangindent=3em \hangafter=1
$D^2= 8.0 = 
8$

\vskip 0.7ex
\hangindent=3em \hangafter=1
$T = ( 0,
\frac{1}{2},
\frac{1}{4},
\frac{3}{4},
\frac{3}{16},
\frac{15}{16} )
$,

\vskip 0.7ex
\hangindent=3em \hangafter=1
$S$ = ($ 1$,
$ 1$,
$ 1$,
$ 1$,
$ \sqrt{2}$,
$ \sqrt{2}$;\ \ 
$ 1$,
$ 1$,
$ 1$,
$ -\sqrt{2}$,
$ -\sqrt{2}$;\ \ 
$ -1$,
$ -1$,
$ -\sqrt{2}$,
$ \sqrt{2}$;\ \ 
$ -1$,
$ \sqrt{2}$,
$ -\sqrt{2}$;\ \ 
$0$,
$0$;\ \ 
$0$)

Factors = $2_{1,2.}^{4,437}\boxtimes 3_{\frac{15}{2},4.}^{16,639} $

\vskip 1ex 
\color{grey}

\noindent(3,3). $6_{\frac{7}{2},8.}^{16,246}$ \irep{39}:\ \ 
$d_i$ = ($1.0$,
$1.0$,
$1.0$,
$1.0$,
$1.414$,
$1.414$) 

\vskip 0.7ex
\hangindent=3em \hangafter=1
$D^2= 8.0 = 
8$

\vskip 0.7ex
\hangindent=3em \hangafter=1
$T = ( 0,
\frac{1}{2},
\frac{1}{4},
\frac{3}{4},
\frac{5}{16},
\frac{9}{16} )
$,

\vskip 0.7ex
\hangindent=3em \hangafter=1
$S$ = ($ 1$,
$ 1$,
$ 1$,
$ 1$,
$ \sqrt{2}$,
$ \sqrt{2}$;\ \ 
$ 1$,
$ 1$,
$ 1$,
$ -\sqrt{2}$,
$ -\sqrt{2}$;\ \ 
$ -1$,
$ -1$,
$ \sqrt{2}$,
$ -\sqrt{2}$;\ \ 
$ -1$,
$ -\sqrt{2}$,
$ \sqrt{2}$;\ \ 
$0$,
$0$;\ \ 
$0$)

Factors = $2_{7,2.}^{4,625}\boxtimes 3_{\frac{9}{2},4.}^{16,156} $

\vskip 1ex 
\color{grey}

\noindent(3,4). $6_{\frac{9}{2},8.}^{16,107}$ \irep{39}:\ \ 
$d_i$ = ($1.0$,
$1.0$,
$1.0$,
$1.0$,
$1.414$,
$1.414$) 

\vskip 0.7ex
\hangindent=3em \hangafter=1
$D^2= 8.0 = 
8$

\vskip 0.7ex
\hangindent=3em \hangafter=1
$T = ( 0,
\frac{1}{2},
\frac{1}{4},
\frac{3}{4},
\frac{7}{16},
\frac{11}{16} )
$,

\vskip 0.7ex
\hangindent=3em \hangafter=1
$S$ = ($ 1$,
$ 1$,
$ 1$,
$ 1$,
$ \sqrt{2}$,
$ \sqrt{2}$;\ \ 
$ 1$,
$ 1$,
$ 1$,
$ -\sqrt{2}$,
$ -\sqrt{2}$;\ \ 
$ -1$,
$ -1$,
$ \sqrt{2}$,
$ -\sqrt{2}$;\ \ 
$ -1$,
$ -\sqrt{2}$,
$ \sqrt{2}$;\ \ 
$0$,
$0$;\ \ 
$0$)

Factors = $2_{1,2.}^{4,437}\boxtimes 3_{\frac{7}{2},4.}^{16,332} $

\vskip 1ex 
\color{grey}

\noindent(3,5). $6_{\frac{3}{2},8.}^{16,438}$ \irep{39}:\ \ 
$d_i$ = ($1.0$,
$1.0$,
$1.0$,
$1.0$,
$-1.414$,
$-1.414$) 

\vskip 0.7ex
\hangindent=3em \hangafter=1
$D^2= 8.0 = 
8$

\vskip 0.7ex
\hangindent=3em \hangafter=1
$T = ( 0,
\frac{1}{2},
\frac{1}{4},
\frac{3}{4},
\frac{1}{16},
\frac{5}{16} )
$,

\vskip 0.7ex
\hangindent=3em \hangafter=1
$S$ = ($ 1$,
$ 1$,
$ 1$,
$ 1$,
$ -\sqrt{2}$,
$ -\sqrt{2}$;\ \ 
$ 1$,
$ 1$,
$ 1$,
$ \sqrt{2}$,
$ \sqrt{2}$;\ \ 
$ -1$,
$ -1$,
$ -\sqrt{2}$,
$ \sqrt{2}$;\ \ 
$ -1$,
$ \sqrt{2}$,
$ -\sqrt{2}$;\ \ 
$0$,
$0$;\ \ 
$0$)

Factors = $2_{7,2.}^{4,625}\boxtimes 3_{\frac{5}{2},4.}^{16,345} $

Pseudo-unitary $\sim$  
$6_{\frac{11}{2},8.}^{16,548}$

\vskip 1ex 
\color{grey}

\noindent(3,6). $6_{\frac{5}{2},8.}^{16,261}$ \irep{39}:\ \ 
$d_i$ = ($1.0$,
$1.0$,
$1.0$,
$1.0$,
$-1.414$,
$-1.414$) 

\vskip 0.7ex
\hangindent=3em \hangafter=1
$D^2= 8.0 = 
8$

\vskip 0.7ex
\hangindent=3em \hangafter=1
$T = ( 0,
\frac{1}{2},
\frac{1}{4},
\frac{3}{4},
\frac{3}{16},
\frac{7}{16} )
$,

\vskip 0.7ex
\hangindent=3em \hangafter=1
$S$ = ($ 1$,
$ 1$,
$ 1$,
$ 1$,
$ -\sqrt{2}$,
$ -\sqrt{2}$;\ \ 
$ 1$,
$ 1$,
$ 1$,
$ \sqrt{2}$,
$ \sqrt{2}$;\ \ 
$ -1$,
$ -1$,
$ -\sqrt{2}$,
$ \sqrt{2}$;\ \ 
$ -1$,
$ \sqrt{2}$,
$ -\sqrt{2}$;\ \ 
$0$,
$0$;\ \ 
$0$)

Factors = $2_{1,2.}^{4,437}\boxtimes 3_{\frac{3}{2},4.}^{16,538} $

Pseudo-unitary $\sim$  
$6_{\frac{13}{2},8.}^{16,107}$

\vskip 1ex 
\color{grey}

\noindent(3,7). $6_{\frac{11}{2},8.}^{16,798}$ \irep{39}:\ \ 
$d_i$ = ($1.0$,
$1.0$,
$1.0$,
$1.0$,
$-1.414$,
$-1.414$) 

\vskip 0.7ex
\hangindent=3em \hangafter=1
$D^2= 8.0 = 
8$

\vskip 0.7ex
\hangindent=3em \hangafter=1
$T = ( 0,
\frac{1}{2},
\frac{1}{4},
\frac{3}{4},
\frac{9}{16},
\frac{13}{16} )
$,

\vskip 0.7ex
\hangindent=3em \hangafter=1
$S$ = ($ 1$,
$ 1$,
$ 1$,
$ 1$,
$ -\sqrt{2}$,
$ -\sqrt{2}$;\ \ 
$ 1$,
$ 1$,
$ 1$,
$ \sqrt{2}$,
$ \sqrt{2}$;\ \ 
$ -1$,
$ -1$,
$ -\sqrt{2}$,
$ \sqrt{2}$;\ \ 
$ -1$,
$ \sqrt{2}$,
$ -\sqrt{2}$;\ \ 
$0$,
$0$;\ \ 
$0$)

Factors = $2_{7,2.}^{4,625}\boxtimes 3_{\frac{13}{2},4.}^{16,830} $

Pseudo-unitary $\sim$  
$6_{\frac{3}{2},8.}^{16,688}$

\vskip 1ex 
\color{grey}

\noindent(3,8). $6_{\frac{13}{2},8.}^{16,132}$ \irep{39}:\ \ 
$d_i$ = ($1.0$,
$1.0$,
$1.0$,
$1.0$,
$-1.414$,
$-1.414$) 

\vskip 0.7ex
\hangindent=3em \hangafter=1
$D^2= 8.0 = 
8$

\vskip 0.7ex
\hangindent=3em \hangafter=1
$T = ( 0,
\frac{1}{2},
\frac{1}{4},
\frac{3}{4},
\frac{11}{16},
\frac{15}{16} )
$,

\vskip 0.7ex
\hangindent=3em \hangafter=1
$S$ = ($ 1$,
$ 1$,
$ 1$,
$ 1$,
$ -\sqrt{2}$,
$ -\sqrt{2}$;\ \ 
$ 1$,
$ 1$,
$ 1$,
$ \sqrt{2}$,
$ \sqrt{2}$;\ \ 
$ -1$,
$ -1$,
$ -\sqrt{2}$,
$ \sqrt{2}$;\ \ 
$ -1$,
$ \sqrt{2}$,
$ -\sqrt{2}$;\ \ 
$0$,
$0$;\ \ 
$0$)

Factors = $2_{1,2.}^{4,437}\boxtimes 3_{\frac{11}{2},4.}^{16,564} $

Pseudo-unitary $\sim$  
$6_{\frac{5}{2},8.}^{16,511}$

\vskip 1ex 
\black

\noindent(4,1). $6_{\frac{24}{5},10.85}^{15,257}$ \irep{38}:\ \ 
$d_i$ = ($1.0$,
$1.0$,
$1.0$,
$1.618$,
$1.618$,
$1.618$) 

\vskip 0.7ex
\hangindent=3em \hangafter=1
$D^2= 10.854 = 
\frac{15+3\sqrt{5}}{2}$

\vskip 0.7ex
\hangindent=3em \hangafter=1
$T = ( 0,
\frac{1}{3},
\frac{1}{3},
\frac{2}{5},
\frac{11}{15},
\frac{11}{15} )
$,

\vskip 0.7ex
\hangindent=3em \hangafter=1
$S$ = ($ 1$,
$ 1$,
$ 1$,
$ \frac{1+\sqrt{5}}{2}$,
$ \frac{1+\sqrt{5}}{2}$,
$ \frac{1+\sqrt{5}}{2}$;\ \ 
$ \zeta_{3}^{1}$,
$ -\zeta_{6}^{1}$,
$ \frac{1+\sqrt{5}}{2}$,
$ -\frac{1+\sqrt{5}}{2}\zeta_{6}^{1}$,
$ \frac{1+\sqrt{5}}{2}\zeta_{3}^{1}$;\ \ 
$ \zeta_{3}^{1}$,
$ \frac{1+\sqrt{5}}{2}$,
$ \frac{1+\sqrt{5}}{2}\zeta_{3}^{1}$,
$ -\frac{1+\sqrt{5}}{2}\zeta_{6}^{1}$;\ \ 
$ -1$,
$ -1$,
$ -1$;\ \ 
$ -\zeta_{3}^{1}$,
$ \zeta_{6}^{1}$;\ \ 
$ -\zeta_{3}^{1}$)

Factors = $2_{\frac{14}{5},3.618}^{5,395}\boxtimes 3_{2,3.}^{3,527} $

\vskip 1ex 
\color{grey}

\noindent(4,2). $6_{\frac{36}{5},10.85}^{15,166}$ \irep{38}:\ \ 
$d_i$ = ($1.0$,
$1.0$,
$1.0$,
$1.618$,
$1.618$,
$1.618$) 

\vskip 0.7ex
\hangindent=3em \hangafter=1
$D^2= 10.854 = 
\frac{15+3\sqrt{5}}{2}$

\vskip 0.7ex
\hangindent=3em \hangafter=1
$T = ( 0,
\frac{1}{3},
\frac{1}{3},
\frac{3}{5},
\frac{14}{15},
\frac{14}{15} )
$,

\vskip 0.7ex
\hangindent=3em \hangafter=1
$S$ = ($ 1$,
$ 1$,
$ 1$,
$ \frac{1+\sqrt{5}}{2}$,
$ \frac{1+\sqrt{5}}{2}$,
$ \frac{1+\sqrt{5}}{2}$;\ \ 
$ \zeta_{3}^{1}$,
$ -\zeta_{6}^{1}$,
$ \frac{1+\sqrt{5}}{2}$,
$ -\frac{1+\sqrt{5}}{2}\zeta_{6}^{1}$,
$ \frac{1+\sqrt{5}}{2}\zeta_{3}^{1}$;\ \ 
$ \zeta_{3}^{1}$,
$ \frac{1+\sqrt{5}}{2}$,
$ \frac{1+\sqrt{5}}{2}\zeta_{3}^{1}$,
$ -\frac{1+\sqrt{5}}{2}\zeta_{6}^{1}$;\ \ 
$ -1$,
$ -1$,
$ -1$;\ \ 
$ -\zeta_{3}^{1}$,
$ \zeta_{6}^{1}$;\ \ 
$ -\zeta_{3}^{1}$)

Factors = $2_{\frac{26}{5},3.618}^{5,720}\boxtimes 3_{2,3.}^{3,527} $

\vskip 1ex 
\color{grey}

\noindent(4,3). $6_{\frac{4}{5},10.85}^{15,801}$ \irep{38}:\ \ 
$d_i$ = ($1.0$,
$1.0$,
$1.0$,
$1.618$,
$1.618$,
$1.618$) 

\vskip 0.7ex
\hangindent=3em \hangafter=1
$D^2= 10.854 = 
\frac{15+3\sqrt{5}}{2}$

\vskip 0.7ex
\hangindent=3em \hangafter=1
$T = ( 0,
\frac{2}{3},
\frac{2}{3},
\frac{2}{5},
\frac{1}{15},
\frac{1}{15} )
$,

\vskip 0.7ex
\hangindent=3em \hangafter=1
$S$ = ($ 1$,
$ 1$,
$ 1$,
$ \frac{1+\sqrt{5}}{2}$,
$ \frac{1+\sqrt{5}}{2}$,
$ \frac{1+\sqrt{5}}{2}$;\ \ 
$ -\zeta_{6}^{1}$,
$ \zeta_{3}^{1}$,
$ \frac{1+\sqrt{5}}{2}$,
$ -\frac{1+\sqrt{5}}{2}\zeta_{6}^{1}$,
$ \frac{1+\sqrt{5}}{2}\zeta_{3}^{1}$;\ \ 
$ -\zeta_{6}^{1}$,
$ \frac{1+\sqrt{5}}{2}$,
$ \frac{1+\sqrt{5}}{2}\zeta_{3}^{1}$,
$ -\frac{1+\sqrt{5}}{2}\zeta_{6}^{1}$;\ \ 
$ -1$,
$ -1$,
$ -1$;\ \ 
$ \zeta_{6}^{1}$,
$ -\zeta_{3}^{1}$;\ \ 
$ \zeta_{6}^{1}$)

Factors = $2_{\frac{14}{5},3.618}^{5,395}\boxtimes 3_{6,3.}^{3,138} $

\vskip 1ex 
\color{grey}

\noindent(4,4). $6_{\frac{16}{5},10.85}^{15,262}$ \irep{38}:\ \ 
$d_i$ = ($1.0$,
$1.0$,
$1.0$,
$1.618$,
$1.618$,
$1.618$) 

\vskip 0.7ex
\hangindent=3em \hangafter=1
$D^2= 10.854 = 
\frac{15+3\sqrt{5}}{2}$

\vskip 0.7ex
\hangindent=3em \hangafter=1
$T = ( 0,
\frac{2}{3},
\frac{2}{3},
\frac{3}{5},
\frac{4}{15},
\frac{4}{15} )
$,

\vskip 0.7ex
\hangindent=3em \hangafter=1
$S$ = ($ 1$,
$ 1$,
$ 1$,
$ \frac{1+\sqrt{5}}{2}$,
$ \frac{1+\sqrt{5}}{2}$,
$ \frac{1+\sqrt{5}}{2}$;\ \ 
$ -\zeta_{6}^{1}$,
$ \zeta_{3}^{1}$,
$ \frac{1+\sqrt{5}}{2}$,
$ -\frac{1+\sqrt{5}}{2}\zeta_{6}^{1}$,
$ \frac{1+\sqrt{5}}{2}\zeta_{3}^{1}$;\ \ 
$ -\zeta_{6}^{1}$,
$ \frac{1+\sqrt{5}}{2}$,
$ \frac{1+\sqrt{5}}{2}\zeta_{3}^{1}$,
$ -\frac{1+\sqrt{5}}{2}\zeta_{6}^{1}$;\ \ 
$ -1$,
$ -1$,
$ -1$;\ \ 
$ \zeta_{6}^{1}$,
$ -\zeta_{3}^{1}$;\ \ 
$ \zeta_{6}^{1}$)

Factors = $2_{\frac{26}{5},3.618}^{5,720}\boxtimes 3_{6,3.}^{3,138} $

\vskip 1ex 
\color{grey}

\noindent(4,5). $6_{\frac{12}{5},4.145}^{15,440}$ \irep{38}:\ \ 
$d_i$ = ($1.0$,
$1.0$,
$1.0$,
$-0.618$,
$-0.618$,
$-0.618$) 

\vskip 0.7ex
\hangindent=3em \hangafter=1
$D^2= 4.145 = 
\frac{15-3\sqrt{5}}{2}$

\vskip 0.7ex
\hangindent=3em \hangafter=1
$T = ( 0,
\frac{1}{3},
\frac{1}{3},
\frac{1}{5},
\frac{8}{15},
\frac{8}{15} )
$,

\vskip 0.7ex
\hangindent=3em \hangafter=1
$S$ = ($ 1$,
$ 1$,
$ 1$,
$ \frac{1-\sqrt{5}}{2}$,
$ \frac{1-\sqrt{5}}{2}$,
$ \frac{1-\sqrt{5}}{2}$;\ \ 
$ \zeta_{3}^{1}$,
$ -\zeta_{6}^{1}$,
$ \frac{1-\sqrt{5}}{2}$,
$ \frac{1-\sqrt{5}}{2}\zeta_{3}^{1}$,
$ -\frac{1-\sqrt{5}}{2}\zeta_{6}^{1}$;\ \ 
$ \zeta_{3}^{1}$,
$ \frac{1-\sqrt{5}}{2}$,
$ -\frac{1-\sqrt{5}}{2}\zeta_{6}^{1}$,
$ \frac{1-\sqrt{5}}{2}\zeta_{3}^{1}$;\ \ 
$ -1$,
$ -1$,
$ -1$;\ \ 
$ -\zeta_{3}^{1}$,
$ \zeta_{6}^{1}$;\ \ 
$ -\zeta_{3}^{1}$)

Factors = $2_{\frac{2}{5},1.381}^{5,120}\boxtimes 3_{2,3.}^{3,527} $

Not pseudo-unitary. 

\vskip 1ex 
\color{grey}

\noindent(4,6). $6_{\frac{8}{5},4.145}^{15,481}$ \irep{38}:\ \ 
$d_i$ = ($1.0$,
$1.0$,
$1.0$,
$-0.618$,
$-0.618$,
$-0.618$) 

\vskip 0.7ex
\hangindent=3em \hangafter=1
$D^2= 4.145 = 
\frac{15-3\sqrt{5}}{2}$

\vskip 0.7ex
\hangindent=3em \hangafter=1
$T = ( 0,
\frac{1}{3},
\frac{1}{3},
\frac{4}{5},
\frac{2}{15},
\frac{2}{15} )
$,

\vskip 0.7ex
\hangindent=3em \hangafter=1
$S$ = ($ 1$,
$ 1$,
$ 1$,
$ \frac{1-\sqrt{5}}{2}$,
$ \frac{1-\sqrt{5}}{2}$,
$ \frac{1-\sqrt{5}}{2}$;\ \ 
$ \zeta_{3}^{1}$,
$ -\zeta_{6}^{1}$,
$ \frac{1-\sqrt{5}}{2}$,
$ \frac{1-\sqrt{5}}{2}\zeta_{3}^{1}$,
$ -\frac{1-\sqrt{5}}{2}\zeta_{6}^{1}$;\ \ 
$ \zeta_{3}^{1}$,
$ \frac{1-\sqrt{5}}{2}$,
$ -\frac{1-\sqrt{5}}{2}\zeta_{6}^{1}$,
$ \frac{1-\sqrt{5}}{2}\zeta_{3}^{1}$;\ \ 
$ -1$,
$ -1$,
$ -1$;\ \ 
$ -\zeta_{3}^{1}$,
$ \zeta_{6}^{1}$;\ \ 
$ -\zeta_{3}^{1}$)

Factors = $2_{\frac{38}{5},1.381}^{5,491}\boxtimes 3_{2,3.}^{3,527} $

Not pseudo-unitary. 

\vskip 1ex 
\color{grey}

\noindent(4,7). $6_{\frac{32}{5},4.145}^{15,350}$ \irep{38}:\ \ 
$d_i$ = ($1.0$,
$1.0$,
$1.0$,
$-0.618$,
$-0.618$,
$-0.618$) 

\vskip 0.7ex
\hangindent=3em \hangafter=1
$D^2= 4.145 = 
\frac{15-3\sqrt{5}}{2}$

\vskip 0.7ex
\hangindent=3em \hangafter=1
$T = ( 0,
\frac{2}{3},
\frac{2}{3},
\frac{1}{5},
\frac{13}{15},
\frac{13}{15} )
$,

\vskip 0.7ex
\hangindent=3em \hangafter=1
$S$ = ($ 1$,
$ 1$,
$ 1$,
$ \frac{1-\sqrt{5}}{2}$,
$ \frac{1-\sqrt{5}}{2}$,
$ \frac{1-\sqrt{5}}{2}$;\ \ 
$ -\zeta_{6}^{1}$,
$ \zeta_{3}^{1}$,
$ \frac{1-\sqrt{5}}{2}$,
$ \frac{1-\sqrt{5}}{2}\zeta_{3}^{1}$,
$ -\frac{1-\sqrt{5}}{2}\zeta_{6}^{1}$;\ \ 
$ -\zeta_{6}^{1}$,
$ \frac{1-\sqrt{5}}{2}$,
$ -\frac{1-\sqrt{5}}{2}\zeta_{6}^{1}$,
$ \frac{1-\sqrt{5}}{2}\zeta_{3}^{1}$;\ \ 
$ -1$,
$ -1$,
$ -1$;\ \ 
$ \zeta_{6}^{1}$,
$ -\zeta_{3}^{1}$;\ \ 
$ \zeta_{6}^{1}$)

Factors = $2_{\frac{2}{5},1.381}^{5,120}\boxtimes 3_{6,3.}^{3,138} $

Not pseudo-unitary. 

\vskip 1ex 
\color{grey}

\noindent(4,8). $6_{\frac{28}{5},4.145}^{15,623}$ \irep{38}:\ \ 
$d_i$ = ($1.0$,
$1.0$,
$1.0$,
$-0.618$,
$-0.618$,
$-0.618$) 

\vskip 0.7ex
\hangindent=3em \hangafter=1
$D^2= 4.145 = 
\frac{15-3\sqrt{5}}{2}$

\vskip 0.7ex
\hangindent=3em \hangafter=1
$T = ( 0,
\frac{2}{3},
\frac{2}{3},
\frac{4}{5},
\frac{7}{15},
\frac{7}{15} )
$,

\vskip 0.7ex
\hangindent=3em \hangafter=1
$S$ = ($ 1$,
$ 1$,
$ 1$,
$ \frac{1-\sqrt{5}}{2}$,
$ \frac{1-\sqrt{5}}{2}$,
$ \frac{1-\sqrt{5}}{2}$;\ \ 
$ -\zeta_{6}^{1}$,
$ \zeta_{3}^{1}$,
$ \frac{1-\sqrt{5}}{2}$,
$ \frac{1-\sqrt{5}}{2}\zeta_{3}^{1}$,
$ -\frac{1-\sqrt{5}}{2}\zeta_{6}^{1}$;\ \ 
$ -\zeta_{6}^{1}$,
$ \frac{1-\sqrt{5}}{2}$,
$ -\frac{1-\sqrt{5}}{2}\zeta_{6}^{1}$,
$ \frac{1-\sqrt{5}}{2}\zeta_{3}^{1}$;\ \ 
$ -1$,
$ -1$,
$ -1$;\ \ 
$ \zeta_{6}^{1}$,
$ -\zeta_{3}^{1}$;\ \ 
$ \zeta_{6}^{1}$)

Factors = $2_{\frac{38}{5},1.381}^{5,491}\boxtimes 3_{6,3.}^{3,138} $

Not pseudo-unitary. 

\vskip 1ex 
\black

\noindent(5,1). $6_{\frac{33}{10},14.47}^{80,798}$ \irep{48}:\ \ 
$d_i$ = ($1.0$,
$1.0$,
$1.414$,
$1.618$,
$1.618$,
$2.288$) 

\vskip 0.7ex
\hangindent=3em \hangafter=1
$D^2= 14.472 = 
10+2\sqrt{5}$

\vskip 0.7ex
\hangindent=3em \hangafter=1
$T = ( 0,
\frac{1}{2},
\frac{1}{16},
\frac{2}{5},
\frac{9}{10},
\frac{37}{80} )
$,

\vskip 0.7ex
\hangindent=3em \hangafter=1
$S$ = ($ 1$,
$ 1$,
$ \sqrt{2}$,
$ \frac{1+\sqrt{5}}{2}$,
$ \frac{1+\sqrt{5}}{2}$,
$ \frac{5+\sqrt{5}}{\sqrt{10}}$;\ \ 
$ 1$,
$ -\sqrt{2}$,
$ \frac{1+\sqrt{5}}{2}$,
$ \frac{1+\sqrt{5}}{2}$,
$ \frac{-5-\sqrt{5}}{\sqrt{10}}$;\ \ 
$0$,
$ \frac{5+\sqrt{5}}{\sqrt{10}}$,
$ \frac{-5-\sqrt{5}}{\sqrt{10}}$,
$0$;\ \ 
$ -1$,
$ -1$,
$ -\sqrt{2}$;\ \ 
$ -1$,
$ \sqrt{2}$;\ \ 
$0$)

Factors = $2_{\frac{14}{5},3.618}^{5,395}\boxtimes 3_{\frac{1}{2},4.}^{16,598} $

\vskip 1ex 
\color{grey}

\noindent(5,2). $6_{\frac{57}{10},14.47}^{80,376}$ \irep{48}:\ \ 
$d_i$ = ($1.0$,
$1.0$,
$1.414$,
$1.618$,
$1.618$,
$2.288$) 

\vskip 0.7ex
\hangindent=3em \hangafter=1
$D^2= 14.472 = 
10+2\sqrt{5}$

\vskip 0.7ex
\hangindent=3em \hangafter=1
$T = ( 0,
\frac{1}{2},
\frac{1}{16},
\frac{3}{5},
\frac{1}{10},
\frac{53}{80} )
$,

\vskip 0.7ex
\hangindent=3em \hangafter=1
$S$ = ($ 1$,
$ 1$,
$ \sqrt{2}$,
$ \frac{1+\sqrt{5}}{2}$,
$ \frac{1+\sqrt{5}}{2}$,
$ \frac{5+\sqrt{5}}{\sqrt{10}}$;\ \ 
$ 1$,
$ -\sqrt{2}$,
$ \frac{1+\sqrt{5}}{2}$,
$ \frac{1+\sqrt{5}}{2}$,
$ \frac{-5-\sqrt{5}}{\sqrt{10}}$;\ \ 
$0$,
$ \frac{5+\sqrt{5}}{\sqrt{10}}$,
$ \frac{-5-\sqrt{5}}{\sqrt{10}}$,
$0$;\ \ 
$ -1$,
$ -1$,
$ -\sqrt{2}$;\ \ 
$ -1$,
$ \sqrt{2}$;\ \ 
$0$)

Factors = $2_{\frac{26}{5},3.618}^{5,720}\boxtimes 3_{\frac{1}{2},4.}^{16,598} $

\vskip 1ex 
\color{grey}

\noindent(5,3). $6_{\frac{63}{10},14.47}^{80,146}$ \irep{48}:\ \ 
$d_i$ = ($1.0$,
$1.0$,
$1.414$,
$1.618$,
$1.618$,
$2.288$) 

\vskip 0.7ex
\hangindent=3em \hangafter=1
$D^2= 14.472 = 
10+2\sqrt{5}$

\vskip 0.7ex
\hangindent=3em \hangafter=1
$T = ( 0,
\frac{1}{2},
\frac{7}{16},
\frac{2}{5},
\frac{9}{10},
\frac{67}{80} )
$,

\vskip 0.7ex
\hangindent=3em \hangafter=1
$S$ = ($ 1$,
$ 1$,
$ \sqrt{2}$,
$ \frac{1+\sqrt{5}}{2}$,
$ \frac{1+\sqrt{5}}{2}$,
$ \frac{5+\sqrt{5}}{\sqrt{10}}$;\ \ 
$ 1$,
$ -\sqrt{2}$,
$ \frac{1+\sqrt{5}}{2}$,
$ \frac{1+\sqrt{5}}{2}$,
$ \frac{-5-\sqrt{5}}{\sqrt{10}}$;\ \ 
$0$,
$ \frac{5+\sqrt{5}}{\sqrt{10}}$,
$ \frac{-5-\sqrt{5}}{\sqrt{10}}$,
$0$;\ \ 
$ -1$,
$ -1$,
$ -\sqrt{2}$;\ \ 
$ -1$,
$ \sqrt{2}$;\ \ 
$0$)

Factors = $2_{\frac{14}{5},3.618}^{5,395}\boxtimes 3_{\frac{7}{2},4.}^{16,332} $

\vskip 1ex 
\color{grey}

\noindent(5,4). $6_{\frac{7}{10},14.47}^{80,111}$ \irep{48}:\ \ 
$d_i$ = ($1.0$,
$1.0$,
$1.414$,
$1.618$,
$1.618$,
$2.288$) 

\vskip 0.7ex
\hangindent=3em \hangafter=1
$D^2= 14.472 = 
10+2\sqrt{5}$

\vskip 0.7ex
\hangindent=3em \hangafter=1
$T = ( 0,
\frac{1}{2},
\frac{7}{16},
\frac{3}{5},
\frac{1}{10},
\frac{3}{80} )
$,

\vskip 0.7ex
\hangindent=3em \hangafter=1
$S$ = ($ 1$,
$ 1$,
$ \sqrt{2}$,
$ \frac{1+\sqrt{5}}{2}$,
$ \frac{1+\sqrt{5}}{2}$,
$ \frac{5+\sqrt{5}}{\sqrt{10}}$;\ \ 
$ 1$,
$ -\sqrt{2}$,
$ \frac{1+\sqrt{5}}{2}$,
$ \frac{1+\sqrt{5}}{2}$,
$ \frac{-5-\sqrt{5}}{\sqrt{10}}$;\ \ 
$0$,
$ \frac{5+\sqrt{5}}{\sqrt{10}}$,
$ \frac{-5-\sqrt{5}}{\sqrt{10}}$,
$0$;\ \ 
$ -1$,
$ -1$,
$ -\sqrt{2}$;\ \ 
$ -1$,
$ \sqrt{2}$;\ \ 
$0$)

Factors = $2_{\frac{26}{5},3.618}^{5,720}\boxtimes 3_{\frac{7}{2},4.}^{16,332} $

\vskip 1ex 
\color{grey}

\noindent(5,5). $6_{\frac{73}{10},14.47}^{80,215}$ \irep{48}:\ \ 
$d_i$ = ($1.0$,
$1.0$,
$1.414$,
$1.618$,
$1.618$,
$2.288$) 

\vskip 0.7ex
\hangindent=3em \hangafter=1
$D^2= 14.472 = 
10+2\sqrt{5}$

\vskip 0.7ex
\hangindent=3em \hangafter=1
$T = ( 0,
\frac{1}{2},
\frac{9}{16},
\frac{2}{5},
\frac{9}{10},
\frac{77}{80} )
$,

\vskip 0.7ex
\hangindent=3em \hangafter=1
$S$ = ($ 1$,
$ 1$,
$ \sqrt{2}$,
$ \frac{1+\sqrt{5}}{2}$,
$ \frac{1+\sqrt{5}}{2}$,
$ \frac{5+\sqrt{5}}{\sqrt{10}}$;\ \ 
$ 1$,
$ -\sqrt{2}$,
$ \frac{1+\sqrt{5}}{2}$,
$ \frac{1+\sqrt{5}}{2}$,
$ \frac{-5-\sqrt{5}}{\sqrt{10}}$;\ \ 
$0$,
$ \frac{5+\sqrt{5}}{\sqrt{10}}$,
$ \frac{-5-\sqrt{5}}{\sqrt{10}}$,
$0$;\ \ 
$ -1$,
$ -1$,
$ -\sqrt{2}$;\ \ 
$ -1$,
$ \sqrt{2}$;\ \ 
$0$)

Factors = $2_{\frac{14}{5},3.618}^{5,395}\boxtimes 3_{\frac{9}{2},4.}^{16,156} $

\vskip 1ex 
\color{grey}

\noindent(5,6). $6_{\frac{17}{10},14.47}^{80,878}$ \irep{48}:\ \ 
$d_i$ = ($1.0$,
$1.0$,
$1.414$,
$1.618$,
$1.618$,
$2.288$) 

\vskip 0.7ex
\hangindent=3em \hangafter=1
$D^2= 14.472 = 
10+2\sqrt{5}$

\vskip 0.7ex
\hangindent=3em \hangafter=1
$T = ( 0,
\frac{1}{2},
\frac{9}{16},
\frac{3}{5},
\frac{1}{10},
\frac{13}{80} )
$,

\vskip 0.7ex
\hangindent=3em \hangafter=1
$S$ = ($ 1$,
$ 1$,
$ \sqrt{2}$,
$ \frac{1+\sqrt{5}}{2}$,
$ \frac{1+\sqrt{5}}{2}$,
$ \frac{5+\sqrt{5}}{\sqrt{10}}$;\ \ 
$ 1$,
$ -\sqrt{2}$,
$ \frac{1+\sqrt{5}}{2}$,
$ \frac{1+\sqrt{5}}{2}$,
$ \frac{-5-\sqrt{5}}{\sqrt{10}}$;\ \ 
$0$,
$ \frac{5+\sqrt{5}}{\sqrt{10}}$,
$ \frac{-5-\sqrt{5}}{\sqrt{10}}$,
$0$;\ \ 
$ -1$,
$ -1$,
$ -\sqrt{2}$;\ \ 
$ -1$,
$ \sqrt{2}$;\ \ 
$0$)

Factors = $2_{\frac{26}{5},3.618}^{5,720}\boxtimes 3_{\frac{9}{2},4.}^{16,156} $

\vskip 1ex 
\color{grey}

\noindent(5,7). $6_{\frac{23}{10},14.47}^{80,108}$ \irep{48}:\ \ 
$d_i$ = ($1.0$,
$1.0$,
$1.414$,
$1.618$,
$1.618$,
$2.288$) 

\vskip 0.7ex
\hangindent=3em \hangafter=1
$D^2= 14.472 = 
10+2\sqrt{5}$

\vskip 0.7ex
\hangindent=3em \hangafter=1
$T = ( 0,
\frac{1}{2},
\frac{15}{16},
\frac{2}{5},
\frac{9}{10},
\frac{27}{80} )
$,

\vskip 0.7ex
\hangindent=3em \hangafter=1
$S$ = ($ 1$,
$ 1$,
$ \sqrt{2}$,
$ \frac{1+\sqrt{5}}{2}$,
$ \frac{1+\sqrt{5}}{2}$,
$ \frac{5+\sqrt{5}}{\sqrt{10}}$;\ \ 
$ 1$,
$ -\sqrt{2}$,
$ \frac{1+\sqrt{5}}{2}$,
$ \frac{1+\sqrt{5}}{2}$,
$ \frac{-5-\sqrt{5}}{\sqrt{10}}$;\ \ 
$0$,
$ \frac{5+\sqrt{5}}{\sqrt{10}}$,
$ \frac{-5-\sqrt{5}}{\sqrt{10}}$,
$0$;\ \ 
$ -1$,
$ -1$,
$ -\sqrt{2}$;\ \ 
$ -1$,
$ \sqrt{2}$;\ \ 
$0$)

Factors = $2_{\frac{14}{5},3.618}^{5,395}\boxtimes 3_{\frac{15}{2},4.}^{16,639} $

\vskip 1ex 
\color{grey}

\noindent(5,8). $6_{\frac{47}{10},14.47}^{80,518}$ \irep{48}:\ \ 
$d_i$ = ($1.0$,
$1.0$,
$1.414$,
$1.618$,
$1.618$,
$2.288$) 

\vskip 0.7ex
\hangindent=3em \hangafter=1
$D^2= 14.472 = 
10+2\sqrt{5}$

\vskip 0.7ex
\hangindent=3em \hangafter=1
$T = ( 0,
\frac{1}{2},
\frac{15}{16},
\frac{3}{5},
\frac{1}{10},
\frac{43}{80} )
$,

\vskip 0.7ex
\hangindent=3em \hangafter=1
$S$ = ($ 1$,
$ 1$,
$ \sqrt{2}$,
$ \frac{1+\sqrt{5}}{2}$,
$ \frac{1+\sqrt{5}}{2}$,
$ \frac{5+\sqrt{5}}{\sqrt{10}}$;\ \ 
$ 1$,
$ -\sqrt{2}$,
$ \frac{1+\sqrt{5}}{2}$,
$ \frac{1+\sqrt{5}}{2}$,
$ \frac{-5-\sqrt{5}}{\sqrt{10}}$;\ \ 
$0$,
$ \frac{5+\sqrt{5}}{\sqrt{10}}$,
$ \frac{-5-\sqrt{5}}{\sqrt{10}}$,
$0$;\ \ 
$ -1$,
$ -1$,
$ -\sqrt{2}$;\ \ 
$ -1$,
$ \sqrt{2}$;\ \ 
$0$)

Factors = $2_{\frac{26}{5},3.618}^{5,720}\boxtimes 3_{\frac{15}{2},4.}^{16,639} $

\vskip 1ex 
\color{grey}

\noindent(5,9). $6_{\frac{43}{10},14.47}^{80,596}$ \irep{48}:\ \ 
$d_i$ = ($1.0$,
$1.0$,
$1.618$,
$1.618$,
$-1.414$,
$-2.288$) 

\vskip 0.7ex
\hangindent=3em \hangafter=1
$D^2= 14.472 = 
10+2\sqrt{5}$

\vskip 0.7ex
\hangindent=3em \hangafter=1
$T = ( 0,
\frac{1}{2},
\frac{2}{5},
\frac{9}{10},
\frac{3}{16},
\frac{47}{80} )
$,

\vskip 0.7ex
\hangindent=3em \hangafter=1
$S$ = ($ 1$,
$ 1$,
$ \frac{1+\sqrt{5}}{2}$,
$ \frac{1+\sqrt{5}}{2}$,
$ -\sqrt{2}$,
$ \frac{-5-\sqrt{5}}{\sqrt{10}}$;\ \ 
$ 1$,
$ \frac{1+\sqrt{5}}{2}$,
$ \frac{1+\sqrt{5}}{2}$,
$ \sqrt{2}$,
$ \frac{5+\sqrt{5}}{\sqrt{10}}$;\ \ 
$ -1$,
$ -1$,
$ \frac{-5-\sqrt{5}}{\sqrt{10}}$,
$ \sqrt{2}$;\ \ 
$ -1$,
$ \frac{5+\sqrt{5}}{\sqrt{10}}$,
$ -\sqrt{2}$;\ \ 
$0$,
$0$;\ \ 
$0$)

Factors = $2_{\frac{14}{5},3.618}^{5,395}\boxtimes 3_{\frac{3}{2},4.}^{16,538} $

Pseudo-unitary $\sim$  
$6_{\frac{3}{10},14.47}^{80,270}$

\vskip 1ex 
\color{grey}

\noindent(5,10). $6_{\frac{53}{10},14.47}^{80,136}$ \irep{48}:\ \ 
$d_i$ = ($1.0$,
$1.0$,
$1.618$,
$1.618$,
$-1.414$,
$-2.288$) 

\vskip 0.7ex
\hangindent=3em \hangafter=1
$D^2= 14.472 = 
10+2\sqrt{5}$

\vskip 0.7ex
\hangindent=3em \hangafter=1
$T = ( 0,
\frac{1}{2},
\frac{2}{5},
\frac{9}{10},
\frac{5}{16},
\frac{57}{80} )
$,

\vskip 0.7ex
\hangindent=3em \hangafter=1
$S$ = ($ 1$,
$ 1$,
$ \frac{1+\sqrt{5}}{2}$,
$ \frac{1+\sqrt{5}}{2}$,
$ -\sqrt{2}$,
$ \frac{-5-\sqrt{5}}{\sqrt{10}}$;\ \ 
$ 1$,
$ \frac{1+\sqrt{5}}{2}$,
$ \frac{1+\sqrt{5}}{2}$,
$ \sqrt{2}$,
$ \frac{5+\sqrt{5}}{\sqrt{10}}$;\ \ 
$ -1$,
$ -1$,
$ \frac{-5-\sqrt{5}}{\sqrt{10}}$,
$ \sqrt{2}$;\ \ 
$ -1$,
$ \frac{5+\sqrt{5}}{\sqrt{10}}$,
$ -\sqrt{2}$;\ \ 
$0$,
$0$;\ \ 
$0$)

Factors = $2_{\frac{14}{5},3.618}^{5,395}\boxtimes 3_{\frac{5}{2},4.}^{16,345} $

Pseudo-unitary $\sim$  
$6_{\frac{13}{10},14.47}^{80,621}$

\vskip 1ex 
\color{grey}

\noindent(5,11). $6_{\frac{3}{10},14.47}^{80,750}$ \irep{48}:\ \ 
$d_i$ = ($1.0$,
$1.0$,
$1.618$,
$1.618$,
$-1.414$,
$-2.288$) 

\vskip 0.7ex
\hangindent=3em \hangafter=1
$D^2= 14.472 = 
10+2\sqrt{5}$

\vskip 0.7ex
\hangindent=3em \hangafter=1
$T = ( 0,
\frac{1}{2},
\frac{2}{5},
\frac{9}{10},
\frac{11}{16},
\frac{7}{80} )
$,

\vskip 0.7ex
\hangindent=3em \hangafter=1
$S$ = ($ 1$,
$ 1$,
$ \frac{1+\sqrt{5}}{2}$,
$ \frac{1+\sqrt{5}}{2}$,
$ -\sqrt{2}$,
$ \frac{-5-\sqrt{5}}{\sqrt{10}}$;\ \ 
$ 1$,
$ \frac{1+\sqrt{5}}{2}$,
$ \frac{1+\sqrt{5}}{2}$,
$ \sqrt{2}$,
$ \frac{5+\sqrt{5}}{\sqrt{10}}$;\ \ 
$ -1$,
$ -1$,
$ \frac{-5-\sqrt{5}}{\sqrt{10}}$,
$ \sqrt{2}$;\ \ 
$ -1$,
$ \frac{5+\sqrt{5}}{\sqrt{10}}$,
$ -\sqrt{2}$;\ \ 
$0$,
$0$;\ \ 
$0$)

Factors = $2_{\frac{14}{5},3.618}^{5,395}\boxtimes 3_{\frac{11}{2},4.}^{16,564} $

Pseudo-unitary $\sim$  
$6_{\frac{43}{10},14.47}^{80,424}$

\vskip 1ex 
\color{grey}

\noindent(5,12). $6_{\frac{13}{10},14.47}^{80,399}$ \irep{48}:\ \ 
$d_i$ = ($1.0$,
$1.0$,
$1.618$,
$1.618$,
$-1.414$,
$-2.288$) 

\vskip 0.7ex
\hangindent=3em \hangafter=1
$D^2= 14.472 = 
10+2\sqrt{5}$

\vskip 0.7ex
\hangindent=3em \hangafter=1
$T = ( 0,
\frac{1}{2},
\frac{2}{5},
\frac{9}{10},
\frac{13}{16},
\frac{17}{80} )
$,

\vskip 0.7ex
\hangindent=3em \hangafter=1
$S$ = ($ 1$,
$ 1$,
$ \frac{1+\sqrt{5}}{2}$,
$ \frac{1+\sqrt{5}}{2}$,
$ -\sqrt{2}$,
$ \frac{-5-\sqrt{5}}{\sqrt{10}}$;\ \ 
$ 1$,
$ \frac{1+\sqrt{5}}{2}$,
$ \frac{1+\sqrt{5}}{2}$,
$ \sqrt{2}$,
$ \frac{5+\sqrt{5}}{\sqrt{10}}$;\ \ 
$ -1$,
$ -1$,
$ \frac{-5-\sqrt{5}}{\sqrt{10}}$,
$ \sqrt{2}$;\ \ 
$ -1$,
$ \frac{5+\sqrt{5}}{\sqrt{10}}$,
$ -\sqrt{2}$;\ \ 
$0$,
$0$;\ \ 
$0$)

Factors = $2_{\frac{14}{5},3.618}^{5,395}\boxtimes 3_{\frac{13}{2},4.}^{16,830} $

Pseudo-unitary $\sim$  
$6_{\frac{53}{10},14.47}^{80,884}$

\vskip 1ex 
\color{grey}

\noindent(5,13). $6_{\frac{67}{10},14.47}^{80,100}$ \irep{48}:\ \ 
$d_i$ = ($1.0$,
$1.0$,
$1.618$,
$1.618$,
$-1.414$,
$-2.288$) 

\vskip 0.7ex
\hangindent=3em \hangafter=1
$D^2= 14.472 = 
10+2\sqrt{5}$

\vskip 0.7ex
\hangindent=3em \hangafter=1
$T = ( 0,
\frac{1}{2},
\frac{3}{5},
\frac{1}{10},
\frac{3}{16},
\frac{63}{80} )
$,

\vskip 0.7ex
\hangindent=3em \hangafter=1
$S$ = ($ 1$,
$ 1$,
$ \frac{1+\sqrt{5}}{2}$,
$ \frac{1+\sqrt{5}}{2}$,
$ -\sqrt{2}$,
$ \frac{-5-\sqrt{5}}{\sqrt{10}}$;\ \ 
$ 1$,
$ \frac{1+\sqrt{5}}{2}$,
$ \frac{1+\sqrt{5}}{2}$,
$ \sqrt{2}$,
$ \frac{5+\sqrt{5}}{\sqrt{10}}$;\ \ 
$ -1$,
$ -1$,
$ \frac{-5-\sqrt{5}}{\sqrt{10}}$,
$ \sqrt{2}$;\ \ 
$ -1$,
$ \frac{5+\sqrt{5}}{\sqrt{10}}$,
$ -\sqrt{2}$;\ \ 
$0$,
$0$;\ \ 
$0$)

Factors = $2_{\frac{26}{5},3.618}^{5,720}\boxtimes 3_{\frac{3}{2},4.}^{16,538} $

Pseudo-unitary $\sim$  
$6_{\frac{27}{10},14.47}^{80,528}$

\vskip 1ex 
\color{grey}

\noindent(5,14). $6_{\frac{77}{10},14.47}^{80,157}$ \irep{48}:\ \ 
$d_i$ = ($1.0$,
$1.0$,
$1.618$,
$1.618$,
$-1.414$,
$-2.288$) 

\vskip 0.7ex
\hangindent=3em \hangafter=1
$D^2= 14.472 = 
10+2\sqrt{5}$

\vskip 0.7ex
\hangindent=3em \hangafter=1
$T = ( 0,
\frac{1}{2},
\frac{3}{5},
\frac{1}{10},
\frac{5}{16},
\frac{73}{80} )
$,

\vskip 0.7ex
\hangindent=3em \hangafter=1
$S$ = ($ 1$,
$ 1$,
$ \frac{1+\sqrt{5}}{2}$,
$ \frac{1+\sqrt{5}}{2}$,
$ -\sqrt{2}$,
$ \frac{-5-\sqrt{5}}{\sqrt{10}}$;\ \ 
$ 1$,
$ \frac{1+\sqrt{5}}{2}$,
$ \frac{1+\sqrt{5}}{2}$,
$ \sqrt{2}$,
$ \frac{5+\sqrt{5}}{\sqrt{10}}$;\ \ 
$ -1$,
$ -1$,
$ \frac{-5-\sqrt{5}}{\sqrt{10}}$,
$ \sqrt{2}$;\ \ 
$ -1$,
$ \frac{5+\sqrt{5}}{\sqrt{10}}$,
$ -\sqrt{2}$;\ \ 
$0$,
$0$;\ \ 
$0$)

Factors = $2_{\frac{26}{5},3.618}^{5,720}\boxtimes 3_{\frac{5}{2},4.}^{16,345} $

Pseudo-unitary $\sim$  
$6_{\frac{37}{10},14.47}^{80,629}$

\vskip 1ex 
\color{grey}

\noindent(5,15). $6_{\frac{27}{10},14.47}^{80,392}$ \irep{48}:\ \ 
$d_i$ = ($1.0$,
$1.0$,
$1.618$,
$1.618$,
$-1.414$,
$-2.288$) 

\vskip 0.7ex
\hangindent=3em \hangafter=1
$D^2= 14.472 = 
10+2\sqrt{5}$

\vskip 0.7ex
\hangindent=3em \hangafter=1
$T = ( 0,
\frac{1}{2},
\frac{3}{5},
\frac{1}{10},
\frac{11}{16},
\frac{23}{80} )
$,

\vskip 0.7ex
\hangindent=3em \hangafter=1
$S$ = ($ 1$,
$ 1$,
$ \frac{1+\sqrt{5}}{2}$,
$ \frac{1+\sqrt{5}}{2}$,
$ -\sqrt{2}$,
$ \frac{-5-\sqrt{5}}{\sqrt{10}}$;\ \ 
$ 1$,
$ \frac{1+\sqrt{5}}{2}$,
$ \frac{1+\sqrt{5}}{2}$,
$ \sqrt{2}$,
$ \frac{5+\sqrt{5}}{\sqrt{10}}$;\ \ 
$ -1$,
$ -1$,
$ \frac{-5-\sqrt{5}}{\sqrt{10}}$,
$ \sqrt{2}$;\ \ 
$ -1$,
$ \frac{5+\sqrt{5}}{\sqrt{10}}$,
$ -\sqrt{2}$;\ \ 
$0$,
$0$;\ \ 
$0$)

Factors = $2_{\frac{26}{5},3.618}^{5,720}\boxtimes 3_{\frac{11}{2},4.}^{16,564} $

Pseudo-unitary $\sim$  
$6_{\frac{67}{10},14.47}^{80,828}$

\vskip 1ex 
\color{grey}

\noindent(5,16). $6_{\frac{37}{10},14.47}^{80,857}$ \irep{48}:\ \ 
$d_i$ = ($1.0$,
$1.0$,
$1.618$,
$1.618$,
$-1.414$,
$-2.288$) 

\vskip 0.7ex
\hangindent=3em \hangafter=1
$D^2= 14.472 = 
10+2\sqrt{5}$

\vskip 0.7ex
\hangindent=3em \hangafter=1
$T = ( 0,
\frac{1}{2},
\frac{3}{5},
\frac{1}{10},
\frac{13}{16},
\frac{33}{80} )
$,

\vskip 0.7ex
\hangindent=3em \hangafter=1
$S$ = ($ 1$,
$ 1$,
$ \frac{1+\sqrt{5}}{2}$,
$ \frac{1+\sqrt{5}}{2}$,
$ -\sqrt{2}$,
$ \frac{-5-\sqrt{5}}{\sqrt{10}}$;\ \ 
$ 1$,
$ \frac{1+\sqrt{5}}{2}$,
$ \frac{1+\sqrt{5}}{2}$,
$ \sqrt{2}$,
$ \frac{5+\sqrt{5}}{\sqrt{10}}$;\ \ 
$ -1$,
$ -1$,
$ \frac{-5-\sqrt{5}}{\sqrt{10}}$,
$ \sqrt{2}$;\ \ 
$ -1$,
$ \frac{5+\sqrt{5}}{\sqrt{10}}$,
$ -\sqrt{2}$;\ \ 
$0$,
$0$;\ \ 
$0$)

Factors = $2_{\frac{26}{5},3.618}^{5,720}\boxtimes 3_{\frac{13}{2},4.}^{16,830} $

Pseudo-unitary $\sim$  
$6_{\frac{77}{10},14.47}^{80,657}$

\vskip 1ex 
\color{grey}

\noindent(5,17). $6_{\frac{69}{10},5.527}^{80,106}$ \irep{48}:\ \ 
$d_i$ = ($1.0$,
$0.874$,
$1.0$,
$-0.618$,
$-0.618$,
$-1.414$) 

\vskip 0.7ex
\hangindent=3em \hangafter=1
$D^2= 5.527 = 
10-2\sqrt{5}$

\vskip 0.7ex
\hangindent=3em \hangafter=1
$T = ( 0,
\frac{1}{80},
\frac{1}{2},
\frac{1}{5},
\frac{7}{10},
\frac{13}{16} )
$,

\vskip 0.7ex
\hangindent=3em \hangafter=1
$S$ = ($ 1$,
$ \frac{5-\sqrt{5}}{\sqrt{10}}$,
$ 1$,
$ \frac{1-\sqrt{5}}{2}$,
$ \frac{1-\sqrt{5}}{2}$,
$ -\sqrt{2}$;\ \ 
$0$,
$ \frac{-5+\sqrt{5}}{\sqrt{10}}$,
$ \sqrt{2}$,
$ -\sqrt{2}$,
$0$;\ \ 
$ 1$,
$ \frac{1-\sqrt{5}}{2}$,
$ \frac{1-\sqrt{5}}{2}$,
$ \sqrt{2}$;\ \ 
$ -1$,
$ -1$,
$ \frac{5-\sqrt{5}}{\sqrt{10}}$;\ \ 
$ -1$,
$ \frac{-5+\sqrt{5}}{\sqrt{10}}$;\ \ 
$0$)

Factors = $2_{\frac{2}{5},1.381}^{5,120}\boxtimes 3_{\frac{13}{2},4.}^{16,830} $

Not pseudo-unitary. 

\vskip 1ex 
\color{grey}

\noindent(5,18). $6_{\frac{21}{10},5.527}^{80,384}$ \irep{48}:\ \ 
$d_i$ = ($1.0$,
$0.874$,
$1.0$,
$-0.618$,
$-0.618$,
$-1.414$) 

\vskip 0.7ex
\hangindent=3em \hangafter=1
$D^2= 5.527 = 
10-2\sqrt{5}$

\vskip 0.7ex
\hangindent=3em \hangafter=1
$T = ( 0,
\frac{9}{80},
\frac{1}{2},
\frac{4}{5},
\frac{3}{10},
\frac{5}{16} )
$,

\vskip 0.7ex
\hangindent=3em \hangafter=1
$S$ = ($ 1$,
$ \frac{5-\sqrt{5}}{\sqrt{10}}$,
$ 1$,
$ \frac{1-\sqrt{5}}{2}$,
$ \frac{1-\sqrt{5}}{2}$,
$ -\sqrt{2}$;\ \ 
$0$,
$ \frac{-5+\sqrt{5}}{\sqrt{10}}$,
$ \sqrt{2}$,
$ -\sqrt{2}$,
$0$;\ \ 
$ 1$,
$ \frac{1-\sqrt{5}}{2}$,
$ \frac{1-\sqrt{5}}{2}$,
$ \sqrt{2}$;\ \ 
$ -1$,
$ -1$,
$ \frac{5-\sqrt{5}}{\sqrt{10}}$;\ \ 
$ -1$,
$ \frac{-5+\sqrt{5}}{\sqrt{10}}$;\ \ 
$0$)

Factors = $2_{\frac{38}{5},1.381}^{5,491}\boxtimes 3_{\frac{5}{2},4.}^{16,345} $

Not pseudo-unitary. 

\vskip 1ex 
\color{grey}

\noindent(5,19). $6_{\frac{19}{10},5.527}^{80,519}$ \irep{48}:\ \ 
$d_i$ = ($1.0$,
$0.874$,
$1.0$,
$-0.618$,
$-0.618$,
$-1.414$) 

\vskip 0.7ex
\hangindent=3em \hangafter=1
$D^2= 5.527 = 
10-2\sqrt{5}$

\vskip 0.7ex
\hangindent=3em \hangafter=1
$T = ( 0,
\frac{31}{80},
\frac{1}{2},
\frac{1}{5},
\frac{7}{10},
\frac{3}{16} )
$,

\vskip 0.7ex
\hangindent=3em \hangafter=1
$S$ = ($ 1$,
$ \frac{5-\sqrt{5}}{\sqrt{10}}$,
$ 1$,
$ \frac{1-\sqrt{5}}{2}$,
$ \frac{1-\sqrt{5}}{2}$,
$ -\sqrt{2}$;\ \ 
$0$,
$ \frac{-5+\sqrt{5}}{\sqrt{10}}$,
$ \sqrt{2}$,
$ -\sqrt{2}$,
$0$;\ \ 
$ 1$,
$ \frac{1-\sqrt{5}}{2}$,
$ \frac{1-\sqrt{5}}{2}$,
$ \sqrt{2}$;\ \ 
$ -1$,
$ -1$,
$ \frac{5-\sqrt{5}}{\sqrt{10}}$;\ \ 
$ -1$,
$ \frac{-5+\sqrt{5}}{\sqrt{10}}$;\ \ 
$0$)

Factors = $2_{\frac{2}{5},1.381}^{5,120}\boxtimes 3_{\frac{3}{2},4.}^{16,538} $

Not pseudo-unitary. 

\vskip 1ex 
\color{grey}

\noindent(5,20). $6_{\frac{51}{10},5.527}^{80,717}$ \irep{48}:\ \ 
$d_i$ = ($1.0$,
$0.874$,
$1.0$,
$-0.618$,
$-0.618$,
$-1.414$) 

\vskip 0.7ex
\hangindent=3em \hangafter=1
$D^2= 5.527 = 
10-2\sqrt{5}$

\vskip 0.7ex
\hangindent=3em \hangafter=1
$T = ( 0,
\frac{39}{80},
\frac{1}{2},
\frac{4}{5},
\frac{3}{10},
\frac{11}{16} )
$,

\vskip 0.7ex
\hangindent=3em \hangafter=1
$S$ = ($ 1$,
$ \frac{5-\sqrt{5}}{\sqrt{10}}$,
$ 1$,
$ \frac{1-\sqrt{5}}{2}$,
$ \frac{1-\sqrt{5}}{2}$,
$ -\sqrt{2}$;\ \ 
$0$,
$ \frac{-5+\sqrt{5}}{\sqrt{10}}$,
$ \sqrt{2}$,
$ -\sqrt{2}$,
$0$;\ \ 
$ 1$,
$ \frac{1-\sqrt{5}}{2}$,
$ \frac{1-\sqrt{5}}{2}$,
$ \sqrt{2}$;\ \ 
$ -1$,
$ -1$,
$ \frac{5-\sqrt{5}}{\sqrt{10}}$;\ \ 
$ -1$,
$ \frac{-5+\sqrt{5}}{\sqrt{10}}$;\ \ 
$0$)

Factors = $2_{\frac{38}{5},1.381}^{5,491}\boxtimes 3_{\frac{11}{2},4.}^{16,564} $

Not pseudo-unitary. 

\vskip 1ex 
\color{grey}

\noindent(5,21). $6_{\frac{9}{10},5.527}^{80,458}$ \irep{48}:\ \ 
$d_i$ = ($1.0$,
$1.0$,
$1.414$,
$-0.618$,
$-0.618$,
$-0.874$) 

\vskip 0.7ex
\hangindent=3em \hangafter=1
$D^2= 5.527 = 
10-2\sqrt{5}$

\vskip 0.7ex
\hangindent=3em \hangafter=1
$T = ( 0,
\frac{1}{2},
\frac{1}{16},
\frac{1}{5},
\frac{7}{10},
\frac{21}{80} )
$,

\vskip 0.7ex
\hangindent=3em \hangafter=1
$S$ = ($ 1$,
$ 1$,
$ \sqrt{2}$,
$ \frac{1-\sqrt{5}}{2}$,
$ \frac{1-\sqrt{5}}{2}$,
$ \frac{-5+\sqrt{5}}{\sqrt{10}}$;\ \ 
$ 1$,
$ -\sqrt{2}$,
$ \frac{1-\sqrt{5}}{2}$,
$ \frac{1-\sqrt{5}}{2}$,
$ \frac{5-\sqrt{5}}{\sqrt{10}}$;\ \ 
$0$,
$ \frac{-5+\sqrt{5}}{\sqrt{10}}$,
$ \frac{5-\sqrt{5}}{\sqrt{10}}$,
$0$;\ \ 
$ -1$,
$ -1$,
$ -\sqrt{2}$;\ \ 
$ -1$,
$ \sqrt{2}$;\ \ 
$0$)

Factors = $2_{\frac{2}{5},1.381}^{5,120}\boxtimes 3_{\frac{1}{2},4.}^{16,598} $

Not pseudo-unitary. 

\vskip 1ex 
\color{grey}

\noindent(5,22). $6_{\frac{1}{10},5.527}^{80,117}$ \irep{48}:\ \ 
$d_i$ = ($1.0$,
$1.0$,
$1.414$,
$-0.618$,
$-0.618$,
$-0.874$) 

\vskip 0.7ex
\hangindent=3em \hangafter=1
$D^2= 5.527 = 
10-2\sqrt{5}$

\vskip 0.7ex
\hangindent=3em \hangafter=1
$T = ( 0,
\frac{1}{2},
\frac{1}{16},
\frac{4}{5},
\frac{3}{10},
\frac{69}{80} )
$,

\vskip 0.7ex
\hangindent=3em \hangafter=1
$S$ = ($ 1$,
$ 1$,
$ \sqrt{2}$,
$ \frac{1-\sqrt{5}}{2}$,
$ \frac{1-\sqrt{5}}{2}$,
$ \frac{-5+\sqrt{5}}{\sqrt{10}}$;\ \ 
$ 1$,
$ -\sqrt{2}$,
$ \frac{1-\sqrt{5}}{2}$,
$ \frac{1-\sqrt{5}}{2}$,
$ \frac{5-\sqrt{5}}{\sqrt{10}}$;\ \ 
$0$,
$ \frac{-5+\sqrt{5}}{\sqrt{10}}$,
$ \frac{5-\sqrt{5}}{\sqrt{10}}$,
$0$;\ \ 
$ -1$,
$ -1$,
$ -\sqrt{2}$;\ \ 
$ -1$,
$ \sqrt{2}$;\ \ 
$0$)

Factors = $2_{\frac{38}{5},1.381}^{5,491}\boxtimes 3_{\frac{1}{2},4.}^{16,598} $

Not pseudo-unitary. 

\vskip 1ex 
\color{grey}

\noindent(5,23). $6_{\frac{39}{10},5.527}^{80,488}$ \irep{48}:\ \ 
$d_i$ = ($1.0$,
$1.0$,
$1.414$,
$-0.618$,
$-0.618$,
$-0.874$) 

\vskip 0.7ex
\hangindent=3em \hangafter=1
$D^2= 5.527 = 
10-2\sqrt{5}$

\vskip 0.7ex
\hangindent=3em \hangafter=1
$T = ( 0,
\frac{1}{2},
\frac{7}{16},
\frac{1}{5},
\frac{7}{10},
\frac{51}{80} )
$,

\vskip 0.7ex
\hangindent=3em \hangafter=1
$S$ = ($ 1$,
$ 1$,
$ \sqrt{2}$,
$ \frac{1-\sqrt{5}}{2}$,
$ \frac{1-\sqrt{5}}{2}$,
$ \frac{-5+\sqrt{5}}{\sqrt{10}}$;\ \ 
$ 1$,
$ -\sqrt{2}$,
$ \frac{1-\sqrt{5}}{2}$,
$ \frac{1-\sqrt{5}}{2}$,
$ \frac{5-\sqrt{5}}{\sqrt{10}}$;\ \ 
$0$,
$ \frac{-5+\sqrt{5}}{\sqrt{10}}$,
$ \frac{5-\sqrt{5}}{\sqrt{10}}$,
$0$;\ \ 
$ -1$,
$ -1$,
$ -\sqrt{2}$;\ \ 
$ -1$,
$ \sqrt{2}$;\ \ 
$0$)

Factors = $2_{\frac{2}{5},1.381}^{5,120}\boxtimes 3_{\frac{7}{2},4.}^{16,332} $

Not pseudo-unitary. 

\vskip 1ex 
\color{grey}

\noindent(5,24). $6_{\frac{31}{10},5.527}^{80,305}$ \irep{48}:\ \ 
$d_i$ = ($1.0$,
$1.0$,
$1.414$,
$-0.618$,
$-0.618$,
$-0.874$) 

\vskip 0.7ex
\hangindent=3em \hangafter=1
$D^2= 5.527 = 
10-2\sqrt{5}$

\vskip 0.7ex
\hangindent=3em \hangafter=1
$T = ( 0,
\frac{1}{2},
\frac{7}{16},
\frac{4}{5},
\frac{3}{10},
\frac{19}{80} )
$,

\vskip 0.7ex
\hangindent=3em \hangafter=1
$S$ = ($ 1$,
$ 1$,
$ \sqrt{2}$,
$ \frac{1-\sqrt{5}}{2}$,
$ \frac{1-\sqrt{5}}{2}$,
$ \frac{-5+\sqrt{5}}{\sqrt{10}}$;\ \ 
$ 1$,
$ -\sqrt{2}$,
$ \frac{1-\sqrt{5}}{2}$,
$ \frac{1-\sqrt{5}}{2}$,
$ \frac{5-\sqrt{5}}{\sqrt{10}}$;\ \ 
$0$,
$ \frac{-5+\sqrt{5}}{\sqrt{10}}$,
$ \frac{5-\sqrt{5}}{\sqrt{10}}$,
$0$;\ \ 
$ -1$,
$ -1$,
$ -\sqrt{2}$;\ \ 
$ -1$,
$ \sqrt{2}$;\ \ 
$0$)

Factors = $2_{\frac{38}{5},1.381}^{5,491}\boxtimes 3_{\frac{7}{2},4.}^{16,332} $

Not pseudo-unitary. 

\vskip 1ex 
\color{grey}

\noindent(5,25). $6_{\frac{49}{10},5.527}^{80,464}$ \irep{48}:\ \ 
$d_i$ = ($1.0$,
$1.0$,
$1.414$,
$-0.618$,
$-0.618$,
$-0.874$) 

\vskip 0.7ex
\hangindent=3em \hangafter=1
$D^2= 5.527 = 
10-2\sqrt{5}$

\vskip 0.7ex
\hangindent=3em \hangafter=1
$T = ( 0,
\frac{1}{2},
\frac{9}{16},
\frac{1}{5},
\frac{7}{10},
\frac{61}{80} )
$,

\vskip 0.7ex
\hangindent=3em \hangafter=1
$S$ = ($ 1$,
$ 1$,
$ \sqrt{2}$,
$ \frac{1-\sqrt{5}}{2}$,
$ \frac{1-\sqrt{5}}{2}$,
$ \frac{-5+\sqrt{5}}{\sqrt{10}}$;\ \ 
$ 1$,
$ -\sqrt{2}$,
$ \frac{1-\sqrt{5}}{2}$,
$ \frac{1-\sqrt{5}}{2}$,
$ \frac{5-\sqrt{5}}{\sqrt{10}}$;\ \ 
$0$,
$ \frac{-5+\sqrt{5}}{\sqrt{10}}$,
$ \frac{5-\sqrt{5}}{\sqrt{10}}$,
$0$;\ \ 
$ -1$,
$ -1$,
$ -\sqrt{2}$;\ \ 
$ -1$,
$ \sqrt{2}$;\ \ 
$0$)

Factors = $2_{\frac{2}{5},1.381}^{5,120}\boxtimes 3_{\frac{9}{2},4.}^{16,156} $

Not pseudo-unitary. 

\vskip 1ex 
\color{grey}

\noindent(5,26). $6_{\frac{41}{10},5.527}^{80,194}$ \irep{48}:\ \ 
$d_i$ = ($1.0$,
$1.0$,
$1.414$,
$-0.618$,
$-0.618$,
$-0.874$) 

\vskip 0.7ex
\hangindent=3em \hangafter=1
$D^2= 5.527 = 
10-2\sqrt{5}$

\vskip 0.7ex
\hangindent=3em \hangafter=1
$T = ( 0,
\frac{1}{2},
\frac{9}{16},
\frac{4}{5},
\frac{3}{10},
\frac{29}{80} )
$,

\vskip 0.7ex
\hangindent=3em \hangafter=1
$S$ = ($ 1$,
$ 1$,
$ \sqrt{2}$,
$ \frac{1-\sqrt{5}}{2}$,
$ \frac{1-\sqrt{5}}{2}$,
$ \frac{-5+\sqrt{5}}{\sqrt{10}}$;\ \ 
$ 1$,
$ -\sqrt{2}$,
$ \frac{1-\sqrt{5}}{2}$,
$ \frac{1-\sqrt{5}}{2}$,
$ \frac{5-\sqrt{5}}{\sqrt{10}}$;\ \ 
$0$,
$ \frac{-5+\sqrt{5}}{\sqrt{10}}$,
$ \frac{5-\sqrt{5}}{\sqrt{10}}$,
$0$;\ \ 
$ -1$,
$ -1$,
$ -\sqrt{2}$;\ \ 
$ -1$,
$ \sqrt{2}$;\ \ 
$0$)

Factors = $2_{\frac{38}{5},1.381}^{5,491}\boxtimes 3_{\frac{9}{2},4.}^{16,156} $

Not pseudo-unitary. 

\vskip 1ex 
\color{grey}

\noindent(5,27). $6_{\frac{79}{10},5.527}^{80,822}$ \irep{48}:\ \ 
$d_i$ = ($1.0$,
$1.0$,
$1.414$,
$-0.618$,
$-0.618$,
$-0.874$) 

\vskip 0.7ex
\hangindent=3em \hangafter=1
$D^2= 5.527 = 
10-2\sqrt{5}$

\vskip 0.7ex
\hangindent=3em \hangafter=1
$T = ( 0,
\frac{1}{2},
\frac{15}{16},
\frac{1}{5},
\frac{7}{10},
\frac{11}{80} )
$,

\vskip 0.7ex
\hangindent=3em \hangafter=1
$S$ = ($ 1$,
$ 1$,
$ \sqrt{2}$,
$ \frac{1-\sqrt{5}}{2}$,
$ \frac{1-\sqrt{5}}{2}$,
$ \frac{-5+\sqrt{5}}{\sqrt{10}}$;\ \ 
$ 1$,
$ -\sqrt{2}$,
$ \frac{1-\sqrt{5}}{2}$,
$ \frac{1-\sqrt{5}}{2}$,
$ \frac{5-\sqrt{5}}{\sqrt{10}}$;\ \ 
$0$,
$ \frac{-5+\sqrt{5}}{\sqrt{10}}$,
$ \frac{5-\sqrt{5}}{\sqrt{10}}$,
$0$;\ \ 
$ -1$,
$ -1$,
$ -\sqrt{2}$;\ \ 
$ -1$,
$ \sqrt{2}$;\ \ 
$0$)

Factors = $2_{\frac{2}{5},1.381}^{5,120}\boxtimes 3_{\frac{15}{2},4.}^{16,639} $

Not pseudo-unitary. 

\vskip 1ex 
\color{grey}

\noindent(5,28). $6_{\frac{71}{10},5.527}^{80,240}$ \irep{48}:\ \ 
$d_i$ = ($1.0$,
$1.0$,
$1.414$,
$-0.618$,
$-0.618$,
$-0.874$) 

\vskip 0.7ex
\hangindent=3em \hangafter=1
$D^2= 5.527 = 
10-2\sqrt{5}$

\vskip 0.7ex
\hangindent=3em \hangafter=1
$T = ( 0,
\frac{1}{2},
\frac{15}{16},
\frac{4}{5},
\frac{3}{10},
\frac{59}{80} )
$,

\vskip 0.7ex
\hangindent=3em \hangafter=1
$S$ = ($ 1$,
$ 1$,
$ \sqrt{2}$,
$ \frac{1-\sqrt{5}}{2}$,
$ \frac{1-\sqrt{5}}{2}$,
$ \frac{-5+\sqrt{5}}{\sqrt{10}}$;\ \ 
$ 1$,
$ -\sqrt{2}$,
$ \frac{1-\sqrt{5}}{2}$,
$ \frac{1-\sqrt{5}}{2}$,
$ \frac{5-\sqrt{5}}{\sqrt{10}}$;\ \ 
$0$,
$ \frac{-5+\sqrt{5}}{\sqrt{10}}$,
$ \frac{5-\sqrt{5}}{\sqrt{10}}$,
$0$;\ \ 
$ -1$,
$ -1$,
$ -\sqrt{2}$;\ \ 
$ -1$,
$ \sqrt{2}$;\ \ 
$0$)

Factors = $2_{\frac{38}{5},1.381}^{5,491}\boxtimes 3_{\frac{15}{2},4.}^{16,639} $

Not pseudo-unitary. 

\vskip 1ex 
\color{grey}

\noindent(5,29). $6_{\frac{29}{10},5.527}^{80,420}$ \irep{48}:\ \ 
$d_i$ = ($1.0$,
$0.874$,
$1.0$,
$-0.618$,
$-0.618$,
$-1.414$) 

\vskip 0.7ex
\hangindent=3em \hangafter=1
$D^2= 5.527 = 
10-2\sqrt{5}$

\vskip 0.7ex
\hangindent=3em \hangafter=1
$T = ( 0,
\frac{41}{80},
\frac{1}{2},
\frac{1}{5},
\frac{7}{10},
\frac{5}{16} )
$,

\vskip 0.7ex
\hangindent=3em \hangafter=1
$S$ = ($ 1$,
$ \frac{5-\sqrt{5}}{\sqrt{10}}$,
$ 1$,
$ \frac{1-\sqrt{5}}{2}$,
$ \frac{1-\sqrt{5}}{2}$,
$ -\sqrt{2}$;\ \ 
$0$,
$ \frac{-5+\sqrt{5}}{\sqrt{10}}$,
$ \sqrt{2}$,
$ -\sqrt{2}$,
$0$;\ \ 
$ 1$,
$ \frac{1-\sqrt{5}}{2}$,
$ \frac{1-\sqrt{5}}{2}$,
$ \sqrt{2}$;\ \ 
$ -1$,
$ -1$,
$ \frac{5-\sqrt{5}}{\sqrt{10}}$;\ \ 
$ -1$,
$ \frac{-5+\sqrt{5}}{\sqrt{10}}$;\ \ 
$0$)

Factors = $2_{\frac{2}{5},1.381}^{5,120}\boxtimes 3_{\frac{5}{2},4.}^{16,345} $

Not pseudo-unitary. 

\vskip 1ex 
\color{grey}

\noindent(5,30). $6_{\frac{61}{10},5.527}^{80,862}$ \irep{48}:\ \ 
$d_i$ = ($1.0$,
$0.874$,
$1.0$,
$-0.618$,
$-0.618$,
$-1.414$) 

\vskip 0.7ex
\hangindent=3em \hangafter=1
$D^2= 5.527 = 
10-2\sqrt{5}$

\vskip 0.7ex
\hangindent=3em \hangafter=1
$T = ( 0,
\frac{49}{80},
\frac{1}{2},
\frac{4}{5},
\frac{3}{10},
\frac{13}{16} )
$,

\vskip 0.7ex
\hangindent=3em \hangafter=1
$S$ = ($ 1$,
$ \frac{5-\sqrt{5}}{\sqrt{10}}$,
$ 1$,
$ \frac{1-\sqrt{5}}{2}$,
$ \frac{1-\sqrt{5}}{2}$,
$ -\sqrt{2}$;\ \ 
$0$,
$ \frac{-5+\sqrt{5}}{\sqrt{10}}$,
$ \sqrt{2}$,
$ -\sqrt{2}$,
$0$;\ \ 
$ 1$,
$ \frac{1-\sqrt{5}}{2}$,
$ \frac{1-\sqrt{5}}{2}$,
$ \sqrt{2}$;\ \ 
$ -1$,
$ -1$,
$ \frac{5-\sqrt{5}}{\sqrt{10}}$;\ \ 
$ -1$,
$ \frac{-5+\sqrt{5}}{\sqrt{10}}$;\ \ 
$0$)

Factors = $2_{\frac{38}{5},1.381}^{5,491}\boxtimes 3_{\frac{13}{2},4.}^{16,830} $

Not pseudo-unitary. 

\vskip 1ex 
\color{grey}

\noindent(5,31). $6_{\frac{59}{10},5.527}^{80,113}$ \irep{48}:\ \ 
$d_i$ = ($1.0$,
$0.874$,
$1.0$,
$-0.618$,
$-0.618$,
$-1.414$) 

\vskip 0.7ex
\hangindent=3em \hangafter=1
$D^2= 5.527 = 
10-2\sqrt{5}$

\vskip 0.7ex
\hangindent=3em \hangafter=1
$T = ( 0,
\frac{71}{80},
\frac{1}{2},
\frac{1}{5},
\frac{7}{10},
\frac{11}{16} )
$,

\vskip 0.7ex
\hangindent=3em \hangafter=1
$S$ = ($ 1$,
$ \frac{5-\sqrt{5}}{\sqrt{10}}$,
$ 1$,
$ \frac{1-\sqrt{5}}{2}$,
$ \frac{1-\sqrt{5}}{2}$,
$ -\sqrt{2}$;\ \ 
$0$,
$ \frac{-5+\sqrt{5}}{\sqrt{10}}$,
$ \sqrt{2}$,
$ -\sqrt{2}$,
$0$;\ \ 
$ 1$,
$ \frac{1-\sqrt{5}}{2}$,
$ \frac{1-\sqrt{5}}{2}$,
$ \sqrt{2}$;\ \ 
$ -1$,
$ -1$,
$ \frac{5-\sqrt{5}}{\sqrt{10}}$;\ \ 
$ -1$,
$ \frac{-5+\sqrt{5}}{\sqrt{10}}$;\ \ 
$0$)

Factors = $2_{\frac{2}{5},1.381}^{5,120}\boxtimes 3_{\frac{11}{2},4.}^{16,564} $

Not pseudo-unitary. 

\vskip 1ex 
\color{grey}

\noindent(5,32). $6_{\frac{11}{10},5.527}^{80,545}$ \irep{48}:\ \ 
$d_i$ = ($1.0$,
$0.874$,
$1.0$,
$-0.618$,
$-0.618$,
$-1.414$) 

\vskip 0.7ex
\hangindent=3em \hangafter=1
$D^2= 5.527 = 
10-2\sqrt{5}$

\vskip 0.7ex
\hangindent=3em \hangafter=1
$T = ( 0,
\frac{79}{80},
\frac{1}{2},
\frac{4}{5},
\frac{3}{10},
\frac{3}{16} )
$,

\vskip 0.7ex
\hangindent=3em \hangafter=1
$S$ = ($ 1$,
$ \frac{5-\sqrt{5}}{\sqrt{10}}$,
$ 1$,
$ \frac{1-\sqrt{5}}{2}$,
$ \frac{1-\sqrt{5}}{2}$,
$ -\sqrt{2}$;\ \ 
$0$,
$ \frac{-5+\sqrt{5}}{\sqrt{10}}$,
$ \sqrt{2}$,
$ -\sqrt{2}$,
$0$;\ \ 
$ 1$,
$ \frac{1-\sqrt{5}}{2}$,
$ \frac{1-\sqrt{5}}{2}$,
$ \sqrt{2}$;\ \ 
$ -1$,
$ -1$,
$ \frac{5-\sqrt{5}}{\sqrt{10}}$;\ \ 
$ -1$,
$ \frac{-5+\sqrt{5}}{\sqrt{10}}$;\ \ 
$0$)

Factors = $2_{\frac{38}{5},1.381}^{5,491}\boxtimes 3_{\frac{3}{2},4.}^{16,538} $

Not pseudo-unitary. 

\vskip 1ex 
\black

\noindent(6,1). $6_{\frac{43}{10},14.47}^{80,424}$ \irep{48}:\ \ 
$d_i$ = ($1.0$,
$1.0$,
$1.414$,
$1.618$,
$1.618$,
$2.288$) 

\vskip 0.7ex
\hangindent=3em \hangafter=1
$D^2= 14.472 = 
10+2\sqrt{5}$

\vskip 0.7ex
\hangindent=3em \hangafter=1
$T = ( 0,
\frac{1}{2},
\frac{3}{16},
\frac{2}{5},
\frac{9}{10},
\frac{47}{80} )
$,

\vskip 0.7ex
\hangindent=3em \hangafter=1
$S$ = ($ 1$,
$ 1$,
$ \sqrt{2}$,
$ \frac{1+\sqrt{5}}{2}$,
$ \frac{1+\sqrt{5}}{2}$,
$ \frac{5+\sqrt{5}}{\sqrt{10}}$;\ \ 
$ 1$,
$ -\sqrt{2}$,
$ \frac{1+\sqrt{5}}{2}$,
$ \frac{1+\sqrt{5}}{2}$,
$ \frac{-5-\sqrt{5}}{\sqrt{10}}$;\ \ 
$0$,
$ \frac{5+\sqrt{5}}{\sqrt{10}}$,
$ \frac{-5-\sqrt{5}}{\sqrt{10}}$,
$0$;\ \ 
$ -1$,
$ -1$,
$ -\sqrt{2}$;\ \ 
$ -1$,
$ \sqrt{2}$;\ \ 
$0$)

Factors = $2_{\frac{14}{5},3.618}^{5,395}\boxtimes 3_{\frac{3}{2},4.}^{16,553} $

\vskip 1ex 
\color{grey}

\noindent(6,2). $6_{\frac{67}{10},14.47}^{80,828}$ \irep{48}:\ \ 
$d_i$ = ($1.0$,
$1.0$,
$1.414$,
$1.618$,
$1.618$,
$2.288$) 

\vskip 0.7ex
\hangindent=3em \hangafter=1
$D^2= 14.472 = 
10+2\sqrt{5}$

\vskip 0.7ex
\hangindent=3em \hangafter=1
$T = ( 0,
\frac{1}{2},
\frac{3}{16},
\frac{3}{5},
\frac{1}{10},
\frac{63}{80} )
$,

\vskip 0.7ex
\hangindent=3em \hangafter=1
$S$ = ($ 1$,
$ 1$,
$ \sqrt{2}$,
$ \frac{1+\sqrt{5}}{2}$,
$ \frac{1+\sqrt{5}}{2}$,
$ \frac{5+\sqrt{5}}{\sqrt{10}}$;\ \ 
$ 1$,
$ -\sqrt{2}$,
$ \frac{1+\sqrt{5}}{2}$,
$ \frac{1+\sqrt{5}}{2}$,
$ \frac{-5-\sqrt{5}}{\sqrt{10}}$;\ \ 
$0$,
$ \frac{5+\sqrt{5}}{\sqrt{10}}$,
$ \frac{-5-\sqrt{5}}{\sqrt{10}}$,
$0$;\ \ 
$ -1$,
$ -1$,
$ -\sqrt{2}$;\ \ 
$ -1$,
$ \sqrt{2}$;\ \ 
$0$)

Factors = $2_{\frac{26}{5},3.618}^{5,720}\boxtimes 3_{\frac{3}{2},4.}^{16,553} $

\vskip 1ex 
\color{grey}

\noindent(6,3). $6_{\frac{53}{10},14.47}^{80,884}$ \irep{48}:\ \ 
$d_i$ = ($1.0$,
$1.0$,
$1.414$,
$1.618$,
$1.618$,
$2.288$) 

\vskip 0.7ex
\hangindent=3em \hangafter=1
$D^2= 14.472 = 
10+2\sqrt{5}$

\vskip 0.7ex
\hangindent=3em \hangafter=1
$T = ( 0,
\frac{1}{2},
\frac{5}{16},
\frac{2}{5},
\frac{9}{10},
\frac{57}{80} )
$,

\vskip 0.7ex
\hangindent=3em \hangafter=1
$S$ = ($ 1$,
$ 1$,
$ \sqrt{2}$,
$ \frac{1+\sqrt{5}}{2}$,
$ \frac{1+\sqrt{5}}{2}$,
$ \frac{5+\sqrt{5}}{\sqrt{10}}$;\ \ 
$ 1$,
$ -\sqrt{2}$,
$ \frac{1+\sqrt{5}}{2}$,
$ \frac{1+\sqrt{5}}{2}$,
$ \frac{-5-\sqrt{5}}{\sqrt{10}}$;\ \ 
$0$,
$ \frac{5+\sqrt{5}}{\sqrt{10}}$,
$ \frac{-5-\sqrt{5}}{\sqrt{10}}$,
$0$;\ \ 
$ -1$,
$ -1$,
$ -\sqrt{2}$;\ \ 
$ -1$,
$ \sqrt{2}$;\ \ 
$0$)

Factors = $2_{\frac{14}{5},3.618}^{5,395}\boxtimes 3_{\frac{5}{2},4.}^{16,465} $

\vskip 1ex 
\color{grey}

\noindent(6,4). $6_{\frac{77}{10},14.47}^{80,657}$ \irep{48}:\ \ 
$d_i$ = ($1.0$,
$1.0$,
$1.414$,
$1.618$,
$1.618$,
$2.288$) 

\vskip 0.7ex
\hangindent=3em \hangafter=1
$D^2= 14.472 = 
10+2\sqrt{5}$

\vskip 0.7ex
\hangindent=3em \hangafter=1
$T = ( 0,
\frac{1}{2},
\frac{5}{16},
\frac{3}{5},
\frac{1}{10},
\frac{73}{80} )
$,

\vskip 0.7ex
\hangindent=3em \hangafter=1
$S$ = ($ 1$,
$ 1$,
$ \sqrt{2}$,
$ \frac{1+\sqrt{5}}{2}$,
$ \frac{1+\sqrt{5}}{2}$,
$ \frac{5+\sqrt{5}}{\sqrt{10}}$;\ \ 
$ 1$,
$ -\sqrt{2}$,
$ \frac{1+\sqrt{5}}{2}$,
$ \frac{1+\sqrt{5}}{2}$,
$ \frac{-5-\sqrt{5}}{\sqrt{10}}$;\ \ 
$0$,
$ \frac{5+\sqrt{5}}{\sqrt{10}}$,
$ \frac{-5-\sqrt{5}}{\sqrt{10}}$,
$0$;\ \ 
$ -1$,
$ -1$,
$ -\sqrt{2}$;\ \ 
$ -1$,
$ \sqrt{2}$;\ \ 
$0$)

Factors = $2_{\frac{26}{5},3.618}^{5,720}\boxtimes 3_{\frac{5}{2},4.}^{16,465} $

\vskip 1ex 
\color{grey}

\noindent(6,5). $6_{\frac{3}{10},14.47}^{80,270}$ \irep{48}:\ \ 
$d_i$ = ($1.0$,
$1.0$,
$1.414$,
$1.618$,
$1.618$,
$2.288$) 

\vskip 0.7ex
\hangindent=3em \hangafter=1
$D^2= 14.472 = 
10+2\sqrt{5}$

\vskip 0.7ex
\hangindent=3em \hangafter=1
$T = ( 0,
\frac{1}{2},
\frac{11}{16},
\frac{2}{5},
\frac{9}{10},
\frac{7}{80} )
$,

\vskip 0.7ex
\hangindent=3em \hangafter=1
$S$ = ($ 1$,
$ 1$,
$ \sqrt{2}$,
$ \frac{1+\sqrt{5}}{2}$,
$ \frac{1+\sqrt{5}}{2}$,
$ \frac{5+\sqrt{5}}{\sqrt{10}}$;\ \ 
$ 1$,
$ -\sqrt{2}$,
$ \frac{1+\sqrt{5}}{2}$,
$ \frac{1+\sqrt{5}}{2}$,
$ \frac{-5-\sqrt{5}}{\sqrt{10}}$;\ \ 
$0$,
$ \frac{5+\sqrt{5}}{\sqrt{10}}$,
$ \frac{-5-\sqrt{5}}{\sqrt{10}}$,
$0$;\ \ 
$ -1$,
$ -1$,
$ -\sqrt{2}$;\ \ 
$ -1$,
$ \sqrt{2}$;\ \ 
$0$)

Factors = $2_{\frac{14}{5},3.618}^{5,395}\boxtimes 3_{\frac{11}{2},4.}^{16,648} $

\vskip 1ex 
\color{grey}

\noindent(6,6). $6_{\frac{27}{10},14.47}^{80,528}$ \irep{48}:\ \ 
$d_i$ = ($1.0$,
$1.0$,
$1.414$,
$1.618$,
$1.618$,
$2.288$) 

\vskip 0.7ex
\hangindent=3em \hangafter=1
$D^2= 14.472 = 
10+2\sqrt{5}$

\vskip 0.7ex
\hangindent=3em \hangafter=1
$T = ( 0,
\frac{1}{2},
\frac{11}{16},
\frac{3}{5},
\frac{1}{10},
\frac{23}{80} )
$,

\vskip 0.7ex
\hangindent=3em \hangafter=1
$S$ = ($ 1$,
$ 1$,
$ \sqrt{2}$,
$ \frac{1+\sqrt{5}}{2}$,
$ \frac{1+\sqrt{5}}{2}$,
$ \frac{5+\sqrt{5}}{\sqrt{10}}$;\ \ 
$ 1$,
$ -\sqrt{2}$,
$ \frac{1+\sqrt{5}}{2}$,
$ \frac{1+\sqrt{5}}{2}$,
$ \frac{-5-\sqrt{5}}{\sqrt{10}}$;\ \ 
$0$,
$ \frac{5+\sqrt{5}}{\sqrt{10}}$,
$ \frac{-5-\sqrt{5}}{\sqrt{10}}$,
$0$;\ \ 
$ -1$,
$ -1$,
$ -\sqrt{2}$;\ \ 
$ -1$,
$ \sqrt{2}$;\ \ 
$0$)

Factors = $2_{\frac{26}{5},3.618}^{5,720}\boxtimes 3_{\frac{11}{2},4.}^{16,648} $

\vskip 1ex 
\color{grey}

\noindent(6,7). $6_{\frac{13}{10},14.47}^{80,621}$ \irep{48}:\ \ 
$d_i$ = ($1.0$,
$1.0$,
$1.414$,
$1.618$,
$1.618$,
$2.288$) 

\vskip 0.7ex
\hangindent=3em \hangafter=1
$D^2= 14.472 = 
10+2\sqrt{5}$

\vskip 0.7ex
\hangindent=3em \hangafter=1
$T = ( 0,
\frac{1}{2},
\frac{13}{16},
\frac{2}{5},
\frac{9}{10},
\frac{17}{80} )
$,

\vskip 0.7ex
\hangindent=3em \hangafter=1
$S$ = ($ 1$,
$ 1$,
$ \sqrt{2}$,
$ \frac{1+\sqrt{5}}{2}$,
$ \frac{1+\sqrt{5}}{2}$,
$ \frac{5+\sqrt{5}}{\sqrt{10}}$;\ \ 
$ 1$,
$ -\sqrt{2}$,
$ \frac{1+\sqrt{5}}{2}$,
$ \frac{1+\sqrt{5}}{2}$,
$ \frac{-5-\sqrt{5}}{\sqrt{10}}$;\ \ 
$0$,
$ \frac{5+\sqrt{5}}{\sqrt{10}}$,
$ \frac{-5-\sqrt{5}}{\sqrt{10}}$,
$0$;\ \ 
$ -1$,
$ -1$,
$ -\sqrt{2}$;\ \ 
$ -1$,
$ \sqrt{2}$;\ \ 
$0$)

Factors = $2_{\frac{14}{5},3.618}^{5,395}\boxtimes 3_{\frac{13}{2},4.}^{16,330} $

\vskip 1ex 
\color{grey}

\noindent(6,8). $6_{\frac{37}{10},14.47}^{80,629}$ \irep{48}:\ \ 
$d_i$ = ($1.0$,
$1.0$,
$1.414$,
$1.618$,
$1.618$,
$2.288$) 

\vskip 0.7ex
\hangindent=3em \hangafter=1
$D^2= 14.472 = 
10+2\sqrt{5}$

\vskip 0.7ex
\hangindent=3em \hangafter=1
$T = ( 0,
\frac{1}{2},
\frac{13}{16},
\frac{3}{5},
\frac{1}{10},
\frac{33}{80} )
$,

\vskip 0.7ex
\hangindent=3em \hangafter=1
$S$ = ($ 1$,
$ 1$,
$ \sqrt{2}$,
$ \frac{1+\sqrt{5}}{2}$,
$ \frac{1+\sqrt{5}}{2}$,
$ \frac{5+\sqrt{5}}{\sqrt{10}}$;\ \ 
$ 1$,
$ -\sqrt{2}$,
$ \frac{1+\sqrt{5}}{2}$,
$ \frac{1+\sqrt{5}}{2}$,
$ \frac{-5-\sqrt{5}}{\sqrt{10}}$;\ \ 
$0$,
$ \frac{5+\sqrt{5}}{\sqrt{10}}$,
$ \frac{-5-\sqrt{5}}{\sqrt{10}}$,
$0$;\ \ 
$ -1$,
$ -1$,
$ -\sqrt{2}$;\ \ 
$ -1$,
$ \sqrt{2}$;\ \ 
$0$)

Factors = $2_{\frac{26}{5},3.618}^{5,720}\boxtimes 3_{\frac{13}{2},4.}^{16,330} $

\vskip 1ex 
\color{grey}

\noindent(6,9). $6_{\frac{33}{10},14.47}^{80,941}$ \irep{48}:\ \ 
$d_i$ = ($1.0$,
$1.0$,
$1.618$,
$1.618$,
$-1.414$,
$-2.288$) 

\vskip 0.7ex
\hangindent=3em \hangafter=1
$D^2= 14.472 = 
10+2\sqrt{5}$

\vskip 0.7ex
\hangindent=3em \hangafter=1
$T = ( 0,
\frac{1}{2},
\frac{2}{5},
\frac{9}{10},
\frac{1}{16},
\frac{37}{80} )
$,

\vskip 0.7ex
\hangindent=3em \hangafter=1
$S$ = ($ 1$,
$ 1$,
$ \frac{1+\sqrt{5}}{2}$,
$ \frac{1+\sqrt{5}}{2}$,
$ -\sqrt{2}$,
$ \frac{-5-\sqrt{5}}{\sqrt{10}}$;\ \ 
$ 1$,
$ \frac{1+\sqrt{5}}{2}$,
$ \frac{1+\sqrt{5}}{2}$,
$ \sqrt{2}$,
$ \frac{5+\sqrt{5}}{\sqrt{10}}$;\ \ 
$ -1$,
$ -1$,
$ \frac{-5-\sqrt{5}}{\sqrt{10}}$,
$ \sqrt{2}$;\ \ 
$ -1$,
$ \frac{5+\sqrt{5}}{\sqrt{10}}$,
$ -\sqrt{2}$;\ \ 
$0$,
$0$;\ \ 
$0$)

Factors = $2_{\frac{14}{5},3.618}^{5,395}\boxtimes 3_{\frac{1}{2},4.}^{16,980} $

Pseudo-unitary $\sim$  
$6_{\frac{73}{10},14.47}^{80,215}$

\vskip 1ex 
\color{grey}

\noindent(6,10). $6_{\frac{63}{10},14.47}^{80,439}$ \irep{48}:\ \ 
$d_i$ = ($1.0$,
$1.0$,
$1.618$,
$1.618$,
$-1.414$,
$-2.288$) 

\vskip 0.7ex
\hangindent=3em \hangafter=1
$D^2= 14.472 = 
10+2\sqrt{5}$

\vskip 0.7ex
\hangindent=3em \hangafter=1
$T = ( 0,
\frac{1}{2},
\frac{2}{5},
\frac{9}{10},
\frac{7}{16},
\frac{67}{80} )
$,

\vskip 0.7ex
\hangindent=3em \hangafter=1
$S$ = ($ 1$,
$ 1$,
$ \frac{1+\sqrt{5}}{2}$,
$ \frac{1+\sqrt{5}}{2}$,
$ -\sqrt{2}$,
$ \frac{-5-\sqrt{5}}{\sqrt{10}}$;\ \ 
$ 1$,
$ \frac{1+\sqrt{5}}{2}$,
$ \frac{1+\sqrt{5}}{2}$,
$ \sqrt{2}$,
$ \frac{5+\sqrt{5}}{\sqrt{10}}$;\ \ 
$ -1$,
$ -1$,
$ \frac{-5-\sqrt{5}}{\sqrt{10}}$,
$ \sqrt{2}$;\ \ 
$ -1$,
$ \frac{5+\sqrt{5}}{\sqrt{10}}$,
$ -\sqrt{2}$;\ \ 
$0$,
$0$;\ \ 
$0$)

Factors = $2_{\frac{14}{5},3.618}^{5,395}\boxtimes 3_{\frac{7}{2},4.}^{16,167} $

Pseudo-unitary $\sim$  
$6_{\frac{23}{10},14.47}^{80,108}$

\vskip 1ex 
\color{grey}

\noindent(6,11). $6_{\frac{73}{10},14.47}^{80,113}$ \irep{48}:\ \ 
$d_i$ = ($1.0$,
$1.0$,
$1.618$,
$1.618$,
$-1.414$,
$-2.288$) 

\vskip 0.7ex
\hangindent=3em \hangafter=1
$D^2= 14.472 = 
10+2\sqrt{5}$

\vskip 0.7ex
\hangindent=3em \hangafter=1
$T = ( 0,
\frac{1}{2},
\frac{2}{5},
\frac{9}{10},
\frac{9}{16},
\frac{77}{80} )
$,

\vskip 0.7ex
\hangindent=3em \hangafter=1
$S$ = ($ 1$,
$ 1$,
$ \frac{1+\sqrt{5}}{2}$,
$ \frac{1+\sqrt{5}}{2}$,
$ -\sqrt{2}$,
$ \frac{-5-\sqrt{5}}{\sqrt{10}}$;\ \ 
$ 1$,
$ \frac{1+\sqrt{5}}{2}$,
$ \frac{1+\sqrt{5}}{2}$,
$ \sqrt{2}$,
$ \frac{5+\sqrt{5}}{\sqrt{10}}$;\ \ 
$ -1$,
$ -1$,
$ \frac{-5-\sqrt{5}}{\sqrt{10}}$,
$ \sqrt{2}$;\ \ 
$ -1$,
$ \frac{5+\sqrt{5}}{\sqrt{10}}$,
$ -\sqrt{2}$;\ \ 
$0$,
$0$;\ \ 
$0$)

Factors = $2_{\frac{14}{5},3.618}^{5,395}\boxtimes 3_{\frac{9}{2},4.}^{16,343} $

Pseudo-unitary $\sim$  
$6_{\frac{33}{10},14.47}^{80,798}$

\vskip 1ex 
\color{grey}

\noindent(6,12). $6_{\frac{23}{10},14.47}^{80,675}$ \irep{48}:\ \ 
$d_i$ = ($1.0$,
$1.0$,
$1.618$,
$1.618$,
$-1.414$,
$-2.288$) 

\vskip 0.7ex
\hangindent=3em \hangafter=1
$D^2= 14.472 = 
10+2\sqrt{5}$

\vskip 0.7ex
\hangindent=3em \hangafter=1
$T = ( 0,
\frac{1}{2},
\frac{2}{5},
\frac{9}{10},
\frac{15}{16},
\frac{27}{80} )
$,

\vskip 0.7ex
\hangindent=3em \hangafter=1
$S$ = ($ 1$,
$ 1$,
$ \frac{1+\sqrt{5}}{2}$,
$ \frac{1+\sqrt{5}}{2}$,
$ -\sqrt{2}$,
$ \frac{-5-\sqrt{5}}{\sqrt{10}}$;\ \ 
$ 1$,
$ \frac{1+\sqrt{5}}{2}$,
$ \frac{1+\sqrt{5}}{2}$,
$ \sqrt{2}$,
$ \frac{5+\sqrt{5}}{\sqrt{10}}$;\ \ 
$ -1$,
$ -1$,
$ \frac{-5-\sqrt{5}}{\sqrt{10}}$,
$ \sqrt{2}$;\ \ 
$ -1$,
$ \frac{5+\sqrt{5}}{\sqrt{10}}$,
$ -\sqrt{2}$;\ \ 
$0$,
$0$;\ \ 
$0$)

Factors = $2_{\frac{14}{5},3.618}^{5,395}\boxtimes 3_{\frac{15}{2},4.}^{16,113} $

Pseudo-unitary $\sim$  
$6_{\frac{63}{10},14.47}^{80,146}$

\vskip 1ex 
\color{grey}

\noindent(6,13). $6_{\frac{57}{10},14.47}^{80,544}$ \irep{48}:\ \ 
$d_i$ = ($1.0$,
$1.0$,
$1.618$,
$1.618$,
$-1.414$,
$-2.288$) 

\vskip 0.7ex
\hangindent=3em \hangafter=1
$D^2= 14.472 = 
10+2\sqrt{5}$

\vskip 0.7ex
\hangindent=3em \hangafter=1
$T = ( 0,
\frac{1}{2},
\frac{3}{5},
\frac{1}{10},
\frac{1}{16},
\frac{53}{80} )
$,

\vskip 0.7ex
\hangindent=3em \hangafter=1
$S$ = ($ 1$,
$ 1$,
$ \frac{1+\sqrt{5}}{2}$,
$ \frac{1+\sqrt{5}}{2}$,
$ -\sqrt{2}$,
$ \frac{-5-\sqrt{5}}{\sqrt{10}}$;\ \ 
$ 1$,
$ \frac{1+\sqrt{5}}{2}$,
$ \frac{1+\sqrt{5}}{2}$,
$ \sqrt{2}$,
$ \frac{5+\sqrt{5}}{\sqrt{10}}$;\ \ 
$ -1$,
$ -1$,
$ \frac{-5-\sqrt{5}}{\sqrt{10}}$,
$ \sqrt{2}$;\ \ 
$ -1$,
$ \frac{5+\sqrt{5}}{\sqrt{10}}$,
$ -\sqrt{2}$;\ \ 
$0$,
$0$;\ \ 
$0$)

Factors = $2_{\frac{26}{5},3.618}^{5,720}\boxtimes 3_{\frac{1}{2},4.}^{16,980} $

Pseudo-unitary $\sim$  
$6_{\frac{17}{10},14.47}^{80,878}$

\vskip 1ex 
\color{grey}

\noindent(6,14). $6_{\frac{7}{10},14.47}^{80,191}$ \irep{48}:\ \ 
$d_i$ = ($1.0$,
$1.0$,
$1.618$,
$1.618$,
$-1.414$,
$-2.288$) 

\vskip 0.7ex
\hangindent=3em \hangafter=1
$D^2= 14.472 = 
10+2\sqrt{5}$

\vskip 0.7ex
\hangindent=3em \hangafter=1
$T = ( 0,
\frac{1}{2},
\frac{3}{5},
\frac{1}{10},
\frac{7}{16},
\frac{3}{80} )
$,

\vskip 0.7ex
\hangindent=3em \hangafter=1
$S$ = ($ 1$,
$ 1$,
$ \frac{1+\sqrt{5}}{2}$,
$ \frac{1+\sqrt{5}}{2}$,
$ -\sqrt{2}$,
$ \frac{-5-\sqrt{5}}{\sqrt{10}}$;\ \ 
$ 1$,
$ \frac{1+\sqrt{5}}{2}$,
$ \frac{1+\sqrt{5}}{2}$,
$ \sqrt{2}$,
$ \frac{5+\sqrt{5}}{\sqrt{10}}$;\ \ 
$ -1$,
$ -1$,
$ \frac{-5-\sqrt{5}}{\sqrt{10}}$,
$ \sqrt{2}$;\ \ 
$ -1$,
$ \frac{5+\sqrt{5}}{\sqrt{10}}$,
$ -\sqrt{2}$;\ \ 
$0$,
$0$;\ \ 
$0$)

Factors = $2_{\frac{26}{5},3.618}^{5,720}\boxtimes 3_{\frac{7}{2},4.}^{16,167} $

Pseudo-unitary $\sim$  
$6_{\frac{47}{10},14.47}^{80,518}$

\vskip 1ex 
\color{grey}

\noindent(6,15). $6_{\frac{17}{10},14.47}^{80,426}$ \irep{48}:\ \ 
$d_i$ = ($1.0$,
$1.0$,
$1.618$,
$1.618$,
$-1.414$,
$-2.288$) 

\vskip 0.7ex
\hangindent=3em \hangafter=1
$D^2= 14.472 = 
10+2\sqrt{5}$

\vskip 0.7ex
\hangindent=3em \hangafter=1
$T = ( 0,
\frac{1}{2},
\frac{3}{5},
\frac{1}{10},
\frac{9}{16},
\frac{13}{80} )
$,

\vskip 0.7ex
\hangindent=3em \hangafter=1
$S$ = ($ 1$,
$ 1$,
$ \frac{1+\sqrt{5}}{2}$,
$ \frac{1+\sqrt{5}}{2}$,
$ -\sqrt{2}$,
$ \frac{-5-\sqrt{5}}{\sqrt{10}}$;\ \ 
$ 1$,
$ \frac{1+\sqrt{5}}{2}$,
$ \frac{1+\sqrt{5}}{2}$,
$ \sqrt{2}$,
$ \frac{5+\sqrt{5}}{\sqrt{10}}$;\ \ 
$ -1$,
$ -1$,
$ \frac{-5-\sqrt{5}}{\sqrt{10}}$,
$ \sqrt{2}$;\ \ 
$ -1$,
$ \frac{5+\sqrt{5}}{\sqrt{10}}$,
$ -\sqrt{2}$;\ \ 
$0$,
$0$;\ \ 
$0$)

Factors = $2_{\frac{26}{5},3.618}^{5,720}\boxtimes 3_{\frac{9}{2},4.}^{16,343} $

Pseudo-unitary $\sim$  
$6_{\frac{57}{10},14.47}^{80,376}$

\vskip 1ex 
\color{grey}

\noindent(6,16). $6_{\frac{47}{10},14.47}^{80,143}$ \irep{48}:\ \ 
$d_i$ = ($1.0$,
$1.0$,
$1.618$,
$1.618$,
$-1.414$,
$-2.288$) 

\vskip 0.7ex
\hangindent=3em \hangafter=1
$D^2= 14.472 = 
10+2\sqrt{5}$

\vskip 0.7ex
\hangindent=3em \hangafter=1
$T = ( 0,
\frac{1}{2},
\frac{3}{5},
\frac{1}{10},
\frac{15}{16},
\frac{43}{80} )
$,

\vskip 0.7ex
\hangindent=3em \hangafter=1
$S$ = ($ 1$,
$ 1$,
$ \frac{1+\sqrt{5}}{2}$,
$ \frac{1+\sqrt{5}}{2}$,
$ -\sqrt{2}$,
$ \frac{-5-\sqrt{5}}{\sqrt{10}}$;\ \ 
$ 1$,
$ \frac{1+\sqrt{5}}{2}$,
$ \frac{1+\sqrt{5}}{2}$,
$ \sqrt{2}$,
$ \frac{5+\sqrt{5}}{\sqrt{10}}$;\ \ 
$ -1$,
$ -1$,
$ \frac{-5-\sqrt{5}}{\sqrt{10}}$,
$ \sqrt{2}$;\ \ 
$ -1$,
$ \frac{5+\sqrt{5}}{\sqrt{10}}$,
$ -\sqrt{2}$;\ \ 
$0$,
$0$;\ \ 
$0$)

Factors = $2_{\frac{26}{5},3.618}^{5,720}\boxtimes 3_{\frac{15}{2},4.}^{16,113} $

Pseudo-unitary $\sim$  
$6_{\frac{7}{10},14.47}^{80,111}$

\vskip 1ex 
\color{grey}

\noindent(6,17). $6_{\frac{79}{10},5.527}^{80,135}$ \irep{48}:\ \ 
$d_i$ = ($1.0$,
$0.874$,
$1.0$,
$-0.618$,
$-0.618$,
$-1.414$) 

\vskip 0.7ex
\hangindent=3em \hangafter=1
$D^2= 5.527 = 
10-2\sqrt{5}$

\vskip 0.7ex
\hangindent=3em \hangafter=1
$T = ( 0,
\frac{11}{80},
\frac{1}{2},
\frac{1}{5},
\frac{7}{10},
\frac{15}{16} )
$,

\vskip 0.7ex
\hangindent=3em \hangafter=1
$S$ = ($ 1$,
$ \frac{5-\sqrt{5}}{\sqrt{10}}$,
$ 1$,
$ \frac{1-\sqrt{5}}{2}$,
$ \frac{1-\sqrt{5}}{2}$,
$ -\sqrt{2}$;\ \ 
$0$,
$ \frac{-5+\sqrt{5}}{\sqrt{10}}$,
$ \sqrt{2}$,
$ -\sqrt{2}$,
$0$;\ \ 
$ 1$,
$ \frac{1-\sqrt{5}}{2}$,
$ \frac{1-\sqrt{5}}{2}$,
$ \sqrt{2}$;\ \ 
$ -1$,
$ -1$,
$ \frac{5-\sqrt{5}}{\sqrt{10}}$;\ \ 
$ -1$,
$ \frac{-5+\sqrt{5}}{\sqrt{10}}$;\ \ 
$0$)

Factors = $2_{\frac{2}{5},1.381}^{5,120}\boxtimes 3_{\frac{15}{2},4.}^{16,113} $

Not pseudo-unitary. 

\vskip 1ex 
\color{grey}

\noindent(6,18). $6_{\frac{31}{10},5.527}^{80,478}$ \irep{48}:\ \ 
$d_i$ = ($1.0$,
$0.874$,
$1.0$,
$-0.618$,
$-0.618$,
$-1.414$) 

\vskip 0.7ex
\hangindent=3em \hangafter=1
$D^2= 5.527 = 
10-2\sqrt{5}$

\vskip 0.7ex
\hangindent=3em \hangafter=1
$T = ( 0,
\frac{19}{80},
\frac{1}{2},
\frac{4}{5},
\frac{3}{10},
\frac{7}{16} )
$,

\vskip 0.7ex
\hangindent=3em \hangafter=1
$S$ = ($ 1$,
$ \frac{5-\sqrt{5}}{\sqrt{10}}$,
$ 1$,
$ \frac{1-\sqrt{5}}{2}$,
$ \frac{1-\sqrt{5}}{2}$,
$ -\sqrt{2}$;\ \ 
$0$,
$ \frac{-5+\sqrt{5}}{\sqrt{10}}$,
$ \sqrt{2}$,
$ -\sqrt{2}$,
$0$;\ \ 
$ 1$,
$ \frac{1-\sqrt{5}}{2}$,
$ \frac{1-\sqrt{5}}{2}$,
$ \sqrt{2}$;\ \ 
$ -1$,
$ -1$,
$ \frac{5-\sqrt{5}}{\sqrt{10}}$;\ \ 
$ -1$,
$ \frac{-5+\sqrt{5}}{\sqrt{10}}$;\ \ 
$0$)

Factors = $2_{\frac{38}{5},1.381}^{5,491}\boxtimes 3_{\frac{7}{2},4.}^{16,167} $

Not pseudo-unitary. 

\vskip 1ex 
\color{grey}

\noindent(6,19). $6_{\frac{9}{10},5.527}^{80,787}$ \irep{48}:\ \ 
$d_i$ = ($1.0$,
$0.874$,
$1.0$,
$-0.618$,
$-0.618$,
$-1.414$) 

\vskip 0.7ex
\hangindent=3em \hangafter=1
$D^2= 5.527 = 
10-2\sqrt{5}$

\vskip 0.7ex
\hangindent=3em \hangafter=1
$T = ( 0,
\frac{21}{80},
\frac{1}{2},
\frac{1}{5},
\frac{7}{10},
\frac{1}{16} )
$,

\vskip 0.7ex
\hangindent=3em \hangafter=1
$S$ = ($ 1$,
$ \frac{5-\sqrt{5}}{\sqrt{10}}$,
$ 1$,
$ \frac{1-\sqrt{5}}{2}$,
$ \frac{1-\sqrt{5}}{2}$,
$ -\sqrt{2}$;\ \ 
$0$,
$ \frac{-5+\sqrt{5}}{\sqrt{10}}$,
$ \sqrt{2}$,
$ -\sqrt{2}$,
$0$;\ \ 
$ 1$,
$ \frac{1-\sqrt{5}}{2}$,
$ \frac{1-\sqrt{5}}{2}$,
$ \sqrt{2}$;\ \ 
$ -1$,
$ -1$,
$ \frac{5-\sqrt{5}}{\sqrt{10}}$;\ \ 
$ -1$,
$ \frac{-5+\sqrt{5}}{\sqrt{10}}$;\ \ 
$0$)

Factors = $2_{\frac{2}{5},1.381}^{5,120}\boxtimes 3_{\frac{1}{2},4.}^{16,980} $

Not pseudo-unitary. 

\vskip 1ex 
\color{grey}

\noindent(6,20). $6_{\frac{41}{10},5.527}^{80,589}$ \irep{48}:\ \ 
$d_i$ = ($1.0$,
$0.874$,
$1.0$,
$-0.618$,
$-0.618$,
$-1.414$) 

\vskip 0.7ex
\hangindent=3em \hangafter=1
$D^2= 5.527 = 
10-2\sqrt{5}$

\vskip 0.7ex
\hangindent=3em \hangafter=1
$T = ( 0,
\frac{29}{80},
\frac{1}{2},
\frac{4}{5},
\frac{3}{10},
\frac{9}{16} )
$,

\vskip 0.7ex
\hangindent=3em \hangafter=1
$S$ = ($ 1$,
$ \frac{5-\sqrt{5}}{\sqrt{10}}$,
$ 1$,
$ \frac{1-\sqrt{5}}{2}$,
$ \frac{1-\sqrt{5}}{2}$,
$ -\sqrt{2}$;\ \ 
$0$,
$ \frac{-5+\sqrt{5}}{\sqrt{10}}$,
$ \sqrt{2}$,
$ -\sqrt{2}$,
$0$;\ \ 
$ 1$,
$ \frac{1-\sqrt{5}}{2}$,
$ \frac{1-\sqrt{5}}{2}$,
$ \sqrt{2}$;\ \ 
$ -1$,
$ -1$,
$ \frac{5-\sqrt{5}}{\sqrt{10}}$;\ \ 
$ -1$,
$ \frac{-5+\sqrt{5}}{\sqrt{10}}$;\ \ 
$0$)

Factors = $2_{\frac{38}{5},1.381}^{5,491}\boxtimes 3_{\frac{9}{2},4.}^{16,343} $

Not pseudo-unitary. 

\vskip 1ex 
\color{grey}

\noindent(6,21). $6_{\frac{19}{10},5.527}^{80,485}$ \irep{48}:\ \ 
$d_i$ = ($1.0$,
$1.0$,
$1.414$,
$-0.618$,
$-0.618$,
$-0.874$) 

\vskip 0.7ex
\hangindent=3em \hangafter=1
$D^2= 5.527 = 
10-2\sqrt{5}$

\vskip 0.7ex
\hangindent=3em \hangafter=1
$T = ( 0,
\frac{1}{2},
\frac{3}{16},
\frac{1}{5},
\frac{7}{10},
\frac{31}{80} )
$,

\vskip 0.7ex
\hangindent=3em \hangafter=1
$S$ = ($ 1$,
$ 1$,
$ \sqrt{2}$,
$ \frac{1-\sqrt{5}}{2}$,
$ \frac{1-\sqrt{5}}{2}$,
$ \frac{-5+\sqrt{5}}{\sqrt{10}}$;\ \ 
$ 1$,
$ -\sqrt{2}$,
$ \frac{1-\sqrt{5}}{2}$,
$ \frac{1-\sqrt{5}}{2}$,
$ \frac{5-\sqrt{5}}{\sqrt{10}}$;\ \ 
$0$,
$ \frac{-5+\sqrt{5}}{\sqrt{10}}$,
$ \frac{5-\sqrt{5}}{\sqrt{10}}$,
$0$;\ \ 
$ -1$,
$ -1$,
$ -\sqrt{2}$;\ \ 
$ -1$,
$ \sqrt{2}$;\ \ 
$0$)

Factors = $2_{\frac{2}{5},1.381}^{5,120}\boxtimes 3_{\frac{3}{2},4.}^{16,553} $

Not pseudo-unitary. 

\vskip 1ex 
\color{grey}

\noindent(6,22). $6_{\frac{11}{10},5.527}^{80,132}$ \irep{48}:\ \ 
$d_i$ = ($1.0$,
$1.0$,
$1.414$,
$-0.618$,
$-0.618$,
$-0.874$) 

\vskip 0.7ex
\hangindent=3em \hangafter=1
$D^2= 5.527 = 
10-2\sqrt{5}$

\vskip 0.7ex
\hangindent=3em \hangafter=1
$T = ( 0,
\frac{1}{2},
\frac{3}{16},
\frac{4}{5},
\frac{3}{10},
\frac{79}{80} )
$,

\vskip 0.7ex
\hangindent=3em \hangafter=1
$S$ = ($ 1$,
$ 1$,
$ \sqrt{2}$,
$ \frac{1-\sqrt{5}}{2}$,
$ \frac{1-\sqrt{5}}{2}$,
$ \frac{-5+\sqrt{5}}{\sqrt{10}}$;\ \ 
$ 1$,
$ -\sqrt{2}$,
$ \frac{1-\sqrt{5}}{2}$,
$ \frac{1-\sqrt{5}}{2}$,
$ \frac{5-\sqrt{5}}{\sqrt{10}}$;\ \ 
$0$,
$ \frac{-5+\sqrt{5}}{\sqrt{10}}$,
$ \frac{5-\sqrt{5}}{\sqrt{10}}$,
$0$;\ \ 
$ -1$,
$ -1$,
$ -\sqrt{2}$;\ \ 
$ -1$,
$ \sqrt{2}$;\ \ 
$0$)

Factors = $2_{\frac{38}{5},1.381}^{5,491}\boxtimes 3_{\frac{3}{2},4.}^{16,553} $

Not pseudo-unitary. 

\vskip 1ex 
\color{grey}

\noindent(6,23). $6_{\frac{29}{10},5.527}^{80,495}$ \irep{48}:\ \ 
$d_i$ = ($1.0$,
$1.0$,
$1.414$,
$-0.618$,
$-0.618$,
$-0.874$) 

\vskip 0.7ex
\hangindent=3em \hangafter=1
$D^2= 5.527 = 
10-2\sqrt{5}$

\vskip 0.7ex
\hangindent=3em \hangafter=1
$T = ( 0,
\frac{1}{2},
\frac{5}{16},
\frac{1}{5},
\frac{7}{10},
\frac{41}{80} )
$,

\vskip 0.7ex
\hangindent=3em \hangafter=1
$S$ = ($ 1$,
$ 1$,
$ \sqrt{2}$,
$ \frac{1-\sqrt{5}}{2}$,
$ \frac{1-\sqrt{5}}{2}$,
$ \frac{-5+\sqrt{5}}{\sqrt{10}}$;\ \ 
$ 1$,
$ -\sqrt{2}$,
$ \frac{1-\sqrt{5}}{2}$,
$ \frac{1-\sqrt{5}}{2}$,
$ \frac{5-\sqrt{5}}{\sqrt{10}}$;\ \ 
$0$,
$ \frac{-5+\sqrt{5}}{\sqrt{10}}$,
$ \frac{5-\sqrt{5}}{\sqrt{10}}$,
$0$;\ \ 
$ -1$,
$ -1$,
$ -\sqrt{2}$;\ \ 
$ -1$,
$ \sqrt{2}$;\ \ 
$0$)

Factors = $2_{\frac{2}{5},1.381}^{5,120}\boxtimes 3_{\frac{5}{2},4.}^{16,465} $

Not pseudo-unitary. 

\vskip 1ex 
\color{grey}

\noindent(6,24). $6_{\frac{21}{10},5.527}^{80,400}$ \irep{48}:\ \ 
$d_i$ = ($1.0$,
$1.0$,
$1.414$,
$-0.618$,
$-0.618$,
$-0.874$) 

\vskip 0.7ex
\hangindent=3em \hangafter=1
$D^2= 5.527 = 
10-2\sqrt{5}$

\vskip 0.7ex
\hangindent=3em \hangafter=1
$T = ( 0,
\frac{1}{2},
\frac{5}{16},
\frac{4}{5},
\frac{3}{10},
\frac{9}{80} )
$,

\vskip 0.7ex
\hangindent=3em \hangafter=1
$S$ = ($ 1$,
$ 1$,
$ \sqrt{2}$,
$ \frac{1-\sqrt{5}}{2}$,
$ \frac{1-\sqrt{5}}{2}$,
$ \frac{-5+\sqrt{5}}{\sqrt{10}}$;\ \ 
$ 1$,
$ -\sqrt{2}$,
$ \frac{1-\sqrt{5}}{2}$,
$ \frac{1-\sqrt{5}}{2}$,
$ \frac{5-\sqrt{5}}{\sqrt{10}}$;\ \ 
$0$,
$ \frac{-5+\sqrt{5}}{\sqrt{10}}$,
$ \frac{5-\sqrt{5}}{\sqrt{10}}$,
$0$;\ \ 
$ -1$,
$ -1$,
$ -\sqrt{2}$;\ \ 
$ -1$,
$ \sqrt{2}$;\ \ 
$0$)

Factors = $2_{\frac{38}{5},1.381}^{5,491}\boxtimes 3_{\frac{5}{2},4.}^{16,465} $

Not pseudo-unitary. 

\vskip 1ex 
\color{grey}

\noindent(6,25). $6_{\frac{59}{10},5.527}^{80,423}$ \irep{48}:\ \ 
$d_i$ = ($1.0$,
$1.0$,
$1.414$,
$-0.618$,
$-0.618$,
$-0.874$) 

\vskip 0.7ex
\hangindent=3em \hangafter=1
$D^2= 5.527 = 
10-2\sqrt{5}$

\vskip 0.7ex
\hangindent=3em \hangafter=1
$T = ( 0,
\frac{1}{2},
\frac{11}{16},
\frac{1}{5},
\frac{7}{10},
\frac{71}{80} )
$,

\vskip 0.7ex
\hangindent=3em \hangafter=1
$S$ = ($ 1$,
$ 1$,
$ \sqrt{2}$,
$ \frac{1-\sqrt{5}}{2}$,
$ \frac{1-\sqrt{5}}{2}$,
$ \frac{-5+\sqrt{5}}{\sqrt{10}}$;\ \ 
$ 1$,
$ -\sqrt{2}$,
$ \frac{1-\sqrt{5}}{2}$,
$ \frac{1-\sqrt{5}}{2}$,
$ \frac{5-\sqrt{5}}{\sqrt{10}}$;\ \ 
$0$,
$ \frac{-5+\sqrt{5}}{\sqrt{10}}$,
$ \frac{5-\sqrt{5}}{\sqrt{10}}$,
$0$;\ \ 
$ -1$,
$ -1$,
$ -\sqrt{2}$;\ \ 
$ -1$,
$ \sqrt{2}$;\ \ 
$0$)

Factors = $2_{\frac{2}{5},1.381}^{5,120}\boxtimes 3_{\frac{11}{2},4.}^{16,648} $

Not pseudo-unitary. 

\vskip 1ex 
\color{grey}

\noindent(6,26). $6_{\frac{51}{10},5.527}^{80,664}$ \irep{48}:\ \ 
$d_i$ = ($1.0$,
$1.0$,
$1.414$,
$-0.618$,
$-0.618$,
$-0.874$) 

\vskip 0.7ex
\hangindent=3em \hangafter=1
$D^2= 5.527 = 
10-2\sqrt{5}$

\vskip 0.7ex
\hangindent=3em \hangafter=1
$T = ( 0,
\frac{1}{2},
\frac{11}{16},
\frac{4}{5},
\frac{3}{10},
\frac{39}{80} )
$,

\vskip 0.7ex
\hangindent=3em \hangafter=1
$S$ = ($ 1$,
$ 1$,
$ \sqrt{2}$,
$ \frac{1-\sqrt{5}}{2}$,
$ \frac{1-\sqrt{5}}{2}$,
$ \frac{-5+\sqrt{5}}{\sqrt{10}}$;\ \ 
$ 1$,
$ -\sqrt{2}$,
$ \frac{1-\sqrt{5}}{2}$,
$ \frac{1-\sqrt{5}}{2}$,
$ \frac{5-\sqrt{5}}{\sqrt{10}}$;\ \ 
$0$,
$ \frac{-5+\sqrt{5}}{\sqrt{10}}$,
$ \frac{5-\sqrt{5}}{\sqrt{10}}$,
$0$;\ \ 
$ -1$,
$ -1$,
$ -\sqrt{2}$;\ \ 
$ -1$,
$ \sqrt{2}$;\ \ 
$0$)

Factors = $2_{\frac{38}{5},1.381}^{5,491}\boxtimes 3_{\frac{11}{2},4.}^{16,648} $

Not pseudo-unitary. 

\vskip 1ex 
\color{grey}

\noindent(6,27). $6_{\frac{69}{10},5.527}^{80,529}$ \irep{48}:\ \ 
$d_i$ = ($1.0$,
$1.0$,
$1.414$,
$-0.618$,
$-0.618$,
$-0.874$) 

\vskip 0.7ex
\hangindent=3em \hangafter=1
$D^2= 5.527 = 
10-2\sqrt{5}$

\vskip 0.7ex
\hangindent=3em \hangafter=1
$T = ( 0,
\frac{1}{2},
\frac{13}{16},
\frac{1}{5},
\frac{7}{10},
\frac{1}{80} )
$,

\vskip 0.7ex
\hangindent=3em \hangafter=1
$S$ = ($ 1$,
$ 1$,
$ \sqrt{2}$,
$ \frac{1-\sqrt{5}}{2}$,
$ \frac{1-\sqrt{5}}{2}$,
$ \frac{-5+\sqrt{5}}{\sqrt{10}}$;\ \ 
$ 1$,
$ -\sqrt{2}$,
$ \frac{1-\sqrt{5}}{2}$,
$ \frac{1-\sqrt{5}}{2}$,
$ \frac{5-\sqrt{5}}{\sqrt{10}}$;\ \ 
$0$,
$ \frac{-5+\sqrt{5}}{\sqrt{10}}$,
$ \frac{5-\sqrt{5}}{\sqrt{10}}$,
$0$;\ \ 
$ -1$,
$ -1$,
$ -\sqrt{2}$;\ \ 
$ -1$,
$ \sqrt{2}$;\ \ 
$0$)

Factors = $2_{\frac{2}{5},1.381}^{5,120}\boxtimes 3_{\frac{13}{2},4.}^{16,330} $

Not pseudo-unitary. 

\vskip 1ex 
\color{grey}

\noindent(6,28). $6_{\frac{61}{10},5.527}^{80,785}$ \irep{48}:\ \ 
$d_i$ = ($1.0$,
$1.0$,
$1.414$,
$-0.618$,
$-0.618$,
$-0.874$) 

\vskip 0.7ex
\hangindent=3em \hangafter=1
$D^2= 5.527 = 
10-2\sqrt{5}$

\vskip 0.7ex
\hangindent=3em \hangafter=1
$T = ( 0,
\frac{1}{2},
\frac{13}{16},
\frac{4}{5},
\frac{3}{10},
\frac{49}{80} )
$,

\vskip 0.7ex
\hangindent=3em \hangafter=1
$S$ = ($ 1$,
$ 1$,
$ \sqrt{2}$,
$ \frac{1-\sqrt{5}}{2}$,
$ \frac{1-\sqrt{5}}{2}$,
$ \frac{-5+\sqrt{5}}{\sqrt{10}}$;\ \ 
$ 1$,
$ -\sqrt{2}$,
$ \frac{1-\sqrt{5}}{2}$,
$ \frac{1-\sqrt{5}}{2}$,
$ \frac{5-\sqrt{5}}{\sqrt{10}}$;\ \ 
$0$,
$ \frac{-5+\sqrt{5}}{\sqrt{10}}$,
$ \frac{5-\sqrt{5}}{\sqrt{10}}$,
$0$;\ \ 
$ -1$,
$ -1$,
$ -\sqrt{2}$;\ \ 
$ -1$,
$ \sqrt{2}$;\ \ 
$0$)

Factors = $2_{\frac{38}{5},1.381}^{5,491}\boxtimes 3_{\frac{13}{2},4.}^{16,330} $

Not pseudo-unitary. 

\vskip 1ex 
\color{grey}

\noindent(6,29). $6_{\frac{39}{10},5.527}^{80,489}$ \irep{48}:\ \ 
$d_i$ = ($1.0$,
$0.874$,
$1.0$,
$-0.618$,
$-0.618$,
$-1.414$) 

\vskip 0.7ex
\hangindent=3em \hangafter=1
$D^2= 5.527 = 
10-2\sqrt{5}$

\vskip 0.7ex
\hangindent=3em \hangafter=1
$T = ( 0,
\frac{51}{80},
\frac{1}{2},
\frac{1}{5},
\frac{7}{10},
\frac{7}{16} )
$,

\vskip 0.7ex
\hangindent=3em \hangafter=1
$S$ = ($ 1$,
$ \frac{5-\sqrt{5}}{\sqrt{10}}$,
$ 1$,
$ \frac{1-\sqrt{5}}{2}$,
$ \frac{1-\sqrt{5}}{2}$,
$ -\sqrt{2}$;\ \ 
$0$,
$ \frac{-5+\sqrt{5}}{\sqrt{10}}$,
$ \sqrt{2}$,
$ -\sqrt{2}$,
$0$;\ \ 
$ 1$,
$ \frac{1-\sqrt{5}}{2}$,
$ \frac{1-\sqrt{5}}{2}$,
$ \sqrt{2}$;\ \ 
$ -1$,
$ -1$,
$ \frac{5-\sqrt{5}}{\sqrt{10}}$;\ \ 
$ -1$,
$ \frac{-5+\sqrt{5}}{\sqrt{10}}$;\ \ 
$0$)

Factors = $2_{\frac{2}{5},1.381}^{5,120}\boxtimes 3_{\frac{7}{2},4.}^{16,167} $

Not pseudo-unitary. 

\vskip 1ex 
\color{grey}

\noindent(6,30). $6_{\frac{71}{10},5.527}^{80,102}$ \irep{48}:\ \ 
$d_i$ = ($1.0$,
$0.874$,
$1.0$,
$-0.618$,
$-0.618$,
$-1.414$) 

\vskip 0.7ex
\hangindent=3em \hangafter=1
$D^2= 5.527 = 
10-2\sqrt{5}$

\vskip 0.7ex
\hangindent=3em \hangafter=1
$T = ( 0,
\frac{59}{80},
\frac{1}{2},
\frac{4}{5},
\frac{3}{10},
\frac{15}{16} )
$,

\vskip 0.7ex
\hangindent=3em \hangafter=1
$S$ = ($ 1$,
$ \frac{5-\sqrt{5}}{\sqrt{10}}$,
$ 1$,
$ \frac{1-\sqrt{5}}{2}$,
$ \frac{1-\sqrt{5}}{2}$,
$ -\sqrt{2}$;\ \ 
$0$,
$ \frac{-5+\sqrt{5}}{\sqrt{10}}$,
$ \sqrt{2}$,
$ -\sqrt{2}$,
$0$;\ \ 
$ 1$,
$ \frac{1-\sqrt{5}}{2}$,
$ \frac{1-\sqrt{5}}{2}$,
$ \sqrt{2}$;\ \ 
$ -1$,
$ -1$,
$ \frac{5-\sqrt{5}}{\sqrt{10}}$;\ \ 
$ -1$,
$ \frac{-5+\sqrt{5}}{\sqrt{10}}$;\ \ 
$0$)

Factors = $2_{\frac{38}{5},1.381}^{5,491}\boxtimes 3_{\frac{15}{2},4.}^{16,113} $

Not pseudo-unitary. 

\vskip 1ex 
\color{grey}

\noindent(6,31). $6_{\frac{49}{10},5.527}^{80,727}$ \irep{48}:\ \ 
$d_i$ = ($1.0$,
$0.874$,
$1.0$,
$-0.618$,
$-0.618$,
$-1.414$) 

\vskip 0.7ex
\hangindent=3em \hangafter=1
$D^2= 5.527 = 
10-2\sqrt{5}$

\vskip 0.7ex
\hangindent=3em \hangafter=1
$T = ( 0,
\frac{61}{80},
\frac{1}{2},
\frac{1}{5},
\frac{7}{10},
\frac{9}{16} )
$,

\vskip 0.7ex
\hangindent=3em \hangafter=1
$S$ = ($ 1$,
$ \frac{5-\sqrt{5}}{\sqrt{10}}$,
$ 1$,
$ \frac{1-\sqrt{5}}{2}$,
$ \frac{1-\sqrt{5}}{2}$,
$ -\sqrt{2}$;\ \ 
$0$,
$ \frac{-5+\sqrt{5}}{\sqrt{10}}$,
$ \sqrt{2}$,
$ -\sqrt{2}$,
$0$;\ \ 
$ 1$,
$ \frac{1-\sqrt{5}}{2}$,
$ \frac{1-\sqrt{5}}{2}$,
$ \sqrt{2}$;\ \ 
$ -1$,
$ -1$,
$ \frac{5-\sqrt{5}}{\sqrt{10}}$;\ \ 
$ -1$,
$ \frac{-5+\sqrt{5}}{\sqrt{10}}$;\ \ 
$0$)

Factors = $2_{\frac{2}{5},1.381}^{5,120}\boxtimes 3_{\frac{9}{2},4.}^{16,343} $

Not pseudo-unitary. 

\vskip 1ex 
\color{grey}

\noindent(6,32). $6_{\frac{1}{10},5.527}^{80,387}$ \irep{48}:\ \ 
$d_i$ = ($1.0$,
$0.874$,
$1.0$,
$-0.618$,
$-0.618$,
$-1.414$) 

\vskip 0.7ex
\hangindent=3em \hangafter=1
$D^2= 5.527 = 
10-2\sqrt{5}$

\vskip 0.7ex
\hangindent=3em \hangafter=1
$T = ( 0,
\frac{69}{80},
\frac{1}{2},
\frac{4}{5},
\frac{3}{10},
\frac{1}{16} )
$,

\vskip 0.7ex
\hangindent=3em \hangafter=1
$S$ = ($ 1$,
$ \frac{5-\sqrt{5}}{\sqrt{10}}$,
$ 1$,
$ \frac{1-\sqrt{5}}{2}$,
$ \frac{1-\sqrt{5}}{2}$,
$ -\sqrt{2}$;\ \ 
$0$,
$ \frac{-5+\sqrt{5}}{\sqrt{10}}$,
$ \sqrt{2}$,
$ -\sqrt{2}$,
$0$;\ \ 
$ 1$,
$ \frac{1-\sqrt{5}}{2}$,
$ \frac{1-\sqrt{5}}{2}$,
$ \sqrt{2}$;\ \ 
$ -1$,
$ -1$,
$ \frac{5-\sqrt{5}}{\sqrt{10}}$;\ \ 
$ -1$,
$ \frac{-5+\sqrt{5}}{\sqrt{10}}$;\ \ 
$0$)

Factors = $2_{\frac{38}{5},1.381}^{5,491}\boxtimes 3_{\frac{1}{2},4.}^{16,980} $

Not pseudo-unitary. 

\vskip 1ex 
\black

\noindent(7,1). $6_{\frac{55}{7},18.59}^{28,108}$ \irep{46}:\ \ 
$d_i$ = ($1.0$,
$1.0$,
$1.801$,
$1.801$,
$2.246$,
$2.246$) 

\vskip 0.7ex
\hangindent=3em \hangafter=1
$D^2= 18.591 = 
12+6c^{1}_{7}
+2c^{2}_{7}
$

\vskip 0.7ex
\hangindent=3em \hangafter=1
$T = ( 0,
\frac{1}{4},
\frac{1}{7},
\frac{11}{28},
\frac{5}{7},
\frac{27}{28} )
$,

\vskip 0.7ex
\hangindent=3em \hangafter=1
$S$ = ($ 1$,
$ 1$,
$ -c_{7}^{3}$,
$ -c_{7}^{3}$,
$ \xi_{7}^{3}$,
$ \xi_{7}^{3}$;\ \ 
$ -1$,
$ -c_{7}^{3}$,
$ c_{7}^{3}$,
$ \xi_{7}^{3}$,
$ -\xi_{7}^{3}$;\ \ 
$ -\xi_{7}^{3}$,
$ -\xi_{7}^{3}$,
$ 1$,
$ 1$;\ \ 
$ \xi_{7}^{3}$,
$ 1$,
$ -1$;\ \ 
$ c_{7}^{3}$,
$ c_{7}^{3}$;\ \ 
$ -c_{7}^{3}$)

Factors = $2_{1,2.}^{4,437}\boxtimes 3_{\frac{48}{7},9.295}^{7,790} $

\vskip 1ex 
\color{grey}

\noindent(7,2). $6_{\frac{15}{7},18.59}^{28,289}$ \irep{46}:\ \ 
$d_i$ = ($1.0$,
$1.0$,
$1.801$,
$1.801$,
$2.246$,
$2.246$) 

\vskip 0.7ex
\hangindent=3em \hangafter=1
$D^2= 18.591 = 
12+6c^{1}_{7}
+2c^{2}_{7}
$

\vskip 0.7ex
\hangindent=3em \hangafter=1
$T = ( 0,
\frac{1}{4},
\frac{6}{7},
\frac{3}{28},
\frac{2}{7},
\frac{15}{28} )
$,

\vskip 0.7ex
\hangindent=3em \hangafter=1
$S$ = ($ 1$,
$ 1$,
$ -c_{7}^{3}$,
$ -c_{7}^{3}$,
$ \xi_{7}^{3}$,
$ \xi_{7}^{3}$;\ \ 
$ -1$,
$ -c_{7}^{3}$,
$ c_{7}^{3}$,
$ \xi_{7}^{3}$,
$ -\xi_{7}^{3}$;\ \ 
$ -\xi_{7}^{3}$,
$ -\xi_{7}^{3}$,
$ 1$,
$ 1$;\ \ 
$ \xi_{7}^{3}$,
$ 1$,
$ -1$;\ \ 
$ c_{7}^{3}$,
$ c_{7}^{3}$;\ \ 
$ -c_{7}^{3}$)

Factors = $2_{1,2.}^{4,437}\boxtimes 3_{\frac{8}{7},9.295}^{7,245} $

\vskip 1ex 
\color{grey}

\noindent(7,3). $6_{\frac{41}{7},18.59}^{28,114}$ \irep{46}:\ \ 
$d_i$ = ($1.0$,
$1.0$,
$1.801$,
$1.801$,
$2.246$,
$2.246$) 

\vskip 0.7ex
\hangindent=3em \hangafter=1
$D^2= 18.591 = 
12+6c^{1}_{7}
+2c^{2}_{7}
$

\vskip 0.7ex
\hangindent=3em \hangafter=1
$T = ( 0,
\frac{3}{4},
\frac{1}{7},
\frac{25}{28},
\frac{5}{7},
\frac{13}{28} )
$,

\vskip 0.7ex
\hangindent=3em \hangafter=1
$S$ = ($ 1$,
$ 1$,
$ -c_{7}^{3}$,
$ -c_{7}^{3}$,
$ \xi_{7}^{3}$,
$ \xi_{7}^{3}$;\ \ 
$ -1$,
$ -c_{7}^{3}$,
$ c_{7}^{3}$,
$ \xi_{7}^{3}$,
$ -\xi_{7}^{3}$;\ \ 
$ -\xi_{7}^{3}$,
$ -\xi_{7}^{3}$,
$ 1$,
$ 1$;\ \ 
$ \xi_{7}^{3}$,
$ 1$,
$ -1$;\ \ 
$ c_{7}^{3}$,
$ c_{7}^{3}$;\ \ 
$ -c_{7}^{3}$)

Factors = $2_{7,2.}^{4,625}\boxtimes 3_{\frac{48}{7},9.295}^{7,790} $

\vskip 1ex 
\color{grey}

\noindent(7,4). $6_{\frac{1}{7},18.59}^{28,212}$ \irep{46}:\ \ 
$d_i$ = ($1.0$,
$1.0$,
$1.801$,
$1.801$,
$2.246$,
$2.246$) 

\vskip 0.7ex
\hangindent=3em \hangafter=1
$D^2= 18.591 = 
12+6c^{1}_{7}
+2c^{2}_{7}
$

\vskip 0.7ex
\hangindent=3em \hangafter=1
$T = ( 0,
\frac{3}{4},
\frac{6}{7},
\frac{17}{28},
\frac{2}{7},
\frac{1}{28} )
$,

\vskip 0.7ex
\hangindent=3em \hangafter=1
$S$ = ($ 1$,
$ 1$,
$ -c_{7}^{3}$,
$ -c_{7}^{3}$,
$ \xi_{7}^{3}$,
$ \xi_{7}^{3}$;\ \ 
$ -1$,
$ -c_{7}^{3}$,
$ c_{7}^{3}$,
$ \xi_{7}^{3}$,
$ -\xi_{7}^{3}$;\ \ 
$ -\xi_{7}^{3}$,
$ -\xi_{7}^{3}$,
$ 1$,
$ 1$;\ \ 
$ \xi_{7}^{3}$,
$ 1$,
$ -1$;\ \ 
$ c_{7}^{3}$,
$ c_{7}^{3}$;\ \ 
$ -c_{7}^{3}$)

Factors = $2_{7,2.}^{4,625}\boxtimes 3_{\frac{8}{7},9.295}^{7,245} $

\vskip 1ex 
\color{grey}

\noindent(7,5). $6_{\frac{5}{7},5.725}^{28,424}$ \irep{46}:\ \ 
$d_i$ = ($1.0$,
$0.554$,
$0.554$,
$1.0$,
$-1.246$,
$-1.246$) 

\vskip 0.7ex
\hangindent=3em \hangafter=1
$D^2= 5.725 = 
10-2  c^{1}_{7}
+4c^{2}_{7}
$

\vskip 0.7ex
\hangindent=3em \hangafter=1
$T = ( 0,
\frac{3}{7},
\frac{5}{28},
\frac{3}{4},
\frac{2}{7},
\frac{1}{28} )
$,

\vskip 0.7ex
\hangindent=3em \hangafter=1
$S$ = ($ 1$,
$ 1+c^{2}_{7}
$,
$ 1+c^{2}_{7}
$,
$ 1$,
$ -c_{7}^{1}$,
$ -c_{7}^{1}$;\ \ 
$ c_{7}^{1}$,
$ c_{7}^{1}$,
$ 1+c^{2}_{7}
$,
$ 1$,
$ 1$;\ \ 
$ -c_{7}^{1}$,
$ -1-c^{2}_{7}
$,
$ 1$,
$ -1$;\ \ 
$ -1$,
$ -c_{7}^{1}$,
$ c_{7}^{1}$;\ \ 
$ -1-c^{2}_{7}
$,
$ -1-c^{2}_{7}
$;\ \ 
$ 1+c^{2}_{7}
$)

Factors = $2_{7,2.}^{4,625}\boxtimes 3_{\frac{12}{7},2.862}^{7,768} $

Not pseudo-unitary. 

\vskip 1ex 
\color{grey}

\noindent(7,6). $6_{\frac{19}{7},5.725}^{28,193}$ \irep{46}:\ \ 
$d_i$ = ($1.0$,
$0.554$,
$0.554$,
$1.0$,
$-1.246$,
$-1.246$) 

\vskip 0.7ex
\hangindent=3em \hangafter=1
$D^2= 5.725 = 
10-2  c^{1}_{7}
+4c^{2}_{7}
$

\vskip 0.7ex
\hangindent=3em \hangafter=1
$T = ( 0,
\frac{3}{7},
\frac{19}{28},
\frac{1}{4},
\frac{2}{7},
\frac{15}{28} )
$,

\vskip 0.7ex
\hangindent=3em \hangafter=1
$S$ = ($ 1$,
$ 1+c^{2}_{7}
$,
$ 1+c^{2}_{7}
$,
$ 1$,
$ -c_{7}^{1}$,
$ -c_{7}^{1}$;\ \ 
$ c_{7}^{1}$,
$ c_{7}^{1}$,
$ 1+c^{2}_{7}
$,
$ 1$,
$ 1$;\ \ 
$ -c_{7}^{1}$,
$ -1-c^{2}_{7}
$,
$ 1$,
$ -1$;\ \ 
$ -1$,
$ -c_{7}^{1}$,
$ c_{7}^{1}$;\ \ 
$ -1-c^{2}_{7}
$,
$ -1-c^{2}_{7}
$;\ \ 
$ 1+c^{2}_{7}
$)

Factors = $2_{1,2.}^{4,437}\boxtimes 3_{\frac{12}{7},2.862}^{7,768} $

Not pseudo-unitary. 

\vskip 1ex 
\color{grey}

\noindent(7,7). $6_{\frac{37}{7},5.725}^{28,257}$ \irep{46}:\ \ 
$d_i$ = ($1.0$,
$0.554$,
$0.554$,
$1.0$,
$-1.246$,
$-1.246$) 

\vskip 0.7ex
\hangindent=3em \hangafter=1
$D^2= 5.725 = 
10-2  c^{1}_{7}
+4c^{2}_{7}
$

\vskip 0.7ex
\hangindent=3em \hangafter=1
$T = ( 0,
\frac{4}{7},
\frac{9}{28},
\frac{3}{4},
\frac{5}{7},
\frac{13}{28} )
$,

\vskip 0.7ex
\hangindent=3em \hangafter=1
$S$ = ($ 1$,
$ 1+c^{2}_{7}
$,
$ 1+c^{2}_{7}
$,
$ 1$,
$ -c_{7}^{1}$,
$ -c_{7}^{1}$;\ \ 
$ c_{7}^{1}$,
$ c_{7}^{1}$,
$ 1+c^{2}_{7}
$,
$ 1$,
$ 1$;\ \ 
$ -c_{7}^{1}$,
$ -1-c^{2}_{7}
$,
$ 1$,
$ -1$;\ \ 
$ -1$,
$ -c_{7}^{1}$,
$ c_{7}^{1}$;\ \ 
$ -1-c^{2}_{7}
$,
$ -1-c^{2}_{7}
$;\ \ 
$ 1+c^{2}_{7}
$)

Factors = $2_{7,2.}^{4,625}\boxtimes 3_{\frac{44}{7},2.862}^{7,531} $

Not pseudo-unitary. 

\vskip 1ex 
\color{grey}

\noindent(7,8). $6_{\frac{51}{7},5.725}^{28,133}$ \irep{46}:\ \ 
$d_i$ = ($1.0$,
$0.554$,
$0.554$,
$1.0$,
$-1.246$,
$-1.246$) 

\vskip 0.7ex
\hangindent=3em \hangafter=1
$D^2= 5.725 = 
10-2  c^{1}_{7}
+4c^{2}_{7}
$

\vskip 0.7ex
\hangindent=3em \hangafter=1
$T = ( 0,
\frac{4}{7},
\frac{23}{28},
\frac{1}{4},
\frac{5}{7},
\frac{27}{28} )
$,

\vskip 0.7ex
\hangindent=3em \hangafter=1
$S$ = ($ 1$,
$ 1+c^{2}_{7}
$,
$ 1+c^{2}_{7}
$,
$ 1$,
$ -c_{7}^{1}$,
$ -c_{7}^{1}$;\ \ 
$ c_{7}^{1}$,
$ c_{7}^{1}$,
$ 1+c^{2}_{7}
$,
$ 1$,
$ 1$;\ \ 
$ -c_{7}^{1}$,
$ -1-c^{2}_{7}
$,
$ 1$,
$ -1$;\ \ 
$ -1$,
$ -c_{7}^{1}$,
$ c_{7}^{1}$;\ \ 
$ -1-c^{2}_{7}
$,
$ -1-c^{2}_{7}
$;\ \ 
$ 1+c^{2}_{7}
$)

Factors = $2_{1,2.}^{4,437}\boxtimes 3_{\frac{44}{7},2.862}^{7,531} $

Not pseudo-unitary. 

\vskip 1ex 
\color{grey}

\noindent(7,9). $6_{\frac{53}{7},3.682}^{28,318}$ \irep{46}:\ \ 
$d_i$ = ($1.0$,
$0.445$,
$0.445$,
$1.0$,
$-0.801$,
$-0.801$) 

\vskip 0.7ex
\hangindent=3em \hangafter=1
$D^2= 3.682 = 
6-4  c^{1}_{7}
-6  c^{2}_{7}
$

\vskip 0.7ex
\hangindent=3em \hangafter=1
$T = ( 0,
\frac{3}{7},
\frac{5}{28},
\frac{3}{4},
\frac{1}{7},
\frac{25}{28} )
$,

\vskip 0.7ex
\hangindent=3em \hangafter=1
$S$ = ($ 1$,
$ -c_{7}^{2}$,
$ -c_{7}^{2}$,
$ 1$,
$ -c^{1}_{7}
-c^{2}_{7}
$,
$ -c^{1}_{7}
-c^{2}_{7}
$;\ \ 
$ c^{1}_{7}
+c^{2}_{7}
$,
$ c^{1}_{7}
+c^{2}_{7}
$,
$ -c_{7}^{2}$,
$ 1$,
$ 1$;\ \ 
$ -c^{1}_{7}
-c^{2}_{7}
$,
$ c_{7}^{2}$,
$ 1$,
$ -1$;\ \ 
$ -1$,
$ -c^{1}_{7}
-c^{2}_{7}
$,
$ c^{1}_{7}
+c^{2}_{7}
$;\ \ 
$ c_{7}^{2}$,
$ c_{7}^{2}$;\ \ 
$ -c_{7}^{2}$)

Factors = $2_{7,2.}^{4,625}\boxtimes 3_{\frac{4}{7},1.841}^{7,953} $

Not pseudo-unitary. 

\vskip 1ex 
\color{grey}

\noindent(7,10). $6_{\frac{11}{7},3.682}^{28,112}$ \irep{46}:\ \ 
$d_i$ = ($1.0$,
$0.445$,
$0.445$,
$1.0$,
$-0.801$,
$-0.801$) 

\vskip 0.7ex
\hangindent=3em \hangafter=1
$D^2= 3.682 = 
6-4  c^{1}_{7}
-6  c^{2}_{7}
$

\vskip 0.7ex
\hangindent=3em \hangafter=1
$T = ( 0,
\frac{3}{7},
\frac{19}{28},
\frac{1}{4},
\frac{1}{7},
\frac{11}{28} )
$,

\vskip 0.7ex
\hangindent=3em \hangafter=1
$S$ = ($ 1$,
$ -c_{7}^{2}$,
$ -c_{7}^{2}$,
$ 1$,
$ -c^{1}_{7}
-c^{2}_{7}
$,
$ -c^{1}_{7}
-c^{2}_{7}
$;\ \ 
$ c^{1}_{7}
+c^{2}_{7}
$,
$ c^{1}_{7}
+c^{2}_{7}
$,
$ -c_{7}^{2}$,
$ 1$,
$ 1$;\ \ 
$ -c^{1}_{7}
-c^{2}_{7}
$,
$ c_{7}^{2}$,
$ 1$,
$ -1$;\ \ 
$ -1$,
$ -c^{1}_{7}
-c^{2}_{7}
$,
$ c^{1}_{7}
+c^{2}_{7}
$;\ \ 
$ c_{7}^{2}$,
$ c_{7}^{2}$;\ \ 
$ -c_{7}^{2}$)

Factors = $2_{1,2.}^{4,437}\boxtimes 3_{\frac{4}{7},1.841}^{7,953} $

Not pseudo-unitary. 

\vskip 1ex 
\color{grey}

\noindent(7,11). $6_{\frac{45}{7},3.682}^{28,452}$ \irep{46}:\ \ 
$d_i$ = ($1.0$,
$0.445$,
$0.445$,
$1.0$,
$-0.801$,
$-0.801$) 

\vskip 0.7ex
\hangindent=3em \hangafter=1
$D^2= 3.682 = 
6-4  c^{1}_{7}
-6  c^{2}_{7}
$

\vskip 0.7ex
\hangindent=3em \hangafter=1
$T = ( 0,
\frac{4}{7},
\frac{9}{28},
\frac{3}{4},
\frac{6}{7},
\frac{17}{28} )
$,

\vskip 0.7ex
\hangindent=3em \hangafter=1
$S$ = ($ 1$,
$ -c_{7}^{2}$,
$ -c_{7}^{2}$,
$ 1$,
$ -c^{1}_{7}
-c^{2}_{7}
$,
$ -c^{1}_{7}
-c^{2}_{7}
$;\ \ 
$ c^{1}_{7}
+c^{2}_{7}
$,
$ c^{1}_{7}
+c^{2}_{7}
$,
$ -c_{7}^{2}$,
$ 1$,
$ 1$;\ \ 
$ -c^{1}_{7}
-c^{2}_{7}
$,
$ c_{7}^{2}$,
$ 1$,
$ -1$;\ \ 
$ -1$,
$ -c^{1}_{7}
-c^{2}_{7}
$,
$ c^{1}_{7}
+c^{2}_{7}
$;\ \ 
$ c_{7}^{2}$,
$ c_{7}^{2}$;\ \ 
$ -c_{7}^{2}$)

Factors = $2_{7,2.}^{4,625}\boxtimes 3_{\frac{52}{7},1.841}^{7,604} $

Not pseudo-unitary. 

\vskip 1ex 
\color{grey}

\noindent(7,12). $6_{\frac{3}{7},3.682}^{28,411}$ \irep{46}:\ \ 
$d_i$ = ($1.0$,
$0.445$,
$0.445$,
$1.0$,
$-0.801$,
$-0.801$) 

\vskip 0.7ex
\hangindent=3em \hangafter=1
$D^2= 3.682 = 
6-4  c^{1}_{7}
-6  c^{2}_{7}
$

\vskip 0.7ex
\hangindent=3em \hangafter=1
$T = ( 0,
\frac{4}{7},
\frac{23}{28},
\frac{1}{4},
\frac{6}{7},
\frac{3}{28} )
$,

\vskip 0.7ex
\hangindent=3em \hangafter=1
$S$ = ($ 1$,
$ -c_{7}^{2}$,
$ -c_{7}^{2}$,
$ 1$,
$ -c^{1}_{7}
-c^{2}_{7}
$,
$ -c^{1}_{7}
-c^{2}_{7}
$;\ \ 
$ c^{1}_{7}
+c^{2}_{7}
$,
$ c^{1}_{7}
+c^{2}_{7}
$,
$ -c_{7}^{2}$,
$ 1$,
$ 1$;\ \ 
$ -c^{1}_{7}
-c^{2}_{7}
$,
$ c_{7}^{2}$,
$ 1$,
$ -1$;\ \ 
$ -1$,
$ -c^{1}_{7}
-c^{2}_{7}
$,
$ c^{1}_{7}
+c^{2}_{7}
$;\ \ 
$ c_{7}^{2}$,
$ c_{7}^{2}$;\ \ 
$ -c_{7}^{2}$)

Factors = $2_{1,2.}^{4,437}\boxtimes 3_{\frac{52}{7},1.841}^{7,604} $

Not pseudo-unitary. 

\vskip 1ex 
\black

\noindent(8,1). $6_{0,20.}^{10,699}$ \irep{27}:\ \ 
$d_i$ = ($1.0$,
$1.0$,
$2.0$,
$2.0$,
$2.236$,
$2.236$) 

\vskip 0.7ex
\hangindent=3em \hangafter=1
$D^2= 20.0 = 
20$

\vskip 0.7ex
\hangindent=3em \hangafter=1
$T = ( 0,
0,
\frac{1}{5},
\frac{4}{5},
0,
\frac{1}{2} )
$,

\vskip 0.7ex
\hangindent=3em \hangafter=1
$S$ = ($ 1$,
$ 1$,
$ 2$,
$ 2$,
$ \sqrt{5}$,
$ \sqrt{5}$;\ \ 
$ 1$,
$ 2$,
$ 2$,
$ -\sqrt{5}$,
$ -\sqrt{5}$;\ \ 
$ -1-\sqrt{5}$,
$ -1+\sqrt{5}$,
$0$,
$0$;\ \ 
$ -1-\sqrt{5}$,
$0$,
$0$;\ \ 
$ \sqrt{5}$,
$ -\sqrt{5}$;\ \ 
$ \sqrt{5}$)

Prime. 

\vskip 1ex 
\color{grey}

\noindent(8,2). $6_{4,20.}^{10,990}$ \irep{27}:\ \ 
$d_i$ = ($1.0$,
$1.0$,
$2.0$,
$2.0$,
$-2.236$,
$-2.236$) 

\vskip 0.7ex
\hangindent=3em \hangafter=1
$D^2= 20.0 = 
20$

\vskip 0.7ex
\hangindent=3em \hangafter=1
$T = ( 0,
0,
\frac{2}{5},
\frac{3}{5},
0,
\frac{1}{2} )
$,

\vskip 0.7ex
\hangindent=3em \hangafter=1
$S$ = ($ 1$,
$ 1$,
$ 2$,
$ 2$,
$ -\sqrt{5}$,
$ -\sqrt{5}$;\ \ 
$ 1$,
$ 2$,
$ 2$,
$ \sqrt{5}$,
$ \sqrt{5}$;\ \ 
$ -1+\sqrt{5}$,
$ -1-\sqrt{5}$,
$0$,
$0$;\ \ 
$ -1+\sqrt{5}$,
$0$,
$0$;\ \ 
$ -\sqrt{5}$,
$ \sqrt{5}$;\ \ 
$ -\sqrt{5}$)

Prime. 

Pseudo-unitary $\sim$  
$6_{4,20.}^{10,101}$

\vskip 1ex 
\black

\noindent(9,1). $6_{4,20.}^{10,101}$ \irep{27}:\ \ 
$d_i$ = ($1.0$,
$1.0$,
$2.0$,
$2.0$,
$2.236$,
$2.236$) 

\vskip 0.7ex
\hangindent=3em \hangafter=1
$D^2= 20.0 = 
20$

\vskip 0.7ex
\hangindent=3em \hangafter=1
$T = ( 0,
0,
\frac{2}{5},
\frac{3}{5},
0,
\frac{1}{2} )
$,

\vskip 0.7ex
\hangindent=3em \hangafter=1
$S$ = ($ 1$,
$ 1$,
$ 2$,
$ 2$,
$ \sqrt{5}$,
$ \sqrt{5}$;\ \ 
$ 1$,
$ 2$,
$ 2$,
$ -\sqrt{5}$,
$ -\sqrt{5}$;\ \ 
$ -1+\sqrt{5}$,
$ -1-\sqrt{5}$,
$0$,
$0$;\ \ 
$ -1+\sqrt{5}$,
$0$,
$0$;\ \ 
$ -\sqrt{5}$,
$ \sqrt{5}$;\ \ 
$ -\sqrt{5}$)

Prime. 

\vskip 1ex 
\color{grey}

\noindent(9,2). $6_{0,20.}^{10,419}$ \irep{27}:\ \ 
$d_i$ = ($1.0$,
$1.0$,
$2.0$,
$2.0$,
$-2.236$,
$-2.236$) 

\vskip 0.7ex
\hangindent=3em \hangafter=1
$D^2= 20.0 = 
20$

\vskip 0.7ex
\hangindent=3em \hangafter=1
$T = ( 0,
0,
\frac{1}{5},
\frac{4}{5},
0,
\frac{1}{2} )
$,

\vskip 0.7ex
\hangindent=3em \hangafter=1
$S$ = ($ 1$,
$ 1$,
$ 2$,
$ 2$,
$ -\sqrt{5}$,
$ -\sqrt{5}$;\ \ 
$ 1$,
$ 2$,
$ 2$,
$ \sqrt{5}$,
$ \sqrt{5}$;\ \ 
$ -1-\sqrt{5}$,
$ -1+\sqrt{5}$,
$0$,
$0$;\ \ 
$ -1-\sqrt{5}$,
$0$,
$0$;\ \ 
$ \sqrt{5}$,
$ -\sqrt{5}$;\ \ 
$ \sqrt{5}$)

Prime. 

Pseudo-unitary $\sim$  
$6_{0,20.}^{10,699}$

\vskip 1ex 
\black

\noindent(10,1). $6_{0,20.}^{20,139}$ \irep{42}:\ \ 
$d_i$ = ($1.0$,
$1.0$,
$2.0$,
$2.0$,
$2.236$,
$2.236$) 

\vskip 0.7ex
\hangindent=3em \hangafter=1
$D^2= 20.0 = 
20$

\vskip 0.7ex
\hangindent=3em \hangafter=1
$T = ( 0,
0,
\frac{1}{5},
\frac{4}{5},
\frac{1}{4},
\frac{3}{4} )
$,

\vskip 0.7ex
\hangindent=3em \hangafter=1
$S$ = ($ 1$,
$ 1$,
$ 2$,
$ 2$,
$ \sqrt{5}$,
$ \sqrt{5}$;\ \ 
$ 1$,
$ 2$,
$ 2$,
$ -\sqrt{5}$,
$ -\sqrt{5}$;\ \ 
$ -1-\sqrt{5}$,
$ -1+\sqrt{5}$,
$0$,
$0$;\ \ 
$ -1-\sqrt{5}$,
$0$,
$0$;\ \ 
$ -\sqrt{5}$,
$ \sqrt{5}$;\ \ 
$ -\sqrt{5}$)

Prime. 

\vskip 1ex 
\color{grey}

\noindent(10,2). $6_{4,20.}^{20,739}$ \irep{42}:\ \ 
$d_i$ = ($1.0$,
$1.0$,
$2.0$,
$2.0$,
$-2.236$,
$-2.236$) 

\vskip 0.7ex
\hangindent=3em \hangafter=1
$D^2= 20.0 = 
20$

\vskip 0.7ex
\hangindent=3em \hangafter=1
$T = ( 0,
0,
\frac{2}{5},
\frac{3}{5},
\frac{1}{4},
\frac{3}{4} )
$,

\vskip 0.7ex
\hangindent=3em \hangafter=1
$S$ = ($ 1$,
$ 1$,
$ 2$,
$ 2$,
$ -\sqrt{5}$,
$ -\sqrt{5}$;\ \ 
$ 1$,
$ 2$,
$ 2$,
$ \sqrt{5}$,
$ \sqrt{5}$;\ \ 
$ -1+\sqrt{5}$,
$ -1-\sqrt{5}$,
$0$,
$0$;\ \ 
$ -1+\sqrt{5}$,
$0$,
$0$;\ \ 
$ \sqrt{5}$,
$ -\sqrt{5}$;\ \ 
$ \sqrt{5}$)

Prime. 

Pseudo-unitary $\sim$  
$6_{4,20.}^{20,180}$

\vskip 1ex 
\black

\noindent(11,1). $6_{4,20.}^{20,180}$ \irep{42}:\ \ 
$d_i$ = ($1.0$,
$1.0$,
$2.0$,
$2.0$,
$2.236$,
$2.236$) 

\vskip 0.7ex
\hangindent=3em \hangafter=1
$D^2= 20.0 = 
20$

\vskip 0.7ex
\hangindent=3em \hangafter=1
$T = ( 0,
0,
\frac{2}{5},
\frac{3}{5},
\frac{1}{4},
\frac{3}{4} )
$,

\vskip 0.7ex
\hangindent=3em \hangafter=1
$S$ = ($ 1$,
$ 1$,
$ 2$,
$ 2$,
$ \sqrt{5}$,
$ \sqrt{5}$;\ \ 
$ 1$,
$ 2$,
$ 2$,
$ -\sqrt{5}$,
$ -\sqrt{5}$;\ \ 
$ -1+\sqrt{5}$,
$ -1-\sqrt{5}$,
$0$,
$0$;\ \ 
$ -1+\sqrt{5}$,
$0$,
$0$;\ \ 
$ \sqrt{5}$,
$ -\sqrt{5}$;\ \ 
$ \sqrt{5}$)

Prime. 

\vskip 1ex 
\color{grey}

\noindent(11,2). $6_{0,20.}^{20,419}$ \irep{42}:\ \ 
$d_i$ = ($1.0$,
$1.0$,
$2.0$,
$2.0$,
$-2.236$,
$-2.236$) 

\vskip 0.7ex
\hangindent=3em \hangafter=1
$D^2= 20.0 = 
20$

\vskip 0.7ex
\hangindent=3em \hangafter=1
$T = ( 0,
0,
\frac{1}{5},
\frac{4}{5},
\frac{1}{4},
\frac{3}{4} )
$,

\vskip 0.7ex
\hangindent=3em \hangafter=1
$S$ = ($ 1$,
$ 1$,
$ 2$,
$ 2$,
$ -\sqrt{5}$,
$ -\sqrt{5}$;\ \ 
$ 1$,
$ 2$,
$ 2$,
$ \sqrt{5}$,
$ \sqrt{5}$;\ \ 
$ -1-\sqrt{5}$,
$ -1+\sqrt{5}$,
$0$,
$0$;\ \ 
$ -1-\sqrt{5}$,
$0$,
$0$;\ \ 
$ -\sqrt{5}$,
$ \sqrt{5}$;\ \ 
$ -\sqrt{5}$)

Prime. 

Pseudo-unitary $\sim$  
$6_{0,20.}^{20,139}$

\vskip 1ex 
\black

\noindent(12,1). $6_{\frac{58}{35},33.63}^{35,955}$ \irep{47}:\ \ 
$d_i$ = ($1.0$,
$1.618$,
$1.801$,
$2.246$,
$2.915$,
$3.635$) 

\vskip 0.7ex
\hangindent=3em \hangafter=1
$D^2= 33.632 = 
15+3c^{1}_{35}
+2c^{4}_{35}
+6c^{5}_{35}
+3c^{6}_{35}
+3c^{7}_{35}
+2c^{10}_{35}
+2c^{11}_{35}
$

\vskip 0.7ex
\hangindent=3em \hangafter=1
$T = ( 0,
\frac{2}{5},
\frac{1}{7},
\frac{5}{7},
\frac{19}{35},
\frac{4}{35} )
$,

\vskip 0.7ex
\hangindent=3em \hangafter=1
$S$ = ($ 1$,
$ \frac{1+\sqrt{5}}{2}$,
$ -c_{7}^{3}$,
$ \xi_{7}^{3}$,
$ c^{1}_{35}
+c^{6}_{35}
$,
$ c^{1}_{35}
+c^{4}_{35}
+c^{6}_{35}
+c^{11}_{35}
$;\ \ 
$ -1$,
$ c^{1}_{35}
+c^{6}_{35}
$,
$ c^{1}_{35}
+c^{4}_{35}
+c^{6}_{35}
+c^{11}_{35}
$,
$ c_{7}^{3}$,
$ -\xi_{7}^{3}$;\ \ 
$ -\xi_{7}^{3}$,
$ 1$,
$ -c^{1}_{35}
-c^{4}_{35}
-c^{6}_{35}
-c^{11}_{35}
$,
$ \frac{1+\sqrt{5}}{2}$;\ \ 
$ c_{7}^{3}$,
$ \frac{1+\sqrt{5}}{2}$,
$ -c^{1}_{35}
-c^{6}_{35}
$;\ \ 
$ \xi_{7}^{3}$,
$ -1$;\ \ 
$ -c_{7}^{3}$)

Factors = $2_{\frac{14}{5},3.618}^{5,395}\boxtimes 3_{\frac{48}{7},9.295}^{7,790} $

\vskip 1ex 
\color{grey}

\noindent(12,2). $6_{\frac{138}{35},33.63}^{35,363}$ \irep{47}:\ \ 
$d_i$ = ($1.0$,
$1.618$,
$1.801$,
$2.246$,
$2.915$,
$3.635$) 

\vskip 0.7ex
\hangindent=3em \hangafter=1
$D^2= 33.632 = 
15+3c^{1}_{35}
+2c^{4}_{35}
+6c^{5}_{35}
+3c^{6}_{35}
+3c^{7}_{35}
+2c^{10}_{35}
+2c^{11}_{35}
$

\vskip 0.7ex
\hangindent=3em \hangafter=1
$T = ( 0,
\frac{2}{5},
\frac{6}{7},
\frac{2}{7},
\frac{9}{35},
\frac{24}{35} )
$,

\vskip 0.7ex
\hangindent=3em \hangafter=1
$S$ = ($ 1$,
$ \frac{1+\sqrt{5}}{2}$,
$ -c_{7}^{3}$,
$ \xi_{7}^{3}$,
$ c^{1}_{35}
+c^{6}_{35}
$,
$ c^{1}_{35}
+c^{4}_{35}
+c^{6}_{35}
+c^{11}_{35}
$;\ \ 
$ -1$,
$ c^{1}_{35}
+c^{6}_{35}
$,
$ c^{1}_{35}
+c^{4}_{35}
+c^{6}_{35}
+c^{11}_{35}
$,
$ c_{7}^{3}$,
$ -\xi_{7}^{3}$;\ \ 
$ -\xi_{7}^{3}$,
$ 1$,
$ -c^{1}_{35}
-c^{4}_{35}
-c^{6}_{35}
-c^{11}_{35}
$,
$ \frac{1+\sqrt{5}}{2}$;\ \ 
$ c_{7}^{3}$,
$ \frac{1+\sqrt{5}}{2}$,
$ -c^{1}_{35}
-c^{6}_{35}
$;\ \ 
$ \xi_{7}^{3}$,
$ -1$;\ \ 
$ -c_{7}^{3}$)

Factors = $2_{\frac{14}{5},3.618}^{5,395}\boxtimes 3_{\frac{8}{7},9.295}^{7,245} $

\vskip 1ex 
\color{grey}

\noindent(12,3). $6_{\frac{142}{35},33.63}^{35,429}$ \irep{47}:\ \ 
$d_i$ = ($1.0$,
$1.618$,
$1.801$,
$2.246$,
$2.915$,
$3.635$) 

\vskip 0.7ex
\hangindent=3em \hangafter=1
$D^2= 33.632 = 
15+3c^{1}_{35}
+2c^{4}_{35}
+6c^{5}_{35}
+3c^{6}_{35}
+3c^{7}_{35}
+2c^{10}_{35}
+2c^{11}_{35}
$

\vskip 0.7ex
\hangindent=3em \hangafter=1
$T = ( 0,
\frac{3}{5},
\frac{1}{7},
\frac{5}{7},
\frac{26}{35},
\frac{11}{35} )
$,

\vskip 0.7ex
\hangindent=3em \hangafter=1
$S$ = ($ 1$,
$ \frac{1+\sqrt{5}}{2}$,
$ -c_{7}^{3}$,
$ \xi_{7}^{3}$,
$ c^{1}_{35}
+c^{6}_{35}
$,
$ c^{1}_{35}
+c^{4}_{35}
+c^{6}_{35}
+c^{11}_{35}
$;\ \ 
$ -1$,
$ c^{1}_{35}
+c^{6}_{35}
$,
$ c^{1}_{35}
+c^{4}_{35}
+c^{6}_{35}
+c^{11}_{35}
$,
$ c_{7}^{3}$,
$ -\xi_{7}^{3}$;\ \ 
$ -\xi_{7}^{3}$,
$ 1$,
$ -c^{1}_{35}
-c^{4}_{35}
-c^{6}_{35}
-c^{11}_{35}
$,
$ \frac{1+\sqrt{5}}{2}$;\ \ 
$ c_{7}^{3}$,
$ \frac{1+\sqrt{5}}{2}$,
$ -c^{1}_{35}
-c^{6}_{35}
$;\ \ 
$ \xi_{7}^{3}$,
$ -1$;\ \ 
$ -c_{7}^{3}$)

Factors = $2_{\frac{26}{5},3.618}^{5,720}\boxtimes 3_{\frac{48}{7},9.295}^{7,790} $

\vskip 1ex 
\color{grey}

\noindent(12,4). $6_{\frac{222}{35},33.63}^{35,224}$ \irep{47}:\ \ 
$d_i$ = ($1.0$,
$1.618$,
$1.801$,
$2.246$,
$2.915$,
$3.635$) 

\vskip 0.7ex
\hangindent=3em \hangafter=1
$D^2= 33.632 = 
15+3c^{1}_{35}
+2c^{4}_{35}
+6c^{5}_{35}
+3c^{6}_{35}
+3c^{7}_{35}
+2c^{10}_{35}
+2c^{11}_{35}
$

\vskip 0.7ex
\hangindent=3em \hangafter=1
$T = ( 0,
\frac{3}{5},
\frac{6}{7},
\frac{2}{7},
\frac{16}{35},
\frac{31}{35} )
$,

\vskip 0.7ex
\hangindent=3em \hangafter=1
$S$ = ($ 1$,
$ \frac{1+\sqrt{5}}{2}$,
$ -c_{7}^{3}$,
$ \xi_{7}^{3}$,
$ c^{1}_{35}
+c^{6}_{35}
$,
$ c^{1}_{35}
+c^{4}_{35}
+c^{6}_{35}
+c^{11}_{35}
$;\ \ 
$ -1$,
$ c^{1}_{35}
+c^{6}_{35}
$,
$ c^{1}_{35}
+c^{4}_{35}
+c^{6}_{35}
+c^{11}_{35}
$,
$ c_{7}^{3}$,
$ -\xi_{7}^{3}$;\ \ 
$ -\xi_{7}^{3}$,
$ 1$,
$ -c^{1}_{35}
-c^{4}_{35}
-c^{6}_{35}
-c^{11}_{35}
$,
$ \frac{1+\sqrt{5}}{2}$;\ \ 
$ c_{7}^{3}$,
$ \frac{1+\sqrt{5}}{2}$,
$ -c^{1}_{35}
-c^{6}_{35}
$;\ \ 
$ \xi_{7}^{3}$,
$ -1$;\ \ 
$ -c_{7}^{3}$)

Factors = $2_{\frac{26}{5},3.618}^{5,720}\boxtimes 3_{\frac{8}{7},9.295}^{7,245} $

\vskip 1ex 
\color{grey}

\noindent(12,5). $6_{\frac{254}{35},12.84}^{35,615}$ \irep{47}:\ \ 
$d_i$ = ($1.0$,
$1.801$,
$2.246$,
$-0.618$,
$-1.113$,
$-1.388$) 

\vskip 0.7ex
\hangindent=3em \hangafter=1
$D^2= 12.846 = 
15-3  c^{1}_{35}
-2  c^{4}_{35}
+9c^{5}_{35}
-3  c^{6}_{35}
-3  c^{7}_{35}
+3c^{10}_{35}
-2  c^{11}_{35}
$

\vskip 0.7ex
\hangindent=3em \hangafter=1
$T = ( 0,
\frac{1}{7},
\frac{5}{7},
\frac{1}{5},
\frac{12}{35},
\frac{32}{35} )
$,

\vskip 0.7ex
\hangindent=3em \hangafter=1
$S$ = ($ 1$,
$ -c_{7}^{3}$,
$ \xi_{7}^{3}$,
$ \frac{1-\sqrt{5}}{2}$,
$ 1-c^{1}_{35}
+c^{5}_{35}
-c^{6}_{35}
+c^{10}_{35}
$,
$ 1-c^{1}_{35}
-c^{4}_{35}
+c^{5}_{35}
-c^{6}_{35}
-c^{11}_{35}
$;\ \ 
$ -\xi_{7}^{3}$,
$ 1$,
$ 1-c^{1}_{35}
+c^{5}_{35}
-c^{6}_{35}
+c^{10}_{35}
$,
$ -1+c^{1}_{35}
+c^{4}_{35}
-c^{5}_{35}
+c^{6}_{35}
+c^{11}_{35}
$,
$ \frac{1-\sqrt{5}}{2}$;\ \ 
$ c_{7}^{3}$,
$ 1-c^{1}_{35}
-c^{4}_{35}
+c^{5}_{35}
-c^{6}_{35}
-c^{11}_{35}
$,
$ \frac{1-\sqrt{5}}{2}$,
$ -1+c^{1}_{35}
-c^{5}_{35}
+c^{6}_{35}
-c^{10}_{35}
$;\ \ 
$ -1$,
$ c_{7}^{3}$,
$ -\xi_{7}^{3}$;\ \ 
$ \xi_{7}^{3}$,
$ -1$;\ \ 
$ -c_{7}^{3}$)

Factors = $2_{\frac{2}{5},1.381}^{5,120}\boxtimes 3_{\frac{48}{7},9.295}^{7,790} $

Not pseudo-unitary. 

\vskip 1ex 
\color{grey}

\noindent(12,6). $6_{\frac{226}{35},12.84}^{35,105}$ \irep{47}:\ \ 
$d_i$ = ($1.0$,
$1.801$,
$2.246$,
$-0.618$,
$-1.113$,
$-1.388$) 

\vskip 0.7ex
\hangindent=3em \hangafter=1
$D^2= 12.846 = 
15-3  c^{1}_{35}
-2  c^{4}_{35}
+9c^{5}_{35}
-3  c^{6}_{35}
-3  c^{7}_{35}
+3c^{10}_{35}
-2  c^{11}_{35}
$

\vskip 0.7ex
\hangindent=3em \hangafter=1
$T = ( 0,
\frac{1}{7},
\frac{5}{7},
\frac{4}{5},
\frac{33}{35},
\frac{18}{35} )
$,

\vskip 0.7ex
\hangindent=3em \hangafter=1
$S$ = ($ 1$,
$ -c_{7}^{3}$,
$ \xi_{7}^{3}$,
$ \frac{1-\sqrt{5}}{2}$,
$ 1-c^{1}_{35}
+c^{5}_{35}
-c^{6}_{35}
+c^{10}_{35}
$,
$ 1-c^{1}_{35}
-c^{4}_{35}
+c^{5}_{35}
-c^{6}_{35}
-c^{11}_{35}
$;\ \ 
$ -\xi_{7}^{3}$,
$ 1$,
$ 1-c^{1}_{35}
+c^{5}_{35}
-c^{6}_{35}
+c^{10}_{35}
$,
$ -1+c^{1}_{35}
+c^{4}_{35}
-c^{5}_{35}
+c^{6}_{35}
+c^{11}_{35}
$,
$ \frac{1-\sqrt{5}}{2}$;\ \ 
$ c_{7}^{3}$,
$ 1-c^{1}_{35}
-c^{4}_{35}
+c^{5}_{35}
-c^{6}_{35}
-c^{11}_{35}
$,
$ \frac{1-\sqrt{5}}{2}$,
$ -1+c^{1}_{35}
-c^{5}_{35}
+c^{6}_{35}
-c^{10}_{35}
$;\ \ 
$ -1$,
$ c_{7}^{3}$,
$ -\xi_{7}^{3}$;\ \ 
$ \xi_{7}^{3}$,
$ -1$;\ \ 
$ -c_{7}^{3}$)

Factors = $2_{\frac{38}{5},1.381}^{5,491}\boxtimes 3_{\frac{48}{7},9.295}^{7,790} $

Not pseudo-unitary. 

\vskip 1ex 
\color{grey}

\noindent(12,7). $6_{\frac{54}{35},12.84}^{35,669}$ \irep{47}:\ \ 
$d_i$ = ($1.0$,
$1.801$,
$2.246$,
$-0.618$,
$-1.113$,
$-1.388$) 

\vskip 0.7ex
\hangindent=3em \hangafter=1
$D^2= 12.846 = 
15-3  c^{1}_{35}
-2  c^{4}_{35}
+9c^{5}_{35}
-3  c^{6}_{35}
-3  c^{7}_{35}
+3c^{10}_{35}
-2  c^{11}_{35}
$

\vskip 0.7ex
\hangindent=3em \hangafter=1
$T = ( 0,
\frac{6}{7},
\frac{2}{7},
\frac{1}{5},
\frac{2}{35},
\frac{17}{35} )
$,

\vskip 0.7ex
\hangindent=3em \hangafter=1
$S$ = ($ 1$,
$ -c_{7}^{3}$,
$ \xi_{7}^{3}$,
$ \frac{1-\sqrt{5}}{2}$,
$ 1-c^{1}_{35}
+c^{5}_{35}
-c^{6}_{35}
+c^{10}_{35}
$,
$ 1-c^{1}_{35}
-c^{4}_{35}
+c^{5}_{35}
-c^{6}_{35}
-c^{11}_{35}
$;\ \ 
$ -\xi_{7}^{3}$,
$ 1$,
$ 1-c^{1}_{35}
+c^{5}_{35}
-c^{6}_{35}
+c^{10}_{35}
$,
$ -1+c^{1}_{35}
+c^{4}_{35}
-c^{5}_{35}
+c^{6}_{35}
+c^{11}_{35}
$,
$ \frac{1-\sqrt{5}}{2}$;\ \ 
$ c_{7}^{3}$,
$ 1-c^{1}_{35}
-c^{4}_{35}
+c^{5}_{35}
-c^{6}_{35}
-c^{11}_{35}
$,
$ \frac{1-\sqrt{5}}{2}$,
$ -1+c^{1}_{35}
-c^{5}_{35}
+c^{6}_{35}
-c^{10}_{35}
$;\ \ 
$ -1$,
$ c_{7}^{3}$,
$ -\xi_{7}^{3}$;\ \ 
$ \xi_{7}^{3}$,
$ -1$;\ \ 
$ -c_{7}^{3}$)

Factors = $2_{\frac{2}{5},1.381}^{5,120}\boxtimes 3_{\frac{8}{7},9.295}^{7,245} $

Not pseudo-unitary. 

\vskip 1ex 
\color{grey}

\noindent(12,8). $6_{\frac{26}{35},12.84}^{35,138}$ \irep{47}:\ \ 
$d_i$ = ($1.0$,
$1.801$,
$2.246$,
$-0.618$,
$-1.113$,
$-1.388$) 

\vskip 0.7ex
\hangindent=3em \hangafter=1
$D^2= 12.846 = 
15-3  c^{1}_{35}
-2  c^{4}_{35}
+9c^{5}_{35}
-3  c^{6}_{35}
-3  c^{7}_{35}
+3c^{10}_{35}
-2  c^{11}_{35}
$

\vskip 0.7ex
\hangindent=3em \hangafter=1
$T = ( 0,
\frac{6}{7},
\frac{2}{7},
\frac{4}{5},
\frac{23}{35},
\frac{3}{35} )
$,

\vskip 0.7ex
\hangindent=3em \hangafter=1
$S$ = ($ 1$,
$ -c_{7}^{3}$,
$ \xi_{7}^{3}$,
$ \frac{1-\sqrt{5}}{2}$,
$ 1-c^{1}_{35}
+c^{5}_{35}
-c^{6}_{35}
+c^{10}_{35}
$,
$ 1-c^{1}_{35}
-c^{4}_{35}
+c^{5}_{35}
-c^{6}_{35}
-c^{11}_{35}
$;\ \ 
$ -\xi_{7}^{3}$,
$ 1$,
$ 1-c^{1}_{35}
+c^{5}_{35}
-c^{6}_{35}
+c^{10}_{35}
$,
$ -1+c^{1}_{35}
+c^{4}_{35}
-c^{5}_{35}
+c^{6}_{35}
+c^{11}_{35}
$,
$ \frac{1-\sqrt{5}}{2}$;\ \ 
$ c_{7}^{3}$,
$ 1-c^{1}_{35}
-c^{4}_{35}
+c^{5}_{35}
-c^{6}_{35}
-c^{11}_{35}
$,
$ \frac{1-\sqrt{5}}{2}$,
$ -1+c^{1}_{35}
-c^{5}_{35}
+c^{6}_{35}
-c^{10}_{35}
$;\ \ 
$ -1$,
$ c_{7}^{3}$,
$ -\xi_{7}^{3}$;\ \ 
$ \xi_{7}^{3}$,
$ -1$;\ \ 
$ -c_{7}^{3}$)

Factors = $2_{\frac{38}{5},1.381}^{5,491}\boxtimes 3_{\frac{8}{7},9.295}^{7,245} $

Not pseudo-unitary. 

\vskip 1ex 
\color{grey}

\noindent(12,9). $6_{\frac{242}{35},10.35}^{35,120}$ \irep{47}:\ \ 
$d_i$ = ($1.0$,
$0.554$,
$0.897$,
$1.618$,
$-1.246$,
$-2.17$) 

\vskip 0.7ex
\hangindent=3em \hangafter=1
$D^2= 10.358 = 
\frac{35}{\sqrt{35}\mathrm{i}}s^{1}_{35}
+\frac{35}{\sqrt{35}\mathrm{i}}s^{4}_{35}
-\frac{35}{\sqrt{35}\mathrm{i}}  s^{6}_{35}
+\frac{35}{\sqrt{35}\mathrm{i}}s^{11}_{35}
$

\vskip 0.7ex
\hangindent=3em \hangafter=1
$T = ( 0,
\frac{3}{7},
\frac{1}{35},
\frac{3}{5},
\frac{2}{7},
\frac{31}{35} )
$,

\vskip 0.7ex
\hangindent=3em \hangafter=1
$S$ = ($ 1$,
$ 1+c^{2}_{7}
$,
$ 1-c^{4}_{35}
+c^{7}_{35}
-c^{11}_{35}
$,
$ \frac{1+\sqrt{5}}{2}$,
$ -c_{7}^{1}$,
$ 1-c^{1}_{35}
-c^{4}_{35}
-c^{6}_{35}
+c^{7}_{35}
-c^{11}_{35}
$;\ \ 
$ c_{7}^{1}$,
$ -1+c^{1}_{35}
+c^{4}_{35}
+c^{6}_{35}
-c^{7}_{35}
+c^{11}_{35}
$,
$ 1-c^{4}_{35}
+c^{7}_{35}
-c^{11}_{35}
$,
$ 1$,
$ \frac{1+\sqrt{5}}{2}$;\ \ 
$ -c_{7}^{1}$,
$ -1-c^{2}_{7}
$,
$ \frac{1+\sqrt{5}}{2}$,
$ -1$;\ \ 
$ -1$,
$ 1-c^{1}_{35}
-c^{4}_{35}
-c^{6}_{35}
+c^{7}_{35}
-c^{11}_{35}
$,
$ c_{7}^{1}$;\ \ 
$ -1-c^{2}_{7}
$,
$ -1+c^{4}_{35}
-c^{7}_{35}
+c^{11}_{35}
$;\ \ 
$ 1+c^{2}_{7}
$)

Factors = $2_{\frac{26}{5},3.618}^{5,720}\boxtimes 3_{\frac{12}{7},2.862}^{7,768} $

Not pseudo-unitary. 

\vskip 1ex 
\color{grey}

\noindent(12,10). $6_{\frac{158}{35},10.35}^{35,274}$ \irep{47}:\ \ 
$d_i$ = ($1.0$,
$0.554$,
$0.897$,
$1.618$,
$-1.246$,
$-2.17$) 

\vskip 0.7ex
\hangindent=3em \hangafter=1
$D^2= 10.358 = 
\frac{35}{\sqrt{35}\mathrm{i}}s^{1}_{35}
+\frac{35}{\sqrt{35}\mathrm{i}}s^{4}_{35}
-\frac{35}{\sqrt{35}\mathrm{i}}  s^{6}_{35}
+\frac{35}{\sqrt{35}\mathrm{i}}s^{11}_{35}
$

\vskip 0.7ex
\hangindent=3em \hangafter=1
$T = ( 0,
\frac{3}{7},
\frac{29}{35},
\frac{2}{5},
\frac{2}{7},
\frac{24}{35} )
$,

\vskip 0.7ex
\hangindent=3em \hangafter=1
$S$ = ($ 1$,
$ 1+c^{2}_{7}
$,
$ 1-c^{4}_{35}
+c^{7}_{35}
-c^{11}_{35}
$,
$ \frac{1+\sqrt{5}}{2}$,
$ -c_{7}^{1}$,
$ 1-c^{1}_{35}
-c^{4}_{35}
-c^{6}_{35}
+c^{7}_{35}
-c^{11}_{35}
$;\ \ 
$ c_{7}^{1}$,
$ -1+c^{1}_{35}
+c^{4}_{35}
+c^{6}_{35}
-c^{7}_{35}
+c^{11}_{35}
$,
$ 1-c^{4}_{35}
+c^{7}_{35}
-c^{11}_{35}
$,
$ 1$,
$ \frac{1+\sqrt{5}}{2}$;\ \ 
$ -c_{7}^{1}$,
$ -1-c^{2}_{7}
$,
$ \frac{1+\sqrt{5}}{2}$,
$ -1$;\ \ 
$ -1$,
$ 1-c^{1}_{35}
-c^{4}_{35}
-c^{6}_{35}
+c^{7}_{35}
-c^{11}_{35}
$,
$ c_{7}^{1}$;\ \ 
$ -1-c^{2}_{7}
$,
$ -1+c^{4}_{35}
-c^{7}_{35}
+c^{11}_{35}
$;\ \ 
$ 1+c^{2}_{7}
$)

Factors = $2_{\frac{14}{5},3.618}^{5,395}\boxtimes 3_{\frac{12}{7},2.862}^{7,768} $

Not pseudo-unitary. 

\vskip 1ex 
\color{grey}

\noindent(12,11). $6_{\frac{122}{35},10.35}^{35,246}$ \irep{47}:\ \ 
$d_i$ = ($1.0$,
$0.554$,
$0.897$,
$1.618$,
$-1.246$,
$-2.17$) 

\vskip 0.7ex
\hangindent=3em \hangafter=1
$D^2= 10.358 = 
\frac{35}{\sqrt{35}\mathrm{i}}s^{1}_{35}
+\frac{35}{\sqrt{35}\mathrm{i}}s^{4}_{35}
-\frac{35}{\sqrt{35}\mathrm{i}}  s^{6}_{35}
+\frac{35}{\sqrt{35}\mathrm{i}}s^{11}_{35}
$

\vskip 0.7ex
\hangindent=3em \hangafter=1
$T = ( 0,
\frac{4}{7},
\frac{6}{35},
\frac{3}{5},
\frac{5}{7},
\frac{11}{35} )
$,

\vskip 0.7ex
\hangindent=3em \hangafter=1
$S$ = ($ 1$,
$ 1+c^{2}_{7}
$,
$ 1-c^{4}_{35}
+c^{7}_{35}
-c^{11}_{35}
$,
$ \frac{1+\sqrt{5}}{2}$,
$ -c_{7}^{1}$,
$ 1-c^{1}_{35}
-c^{4}_{35}
-c^{6}_{35}
+c^{7}_{35}
-c^{11}_{35}
$;\ \ 
$ c_{7}^{1}$,
$ -1+c^{1}_{35}
+c^{4}_{35}
+c^{6}_{35}
-c^{7}_{35}
+c^{11}_{35}
$,
$ 1-c^{4}_{35}
+c^{7}_{35}
-c^{11}_{35}
$,
$ 1$,
$ \frac{1+\sqrt{5}}{2}$;\ \ 
$ -c_{7}^{1}$,
$ -1-c^{2}_{7}
$,
$ \frac{1+\sqrt{5}}{2}$,
$ -1$;\ \ 
$ -1$,
$ 1-c^{1}_{35}
-c^{4}_{35}
-c^{6}_{35}
+c^{7}_{35}
-c^{11}_{35}
$,
$ c_{7}^{1}$;\ \ 
$ -1-c^{2}_{7}
$,
$ -1+c^{4}_{35}
-c^{7}_{35}
+c^{11}_{35}
$;\ \ 
$ 1+c^{2}_{7}
$)

Factors = $2_{\frac{26}{5},3.618}^{5,720}\boxtimes 3_{\frac{44}{7},2.862}^{7,531} $

Not pseudo-unitary. 

\vskip 1ex 
\color{grey}

\noindent(12,12). $6_{\frac{38}{35},10.35}^{35,423}$ \irep{47}:\ \ 
$d_i$ = ($1.0$,
$0.554$,
$0.897$,
$1.618$,
$-1.246$,
$-2.17$) 

\vskip 0.7ex
\hangindent=3em \hangafter=1
$D^2= 10.358 = 
\frac{35}{\sqrt{35}\mathrm{i}}s^{1}_{35}
+\frac{35}{\sqrt{35}\mathrm{i}}s^{4}_{35}
-\frac{35}{\sqrt{35}\mathrm{i}}  s^{6}_{35}
+\frac{35}{\sqrt{35}\mathrm{i}}s^{11}_{35}
$

\vskip 0.7ex
\hangindent=3em \hangafter=1
$T = ( 0,
\frac{4}{7},
\frac{34}{35},
\frac{2}{5},
\frac{5}{7},
\frac{4}{35} )
$,

\vskip 0.7ex
\hangindent=3em \hangafter=1
$S$ = ($ 1$,
$ 1+c^{2}_{7}
$,
$ 1-c^{4}_{35}
+c^{7}_{35}
-c^{11}_{35}
$,
$ \frac{1+\sqrt{5}}{2}$,
$ -c_{7}^{1}$,
$ 1-c^{1}_{35}
-c^{4}_{35}
-c^{6}_{35}
+c^{7}_{35}
-c^{11}_{35}
$;\ \ 
$ c_{7}^{1}$,
$ -1+c^{1}_{35}
+c^{4}_{35}
+c^{6}_{35}
-c^{7}_{35}
+c^{11}_{35}
$,
$ 1-c^{4}_{35}
+c^{7}_{35}
-c^{11}_{35}
$,
$ 1$,
$ \frac{1+\sqrt{5}}{2}$;\ \ 
$ -c_{7}^{1}$,
$ -1-c^{2}_{7}
$,
$ \frac{1+\sqrt{5}}{2}$,
$ -1$;\ \ 
$ -1$,
$ 1-c^{1}_{35}
-c^{4}_{35}
-c^{6}_{35}
+c^{7}_{35}
-c^{11}_{35}
$,
$ c_{7}^{1}$;\ \ 
$ -1-c^{2}_{7}
$,
$ -1+c^{4}_{35}
-c^{7}_{35}
+c^{11}_{35}
$;\ \ 
$ 1+c^{2}_{7}
$)

Factors = $2_{\frac{14}{5},3.618}^{5,395}\boxtimes 3_{\frac{44}{7},2.862}^{7,531} $

Not pseudo-unitary. 

\vskip 1ex 
\color{grey}

\noindent(12,13). $6_{\frac{202}{35},6.661}^{35,488}$ \irep{47}:\ \ 
$d_i$ = ($1.0$,
$0.445$,
$0.720$,
$1.618$,
$-0.801$,
$-1.297$) 

\vskip 0.7ex
\hangindent=3em \hangafter=1
$D^2= 6.661 = 
11-2  c^{1}_{35}
+c^{4}_{35}
-4  c^{5}_{35}
-2  c^{6}_{35}
+5c^{7}_{35}
-6  c^{10}_{35}
+c^{11}_{35}
$

\vskip 0.7ex
\hangindent=3em \hangafter=1
$T = ( 0,
\frac{3}{7},
\frac{1}{35},
\frac{3}{5},
\frac{1}{7},
\frac{26}{35} )
$,

\vskip 0.7ex
\hangindent=3em \hangafter=1
$S$ = ($ 1$,
$ -c_{7}^{2}$,
$ c^{4}_{35}
+c^{11}_{35}
$,
$ \frac{1+\sqrt{5}}{2}$,
$ -c^{1}_{7}
-c^{2}_{7}
$,
$ 1-c^{1}_{35}
-c^{6}_{35}
+c^{7}_{35}
$;\ \ 
$ c^{1}_{7}
+c^{2}_{7}
$,
$ -1+c^{1}_{35}
+c^{6}_{35}
-c^{7}_{35}
$,
$ c^{4}_{35}
+c^{11}_{35}
$,
$ 1$,
$ \frac{1+\sqrt{5}}{2}$;\ \ 
$ -c^{1}_{7}
-c^{2}_{7}
$,
$ c_{7}^{2}$,
$ \frac{1+\sqrt{5}}{2}$,
$ -1$;\ \ 
$ -1$,
$ 1-c^{1}_{35}
-c^{6}_{35}
+c^{7}_{35}
$,
$ c^{1}_{7}
+c^{2}_{7}
$;\ \ 
$ c_{7}^{2}$,
$ -c^{4}_{35}
-c^{11}_{35}
$;\ \ 
$ -c_{7}^{2}$)

Factors = $2_{\frac{26}{5},3.618}^{5,720}\boxtimes 3_{\frac{4}{7},1.841}^{7,953} $

Not pseudo-unitary. 

\vskip 1ex 
\color{grey}

\noindent(12,14). $6_{\frac{118}{35},6.661}^{35,153}$ \irep{47}:\ \ 
$d_i$ = ($1.0$,
$0.445$,
$0.720$,
$1.618$,
$-0.801$,
$-1.297$) 

\vskip 0.7ex
\hangindent=3em \hangafter=1
$D^2= 6.661 = 
11-2  c^{1}_{35}
+c^{4}_{35}
-4  c^{5}_{35}
-2  c^{6}_{35}
+5c^{7}_{35}
-6  c^{10}_{35}
+c^{11}_{35}
$

\vskip 0.7ex
\hangindent=3em \hangafter=1
$T = ( 0,
\frac{3}{7},
\frac{29}{35},
\frac{2}{5},
\frac{1}{7},
\frac{19}{35} )
$,

\vskip 0.7ex
\hangindent=3em \hangafter=1
$S$ = ($ 1$,
$ -c_{7}^{2}$,
$ c^{4}_{35}
+c^{11}_{35}
$,
$ \frac{1+\sqrt{5}}{2}$,
$ -c^{1}_{7}
-c^{2}_{7}
$,
$ 1-c^{1}_{35}
-c^{6}_{35}
+c^{7}_{35}
$;\ \ 
$ c^{1}_{7}
+c^{2}_{7}
$,
$ -1+c^{1}_{35}
+c^{6}_{35}
-c^{7}_{35}
$,
$ c^{4}_{35}
+c^{11}_{35}
$,
$ 1$,
$ \frac{1+\sqrt{5}}{2}$;\ \ 
$ -c^{1}_{7}
-c^{2}_{7}
$,
$ c_{7}^{2}$,
$ \frac{1+\sqrt{5}}{2}$,
$ -1$;\ \ 
$ -1$,
$ 1-c^{1}_{35}
-c^{6}_{35}
+c^{7}_{35}
$,
$ c^{1}_{7}
+c^{2}_{7}
$;\ \ 
$ c_{7}^{2}$,
$ -c^{4}_{35}
-c^{11}_{35}
$;\ \ 
$ -c_{7}^{2}$)

Factors = $2_{\frac{14}{5},3.618}^{5,395}\boxtimes 3_{\frac{4}{7},1.841}^{7,953} $

Not pseudo-unitary. 

\vskip 1ex 
\color{grey}

\noindent(12,15). $6_{\frac{162}{35},6.661}^{35,532}$ \irep{47}:\ \ 
$d_i$ = ($1.0$,
$0.445$,
$0.720$,
$1.618$,
$-0.801$,
$-1.297$) 

\vskip 0.7ex
\hangindent=3em \hangafter=1
$D^2= 6.661 = 
11-2  c^{1}_{35}
+c^{4}_{35}
-4  c^{5}_{35}
-2  c^{6}_{35}
+5c^{7}_{35}
-6  c^{10}_{35}
+c^{11}_{35}
$

\vskip 0.7ex
\hangindent=3em \hangafter=1
$T = ( 0,
\frac{4}{7},
\frac{6}{35},
\frac{3}{5},
\frac{6}{7},
\frac{16}{35} )
$,

\vskip 0.7ex
\hangindent=3em \hangafter=1
$S$ = ($ 1$,
$ -c_{7}^{2}$,
$ c^{4}_{35}
+c^{11}_{35}
$,
$ \frac{1+\sqrt{5}}{2}$,
$ -c^{1}_{7}
-c^{2}_{7}
$,
$ 1-c^{1}_{35}
-c^{6}_{35}
+c^{7}_{35}
$;\ \ 
$ c^{1}_{7}
+c^{2}_{7}
$,
$ -1+c^{1}_{35}
+c^{6}_{35}
-c^{7}_{35}
$,
$ c^{4}_{35}
+c^{11}_{35}
$,
$ 1$,
$ \frac{1+\sqrt{5}}{2}$;\ \ 
$ -c^{1}_{7}
-c^{2}_{7}
$,
$ c_{7}^{2}$,
$ \frac{1+\sqrt{5}}{2}$,
$ -1$;\ \ 
$ -1$,
$ 1-c^{1}_{35}
-c^{6}_{35}
+c^{7}_{35}
$,
$ c^{1}_{7}
+c^{2}_{7}
$;\ \ 
$ c_{7}^{2}$,
$ -c^{4}_{35}
-c^{11}_{35}
$;\ \ 
$ -c_{7}^{2}$)

Factors = $2_{\frac{26}{5},3.618}^{5,720}\boxtimes 3_{\frac{52}{7},1.841}^{7,604} $

Not pseudo-unitary. 

\vskip 1ex 
\color{grey}

\noindent(12,16). $6_{\frac{78}{35},6.661}^{35,121}$ \irep{47}:\ \ 
$d_i$ = ($1.0$,
$0.445$,
$0.720$,
$1.618$,
$-0.801$,
$-1.297$) 

\vskip 0.7ex
\hangindent=3em \hangafter=1
$D^2= 6.661 = 
11-2  c^{1}_{35}
+c^{4}_{35}
-4  c^{5}_{35}
-2  c^{6}_{35}
+5c^{7}_{35}
-6  c^{10}_{35}
+c^{11}_{35}
$

\vskip 0.7ex
\hangindent=3em \hangafter=1
$T = ( 0,
\frac{4}{7},
\frac{34}{35},
\frac{2}{5},
\frac{6}{7},
\frac{9}{35} )
$,

\vskip 0.7ex
\hangindent=3em \hangafter=1
$S$ = ($ 1$,
$ -c_{7}^{2}$,
$ c^{4}_{35}
+c^{11}_{35}
$,
$ \frac{1+\sqrt{5}}{2}$,
$ -c^{1}_{7}
-c^{2}_{7}
$,
$ 1-c^{1}_{35}
-c^{6}_{35}
+c^{7}_{35}
$;\ \ 
$ c^{1}_{7}
+c^{2}_{7}
$,
$ -1+c^{1}_{35}
+c^{6}_{35}
-c^{7}_{35}
$,
$ c^{4}_{35}
+c^{11}_{35}
$,
$ 1$,
$ \frac{1+\sqrt{5}}{2}$;\ \ 
$ -c^{1}_{7}
-c^{2}_{7}
$,
$ c_{7}^{2}$,
$ \frac{1+\sqrt{5}}{2}$,
$ -1$;\ \ 
$ -1$,
$ 1-c^{1}_{35}
-c^{6}_{35}
+c^{7}_{35}
$,
$ c^{1}_{7}
+c^{2}_{7}
$;\ \ 
$ c_{7}^{2}$,
$ -c^{4}_{35}
-c^{11}_{35}
$;\ \ 
$ -c_{7}^{2}$)

Factors = $2_{\frac{14}{5},3.618}^{5,395}\boxtimes 3_{\frac{52}{7},1.841}^{7,604} $

Not pseudo-unitary. 

\vskip 1ex 
\color{grey}

\noindent(12,17). $6_{\frac{46}{35},3.956}^{35,437}$ \irep{47}:\ \ 
$d_i$ = ($1.0$,
$0.554$,
$0.770$,
$-0.342$,
$-0.618$,
$-1.246$) 

\vskip 0.7ex
\hangindent=3em \hangafter=1
$D^2= 3.956 = 
9+c^{1}_{35}
+3c^{4}_{35}
-3  c^{5}_{35}
+c^{6}_{35}
-6  c^{7}_{35}
+6c^{10}_{35}
+3c^{11}_{35}
$

\vskip 0.7ex
\hangindent=3em \hangafter=1
$T = ( 0,
\frac{3}{7},
\frac{3}{35},
\frac{8}{35},
\frac{4}{5},
\frac{2}{7} )
$,

\vskip 0.7ex
\hangindent=3em \hangafter=1
$S$ = ($ 1$,
$ 1+c^{2}_{7}
$,
$ -1+c^{1}_{35}
+c^{4}_{35}
-c^{5}_{35}
+c^{6}_{35}
-c^{7}_{35}
+c^{11}_{35}
$,
$ c^{4}_{35}
-c^{7}_{35}
+c^{10}_{35}
+c^{11}_{35}
$,
$ \frac{1-\sqrt{5}}{2}$,
$ -c_{7}^{1}$;\ \ 
$ c_{7}^{1}$,
$ \frac{1-\sqrt{5}}{2}$,
$ 1-c^{1}_{35}
-c^{4}_{35}
+c^{5}_{35}
-c^{6}_{35}
+c^{7}_{35}
-c^{11}_{35}
$,
$ c^{4}_{35}
-c^{7}_{35}
+c^{10}_{35}
+c^{11}_{35}
$,
$ 1$;\ \ 
$ 1+c^{2}_{7}
$,
$ -1$,
$ c_{7}^{1}$,
$ -c^{4}_{35}
+c^{7}_{35}
-c^{10}_{35}
-c^{11}_{35}
$;\ \ 
$ -c_{7}^{1}$,
$ -1-c^{2}_{7}
$,
$ \frac{1-\sqrt{5}}{2}$;\ \ 
$ -1$,
$ -1+c^{1}_{35}
+c^{4}_{35}
-c^{5}_{35}
+c^{6}_{35}
-c^{7}_{35}
+c^{11}_{35}
$;\ \ 
$ -1-c^{2}_{7}
$)

Factors = $2_{\frac{38}{5},1.381}^{5,491}\boxtimes 3_{\frac{12}{7},2.862}^{7,768} $

Not pseudo-unitary. 

\vskip 1ex 
\color{grey}

\noindent(12,18). $6_{\frac{74}{35},3.956}^{35,769}$ \irep{47}:\ \ 
$d_i$ = ($1.0$,
$0.554$,
$0.770$,
$-0.342$,
$-0.618$,
$-1.246$) 

\vskip 0.7ex
\hangindent=3em \hangafter=1
$D^2= 3.956 = 
9+c^{1}_{35}
+3c^{4}_{35}
-3  c^{5}_{35}
+c^{6}_{35}
-6  c^{7}_{35}
+6c^{10}_{35}
+3c^{11}_{35}
$

\vskip 0.7ex
\hangindent=3em \hangafter=1
$T = ( 0,
\frac{3}{7},
\frac{17}{35},
\frac{22}{35},
\frac{1}{5},
\frac{2}{7} )
$,

\vskip 0.7ex
\hangindent=3em \hangafter=1
$S$ = ($ 1$,
$ 1+c^{2}_{7}
$,
$ -1+c^{1}_{35}
+c^{4}_{35}
-c^{5}_{35}
+c^{6}_{35}
-c^{7}_{35}
+c^{11}_{35}
$,
$ c^{4}_{35}
-c^{7}_{35}
+c^{10}_{35}
+c^{11}_{35}
$,
$ \frac{1-\sqrt{5}}{2}$,
$ -c_{7}^{1}$;\ \ 
$ c_{7}^{1}$,
$ \frac{1-\sqrt{5}}{2}$,
$ 1-c^{1}_{35}
-c^{4}_{35}
+c^{5}_{35}
-c^{6}_{35}
+c^{7}_{35}
-c^{11}_{35}
$,
$ c^{4}_{35}
-c^{7}_{35}
+c^{10}_{35}
+c^{11}_{35}
$,
$ 1$;\ \ 
$ 1+c^{2}_{7}
$,
$ -1$,
$ c_{7}^{1}$,
$ -c^{4}_{35}
+c^{7}_{35}
-c^{10}_{35}
-c^{11}_{35}
$;\ \ 
$ -c_{7}^{1}$,
$ -1-c^{2}_{7}
$,
$ \frac{1-\sqrt{5}}{2}$;\ \ 
$ -1$,
$ -1+c^{1}_{35}
+c^{4}_{35}
-c^{5}_{35}
+c^{6}_{35}
-c^{7}_{35}
+c^{11}_{35}
$;\ \ 
$ -1-c^{2}_{7}
$)

Factors = $2_{\frac{2}{5},1.381}^{5,120}\boxtimes 3_{\frac{12}{7},2.862}^{7,768} $

Not pseudo-unitary. 

\vskip 1ex 
\color{grey}

\noindent(12,19). $6_{\frac{206}{35},3.956}^{35,723}$ \irep{47}:\ \ 
$d_i$ = ($1.0$,
$0.554$,
$0.770$,
$-0.342$,
$-0.618$,
$-1.246$) 

\vskip 0.7ex
\hangindent=3em \hangafter=1
$D^2= 3.956 = 
9+c^{1}_{35}
+3c^{4}_{35}
-3  c^{5}_{35}
+c^{6}_{35}
-6  c^{7}_{35}
+6c^{10}_{35}
+3c^{11}_{35}
$

\vskip 0.7ex
\hangindent=3em \hangafter=1
$T = ( 0,
\frac{4}{7},
\frac{18}{35},
\frac{13}{35},
\frac{4}{5},
\frac{5}{7} )
$,

\vskip 0.7ex
\hangindent=3em \hangafter=1
$S$ = ($ 1$,
$ 1+c^{2}_{7}
$,
$ -1+c^{1}_{35}
+c^{4}_{35}
-c^{5}_{35}
+c^{6}_{35}
-c^{7}_{35}
+c^{11}_{35}
$,
$ c^{4}_{35}
-c^{7}_{35}
+c^{10}_{35}
+c^{11}_{35}
$,
$ \frac{1-\sqrt{5}}{2}$,
$ -c_{7}^{1}$;\ \ 
$ c_{7}^{1}$,
$ \frac{1-\sqrt{5}}{2}$,
$ 1-c^{1}_{35}
-c^{4}_{35}
+c^{5}_{35}
-c^{6}_{35}
+c^{7}_{35}
-c^{11}_{35}
$,
$ c^{4}_{35}
-c^{7}_{35}
+c^{10}_{35}
+c^{11}_{35}
$,
$ 1$;\ \ 
$ 1+c^{2}_{7}
$,
$ -1$,
$ c_{7}^{1}$,
$ -c^{4}_{35}
+c^{7}_{35}
-c^{10}_{35}
-c^{11}_{35}
$;\ \ 
$ -c_{7}^{1}$,
$ -1-c^{2}_{7}
$,
$ \frac{1-\sqrt{5}}{2}$;\ \ 
$ -1$,
$ -1+c^{1}_{35}
+c^{4}_{35}
-c^{5}_{35}
+c^{6}_{35}
-c^{7}_{35}
+c^{11}_{35}
$;\ \ 
$ -1-c^{2}_{7}
$)

Factors = $2_{\frac{38}{5},1.381}^{5,491}\boxtimes 3_{\frac{44}{7},2.862}^{7,531} $

Not pseudo-unitary. 

\vskip 1ex 
\color{grey}

\noindent(12,20). $6_{\frac{234}{35},3.956}^{35,690}$ \irep{47}:\ \ 
$d_i$ = ($1.0$,
$0.554$,
$0.770$,
$-0.342$,
$-0.618$,
$-1.246$) 

\vskip 0.7ex
\hangindent=3em \hangafter=1
$D^2= 3.956 = 
9+c^{1}_{35}
+3c^{4}_{35}
-3  c^{5}_{35}
+c^{6}_{35}
-6  c^{7}_{35}
+6c^{10}_{35}
+3c^{11}_{35}
$

\vskip 0.7ex
\hangindent=3em \hangafter=1
$T = ( 0,
\frac{4}{7},
\frac{32}{35},
\frac{27}{35},
\frac{1}{5},
\frac{5}{7} )
$,

\vskip 0.7ex
\hangindent=3em \hangafter=1
$S$ = ($ 1$,
$ 1+c^{2}_{7}
$,
$ -1+c^{1}_{35}
+c^{4}_{35}
-c^{5}_{35}
+c^{6}_{35}
-c^{7}_{35}
+c^{11}_{35}
$,
$ c^{4}_{35}
-c^{7}_{35}
+c^{10}_{35}
+c^{11}_{35}
$,
$ \frac{1-\sqrt{5}}{2}$,
$ -c_{7}^{1}$;\ \ 
$ c_{7}^{1}$,
$ \frac{1-\sqrt{5}}{2}$,
$ 1-c^{1}_{35}
-c^{4}_{35}
+c^{5}_{35}
-c^{6}_{35}
+c^{7}_{35}
-c^{11}_{35}
$,
$ c^{4}_{35}
-c^{7}_{35}
+c^{10}_{35}
+c^{11}_{35}
$,
$ 1$;\ \ 
$ 1+c^{2}_{7}
$,
$ -1$,
$ c_{7}^{1}$,
$ -c^{4}_{35}
+c^{7}_{35}
-c^{10}_{35}
-c^{11}_{35}
$;\ \ 
$ -c_{7}^{1}$,
$ -1-c^{2}_{7}
$,
$ \frac{1-\sqrt{5}}{2}$;\ \ 
$ -1$,
$ -1+c^{1}_{35}
+c^{4}_{35}
-c^{5}_{35}
+c^{6}_{35}
-c^{7}_{35}
+c^{11}_{35}
$;\ \ 
$ -1-c^{2}_{7}
$)

Factors = $2_{\frac{2}{5},1.381}^{5,120}\boxtimes 3_{\frac{44}{7},2.862}^{7,531} $

Not pseudo-unitary. 

\vskip 1ex 
\color{grey}

\noindent(12,21). $6_{\frac{34}{35},2.544}^{35,711}$ \irep{47}:\ \ 
$d_i$ = ($1.0$,
$0.445$,
$0.495$,
$-0.275$,
$-0.618$,
$-0.801$) 

\vskip 0.7ex
\hangindent=3em \hangafter=1
$D^2= 2.544 = 
4+2c^{1}_{35}
-c^{4}_{35}
-6  c^{5}_{35}
+2c^{6}_{35}
-5  c^{7}_{35}
-9  c^{10}_{35}
-c^{11}_{35}
$

\vskip 0.7ex
\hangindent=3em \hangafter=1
$T = ( 0,
\frac{3}{7},
\frac{12}{35},
\frac{22}{35},
\frac{1}{5},
\frac{1}{7} )
$,

\vskip 0.7ex
\hangindent=3em \hangafter=1
$S$ = ($ 1$,
$ -c_{7}^{2}$,
$ -1+c^{1}_{35}
-c^{5}_{35}
+c^{6}_{35}
-c^{7}_{35}
-c^{10}_{35}
$,
$ -c^{4}_{35}
-c^{10}_{35}
-c^{11}_{35}
$,
$ \frac{1-\sqrt{5}}{2}$,
$ -c^{1}_{7}
-c^{2}_{7}
$;\ \ 
$ c^{1}_{7}
+c^{2}_{7}
$,
$ \frac{1-\sqrt{5}}{2}$,
$ 1-c^{1}_{35}
+c^{5}_{35}
-c^{6}_{35}
+c^{7}_{35}
+c^{10}_{35}
$,
$ -c^{4}_{35}
-c^{10}_{35}
-c^{11}_{35}
$,
$ 1$;\ \ 
$ -c_{7}^{2}$,
$ -1$,
$ c^{1}_{7}
+c^{2}_{7}
$,
$ c^{4}_{35}
+c^{10}_{35}
+c^{11}_{35}
$;\ \ 
$ -c^{1}_{7}
-c^{2}_{7}
$,
$ c_{7}^{2}$,
$ \frac{1-\sqrt{5}}{2}$;\ \ 
$ -1$,
$ -1+c^{1}_{35}
-c^{5}_{35}
+c^{6}_{35}
-c^{7}_{35}
-c^{10}_{35}
$;\ \ 
$ c_{7}^{2}$)

Factors = $2_{\frac{2}{5},1.381}^{5,120}\boxtimes 3_{\frac{4}{7},1.841}^{7,953} $

Not pseudo-unitary. 

\vskip 1ex 
\color{grey}

\noindent(12,22). $6_{\frac{6}{35},2.544}^{35,346}$ \irep{47}:\ \ 
$d_i$ = ($1.0$,
$0.445$,
$0.495$,
$-0.275$,
$-0.618$,
$-0.801$) 

\vskip 0.7ex
\hangindent=3em \hangafter=1
$D^2= 2.544 = 
4+2c^{1}_{35}
-c^{4}_{35}
-6  c^{5}_{35}
+2c^{6}_{35}
-5  c^{7}_{35}
-9  c^{10}_{35}
-c^{11}_{35}
$

\vskip 0.7ex
\hangindent=3em \hangafter=1
$T = ( 0,
\frac{3}{7},
\frac{33}{35},
\frac{8}{35},
\frac{4}{5},
\frac{1}{7} )
$,

\vskip 0.7ex
\hangindent=3em \hangafter=1
$S$ = ($ 1$,
$ -c_{7}^{2}$,
$ -1+c^{1}_{35}
-c^{5}_{35}
+c^{6}_{35}
-c^{7}_{35}
-c^{10}_{35}
$,
$ -c^{4}_{35}
-c^{10}_{35}
-c^{11}_{35}
$,
$ \frac{1-\sqrt{5}}{2}$,
$ -c^{1}_{7}
-c^{2}_{7}
$;\ \ 
$ c^{1}_{7}
+c^{2}_{7}
$,
$ \frac{1-\sqrt{5}}{2}$,
$ 1-c^{1}_{35}
+c^{5}_{35}
-c^{6}_{35}
+c^{7}_{35}
+c^{10}_{35}
$,
$ -c^{4}_{35}
-c^{10}_{35}
-c^{11}_{35}
$,
$ 1$;\ \ 
$ -c_{7}^{2}$,
$ -1$,
$ c^{1}_{7}
+c^{2}_{7}
$,
$ c^{4}_{35}
+c^{10}_{35}
+c^{11}_{35}
$;\ \ 
$ -c^{1}_{7}
-c^{2}_{7}
$,
$ c_{7}^{2}$,
$ \frac{1-\sqrt{5}}{2}$;\ \ 
$ -1$,
$ -1+c^{1}_{35}
-c^{5}_{35}
+c^{6}_{35}
-c^{7}_{35}
-c^{10}_{35}
$;\ \ 
$ c_{7}^{2}$)

Factors = $2_{\frac{38}{5},1.381}^{5,491}\boxtimes 3_{\frac{4}{7},1.841}^{7,953} $

Not pseudo-unitary. 

\vskip 1ex 
\color{grey}

\noindent(12,23). $6_{\frac{274}{35},2.544}^{35,692}$ \irep{47}:\ \ 
$d_i$ = ($1.0$,
$0.445$,
$0.495$,
$-0.275$,
$-0.618$,
$-0.801$) 

\vskip 0.7ex
\hangindent=3em \hangafter=1
$D^2= 2.544 = 
4+2c^{1}_{35}
-c^{4}_{35}
-6  c^{5}_{35}
+2c^{6}_{35}
-5  c^{7}_{35}
-9  c^{10}_{35}
-c^{11}_{35}
$

\vskip 0.7ex
\hangindent=3em \hangafter=1
$T = ( 0,
\frac{4}{7},
\frac{2}{35},
\frac{27}{35},
\frac{1}{5},
\frac{6}{7} )
$,

\vskip 0.7ex
\hangindent=3em \hangafter=1
$S$ = ($ 1$,
$ -c_{7}^{2}$,
$ -1+c^{1}_{35}
-c^{5}_{35}
+c^{6}_{35}
-c^{7}_{35}
-c^{10}_{35}
$,
$ -c^{4}_{35}
-c^{10}_{35}
-c^{11}_{35}
$,
$ \frac{1-\sqrt{5}}{2}$,
$ -c^{1}_{7}
-c^{2}_{7}
$;\ \ 
$ c^{1}_{7}
+c^{2}_{7}
$,
$ \frac{1-\sqrt{5}}{2}$,
$ 1-c^{1}_{35}
+c^{5}_{35}
-c^{6}_{35}
+c^{7}_{35}
+c^{10}_{35}
$,
$ -c^{4}_{35}
-c^{10}_{35}
-c^{11}_{35}
$,
$ 1$;\ \ 
$ -c_{7}^{2}$,
$ -1$,
$ c^{1}_{7}
+c^{2}_{7}
$,
$ c^{4}_{35}
+c^{10}_{35}
+c^{11}_{35}
$;\ \ 
$ -c^{1}_{7}
-c^{2}_{7}
$,
$ c_{7}^{2}$,
$ \frac{1-\sqrt{5}}{2}$;\ \ 
$ -1$,
$ -1+c^{1}_{35}
-c^{5}_{35}
+c^{6}_{35}
-c^{7}_{35}
-c^{10}_{35}
$;\ \ 
$ c_{7}^{2}$)

Factors = $2_{\frac{2}{5},1.381}^{5,120}\boxtimes 3_{\frac{52}{7},1.841}^{7,604} $

Not pseudo-unitary. 

\vskip 1ex 
\color{grey}

\noindent(12,24). $6_{\frac{246}{35},2.544}^{35,724}$ \irep{47}:\ \ 
$d_i$ = ($1.0$,
$0.445$,
$0.495$,
$-0.275$,
$-0.618$,
$-0.801$) 

\vskip 0.7ex
\hangindent=3em \hangafter=1
$D^2= 2.544 = 
4+2c^{1}_{35}
-c^{4}_{35}
-6  c^{5}_{35}
+2c^{6}_{35}
-5  c^{7}_{35}
-9  c^{10}_{35}
-c^{11}_{35}
$

\vskip 0.7ex
\hangindent=3em \hangafter=1
$T = ( 0,
\frac{4}{7},
\frac{23}{35},
\frac{13}{35},
\frac{4}{5},
\frac{6}{7} )
$,

\vskip 0.7ex
\hangindent=3em \hangafter=1
$S$ = ($ 1$,
$ -c_{7}^{2}$,
$ -1+c^{1}_{35}
-c^{5}_{35}
+c^{6}_{35}
-c^{7}_{35}
-c^{10}_{35}
$,
$ -c^{4}_{35}
-c^{10}_{35}
-c^{11}_{35}
$,
$ \frac{1-\sqrt{5}}{2}$,
$ -c^{1}_{7}
-c^{2}_{7}
$;\ \ 
$ c^{1}_{7}
+c^{2}_{7}
$,
$ \frac{1-\sqrt{5}}{2}$,
$ 1-c^{1}_{35}
+c^{5}_{35}
-c^{6}_{35}
+c^{7}_{35}
+c^{10}_{35}
$,
$ -c^{4}_{35}
-c^{10}_{35}
-c^{11}_{35}
$,
$ 1$;\ \ 
$ -c_{7}^{2}$,
$ -1$,
$ c^{1}_{7}
+c^{2}_{7}
$,
$ c^{4}_{35}
+c^{10}_{35}
+c^{11}_{35}
$;\ \ 
$ -c^{1}_{7}
-c^{2}_{7}
$,
$ c_{7}^{2}$,
$ \frac{1-\sqrt{5}}{2}$;\ \ 
$ -1$,
$ -1+c^{1}_{35}
-c^{5}_{35}
+c^{6}_{35}
-c^{7}_{35}
-c^{10}_{35}
$;\ \ 
$ c_{7}^{2}$)

Factors = $2_{\frac{38}{5},1.381}^{5,491}\boxtimes 3_{\frac{52}{7},1.841}^{7,604} $

Not pseudo-unitary. 

\vskip 1ex 
\black

\noindent(13,1). $6_{\frac{46}{13},56.74}^{13,131}$ \irep{35}:\ \ 
$d_i$ = ($1.0$,
$1.941$,
$2.770$,
$3.438$,
$3.907$,
$4.148$) 

\vskip 0.7ex
\hangindent=3em \hangafter=1
$D^2= 56.746 = 
21+15c^{1}_{13}
+10c^{2}_{13}
+6c^{3}_{13}
+3c^{4}_{13}
+c^{5}_{13}
$

\vskip 0.7ex
\hangindent=3em \hangafter=1
$T = ( 0,
\frac{4}{13},
\frac{2}{13},
\frac{7}{13},
\frac{6}{13},
\frac{12}{13} )
$,

\vskip 0.7ex
\hangindent=3em \hangafter=1
$S$ = ($ 1$,
$ -c_{13}^{6}$,
$ \xi_{13}^{3}$,
$ \xi_{13}^{9}$,
$ \xi_{13}^{5}$,
$ \xi_{13}^{7}$;\ \ 
$ -\xi_{13}^{9}$,
$ \xi_{13}^{7}$,
$ -\xi_{13}^{5}$,
$ \xi_{13}^{3}$,
$ -1$;\ \ 
$ \xi_{13}^{9}$,
$ 1$,
$ c_{13}^{6}$,
$ -\xi_{13}^{5}$;\ \ 
$ \xi_{13}^{3}$,
$ -\xi_{13}^{7}$,
$ -c_{13}^{6}$;\ \ 
$ -1$,
$ \xi_{13}^{9}$;\ \ 
$ -\xi_{13}^{3}$)

Prime. 

\vskip 1ex 
\color{grey}

\noindent(13,2). $6_{\frac{58}{13},56.74}^{13,502}$ \irep{35}:\ \ 
$d_i$ = ($1.0$,
$1.941$,
$2.770$,
$3.438$,
$3.907$,
$4.148$) 

\vskip 0.7ex
\hangindent=3em \hangafter=1
$D^2= 56.746 = 
21+15c^{1}_{13}
+10c^{2}_{13}
+6c^{3}_{13}
+3c^{4}_{13}
+c^{5}_{13}
$

\vskip 0.7ex
\hangindent=3em \hangafter=1
$T = ( 0,
\frac{9}{13},
\frac{11}{13},
\frac{6}{13},
\frac{7}{13},
\frac{1}{13} )
$,

\vskip 0.7ex
\hangindent=3em \hangafter=1
$S$ = ($ 1$,
$ -c_{13}^{6}$,
$ \xi_{13}^{3}$,
$ \xi_{13}^{9}$,
$ \xi_{13}^{5}$,
$ \xi_{13}^{7}$;\ \ 
$ -\xi_{13}^{9}$,
$ \xi_{13}^{7}$,
$ -\xi_{13}^{5}$,
$ \xi_{13}^{3}$,
$ -1$;\ \ 
$ \xi_{13}^{9}$,
$ 1$,
$ c_{13}^{6}$,
$ -\xi_{13}^{5}$;\ \ 
$ \xi_{13}^{3}$,
$ -\xi_{13}^{7}$,
$ -c_{13}^{6}$;\ \ 
$ -1$,
$ \xi_{13}^{9}$;\ \ 
$ -\xi_{13}^{3}$)

Prime. 

\vskip 1ex 
\color{grey}

\noindent(13,3). $6_{\frac{90}{13},15.04}^{13,102}$ \irep{35}:\ \ 
$d_i$ = ($1.0$,
$1.426$,
$2.136$,
$-0.514$,
$-1.770$,
$-2.11$) 

\vskip 0.7ex
\hangindent=3em \hangafter=1
$D^2= 15.48 = 
15-6  c^{1}_{13}
+9c^{2}_{13}
-5  c^{3}_{13}
+4c^{4}_{13}
-3  c^{5}_{13}
$

\vskip 0.7ex
\hangindent=3em \hangafter=1
$T = ( 0,
\frac{1}{13},
\frac{9}{13},
\frac{2}{13},
\frac{5}{13},
\frac{12}{13} )
$,

\vskip 0.7ex
\hangindent=3em \hangafter=1
$S$ = ($ 1$,
$ 1+c^{2}_{13}
+c^{4}_{13}
$,
$ 1+c^{2}_{13}
$,
$ -c^{1}_{13}
-c^{3}_{13}
-c^{5}_{13}
$,
$ -c_{13}^{1}$,
$ -c^{1}_{13}
-c^{3}_{13}
$;\ \ 
$ -1$,
$ c_{13}^{1}$,
$ -c^{1}_{13}
-c^{3}_{13}
$,
$ 1+c^{2}_{13}
$,
$ c^{1}_{13}
+c^{3}_{13}
+c^{5}_{13}
$;\ \ 
$ -c^{1}_{13}
-c^{3}_{13}
$,
$ -1-c^{2}_{13}
-c^{4}_{13}
$,
$ -c^{1}_{13}
-c^{3}_{13}
-c^{5}_{13}
$,
$ 1$;\ \ 
$ -1-c^{2}_{13}
$,
$ -1$,
$ -c_{13}^{1}$;\ \ 
$ c^{1}_{13}
+c^{3}_{13}
$,
$ -1-c^{2}_{13}
-c^{4}_{13}
$;\ \ 
$ 1+c^{2}_{13}
$)

Prime. 

Not pseudo-unitary. 

\vskip 1ex 
\color{grey}

\noindent(13,4). $6_{\frac{14}{13},15.04}^{13,300}$ \irep{35}:\ \ 
$d_i$ = ($1.0$,
$1.426$,
$2.136$,
$-0.514$,
$-1.770$,
$-2.11$) 

\vskip 0.7ex
\hangindent=3em \hangafter=1
$D^2= 15.48 = 
15-6  c^{1}_{13}
+9c^{2}_{13}
-5  c^{3}_{13}
+4c^{4}_{13}
-3  c^{5}_{13}
$

\vskip 0.7ex
\hangindent=3em \hangafter=1
$T = ( 0,
\frac{12}{13},
\frac{4}{13},
\frac{11}{13},
\frac{8}{13},
\frac{1}{13} )
$,

\vskip 0.7ex
\hangindent=3em \hangafter=1
$S$ = ($ 1$,
$ 1+c^{2}_{13}
+c^{4}_{13}
$,
$ 1+c^{2}_{13}
$,
$ -c^{1}_{13}
-c^{3}_{13}
-c^{5}_{13}
$,
$ -c_{13}^{1}$,
$ -c^{1}_{13}
-c^{3}_{13}
$;\ \ 
$ -1$,
$ c_{13}^{1}$,
$ -c^{1}_{13}
-c^{3}_{13}
$,
$ 1+c^{2}_{13}
$,
$ c^{1}_{13}
+c^{3}_{13}
+c^{5}_{13}
$;\ \ 
$ -c^{1}_{13}
-c^{3}_{13}
$,
$ -1-c^{2}_{13}
-c^{4}_{13}
$,
$ -c^{1}_{13}
-c^{3}_{13}
-c^{5}_{13}
$,
$ 1$;\ \ 
$ -1-c^{2}_{13}
$,
$ -1$,
$ -c_{13}^{1}$;\ \ 
$ c^{1}_{13}
+c^{3}_{13}
$,
$ -1-c^{2}_{13}
-c^{4}_{13}
$;\ \ 
$ 1+c^{2}_{13}
$)

Prime. 

Not pseudo-unitary. 

\vskip 1ex 
\color{grey}

\noindent(13,5). $6_{\frac{18}{13},7.390}^{13,415}$ \irep{35}:\ \ 
$d_i$ = ($1.0$,
$0.360$,
$1.241$,
$1.497$,
$-0.700$,
$-1.410$) 

\vskip 0.7ex
\hangindent=3em \hangafter=1
$D^2= 7.390 = 
11-7  c^{1}_{13}
-9  c^{2}_{13}
+5c^{3}_{13}
-4  c^{4}_{13}
-10  c^{5}_{13}
$

\vskip 0.7ex
\hangindent=3em \hangafter=1
$T = ( 0,
\frac{5}{13},
\frac{7}{13},
\frac{1}{13},
\frac{8}{13},
\frac{3}{13} )
$,

\vskip 0.7ex
\hangindent=3em \hangafter=1
$S$ = ($ 1$,
$ -c^{2}_{13}
-c^{5}_{13}
$,
$ 1+c^{3}_{13}
$,
$ -c_{13}^{5}$,
$ -c^{1}_{13}
-c^{2}_{13}
-c^{4}_{13}
-c^{5}_{13}
$,
$ -c^{1}_{13}
-c^{2}_{13}
-c^{5}_{13}
$;\ \ 
$ 1+c^{3}_{13}
$,
$ 1$,
$ c^{1}_{13}
+c^{2}_{13}
+c^{4}_{13}
+c^{5}_{13}
$,
$ c^{1}_{13}
+c^{2}_{13}
+c^{5}_{13}
$,
$ -c_{13}^{5}$;\ \ 
$ -c^{2}_{13}
-c^{5}_{13}
$,
$ -c^{1}_{13}
-c^{2}_{13}
-c^{5}_{13}
$,
$ c_{13}^{5}$,
$ c^{1}_{13}
+c^{2}_{13}
+c^{4}_{13}
+c^{5}_{13}
$;\ \ 
$ c^{2}_{13}
+c^{5}_{13}
$,
$ 1+c^{3}_{13}
$,
$ -1$;\ \ 
$ -1$,
$ -c^{2}_{13}
-c^{5}_{13}
$;\ \ 
$ -1-c^{3}_{13}
$)

Prime. 

Not pseudo-unitary. 

\vskip 1ex 
\color{grey}

\noindent(13,6). $6_{\frac{86}{13},7.390}^{13,241}$ \irep{35}:\ \ 
$d_i$ = ($1.0$,
$0.360$,
$1.241$,
$1.497$,
$-0.700$,
$-1.410$) 

\vskip 0.7ex
\hangindent=3em \hangafter=1
$D^2= 7.390 = 
11-7  c^{1}_{13}
-9  c^{2}_{13}
+5c^{3}_{13}
-4  c^{4}_{13}
-10  c^{5}_{13}
$

\vskip 0.7ex
\hangindent=3em \hangafter=1
$T = ( 0,
\frac{8}{13},
\frac{6}{13},
\frac{12}{13},
\frac{5}{13},
\frac{10}{13} )
$,

\vskip 0.7ex
\hangindent=3em \hangafter=1
$S$ = ($ 1$,
$ -c^{2}_{13}
-c^{5}_{13}
$,
$ 1+c^{3}_{13}
$,
$ -c_{13}^{5}$,
$ -c^{1}_{13}
-c^{2}_{13}
-c^{4}_{13}
-c^{5}_{13}
$,
$ -c^{1}_{13}
-c^{2}_{13}
-c^{5}_{13}
$;\ \ 
$ 1+c^{3}_{13}
$,
$ 1$,
$ c^{1}_{13}
+c^{2}_{13}
+c^{4}_{13}
+c^{5}_{13}
$,
$ c^{1}_{13}
+c^{2}_{13}
+c^{5}_{13}
$,
$ -c_{13}^{5}$;\ \ 
$ -c^{2}_{13}
-c^{5}_{13}
$,
$ -c^{1}_{13}
-c^{2}_{13}
-c^{5}_{13}
$,
$ c_{13}^{5}$,
$ c^{1}_{13}
+c^{2}_{13}
+c^{4}_{13}
+c^{5}_{13}
$;\ \ 
$ c^{2}_{13}
+c^{5}_{13}
$,
$ 1+c^{3}_{13}
$,
$ -1$;\ \ 
$ -1$,
$ -c^{2}_{13}
-c^{5}_{13}
$;\ \ 
$ -1-c^{3}_{13}
$)

Prime. 

Not pseudo-unitary. 

\vskip 1ex 
\color{grey}

\noindent(13,7). $6_{\frac{102}{13},4.798}^{13,107}$ \irep{35}:\ \ 
$d_i$ = ($1.0$,
$0.290$,
$0.564$,
$0.805$,
$-1.136$,
$-1.206$) 

\vskip 0.7ex
\hangindent=3em \hangafter=1
$D^2= 4.798 = 
20+5c^{1}_{13}
-c^{2}_{13}
+2c^{3}_{13}
+14c^{4}_{13}
+9c^{5}_{13}
$

\vskip 0.7ex
\hangindent=3em \hangafter=1
$T = ( 0,
\frac{5}{13},
\frac{4}{13},
\frac{11}{13},
\frac{10}{13},
\frac{2}{13} )
$,

\vskip 0.7ex
\hangindent=3em \hangafter=1
$S$ = ($ 1$,
$ 1+c^{4}_{13}
$,
$ 1+c^{1}_{13}
+c^{4}_{13}
+c^{5}_{13}
$,
$ 1+c^{1}_{13}
+c^{3}_{13}
+c^{4}_{13}
+c^{5}_{13}
$,
$ -c_{13}^{2}$,
$ 1+c^{4}_{13}
+c^{5}_{13}
$;\ \ 
$ 1+c^{1}_{13}
+c^{3}_{13}
+c^{4}_{13}
+c^{5}_{13}
$,
$ -1-c^{4}_{13}
-c^{5}_{13}
$,
$ 1$,
$ 1+c^{1}_{13}
+c^{4}_{13}
+c^{5}_{13}
$,
$ c_{13}^{2}$;\ \ 
$ -1-c^{4}_{13}
$,
$ -c_{13}^{2}$,
$ -1$,
$ 1+c^{1}_{13}
+c^{3}_{13}
+c^{4}_{13}
+c^{5}_{13}
$;\ \ 
$ 1+c^{4}_{13}
$,
$ -1-c^{4}_{13}
-c^{5}_{13}
$,
$ -1-c^{1}_{13}
-c^{4}_{13}
-c^{5}_{13}
$;\ \ 
$ -1-c^{1}_{13}
-c^{3}_{13}
-c^{4}_{13}
-c^{5}_{13}
$,
$ 1+c^{4}_{13}
$;\ \ 
$ -1$)

Prime. 

Not pseudo-unitary. 

\vskip 1ex 
\color{grey}

\noindent(13,8). $6_{\frac{2}{13},4.798}^{13,604}$ \irep{35}:\ \ 
$d_i$ = ($1.0$,
$0.290$,
$0.564$,
$0.805$,
$-1.136$,
$-1.206$) 

\vskip 0.7ex
\hangindent=3em \hangafter=1
$D^2= 4.798 = 
20+5c^{1}_{13}
-c^{2}_{13}
+2c^{3}_{13}
+14c^{4}_{13}
+9c^{5}_{13}
$

\vskip 0.7ex
\hangindent=3em \hangafter=1
$T = ( 0,
\frac{8}{13},
\frac{9}{13},
\frac{2}{13},
\frac{3}{13},
\frac{11}{13} )
$,

\vskip 0.7ex
\hangindent=3em \hangafter=1
$S$ = ($ 1$,
$ 1+c^{4}_{13}
$,
$ 1+c^{1}_{13}
+c^{4}_{13}
+c^{5}_{13}
$,
$ 1+c^{1}_{13}
+c^{3}_{13}
+c^{4}_{13}
+c^{5}_{13}
$,
$ -c_{13}^{2}$,
$ 1+c^{4}_{13}
+c^{5}_{13}
$;\ \ 
$ 1+c^{1}_{13}
+c^{3}_{13}
+c^{4}_{13}
+c^{5}_{13}
$,
$ -1-c^{4}_{13}
-c^{5}_{13}
$,
$ 1$,
$ 1+c^{1}_{13}
+c^{4}_{13}
+c^{5}_{13}
$,
$ c_{13}^{2}$;\ \ 
$ -1-c^{4}_{13}
$,
$ -c_{13}^{2}$,
$ -1$,
$ 1+c^{1}_{13}
+c^{3}_{13}
+c^{4}_{13}
+c^{5}_{13}
$;\ \ 
$ 1+c^{4}_{13}
$,
$ -1-c^{4}_{13}
-c^{5}_{13}
$,
$ -1-c^{1}_{13}
-c^{4}_{13}
-c^{5}_{13}
$;\ \ 
$ -1-c^{1}_{13}
-c^{3}_{13}
-c^{4}_{13}
-c^{5}_{13}
$,
$ 1+c^{4}_{13}
$;\ \ 
$ -1$)

Prime. 

Not pseudo-unitary. 

\vskip 1ex 
\color{grey}

\noindent(13,9). $6_{\frac{30}{13},3.717}^{13,162}$ \irep{35}:\ \ 
$d_i$ = ($1.0$,
$0.709$,
$0.880$,
$-0.255$,
$-0.497$,
$-1.61$) 

\vskip 0.7ex
\hangindent=3em \hangafter=1
$D^2= 3.717 = 
18-2  c^{1}_{13}
+3c^{2}_{13}
+7c^{3}_{13}
-3  c^{4}_{13}
+12c^{5}_{13}
$

\vskip 0.7ex
\hangindent=3em \hangafter=1
$T = ( 0,
\frac{6}{13},
\frac{5}{13},
\frac{9}{13},
\frac{3}{13},
\frac{4}{13} )
$,

\vskip 0.7ex
\hangindent=3em \hangafter=1
$S$ = ($ 1$,
$ -c_{13}^{4}$,
$ 1+c^{2}_{13}
+c^{3}_{13}
+c^{5}_{13}
$,
$ 1+c^{3}_{13}
+c^{5}_{13}
$,
$ 1+c^{5}_{13}
$,
$ -c^{1}_{13}
-c^{4}_{13}
$;\ \ 
$ c^{1}_{13}
+c^{4}_{13}
$,
$ -1$,
$ 1+c^{5}_{13}
$,
$ 1+c^{2}_{13}
+c^{3}_{13}
+c^{5}_{13}
$,
$ -1-c^{3}_{13}
-c^{5}_{13}
$;\ \ 
$ -1-c^{5}_{13}
$,
$ -c^{1}_{13}
-c^{4}_{13}
$,
$ -1-c^{3}_{13}
-c^{5}_{13}
$,
$ -c_{13}^{4}$;\ \ 
$ -1$,
$ c_{13}^{4}$,
$ -1-c^{2}_{13}
-c^{3}_{13}
-c^{5}_{13}
$;\ \ 
$ -c^{1}_{13}
-c^{4}_{13}
$,
$ 1$;\ \ 
$ 1+c^{5}_{13}
$)

Prime. 

Not pseudo-unitary. 

\vskip 1ex 
\color{grey}

\noindent(13,10). $6_{\frac{74}{13},3.717}^{13,481}$ \irep{35}:\ \ 
$d_i$ = ($1.0$,
$0.709$,
$0.880$,
$-0.255$,
$-0.497$,
$-1.61$) 

\vskip 0.7ex
\hangindent=3em \hangafter=1
$D^2= 3.717 = 
18-2  c^{1}_{13}
+3c^{2}_{13}
+7c^{3}_{13}
-3  c^{4}_{13}
+12c^{5}_{13}
$

\vskip 0.7ex
\hangindent=3em \hangafter=1
$T = ( 0,
\frac{7}{13},
\frac{8}{13},
\frac{4}{13},
\frac{10}{13},
\frac{9}{13} )
$,

\vskip 0.7ex
\hangindent=3em \hangafter=1
$S$ = ($ 1$,
$ -c_{13}^{4}$,
$ 1+c^{2}_{13}
+c^{3}_{13}
+c^{5}_{13}
$,
$ 1+c^{3}_{13}
+c^{5}_{13}
$,
$ 1+c^{5}_{13}
$,
$ -c^{1}_{13}
-c^{4}_{13}
$;\ \ 
$ c^{1}_{13}
+c^{4}_{13}
$,
$ -1$,
$ 1+c^{5}_{13}
$,
$ 1+c^{2}_{13}
+c^{3}_{13}
+c^{5}_{13}
$,
$ -1-c^{3}_{13}
-c^{5}_{13}
$;\ \ 
$ -1-c^{5}_{13}
$,
$ -c^{1}_{13}
-c^{4}_{13}
$,
$ -1-c^{3}_{13}
-c^{5}_{13}
$,
$ -c_{13}^{4}$;\ \ 
$ -1$,
$ c_{13}^{4}$,
$ -1-c^{2}_{13}
-c^{3}_{13}
-c^{5}_{13}
$;\ \ 
$ -c^{1}_{13}
-c^{4}_{13}
$,
$ 1$;\ \ 
$ 1+c^{5}_{13}
$)

Prime. 

Not pseudo-unitary. 

\vskip 1ex 
\color{grey}

\noindent(13,11). $6_{\frac{94}{13},3.297}^{13,764}$ \irep{35}:\ \ 
$d_i$ = ($1.0$,
$0.468$,
$0.829$,
$-0.241$,
$-0.667$,
$-0.941$) 

\vskip 0.7ex
\hangindent=3em \hangafter=1
$D^2= 3.297 = 
6-5  c^{1}_{13}
-12  c^{2}_{13}
-15  c^{3}_{13}
-14  c^{4}_{13}
-9  c^{5}_{13}
$

\vskip 0.7ex
\hangindent=3em \hangafter=1
$T = ( 0,
\frac{3}{13},
\frac{10}{13},
\frac{11}{13},
\frac{7}{13},
\frac{12}{13} )
$,

\vskip 0.7ex
\hangindent=3em \hangafter=1
$S$ = ($ 1$,
$ -c^{3}_{13}
-c^{4}_{13}
$,
$ -c^{2}_{13}
-c^{3}_{13}
-c^{4}_{13}
-c^{5}_{13}
$,
$ -c_{13}^{3}$,
$ -c^{2}_{13}
-c^{3}_{13}
-c^{4}_{13}
$,
$ -c^{1}_{13}
-c^{2}_{13}
-c^{3}_{13}
-c^{4}_{13}
-c^{5}_{13}
$;\ \ 
$ -c^{1}_{13}
-c^{2}_{13}
-c^{3}_{13}
-c^{4}_{13}
-c^{5}_{13}
$,
$ c^{2}_{13}
+c^{3}_{13}
+c^{4}_{13}
$,
$ c^{2}_{13}
+c^{3}_{13}
+c^{4}_{13}
+c^{5}_{13}
$,
$ -c_{13}^{3}$,
$ 1$;\ \ 
$ -1$,
$ -c^{1}_{13}
-c^{2}_{13}
-c^{3}_{13}
-c^{4}_{13}
-c^{5}_{13}
$,
$ -c^{3}_{13}
-c^{4}_{13}
$,
$ c_{13}^{3}$;\ \ 
$ c^{3}_{13}
+c^{4}_{13}
$,
$ -1$,
$ -c^{2}_{13}
-c^{3}_{13}
-c^{4}_{13}
$;\ \ 
$ c^{1}_{13}
+c^{2}_{13}
+c^{3}_{13}
+c^{4}_{13}
+c^{5}_{13}
$,
$ c^{2}_{13}
+c^{3}_{13}
+c^{4}_{13}
+c^{5}_{13}
$;\ \ 
$ -c^{3}_{13}
-c^{4}_{13}
$)

Prime. 

Not pseudo-unitary. 

\vskip 1ex 
\color{grey}

\noindent(13,12). $6_{\frac{10}{13},3.297}^{13,560}$ \irep{35}:\ \ 
$d_i$ = ($1.0$,
$0.468$,
$0.829$,
$-0.241$,
$-0.667$,
$-0.941$) 

\vskip 0.7ex
\hangindent=3em \hangafter=1
$D^2= 3.297 = 
6-5  c^{1}_{13}
-12  c^{2}_{13}
-15  c^{3}_{13}
-14  c^{4}_{13}
-9  c^{5}_{13}
$

\vskip 0.7ex
\hangindent=3em \hangafter=1
$T = ( 0,
\frac{10}{13},
\frac{3}{13},
\frac{2}{13},
\frac{6}{13},
\frac{1}{13} )
$,

\vskip 0.7ex
\hangindent=3em \hangafter=1
$S$ = ($ 1$,
$ -c^{3}_{13}
-c^{4}_{13}
$,
$ -c^{2}_{13}
-c^{3}_{13}
-c^{4}_{13}
-c^{5}_{13}
$,
$ -c_{13}^{3}$,
$ -c^{2}_{13}
-c^{3}_{13}
-c^{4}_{13}
$,
$ -c^{1}_{13}
-c^{2}_{13}
-c^{3}_{13}
-c^{4}_{13}
-c^{5}_{13}
$;\ \ 
$ -c^{1}_{13}
-c^{2}_{13}
-c^{3}_{13}
-c^{4}_{13}
-c^{5}_{13}
$,
$ c^{2}_{13}
+c^{3}_{13}
+c^{4}_{13}
$,
$ c^{2}_{13}
+c^{3}_{13}
+c^{4}_{13}
+c^{5}_{13}
$,
$ -c_{13}^{3}$,
$ 1$;\ \ 
$ -1$,
$ -c^{1}_{13}
-c^{2}_{13}
-c^{3}_{13}
-c^{4}_{13}
-c^{5}_{13}
$,
$ -c^{3}_{13}
-c^{4}_{13}
$,
$ c_{13}^{3}$;\ \ 
$ c^{3}_{13}
+c^{4}_{13}
$,
$ -1$,
$ -c^{2}_{13}
-c^{3}_{13}
-c^{4}_{13}
$;\ \ 
$ c^{1}_{13}
+c^{2}_{13}
+c^{3}_{13}
+c^{4}_{13}
+c^{5}_{13}
$,
$ c^{2}_{13}
+c^{3}_{13}
+c^{4}_{13}
+c^{5}_{13}
$;\ \ 
$ -c^{3}_{13}
-c^{4}_{13}
$)

Prime. 

Not pseudo-unitary. 

\vskip 1ex 
\black

\noindent(14,1). $6_{\frac{8}{3},74.61}^{9,186}$ \irep{26}:\ \ 
$d_i$ = ($1.0$,
$2.879$,
$2.879$,
$2.879$,
$4.411$,
$5.411$) 

\vskip 0.7ex
\hangindent=3em \hangafter=1
$D^2= 74.617 = 
27+27c^{1}_{9}
+18c^{2}_{9}
$

\vskip 0.7ex
\hangindent=3em \hangafter=1
$T = ( 0,
\frac{1}{9},
\frac{1}{9},
\frac{1}{9},
\frac{1}{3},
\frac{2}{3} )
$,

\vskip 0.7ex
\hangindent=3em \hangafter=1
$S$ = ($ 1$,
$ \xi_{9}^{5}$,
$ \xi_{9}^{5}$,
$ \xi_{9}^{5}$,
$ 1+2c^{1}_{9}
+c^{2}_{9}
$,
$ 2+2c^{1}_{9}
+c^{2}_{9}
$;\ \ 
$ 2\xi_{9}^{5}$,
$ -\xi_{9}^{5}$,
$ -\xi_{9}^{5}$,
$ \xi_{9}^{5}$,
$ -\xi_{9}^{5}$;\ \ 
$ 2\xi_{9}^{5}$,
$ -\xi_{9}^{5}$,
$ \xi_{9}^{5}$,
$ -\xi_{9}^{5}$;\ \ 
$ 2\xi_{9}^{5}$,
$ \xi_{9}^{5}$,
$ -\xi_{9}^{5}$;\ \ 
$ -2-2  c^{1}_{9}
-c^{2}_{9}
$,
$ -1$;\ \ 
$ 1+2c^{1}_{9}
+c^{2}_{9}
$)

Prime. 

\vskip 1ex 
\color{grey}

\noindent(14,2). $6_{\frac{16}{3},74.61}^{9,452}$ \irep{26}:\ \ 
$d_i$ = ($1.0$,
$2.879$,
$2.879$,
$2.879$,
$4.411$,
$5.411$) 

\vskip 0.7ex
\hangindent=3em \hangafter=1
$D^2= 74.617 = 
27+27c^{1}_{9}
+18c^{2}_{9}
$

\vskip 0.7ex
\hangindent=3em \hangafter=1
$T = ( 0,
\frac{8}{9},
\frac{8}{9},
\frac{8}{9},
\frac{2}{3},
\frac{1}{3} )
$,

\vskip 0.7ex
\hangindent=3em \hangafter=1
$S$ = ($ 1$,
$ \xi_{9}^{5}$,
$ \xi_{9}^{5}$,
$ \xi_{9}^{5}$,
$ 1+2c^{1}_{9}
+c^{2}_{9}
$,
$ 2+2c^{1}_{9}
+c^{2}_{9}
$;\ \ 
$ 2\xi_{9}^{5}$,
$ -\xi_{9}^{5}$,
$ -\xi_{9}^{5}$,
$ \xi_{9}^{5}$,
$ -\xi_{9}^{5}$;\ \ 
$ 2\xi_{9}^{5}$,
$ -\xi_{9}^{5}$,
$ \xi_{9}^{5}$,
$ -\xi_{9}^{5}$;\ \ 
$ 2\xi_{9}^{5}$,
$ \xi_{9}^{5}$,
$ -\xi_{9}^{5}$;\ \ 
$ -2-2  c^{1}_{9}
-c^{2}_{9}
$,
$ -1$;\ \ 
$ 1+2c^{1}_{9}
+c^{2}_{9}
$)

Prime. 

\vskip 1ex 
\color{grey}

\noindent(14,3). $6_{\frac{8}{3},3.834}^{9,226}$ \irep{26}:\ \ 
$d_i$ = ($1.0$,
$0.652$,
$0.652$,
$0.652$,
$-0.226$,
$-1.226$) 

\vskip 0.7ex
\hangindent=3em \hangafter=1
$D^2= 3.834 = 
27-9  c^{1}_{9}
-27  c^{2}_{9}
$

\vskip 0.7ex
\hangindent=3em \hangafter=1
$T = ( 0,
\frac{4}{9},
\frac{4}{9},
\frac{4}{9},
\frac{2}{3},
\frac{1}{3} )
$,

\vskip 0.7ex
\hangindent=3em \hangafter=1
$S$ = ($ 1$,
$ 1-c^{2}_{9}
$,
$ 1-c^{2}_{9}
$,
$ 1-c^{2}_{9}
$,
$ 2-c^{1}_{9}
-2  c^{2}_{9}
$,
$ 1-c^{1}_{9}
-2  c^{2}_{9}
$;\ \ 
$ 2-2  c^{2}_{9}
$,
$ -1+c^{2}_{9}
$,
$ -1+c^{2}_{9}
$,
$ -1+c^{2}_{9}
$,
$ 1-c^{2}_{9}
$;\ \ 
$ 2-2  c^{2}_{9}
$,
$ -1+c^{2}_{9}
$,
$ -1+c^{2}_{9}
$,
$ 1-c^{2}_{9}
$;\ \ 
$ 2-2  c^{2}_{9}
$,
$ -1+c^{2}_{9}
$,
$ 1-c^{2}_{9}
$;\ \ 
$ 1-c^{1}_{9}
-2  c^{2}_{9}
$,
$ -1$;\ \ 
$ -2+c^{1}_{9}
+2c^{2}_{9}
$)

Prime. 

Not pseudo-unitary. 

\vskip 1ex 
\color{grey}

\noindent(14,4). $6_{\frac{16}{3},3.834}^{9,342}$ \irep{26}:\ \ 
$d_i$ = ($1.0$,
$0.652$,
$0.652$,
$0.652$,
$-0.226$,
$-1.226$) 

\vskip 0.7ex
\hangindent=3em \hangafter=1
$D^2= 3.834 = 
27-9  c^{1}_{9}
-27  c^{2}_{9}
$

\vskip 0.7ex
\hangindent=3em \hangafter=1
$T = ( 0,
\frac{5}{9},
\frac{5}{9},
\frac{5}{9},
\frac{1}{3},
\frac{2}{3} )
$,

\vskip 0.7ex
\hangindent=3em \hangafter=1
$S$ = ($ 1$,
$ 1-c^{2}_{9}
$,
$ 1-c^{2}_{9}
$,
$ 1-c^{2}_{9}
$,
$ 2-c^{1}_{9}
-2  c^{2}_{9}
$,
$ 1-c^{1}_{9}
-2  c^{2}_{9}
$;\ \ 
$ 2-2  c^{2}_{9}
$,
$ -1+c^{2}_{9}
$,
$ -1+c^{2}_{9}
$,
$ -1+c^{2}_{9}
$,
$ 1-c^{2}_{9}
$;\ \ 
$ 2-2  c^{2}_{9}
$,
$ -1+c^{2}_{9}
$,
$ -1+c^{2}_{9}
$,
$ 1-c^{2}_{9}
$;\ \ 
$ 2-2  c^{2}_{9}
$,
$ -1+c^{2}_{9}
$,
$ 1-c^{2}_{9}
$;\ \ 
$ 1-c^{1}_{9}
-2  c^{2}_{9}
$,
$ -1$;\ \ 
$ -2+c^{1}_{9}
+2c^{2}_{9}
$)

Prime. 

Not pseudo-unitary. 

\vskip 1ex 
\color{grey}

\noindent(14,5). $6_{\frac{4}{3},2.548}^{9,789}$ \irep{26}:\ \ 
$d_i$ = ($1.0$,
$0.815$,
$-0.184$,
$-0.532$,
$-0.532$,
$-0.532$) 

\vskip 0.7ex
\hangindent=3em \hangafter=1
$D^2= 2.548 = 
27-18  c^{1}_{9}
+9c^{2}_{9}
$

\vskip 0.7ex
\hangindent=3em \hangafter=1
$T = ( 0,
\frac{1}{3},
\frac{2}{3},
\frac{2}{9},
\frac{2}{9},
\frac{2}{9} )
$,

\vskip 0.7ex
\hangindent=3em \hangafter=1
$S$ = ($ 1$,
$ 2-c^{1}_{9}
+c^{2}_{9}
$,
$ 1-c^{1}_{9}
+c^{2}_{9}
$,
$ 1-c^{1}_{9}
$,
$ 1-c^{1}_{9}
$,
$ 1-c^{1}_{9}
$;\ \ 
$ 1-c^{1}_{9}
+c^{2}_{9}
$,
$ -1$,
$ -1+c^{1}_{9}
$,
$ -1+c^{1}_{9}
$,
$ -1+c^{1}_{9}
$;\ \ 
$ -2+c^{1}_{9}
-c^{2}_{9}
$,
$ 1-c^{1}_{9}
$,
$ 1-c^{1}_{9}
$,
$ 1-c^{1}_{9}
$;\ \ 
$ 2-2  c^{1}_{9}
$,
$ -1+c^{1}_{9}
$,
$ -1+c^{1}_{9}
$;\ \ 
$ 2-2  c^{1}_{9}
$,
$ -1+c^{1}_{9}
$;\ \ 
$ 2-2  c^{1}_{9}
$)

Prime. 

Not pseudo-unitary. 

\vskip 1ex 
\color{grey}

\noindent(14,6). $6_{\frac{20}{3},2.548}^{9,632}$ \irep{26}:\ \ 
$d_i$ = ($1.0$,
$0.815$,
$-0.184$,
$-0.532$,
$-0.532$,
$-0.532$) 

\vskip 0.7ex
\hangindent=3em \hangafter=1
$D^2= 2.548 = 
27-18  c^{1}_{9}
+9c^{2}_{9}
$

\vskip 0.7ex
\hangindent=3em \hangafter=1
$T = ( 0,
\frac{2}{3},
\frac{1}{3},
\frac{7}{9},
\frac{7}{9},
\frac{7}{9} )
$,

\vskip 0.7ex
\hangindent=3em \hangafter=1
$S$ = ($ 1$,
$ 2-c^{1}_{9}
+c^{2}_{9}
$,
$ 1-c^{1}_{9}
+c^{2}_{9}
$,
$ 1-c^{1}_{9}
$,
$ 1-c^{1}_{9}
$,
$ 1-c^{1}_{9}
$;\ \ 
$ 1-c^{1}_{9}
+c^{2}_{9}
$,
$ -1$,
$ -1+c^{1}_{9}
$,
$ -1+c^{1}_{9}
$,
$ -1+c^{1}_{9}
$;\ \ 
$ -2+c^{1}_{9}
-c^{2}_{9}
$,
$ 1-c^{1}_{9}
$,
$ 1-c^{1}_{9}
$,
$ 1-c^{1}_{9}
$;\ \ 
$ 2-2  c^{1}_{9}
$,
$ -1+c^{1}_{9}
$,
$ -1+c^{1}_{9}
$;\ \ 
$ 2-2  c^{1}_{9}
$,
$ -1+c^{1}_{9}
$;\ \ 
$ 2-2  c^{1}_{9}
$)

Prime. 

Not pseudo-unitary. 

\vskip 1ex 
\black

\noindent(15,1). $6_{6,100.6}^{21,154}$ \irep{43}:\ \ 
$d_i$ = ($1.0$,
$3.791$,
$3.791$,
$3.791$,
$4.791$,
$5.791$) 

\vskip 0.7ex
\hangindent=3em \hangafter=1
$D^2= 100.617 = 
\frac{105+21\sqrt{21}}{2}$

\vskip 0.7ex
\hangindent=3em \hangafter=1
$T = ( 0,
\frac{1}{7},
\frac{2}{7},
\frac{4}{7},
0,
\frac{2}{3} )
$,

\vskip 0.7ex
\hangindent=3em \hangafter=1
$S$ = ($ 1$,
$ \frac{3+\sqrt{21}}{2}$,
$ \frac{3+\sqrt{21}}{2}$,
$ \frac{3+\sqrt{21}}{2}$,
$ \frac{5+\sqrt{21}}{2}$,
$ \frac{7+\sqrt{21}}{2}$;\ \ 
$ 2-c^{1}_{21}
-2  c^{2}_{21}
+3c^{3}_{21}
+2c^{4}_{21}
-2  c^{5}_{21}
$,
$ -c^{2}_{21}
-2  c^{3}_{21}
-c^{4}_{21}
+c^{5}_{21}
$,
$ -1+2c^{1}_{21}
+3c^{2}_{21}
-c^{3}_{21}
+2c^{5}_{21}
$,
$ -\frac{3+\sqrt{21}}{2}$,
$0$;\ \ 
$ -1+2c^{1}_{21}
+3c^{2}_{21}
-c^{3}_{21}
+2c^{5}_{21}
$,
$ 2-c^{1}_{21}
-2  c^{2}_{21}
+3c^{3}_{21}
+2c^{4}_{21}
-2  c^{5}_{21}
$,
$ -\frac{3+\sqrt{21}}{2}$,
$0$;\ \ 
$ -c^{2}_{21}
-2  c^{3}_{21}
-c^{4}_{21}
+c^{5}_{21}
$,
$ -\frac{3+\sqrt{21}}{2}$,
$0$;\ \ 
$ 1$,
$ \frac{7+\sqrt{21}}{2}$;\ \ 
$ -\frac{7+\sqrt{21}}{2}$)

Prime. 

\vskip 1ex 
\color{grey}

\noindent(15,2). $6_{2,100.6}^{21,320}$ \irep{43}:\ \ 
$d_i$ = ($1.0$,
$3.791$,
$3.791$,
$3.791$,
$4.791$,
$5.791$) 

\vskip 0.7ex
\hangindent=3em \hangafter=1
$D^2= 100.617 = 
\frac{105+21\sqrt{21}}{2}$

\vskip 0.7ex
\hangindent=3em \hangafter=1
$T = ( 0,
\frac{3}{7},
\frac{5}{7},
\frac{6}{7},
0,
\frac{1}{3} )
$,

\vskip 0.7ex
\hangindent=3em \hangafter=1
$S$ = ($ 1$,
$ \frac{3+\sqrt{21}}{2}$,
$ \frac{3+\sqrt{21}}{2}$,
$ \frac{3+\sqrt{21}}{2}$,
$ \frac{5+\sqrt{21}}{2}$,
$ \frac{7+\sqrt{21}}{2}$;\ \ 
$ -c^{2}_{21}
-2  c^{3}_{21}
-c^{4}_{21}
+c^{5}_{21}
$,
$ 2-c^{1}_{21}
-2  c^{2}_{21}
+3c^{3}_{21}
+2c^{4}_{21}
-2  c^{5}_{21}
$,
$ -1+2c^{1}_{21}
+3c^{2}_{21}
-c^{3}_{21}
+2c^{5}_{21}
$,
$ -\frac{3+\sqrt{21}}{2}$,
$0$;\ \ 
$ -1+2c^{1}_{21}
+3c^{2}_{21}
-c^{3}_{21}
+2c^{5}_{21}
$,
$ -c^{2}_{21}
-2  c^{3}_{21}
-c^{4}_{21}
+c^{5}_{21}
$,
$ -\frac{3+\sqrt{21}}{2}$,
$0$;\ \ 
$ 2-c^{1}_{21}
-2  c^{2}_{21}
+3c^{3}_{21}
+2c^{4}_{21}
-2  c^{5}_{21}
$,
$ -\frac{3+\sqrt{21}}{2}$,
$0$;\ \ 
$ 1$,
$ \frac{7+\sqrt{21}}{2}$;\ \ 
$ -\frac{7+\sqrt{21}}{2}$)

Prime. 

\vskip 1ex 
\color{grey}

\noindent(15,3). $6_{2,4.382}^{21,215}$ \irep{43}:\ \ 
$d_i$ = ($1.0$,
$0.208$,
$1.208$,
$-0.791$,
$-0.791$,
$-0.791$) 

\vskip 0.7ex
\hangindent=3em \hangafter=1
$D^2= 4.382 = 
\frac{105-21\sqrt{21}}{2}$

\vskip 0.7ex
\hangindent=3em \hangafter=1
$T = ( 0,
0,
\frac{1}{3},
\frac{1}{7},
\frac{2}{7},
\frac{4}{7} )
$,

\vskip 0.7ex
\hangindent=3em \hangafter=1
$S$ = ($ 1$,
$ \frac{5-\sqrt{21}}{2}$,
$ \frac{7-\sqrt{21}}{2}$,
$ \frac{3-\sqrt{21}}{2}$,
$ \frac{3-\sqrt{21}}{2}$,
$ \frac{3-\sqrt{21}}{2}$;\ \ 
$ 1$,
$ \frac{7-\sqrt{21}}{2}$,
$ -\frac{3-\sqrt{21}}{2}$,
$ -\frac{3-\sqrt{21}}{2}$,
$ -\frac{3-\sqrt{21}}{2}$;\ \ 
$ -\frac{7-\sqrt{21}}{2}$,
$0$,
$0$,
$0$;\ \ 
$ 1+c^{1}_{21}
-c^{2}_{21}
-2  c^{4}_{21}
-c^{5}_{21}
$,
$ c^{2}_{21}
-c^{3}_{21}
+c^{4}_{21}
-c^{5}_{21}
$,
$ 1-2  c^{1}_{21}
+c^{3}_{21}
+c^{5}_{21}
$;\ \ 
$ 1-2  c^{1}_{21}
+c^{3}_{21}
+c^{5}_{21}
$,
$ 1+c^{1}_{21}
-c^{2}_{21}
-2  c^{4}_{21}
-c^{5}_{21}
$;\ \ 
$ c^{2}_{21}
-c^{3}_{21}
+c^{4}_{21}
-c^{5}_{21}
$)

Prime. 

Not pseudo-unitary. 

\vskip 1ex 
\color{grey}

\noindent(15,4). $6_{6,4.382}^{21,604}$ \irep{43}:\ \ 
$d_i$ = ($1.0$,
$0.208$,
$1.208$,
$-0.791$,
$-0.791$,
$-0.791$) 

\vskip 0.7ex
\hangindent=3em \hangafter=1
$D^2= 4.382 = 
\frac{105-21\sqrt{21}}{2}$

\vskip 0.7ex
\hangindent=3em \hangafter=1
$T = ( 0,
0,
\frac{2}{3},
\frac{3}{7},
\frac{5}{7},
\frac{6}{7} )
$,

\vskip 0.7ex
\hangindent=3em \hangafter=1
$S$ = ($ 1$,
$ \frac{5-\sqrt{21}}{2}$,
$ \frac{7-\sqrt{21}}{2}$,
$ \frac{3-\sqrt{21}}{2}$,
$ \frac{3-\sqrt{21}}{2}$,
$ \frac{3-\sqrt{21}}{2}$;\ \ 
$ 1$,
$ \frac{7-\sqrt{21}}{2}$,
$ -\frac{3-\sqrt{21}}{2}$,
$ -\frac{3-\sqrt{21}}{2}$,
$ -\frac{3-\sqrt{21}}{2}$;\ \ 
$ -\frac{7-\sqrt{21}}{2}$,
$0$,
$0$,
$0$;\ \ 
$ c^{2}_{21}
-c^{3}_{21}
+c^{4}_{21}
-c^{5}_{21}
$,
$ 1+c^{1}_{21}
-c^{2}_{21}
-2  c^{4}_{21}
-c^{5}_{21}
$,
$ 1-2  c^{1}_{21}
+c^{3}_{21}
+c^{5}_{21}
$;\ \ 
$ 1-2  c^{1}_{21}
+c^{3}_{21}
+c^{5}_{21}
$,
$ c^{2}_{21}
-c^{3}_{21}
+c^{4}_{21}
-c^{5}_{21}
$;\ \ 
$ 1+c^{1}_{21}
-c^{2}_{21}
-2  c^{4}_{21}
-c^{5}_{21}
$)

Prime. 

Not pseudo-unitary. 

\vskip 1ex 
\color{blue}

\noindent(16,1). $6_{\frac{55}{7},18.59}^{28,124}$ \irep{46}:\ \ 
$d_i$ = ($1.0$,
$1.801$,
$2.246$,
$-1.0$,
$-1.801$,
$-2.246$) 

\vskip 0.7ex
\hangindent=3em \hangafter=1
$D^2= 18.591 = 
12+6c^{1}_{7}
+2c^{2}_{7}
$

\vskip 0.7ex
\hangindent=3em \hangafter=1
$T = ( 0,
\frac{1}{7},
\frac{5}{7},
\frac{1}{4},
\frac{11}{28},
\frac{27}{28} )
$,

\vskip 0.7ex
\hangindent=3em \hangafter=1
$S$ = ($ 1$,
$ -c_{7}^{3}$,
$ \xi_{7}^{3}$,
$ -1$,
$ c_{7}^{3}$,
$ -\xi_{7}^{3}$;\ \ 
$ -\xi_{7}^{3}$,
$ 1$,
$ c_{7}^{3}$,
$ \xi_{7}^{3}$,
$ -1$;\ \ 
$ c_{7}^{3}$,
$ -\xi_{7}^{3}$,
$ -1$,
$ -c_{7}^{3}$;\ \ 
$ -1$,
$ c_{7}^{3}$,
$ -\xi_{7}^{3}$;\ \ 
$ \xi_{7}^{3}$,
$ -1$;\ \ 
$ -c_{7}^{3}$)

Factors = $2_{1,2.}^{4,625}\boxtimes 3_{\frac{48}{7},9.295}^{7,790} $

Pseudo-unitary $\sim$  
$6_{\frac{41}{7},18.59}^{28,114}$

\vskip 1ex 
\color{grey}

\noindent(16,2). $6_{\frac{41}{7},18.59}^{28,130}$ \irep{46}:\ \ 
$d_i$ = ($1.0$,
$1.801$,
$2.246$,
$-1.0$,
$-1.801$,
$-2.246$) 

\vskip 0.7ex
\hangindent=3em \hangafter=1
$D^2= 18.591 = 
12+6c^{1}_{7}
+2c^{2}_{7}
$

\vskip 0.7ex
\hangindent=3em \hangafter=1
$T = ( 0,
\frac{1}{7},
\frac{5}{7},
\frac{3}{4},
\frac{25}{28},
\frac{13}{28} )
$,

\vskip 0.7ex
\hangindent=3em \hangafter=1
$S$ = ($ 1$,
$ -c_{7}^{3}$,
$ \xi_{7}^{3}$,
$ -1$,
$ c_{7}^{3}$,
$ -\xi_{7}^{3}$;\ \ 
$ -\xi_{7}^{3}$,
$ 1$,
$ c_{7}^{3}$,
$ \xi_{7}^{3}$,
$ -1$;\ \ 
$ c_{7}^{3}$,
$ -\xi_{7}^{3}$,
$ -1$,
$ -c_{7}^{3}$;\ \ 
$ -1$,
$ c_{7}^{3}$,
$ -\xi_{7}^{3}$;\ \ 
$ \xi_{7}^{3}$,
$ -1$;\ \ 
$ -c_{7}^{3}$)

Factors = $2_{7,2.}^{4,562}\boxtimes 3_{\frac{48}{7},9.295}^{7,790} $

Pseudo-unitary $\sim$  
$6_{\frac{55}{7},18.59}^{28,108}$

\vskip 1ex 
\color{grey}

\noindent(16,3). $6_{\frac{15}{7},18.59}^{28,779}$ \irep{46}:\ \ 
$d_i$ = ($1.0$,
$1.801$,
$2.246$,
$-1.0$,
$-1.801$,
$-2.246$) 

\vskip 0.7ex
\hangindent=3em \hangafter=1
$D^2= 18.591 = 
12+6c^{1}_{7}
+2c^{2}_{7}
$

\vskip 0.7ex
\hangindent=3em \hangafter=1
$T = ( 0,
\frac{6}{7},
\frac{2}{7},
\frac{1}{4},
\frac{3}{28},
\frac{15}{28} )
$,

\vskip 0.7ex
\hangindent=3em \hangafter=1
$S$ = ($ 1$,
$ -c_{7}^{3}$,
$ \xi_{7}^{3}$,
$ -1$,
$ c_{7}^{3}$,
$ -\xi_{7}^{3}$;\ \ 
$ -\xi_{7}^{3}$,
$ 1$,
$ c_{7}^{3}$,
$ \xi_{7}^{3}$,
$ -1$;\ \ 
$ c_{7}^{3}$,
$ -\xi_{7}^{3}$,
$ -1$,
$ -c_{7}^{3}$;\ \ 
$ -1$,
$ c_{7}^{3}$,
$ -\xi_{7}^{3}$;\ \ 
$ \xi_{7}^{3}$,
$ -1$;\ \ 
$ -c_{7}^{3}$)

Factors = $2_{1,2.}^{4,625}\boxtimes 3_{\frac{8}{7},9.295}^{7,245} $

Pseudo-unitary $\sim$  
$6_{\frac{1}{7},18.59}^{28,212}$

\vskip 1ex 
\color{grey}

\noindent(16,4). $6_{\frac{1}{7},18.59}^{28,277}$ \irep{46}:\ \ 
$d_i$ = ($1.0$,
$1.801$,
$2.246$,
$-1.0$,
$-1.801$,
$-2.246$) 

\vskip 0.7ex
\hangindent=3em \hangafter=1
$D^2= 18.591 = 
12+6c^{1}_{7}
+2c^{2}_{7}
$

\vskip 0.7ex
\hangindent=3em \hangafter=1
$T = ( 0,
\frac{6}{7},
\frac{2}{7},
\frac{3}{4},
\frac{17}{28},
\frac{1}{28} )
$,

\vskip 0.7ex
\hangindent=3em \hangafter=1
$S$ = ($ 1$,
$ -c_{7}^{3}$,
$ \xi_{7}^{3}$,
$ -1$,
$ c_{7}^{3}$,
$ -\xi_{7}^{3}$;\ \ 
$ -\xi_{7}^{3}$,
$ 1$,
$ c_{7}^{3}$,
$ \xi_{7}^{3}$,
$ -1$;\ \ 
$ c_{7}^{3}$,
$ -\xi_{7}^{3}$,
$ -1$,
$ -c_{7}^{3}$;\ \ 
$ -1$,
$ c_{7}^{3}$,
$ -\xi_{7}^{3}$;\ \ 
$ \xi_{7}^{3}$,
$ -1$;\ \ 
$ -c_{7}^{3}$)

Factors = $2_{7,2.}^{4,562}\boxtimes 3_{\frac{8}{7},9.295}^{7,245} $

Pseudo-unitary $\sim$  
$6_{\frac{15}{7},18.59}^{28,289}$

\vskip 1ex 
\color{grey}

\noindent(16,5). $6_{\frac{5}{7},5.725}^{28,578}$ \irep{46}:\ \ 
$d_i$ = ($1.0$,
$0.554$,
$1.246$,
$-0.554$,
$-1.0$,
$-1.246$) 

\vskip 0.7ex
\hangindent=3em \hangafter=1
$D^2= 5.725 = 
10-2  c^{1}_{7}
+4c^{2}_{7}
$

\vskip 0.7ex
\hangindent=3em \hangafter=1
$T = ( 0,
\frac{3}{7},
\frac{1}{28},
\frac{5}{28},
\frac{3}{4},
\frac{2}{7} )
$,

\vskip 0.7ex
\hangindent=3em \hangafter=1
$S$ = ($ 1$,
$ 1+c^{2}_{7}
$,
$ c_{7}^{1}$,
$ -1-c^{2}_{7}
$,
$ -1$,
$ -c_{7}^{1}$;\ \ 
$ c_{7}^{1}$,
$ -1$,
$ -c_{7}^{1}$,
$ -1-c^{2}_{7}
$,
$ 1$;\ \ 
$ 1+c^{2}_{7}
$,
$ -1$,
$ c_{7}^{1}$,
$ 1+c^{2}_{7}
$;\ \ 
$ -c_{7}^{1}$,
$ -1-c^{2}_{7}
$,
$ -1$;\ \ 
$ -1$,
$ c_{7}^{1}$;\ \ 
$ -1-c^{2}_{7}
$)

Factors = $2_{7,2.}^{4,562}\boxtimes 3_{\frac{12}{7},2.862}^{7,768} $

Not pseudo-unitary. 

\vskip 1ex 
\color{grey}

\noindent(16,6). $6_{\frac{19}{7},5.725}^{28,399}$ \irep{46}:\ \ 
$d_i$ = ($1.0$,
$0.554$,
$1.246$,
$-0.554$,
$-1.0$,
$-1.246$) 

\vskip 0.7ex
\hangindent=3em \hangafter=1
$D^2= 5.725 = 
10-2  c^{1}_{7}
+4c^{2}_{7}
$

\vskip 0.7ex
\hangindent=3em \hangafter=1
$T = ( 0,
\frac{3}{7},
\frac{15}{28},
\frac{19}{28},
\frac{1}{4},
\frac{2}{7} )
$,

\vskip 0.7ex
\hangindent=3em \hangafter=1
$S$ = ($ 1$,
$ 1+c^{2}_{7}
$,
$ c_{7}^{1}$,
$ -1-c^{2}_{7}
$,
$ -1$,
$ -c_{7}^{1}$;\ \ 
$ c_{7}^{1}$,
$ -1$,
$ -c_{7}^{1}$,
$ -1-c^{2}_{7}
$,
$ 1$;\ \ 
$ 1+c^{2}_{7}
$,
$ -1$,
$ c_{7}^{1}$,
$ 1+c^{2}_{7}
$;\ \ 
$ -c_{7}^{1}$,
$ -1-c^{2}_{7}
$,
$ -1$;\ \ 
$ -1$,
$ c_{7}^{1}$;\ \ 
$ -1-c^{2}_{7}
$)

Factors = $2_{1,2.}^{4,625}\boxtimes 3_{\frac{12}{7},2.862}^{7,768} $

Not pseudo-unitary. 

\vskip 1ex 
\color{grey}

\noindent(16,7). $6_{\frac{37}{7},5.725}^{28,806}$ \irep{46}:\ \ 
$d_i$ = ($1.0$,
$0.554$,
$1.246$,
$-0.554$,
$-1.0$,
$-1.246$) 

\vskip 0.7ex
\hangindent=3em \hangafter=1
$D^2= 5.725 = 
10-2  c^{1}_{7}
+4c^{2}_{7}
$

\vskip 0.7ex
\hangindent=3em \hangafter=1
$T = ( 0,
\frac{4}{7},
\frac{13}{28},
\frac{9}{28},
\frac{3}{4},
\frac{5}{7} )
$,

\vskip 0.7ex
\hangindent=3em \hangafter=1
$S$ = ($ 1$,
$ 1+c^{2}_{7}
$,
$ c_{7}^{1}$,
$ -1-c^{2}_{7}
$,
$ -1$,
$ -c_{7}^{1}$;\ \ 
$ c_{7}^{1}$,
$ -1$,
$ -c_{7}^{1}$,
$ -1-c^{2}_{7}
$,
$ 1$;\ \ 
$ 1+c^{2}_{7}
$,
$ -1$,
$ c_{7}^{1}$,
$ 1+c^{2}_{7}
$;\ \ 
$ -c_{7}^{1}$,
$ -1-c^{2}_{7}
$,
$ -1$;\ \ 
$ -1$,
$ c_{7}^{1}$;\ \ 
$ -1-c^{2}_{7}
$)

Factors = $2_{7,2.}^{4,562}\boxtimes 3_{\frac{44}{7},2.862}^{7,531} $

Not pseudo-unitary. 

\vskip 1ex 
\color{grey}

\noindent(16,8). $6_{\frac{51}{7},5.725}^{28,267}$ \irep{46}:\ \ 
$d_i$ = ($1.0$,
$0.554$,
$1.246$,
$-0.554$,
$-1.0$,
$-1.246$) 

\vskip 0.7ex
\hangindent=3em \hangafter=1
$D^2= 5.725 = 
10-2  c^{1}_{7}
+4c^{2}_{7}
$

\vskip 0.7ex
\hangindent=3em \hangafter=1
$T = ( 0,
\frac{4}{7},
\frac{27}{28},
\frac{23}{28},
\frac{1}{4},
\frac{5}{7} )
$,

\vskip 0.7ex
\hangindent=3em \hangafter=1
$S$ = ($ 1$,
$ 1+c^{2}_{7}
$,
$ c_{7}^{1}$,
$ -1-c^{2}_{7}
$,
$ -1$,
$ -c_{7}^{1}$;\ \ 
$ c_{7}^{1}$,
$ -1$,
$ -c_{7}^{1}$,
$ -1-c^{2}_{7}
$,
$ 1$;\ \ 
$ 1+c^{2}_{7}
$,
$ -1$,
$ c_{7}^{1}$,
$ 1+c^{2}_{7}
$;\ \ 
$ -c_{7}^{1}$,
$ -1-c^{2}_{7}
$,
$ -1$;\ \ 
$ -1$,
$ c_{7}^{1}$;\ \ 
$ -1-c^{2}_{7}
$)

Factors = $2_{1,2.}^{4,625}\boxtimes 3_{\frac{44}{7},2.862}^{7,531} $

Not pseudo-unitary. 

\vskip 1ex 
\color{grey}

\noindent(16,9). $6_{\frac{11}{7},3.682}^{28,782}$ \irep{46}:\ \ 
$d_i$ = ($1.0$,
$0.445$,
$0.801$,
$-0.445$,
$-0.801$,
$-1.0$) 

\vskip 0.7ex
\hangindent=3em \hangafter=1
$D^2= 3.682 = 
6-4  c^{1}_{7}
-6  c^{2}_{7}
$

\vskip 0.7ex
\hangindent=3em \hangafter=1
$T = ( 0,
\frac{3}{7},
\frac{11}{28},
\frac{19}{28},
\frac{1}{7},
\frac{1}{4} )
$,

\vskip 0.7ex
\hangindent=3em \hangafter=1
$S$ = ($ 1$,
$ -c_{7}^{2}$,
$ c^{1}_{7}
+c^{2}_{7}
$,
$ c_{7}^{2}$,
$ -c^{1}_{7}
-c^{2}_{7}
$,
$ -1$;\ \ 
$ c^{1}_{7}
+c^{2}_{7}
$,
$ -1$,
$ -c^{1}_{7}
-c^{2}_{7}
$,
$ 1$,
$ c_{7}^{2}$;\ \ 
$ -c_{7}^{2}$,
$ -1$,
$ -c_{7}^{2}$,
$ c^{1}_{7}
+c^{2}_{7}
$;\ \ 
$ -c^{1}_{7}
-c^{2}_{7}
$,
$ -1$,
$ c_{7}^{2}$;\ \ 
$ c_{7}^{2}$,
$ c^{1}_{7}
+c^{2}_{7}
$;\ \ 
$ -1$)

Factors = $2_{1,2.}^{4,625}\boxtimes 3_{\frac{4}{7},1.841}^{7,953} $

Not pseudo-unitary. 

\vskip 1ex 
\color{grey}

\noindent(16,10). $6_{\frac{53}{7},3.682}^{28,127}$ \irep{46}:\ \ 
$d_i$ = ($1.0$,
$0.445$,
$0.801$,
$-0.445$,
$-0.801$,
$-1.0$) 

\vskip 0.7ex
\hangindent=3em \hangafter=1
$D^2= 3.682 = 
6-4  c^{1}_{7}
-6  c^{2}_{7}
$

\vskip 0.7ex
\hangindent=3em \hangafter=1
$T = ( 0,
\frac{3}{7},
\frac{25}{28},
\frac{5}{28},
\frac{1}{7},
\frac{3}{4} )
$,

\vskip 0.7ex
\hangindent=3em \hangafter=1
$S$ = ($ 1$,
$ -c_{7}^{2}$,
$ c^{1}_{7}
+c^{2}_{7}
$,
$ c_{7}^{2}$,
$ -c^{1}_{7}
-c^{2}_{7}
$,
$ -1$;\ \ 
$ c^{1}_{7}
+c^{2}_{7}
$,
$ -1$,
$ -c^{1}_{7}
-c^{2}_{7}
$,
$ 1$,
$ c_{7}^{2}$;\ \ 
$ -c_{7}^{2}$,
$ -1$,
$ -c_{7}^{2}$,
$ c^{1}_{7}
+c^{2}_{7}
$;\ \ 
$ -c^{1}_{7}
-c^{2}_{7}
$,
$ -1$,
$ c_{7}^{2}$;\ \ 
$ c_{7}^{2}$,
$ c^{1}_{7}
+c^{2}_{7}
$;\ \ 
$ -1$)

Factors = $2_{7,2.}^{4,562}\boxtimes 3_{\frac{4}{7},1.841}^{7,953} $

Not pseudo-unitary. 

\vskip 1ex 
\color{grey}

\noindent(16,11). $6_{\frac{3}{7},3.682}^{28,797}$ \irep{46}:\ \ 
$d_i$ = ($1.0$,
$0.445$,
$0.801$,
$-0.445$,
$-0.801$,
$-1.0$) 

\vskip 0.7ex
\hangindent=3em \hangafter=1
$D^2= 3.682 = 
6-4  c^{1}_{7}
-6  c^{2}_{7}
$

\vskip 0.7ex
\hangindent=3em \hangafter=1
$T = ( 0,
\frac{4}{7},
\frac{3}{28},
\frac{23}{28},
\frac{6}{7},
\frac{1}{4} )
$,

\vskip 0.7ex
\hangindent=3em \hangafter=1
$S$ = ($ 1$,
$ -c_{7}^{2}$,
$ c^{1}_{7}
+c^{2}_{7}
$,
$ c_{7}^{2}$,
$ -c^{1}_{7}
-c^{2}_{7}
$,
$ -1$;\ \ 
$ c^{1}_{7}
+c^{2}_{7}
$,
$ -1$,
$ -c^{1}_{7}
-c^{2}_{7}
$,
$ 1$,
$ c_{7}^{2}$;\ \ 
$ -c_{7}^{2}$,
$ -1$,
$ -c_{7}^{2}$,
$ c^{1}_{7}
+c^{2}_{7}
$;\ \ 
$ -c^{1}_{7}
-c^{2}_{7}
$,
$ -1$,
$ c_{7}^{2}$;\ \ 
$ c_{7}^{2}$,
$ c^{1}_{7}
+c^{2}_{7}
$;\ \ 
$ -1$)

Factors = $2_{1,2.}^{4,625}\boxtimes 3_{\frac{52}{7},1.841}^{7,604} $

Not pseudo-unitary. 

\vskip 1ex 
\color{grey}

\noindent(16,12). $6_{\frac{45}{7},3.682}^{28,756}$ \irep{46}:\ \ 
$d_i$ = ($1.0$,
$0.445$,
$0.801$,
$-0.445$,
$-0.801$,
$-1.0$) 

\vskip 0.7ex
\hangindent=3em \hangafter=1
$D^2= 3.682 = 
6-4  c^{1}_{7}
-6  c^{2}_{7}
$

\vskip 0.7ex
\hangindent=3em \hangafter=1
$T = ( 0,
\frac{4}{7},
\frac{17}{28},
\frac{9}{28},
\frac{6}{7},
\frac{3}{4} )
$,

\vskip 0.7ex
\hangindent=3em \hangafter=1
$S$ = ($ 1$,
$ -c_{7}^{2}$,
$ c^{1}_{7}
+c^{2}_{7}
$,
$ c_{7}^{2}$,
$ -c^{1}_{7}
-c^{2}_{7}
$,
$ -1$;\ \ 
$ c^{1}_{7}
+c^{2}_{7}
$,
$ -1$,
$ -c^{1}_{7}
-c^{2}_{7}
$,
$ 1$,
$ c_{7}^{2}$;\ \ 
$ -c_{7}^{2}$,
$ -1$,
$ -c_{7}^{2}$,
$ c^{1}_{7}
+c^{2}_{7}
$;\ \ 
$ -c^{1}_{7}
-c^{2}_{7}
$,
$ -1$,
$ c_{7}^{2}$;\ \ 
$ c_{7}^{2}$,
$ c^{1}_{7}
+c^{2}_{7}
$;\ \ 
$ -1$)

Factors = $2_{7,2.}^{4,562}\boxtimes 3_{\frac{52}{7},1.841}^{7,604} $

Not pseudo-unitary. 

\vskip 1ex 
\color{blue}

\noindent(17,1). $6_{\frac{16}{3},9.}^{9,746}$ \irep{25}:\ \ 
$d_i$ = ($1.0$,
$0.347$,
$1.0$,
$1.532$,
$-1.0$,
$-1.879$) 

\vskip 0.7ex
\hangindent=3em \hangafter=1
$D^2= 9.0 = 
9$

\vskip 0.7ex
\hangindent=3em \hangafter=1
$T = ( 0,
\frac{1}{9},
\frac{2}{3},
\frac{4}{9},
\frac{1}{3},
\frac{7}{9} )
$,

\vskip 0.7ex
\hangindent=3em \hangafter=1
$S$ = ($ 1$,
$ c_{9}^{2}$,
$ 1$,
$ c_{9}^{1}$,
$ -1$,
$ c_{9}^{4}$;\ \ 
$ 1$,
$ c_{9}^{1}$,
$ 1$,
$ -c_{9}^{4}$,
$ 1$;\ \ 
$ 1$,
$ c_{9}^{4}$,
$ -1$,
$ c_{9}^{2}$;\ \ 
$ 1$,
$ -c_{9}^{2}$,
$ 1$;\ \ 
$ 1$,
$ -c_{9}^{1}$;\ \ 
$ 1$)

Prime. 

Not pseudo-unitary. 

\vskip 1ex 
\color{grey}

\noindent(17,2). $6_{\frac{8}{3},9.}^{9,711}$ \irep{25}:\ \ 
$d_i$ = ($1.0$,
$0.347$,
$1.0$,
$1.532$,
$-1.0$,
$-1.879$) 

\vskip 0.7ex
\hangindent=3em \hangafter=1
$D^2= 9.0 = 
9$

\vskip 0.7ex
\hangindent=3em \hangafter=1
$T = ( 0,
\frac{8}{9},
\frac{1}{3},
\frac{5}{9},
\frac{2}{3},
\frac{2}{9} )
$,

\vskip 0.7ex
\hangindent=3em \hangafter=1
$S$ = ($ 1$,
$ c_{9}^{2}$,
$ 1$,
$ c_{9}^{1}$,
$ -1$,
$ c_{9}^{4}$;\ \ 
$ 1$,
$ c_{9}^{1}$,
$ 1$,
$ -c_{9}^{4}$,
$ 1$;\ \ 
$ 1$,
$ c_{9}^{4}$,
$ -1$,
$ c_{9}^{2}$;\ \ 
$ 1$,
$ -c_{9}^{2}$,
$ 1$;\ \ 
$ 1$,
$ -c_{9}^{1}$;\ \ 
$ 1$)

Prime. 

Not pseudo-unitary. 

\vskip 1ex 
\color{blue}

\noindent(18,1). $6_{\frac{3}{2},8.}^{16,507}$ \irep{39}:\ \ 
$d_i$ = ($1.0$,
$1.0$,
$1.414$,
$-1.0$,
$-1.0$,
$-1.414$) 

\vskip 0.7ex
\hangindent=3em \hangafter=1
$D^2= 8.0 = 
8$

\vskip 0.7ex
\hangindent=3em \hangafter=1
$T = ( 0,
\frac{1}{2},
\frac{1}{16},
\frac{1}{4},
\frac{3}{4},
\frac{5}{16} )
$,

\vskip 0.7ex
\hangindent=3em \hangafter=1
$S$ = ($ 1$,
$ 1$,
$ \sqrt{2}$,
$ -1$,
$ -1$,
$ -\sqrt{2}$;\ \ 
$ 1$,
$ -\sqrt{2}$,
$ -1$,
$ -1$,
$ \sqrt{2}$;\ \ 
$0$,
$ -\sqrt{2}$,
$ \sqrt{2}$,
$0$;\ \ 
$ -1$,
$ -1$,
$ -\sqrt{2}$;\ \ 
$ -1$,
$ \sqrt{2}$;\ \ 
$0$)

Factors = $2_{1,2.}^{4,625}\boxtimes 3_{\frac{1}{2},4.}^{16,598} $

Pseudo-unitary $\sim$  
$6_{\frac{15}{2},8.}^{16,107}$

\vskip 1ex 
\color{grey}

\noindent(18,2). $6_{\frac{15}{2},8.}^{16,130}$ \irep{39}:\ \ 
$d_i$ = ($1.0$,
$1.0$,
$1.414$,
$-1.0$,
$-1.0$,
$-1.414$) 

\vskip 0.7ex
\hangindent=3em \hangafter=1
$D^2= 8.0 = 
8$

\vskip 0.7ex
\hangindent=3em \hangafter=1
$T = ( 0,
\frac{1}{2},
\frac{1}{16},
\frac{1}{4},
\frac{3}{4},
\frac{13}{16} )
$,

\vskip 0.7ex
\hangindent=3em \hangafter=1
$S$ = ($ 1$,
$ 1$,
$ \sqrt{2}$,
$ -1$,
$ -1$,
$ -\sqrt{2}$;\ \ 
$ 1$,
$ -\sqrt{2}$,
$ -1$,
$ -1$,
$ \sqrt{2}$;\ \ 
$0$,
$ \sqrt{2}$,
$ -\sqrt{2}$,
$0$;\ \ 
$ -1$,
$ -1$,
$ \sqrt{2}$;\ \ 
$ -1$,
$ -\sqrt{2}$;\ \ 
$0$)

Factors = $2_{1,2.}^{4,625}\boxtimes 3_{\frac{13}{2},4.}^{16,830} $

Pseudo-unitary $\sim$  
$6_{\frac{3}{2},8.}^{16,688}$

\vskip 1ex 
\color{grey}

\noindent(18,3). $6_{\frac{5}{2},8.}^{16,154}$ \irep{39}:\ \ 
$d_i$ = ($1.0$,
$1.0$,
$1.414$,
$-1.0$,
$-1.0$,
$-1.414$) 

\vskip 0.7ex
\hangindent=3em \hangafter=1
$D^2= 8.0 = 
8$

\vskip 0.7ex
\hangindent=3em \hangafter=1
$T = ( 0,
\frac{1}{2},
\frac{7}{16},
\frac{1}{4},
\frac{3}{4},
\frac{3}{16} )
$,

\vskip 0.7ex
\hangindent=3em \hangafter=1
$S$ = ($ 1$,
$ 1$,
$ \sqrt{2}$,
$ -1$,
$ -1$,
$ -\sqrt{2}$;\ \ 
$ 1$,
$ -\sqrt{2}$,
$ -1$,
$ -1$,
$ \sqrt{2}$;\ \ 
$0$,
$ \sqrt{2}$,
$ -\sqrt{2}$,
$0$;\ \ 
$ -1$,
$ -1$,
$ \sqrt{2}$;\ \ 
$ -1$,
$ -\sqrt{2}$;\ \ 
$0$)

Factors = $2_{1,2.}^{4,625}\boxtimes 3_{\frac{3}{2},4.}^{16,538} $

Pseudo-unitary $\sim$  
$6_{\frac{9}{2},8.}^{16,107}$

\vskip 1ex 
\color{grey}

\noindent(18,4). $6_{\frac{9}{2},8.}^{16,772}$ \irep{39}:\ \ 
$d_i$ = ($1.0$,
$1.0$,
$1.414$,
$-1.0$,
$-1.0$,
$-1.414$) 

\vskip 0.7ex
\hangindent=3em \hangafter=1
$D^2= 8.0 = 
8$

\vskip 0.7ex
\hangindent=3em \hangafter=1
$T = ( 0,
\frac{1}{2},
\frac{7}{16},
\frac{1}{4},
\frac{3}{4},
\frac{11}{16} )
$,

\vskip 0.7ex
\hangindent=3em \hangafter=1
$S$ = ($ 1$,
$ 1$,
$ \sqrt{2}$,
$ -1$,
$ -1$,
$ -\sqrt{2}$;\ \ 
$ 1$,
$ -\sqrt{2}$,
$ -1$,
$ -1$,
$ \sqrt{2}$;\ \ 
$0$,
$ -\sqrt{2}$,
$ \sqrt{2}$,
$0$;\ \ 
$ -1$,
$ -1$,
$ -\sqrt{2}$;\ \ 
$ -1$,
$ \sqrt{2}$;\ \ 
$0$)

Factors = $2_{1,2.}^{4,625}\boxtimes 3_{\frac{7}{2},4.}^{16,332} $

Pseudo-unitary $\sim$  
$6_{\frac{5}{2},8.}^{16,511}$

\vskip 1ex 
\color{grey}

\noindent(18,5). $6_{\frac{7}{2},8.}^{16,656}$ \irep{39}:\ \ 
$d_i$ = ($1.0$,
$1.0$,
$1.414$,
$-1.0$,
$-1.0$,
$-1.414$) 

\vskip 0.7ex
\hangindent=3em \hangafter=1
$D^2= 8.0 = 
8$

\vskip 0.7ex
\hangindent=3em \hangafter=1
$T = ( 0,
\frac{1}{2},
\frac{9}{16},
\frac{1}{4},
\frac{3}{4},
\frac{5}{16} )
$,

\vskip 0.7ex
\hangindent=3em \hangafter=1
$S$ = ($ 1$,
$ 1$,
$ \sqrt{2}$,
$ -1$,
$ -1$,
$ -\sqrt{2}$;\ \ 
$ 1$,
$ -\sqrt{2}$,
$ -1$,
$ -1$,
$ \sqrt{2}$;\ \ 
$0$,
$ \sqrt{2}$,
$ -\sqrt{2}$,
$0$;\ \ 
$ -1$,
$ -1$,
$ \sqrt{2}$;\ \ 
$ -1$,
$ -\sqrt{2}$;\ \ 
$0$)

Factors = $2_{1,2.}^{4,625}\boxtimes 3_{\frac{5}{2},4.}^{16,345} $

Pseudo-unitary $\sim$  
$6_{\frac{11}{2},8.}^{16,548}$

\vskip 1ex 
\color{grey}

\noindent(18,6). $6_{\frac{11}{2},8.}^{16,861}$ \irep{39}:\ \ 
$d_i$ = ($1.0$,
$1.0$,
$1.414$,
$-1.0$,
$-1.0$,
$-1.414$) 

\vskip 0.7ex
\hangindent=3em \hangafter=1
$D^2= 8.0 = 
8$

\vskip 0.7ex
\hangindent=3em \hangafter=1
$T = ( 0,
\frac{1}{2},
\frac{9}{16},
\frac{1}{4},
\frac{3}{4},
\frac{13}{16} )
$,

\vskip 0.7ex
\hangindent=3em \hangafter=1
$S$ = ($ 1$,
$ 1$,
$ \sqrt{2}$,
$ -1$,
$ -1$,
$ -\sqrt{2}$;\ \ 
$ 1$,
$ -\sqrt{2}$,
$ -1$,
$ -1$,
$ \sqrt{2}$;\ \ 
$0$,
$ -\sqrt{2}$,
$ \sqrt{2}$,
$0$;\ \ 
$ -1$,
$ -1$,
$ -\sqrt{2}$;\ \ 
$ -1$,
$ \sqrt{2}$;\ \ 
$0$)

Factors = $2_{1,2.}^{4,625}\boxtimes 3_{\frac{9}{2},4.}^{16,156} $

Pseudo-unitary $\sim$  
$6_{\frac{7}{2},8.}^{16,246}$

\vskip 1ex 
\color{grey}

\noindent(18,7). $6_{\frac{1}{2},8.}^{16,818}$ \irep{39}:\ \ 
$d_i$ = ($1.0$,
$1.0$,
$1.414$,
$-1.0$,
$-1.0$,
$-1.414$) 

\vskip 0.7ex
\hangindent=3em \hangafter=1
$D^2= 8.0 = 
8$

\vskip 0.7ex
\hangindent=3em \hangafter=1
$T = ( 0,
\frac{1}{2},
\frac{15}{16},
\frac{1}{4},
\frac{3}{4},
\frac{3}{16} )
$,

\vskip 0.7ex
\hangindent=3em \hangafter=1
$S$ = ($ 1$,
$ 1$,
$ \sqrt{2}$,
$ -1$,
$ -1$,
$ -\sqrt{2}$;\ \ 
$ 1$,
$ -\sqrt{2}$,
$ -1$,
$ -1$,
$ \sqrt{2}$;\ \ 
$0$,
$ -\sqrt{2}$,
$ \sqrt{2}$,
$0$;\ \ 
$ -1$,
$ -1$,
$ -\sqrt{2}$;\ \ 
$ -1$,
$ \sqrt{2}$;\ \ 
$0$)

Factors = $2_{1,2.}^{4,625}\boxtimes 3_{\frac{15}{2},4.}^{16,639} $

Pseudo-unitary $\sim$  
$6_{\frac{13}{2},8.}^{16,107}$

\vskip 1ex 
\color{grey}

\noindent(18,8). $6_{\frac{13}{2},8.}^{16,199}$ \irep{39}:\ \ 
$d_i$ = ($1.0$,
$1.0$,
$1.414$,
$-1.0$,
$-1.0$,
$-1.414$) 

\vskip 0.7ex
\hangindent=3em \hangafter=1
$D^2= 8.0 = 
8$

\vskip 0.7ex
\hangindent=3em \hangafter=1
$T = ( 0,
\frac{1}{2},
\frac{15}{16},
\frac{1}{4},
\frac{3}{4},
\frac{11}{16} )
$,

\vskip 0.7ex
\hangindent=3em \hangafter=1
$S$ = ($ 1$,
$ 1$,
$ \sqrt{2}$,
$ -1$,
$ -1$,
$ -\sqrt{2}$;\ \ 
$ 1$,
$ -\sqrt{2}$,
$ -1$,
$ -1$,
$ \sqrt{2}$;\ \ 
$0$,
$ \sqrt{2}$,
$ -\sqrt{2}$,
$0$;\ \ 
$ -1$,
$ -1$,
$ \sqrt{2}$;\ \ 
$ -1$,
$ -\sqrt{2}$;\ \ 
$0$)

Factors = $2_{1,2.}^{4,625}\boxtimes 3_{\frac{11}{2},4.}^{16,564} $

Pseudo-unitary $\sim$  
$6_{\frac{1}{2},8.}^{16,460}$

\vskip 1ex 
\color{blue}

\noindent(19,1). $6_{\frac{5}{2},8.}^{16,595}$ \irep{39}:\ \ 
$d_i$ = ($1.0$,
$1.0$,
$1.414$,
$-1.0$,
$-1.0$,
$-1.414$) 

\vskip 0.7ex
\hangindent=3em \hangafter=1
$D^2= 8.0 = 
8$

\vskip 0.7ex
\hangindent=3em \hangafter=1
$T = ( 0,
\frac{1}{2},
\frac{3}{16},
\frac{1}{4},
\frac{3}{4},
\frac{7}{16} )
$,

\vskip 0.7ex
\hangindent=3em \hangafter=1
$S$ = ($ 1$,
$ 1$,
$ \sqrt{2}$,
$ -1$,
$ -1$,
$ -\sqrt{2}$;\ \ 
$ 1$,
$ -\sqrt{2}$,
$ -1$,
$ -1$,
$ \sqrt{2}$;\ \ 
$0$,
$ -\sqrt{2}$,
$ \sqrt{2}$,
$0$;\ \ 
$ -1$,
$ -1$,
$ -\sqrt{2}$;\ \ 
$ -1$,
$ \sqrt{2}$;\ \ 
$0$)

Factors = $2_{1,2.}^{4,625}\boxtimes 3_{\frac{3}{2},4.}^{16,553} $

Pseudo-unitary $\sim$  
$6_{\frac{1}{2},8.}^{16,460}$

\vskip 1ex 
\color{grey}

\noindent(19,2). $6_{\frac{1}{2},8.}^{16,156}$ \irep{39}:\ \ 
$d_i$ = ($1.0$,
$1.0$,
$1.414$,
$-1.0$,
$-1.0$,
$-1.414$) 

\vskip 0.7ex
\hangindent=3em \hangafter=1
$D^2= 8.0 = 
8$

\vskip 0.7ex
\hangindent=3em \hangafter=1
$T = ( 0,
\frac{1}{2},
\frac{3}{16},
\frac{1}{4},
\frac{3}{4},
\frac{15}{16} )
$,

\vskip 0.7ex
\hangindent=3em \hangafter=1
$S$ = ($ 1$,
$ 1$,
$ \sqrt{2}$,
$ -1$,
$ -1$,
$ -\sqrt{2}$;\ \ 
$ 1$,
$ -\sqrt{2}$,
$ -1$,
$ -1$,
$ \sqrt{2}$;\ \ 
$0$,
$ \sqrt{2}$,
$ -\sqrt{2}$,
$0$;\ \ 
$ -1$,
$ -1$,
$ \sqrt{2}$;\ \ 
$ -1$,
$ -\sqrt{2}$;\ \ 
$0$)

Factors = $2_{1,2.}^{4,625}\boxtimes 3_{\frac{15}{2},4.}^{16,113} $

Pseudo-unitary $\sim$  
$6_{\frac{5}{2},8.}^{16,511}$

\vskip 1ex 
\color{grey}

\noindent(19,3). $6_{\frac{3}{2},8.}^{16,242}$ \irep{39}:\ \ 
$d_i$ = ($1.0$,
$1.0$,
$1.414$,
$-1.0$,
$-1.0$,
$-1.414$) 

\vskip 0.7ex
\hangindent=3em \hangafter=1
$D^2= 8.0 = 
8$

\vskip 0.7ex
\hangindent=3em \hangafter=1
$T = ( 0,
\frac{1}{2},
\frac{5}{16},
\frac{1}{4},
\frac{3}{4},
\frac{1}{16} )
$,

\vskip 0.7ex
\hangindent=3em \hangafter=1
$S$ = ($ 1$,
$ 1$,
$ \sqrt{2}$,
$ -1$,
$ -1$,
$ -\sqrt{2}$;\ \ 
$ 1$,
$ -\sqrt{2}$,
$ -1$,
$ -1$,
$ \sqrt{2}$;\ \ 
$0$,
$ \sqrt{2}$,
$ -\sqrt{2}$,
$0$;\ \ 
$ -1$,
$ -1$,
$ \sqrt{2}$;\ \ 
$ -1$,
$ -\sqrt{2}$;\ \ 
$0$)

Factors = $2_{1,2.}^{4,625}\boxtimes 3_{\frac{1}{2},4.}^{16,980} $

Pseudo-unitary $\sim$  
$6_{\frac{7}{2},8.}^{16,246}$

\vskip 1ex 
\color{grey}

\noindent(19,4). $6_{\frac{7}{2},8.}^{16,684}$ \irep{39}:\ \ 
$d_i$ = ($1.0$,
$1.0$,
$1.414$,
$-1.0$,
$-1.0$,
$-1.414$) 

\vskip 0.7ex
\hangindent=3em \hangafter=1
$D^2= 8.0 = 
8$

\vskip 0.7ex
\hangindent=3em \hangafter=1
$T = ( 0,
\frac{1}{2},
\frac{5}{16},
\frac{1}{4},
\frac{3}{4},
\frac{9}{16} )
$,

\vskip 0.7ex
\hangindent=3em \hangafter=1
$S$ = ($ 1$,
$ 1$,
$ \sqrt{2}$,
$ -1$,
$ -1$,
$ -\sqrt{2}$;\ \ 
$ 1$,
$ -\sqrt{2}$,
$ -1$,
$ -1$,
$ \sqrt{2}$;\ \ 
$0$,
$ -\sqrt{2}$,
$ \sqrt{2}$,
$0$;\ \ 
$ -1$,
$ -1$,
$ -\sqrt{2}$;\ \ 
$ -1$,
$ \sqrt{2}$;\ \ 
$0$)

Factors = $2_{1,2.}^{4,625}\boxtimes 3_{\frac{5}{2},4.}^{16,465} $

Pseudo-unitary $\sim$  
$6_{\frac{3}{2},8.}^{16,688}$

\vskip 1ex 
\color{grey}

\noindent(19,5). $6_{\frac{9}{2},8.}^{16,227}$ \irep{39}:\ \ 
$d_i$ = ($1.0$,
$1.0$,
$1.414$,
$-1.0$,
$-1.0$,
$-1.414$) 

\vskip 0.7ex
\hangindent=3em \hangafter=1
$D^2= 8.0 = 
8$

\vskip 0.7ex
\hangindent=3em \hangafter=1
$T = ( 0,
\frac{1}{2},
\frac{11}{16},
\frac{1}{4},
\frac{3}{4},
\frac{7}{16} )
$,

\vskip 0.7ex
\hangindent=3em \hangafter=1
$S$ = ($ 1$,
$ 1$,
$ \sqrt{2}$,
$ -1$,
$ -1$,
$ -\sqrt{2}$;\ \ 
$ 1$,
$ -\sqrt{2}$,
$ -1$,
$ -1$,
$ \sqrt{2}$;\ \ 
$0$,
$ \sqrt{2}$,
$ -\sqrt{2}$,
$0$;\ \ 
$ -1$,
$ -1$,
$ \sqrt{2}$;\ \ 
$ -1$,
$ -\sqrt{2}$;\ \ 
$0$)

Factors = $2_{1,2.}^{4,625}\boxtimes 3_{\frac{7}{2},4.}^{16,167} $

Pseudo-unitary $\sim$  
$6_{\frac{13}{2},8.}^{16,107}$

\vskip 1ex 
\color{grey}

\noindent(19,6). $6_{\frac{13}{2},8.}^{16,949}$ \irep{39}:\ \ 
$d_i$ = ($1.0$,
$1.0$,
$1.414$,
$-1.0$,
$-1.0$,
$-1.414$) 

\vskip 0.7ex
\hangindent=3em \hangafter=1
$D^2= 8.0 = 
8$

\vskip 0.7ex
\hangindent=3em \hangafter=1
$T = ( 0,
\frac{1}{2},
\frac{11}{16},
\frac{1}{4},
\frac{3}{4},
\frac{15}{16} )
$,

\vskip 0.7ex
\hangindent=3em \hangafter=1
$S$ = ($ 1$,
$ 1$,
$ \sqrt{2}$,
$ -1$,
$ -1$,
$ -\sqrt{2}$;\ \ 
$ 1$,
$ -\sqrt{2}$,
$ -1$,
$ -1$,
$ \sqrt{2}$;\ \ 
$0$,
$ -\sqrt{2}$,
$ \sqrt{2}$,
$0$;\ \ 
$ -1$,
$ -1$,
$ -\sqrt{2}$;\ \ 
$ -1$,
$ \sqrt{2}$;\ \ 
$0$)

Factors = $2_{1,2.}^{4,625}\boxtimes 3_{\frac{11}{2},4.}^{16,648} $

Pseudo-unitary $\sim$  
$6_{\frac{9}{2},8.}^{16,107}$

\vskip 1ex 
\color{grey}

\noindent(19,7). $6_{\frac{15}{2},8.}^{16,553}$ \irep{39}:\ \ 
$d_i$ = ($1.0$,
$1.0$,
$1.414$,
$-1.0$,
$-1.0$,
$-1.414$) 

\vskip 0.7ex
\hangindent=3em \hangafter=1
$D^2= 8.0 = 
8$

\vskip 0.7ex
\hangindent=3em \hangafter=1
$T = ( 0,
\frac{1}{2},
\frac{13}{16},
\frac{1}{4},
\frac{3}{4},
\frac{1}{16} )
$,

\vskip 0.7ex
\hangindent=3em \hangafter=1
$S$ = ($ 1$,
$ 1$,
$ \sqrt{2}$,
$ -1$,
$ -1$,
$ -\sqrt{2}$;\ \ 
$ 1$,
$ -\sqrt{2}$,
$ -1$,
$ -1$,
$ \sqrt{2}$;\ \ 
$0$,
$ -\sqrt{2}$,
$ \sqrt{2}$,
$0$;\ \ 
$ -1$,
$ -1$,
$ -\sqrt{2}$;\ \ 
$ -1$,
$ \sqrt{2}$;\ \ 
$0$)

Factors = $2_{1,2.}^{4,625}\boxtimes 3_{\frac{13}{2},4.}^{16,330} $

Pseudo-unitary $\sim$  
$6_{\frac{11}{2},8.}^{16,548}$

\vskip 1ex 
\color{grey}

\noindent(19,8). $6_{\frac{11}{2},8.}^{16,111}$ \irep{39}:\ \ 
$d_i$ = ($1.0$,
$1.0$,
$1.414$,
$-1.0$,
$-1.0$,
$-1.414$) 

\vskip 0.7ex
\hangindent=3em \hangafter=1
$D^2= 8.0 = 
8$

\vskip 0.7ex
\hangindent=3em \hangafter=1
$T = ( 0,
\frac{1}{2},
\frac{13}{16},
\frac{1}{4},
\frac{3}{4},
\frac{9}{16} )
$,

\vskip 0.7ex
\hangindent=3em \hangafter=1
$S$ = ($ 1$,
$ 1$,
$ \sqrt{2}$,
$ -1$,
$ -1$,
$ -\sqrt{2}$;\ \ 
$ 1$,
$ -\sqrt{2}$,
$ -1$,
$ -1$,
$ \sqrt{2}$;\ \ 
$0$,
$ \sqrt{2}$,
$ -\sqrt{2}$,
$0$;\ \ 
$ -1$,
$ -1$,
$ \sqrt{2}$;\ \ 
$ -1$,
$ -\sqrt{2}$;\ \ 
$0$)

Factors = $2_{1,2.}^{4,625}\boxtimes 3_{\frac{9}{2},4.}^{16,343} $

Pseudo-unitary $\sim$  
$6_{\frac{15}{2},8.}^{16,107}$

\vskip 1ex 
\color{blue}

\noindent(20,1). $6_{3,6.}^{12,520}$ \irep{34}:\ \ 
$d_i$ = ($1.0$,
$1.0$,
$1.0$,
$-1.0$,
$-1.0$,
$-1.0$) 

\vskip 0.7ex
\hangindent=3em \hangafter=1
$D^2= 6.0 = 
6$

\vskip 0.7ex
\hangindent=3em \hangafter=1
$T = ( 0,
\frac{1}{3},
\frac{1}{3},
\frac{1}{4},
\frac{7}{12},
\frac{7}{12} )
$,

\vskip 0.7ex
\hangindent=3em \hangafter=1
$S$ = ($ 1$,
$ 1$,
$ 1$,
$ -1$,
$ -1$,
$ -1$;\ \ 
$ \zeta_{3}^{1}$,
$ -\zeta_{6}^{1}$,
$ -1$,
$ -\zeta_{3}^{1}$,
$ \zeta_{6}^{1}$;\ \ 
$ \zeta_{3}^{1}$,
$ -1$,
$ \zeta_{6}^{1}$,
$ -\zeta_{3}^{1}$;\ \ 
$ -1$,
$ -1$,
$ -1$;\ \ 
$ -\zeta_{3}^{1}$,
$ \zeta_{6}^{1}$;\ \ 
$ -\zeta_{3}^{1}$)

Factors = $2_{1,2.}^{4,625}\boxtimes 3_{2,3.}^{3,527} $

Pseudo-unitary $\sim$  
$6_{1,6.}^{12,701}$

\vskip 1ex 
\color{grey}

\noindent(20,2). $6_{1,6.}^{12,354}$ \irep{34}:\ \ 
$d_i$ = ($1.0$,
$1.0$,
$1.0$,
$-1.0$,
$-1.0$,
$-1.0$) 

\vskip 0.7ex
\hangindent=3em \hangafter=1
$D^2= 6.0 = 
6$

\vskip 0.7ex
\hangindent=3em \hangafter=1
$T = ( 0,
\frac{1}{3},
\frac{1}{3},
\frac{3}{4},
\frac{1}{12},
\frac{1}{12} )
$,

\vskip 0.7ex
\hangindent=3em \hangafter=1
$S$ = ($ 1$,
$ 1$,
$ 1$,
$ -1$,
$ -1$,
$ -1$;\ \ 
$ \zeta_{3}^{1}$,
$ -\zeta_{6}^{1}$,
$ -1$,
$ -\zeta_{3}^{1}$,
$ \zeta_{6}^{1}$;\ \ 
$ \zeta_{3}^{1}$,
$ -1$,
$ \zeta_{6}^{1}$,
$ -\zeta_{3}^{1}$;\ \ 
$ -1$,
$ -1$,
$ -1$;\ \ 
$ -\zeta_{3}^{1}$,
$ \zeta_{6}^{1}$;\ \ 
$ -\zeta_{3}^{1}$)

Factors = $2_{7,2.}^{4,562}\boxtimes 3_{2,3.}^{3,527} $

Pseudo-unitary $\sim$  
$6_{3,6.}^{12,534}$

\vskip 1ex 
\color{grey}

\noindent(20,3). $6_{7,6.}^{12,854}$ \irep{34}:\ \ 
$d_i$ = ($1.0$,
$1.0$,
$1.0$,
$-1.0$,
$-1.0$,
$-1.0$) 

\vskip 0.7ex
\hangindent=3em \hangafter=1
$D^2= 6.0 = 
6$

\vskip 0.7ex
\hangindent=3em \hangafter=1
$T = ( 0,
\frac{2}{3},
\frac{2}{3},
\frac{1}{4},
\frac{11}{12},
\frac{11}{12} )
$,

\vskip 0.7ex
\hangindent=3em \hangafter=1
$S$ = ($ 1$,
$ 1$,
$ 1$,
$ -1$,
$ -1$,
$ -1$;\ \ 
$ -\zeta_{6}^{1}$,
$ \zeta_{3}^{1}$,
$ -1$,
$ -\zeta_{3}^{1}$,
$ \zeta_{6}^{1}$;\ \ 
$ -\zeta_{6}^{1}$,
$ -1$,
$ \zeta_{6}^{1}$,
$ -\zeta_{3}^{1}$;\ \ 
$ -1$,
$ -1$,
$ -1$;\ \ 
$ \zeta_{6}^{1}$,
$ -\zeta_{3}^{1}$;\ \ 
$ \zeta_{6}^{1}$)

Factors = $2_{1,2.}^{4,625}\boxtimes 3_{6,3.}^{3,138} $

Pseudo-unitary $\sim$  
$6_{5,6.}^{12,298}$

\vskip 1ex 
\color{grey}

\noindent(20,4). $6_{5,6.}^{12,208}$ \irep{34}:\ \ 
$d_i$ = ($1.0$,
$1.0$,
$1.0$,
$-1.0$,
$-1.0$,
$-1.0$) 

\vskip 0.7ex
\hangindent=3em \hangafter=1
$D^2= 6.0 = 
6$

\vskip 0.7ex
\hangindent=3em \hangafter=1
$T = ( 0,
\frac{2}{3},
\frac{2}{3},
\frac{3}{4},
\frac{5}{12},
\frac{5}{12} )
$,

\vskip 0.7ex
\hangindent=3em \hangafter=1
$S$ = ($ 1$,
$ 1$,
$ 1$,
$ -1$,
$ -1$,
$ -1$;\ \ 
$ -\zeta_{6}^{1}$,
$ \zeta_{3}^{1}$,
$ -1$,
$ -\zeta_{3}^{1}$,
$ \zeta_{6}^{1}$;\ \ 
$ -\zeta_{6}^{1}$,
$ -1$,
$ \zeta_{6}^{1}$,
$ -\zeta_{3}^{1}$;\ \ 
$ -1$,
$ -1$,
$ -1$;\ \ 
$ \zeta_{6}^{1}$,
$ -\zeta_{3}^{1}$;\ \ 
$ \zeta_{6}^{1}$)

Factors = $2_{7,2.}^{4,562}\boxtimes 3_{6,3.}^{3,138} $

Pseudo-unitary $\sim$  
$6_{7,6.}^{12,113}$

\vskip 1ex 

}

\subsection{Rank 7}

{\small
\black

\noindent(1,1). $7_{2,7.}^{7,892}$ \irep{48}:\ \ 
$d_i$ = ($1.0$,
$1.0$,
$1.0$,
$1.0$,
$1.0$,
$1.0$,
$1.0$) 

\vskip 0.7ex
\hangindent=3em \hangafter=1
$D^2= 7.0 = 
7$

\vskip 0.7ex
\hangindent=3em \hangafter=1
$T = ( 0,
\frac{1}{7},
\frac{1}{7},
\frac{2}{7},
\frac{2}{7},
\frac{4}{7},
\frac{4}{7} )
$,

\vskip 0.7ex
\hangindent=3em \hangafter=1
$S$ = ($ 1$,
$ 1$,
$ 1$,
$ 1$,
$ 1$,
$ 1$,
$ 1$;\ \ 
$ -\zeta_{14}^{3}$,
$ \zeta_{7}^{2}$,
$ -\zeta_{14}^{5}$,
$ \zeta_{7}^{1}$,
$ -\zeta_{14}^{1}$,
$ \zeta_{7}^{3}$;\ \ 
$ -\zeta_{14}^{3}$,
$ \zeta_{7}^{1}$,
$ -\zeta_{14}^{5}$,
$ \zeta_{7}^{3}$,
$ -\zeta_{14}^{1}$;\ \ 
$ \zeta_{7}^{3}$,
$ -\zeta_{14}^{1}$,
$ \zeta_{7}^{2}$,
$ -\zeta_{14}^{3}$;\ \ 
$ \zeta_{7}^{3}$,
$ -\zeta_{14}^{3}$,
$ \zeta_{7}^{2}$;\ \ 
$ -\zeta_{14}^{5}$,
$ \zeta_{7}^{1}$;\ \ 
$ -\zeta_{14}^{5}$)

Prime. 

\vskip 1ex 
\color{grey}

\noindent(1,2). $7_{6,7.}^{7,110}$ \irep{48}:\ \ 
$d_i$ = ($1.0$,
$1.0$,
$1.0$,
$1.0$,
$1.0$,
$1.0$,
$1.0$) 

\vskip 0.7ex
\hangindent=3em \hangafter=1
$D^2= 7.0 = 
7$

\vskip 0.7ex
\hangindent=3em \hangafter=1
$T = ( 0,
\frac{3}{7},
\frac{3}{7},
\frac{5}{7},
\frac{5}{7},
\frac{6}{7},
\frac{6}{7} )
$,

\vskip 0.7ex
\hangindent=3em \hangafter=1
$S$ = ($ 1$,
$ 1$,
$ 1$,
$ 1$,
$ 1$,
$ 1$,
$ 1$;\ \ 
$ \zeta_{7}^{1}$,
$ -\zeta_{14}^{5}$,
$ -\zeta_{14}^{3}$,
$ \zeta_{7}^{2}$,
$ -\zeta_{14}^{1}$,
$ \zeta_{7}^{3}$;\ \ 
$ \zeta_{7}^{1}$,
$ \zeta_{7}^{2}$,
$ -\zeta_{14}^{3}$,
$ \zeta_{7}^{3}$,
$ -\zeta_{14}^{1}$;\ \ 
$ -\zeta_{14}^{1}$,
$ \zeta_{7}^{3}$,
$ -\zeta_{14}^{5}$,
$ \zeta_{7}^{1}$;\ \ 
$ -\zeta_{14}^{1}$,
$ \zeta_{7}^{1}$,
$ -\zeta_{14}^{5}$;\ \ 
$ \zeta_{7}^{2}$,
$ -\zeta_{14}^{3}$;\ \ 
$ \zeta_{7}^{2}$)

Prime. 

\vskip 1ex 
\black

\noindent(2,1). $7_{\frac{27}{4},27.31}^{32,396}$ \irep{96}:\ \ 
$d_i$ = ($1.0$,
$1.0$,
$1.847$,
$1.847$,
$2.414$,
$2.414$,
$2.613$) 

\vskip 0.7ex
\hangindent=3em \hangafter=1
$D^2= 27.313 = 
16+8\sqrt{2}$

\vskip 0.7ex
\hangindent=3em \hangafter=1
$T = ( 0,
\frac{1}{2},
\frac{1}{32},
\frac{1}{32},
\frac{1}{4},
\frac{3}{4},
\frac{21}{32} )
$,

\vskip 0.7ex
\hangindent=3em \hangafter=1
$S$ = ($ 1$,
$ 1$,
$ c_{16}^{1}$,
$ c_{16}^{1}$,
$ 1+\sqrt{2}$,
$ 1+\sqrt{2}$,
$ c^{1}_{16}
+c^{3}_{16}
$;\ \ 
$ 1$,
$ -c_{16}^{1}$,
$ -c_{16}^{1}$,
$ 1+\sqrt{2}$,
$ 1+\sqrt{2}$,
$ -c^{1}_{16}
-c^{3}_{16}
$;\ \ 
$(-c^{1}_{16}
-c^{3}_{16}
)\mathrm{i}$,
$(c^{1}_{16}
+c^{3}_{16}
)\mathrm{i}$,
$ -c_{16}^{1}$,
$ c_{16}^{1}$,
$0$;\ \ 
$(-c^{1}_{16}
-c^{3}_{16}
)\mathrm{i}$,
$ -c_{16}^{1}$,
$ c_{16}^{1}$,
$0$;\ \ 
$ -1$,
$ -1$,
$ c^{1}_{16}
+c^{3}_{16}
$;\ \ 
$ -1$,
$ -c^{1}_{16}
-c^{3}_{16}
$;\ \ 
$0$)

Prime. 

\vskip 1ex 
\color{grey}

\noindent(2,2). $7_{\frac{21}{4},27.31}^{32,114}$ \irep{96}:\ \ 
$d_i$ = ($1.0$,
$1.0$,
$1.847$,
$1.847$,
$2.414$,
$2.414$,
$2.613$) 

\vskip 0.7ex
\hangindent=3em \hangafter=1
$D^2= 27.313 = 
16+8\sqrt{2}$

\vskip 0.7ex
\hangindent=3em \hangafter=1
$T = ( 0,
\frac{1}{2},
\frac{15}{32},
\frac{15}{32},
\frac{1}{4},
\frac{3}{4},
\frac{27}{32} )
$,

\vskip 0.7ex
\hangindent=3em \hangafter=1
$S$ = ($ 1$,
$ 1$,
$ c_{16}^{1}$,
$ c_{16}^{1}$,
$ 1+\sqrt{2}$,
$ 1+\sqrt{2}$,
$ c^{1}_{16}
+c^{3}_{16}
$;\ \ 
$ 1$,
$ -c_{16}^{1}$,
$ -c_{16}^{1}$,
$ 1+\sqrt{2}$,
$ 1+\sqrt{2}$,
$ -c^{1}_{16}
-c^{3}_{16}
$;\ \ 
$(c^{1}_{16}
+c^{3}_{16}
)\mathrm{i}$,
$(-c^{1}_{16}
-c^{3}_{16}
)\mathrm{i}$,
$ c_{16}^{1}$,
$ -c_{16}^{1}$,
$0$;\ \ 
$(c^{1}_{16}
+c^{3}_{16}
)\mathrm{i}$,
$ c_{16}^{1}$,
$ -c_{16}^{1}$,
$0$;\ \ 
$ -1$,
$ -1$,
$ -c^{1}_{16}
-c^{3}_{16}
$;\ \ 
$ -1$,
$ c^{1}_{16}
+c^{3}_{16}
$;\ \ 
$0$)

Prime. 

\vskip 1ex 
\color{grey}

\noindent(2,3). $7_{\frac{11}{4},27.31}^{32,418}$ \irep{96}:\ \ 
$d_i$ = ($1.0$,
$1.0$,
$1.847$,
$1.847$,
$2.414$,
$2.414$,
$2.613$) 

\vskip 0.7ex
\hangindent=3em \hangafter=1
$D^2= 27.313 = 
16+8\sqrt{2}$

\vskip 0.7ex
\hangindent=3em \hangafter=1
$T = ( 0,
\frac{1}{2},
\frac{17}{32},
\frac{17}{32},
\frac{1}{4},
\frac{3}{4},
\frac{5}{32} )
$,

\vskip 0.7ex
\hangindent=3em \hangafter=1
$S$ = ($ 1$,
$ 1$,
$ c_{16}^{1}$,
$ c_{16}^{1}$,
$ 1+\sqrt{2}$,
$ 1+\sqrt{2}$,
$ c^{1}_{16}
+c^{3}_{16}
$;\ \ 
$ 1$,
$ -c_{16}^{1}$,
$ -c_{16}^{1}$,
$ 1+\sqrt{2}$,
$ 1+\sqrt{2}$,
$ -c^{1}_{16}
-c^{3}_{16}
$;\ \ 
$(-c^{1}_{16}
-c^{3}_{16}
)\mathrm{i}$,
$(c^{1}_{16}
+c^{3}_{16}
)\mathrm{i}$,
$ -c_{16}^{1}$,
$ c_{16}^{1}$,
$0$;\ \ 
$(-c^{1}_{16}
-c^{3}_{16}
)\mathrm{i}$,
$ -c_{16}^{1}$,
$ c_{16}^{1}$,
$0$;\ \ 
$ -1$,
$ -1$,
$ c^{1}_{16}
+c^{3}_{16}
$;\ \ 
$ -1$,
$ -c^{1}_{16}
-c^{3}_{16}
$;\ \ 
$0$)

Prime. 

\vskip 1ex 
\color{grey}

\noindent(2,4). $7_{\frac{5}{4},27.31}^{32,225}$ \irep{96}:\ \ 
$d_i$ = ($1.0$,
$1.0$,
$1.847$,
$1.847$,
$2.414$,
$2.414$,
$2.613$) 

\vskip 0.7ex
\hangindent=3em \hangafter=1
$D^2= 27.313 = 
16+8\sqrt{2}$

\vskip 0.7ex
\hangindent=3em \hangafter=1
$T = ( 0,
\frac{1}{2},
\frac{31}{32},
\frac{31}{32},
\frac{1}{4},
\frac{3}{4},
\frac{11}{32} )
$,

\vskip 0.7ex
\hangindent=3em \hangafter=1
$S$ = ($ 1$,
$ 1$,
$ c_{16}^{1}$,
$ c_{16}^{1}$,
$ 1+\sqrt{2}$,
$ 1+\sqrt{2}$,
$ c^{1}_{16}
+c^{3}_{16}
$;\ \ 
$ 1$,
$ -c_{16}^{1}$,
$ -c_{16}^{1}$,
$ 1+\sqrt{2}$,
$ 1+\sqrt{2}$,
$ -c^{1}_{16}
-c^{3}_{16}
$;\ \ 
$(c^{1}_{16}
+c^{3}_{16}
)\mathrm{i}$,
$(-c^{1}_{16}
-c^{3}_{16}
)\mathrm{i}$,
$ c_{16}^{1}$,
$ -c_{16}^{1}$,
$0$;\ \ 
$(c^{1}_{16}
+c^{3}_{16}
)\mathrm{i}$,
$ c_{16}^{1}$,
$ -c_{16}^{1}$,
$0$;\ \ 
$ -1$,
$ -1$,
$ -c^{1}_{16}
-c^{3}_{16}
$;\ \ 
$ -1$,
$ c^{1}_{16}
+c^{3}_{16}
$;\ \ 
$0$)

Prime. 

\vskip 1ex 
\color{grey}

\noindent(2,5). $7_{\frac{13}{4},27.31}^{32,530}$ \irep{96}:\ \ 
$d_i$ = ($1.0$,
$1.0$,
$2.414$,
$2.414$,
$-1.847$,
$-1.847$,
$-2.613$) 

\vskip 0.7ex
\hangindent=3em \hangafter=1
$D^2= 27.313 = 
16+8\sqrt{2}$

\vskip 0.7ex
\hangindent=3em \hangafter=1
$T = ( 0,
\frac{1}{2},
\frac{1}{4},
\frac{3}{4},
\frac{7}{32},
\frac{7}{32},
\frac{19}{32} )
$,

\vskip 0.7ex
\hangindent=3em \hangafter=1
$S$ = ($ 1$,
$ 1$,
$ 1+\sqrt{2}$,
$ 1+\sqrt{2}$,
$ -c_{16}^{1}$,
$ -c_{16}^{1}$,
$ -c^{1}_{16}
-c^{3}_{16}
$;\ \ 
$ 1$,
$ 1+\sqrt{2}$,
$ 1+\sqrt{2}$,
$ c_{16}^{1}$,
$ c_{16}^{1}$,
$ c^{1}_{16}
+c^{3}_{16}
$;\ \ 
$ -1$,
$ -1$,
$ -c_{16}^{1}$,
$ -c_{16}^{1}$,
$ c^{1}_{16}
+c^{3}_{16}
$;\ \ 
$ -1$,
$ c_{16}^{1}$,
$ c_{16}^{1}$,
$ -c^{1}_{16}
-c^{3}_{16}
$;\ \ 
$(-c^{1}_{16}
-c^{3}_{16}
)\mathrm{i}$,
$(c^{1}_{16}
+c^{3}_{16}
)\mathrm{i}$,
$0$;\ \ 
$(-c^{1}_{16}
-c^{3}_{16}
)\mathrm{i}$,
$0$;\ \ 
$0$)

Prime. 

Pseudo-unitary $\sim$  
$7_{\frac{29}{4},27.31}^{32,406}$

\vskip 1ex 
\color{grey}

\noindent(2,6). $7_{\frac{3}{4},27.31}^{32,809}$ \irep{96}:\ \ 
$d_i$ = ($1.0$,
$1.0$,
$2.414$,
$2.414$,
$-1.847$,
$-1.847$,
$-2.613$) 

\vskip 0.7ex
\hangindent=3em \hangafter=1
$D^2= 27.313 = 
16+8\sqrt{2}$

\vskip 0.7ex
\hangindent=3em \hangafter=1
$T = ( 0,
\frac{1}{2},
\frac{1}{4},
\frac{3}{4},
\frac{9}{32},
\frac{9}{32},
\frac{29}{32} )
$,

\vskip 0.7ex
\hangindent=3em \hangafter=1
$S$ = ($ 1$,
$ 1$,
$ 1+\sqrt{2}$,
$ 1+\sqrt{2}$,
$ -c_{16}^{1}$,
$ -c_{16}^{1}$,
$ -c^{1}_{16}
-c^{3}_{16}
$;\ \ 
$ 1$,
$ 1+\sqrt{2}$,
$ 1+\sqrt{2}$,
$ c_{16}^{1}$,
$ c_{16}^{1}$,
$ c^{1}_{16}
+c^{3}_{16}
$;\ \ 
$ -1$,
$ -1$,
$ c_{16}^{1}$,
$ c_{16}^{1}$,
$ -c^{1}_{16}
-c^{3}_{16}
$;\ \ 
$ -1$,
$ -c_{16}^{1}$,
$ -c_{16}^{1}$,
$ c^{1}_{16}
+c^{3}_{16}
$;\ \ 
$(c^{1}_{16}
+c^{3}_{16}
)\mathrm{i}$,
$(-c^{1}_{16}
-c^{3}_{16}
)\mathrm{i}$,
$0$;\ \ 
$(c^{1}_{16}
+c^{3}_{16}
)\mathrm{i}$,
$0$;\ \ 
$0$)

Prime. 

Pseudo-unitary $\sim$  
$7_{\frac{19}{4},27.31}^{32,116}$

\vskip 1ex 
\color{grey}

\noindent(2,7). $7_{\frac{29}{4},27.31}^{32,303}$ \irep{96}:\ \ 
$d_i$ = ($1.0$,
$1.0$,
$2.414$,
$2.414$,
$-1.847$,
$-1.847$,
$-2.613$) 

\vskip 0.7ex
\hangindent=3em \hangafter=1
$D^2= 27.313 = 
16+8\sqrt{2}$

\vskip 0.7ex
\hangindent=3em \hangafter=1
$T = ( 0,
\frac{1}{2},
\frac{1}{4},
\frac{3}{4},
\frac{23}{32},
\frac{23}{32},
\frac{3}{32} )
$,

\vskip 0.7ex
\hangindent=3em \hangafter=1
$S$ = ($ 1$,
$ 1$,
$ 1+\sqrt{2}$,
$ 1+\sqrt{2}$,
$ -c_{16}^{1}$,
$ -c_{16}^{1}$,
$ -c^{1}_{16}
-c^{3}_{16}
$;\ \ 
$ 1$,
$ 1+\sqrt{2}$,
$ 1+\sqrt{2}$,
$ c_{16}^{1}$,
$ c_{16}^{1}$,
$ c^{1}_{16}
+c^{3}_{16}
$;\ \ 
$ -1$,
$ -1$,
$ -c_{16}^{1}$,
$ -c_{16}^{1}$,
$ c^{1}_{16}
+c^{3}_{16}
$;\ \ 
$ -1$,
$ c_{16}^{1}$,
$ c_{16}^{1}$,
$ -c^{1}_{16}
-c^{3}_{16}
$;\ \ 
$(-c^{1}_{16}
-c^{3}_{16}
)\mathrm{i}$,
$(c^{1}_{16}
+c^{3}_{16}
)\mathrm{i}$,
$0$;\ \ 
$(-c^{1}_{16}
-c^{3}_{16}
)\mathrm{i}$,
$0$;\ \ 
$0$)

Prime. 

Pseudo-unitary $\sim$  
$7_{\frac{13}{4},27.31}^{32,427}$

\vskip 1ex 
\color{grey}

\noindent(2,8). $7_{\frac{19}{4},27.31}^{32,105}$ \irep{96}:\ \ 
$d_i$ = ($1.0$,
$1.0$,
$2.414$,
$2.414$,
$-1.847$,
$-1.847$,
$-2.613$) 

\vskip 0.7ex
\hangindent=3em \hangafter=1
$D^2= 27.313 = 
16+8\sqrt{2}$

\vskip 0.7ex
\hangindent=3em \hangafter=1
$T = ( 0,
\frac{1}{2},
\frac{1}{4},
\frac{3}{4},
\frac{25}{32},
\frac{25}{32},
\frac{13}{32} )
$,

\vskip 0.7ex
\hangindent=3em \hangafter=1
$S$ = ($ 1$,
$ 1$,
$ 1+\sqrt{2}$,
$ 1+\sqrt{2}$,
$ -c_{16}^{1}$,
$ -c_{16}^{1}$,
$ -c^{1}_{16}
-c^{3}_{16}
$;\ \ 
$ 1$,
$ 1+\sqrt{2}$,
$ 1+\sqrt{2}$,
$ c_{16}^{1}$,
$ c_{16}^{1}$,
$ c^{1}_{16}
+c^{3}_{16}
$;\ \ 
$ -1$,
$ -1$,
$ c_{16}^{1}$,
$ c_{16}^{1}$,
$ -c^{1}_{16}
-c^{3}_{16}
$;\ \ 
$ -1$,
$ -c_{16}^{1}$,
$ -c_{16}^{1}$,
$ c^{1}_{16}
+c^{3}_{16}
$;\ \ 
$(c^{1}_{16}
+c^{3}_{16}
)\mathrm{i}$,
$(-c^{1}_{16}
-c^{3}_{16}
)\mathrm{i}$,
$0$;\ \ 
$(c^{1}_{16}
+c^{3}_{16}
)\mathrm{i}$,
$0$;\ \ 
$0$)

Prime. 

Pseudo-unitary $\sim$  
$7_{\frac{3}{4},27.31}^{32,913}$

\vskip 1ex 
\color{grey}

\noindent(2,9). $7_{\frac{1}{4},4.686}^{32,141}$ \irep{96}:\ \ 
$d_i$ = ($1.0$,
$0.765$,
$0.765$,
$1.0$,
$-0.414$,
$-0.414$,
$-1.82$) 

\vskip 0.7ex
\hangindent=3em \hangafter=1
$D^2= 4.686 = 
16-8\sqrt{2}$

\vskip 0.7ex
\hangindent=3em \hangafter=1
$T = ( 0,
\frac{3}{32},
\frac{3}{32},
\frac{1}{2},
\frac{1}{4},
\frac{3}{4},
\frac{31}{32} )
$,

\vskip 0.7ex
\hangindent=3em \hangafter=1
$S$ = ($ 1$,
$ c_{16}^{3}$,
$ c_{16}^{3}$,
$ 1$,
$ 1-\sqrt{2}$,
$ 1-\sqrt{2}$,
$ -c^{1}_{16}
+c^{3}_{16}
$;\ \ 
$(-c^{1}_{16}
+c^{3}_{16}
)\mathrm{i}$,
$(c^{1}_{16}
-c^{3}_{16}
)\mathrm{i}$,
$ -c_{16}^{3}$,
$ c_{16}^{3}$,
$ -c_{16}^{3}$,
$0$;\ \ 
$(-c^{1}_{16}
+c^{3}_{16}
)\mathrm{i}$,
$ -c_{16}^{3}$,
$ c_{16}^{3}$,
$ -c_{16}^{3}$,
$0$;\ \ 
$ 1$,
$ 1-\sqrt{2}$,
$ 1-\sqrt{2}$,
$ c^{1}_{16}
-c^{3}_{16}
$;\ \ 
$ -1$,
$ -1$,
$ c^{1}_{16}
-c^{3}_{16}
$;\ \ 
$ -1$,
$ -c^{1}_{16}
+c^{3}_{16}
$;\ \ 
$0$)

Prime. 

Not pseudo-unitary. 

\vskip 1ex 
\color{grey}

\noindent(2,10). $7_{\frac{15}{4},4.686}^{32,466}$ \irep{96}:\ \ 
$d_i$ = ($1.0$,
$0.765$,
$0.765$,
$1.0$,
$-0.414$,
$-0.414$,
$-1.82$) 

\vskip 0.7ex
\hangindent=3em \hangafter=1
$D^2= 4.686 = 
16-8\sqrt{2}$

\vskip 0.7ex
\hangindent=3em \hangafter=1
$T = ( 0,
\frac{13}{32},
\frac{13}{32},
\frac{1}{2},
\frac{1}{4},
\frac{3}{4},
\frac{17}{32} )
$,

\vskip 0.7ex
\hangindent=3em \hangafter=1
$S$ = ($ 1$,
$ c_{16}^{3}$,
$ c_{16}^{3}$,
$ 1$,
$ 1-\sqrt{2}$,
$ 1-\sqrt{2}$,
$ -c^{1}_{16}
+c^{3}_{16}
$;\ \ 
$(c^{1}_{16}
-c^{3}_{16}
)\mathrm{i}$,
$(-c^{1}_{16}
+c^{3}_{16}
)\mathrm{i}$,
$ -c_{16}^{3}$,
$ -c_{16}^{3}$,
$ c_{16}^{3}$,
$0$;\ \ 
$(c^{1}_{16}
-c^{3}_{16}
)\mathrm{i}$,
$ -c_{16}^{3}$,
$ -c_{16}^{3}$,
$ c_{16}^{3}$,
$0$;\ \ 
$ 1$,
$ 1-\sqrt{2}$,
$ 1-\sqrt{2}$,
$ c^{1}_{16}
-c^{3}_{16}
$;\ \ 
$ -1$,
$ -1$,
$ -c^{1}_{16}
+c^{3}_{16}
$;\ \ 
$ -1$,
$ c^{1}_{16}
-c^{3}_{16}
$;\ \ 
$0$)

Prime. 

Not pseudo-unitary. 

\vskip 1ex 
\color{grey}

\noindent(2,11). $7_{\frac{9}{4},4.686}^{32,318}$ \irep{96}:\ \ 
$d_i$ = ($1.0$,
$1.0$,
$1.82$,
$-0.414$,
$-0.414$,
$-0.765$,
$-0.765$) 

\vskip 0.7ex
\hangindent=3em \hangafter=1
$D^2= 4.686 = 
16-8\sqrt{2}$

\vskip 0.7ex
\hangindent=3em \hangafter=1
$T = ( 0,
\frac{1}{2},
\frac{7}{32},
\frac{1}{4},
\frac{3}{4},
\frac{11}{32},
\frac{11}{32} )
$,

\vskip 0.7ex
\hangindent=3em \hangafter=1
$S$ = ($ 1$,
$ 1$,
$ c^{1}_{16}
-c^{3}_{16}
$,
$ 1-\sqrt{2}$,
$ 1-\sqrt{2}$,
$ -c_{16}^{3}$,
$ -c_{16}^{3}$;\ \ 
$ 1$,
$ -c^{1}_{16}
+c^{3}_{16}
$,
$ 1-\sqrt{2}$,
$ 1-\sqrt{2}$,
$ c_{16}^{3}$,
$ c_{16}^{3}$;\ \ 
$0$,
$ -c^{1}_{16}
+c^{3}_{16}
$,
$ c^{1}_{16}
-c^{3}_{16}
$,
$0$,
$0$;\ \ 
$ -1$,
$ -1$,
$ -c_{16}^{3}$,
$ -c_{16}^{3}$;\ \ 
$ -1$,
$ c_{16}^{3}$,
$ c_{16}^{3}$;\ \ 
$(c^{1}_{16}
-c^{3}_{16}
)\mathrm{i}$,
$(-c^{1}_{16}
+c^{3}_{16}
)\mathrm{i}$;\ \ 
$(c^{1}_{16}
-c^{3}_{16}
)\mathrm{i}$)

Prime. 

Not pseudo-unitary. 

\vskip 1ex 
\color{grey}

\noindent(2,12). $7_{\frac{7}{4},4.686}^{32,141}$ \irep{96}:\ \ 
$d_i$ = ($1.0$,
$1.0$,
$1.82$,
$-0.414$,
$-0.414$,
$-0.765$,
$-0.765$) 

\vskip 0.7ex
\hangindent=3em \hangafter=1
$D^2= 4.686 = 
16-8\sqrt{2}$

\vskip 0.7ex
\hangindent=3em \hangafter=1
$T = ( 0,
\frac{1}{2},
\frac{9}{32},
\frac{1}{4},
\frac{3}{4},
\frac{5}{32},
\frac{5}{32} )
$,

\vskip 0.7ex
\hangindent=3em \hangafter=1
$S$ = ($ 1$,
$ 1$,
$ c^{1}_{16}
-c^{3}_{16}
$,
$ 1-\sqrt{2}$,
$ 1-\sqrt{2}$,
$ -c_{16}^{3}$,
$ -c_{16}^{3}$;\ \ 
$ 1$,
$ -c^{1}_{16}
+c^{3}_{16}
$,
$ 1-\sqrt{2}$,
$ 1-\sqrt{2}$,
$ c_{16}^{3}$,
$ c_{16}^{3}$;\ \ 
$0$,
$ c^{1}_{16}
-c^{3}_{16}
$,
$ -c^{1}_{16}
+c^{3}_{16}
$,
$0$,
$0$;\ \ 
$ -1$,
$ -1$,
$ c_{16}^{3}$,
$ c_{16}^{3}$;\ \ 
$ -1$,
$ -c_{16}^{3}$,
$ -c_{16}^{3}$;\ \ 
$(-c^{1}_{16}
+c^{3}_{16}
)\mathrm{i}$,
$(c^{1}_{16}
-c^{3}_{16}
)\mathrm{i}$;\ \ 
$(-c^{1}_{16}
+c^{3}_{16}
)\mathrm{i}$)

Prime. 

Not pseudo-unitary. 

\vskip 1ex 
\color{grey}

\noindent(2,13). $7_{\frac{25}{4},4.686}^{32,720}$ \irep{96}:\ \ 
$d_i$ = ($1.0$,
$1.0$,
$1.82$,
$-0.414$,
$-0.414$,
$-0.765$,
$-0.765$) 

\vskip 0.7ex
\hangindent=3em \hangafter=1
$D^2= 4.686 = 
16-8\sqrt{2}$

\vskip 0.7ex
\hangindent=3em \hangafter=1
$T = ( 0,
\frac{1}{2},
\frac{23}{32},
\frac{1}{4},
\frac{3}{4},
\frac{27}{32},
\frac{27}{32} )
$,

\vskip 0.7ex
\hangindent=3em \hangafter=1
$S$ = ($ 1$,
$ 1$,
$ c^{1}_{16}
-c^{3}_{16}
$,
$ 1-\sqrt{2}$,
$ 1-\sqrt{2}$,
$ -c_{16}^{3}$,
$ -c_{16}^{3}$;\ \ 
$ 1$,
$ -c^{1}_{16}
+c^{3}_{16}
$,
$ 1-\sqrt{2}$,
$ 1-\sqrt{2}$,
$ c_{16}^{3}$,
$ c_{16}^{3}$;\ \ 
$0$,
$ -c^{1}_{16}
+c^{3}_{16}
$,
$ c^{1}_{16}
-c^{3}_{16}
$,
$0$,
$0$;\ \ 
$ -1$,
$ -1$,
$ -c_{16}^{3}$,
$ -c_{16}^{3}$;\ \ 
$ -1$,
$ c_{16}^{3}$,
$ c_{16}^{3}$;\ \ 
$(c^{1}_{16}
-c^{3}_{16}
)\mathrm{i}$,
$(-c^{1}_{16}
+c^{3}_{16}
)\mathrm{i}$;\ \ 
$(c^{1}_{16}
-c^{3}_{16}
)\mathrm{i}$)

Prime. 

Not pseudo-unitary. 

\vskip 1ex 
\color{grey}

\noindent(2,14). $7_{\frac{23}{4},4.686}^{32,188}$ \irep{96}:\ \ 
$d_i$ = ($1.0$,
$1.0$,
$1.82$,
$-0.414$,
$-0.414$,
$-0.765$,
$-0.765$) 

\vskip 0.7ex
\hangindent=3em \hangafter=1
$D^2= 4.686 = 
16-8\sqrt{2}$

\vskip 0.7ex
\hangindent=3em \hangafter=1
$T = ( 0,
\frac{1}{2},
\frac{25}{32},
\frac{1}{4},
\frac{3}{4},
\frac{21}{32},
\frac{21}{32} )
$,

\vskip 0.7ex
\hangindent=3em \hangafter=1
$S$ = ($ 1$,
$ 1$,
$ c^{1}_{16}
-c^{3}_{16}
$,
$ 1-\sqrt{2}$,
$ 1-\sqrt{2}$,
$ -c_{16}^{3}$,
$ -c_{16}^{3}$;\ \ 
$ 1$,
$ -c^{1}_{16}
+c^{3}_{16}
$,
$ 1-\sqrt{2}$,
$ 1-\sqrt{2}$,
$ c_{16}^{3}$,
$ c_{16}^{3}$;\ \ 
$0$,
$ c^{1}_{16}
-c^{3}_{16}
$,
$ -c^{1}_{16}
+c^{3}_{16}
$,
$0$,
$0$;\ \ 
$ -1$,
$ -1$,
$ c_{16}^{3}$,
$ c_{16}^{3}$;\ \ 
$ -1$,
$ -c_{16}^{3}$,
$ -c_{16}^{3}$;\ \ 
$(-c^{1}_{16}
+c^{3}_{16}
)\mathrm{i}$,
$(c^{1}_{16}
-c^{3}_{16}
)\mathrm{i}$;\ \ 
$(-c^{1}_{16}
+c^{3}_{16}
)\mathrm{i}$)

Prime. 

Not pseudo-unitary. 

\vskip 1ex 
\color{grey}

\noindent(2,15). $7_{\frac{17}{4},4.686}^{32,112}$ \irep{96}:\ \ 
$d_i$ = ($1.0$,
$0.765$,
$0.765$,
$1.0$,
$-0.414$,
$-0.414$,
$-1.82$) 

\vskip 0.7ex
\hangindent=3em \hangafter=1
$D^2= 4.686 = 
16-8\sqrt{2}$

\vskip 0.7ex
\hangindent=3em \hangafter=1
$T = ( 0,
\frac{19}{32},
\frac{19}{32},
\frac{1}{2},
\frac{1}{4},
\frac{3}{4},
\frac{15}{32} )
$,

\vskip 0.7ex
\hangindent=3em \hangafter=1
$S$ = ($ 1$,
$ c_{16}^{3}$,
$ c_{16}^{3}$,
$ 1$,
$ 1-\sqrt{2}$,
$ 1-\sqrt{2}$,
$ -c^{1}_{16}
+c^{3}_{16}
$;\ \ 
$(-c^{1}_{16}
+c^{3}_{16}
)\mathrm{i}$,
$(c^{1}_{16}
-c^{3}_{16}
)\mathrm{i}$,
$ -c_{16}^{3}$,
$ c_{16}^{3}$,
$ -c_{16}^{3}$,
$0$;\ \ 
$(-c^{1}_{16}
+c^{3}_{16}
)\mathrm{i}$,
$ -c_{16}^{3}$,
$ c_{16}^{3}$,
$ -c_{16}^{3}$,
$0$;\ \ 
$ 1$,
$ 1-\sqrt{2}$,
$ 1-\sqrt{2}$,
$ c^{1}_{16}
-c^{3}_{16}
$;\ \ 
$ -1$,
$ -1$,
$ c^{1}_{16}
-c^{3}_{16}
$;\ \ 
$ -1$,
$ -c^{1}_{16}
+c^{3}_{16}
$;\ \ 
$0$)

Prime. 

Not pseudo-unitary. 

\vskip 1ex 
\color{grey}

\noindent(2,16). $7_{\frac{31}{4},4.686}^{32,842}$ \irep{96}:\ \ 
$d_i$ = ($1.0$,
$0.765$,
$0.765$,
$1.0$,
$-0.414$,
$-0.414$,
$-1.82$) 

\vskip 0.7ex
\hangindent=3em \hangafter=1
$D^2= 4.686 = 
16-8\sqrt{2}$

\vskip 0.7ex
\hangindent=3em \hangafter=1
$T = ( 0,
\frac{29}{32},
\frac{29}{32},
\frac{1}{2},
\frac{1}{4},
\frac{3}{4},
\frac{1}{32} )
$,

\vskip 0.7ex
\hangindent=3em \hangafter=1
$S$ = ($ 1$,
$ c_{16}^{3}$,
$ c_{16}^{3}$,
$ 1$,
$ 1-\sqrt{2}$,
$ 1-\sqrt{2}$,
$ -c^{1}_{16}
+c^{3}_{16}
$;\ \ 
$(c^{1}_{16}
-c^{3}_{16}
)\mathrm{i}$,
$(-c^{1}_{16}
+c^{3}_{16}
)\mathrm{i}$,
$ -c_{16}^{3}$,
$ -c_{16}^{3}$,
$ c_{16}^{3}$,
$0$;\ \ 
$(c^{1}_{16}
-c^{3}_{16}
)\mathrm{i}$,
$ -c_{16}^{3}$,
$ -c_{16}^{3}$,
$ c_{16}^{3}$,
$0$;\ \ 
$ 1$,
$ 1-\sqrt{2}$,
$ 1-\sqrt{2}$,
$ c^{1}_{16}
-c^{3}_{16}
$;\ \ 
$ -1$,
$ -1$,
$ -c^{1}_{16}
+c^{3}_{16}
$;\ \ 
$ -1$,
$ c^{1}_{16}
-c^{3}_{16}
$;\ \ 
$0$)

Prime. 

Not pseudo-unitary. 

\vskip 1ex 
\black

\noindent(3,1). $7_{\frac{9}{4},27.31}^{32,918}$ \irep{97}:\ \ 
$d_i$ = ($1.0$,
$1.0$,
$1.847$,
$1.847$,
$2.414$,
$2.414$,
$2.613$) 

\vskip 0.7ex
\hangindent=3em \hangafter=1
$D^2= 27.313 = 
16+8\sqrt{2}$

\vskip 0.7ex
\hangindent=3em \hangafter=1
$T = ( 0,
\frac{1}{2},
\frac{3}{32},
\frac{3}{32},
\frac{1}{4},
\frac{3}{4},
\frac{15}{32} )
$,

\vskip 0.7ex
\hangindent=3em \hangafter=1
$S$ = ($ 1$,
$ 1$,
$ c_{16}^{1}$,
$ c_{16}^{1}$,
$ 1+\sqrt{2}$,
$ 1+\sqrt{2}$,
$ c^{1}_{16}
+c^{3}_{16}
$;\ \ 
$ 1$,
$ -c_{16}^{1}$,
$ -c_{16}^{1}$,
$ 1+\sqrt{2}$,
$ 1+\sqrt{2}$,
$ -c^{1}_{16}
-c^{3}_{16}
$;\ \ 
$ c^{1}_{16}
+c^{3}_{16}
$,
$ -c^{1}_{16}
-c^{3}_{16}
$,
$ c_{16}^{1}$,
$ -c_{16}^{1}$,
$0$;\ \ 
$ c^{1}_{16}
+c^{3}_{16}
$,
$ c_{16}^{1}$,
$ -c_{16}^{1}$,
$0$;\ \ 
$ -1$,
$ -1$,
$ -c^{1}_{16}
-c^{3}_{16}
$;\ \ 
$ -1$,
$ c^{1}_{16}
+c^{3}_{16}
$;\ \ 
$0$)

Prime. 

\vskip 1ex 
\color{grey}

\noindent(3,2). $7_{\frac{7}{4},27.31}^{32,912}$ \irep{97}:\ \ 
$d_i$ = ($1.0$,
$1.0$,
$1.847$,
$1.847$,
$2.414$,
$2.414$,
$2.613$) 

\vskip 0.7ex
\hangindent=3em \hangafter=1
$D^2= 27.313 = 
16+8\sqrt{2}$

\vskip 0.7ex
\hangindent=3em \hangafter=1
$T = ( 0,
\frac{1}{2},
\frac{13}{32},
\frac{13}{32},
\frac{1}{4},
\frac{3}{4},
\frac{1}{32} )
$,

\vskip 0.7ex
\hangindent=3em \hangafter=1
$S$ = ($ 1$,
$ 1$,
$ c_{16}^{1}$,
$ c_{16}^{1}$,
$ 1+\sqrt{2}$,
$ 1+\sqrt{2}$,
$ c^{1}_{16}
+c^{3}_{16}
$;\ \ 
$ 1$,
$ -c_{16}^{1}$,
$ -c_{16}^{1}$,
$ 1+\sqrt{2}$,
$ 1+\sqrt{2}$,
$ -c^{1}_{16}
-c^{3}_{16}
$;\ \ 
$ c^{1}_{16}
+c^{3}_{16}
$,
$ -c^{1}_{16}
-c^{3}_{16}
$,
$ -c_{16}^{1}$,
$ c_{16}^{1}$,
$0$;\ \ 
$ c^{1}_{16}
+c^{3}_{16}
$,
$ -c_{16}^{1}$,
$ c_{16}^{1}$,
$0$;\ \ 
$ -1$,
$ -1$,
$ c^{1}_{16}
+c^{3}_{16}
$;\ \ 
$ -1$,
$ -c^{1}_{16}
-c^{3}_{16}
$;\ \ 
$0$)

Prime. 

\vskip 1ex 
\color{grey}

\noindent(3,3). $7_{\frac{25}{4},27.31}^{32,222}$ \irep{97}:\ \ 
$d_i$ = ($1.0$,
$1.0$,
$1.847$,
$1.847$,
$2.414$,
$2.414$,
$2.613$) 

\vskip 0.7ex
\hangindent=3em \hangafter=1
$D^2= 27.313 = 
16+8\sqrt{2}$

\vskip 0.7ex
\hangindent=3em \hangafter=1
$T = ( 0,
\frac{1}{2},
\frac{19}{32},
\frac{19}{32},
\frac{1}{4},
\frac{3}{4},
\frac{31}{32} )
$,

\vskip 0.7ex
\hangindent=3em \hangafter=1
$S$ = ($ 1$,
$ 1$,
$ c_{16}^{1}$,
$ c_{16}^{1}$,
$ 1+\sqrt{2}$,
$ 1+\sqrt{2}$,
$ c^{1}_{16}
+c^{3}_{16}
$;\ \ 
$ 1$,
$ -c_{16}^{1}$,
$ -c_{16}^{1}$,
$ 1+\sqrt{2}$,
$ 1+\sqrt{2}$,
$ -c^{1}_{16}
-c^{3}_{16}
$;\ \ 
$ c^{1}_{16}
+c^{3}_{16}
$,
$ -c^{1}_{16}
-c^{3}_{16}
$,
$ c_{16}^{1}$,
$ -c_{16}^{1}$,
$0$;\ \ 
$ c^{1}_{16}
+c^{3}_{16}
$,
$ c_{16}^{1}$,
$ -c_{16}^{1}$,
$0$;\ \ 
$ -1$,
$ -1$,
$ -c^{1}_{16}
-c^{3}_{16}
$;\ \ 
$ -1$,
$ c^{1}_{16}
+c^{3}_{16}
$;\ \ 
$0$)

Prime. 

\vskip 1ex 
\color{grey}

\noindent(3,4). $7_{\frac{23}{4},27.31}^{32,224}$ \irep{97}:\ \ 
$d_i$ = ($1.0$,
$1.0$,
$1.847$,
$1.847$,
$2.414$,
$2.414$,
$2.613$) 

\vskip 0.7ex
\hangindent=3em \hangafter=1
$D^2= 27.313 = 
16+8\sqrt{2}$

\vskip 0.7ex
\hangindent=3em \hangafter=1
$T = ( 0,
\frac{1}{2},
\frac{29}{32},
\frac{29}{32},
\frac{1}{4},
\frac{3}{4},
\frac{17}{32} )
$,

\vskip 0.7ex
\hangindent=3em \hangafter=1
$S$ = ($ 1$,
$ 1$,
$ c_{16}^{1}$,
$ c_{16}^{1}$,
$ 1+\sqrt{2}$,
$ 1+\sqrt{2}$,
$ c^{1}_{16}
+c^{3}_{16}
$;\ \ 
$ 1$,
$ -c_{16}^{1}$,
$ -c_{16}^{1}$,
$ 1+\sqrt{2}$,
$ 1+\sqrt{2}$,
$ -c^{1}_{16}
-c^{3}_{16}
$;\ \ 
$ c^{1}_{16}
+c^{3}_{16}
$,
$ -c^{1}_{16}
-c^{3}_{16}
$,
$ -c_{16}^{1}$,
$ c_{16}^{1}$,
$0$;\ \ 
$ c^{1}_{16}
+c^{3}_{16}
$,
$ -c_{16}^{1}$,
$ c_{16}^{1}$,
$0$;\ \ 
$ -1$,
$ -1$,
$ c^{1}_{16}
+c^{3}_{16}
$;\ \ 
$ -1$,
$ -c^{1}_{16}
-c^{3}_{16}
$;\ \ 
$0$)

Prime. 

\vskip 1ex 
\color{grey}

\noindent(3,5). $7_{\frac{31}{4},27.31}^{32,562}$ \irep{97}:\ \ 
$d_i$ = ($1.0$,
$1.0$,
$2.414$,
$2.414$,
$-1.847$,
$-1.847$,
$-2.613$) 

\vskip 0.7ex
\hangindent=3em \hangafter=1
$D^2= 27.313 = 
16+8\sqrt{2}$

\vskip 0.7ex
\hangindent=3em \hangafter=1
$T = ( 0,
\frac{1}{2},
\frac{1}{4},
\frac{3}{4},
\frac{5}{32},
\frac{5}{32},
\frac{25}{32} )
$,

\vskip 0.7ex
\hangindent=3em \hangafter=1
$S$ = ($ 1$,
$ 1$,
$ 1+\sqrt{2}$,
$ 1+\sqrt{2}$,
$ -c_{16}^{1}$,
$ -c_{16}^{1}$,
$ -c^{1}_{16}
-c^{3}_{16}
$;\ \ 
$ 1$,
$ 1+\sqrt{2}$,
$ 1+\sqrt{2}$,
$ c_{16}^{1}$,
$ c_{16}^{1}$,
$ c^{1}_{16}
+c^{3}_{16}
$;\ \ 
$ -1$,
$ -1$,
$ c_{16}^{1}$,
$ c_{16}^{1}$,
$ -c^{1}_{16}
-c^{3}_{16}
$;\ \ 
$ -1$,
$ -c_{16}^{1}$,
$ -c_{16}^{1}$,
$ c^{1}_{16}
+c^{3}_{16}
$;\ \ 
$ -c^{1}_{16}
-c^{3}_{16}
$,
$ c^{1}_{16}
+c^{3}_{16}
$,
$0$;\ \ 
$ -c^{1}_{16}
-c^{3}_{16}
$,
$0$;\ \ 
$0$)

Prime. 

Pseudo-unitary $\sim$  
$7_{\frac{15}{4},27.31}^{32,272}$

\vskip 1ex 
\color{grey}

\noindent(3,6). $7_{\frac{17}{4},27.31}^{32,157}$ \irep{97}:\ \ 
$d_i$ = ($1.0$,
$1.0$,
$2.414$,
$2.414$,
$-1.847$,
$-1.847$,
$-2.613$) 

\vskip 0.7ex
\hangindent=3em \hangafter=1
$D^2= 27.313 = 
16+8\sqrt{2}$

\vskip 0.7ex
\hangindent=3em \hangafter=1
$T = ( 0,
\frac{1}{2},
\frac{1}{4},
\frac{3}{4},
\frac{11}{32},
\frac{11}{32},
\frac{23}{32} )
$,

\vskip 0.7ex
\hangindent=3em \hangafter=1
$S$ = ($ 1$,
$ 1$,
$ 1+\sqrt{2}$,
$ 1+\sqrt{2}$,
$ -c_{16}^{1}$,
$ -c_{16}^{1}$,
$ -c^{1}_{16}
-c^{3}_{16}
$;\ \ 
$ 1$,
$ 1+\sqrt{2}$,
$ 1+\sqrt{2}$,
$ c_{16}^{1}$,
$ c_{16}^{1}$,
$ c^{1}_{16}
+c^{3}_{16}
$;\ \ 
$ -1$,
$ -1$,
$ -c_{16}^{1}$,
$ -c_{16}^{1}$,
$ c^{1}_{16}
+c^{3}_{16}
$;\ \ 
$ -1$,
$ c_{16}^{1}$,
$ c_{16}^{1}$,
$ -c^{1}_{16}
-c^{3}_{16}
$;\ \ 
$ -c^{1}_{16}
-c^{3}_{16}
$,
$ c^{1}_{16}
+c^{3}_{16}
$,
$0$;\ \ 
$ -c^{1}_{16}
-c^{3}_{16}
$,
$0$;\ \ 
$0$)

Prime. 

Pseudo-unitary $\sim$  
$7_{\frac{1}{4},27.31}^{32,123}$

\vskip 1ex 
\color{grey}

\noindent(3,7). $7_{\frac{15}{4},27.31}^{32,169}$ \irep{97}:\ \ 
$d_i$ = ($1.0$,
$1.0$,
$2.414$,
$2.414$,
$-1.847$,
$-1.847$,
$-2.613$) 

\vskip 0.7ex
\hangindent=3em \hangafter=1
$D^2= 27.313 = 
16+8\sqrt{2}$

\vskip 0.7ex
\hangindent=3em \hangafter=1
$T = ( 0,
\frac{1}{2},
\frac{1}{4},
\frac{3}{4},
\frac{21}{32},
\frac{21}{32},
\frac{9}{32} )
$,

\vskip 0.7ex
\hangindent=3em \hangafter=1
$S$ = ($ 1$,
$ 1$,
$ 1+\sqrt{2}$,
$ 1+\sqrt{2}$,
$ -c_{16}^{1}$,
$ -c_{16}^{1}$,
$ -c^{1}_{16}
-c^{3}_{16}
$;\ \ 
$ 1$,
$ 1+\sqrt{2}$,
$ 1+\sqrt{2}$,
$ c_{16}^{1}$,
$ c_{16}^{1}$,
$ c^{1}_{16}
+c^{3}_{16}
$;\ \ 
$ -1$,
$ -1$,
$ c_{16}^{1}$,
$ c_{16}^{1}$,
$ -c^{1}_{16}
-c^{3}_{16}
$;\ \ 
$ -1$,
$ -c_{16}^{1}$,
$ -c_{16}^{1}$,
$ c^{1}_{16}
+c^{3}_{16}
$;\ \ 
$ -c^{1}_{16}
-c^{3}_{16}
$,
$ c^{1}_{16}
+c^{3}_{16}
$,
$0$;\ \ 
$ -c^{1}_{16}
-c^{3}_{16}
$,
$0$;\ \ 
$0$)

Prime. 

Pseudo-unitary $\sim$  
$7_{\frac{31}{4},27.31}^{32,159}$

\vskip 1ex 
\color{grey}

\noindent(3,8). $7_{\frac{1}{4},27.31}^{32,112}$ \irep{97}:\ \ 
$d_i$ = ($1.0$,
$1.0$,
$2.414$,
$2.414$,
$-1.847$,
$-1.847$,
$-2.613$) 

\vskip 0.7ex
\hangindent=3em \hangafter=1
$D^2= 27.313 = 
16+8\sqrt{2}$

\vskip 0.7ex
\hangindent=3em \hangafter=1
$T = ( 0,
\frac{1}{2},
\frac{1}{4},
\frac{3}{4},
\frac{27}{32},
\frac{27}{32},
\frac{7}{32} )
$,

\vskip 0.7ex
\hangindent=3em \hangafter=1
$S$ = ($ 1$,
$ 1$,
$ 1+\sqrt{2}$,
$ 1+\sqrt{2}$,
$ -c_{16}^{1}$,
$ -c_{16}^{1}$,
$ -c^{1}_{16}
-c^{3}_{16}
$;\ \ 
$ 1$,
$ 1+\sqrt{2}$,
$ 1+\sqrt{2}$,
$ c_{16}^{1}$,
$ c_{16}^{1}$,
$ c^{1}_{16}
+c^{3}_{16}
$;\ \ 
$ -1$,
$ -1$,
$ -c_{16}^{1}$,
$ -c_{16}^{1}$,
$ c^{1}_{16}
+c^{3}_{16}
$;\ \ 
$ -1$,
$ c_{16}^{1}$,
$ c_{16}^{1}$,
$ -c^{1}_{16}
-c^{3}_{16}
$;\ \ 
$ -c^{1}_{16}
-c^{3}_{16}
$,
$ c^{1}_{16}
+c^{3}_{16}
$,
$0$;\ \ 
$ -c^{1}_{16}
-c^{3}_{16}
$,
$0$;\ \ 
$0$)

Prime. 

Pseudo-unitary $\sim$  
$7_{\frac{17}{4},27.31}^{32,261}$

\vskip 1ex 
\color{grey}

\noindent(3,9). $7_{\frac{5}{4},4.686}^{32,350}$ \irep{97}:\ \ 
$d_i$ = ($1.0$,
$0.765$,
$0.765$,
$1.0$,
$-0.414$,
$-0.414$,
$-1.82$) 

\vskip 0.7ex
\hangindent=3em \hangafter=1
$D^2= 4.686 = 
16-8\sqrt{2}$

\vskip 0.7ex
\hangindent=3em \hangafter=1
$T = ( 0,
\frac{7}{32},
\frac{7}{32},
\frac{1}{2},
\frac{1}{4},
\frac{3}{4},
\frac{3}{32} )
$,

\vskip 0.7ex
\hangindent=3em \hangafter=1
$S$ = ($ 1$,
$ c_{16}^{3}$,
$ c_{16}^{3}$,
$ 1$,
$ 1-\sqrt{2}$,
$ 1-\sqrt{2}$,
$ -c^{1}_{16}
+c^{3}_{16}
$;\ \ 
$ -c^{1}_{16}
+c^{3}_{16}
$,
$ c^{1}_{16}
-c^{3}_{16}
$,
$ -c_{16}^{3}$,
$ c_{16}^{3}$,
$ -c_{16}^{3}$,
$0$;\ \ 
$ -c^{1}_{16}
+c^{3}_{16}
$,
$ -c_{16}^{3}$,
$ c_{16}^{3}$,
$ -c_{16}^{3}$,
$0$;\ \ 
$ 1$,
$ 1-\sqrt{2}$,
$ 1-\sqrt{2}$,
$ c^{1}_{16}
-c^{3}_{16}
$;\ \ 
$ -1$,
$ -1$,
$ c^{1}_{16}
-c^{3}_{16}
$;\ \ 
$ -1$,
$ -c^{1}_{16}
+c^{3}_{16}
$;\ \ 
$0$)

Prime. 

Not pseudo-unitary. 

\vskip 1ex 
\color{grey}

\noindent(3,10). $7_{\frac{11}{4},4.686}^{32,471}$ \irep{97}:\ \ 
$d_i$ = ($1.0$,
$0.765$,
$0.765$,
$1.0$,
$-0.414$,
$-0.414$,
$-1.82$) 

\vskip 0.7ex
\hangindent=3em \hangafter=1
$D^2= 4.686 = 
16-8\sqrt{2}$

\vskip 0.7ex
\hangindent=3em \hangafter=1
$T = ( 0,
\frac{9}{32},
\frac{9}{32},
\frac{1}{2},
\frac{1}{4},
\frac{3}{4},
\frac{13}{32} )
$,

\vskip 0.7ex
\hangindent=3em \hangafter=1
$S$ = ($ 1$,
$ c_{16}^{3}$,
$ c_{16}^{3}$,
$ 1$,
$ 1-\sqrt{2}$,
$ 1-\sqrt{2}$,
$ -c^{1}_{16}
+c^{3}_{16}
$;\ \ 
$ -c^{1}_{16}
+c^{3}_{16}
$,
$ c^{1}_{16}
-c^{3}_{16}
$,
$ -c_{16}^{3}$,
$ -c_{16}^{3}$,
$ c_{16}^{3}$,
$0$;\ \ 
$ -c^{1}_{16}
+c^{3}_{16}
$,
$ -c_{16}^{3}$,
$ -c_{16}^{3}$,
$ c_{16}^{3}$,
$0$;\ \ 
$ 1$,
$ 1-\sqrt{2}$,
$ 1-\sqrt{2}$,
$ c^{1}_{16}
-c^{3}_{16}
$;\ \ 
$ -1$,
$ -1$,
$ -c^{1}_{16}
+c^{3}_{16}
$;\ \ 
$ -1$,
$ c^{1}_{16}
-c^{3}_{16}
$;\ \ 
$0$)

Prime. 

Not pseudo-unitary. 

\vskip 1ex 
\color{grey}

\noindent(3,11). $7_{\frac{3}{4},4.686}^{32,164}$ \irep{97}:\ \ 
$d_i$ = ($1.0$,
$1.0$,
$1.82$,
$-0.414$,
$-0.414$,
$-0.765$,
$-0.765$) 

\vskip 0.7ex
\hangindent=3em \hangafter=1
$D^2= 4.686 = 
16-8\sqrt{2}$

\vskip 0.7ex
\hangindent=3em \hangafter=1
$T = ( 0,
\frac{1}{2},
\frac{5}{32},
\frac{1}{4},
\frac{3}{4},
\frac{1}{32},
\frac{1}{32} )
$,

\vskip 0.7ex
\hangindent=3em \hangafter=1
$S$ = ($ 1$,
$ 1$,
$ c^{1}_{16}
-c^{3}_{16}
$,
$ 1-\sqrt{2}$,
$ 1-\sqrt{2}$,
$ -c_{16}^{3}$,
$ -c_{16}^{3}$;\ \ 
$ 1$,
$ -c^{1}_{16}
+c^{3}_{16}
$,
$ 1-\sqrt{2}$,
$ 1-\sqrt{2}$,
$ c_{16}^{3}$,
$ c_{16}^{3}$;\ \ 
$0$,
$ c^{1}_{16}
-c^{3}_{16}
$,
$ -c^{1}_{16}
+c^{3}_{16}
$,
$0$,
$0$;\ \ 
$ -1$,
$ -1$,
$ c_{16}^{3}$,
$ c_{16}^{3}$;\ \ 
$ -1$,
$ -c_{16}^{3}$,
$ -c_{16}^{3}$;\ \ 
$ c^{1}_{16}
-c^{3}_{16}
$,
$ -c^{1}_{16}
+c^{3}_{16}
$;\ \ 
$ c^{1}_{16}
-c^{3}_{16}
$)

Prime. 

Not pseudo-unitary. 

\vskip 1ex 
\color{grey}

\noindent(3,12). $7_{\frac{13}{4},4.686}^{32,398}$ \irep{97}:\ \ 
$d_i$ = ($1.0$,
$1.0$,
$1.82$,
$-0.414$,
$-0.414$,
$-0.765$,
$-0.765$) 

\vskip 0.7ex
\hangindent=3em \hangafter=1
$D^2= 4.686 = 
16-8\sqrt{2}$

\vskip 0.7ex
\hangindent=3em \hangafter=1
$T = ( 0,
\frac{1}{2},
\frac{11}{32},
\frac{1}{4},
\frac{3}{4},
\frac{15}{32},
\frac{15}{32} )
$,

\vskip 0.7ex
\hangindent=3em \hangafter=1
$S$ = ($ 1$,
$ 1$,
$ c^{1}_{16}
-c^{3}_{16}
$,
$ 1-\sqrt{2}$,
$ 1-\sqrt{2}$,
$ -c_{16}^{3}$,
$ -c_{16}^{3}$;\ \ 
$ 1$,
$ -c^{1}_{16}
+c^{3}_{16}
$,
$ 1-\sqrt{2}$,
$ 1-\sqrt{2}$,
$ c_{16}^{3}$,
$ c_{16}^{3}$;\ \ 
$0$,
$ -c^{1}_{16}
+c^{3}_{16}
$,
$ c^{1}_{16}
-c^{3}_{16}
$,
$0$,
$0$;\ \ 
$ -1$,
$ -1$,
$ -c_{16}^{3}$,
$ -c_{16}^{3}$;\ \ 
$ -1$,
$ c_{16}^{3}$,
$ c_{16}^{3}$;\ \ 
$ c^{1}_{16}
-c^{3}_{16}
$,
$ -c^{1}_{16}
+c^{3}_{16}
$;\ \ 
$ c^{1}_{16}
-c^{3}_{16}
$)

Prime. 

Not pseudo-unitary. 

\vskip 1ex 
\color{grey}

\noindent(3,13). $7_{\frac{19}{4},4.686}^{32,155}$ \irep{97}:\ \ 
$d_i$ = ($1.0$,
$1.0$,
$1.82$,
$-0.414$,
$-0.414$,
$-0.765$,
$-0.765$) 

\vskip 0.7ex
\hangindent=3em \hangafter=1
$D^2= 4.686 = 
16-8\sqrt{2}$

\vskip 0.7ex
\hangindent=3em \hangafter=1
$T = ( 0,
\frac{1}{2},
\frac{21}{32},
\frac{1}{4},
\frac{3}{4},
\frac{17}{32},
\frac{17}{32} )
$,

\vskip 0.7ex
\hangindent=3em \hangafter=1
$S$ = ($ 1$,
$ 1$,
$ c^{1}_{16}
-c^{3}_{16}
$,
$ 1-\sqrt{2}$,
$ 1-\sqrt{2}$,
$ -c_{16}^{3}$,
$ -c_{16}^{3}$;\ \ 
$ 1$,
$ -c^{1}_{16}
+c^{3}_{16}
$,
$ 1-\sqrt{2}$,
$ 1-\sqrt{2}$,
$ c_{16}^{3}$,
$ c_{16}^{3}$;\ \ 
$0$,
$ c^{1}_{16}
-c^{3}_{16}
$,
$ -c^{1}_{16}
+c^{3}_{16}
$,
$0$,
$0$;\ \ 
$ -1$,
$ -1$,
$ c_{16}^{3}$,
$ c_{16}^{3}$;\ \ 
$ -1$,
$ -c_{16}^{3}$,
$ -c_{16}^{3}$;\ \ 
$ c^{1}_{16}
-c^{3}_{16}
$,
$ -c^{1}_{16}
+c^{3}_{16}
$;\ \ 
$ c^{1}_{16}
-c^{3}_{16}
$)

Prime. 

Not pseudo-unitary. 

\vskip 1ex 
\color{grey}

\noindent(3,14). $7_{\frac{29}{4},4.686}^{32,855}$ \irep{97}:\ \ 
$d_i$ = ($1.0$,
$1.0$,
$1.82$,
$-0.414$,
$-0.414$,
$-0.765$,
$-0.765$) 

\vskip 0.7ex
\hangindent=3em \hangafter=1
$D^2= 4.686 = 
16-8\sqrt{2}$

\vskip 0.7ex
\hangindent=3em \hangafter=1
$T = ( 0,
\frac{1}{2},
\frac{27}{32},
\frac{1}{4},
\frac{3}{4},
\frac{31}{32},
\frac{31}{32} )
$,

\vskip 0.7ex
\hangindent=3em \hangafter=1
$S$ = ($ 1$,
$ 1$,
$ c^{1}_{16}
-c^{3}_{16}
$,
$ 1-\sqrt{2}$,
$ 1-\sqrt{2}$,
$ -c_{16}^{3}$,
$ -c_{16}^{3}$;\ \ 
$ 1$,
$ -c^{1}_{16}
+c^{3}_{16}
$,
$ 1-\sqrt{2}$,
$ 1-\sqrt{2}$,
$ c_{16}^{3}$,
$ c_{16}^{3}$;\ \ 
$0$,
$ -c^{1}_{16}
+c^{3}_{16}
$,
$ c^{1}_{16}
-c^{3}_{16}
$,
$0$,
$0$;\ \ 
$ -1$,
$ -1$,
$ -c_{16}^{3}$,
$ -c_{16}^{3}$;\ \ 
$ -1$,
$ c_{16}^{3}$,
$ c_{16}^{3}$;\ \ 
$ c^{1}_{16}
-c^{3}_{16}
$,
$ -c^{1}_{16}
+c^{3}_{16}
$;\ \ 
$ c^{1}_{16}
-c^{3}_{16}
$)

Prime. 

Not pseudo-unitary. 

\vskip 1ex 
\color{grey}

\noindent(3,15). $7_{\frac{21}{4},4.686}^{32,466}$ \irep{97}:\ \ 
$d_i$ = ($1.0$,
$0.765$,
$0.765$,
$1.0$,
$-0.414$,
$-0.414$,
$-1.82$) 

\vskip 0.7ex
\hangindent=3em \hangafter=1
$D^2= 4.686 = 
16-8\sqrt{2}$

\vskip 0.7ex
\hangindent=3em \hangafter=1
$T = ( 0,
\frac{23}{32},
\frac{23}{32},
\frac{1}{2},
\frac{1}{4},
\frac{3}{4},
\frac{19}{32} )
$,

\vskip 0.7ex
\hangindent=3em \hangafter=1
$S$ = ($ 1$,
$ c_{16}^{3}$,
$ c_{16}^{3}$,
$ 1$,
$ 1-\sqrt{2}$,
$ 1-\sqrt{2}$,
$ -c^{1}_{16}
+c^{3}_{16}
$;\ \ 
$ -c^{1}_{16}
+c^{3}_{16}
$,
$ c^{1}_{16}
-c^{3}_{16}
$,
$ -c_{16}^{3}$,
$ c_{16}^{3}$,
$ -c_{16}^{3}$,
$0$;\ \ 
$ -c^{1}_{16}
+c^{3}_{16}
$,
$ -c_{16}^{3}$,
$ c_{16}^{3}$,
$ -c_{16}^{3}$,
$0$;\ \ 
$ 1$,
$ 1-\sqrt{2}$,
$ 1-\sqrt{2}$,
$ c^{1}_{16}
-c^{3}_{16}
$;\ \ 
$ -1$,
$ -1$,
$ c^{1}_{16}
-c^{3}_{16}
$;\ \ 
$ -1$,
$ -c^{1}_{16}
+c^{3}_{16}
$;\ \ 
$0$)

Prime. 

Not pseudo-unitary. 

\vskip 1ex 
\color{grey}

\noindent(3,16). $7_{\frac{27}{4},4.686}^{32,368}$ \irep{97}:\ \ 
$d_i$ = ($1.0$,
$0.765$,
$0.765$,
$1.0$,
$-0.414$,
$-0.414$,
$-1.82$) 

\vskip 0.7ex
\hangindent=3em \hangafter=1
$D^2= 4.686 = 
16-8\sqrt{2}$

\vskip 0.7ex
\hangindent=3em \hangafter=1
$T = ( 0,
\frac{25}{32},
\frac{25}{32},
\frac{1}{2},
\frac{1}{4},
\frac{3}{4},
\frac{29}{32} )
$,

\vskip 0.7ex
\hangindent=3em \hangafter=1
$S$ = ($ 1$,
$ c_{16}^{3}$,
$ c_{16}^{3}$,
$ 1$,
$ 1-\sqrt{2}$,
$ 1-\sqrt{2}$,
$ -c^{1}_{16}
+c^{3}_{16}
$;\ \ 
$ -c^{1}_{16}
+c^{3}_{16}
$,
$ c^{1}_{16}
-c^{3}_{16}
$,
$ -c_{16}^{3}$,
$ -c_{16}^{3}$,
$ c_{16}^{3}$,
$0$;\ \ 
$ -c^{1}_{16}
+c^{3}_{16}
$,
$ -c_{16}^{3}$,
$ -c_{16}^{3}$,
$ c_{16}^{3}$,
$0$;\ \ 
$ 1$,
$ 1-\sqrt{2}$,
$ 1-\sqrt{2}$,
$ c^{1}_{16}
-c^{3}_{16}
$;\ \ 
$ -1$,
$ -1$,
$ -c^{1}_{16}
+c^{3}_{16}
$;\ \ 
$ -1$,
$ c^{1}_{16}
-c^{3}_{16}
$;\ \ 
$0$)

Prime. 

Not pseudo-unitary. 

\vskip 1ex 
\black

\noindent(4,1). $7_{\frac{31}{4},27.31}^{32,159}$ \irep{97}:\ \ 
$d_i$ = ($1.0$,
$1.0$,
$1.847$,
$1.847$,
$2.414$,
$2.414$,
$2.613$) 

\vskip 0.7ex
\hangindent=3em \hangafter=1
$D^2= 27.313 = 
16+8\sqrt{2}$

\vskip 0.7ex
\hangindent=3em \hangafter=1
$T = ( 0,
\frac{1}{2},
\frac{5}{32},
\frac{5}{32},
\frac{1}{4},
\frac{3}{4},
\frac{25}{32} )
$,

\vskip 0.7ex
\hangindent=3em \hangafter=1
$S$ = ($ 1$,
$ 1$,
$ c_{16}^{1}$,
$ c_{16}^{1}$,
$ 1+\sqrt{2}$,
$ 1+\sqrt{2}$,
$ c^{1}_{16}
+c^{3}_{16}
$;\ \ 
$ 1$,
$ -c_{16}^{1}$,
$ -c_{16}^{1}$,
$ 1+\sqrt{2}$,
$ 1+\sqrt{2}$,
$ -c^{1}_{16}
-c^{3}_{16}
$;\ \ 
$ -c^{1}_{16}
-c^{3}_{16}
$,
$ c^{1}_{16}
+c^{3}_{16}
$,
$ -c_{16}^{1}$,
$ c_{16}^{1}$,
$0$;\ \ 
$ -c^{1}_{16}
-c^{3}_{16}
$,
$ -c_{16}^{1}$,
$ c_{16}^{1}$,
$0$;\ \ 
$ -1$,
$ -1$,
$ c^{1}_{16}
+c^{3}_{16}
$;\ \ 
$ -1$,
$ -c^{1}_{16}
-c^{3}_{16}
$;\ \ 
$0$)

Prime. 

\vskip 1ex 
\color{grey}

\noindent(4,2). $7_{\frac{17}{4},27.31}^{32,261}$ \irep{97}:\ \ 
$d_i$ = ($1.0$,
$1.0$,
$1.847$,
$1.847$,
$2.414$,
$2.414$,
$2.613$) 

\vskip 0.7ex
\hangindent=3em \hangafter=1
$D^2= 27.313 = 
16+8\sqrt{2}$

\vskip 0.7ex
\hangindent=3em \hangafter=1
$T = ( 0,
\frac{1}{2},
\frac{11}{32},
\frac{11}{32},
\frac{1}{4},
\frac{3}{4},
\frac{23}{32} )
$,

\vskip 0.7ex
\hangindent=3em \hangafter=1
$S$ = ($ 1$,
$ 1$,
$ c_{16}^{1}$,
$ c_{16}^{1}$,
$ 1+\sqrt{2}$,
$ 1+\sqrt{2}$,
$ c^{1}_{16}
+c^{3}_{16}
$;\ \ 
$ 1$,
$ -c_{16}^{1}$,
$ -c_{16}^{1}$,
$ 1+\sqrt{2}$,
$ 1+\sqrt{2}$,
$ -c^{1}_{16}
-c^{3}_{16}
$;\ \ 
$ -c^{1}_{16}
-c^{3}_{16}
$,
$ c^{1}_{16}
+c^{3}_{16}
$,
$ c_{16}^{1}$,
$ -c_{16}^{1}$,
$0$;\ \ 
$ -c^{1}_{16}
-c^{3}_{16}
$,
$ c_{16}^{1}$,
$ -c_{16}^{1}$,
$0$;\ \ 
$ -1$,
$ -1$,
$ -c^{1}_{16}
-c^{3}_{16}
$;\ \ 
$ -1$,
$ c^{1}_{16}
+c^{3}_{16}
$;\ \ 
$0$)

Prime. 

\vskip 1ex 
\color{grey}

\noindent(4,3). $7_{\frac{15}{4},27.31}^{32,272}$ \irep{97}:\ \ 
$d_i$ = ($1.0$,
$1.0$,
$1.847$,
$1.847$,
$2.414$,
$2.414$,
$2.613$) 

\vskip 0.7ex
\hangindent=3em \hangafter=1
$D^2= 27.313 = 
16+8\sqrt{2}$

\vskip 0.7ex
\hangindent=3em \hangafter=1
$T = ( 0,
\frac{1}{2},
\frac{21}{32},
\frac{21}{32},
\frac{1}{4},
\frac{3}{4},
\frac{9}{32} )
$,

\vskip 0.7ex
\hangindent=3em \hangafter=1
$S$ = ($ 1$,
$ 1$,
$ c_{16}^{1}$,
$ c_{16}^{1}$,
$ 1+\sqrt{2}$,
$ 1+\sqrt{2}$,
$ c^{1}_{16}
+c^{3}_{16}
$;\ \ 
$ 1$,
$ -c_{16}^{1}$,
$ -c_{16}^{1}$,
$ 1+\sqrt{2}$,
$ 1+\sqrt{2}$,
$ -c^{1}_{16}
-c^{3}_{16}
$;\ \ 
$ -c^{1}_{16}
-c^{3}_{16}
$,
$ c^{1}_{16}
+c^{3}_{16}
$,
$ -c_{16}^{1}$,
$ c_{16}^{1}$,
$0$;\ \ 
$ -c^{1}_{16}
-c^{3}_{16}
$,
$ -c_{16}^{1}$,
$ c_{16}^{1}$,
$0$;\ \ 
$ -1$,
$ -1$,
$ c^{1}_{16}
+c^{3}_{16}
$;\ \ 
$ -1$,
$ -c^{1}_{16}
-c^{3}_{16}
$;\ \ 
$0$)

Prime. 

\vskip 1ex 
\color{grey}

\noindent(4,4). $7_{\frac{1}{4},27.31}^{32,123}$ \irep{97}:\ \ 
$d_i$ = ($1.0$,
$1.0$,
$1.847$,
$1.847$,
$2.414$,
$2.414$,
$2.613$) 

\vskip 0.7ex
\hangindent=3em \hangafter=1
$D^2= 27.313 = 
16+8\sqrt{2}$

\vskip 0.7ex
\hangindent=3em \hangafter=1
$T = ( 0,
\frac{1}{2},
\frac{27}{32},
\frac{27}{32},
\frac{1}{4},
\frac{3}{4},
\frac{7}{32} )
$,

\vskip 0.7ex
\hangindent=3em \hangafter=1
$S$ = ($ 1$,
$ 1$,
$ c_{16}^{1}$,
$ c_{16}^{1}$,
$ 1+\sqrt{2}$,
$ 1+\sqrt{2}$,
$ c^{1}_{16}
+c^{3}_{16}
$;\ \ 
$ 1$,
$ -c_{16}^{1}$,
$ -c_{16}^{1}$,
$ 1+\sqrt{2}$,
$ 1+\sqrt{2}$,
$ -c^{1}_{16}
-c^{3}_{16}
$;\ \ 
$ -c^{1}_{16}
-c^{3}_{16}
$,
$ c^{1}_{16}
+c^{3}_{16}
$,
$ c_{16}^{1}$,
$ -c_{16}^{1}$,
$0$;\ \ 
$ -c^{1}_{16}
-c^{3}_{16}
$,
$ c_{16}^{1}$,
$ -c_{16}^{1}$,
$0$;\ \ 
$ -1$,
$ -1$,
$ -c^{1}_{16}
-c^{3}_{16}
$;\ \ 
$ -1$,
$ c^{1}_{16}
+c^{3}_{16}
$;\ \ 
$0$)

Prime. 

\vskip 1ex 
\color{grey}

\noindent(4,5). $7_{\frac{9}{4},27.31}^{32,102}$ \irep{97}:\ \ 
$d_i$ = ($1.0$,
$1.0$,
$2.414$,
$2.414$,
$-1.847$,
$-1.847$,
$-2.613$) 

\vskip 0.7ex
\hangindent=3em \hangafter=1
$D^2= 27.313 = 
16+8\sqrt{2}$

\vskip 0.7ex
\hangindent=3em \hangafter=1
$T = ( 0,
\frac{1}{2},
\frac{1}{4},
\frac{3}{4},
\frac{3}{32},
\frac{3}{32},
\frac{15}{32} )
$,

\vskip 0.7ex
\hangindent=3em \hangafter=1
$S$ = ($ 1$,
$ 1$,
$ 1+\sqrt{2}$,
$ 1+\sqrt{2}$,
$ -c_{16}^{1}$,
$ -c_{16}^{1}$,
$ -c^{1}_{16}
-c^{3}_{16}
$;\ \ 
$ 1$,
$ 1+\sqrt{2}$,
$ 1+\sqrt{2}$,
$ c_{16}^{1}$,
$ c_{16}^{1}$,
$ c^{1}_{16}
+c^{3}_{16}
$;\ \ 
$ -1$,
$ -1$,
$ -c_{16}^{1}$,
$ -c_{16}^{1}$,
$ c^{1}_{16}
+c^{3}_{16}
$;\ \ 
$ -1$,
$ c_{16}^{1}$,
$ c_{16}^{1}$,
$ -c^{1}_{16}
-c^{3}_{16}
$;\ \ 
$ c^{1}_{16}
+c^{3}_{16}
$,
$ -c^{1}_{16}
-c^{3}_{16}
$,
$0$;\ \ 
$ c^{1}_{16}
+c^{3}_{16}
$,
$0$;\ \ 
$0$)

Prime. 

Pseudo-unitary $\sim$  
$7_{\frac{25}{4},27.31}^{32,222}$

\vskip 1ex 
\color{grey}

\noindent(4,6). $7_{\frac{7}{4},27.31}^{32,101}$ \irep{97}:\ \ 
$d_i$ = ($1.0$,
$1.0$,
$2.414$,
$2.414$,
$-1.847$,
$-1.847$,
$-2.613$) 

\vskip 0.7ex
\hangindent=3em \hangafter=1
$D^2= 27.313 = 
16+8\sqrt{2}$

\vskip 0.7ex
\hangindent=3em \hangafter=1
$T = ( 0,
\frac{1}{2},
\frac{1}{4},
\frac{3}{4},
\frac{13}{32},
\frac{13}{32},
\frac{1}{32} )
$,

\vskip 0.7ex
\hangindent=3em \hangafter=1
$S$ = ($ 1$,
$ 1$,
$ 1+\sqrt{2}$,
$ 1+\sqrt{2}$,
$ -c_{16}^{1}$,
$ -c_{16}^{1}$,
$ -c^{1}_{16}
-c^{3}_{16}
$;\ \ 
$ 1$,
$ 1+\sqrt{2}$,
$ 1+\sqrt{2}$,
$ c_{16}^{1}$,
$ c_{16}^{1}$,
$ c^{1}_{16}
+c^{3}_{16}
$;\ \ 
$ -1$,
$ -1$,
$ c_{16}^{1}$,
$ c_{16}^{1}$,
$ -c^{1}_{16}
-c^{3}_{16}
$;\ \ 
$ -1$,
$ -c_{16}^{1}$,
$ -c_{16}^{1}$,
$ c^{1}_{16}
+c^{3}_{16}
$;\ \ 
$ c^{1}_{16}
+c^{3}_{16}
$,
$ -c^{1}_{16}
-c^{3}_{16}
$,
$0$;\ \ 
$ c^{1}_{16}
+c^{3}_{16}
$,
$0$;\ \ 
$0$)

Prime. 

Pseudo-unitary $\sim$  
$7_{\frac{23}{4},27.31}^{32,224}$

\vskip 1ex 
\color{grey}

\noindent(4,7). $7_{\frac{25}{4},27.31}^{32,212}$ \irep{97}:\ \ 
$d_i$ = ($1.0$,
$1.0$,
$2.414$,
$2.414$,
$-1.847$,
$-1.847$,
$-2.613$) 

\vskip 0.7ex
\hangindent=3em \hangafter=1
$D^2= 27.313 = 
16+8\sqrt{2}$

\vskip 0.7ex
\hangindent=3em \hangafter=1
$T = ( 0,
\frac{1}{2},
\frac{1}{4},
\frac{3}{4},
\frac{19}{32},
\frac{19}{32},
\frac{31}{32} )
$,

\vskip 0.7ex
\hangindent=3em \hangafter=1
$S$ = ($ 1$,
$ 1$,
$ 1+\sqrt{2}$,
$ 1+\sqrt{2}$,
$ -c_{16}^{1}$,
$ -c_{16}^{1}$,
$ -c^{1}_{16}
-c^{3}_{16}
$;\ \ 
$ 1$,
$ 1+\sqrt{2}$,
$ 1+\sqrt{2}$,
$ c_{16}^{1}$,
$ c_{16}^{1}$,
$ c^{1}_{16}
+c^{3}_{16}
$;\ \ 
$ -1$,
$ -1$,
$ -c_{16}^{1}$,
$ -c_{16}^{1}$,
$ c^{1}_{16}
+c^{3}_{16}
$;\ \ 
$ -1$,
$ c_{16}^{1}$,
$ c_{16}^{1}$,
$ -c^{1}_{16}
-c^{3}_{16}
$;\ \ 
$ c^{1}_{16}
+c^{3}_{16}
$,
$ -c^{1}_{16}
-c^{3}_{16}
$,
$0$;\ \ 
$ c^{1}_{16}
+c^{3}_{16}
$,
$0$;\ \ 
$0$)

Prime. 

Pseudo-unitary $\sim$  
$7_{\frac{9}{4},27.31}^{32,918}$

\vskip 1ex 
\color{grey}

\noindent(4,8). $7_{\frac{23}{4},27.31}^{32,214}$ \irep{97}:\ \ 
$d_i$ = ($1.0$,
$1.0$,
$2.414$,
$2.414$,
$-1.847$,
$-1.847$,
$-2.613$) 

\vskip 0.7ex
\hangindent=3em \hangafter=1
$D^2= 27.313 = 
16+8\sqrt{2}$

\vskip 0.7ex
\hangindent=3em \hangafter=1
$T = ( 0,
\frac{1}{2},
\frac{1}{4},
\frac{3}{4},
\frac{29}{32},
\frac{29}{32},
\frac{17}{32} )
$,

\vskip 0.7ex
\hangindent=3em \hangafter=1
$S$ = ($ 1$,
$ 1$,
$ 1+\sqrt{2}$,
$ 1+\sqrt{2}$,
$ -c_{16}^{1}$,
$ -c_{16}^{1}$,
$ -c^{1}_{16}
-c^{3}_{16}
$;\ \ 
$ 1$,
$ 1+\sqrt{2}$,
$ 1+\sqrt{2}$,
$ c_{16}^{1}$,
$ c_{16}^{1}$,
$ c^{1}_{16}
+c^{3}_{16}
$;\ \ 
$ -1$,
$ -1$,
$ c_{16}^{1}$,
$ c_{16}^{1}$,
$ -c^{1}_{16}
-c^{3}_{16}
$;\ \ 
$ -1$,
$ -c_{16}^{1}$,
$ -c_{16}^{1}$,
$ c^{1}_{16}
+c^{3}_{16}
$;\ \ 
$ c^{1}_{16}
+c^{3}_{16}
$,
$ -c^{1}_{16}
-c^{3}_{16}
$,
$0$;\ \ 
$ c^{1}_{16}
+c^{3}_{16}
$,
$0$;\ \ 
$0$)

Prime. 

Pseudo-unitary $\sim$  
$7_{\frac{7}{4},27.31}^{32,912}$

\vskip 1ex 
\color{grey}

\noindent(4,9). $7_{\frac{3}{4},4.686}^{32,438}$ \irep{97}:\ \ 
$d_i$ = ($1.0$,
$0.765$,
$0.765$,
$1.0$,
$-0.414$,
$-0.414$,
$-1.82$) 

\vskip 0.7ex
\hangindent=3em \hangafter=1
$D^2= 4.686 = 
16-8\sqrt{2}$

\vskip 0.7ex
\hangindent=3em \hangafter=1
$T = ( 0,
\frac{1}{32},
\frac{1}{32},
\frac{1}{2},
\frac{1}{4},
\frac{3}{4},
\frac{5}{32} )
$,

\vskip 0.7ex
\hangindent=3em \hangafter=1
$S$ = ($ 1$,
$ c_{16}^{3}$,
$ c_{16}^{3}$,
$ 1$,
$ 1-\sqrt{2}$,
$ 1-\sqrt{2}$,
$ -c^{1}_{16}
+c^{3}_{16}
$;\ \ 
$ c^{1}_{16}
-c^{3}_{16}
$,
$ -c^{1}_{16}
+c^{3}_{16}
$,
$ -c_{16}^{3}$,
$ -c_{16}^{3}$,
$ c_{16}^{3}$,
$0$;\ \ 
$ c^{1}_{16}
-c^{3}_{16}
$,
$ -c_{16}^{3}$,
$ -c_{16}^{3}$,
$ c_{16}^{3}$,
$0$;\ \ 
$ 1$,
$ 1-\sqrt{2}$,
$ 1-\sqrt{2}$,
$ c^{1}_{16}
-c^{3}_{16}
$;\ \ 
$ -1$,
$ -1$,
$ -c^{1}_{16}
+c^{3}_{16}
$;\ \ 
$ -1$,
$ c^{1}_{16}
-c^{3}_{16}
$;\ \ 
$0$)

Prime. 

Not pseudo-unitary. 

\vskip 1ex 
\color{grey}

\noindent(4,10). $7_{\frac{13}{4},4.686}^{32,205}$ \irep{97}:\ \ 
$d_i$ = ($1.0$,
$0.765$,
$0.765$,
$1.0$,
$-0.414$,
$-0.414$,
$-1.82$) 

\vskip 0.7ex
\hangindent=3em \hangafter=1
$D^2= 4.686 = 
16-8\sqrt{2}$

\vskip 0.7ex
\hangindent=3em \hangafter=1
$T = ( 0,
\frac{15}{32},
\frac{15}{32},
\frac{1}{2},
\frac{1}{4},
\frac{3}{4},
\frac{11}{32} )
$,

\vskip 0.7ex
\hangindent=3em \hangafter=1
$S$ = ($ 1$,
$ c_{16}^{3}$,
$ c_{16}^{3}$,
$ 1$,
$ 1-\sqrt{2}$,
$ 1-\sqrt{2}$,
$ -c^{1}_{16}
+c^{3}_{16}
$;\ \ 
$ c^{1}_{16}
-c^{3}_{16}
$,
$ -c^{1}_{16}
+c^{3}_{16}
$,
$ -c_{16}^{3}$,
$ c_{16}^{3}$,
$ -c_{16}^{3}$,
$0$;\ \ 
$ c^{1}_{16}
-c^{3}_{16}
$,
$ -c_{16}^{3}$,
$ c_{16}^{3}$,
$ -c_{16}^{3}$,
$0$;\ \ 
$ 1$,
$ 1-\sqrt{2}$,
$ 1-\sqrt{2}$,
$ c^{1}_{16}
-c^{3}_{16}
$;\ \ 
$ -1$,
$ -1$,
$ c^{1}_{16}
-c^{3}_{16}
$;\ \ 
$ -1$,
$ -c^{1}_{16}
+c^{3}_{16}
$;\ \ 
$0$)

Prime. 

Not pseudo-unitary. 

\vskip 1ex 
\color{grey}

\noindent(4,11). $7_{\frac{5}{4},4.686}^{32,253}$ \irep{97}:\ \ 
$d_i$ = ($1.0$,
$1.0$,
$1.82$,
$-0.414$,
$-0.414$,
$-0.765$,
$-0.765$) 

\vskip 0.7ex
\hangindent=3em \hangafter=1
$D^2= 4.686 = 
16-8\sqrt{2}$

\vskip 0.7ex
\hangindent=3em \hangafter=1
$T = ( 0,
\frac{1}{2},
\frac{3}{32},
\frac{1}{4},
\frac{3}{4},
\frac{7}{32},
\frac{7}{32} )
$,

\vskip 0.7ex
\hangindent=3em \hangafter=1
$S$ = ($ 1$,
$ 1$,
$ c^{1}_{16}
-c^{3}_{16}
$,
$ 1-\sqrt{2}$,
$ 1-\sqrt{2}$,
$ -c_{16}^{3}$,
$ -c_{16}^{3}$;\ \ 
$ 1$,
$ -c^{1}_{16}
+c^{3}_{16}
$,
$ 1-\sqrt{2}$,
$ 1-\sqrt{2}$,
$ c_{16}^{3}$,
$ c_{16}^{3}$;\ \ 
$0$,
$ -c^{1}_{16}
+c^{3}_{16}
$,
$ c^{1}_{16}
-c^{3}_{16}
$,
$0$,
$0$;\ \ 
$ -1$,
$ -1$,
$ -c_{16}^{3}$,
$ -c_{16}^{3}$;\ \ 
$ -1$,
$ c_{16}^{3}$,
$ c_{16}^{3}$;\ \ 
$ -c^{1}_{16}
+c^{3}_{16}
$,
$ c^{1}_{16}
-c^{3}_{16}
$;\ \ 
$ -c^{1}_{16}
+c^{3}_{16}
$)

Prime. 

Not pseudo-unitary. 

\vskip 1ex 
\color{grey}

\noindent(4,12). $7_{\frac{11}{4},4.686}^{32,132}$ \irep{97}:\ \ 
$d_i$ = ($1.0$,
$1.0$,
$1.82$,
$-0.414$,
$-0.414$,
$-0.765$,
$-0.765$) 

\vskip 0.7ex
\hangindent=3em \hangafter=1
$D^2= 4.686 = 
16-8\sqrt{2}$

\vskip 0.7ex
\hangindent=3em \hangafter=1
$T = ( 0,
\frac{1}{2},
\frac{13}{32},
\frac{1}{4},
\frac{3}{4},
\frac{9}{32},
\frac{9}{32} )
$,

\vskip 0.7ex
\hangindent=3em \hangafter=1
$S$ = ($ 1$,
$ 1$,
$ c^{1}_{16}
-c^{3}_{16}
$,
$ 1-\sqrt{2}$,
$ 1-\sqrt{2}$,
$ -c_{16}^{3}$,
$ -c_{16}^{3}$;\ \ 
$ 1$,
$ -c^{1}_{16}
+c^{3}_{16}
$,
$ 1-\sqrt{2}$,
$ 1-\sqrt{2}$,
$ c_{16}^{3}$,
$ c_{16}^{3}$;\ \ 
$0$,
$ c^{1}_{16}
-c^{3}_{16}
$,
$ -c^{1}_{16}
+c^{3}_{16}
$,
$0$,
$0$;\ \ 
$ -1$,
$ -1$,
$ c_{16}^{3}$,
$ c_{16}^{3}$;\ \ 
$ -1$,
$ -c_{16}^{3}$,
$ -c_{16}^{3}$;\ \ 
$ -c^{1}_{16}
+c^{3}_{16}
$,
$ c^{1}_{16}
-c^{3}_{16}
$;\ \ 
$ -c^{1}_{16}
+c^{3}_{16}
$)

Prime. 

Not pseudo-unitary. 

\vskip 1ex 
\color{grey}

\noindent(4,13). $7_{\frac{21}{4},4.686}^{32,598}$ \irep{97}:\ \ 
$d_i$ = ($1.0$,
$1.0$,
$1.82$,
$-0.414$,
$-0.414$,
$-0.765$,
$-0.765$) 

\vskip 0.7ex
\hangindent=3em \hangafter=1
$D^2= 4.686 = 
16-8\sqrt{2}$

\vskip 0.7ex
\hangindent=3em \hangafter=1
$T = ( 0,
\frac{1}{2},
\frac{19}{32},
\frac{1}{4},
\frac{3}{4},
\frac{23}{32},
\frac{23}{32} )
$,

\vskip 0.7ex
\hangindent=3em \hangafter=1
$S$ = ($ 1$,
$ 1$,
$ c^{1}_{16}
-c^{3}_{16}
$,
$ 1-\sqrt{2}$,
$ 1-\sqrt{2}$,
$ -c_{16}^{3}$,
$ -c_{16}^{3}$;\ \ 
$ 1$,
$ -c^{1}_{16}
+c^{3}_{16}
$,
$ 1-\sqrt{2}$,
$ 1-\sqrt{2}$,
$ c_{16}^{3}$,
$ c_{16}^{3}$;\ \ 
$0$,
$ -c^{1}_{16}
+c^{3}_{16}
$,
$ c^{1}_{16}
-c^{3}_{16}
$,
$0$,
$0$;\ \ 
$ -1$,
$ -1$,
$ -c_{16}^{3}$,
$ -c_{16}^{3}$;\ \ 
$ -1$,
$ c_{16}^{3}$,
$ c_{16}^{3}$;\ \ 
$ -c^{1}_{16}
+c^{3}_{16}
$,
$ c^{1}_{16}
-c^{3}_{16}
$;\ \ 
$ -c^{1}_{16}
+c^{3}_{16}
$)

Prime. 

Not pseudo-unitary. 

\vskip 1ex 
\color{grey}

\noindent(4,14). $7_{\frac{27}{4},4.686}^{32,235}$ \irep{97}:\ \ 
$d_i$ = ($1.0$,
$1.0$,
$1.82$,
$-0.414$,
$-0.414$,
$-0.765$,
$-0.765$) 

\vskip 0.7ex
\hangindent=3em \hangafter=1
$D^2= 4.686 = 
16-8\sqrt{2}$

\vskip 0.7ex
\hangindent=3em \hangafter=1
$T = ( 0,
\frac{1}{2},
\frac{29}{32},
\frac{1}{4},
\frac{3}{4},
\frac{25}{32},
\frac{25}{32} )
$,

\vskip 0.7ex
\hangindent=3em \hangafter=1
$S$ = ($ 1$,
$ 1$,
$ c^{1}_{16}
-c^{3}_{16}
$,
$ 1-\sqrt{2}$,
$ 1-\sqrt{2}$,
$ -c_{16}^{3}$,
$ -c_{16}^{3}$;\ \ 
$ 1$,
$ -c^{1}_{16}
+c^{3}_{16}
$,
$ 1-\sqrt{2}$,
$ 1-\sqrt{2}$,
$ c_{16}^{3}$,
$ c_{16}^{3}$;\ \ 
$0$,
$ c^{1}_{16}
-c^{3}_{16}
$,
$ -c^{1}_{16}
+c^{3}_{16}
$,
$0$,
$0$;\ \ 
$ -1$,
$ -1$,
$ c_{16}^{3}$,
$ c_{16}^{3}$;\ \ 
$ -1$,
$ -c_{16}^{3}$,
$ -c_{16}^{3}$;\ \ 
$ -c^{1}_{16}
+c^{3}_{16}
$,
$ c^{1}_{16}
-c^{3}_{16}
$;\ \ 
$ -c^{1}_{16}
+c^{3}_{16}
$)

Prime. 

Not pseudo-unitary. 

\vskip 1ex 
\color{grey}

\noindent(4,15). $7_{\frac{19}{4},4.686}^{32,447}$ \irep{97}:\ \ 
$d_i$ = ($1.0$,
$0.765$,
$0.765$,
$1.0$,
$-0.414$,
$-0.414$,
$-1.82$) 

\vskip 0.7ex
\hangindent=3em \hangafter=1
$D^2= 4.686 = 
16-8\sqrt{2}$

\vskip 0.7ex
\hangindent=3em \hangafter=1
$T = ( 0,
\frac{17}{32},
\frac{17}{32},
\frac{1}{2},
\frac{1}{4},
\frac{3}{4},
\frac{21}{32} )
$,

\vskip 0.7ex
\hangindent=3em \hangafter=1
$S$ = ($ 1$,
$ c_{16}^{3}$,
$ c_{16}^{3}$,
$ 1$,
$ 1-\sqrt{2}$,
$ 1-\sqrt{2}$,
$ -c^{1}_{16}
+c^{3}_{16}
$;\ \ 
$ c^{1}_{16}
-c^{3}_{16}
$,
$ -c^{1}_{16}
+c^{3}_{16}
$,
$ -c_{16}^{3}$,
$ -c_{16}^{3}$,
$ c_{16}^{3}$,
$0$;\ \ 
$ c^{1}_{16}
-c^{3}_{16}
$,
$ -c_{16}^{3}$,
$ -c_{16}^{3}$,
$ c_{16}^{3}$,
$0$;\ \ 
$ 1$,
$ 1-\sqrt{2}$,
$ 1-\sqrt{2}$,
$ c^{1}_{16}
-c^{3}_{16}
$;\ \ 
$ -1$,
$ -1$,
$ -c^{1}_{16}
+c^{3}_{16}
$;\ \ 
$ -1$,
$ c^{1}_{16}
-c^{3}_{16}
$;\ \ 
$0$)

Prime. 

Not pseudo-unitary. 

\vskip 1ex 
\color{grey}

\noindent(4,16). $7_{\frac{29}{4},4.686}^{32,252}$ \irep{97}:\ \ 
$d_i$ = ($1.0$,
$0.765$,
$0.765$,
$1.0$,
$-0.414$,
$-0.414$,
$-1.82$) 

\vskip 0.7ex
\hangindent=3em \hangafter=1
$D^2= 4.686 = 
16-8\sqrt{2}$

\vskip 0.7ex
\hangindent=3em \hangafter=1
$T = ( 0,
\frac{31}{32},
\frac{31}{32},
\frac{1}{2},
\frac{1}{4},
\frac{3}{4},
\frac{27}{32} )
$,

\vskip 0.7ex
\hangindent=3em \hangafter=1
$S$ = ($ 1$,
$ c_{16}^{3}$,
$ c_{16}^{3}$,
$ 1$,
$ 1-\sqrt{2}$,
$ 1-\sqrt{2}$,
$ -c^{1}_{16}
+c^{3}_{16}
$;\ \ 
$ c^{1}_{16}
-c^{3}_{16}
$,
$ -c^{1}_{16}
+c^{3}_{16}
$,
$ -c_{16}^{3}$,
$ c_{16}^{3}$,
$ -c_{16}^{3}$,
$0$;\ \ 
$ c^{1}_{16}
-c^{3}_{16}
$,
$ -c_{16}^{3}$,
$ c_{16}^{3}$,
$ -c_{16}^{3}$,
$0$;\ \ 
$ 1$,
$ 1-\sqrt{2}$,
$ 1-\sqrt{2}$,
$ c^{1}_{16}
-c^{3}_{16}
$;\ \ 
$ -1$,
$ -1$,
$ c^{1}_{16}
-c^{3}_{16}
$;\ \ 
$ -1$,
$ -c^{1}_{16}
+c^{3}_{16}
$;\ \ 
$0$)

Prime. 

Not pseudo-unitary. 

\vskip 1ex 
\black

\noindent(5,1). $7_{\frac{13}{4},27.31}^{32,427}$ \irep{96}:\ \ 
$d_i$ = ($1.0$,
$1.0$,
$1.847$,
$1.847$,
$2.414$,
$2.414$,
$2.613$) 

\vskip 0.7ex
\hangindent=3em \hangafter=1
$D^2= 27.313 = 
16+8\sqrt{2}$

\vskip 0.7ex
\hangindent=3em \hangafter=1
$T = ( 0,
\frac{1}{2},
\frac{7}{32},
\frac{7}{32},
\frac{1}{4},
\frac{3}{4},
\frac{19}{32} )
$,

\vskip 0.7ex
\hangindent=3em \hangafter=1
$S$ = ($ 1$,
$ 1$,
$ c_{16}^{1}$,
$ c_{16}^{1}$,
$ 1+\sqrt{2}$,
$ 1+\sqrt{2}$,
$ c^{1}_{16}
+c^{3}_{16}
$;\ \ 
$ 1$,
$ -c_{16}^{1}$,
$ -c_{16}^{1}$,
$ 1+\sqrt{2}$,
$ 1+\sqrt{2}$,
$ -c^{1}_{16}
-c^{3}_{16}
$;\ \ 
$(-c^{1}_{16}
-c^{3}_{16}
)\mathrm{i}$,
$(c^{1}_{16}
+c^{3}_{16}
)\mathrm{i}$,
$ c_{16}^{1}$,
$ -c_{16}^{1}$,
$0$;\ \ 
$(-c^{1}_{16}
-c^{3}_{16}
)\mathrm{i}$,
$ c_{16}^{1}$,
$ -c_{16}^{1}$,
$0$;\ \ 
$ -1$,
$ -1$,
$ -c^{1}_{16}
-c^{3}_{16}
$;\ \ 
$ -1$,
$ c^{1}_{16}
+c^{3}_{16}
$;\ \ 
$0$)

Prime. 

\vskip 1ex 
\color{grey}

\noindent(5,2). $7_{\frac{3}{4},27.31}^{32,913}$ \irep{96}:\ \ 
$d_i$ = ($1.0$,
$1.0$,
$1.847$,
$1.847$,
$2.414$,
$2.414$,
$2.613$) 

\vskip 0.7ex
\hangindent=3em \hangafter=1
$D^2= 27.313 = 
16+8\sqrt{2}$

\vskip 0.7ex
\hangindent=3em \hangafter=1
$T = ( 0,
\frac{1}{2},
\frac{9}{32},
\frac{9}{32},
\frac{1}{4},
\frac{3}{4},
\frac{29}{32} )
$,

\vskip 0.7ex
\hangindent=3em \hangafter=1
$S$ = ($ 1$,
$ 1$,
$ c_{16}^{1}$,
$ c_{16}^{1}$,
$ 1+\sqrt{2}$,
$ 1+\sqrt{2}$,
$ c^{1}_{16}
+c^{3}_{16}
$;\ \ 
$ 1$,
$ -c_{16}^{1}$,
$ -c_{16}^{1}$,
$ 1+\sqrt{2}$,
$ 1+\sqrt{2}$,
$ -c^{1}_{16}
-c^{3}_{16}
$;\ \ 
$(c^{1}_{16}
+c^{3}_{16}
)\mathrm{i}$,
$(-c^{1}_{16}
-c^{3}_{16}
)\mathrm{i}$,
$ -c_{16}^{1}$,
$ c_{16}^{1}$,
$0$;\ \ 
$(c^{1}_{16}
+c^{3}_{16}
)\mathrm{i}$,
$ -c_{16}^{1}$,
$ c_{16}^{1}$,
$0$;\ \ 
$ -1$,
$ -1$,
$ c^{1}_{16}
+c^{3}_{16}
$;\ \ 
$ -1$,
$ -c^{1}_{16}
-c^{3}_{16}
$;\ \ 
$0$)

Prime. 

\vskip 1ex 
\color{grey}

\noindent(5,3). $7_{\frac{29}{4},27.31}^{32,406}$ \irep{96}:\ \ 
$d_i$ = ($1.0$,
$1.0$,
$1.847$,
$1.847$,
$2.414$,
$2.414$,
$2.613$) 

\vskip 0.7ex
\hangindent=3em \hangafter=1
$D^2= 27.313 = 
16+8\sqrt{2}$

\vskip 0.7ex
\hangindent=3em \hangafter=1
$T = ( 0,
\frac{1}{2},
\frac{23}{32},
\frac{23}{32},
\frac{1}{4},
\frac{3}{4},
\frac{3}{32} )
$,

\vskip 0.7ex
\hangindent=3em \hangafter=1
$S$ = ($ 1$,
$ 1$,
$ c_{16}^{1}$,
$ c_{16}^{1}$,
$ 1+\sqrt{2}$,
$ 1+\sqrt{2}$,
$ c^{1}_{16}
+c^{3}_{16}
$;\ \ 
$ 1$,
$ -c_{16}^{1}$,
$ -c_{16}^{1}$,
$ 1+\sqrt{2}$,
$ 1+\sqrt{2}$,
$ -c^{1}_{16}
-c^{3}_{16}
$;\ \ 
$(-c^{1}_{16}
-c^{3}_{16}
)\mathrm{i}$,
$(c^{1}_{16}
+c^{3}_{16}
)\mathrm{i}$,
$ c_{16}^{1}$,
$ -c_{16}^{1}$,
$0$;\ \ 
$(-c^{1}_{16}
-c^{3}_{16}
)\mathrm{i}$,
$ c_{16}^{1}$,
$ -c_{16}^{1}$,
$0$;\ \ 
$ -1$,
$ -1$,
$ -c^{1}_{16}
-c^{3}_{16}
$;\ \ 
$ -1$,
$ c^{1}_{16}
+c^{3}_{16}
$;\ \ 
$0$)

Prime. 

\vskip 1ex 
\color{grey}

\noindent(5,4). $7_{\frac{19}{4},27.31}^{32,116}$ \irep{96}:\ \ 
$d_i$ = ($1.0$,
$1.0$,
$1.847$,
$1.847$,
$2.414$,
$2.414$,
$2.613$) 

\vskip 0.7ex
\hangindent=3em \hangafter=1
$D^2= 27.313 = 
16+8\sqrt{2}$

\vskip 0.7ex
\hangindent=3em \hangafter=1
$T = ( 0,
\frac{1}{2},
\frac{25}{32},
\frac{25}{32},
\frac{1}{4},
\frac{3}{4},
\frac{13}{32} )
$,

\vskip 0.7ex
\hangindent=3em \hangafter=1
$S$ = ($ 1$,
$ 1$,
$ c_{16}^{1}$,
$ c_{16}^{1}$,
$ 1+\sqrt{2}$,
$ 1+\sqrt{2}$,
$ c^{1}_{16}
+c^{3}_{16}
$;\ \ 
$ 1$,
$ -c_{16}^{1}$,
$ -c_{16}^{1}$,
$ 1+\sqrt{2}$,
$ 1+\sqrt{2}$,
$ -c^{1}_{16}
-c^{3}_{16}
$;\ \ 
$(c^{1}_{16}
+c^{3}_{16}
)\mathrm{i}$,
$(-c^{1}_{16}
-c^{3}_{16}
)\mathrm{i}$,
$ -c_{16}^{1}$,
$ c_{16}^{1}$,
$0$;\ \ 
$(c^{1}_{16}
+c^{3}_{16}
)\mathrm{i}$,
$ -c_{16}^{1}$,
$ c_{16}^{1}$,
$0$;\ \ 
$ -1$,
$ -1$,
$ c^{1}_{16}
+c^{3}_{16}
$;\ \ 
$ -1$,
$ -c^{1}_{16}
-c^{3}_{16}
$;\ \ 
$0$)

Prime. 

\vskip 1ex 
\color{grey}

\noindent(5,5). $7_{\frac{27}{4},27.31}^{32,499}$ \irep{96}:\ \ 
$d_i$ = ($1.0$,
$1.0$,
$2.414$,
$2.414$,
$-1.847$,
$-1.847$,
$-2.613$) 

\vskip 0.7ex
\hangindent=3em \hangafter=1
$D^2= 27.313 = 
16+8\sqrt{2}$

\vskip 0.7ex
\hangindent=3em \hangafter=1
$T = ( 0,
\frac{1}{2},
\frac{1}{4},
\frac{3}{4},
\frac{1}{32},
\frac{1}{32},
\frac{21}{32} )
$,

\vskip 0.7ex
\hangindent=3em \hangafter=1
$S$ = ($ 1$,
$ 1$,
$ 1+\sqrt{2}$,
$ 1+\sqrt{2}$,
$ -c_{16}^{1}$,
$ -c_{16}^{1}$,
$ -c^{1}_{16}
-c^{3}_{16}
$;\ \ 
$ 1$,
$ 1+\sqrt{2}$,
$ 1+\sqrt{2}$,
$ c_{16}^{1}$,
$ c_{16}^{1}$,
$ c^{1}_{16}
+c^{3}_{16}
$;\ \ 
$ -1$,
$ -1$,
$ c_{16}^{1}$,
$ c_{16}^{1}$,
$ -c^{1}_{16}
-c^{3}_{16}
$;\ \ 
$ -1$,
$ -c_{16}^{1}$,
$ -c_{16}^{1}$,
$ c^{1}_{16}
+c^{3}_{16}
$;\ \ 
$(-c^{1}_{16}
-c^{3}_{16}
)\mathrm{i}$,
$(c^{1}_{16}
+c^{3}_{16}
)\mathrm{i}$,
$0$;\ \ 
$(-c^{1}_{16}
-c^{3}_{16}
)\mathrm{i}$,
$0$;\ \ 
$0$)

Prime. 

Pseudo-unitary $\sim$  
$7_{\frac{11}{4},27.31}^{32,418}$

\vskip 1ex 
\color{grey}

\noindent(5,6). $7_{\frac{21}{4},27.31}^{32,104}$ \irep{96}:\ \ 
$d_i$ = ($1.0$,
$1.0$,
$2.414$,
$2.414$,
$-1.847$,
$-1.847$,
$-2.613$) 

\vskip 0.7ex
\hangindent=3em \hangafter=1
$D^2= 27.313 = 
16+8\sqrt{2}$

\vskip 0.7ex
\hangindent=3em \hangafter=1
$T = ( 0,
\frac{1}{2},
\frac{1}{4},
\frac{3}{4},
\frac{15}{32},
\frac{15}{32},
\frac{27}{32} )
$,

\vskip 0.7ex
\hangindent=3em \hangafter=1
$S$ = ($ 1$,
$ 1$,
$ 1+\sqrt{2}$,
$ 1+\sqrt{2}$,
$ -c_{16}^{1}$,
$ -c_{16}^{1}$,
$ -c^{1}_{16}
-c^{3}_{16}
$;\ \ 
$ 1$,
$ 1+\sqrt{2}$,
$ 1+\sqrt{2}$,
$ c_{16}^{1}$,
$ c_{16}^{1}$,
$ c^{1}_{16}
+c^{3}_{16}
$;\ \ 
$ -1$,
$ -1$,
$ -c_{16}^{1}$,
$ -c_{16}^{1}$,
$ c^{1}_{16}
+c^{3}_{16}
$;\ \ 
$ -1$,
$ c_{16}^{1}$,
$ c_{16}^{1}$,
$ -c^{1}_{16}
-c^{3}_{16}
$;\ \ 
$(c^{1}_{16}
+c^{3}_{16}
)\mathrm{i}$,
$(-c^{1}_{16}
-c^{3}_{16}
)\mathrm{i}$,
$0$;\ \ 
$(c^{1}_{16}
+c^{3}_{16}
)\mathrm{i}$,
$0$;\ \ 
$0$)

Prime. 

Pseudo-unitary $\sim$  
$7_{\frac{5}{4},27.31}^{32,225}$

\vskip 1ex 
\color{grey}

\noindent(5,7). $7_{\frac{11}{4},27.31}^{32,522}$ \irep{96}:\ \ 
$d_i$ = ($1.0$,
$1.0$,
$2.414$,
$2.414$,
$-1.847$,
$-1.847$,
$-2.613$) 

\vskip 0.7ex
\hangindent=3em \hangafter=1
$D^2= 27.313 = 
16+8\sqrt{2}$

\vskip 0.7ex
\hangindent=3em \hangafter=1
$T = ( 0,
\frac{1}{2},
\frac{1}{4},
\frac{3}{4},
\frac{17}{32},
\frac{17}{32},
\frac{5}{32} )
$,

\vskip 0.7ex
\hangindent=3em \hangafter=1
$S$ = ($ 1$,
$ 1$,
$ 1+\sqrt{2}$,
$ 1+\sqrt{2}$,
$ -c_{16}^{1}$,
$ -c_{16}^{1}$,
$ -c^{1}_{16}
-c^{3}_{16}
$;\ \ 
$ 1$,
$ 1+\sqrt{2}$,
$ 1+\sqrt{2}$,
$ c_{16}^{1}$,
$ c_{16}^{1}$,
$ c^{1}_{16}
+c^{3}_{16}
$;\ \ 
$ -1$,
$ -1$,
$ c_{16}^{1}$,
$ c_{16}^{1}$,
$ -c^{1}_{16}
-c^{3}_{16}
$;\ \ 
$ -1$,
$ -c_{16}^{1}$,
$ -c_{16}^{1}$,
$ c^{1}_{16}
+c^{3}_{16}
$;\ \ 
$(-c^{1}_{16}
-c^{3}_{16}
)\mathrm{i}$,
$(c^{1}_{16}
+c^{3}_{16}
)\mathrm{i}$,
$0$;\ \ 
$(-c^{1}_{16}
-c^{3}_{16}
)\mathrm{i}$,
$0$;\ \ 
$0$)

Prime. 

Pseudo-unitary $\sim$  
$7_{\frac{27}{4},27.31}^{32,396}$

\vskip 1ex 
\color{grey}

\noindent(5,8). $7_{\frac{5}{4},27.31}^{32,214}$ \irep{96}:\ \ 
$d_i$ = ($1.0$,
$1.0$,
$2.414$,
$2.414$,
$-1.847$,
$-1.847$,
$-2.613$) 

\vskip 0.7ex
\hangindent=3em \hangafter=1
$D^2= 27.313 = 
16+8\sqrt{2}$

\vskip 0.7ex
\hangindent=3em \hangafter=1
$T = ( 0,
\frac{1}{2},
\frac{1}{4},
\frac{3}{4},
\frac{31}{32},
\frac{31}{32},
\frac{11}{32} )
$,

\vskip 0.7ex
\hangindent=3em \hangafter=1
$S$ = ($ 1$,
$ 1$,
$ 1+\sqrt{2}$,
$ 1+\sqrt{2}$,
$ -c_{16}^{1}$,
$ -c_{16}^{1}$,
$ -c^{1}_{16}
-c^{3}_{16}
$;\ \ 
$ 1$,
$ 1+\sqrt{2}$,
$ 1+\sqrt{2}$,
$ c_{16}^{1}$,
$ c_{16}^{1}$,
$ c^{1}_{16}
+c^{3}_{16}
$;\ \ 
$ -1$,
$ -1$,
$ -c_{16}^{1}$,
$ -c_{16}^{1}$,
$ c^{1}_{16}
+c^{3}_{16}
$;\ \ 
$ -1$,
$ c_{16}^{1}$,
$ c_{16}^{1}$,
$ -c^{1}_{16}
-c^{3}_{16}
$;\ \ 
$(c^{1}_{16}
+c^{3}_{16}
)\mathrm{i}$,
$(-c^{1}_{16}
-c^{3}_{16}
)\mathrm{i}$,
$0$;\ \ 
$(c^{1}_{16}
+c^{3}_{16}
)\mathrm{i}$,
$0$;\ \ 
$0$)

Prime. 

Pseudo-unitary $\sim$  
$7_{\frac{21}{4},27.31}^{32,114}$

\vskip 1ex 
\color{grey}

\noindent(5,9). $7_{\frac{7}{4},4.686}^{32,462}$ \irep{96}:\ \ 
$d_i$ = ($1.0$,
$0.765$,
$0.765$,
$1.0$,
$-0.414$,
$-0.414$,
$-1.82$) 

\vskip 0.7ex
\hangindent=3em \hangafter=1
$D^2= 4.686 = 
16-8\sqrt{2}$

\vskip 0.7ex
\hangindent=3em \hangafter=1
$T = ( 0,
\frac{5}{32},
\frac{5}{32},
\frac{1}{2},
\frac{1}{4},
\frac{3}{4},
\frac{9}{32} )
$,

\vskip 0.7ex
\hangindent=3em \hangafter=1
$S$ = ($ 1$,
$ c_{16}^{3}$,
$ c_{16}^{3}$,
$ 1$,
$ 1-\sqrt{2}$,
$ 1-\sqrt{2}$,
$ -c^{1}_{16}
+c^{3}_{16}
$;\ \ 
$(-c^{1}_{16}
+c^{3}_{16}
)\mathrm{i}$,
$(c^{1}_{16}
-c^{3}_{16}
)\mathrm{i}$,
$ -c_{16}^{3}$,
$ -c_{16}^{3}$,
$ c_{16}^{3}$,
$0$;\ \ 
$(-c^{1}_{16}
+c^{3}_{16}
)\mathrm{i}$,
$ -c_{16}^{3}$,
$ -c_{16}^{3}$,
$ c_{16}^{3}$,
$0$;\ \ 
$ 1$,
$ 1-\sqrt{2}$,
$ 1-\sqrt{2}$,
$ c^{1}_{16}
-c^{3}_{16}
$;\ \ 
$ -1$,
$ -1$,
$ -c^{1}_{16}
+c^{3}_{16}
$;\ \ 
$ -1$,
$ c^{1}_{16}
-c^{3}_{16}
$;\ \ 
$0$)

Prime. 

Not pseudo-unitary. 

\vskip 1ex 
\color{grey}

\noindent(5,10). $7_{\frac{9}{4},4.686}^{32,284}$ \irep{96}:\ \ 
$d_i$ = ($1.0$,
$0.765$,
$0.765$,
$1.0$,
$-0.414$,
$-0.414$,
$-1.82$) 

\vskip 0.7ex
\hangindent=3em \hangafter=1
$D^2= 4.686 = 
16-8\sqrt{2}$

\vskip 0.7ex
\hangindent=3em \hangafter=1
$T = ( 0,
\frac{11}{32},
\frac{11}{32},
\frac{1}{2},
\frac{1}{4},
\frac{3}{4},
\frac{7}{32} )
$,

\vskip 0.7ex
\hangindent=3em \hangafter=1
$S$ = ($ 1$,
$ c_{16}^{3}$,
$ c_{16}^{3}$,
$ 1$,
$ 1-\sqrt{2}$,
$ 1-\sqrt{2}$,
$ -c^{1}_{16}
+c^{3}_{16}
$;\ \ 
$(c^{1}_{16}
-c^{3}_{16}
)\mathrm{i}$,
$(-c^{1}_{16}
+c^{3}_{16}
)\mathrm{i}$,
$ -c_{16}^{3}$,
$ c_{16}^{3}$,
$ -c_{16}^{3}$,
$0$;\ \ 
$(c^{1}_{16}
-c^{3}_{16}
)\mathrm{i}$,
$ -c_{16}^{3}$,
$ c_{16}^{3}$,
$ -c_{16}^{3}$,
$0$;\ \ 
$ 1$,
$ 1-\sqrt{2}$,
$ 1-\sqrt{2}$,
$ c^{1}_{16}
-c^{3}_{16}
$;\ \ 
$ -1$,
$ -1$,
$ c^{1}_{16}
-c^{3}_{16}
$;\ \ 
$ -1$,
$ -c^{1}_{16}
+c^{3}_{16}
$;\ \ 
$0$)

Prime. 

Not pseudo-unitary. 

\vskip 1ex 
\color{grey}

\noindent(5,11). $7_{\frac{31}{4},4.686}^{32,144}$ \irep{96}:\ \ 
$d_i$ = ($1.0$,
$1.0$,
$1.82$,
$-0.414$,
$-0.414$,
$-0.765$,
$-0.765$) 

\vskip 0.7ex
\hangindent=3em \hangafter=1
$D^2= 4.686 = 
16-8\sqrt{2}$

\vskip 0.7ex
\hangindent=3em \hangafter=1
$T = ( 0,
\frac{1}{2},
\frac{1}{32},
\frac{1}{4},
\frac{3}{4},
\frac{29}{32},
\frac{29}{32} )
$,

\vskip 0.7ex
\hangindent=3em \hangafter=1
$S$ = ($ 1$,
$ 1$,
$ c^{1}_{16}
-c^{3}_{16}
$,
$ 1-\sqrt{2}$,
$ 1-\sqrt{2}$,
$ -c_{16}^{3}$,
$ -c_{16}^{3}$;\ \ 
$ 1$,
$ -c^{1}_{16}
+c^{3}_{16}
$,
$ 1-\sqrt{2}$,
$ 1-\sqrt{2}$,
$ c_{16}^{3}$,
$ c_{16}^{3}$;\ \ 
$0$,
$ c^{1}_{16}
-c^{3}_{16}
$,
$ -c^{1}_{16}
+c^{3}_{16}
$,
$0$,
$0$;\ \ 
$ -1$,
$ -1$,
$ c_{16}^{3}$,
$ c_{16}^{3}$;\ \ 
$ -1$,
$ -c_{16}^{3}$,
$ -c_{16}^{3}$;\ \ 
$(c^{1}_{16}
-c^{3}_{16}
)\mathrm{i}$,
$(-c^{1}_{16}
+c^{3}_{16}
)\mathrm{i}$;\ \ 
$(c^{1}_{16}
-c^{3}_{16}
)\mathrm{i}$)

Prime. 

Not pseudo-unitary. 

\vskip 1ex 
\color{grey}

\noindent(5,12). $7_{\frac{17}{4},4.686}^{32,491}$ \irep{96}:\ \ 
$d_i$ = ($1.0$,
$1.0$,
$1.82$,
$-0.414$,
$-0.414$,
$-0.765$,
$-0.765$) 

\vskip 0.7ex
\hangindent=3em \hangafter=1
$D^2= 4.686 = 
16-8\sqrt{2}$

\vskip 0.7ex
\hangindent=3em \hangafter=1
$T = ( 0,
\frac{1}{2},
\frac{15}{32},
\frac{1}{4},
\frac{3}{4},
\frac{19}{32},
\frac{19}{32} )
$,

\vskip 0.7ex
\hangindent=3em \hangafter=1
$S$ = ($ 1$,
$ 1$,
$ c^{1}_{16}
-c^{3}_{16}
$,
$ 1-\sqrt{2}$,
$ 1-\sqrt{2}$,
$ -c_{16}^{3}$,
$ -c_{16}^{3}$;\ \ 
$ 1$,
$ -c^{1}_{16}
+c^{3}_{16}
$,
$ 1-\sqrt{2}$,
$ 1-\sqrt{2}$,
$ c_{16}^{3}$,
$ c_{16}^{3}$;\ \ 
$0$,
$ -c^{1}_{16}
+c^{3}_{16}
$,
$ c^{1}_{16}
-c^{3}_{16}
$,
$0$,
$0$;\ \ 
$ -1$,
$ -1$,
$ -c_{16}^{3}$,
$ -c_{16}^{3}$;\ \ 
$ -1$,
$ c_{16}^{3}$,
$ c_{16}^{3}$;\ \ 
$(-c^{1}_{16}
+c^{3}_{16}
)\mathrm{i}$,
$(c^{1}_{16}
-c^{3}_{16}
)\mathrm{i}$;\ \ 
$(-c^{1}_{16}
+c^{3}_{16}
)\mathrm{i}$)

Prime. 

Not pseudo-unitary. 

\vskip 1ex 
\color{grey}

\noindent(5,13). $7_{\frac{15}{4},4.686}^{32,136}$ \irep{96}:\ \ 
$d_i$ = ($1.0$,
$1.0$,
$1.82$,
$-0.414$,
$-0.414$,
$-0.765$,
$-0.765$) 

\vskip 0.7ex
\hangindent=3em \hangafter=1
$D^2= 4.686 = 
16-8\sqrt{2}$

\vskip 0.7ex
\hangindent=3em \hangafter=1
$T = ( 0,
\frac{1}{2},
\frac{17}{32},
\frac{1}{4},
\frac{3}{4},
\frac{13}{32},
\frac{13}{32} )
$,

\vskip 0.7ex
\hangindent=3em \hangafter=1
$S$ = ($ 1$,
$ 1$,
$ c^{1}_{16}
-c^{3}_{16}
$,
$ 1-\sqrt{2}$,
$ 1-\sqrt{2}$,
$ -c_{16}^{3}$,
$ -c_{16}^{3}$;\ \ 
$ 1$,
$ -c^{1}_{16}
+c^{3}_{16}
$,
$ 1-\sqrt{2}$,
$ 1-\sqrt{2}$,
$ c_{16}^{3}$,
$ c_{16}^{3}$;\ \ 
$0$,
$ c^{1}_{16}
-c^{3}_{16}
$,
$ -c^{1}_{16}
+c^{3}_{16}
$,
$0$,
$0$;\ \ 
$ -1$,
$ -1$,
$ c_{16}^{3}$,
$ c_{16}^{3}$;\ \ 
$ -1$,
$ -c_{16}^{3}$,
$ -c_{16}^{3}$;\ \ 
$(c^{1}_{16}
-c^{3}_{16}
)\mathrm{i}$,
$(-c^{1}_{16}
+c^{3}_{16}
)\mathrm{i}$;\ \ 
$(c^{1}_{16}
-c^{3}_{16}
)\mathrm{i}$)

Prime. 

Not pseudo-unitary. 

\vskip 1ex 
\color{grey}

\noindent(5,14). $7_{\frac{1}{4},4.686}^{32,812}$ \irep{96}:\ \ 
$d_i$ = ($1.0$,
$1.0$,
$1.82$,
$-0.414$,
$-0.414$,
$-0.765$,
$-0.765$) 

\vskip 0.7ex
\hangindent=3em \hangafter=1
$D^2= 4.686 = 
16-8\sqrt{2}$

\vskip 0.7ex
\hangindent=3em \hangafter=1
$T = ( 0,
\frac{1}{2},
\frac{31}{32},
\frac{1}{4},
\frac{3}{4},
\frac{3}{32},
\frac{3}{32} )
$,

\vskip 0.7ex
\hangindent=3em \hangafter=1
$S$ = ($ 1$,
$ 1$,
$ c^{1}_{16}
-c^{3}_{16}
$,
$ 1-\sqrt{2}$,
$ 1-\sqrt{2}$,
$ -c_{16}^{3}$,
$ -c_{16}^{3}$;\ \ 
$ 1$,
$ -c^{1}_{16}
+c^{3}_{16}
$,
$ 1-\sqrt{2}$,
$ 1-\sqrt{2}$,
$ c_{16}^{3}$,
$ c_{16}^{3}$;\ \ 
$0$,
$ -c^{1}_{16}
+c^{3}_{16}
$,
$ c^{1}_{16}
-c^{3}_{16}
$,
$0$,
$0$;\ \ 
$ -1$,
$ -1$,
$ -c_{16}^{3}$,
$ -c_{16}^{3}$;\ \ 
$ -1$,
$ c_{16}^{3}$,
$ c_{16}^{3}$;\ \ 
$(-c^{1}_{16}
+c^{3}_{16}
)\mathrm{i}$,
$(c^{1}_{16}
-c^{3}_{16}
)\mathrm{i}$;\ \ 
$(-c^{1}_{16}
+c^{3}_{16}
)\mathrm{i}$)

Prime. 

Not pseudo-unitary. 

\vskip 1ex 
\color{grey}

\noindent(5,15). $7_{\frac{23}{4},4.686}^{32,415}$ \irep{96}:\ \ 
$d_i$ = ($1.0$,
$0.765$,
$0.765$,
$1.0$,
$-0.414$,
$-0.414$,
$-1.82$) 

\vskip 0.7ex
\hangindent=3em \hangafter=1
$D^2= 4.686 = 
16-8\sqrt{2}$

\vskip 0.7ex
\hangindent=3em \hangafter=1
$T = ( 0,
\frac{21}{32},
\frac{21}{32},
\frac{1}{2},
\frac{1}{4},
\frac{3}{4},
\frac{25}{32} )
$,

\vskip 0.7ex
\hangindent=3em \hangafter=1
$S$ = ($ 1$,
$ c_{16}^{3}$,
$ c_{16}^{3}$,
$ 1$,
$ 1-\sqrt{2}$,
$ 1-\sqrt{2}$,
$ -c^{1}_{16}
+c^{3}_{16}
$;\ \ 
$(-c^{1}_{16}
+c^{3}_{16}
)\mathrm{i}$,
$(c^{1}_{16}
-c^{3}_{16}
)\mathrm{i}$,
$ -c_{16}^{3}$,
$ -c_{16}^{3}$,
$ c_{16}^{3}$,
$0$;\ \ 
$(-c^{1}_{16}
+c^{3}_{16}
)\mathrm{i}$,
$ -c_{16}^{3}$,
$ -c_{16}^{3}$,
$ c_{16}^{3}$,
$0$;\ \ 
$ 1$,
$ 1-\sqrt{2}$,
$ 1-\sqrt{2}$,
$ c^{1}_{16}
-c^{3}_{16}
$;\ \ 
$ -1$,
$ -1$,
$ -c^{1}_{16}
+c^{3}_{16}
$;\ \ 
$ -1$,
$ c^{1}_{16}
-c^{3}_{16}
$;\ \ 
$0$)

Prime. 

Not pseudo-unitary. 

\vskip 1ex 
\color{grey}

\noindent(5,16). $7_{\frac{25}{4},4.686}^{32,116}$ \irep{96}:\ \ 
$d_i$ = ($1.0$,
$0.765$,
$0.765$,
$1.0$,
$-0.414$,
$-0.414$,
$-1.82$) 

\vskip 0.7ex
\hangindent=3em \hangafter=1
$D^2= 4.686 = 
16-8\sqrt{2}$

\vskip 0.7ex
\hangindent=3em \hangafter=1
$T = ( 0,
\frac{27}{32},
\frac{27}{32},
\frac{1}{2},
\frac{1}{4},
\frac{3}{4},
\frac{23}{32} )
$,

\vskip 0.7ex
\hangindent=3em \hangafter=1
$S$ = ($ 1$,
$ c_{16}^{3}$,
$ c_{16}^{3}$,
$ 1$,
$ 1-\sqrt{2}$,
$ 1-\sqrt{2}$,
$ -c^{1}_{16}
+c^{3}_{16}
$;\ \ 
$(c^{1}_{16}
-c^{3}_{16}
)\mathrm{i}$,
$(-c^{1}_{16}
+c^{3}_{16}
)\mathrm{i}$,
$ -c_{16}^{3}$,
$ c_{16}^{3}$,
$ -c_{16}^{3}$,
$0$;\ \ 
$(c^{1}_{16}
-c^{3}_{16}
)\mathrm{i}$,
$ -c_{16}^{3}$,
$ c_{16}^{3}$,
$ -c_{16}^{3}$,
$0$;\ \ 
$ 1$,
$ 1-\sqrt{2}$,
$ 1-\sqrt{2}$,
$ c^{1}_{16}
-c^{3}_{16}
$;\ \ 
$ -1$,
$ -1$,
$ c^{1}_{16}
-c^{3}_{16}
$;\ \ 
$ -1$,
$ -c^{1}_{16}
+c^{3}_{16}
$;\ \ 
$0$)

Prime. 

Not pseudo-unitary. 

\vskip 1ex 
\black

\noindent(6,1). $7_{2,28.}^{56,139}$ \irep{99}:\ \ 
$d_i$ = ($1.0$,
$1.0$,
$2.0$,
$2.0$,
$2.0$,
$2.645$,
$2.645$) 

\vskip 0.7ex
\hangindent=3em \hangafter=1
$D^2= 28.0 = 
28$

\vskip 0.7ex
\hangindent=3em \hangafter=1
$T = ( 0,
0,
\frac{1}{7},
\frac{2}{7},
\frac{4}{7},
\frac{1}{8},
\frac{5}{8} )
$,

\vskip 0.7ex
\hangindent=3em \hangafter=1
$S$ = ($ 1$,
$ 1$,
$ 2$,
$ 2$,
$ 2$,
$ \sqrt{7}$,
$ \sqrt{7}$;\ \ 
$ 1$,
$ 2$,
$ 2$,
$ 2$,
$ -\sqrt{7}$,
$ -\sqrt{7}$;\ \ 
$ 2c_{7}^{2}$,
$ 2c_{7}^{1}$,
$ 2c_{7}^{3}$,
$0$,
$0$;\ \ 
$ 2c_{7}^{3}$,
$ 2c_{7}^{2}$,
$0$,
$0$;\ \ 
$ 2c_{7}^{1}$,
$0$,
$0$;\ \ 
$ \sqrt{7}$,
$ -\sqrt{7}$;\ \ 
$ \sqrt{7}$)

Prime. 

\vskip 1ex 
\color{grey}

\noindent(6,2). $7_{6,28.}^{56,193}$ \irep{99}:\ \ 
$d_i$ = ($1.0$,
$1.0$,
$2.0$,
$2.0$,
$2.0$,
$2.645$,
$2.645$) 

\vskip 0.7ex
\hangindent=3em \hangafter=1
$D^2= 28.0 = 
28$

\vskip 0.7ex
\hangindent=3em \hangafter=1
$T = ( 0,
0,
\frac{3}{7},
\frac{5}{7},
\frac{6}{7},
\frac{3}{8},
\frac{7}{8} )
$,

\vskip 0.7ex
\hangindent=3em \hangafter=1
$S$ = ($ 1$,
$ 1$,
$ 2$,
$ 2$,
$ 2$,
$ \sqrt{7}$,
$ \sqrt{7}$;\ \ 
$ 1$,
$ 2$,
$ 2$,
$ 2$,
$ -\sqrt{7}$,
$ -\sqrt{7}$;\ \ 
$ 2c_{7}^{1}$,
$ 2c_{7}^{2}$,
$ 2c_{7}^{3}$,
$0$,
$0$;\ \ 
$ 2c_{7}^{3}$,
$ 2c_{7}^{1}$,
$0$,
$0$;\ \ 
$ 2c_{7}^{2}$,
$0$,
$0$;\ \ 
$ \sqrt{7}$,
$ -\sqrt{7}$;\ \ 
$ \sqrt{7}$)

Prime. 

\vskip 1ex 
\color{grey}

\noindent(6,3). $7_{2,28.}^{56,221}$ \irep{99}:\ \ 
$d_i$ = ($1.0$,
$1.0$,
$2.0$,
$2.0$,
$2.0$,
$-2.645$,
$-2.645$) 

\vskip 0.7ex
\hangindent=3em \hangafter=1
$D^2= 28.0 = 
28$

\vskip 0.7ex
\hangindent=3em \hangafter=1
$T = ( 0,
0,
\frac{1}{7},
\frac{2}{7},
\frac{4}{7},
\frac{3}{8},
\frac{7}{8} )
$,

\vskip 0.7ex
\hangindent=3em \hangafter=1
$S$ = ($ 1$,
$ 1$,
$ 2$,
$ 2$,
$ 2$,
$ -\sqrt{7}$,
$ -\sqrt{7}$;\ \ 
$ 1$,
$ 2$,
$ 2$,
$ 2$,
$ \sqrt{7}$,
$ \sqrt{7}$;\ \ 
$ 2c_{7}^{2}$,
$ 2c_{7}^{1}$,
$ 2c_{7}^{3}$,
$0$,
$0$;\ \ 
$ 2c_{7}^{3}$,
$ 2c_{7}^{2}$,
$0$,
$0$;\ \ 
$ 2c_{7}^{1}$,
$0$,
$0$;\ \ 
$ -\sqrt{7}$,
$ \sqrt{7}$;\ \ 
$ -\sqrt{7}$)

Prime. 

Pseudo-unitary $\sim$  
$7_{2,28.}^{56,680}$

\vskip 1ex 
\color{grey}

\noindent(6,4). $7_{6,28.}^{56,110}$ \irep{99}:\ \ 
$d_i$ = ($1.0$,
$1.0$,
$2.0$,
$2.0$,
$2.0$,
$-2.645$,
$-2.645$) 

\vskip 0.7ex
\hangindent=3em \hangafter=1
$D^2= 28.0 = 
28$

\vskip 0.7ex
\hangindent=3em \hangafter=1
$T = ( 0,
0,
\frac{3}{7},
\frac{5}{7},
\frac{6}{7},
\frac{1}{8},
\frac{5}{8} )
$,

\vskip 0.7ex
\hangindent=3em \hangafter=1
$S$ = ($ 1$,
$ 1$,
$ 2$,
$ 2$,
$ 2$,
$ -\sqrt{7}$,
$ -\sqrt{7}$;\ \ 
$ 1$,
$ 2$,
$ 2$,
$ 2$,
$ \sqrt{7}$,
$ \sqrt{7}$;\ \ 
$ 2c_{7}^{1}$,
$ 2c_{7}^{2}$,
$ 2c_{7}^{3}$,
$0$,
$0$;\ \ 
$ 2c_{7}^{3}$,
$ 2c_{7}^{1}$,
$0$,
$0$;\ \ 
$ 2c_{7}^{2}$,
$0$,
$0$;\ \ 
$ -\sqrt{7}$,
$ \sqrt{7}$;\ \ 
$ -\sqrt{7}$)

Prime. 

Pseudo-unitary $\sim$  
$7_{6,28.}^{56,609}$

\vskip 1ex 
\black

\noindent(7,1). $7_{2,28.}^{56,680}$ \irep{99}:\ \ 
$d_i$ = ($1.0$,
$1.0$,
$2.0$,
$2.0$,
$2.0$,
$2.645$,
$2.645$) 

\vskip 0.7ex
\hangindent=3em \hangafter=1
$D^2= 28.0 = 
28$

\vskip 0.7ex
\hangindent=3em \hangafter=1
$T = ( 0,
0,
\frac{1}{7},
\frac{2}{7},
\frac{4}{7},
\frac{3}{8},
\frac{7}{8} )
$,

\vskip 0.7ex
\hangindent=3em \hangafter=1
$S$ = ($ 1$,
$ 1$,
$ 2$,
$ 2$,
$ 2$,
$ \sqrt{7}$,
$ \sqrt{7}$;\ \ 
$ 1$,
$ 2$,
$ 2$,
$ 2$,
$ -\sqrt{7}$,
$ -\sqrt{7}$;\ \ 
$ 2c_{7}^{2}$,
$ 2c_{7}^{1}$,
$ 2c_{7}^{3}$,
$0$,
$0$;\ \ 
$ 2c_{7}^{3}$,
$ 2c_{7}^{2}$,
$0$,
$0$;\ \ 
$ 2c_{7}^{1}$,
$0$,
$0$;\ \ 
$ -\sqrt{7}$,
$ \sqrt{7}$;\ \ 
$ -\sqrt{7}$)

Prime. 

\vskip 1ex 
\color{grey}

\noindent(7,2). $7_{6,28.}^{56,609}$ \irep{99}:\ \ 
$d_i$ = ($1.0$,
$1.0$,
$2.0$,
$2.0$,
$2.0$,
$2.645$,
$2.645$) 

\vskip 0.7ex
\hangindent=3em \hangafter=1
$D^2= 28.0 = 
28$

\vskip 0.7ex
\hangindent=3em \hangafter=1
$T = ( 0,
0,
\frac{3}{7},
\frac{5}{7},
\frac{6}{7},
\frac{1}{8},
\frac{5}{8} )
$,

\vskip 0.7ex
\hangindent=3em \hangafter=1
$S$ = ($ 1$,
$ 1$,
$ 2$,
$ 2$,
$ 2$,
$ \sqrt{7}$,
$ \sqrt{7}$;\ \ 
$ 1$,
$ 2$,
$ 2$,
$ 2$,
$ -\sqrt{7}$,
$ -\sqrt{7}$;\ \ 
$ 2c_{7}^{1}$,
$ 2c_{7}^{2}$,
$ 2c_{7}^{3}$,
$0$,
$0$;\ \ 
$ 2c_{7}^{3}$,
$ 2c_{7}^{1}$,
$0$,
$0$;\ \ 
$ 2c_{7}^{2}$,
$0$,
$0$;\ \ 
$ -\sqrt{7}$,
$ \sqrt{7}$;\ \ 
$ -\sqrt{7}$)

Prime. 

\vskip 1ex 
\color{grey}

\noindent(7,3). $7_{2,28.}^{56,894}$ \irep{99}:\ \ 
$d_i$ = ($1.0$,
$1.0$,
$2.0$,
$2.0$,
$2.0$,
$-2.645$,
$-2.645$) 

\vskip 0.7ex
\hangindent=3em \hangafter=1
$D^2= 28.0 = 
28$

\vskip 0.7ex
\hangindent=3em \hangafter=1
$T = ( 0,
0,
\frac{1}{7},
\frac{2}{7},
\frac{4}{7},
\frac{1}{8},
\frac{5}{8} )
$,

\vskip 0.7ex
\hangindent=3em \hangafter=1
$S$ = ($ 1$,
$ 1$,
$ 2$,
$ 2$,
$ 2$,
$ -\sqrt{7}$,
$ -\sqrt{7}$;\ \ 
$ 1$,
$ 2$,
$ 2$,
$ 2$,
$ \sqrt{7}$,
$ \sqrt{7}$;\ \ 
$ 2c_{7}^{2}$,
$ 2c_{7}^{1}$,
$ 2c_{7}^{3}$,
$0$,
$0$;\ \ 
$ 2c_{7}^{3}$,
$ 2c_{7}^{2}$,
$0$,
$0$;\ \ 
$ 2c_{7}^{1}$,
$0$,
$0$;\ \ 
$ \sqrt{7}$,
$ -\sqrt{7}$;\ \ 
$ \sqrt{7}$)

Prime. 

Pseudo-unitary $\sim$  
$7_{2,28.}^{56,139}$

\vskip 1ex 
\color{grey}

\noindent(7,4). $7_{6,28.}^{56,217}$ \irep{99}:\ \ 
$d_i$ = ($1.0$,
$1.0$,
$2.0$,
$2.0$,
$2.0$,
$-2.645$,
$-2.645$) 

\vskip 0.7ex
\hangindent=3em \hangafter=1
$D^2= 28.0 = 
28$

\vskip 0.7ex
\hangindent=3em \hangafter=1
$T = ( 0,
0,
\frac{3}{7},
\frac{5}{7},
\frac{6}{7},
\frac{3}{8},
\frac{7}{8} )
$,

\vskip 0.7ex
\hangindent=3em \hangafter=1
$S$ = ($ 1$,
$ 1$,
$ 2$,
$ 2$,
$ 2$,
$ -\sqrt{7}$,
$ -\sqrt{7}$;\ \ 
$ 1$,
$ 2$,
$ 2$,
$ 2$,
$ \sqrt{7}$,
$ \sqrt{7}$;\ \ 
$ 2c_{7}^{1}$,
$ 2c_{7}^{2}$,
$ 2c_{7}^{3}$,
$0$,
$0$;\ \ 
$ 2c_{7}^{3}$,
$ 2c_{7}^{1}$,
$0$,
$0$;\ \ 
$ 2c_{7}^{2}$,
$0$,
$0$;\ \ 
$ \sqrt{7}$,
$ -\sqrt{7}$;\ \ 
$ \sqrt{7}$)

Prime. 

Pseudo-unitary $\sim$  
$7_{6,28.}^{56,193}$

\vskip 1ex 
\black

\noindent(8,1). $7_{\frac{32}{5},86.75}^{15,205}$ \irep{80}:\ \ 
$d_i$ = ($1.0$,
$1.956$,
$2.827$,
$3.574$,
$4.165$,
$4.574$,
$4.783$) 

\vskip 0.7ex
\hangindent=3em \hangafter=1
$D^2= 86.750 = 
30+15c^{1}_{15}
+15c^{2}_{15}
+15c^{3}_{15}
$

\vskip 0.7ex
\hangindent=3em \hangafter=1
$T = ( 0,
\frac{1}{5},
\frac{13}{15},
0,
\frac{3}{5},
\frac{2}{3},
\frac{1}{5} )
$,

\vskip 0.7ex
\hangindent=3em \hangafter=1
$S$ = ($ 1$,
$ -c_{15}^{7}$,
$ \xi_{15}^{3}$,
$ \xi_{15}^{11}$,
$ \xi_{15}^{5}$,
$ \xi_{15}^{9}$,
$ \xi_{15}^{7}$;\ \ 
$ -\xi_{15}^{11}$,
$ \xi_{15}^{9}$,
$ -\xi_{15}^{7}$,
$ \xi_{15}^{5}$,
$ -\xi_{15}^{3}$,
$ 1$;\ \ 
$ \xi_{15}^{9}$,
$ \xi_{15}^{3}$,
$0$,
$ -\xi_{15}^{3}$,
$ -\xi_{15}^{9}$;\ \ 
$ 1$,
$ -\xi_{15}^{5}$,
$ \xi_{15}^{9}$,
$ c_{15}^{7}$;\ \ 
$ -\xi_{15}^{5}$,
$0$,
$ \xi_{15}^{5}$;\ \ 
$ -\xi_{15}^{9}$,
$ \xi_{15}^{3}$;\ \ 
$ -\xi_{15}^{11}$)

Prime. 

\vskip 1ex 
\color{grey}

\noindent(8,2). $7_{\frac{8}{5},86.75}^{15,181}$ \irep{80}:\ \ 
$d_i$ = ($1.0$,
$1.956$,
$2.827$,
$3.574$,
$4.165$,
$4.574$,
$4.783$) 

\vskip 0.7ex
\hangindent=3em \hangafter=1
$D^2= 86.750 = 
30+15c^{1}_{15}
+15c^{2}_{15}
+15c^{3}_{15}
$

\vskip 0.7ex
\hangindent=3em \hangafter=1
$T = ( 0,
\frac{4}{5},
\frac{2}{15},
0,
\frac{2}{5},
\frac{1}{3},
\frac{4}{5} )
$,

\vskip 0.7ex
\hangindent=3em \hangafter=1
$S$ = ($ 1$,
$ -c_{15}^{7}$,
$ \xi_{15}^{3}$,
$ \xi_{15}^{11}$,
$ \xi_{15}^{5}$,
$ \xi_{15}^{9}$,
$ \xi_{15}^{7}$;\ \ 
$ -\xi_{15}^{11}$,
$ \xi_{15}^{9}$,
$ -\xi_{15}^{7}$,
$ \xi_{15}^{5}$,
$ -\xi_{15}^{3}$,
$ 1$;\ \ 
$ \xi_{15}^{9}$,
$ \xi_{15}^{3}$,
$0$,
$ -\xi_{15}^{3}$,
$ -\xi_{15}^{9}$;\ \ 
$ 1$,
$ -\xi_{15}^{5}$,
$ \xi_{15}^{9}$,
$ c_{15}^{7}$;\ \ 
$ -\xi_{15}^{5}$,
$0$,
$ \xi_{15}^{5}$;\ \ 
$ -\xi_{15}^{9}$,
$ \xi_{15}^{3}$;\ \ 
$ -\xi_{15}^{11}$)

Prime. 

\vskip 1ex 
\color{grey}

\noindent(8,3). $7_{\frac{4}{5},22.66}^{15,906}$ \irep{80}:\ \ 
$d_i$ = ($1.0$,
$0.511$,
$2.129$,
$2.338$,
$-1.445$,
$-1.827$,
$-2.445$) 

\vskip 0.7ex
\hangindent=3em \hangafter=1
$D^2= 22.667 = 
30-15  c^{1}_{15}
+15c^{2}_{15}
$

\vskip 0.7ex
\hangindent=3em \hangafter=1
$T = ( 0,
\frac{2}{5},
\frac{1}{5},
\frac{11}{15},
\frac{1}{3},
\frac{2}{5},
0 )
$,

\vskip 0.7ex
\hangindent=3em \hangafter=1
$S$ = ($ 1$,
$ 1-c^{1}_{15}
+c^{2}_{15}
$,
$ 2-c^{1}_{15}
+c^{2}_{15}
+c^{3}_{15}
$,
$ 1+c^{2}_{15}
$,
$ 1-c^{1}_{15}
-c^{3}_{15}
$,
$ -c_{15}^{1}$,
$ -c^{1}_{15}
-c^{3}_{15}
$;\ \ 
$ c^{1}_{15}
+c^{3}_{15}
$,
$ 2-c^{1}_{15}
+c^{2}_{15}
+c^{3}_{15}
$,
$ -1+c^{1}_{15}
+c^{3}_{15}
$,
$ 1+c^{2}_{15}
$,
$ 1$,
$ c_{15}^{1}$;\ \ 
$ -2+c^{1}_{15}
-c^{2}_{15}
-c^{3}_{15}
$,
$0$,
$0$,
$ 2-c^{1}_{15}
+c^{2}_{15}
+c^{3}_{15}
$,
$ -2+c^{1}_{15}
-c^{2}_{15}
-c^{3}_{15}
$;\ \ 
$ 1-c^{1}_{15}
-c^{3}_{15}
$,
$ -1-c^{2}_{15}
$,
$ 1-c^{1}_{15}
-c^{3}_{15}
$,
$ 1+c^{2}_{15}
$;\ \ 
$ -1+c^{1}_{15}
+c^{3}_{15}
$,
$ -1-c^{2}_{15}
$,
$ 1-c^{1}_{15}
-c^{3}_{15}
$;\ \ 
$ c^{1}_{15}
+c^{3}_{15}
$,
$ -1+c^{1}_{15}
-c^{2}_{15}
$;\ \ 
$ 1$)

Prime. 

Not pseudo-unitary. 

\vskip 1ex 
\color{grey}

\noindent(8,4). $7_{\frac{36}{5},22.66}^{15,347}$ \irep{80}:\ \ 
$d_i$ = ($1.0$,
$0.511$,
$2.129$,
$2.338$,
$-1.445$,
$-1.827$,
$-2.445$) 

\vskip 0.7ex
\hangindent=3em \hangafter=1
$D^2= 22.667 = 
30-15  c^{1}_{15}
+15c^{2}_{15}
$

\vskip 0.7ex
\hangindent=3em \hangafter=1
$T = ( 0,
\frac{3}{5},
\frac{4}{5},
\frac{4}{15},
\frac{2}{3},
\frac{3}{5},
0 )
$,

\vskip 0.7ex
\hangindent=3em \hangafter=1
$S$ = ($ 1$,
$ 1-c^{1}_{15}
+c^{2}_{15}
$,
$ 2-c^{1}_{15}
+c^{2}_{15}
+c^{3}_{15}
$,
$ 1+c^{2}_{15}
$,
$ 1-c^{1}_{15}
-c^{3}_{15}
$,
$ -c_{15}^{1}$,
$ -c^{1}_{15}
-c^{3}_{15}
$;\ \ 
$ c^{1}_{15}
+c^{3}_{15}
$,
$ 2-c^{1}_{15}
+c^{2}_{15}
+c^{3}_{15}
$,
$ -1+c^{1}_{15}
+c^{3}_{15}
$,
$ 1+c^{2}_{15}
$,
$ 1$,
$ c_{15}^{1}$;\ \ 
$ -2+c^{1}_{15}
-c^{2}_{15}
-c^{3}_{15}
$,
$0$,
$0$,
$ 2-c^{1}_{15}
+c^{2}_{15}
+c^{3}_{15}
$,
$ -2+c^{1}_{15}
-c^{2}_{15}
-c^{3}_{15}
$;\ \ 
$ 1-c^{1}_{15}
-c^{3}_{15}
$,
$ -1-c^{2}_{15}
$,
$ 1-c^{1}_{15}
-c^{3}_{15}
$,
$ 1+c^{2}_{15}
$;\ \ 
$ -1+c^{1}_{15}
+c^{3}_{15}
$,
$ -1-c^{2}_{15}
$,
$ 1-c^{1}_{15}
-c^{3}_{15}
$;\ \ 
$ c^{1}_{15}
+c^{3}_{15}
$,
$ -1+c^{1}_{15}
-c^{2}_{15}
$;\ \ 
$ 1$)

Prime. 

Not pseudo-unitary. 

\vskip 1ex 
\color{grey}

\noindent(8,5). $7_{\frac{28}{5},6.790}^{15,727}$ \irep{80}:\ \ 
$d_i$ = ($1.0$,
$0.279$,
$0.790$,
$1.279$,
$-0.547$,
$-1.165$,
$-1.338$) 

\vskip 0.7ex
\hangindent=3em \hangafter=1
$D^2= 6.790 = 
45-15  c^{1}_{15}
-15  c^{2}_{15}
+15c^{3}_{15}
$

\vskip 0.7ex
\hangindent=3em \hangafter=1
$T = ( 0,
0,
\frac{7}{15},
\frac{2}{3},
\frac{4}{5},
\frac{2}{5},
\frac{4}{5} )
$,

\vskip 0.7ex
\hangindent=3em \hangafter=1
$S$ = ($ 1$,
$ 1-c^{2}_{15}
+c^{3}_{15}
$,
$ 2-c^{1}_{15}
+c^{3}_{15}
$,
$ 2-c^{2}_{15}
+c^{3}_{15}
$,
$ 2-c^{1}_{15}
-c^{2}_{15}
+c^{3}_{15}
$,
$ 2-c^{1}_{15}
-c^{2}_{15}
$,
$ -c_{15}^{2}$;\ \ 
$ 1$,
$ 2-c^{1}_{15}
+c^{3}_{15}
$,
$ 2-c^{2}_{15}
+c^{3}_{15}
$,
$ c_{15}^{2}$,
$ -2+c^{1}_{15}
+c^{2}_{15}
$,
$ -2+c^{1}_{15}
+c^{2}_{15}
-c^{3}_{15}
$;\ \ 
$ 2-c^{2}_{15}
+c^{3}_{15}
$,
$ -2+c^{1}_{15}
-c^{3}_{15}
$,
$ -2+c^{2}_{15}
-c^{3}_{15}
$,
$0$,
$ 2-c^{2}_{15}
+c^{3}_{15}
$;\ \ 
$ -2+c^{2}_{15}
-c^{3}_{15}
$,
$ 2-c^{1}_{15}
+c^{3}_{15}
$,
$0$,
$ -2+c^{1}_{15}
-c^{3}_{15}
$;\ \ 
$ -1+c^{2}_{15}
-c^{3}_{15}
$,
$ 2-c^{1}_{15}
-c^{2}_{15}
$,
$ 1$;\ \ 
$ -2+c^{1}_{15}
+c^{2}_{15}
$,
$ 2-c^{1}_{15}
-c^{2}_{15}
$;\ \ 
$ -1+c^{2}_{15}
-c^{3}_{15}
$)

Prime. 

Not pseudo-unitary. 

\vskip 1ex 
\color{grey}

\noindent(8,6). $7_{\frac{12}{5},6.790}^{15,202}$ \irep{80}:\ \ 
$d_i$ = ($1.0$,
$0.279$,
$0.790$,
$1.279$,
$-0.547$,
$-1.165$,
$-1.338$) 

\vskip 0.7ex
\hangindent=3em \hangafter=1
$D^2= 6.790 = 
45-15  c^{1}_{15}
-15  c^{2}_{15}
+15c^{3}_{15}
$

\vskip 0.7ex
\hangindent=3em \hangafter=1
$T = ( 0,
0,
\frac{8}{15},
\frac{1}{3},
\frac{1}{5},
\frac{3}{5},
\frac{1}{5} )
$,

\vskip 0.7ex
\hangindent=3em \hangafter=1
$S$ = ($ 1$,
$ 1-c^{2}_{15}
+c^{3}_{15}
$,
$ 2-c^{1}_{15}
+c^{3}_{15}
$,
$ 2-c^{2}_{15}
+c^{3}_{15}
$,
$ 2-c^{1}_{15}
-c^{2}_{15}
+c^{3}_{15}
$,
$ 2-c^{1}_{15}
-c^{2}_{15}
$,
$ -c_{15}^{2}$;\ \ 
$ 1$,
$ 2-c^{1}_{15}
+c^{3}_{15}
$,
$ 2-c^{2}_{15}
+c^{3}_{15}
$,
$ c_{15}^{2}$,
$ -2+c^{1}_{15}
+c^{2}_{15}
$,
$ -2+c^{1}_{15}
+c^{2}_{15}
-c^{3}_{15}
$;\ \ 
$ 2-c^{2}_{15}
+c^{3}_{15}
$,
$ -2+c^{1}_{15}
-c^{3}_{15}
$,
$ -2+c^{2}_{15}
-c^{3}_{15}
$,
$0$,
$ 2-c^{2}_{15}
+c^{3}_{15}
$;\ \ 
$ -2+c^{2}_{15}
-c^{3}_{15}
$,
$ 2-c^{1}_{15}
+c^{3}_{15}
$,
$0$,
$ -2+c^{1}_{15}
-c^{3}_{15}
$;\ \ 
$ -1+c^{2}_{15}
-c^{3}_{15}
$,
$ 2-c^{1}_{15}
-c^{2}_{15}
$,
$ 1$;\ \ 
$ -2+c^{1}_{15}
+c^{2}_{15}
$,
$ 2-c^{1}_{15}
-c^{2}_{15}
$;\ \ 
$ -1+c^{2}_{15}
-c^{3}_{15}
$)

Prime. 

Not pseudo-unitary. 

\vskip 1ex 
\color{grey}

\noindent(8,7). $7_{\frac{4}{5},3.791}^{15,675}$ \irep{80}:\ \ 
$d_i$ = ($1.0$,
$0.209$,
$0.591$,
$0.870$,
$-0.408$,
$-0.747$,
$-0.956$) 

\vskip 0.7ex
\hangindent=3em \hangafter=1
$D^2= 3.791 = 
15+15c^{1}_{15}
-15  c^{2}_{15}
-30  c^{3}_{15}
$

\vskip 0.7ex
\hangindent=3em \hangafter=1
$T = ( 0,
\frac{2}{5},
\frac{2}{3},
\frac{1}{5},
0,
\frac{2}{5},
\frac{1}{15} )
$,

\vskip 0.7ex
\hangindent=3em \hangafter=1
$S$ = ($ 1$,
$ -c_{15}^{4}$,
$ c^{1}_{15}
-2  c^{3}_{15}
$,
$ 1+c^{1}_{15}
-c^{2}_{15}
-c^{3}_{15}
$,
$ -1+c^{1}_{15}
-2  c^{3}_{15}
$,
$ c^{1}_{15}
-c^{2}_{15}
-2  c^{3}_{15}
$,
$ 1-c^{2}_{15}
-c^{3}_{15}
$;\ \ 
$ 1-c^{1}_{15}
+2c^{3}_{15}
$,
$ -1+c^{2}_{15}
+c^{3}_{15}
$,
$ 1+c^{1}_{15}
-c^{2}_{15}
-c^{3}_{15}
$,
$ -c^{1}_{15}
+c^{2}_{15}
+2c^{3}_{15}
$,
$ 1$,
$ c^{1}_{15}
-2  c^{3}_{15}
$;\ \ 
$ -c^{1}_{15}
+2c^{3}_{15}
$,
$0$,
$ c^{1}_{15}
-2  c^{3}_{15}
$,
$ 1-c^{2}_{15}
-c^{3}_{15}
$,
$ -1+c^{2}_{15}
+c^{3}_{15}
$;\ \ 
$ -1-c^{1}_{15}
+c^{2}_{15}
+c^{3}_{15}
$,
$ -1-c^{1}_{15}
+c^{2}_{15}
+c^{3}_{15}
$,
$ 1+c^{1}_{15}
-c^{2}_{15}
-c^{3}_{15}
$,
$0$;\ \ 
$ 1$,
$ c_{15}^{4}$,
$ 1-c^{2}_{15}
-c^{3}_{15}
$;\ \ 
$ 1-c^{1}_{15}
+2c^{3}_{15}
$,
$ -c^{1}_{15}
+2c^{3}_{15}
$;\ \ 
$ c^{1}_{15}
-2  c^{3}_{15}
$)

Prime. 

Not pseudo-unitary. 

\vskip 1ex 
\color{grey}

\noindent(8,8). $7_{\frac{36}{5},3.791}^{15,543}$ \irep{80}:\ \ 
$d_i$ = ($1.0$,
$0.209$,
$0.591$,
$0.870$,
$-0.408$,
$-0.747$,
$-0.956$) 

\vskip 0.7ex
\hangindent=3em \hangafter=1
$D^2= 3.791 = 
15+15c^{1}_{15}
-15  c^{2}_{15}
-30  c^{3}_{15}
$

\vskip 0.7ex
\hangindent=3em \hangafter=1
$T = ( 0,
\frac{3}{5},
\frac{1}{3},
\frac{4}{5},
0,
\frac{3}{5},
\frac{14}{15} )
$,

\vskip 0.7ex
\hangindent=3em \hangafter=1
$S$ = ($ 1$,
$ -c_{15}^{4}$,
$ c^{1}_{15}
-2  c^{3}_{15}
$,
$ 1+c^{1}_{15}
-c^{2}_{15}
-c^{3}_{15}
$,
$ -1+c^{1}_{15}
-2  c^{3}_{15}
$,
$ c^{1}_{15}
-c^{2}_{15}
-2  c^{3}_{15}
$,
$ 1-c^{2}_{15}
-c^{3}_{15}
$;\ \ 
$ 1-c^{1}_{15}
+2c^{3}_{15}
$,
$ -1+c^{2}_{15}
+c^{3}_{15}
$,
$ 1+c^{1}_{15}
-c^{2}_{15}
-c^{3}_{15}
$,
$ -c^{1}_{15}
+c^{2}_{15}
+2c^{3}_{15}
$,
$ 1$,
$ c^{1}_{15}
-2  c^{3}_{15}
$;\ \ 
$ -c^{1}_{15}
+2c^{3}_{15}
$,
$0$,
$ c^{1}_{15}
-2  c^{3}_{15}
$,
$ 1-c^{2}_{15}
-c^{3}_{15}
$,
$ -1+c^{2}_{15}
+c^{3}_{15}
$;\ \ 
$ -1-c^{1}_{15}
+c^{2}_{15}
+c^{3}_{15}
$,
$ -1-c^{1}_{15}
+c^{2}_{15}
+c^{3}_{15}
$,
$ 1+c^{1}_{15}
-c^{2}_{15}
-c^{3}_{15}
$,
$0$;\ \ 
$ 1$,
$ c_{15}^{4}$,
$ 1-c^{2}_{15}
-c^{3}_{15}
$;\ \ 
$ 1-c^{1}_{15}
+2c^{3}_{15}
$,
$ -c^{1}_{15}
+2c^{3}_{15}
$;\ \ 
$ c^{1}_{15}
-2  c^{3}_{15}
$)

Prime. 

Not pseudo-unitary. 

\vskip 1ex 
\black

\noindent(9,1). $7_{1,93.25}^{8,230}$ \irep{54}:\ \ 
$d_i$ = ($1.0$,
$2.414$,
$2.414$,
$3.414$,
$3.414$,
$4.828$,
$5.828$) 

\vskip 0.7ex
\hangindent=3em \hangafter=1
$D^2= 93.254 = 
48+32\sqrt{2}$

\vskip 0.7ex
\hangindent=3em \hangafter=1
$T = ( 0,
\frac{1}{2},
\frac{1}{2},
\frac{1}{4},
\frac{1}{4},
\frac{5}{8},
0 )
$,

\vskip 0.7ex
\hangindent=3em \hangafter=1
$S$ = ($ 1$,
$ 1+\sqrt{2}$,
$ 1+\sqrt{2}$,
$ 2+\sqrt{2}$,
$ 2+\sqrt{2}$,
$ 2+2\sqrt{2}$,
$ 3+2\sqrt{2}$;\ \ 
$ -1-2  \zeta^{1}_{8}
-2  \zeta^{2}_{8}
$,
$ -1-2  \zeta^{-1}_{8}
+2\zeta^{2}_{8}
$,
$(-2-\sqrt{2})\mathrm{i}$,
$(2+\sqrt{2})\mathrm{i}$,
$ 2+2\sqrt{2}$,
$ -1-\sqrt{2}$;\ \ 
$ -1-2  \zeta^{1}_{8}
-2  \zeta^{2}_{8}
$,
$(2+\sqrt{2})\mathrm{i}$,
$(-2-\sqrt{2})\mathrm{i}$,
$ 2+2\sqrt{2}$,
$ -1-\sqrt{2}$;\ \ 
$ (2+2\sqrt{2})\zeta_{8}^{3}$,
$ (-2-2\sqrt{2})\zeta_{8}^{1}$,
$0$,
$ 2+\sqrt{2}$;\ \ 
$ (2+2\sqrt{2})\zeta_{8}^{3}$,
$0$,
$ 2+\sqrt{2}$;\ \ 
$0$,
$ -2-2\sqrt{2}$;\ \ 
$ 1$)

Prime. 

\vskip 1ex 
\color{grey}

\noindent(9,2). $7_{7,93.25}^{8,101}$ \irep{54}:\ \ 
$d_i$ = ($1.0$,
$2.414$,
$2.414$,
$3.414$,
$3.414$,
$4.828$,
$5.828$) 

\vskip 0.7ex
\hangindent=3em \hangafter=1
$D^2= 93.254 = 
48+32\sqrt{2}$

\vskip 0.7ex
\hangindent=3em \hangafter=1
$T = ( 0,
\frac{1}{2},
\frac{1}{2},
\frac{3}{4},
\frac{3}{4},
\frac{3}{8},
0 )
$,

\vskip 0.7ex
\hangindent=3em \hangafter=1
$S$ = ($ 1$,
$ 1+\sqrt{2}$,
$ 1+\sqrt{2}$,
$ 2+\sqrt{2}$,
$ 2+\sqrt{2}$,
$ 2+2\sqrt{2}$,
$ 3+2\sqrt{2}$;\ \ 
$ -1-2  \zeta^{-1}_{8}
+2\zeta^{2}_{8}
$,
$ -1-2  \zeta^{1}_{8}
-2  \zeta^{2}_{8}
$,
$(-2-\sqrt{2})\mathrm{i}$,
$(2+\sqrt{2})\mathrm{i}$,
$ 2+2\sqrt{2}$,
$ -1-\sqrt{2}$;\ \ 
$ -1-2  \zeta^{-1}_{8}
+2\zeta^{2}_{8}
$,
$(2+\sqrt{2})\mathrm{i}$,
$(-2-\sqrt{2})\mathrm{i}$,
$ 2+2\sqrt{2}$,
$ -1-\sqrt{2}$;\ \ 
$ (-2-2\sqrt{2})\zeta_{8}^{1}$,
$ (2+2\sqrt{2})\zeta_{8}^{3}$,
$0$,
$ 2+\sqrt{2}$;\ \ 
$ (-2-2\sqrt{2})\zeta_{8}^{1}$,
$0$,
$ 2+\sqrt{2}$;\ \ 
$0$,
$ -2-2\sqrt{2}$;\ \ 
$ 1$)

Prime. 

\vskip 1ex 
\color{grey}

\noindent(9,3). $7_{1,2.745}^{8,318}$ \irep{54}:\ \ 
$d_i$ = ($1.0$,
$0.171$,
$0.585$,
$0.585$,
$-0.414$,
$-0.414$,
$-0.828$) 

\vskip 0.7ex
\hangindent=3em \hangafter=1
$D^2= 2.745 = 
48-32\sqrt{2}$

\vskip 0.7ex
\hangindent=3em \hangafter=1
$T = ( 0,
0,
\frac{1}{4},
\frac{1}{4},
\frac{1}{2},
\frac{1}{2},
\frac{1}{8} )
$,

\vskip 0.7ex
\hangindent=3em \hangafter=1
$S$ = ($ 1$,
$ 3-2\sqrt{2}$,
$ 2-\sqrt{2}$,
$ 2-\sqrt{2}$,
$ 1-\sqrt{2}$,
$ 1-\sqrt{2}$,
$ 2-2\sqrt{2}$;\ \ 
$ 1$,
$ 2-\sqrt{2}$,
$ 2-\sqrt{2}$,
$ -1+\sqrt{2}$,
$ -1+\sqrt{2}$,
$ -2+2\sqrt{2}$;\ \ 
$ (-2+2\sqrt{2})\zeta_{8}^{3}$,
$ (2-2\sqrt{2})\zeta_{8}^{1}$,
$(2-\sqrt{2})\mathrm{i}$,
$(-2+\sqrt{2})\mathrm{i}$,
$0$;\ \ 
$ (-2+2\sqrt{2})\zeta_{8}^{3}$,
$(-2+\sqrt{2})\mathrm{i}$,
$(2-\sqrt{2})\mathrm{i}$,
$0$;\ \ 
$ -1+2\zeta^{1}_{8}
-2  \zeta^{2}_{8}
$,
$ -1+2\zeta^{-1}_{8}
+2\zeta^{2}_{8}
$,
$ 2-2\sqrt{2}$;\ \ 
$ -1+2\zeta^{1}_{8}
-2  \zeta^{2}_{8}
$,
$ 2-2\sqrt{2}$;\ \ 
$0$)

Prime. 

Not pseudo-unitary. 

\vskip 1ex 
\color{grey}

\noindent(9,4). $7_{7,2.745}^{8,353}$ \irep{54}:\ \ 
$d_i$ = ($1.0$,
$0.171$,
$0.585$,
$0.585$,
$-0.414$,
$-0.414$,
$-0.828$) 

\vskip 0.7ex
\hangindent=3em \hangafter=1
$D^2= 2.745 = 
48-32\sqrt{2}$

\vskip 0.7ex
\hangindent=3em \hangafter=1
$T = ( 0,
0,
\frac{3}{4},
\frac{3}{4},
\frac{1}{2},
\frac{1}{2},
\frac{7}{8} )
$,

\vskip 0.7ex
\hangindent=3em \hangafter=1
$S$ = ($ 1$,
$ 3-2\sqrt{2}$,
$ 2-\sqrt{2}$,
$ 2-\sqrt{2}$,
$ 1-\sqrt{2}$,
$ 1-\sqrt{2}$,
$ 2-2\sqrt{2}$;\ \ 
$ 1$,
$ 2-\sqrt{2}$,
$ 2-\sqrt{2}$,
$ -1+\sqrt{2}$,
$ -1+\sqrt{2}$,
$ -2+2\sqrt{2}$;\ \ 
$ (2-2\sqrt{2})\zeta_{8}^{1}$,
$ (-2+2\sqrt{2})\zeta_{8}^{3}$,
$(2-\sqrt{2})\mathrm{i}$,
$(-2+\sqrt{2})\mathrm{i}$,
$0$;\ \ 
$ (2-2\sqrt{2})\zeta_{8}^{1}$,
$(-2+\sqrt{2})\mathrm{i}$,
$(2-\sqrt{2})\mathrm{i}$,
$0$;\ \ 
$ -1+2\zeta^{-1}_{8}
+2\zeta^{2}_{8}
$,
$ -1+2\zeta^{1}_{8}
-2  \zeta^{2}_{8}
$,
$ 2-2\sqrt{2}$;\ \ 
$ -1+2\zeta^{-1}_{8}
+2\zeta^{2}_{8}
$,
$ 2-2\sqrt{2}$;\ \ 
$0$)

Prime. 

Not pseudo-unitary. 

\vskip 1ex 
\black

\noindent(10,1). $7_{\frac{30}{11},135.7}^{11,157}$ \irep{66}:\ \ 
$d_i$ = ($1.0$,
$2.918$,
$3.513$,
$3.513$,
$4.601$,
$5.911$,
$6.742$) 

\vskip 0.7ex
\hangindent=3em \hangafter=1
$D^2= 135.778 = 
55+44c^{1}_{11}
+33c^{2}_{11}
+22c^{3}_{11}
+11c^{4}_{11}
$

\vskip 0.7ex
\hangindent=3em \hangafter=1
$T = ( 0,
\frac{1}{11},
\frac{4}{11},
\frac{4}{11},
\frac{3}{11},
\frac{6}{11},
\frac{10}{11} )
$,

\vskip 0.7ex
\hangindent=3em \hangafter=1
$S$ = ($ 1$,
$ 2+c^{1}_{11}
+c^{2}_{11}
+c^{3}_{11}
+c^{4}_{11}
$,
$ \xi_{11}^{5}$,
$ \xi_{11}^{5}$,
$ 2+2c^{1}_{11}
+c^{2}_{11}
+c^{3}_{11}
+c^{4}_{11}
$,
$ 2+2c^{1}_{11}
+c^{2}_{11}
+c^{3}_{11}
$,
$ 2+2c^{1}_{11}
+2c^{2}_{11}
+c^{3}_{11}
$;\ \ 
$ 2+2c^{1}_{11}
+2c^{2}_{11}
+c^{3}_{11}
$,
$ -\xi_{11}^{5}$,
$ -\xi_{11}^{5}$,
$ 2+2c^{1}_{11}
+c^{2}_{11}
+c^{3}_{11}
$,
$ 1$,
$ -2-2  c^{1}_{11}
-c^{2}_{11}
-c^{3}_{11}
-c^{4}_{11}
$;\ \ 
$ s^{2}_{11}
+2\zeta^{3}_{11}
-\zeta^{-3}_{11}
+\zeta^{4}_{11}
+\zeta^{5}_{11}
$,
$ -1-c^{1}_{11}
-2  \zeta^{2}_{11}
-2  \zeta^{3}_{11}
+\zeta^{-3}_{11}
-\zeta^{4}_{11}
-\zeta^{5}_{11}
$,
$ \xi_{11}^{5}$,
$ -\xi_{11}^{5}$,
$ \xi_{11}^{5}$;\ \ 
$ s^{2}_{11}
+2\zeta^{3}_{11}
-\zeta^{-3}_{11}
+\zeta^{4}_{11}
+\zeta^{5}_{11}
$,
$ \xi_{11}^{5}$,
$ -\xi_{11}^{5}$,
$ \xi_{11}^{5}$;\ \ 
$ -2-c^{1}_{11}
-c^{2}_{11}
-c^{3}_{11}
-c^{4}_{11}
$,
$ -2-2  c^{1}_{11}
-2  c^{2}_{11}
-c^{3}_{11}
$,
$ 1$;\ \ 
$ 2+2c^{1}_{11}
+c^{2}_{11}
+c^{3}_{11}
+c^{4}_{11}
$,
$ 2+c^{1}_{11}
+c^{2}_{11}
+c^{3}_{11}
+c^{4}_{11}
$;\ \ 
$ -2-2  c^{1}_{11}
-c^{2}_{11}
-c^{3}_{11}
$)

Prime. 

\vskip 1ex 
\color{grey}

\noindent(10,2). $7_{\frac{58}{11},135.7}^{11,191}$ \irep{66}:\ \ 
$d_i$ = ($1.0$,
$2.918$,
$3.513$,
$3.513$,
$4.601$,
$5.911$,
$6.742$) 

\vskip 0.7ex
\hangindent=3em \hangafter=1
$D^2= 135.778 = 
55+44c^{1}_{11}
+33c^{2}_{11}
+22c^{3}_{11}
+11c^{4}_{11}
$

\vskip 0.7ex
\hangindent=3em \hangafter=1
$T = ( 0,
\frac{10}{11},
\frac{7}{11},
\frac{7}{11},
\frac{8}{11},
\frac{5}{11},
\frac{1}{11} )
$,

\vskip 0.7ex
\hangindent=3em \hangafter=1
$S$ = ($ 1$,
$ 2+c^{1}_{11}
+c^{2}_{11}
+c^{3}_{11}
+c^{4}_{11}
$,
$ \xi_{11}^{5}$,
$ \xi_{11}^{5}$,
$ 2+2c^{1}_{11}
+c^{2}_{11}
+c^{3}_{11}
+c^{4}_{11}
$,
$ 2+2c^{1}_{11}
+c^{2}_{11}
+c^{3}_{11}
$,
$ 2+2c^{1}_{11}
+2c^{2}_{11}
+c^{3}_{11}
$;\ \ 
$ 2+2c^{1}_{11}
+2c^{2}_{11}
+c^{3}_{11}
$,
$ -\xi_{11}^{5}$,
$ -\xi_{11}^{5}$,
$ 2+2c^{1}_{11}
+c^{2}_{11}
+c^{3}_{11}
$,
$ 1$,
$ -2-2  c^{1}_{11}
-c^{2}_{11}
-c^{3}_{11}
-c^{4}_{11}
$;\ \ 
$ -1-c^{1}_{11}
-2  \zeta^{2}_{11}
-2  \zeta^{3}_{11}
+\zeta^{-3}_{11}
-\zeta^{4}_{11}
-\zeta^{5}_{11}
$,
$ s^{2}_{11}
+2\zeta^{3}_{11}
-\zeta^{-3}_{11}
+\zeta^{4}_{11}
+\zeta^{5}_{11}
$,
$ \xi_{11}^{5}$,
$ -\xi_{11}^{5}$,
$ \xi_{11}^{5}$;\ \ 
$ -1-c^{1}_{11}
-2  \zeta^{2}_{11}
-2  \zeta^{3}_{11}
+\zeta^{-3}_{11}
-\zeta^{4}_{11}
-\zeta^{5}_{11}
$,
$ \xi_{11}^{5}$,
$ -\xi_{11}^{5}$,
$ \xi_{11}^{5}$;\ \ 
$ -2-c^{1}_{11}
-c^{2}_{11}
-c^{3}_{11}
-c^{4}_{11}
$,
$ -2-2  c^{1}_{11}
-2  c^{2}_{11}
-c^{3}_{11}
$,
$ 1$;\ \ 
$ 2+2c^{1}_{11}
+c^{2}_{11}
+c^{3}_{11}
+c^{4}_{11}
$,
$ 2+c^{1}_{11}
+c^{2}_{11}
+c^{3}_{11}
+c^{4}_{11}
$;\ \ 
$ -2-2  c^{1}_{11}
-c^{2}_{11}
-c^{3}_{11}
$)

Prime. 

\vskip 1ex 
\color{grey}

\noindent(10,3). $7_{\frac{86}{11},15.93}^{11,119}$ \irep{66}:\ \ 
$d_i$ = ($1.0$,
$0.342$,
$2.25$,
$2.309$,
$-1.203$,
$-1.203$,
$-1.576$) 

\vskip 0.7ex
\hangindent=3em \hangafter=1
$D^2= 15.935 = 
22-22  c^{1}_{11}
-11  c^{2}_{11}
+11c^{3}_{11}
-33  c^{4}_{11}
$

\vskip 0.7ex
\hangindent=3em \hangafter=1
$T = ( 0,
\frac{4}{11},
\frac{2}{11},
\frac{8}{11},
\frac{10}{11},
\frac{10}{11},
\frac{3}{11} )
$,

\vskip 0.7ex
\hangindent=3em \hangafter=1
$S$ = ($ 1$,
$ 1-c^{1}_{11}
+c^{3}_{11}
-c^{4}_{11}
$,
$ 1+c^{3}_{11}
-c^{4}_{11}
$,
$ 1-c^{4}_{11}
$,
$ -c^{1}_{11}
-c^{2}_{11}
-c^{4}_{11}
$,
$ -c^{1}_{11}
-c^{2}_{11}
-c^{4}_{11}
$,
$ -2  c^{1}_{11}
-c^{2}_{11}
-2  c^{4}_{11}
$;\ \ 
$ 1+c^{3}_{11}
-c^{4}_{11}
$,
$ 2c^{1}_{11}
+c^{2}_{11}
+2c^{4}_{11}
$,
$ 1$,
$ c^{1}_{11}
+c^{2}_{11}
+c^{4}_{11}
$,
$ c^{1}_{11}
+c^{2}_{11}
+c^{4}_{11}
$,
$ 1-c^{4}_{11}
$;\ \ 
$ -1+c^{4}_{11}
$,
$ 1-c^{1}_{11}
+c^{3}_{11}
-c^{4}_{11}
$,
$ -c^{1}_{11}
-c^{2}_{11}
-c^{4}_{11}
$,
$ -c^{1}_{11}
-c^{2}_{11}
-c^{4}_{11}
$,
$ 1$;\ \ 
$ -2  c^{1}_{11}
-c^{2}_{11}
-2  c^{4}_{11}
$,
$ c^{1}_{11}
+c^{2}_{11}
+c^{4}_{11}
$,
$ c^{1}_{11}
+c^{2}_{11}
+c^{4}_{11}
$,
$ -1-c^{3}_{11}
+c^{4}_{11}
$;\ \ 
$ 1+\zeta^{1}_{11}
+2\zeta^{-1}_{11}
+3\zeta^{2}_{11}
+c^{3}_{11}
+\zeta^{4}_{11}
+2\zeta^{-4}_{11}
+2\zeta^{5}_{11}
$,
$ -1-\zeta^{-1}_{11}
-2  \zeta^{2}_{11}
+\zeta^{-2}_{11}
-c^{3}_{11}
-\zeta^{-4}_{11}
-2  \zeta^{5}_{11}
$,
$ -c^{1}_{11}
-c^{2}_{11}
-c^{4}_{11}
$;\ \ 
$ 1+\zeta^{1}_{11}
+2\zeta^{-1}_{11}
+3\zeta^{2}_{11}
+c^{3}_{11}
+\zeta^{4}_{11}
+2\zeta^{-4}_{11}
+2\zeta^{5}_{11}
$,
$ -c^{1}_{11}
-c^{2}_{11}
-c^{4}_{11}
$;\ \ 
$ -1+c^{1}_{11}
-c^{3}_{11}
+c^{4}_{11}
$)

Prime. 

Not pseudo-unitary. 

\vskip 1ex 
\color{grey}

\noindent(10,4). $7_{\frac{2}{11},15.93}^{11,389}$ \irep{66}:\ \ 
$d_i$ = ($1.0$,
$0.342$,
$2.25$,
$2.309$,
$-1.203$,
$-1.203$,
$-1.576$) 

\vskip 0.7ex
\hangindent=3em \hangafter=1
$D^2= 15.935 = 
22-22  c^{1}_{11}
-11  c^{2}_{11}
+11c^{3}_{11}
-33  c^{4}_{11}
$

\vskip 0.7ex
\hangindent=3em \hangafter=1
$T = ( 0,
\frac{7}{11},
\frac{9}{11},
\frac{3}{11},
\frac{1}{11},
\frac{1}{11},
\frac{8}{11} )
$,

\vskip 0.7ex
\hangindent=3em \hangafter=1
$S$ = ($ 1$,
$ 1-c^{1}_{11}
+c^{3}_{11}
-c^{4}_{11}
$,
$ 1+c^{3}_{11}
-c^{4}_{11}
$,
$ 1-c^{4}_{11}
$,
$ -c^{1}_{11}
-c^{2}_{11}
-c^{4}_{11}
$,
$ -c^{1}_{11}
-c^{2}_{11}
-c^{4}_{11}
$,
$ -2  c^{1}_{11}
-c^{2}_{11}
-2  c^{4}_{11}
$;\ \ 
$ 1+c^{3}_{11}
-c^{4}_{11}
$,
$ 2c^{1}_{11}
+c^{2}_{11}
+2c^{4}_{11}
$,
$ 1$,
$ c^{1}_{11}
+c^{2}_{11}
+c^{4}_{11}
$,
$ c^{1}_{11}
+c^{2}_{11}
+c^{4}_{11}
$,
$ 1-c^{4}_{11}
$;\ \ 
$ -1+c^{4}_{11}
$,
$ 1-c^{1}_{11}
+c^{3}_{11}
-c^{4}_{11}
$,
$ -c^{1}_{11}
-c^{2}_{11}
-c^{4}_{11}
$,
$ -c^{1}_{11}
-c^{2}_{11}
-c^{4}_{11}
$,
$ 1$;\ \ 
$ -2  c^{1}_{11}
-c^{2}_{11}
-2  c^{4}_{11}
$,
$ c^{1}_{11}
+c^{2}_{11}
+c^{4}_{11}
$,
$ c^{1}_{11}
+c^{2}_{11}
+c^{4}_{11}
$,
$ -1-c^{3}_{11}
+c^{4}_{11}
$;\ \ 
$ -1-\zeta^{-1}_{11}
-2  \zeta^{2}_{11}
+\zeta^{-2}_{11}
-c^{3}_{11}
-\zeta^{-4}_{11}
-2  \zeta^{5}_{11}
$,
$ 1+\zeta^{1}_{11}
+2\zeta^{-1}_{11}
+3\zeta^{2}_{11}
+c^{3}_{11}
+\zeta^{4}_{11}
+2\zeta^{-4}_{11}
+2\zeta^{5}_{11}
$,
$ -c^{1}_{11}
-c^{2}_{11}
-c^{4}_{11}
$;\ \ 
$ -1-\zeta^{-1}_{11}
-2  \zeta^{2}_{11}
+\zeta^{-2}_{11}
-c^{3}_{11}
-\zeta^{-4}_{11}
-2  \zeta^{5}_{11}
$,
$ -c^{1}_{11}
-c^{2}_{11}
-c^{4}_{11}
$;\ \ 
$ -1+c^{1}_{11}
-c^{3}_{11}
+c^{4}_{11}
$)

Prime. 

Not pseudo-unitary. 

\vskip 1ex 
\color{grey}

\noindent(10,5). $7_{\frac{26}{11},6.412}^{11,370}$ \irep{66}:\ \ 
$d_i$ = ($1.0$,
$0.217$,
$0.763$,
$0.763$,
$1.284$,
$-0.634$,
$-1.465$) 

\vskip 0.7ex
\hangindent=3em \hangafter=1
$D^2= 6.412 = 
11-11  c^{1}_{11}
-33  c^{2}_{11}
-44  c^{3}_{11}
-22  c^{4}_{11}
$

\vskip 0.7ex
\hangindent=3em \hangafter=1
$T = ( 0,
\frac{5}{11},
\frac{2}{11},
\frac{2}{11},
\frac{6}{11},
\frac{7}{11},
\frac{3}{11} )
$,

\vskip 0.7ex
\hangindent=3em \hangafter=1
$S$ = ($ 1$,
$ -2  c^{2}_{11}
-2  c^{3}_{11}
-c^{4}_{11}
$,
$ -c^{2}_{11}
-c^{3}_{11}
-c^{4}_{11}
$,
$ -c^{2}_{11}
-c^{3}_{11}
-c^{4}_{11}
$,
$ 1-c^{3}_{11}
$,
$ -c^{1}_{11}
-c^{2}_{11}
-2  c^{3}_{11}
-c^{4}_{11}
$,
$ -c^{1}_{11}
-2  c^{2}_{11}
-2  c^{3}_{11}
-c^{4}_{11}
$;\ \ 
$ c^{1}_{11}
+2c^{2}_{11}
+2c^{3}_{11}
+c^{4}_{11}
$,
$ -c^{2}_{11}
-c^{3}_{11}
-c^{4}_{11}
$,
$ -c^{2}_{11}
-c^{3}_{11}
-c^{4}_{11}
$,
$ c^{1}_{11}
+c^{2}_{11}
+2c^{3}_{11}
+c^{4}_{11}
$,
$ 1$,
$ 1-c^{3}_{11}
$;\ \ 
$ s^{1}_{11}
+\zeta^{2}_{11}
+\zeta^{-3}_{11}
-\zeta^{4}_{11}
+2\zeta^{-4}_{11}
$,
$ -s^{1}_{11}
+\zeta^{-2}_{11}
+\zeta^{3}_{11}
+2\zeta^{4}_{11}
-\zeta^{-4}_{11}
$,
$ c^{2}_{11}
+c^{3}_{11}
+c^{4}_{11}
$,
$ -c^{2}_{11}
-c^{3}_{11}
-c^{4}_{11}
$,
$ c^{2}_{11}
+c^{3}_{11}
+c^{4}_{11}
$;\ \ 
$ s^{1}_{11}
+\zeta^{2}_{11}
+\zeta^{-3}_{11}
-\zeta^{4}_{11}
+2\zeta^{-4}_{11}
$,
$ c^{2}_{11}
+c^{3}_{11}
+c^{4}_{11}
$,
$ -c^{2}_{11}
-c^{3}_{11}
-c^{4}_{11}
$,
$ c^{2}_{11}
+c^{3}_{11}
+c^{4}_{11}
$;\ \ 
$ -2  c^{2}_{11}
-2  c^{3}_{11}
-c^{4}_{11}
$,
$ -c^{1}_{11}
-2  c^{2}_{11}
-2  c^{3}_{11}
-c^{4}_{11}
$,
$ 1$;\ \ 
$ -1+c^{3}_{11}
$,
$ 2c^{2}_{11}
+2c^{3}_{11}
+c^{4}_{11}
$;\ \ 
$ -c^{1}_{11}
-c^{2}_{11}
-2  c^{3}_{11}
-c^{4}_{11}
$)

Prime. 

Not pseudo-unitary. 

\vskip 1ex 
\color{grey}

\noindent(10,6). $7_{\frac{62}{11},6.412}^{11,110}$ \irep{66}:\ \ 
$d_i$ = ($1.0$,
$0.217$,
$0.763$,
$0.763$,
$1.284$,
$-0.634$,
$-1.465$) 

\vskip 0.7ex
\hangindent=3em \hangafter=1
$D^2= 6.412 = 
11-11  c^{1}_{11}
-33  c^{2}_{11}
-44  c^{3}_{11}
-22  c^{4}_{11}
$

\vskip 0.7ex
\hangindent=3em \hangafter=1
$T = ( 0,
\frac{6}{11},
\frac{9}{11},
\frac{9}{11},
\frac{5}{11},
\frac{4}{11},
\frac{8}{11} )
$,

\vskip 0.7ex
\hangindent=3em \hangafter=1
$S$ = ($ 1$,
$ -2  c^{2}_{11}
-2  c^{3}_{11}
-c^{4}_{11}
$,
$ -c^{2}_{11}
-c^{3}_{11}
-c^{4}_{11}
$,
$ -c^{2}_{11}
-c^{3}_{11}
-c^{4}_{11}
$,
$ 1-c^{3}_{11}
$,
$ -c^{1}_{11}
-c^{2}_{11}
-2  c^{3}_{11}
-c^{4}_{11}
$,
$ -c^{1}_{11}
-2  c^{2}_{11}
-2  c^{3}_{11}
-c^{4}_{11}
$;\ \ 
$ c^{1}_{11}
+2c^{2}_{11}
+2c^{3}_{11}
+c^{4}_{11}
$,
$ -c^{2}_{11}
-c^{3}_{11}
-c^{4}_{11}
$,
$ -c^{2}_{11}
-c^{3}_{11}
-c^{4}_{11}
$,
$ c^{1}_{11}
+c^{2}_{11}
+2c^{3}_{11}
+c^{4}_{11}
$,
$ 1$,
$ 1-c^{3}_{11}
$;\ \ 
$ -s^{1}_{11}
+\zeta^{-2}_{11}
+\zeta^{3}_{11}
+2\zeta^{4}_{11}
-\zeta^{-4}_{11}
$,
$ s^{1}_{11}
+\zeta^{2}_{11}
+\zeta^{-3}_{11}
-\zeta^{4}_{11}
+2\zeta^{-4}_{11}
$,
$ c^{2}_{11}
+c^{3}_{11}
+c^{4}_{11}
$,
$ -c^{2}_{11}
-c^{3}_{11}
-c^{4}_{11}
$,
$ c^{2}_{11}
+c^{3}_{11}
+c^{4}_{11}
$;\ \ 
$ -s^{1}_{11}
+\zeta^{-2}_{11}
+\zeta^{3}_{11}
+2\zeta^{4}_{11}
-\zeta^{-4}_{11}
$,
$ c^{2}_{11}
+c^{3}_{11}
+c^{4}_{11}
$,
$ -c^{2}_{11}
-c^{3}_{11}
-c^{4}_{11}
$,
$ c^{2}_{11}
+c^{3}_{11}
+c^{4}_{11}
$;\ \ 
$ -2  c^{2}_{11}
-2  c^{3}_{11}
-c^{4}_{11}
$,
$ -c^{1}_{11}
-2  c^{2}_{11}
-2  c^{3}_{11}
-c^{4}_{11}
$,
$ 1$;\ \ 
$ -1+c^{3}_{11}
$,
$ 2c^{2}_{11}
+2c^{3}_{11}
+c^{4}_{11}
$;\ \ 
$ -c^{1}_{11}
-c^{2}_{11}
-2  c^{3}_{11}
-c^{4}_{11}
$)

Prime. 

Not pseudo-unitary. 

\vskip 1ex 
\color{grey}

\noindent(10,7). $7_{\frac{10}{11},3.885}^{11,349}$ \irep{66}:\ \ 
$d_i$ = ($1.0$,
$0.169$,
$0.493$,
$0.778$,
$-0.594$,
$-0.594$,
$-1.140$) 

\vskip 0.7ex
\hangindent=3em \hangafter=1
$D^2= 3.885 = 
44+11c^{1}_{11}
-11  c^{2}_{11}
+22c^{3}_{11}
+33c^{4}_{11}
$

\vskip 0.7ex
\hangindent=3em \hangafter=1
$T = ( 0,
\frac{4}{11},
\frac{7}{11},
\frac{2}{11},
\frac{5}{11},
\frac{5}{11},
\frac{1}{11} )
$,

\vskip 0.7ex
\hangindent=3em \hangafter=1
$S$ = ($ 1$,
$ 1-c^{2}_{11}
$,
$ 2+c^{1}_{11}
+2c^{3}_{11}
+2c^{4}_{11}
$,
$ 2+c^{1}_{11}
+c^{3}_{11}
+2c^{4}_{11}
$,
$ 1+c^{3}_{11}
+c^{4}_{11}
$,
$ 1+c^{3}_{11}
+c^{4}_{11}
$,
$ 1-c^{2}_{11}
+c^{4}_{11}
$;\ \ 
$ 2+c^{1}_{11}
+2c^{3}_{11}
+2c^{4}_{11}
$,
$ -1+c^{2}_{11}
-c^{4}_{11}
$,
$ 1$,
$ -1-c^{3}_{11}
-c^{4}_{11}
$,
$ -1-c^{3}_{11}
-c^{4}_{11}
$,
$ 2+c^{1}_{11}
+c^{3}_{11}
+2c^{4}_{11}
$;\ \ 
$ -2-c^{1}_{11}
-c^{3}_{11}
-2  c^{4}_{11}
$,
$ 1-c^{2}_{11}
$,
$ 1+c^{3}_{11}
+c^{4}_{11}
$,
$ 1+c^{3}_{11}
+c^{4}_{11}
$,
$ 1$;\ \ 
$ 1-c^{2}_{11}
+c^{4}_{11}
$,
$ -1-c^{3}_{11}
-c^{4}_{11}
$,
$ -1-c^{3}_{11}
-c^{4}_{11}
$,
$ -2-c^{1}_{11}
-2  c^{3}_{11}
-2  c^{4}_{11}
$;\ \ 
$ 2\zeta^{1}_{11}
-\zeta^{-1}_{11}
+\zeta^{-2}_{11}
-s^{3}_{11}
+\zeta^{5}_{11}
$,
$ -1-2  \zeta^{1}_{11}
+\zeta^{-1}_{11}
-\zeta^{-2}_{11}
-2  \zeta^{-3}_{11}
-c^{4}_{11}
-\zeta^{5}_{11}
$,
$ 1+c^{3}_{11}
+c^{4}_{11}
$;\ \ 
$ 2\zeta^{1}_{11}
-\zeta^{-1}_{11}
+\zeta^{-2}_{11}
-s^{3}_{11}
+\zeta^{5}_{11}
$,
$ 1+c^{3}_{11}
+c^{4}_{11}
$;\ \ 
$ -1+c^{2}_{11}
$)

Prime. 

Not pseudo-unitary. 

\vskip 1ex 
\color{grey}

\noindent(10,8). $7_{\frac{78}{11},3.885}^{11,669}$ \irep{66}:\ \ 
$d_i$ = ($1.0$,
$0.169$,
$0.493$,
$0.778$,
$-0.594$,
$-0.594$,
$-1.140$) 

\vskip 0.7ex
\hangindent=3em \hangafter=1
$D^2= 3.885 = 
44+11c^{1}_{11}
-11  c^{2}_{11}
+22c^{3}_{11}
+33c^{4}_{11}
$

\vskip 0.7ex
\hangindent=3em \hangafter=1
$T = ( 0,
\frac{7}{11},
\frac{4}{11},
\frac{9}{11},
\frac{6}{11},
\frac{6}{11},
\frac{10}{11} )
$,

\vskip 0.7ex
\hangindent=3em \hangafter=1
$S$ = ($ 1$,
$ 1-c^{2}_{11}
$,
$ 2+c^{1}_{11}
+2c^{3}_{11}
+2c^{4}_{11}
$,
$ 2+c^{1}_{11}
+c^{3}_{11}
+2c^{4}_{11}
$,
$ 1+c^{3}_{11}
+c^{4}_{11}
$,
$ 1+c^{3}_{11}
+c^{4}_{11}
$,
$ 1-c^{2}_{11}
+c^{4}_{11}
$;\ \ 
$ 2+c^{1}_{11}
+2c^{3}_{11}
+2c^{4}_{11}
$,
$ -1+c^{2}_{11}
-c^{4}_{11}
$,
$ 1$,
$ -1-c^{3}_{11}
-c^{4}_{11}
$,
$ -1-c^{3}_{11}
-c^{4}_{11}
$,
$ 2+c^{1}_{11}
+c^{3}_{11}
+2c^{4}_{11}
$;\ \ 
$ -2-c^{1}_{11}
-c^{3}_{11}
-2  c^{4}_{11}
$,
$ 1-c^{2}_{11}
$,
$ 1+c^{3}_{11}
+c^{4}_{11}
$,
$ 1+c^{3}_{11}
+c^{4}_{11}
$,
$ 1$;\ \ 
$ 1-c^{2}_{11}
+c^{4}_{11}
$,
$ -1-c^{3}_{11}
-c^{4}_{11}
$,
$ -1-c^{3}_{11}
-c^{4}_{11}
$,
$ -2-c^{1}_{11}
-2  c^{3}_{11}
-2  c^{4}_{11}
$;\ \ 
$ -1-2  \zeta^{1}_{11}
+\zeta^{-1}_{11}
-\zeta^{-2}_{11}
-2  \zeta^{-3}_{11}
-c^{4}_{11}
-\zeta^{5}_{11}
$,
$ 2\zeta^{1}_{11}
-\zeta^{-1}_{11}
+\zeta^{-2}_{11}
-s^{3}_{11}
+\zeta^{5}_{11}
$,
$ 1+c^{3}_{11}
+c^{4}_{11}
$;\ \ 
$ -1-2  \zeta^{1}_{11}
+\zeta^{-1}_{11}
-\zeta^{-2}_{11}
-2  \zeta^{-3}_{11}
-c^{4}_{11}
-\zeta^{5}_{11}
$,
$ 1+c^{3}_{11}
+c^{4}_{11}
$;\ \ 
$ -1+c^{2}_{11}
$)

Prime. 

Not pseudo-unitary. 

\vskip 1ex 
\color{grey}

\noindent(10,9). $7_{\frac{6}{11},2.987}^{11,285}$ \irep{66}:\ \ 
$d_i$ = ($1.0$,
$0.148$,
$0.432$,
$0.521$,
$0.521$,
$-0.682$,
$-0.876$) 

\vskip 0.7ex
\hangindent=3em \hangafter=1
$D^2= 2.987 = 
33-22  c^{1}_{11}
+22c^{2}_{11}
-11  c^{3}_{11}
+11c^{4}_{11}
$

\vskip 0.7ex
\hangindent=3em \hangafter=1
$T = ( 0,
\frac{5}{11},
\frac{10}{11},
\frac{3}{11},
\frac{3}{11},
\frac{9}{11},
\frac{2}{11} )
$,

\vskip 0.7ex
\hangindent=3em \hangafter=1
$S$ = ($ 1$,
$ 1-c^{1}_{11}
+c^{2}_{11}
$,
$ 1-c^{1}_{11}
+c^{2}_{11}
-c^{3}_{11}
$,
$ 1+c^{2}_{11}
+c^{4}_{11}
$,
$ 1+c^{2}_{11}
+c^{4}_{11}
$,
$ 1-c^{1}_{11}
$,
$ 1-c^{1}_{11}
+c^{2}_{11}
-c^{3}_{11}
+c^{4}_{11}
$;\ \ 
$ -1+c^{1}_{11}
$,
$ -1+c^{1}_{11}
-c^{2}_{11}
+c^{3}_{11}
-c^{4}_{11}
$,
$ 1+c^{2}_{11}
+c^{4}_{11}
$,
$ 1+c^{2}_{11}
+c^{4}_{11}
$,
$ 1-c^{1}_{11}
+c^{2}_{11}
-c^{3}_{11}
$,
$ 1$;\ \ 
$ 1-c^{1}_{11}
+c^{2}_{11}
$,
$ -1-c^{2}_{11}
-c^{4}_{11}
$,
$ -1-c^{2}_{11}
-c^{4}_{11}
$,
$ 1$,
$ 1-c^{1}_{11}
$;\ \ 
$ 1+2\zeta^{1}_{11}
+\zeta^{-1}_{11}
+c^{2}_{11}
+2\zeta^{3}_{11}
+\zeta^{-3}_{11}
+2\zeta^{-4}_{11}
+3\zeta^{5}_{11}
$,
$ -2-2  \zeta^{1}_{11}
-\zeta^{-1}_{11}
-2  c^{2}_{11}
-2  \zeta^{3}_{11}
-\zeta^{-3}_{11}
-\zeta^{4}_{11}
-3  \zeta^{-4}_{11}
-3  \zeta^{5}_{11}
$,
$ -1-c^{2}_{11}
-c^{4}_{11}
$,
$ 1+c^{2}_{11}
+c^{4}_{11}
$;\ \ 
$ 1+2\zeta^{1}_{11}
+\zeta^{-1}_{11}
+c^{2}_{11}
+2\zeta^{3}_{11}
+\zeta^{-3}_{11}
+2\zeta^{-4}_{11}
+3\zeta^{5}_{11}
$,
$ -1-c^{2}_{11}
-c^{4}_{11}
$,
$ 1+c^{2}_{11}
+c^{4}_{11}
$;\ \ 
$ 1-c^{1}_{11}
+c^{2}_{11}
-c^{3}_{11}
+c^{4}_{11}
$,
$ -1+c^{1}_{11}
-c^{2}_{11}
$;\ \ 
$ -1+c^{1}_{11}
-c^{2}_{11}
+c^{3}_{11}
$)

Prime. 

Not pseudo-unitary. 

\vskip 1ex 
\color{grey}

\noindent(10,10). $7_{\frac{82}{11},2.987}^{11,276}$ \irep{66}:\ \ 
$d_i$ = ($1.0$,
$0.148$,
$0.432$,
$0.521$,
$0.521$,
$-0.682$,
$-0.876$) 

\vskip 0.7ex
\hangindent=3em \hangafter=1
$D^2= 2.987 = 
33-22  c^{1}_{11}
+22c^{2}_{11}
-11  c^{3}_{11}
+11c^{4}_{11}
$

\vskip 0.7ex
\hangindent=3em \hangafter=1
$T = ( 0,
\frac{6}{11},
\frac{1}{11},
\frac{8}{11},
\frac{8}{11},
\frac{2}{11},
\frac{9}{11} )
$,

\vskip 0.7ex
\hangindent=3em \hangafter=1
$S$ = ($ 1$,
$ 1-c^{1}_{11}
+c^{2}_{11}
$,
$ 1-c^{1}_{11}
+c^{2}_{11}
-c^{3}_{11}
$,
$ 1+c^{2}_{11}
+c^{4}_{11}
$,
$ 1+c^{2}_{11}
+c^{4}_{11}
$,
$ 1-c^{1}_{11}
$,
$ 1-c^{1}_{11}
+c^{2}_{11}
-c^{3}_{11}
+c^{4}_{11}
$;\ \ 
$ -1+c^{1}_{11}
$,
$ -1+c^{1}_{11}
-c^{2}_{11}
+c^{3}_{11}
-c^{4}_{11}
$,
$ 1+c^{2}_{11}
+c^{4}_{11}
$,
$ 1+c^{2}_{11}
+c^{4}_{11}
$,
$ 1-c^{1}_{11}
+c^{2}_{11}
-c^{3}_{11}
$,
$ 1$;\ \ 
$ 1-c^{1}_{11}
+c^{2}_{11}
$,
$ -1-c^{2}_{11}
-c^{4}_{11}
$,
$ -1-c^{2}_{11}
-c^{4}_{11}
$,
$ 1$,
$ 1-c^{1}_{11}
$;\ \ 
$ -2-2  \zeta^{1}_{11}
-\zeta^{-1}_{11}
-2  c^{2}_{11}
-2  \zeta^{3}_{11}
-\zeta^{-3}_{11}
-\zeta^{4}_{11}
-3  \zeta^{-4}_{11}
-3  \zeta^{5}_{11}
$,
$ 1+2\zeta^{1}_{11}
+\zeta^{-1}_{11}
+c^{2}_{11}
+2\zeta^{3}_{11}
+\zeta^{-3}_{11}
+2\zeta^{-4}_{11}
+3\zeta^{5}_{11}
$,
$ -1-c^{2}_{11}
-c^{4}_{11}
$,
$ 1+c^{2}_{11}
+c^{4}_{11}
$;\ \ 
$ -2-2  \zeta^{1}_{11}
-\zeta^{-1}_{11}
-2  c^{2}_{11}
-2  \zeta^{3}_{11}
-\zeta^{-3}_{11}
-\zeta^{4}_{11}
-3  \zeta^{-4}_{11}
-3  \zeta^{5}_{11}
$,
$ -1-c^{2}_{11}
-c^{4}_{11}
$,
$ 1+c^{2}_{11}
+c^{4}_{11}
$;\ \ 
$ 1-c^{1}_{11}
+c^{2}_{11}
-c^{3}_{11}
+c^{4}_{11}
$,
$ -1+c^{1}_{11}
-c^{2}_{11}
$;\ \ 
$ -1+c^{1}_{11}
-c^{2}_{11}
+c^{3}_{11}
$)

Prime. 

Not pseudo-unitary. 

\vskip 1ex 

}

\subsection{Rank 8}

{\small
\input{SsL8lng_}
}

\subsection{Rank 9 }

{\small
\input{SsL9lng_}
}

\subsection{Rank 10 }

{\small
\input{SsL10lng_}
}

\subsection{Rank 11}

{\small
\input{SsL11lng_}
}

\subsection{Rank 12}

{\small
\input{SsL12lng_}
}

\bibliography{all,publst,ref}

\end{document}